\documentclass[12pt,american,english]{article}
\usepackage[LGR,T1]{fontenc}
\usepackage[latin9]{inputenc}
\usepackage{color}
\usepackage{float}
\usepackage{booktabs}
\usepackage{varwidth}
\usepackage{amsmath}
\usepackage{amsthm}
\usepackage{amssymb}
\usepackage{graphicx}
\usepackage[a4paper]{geometry}
\makeatletter

\DeclareRobustCommand{\greektext}{%
  \fontencoding{LGR}\selectfont\def\encodingdefault{LGR}}
\DeclareRobustCommand{\textgreek}[1]{\leavevmode{\greektext #1}}

\providecommand{\tabularnewline}{\\}
\newenvironment{cellvarwidth}[1][t]
    {\begin{varwidth}[#1]{\linewidth}}
    {\@finalstrut\@arstrutbox\end{varwidth}}

\usepackage{babel}
\usepackage{float}
\usepackage{amsmath}
\usepackage{amssymb}
\usepackage{amsthm}
\usepackage{graphicx}
\usepackage{setspace}
\usepackage{color}
\usepackage{pdflscape}
\usepackage{graphicx}
\usepackage{placeins}
\usepackage{caption}
\usepackage{subcaption}
\usepackage{graphics}
\usepackage{booktabs}
\usepackage{multicol}
\usepackage{multirow}

\graphicspath{ {./Figures/} }

\allowdisplaybreaks
\providecommand{\tabularnewline}{\\}
\theoremstyle{remark}

\newtheorem{remark}{Remark}\newtheorem{theorem}{Theorem}\newtheorem{lemma}{Lemma}\newtheorem{corollary}{Corollary}\newtheorem{proposition}{Proposition}

\newcounter{FSalphaM}

\newcounter{SSalphaM}

\theoremstyle{plain}
\theoremstyle{plain}

\providecommand{\lemmaname}{Lemma}

\providecommand{\theoremname}{Theorem}

\newcommand{\be}{\begin{eqnarray}}
\newcommand{\ee}{\end{eqnarray}}
\newcommand{\bestar}{\begin{eqnarray*}}
\newcommand{\eestar}{\end{eqnarray*}}

\newcommand{\ignore}[1]{}
{} \theoremstyle{plain}

\def\>{\geq}
\def\<{\leq}

\usepackage{lmodern}

\makeatletter
\renewcommand{\maketitle}{%
  \begin{center}
    \vspace*{-1.35cm}
    {\LARGE\bfseries \@title\par}%
    \vspace{1.30em}
    {\normalsize \@author\par}%
    \vspace{0.55em}%
    {\small \@date\par}%
  \end{center}%
  \vspace{0.55em}%
}

\renewcommand{\@seccntformat}[1]{\csname the#1\endcsname.\quad}
\renewcommand\section{\@startsection{section}{1}{\z@}%
  {-2.25ex \@plus -.6ex \@minus -.2ex}%
  {0.95ex \@plus .2ex}%
  {\normalfont\normalsize\bfseries}}
\renewcommand\subsection{\@startsection{subsection}{2}{\z@}%
  {-1.85ex \@plus -.5ex \@minus -.2ex}%
  {0.75ex \@plus .2ex}%
  {\normalfont\normalsize\bfseries}}
\renewcommand\subsubsection{\@startsection{subsubsection}{3}{\z@}%
  {-1.55ex \@plus -.4ex \@minus -.2ex}%
  {0.60ex \@plus .15ex}%
  {\normalfont\normalsize\itshape}}
\makeatother

\usepackage{adjustbox}
\let\HTtextcolor\textcolor
\newcommand{\HTpanel}[1]{\vspace{2.6ex}\noindent\textbf{#1}\par\nobreak\vspace{1.3ex}}
\newcommand{\HTmemorybegin}{\begingroup\footnotesize\setlength{\tabcolsep}{3.2pt}\renewcommand{\arraystretch}{0.92}}
\let\HTmemoryend\endgroup

\makeatother

\usepackage{babel}
\begin{document}
\title{Conventional Inference under Nonstationary Heavy Tailed Fractional
Processes \& Arrays}
\author{%
\begin{tabular}{c}
{\large Ioannis Kasparis}\tabularnewline
\begin{cellvarwidth}[t]
\centering
{\small\textit{Department of Economics}}\\
{\small\textit{University of Cyprus}}
\end{cellvarwidth}\tabularnewline
\end{tabular}}
\date{October 2026}
\maketitle
\begin{abstract}
Existing regression inference methods for nonstationary heavy tailed
data are largely limited to unit roots tests. Functional central limit
theorems (FCLTs), which are commonly used for the asymptotic analysis
of nonstationary models, are inapplicable to nonstationary data characterized
by a memory parameter $d<1$ and a stability parameter $\alpha<2$.
Furthermore, it is well known that inference under nonstationarity
is often non-conventional, and this issue is exacerbated in the presence
of heavy-tailed data. The interaction of long memory and heavy tails
significantly enlarges the nonstationary region---i.e., the covariate
model space---where conventional inference methods typically fail. 

To address these challenges, we propose a new limit theory and methods
that facilitate both workable asymptotics and conventional inference
in regressions involving general nonstationary heavy-tailed covariate
processes. Specifically, our approach develops a limit theory for
heavy-tailed weakly nonstationary processes (hereafter, HT-WNPs),
which lie on the threshold of nonstationarity (cf. Duffy and Kasparis,
2021, Annals of Statistics). We demonstrate that the proposed methods---including
IVX-based and related instrumental-variable procedures---facilitate
conventional inference for a broad class of stationary and nonstationary
heavy-tailed covariates. This class includes nearly and mildly integrated
arrays, as well as nonstationary fractional processes with memory
parameter $1-1/\alpha<d<2-1/\alpha$, where $\alpha\in(1,2]$. Our
analysis is applicable to both parametric and kernel regression models.

The finite-sample performance of the proposed methods is evaluated
through a comprehensive simulation study. Additionally, we provide
an empirical application to the predictability of stock returns using
long-memory heavy-tailed predictors.
\end{abstract}
\vspace{-0.5em}

{\small\textit{Keywords. }}{\small Domains of attraction; fractional
process; heavy-tailed process; kernel regression; IVX; mildly integrated
process; nearly integrated process; stable distribution; weakly nonstationary
process}{\small\par}

\section{Introduction}\label{sec:Introduction}

Heavy-tailed (HT) processes are prevalent across various fields, including
economics, finance, urban studies, actuarial science, telecommunications,
meteorology, and epidemiology. For instance, Cohen, Davis, and Samorodnitsky
(2022) provide evidence that COVID-19 death counts in the U.S. exhibit
infinite variance. Additional examples can be found in the works of
Embrechts, Klüppelberg, and Mikosch (1997), Alder, Feldman, and Taqqu
(1998), and Ibragimov, Ibragimov, and Walden (2015). While significant
early research has focused on estimating stationary and unit root
HT autoregressions, considerably less attention has been paid to regressions
that incorporate general HT covariates, which may display either stationary
or nonstationary characteristics with varying persistence (e.g., memory).

In this study, we primarily investigate inference in parametric regressions
of the form 
\begin{equation}
y_{t}=\mu+\beta f(x_{t-1})+u_{t},\text{ }t=1,...,n,\label{Model0}
\end{equation}
where $x_{t}$ is some times series process that is potentially nonstationary
and heavy tailed, $u_{t}$ a finite variance m.d. sequence, and $\left[\mu,\beta\right]$
unknown regressions parameters. We also provide results for the case
where $y_{t}$ is given by a nonparametric regression of the form
\begin{equation}
y_{t}=m(x_{t-1})+u_{t},\text{ }t=1,...,n,\label{Model02}
\end{equation}
with $m(.)$ being some unknown regression function estimated by kernel
methods are also provided. Our goal is to develop methods that enable
conventional inference regardless of the covariate's stationarity
properties and the existence of moments. To address these challenges,
we focus on the asymptotic analysis of Heavy-Tailed Weakly Nonstationary
Processes (HT-WNPs), a generalization of processes recently examined
by Duffy and Kasparis (2021), encompassing processes at the boundary
of nonstationarity.

In this study, WNPs may arise either from the model itself or from
the estimation method. In the former case, the covariate ($x_{t}$)
is assumed to belong to a particular class of WNPs. A prominent example
of a method-specific WNP is provided by the so-called IVX methodology;
see, inter alia, Phillips and Magdalinos (2009) and Kostakis, Magdalinos,
and Stamatogiannis (2015) - KMS2015, hereafter. IVX instruments are
constructed by applying an appropriate linear filter to the covariate
($x_{t}$), thereby generating a WNP that falls within the asymptotic
framework developed in this study.

We consider three estimators: (i) a kernel regression estimator for
the unknown function ($m$) in (\ref{Model0}); (ii) a weighted least-squares
(WIV) estimator, comparable to that studied by Samorodnitsky, Kurz-Kim,
and Rachev (2007); and (iii) a weighted version of the IVX estimator,
referred to as the W-IVX estimator. The kernel regression estimator
is used for both estimation of and inference on the unknown regression
function ($m$). We show that it supports conventional inference for
a broad class of stationary and nonstationary covariates. The latter
two estimators are used for estimation of and inference on the slope
parameter in (\ref{Model0}). The WIV estimator is primarily intended
for point estimation, whereas the W-IVX estimator is designed principally
for inference.

More specifically, the WIV estimator attains parametric convergence
rates but does not always support conventional inference when the
regressor is nonstationary. By contrast, under nonstationarity, the
W-IVX estimator generally converges at semiparametric rates but yields
conventional inference for a broad class of stationary and nonstationary
covariates. Consequently, practitioners may conduct W-IVX-based inference
without precise knowledge of the persistence properties of $x_{t}$.
Finally, when $x_{t}$ is stationary, we show that the WIV and W-IVX
estimators are asymptotically equivalent.

The interplay between heavy tails and persistence introduces inferential
complications that go beyond those arising in the conventional near-unit-root
framework. Existing analyses of highly persistent predictive regressors,
such as Campbell and Yogo (2006), Mikusheva (2007), and Phillips (2014),
maintain finite-variance, short-memory innovations. Within that framework,
nonstandard inference arises from local-to-unity persistence, with
the localizing parameter $c$ constituting the relevant nuisance parameter.
Relaxing the finite-variance and short-memory assumptions introduces
two additional complications. First, heavy tails enlarge the nonstationary
region. For fractional processes, in particular, the conventional
finite-variance stationarity boundary $d=1/2$ is shifted to
\[
d=1-\frac{1}{\alpha},\;\;0<\alpha\leq2,
\]
so that processes that would be classified as stationary under finite
variance may become nonstationary under sufficiently heavy tails.
Second, once fractional dependence and heavy-tailed innovations are
allowed jointly, the limiting distributions can depend not only on
the persistence parameter but also on additional nuisance parameters
governing memory and the tail behavior of the innovations, including
the stability index and skewness. The appearance of these additional
nuisance parameters is specific to the more general framework considered
here and constitutes an important departure from the existing near-unit-root
literature.

Therefore, it likely that data sets exhibiting long memory and heavy
tail characteristics are more accurately described by nonstationary
processes, where conventional, nuisance parameter-free inference generally
does not hold. In economics, which is the focus of our empirical application,
key examples include growth rates of macroeconomic and financial variables
such as inflation (Charemza, Hristova, and Burridge, 2006), exchange
rate returns (Falk and Wang, 2003; Nicolau and Rodrigues, 2019), and
financial return volatility (Kim and Meddahi, 2020; Kim, Meddahi,
and Yamashita, 2021). Additionally, economic activity metrics connected
to extreme events or risk distress are relevant, as noted by Barigozzi,
Cavaliere, and Trapani (2021).

The limit distributions of the estimators under consideration are
also derived for cases where regression error ($u_{t}$) is a heavy-tailed
(HT) term. When the regression error is heavy-tailed, both the instrumental
variable (IV) and kernel methods may produce inconsistent estimators
due to a weak signal-to-noise ratio. However, when the signal-to-noise
ratio is substantial, the estimators are consistent; nonetheless,
convergence rates are slower, and the limit distributions include
nuisance parameters (e.g., tail and skewness parameters) related to
the regression error distribution, which complicates the proposed
tests. One potential solution to this issue is to combine the proposed
methods with least absolute deviations (LAD) or quantile regression.
Alternatively, bootstrap methods may be employed to address the nuisance
parameter problem (see Hall and Yao, 2003; Cavaliere, Georgiev, and
Taylor, 2016, 2018). However, the latter approach is less appealing
for the models considered here, as those studies focus on bootstrap
methods in the context of autoregressive time series models, where
there is no deterioration in the signal-to-noise ratio in cases of
heavy-tailed regression errors. To keep this work a manageable length,
we will leave the exploration of quantile methods for future research.

The literature on nonstationary HT regression models has predominantly
concentrated on random walks i.e. fractional processes with $d=1$
-see e.g. Chan and Tran (1989), Knight (1989), Phillips (1991), Caner
(1997, 1998), Hallin, van den Akker and Werker (2011), Cavaliere,
Georgiev and Taylor (2018). The more recent work of Barigozzi, Cavaliere\ and
Trapani (2021) develops cointegration rank test for systems with HT
unit root processes. The asymptotic theory in the aforementioned studies
a relies on FCLTs for HT processes (e.g. Astrauskas, 1983; Avram and
Taqqu, 1992). For fractional parameter $d<1$ and tail parameter $\alpha<2$
however, FCLTs do not apply.\footnote{In this case partial sum processes converge in terms of finite dimensional
distributions to fractonal stable motions that exhibit discontinuities
of the second kind (see e.g. Astrauskas, 1983). In fact, it can be
demonstated that a fractional stable motion has paths that $a.s.$
unbounded on any interval of positive length, when $d<1$ and $\alpha<2$
-see Samorodnisky and Taqqu (1994) Corollary 10.2.4 and Example 10.2.5.} This technical difficulty possibly explains why there has been little
attention to HT fractional processes that exhibit weaker signal than
unit roots. Moreover, just like the stationary case limit distribution
dependent on tail parameters. In particular, for simple unit root
tests, limit distributions are determined by stochastic integrals
depending on the tail (nuisance) parameter. Hallin et al. (2011) propose
a distribution-free rank test for the unit root hypothesis while Georgiev
and Taylor (2018) develop pivotal unit root tests based on bootstrap
methods. The more recent work of Barigozzi, Cavaliere and Trapani
(2021) develops a randomized test for common stochastic trends that
is free of the tail parameter and at the same time circumvents the
utilization of FCLTs. A comprehensive literature review on methods
for heavy tailed time series models can be found in the aforementioned
work.

It is important to note that developing pivotal tests within a nonstationary
framework is a challenging endeavor, even when second moments exist.
In nonstationary contexts, the limit distributions of various estimators
are non-standard---typically determined by stochastic integrals---and
depend on nuisance parameters, such as those related to near-unit
roots and memory. When infinite variance is present, tail parameters
also emerge in the limit distributions, adding an additional layer
of complexity to an already intricate problem. Conventional estimators,
such as ordinary least squares (OLS), do not follow the standard distributions
under typical assumptions of finite second moments (see Chan and Wei,
1987; Phillips, 1987a). Furthermore, methods that provide valid and
conventional inference under nonstationarity (e.g., Phillips and Hansen,
1990; Johansen, 1991; Phillips, 1995; Robinson and Hualde, 2003; Johansen
and Nielsen, 2012a) lack robustness to local deviations from the (fractional)
unit root paradigm (e.g., Phillips, 1987b; Chan and Wei, 1988; Buchmann
and Chan, 2007). In general, inferential methods that are valid for
nonstationary models may not apply to stationary models, and vice
versa.

A number of relatively recent studies aim to bridge the gap between
stationary and nonstationary environments under the assumption of
finite second moments. One strand of this literature focuses on near-unit-root
forms of nonstationarity. Notably, Campbell and Yogo (2006) and Mikusheva
(2007), among others (see Phillips, 2014, for a review), develop conservative
methods of inference. A drawback of these methods is that the associated
test statistics can be complicated to construct and their critical
values are not readily available from standard statistical tables.
Moreover, although these procedures robustify inference to local deviations
from a unit root, they do not accommodate fractional dependence or
heavy-tailed innovations. Under fractional or heavy-tailed data, the
corresponding limiting distributions generally acquire additional
nuisance parameters associated with memory and the underlying stable
law, rendering the resulting tests non-pivotal.

An alternative approach is provided by Phillips and Magdalinos (2009)
and KMS2015, who develop IV-type estimators based on specific prefiltering
of the covariates, commonly referred to as IVX, that facilitate conventional
inference, including t-statistics with N(0,1) limiting distributions.
Existing IVX theory, however, is likewise developed under finite-variance
assumptions and for covariates belonging to conventional persistence
classes, such as stationary short-memory linear processes, mildly
integrated processes, and nearly integrated processes. Thus, neither
fractional dependence nor heavy-tailed covariates are accommodated
by the existing IVX results. Extending IVX inference to such environments
is therefore nontrivial and constitutes part of the contribution of
the present paper.

IV methods such as IVX are nevertheless particularly attractive in
the present framework. Because inference is based on studentized estimators,
nuisance parameters associated with tail thickness, skewness, and
memory can potentially scale out of the limiting distributions, thereby
restoring pivotal inference. In addition, the resulting test statistics
are straightforward to construct and their critical values are available
from standard statistical tables. For additional work on IVX methods,
see, e.g., Yang, Long, Peng and Cai (2019), Demetrescu, Georgiev,
Rodrigues and Taylor (2022), Magdalinos and Petrova (2022), and the
references therein.

Another strand of the literature develops methods for inference in
fractional systems under finite-variance conditions; see, for example,
Robinson and Hualde (2003), Hualde and Robinson (2010), and Andersen
and Varneskov (2021). These frameworks are more general than the one
considered in the present paper in several important respects. In
particular, they allow for systems of equations and can accommodate
regression disturbances that themselves exhibit long memory. Their
validity, however, relies on assumptions concerning fractional integration
and the existence of finite second moments. Some of these methods
bear a close resemblance to fully modified least squares (FMLS; Phillips,
1995), in that they exploit transformations of the covariates to recover
short-memory innovations and correct for endogeneity. Whereas FMLS
relies on ordinary differencing, Robinson and Hualde (2003), for example,
employ fractional differencing together with GLS-type orthogonalization
to accommodate fractional cointegrating relationships. However, when
the covariate is a nearly integrated array, fractional differencing
no longer recovers the underlying innovations exactly, leaving an
additional component that depends on the localizing parameter. This
component can affect the first-order asymptotic behavior of the estimators
and compromise the validity of conventional inference. Furthermore,
the existing results do not accommodate heavy-tailed innovations with
infinite variance. Thus, although the fractional-cointegration literature
permits substantially richer multivariate structures and more general
dependence in the regression disturbances than considered here, it
does not provide inference that is simultaneously robust to fractional
dependence, near integration, and heavy tails.

Finally, a different approach is developed by Wang and Phillips (2012)
and Duffy and Kasparis (2021), who investigate inferential procedures
based on kernel methods. An important feature of this approach is
that suitably normalized kernel functionals exhibit remarkably similar
limiting behavior across very different persistence regimes. This
facilitates conventional inference for regression models with covariates
ranging from stationary processes to nonstationary fractional processes
and nearly integrated arrays, without requiring the practitioner to
know the precise persistence class of the covariate. In this sense,
kernel methods provide an effective means of conducting inference
while remaining agnostic about the degree and form of persistence
in the data. The results developed in the present paper show that
this robustness extends further to heavy-tailed environments. In particular,
kernel-based procedures can remain pivotal despite the additional
nuisance parameters associated with tail thickness, skewness, and
memory. Hence, kernel methods provide a route to inference that is
robust across stationary and nonstationary regimes, fractional and
near-integrated forms of persistence, and, as established here, heavy-tailed
environments.

The remainder of the paper is organised as follows. Section \ref{sec:Models-and-Main}
presents the models considered in this study. Section \ref{sec:Limit-Theory-WNPs}
develops basic limit theory for functionals of HT-WNPs relevant to
both parametric and nonparametric estimation, while Section \ref{sec:Limit-Theory-SNPs}
establishes corresponding results for functionals of HT-SNPs. Section
\ref{sec:Kernel-Regression} considers estimation and inference in
nonparametric regression models and shows that the resulting limit
distributions are conventional for a broad class of heavy-tailed regressors,
including stationary processes, WNPs, and SNPs. Section \ref{sec:IV-Methods}
considers estimation in parametric regression models, whereas Section
\ref{sec:IVX-Inference} develops robust IVX-based methods for parametric
inference. Section \ref{sec:Simulations} presents simulation results
examining the finite-sample performance of the proposed methods. Finally,
Section \ref{sec:Application} provides an empirical application to
the predictability of stock returns using persistent and potentially
heavy-tailed predictors.

\paragraph{Notation}

For two deterministic sequences $\{a_{n}\}$ and $\{b_{n}\}$, $a_{n}\sim b_{n}$,
$a_{n}\asymp b_{n}$, and $a_{n}\preceq b_{n}$ denote $\lim_{n\rightarrow\infty}a_{n}/b_{n}=1,$
$\lim_{n\rightarrow\infty}a_{n}/b_{n}=C$, and $\lim_{n\rightarrow\infty}a_{n}/b_{n}\leq C\in(0,\infty)$
respectively. For two random variables $X$ and $Y$, $X\sim Y$ and
$X\overset{d}{=}Y$ denote distributional equality between $X$ and
$Y$. For a real number $x$, $\lfloor x\rfloor$ denotes its integer
part. For $a,b\in\mathbb{R}$, $\sum^{b}_{j=a}(.)_{j}$ will denote
the sum $\sum^{\lfloor b\rfloor}_{j=\lfloor a\rfloor}(.)_{j}$. $1\left\{ A\right\} $
denotes the indicator function for the set $A$. As usual, $\overline{\mathbb{R}}$
is the extended real line. $\sigma(X)$ denotes the sigma algebra
generated by random variable $X$. \textit{SV }and \textit{RV} stands
for slowly varying and regularly varying function respectively. In
addition, $RV_{p}$ denotes a regularly varying function of index
$p\in\mathbb{R}$ (see Bingham et al. 1987; p. 6 and p. 18 respectively).
For $x,a\in\mathbb{R}$, $x^{\left\langle a\right\rangle }$ stands
for $\mathrm{sgn}(x)\left|x\right|^{a}$ (signed polynomial) with
$\mathrm{sgn}(.)$ being the sign function. As usual $v_{t}\sim I(\delta)$
stands for a process that is integrated of order $\delta$ (or equivalently
has memory $\delta$). All limits are taken as $n\rightarrow\infty$
unless otherwise indicated, and $:=$ denotes definitional equality.

\section{Models and Main Assumptions}\label{sec:Models-and-Main}

Next, we clearly state the set of assumptions that define the specifications
of (\ref{Model0}) and (\ref{Model02}). As previously mentioned,
the regression error $u_{t}$ is assumed to form a martingale difference
sequence with respect to an appropriate filtration and to possess
finite conditional variance---this is the main assumption. Additionally,
the possibility that $u_{t}$ is an adapted heavy tailed sequence
is also considered. For the regressor process $x_{t}$ three classes
are distinguished:
\begin{itemize}
\item \textbf{HT-SP}: HT Stationary Process.
\item \textbf{HT-WNP}: HT Weakly Nonstationary Process.
\item \textbf{HT-SNP}: HT Strongly Nonstationary Process.
\end{itemize}
Among these, the HT-SP class plays a particularly central role. It
not only defines a class of covariate processes but also serves as
a foundational building block for the other two. In fact, certain
HT-WNP and HT-SNP processes can be constructed as linear filters of
HT-SP processes. The first and the last class essentially comprise
stationary and nonstationary fractional processes. These processes
have been extensively studied under finite variance assumptions see
e.g. Robinson and Hualde (2003), Phillips and Shimotsu (2004), Shimotsu
and Phillips (2005), Christensen and Nielsen (2006), Hualde and Robinson
(2010), Johansen and Nielsen (2012a) among others. The second class
comprises of fractional processes on nonstationarity boundary (i.e.
$d=1-1/\alpha$) and mildly integrated (MI, hereafter) arrays like
those introduced by Giraitis and Phillips (2006), Phillips and Magdalinos
(2007, 2009). A unified asymptotic framework for WNPs was developed
recently by Duffy and Kasparis (2021) (DK2021, hereafter) under finite
variance assumptions. It is important to emphasize that in this work,
MI processes are not only considered as potential models for the covariate
process but also as instrumental variables. The so-called IVX method
employs MI processes to construct t- and F-tests that are robust to
the degree of dependence in the covariates. IVX instruments are obtained
by applying a specific linear filter---referred to as the IVX filter---to
the covariates. Within the current framework, applying the IVX filter
to the HT-SNP class yields an HT-MI process, which can then be used
as an instrument.

\subsection{Covariates }

The first assumption below specifies the innovations that drive the
covariate term $x_{t}$.\bigskip{}

\noindent\textbf{Assumption INN (innovations). }

\noindent\textbf{(i)} $\left\{ \xi_{j}\right\} _{j\in\mathbb{Z}}$\textit{
is i.i.d. in the domain of attraction of a strictly stable law with
stability parameter }$\alpha\in(0,2]$\textit{ such that as }$\lambda\rightarrow0$
\begin{equation}
\ln\mathbf{E(}e^{i\lambda\xi_{1}})=-\sigma^{\alpha}\left\vert \lambda\right\vert ^{\alpha}H\left(\left\vert \lambda\right\vert ^{-1}\right)\left[1+i\beta^{*}\text{sgn}(\lambda)\tan\left(\frac{\pi\alpha}{2}\right)\right]\left[1+o(1)\right],\label{INN}
\end{equation}
w\textit{here }$H$\textit{ is slowly varying at infinity (SV hereafter),
}$\sigma^{\alpha}$ ($\sigma>0)$\textit{ is the scale parameter,
}$\left\vert \beta\right\vert \leq1$ \textit{the skewness parameter,
and $\beta^{*}=\beta I_{\left\{ \alpha\neq1\right\} }.$}

\noindent\textbf{(ii)}\textit{ For }$\alpha=2$,\textit{ and $x>0$
we have }$1-P\left(\xi_{1}\leq x\right)+P\left(\xi_{1}\leq-x\right)=x^{-2}h(x)\left[1+o(1)\right]$\textit{,
with }$h(x)$\textit{ SV.}

\noindent\textbf{(iii)} $\xi_{1}$\textit{ has an absolutely continuous
distribution and its characteristic function} $\psi_{\xi}\left(\lambda\right)$\textit{
satisfies} $\int_{\mathbb{R}}\left\vert \lambda\right\vert \left\vert \psi_{\xi}\left(\lambda\right)\right\vert ^{\theta}d\lambda<\infty$\textit{
for some }$\theta\geq1$\textit{.}

\bigskip{}

\noindent\textbf{Notation}. If the characteristic function of some
random variable $\xi$ admits the representation of (\ref{INN}),
this will be denoted $\xi\sim\tilde{S}_{\alpha}\left\{ \sigma^{\alpha},\beta,H\right\} $.\bigskip{}

\begin{remark} \textbf{\label{Rmk1}}

\textbf{(a)} Following Astrauskas (1983) and Jeganathan (2004, 2008),
we assume that for $\alpha\neq1$, $\xi_{1}$ has location parameter
zero (i.e. strictly stable), and that for $\alpha=1$, $\xi_{1}$
is symmetric (i.e. location and skewness are both zero).

\textbf{(b) }In general, $\xi_{1}$ is in the domain of attraction
of an $\alpha$-stable law with $\alpha\in(0,2)$, if and only if\footnote{See e.g. Ibragimov and Linnik (1971) or Aaronson and Denker (1998).}
there are constants $c_{1},c_{2}\geq0$, $c_{1}+c_{2}>0$, and $h$
slowly varying such that as $0<x\rightarrow\infty$ 
\[
P\left(\xi_{1}>x\right)=c_{1}x^{-\alpha}h(x)\left[1+o(1)\right]
\]
\[
P\left(\xi_{1}\leq-x\right)=c_{2}x^{-\alpha}h(x)\left[1+o(1)\right].
\]

\textbf{(c) }If $Var(\xi_{1})<\infty$ or \textbf{Assumption INN(ii)}
holds, then $\xi_{1}$ is in the domain of attraction of normal law
(i.e. (\ref{INN}) holds with $\alpha=2$ -see Ibragimov and Linnik,
1971; p. 79, 84). Related results under finite variance are provided
by DK2021.\footnote{DK2021 consider limit theory for mildly integrated process driven
by short memory errors. In this work also allow for fractional errors
that are either of positive memory or negative memory. 

{}} In this work for $\alpha=2$ we assume that $\xi_{1}$ has $RV_{2}$
tails which allows for infinite variance.\footnote{Under \textbf{Assumption INN(ii)} we have $Var(\xi_{1}I_{\left\{ \left|\xi_{1}\right|<x\right\} })\sim2\int^{x}_{0}s^{-1}h(s)ds\sim H(x).$
Hence, $\xi_{1}$ has a finite variance if and only if $H$ converges.} 

\textbf{(d) }The connection between $h$ and the SV of (\ref{INN})
is the following (see Ibragimov and Linnik, 1971; proof of Thm 2.6.5)
\[
H(x)=\left\{ \begin{array}{l}
h(x),\text{ }\alpha\in(0,2)\\
-\int^{x}_{0}s^{2}d\left(\frac{h(s)}{s^{2}}\right),\text{ }\alpha=2
\end{array}\right..
\]
Further, when the skewness parameter is zero we have $c_{1}=c_{2}$.

\textbf{(e)} In view of Theorem 1.8.2. in BGT 1987, $h(x)$ can be
chosen to be smooth SV i.e. smooth regularly varying of order zero
(see e.g. BGT 1987, p. 44). \end{remark}

\bigskip{}

\subsubsection{Stationary and Strongly Nonstationary HT Processes}

We begin our analysis with the first and last classes, which have
received the most attention in the literature. First, we introduce
a set of conditions that define the stationary process. As noted earlier,
the stationary process serves as a foundational building block for
both HT-WNP and HT-SNP classes. For this reason, we adopt a distinct
notation---denoting the stationary process by $v_{t}$ rather than
$x_{t}$. However, when the covariate belongs to the HT-SP class,
we may write $x_{t}=v_{t}$. Let $v_{t}$ be a linear process of the
form
\[
v_{t}=\sum^{\infty}_{i=0}c_{i}\xi_{t-i}.
\]
The following assumption, fully specifies $v_{t}.$

\bigskip{}

\noindent\textbf{Assumption SP (stationary process). }$v_{t}$\textit{
is generated by (\ref{model01}), and one of the following conditions
holds:}
\begin{itemize}
\item[\textbf{SM}] (short memory) $\sum^{\infty}_{s=0}c_{s}\neq0$.\emph{ Moreover,}
\textit{either: (i) }$\sum^{\infty}_{s=0}\left\vert c_{s}\right\vert ^{\nu}<\infty$
\textit{for some }$\nu:=\tau\wedge1$\textit{, }$\tau\in(0,\alpha)$\textit{,
}$\alpha\in(0,2]$\textit{, and for $\alpha\in(1,2]$ $\int^{x}_{1}h(u)/udu$
is bounded; or (ii) }$\sum^{\infty}_{s=0}\left\vert c_{s}\right\vert ^{\nu}<\infty$
\textit{for some }$\nu\in(0,\alpha\wedge1)$\textit{, }$\alpha\in(0,2]$.
\item[\textbf{LM}] (long memory) $c_{j}=m_{0}\cdot j^{-m}\left(1+O(j^{-\nu})\right)$,
\textit{as $j\rightarrow\infty$, for some $m_{0}\neq0$, }$0<\nu<1-m$\textit{,}
$m\in\left(1/\alpha,1\right)$, $\alpha\in(1,2]$;
\item[\textbf{NM}] (negative memory) $c_{j}=m_{0}\cdot j^{-m}\left(1+O(j^{-\nu})\right)$,
\textit{as $j\rightarrow\infty$, for some $m_{0}\neq0$, }$0<\nu$\textit{,}
$m\in(1,1+1/\alpha)$, $\alpha\in(1,2]$\textit{. In addition, }$\sum^{\infty}_{i=0}c_{i}=0$\textit{,
}$v_{t}=\sum^{t}_{i=0}c_{i}\xi_{t-i}$\textit{, and }$H(x)\rightarrow1$\textit{
as }$x\rightarrow\infty$; 
\end{itemize}
\bigskip{}

\begin{remark} \textbf{\label{Rmk5a}} Under \textbf{Assumption SP-SM},
$v_{t}$ is a short memory errors linear process. The summability
requirement under condition \textbf{SM(i)} ensures that linear process
$v_{t}$ is well defined -see e.g. Avram and Taqqu (1992), Jeganathan
(2004). \textbf{SM(i)} entails an additional restriction on the SV
function $h$, when $\alpha>1$. Condition \textbf{SM(ii)}, provides
an alternative requirement with a slightly stronger condition on the
linear process coefficients, but without any restrictions on $h$.
Most of the technical results provided can be proved under \textbf{SM(i)}
without any restrictions on $h$. Nevertheless, for the purposes of
\textbf{Lemma \ref{MartingApprox}} some additional restrictions are
required. The aforementioned lemma shows that for any $g$\textit{
}bounded and Lipschitz an appropriate normalized WNP (MI) process
($x_{t}$) driven by \textbf{SM} errors, satisfies
\[
\frac{1}{n}\sum^{n}_{t=1}g\left(\gamma^{-1}_{n}x_{t}\right)=\frac{1}{n}\sum^{n}_{t=1}\mathbf{E}\left[g\left(\gamma^{-1}_{n}x_{t}\right)\mid\mathcal{F}_{0}\right]+o_{p}(1).
\]
\end{remark}

\bigskip{}

\begin{remark} \textbf{\label{Rmk5b}} Under \textbf{Assumption SP-LM},
$v_{t}$ is a stationary long memory error $v_{t}$. In particular,
our formulation allows for $v_{t}\sim I(\delta)$ with $0<\delta<1-1/\alpha$,
$\alpha\in(1,2]$. The explicit relation between the memory parameter
$\delta$ and the exponent $m$ is $\delta=1-m$, see e.g., Johansen
and Nielsen (2012b; Definition 2). We remark that the restriction
on the stability parameter ensures that $v_{t}$ is well defined -see
e.g. Samorodnitsky and Taqqu (1994), Theorem 7.13.2. \end{remark}
Most of the technical results presented here can be derived under
the weaker assumption $c_{j}\sim m_{0}\cdot j^{-m}$. However, for
the purposes of \textbf{Lemma \ref{WCon}}, we impose the slightly
stronger second-order condition $c_{j}=m_{0}\cdot j^{-m}\left(1+O(j^{-\nu})\right)$,\textit{
}$0<\nu<1-m$. This condition remains weaker than the coefficient
approximation typically satisfied by parametric ARFIMA models. In
particular, for such models one has $c_{j}=m_{0}\cdot j^{-m}\left(1+O(j^{-1})\right)$;
see, e.g., Giraitis and Surgailis (2012; Theorem 7.2.1). 

\bigskip{}

\begin{remark} \textbf{\label{Rmk5c}} Under \textbf{Assumption SP-NM},
$v_{t}$ exhibits negative memory. In particular, $v_{t}$ has a memory
parameter $-1/\alpha<\delta<0$, $\alpha\in(1,2]$. In addition, we
assume that $v_{t}$ is a truncated linear of the form $\sum^{\infty}_{i=0}c_{i}\xi_{t-i}I_{\left\{ t-i\geq0\right\} }$,
which corresponds to a type II fractional process -see e.g. Robinson
and Hualde (2003), Phillips and Shimotsu (2004). The latter restriction
is imposed to simplify derivations relating a technical result that
pertains to the WNPs case (i.e. \textbf{Lemma \ref{lem:FSalphaM}}),
and can be possibly relaxed. Similarly to the \textbf{SP-LM} case,
we impose a second-order condition on the coefficients $c_{j}$ that
is satisfied by parametric ARFIMA specifications. Finally, we note
that the restriction on the stability parameter (i.e. $\alpha>1$)
postulates that the innovations have a mean. It is possible to formulate
a negative memory linear process of order $-1/\alpha<\delta<1-1/\alpha$
with $\alpha\in(0,1]$ by restricting the exponent $m\in(1/\alpha,1+1/\alpha)$.
We leave this case for future research. \end{remark}

\bigskip{}

Next, we specify the two nonstationary classes. In general, these
processes can be represented as arrays of the form $\left\{ x_{t}\left(n\right)\right\} ^{n}_{t=1}$,
$n\in\mathbb{N}$, where each element is defined as
\begin{equation}
x_{t}:=x_{t}\left(n\right)=\sum^{t-1}_{j=0}\phi_{j}\left(n\right)v_{t-j},\qquad\text{with}\qquad v_{t}=\sum^{\infty}_{i=0}c_{i}\xi_{t-i}.\label{model01}
\end{equation}
Moreover, for convenience, we assume that $x_{0}=0$. This assumption
can be relaxed through an alternative formulation, as in KMS2015.
In most cases, the process $v_{t}$ is assumed to satisfy \textbf{Assumption
SP }, except in the boundary fractional case with $d=1-1/\alpha$
which is treated separately under \textbf{Assumption BFR} (see below).
\textbf{Assumption SNP} next provides an explicit definition of the
HT-SNP class that essentially comprises nonstationary heavy tailed
fractional arrays. \bigskip{}

\noindent\textbf{Assumption SNP (strongly nonstationary process).
}$x_{t}$\textit{ is generated by (\ref{model01}) with} $\phi_{j}(n)=1-\frac{c}{n}$\textit{,
$c\in\mathbb{R}$ and }$v_{t}$\textit{ satisfies }\textbf{\textit{Assumption
SP}}\textit{.}

\bigskip{}

\begin{remark} \textbf{\label{Rmk5d} }Notice that under \textbf{Assumption
SNP} $x_{t}$ is a nearly integrated (\textbf{NI} hereafter) driven
by stationary fractional errors ($v_{t}$) of memory $\delta$. In
the specific case in which the near-to-unity parameter satisfies $c=0$,
$x_{t}$ reduces to a pure fractional process with memory parameter
$d=1+\delta$. It is worth noting that even if $c\neq0$, the effective
memory parameter of $x_{t}$ remains $d$. By effective memory, we
mean the order of magnitude of the normalising sequence required for
the appropriately rescaled $x_{n}$ to converge weakly to a nondegenerate
limit. Under our assumptions the range of the memory ($d$) and stability
($\alpha$) parameters is
\[
d\in\left\{ \begin{array}{l}
\left\{ 1\right\} \text{, }\alpha\in(0,2]\\
\left(1-1/\alpha,1\right)\cup\left(1,2-1/\alpha\right)\text{, }\alpha\in(1,2]
\end{array}\right..
\]
Notice that for the specific case where $\alpha=2$, this specializes
to the usual nonstationary fractional class $1/2<d<3/2$ that has
been widely considered in the literature -see e.g. Robinson and Hualde
(2003), Phillips and Shimotsu (2004) among others. \end{remark}

\bigskip{}

Hereafter, we will refer to $x_{t}$ that satisfies \textbf{Assumption
SNP }as \textbf{HT-SNP}. 

\subsubsection{Processes on the Nonstationarity Boundary: HT-WNPs}

WNPs play a central role in this work, with MI processes being particularly
important for the development of robust inferential procedures. Accordingly,
we provide a dedicated treatment of this class in the present subsection.
Following DK2021, we consider two types of WNPs. The first class consists
of boundary fractional processes (hereafter BFR), i.e., processes
that lie exactly on the boundary of nonstationarity. In the current
framework, this corresponds to a memory parameter of $d=1-1/\alpha$,
which represents the lower bound of the memory parameter range for
nonstationarity. Under the assumption of finite variance (i.e. $\alpha=2$),
the nonstationarity threshold becomes $d=1/2$ recovering the results
of DK2021 as special case. The second class of WNPs comprises MI arrays,
i.e., autoregressive processes with coefficients that drift slowly
toward unity---at a rate slower than $n$. Both types can be incorporated
into the specification of (\ref{model01}), with the latter class
relying on \textbf{Assumption SP} to define $v_{t}$. For the former
type, i.e. BFR, a different formulation of $v_{t}$ is required. The
following assumption characterizes the BFR class in detail.

\bigskip{}

\noindent\textbf{Assumption BFR (boundary fractional process). }$x_{t}\left(n\right)$\textit{
is generated by (\ref{model01}), as }$x\rightarrow\infty$
\[
L(x):=\int^{x}_{1}\frac{H(u)}{u}du\rightarrow\infty
\]
\textit{and either:}

\bigskip{}

\noindent\textbf{BFR1 }$\phi_{j}=1\ \forall j\geq0$\textit{, }$c_{s}\sim s^{-\left(1+\frac{1}{\alpha}\right)}$\textit{,
and }$\sum^{\infty}_{s=0}c_{s}=0$\textit{; or}

\bigskip{}

\noindent\textbf{BFR2 }$\phi_{j}\sim j^{-1/\alpha}$\textit{, }$\phi_{0}\neq0$\textit{,
}\textbf{\textit{Assumption SM}}\textit{ holds with,} $c_{\infty}:=\sum^{\infty}_{s=0}c_{s}\neq0$\textit{.
Further, one of the following holds:}
\begin{itemize}
\item[(a)] \begin{flushleft}
$\lim\sup_{j\rightarrow\infty}\bar{c}_{j}\left\vert c_{j}\right\vert <\infty$,
\textit{where }$\bar{c}_{j}:=\varsigma(j)j^{\frac{1}{\alpha\wedge1}}$\textit{
with }$\varsigma(x)$\textit{ SV such that}
\par\end{flushleft}

\end{itemize}
$\begin{array}{l}
\lim_{x\rightarrow\infty}\varsigma(x)^{-1}=0,\:\mathit{if}\:\alpha<1,\\
\lim_{x\rightarrow\infty}\varsigma(x)^{-1}\ln(x)=0,\:\mathit{if}\:\alpha=1,\\
\varsigma(x)=1,\:\mathit{if}\:\alpha\in(1,2];
\end{array}$
\begin{flushleft}
\textit{or}
\par\end{flushleft}
\begin{itemize}
\item[(b)] $\sum^{\infty}_{j=1}j^{1/\alpha}\left\vert c_{j}\right\vert <\infty$\textit{.}\bigskip{}
\end{itemize}
\begin{remark} \textbf{\label{Rmk2}}

(a) The condition $L(x)\rightarrow\infty$ effectively requires $x_{t}$
being nonstationary with $x_{\lfloor nr\rfloor}$ converging weakly
upon standardization. The function $L(x)$ defines a normalizing sequence
$\gamma_{n}\rightarrow\infty$ required for the term $\gamma^{-1}_{n}x_{n}$
to converge weakly under \textbf{Assumption BFR} -see \textbf{Definition
1} and \textbf{Lemma \ref{WCon}} below.

(b) $L(x)$ is SV (see e.g. BGT 1987; Proposition 1.5.9a).

(c) Notice for some $x_{0}>0$ large enough, $\int^{x}_{x_{0}}\frac{H(u)}{u}du$,
is well defined (see e.g. BGT 1987; Lemma 1.3.2). For simplicity in
\textbf{Assumption BFR} we assume $x_{0}=1$. \end{remark}

\begin{remark} \textbf{\label{Rmk3}} Conditions (a) and (b) of \textbf{Assumption
BFR2} are technical restrictions on the short memory component ($v_{t}$)
that ensure
\[
\sum^{s}_{j=0}\phi_{j}c_{s-j}\sim s^{-1/\alpha}\sum^{\infty}_{j=0}c_{j}\text{.}
\]
\end{remark}

\bigskip{}

\begin{remark} \textbf{\label{Rmk5}} Under \textbf{BFR1} and \textbf{BFR2}
we have type I and type II fractional $I\left(1-1/\alpha\right)$,
$\alpha\in(0,2]$ processes respectively. The parametric counterpart
of these processes is an ARFIMA of memory parameter $d=1-1/\alpha$.
Clearly, when $\xi_{t}$ possesses second moment, $\alpha=2$ and
therefore we get $d=1/2$ as a special case (cf. DK2021). \end{remark}

\bigskip{}

\noindent We now introduce the MI class. Similar to HT-SNPs, MI processes
can be modeled using (\ref{model01}) under\textbf{ Assumption SP}.
The key distinction between these two classes lies in the specification
of the parameter $\phi_{j}(n)$. For HT-SNPs, $\phi_{j}(n)$, rendering
the process a near unit root model applied to the linear process errors
described in \textbf{Assumption SP}. In contrast, for MI processes,
$\phi_{j}(n)$ is a sequence that drifts slowly toward unity, as shown
in the equation below. Let 
\begin{equation}
\rho_{n}=1-\kappa^{-1}_{n},\label{rho_n}
\end{equation}
where $\kappa_{n}>0$ such that 
\[
\kappa_{n}=n^{\alpha_{\kappa}}\ell_{\kappa}(n),
\]
with $\ell_{\kappa}$ SV and $\alpha_{\kappa}\in[0,1)$, $\kappa_{n}\rightarrow\infty$
and $\sup_{n\geq1,1\leq t\leq n}\sup\kappa_{t}/\kappa_{n}<\infty$.

\bigskip{}

\noindent\textbf{Assumption MI (mildly integrated process). }$x_{t}\left(n\right)$\textit{
is generated by (\ref{model01}) with} $\phi_{j}(n)=\rho^{j}_{n}$\textit{
and }$v_{t}$\textit{ satisfies }\textbf{\textit{Assumption SP}}\textit{.}\bigskip{}

\begin{remark} \textbf{\label{Rmk5e}} \textbf{Assumption MI} encompasses
three types of MI processes, extending the framework of DK2021, who
consider only MI processes driven by short-memory errors. This generalization
is crucial for the development of IVX methods applicable to HT-SNP
covariates. 

\end{remark}

\subsection{Regression Errors}

Next, specify the regression error $u_{t}$ of (\ref{Model0}) and
(\ref{Model02}). We consider two possibilities for the regression
error. First, we assume that $u_{t}$ is (together with an appropriate
filtration) forms a martingale difference error with finite conditional
variance --see \textbf{Assumption MDE} below. We show that under
this assumption the IVX method proposed here delivers conventional
inference. We also consider the possibility of heavy-tailed regression
errors (Assumption HTE below) in a number of illustrative examples
involving the kernel and WIV estimators. These examples highlight
the importance of the noise-to-signal ratio for the asymptotic behavior
of the estimators. In particular, when this ratio is sufficiently
large, the estimators may become inconsistent. Under more favorable
noise-to-signal conditions, consistency may be preserved, although
convergence generally occurs at reduced rates. Moreover, the resulting
limiting distributions typically depend on additional nuisance parameters
associated with the heavy-tailed regression errors, rendering conventional
regression tests non-pivotal. These findings illustrate the challenges
of extending conventional inference to models with heavy-tailed regression
disturbances. A potentially promising approach is to combine the methods
developed here with LAD or quantile regression techniques. We leave
this development for future research.\bigskip{}

\noindent\textbf{Assumption MDE (martingale difference error).}
\begin{itemize}
\item[\textbf{(i)}] \emph{$\left\{ x_{t},u_{t}\right\} _{t\geq1}$ are $\mathcal{F}_{t}$-measurable
with $\mathcal{F}_{t}=\mathcal{F}_{0}\vee\sigma\left(u_{1},...,u_{t},\xi_{1},...,\xi_{t}\right)$,
$t\geq1$, and $\mathcal{F}_{0}=\sigma\left(...,\xi_{-1},\xi_{0}\right)$;}
\item[\textbf{(ii)}] \emph{$\left\{ u_{t},\mathcal{F}_{t}\right\} $ is martingale difference
such that $\mathbf{E}\left(u^{2}_{t}\mid\mathcal{F}_{t-1}\right)=\sigma^{2}_{u}<\infty$
$a.s$;}
\item[\textbf{(iii)}] \emph{For any deterministic $0<A_{n}\rightarrow\infty$, $\max_{1\leq t\leq n}\mathbf{E}\left(u^{2}_{t}I\left\{ \left\vert u_{t}\right\vert \geq A_{n}\right\} \mid\mathcal{F}_{t-1}\right)=o_{p}(1)$.}
\end{itemize}
\bigskip{}

The requirements of \textbf{Assumption MDE} are standard in the literature
see e.g. Park and Phillips (2001), Wang and Phillips (2009), Wang
(2014), DK2021. The following alternative assumption relaxes these
requirements by removing the need for a finite conditional variance.

\bigskip{}

\noindent\textbf{Assumption HTE (heavy-tailed error).}
\begin{itemize}
\item[(i)] \textbf{}$\left\{ x_{t},u_{t}\right\} _{t\geq1}$\textit{ are }$\mathcal{F}_{t}$\textit{-measurable
with }$\mathcal{F}_{t}=\mathcal{F}_{0}\vee\sigma\left(u_{1},...,u_{t},\xi_{1},...,\xi_{t}\right)$\textit{,
}$t\geq1$\textit{, and }$\mathcal{F}_{0}=\sigma\left(...,\xi_{-1},\xi_{0}\right)$;
\item[(ii)] \textbf{}\textit{ Let }$\alpha_{u}\in(0,2]$, $\sigma_{u}>0$, $\beta_{u}\in[-1,1]$,
$\beta^{*}_{u}=\beta_{u}I_{\left\{ \alpha_{u}\neq1\right\} }$\textit{
and }$H_{u}$ \textit{SV. The conditional characteristic function
satisfies}
\[
\ln\mathbf{E(}e^{i\lambda u_{t}}\mid\mathcal{F}_{t-1})=-\sigma^{\alpha_{u}}_{u}\left\vert \lambda\right\vert ^{\alpha_{u}}H_{u}\left(\left\vert \lambda\right\vert ^{-1}\right)\left[1+i\beta^{*}_{u}\text{sgn}(\lambda)\tan\left(\frac{\pi\alpha_{u}}{2}\right)\right]\left[1+R_{t}(\lambda)\right],
\]
\textit{where for all $\varepsilon>0$,
\[
\lim_{\delta\downarrow0}\limsup_{n\rightarrow\infty}P\left(\sup_{\left|\lambda\right|\leq\delta}\sup_{1\leq t\leq n}\left|R_{t}(\lambda)\right|>\varepsilon\right)=0.
\]
}
\item[\textbf{(iii)}] $\lim_{x\rightarrow\infty}H_{u}(x)=1$\textit{.} 
\end{itemize}
\bigskip{}
Condition (ii) postulates a conditional characteristic-function expansion
corresponding to a regression error in the domain of attraction of
an $\alpha_{u}$-stable law, i.e.\textit{ }$u_{t}\mid\mathcal{F}_{t-1}\sim S_{\alpha_{u}}\left\{ \sigma^{\alpha_{u}}_{u},\beta_{u},H_{u}\right\} $
with probability approaching one. The remainder condition is uniform
over the triangular array $1\le t\le n$. In particular, the stronger
condition \textit{$\sup_{t\geq1}\left|R_{t}(\lambda)\right|=o_{a.s.}(1)$},
as\textit{ $\lambda\rightarrow0$} is sufficient, since it implies
\textit{$\sup_{\left|\lambda\right|\leq\delta}\sup_{t\geq1}\left|R_{t}(\lambda)\right|=o_{p}(1)$},
as \textit{$\delta\rightarrow0,$} and hence also the requirement
of condition (ii) above.

\section{Limit Theory for Functionals of HT-WNPs}\label{sec:Limit-Theory-WNPs}
\begin{quote}
This section develops a foundational limit theory for functionals
of HT-WNPs. We begin by introducing the necessary notation and definitions.
In view of (\ref{model01}), we can write 
\begin{equation}
x_{t}\left(n\right)=\sum^{\infty}_{k=0}a_{k,t}\left(n\right)\xi_{t-k},\label{linproc}
\end{equation}
where (suppressing the dependence of these quantities on $n$ for
the sake a readability, as we shall do freely below) 
\begin{equation}
a_{k,t}=\sum^{(t-1)\wedge k}_{j=0}\phi_{j}c_{k-j}=\begin{cases}
\sum^{k}_{j=0}\phi_{j}c_{k-j}=:a_{k} & \text{if }0\leq k\leq t-1,\\
\sum^{t-1}_{j=0}\phi_{j}c_{k-j}=:a^{-}_{k,t} & \text{if }k\geq t.
\end{cases}\label{lincoefs}
\end{equation}
Further, decompose 
\begin{equation}
x_{t}(n)=\sum^{t-1}_{k=0}a_{k}\left(n\right)\xi_{t-k}+\sum^{\infty}_{k=t}a^{-}_{k,t}\left(n\right)\xi_{t-k}=:x^{+}_{t}\left(n\right)+x^{-}_{t}\left(n\right).\label{linprocsplit}
\end{equation}
\bigskip{}
\end{quote}
After appropriate normalization --denoted $\gamma_{n}$-- the process
$x_{\left\lfloor nr\right\rfloor }(n)$, $r>0$ converges weakly and
possesses finite sub-$\alpha$ moments. Both properties are essential
for establishing the basic limit theory. The following definition
formalizes the normalizing sequence $\gamma_{n}$.\bigskip{}

\noindent\textbf{Definition 1} \textbf{(normalizing sequence $\gamma_{n}$)}
\begin{itemize}
\item[(i)] \textbf{ Under Assumption BFR}: $\gamma(x)>0$\textit{ is SV such
that $\gamma_{n}:=\gamma(n)$ satisfies}
\[
\gamma_{n}\sim L(n^{1/\alpha}\gamma_{n})^{1/\alpha}.
\]
\item[(ii)] \textbf{ Under Assumption MI}: $\gamma(x)>0$\textit{ and $\gamma_{n}:=\gamma(n)$
satisfies}
\[
\gamma_{n}\sim\kappa^{b/\alpha}_{n}\ell_{1}\left(\kappa_{n}\right)^{1/\alpha},
\]
\textit{where }$\ell_{1}>0$\textit{ is SV} \textit{such that}
\[
\ell_{1}\left(x\right)\sim H\left(x^{1/\alpha}\ell_{1}\left(x\right)^{1/\alpha}\right),
\]
\textit{and}
\[
b:=\left\{ \begin{array}{l}
1\text{, \textit{under }\textbf{SM}, }\\
1+\alpha-\alpha m\text{, \textit{under }\textbf{LM} \textit{\& } \textbf{NM}.}
\end{array}\right.
\]
\end{itemize}
\bigskip{}

\begin{remark} The existence of the SV sequences in \textbf{Definition
1} follows from asymptotic inversion of regularly varying functions,
see e.g., Arvanitis and Louka (2024).\textit{ }Moreover, notice that
by definition the following relation holds.
\[
b\in\left\{ \begin{array}{l}
\left(1,\alpha\right)\text{, under \textbf{LM}, }\alpha\in(1,2]\\
\left\{ 1\right\} \text{, under \textbf{SM}, }\alpha\in(0,2]\\
(0,1)\text{, under \textbf{NM}, }\alpha\in(1,2]
\end{array}\right..
\]

This implies that under \textbf{Assumption MI} $\gamma_{n}$ is regularly
varying in $\kappa_{n}$ of order

\[
b/\alpha\in\left\{ \begin{array}{l}
\left(1/\alpha,1\right)\text{, under \textbf{LM}, }\alpha\in(1,2]\\
\left\{ 1/\alpha\right\} \text{, under \textbf{SM}, }\alpha\in(0,2]\\
(0,1/\alpha)\text{, under \textbf{NM}, }\alpha\in(1,2]
\end{array}\right..
\]
It is worth noting that when the MI is obtained from IVX filtering
of an $I(d)$, $1-1/\alpha<d<2-1/\alpha$ process, $b/\alpha$ corresponds
to the \emph{Hurst} ($\mathcal{H}$) index of underlying process i.e.
\[
b/\alpha=\mathcal{H}=d-1+1/\alpha.
\]
 Hence, $\gamma_{n}$ regularly varying in $\kappa_{n}$ of index
$\mathcal{H}$.

\end{remark}

\bigskip{}

\noindent\refstepcounter{FSalphaM}\label{lem:FSalphaM}\textbf{Lemma
FS}$_{\alpha}$\textbf{M (finite sub-}$\alpha$\textbf{ moments).}
\textit{Suppose that }\textbf{\textit{Assumption INN}}\textit{ holds
and }$\alpha^{\prime}$\textit{ is an arbitrary number in }$\left(0,\alpha\right)$\textit{.
The following hold:}\bigskip{}

\textbf{(a)} \textit{Under }\textbf{\textit{Assumption BFR}}\textit{
and }\textbf{\textit{MI}}\textit{,}
\[
\sup_{n\geq1,1\leq t\leq n}\mathbf{E}\left\vert \gamma^{-1}_{n}x_{t}\right\vert ^{\alpha^{\prime}}<\infty;
\]

\textbf{(b)} \textit{Under }\textbf{\textit{Assumption BFR1}}\textit{
for,} $0<s<r\leq1$\emph{ as }$n\rightarrow\infty$\emph{,} 
\[
\mathbf{E}\left\vert \gamma^{-1}_{n}\left(x^{-}_{\lfloor nr\rfloor}-x^{-}_{\lfloor ns\rfloor}\right)\right\vert ^{\alpha^{\prime}}=o(1);
\]

\textbf{(c)} \textit{Under }\textbf{\textit{Assumption BFR2}}\textit{
and }\textbf{\textit{MI}}\textit{, for all }$0<r\leq1$\textit{ as
}$n\rightarrow\infty$\textit{,}
\[
\mathbf{E}\left\vert \gamma^{-1}_{n}x^{-}_{\lfloor nr\rfloor}\right\vert ^{\alpha^{\prime}}=o(1).
\]

\bigskip{}
Using the same notation as in \textbf{Assumption INN}, we have the
following definition.\bigskip{}

\noindent\textbf{Definition 2} \textbf{(stable distribution) }\textit{A
random variable }$\xi$\textit{ follows an }$\alpha$\textit{-stable
distribution }$S_{\alpha}\left\{ \sigma^{\alpha},\beta\right\} $\textit{
if,}
\begin{equation}
\mathbf{E(}e^{i\lambda\xi})=\exp\left\{ -\sigma^{\alpha}\left\vert \lambda\right\vert ^{\alpha}\left[1+i\beta^{*}\text{sgn}(\lambda)\tan\left(\frac{\pi\alpha}{2}\right)\right]\right\} ,\:\beta^{*}=\beta I_{\left\{ \alpha\neq1\right\} }.\label{SD}
\end{equation}

\bigskip{}

Next, we show that the ``partial sum''\
process $\gamma^{-1}_{n}x_{\lfloor nr\rfloor}$, $0<r\leq1$\
converges weakly. Recall that $c_{\infty}:=\sum^{\infty}_{s=0}c_{s}$
when the sum exists. We have the following result.

\bigskip{}

\begin{lemma} \label{WCon}\textbf{ (weak convergence)}. \textit{Suppose
that }\textbf{\textit{Assumption INN}}\textit{, and either }\textbf{\textit{Assumption
BFR}}\textit{ or }\textbf{\textit{MI}}\textit{ hold. Then for all
}$0<r\leq1$\textit{, we have} \end{lemma}

\[
\gamma^{-1}_{n}\left[x^{+}_{\lfloor nr\rfloor},x^{-}_{\lfloor nr\rfloor}\right]\overset{d}{\rightarrow}\left[X^{+},X^{-}\right],
\]
\textit{where under }\textbf{\textit{Assumption BFR1}}\textit{ }$X^{+},X^{-}\sim iid$\textit{,
whilst under }\textbf{\textit{Assumption BFR2}}\textit{ and }\textbf{\textit{MI}}\textit{
}$X^{-}=0$\textit{. In addition,}

\bigskip{}

\begin{itemize}
\item[\textbf{(i)}] \textit{Under }\textbf{\textit{Assumption BFR1}}
\[
X^{+}\sim S_{\alpha}\left\{ \sigma^{\alpha}\alpha^{\alpha+1},\beta^{*}\right\} .
\]
\item[\textbf{(ii)}] \textit{Under }\textbf{\textit{Assumption BFR2}}\textit{ }
\[
X^{+}\sim S_{\alpha}\left\{ \sigma^{\alpha}\alpha\left\vert c_{\infty}\right\vert ^{\alpha},\beta^{*}\text{sgn}\left(c_{\infty}\right)\right\} .
\]
\item[\textbf{(iii)}] \textit{Under }\textbf{\textit{Assumption MI-SM}}\textit{ }
\[
X^{+}\sim S_{\alpha}\left\{ \sigma^{\alpha}\alpha^{-1}\left\vert c_{\infty}\right\vert ^{\alpha},\beta^{*}\text{sgn}\left(c_{\infty}\right)\right\} .
\]
\item[\textbf{(iv)}] \textit{Under }\textbf{\textit{Assumption MI-LM}}\textit{ }
\[
X^{+}\sim S_{\alpha}\left\{ \sigma^{\alpha}\left\vert m_{0}\right\vert ^{\alpha}\mathbb{V}_{1},\beta^{*}\text{sgn}\left(m_{0}\right)\right\} ,
\]
\textit{where }
\[
\mathbb{V}_{1}:=\int^{\infty}_{x=0}e^{-x}\left(\int^{x}_{y=0}e^{y}y^{-m}dy\right)^{\alpha}dx.
\]
\item[\textbf{(v)}] \textit{Under }\textbf{\textit{Assumption MI-NM}}
\[
X^{+}\sim S_{\alpha}\left\{ \sigma^{\alpha}\left\vert m_{0}\right\vert ^{\alpha}\mathbb{V}_{2},-\beta^{*}\mathbb{V}_{3}\text{sgn}\left(m_{0}\right)\right\} ,
\]
\textit{where }
\[
\mathbb{V}_{2}:=\left\vert 1-m\right\vert ^{-\alpha}\int^{\infty}_{x=0}\underset{=:p(x)}{\left|\right.\underbrace{x^{1-m}-e^{-x}\int^{x}_{y=0}e^{y}y^{1-m}dy}\left.\right|^{\alpha}}dx\mathit{,}
\]
\textit{and }
\[
\mathbb{V}_{3}=\mathbb{V}^{-1}_{2}\cdot\left\vert 1-m\right\vert ^{-\alpha}\int^{\infty}_{0}p(x)^{\left\langle \alpha\right\rangle }dx.
\]
\end{itemize}
\bigskip{}

We remark that if $\alpha=2$ (i.e. Gaussian case), by definition
all skewness parameters equal zero. In addition, $\mathbb{V}_{1}$
reduces to 
\[
\mathbb{V}_{1}=\Gamma\left[2(1-m)\right]\cdot B\left[1-m,2m-1\right],
\]
where $\Gamma$ and $B$ are Gamma and Beta functions respectively.\bigskip{}

Next, we introduce some notation. Hereafter, $\Phi_{X^{+}}$ is the
density of $X^{+}.$\bigskip{}

\begin{theorem} \label{MainLL}\textit{Suppose that }$f:\mathbb{R}\rightarrow\mathbb{R}$\textit{\
is locally Lebesgue integrable. Furthermore, suppose:}
\begin{itemize}
\item[\textbf{(i)}]  $x_{t}$ \textit{satisfies}\textbf{\textit{ Assumption INN}}\textit{
and either }\textbf{\textit{BFR}}\textit{ or }\textbf{\textit{MI}}\textit{
hold;}
\item[\textbf{(ii)}] \textit{ there is some }$\lambda\in\left(0,\alpha\right)$\textit{
such that }$\lvert f(x)\rvert=O\left(\left\vert x\right\vert ^{\lambda}\right)$,
\textit{as }$\left\vert x\right\vert \rightarrow\infty$;
\item[\textbf{(iii)} ] \textit{ (a) either for each $t'$ fixed }$\frac{1}{n}\sum^{t'}_{t=1}f(\gamma^{-1}_{n}x_{t})=o_{p}(1)$\textit{,
or (b) }\textbf{\textit{Assumption INN}}\textit{ holds with $\theta=1$. }
\end{itemize}
\noindent\textit{Then,} \textit{as $n\rightarrow\infty$,} 
\begin{equation}
\frac{1}{n}\sum^{n}_{t=1}f(\gamma^{-1}_{n}x_{t})\overset{d}{\rightarrow}\int_{\mathbb{R}}f(x+X^{-})\Phi_{X^{+}}(x)dx,\label{Corollary2.4}
\end{equation}
\textit{where }$X^{+}$, $X^{-}$ are as in \textbf{Lemma \ref{WCon}}\textit{.} 

\noindent\textit{Furthermore, suppose that }$\tilde{f}:\mathbb{R}\rightarrow\overline{\mathbb{R}}$\textit{
is such that }$\tilde{f}(x)=f(x)$\textit{ almost everywhere and (a)
either for each $t'$ fixed }$\frac{1}{n}\sum^{t'}_{t=1}\tilde{f}(\gamma^{-1}_{n}x_{t})=o_{p}(1)$\textit{,
or (b) }\textbf{\textit{Assumption INN}}\textit{ holds with $\theta=1$.
Then (\ref{Corollary2.4}) holds with }$f$\textit{ replaced by }$\tilde{f}$.
\end{theorem}

\bigskip{}

\noindent\textbf{Assumption K} (kernel). $K:\mathbb{R}\rightarrow\overline{\mathbb{R}}$
\textit{is such that $K$ and $K^{2}$ are Lebesgue integrable, with
}$\int_{\mathbb{R}}K(x)dx\neq0$ .\bigskip{}

\begin{theorem} \label{MainLL_kernel} \textit{Suppose that }
\begin{itemize}
\item[(i)] $x_{t}$ \textit{satisfies}\textbf{\textit{ Assumption INN}}\textit{
and either }\textbf{\textit{BFR}}\textit{ or }\textbf{\textit{MI}}\textit{
hold;}
\item[(ii)] $\{h_{n}\}$ \textit{a positive sequence with }$\gamma^{-1}_{n}h_{n}+\gamma_{n}(nh_{n})^{-1}\rightarrow0$\textit{;}
\item[(iii)] $K$ \textit{satisfies}\textbf{\textit{ Assumption K}}\textit{;}
\item[(iv)] \textit{ for each $t'$ fixed }$\frac{\gamma_{n}}{h_{n}n}\sum^{t'}_{t=1}K\left(\frac{x_{t}-x}{h_{n}}\right)=o_{p}(1)$\textit{.}
\end{itemize}
\noindent\textit{Then,} \textit{as $n\rightarrow\infty$,}
\begin{equation}
\frac{\gamma_{n}}{h_{n}n}\sum^{n}_{t=1}K\left(\frac{x_{t}-x}{h_{n}}\right)\overset{d}{\rightarrow}\Phi_{X^{+}}(-X^{-})\int_{\mathbb{R}}K(x)dx,
\end{equation}
\textit{where }$X^{+}$, $X^{-}$\textit{are as in }\textbf{\textit{Lemma
\ref{WCon}}}\textit{.}

\noindent\end{theorem}

\bigskip{}

We note that condition (iv) in the theorem above is very mild and
satisfied for all bounded kernels.

To simplify notation set 
\[
\varrho(x):=\left\{ \begin{array}{l}
\Phi_{X^{+}}(x-X^{-})\text{, \textbf{under BFR1}}\\
\Phi_{X^{+}}(x)\text{, \textbf{under BFR2, MI}}
\end{array}\right.
\]
Therefore, we can write the limit term of (\ref{Corollary2.4}) as
\[
\int_{\mathbb{R}}f(x+X^{-})\Phi_{X^{+}}(x)dx=\int_{\mathbb{R}}f(x)\varrho(x)dx.
\]

\bigskip{}

\begin{remark} \textbf{\label{Rmk6}} 
\begin{itemize}
\item[(a)]  In view of Condition (ii), Condition (iv) is satisfied when the
kernel functional is bounded. Under some additional requirements on
the coefficients $a_{k,t}\left(n\right)$, it can be shown that \textbf{Conditions
HL7 }and\textbf{ HL9} (Appendix) hold with $t_{0}=1$. In this case
condition (iv) can be dispensed with. 
\item[(b)] Note that the limit term above can be also written as $\varrho(0)\int_{\mathbb{R}}K(x)dx$.
\item[(c)]  It follows from the above that for $h_{n}=1$, 
\end{itemize}
\[
\frac{\gamma_{n}}{n}\sum^{n}_{t=1}K\left(x_{t}\right)\overset{d}{\rightarrow}\Phi_{X^{+}}(-X^{-})\int_{\mathbb{R}}K(x)dx.
\]
\end{remark}

\bigskip{}

\begin{remark} \textbf{\label{Rmk6aa}} \textbf{(comparison with
the stationary case)}.

\textbf{(a)} (finite moment functionals) Consider the fractional linear
process
\[
x_{t}=\sum^{\infty}_{j=0}c_{j}\xi_{t-j},\text{ }c_{j}\sim cons.\cdot j^{d-1}\text{, }d\in\left(-\infty,1-1/\alpha\right)\text{, }\alpha\in(0,2]\text{,}
\]
with $\xi_{t}\sim id.S_{\alpha}\left\{ \sigma^{\alpha},\beta,\mu\right\} $
i.e. $\alpha$-stable with scale, skewness and location parameters
$\sigma^{\alpha},\beta,\mu$ respectively (see e.g. Samorodnitsky
and Taqqu, 1994; p. 381). If in addition $\mathbf{E}\left\vert f\left(x_{t}\right)\right\vert <\infty$,
by ergodicity we have (e.g. Shiryaev, 1996; p. 413)
\[
\frac{1}{n}\sum^{n}_{t=1}f(x_{t})\overset{L_{1}}{\rightarrow}\mathbf{E}f\left(x_{t}\right).
\]
Further, if $\xi_{t}$ has a bounded density the ergodic theorem yields
\[
\frac{1}{h_{n}n}\sum^{n}_{t=1}K\left(\frac{x_{t}-x}{h_{n}}\right)\overset{p}{\rightarrow}\Phi_{x_{1}}(x)\int_{\mathbb{R}}K(x)dx,
\]
with $\Phi_{x_{1}}(.)$ being the density of $x_{1}$ (cf. Wu and
Mielniczuk, 2002; p. 1451).

\textbf{(b)} (infinite moment functionals) Davis and Resnick (1985a,b)
provide limit theory sample covariance functionals of stationary linear
processes of the form $x_{t}=\sum^{\infty}_{j=0}c_{j}\xi_{t-j}$,
with $\sum^{\infty}_{j=0}\left\vert c_{j}\right\vert ^{\tau\wedge1}<\infty$,
$\tau\in(0,\alpha)$ and $\xi_{t}$ are iid possessing regularly varying
tails $RV_{-\alpha}$ $\alpha\in(0,2]$ -see e.g. Davis and Resnick
(1985a; p. 183). Under some additional regularity conditions on $\{c_{j}\}$,
we have the following limit result that is relevant for regression
analysis
\[
\frac{1}{a^{2}_{n}}\sum^{n}_{t=1}x^{2}_{t}\overset{d}{\rightarrow}\sum^{\infty}_{j=0}c^{2}_{j}\cdot X_{0},
\]
with $a_{n}$ being an appropriate diverging sequence and $X_{0}$
an $\alpha/2$-stable variate that corresponds to the weak limit of
$a^{-2}_{n}\sum^{n}_{t=1}\xi^{2}_{t}$ -see e.g. Davis and Resnick
(1985b; p. 271). \end{remark}

\bigskip{}

\begin{remark} \textbf{\label{Rmk6a}} \textbf{(comparison with the
SNP case)}. Theorem \ref{MainLL_kernel} is complementary to the limit
theory of Jeganathan (2004, 2008) for kernel functionals of HT-SNPs.
The aforementioned papers consider type I fractional $x_{t}$\ with
$1-1/\alpha<d<2-1/\alpha$. Note that the requirement on the memory
range implies that $\alpha>1$, when $d\geq1$. In particular, Suppose
that $x_{t}$ satisfies \textbf{Assumption SNP} with $c=0$ (near-to-unity
parameter). It is known (e.g. Astrauskas, 1983; Kasahara and Maejima,
1988,) that for an appropriate $\gamma^{\ast}_{n}\rightarrow\infty$
sequence\footnote{The sequence $\gamma^{\ast}_{n}$ is regularly varying in $n$ of
order $\mathcal{H}=d-1+1/\alpha$. Notice that $\mathcal{H}$ is the
Hurst index of $\Lambda_{\alpha,d}(r)$ -see e.g. Jeganathan (2004).}, the following finite dimensional distribution result holds
\begin{equation}
\gamma^{\ast-1}_{n}x_{\lfloor nr\rfloor}\stackrel{f.d.d.}{\rightarrow}\Lambda_{\alpha,d}(r),\text{ }r\in[0,1],\label{SNP2}
\end{equation}
with $\Lambda_{\alpha,d}$ being the Linear Fractional Stable Motion
(LFSM hereafter) 
\[
\Lambda_{\alpha,d}(r):=\left\{ \begin{array}{l}
\sigma_{0}\int_{\mathbb{R}}\left[\left(r-u\right)^{d-1}_{+}-(-u)^{d-1}_{+}\right]dZ_{\alpha}(u),\text{ }d\neq1\\
\sigma_{0}Z_{\alpha}(r),\text{ }d=1
\end{array}\right.,
\]
where $(x)_{+}:=xI_{\{x>0\}}$, $\sigma_{0}>0$, and $\left\{ Z_{\alpha}(r)\right\} _{r>0}$
is an $\alpha$-stable Lévy process. The limit of (\ref{SNP2}) can
be strengthened to weak convergence in the $D[0,1]$ space, if $\alpha=2$
or $\alpha<2$ with $d\geq1$ (Astrauskas, 1983; Avram and Taqqu,
1992). Note that for $\alpha=2$ $\Lambda_{\alpha,d}(r)$ reduces
to a fractional BM with $1/2<d<3/2$.

Jeganathan (2004, 2008) shows that, for integrable $K$ and $\frac{\gamma^{\ast}_{n}}{h_{n}n}\rightarrow0$,
\begin{equation}
\frac{\gamma^{\ast}_{n}}{h_{n}n}\sum^{n}_{t=1}K\left(\frac{x_{t}-x}{h_{n}}\right)\overset{d}{\rightarrow}L_{\Lambda_{\alpha,d}}(0,1)\int_{\mathbb{R}}K(x)dx,\label{SNP3}
\end{equation}
where $L_{\Lambda_{\alpha,d}}(0,1)$ denotes the local time of $\Lambda_{\alpha,d}$
at the origin up to time 1. It should be emphasized that this limit
result holds under (\ref{SNP2}) and does not require convergence
of $\gamma^{\ast-1}_{n}x_{\lfloor nr\rfloor}$ in the $D[0,1]$ space.
In fact, under \textbf{Assumption SNP}, with $c=0$, and for any locally
integrable $f:\mathbb{R\rightarrow\bar{R}}$ that satisfies certain
polynomial growth conditions, it is possible to show that (see \textbf{Theorem
\ref{SNP-1}} in the following section)
\begin{equation}
\frac{1}{n}\sum^{n}_{t=1}f\left(\frac{x_{t}}{\gamma^{\ast}_{n}}\right)\overset{d}{\rightarrow}\int^{1}_{0}f\left(\Lambda_{\alpha,d}(r)\right)dr,\label{SNP5}
\end{equation}
where $\gamma^{\ast}_{n}$ and $\Lambda_{\alpha,d}(r)$ are given
in (\ref{SNP2}). It is worth noting that (\ref{SNP2}) does not require
an FCLT and is valid even if $d<1$ or $\alpha<2$.

Preliminary results obtained by Kasparis and Wang (2026) suggest that
(\ref{SNP5}) extends to the case $c\neq0$ with $\Lambda_{\alpha,d}(r)$
replaced by the fractional $\alpha$-stable Ornstein--Uhlenbeck process
\[
J_{\alpha,c,d}(r):=\Lambda_{\alpha,d}(r)+c\int^{r}_{0}e^{c\left(r-u\right)}\Lambda_{\alpha,d}(u)du.
\]
A rigorous treatment of this generalization is left for future work.
\end{remark}

\bigskip{}

For the purposes of our regression analysis we need some additional\ limit
\
results involving sample covariance terms of $f(\gamma^{-1}_{n}x_{t})$
or $K\left(\frac{x_{t}-x}{h_{n}}\right)$ and some term $u_{t}$ corresponding
to regression error -cf. (\ref{Model0}). We will primarily assume
that the regression error is a finite variance martingale difference.
Some limit results are also provided for the case where $u_{t}$ is
an iid process in the domain of attraction of an alpha-stable law.
In the latter case the proposed estimators do not yield pivotal statistical
tests due nuisance tail parameters related to the regression error.
Usable inference in the later case can developed using analogous methods
to those considered by Zhu and Ling (2012) or Cavaliere, Georgiev
and Taylor (2016, 2018). We leave explorations in this direction for
future work.

\bigskip{}

Hereafter, in order to simplify conditions we assume that \textbf{Assumption
INN} holds with $\theta=1$. As demonstrated in \textbf{Theorem \ref{MainLL}}
and \textbf{Theorem \ref{MainLL_kernel}}, this requirement can be
relaxed if finite sums of HT-WNPs functionals are asymptotically negligible. 

The following two results generalize \textbf{Theorems \ref{MainLL}}
and \textbf{\ref{MainLL_kernel}}. In particular, \textbf{Theorems
\ref{MCLT-LL}} and \textbf{\ref{MCLT-LLker}} below are relevant
for parametric and nonparametric analysis respectively.\bigskip{}

\begin{theorem} \label{MCLT-LL}\textit{Suppose that }$f^{2},g:\mathbb{R}\rightarrow\mathbb{R}$\textit{
are locally Lebesgue integrable. In addition, suppose that:}
\begin{itemize}
\item[\textbf{(i)}] \textbf{\textit{ Assumption INN}}\textit{ holds with $\theta=1$,
and }$x_{t}$ \textit{satisfies either }\textbf{\textit{Assumption}}\textit{
}\textbf{\textit{BFR}}\textit{ or }\textbf{\textit{MI}}\textit{;}
\item[\textbf{(ii)}] \textit{ there is some }$\lambda\in\left(0,\alpha\right)$\textit{
such that $f(x)^{2}$, }$\lvert g(x)\rvert=O\left(\left\vert x\right\vert ^{\lambda}\right)$\textit{,
as }$\left\vert x\right\vert \rightarrow\infty$.
\item[\textbf{(iii)}]  $x_{t}$\emph{,} $u_{t}$ \textit{satisfy }\textbf{\textit{Assumption
MDE}}\textbf{.} 
\end{itemize}
\noindent\textit{Then,} \textit{as $n\rightarrow\infty$,}
\begin{equation}
\left[\frac{1}{\sqrt{n}}\sum^{n}_{t=1}f(\gamma^{-1}_{n}x_{t-1})u_{t},\text{ }\frac{1}{n}\sum^{n}_{t=1}f(\gamma^{-1}_{n}x_{t-1})^{2},\text{ }\frac{1}{n}\sum^{n}_{t=1}g(\gamma^{-1}_{n}x_{t-1})\right]
\end{equation}
\[
\overset{d}{\rightarrow}\left[MN\left[0,\sigma^{2}_{u}\int_{\mathbb{R}}f(x)^{2}\varrho(x)dx\right],\text{ }\int_{\mathbb{R}}f(x)^{2}\varrho(x)dx,\text{ }\int_{\mathbb{R}}g(x)\varrho(x)dx\right]
\]
\textit{where }$X^{+}$\textit{,} $X^{-}$ \textit{are defined in
}\textbf{\textit{Lemma \ref{WCon}}}\textit{.}\end{theorem}

\bigskip{}
\textbf{Remark}. Under conditions analogous to those imposed in \textbf{Theorem
\ref{MainLL}}, it can be shown that\emph{ }\textit{\emph{(\ref{Corollary2.4})
continues to hold when $f$ and $g$ are extended-real-valued functions.}}

\bigskip{}

\begin{theorem} \label{MCLT-LLker} \textit{Suppose that:}
\begin{itemize}
\item[\textbf{(i)}] \textbf{\textit{Assumption INN}}\textit{ holds with $\theta=1$,
and }$x_{t}$ \textit{satisfies either }\textbf{\textit{Assumption}}\textit{
}\textbf{\textit{BFR}}\textit{ or }\textbf{\textit{MI}}\textit{;}
\item[\textbf{(ii)}] $x_{t}$\emph{, }$u_{t}$\textit{ satisfy }\textbf{\textit{Assumption
MDE}}\textit{;}
\item[\textbf{(iii)}] $\{h_{n}\}$ \textit{is a positive sequence with }$\gamma^{-1}_{n}h_{n}+\gamma_{n}(nh_{n})^{-1}\rightarrow0$\textit{;}
\item[\textbf{(iv)}] $K,$ $K^{2}$\textit{ and }$L$ \textit{satisfy }\textbf{\textit{Assumption
K}}\textit{. }
\end{itemize}
\noindent\textit{Then,} \textit{as $n\rightarrow\infty$,}
\begin{equation}
\left[\sqrt{\frac{\gamma_{n}}{h_{n}n}}\sum^{n}_{t=1}K\left(\frac{x_{t-1}-x}{h_{n}}\right)u_{t},\text{ }\frac{\gamma_{n}}{h_{n}n}\sum^{n}_{t=1}K\left(\frac{x_{t-1}-x}{h_{n}}\right)^{2},\text{ }\frac{\gamma_{n}}{h_{n}n}\sum^{n}_{t=1}L\left(\frac{x_{t-1}-x}{h_{n}}\right)\right]
\end{equation}
\[
\overset{d}{\rightarrow}\left[MN\left[0,\sigma^{2}_{u}\varrho(0)\int_{\mathbb{R}}K(x)^{2}dx\right],\text{ }\varrho(0)\int_{\mathbb{R}}K(x)^{2}dx,\text{ }\varrho(0)\int_{\mathbb{R}}L(x)dx\right].
\]
\end{theorem}

\bigskip{}
It follows directly from \textbf{Theorem \ref{MCLT-LL}} and \textbf{Theorem
\ref{MCLT-LLker}} that, when the regression error $u_{t}$ possesses
a finite second moment, the ``score'' term admits a mixed Gaussian
weak limit. This property underpins the validity of conventional inference.

The next result provides a limit theory for the score terms in the
case where the regression error is an adapted process lying in the
domain of attraction of an $\alpha_{u}$-stable distribution, with
$0<\alpha_{u}\leq2$. To keep the technical exposition simple, we
exclude the case of a \textbf{BFR1}\textit{ }covariate $x_{t}$. 

\begin{theorem} \label{HT-MCLT-LL}\textit{Suppose that:}
\begin{itemize}
\item[\textbf{(i)}] \textbf{\textit{ Assumption INN}}\textit{ holds with $\theta=1$,
and }$x_{t}$ \textit{satisfies either }\textbf{\textit{Assumption}}\textit{
}\textbf{\textit{BFR}}\textit{ or }\textbf{\textit{MI}}\textit{;}
\item[\textbf{(ii)}] $x_{t}$, $u_{t}$ \textit{satisfy} \textbf{\textit{Assumption HTE}}; 
\end{itemize}
\end{theorem}
\begin{itemize}
\item[\textbf{(a)}] \textit{ Suppose that in addition to conditions (i)-(ii) the following
hold:}
\item[\textbf{(iii)}] $\left\vert f\right\vert ^{\alpha_{u}}:\mathbb{R}\rightarrow\mathbb{R}$\textit{
is locally Lebesgue integrable, and $\int_{\mathbb{R}}\left\vert f(x)\right\vert ^{\alpha_{u}}\varrho(x)dx>0$;}
\item[\textbf{(iv)}] \textit{ there is some $\lambda\in\left(0,\alpha\right)$ such that
$\left\vert f(x)\right\vert ^{\alpha_{u}}=O\left(\left\vert x\right\vert ^{\lambda}\right)$,
as }$\left\vert x\right\vert \rightarrow\infty$. 
\end{itemize}
\noindent\textit{Then,} \textit{as $n\rightarrow\infty$,}
\[
n^{-1/\alpha_{u}}\sum^{n}_{t=1}f(\gamma^{-1}_{n}x_{t-1})u_{t}\overset{d}{\rightarrow}S_{\alpha_{u}}\left\{ \sigma^{\alpha_{u}}_{u}\int_{\mathbb{R}}\left\vert f(x)\right\vert ^{\alpha_{u}}\varrho(x)dx,\beta^{\ast}_{u}\frac{\int_{\mathbb{R}}f(x)^{\left\langle \alpha_{u}\right\rangle }\varrho(x)dx}{\int_{\mathbb{R}}\left\vert f(x)\right\vert ^{\alpha_{u}}\varrho(x)dx}\right\} .
\]

\begin{itemize}
\item[\textbf{(b)}] \textit{ Suppose that addition to (i)-(ii) the following hold:}
\item[\textbf{(v)}] $\{h_{n}\}$\textit{is a positive sequence with }$\gamma^{-1}_{n}h_{n}+\gamma_{n}(nh_{n})^{-1}\rightarrow0$\textit{;}
\item[\textbf{(vi)}] $\left\vert K\right\vert ^{\alpha_{u}}$ \textit{satisfies }\textbf{\textit{Assumption
K}}\textit{, and $\int_{\mathbb{R}}\left\vert K(x)\right\vert ^{\alpha_{u}}dx>0$.}
\end{itemize}
\noindent\textit{Then,} \textit{as $n\rightarrow\infty$,}
\[
\left(\frac{\gamma_{n}}{h_{n}n}\right)^{1/\alpha_{u}}\sum^{n}_{t=1}K\left(\frac{x_{t-1}-x}{h_{n}}\right)u_{t}\overset{d}{\rightarrow}S_{\alpha_{u}}\left\{ \sigma^{\alpha_{u}}_{u}\varrho(0)\int_{\mathbb{R}}\left\vert K(x)\right\vert ^{\alpha_{u}}dx,\beta^{\ast}_{u}\frac{\int_{\mathbb{R}}K(x)^{\left\langle \alpha_{u}\right\rangle }dx}{\int_{\mathbb{R}}\left\vert K(x)\right\vert ^{\alpha_{u}}dx}\right\} .
\]

\bigskip{}

Under the assumptions of \textbf{Theorem \ref{HT-MCLT-LL}}, the score
terms possess stable limit distributions. In the special case where
$\alpha_{u}=2$ (i.e. the error is in the domain of attraction of
a normal law), regression estimators and t-statistics admit conventional
limit distributions. By contrast, for $\alpha_{u}<2$, limit distributions
are non-standard, involving nuisance parameters related to the regression
error---namely $\alpha_{u}$, $\beta^{\ast}_{u}$, $\sigma_{u}$.
We conjecture that \textbf{Theorem \ref{HT-MCLT-LL}} can be extended
to settings with a covariate satisfying \textbf{Assumption}\textit{
}\textbf{BFR1}. In such cases, it is natural to anticipate mixed-stable
limit distributions, with $X^{-}$ serving as the mixing variate.
We leave these extensions to future work.

\section{Limit Theory for Functionals of HT-SNPs}\label{sec:Limit-Theory-SNPs}

This section develops limit theory for sample functionals of heavy-tailed
strongly nonstationary processes (\textbf{HT-SNPs}). The results provide
the theoretical foundation for parametric estimation procedures, particularly
the weighted instrumental-variable (WIV) estimator, and for the consistent
estimation of nuisance parameters, including the regression-error
variance ($\sigma^{2}_{u}$), under \textbf{Assumption MDE}. Unlike
the W-IVX procedure, which facilitates conventional inference at the
cost of slower, semiparametric convergence rates, WIV generally achieves
faster convergence rates but does not support conventional inference
in all persistence regimes. Nevertheless, its superior convergence
rates make WIV particularly attractive for point estimation of regression
parameters and associated nuisance quantities, even when the limiting
distributions of the corresponding test statistics are non-pivotal.
Thus, the asymptotic theory developed in this section is relevant
not only for inference but also for establishing the consistency and
convergence rates of parametric estimators.

The main objects of interest are normalized sample averages of the
form $\frac{1}{n}\sum^{n}_{t=1}f\left(\frac{x_{t}}{\gamma^{*}_{n}}\right)$,
where $x_{t}$ is an HT-SNP and $\gamma^{*}_{n}$ is the appropriate
normalizing sequence. The goal is to show that, for a suitable class
of functions $f$, 
\begin{equation}
\frac{1}{n}\sum^{n}_{t=1}f\left(\frac{x_{t}}{\gamma^{*}_{n}}\right)\xrightarrow{d}\int^{1}_{0}f\bigl(\Lambda_{\alpha,d}(r)\bigr)dr,\label{SNP5-1}
\end{equation}
where $\Lambda_{\alpha,d}$ denotes linear fractional stable motion
(LFSM). This result is useful even in regions where a functional central
limit theorem (FCLT) in the usual Skorokhod space is not available.

When $x_{t}$ has finite variance, results such as Marinucci and Robinson
(2000) establish weak convergence of $x_{\lfloor nr\rfloor}/\gamma^{*}_{n}$
in $D[0,1]$ for $1/2<d<3/2$, under moment restrictions on the innovations.
In particular, if the innovations $\xi_{t}$ are i.i.d. with $\mathbf{E}|\xi_{t}|^{q}<\infty$,
the relevant restriction is 
\[
q>2\vee\frac{2}{2d-1}.
\]
Thus, as $d\downarrow1/2$, progressively stronger moment assumptions
are required. At the boundary $d=1/2$, the FCLT no longer applies;
a different weakly nonstationary limit theory is required, as in Duffy
and Kasparis (2021) -see also Johansen and Nielsen (2012b) for impossibility
results.

For HT innovations in the domain of attraction of an $\alpha$-stable
law, with $1<\alpha\leq2$, Astrauskas (1983) obtains weak convergence
of the continuous polygonal interpolation of $x_{\lfloor nr\rfloor}/\gamma^{*}_{n}$
in $C[0,1]$ when 
\[
1<d<2-\frac{1}{\alpha}.
\]
The case $d=1$ corresponds to stable Lévy motion. Avram and Taqqu
(1992) establish weak convergence in $D[0,1]$ for this integrated
short-memory case under additional technical conditions on the short-memory
component.

The remaining nonstationary region is fundamentally different. When
\[
1-\frac{1}{\alpha}<d<1,\:\alpha<2,
\]
the finite-dimensional distributions of $x_{\lfloor nr\rfloor}/\gamma^{*}_{n}$
still converge to those of $\Lambda_{\alpha,d}$; see, for example,
Astrauskas (1983) and Kasahara and Maejima (1988). However, $\Lambda_{\alpha,d}$
is not a $D[0,1]$-valued process in this range. Indeed, Corollary
10.2.4 and Example 10.2.5 of Samorodnitsky and Taqqu (1994) imply
that, for $d<1$ and $\alpha<2$, LFSM has sample paths that are almost
surely unbounded on every interval of positive length. Since every
cadlag function on a compact interval is bounded, weak convergence
in $D[0,1]$ cannot hold. Consequently, the proof of \eqref{SNP5-1}
cannot rely on a continuous mapping argument applied to an FCLT. Instead,
it combines finite-dimensional convergence with suitable moment and
smoothness bounds for the transformed process.

As explained earlier, mildly integrated (MI) processes arise from
applying the IVX filter to strongly nonstationary processes (SNPs).
For the limit theory developed in the previous section, additional
technical restrictions were imposed on the MI/SNP construction. We
now introduce \textbf{Assumptions SP{*}} and \textbf{SNP{*}}, which
are variants of \textbf{Assumptions SP }and\textbf{ SNP}, respectively.
\textbf{Assumption SP{*} }is weaker than \textbf{Assumption SP} and
allows for a more general specification of the stationary component
$v_{t}$, whereas \textbf{Assumption SNP{*}} imposes stronger conditions
than \textbf{Assumption SNP}. The more general specification of $v_{t}$
accommodates a broader class of strongly nonstationary processes.
To keep the paper to a manageable length, we restrict attention to
the case in which the localizing parameter satisfies $c=0$, so that
$\rho_{n}=1$. We begin by specifying $v_{t}$.

\bigskip{}

\noindent\textbf{Assumption SP{*} (stationary process).} \textit{The
process $v_{t}$ is generated by the linear process representation
in \eqref{model01}, and one of the following conditions holds: }
\begin{itemize}
\item[\textbf{SM{*}}] (short memory) $\sum^{\infty}_{s=0}c_{s}\neq0$\emph{,}\textit{ }$\sum^{\infty}_{s=0}\left\vert c_{s}\right\vert ^{\nu}<\infty$
\textit{for some }$\nu:=\tau\wedge1$\textit{, }$\tau\in(0,\alpha)$\textit{,
}$\alpha\in(0,2]$.
\item[\textbf{LM{*}}] (long memory) $c_{j}\sim m_{0}\cdot j^{-m}$\emph{,} \textit{as $j\rightarrow\infty$,
for some $m_{0}\neq0$,} $m\in\left(1/\alpha,1\right)$\emph{,} $\alpha\in(1,2]$;
\item[\textbf{NM{*}}] (negative memory) $c_{j}\sim m_{0}\cdot j^{-m}$\emph{,} \textit{as
$j\rightarrow\infty$, for some $m_{0}\neq0$,} $m\in(1,1+1/\alpha)$\emph{,
}$\alpha\in(1,2]$\textit{. In addition, }$\sum^{\infty}_{i=0}c_{i}=0$. 
\end{itemize}
\bigskip{}
\noindent In view of \textbf{Assumption SP{*}}, we now introduce a weaker version
of \textbf{Assumption SNP.}\bigskip{}
\noindent{}

\noindent\textbf{Assumption SNP{*} (strongly nonstationary process).}
\textit{The process $x_{t}$ is generated by \eqref{model01} with
$\varphi_{j}(n)=1$, for all $j\geq0$, and the stationary component
$v_{t}$ satisfies }\textbf{\textit{Assumption SP}}\textbf{{*}}.

\bigskip{}
\noindent The most notable extension concerns the negative-memory case. Under
Assumption \textbf{SNP{*}-NM}, the process $x_{t}$ is a Type I fractional
process; in particular, it is not truncated. We next define the normalizing
sequence required for the finite-dimensional convergence of $(\gamma^{*}_{n})^{-1}x_{\lfloor nr\rfloor}$.\bigskip{}
\noindent{}

\noindent\textbf{Definition 3 (normalizing sequence $\gamma^{*}_{n}$).}
\textit{Under Assumption }\textbf{\textit{SNP{*}}}\textit{, let $\gamma^{*}(x)>0$
and set $\gamma^{*}_{n}:=\gamma^{*}(n)$, where 
\[
\gamma^{*}_{n}\sim n^{b/\alpha}\ell_{1}(n)^{1/\alpha}.
\]
Here $\ell_{1}>0$ is slowly varying and satisfies 
\[
\ell_{1}(x)\sim H\left(x^{1/\alpha}\ell_{1}(x)^{1/\alpha}\right),
\]
and }
\[
b:=\left\{ \begin{array}{l}
1\text{, \textit{under }\textbf{SM}*, }\\
1+\alpha-\alpha m\text{, \textit{under }\textbf{LM*} \textit{\& } \textbf{NM*}.}
\end{array}\right.
\]
\bigskip{}

It is useful to compare this normalization with the corresponding
one for MI processes. Under \textbf{Assumption SNP{*}}, $\gamma^{*}_{n}$
is regularly varying in $n$ with index $\mathcal{H}=\frac{b}{\alpha}=d-1+\frac{1}{\alpha}$,
where $\mathcal{H}$ is the Hurst, or self-similarity, index of the
limiting LFSM. Under Assumption MI, the analogous normalizing sequence
is regularly varying in the localizing sequence $\kappa_{n}$ with
the same index $\mathcal{H}$. Thus, the SNP normalization uses calendar
time $n$, whereas the MI normalization uses the effective persistence
scale $\kappa_{n}$.

The finite-dimensional convergence 
\begin{equation}
(\gamma^{*}_{n})^{-1}x_{\lfloor nr\rfloor}\stackrel{f.d.d.}{\rightarrow}\Lambda_{\alpha,d}(r),\:r\in[0,1],\label{SNP2-1}
\end{equation}
established, for example by Astrauskas, (1983), and Kashara and Maejima
(1986), together with the following moment-continuity result of \textbf{Lemma
\ref{lem:SSalphaM} }below is sufficient for (\ref{SNP5-1}).\bigskip{}

Before we present the limit theory for nonlinear functionals of HT-SNP
processes, we introduce some class of functions that is relevant to
our asymptotic analysis and provide a preliminary Lemma.

\bigskip{}
\noindent\textbf{Definition LH (local Hölder).} \textit{Let $f:\mathbb{R}\to\mathbb{R}$.
We will say that $f$ is }\textbf{\textit{LH}}\textit{, if for some
$C_{\alpha'}\in(0,\infty)$, 
\begin{equation}
|f(x)-f(y)|\leq C_{\alpha'}\cdot\begin{cases}
|x|^{\alpha'-1}|y-x|+|y-x|^{\alpha'}, & 1<\alpha'<\alpha,\ \alpha>1,\\
|y-x|^{\alpha'}, & 0<\alpha'<\alpha,\ \alpha\leq1.
\end{cases}\label{SNP4-1-1-1}
\end{equation}
}\bigskip{}
\noindent The smoothness restriction in (\ref{SNP4-1-1-1}) is a Hölder-type
condition designed for functions with polynomial growth. The requirement
is mild enough to include many functions used in regression applications,
while being strong enough to facilitate a Riemann-sum approximation
underlying (\ref{SNP5-1}).

\bigskip{}

\begin{remark} Note that (\ref{SNP4-1-1-1}) implies the following
growth condition for $f$:\textit{ For some $C^{\prime}_{\alpha'}\in(0,\infty)$,
\[
|f(x)|\leq C^{\prime}_{\alpha'}\cdot\left(1+|x|^{\alpha'}\right)\:0<\alpha'<\alpha.
\]
}This growth restriction ensures uniform integrability at sub-$\alpha$
orders.

\end{remark} 

\bigskip{}

\noindent\refstepcounter{SSalphaM}\label{lem:SSalphaM}\textbf{Lemma
SS$_{\alpha}$M (smooth sub-$\alpha$ moments).}\textit{ Suppose that: }
\begin{itemize}
\item[$\mathbf{(i)}$] \textit{ $x_{t}$ satisfies}\textbf{\textit{ Assumption INN }}\textit{and
}\textbf{\textit{Assumption SNP{*}}}\textit{; }
\item[$\mathbf{(ii)}$] \textit{ either:}
\begin{itemize}
\item[$\mathbf{(a)}$] \textit{ $f$ is }\textbf{\textit{LH}}\textit{; or}
\item[$\mathbf{(b)}$] \textit{ $f$ is bounded and Lipschitz}
\end{itemize}
\end{itemize}
\textit{Then 
\begin{equation}
\sup_{n\geq1}\sup_{r\in[0,1]}\mathbf{E}\left|f\left(\gamma^{*}_{n}{}^{-1}x_{\lfloor nr\rfloor}\right)\right|<\infty,\label{SNP4-1-1-1-1-2}
\end{equation}
and }
\begin{equation}
\lim_{\delta\downarrow0}\limsup_{n\to\infty}\sup_{\substack{r,s\in[0,1]\\
|r-s|\leq\delta
}
}\mathbf{E}\left|f\left(\gamma^{*}_{n}{}^{-1}x_{\lfloor nr\rfloor}\right)-f\left(\gamma^{*}_{n}{}^{-1}x_{\lfloor ns\rfloor}\right)\right|=0.\label{SNP4-1-1-1-1-3}
\end{equation}
\bigskip{}

It follows from Gikhman and Skorokhod (1969) that (\ref{SNP2-1}),
(\ref{SNP4-1-1-1-1-2}) and (\ref{SNP4-1-1-1-1-3}) are sufficient
for the limit result of (\ref{SNP5-1}). Consequently, (\ref{SNP5-1})
holds whenever \textit{$f$ }\textit{\emph{is }}\textbf{\textit{\emph{LH}}}\textit{\emph{
or bounded and Lipschitz. It can further be shown that, if $f$ locally
integrable and satisfies the polynomial growth condition }}\textit{$\left|f(x)\right|=O\left(\left|x\right|^{\lambda}\right)$
}\textit{\emph{for some}}\textit{ $0<\lambda<\alpha$}\textit{\emph{,
then sample functionals involving $f$ can be approximated by corresponding
sample functionals of a bounded and Lipschitz function, for which
(\ref{SNP5-1}) holds by the preceding argument. This approximation
argument is implicit in the following theorem, which provides the
main limit result for sample functionals of HT-SNP processes.}}

\bigskip{}

\begin{theorem}\label{SNP-1}\textit{ Suppose that: }
\begin{itemize}
\item[$\mathbf{(i)}$] \textit{ $x_{t}$ satisfies}\textbf{\textit{ Assumption INN }}\textit{and
}\textbf{\textit{Assumption SNP{*}}}\textit{; }
\item[$\mathbf{(ii)}$] \textit{ $f:\mathbb{R}\to\mathbb{R}$ is locally integrable;}
\item[$\mathbf{(iii)}$] \textit{ As $\left|x\right|\rightarrow\infty,$ $\left|f(x)\right|=O\left(\left|x\right|^{\lambda}\right)$
for some $0<\lambda<\alpha$;}
\item[(\textbf{iv})] \textit{(a) either for each $t'$ fixed }$\frac{1}{n}\sum^{t'}_{t=1}f(\gamma^{*-1}_{n}x_{t})=o_{p}(1)$\textit{,
or (b) }\textbf{\textit{Assumption INN}}\textit{ holds with $\theta=1$. }
\end{itemize}
\noindent\textit{Then,} \textit{as $n\rightarrow\infty$,} 

\begin{equation}
\frac{1}{n}\sum^{n}_{t=1}f\left(\frac{x_{t}}{\gamma^{*}_{n}}\right)\xrightarrow{d}\int^{1}_{0}f\bigl(\Lambda_{\alpha,d}(r)\bigr)dr.\label{SNP1}
\end{equation}
\textit{Furthermore, suppose that }$\tilde{f}:\mathbb{R}\rightarrow\overline{\mathbb{R}}$\textit{
is such that }$\tilde{f}(x)=f(x)$\textit{ almost everywhere and (a)
either for each $t'$ fixed }$\frac{1}{n}\sum^{t'}_{t=1}\tilde{f}(\gamma^{*-1}_{n}x_{t})=o_{p}(1)$\textit{,
or (b) }\textbf{\textit{Assumption INN}}\textit{ holds with $\theta=1$.
Then (\ref{SNP1}) holds with }$f$\textit{ replaced by }$\tilde{f}$.

\end{theorem}\bigskip{}

\begin{remark} \textbf{Theorem \ref{SNP-1}} requires that the function
$f$ is locally integrable and satisfy a sub-$\alpha$ polynomial
growth rate. This is a technical condition imposed to establish (\ref{SNP1}),
and it accommodates many regression functions used in practice. Under
these conditions, it can also be shown that the limiting random variable
is well defined, in the sense that $\int^{1}_{0}\left|f\bigl(\Lambda_{\alpha,d}(r)\bigr)\right|dr<\infty$
a.s. See Proposition \ref{IntegralExist} in the Appendix and, for
more general conditions, Samorodnisky and Taqqu (1994), Theorem 11.3.2.
A more general version of (\ref{SNP1}) is under development by Kasparis
and Wang (2026). In particular, preliminary arguments indicate that
(\ref{SNP1}) continues to hold when $f$ is locally integrable and
satisfies\textit{ }$\left|f(x)\right|=O\left(\left|x\right|^{\lambda}\right)$,
for some $\lambda>0$ such that $\mathcal{H}>\frac{1}{\alpha}-\frac{1}{\lambda}$
where recall that the Hurst index is $\mathcal{H}=d-1+1/\alpha$.

\end{remark} 

\section{Kernel Regression}\label{sec:Kernel-Regression}

Next, we consider estimation and inference in regression models of
the form
\begin{equation}
y_{t}=m(x_{t-1})+u_{t},\label{model0}
\end{equation}
where $m$ is an unknown regression function and $x_{t},u_{t}$ are
$\mathcal{F}_{t}$-measurable for some appropriate filtration $\mathcal{F}_{t}$,
$t\geq0$ -e.g. \textbf{Assumption MDE}. In particular, we consider
the limit distribution of the kernel regression estimator (local level)
for the function $m(.)$ in (\ref{model0}) for the case where $x_{t}$
is an HT-WNP. This result is complementary to e.g. Wu and Mielniczuk
(2002), who consider kernel regressions with stationary linear process
covariates. In addition, Jeganathan (2004, 2008) develops foundational
limit theory for kernel functionals of HT-SNPs under conditions analogous
to those imposed by \textbf{Assumption SNP{*}} in the present paper.
These results can be adapted with relatively minor modifications to
the asymptotic analysis of kernel regression estimators. In particular,
the aforementioned studies consider nonstationary fractional covariates
of memory $1-1/\alpha<d<2-1/\alpha$, $\alpha\in(0,2]$.\footnote{Notice that for $d=1$ this restrictions are equivalent to $\alpha\in(1,2]$.}

Set $K_{n,t}\left(x\right):=K\left(\frac{x_{t}-x}{h_{n}}\right)$
and consider the Nadaraya-Watson (NW hereafter) kernel regression
estimator
\[
\hat{m}(x)=\frac{\sum^{n}_{t=1}K_{n,t-1}\left(x\right)y_{t}}{\sum^{n}_{t=1}K_{n,t-1}\left(x\right)}.
\]
\bigskip{}

\begin{theorem} \label{MainNonParam}\emph{ Suppose that $\left\{ y_{t}\right\} $
is generated by (\ref{model0}). In addition:}
\begin{itemize}
\item[\textbf{(i)}] \textbf{\emph{ Assumption INN}}\emph{ holds with $\theta=1$ and
$x_{t}$ satisfies either }\textbf{\emph{Assumption BFR}}\emph{ or
}\textbf{\emph{Assumption MI}}\emph{;}
\item[\textbf{(ii)}] \emph{$x_{t}$, $u_{t}$ satisfy }\textbf{\emph{Assumption MDE}}\emph{;}
\item[\textbf{(iii)}] \emph{$K$ and $K^{2}$ satisfy }\textbf{\emph{Assumption K}}\emph{;}
\item[\textbf{(iv)}] \emph{$\{h_{n}\}$ a positive sequence such that $h_{n}+\gamma_{n}(nh_{n})^{-1}\rightarrow0$;}
\item[\textbf{(v)}] \emph{There is some function $\bar{m}(x,z)$ such that $\left\vert m(x+hz)-m(x)\right\vert \leq h^{\nu}\bar{m}(x,z)$,
with $\nu\in(0,1]$ for all $h$ sufficiently small, and $\int_{\mathbb{R}}\bar{m}(x,z)\left|K(z)\right|dz<\infty$
for all $x\in\mathbb{R}$;}
\end{itemize}
\textit{Then, as $n\rightarrow\infty$,}

\noindent\emph{
\[
\hat{m}(x)\overset{p}{\rightarrow}m(x).
\]
 If in addition, $\left(\frac{nh_{n}}{\gamma_{n}}\right)^{1/2}h^{\nu}_{n}\rightarrow0$,
then }
\begin{equation}
\left(\frac{nh_{n}}{\gamma_{n}}\right)^{1/2}\left[\hat{m}(x)-m(x)\right]\overset{d}{\rightarrow}MN\left[0,\frac{\sigma^{2}_{u}}{\varrho(0)\int_{\mathbb{R}}K(z)^{2}dz}\right].\label{NW1}
\end{equation}

\end{theorem}

\bigskip{}

Next, we consider a non parametric t-test for the hypothesis $H_{0}:m(x)=m_{0}(x)$,
where $m_{0}(x)$ is a prespecified function. To construct the test
statistic, we employ the following local estimator of the regression-error
variance, based on nonparametric residuals:
\begin{equation}
\hat{\sigma}^{2}_{u}\left(x\right)=\frac{\sum^{n}_{k=1}\left[y_{t}-\hat{m}(x)\right]^{2}K_{n,t-1}\left(x\right)}{\sum^{n}_{k=1}K_{n,t-1}\left(x\right)},
\end{equation}
can be shown to be consistent for $\sigma^{2}_{u}$ under our assumptions.
Set $\hat{t}(x,m_{0})$
\[
\hat{t}(x,m_{0}):=\frac{\sum^{n}_{k=1}K_{n,t-1}\left(x\right)\left[\hat{m}(x)-m_{0}(x)\right]}{\sqrt{\hat{\sigma}^{2}_{u}\left(x\right)\left[\sum^{n}_{k=1}K^{2}_{n,t-1}\left(x\right)\right]}}.
\]
Under the null hypothesis the t-statistic has approximately standard
normal limit distribution. This is demonstrated by the following result.\bigskip{}

\begin{theorem} \label{NonParam_t-test} \emph{Suppose that in addition
to the conditions of }\textbf{\emph{Theorem \ref{MainNonParam}}}\emph{,
there exists some $p>1$ such that:}
\begin{itemize}
\item[\textbf{(i)}]  $K>0$ \textit{and}\emph{ $K^{p}$ and $K^{2p}$ satisfy }\textbf{\emph{Assumption
K}}\emph{;}
\item[\textbf{(ii)}] \emph{ $\mathrm{sup}_{n\geq1,1\leq t\leq n}\mathbf{E}\left[\left|u_{t}\right|^{2p}\mid\mathcal{F}_{t-1}\right]<\infty$,
$a.s.$,}
\end{itemize}
\emph{Then for each $x$, under $H_{0}$,}
\[
\hat{t}(x,m_{0})\overset{d}{\rightarrow}N(0,1).
\]
\end{theorem}

\bigskip{}

\begin{itemize}
\item \textbf{Comparison with the stationary case}: For the case where $x_{t}$
is a stationary linear process on i.i.d. innovations (possibly heavy
tailed), under some additional technical conditions on the smoothness
of $m$ it can be shown (cf. Wu and Mielniczuk, 2002; p. 1451) that
(\ref{NW1}) holds with $\gamma_{n}=1$ and $\varrho(0)$ replaced
with $\Phi_{x_{1}}(0)$ i.e. the density of $x_{1}$ at the origin. 
\item \textbf{Comparison with the strongly nonstationary case}: For $x_{t}$
fractional $1-1/\alpha<d<2-1/\alpha$, $\alpha\in(0,2]$, under smoothness
requirements similar to those of condition (v) of the Theorem above,
(\ref{NW1}) holds (cf. Jeganathan, 2008; Theorem 4) with $\gamma_{n}=\gamma^{*}_{n}$
and $\varrho(0)$ replaced by $L_{\Lambda_{\alpha,d}}(0,1)$ i.e.
the local time of the fractional stable motion $\Lambda_{\alpha,d}$
--see (\ref{SNP2}).
\end{itemize}
\bigskip{}

These results show that the kernel regression estimator has conventional
limit distribution irrespective of the stationarity properties of
the data and moments properties of the covariate. Therefore, kernel
based inference is conventional.

Finally, we consider the case of heavy-tailed regression errors satisfying
\textbf{Assumption HTE}. We restrict attention to WNP covariates satisfying
either \textbf{Assumption BFR2} or \textbf{Assumption MI}. Before
stating the corresponding limit result, we introduce some notation.
Set 
\[
\tau_{1}:=\frac{\sigma^{\alpha_{u}}_{u}\varrho(0)^{1-\alpha_{u}}\int_{\mathbb{R}}\left\vert K(x)\right\vert ^{\alpha_{u}}dx}{\left(\int_{\mathbb{R}}K(x)dx\right)^{\alpha_{u}}},\quad\tau_{2}:=\frac{\int_{\mathbb{R}}K(x)^{\left\langle \alpha_{u}\right\rangle }dx}{\int_{\mathbb{R}}\left\vert K(x)\right\vert ^{\alpha_{u}}dx}.
\]
These quantities are the scale and skewness parameters of the limit
distribution of the rescaled score term asscociated with the NW estimator.

\bigskip{}

\begin{theorem} \label{MainNonParam-HTE}\emph{ Suppose that $\left\{ y_{t}\right\} $
is generated by (\ref{model0}). In addition:}
\begin{itemize}
\item[\textbf{(i)}] \textbf{\emph{ Assumption INN}}\emph{ holds with $\theta=1$ and
$x_{t}$ satisfies either }\textbf{\emph{Assumption BFR2}}\emph{ or
}\textbf{\emph{Assumption MI}}\emph{;}
\item[\textbf{(ii)}] \emph{$x_{t}$, $u_{t}$ satisfy }\textbf{\emph{Assumption HTE}}\emph{;}
\item[\textbf{(iii)}] \emph{$K$ and $\left|K\right|^{\alpha_{u}}$ satisfy }\textbf{\emph{Assumption
K}}\emph{, and $\int K>0;$ }
\item[\textbf{(iv)}] \emph{$\{h_{n}\}$ a positive sequence such that h$_{n}+$$\gamma_{n}(nh_{n})^{-1}\rightarrow0$;}
\item[\textbf{(v)}] \emph{There is some function $\bar{m}(x,z)$ such that $\left\vert m(x+hz)-m(x)\right\vert \leq h^{\nu}\bar{m}(x,z)$,
with $\nu\in(0,1]$ for all $h$ sufficiently small, and $\int_{\mathbb{R}}\bar{m}(x,z)\left|K(z)\right|dz<\infty$
for all $x\in\mathbb{R}$.}
\end{itemize}
\textbf{\textit{\emph{(a)}}}\textit{ Suppose that $\alpha_{u}>1$.
Then, as $n\rightarrow\infty$,}

\noindent\emph{
\[
\hat{m}(x)\overset{p}{\rightarrow}m(x).
\]
 If in addition $\left(\frac{nh_{n}}{\gamma_{n}}\right)^{1-1/\alpha_{u}}h^{\nu}_{n}\rightarrow0$,
then }
\begin{equation}
\left(\frac{nh_{n}}{\gamma_{n}}\right)^{1-1/\alpha_{u}}\left[\hat{m}(x)-m(x)\right]\overset{d}{\rightarrow}S_{\alpha_{u}}\left\{ \tau_{1},\beta^{\ast}_{u}\tau_{2}\right\} .\label{NW1-HTE-1}
\end{equation}
\textbf{\textit{\emph{(b) }}}\textit{Suppose that $\alpha_{u}=1$.
Then, as $n\rightarrow\infty$,
\begin{equation}
\hat{m}(x)\overset{d}{\rightarrow}m(x)+S_{\alpha_{u}}\left\{ \tau_{1},\beta^{\ast}_{u}\tau_{2}\right\} .\label{NW1-HTE-2}
\end{equation}
}\textbf{\textit{\emph{(c) }}}\textit{Suppose that $\alpha_{u}<1$.
Then, as $n\rightarrow\infty$,
\begin{equation}
\left(\frac{nh_{n}}{\gamma_{n}}\right)^{1-1/\alpha_{u}}\hat{m}(x)\overset{d}{\rightarrow}S_{\alpha_{u}}\left\{ \tau_{1},\beta^{\ast}_{u}\tau_{2}\right\} .\label{NW1-HTE-3}
\end{equation}
}

\end{theorem}

\bigskip{}

\begin{remark} \textbf{Theorem \ref{MainNonParam-HTE}} highlights
important limitations of the NW estimator in the presence of heavy
tailed regression errors. In particular, the estimator is consistent
only when $\alpha_{u}>1$. For $\alpha_{u}=1$, it converges in distribution
to a nondegenerate random pseudo-true limit, whereas for $\alpha_{u}<1$,
it diverges in absolute value. Even when consistency holds, the convergence
rate is slower than the conventional Gaussian rate for $\alpha_{u}\in(1,2)$,
and the limiting distribution is non-Gaussian and depends on unknown
nuisance parameters, including the tail index, scale, and skewness
of the regression errors. Consequently, conventional Gaussian-based
inference is generally invalid, while inference based on the stable
limiting distribution is non-pivotal, complicating the construction
of valid statistical tests.

\end{remark}

\section{IV Methods for Parametric Estimation}\label{sec:IV-Methods}

The focus of this section is the parametric estimation of regressions
of the form
\begin{equation}
y_{t}=\mu+\beta f(x_{t-1})+u_{t}.\label{parametric model IV}
\end{equation}
As in the previous section, $x_{t}$ and $u_{t}$ are assumed to be
$\mathcal{F}_{t}$-measurable with respect to an appropriate filtration
$\mathcal{F}_{t}$, $t\geq0$, and $\left\{ u_{t},\mathcal{F}_{t}\right\} $
forms a martingale difference sequence. The regressions function $f(x)$
is assumed to be known, whereas the parameter vector $\mathbf{b}:=\left[\mu,\beta\right]^{\prime}\in\mathbb{R}^{2}$
is unknown.

We consider an IV estimator that uses reduced-growth transformations,
or weights, of the covariates as instruments. We refer to this estimator
as the weighted IV estimator (WIV hereafter). Comparable estimators
have been considered by Blattberg and Sargent (1971), Samorodnitsky,
Kurz-Kim, and Rachev (2007), and Kim and Meddahi (2020), among others.
Samorodnitsky et al. (2007) employ signed-power functions as weights
and show that, under assumptions more restrictive than those imposed
in the present work, the resulting estimator is efficient. In our
framework, when the covariate is stationary, reduced-growth transformations
ensure the existence of sufficiently high moments for an LLN to apply;
see, for example, Remark \ref{Rmk6aa}. Similarly, such transformations
facilitate the development of our limit theory when the covariate
is either an an HT-WNP or an HT-SNP.

It should be emphasised, however, that the WIV estimator does not
always support conventional inference under strong nonstationarity.
To address this limitation, we consider a more sophisticated IV method
in the subsequent section. The WIV estimator of the slope parameter
in (\ref{Model0}) takes the form 
\begin{equation}
\mathbf{\tilde{\mathbf{b}}}_{WIV}=\left[\sum^{n}_{t=1}\mathbf{w}_{t-1}\mathbf{x}^{\prime}_{t-1}\right]^{-1}\left[\sum^{n}_{t=1}\mathbf{w}_{t-1}y_{t}\right],\label{WIV}
\end{equation}
with $\mathbf{x}^{\prime}_{t}=[1,f(x_{t})]$ and $\mathbf{w}^{\prime}_{t}=[w_{0}(x_{t}),w_{1}(x_{t})]$.
In particular, we assume that the weight functions $w_{i}(.)$, $i=\{0,1\}$
satisfy the following growth
\begin{equation}
\lvert f(x)w_{i}(x)\rvert=O\left(\left\vert x\right\vert ^{\lambda_{i}}\right),\text{ }\lambda_{i}\in(0,\alpha).\label{weight1}
\end{equation}
For example, the\textit{ Cauchy} estimator (see e.g. So and Shin (1999),
Phillips, Park and Chang (2004), Choi, Jacewitz and Park (2016), Kim
and Meddahi (2020) among others) is a special case of the WIV estimator
with $w_{1}(x)=$sgn$(x)$.

\bigskip{}

It follows from (\ref{WIV}) that the WIV estimator involves nonlinear
transformations through both the weight functions and the regression
function $f$. For the case where $x_{t}$ is a nonstationary process,
we consider a class of functions introduced by Park and Phillips (1999,
2001). In particular, we focus on locally integrable functions that
are not globally integrable. Although \textbf{Theorem \ref{MainLL_kernel}}
(with $h_{n}=1$) also enables the analysis to be extended to integrable
transformations, we restrict attention to the nonintegrable case for
brevity and because such transformations arise more frequently in
applications. In particular, we consider functions that are asymptotically
homogeneous, as defined below (see also DK2021).\bigskip{}

\noindent\textbf{Definition 4 (Asymptotically Homogeneous Function)}. 
\begin{itemize}
\item[\textbf{(i)}] \textbf{\emph{ }}\emph{Let $\left\{ x_{t}\right\} $ a random sequence
and $\mu_{n}$ a strictly positive sequence. $g:\mathbb{R}\rightarrow\overline{\mathbb{R}}$
is an }\textbf{\emph{asymptotically homogeneous function}}\emph{ for
$\left\{ x_{t},\mu_{n}\right\} $ (}\textbf{\emph{AHF}}\emph{$\left\{ x_{t},\mu_{n}\right\} $,
hereafter), if for all $\lambda>0$ it admits the decomposition 
\[
g(x)=k_{g}(\lambda)g_{H}(x/\lambda)+R_{g}(x,\lambda),
\]
where $k_{g}:\mathbb{R}^{*}_{+}\rightarrow\mathbb{R}^{*}_{+}$, $g_{H}:\mathbb{R}\rightarrow\overline{\mathbb{R}}$
is locally integrable, $R_{g}:\mathbb{R\times}\mathbb{R}^{*}_{+}\rightarrow\overline{\mathbb{R}}$
and,
\[
R_{g,n}:=\frac{1}{k_{g}(\mu_{n})n}\sum^{n}_{t=1}\mathbf{E}\left|R_{g}(x_{t},\mu_{n})\right|=o(1),
\]
as $n\rightarrow\infty$. $g_{H}$ and $k_{g}$ are the limit homogeneous
function and asymptotic order of $g$ respectively.}
\item[\textbf{(ii)}] \textbf{\emph{ }}\emph{If $R_{g}(x,\lambda)=0$, }$g$\emph{ is a
}\textbf{\emph{homogeneous function}}\emph{ (}\textbf{\emph{HF}}\emph{,
hereafter).}
\item[\textbf{(iii)}] \textbf{\emph{ }}\emph{A matrix $\mathbf{A}$ is }\textbf{\emph{AHF}}\emph{$\left\{ x_{t},\mu_{n}\right\} $,
if its elements are}\textbf{\emph{ AHF}}\emph{$\left\{ x_{t},\mu_{n}\right\} $.
}\bigskip{}
\end{itemize}
\noindent Next, set 
\[
\mathbf{A}(x):=\left[\begin{array}{cc}
w_{0}(x) & \left(w_{0}\cdot f\right)(x)\\
w_{1}(x) & \left(w_{1}\cdot f\right)(x)
\end{array}\right]\text{, }\mathbf{B}(x):=\left[\begin{array}{c}
\left(w_{0}\cdot f\right)(x)\\
\left(w_{1}\cdot f\right)(x)
\end{array}\right]\text{,}
\]
\[
\mathbf{W}(x):=\left[\begin{array}{cc}
w^{2}_{0}(x) & \left(w_{0}\cdot w_{1}\right)(x)\\
\left(w_{0}\cdot w_{1}\right)(x) & w^{2}_{1}(x)
\end{array}\right],
\]
Under nonstationarity, we assume that the elements of the matrices
above are \textbf{AHF}, for some appropriate random array and sequence.
We denote the corresponding limit homogeneous matrix-valued functions
by $\mathbf{A}_{H}(x)$, $\mathbf{B}_{H}(x)$, $\mathbf{W}_{H}(x)$. 

Moreover, we consider following t-statistic for the null hypothesis
$H_{0}:\beta=\beta_{0}$,

\[
\tilde{t}_{WIV}=\frac{\tilde{\beta}_{WIV}-\beta_{0}}{\sqrt{\tilde{\sigma}^{2}_{u}\left[\left[\sum^{n}_{t=1}\mathbf{A}(x_{t-1})\right]^{-1}\sum^{n}_{t=1}\mathbf{W}(x_{t-1})\left[\sum^{n}_{t=1}\mathbf{A^{\prime}}(x_{t-1})\right]^{-1}\right]_{22}}},
\]
where

\begin{equation}
\tilde{\sigma}^{2}_{w}=\frac{\sum^{n}_{k=1}\left[y_{t}-\tilde{m}_{t-1}\right]^{2}w_{u}(x_{t-1})}{\sum^{n}_{k=1}w_{u}(x_{t-1})},\text{ }w_{u}(x)\geq0\label{sigma_u}
\end{equation}
and $\tilde{m}_{t-1}:=\mathbf{x}^{\prime}_{t-1}\mathbf{\tilde{\mathbf{b}}}$.
In addition, to simplify our exposition we assume that the weight
$w_{u}$ is an HF. i.e. for all $\lambda>0$,
\[
w_{u}(\lambda x)=k_{w_{u}}(\lambda)\text{\ensuremath{w_{u}}}(x).
\]
The decomposition above is relevant only when $x_{t}$ is nonstationary.
The following two results characterize the limiting behavior of $\tilde{t}_{WIV}$
is an SP and a WNP, respectively. In both cases inference is asympotically
Gaussian.

\bigskip{}

\begin{theorem} \label{WIV_SP}\emph{Suppose that $y_{t}$ is generated
by (\ref{parametric model IV}). }

\noindent\textbf{(a)} \emph{In addition, suppose that the following
conditions hold:}
\begin{itemize}
\item[\textbf{(i)}]  \textbf{\emph{Assumption INN}}\emph{ holds;}
\item[\textbf{(ii)}]  \emph{$x_{t}$ satisfies }\textbf{\emph{Assumption SP}}\emph{;}
\item[\textbf{(iii)}]  $\mathbf{E}\left\Vert \mathbf{A}(x_{t})\right\Vert $\emph{, $\mathbf{E}\left\Vert \mathbf{W}(x_{t})\right\Vert ^{1+\delta}<\infty$,
for some $\delta>0$;}
\item[\textbf{(iv)}] \textbf{\emph{ Assumption MDE}}\emph{ holds;}
\end{itemize}
\noindent\textit{Then}\emph{ as $n\rightarrow\infty$,}
\[
\sqrt{n}\left(\mathbf{\tilde{\mathbf{b}}}_{WIV}-\mathbf{b}\right)\stackrel{d}{\rightarrow}\left[\mathbf{E}\left(\mathbf{A}(x_{t})\right)\right]^{-1}\cdot N\left[\mathbf{0},\sigma^{2}_{u}\mathbf{E}\left(\mathbf{W}(x_{t})\right)\right].
\]

\noindent\textbf{\textit{\emph{(b)}}}\textit{ }\emph{In addition,
suppose that for some $1<p<2$ the following conditions hold:}
\begin{itemize}
\item[\textbf{(v)}] \emph{ $\mathrm{sup}_{n\geq1,1\leq t\leq n}\mathbf{E}\left[\left|u_{t}\right|^{2p}\mid\mathcal{F}_{t-1}\right]<\infty$
a.s. }
\item[\textbf{(vi)}] \emph{ $\mathbf{E}\left[f(x_{t})^{2}w_{u}(x_{t})\right]$, $\mathbf{E}\left[w_{u}(x_{t})^{p}\right]<\infty$. }
\end{itemize}
\noindent\textit{Then }\emph{under $H_{0}$, as $n\rightarrow\infty$,}
\[
\hat{t}_{WIV}\overset{d}{\rightarrow}N(0,1).
\]
\end{theorem}

\bigskip{}

The following result gives the limit distribution of the WIV estimator
when the covariate is an HT-WNPs process.\bigskip{}

\noindent\begin{theorem} \emph{\label{MainParametric} Suppose that
$\left\{ y_{t}\right\} $ is generated by (\ref{parametric model IV}). }

\noindent\textbf{(a)} \emph{In addition:}
\begin{itemize}
\item[\textbf{(i)}] \textbf{\emph{ }}\textbf{\emph{Assumption INN}}\emph{ holds with
$\theta=1$ and $x_{t}$ satisfies either }\textbf{\emph{Assumption
BFR}}\emph{ or }\textbf{\emph{Assumption MI}}\emph{;}
\item[\textbf{(ii)}]  \emph{$x_{t}$, $u_{t}$ satisfy }\textbf{\emph{Assumption MDE}}\emph{;}
\item[\textbf{(iii)}] $\mathbf{A}(x)$\emph{, $\mathbf{W}(x)$ are }\textbf{\emph{AHF}}\emph{$\left\{ x_{t},\gamma_{n}\right\} $
and their elements are $O\left(\left\vert x\right\vert ^{\lambda}\right)$,
$\lambda\in(0,\alpha)$;}
\item[\textbf{(iv)}] \emph{ The inverse of }$\int_{\mathbb{R}}\mathbf{A}_{H}(x)\varrho(x)dx$
\emph{exists $a.s.$ }
\end{itemize}
\noindent\emph{Then as }$n\rightarrow\infty$,
\[
\sqrt{n}diag\{1,k_{f}(\gamma_{n})\}\left(\mathbf{\tilde{\mathbf{b}}}_{WIV}-\mathbf{b}\right)
\]
\begin{equation}
\overset{d}{\rightarrow}\left[\int_{\mathbb{R}}\mathbf{A}_{H}(x)\varrho(x)dx\right]^{-1}\cdot MN\left[\mathbf{0},\sigma^{2}_{u}\int_{\mathbb{R}}\mathbf{W}_{H}(x)\varrho(x)dx\right].\label{mainparam}
\end{equation}

\noindent\textbf{(b)} \emph{In addition to the conditions (i)-(iv),
suppose there is some $1<p<2$ such that:}
\begin{itemize}
\item[\textbf{(v)}] \emph{$w_{u}$ is an }\textbf{\emph{HF }}\emph{and $fw_{u}$, $f^{2}w_{u}$
are }\textbf{\emph{AHF}}\emph{$\left\{ x_{t},\gamma_{n}\right\} $;}
\item[\textbf{(vi)}] $f^{2}w_{u},\text{\ensuremath{w^{p}_{u}}}$\emph{, are of order }$O\left(\left\vert x\right\vert ^{\lambda}\right)$,
$\lambda\in(0,\alpha)$\emph{;} 
\item[\textbf{(vii)}]  $\mathrm{sup}_{n\geq1,1\leq t\leq n}\mathbf{E}\left[\left|u_{t}\right|^{2p}\mid\mathcal{F}_{t-1}\right]<\infty$
$a.s.$
\end{itemize}
\emph{Then as $n\rightarrow\infty$, }

\[
\tilde{t}_{WIV}\overset{d}{\rightarrow}N(0,1).
\]

\noindent\end{theorem}

\bigskip{}

\noindent\textbf{Example 1} (Cauchy estimator). Suppose that $f(x)=x$
and $w_{1}(x)=$sgn$(x)$, $w_{0}(x)=1$. Then (\ref{mainparam})
holds for all $1<\alpha\leq2$. In particular, the IV estimator for
the slope parameter $\beta$ satisfies
\[
\sqrt{n}\gamma_{n}\left(\tilde{\beta}_{WIV}-\beta\right)\overset{d}{\rightarrow}MN\left[0,\sigma^{2}_{u}\frac{1-\left[\int_{\mathbb{R}}\textrm{sgn}(x)\varrho(x)dx\right]^{2}}{\left[\int_{\mathbb{R}}\left|x\right|\varrho(x)dx-X^{-}\int_{\mathbb{R}}\textrm{sgn}(x)\varrho(x)dx\right]^{2}}\right]
\]

\bigskip{}

\noindent\textbf{Example 2}. Suppose that $f(x)=x$ and $w_{1}(x)=$sgn$(x)\left\vert x\right\vert ^{-b_{1}}$,
$w_{0}(x)=$sgn$(x)\left\vert x\right\vert ^{-b_{0}}$ with $0<2b_{1},2b_{0}<1$,
$b_{1}\neq b_{0}$. Then (\ref{mainparam}) holds for all $1-\left(b_{1}\wedge b_{0}\right)<\alpha\leq2$.

\bigskip{}

The following result establishes the convergence rate of the WIV estimator
when the covariate is an SNP process satisfying \textbf{Assumption
SNP{*}}. In particular, it shows that WIV attains a parametric rate
of convergence and is therefore more suitable for estimation than
IVX. Recall that, under nonstationarity, IVX generally attains only
semiparametric convergence rates, see e.g. Phillips and Magdalinos
(2009), KMS2015. Moreover, this convergence-rate result is used to
establish the consistency of an estimator of $\sigma^{2}_{u}$ based
on WIV residuals. This variance estimator is subsequently employed
in the studentization of IVX-based t-statistics.\bigskip{}

The following result characterizes the limiting behavior of the WIV
estimator when the regression errors satisfy Assumption HTE and the
covariate is a weakly nonstationary process (WNP). To simplify the
exposition, we exclude BFR1 processes and restrict attention to homogeneous
weight functions (\textbf{HF}).

\bigskip{}

\noindent\begin{theorem} \emph{\label{MainParametric-HTE} Suppose
that $\left\{ y_{t}\right\} $ is generated by (\ref{parametric model IV}).
In addition:}
\begin{itemize}
\item[\textbf{(i)}] \textbf{\emph{ }}\textbf{\emph{Assumption INN}}\emph{ holds with
$\theta=1$ and $x_{t}$ satisfies either }\textbf{\emph{Assumption
BFR2}}\emph{ or }\textbf{\emph{Assumption MI}}\emph{;}
\item[\textbf{(ii)}]  \emph{$x_{t}$, $u_{t}$ satisfy }\textbf{\emph{Assumption HTE}}\emph{;}
\item[\textbf{(iii)}]  $\mathbf{A}(x)$\emph{, is }\textbf{\emph{AHF}}\emph{$\left\{ x_{t},\gamma_{n}\right\} $
and its elements are $O\left(\left\vert x\right\vert ^{\lambda}\right)$,
$\lambda\in(0,\alpha)$;}
\item[\textbf{(iv)}]  \emph{for $i=\left\{ 0,1\right\} $, $w_{i}(x)$, $\left|w_{i}(x)\right|^{\alpha_{u}}$,
is }\textbf{\emph{HF}}\emph{ and }$\left|w_{i}(x)\right|^{\alpha_{u}}=O\left(\left|x\right|^{\lambda_{i}}\right)$\textit{,}
\emph{$\lambda_{i}\in(0,\alpha)$;} 
\item[\textbf{(v)}] \emph{ The inverse of }$\int_{\mathbb{R}}\mathbf{A}_{H}(x)\varrho(x)dx$
\emph{exists. }
\end{itemize}
\noindent\emph{Then as }$n\rightarrow\infty$\textit{,}
\begin{equation}
n^{1-1/\alpha_{u}}diag\{1,k_{f}(\gamma_{n})\}\left(\mathbf{\tilde{\mathbf{b}}}_{WIV}-\mathbf{b}\right)\overset{d}{\rightarrow}\left[\int_{\mathbb{R}}\mathbf{A}_{H}(x)\varrho(x)dx\right]^{-1}\cdot\mathbf{S}_{u},\label{MainParametric-HTE-1}
\end{equation}
\textit{where for $\mathbf{v}\in\mathbb{R}^{2}$,
\[
\mathbf{E}\left(e^{i\mathbf{v}^{\prime}\mathbf{S}_{u}}\right)=\exp\left\{ -\sigma^{\alpha_{u}}_{u}\left[\int_{\mathbb{R}}\left|\mathbf{v}^{\prime}\mathbf{w}(x)\right|^{\alpha_{u}}\left[1+i\beta^{\ast}_{u}\textrm{sgn}\left(\mathbf{v}^{\prime}\mathbf{w}(x)\right)\right]\tan\left(\frac{\pi\alpha_{u}}{2}\right)\varrho(x)dx\right]\right\} .
\]
}

\noindent\end{theorem}

\bigskip{}

\begin{remark} \textbf{\label{HTEparam}} The limit result of (\ref{MainParametric-HTE-1})
shows that, in contrast to (\ref{mainparam}), the WIV estimator generally
has a non-Gaussian and non-pivotal stable limiting distribution involving
the nuisance parameters $\sigma_{u},\alpha_{u},\beta^{\ast}_{u}$.
Moreover, for $\alpha_{u}<2$, the convergence rates are slower than
those obtained under \textbf{Assumption MDE}. An important implication
of this result is that the consistency of the WIV estimator depends
on the relative strength of the regression signal and the heavy-tailed
noise. In particular, when $\alpha_{u}\leq1$, the intercept estimator
is inconsistent. For $\alpha_{u}=1$, its estimation error converges
to a nondegenerate stable distribution, whereas for $\alpha_{u}<1$,
the estimation error diverges in absolute value in probability. By
contrast, the slope estimator may remain consistent even when $\alpha_{u}\leq1$.
Specifically, provided its normalized limiting distribution is nondegenerate,
consistency holds if and only if

\[
n^{1-1/\alpha_{u}}k_{f}(\gamma_{n})\rightarrow\infty
\]
Thus, sufficiently strong growth of the regression function can compensate
for the heavy-tailed regression disturbances, preserving consistency
of the slope estimator even when the intercept estimator is inconsistent.

\end{remark}

\bigskip{}

\begin{theorem} \label{MainParametricWLS} \emph{Suppose that $\left\{ y_{t}\right\} $
is generated by (\ref{model0}). }

\noindent\textbf{(a)}\emph{ In addition suppose that the following
hold:}
\begin{itemize}
\item[\textbf{(i)}] \textbf{ }\textbf{\emph{$x_{t}$ }}\emph{satisfies }\textbf{\emph{Assumption
INN}}\emph{ and }\textbf{\emph{Assumption SNP{*}}}\emph{;}
\item[\textbf{(ii)}]  $\mathbf{A}(x)$\emph{, $\mathbf{W}(x)$ are }\textbf{\emph{AHF}}\emph{$\left\{ x_{t},\gamma^{*}_{n}\right\} $
and their elements satisfy the assumptions imposed on $f$ in }\textbf{\emph{Theorem
\ref{SNP-1}}}\emph{;}
\item[\textbf{(iii)}] \textbf{$x_{t}$} \emph{and $u_{t}$ satisfy }\textbf{\emph{Assumption
MDE.}}
\end{itemize}
\emph{Then as }$n\rightarrow\infty$,

\noindent
\[
\sqrt{n}diag\{1,k_{f}(\gamma^{\ast}_{n})\}\left(\mathbf{\tilde{\mathbf{b}}}_{WIV}-\mathbf{b}\right)=O_{p}(1).
\]

\textbf{(b)}\emph{ Suppose that in addition to the conditions (i)-(iii)
there is some $1<p<2$ such that:}
\begin{itemize}
\item[\textbf{(iv)}]  \emph{$w_{u}$ is an }\textbf{\emph{HF}}\emph{ and $fw_{u}$, $f^{2}w_{u}$
are }\textbf{\emph{AHF}}\emph{$\left\{ x_{t},\gamma^{*}_{n}\right\} $;}
\item[\textbf{(v)}]  $w_{u}f^{2}$, $w_{u}f$, $w^{p}_{u}$\emph{ satisfy the assumptions
imposed on $f$ in }\textbf{\emph{Theorem \ref{SNP-1}}}\emph{;}
\item[\textbf{(vi)}] \emph{$\mathrm{sup}_{n\geq1,1\leq t\leq n}\mathbf{E}\left[\left|u_{t}\right|^{2p}\mid\mathcal{F}_{t-1}\right]<\infty$
a.s.}
\end{itemize}
\emph{Then, as }$n\rightarrow\infty$,
\[
\tilde{\sigma}^{2}_{w}\stackrel{p}{\rightarrow}\sigma^{2}_{u}.
\]

\end{theorem}

\bigskip{}

\begin{remark} \textbf{\label{Rmk7}} \textbf{(WIV with SNP Covariates)}.
\textbf{(a)} It least for the case where FCLTs apply, it can demonstrated
directly that the IV estimator has limit distribution determined by
stochastic integrals. Therefore, in this case inference is non-conventional
and test statistics have limit distributions that depend on nuisance
parameters. To see this suppose that weak limit of (\ref{SNP2}) can
be strengthened to weak convergence in the $D[0,1]$ space i.e.
\[
\gamma^{\ast-1}_{n}x_{\lfloor nr\rfloor}\Rightarrow\Lambda_{\alpha,d}(r),\text{ }r\in[0,1].
\]
This is true if $\alpha=2$ or $d\geq1$, see e.g. Astrauskas (1983),
Avram and Taqqu (1992). Arguments similar to those used in Park and
Phillips (2001), Robsinon and Hualde (2003)
\begin{equation}
\sqrt{n}diag\{1,k_{f}(\gamma^{\ast}_{n})\}\left(\mathbf{\tilde{\mathbf{b}}}_{WIV}-\mathbf{b}\right)\overset{d}{\rightarrow}\left[\int^{1}_{0}\mathbf{A}_{H}(\Lambda_{\alpha,d}(r))dr\right]^{-1}\cdot\int^{1}_{0}\mathbf{B}_{H}(\Lambda_{\alpha,d}(r))dU(r),\label{IVSNPs}
\end{equation}
where $U$ is a Brownian Motion relating to the regression error $u_{t}$.
Notice that stochastic integral above mixed normal only under very
restrictive conditions i.e. $\Lambda_{\alpha,d}\perp U$. In practice,
this requires strong exogeneity in the model.

\textbf{(b)} Some preliminary results (cf. Kasparis and Wang (2026))
indicate that for $0<\alpha<2$ the $\mathbf{\tilde{\mathbf{b}}}_{WIV}$
estimator has a mixed normal distribution with the associate score
term having asymptotic distribution
\[
MN\left[\mathbf{0},\sigma^{2}_{u}\int^{1}_{0}\mathbf{B}^{2}_{H}(\Lambda_{\alpha,d}(r))dr\right].
\]
An intuitive explanation is provided by the following result. As shown,
for example, by Jeganathan (2008; eq. (133)), for $0<\alpha<2$, there
is some SV $\varsigma(x)$ satisfying $\varsigma(x)^{\alpha}\sim H\left(x^{1/\alpha}\varsigma(x)\right)$
such that, as $n\rightarrow\infty$,
\[
\left[\frac{1}{n^{1/\alpha}\varsigma(n)}\sum^{\left\lfloor nr\right\rfloor }_{t=1}\xi_{t},\frac{1}{\sqrt{n}}\sum^{\left\lfloor nr\right\rfloor }_{t=1}u_{t}\right]\stackrel{f.d.d.}{\rightarrow}\left[Z_{\alpha}(r),U(r)\right],
\]
where $Z_{\alpha}$ is an $\alpha$-stable Lévy process, $U$ a Brownian
motion and $Z_{\alpha}\bot U$. This asymptotic independence suggests
that, in this setting, the WIV estimator may support conventional
inference. Nevertheless, two qualifications are important:

(i) The WIV estimator is not robust in general. In particular, if
$\xi_{t}$ belongs to the domain of attraction of a normal law, the
resulting inference is nonconventional and is governed by stochastic-integral
limits, as in (\ref{IVSNPs}).

(ii) The simulation results presented in this work show that WIV-based
t-tests exhibit substantial finite-sample size distortions when the
stability parameter satisfies $\alpha>1$. \end{remark}

\section{IVX for Parametric Inference }\label{sec:IVX-Inference}

The results of the previous section demonstrate that the WIV estimator
does not, in general, support conventional inference when the covariate
is an HT-SNP. To address this limitation, we develop a weighted extension
of IVX (W-IVX) that facilitates conventional inference for the slope
parameter in (\ref{parametric model IV}). A central feature of IVX
is the construction of instruments that are either exactly or approximately
mildly integrated (MI), depending on the persistence properties of
the original covariate. This construction reduces the degree of persistence
of the instruments and facilitates asymptotically Gaussian inference.
In the present heavy-tailed framework, however, the MI structure alone
is not sufficient to ensure the desired limiting behavior of the relevant
sample functionals. The weight functions are therefore introduced
to control the growth of functionals of the IVX instruments, ensuring,
in particular, that condition (ii) of \textbf{Theorem \ref{MainLL}}
is satisfied. This allows us to exploit the limit theory developed
in Section \ref{sec:Limit-Theory-WNPs}. and establish conventional
inference for the resulting W-IVX estimator. In this respect, W-IVX
can be viewed as the semiparametric counterpart of WIV, in much the
same way that the IVX estimator of Phillips and Magdalinos (2009)
provides a semiparametric counterpart to OLS in conventional linear
predictive regressions. The distinction reflects a trade-off between
estimation efficiency and robustness of inference: while WIV typically
achieves faster, parametric convergence rates, W-IVX sacrifices some
of this convergence speed in exchange for conventional, asymptotically
pivotal inference across a broader range of persistence regimes.

IVX methods were introduced by Phillips and Magdalinos (2009) and
further developed by KMS2015 and Magdalinos and Petrova (2022), among
others. These studies show that IVX-based tests support conventional
inference for a wide range of regressors. Specifically, the original
IVX results apply to linear models in which the regressors are finite-variance
nearly integrated or mildly integrated processes; see Phillips and
Magdalinos (2009). KMS2015 extend these results to stationary short-memory
regressors, while Magdalinos and Petrova (2022) further generalize
the method to accommodate potentially explosive or mildly explosive
covariates driven by short-memory innovations. Extensive simulation
studies indicate that IVX-based tests exhibit very good finite-sample
performance and outperform competing methods in most scenarios. For
further details, see KMS2015 and Magdalinos and Petrova (2022).

Existing studies, however, do not consider fractionally integrated
or heavy-tailed regressors or nonlinear regression functions. Thus,
the validity of IVX inference with fractionally integrated regressors
has not previously been established. Indeed, Andersen and Varneskov
(2021, Remark 9) describe their LCM methodology as an analogue of
IVX for general multivariate fractionally integrated processes and
appear to suggest that IVX does not, in general, deliver Gaussian
inference in fractional settings:\bigskip{}

\textquotedblleft \emph{...however, {[}for IVX{]} the dependent variable
is restricted from exhibiting long-range dependence. Whereas the IVX
inference is not Gaussian uniformly in the persistence regimes, the
corresponding Wald significance test is robust to the specific form
of persistence.}\textquotedblright \bigskip{}

In this study, we extend the IVX framework by allowing for nonlinearities
in both the regression and instrument functions. We also accommodate
heavy-tailed, SNP, MI and SP fractional regressors. Contrary to the
suggestion that IVX does not support Gaussian inference in fractional
settings, we show that the resulting IVX-based t-statistics have standard
Gaussian limit distributions even when the regressors are both fractionally
integrated and heavy-tailed.

\subsection{IVX Instruments}

We establish the validity of the IVX method when the covariates are
HT-SNPs, HT-MI or HT-SPs. Specifically, given a series $x_{t}$, the
IVX instruments $z_{t}(n)$ are obtained by applying the following
linear filter:
\begin{equation}
z_{t}:=z_{t}(n):=\sum^{t}_{j=1}\rho^{t-j}_{z,n}\left(x_{j}-x_{j-1}\right),\label{IVXz0}
\end{equation}
where $\rho_{z,n}$ is defined as in (\ref{rho_n}) i.e.
\begin{equation}
\rho_{z,n}:=1-\kappa^{-1}_{z,n},\;\kappa_{z,n}:=n^{\alpha_{z,\kappa}}\ell_{z,\kappa}(n),\label{IVXz}
\end{equation}
with $\kappa_{z,n}>0$ such that $\kappa_{z,n}\rightarrow\infty$
and $\sup_{n\geq1,1\leq t\leq n}\sup\kappa_{z,t}/\kappa_{z,n}<\infty$,
$\ell_{z,\kappa}$ SV and $\alpha_{z,\kappa}\in[0,1)$.

Recall that $\rho_{n}$ denotes the model-specific sequence associated
with the autoregressive parameter of $x_{t}$ when $x_{t}$ is a MI
process satisfying \textbf{Assumption MI}, whereas $\rho_{z,n}$ denotes
the analogous method-specific sequence associated with the IVX method.
In this context, the localizing sequence $\kappa_{z,n}$ functions
analogously to a bandwidth parameter and is selected by the practitioner.
The IVX estimator achieves conventional inference under minimal conditions
on $\kappa_{z,n}$, namely $\kappa^{-1}_{z,n}+\kappa_{z,n}n^{-1}\rightarrow0$.

If $x_{t}$ is an HT-SNP or HT-MI process, then $z_{t}$ behaves asymptotically
as an HT-MI process whose underlying autoregressive sequence is either
$\rho_{n}$ or $\rho_{z,n}$. Finally, if $x_{t}$ is a HT-SP, then
$z_{t}$ is asymptotically equivalent to $x_{t}$. These results parallel
those established by KMS2015.

To keep the length of the paper manageable, we impose some restrictions
on the parameter space of $x_{t}$. Specifically, we only consider
array (i.e. NI and MI) covariates driven by short memory errors and
stability parameter $1<\alpha\leq2$. The latter restriction is not
binding when $x_{t}$ is a pure $I(1)$ process, that is, when $\rho_{n}=1$.
The asymptotic behavior of the IVX instruments depends on the underlying
properties of $x_{t}$ and may be summarized as follows:
\begin{itemize}
\item[\textbf{\emph{HT-SNP:}}]  If $x_{t}\sim$\emph{HT-SNP} with $\rho_{n}=1$, then we have 
\[
d\in\left\{ \begin{array}{l}
(1-1/\alpha,1),\\
\left\{ 1\right\} ,\\
(1,2-1/\alpha)
\end{array}\right\} \stackrel{\textrm{IVX filtering}}{\rightarrow}z_{t}\sim\left\{ \begin{array}{l}
\textrm{MI-NM,}\:\textrm{on}\:\rho_{z,n}\\
\textrm{MI-SM},\:\textrm{on}\:\rho_{z,n}\\
\textrm{MI-LM}\:\textrm{on}\:\rho_{z,n}
\end{array}\right..
\]
\item[\textbf{\emph{HT-SNP/MI:}} ] \textbf{(arrays)} If $x_{t}\sim$\emph{HT-SNP-SM} or \emph{HT-MI}-SM,
with $\kappa_{n}/\kappa_{z,n}\rightarrow\infty$, then we have 
\[
x_{t}\sim\textrm{SNP/ MI-SM},\:\textrm{on}\:\rho_{n}\stackrel{\textrm{IVX filtering}}{\rightarrow}z_{t}\stackrel{approx.}{\sim}\textrm{MI-SM},\:\textrm{on}\:\rho_{z,n}.
\]
In addition, if $\kappa_{z,n}/\kappa_{n}\rightarrow\infty$, 
\[
x_{t}\sim\textrm{MI-SM},\:\textrm{on}\:\rho_{n}\stackrel{\textrm{IVX filtering}}{\rightarrow}z_{t}\stackrel{approx.}{\sim}\textrm{MI-SM},\:\textrm{on}\:\rho_{n}.
\]
\item[\textbf{\emph{HT-SP:}}]  If $x_{t}\sim$\emph{HT-SP},
\[
x_{t}\sim\textrm{HT-SP}\stackrel{\textrm{IVX filtering}}{\rightarrow}z_{t}\stackrel{approx.}{\sim}x_{t}.
\]
\end{itemize}
As indicated above, when $x_{t}\sim$\emph{HT-SNP} with $\rho_{n}=1$,
the IVX instruments $z_{t}$ are exactly MI processes. When $x_{t}$
is a nonstationary array, either NI or MI, the IVX instruments behave
asymptotically as MI processes. Finally, when $x_{t}$ is stationary,
the IVX instruments approximate $x_{t}$. The latter two approximation
results are formally stated in the following lemma. Before presenting
the lemma, we introduce an additional normalising sequence that supplements
\textbf{Definition 1}.

\bigskip{}

\noindent\textbf{Definition 1{*}}. \emph{Given a localizing sequence
$\kappa_{z,n}$, under }\textbf{\emph{Assumption HT-SNP}}\emph{ or
}\textbf{\emph{Assumption HT-MI}}\emph{, define $\gamma_{z,n}$ by}
\[
\gamma_{z,n}\sim\kappa^{b/\alpha}_{z,n}\ell_{1}\left(\kappa_{z,n}\right)^{1/\alpha},
\]
\textit{where $b$ and }$\ell_{1}$ \emph{are defined in }\textbf{\emph{Definition
1}}\emph{.}\bigskip{}

\textbf{Lemma \ref{Lemma_kappa_z}} below introduces random arrays
that asymptotically approximate the IVX instruments $z_{t}$ suitable
suitable normalization. The corresponding normalizing sequences, which
characterize the orders of magnitude of these instruments, are specified
in the following definition.

\bigskip{}

\noindent\textbf{Definition 5}. \emph{Suppose that $x_{t}$ satisfies
either }\textbf{\emph{Assumption SNP }}\emph{or }\textbf{\emph{Assumption
MI}}\emph{. Define the normalizing sequence $\mu_{n}$ by}
\[
\mu_{n}:=\left\{ \begin{array}{l}
\gamma_{z,n},\:\textrm{\emph{under \textbf{Assumption SNP} or \textbf{Assumption MI} with} \ensuremath{\kappa_{n}}/\ensuremath{\kappa_{z,n}\rightarrow\infty}},\\
\\\gamma_{n},\:\textrm{\emph{under \textbf{Assumption MI} with} \ensuremath{\kappa_{z,n}}/\ensuremath{\kappa_{n}\rightarrow\infty}}.
\end{array}\right.
\]

\bigskip{}

\noindent\begin{lemma}\textit{ \label{Lemma_kappa_z} Suppose that
$z_{t}$ is defined by (\ref{IVXz0}) and $\rho_{z,n}$ by (\ref{IVXz}).
The following hold:}

\bigskip{}

\begin{itemize}
\item[\textbf{(a)}] \textit{ Suppose that} $x_{t}$\textit{ satisfies }\textbf{\textit{Assumption
SNP-SM}}\textit{ or }\textbf{\textit{Assumption MI-SM}}\textit{ with
}$\kappa_{n}/\kappa_{z,n}\rightarrow\infty$\textit{. Then, for all
$1<\alpha^{\prime}<\alpha$ as }$n\rightarrow\infty$\textit{,}
\[
\frac{1}{\mu_{n}}\sup_{1\leq t\leq n}\mathbf{E}\left\vert z_{t}(n)-\sum^{t}_{j=1}\rho^{t-j}_{z,n}v_{j}\right\vert ^{\alpha^{\prime}}=o(1).
\]
\item[\textbf{(b)}] \textit{ Suppose that} $x_{t}$\textit{ satisfies }\textbf{\textit{Assumption
MI-SM}}\textit{ with }$\kappa_{z,n}/\kappa_{n}\rightarrow\infty$\textit{.
Then, for all $1<\alpha^{\prime}<\alpha$ as }$n\rightarrow\infty$\textit{,}
\[
\frac{1}{\mu_{n}}\sup_{1\leq t\leq n}\mathbf{E}\left\vert z_{t}(n)-x_{t}\right\vert ^{\alpha^{\prime}}=o(1).
\]
\item[\textbf{(c)}] \textit{ Suppose that} $x_{t}=v_{t}$\emph{, where}\textit{ }$v_{t}$\textit{
satisfies }\textbf{\textit{Assumption SP}}\textit{. Then, for all
$0<\alpha^{\prime}<\alpha$ as }$n\rightarrow\infty$\textit{,}
\[
\sup_{1\leq t\leq n}\mathbf{E}\left\vert z_{t}(n)-x_{t}\right\vert ^{\alpha^{\prime}}=o(1).
\]
\end{itemize}
\noindent\end{lemma} \bigskip{}

\noindent\begin{remark} It can be easily seen from \textbf{Lemma
\ref{WCon}} and \textbf{Lemma \ref{Lemma_kappa_z}} that for nonstationary
$x_{t}$, 
\[
z_{n}=\left\{ \begin{array}{l}
O_{p}(\gamma_{z,n}),\;\textrm{if \ensuremath{x_{t}} HT-SNP or HT-MI with}\;\kappa_{n}/\kappa_{z,n}\rightarrow\infty\\
O_{p}(\gamma_{n}),\;\textrm{if \ensuremath{x_{t}} HT-MI with}\;\kappa_{z,n}/\kappa_{n}\rightarrow\infty
\end{array}\right..
\]

\noindent\end{remark}

It follows from the preceding discussion that the IVX instrument is
obtained by applying a linear filter to the covariate $x_{t}$. As
discussed in the previous section, when $x_{t}$ is an SNP, parametric
estimators such as W-IV generally do not have Gaussian limit distributions.
By contrast, IVX filtering attenuates the signal in $x_{t}$ sufficiently
for an MCLT to apply, thereby validating conventional inference. This
signal attenuation is reflected in the relative orders of magnitude
of $x_{n}$ and $z_{n}$ when $x_{t}$ is an HT-SNP. In particular,
under our assumptions, 
\[
\gamma_{z,n},\gamma_{n}=o(\gamma^{\ast}_{n}).
\]
Thus, the normalizing sequence required for the weak convergence of
$x_{n}$ dominates those required for the weak convergence of $z_{n}$
and of the corresponding mildly integrated process. To see this, note
the following.
\begin{itemize}
\item $\gamma^{\ast}_{n}$ is $RV_{\mathcal{H}}$ in $n$ (i.e. regularly
varying of order $\mathcal{H}$) with $\mathcal{H}=d-(1-1/\alpha)$
being the Hurst index -see e.g. Jeganathan (2004, 2008).
\item $\gamma_{n}$ is $RV_{\alpha_{\kappa}\cdot\mathcal{H}}$ in $n$ with,
$\alpha_{\kappa}\in[0,1)$. \footnote{Recall that under \textbf{Assumption MI}, $\gamma_{n}$ is \textit{RV}$_{p}$
in $n$ with $p=\alpha_{\kappa}b/\alpha$, $b=1+\alpha-\alpha m$.
Moreover, the relationship between the memory parameter $\delta$
of $v_{t}$ and the exponent $m$ is $\delta=1-m$, while $d=1+\delta$.
Consequently, $p=\alpha_{\kappa}\left[d-(1-1/\alpha)\right]$.} 
\item $\gamma_{z,n}$ is $RV_{\alpha_{\kappa_{z}}\cdot\mathcal{H}}$ in
$n$ with, $\alpha_{\kappa_{z}}\in[0,1)$.
\end{itemize}
\bigskip{}

Before we present the IVX estimator, we introduce some notation. As
usual, for any sequence $\left\{ a_{t}\right\} ^{n}_{t=1}$ set 
\[
\bar{\alpha}=\frac{\sum^{n}_{t=1}a_{t}w_{0}(z_{t-1})}{\sum^{n}_{t=1}w_{0}(z_{t-1})}
\]
 and $\bar{\alpha}_{t}=a_{t}-\bar{\alpha}$. Note that for $w_{0}=1$,
$\bar{\alpha}$ reduces to the sample mean. Consider the weighted
IVX estimator (W-IVX, herafter) estimator for the slope parameter
of (\ref{Model0})
\begin{equation}
\hat{\beta}_{IVX}=\left[\sum^{n}_{t=1}w_{1}(z_{t-1})\overline{f}\left(x_{t-1}\right)\right]^{-1}\sum^{n}_{t=1}w_{1}(z_{t-1})\overline{y}_{t},\label{IVX1}
\end{equation}
where $z_{t}$ is given in (\ref{IVXz}, and $w$ is some weight function.
If we set $w(x)=x$ (and $f(x)$) we get the IVX estimator of Phillips
and Magdalinos (2009) as special case. As mentioned before, the possibility
of fractional nonstationary covariates has not been considered in
the IVX literature. In addition, all existing results concern covariates
with finite variance. Our limit theory is sufficient for the characterization
of IVX based test-statistics for the case where $x_{t}$ is an HT-SNP.
Consider the IVX t-statistic for the $H_{0}:\beta=\beta_{0}$
\begin{equation}
\hat{t}_{IVX}=\sum^{n}_{t=1}w_{1}(z_{t-1})\overline{f}\left(x_{t-1}\right)\cdot\frac{\hat{\beta}_{IVX}-\beta_{0}}{\sqrt{\hat{\sigma}^{2}\sum^{n}_{t=1}\overline{w}_{1}(z_{t-1})^{2}}},\label{IVX1b}
\end{equation}
where $\hat{\sigma}^{2}$ is some consistent estimator for $\sigma^{2}_{u}$
e.g. eq. (\ref{sigma_u}).

\bigskip{}

We start with the limiting behaviour of the IVX estimator when $x_{t}$
is an HT-SP. KMS2015 show that, for finite-variance short-memory linear
processes, the IVX estimator is asymptotically equivalent to OLS.
In our framework, the parametric counterpart of W-IVX is WIV. This
equivalence is important from a practical perspective, as it implies
that W-IVX attains parametric convergence rates under stationarity.
Consequently, W-IVX-based tests have very good asymptotic power under
stationarity, comparable to that of tests based on parametric estimators.
The following theorem establishes that WIV and W-IVX are asymptotically
equivalent when the covariate is an HT-SP. 

\bigskip{}

\begin{theorem} \label{IVX_SP}\emph{Suppose that $y_{t}$ is generated
by (\ref{parametric model IV}). In addition, suppose that the following
conditions hold:}
\begin{itemize}
\item[\textbf{(i)}]  \textbf{\emph{Assumption INN}}\emph{ holds with $\theta=1$;}
\item[\textbf{(ii)}]  $z_{t}$\emph{ }\textit{is defined by (\ref{IVXz0}) and $\rho_{z,n}$
by (\ref{IVXz}}\emph{);}
\item[\textbf{(iii)}]  \emph{$x_{t}$ satisfies }\textbf{\emph{Assumption SP}}\emph{;}
\item[\textbf{(iv)}] \textbf{\emph{ Assumption MDE}}\emph{ holds;}
\item[\textbf{(v)}] \emph{ }$\mathbf{E}\left\Vert \mathbf{A}(x_{t})\right\Vert $\emph{,
}\textbf{\emph{ }}\emph{$\mathbf{E}\left\Vert \mathbf{W}(x_{t})\right\Vert ^{1+\delta}<\infty$,
for some $\delta>0$;}
\item[\textbf{(vi)}] \emph{ the elements of $\mathbf{A}(x),\mathbf{W}(x)$ are $O\left(\left\vert x\right\vert ^{\lambda}\right)$,
$\lambda\in(0,\alpha)$;}
\item[\textbf{(vii)}]  $\hat{\sigma}^{2}$\emph{ is consistent for $\sigma^{2}_{u}$. }
\end{itemize}
\noindent\textit{The following hold:}
\begin{itemize}
\item[\textbf{(a)}]  \emph{As $n\rightarrow\infty$,
\[
\sqrt{n}\left(\mathbf{\tilde{\mathbf{b}}}_{WIV}-\mathbf{b}\right)=\sqrt{n}\left(\mathbf{\hat{b}}_{IVX}-\mathbf{b}\right)+o_{p}(1).
\]
}
\item[\textbf{(b)}]  \emph{Under $H_{0}$, as $n\rightarrow\infty$,}
\[
\hat{t}_{IVX}\overset{d}{\rightarrow}N(0,1).
\]
\end{itemize}
\end{theorem}

\bigskip{}
The conditions of \textbf{Theorem \ref{IVX_SP}} are comparable to
those of \textbf{Theorem \ref{WIV_SP}}. However, the latter imposes
an additional requirement, namely condition (vi), on the growth rates
of certain nonlinear transformations involving the regression and
weight functions. This technical condition facilitates the approximation
established in \textbf{Lemma \ref{Lemma_kappa_z}} without requiring
any smoothness restrictions on these functions. 

Next, we consider the case in which the covariate is either a nearly
integrated or an MI heavy-tailed array. Here, we restrict attention
to covariates driven by short-memory errors and assume that $1<\alpha\leq2$.
These restrictions are imposed to facilitate the approximation established
in \textbf{Lemma \ref{Lemma_kappa_z}}. Finally, we require the relevant
nonlinear functionals to satisfy\textbf{ Definition LH}. This latter
restriction can be dispensed with if $\mu^{-1}_{t}z_{t}$ possesses
a uniformly bounded density from some index onward. We conjecture
that this property holds under our assumptions, but leave its formal
verification, and hence the corresponding extension of the result,
for future work.

\bigskip{}

\begin{theorem} \label{IVX_SNP-2}\emph{Suppose that $y_{t}$ is
generated by (\ref{parametric model IV}). In addition, suppose that
the following hold:}
\begin{itemize}
\item[\textbf{(i)}]  \textbf{\emph{Assumption INN}}\emph{ holds with $1<\alpha\leq2$;}
\item[\textbf{(ii)}]  $z_{t}$\emph{ }\textit{is defined by (\ref{IVXz0}) and $\rho_{z,n}$
by (\ref{IVXz}}\emph{);}
\item[\textbf{(iii)}]  \emph{$x_{t}$ satisfies either} \textbf{\emph{Assumption SNP-SM}}\emph{,
or }\textbf{\emph{Assumption MI-SM}}\emph{ with $\kappa_{n}/\kappa_{z,n}\rightarrow\infty$
or }\textbf{\emph{Assumption MI-SM}}\emph{ with $\kappa_{z,n}/\kappa_{n}\rightarrow\infty$; }
\item[\textbf{(iv)}] \textbf{\emph{ Assumption MDE}}\emph{ holds;}
\item[\textbf{(v)}] \emph{ $\mathbf{W}(x)$ is }\textbf{\emph{AHF}}\emph{$\left\{ z_{t},\mu_{n}\right\} $
and }\textbf{\emph{LH}}\emph{;}
\item[\textbf{(vi)}]  $\int_{\mathbb{R}}\left[w_{1,H}(x)-a_{H}w_{0,H}\right]^{2}\varrho(x)dx>0$\emph{,
$a_{H}:=\frac{\int_{\mathbb{R}}w_{1,H}(x)\varrho(x)dx}{\int_{\mathbb{R}}w_{0,H}(x)\varrho(x)dx}$;}
\item[\textbf{(vii)}] $\hat{\sigma}^{2}$\emph{ is consistent for $\sigma^{2}_{u}$. }
\end{itemize}
\emph{Then under $H_{0}$, as $n\rightarrow\infty$,}
\[
\hat{t}_{IVX}\overset{d}{\rightarrow}N(0,1).
\]
\end{theorem}

\bigskip{}

Finally, the following limit result establishes that $\hat{t}_{IVX}$
is asymptotically standard normal when the covariate is HT-SNP process
with $\rho_{n}=1$.\bigskip{}

\begin{theorem} \label{IVX_SNP}\emph{Suppose that $y_{t}$ is generated
by (\ref{parametric model IV}). In addition, suppose that the following
conditions hold:}
\begin{itemize}
\item[\textbf{(i)}]  \textbf{\emph{Assumption INN}}\emph{ holds;}
\item[\textbf{(ii)}]  $z_{t}$\emph{ }\textit{is defined by (\ref{IVXz0}) and $\rho_{z,n}$
by (\ref{IVXz}}\emph{);}
\item[\textbf{(iii)}]  \emph{$x_{t}$ satisfies }\textbf{\emph{Assumption SNP}}\emph{ with
$\rho_{n}=1$;}
\item[\textbf{(iv)}] \textbf{\emph{ Assumption MDE}}\emph{ holds;}
\item[\textbf{(v)}] \emph{ $\mathbf{W}(x)$ is }\textbf{\emph{AHF}}\emph{$\left\{ z_{t},\mu_{n}\right\} $
and its elements are $O\left(\left\vert x\right\vert ^{\lambda}\right)$,
$\lambda\in(0,\alpha)$;}
\item[\textbf{(vi)}] \emph{ }$\int_{\mathbb{R}}\left[w_{1,H}(x)-a_{H}w_{0,H}\right]^{2}\varrho(x)dx>0$\emph{,
$a_{H}:=\frac{\int_{\mathbb{R}}w_{1,H}(x)\varrho(x)dx}{\int_{\mathbb{R}}w_{0,H}(x)\varrho(x)dx}$;}
\item[\textbf{(vii)}] $\hat{\sigma}^{2}$\emph{ is consistent for $\sigma^{2}_{u}$. }
\end{itemize}
\emph{Then under $H_{0}$, as $n\rightarrow\infty$,}
\[
\hat{t}_{IVX}\overset{d}{\rightarrow}N(0,1).
\]
\end{theorem}\bigskip{}

\section{Simulations}\label{sec:Simulations}

The tables and figures discussed in this section are collected in
Section \ref{sec:Simulation-Results} at the end of the appendix.

In this section we explore the finite sample performance of the weighted
IVX t-statistic -cf. (\ref{IVX1b}) for the hypothesis
\[
H_{0}:\beta=0
\]
\[
H_{1}:\beta=b\neq0
\]
in the context of the parametric regression of (\ref{Model0}) i.e.
\[
y_{t}=\mu+\beta f(x_{t-1})+u_{t},\:t=1,...,n.
\]
 Data are generated from 20,000 replication paths as follows. The
covariate $x_{t}$ is a type-II fractional process (e.g. Robinson
and Hualde, 2003; Shimotsu and Phillips, 2005) generated from
\[
\left(I-L\right)^{d}x_{t}=v_{t}I\left\{ t>0\right\} ,
\]
 where $v_{t}$ is some innovation process that is specified in detail
below. For of the memory parameter we consider a range of values in
the interval $d\in[d_{*}-0.6,d_{*}+0.6]$ with $d_{*}:=d_{*}(\alpha)$
being the nonstationarity boundary 
\[
d_{*}=1-1/\alpha.
\]
The regression errors are generated from
\begin{equation}
u_{t}=\alpha^{-1}\cdot\varepsilon^{\left\langle \alpha/\nu\right\rangle }_{t},\text{ }\nu=2.1,\label{sims2}
\end{equation}
where $\varepsilon_{t}$ is a secondary heavy tailed innovation specified
below. Recall that $x^{\left\langle a\right\rangle }=\textrm{sgn}(x)\left|x\right|^{a}$.
Finally, we specify the innovation vector $\left[v_{t},\varepsilon_{t}\right]$.
For $t=1,...,n$, $n=\{250,500,1000\}$, we consider three alternative
DGPs: 
\begin{itemize}
\item \textbf{DGP1}: $\left[v_{t},\varepsilon_{t}\right]$ follows a bi-variate
t-distribution with location parameter zero, $\alpha$-degrees of
freedom (tail parameter) and scale matrix 
\end{itemize}
\begin{equation}
\Sigma_{\delta}=\left[\begin{array}{cc}
1 & \delta\\
\delta & 1
\end{array}\right],\:\left|\delta\right|<1.\label{sims1}
\end{equation}

\begin{itemize}
\item \textbf{DGP2}: $\left[v_{t},\varepsilon_{t}\right]$ is sub-Gaussian
with stability parameter $\alpha$, and scale matrix $\Sigma_{\delta}$.
This is equivalent to $\left[v_{t},\varepsilon_{t}\right]$ having
characteristic function $\phi(\mathbf{t})=\exp\left(-\left|\frac{1}{2}\mathbf{t}'\Sigma_{\delta}\mathbf{t}\right|^{\alpha/2}\right)$.
For more details, see Samorodnitsky and Taqqu (1994; p. 78).
\end{itemize}
The following values for the stability parameter are considered 
\[
\alpha=[2,1.99,1.8,1.6,1.4,1.2,1,0.8].
\]
Under \textbf{DGP1}, when $\alpha=2$, $\left[v_{t},\varepsilon_{t}\right]$
is in the domain of attraction of normal law but has infinite variance\footnote{In this case, as $x\rightarrow\infty$ $\mathbf{E}\left(v^{2}_{t}I\left\{ \left|v_{t}\right|\leq x\right\} \right)\sim2\ln(x)$,
see e.g., Abadir and Magnus (2004) }. Under \textbf{DGP2}, when $\alpha=2$, the vector $\left[v_{t},\varepsilon_{t}\right]$
is Gaussian. In all cases, for $\alpha<2$ ($\alpha<1$), $v_{t},\varepsilon_{t}$
are infinite variance (mean) with tails $P\left(\left|v_{t}\right|>x\right),P\left(\left|\varepsilon_{t}\right|>x\right)\sim cons.x^{-\alpha}$.
The formulation of (\ref{sims2}) together with \textbf{DGP1-DGP2}
postulates dependence between $v_{t}$ (the innovations that drive
the covariate) and $u_{t}$ (regression error) are dependent with
each other. In particular, for $\alpha<2$, $v_{t}$ is always infinite-variance
whereas $u_{t}$ is heavy tailed but has finite variance with tail
exponent $\nu=2.1$.

We consider two functional forms for the regression function $f(x)$
along with multiple choices of the weight function $w(x)$ for each
regression function:
\begin{itemize}
\item \textbf{Specification 1}: $\left[f(x),\:w(x)\right]=\left[x,\:x^{\left\langle j\right\rangle }\right],j=\{0,0.5,0.75,1\}$;
\item \textbf{Specification 2:} $\left[f(x),\:w(x)\right]=\left[x^{\left\langle 0.5\right\rangle },\:x^{\left\langle 0.5\right\rangle }\right]$.
\end{itemize}
The first specification corresponds to a linear-in-variables model
with four alternative choices for the weight function. In particular,
three of these weight functions grow at a slower rate than $f$. To
assess the robustness of the proposed methods to more general specifications,
the simulation set-up includes models that are more general than those
permitted under our theoretical assumptions, specifically condition
(v) of Theorem \ref{IVX_SNP}. The second specification involves a
nonlinear regression function with an identical weight function; this
specification is directly relevant for the empirical application presented
in the following section. We note that for $\left[f(x),\:w(x)\right]=\left[x,\:x^{\left\langle 0.5\right\rangle }\right],\left[x^{\left\langle 0.5\right\rangle },\:x^{\left\langle 0.5\right\rangle }\right]$
our theoretical assumptions on the growth rates of the transformation
functions are satisfied whenever the stability parameter of the covariate
$x_{t}$ is $\alpha>1$. Following KMS2015, we set the IVX bandwidth
term to $\kappa^{-1}_{n}=n^{-0.95}$ and employ the studentization
method recommended therein. This approach uses HAC estimates based
on auto-regressive residuals. Whereas KMS2015 obtain regression residuals
using OLS autoregressions, here we instead use WIV autoregressions
(i.e., IV autoregressions that instrument with $w(x_{t-1})$). 

We report size and power results for three methods: 
\begin{itemize}
\item weighted IVX which uses the instrument $w(z_{t-1})$;
\item WIV which uses the instrument $w(x_{t-1})$;
\item OLS, interpreted as IV with instrument with instrument $f(x_{t-1}).$
\end{itemize}
It is worth noting that for Specification 1 with $w(x)=x$ and Specification
2, the IV and the OLS estimators coincide. 

The simulation results under \textbf{DGP1} are nearly identical to
those under \textbf{DGP2}. For this reason, we primarily report results
for \textbf{DGP1}. Tables \ref{tab:sim-dgp1-linear-w05-alpha-2p00}--\ref{tab:sim-dgp1-linear-w05-alpha-0p80}
present empirical size results for t-tests based on IVX, WIV, and
OLS for Specification 1 with weight function $w(x)=x^{\left\langle 0.5\right\rangle }$,
under \textbf{DGP1}. It is can easily seen from Tables \ref{tab:sim-dgp1-linear-w05-alpha-2p00}-\ref{tab:sim-dgp1-linear-w05-alpha-0p80}
that IVX tests exhibit very good size for all configurations even
if endogeneity is extreme (i.e. $\delta=\pm0.99$) and sample sizes
relatively small (e.g. $n=250$). On the other hand the IV and OLS
test exhibit severe size distortions when the covariate is in the
nonstationary region with the former performing slightly better in
general. The IVX tests in some cases are slightly oversized with empirical
size around 8\%. In contrast, the IV and OLS tests display severe
size distortions when the covariate lies in the nonstationary region,
with IV generally performing slightly better than OLS. In some cases,
the IVX tests are mildly oversized, with empirical size around 8\%.
By comparison, for the same configurations the IV and OLS tests exhibit
empirical sizes exceeding 35\%, i.e., more than six times the nominal
level. Interestingly, the IVX tests continue to exhibit very good
size even when Condition iv of Theorem \ref{IVX_SNP} is violated,
suggesting that the IVX t-statistics still converge to a standard
normal distribution, despite the fact that the basic limit theory
(cf. Theorem \ref{MainLL}), does not apply in this setting.

The simulation results under \textbf{DGP1} and \textbf{DGP2} are nearly
identical. We therefore focus primarily on\textbf{ DGP1}. Tables \ref{tab:sim-dgp1-linear-w05-alpha-2p00}--\ref{tab:sim-dgp1-linear-w05-alpha-0p80}
report the empirical sizes of t-tests based on IVX, WIV, and OLS for
Specification 1, with weight function $w(x)=x^{\left\langle 0.5\right\rangle }$.
The IVX tests exhibit very good size control across the configurations
considered, even under extreme endogeneity (i.e. $\delta=\pm0.99$)
and relatively small sample sizes $(\ensuremath{n=250})$. Although
mild over-rejection occurs in some cases, with empirical sizes reaching
approximately 8\%, the corresponding distortions are substantially
smaller than those associated with WIV and OLS. In particular, when
the covariate lies in the nonstationary region, the latter procedures
exhibit severe size distortions, with empirical rejection frequencies
exceeding 35\% in some configurations, compared with a nominal level
of 5\%. Overall, WIV tends to perform slightly better than OLS. Interestingly,
the IVX tests continue to exhibit very good size control even when
condition (v) of \textbf{Theorem \ref{IVX_SNP}} is violated. This
finding suggests that the validity of conventional IVX inference may
extend beyond the sufficient conditions established by our limit theory,
particularly since \textbf{Theorem \ref{MainLL}} does not apply in
these configurations.

Next, Tables \ref{tab:sim-dgp1-linear-weights-alpha-2p00}--\ref{tab:sim-dgp1-linear-weights-alpha-0p80}
examine the effects of alternative weight functions on the empirical
size of IVX tests under \textbf{DGP1}. Specifically, we consider the
following choices for $w(x)$: $w(x)=x^{\left\langle j\right\rangle },j=\{0,0.5,0.75,1\}$.
Overall, size performance is very good across specifications, although
the choice $w(x)=x^{\left\langle 0\right\rangle }$ (the Cauchy estimator
weight function) exhibits somewhat more pronounced over-sizing in
certain cases, with empirical size approaching 9\%. 

Tables \ref{tab:sim-dgp1-nonlinear-alpha-2p00}--\ref{tab:sim-dgp1-nonlinear-alpha-0p80}
report empirical size results for IVX and OLS tests under Specification
2 (i.e., a nonlinear regression function) for \textbf{DGP1}. We note
that in this setting the WIV and OLS estimators coincide. The overall
size performance is comparable to that reported in Tables \ref{tab:sim-dgp1-linear-w05-alpha-2p00}--\ref{tab:sim-dgp1-linear-w05-alpha-0p80}:
the IVX test continues to perform very well across all configurations,
whereas the OLS test exhibits severe over-sizing when the covariate
lies in the nonstationary region.

Finally, Tables \ref{tab:sim-dgp2-linear-weights-alpha-2p00}--\ref{tab:sim-dgp2-linear-weights-alpha-0p80}
report empirical size results for IVX tests under \textbf{DGP2} for
different choices of the weight function. The simulation design is
comparable to that used in Tables \ref{tab:sim-dgp1-linear-weights-alpha-2p00}--\ref{tab:sim-dgp1-linear-weights-alpha-0p80}.
Once again, size performance is very good across all weight-function
choices and across different memory and tail parameter configurations.

Finally, empirical power results for the three tests are presented
in Figures \ref{fig:sim-power-linear-w00}--\ref{fig:sim-sap-nonlinear-f05-w05}.
We focus on the case $n=1000$ and $\delta=-0.95$, which aligns the
Monte Carlo design with that of KMS2015. Regression errors are generated
according to DGP1. The following values for the stability and memory
parameters are considered:
\[
\alpha=[2,1.8,1.6,1.4],\:d\in[d_{*}-0.6,d_{*}+0.6].
\]

Figures \ref{fig:sim-power-linear-w00}--\ref{fig:sim-power-linear-w10}
display empirical power results for IVX-, WIV-, and OLS-based tests
under Specification 1. For nonstationary predictors, WIV and OLS yield
higher empirical power, albeit at the cost of severe over-sizing under
the null hypothesis. The asymptotic power advantage of these tests
over IVX diminishes progressively as the process approaches the stationary
region. For smaller values of the memory parameter, namely $d=d_{*}-0.4,d_{*}-0.6$,
the power performance of all three tests is nearly identical. This
finding is consistent with \textbf{Theorem \ref{IVX_SP}}, which establishes
the asymptotic equivalence of the IVX and WIV estimators under stationary
short-memory covariates.

Figures \ref{fig:sim-sap-linear-w00}--\ref{fig:sim-sap-nonlinear-f05-w05}
report size-adjusted power (SAP) results.\footnote{Under the alternative hypothesis $H_{1}$, SAP is computed as
\[
\mathrm{Empirical\:Power-}\max\mathrm{\{(Empirical\:Size-Nominal\:Size),0\}.}
\]
 } For highly persistent regressors (e.g., $d=d_{*}+0.4,d_{*}+0.6$),
IVX tests outperform WIV and OLS across a wide range of slope parameter
values. Under Specification 1, IVX exhibits particularly strong performance
for the weight functions $w(x)=x^{\left\langle j\right\rangle },j=\{0.75,1\}$.
Consistent with the results shown in Figures \ref{fig:sim-power-linear-w00}--\ref{fig:sim-power-linear-w10},
under stationarity the IVX tests exhibit size performance comparable
to that of the WIV and OLS tests, particularly for the weight functions
$w(x)=x^{\left\langle 0.75\right\rangle },x^{\left\langle 1\right\rangle }$.

\bigskip{}

\section{Application to the Predictability of Stock Returns}\label{sec:Application}

The predictability of stock returns using publicly available information
has been the subject of a vast literature in empirical finance and
econometrics. For reviews and recent developments, see, among others,
Campbell and Yogo (2006), Welch and Goyal (2008), Bollerslev et al.
(2013), Kasparis et al. (2015), KMS2015, Phillips (2015), and Andersen
and Varneskov (2021). The most common approach to testing the so-called
predictability hypothesis relies on t-tests in predictive regressions
of the form
\[
r_{t}=\beta_{0}+\beta_{1}x_{t-1}+u_{t},
\]
where $r_{t}$ denotes stock returns over a given horizon, $x_{t-1}$
is a lagged predictive variable, and $u_{t}$ is a martingale difference
error term. The null hypothesis $H_{0}:\beta_{1}=0$ corresponds to
the absence of predictive power of $x_{t}$ for future returns. 

Testing this hypothesis is challenging due to uncertainty surrounding
the stationarity---and more generally, the persistence properties---of
the predictive variables. The predictors most commonly used in empirical
work include financial variables such as dividend yields, return volatility,
and earnings-price ratios, as well as macroeconomic variables such
as inflation and short-term interest rates. These variables typically
exhibit high persistence, often consistent with nonstationary behavior.
It is well documented, for example, that measures of stock market
risk such as realized variance display long-memory dynamics (see Christensen
and Nielsen, 2006, 2007; Bollerslev et al., 2013; Andersen and Varneskov,
2021). By contrast, the possibility that predictive variables may
also exhibit heavy-tailed behavior has received comparatively little
attention. Recent work by Kim and Meddahi (2020) and Kim, Meddahi,
and Yamashita (2021) shows that return volatility -a widely used predictor-
exhibits thick tails.

In this study, we apply IVX tests to a set of predictive variables
from Welch and Goyal (2024) that appear to exhibit long memory and,
in some cases, heavy-tailed behavior. As demonstrated earlier, the
IVX tests considered here are agnostic to the exact values of the
memory and tail parameters. Nevertheless, to underscore the relevance
of this dataset for the methods employed, we report preliminary estimates
of memory parameters using the Local Whittle and Exact Local Whittle
(ELW) estimators (see Shimotsu and Phillips, 2005), along with Hill
estimates for the tail index.

Shimotsu and Phillips (2005) propose the ELW estimator, which remains
valid under nonstationarity provided that second moments are finite.
The validity of the Hill and Whittle estimators for stationary heavy-tailed
processes has been established by Hill (2010) and Kokoszka and Taqqu
(1996), respectively. More specifically, for stationary heavy-tailed
fractional processes (i.e., $d<1-1/\alpha$), Hill (2010) establishes
consistency and asymptotic normality of the Hill estimator, while
Kokoszka and Taqqu (1996) establish analogous results for the Whittle
estimator.

Given the uncertainty surrounding the persistence and tail properties
of the predictive variables, we conducted supplementary estimations
using both first-differenced and fractionally differenced series,
namely $\Delta x_{t}$ and $\Delta^{\hat{d}_{ELW}}x_{t}$, where $\hat{d}_{ELW}$
denotes the ELW estimate of the memory parameter. These transformations
were employed to assess the robustness of the estimated memory and
tail indices to the persistence properties of the underlying series.
After accounting for the effects of differencing on the memory parameter,
the supplementary estimates were broadly consistent with those obtained
from the original series in levels. Since these exercises yielded
qualitatively similar conclusions, we report only the estimates based
on the original series.

\bigskip{}

\textbf{Rebalancing Memory and Tails in Predictive Regressions}. A
well-known stylized fact in the return predictability literature (see,
e.g., Marmer, 2007; Kasparis et al., 2015; Phillips, 2015) is that
stock returns exhibit markedly different memory properties from those
of commonly used predictive variables such as return volatility, dividend
yield, earnings--price ratios, and inflation. This discrepancy creates
both modeling and statistical challenges. First, differences in memory
between the dependent and independent variables may signal model misspecification.
Second, it is well known that regressing a stationary variable on
a nonstationary one yields coefficient estimates that vanish asymptotically.
Phillips (2015) refers to this phenomenon as misbalanced regression.
As shown in Kasparis et al. (2015) and Hu, Kasparis, and Wang (2004),
such misbalancing can substantially reduce the power of predictability
tests. While memory misbalancing has received some attention in the
literature, the possibility of \textit{tail misbalancing} has, to
our knowledge, remained entirely unexplored.

To address memory misbalancing issues, two broad approaches have been
proposed. The first involves\textit{ linear filtering} of the data,
while the second relies on\textit{ nonlinear transformations}. The
linear filtering approach is advocated by Christensen and Nielsen
(2007), Bollerslev et al. (2013), and Andersen and Varneskov (2021),
who consider regressions involving fractionally differenced variables.
In particular, they study specifications of the form
\[
y_{t}=\hat{\mu}+\hat{\beta}\Delta^{\hat{d}_{x}}x_{t-1}+\hat{u}_{t}\text{ and }\Delta^{\hat{d}_{y}}y_{t}=\hat{\mu}+\hat{\beta}\Delta^{\hat{d}_{x}}x_{t-1}+\hat{u}_{t},
\]
with $\hat{d}_{y}$, $\hat{d}_{x}$ preliminary memory estimators
for $y_{t}$ and $x_{t}$ respectively. 

It is worth noting that while differencing reduces persistence, it
does not alter the tail behavior of the series. On the other hand,
nonlinear filtering may achieve both memory and tail rebalancing.
Marmer (2007), Kasparis (2011), Kasparis et al. (2015), Phillips (2015),
and Hu, Kasparis, and Wang (2024), among others, consider possibly
nonlinear predictive regressions of the form
\[
y_{t}=\mu+\beta f(x_{t-1})+u_{t},
\]
where $u_{t}$ is a martingale difference error term. As noted in
these studies, certain nonlinear transformations---such as reduced-growth,
bounded, or integrable functions---can weaken the signal of highly
persistent predictors, thereby improving the memory balance between
returns and predictors.

To illustrate this effect, suppose that $x_{t}$ is a finite-variance
random walk. Then the following orders of magnitude apply (see Park
and Phillips, 1999, 2001):
\[
\sum^{n}_{t=1}\left\vert f(x_{t})\right\vert =\left\{ \begin{array}{l}
O_{p}(n^{3/2})\text{, }f(x)=x\\
O_{p}(n^{1+b/2})\text{, }f(x)=\left\vert x\right\vert ^{b}\text{, }b>-1\\
O_{p}(n\ln(n))\text{, }f(x)=\ln\left\vert x\right\vert \\
O_{p}(n)\text{, }f(x)=1\left\{ x>0\right\} \\
O_{p}(n^{1/2})\text{, }f(x)\text{ integrable}
\end{array}\right.
\]
These results demonstrate that certain transformations of random walks
can produce sample functionals with weaker signal strength than those
associated with stationary processes, for which one typically has
$\sum^{n}_{t=1}\left\vert f(x_{t})\right\vert =O_{p}(n)$ (see e.g.
Giraitis, 1985). Simulation experiments in Kasparis et al. (2015)
and Phillips (2015) further show that, for specific choices of $f$,
the sample paths of a martingale difference sequence $u_{t}$ and
those of $f(x_{t})+u_{t}$ are almost indistinguishable, even if $x_{t}$
is a highly persistent process. 

In addition to reducing persistence, transformations with reduced
growth rates are readily seen to mitigate heavy-tailed behavior. For
instance, suppose that $X\sim\tilde{S}_{\alpha}$. Then the transformation
$X^{\left\langle 0.5\right\rangle }$ has finite variance whenever
the stability parameter satisfies $\alpha>1$. Moreover, the logarithmic
transformation --widely used in empirical finance-- $\ln\left(1+\left|X\right|\right)$
produces a random variable with exponentially decaying tails. \bigskip{}

\textbf{Empirical Evidence for the Predictability of Stock Returns}.
We primarily focus on monthly returns. For this frequency, we consider
four types of return series and nine predictive variables. Preliminary
results indicate that the most pronounced differences between monthly
and quarterly return predictability arise for the book-to-market (BM)
ratio and stock variance (svar) predictors. In particular, stock variance
exhibits no predictive power at the monthly frequency but displays
strong predictability for quarterly returns. The book-to-market ratio
appears to have predictive power at both frequencies; however, its
predictive strength is substantially greater at the quarterly horizon.
For this reason, we restrict attention to the BM and svar predictors
in the analysis of quarterly returns.

We begin our analysis by reporting preliminary estimates of the memory
and tail parameters. As noted above, the IVX methodology does not
rely on such preliminary estimates. Nevertheless, these estimates
help illustrate the relevance of our methods for the empirical application
considered here. Table \ref{tab:app-memory-tail} reports memory and
tail parameter estimates for the return series and the predictors.
Specifically, we employ the Local Whittle (LW) and Exact Local Whittle
(ELW) estimators using a bandwidth of $n^{0.65}$, where $n$ denotes
the sample size, as recommended by Shimotsu and Phillips (2005). For
the Hill estimator we consider four bandwidth choices: $k_{1}=0.05n$,
$k_{2}=0.1n$, $k_{3}=0.15n$ and $k_{4}=0.2n$ yielding the corresponding
estimates $\hat{\alpha}_{1}$, $\hat{\alpha}_{1}$, $\hat{\alpha}_{3}$,
and $\hat{\alpha}_{4}$. In addition, for each tail estimate we report
the implied nonstationarity boundary for the memory parameter, defined
as $\hat{d}_{i*}=1-1/$$\hat{\alpha}_{i}$, $i=\{1,2,3,4\}$.

Table \ref{tab:app-memory-tail} shows that stock returns exhibit
negligible memory and possess finite second moments. The tail index
estimates range from approximately 2.3 to 2.8 for monthly returns
and are slightly higher for quarterly returns. In contrast, all nine
monthly predictors considered display pronounced long-memory behavior,
with memory parameter estimates ranging from 0.518 to 1.085. Moreover,
for five of the nine predictors the Hill tail index estimates are
consistently below two, indicating infinite variance. In all cases,
the estimated memory parameters exceed the corresponding estimated
nonstationarity thresholds, implying that the predictors exhibit nonstationary
long-memory dynamics.

Tables \ref{tab:app-monthly-linear}--\ref{tab:app-quarterly-nonlinear}
report IVX F-statistic values for tests of the predictability hypothesis
across various predictors, regression specifications, and weight functions.
In particular, Tables \ref{tab:app-monthly-linear} and \ref{tab:app-monthly-nonlinear}
present results for monthly returns. Table \ref{tab:app-monthly-linear}
focuses on a linear regression specification and considers five choices
of the weight function, namely $w(x)=x^{\left\langle i\right\rangle }$,
$i=\{1,0.75,0.5,0.25,0.01\}$.\footnote{We also experimented with the Cauchy estimator weight $w(x)=x^{\left\langle 0\right\rangle }$;
however, this choice led to singularity issues for some predictor
series. To obtain comparable performance while avoiding singular matrices,
we instead employ the weight function $w(x)=x^{\left\langle 0.01\right\rangle }$.} Table \ref{tab:app-monthly-nonlinear} reports IVX F-statistics for
nonlinear specifications in which both the regression and weight functions
are given by $f(x)=w(x)=x^{\left\langle i\right\rangle }$, $i=\{1,0.75,0.5,0.25\}$.
Both tables indicate strong evidence of predictability for ntis and
inflation. In addition, there is some evidence of predictability for
tbl, and to a lesser extent for BM and EP. By contrast, there is no
evidence of predictability for widely used predictors such as DY and
svar (at least at the monthly frequency). Interestingly, the test
statistics for inflation are largest under nonlinear specifications,
suggesting the presence of nonlinear predictive effects. Finally,
Tables \ref{tab:app-quarterly-linear} and \ref{tab:app-quarterly-nonlinear}
report results for quarterly returns. As noted above, the analysis
at this frequency focuses on two predictors only: BM and svar. Both
variables exhibit strong predictive power, with BM significant at
the 1\% level and svar at the 5\% level. Notably, svar displays no
predictive power at the monthly frequency. Moreover, the test statistics
for the svar predictive regression are maximized under nonlinear specifications,
once again suggesting the presence of nonlinearities in the predictive
relationship.

\FloatBarrier

\begin{landscape} 

\begin{table}[H]
\centering \caption{Memory and Tail Estimates (1926--2024)}
\label{tab:app-memory-tail}

\HTmemorybegin

\HTpanel{Panel A: Monthly Returns; 1926:01 -- 2024:12} \begin{adjustbox}{max totalsize={\linewidth}{0.77\textwidth},center}
\begin{tabular}{lrrrrrr}
\textbf{Dependent Variables}  & \textbf{LW}  & \textbf{ELW}  & $\hat{\alpha}_{1}$  & $\hat{\alpha}_{2}$  & $\hat{\alpha}_{3}$  & $\hat{\alpha}_{4}$ \tabularnewline
Returns  & \HTtextcolor[rgb]{.047,.463,.62}{0.051}  & \HTtextcolor[rgb]{.047,.463,.62}{0.054}  & \HTtextcolor[rgb]{1,0,0}{2.697}  & \HTtextcolor[rgb]{1,0,0}{2.690}  & \HTtextcolor[rgb]{1,0,0}{2.614}  & \HTtextcolor[rgb]{1,0,0}{2.365} \tabularnewline
Returns-Rfree  & \HTtextcolor[rgb]{.047,.463,.62}{0.055}  & \HTtextcolor[rgb]{.047,.463,.62}{0.058}  & \HTtextcolor[rgb]{1,0,0}{2.774}  & \HTtextcolor[rgb]{1,0,0}{2.670}  & \HTtextcolor[rgb]{1,0,0}{2.602}  & \HTtextcolor[rgb]{1,0,0}{2.348} \tabularnewline
CRSP\_SPvw  & \HTtextcolor[rgb]{.047,.463,.62}{-0.012}  & \HTtextcolor[rgb]{.047,.463,.62}{-0.011}  & \HTtextcolor[rgb]{1,0,0}{2.658}  & \HTtextcolor[rgb]{1,0,0}{2.885}  & \HTtextcolor[rgb]{1,0,0}{2.772}  & \HTtextcolor[rgb]{1,0,0}{2.567} \tabularnewline
CRSP\_SPvw-Rfree  & \HTtextcolor[rgb]{.047,.463,.62}{0.011}  & \HTtextcolor[rgb]{.047,.463,.62}{0.010}  & \HTtextcolor[rgb]{1,0,0}{2.696}  & \HTtextcolor[rgb]{1,0,0}{2.851}  & \HTtextcolor[rgb]{1,0,0}{2.768}  & \HTtextcolor[rgb]{1,0,0}{2.548} \tabularnewline
\end{tabular}\end{adjustbox}

\HTpanel{Panel B: Monthly Predictors; 1926:01 -- 2024:12} \begin{adjustbox}{max totalsize={\linewidth}{0.77\textwidth},center}
\begin{tabular}{lrrrrrrrrrr}
\textbf{Predictors}  & \textbf{LW}  & \textbf{ELW}  & $\hat{\alpha}_{1}$  & $d_{1*}$  & $\hat{\alpha}_{2}$  & $d_{2*}$  & $\hat{\alpha}_{3}$  & $d_{3*}$  & $\hat{\alpha}_{4}$  & $d_{4*}$ \tabularnewline
BM  & \HTtextcolor[rgb]{.047,.463,.62}{\textbf{0.943}}  & \HTtextcolor[rgb]{.047,.463,.62}{\textbf{0.951}}  & \HTtextcolor[rgb]{1,0,0}{\textbf{1.707}}  & \HTtextcolor[rgb]{.235,.49,.133}{\textbf{0.414}}  & \HTtextcolor[rgb]{1,0,0}{\textbf{1.879}}  & \HTtextcolor[rgb]{.235,.49,.133}{\textbf{0.468}}  & \HTtextcolor[rgb]{1,0,0}{\textbf{1.916}}  & \HTtextcolor[rgb]{.235,.49,.133}{\textbf{0.478}}  & \HTtextcolor[rgb]{1,0,0}{\textbf{1.812}}  & \HTtextcolor[rgb]{.235,.49,.133}{\textbf{0.448}} \tabularnewline
tbl  & \HTtextcolor[rgb]{.047,.463,.62}{\textbf{1.080}}  & \HTtextcolor[rgb]{.047,.463,.62}{\textbf{1.085}}  & \HTtextcolor[rgb]{1,0,0}{\textbf{1.899}}  & \HTtextcolor[rgb]{.235,.49,.133}{\textbf{0.473}}  & \HTtextcolor[rgb]{1,0,0}{\textbf{1.762}}  & \HTtextcolor[rgb]{.235,.49,.133}{\textbf{0.432}}  & \HTtextcolor[rgb]{1,0,0}{\textbf{1.527}}  & \HTtextcolor[rgb]{.235,.49,.133}{\textbf{0.345}}  & \HTtextcolor[rgb]{1,0,0}{\textbf{1.439}}  & \HTtextcolor[rgb]{.235,.49,.133}{\textbf{0.305}} \tabularnewline
lty  & \HTtextcolor[rgb]{.047,.463,.62}{\textbf{1.042}}  & \HTtextcolor[rgb]{.047,.463,.62}{\textbf{1.048}}  & \HTtextcolor[rgb]{1,0,0}{2.628}  & \HTtextcolor[rgb]{.235,.49,.133}{0.500}  & \HTtextcolor[rgb]{1,0,0}{2.391}  & \HTtextcolor[rgb]{.235,.49,.133}{0.500}  & \HTtextcolor[rgb]{1,0,0}{2.228}  & \HTtextcolor[rgb]{.235,.49,.133}{0.500}  & \HTtextcolor[rgb]{1,0,0}{\textbf{1.972}}  & \HTtextcolor[rgb]{.235,.49,.133}{\textbf{0.493}} \tabularnewline
ntis  & \HTtextcolor[rgb]{.047,.463,.62}{\textbf{0.937}}  & \HTtextcolor[rgb]{.047,.463,.62}{\textbf{0.859}}  & \HTtextcolor[rgb]{1,0,0}{2.011}  & \HTtextcolor[rgb]{.235,.49,.133}{0.500}  & \HTtextcolor[rgb]{1,0,0}{\textbf{1.923}}  & \HTtextcolor[rgb]{.235,.49,.133}{\textbf{0.480}}  & \HTtextcolor[rgb]{1,0,0}{\textbf{1.901}}  & \HTtextcolor[rgb]{.235,.49,.133}{\textbf{0.474}}  & \HTtextcolor[rgb]{1,0,0}{\textbf{1.897}}  & \HTtextcolor[rgb]{.235,.49,.133}{\textbf{0.473}} \tabularnewline
Inflation  & \HTtextcolor[rgb]{.047,.463,.62}{\textbf{0.516}}  & \HTtextcolor[rgb]{.047,.463,.62}{\textbf{0.538}}  & \HTtextcolor[rgb]{1,0,0}{2.666}  & \HTtextcolor[rgb]{.235,.49,.133}{0.500}  & \HTtextcolor[rgb]{1,0,0}{2.888}  & \HTtextcolor[rgb]{.235,.49,.133}{0.500}  & \HTtextcolor[rgb]{1,0,0}{2.761}  & \HTtextcolor[rgb]{.235,.49,.133}{0.500}  & \HTtextcolor[rgb]{1,0,0}{2.503}  & \HTtextcolor[rgb]{.235,.49,.133}{0.500} \tabularnewline
svar  & \HTtextcolor[rgb]{.047,.463,.62}{\textbf{0.518}}  & \HTtextcolor[rgb]{.047,.463,.62}{\textbf{0.518}}  & \HTtextcolor[rgb]{1,0,0}{\textbf{1.469}}  & \HTtextcolor[rgb]{.235,.49,.133}{\textbf{0.319}}  & \HTtextcolor[rgb]{1,0,0}{\textbf{1.299}}  & \HTtextcolor[rgb]{.235,.49,.133}{\textbf{0.230}}  & \HTtextcolor[rgb]{1,0,0}{\textbf{1.169}}  & \HTtextcolor[rgb]{.235,.49,.133}{\textbf{0.144}}  & \HTtextcolor[rgb]{1,0,0}{\textbf{1.144}}  & \HTtextcolor[rgb]{.235,.49,.133}{\textbf{0.126}} \tabularnewline
DY  & \HTtextcolor[rgb]{.047,.463,.62}{\textbf{0.886}}  & \HTtextcolor[rgb]{.047,.463,.62}{\textbf{0.906}}  & \HTtextcolor[rgb]{1,0,0}{2.327}  & \HTtextcolor[rgb]{.235,.49,.133}{0.500}  & \HTtextcolor[rgb]{1,0,0}{2.373}  & \HTtextcolor[rgb]{.235,.49,.133}{0.500}  & \HTtextcolor[rgb]{1,0,0}{2.414}  & \HTtextcolor[rgb]{.235,.49,.133}{0.500}  & \HTtextcolor[rgb]{1,0,0}{2.263}  & \HTtextcolor[rgb]{.235,.49,.133}{0.500} \tabularnewline
EP  & \HTtextcolor[rgb]{.047,.463,.62}{\textbf{0.698}}  & \HTtextcolor[rgb]{.047,.463,.62}{\textbf{0.699}}  & \HTtextcolor[rgb]{1,0,0}{2.308}  & \HTtextcolor[rgb]{.235,.49,.133}{0.500}  & \HTtextcolor[rgb]{1,0,0}{2.103}  & \HTtextcolor[rgb]{.235,.49,.133}{0.500}  & \HTtextcolor[rgb]{1,0,0}{2.129}  & \HTtextcolor[rgb]{.235,.49,.133}{0.500}  & \HTtextcolor[rgb]{1,0,0}{2.065}  & \HTtextcolor[rgb]{.235,.49,.133}{0.500} \tabularnewline
Dpay  & \HTtextcolor[rgb]{.047,.463,.62}{\textbf{0.657}}  & \HTtextcolor[rgb]{.047,.463,.62}{\textbf{0.630}}  & \HTtextcolor[rgb]{1,0,0}{\textbf{1.754}}  & \HTtextcolor[rgb]{.235,.49,.133}{\textbf{0.430}}  & \HTtextcolor[rgb]{1,0,0}{\textbf{1.943}}  & \HTtextcolor[rgb]{.235,.49,.133}{\textbf{0.485}}  & \HTtextcolor[rgb]{1,0,0}{\textbf{1.968}}  & \HTtextcolor[rgb]{.235,.49,.133}{\textbf{0.492}}  & \HTtextcolor[rgb]{1,0,0}{\textbf{1.726}}  & \HTtextcolor[rgb]{.235,.49,.133}{\textbf{0.421}} \tabularnewline
\end{tabular}\end{adjustbox}

\HTpanel{Panel C: Quarterly Returns; 1926:1 -- 2024:4} \begin{adjustbox}{max totalsize={\linewidth}{0.77\textwidth},center}
\begin{tabular}{lrrrrrr}
\textbf{Dependent Variables}  & \textbf{LW}  & \textbf{ELW}  & $\hat{\alpha}_{1}$  & $\hat{\alpha}_{2}$  & $\hat{\alpha}_{3}$  & $\hat{\alpha}_{4}$ \tabularnewline
Returns  & \HTtextcolor[rgb]{.082,.376,.51}{-0.154}  & \HTtextcolor[rgb]{.082,.376,.51}{-0.135}  & \HTtextcolor[rgb]{1,0,0}{2.639}  & \HTtextcolor[rgb]{1,0,0}{2.696}  & \HTtextcolor[rgb]{1,0,0}{2.464}  & \HTtextcolor[rgb]{1,0,0}{2.291} \tabularnewline
Returns-Rfree  & \HTtextcolor[rgb]{.082,.376,.51}{-0.140}  & \HTtextcolor[rgb]{.082,.376,.51}{-0.132}  & \HTtextcolor[rgb]{1,0,0}{2.378}  & \HTtextcolor[rgb]{1,0,0}{2.551}  & \HTtextcolor[rgb]{1,0,0}{2.283}  & \HTtextcolor[rgb]{1,0,0}{2.300} \tabularnewline
CRSP\_SPvw  & \HTtextcolor[rgb]{.082,.376,.51}{-0.309}  & \HTtextcolor[rgb]{.082,.376,.51}{-0.312}  & \HTtextcolor[rgb]{1,0,0}{2.743}  & \HTtextcolor[rgb]{1,0,0}{2.940}  & \HTtextcolor[rgb]{1,0,0}{2.905}  & \HTtextcolor[rgb]{1,0,0}{2.625} \tabularnewline
CRSP\_SPvw-Rfree  & \HTtextcolor[rgb]{.082,.376,.51}{-0.199}  & \HTtextcolor[rgb]{.082,.376,.51}{-0.204}  & \HTtextcolor[rgb]{1,0,0}{3.020}  & \HTtextcolor[rgb]{1,0,0}{2.937}  & \HTtextcolor[rgb]{1,0,0}{2.680}  & \HTtextcolor[rgb]{1,0,0}{2.554} \tabularnewline
\end{tabular}\end{adjustbox}

\HTpanel{Panel D: Quarterly Predictors; 1926:1 -- 2024:4} \begin{adjustbox}{max totalsize={\linewidth}{0.77\textwidth},center}
\begin{tabular}{lrrrrrrrrrr}
\textbf{Predictors}  & \textbf{LW}  & \textbf{ELW}  & $\hat{\alpha}_{1}$  & $d_{1*}$  & $\hat{\alpha}_{2}$  & $d_{2*}$  & $\hat{\alpha}_{3}$  & $d_{3*}$  & $\hat{\alpha}_{4}$  & $d_{4*}$ \tabularnewline
BM  & \HTtextcolor[rgb]{.047,.463,.62}{\textbf{0.739}}  & \HTtextcolor[rgb]{.047,.463,.62}{\textbf{0.741}}  & \HTtextcolor[rgb]{1,0,0}{\textbf{1.885}}  & \HTtextcolor[rgb]{.235,.49,.133}{\textbf{0.470}}  & \HTtextcolor[rgb]{1,0,0}{\textbf{1.681}}  & \HTtextcolor[rgb]{.235,.49,.133}{\textbf{0.405}}  & \HTtextcolor[rgb]{1,0,0}{\textbf{1.787}}  & \HTtextcolor[rgb]{.235,.49,.133}{\textbf{0.440}}  & \HTtextcolor[rgb]{1,0,0}{\textbf{1.808}}  & \HTtextcolor[rgb]{.235,.49,.133}{\textbf{0.447}} \tabularnewline
svar  & \HTtextcolor[rgb]{.047,.463,.62}{\textbf{0.388}}  & \HTtextcolor[rgb]{.047,.463,.62}{\textbf{0.394}}  & \HTtextcolor[rgb]{1,0,0}{\textbf{1.369}}  & \HTtextcolor[rgb]{.235,.49,.133}{\textbf{0.269}}  & \HTtextcolor[rgb]{1,0,0}{\textbf{1.323}}  & \HTtextcolor[rgb]{.235,.49,.133}{\textbf{0.244}}  & \HTtextcolor[rgb]{1,0,0}{\textbf{1.267}}  & \HTtextcolor[rgb]{.235,.49,.133}{\textbf{0.211}}  & \HTtextcolor[rgb]{1,0,0}{\textbf{1.204}}  & \HTtextcolor[rgb]{.235,.49,.133}{\textbf{0.169}} \tabularnewline
\end{tabular}\end{adjustbox}

\HTmemoryend 
\end{table}

\end{landscape} \clearpage{}

\begin{table}[H]
\centering \caption{IVX predictability tests; monthly data; linear fits: $f(x)=x$}
\begin{adjustbox}{max totalsize={\linewidth}{0.80\textheight},center}
\begin{tabular}{llllll}
\toprule 
$w(x)$  & $x^{\langle1\rangle}$  & $x^{\langle0.75\rangle}$  & $x^{\langle0.5\rangle}$  & $x^{\langle0.25\rangle}$  & $x^{\langle0.01\rangle}$ \tabularnewline
\midrule 
\textbf{BM}  &  &  &  &  & \tabularnewline
Returns  & 0.283  & 0.053  & 0.000  & 0.005  & 0.004 \tabularnewline
Returns-Rfree  & 0.316  & 0.045  & 0.003  & 0.053  & 0.096 \tabularnewline
CRSP\_SPvw  & \textbf{3.581{*}}  & 2.071  & 1.204  & 0.779  & 0.621 \tabularnewline
CRSP\_SPvw-Rfree  & \textbf{3.686{*}}  & 2.014  & 1.032  & 0.522  & 0.287 \tabularnewline
\midrule 
\textbf{tbl}  &  &  &  &  & \tabularnewline
Returns  & 0.266  & 0.089  & 0  & 0.097  & 0.344 \tabularnewline
Returns-Rfree  & \textbf{4.383{*}{*}}  & \textbf{3.305{*}}  & 2.014  & 0.818  & 0.127 \tabularnewline
CRSP\_SPvw  & 0.171  & 0.04  & 0.006  & 0.163  & 0.445 \tabularnewline
CRSP\_SPvw-Rfree  & \textbf{4.071{*}{*}}  & \textbf{3.047{*}}  & 1.814  & 0.692  & 0.085 \tabularnewline
\midrule 
\textbf{lty}  &  &  &  &  & \tabularnewline
Returns  & 0.321  & 0.471  & 0.631  & 0.733  & 0.707 \tabularnewline
Returns-Rfree  & \textbf{3.173{*}}  & \textbf{3.277{*}}  & \textbf{3.1{*}}  & 2.569  & 1.807 \tabularnewline
CRSP\_SPvw  & 0.056  & 0.093  & 0.137  & 0.174  & 0.187 \tabularnewline
CRSP\_SPvw-Rfree  & 2.087  & 2.025  & 1.773  & 1.345  & 0.872 \tabularnewline
\midrule 
\textbf{Inflation}  &  &  &  &  & \tabularnewline
Returns  & 2.242  & 2.256  & 2.232  & 2.149  & 1.989 \tabularnewline
Returns-Rfree  & \textbf{3.346{*}}  & \textbf{3.609{*}{*}}  & \textbf{3.735{*}}  & \textbf{3.646{*}}  & \textbf{3.265{*}} \tabularnewline
CRSP\_SPvw  & \textbf{2.902{*}}  & \textbf{2.996{*}}  & \textbf{2.922{*}}  & 2.613  & 2.062 \tabularnewline
CRSP\_SPvw-Rfree  & \textbf{4.145{*}{*}}  & \textbf{4.534{*}{*}}  & \textbf{4.618{*}{*}}  & \textbf{4.248{*}{*}}  & \textbf{3.364{*}} \tabularnewline
\midrule 
\textbf{ntis}  &  &  &  &  & \tabularnewline
Returns  & \textbf{6.412{*}{*}}  & \textbf{5.434{*}{*}}  & \textbf{3.704{*}}  & 1.729  & 0.256 \tabularnewline
Returns-Rfree  & \textbf{5.975{*}{*}}  & \textbf{4.920{*}{*}}  & \textbf{3.195{*}}  & 1.34  & 0.115 \tabularnewline
CRSP\_SPvw  & \textbf{6.520{*}{*}}  & \textbf{5.653{*}{*}}  & \textbf{3.971{*}{*}}  & 1.959  & 0.367 \tabularnewline
CRSP\_SPvw-Rfree  & \textbf{6.068{*}{*}}  & \textbf{5.116{*}{*}}  & \textbf{3.431}  & 1.535  & 0.189 \tabularnewline
\midrule 
\textbf{svar}  &  &  &  &  & \tabularnewline
Returns  & 1.443  & 1.478  & 1.03  & 0.26  & 0.015 \tabularnewline
Returns-Rfree  & 0.827  & 0.836  & 0.575  & 0.117  & 0.047 \tabularnewline
CRSP\_SPvw  & 0.035  & 0.061  & 0.12  & 0.3  & 0.704 \tabularnewline
CRSP\_SPvw-Rfree  & 0.23  & 0.299  & 0.362  & 0.509  & 0.866 \tabularnewline
\midrule 
\textbf{DY}  &  &  &  &  & \tabularnewline
Returns  & 0.541  & 0.318  & 0.147  & 0.037  & 0 \tabularnewline
Returns-Rfree  & 0.332  & 0.195  & 0.089  & 0.023  & 0 \tabularnewline
CRSP\_SPvw  & 0.406  & 0.24  & 0.111  & 0.028  & 0 \tabularnewline
CRSP\_SPvw-Rfree  & 0.231  & 0.137  & 0.063  & 0.016  & 0 \tabularnewline
\midrule 
\textbf{EP}  &  &  &  &  & \tabularnewline
Returns  & \textbf{3.415{*}}  & 2.589  & 1.541  & 0.486  & 0 \tabularnewline
Returns-Rfree  & 1.802  & 1.339  & 0.781  & 0.241  & 0 \tabularnewline
CRSP\_SPvw  & 1.979  & 1.496  & 0.885  & 0.278  & 0 \tabularnewline
CRSP\_SPvw-Rfree  & 0.832  & 0.615  & 0.356  & 0.11  & 0 \tabularnewline
\midrule 
\textbf{Dpay}  &  &  &  &  & \tabularnewline
Returns  & 0.9  & 1.592  & 2.122  & 2.123  & 1.353 \tabularnewline
Returns-Rfree  & 0.264  & 0.679  & 1.063  & 1.13  & 0.675 \tabularnewline
CRSP\_SPvw  & 0.09  & 0.263  & 0.369  & 0.266  & 0.027 \tabularnewline
CRSP\_SPvw-Rfree  & 0.02  & 0.003  & 0.027  & 0.01  & 0.036 \tabularnewline
\bottomrule
\end{tabular}\end{adjustbox}\label{tab:app-monthly-linear}
\end{table}

\clearpage{}

\begin{table}[H]
\centering \caption{IVX predictability tests; monthly data; nonlinear fits: $f(x)=w(x)$}
\begin{adjustbox}{max totalsize={\linewidth}{0.80\textheight},center}
\begin{tabular}{lllll}
\toprule 
$f(x)$  & $x^{\langle1\rangle}$  & $x^{\langle0.75\rangle}$  & $x^{\langle0.5\rangle}$  & $x^{\langle0.25\rangle}$ \tabularnewline
\midrule 
\textbf{BM}  &  &  &  & \tabularnewline
Returns  & 0.283  & 0.052  & 0  & 0.005 \tabularnewline
Returns-Rfree  & 0.316  & 0.045  & 0.003  & 0.053 \tabularnewline
CRSP\_SPvw  & \textbf{3.581{*}}  & 2.064  & 1.198  & 0.775 \tabularnewline
CRSP\_SPvw-Rfree  & \textbf{3.686{*}}  & 2.007  & 1.027  & 0.52 \tabularnewline
\midrule 
\textbf{tbl}  &  &  &  & \tabularnewline
Returns  & 0.266  & 0.089  & 0  & 0.097 \tabularnewline
Returns-Rfree  & \textbf{4.383{*}{*}}  & \textbf{3.304{*}}  & 2.014  & 0.818 \tabularnewline
CRSP\_SPvw  & 0.171  & 0.04  & 0.006  & 0.164 \tabularnewline
CRSP\_SPvw-Rfree  & \textbf{4.071{*}{*}}  & \textbf{3.047{*}}  & 1.814  & 0.692 \tabularnewline
\midrule 
\textbf{lty}  &  &  &  & \tabularnewline
Returns  & 0.321  & 0.471  & 0.631  & 0.733 \tabularnewline
Returns-Rfree  & \textbf{3.173{*}}  & \textbf{3.277{*}}  & \textbf{3.1{*}}  & 2.57 \tabularnewline
CRSP\_SPvw  & 0.056  & 0.093  & 0.137  & 0.174 \tabularnewline
CRSP\_SPvw-Rfree  & 2.087  & 2.026  & 1.773  & 1.345 \tabularnewline
\midrule 
\textbf{Inflation}  &  &  &  & \tabularnewline
Returns  & 2.242  & 2.256  & 2.231  & 2.148 \tabularnewline
Returns-Rfree  & \textbf{3.346{*}}  & \textbf{3.609{*}}  & \textbf{3.734{*}}  & \textbf{3.644{*}} \tabularnewline
CRSP\_SPvw  & \textbf{2.902{*}}  & \textbf{2.996{*}}  & \textbf{2.921{*}}  & 2.612 \tabularnewline
CRSP\_SPvw-Rfree  & \textbf{4.145{*}{*}}  & \textbf{4.535{*}{*}}  & \textbf{4.618{*}{*}}  & \textbf{4.247{*}{*}} \tabularnewline
\midrule 
\textbf{ntis}  &  &  &  & \tabularnewline
Returns  & \textbf{6.412{*}{*}}  & \textbf{5.431{*}{*}}  & \textbf{3.701{*}}  & 1.727 \tabularnewline
Returns-Rfree  & \textbf{5.975{*}{*}}  & \textbf{4.918{*}{*}}  & \textbf{3.192{*}}  & 1.339 \tabularnewline
CRSP\_SPvw  & \textbf{6.520{*}{*}}  & \textbf{5.65{*}{*}}  & \textbf{3.967{*}{*}}  & 1.957 \tabularnewline
CRSP\_SPvw-Rfree  & \textbf{6.068{*}{*}}  & \textbf{5.114{*}{*}}  & \textbf{3.428{*}}  & 1.533 \tabularnewline
\midrule 
\textbf{svar}  &  &  &  & \tabularnewline
Returns  & 1.443  & 1.478  & 1.03  & 0.26 \tabularnewline
Returns-Rfree  & 0.827  & 0.836  & 0.575  & 0.117 \tabularnewline
CRSP\_SPvw  & 0.035  & 0.061  & 0.12  & 0.3 \tabularnewline
CRSP\_SPvw-Rfree  & 0.23  & 0.299  & 0.362  & 0.509 \tabularnewline
\midrule 
\textbf{DY}  &  &  &  & \tabularnewline
Returns  & 0.541  & 0.318  & 0.147  & 0.037 \tabularnewline
Returns-Rfree  & 0.332  & 0.195  & 0.089  & 0.023 \tabularnewline
CRSP\_SPvw  & 0.406  & 0.24  & 0.111  & 0.028 \tabularnewline
CRSP\_SPvw-Rfree  & 0.231  & 0.137  & 0.063  & 0.016 \tabularnewline
\midrule 
\textbf{EP}  &  &  &  & \tabularnewline
Returns  & \textbf{3.415{*}}  & 2.581  & 1.527  & 0.478 \tabularnewline
Returns-Rfree  & 1.802  & 1.336  & 0.776  & 0.238 \tabularnewline
CRSP\_SPvw  & 1.979  & 1.488  & 0.874  & 0.272 \tabularnewline
CRSP\_SPvw-Rfree  & 0.832  & 0.612  & 0.352  & 0.108 \tabularnewline
\midrule 
\textbf{Dpay}  &  &  &  & \tabularnewline
Returns  & 0.9  & 1.592  & 2.131  & 2.172 \tabularnewline
Returns-Rfree  & 0.264  & 0.679  & 1.067  & 1.148 \tabularnewline
CRSP\_SPvw  & 0.09  & 0.263  & 0.369  & 0.267 \tabularnewline
CRSP\_SPvw-Rfree  & 0.02  & 0.003  & 0.027  & 0.01 \tabularnewline
\bottomrule
\end{tabular}\end{adjustbox}\label{tab:app-monthly-nonlinear}
\end{table}

\clearpage{}

\begin{table}[H]
\centering \caption{IVX predictability tests; quarterly data; linear fits: $f(x)=x$}
\begin{adjustbox}{max totalsize={\linewidth}{0.80\textheight},center}
\begin{tabular}{llllll}
\toprule 
$w(x)$  & $x^{\langle1\rangle}$  & $x^{\langle0.75\rangle}$  & $x^{\langle0.5\rangle}$  & $x^{\langle0.25\rangle}$  & $x^{\langle0.01\rangle}$ \tabularnewline
\midrule 
\textbf{BM}  &  &  &  &  & \tabularnewline
Returns  & 0.963  & 0.305  & 0.065  & 0.01  & 0.011 \tabularnewline
Returns-Rfree  & 0.97  & 0.263  & 0.025  & 0.004  & 0.02 \tabularnewline
CRSP\_SPvw  & \textbf{7.258{*}{*}{*}}  & \textbf{3.95{*}{*}}  & 2.148  & 1.293  & 1.005 \tabularnewline
CRSP\_SPvw-Rfree  & \textbf{7.261{*}{*}{*}}  & \textbf{3.789{*}}  & 1.877  & 0.939  & 0.558 \tabularnewline
\midrule 
\textbf{svar}  &  &  &  &  & \tabularnewline
Returns  & 0.003  & 0  & 0.002  & 0.012  & 0.003 \tabularnewline
Returns-Rfree  & 0.077  & 0.092  & 0.098  & 0.066  & 0.011 \tabularnewline
CRSP\_SPvw  & \textbf{3.892{*}{*}}  & \textbf{3.943{*}{*}}  & \textbf{3.34{*}}  & 2.072  & 0.881 \tabularnewline
CRSP\_SPvw-Rfree  & \textbf{5.32{*}{*}}  & \textbf{5.333{*}{*}}  & \textbf{4.348{*}{*}}  & 2.495  & 0.959 \tabularnewline
\bottomrule
\end{tabular}\end{adjustbox}\label{tab:app-quarterly-linear}
\end{table}

\clearpage{}

\begin{table}[H]
\centering \caption{IVX predictability tests; quarterly data; nonlinear fits: $f(x)=w(x)$}
\begin{adjustbox}{max totalsize={\linewidth}{0.80\textheight},center}
\begin{tabular}{lllll}
\toprule 
$f(x)$  & $x^{\langle1\rangle}$  & $x^{\langle0.75\rangle}$  & $x^{\langle0.5\rangle}$  & $x^{\langle0.25\rangle}$ \tabularnewline
\midrule 
\textbf{BM}  &  &  &  & \tabularnewline
Returns  & 0.963  & 0.304  & 0.064  & 0.01 \tabularnewline
Returns-Rfree  & 0.97  & 0.262  & 0.025  & 0.004 \tabularnewline
CRSP\_SPvw  & \textbf{7.258{*}{*}{*}}  & \textbf{3.907{*}{*}}  & 2.111  & 1.27 \tabularnewline
CRSP\_SPvw-Rfree  & \textbf{7.261{*}{*}{*}}  & \textbf{3.751{*}{*}}  & 1.847  & 0.923 \tabularnewline
\midrule 
\textbf{svar}  &  &  &  & \tabularnewline
Returns  & 0.003  & 0  & 0.002  & 0.012 \tabularnewline
Returns-Rfree  & 0.077  & 0.092  & 0.098  & 0.066 \tabularnewline
CRSP\_SPvw  & \textbf{3.892{*}{*}}  & \textbf{3.941{*}{*}}  & \textbf{3.337{*}{*}}  & 2.071 \tabularnewline
CRSP\_SPvw-Rfree  & \textbf{5.32{*}{*}}  & \textbf{5.329{*}{*}}  & \textbf{4.342{*}{*}}  & 2.493 \tabularnewline
\bottomrule
\end{tabular}\end{adjustbox}\label{tab:app-quarterly-nonlinear}
\end{table}

\clearpage{}

\FloatBarrier

\newpage
\[
\mathbf{APPENDIX}
\]
\bigskip{}

\section{High Level Conditions}

Throughout, $C,C_{1},C_{2},...\in(0,\infty)$ denote generic constants
which may take different values at each appearance, even within the
same proof. 

The following inequality will be utilized repeatedly. For all $\theta>0$
and some $C_{\theta}\in(0,\infty)$ 
\begin{equation}
\left\vert \left\vert Y\right\vert ^{\theta}-\left\vert X\right\vert ^{\theta}\right\vert \leq C_{\theta}\left\{ \begin{array}{l}
\left(\left\vert X\right\vert ^{\theta-1}\left\vert Y-X\right\vert +\left\vert Y-X\right\vert ^{\theta}\right),\text{ }\theta\geq1\\
\left\vert Y-X\right\vert ^{\theta},\text{ }0<\theta<1
\end{array}\right.\text{.}\label{MVT}
\end{equation}
\bigskip{}

Below, we consider a set of \textbf{high-level (HL)} conditions that
closely resemble those in DK2021 and yield limit theory for locally
integrable and integrable kernel sample functionals of stochastic
arrays $X_{t}(n)$. These conditions are then shown to hold for HT-WNPs
under \textbf{Assumptions BFR }and\textbf{ MI}.

There is one notable difference between the \textbf{HL} conditions
used here and those considered by DK2021. In this work, we assume
that, upon normalization, $X_{t}(n)$ possesses a uniformly bounded
sub-\textgreek{\textalpha} moment, whereas DK2021 assume at least
two moments. Before presenting these conditions, we introduce some
notation. Let
\[
\Omega_{n}\left(\eta\right):=\left\{ \left\{ s,t\right\} \in\mathbb{N}:\left\lfloor \eta n\right\rfloor \leq s\leq\left\lfloor ns\right\rfloor ,\left\lfloor \eta n\right\rfloor +s\leq t\leq n\right\} ,
\]
for each $\eta\in(0,1/2)$.

\bigskip{}

\noindent\textbf{Assumption HL (high-level conditions)}\textit{\ }

\bigskip{}

\noindent\textbf{HL0}\textit{ Let }$\left\{ X_{t}\left(n\right)\right\} ^{n}_{t=1}$\textit{,
}$n\in\mathbb{N}$ \textit{be a random array and }$\left\{ \mathcal{F}_{t}\right\} ^{\infty}_{t=-\infty}$\textit{
a filtration such that }$X_{t}\left(n\right)$\textit{ is }$\mathcal{F}_{t}$\textit{-measurable
for all }$t$\textit{ and }$n$\textit{. }

\bigskip{}

\noindent\textbf{HL1} $X_{t}\left(n\right)=X_{t}\left(n\right)^{+}+X_{t}\left(n\right)^{-}+R_{t}(n)$\textit{,
where }$X_{t}\left(n\right)^{-}$\textit{ is }$\mathcal{F}_{0}$\textit{-measurable,
and }$\sup_{1\leq t\leq n}P\{\gamma^{-1}_{n}\lvert R_{t}(n)\rvert>\epsilon\}\rightarrow0$\textit{
for every }$\epsilon>0$\textit{.}

\bigskip{}

\noindent\textbf{HL2} \textit{Let }$\gamma_{n}>0$\textit{ a deterministic
sequence with }$\gamma_{n}\rightarrow\infty$\textit{. There are random
variables }$X^{+}$\textit{ and }$X^{-}$\textit{, where }$X^{+}$
\textit{has bounded density }$\Phi_{X^{+}}$\textit{ such that: for
every }$\delta\in(0,1)$ and $\left\{ t_{n}\right\} $\textit{ with
}$\lfloor n\delta\rfloor\leq t_{n}\leq n$

\textbf{(a)} $\gamma^{-1}_{n}X_{t_{n}}\left(n\right)^{+}\overset{d}{\rightarrow}X^{+}$\textit{,
conditionally on }$\mathcal{F}_{0}$\textit{ in the sense that for
all bounded and continuous $h:\mathbb{R}\rightarrow\mathbb{R}$ }
\[
\mathbf{E}\left[h\left(\gamma^{-1}_{n}X_{t_{n}}\left(n\right)^{+}\right)\mid\mathcal{F}_{0}\right]\overset{p}{\rightarrow}\int_{\mathbb{R}}h\left(x\right)\Phi_{X^{+}}(x)dx\text{\textit{;} \textit{and}}
\]

\textbf{(b)} $\gamma^{-1}_{n}X_{t_{n}}\left(n\right)^{-}\overset{d}{\rightarrow}X^{-}$,\textit{
and }$\gamma^{-1}_{n}\left[X_{n}\left(n\right)^{-}-X_{t_{n}}\left(n\right)^{-}\right]\overset{p}{\rightarrow}0$\textit{.}

\bigskip{}

\noindent\textbf{HL3} $\gamma^{-1}_{t}X_{t}\left(n\right)$\textit{
has Lebesgue density }$\mathcal{D}_{n,t}\left(x\right)$\textit{ such
that for some $n_{0}\geq t_{0}\geq1$},
\[
\sup_{n\geq n_{0},t_{0}\leq t\leq n}\sup_{x}\mathcal{D}_{n,t}\left(x\right)<\infty.
\]

\bigskip{}

\noindent\textbf{HL4} \textit{For every bounded and Lipschitz continuous
}$g:\mathbb{R}\rightarrow\mathbb{R}$ 
\[
\frac{1}{n}\sum^{n}_{t=1}g\left(\gamma^{-1}_{n}X_{t}\left(n\right)\right)=\frac{1}{n}\sum^{n}_{t=1}\mathbf{E}\left[g\left(\gamma^{-1}_{n}X_{t}\left(n\right)\right)\mid\mathcal{F}_{0}\right]+o_{p}(1).
\]

\noindent\textbf{HL5 }\textit{For all }$\alpha^{\prime}\in\left(0,\alpha\right)$
\[
\sup_{n\geq1,1\leq t\leq n}\mathbf{E}\left\vert \gamma^{-1}_{n}X_{t}\left(n\right)\right\vert ^{\alpha^{\prime}}<\infty\mathit{;}
\]

\noindent\textbf{HL6 }$\sup_{n\geq1}\frac{\gamma_{n}}{n}\sum^{n}_{t=1}\gamma^{-1}_{t}<\infty$\textit{.}

\bigskip{}

\noindent\textbf{HL7}\textit{ Let }$X_{0}(n)=0$\textit{ and }$t>s\geq0$\textit{.
Conditionally on $\mathcal{F}$}$_{s}$\textit{,} $\gamma^{-1}_{t-s}(X_{t}(n)-X_{s}(n))$\textit{
has density $\mathcal{D}$}$_{t,s,n}(x)$\textit{ such that} 
\[
\sup_{n\geq n_{0},0\leq s<t\leq n,t-s\geq t_{0}}\sup_{x}\mathcal{D}_{t,s,n}(x)<\infty.
\]

\noindent\textbf{HL8} \textit{For all }$q_{0},q_{1}>0$\textit{\ }
\[
\lim_{\eta\downarrow0}\lim\sup_{n\rightarrow\infty}\sup_{\left(s,t\right)\in\Omega_{n}\left(\eta\right)}\sup_{\left\vert x\right\vert \leq q_{0}\eta^{q_{1}}}\left\vert \mathcal{D}_{t,s,n}(x)-\mathcal{D}_{t,s,n}(0)\right\vert =0
\]

\noindent\textbf{HL9}

\textbf{(a)} $\lim_{\eta\rightarrow0}\lim_{n\rightarrow\infty}\frac{\gamma_{n}}{n}\sum^{n}_{t=\lfloor\left(1-\eta\right)n\rfloor}\gamma^{-1}_{t}=0;$

\textbf{(b) }$\lim_{\eta\rightarrow0}\lim_{n\rightarrow\infty}\frac{\gamma_{n}}{n}\sum^{\lfloor\eta n\rfloor}_{t=t_{0}}\gamma^{-1}_{t}=0$\textit{;}

\textbf{(c)} $\lim\sup_{n\rightarrow\infty}\frac{\gamma_{n}}{n}\sum^{n}_{t=t_{0}}\gamma^{-1}_{t}<\infty$\textit{;}

\textbf{(d) }\textit{for each }$\eta\in(0,1)$\textit{, there exist
}$l_{0},l_{1}>0$\textit{\  such that }
\[
\mathit{\lim}\inf_{n\rightarrow\infty}\gamma^{-1}_{n}\inf_{\left(s,t\right)\in\Omega_{n}\left(\eta\right)}\gamma_{t-s}\geq\eta^{l_{1}}/l_{0}\text{\textit{;}}
\]

\textbf{(e)} $\sup_{n\geq1}\sup_{1\leq t\leq n}\gamma^{-1}_{n}\gamma_{t}<\infty$.

\bigskip{}

Conditions \textbf{HL0}--\textbf{HL6} facilitate limit theory for
locally integrable functionals of HT-WNP arrays upon normalization.
For kernel functionals, \textbf{HL7--HL9} are required in addition
to these conditions. We note that \textbf{HL7} is a strengthened version
of \textbf{HL3}. Limit theory for locally integrable and kernel functionals
is provided in\textbf{ Theorems \ref{MainThmHL}} and \textbf{\ref{MainKernelHL}}\
below. \bigskip{}

\begin{theorem}\label{MainThmHL}\textit{Suppose that }$f:\mathbb{R}\rightarrow\mathbb{R}$\textit{
is locally Lebesgue integrable, and }\textbf{HL0--HL4}\textit{, }\textbf{HL6}\textit{
hold. Furthermore, suppose }$f$\textit{ is bounded or that the following
hold:}
\begin{itemize}
\item[(i)] $\int_{\mathbb{R}}\left|f\left(x+y\right)\right|\Phi(x)dx<\infty$\textit{,
for all }$y\in\mathcal{Y}$\textit{, and }$P\left(X^{-}\in\mathcal{Y}\right)=1$\textit{;}
\item[(ii)] \textbf{HL5} \textit{holds;}
\item[(iii)] \textit{ for some} $q\in(0,\alpha)$, $\left\vert f\left(x\right)\right\vert =O(\left\vert x\right\vert ^{q})$\textit{,
as }$\left\vert x\right\vert \rightarrow\infty$\textit{;} 
\item[(iv)] \textit{ $\frac{1}{n}\sum^{t_{0}-1}_{t=1}f\left(\gamma^{-1}_{n}X_{t}\left(n\right)\right)=o_{p}(1)$,
with $t_{0}$ as in} \textbf{HL3}\textit{.}
\end{itemize}
\noindent\textit{Then as }$n\rightarrow\infty$\textit{,}
\begin{equation}
\frac{1}{n}\sum^{n}_{t=1}f\left(\gamma^{-1}_{n}X_{t}\left(n\right)\right)\overset{d}{\rightarrow}\int_{\mathbb{R}}f\left(x+X^{-}\right)\Phi_{X^{+}}(x)dx.\label{Thm2.1}
\end{equation}

\end{theorem}

\bigskip{}

\noindent\textbf{Remark}. If $\tilde{f}:\mathbb{R\rightarrow\overline{\mathbb{R}}}$
is such that $\tilde{f}(x)=f(x)$ almost everywhere, and \textit{$\frac{1}{n}\sum^{t_{0}-1}_{t=1}\tilde{f}\left(\gamma^{-1}_{n}X_{t}\left(n\right)\right)=o_{p}(1)$,
with $t_{0}$ as in} \textbf{HL3}\textit{.} then (\ref{Thm2.1}) holds
with $f$ replaced by $\tilde{f}$.

\bigskip{}

\noindent\textbf{Proof of Theorem \ref{MainThmHL}}. The proof follows
from Theorem 3.1 in DK2021. For convenience we provide a summary of
the key steps. Given that $f$ is locally integrable, for each $\varepsilon>0$,
there is a Lipschitz continuous $f_{\varepsilon}$ that approximates
$f$ in $L_{1}$. In particular,
\[
\int_{\left|x\right|\leq\varepsilon^{-1}}\left|f\left(x\right)-f_{\varepsilon}(x)\right|dx<\varepsilon\:\mathrm{and}\;f_{\varepsilon}(x)=0,\left|x\right|>\varepsilon^{-1}.
\]
 It can be shown that for each $\varepsilon>0$ and $0<\delta<1$
\begin{equation}
\frac{1}{n}\sum^{n}_{t=1}f\left(\gamma^{-1}_{n}X_{t}\left(n\right)\right)=G_{\varepsilon}\left[\gamma^{-1}_{n}X_{n}(n)^{-}\right]+R_{n,\varepsilon,\delta}\label{MainThmHL1}
\end{equation}
with $G_{\varepsilon}\left[.\right]=\int_{\mathbb{R}}f_{\varepsilon}\left(x+.\right)dF_{X^{+}}(x)$
continuous and $R_{n,\varepsilon,\delta}=o_{p}(1)$ as $n\rightarrow\infty$
first, then $\delta\rightarrow0$, and then $\varepsilon\rightarrow0$.
Moreover, as $n\rightarrow\infty$ we have
\begin{equation}
G_{\varepsilon}\left[\gamma^{-1}_{n}X_{n}(n)^{-}\right]\overset{d}{\rightarrow}G_{\varepsilon}\left[X^{-}\right],\label{MainThmHL2}
\end{equation}
and
\begin{equation}
G_{\varepsilon}\left[X^{-}\right]=\int_{\mathbb{R}}f\left(x+X^{-}\right)\Phi_{X^{+}}(x)dx+o_{a.s.}(1),\label{MainThmHL3}
\end{equation}
as $\varepsilon\rightarrow0$. Finally, we note that the approximation
with respect to to $\varepsilon$ can be dropped, if $f$ is bounded
and Lipschitz. In particular, (\ref{MainThmHL1}) simplifies to
\begin{equation}
\frac{1}{n}\sum^{n}_{t=1}f\left(\gamma^{-1}_{n}X_{t}\left(n\right)\right)=G\left[\gamma^{-1}_{n}X_{n}(n)^{-}\right]+R_{n,\delta}\label{MainThmHL3b}
\end{equation}
with $R_{n,\delta}=o_{p}(1)$ as $n\rightarrow\infty$ first, then
$\delta\rightarrow0$, where $G[.]=\int_{\mathbb{R}}f\left(x+.\right)\Phi_{X^{+}}(x)dx$.
$\blacksquare$

\bigskip{}

\begin{theorem} \label{MainKernelHL} \textit{Let }$\{h_{n}\}$\textit{
a positive sequence with }$\gamma^{-1}_{n}h_{n}+\gamma_{n}(nh_{n})^{-1}\rightarrow0$.\textit{
In addition, }$\frac{\gamma_{n}}{h_{n}n}\sum^{t_{0}-1}_{t=1}K\left(\frac{X_{t}(n)-x}{h_{n}}\right)=o_{p}(1)$,\textit{
where }$t_{0}$\textit{ is as in} \textbf{HL7}. \textit{Then under
}\textbf{Assumptions K}, \textbf{HL0--HL4} \textit{and }\textbf{HL6--HL9}\textit{,
as }$n\rightarrow\infty$,
\[
\frac{\gamma_{n}}{h_{n}n}\sum^{n}_{t=1}K\left(\frac{X_{t}(n)-x}{h_{n}}\right)\overset{d}{\rightarrow}\Phi_{X^{+}}(-X^{-})\int_{\mathbb{R}}K(u)du.
\]
\end{theorem}

\bigskip{}

\noindent\textbf{Proof of Theorem \ref{MainKernelHL}}. The proof
is identical to that of Theorem 3.3 in DK2021. $\blacksquare$

\section{Proof of Main Limit Theorems}

\noindent\textbf{Proof of Theorem \ref{MainLL}}. In view of \textbf{Theorem
\ref{MainThmHL}}, it suffices to verify conditions (i), (iii), (iv),
of that theorem, together with \textbf{HL0}-\textbf{HL6}. We begin
by verifying \textbf{HL0}-\textbf{HL6}.
\begin{itemize}
\item[\textbf{HL0}-\textbf{1}] : Under \textbf{Assumptions INN}, \textbf{BFR} and \textbf{MI}, \textbf{HL0}
and \textbf{HL1} hold trivially, by setting $X_{t}(n)^{+}=x^{+}_{t}(n)$,
$X_{t}(n)^{-}=x^{-}_{t}(n)$ and $R_{t}(n)=0$ where $x^{+}_{t}(n)$,
$x^{-}_{t}(n)$ are defined in (\ref{linprocsplit}), and by taking
$\mathcal{F}_{0}=\sigma\left(\xi_{0},\xi_{-1},...\right)$.
\item[\textbf{HL2}] : Let $\gamma_{n}$ be as in \textbf{Definition 1}, let $X^{+}$,
$X^{-}$ be as in \textbf{Lemma \ref{WCon}}. By \textbf{Definition
2} ($\alpha$-stable distribution), $X^{+}$ has an integrable characteristic
function and therefore admits a bounded density. Moreover, \textbf{HL2(a)}
and the first part of \textbf{HL2(b)} follow immediately from \textbf{Lemma}
\textbf{\ref{WCon}}. Finally, the second part of \textbf{HL2(b)}
follows from \textbf{Lemma \ref{lem:FSalphaM}(b)}.
\item[\textbf{HL3}] : The result follows immediately from\textbf{ Lemma \ref{UBD}}.
In particular, that lemma shows that $\gamma^{-1}_{n}x\left(n\right)$
and $\gamma^{-1}_{n}x^{+}_{t}\left(n\right)$ have integrable characteristic
functions which is sufficient for the result; see the verification
of\textbf{ HL8} in the proof of \textbf{Theorem \ref{MainLL_kernel}}
below for more details.
\item[\textbf{HL4}] : This follows from \textbf{Lemma \ref{MartingApprox}}.
\item[\textbf{HL5}] : This is a direct consequence of \textbf{Lemma \ref{lem:FSalphaM}(a)}. 
\item[\textbf{HL6}] : It follows from \textbf{Definition 1} and the subsequent remarks
that $\gamma_{n}$ is is regularly varying $n$ with index $0\leq p<1$.
In particular, under \textbf{BFR}, $p=0$, whereas under \textbf{MI,}
$p=\alpha_{\kappa}\mathcal{H}$, where $0\leq\alpha_{\kappa}<1$ and
$0<\mathcal{H}<1$. Consequently, \textbf{HL6} follows from standard
arguments for regularly varying sequences; see, e.g., Lemma B.2 of
DK2021.
\end{itemize}
\bigskip{}

Finally, we verify conditions (i), (iii), (iv) of \textbf{Theorem
\ref{MainThmHL}}.
\begin{itemize}
\item[Condition (i)] : For each $y\in\mathbb{R}$ define $Y=X^{+}+y$. Note that $Y$
is $\alpha$-stable with location parameter $y$ and therefore $\mathbf{E}\left\vert Y\right\vert ^{\alpha^{\prime}}<\infty$
for all $\alpha^{\prime}\in(0,\alpha)$ -see e.g. Samorodnitsky and
Taqqu (1994; Property 1.2.16. p.18). Furthermore, $\Phi_{Y}$, the
density of $Y$, exists and is bounded due to the fact that $Y$ has
an integrable characteristic function. Then for $m>0$ large enough,
\[
\int_{\mathbb{R}}\left\vert f\left(x+y\right)\right\vert \Phi_{X^{+}}(x)dx=\mathbf{E}\left|f(Y)\right|\leq_{(1)}\int_{\left\vert u\right\vert \leq m}\left\vert f(u)\right\vert \Phi_{Y}(u)du+C\mathbf{E}\left\vert Y\right\vert ^{\alpha^{\prime}}
\]
\[
\leq\sup_{u}\Phi_{Y}(u)\int_{\left\vert u\right\vert \leq m}\left\vert f\left(u\right)\right\vert du+C\mathbf{E}\left\vert Y\right\vert ^{\alpha^{\prime}}<\infty,
\]
where $\leq_{(1)}$ follows from condition (ii) of \textbf{Theorem
\ref{MainLL}}.
\item[Condition (iii)] : this is identical to condition (ii) of \textbf{Theorem \ref{MainLL}}.
\item[Condition (iv)] : It can be readily seen that condition (iii.a) of \textbf{Theorem
\ref{MainLL}} implies condition (iv). \\
Moreover, under condition (iii.b), \textbf{HL4} holds with $t_{0}=1$.
To see this, first recall that
\[
x_{t}=\sum^{\infty}_{k=0}a_{k,t}\xi_{t-k},
\]
where
\[
a_{k,t}(n)=\sum^{(t-1)\wedge k}_{j=0}\phi_{j}c_{k-j}=\begin{cases}
\sum^{k}_{j=0}\phi_{j}(n)c_{k-j}=:a_{k} & \text{if }0\leq k\leq t-1,\\
\sum^{t-1}_{j=0}\phi_{j}(n)c_{k-j}=:a^{-}_{k,t} & \text{if }k\geq t,
\end{cases}
\]
and $\phi_{0}(n)\neq0$, for all $n$. Note that, under \textbf{BFR}
$\phi_{0}(n)=\phi_{0}\neq0$ by assumption, whereas under \textbf{MI},
$\phi_{j}(n)=\rho^{j}_{n}$. Hence, in the latter case $\phi_{0}(n)=1$.
In addition, under \textbf{BFR} and \textbf{MI}, there is at least
one coefficient $\left\{ c_{j}\right\} _{j\in\mathbb{Z}_{+}}$ that
is non-zero. Let $k^{*}$ be the minimal index of the non-zero coefficients.
Then for all $t\geq1$,
\begin{equation}
\left|a_{k^{*},t}(n)\right|\geq\textrm{inf}_{n\in\mathbb{N}}\left|\phi_{0}(n)c_{k^{*}}\right|>0.\label{HL4.1}
\end{equation}
Next, note that by \textbf{HL3}, $\gamma^{-1}_{t}x_{t}$ admits a
density $\mathcal{D}_{n,t}$ such that
\[
\textrm{sup}_{n\geq n_{0},t_{0}\leq t\leq n}\textrm{sup}_{x}\mathcal{D}_{n,t}(x)<\infty,
\]
for some $n_{0}\geq1$ and $t_{0}\geq1$. The required result will
follow if we show that for any $t_{0}\geq1$,
\[
\textrm{sup}_{n\geq n_{0},1\leq t\leq t_{0}}\textrm{sup}_{x}\mathcal{D}_{n,t}(x)<\infty.
\]
Indeed, by (\ref{HL4.1}), for all $1\leq t\leq t_{0}$ we get
\[
\int_{\mathbb{R}}\left|\mathbf{E}\left(e^{i\lambda\gamma^{-1}_{t}x_{t}}\right)\right|d\lambda=\int_{\mathbb{R}}\prod^{\infty}_{k=0}\left|\psi_{\xi}\left(\gamma^{-1}_{t}a_{k,t}(n)\lambda\right)\right|d\lambda
\]
\[
\leq\frac{\textrm{sup}_{1\leq t\leq t_{0}}\gamma_{t}}{\textrm{inf}_{1\leq t\leq t_{0},n\in\mathbb{N}}\left|a_{k^{*},t}(n)\right|}\int_{\mathbb{R}}\left|\psi_{\xi}\left(\lambda\right)\right|d\lambda<\infty,
\]
as required. $\blacksquare$
\end{itemize}
\bigskip{}

\noindent\textbf{Proof of Theorem \ref{MainLL_kernel}}. In view
of \textbf{Theorem \ref{MainKernelHL}}, it suffices to show that
Conditions \textbf{HL0--HL4} and \textbf{HL6--HL9} hold. Conditions
\textbf{HL0--HL6} have been already established in the proof of \textbf{Theorem
\ref{MainThmHL}}. We shall show that \textbf{HL7--HL9} hold.
\begin{itemize}
\item[\textbf{HL7}] :\textbf{ }We show that \textbf{HL7} holds with $X_{t}(n)=x_{t}(n)$,
$t\geq1$. In view of (\ref{linprocsplit}) write 
\begin{eqnarray*}
x^{+}_{t}\left(n\right) & = & \sum^{t-1}_{k=0}a_{k}\left(n\right)\xi_{t-k}=\sum^{t}_{l=1}a_{t-l}\left(n\right)\xi_{l}.\\
 &  & \left[l=t-k\rightarrow k=t-l\right]
\end{eqnarray*}
Hence, for $t>s$
\[
x_{t}\left(n\right)-x_{s}\left(n\right)=\sum^{t}_{l=1}a_{t-l}\left(n\right)\xi_{l}-\sum^{s}_{l=1}a_{s-l}\left(n\right)\xi_{l}+x^{-}_{t}\left(n\right)-x^{-}_{s}\left(n\right)
\]
\[
=\sum^{t}_{l=s+1}a_{t-l}\left(n\right)\xi_{l}+\sum^{s}_{l=1}\left[a_{t-l}\left(n\right)-a_{s-l}\left(n\right)\right]\xi_{l}+x^{-}_{t}\left(n\right)-x^{-}_{s}\left(n\right)=:\sum^{t}_{l=s+1}a_{t-l}\left(n\right)\xi_{l}+\mathcal{X}_{t,s}(n).
\]
The first term on the r.h.s. in the last line above is 
\begin{eqnarray*}
\sum^{t}_{l=s+1}a_{t-l}\left(n\right)\xi_{l} & = & \sum^{t-s-1}_{k=0}a_{k}\left(n\right)\xi_{t-k}\overset{d}{=}x^{+}_{t-s}\left(n\right).\\
 &  & \left[l=t-k\rightarrow k=t-l\right]
\end{eqnarray*}
It follows from\textbf{ Lemma \ref{UBD}} that $\gamma^{-1}_{t}x^{+}_{t}\left(n\right)$
possesses a density, $\mathcal{D}^{+}_{t,n}(x)$\ say, such that
$\sup_{n\geq n_{0},t_{0}\leq t\leq n}\sup_{x}\mathcal{D}^{+}_{t,n}(x)<\infty$,
for some $n_{0}\geq t_{0}\geq1$. In view of the above, conditionally
on $\mathcal{F}_{s}$, $\gamma^{-1}_{t-s}\left[x_{t}\left(n\right)-x_{s}\left(n\right)\right]$
has a density $\mathcal{D}_{t,s,n}(x)$ such that $\mathcal{D}_{t,s,n}(x)$
$=\mathcal{D}^{+}_{t-s,n}(x-\mathcal{X}_{t,s}(n))\ $and
\[
\sup_{n\geq n_{0},0\leq s<t\leq n,t-s\geq t_{0}}\sup_{x}\mathcal{D}^{+}_{t-s,n}(x-\mathcal{X}_{t,s}(n))\leq\sup_{n\geq n_{0},t_{0}\leq t\leq n}\sup_{x}\mathcal{D}^{+}_{t,n}(x)<\infty.
\]
\item[\textbf{HL8}] : We begin with some notation. For a function $\varphi(x),$ let
$\tilde{\varphi}(\lambda)$ denote its Fourier transform, defined
by 
\[
\tilde{\varphi}(\lambda)=\frac{1}{\sqrt{2\pi}}\int_{\mathbb{R}}e^{-i\lambda x}\varphi(x)dx.
\]
By Fourier inversion; see, for example, Theorem 5.1 of Lang (1993),
we have
\[
\varphi(-x)=\frac{1}{\sqrt{2\pi}}\int_{\mathbb{R}}e^{-ix\lambda}\tilde{\varphi}(\lambda)d\lambda,
\]
and
\[
\varphi(x)=\frac{1}{\sqrt{2\pi}}\int_{\mathbb{R}}e^{ix\lambda}\tilde{\varphi}(\lambda)d\lambda.
\]
Next, with $\mathcal{D}^{+}_{t,n}(x)$ as defined in the verification
of \textbf{HL7} above, we have 
\[
\int_{\mathbb{R}}\left\vert \lambda\right\vert \left\vert \mathcal{\tilde{D}}^{+}_{t,n}(\lambda)\right\vert d\lambda=\int_{\mathbb{R}}\left\vert \lambda\right\vert \left\vert \int_{\mathbb{R}}e^{-i\lambda x}\mathcal{D}^{+}_{t,n}(x)dx\right\vert d\lambda=\int_{\mathbb{R}}\left\vert \lambda\right\vert \left\vert \mathbf{E}e^{-i\lambda x^{+}_{t}\left(n\right)}\right\vert d\lambda.
\]
Combining the above with \textbf{Lemma \ref{UBD}} we obtain
\[
\sup_{n\geq n_{0},t_{0}\leq t\leq n}\int_{\mathbb{R}}\left\vert \lambda\right\vert \left\vert \mathcal{\tilde{D}}^{+}_{t,n}(\lambda)\right\vert d\lambda=\frac{1}{\sqrt{2\pi}}\sup_{n\geq n_{0},t_{0}\leq t\leq n}\int_{\mathbb{R}}\left\vert \lambda\right\vert \left\vert \mathbf{E}e^{i\lambda x^{+}_{t}\left(n\right)}\right\vert d\lambda<\infty.
\]
Since
\[
\left\vert \lambda x\right\vert \wedge1\leq\left(1+\left\vert \lambda\right\vert \right)\left\vert x\right\vert \wedge\left(1+\left\vert \lambda\right\vert \right)=\left(1+\left\vert \lambda\right\vert \right)\left(\left\vert x\right\vert \wedge1\right),
\]
it follows that
\[
\sup_{\left(s,t\right)\in\Omega_{n}\left(\eta\right)}\sup_{\left\vert x\right\vert \leq q_{0}\eta^{q_{1}}}\left\vert \mathcal{D}_{t,s,n}(x)-\mathcal{D}_{t,s,n}(0)\right\vert 
\]
\[
=\sup_{\left(s,t\right)\in\Omega_{n}\left(\eta\right)}\sup_{\left\vert x\right\vert \leq q_{0}\eta^{q_{1}}}\left\vert \mathcal{D}^{+}_{t-s,n}(x-\mathcal{X}_{t,s}(n))-\mathcal{D}^{+}_{t-s,n}(\mathcal{X}_{t,s}(n))\right\vert 
\]
\[
=\sup_{\left(s,t\right)\in\Omega_{n}\left(\eta\right)}\sup_{\left\vert x\right\vert \leq q_{0}\eta^{q_{1}}}\left\vert \frac{1}{\sqrt{2\pi}}\int_{\mathbb{R}}e^{i\left[x-\mathcal{X}_{t,s}(n))\right]\lambda}\mathcal{\tilde{D}}^{+}_{t-s,n}(\lambda)d\lambda-\frac{1}{\sqrt{2\pi}}\int_{\mathbb{R}}e^{-i\mathcal{X}_{t,s}(n))\lambda}\mathcal{\tilde{D}}^{+}_{t-s,n}(\lambda)d\lambda\right\vert 
\]
\[
\leq\frac{2}{\sqrt{2\pi}}\sup_{\left(s,t\right)\in\Omega_{n}\left(\eta\right)}\sup_{\left\vert x\right\vert \leq q_{0}\eta^{q_{1}}}\int_{\mathbb{R}}\left(\left\vert \lambda x\right\vert \wedge1\right)\left\vert \mathcal{\tilde{D}}^{+}_{t-s,n}(\lambda)\right\vert d\lambda
\]
\[
\leq\frac{2}{\sqrt{2\pi}}\sup_{\left(s,t\right)\in\Omega_{n}\left(\eta\right)}\sup_{\left\vert x\right\vert \leq q_{0}\eta^{q_{1}}}\left(\left\vert x\right\vert \wedge1\right)\int_{\mathbb{R}}\left(1+\left\vert \lambda\right\vert \right)\left\vert \mathcal{\tilde{D}}^{+}_{t-s,n}(\lambda)\right\vert d\lambda
\]
\[
\leq\frac{2}{\sqrt{2\pi}}\sup_{\left\vert x\right\vert \leq q_{0}\eta^{q_{1}}}\left(\left\vert x\right\vert \wedge1\right)\sup_{n\geq n_{0},t_{0}\leq t\leq n}\int_{\mathbb{R}}\left(1+\left\vert \lambda\right\vert \right)\left\vert \mathcal{\tilde{D}}^{+}_{t,n}(\lambda)\right\vert d\lambda
\]
\[
\leq C\sup_{\left\vert x\right\vert \leq q_{0}\eta^{q_{1}}}\left(\left\vert x\right\vert \wedge1\right)\rightarrow0,
\]
as $\eta\rightarrow0$.
\item[\textbf{HL9}] : Since $\gamma_{n}$ is regularly varying, \textbf{HL9} follows
directly from the proof of Theorem 3.4 in DK2021. $\blacksquare$ 
\end{itemize}
\bigskip{}

\noindent\textbf{Proof of Theorem \ref{MCLT-LL}}. By \textbf{Theorem
\ref{MainLL}},
\[
\frac{1}{n}\sum^{n}_{t=1}f(\gamma^{-1}_{n}x_{t-1})^{2}\overset{d}{\rightarrow}\int_{\mathbb{R}}f^{2}(x)\varrho(x)dx.
\]
In view of this, it follows from Wang's (2014) MCLT -cf. Wang (2014;
cf. Thm 2.1)- that 
\[
\left[\frac{1}{\sqrt{n}}\sum^{n}_{t=1}f(\gamma^{-1}_{n}x_{t-1})u_{t},\:\frac{1}{n}\sum^{n}_{t=1}f(\gamma^{-1}_{n}x_{t-1})^{2}\right]
\]
\begin{equation}
\overset{d}{\rightarrow}\left[MN\left[0,\sigma^{2}_{u}\int_{\mathbb{R}}f^{2}(x)\varrho(x)dx\right],\int_{\mathbb{R}}f^{2}(x)\varrho(x)dx\right],\label{MainMCLT1}
\end{equation}
if we show that
\begin{equation}
n^{-1/2}\max_{1\leq t\leq n}\left\vert f(\gamma^{-1}_{n}x_{t})\right\vert =o_{p}(1).\label{MainMCLT2}
\end{equation}
Eq. (\ref{MainMCLT2}) follows from arguments identical to those used
for the proof of Thm 4.1 in DK2021. We note that eq. (2.2) and (2.3)
in Wang (2014) are not satisfied under the conditions of \textbf{Theorem
\ref{MCLT-LL}}; however, they are unnecessary in the present setting.
Indeed, here the mixing variate $X^{-}$ is the weak limit of an $\mathcal{F}_{0}$-measurable
term. Moreover, a careful reading of Wang's (2014) derivations shows
that (\ref{MainMCLT1}) can be strengthened to 
\[
\left[\frac{1}{\sqrt{n}}\sum^{n}_{t=1}f(\gamma^{-1}_{n}x_{t-1})u_{t},\frac{1}{n}\sum^{n}_{t=1}f(\gamma^{-1}_{n}x_{t-1})^{2},\text{ }\gamma^{-1}_{n}x^{-}_{n}\right]
\]
\begin{equation}
\overset{d}{\rightarrow}\left[MN\left[0,\sigma^{2}_{u}\int_{\mathbb{R}}f^{2}(x)\varrho(x)dx\right],\int_{\mathbb{R}}f^{2}(x)\varrho(x)dx,X^{-}\right].\label{MainMCLT3}
\end{equation}
Next, note that by (\ref{MainThmHL1})-(\ref{MainThmHL3}) for any
$0<\delta<1$, $\varepsilon>0$ there is some continuous $G_{\varepsilon}\left[.\right]$
and $R_{n,\varepsilon,\delta}$ -cf. (\ref{MainThmHL1})- such that
\[
\frac{1}{n}\sum^{n}_{t=1}g(\gamma^{-1}_{n}x_{t-1})=G_{\varepsilon}\left[\gamma^{-1}_{n}x^{-}_{n}\right]+R_{n,\varepsilon,\delta},
\]
and $G_{\varepsilon}\left[\gamma^{-1}_{n}x^{-}_{n}\right]\overset{d}{\rightarrow}G_{\varepsilon}\left[X^{-}\right]$,
as $n\rightarrow\infty$, whilst $R_{n,\varepsilon,\delta}=o_{p}(1)$
as as $n\rightarrow\infty$ first, then $\delta\rightarrow0$, and
then $\varepsilon\rightarrow0$. Furthermore, $G_{\varepsilon}\left[X^{-}\right]=\int_{\mathbb{R}}g(x)\varrho(x)dx+o_{a.s.}(1)$
as $\varepsilon\rightarrow0$. In view of this, (\ref{MainMCLT3}),
the continuous mapping theorem, and Billingsley (1968, Thm 4.2) as
$n\rightarrow\infty$ first, then $\delta\rightarrow0$, and then
$\varepsilon\rightarrow0$
\[
\left[\frac{1}{\sqrt{n}}\sum^{n}_{t=1}f(\gamma^{-1}_{n}x_{t-1})u_{t},\text{ }\frac{1}{n}\sum^{n}_{t=1}f(\gamma^{-1}_{n}x_{t-1})^{2},\text{ }\frac{1}{n}\sum^{n}_{t=1}g(\gamma^{-1}_{n}x_{t-1})\right]
\]
\[
=\left[\frac{1}{\sqrt{n}}\sum^{n}_{t=1}f(\gamma^{-1}_{n}x_{t-1})u_{t},\text{ }\frac{1}{n}\sum^{n}_{t=1}f(\gamma^{-1}_{n}x_{t-1})^{2},\text{ }G_{\varepsilon}\left[\gamma^{-1}_{n}x^{-}_{n}\right]+R_{n,\varepsilon,\delta}\right]
\]
\[
=\left[\frac{1}{\sqrt{n}}\sum^{n}_{t=1}f(\gamma^{-1}_{n}x_{t-1})u_{t},\text{ }\frac{1}{n}\sum^{n}_{t=1}f(\gamma^{-1}_{n}x_{t-1})^{2},\text{ }G_{\varepsilon}\left[\gamma^{-1}_{n}x^{-}_{n}\right]\right]+o_{p}(1)\text{\:(}n\rightarrow\infty,\delta\rightarrow0,\varepsilon\rightarrow0\text{)}
\]
\[
\overset{d}{\rightarrow}\left[MN\left[0,\sigma^{2}_{u}\int_{\mathbb{R}}f^{2}(x)\varrho(x)dx\right],\int_{\mathbb{R}}f^{2}(x)\varrho(x)dx,G_{\varepsilon}[X^{-}]\right]\:\text{ (}n\rightarrow\infty\text{)}
\]
\[
=\left[MN\left[0,\sigma^{2}_{u}\int_{\mathbb{R}}f^{2}(x)\varrho(x)dx\right],\int_{\mathbb{R}}f^{2}(x)\varrho(x)dx,\int_{\mathbb{R}}g(x)\varrho(x)dx\right]+o_{a.s.}(1)\:(\varepsilon\rightarrow0\text{),}
\]
as required. $\blacksquare$

\bigskip{}

\noindent\textbf{Proof of Theorem \ref{MCLT-LLker}}. The proof is
similar to that of \textbf{Theorem \ref{MCLT-LL}} above. By \textbf{Theorem
\ref{MainLL_kernel}},
\[
\frac{\gamma_{n}}{h_{n}n}\sum^{n}_{t=1}K\left(\frac{x_{t-1}-x}{h_{n}}\right)\overset{d}{\rightarrow}\varrho(0)\int_{\mathbb{R}}K(u)du.
\]
In view of this convergence and Wang (2014; cf. Thm 2.1), we have
\[
\left[\sqrt{\frac{\gamma_{n}}{h_{n}n}}\sum^{n}_{t=1}K\left(\frac{x_{t-1}-x}{h_{n}}\right)u_{t},\text{ }\frac{\gamma_{n}}{h_{n}n}\sum^{n}_{t=1}K\left(\frac{x_{t-1}-x}{h_{n}}\right)^{2}\right]
\]
\begin{equation}
\overset{d}{\rightarrow}\left[MN\left[0,\sigma^{2}_{u}\varrho(0)\int_{\mathbb{R}}K^{2}(u)du\right],\text{ }\varrho(0)\int_{\mathbb{R}}K^{2}(u)du\right],\label{MCLT-LLker1}
\end{equation}
provided that
\begin{equation}
\sqrt{\frac{\gamma_{n}}{h_{n}n}}\max_{1\leq t\leq n}\left\vert K\left(\frac{x_{t-1}-x}{h_{n}}\right)\right\vert =o_{p}(1).\label{MCLT-LLker02}
\end{equation}
The latter approximation follows from an argument similar to the one
used for the proof of Thm 4.1 in DK2021. We note that the mixing variate
in the mixed Gaussian limit in (\ref{MCLT-LLker1}) is the weak limit
of $\gamma^{-1}_{n}x^{-}_{n}$, which is $\mathcal{F}_{0}$-measurable.
Therefore, the requirements in eq. (2.2) and (2.3) of Wang (2014)
are unnecessary in the present setting. Furthermore, as mentioned
in the previous proof, (\ref{MCLT-LLker1}) can be generalized as
follows:
\[
\left[\sqrt{\frac{\gamma_{n}}{h_{n}n}}\sum^{n}_{t=1}K\left(\frac{x_{t-1}-x}{h_{n}}\right)u_{t},\text{ }\frac{\gamma_{n}}{h_{n}n}\sum^{n}_{t=1}K\left(\frac{x_{t-1}-x}{h_{n}}\right)^{2},\text{ }\gamma^{-1}_{n}x^{-}_{n}\right]
\]
\begin{equation}
\overset{d}{\rightarrow}\left[MN\left[0,\sigma^{2}_{u}\varrho(0)\int_{\mathbb{R}}K^{2}(u)du\right],\text{ }\varrho(0)\int_{\mathbb{R}}K^{2}(u)du,\text{ }X^{-}\right].\label{MCLT-LLker2}
\end{equation}
Next, let $\varphi_{\eta^{2}}(x)$ the density of a $N(0,\eta^{2})$
variate. It can be shown (see e.g. DK2021, Theorem 3.3) that under
our assumptions as $n\rightarrow\infty$ first and the as $\eta\rightarrow0$
\[
\frac{\gamma_{n}}{h_{n}n}\sum^{n}_{t=1}L\left(\frac{x_{t-1}-x}{h_{n}}\right)=\int_{\mathbb{R}}L(u)du\cdot\frac{1}{n}\sum^{n}_{t=1}\varphi_{\eta^{2}}\left(\gamma^{-1}_{n}x_{t-1}\right)+o_{p}(1)
\]
and 
\begin{gather*}
\frac{1}{n}\sum^{n}_{t=1}\varphi_{\eta^{2}}\left(\gamma^{-1}_{n}x_{t}\right)\overset{d}{\rightarrow}\int_{\mathbb{R}}\varphi_{\eta^{2}}\varrho(x)dx,\text{ as }n\rightarrow\infty\\
\overset{a.s.}{\rightarrow}\varrho(0),\text{ as }\eta\rightarrow0.
\end{gather*}
More specifically, by (\ref{MainThmHL3b}), noting that $\varphi_{\eta^{2}}$
is bounded and Lipschitz, there is a remainder $R_{n,\delta}=o_{p}(1)$
as $n\rightarrow\infty$ first followed by $\delta\rightarrow0$ and
a continuous $G_{\eta}\left[.\right]\left(=\int_{\mathbb{R}}\varphi_{\eta^{2}}(x+.)dF_{X^{+}}(x)\right)$
such that 
\[
\frac{1}{n}\sum^{n}_{t=1}\varphi_{\eta^{2}}\left(\gamma^{-1}_{n}x_{t}\right)=G_{\eta}\left[\gamma^{-1}_{n}x^{-}_{n}\right]+R_{n,\delta}
\]
\[
=G_{\eta}\left[\gamma^{-1}_{n}x^{-}_{n}\right]+o_{p}(1)\text{ (}n\rightarrow\infty,\delta\rightarrow0\text{)}
\]
\[
\overset{d}{\rightarrow}G_{\eta}\left[X^{-}\right]\text{ (}n\rightarrow\infty\text{)}
\]
\[
=\int_{\mathbb{R}}\varphi_{\eta^{2}}\varrho(x)dx+o(1)\overset{a.s.}{\rightarrow}\varrho(0)\text{ (}\eta\rightarrow0\text{)}.
\]
In view of the above, (\ref{MCLT-LLker2}), the continuous mapping
theorem, and Billingsley (1968, Thm 4.2) we obtain
\[
\left[\sqrt{\frac{\gamma_{n}}{h_{n}n}}\sum^{n}_{t=1}K\left(\frac{x_{t-1}-x}{h_{n}}\right)u_{t},\text{ }\frac{\gamma_{n}}{h_{n}n}\sum^{n}_{t=1}K\left(\frac{x_{t-1}-x}{h_{n}}\right)^{2},\text{ }\frac{\gamma_{n}}{h_{n}n}\sum^{n}_{t=1}L\left(\frac{x_{t-1}-x}{h_{n}}\right)\right]
\]
\[
=\left[\sqrt{\frac{\gamma_{n}}{h_{n}n}}\sum^{n}_{t=1}K\left(\frac{x_{t-1}-x}{h_{n}}\right)u_{t},\text{ }\frac{\gamma_{n}}{h_{n}n}\sum^{n}_{t=1}K\left(\frac{x_{t-1}-x}{h_{n}}\right)^{2},\text{ }\int_{\mathbb{R}}L(u)du\cdot\frac{1}{n}\sum^{n}_{t=1}\varphi_{\eta^{2}}\left(\gamma^{-1}_{n}x_{t}\right)\right]+o_{p}(1)
\]
\[
\text{ (}n\rightarrow\infty\text{)}
\]
\[
=\left[\sqrt{\frac{\gamma_{n}}{h_{n}n}}\sum^{n}_{t=1}K\left(\frac{x_{t}-x}{h_{n}}\right)u_{t},\text{ }\text{ }\frac{\gamma_{n}}{h_{n}n}\sum^{n}_{t=1}K\left(\frac{x_{t-1}-x}{h_{n}}\right)^{2},\text{ }\int_{\mathbb{R}}L(u)du\cdot G_{\eta}\left[\gamma^{-1}_{n}x^{-}_{n}\right]\right]+o_{p}(1)
\]
\[
\text{ (}n\rightarrow\infty\text{)}
\]
\[
\overset{d}{\rightarrow}\left[MN\left[0,\sigma^{2}_{u}\varrho(0)\int_{\mathbb{R}}K^{2}(u)du\right],\text{ }\varrho(0)\int_{\mathbb{R}}K^{2}(u)du\text{, }\int_{\mathbb{R}}L(u)du\cdot G_{\eta}\left[X^{-}\right]\right]+o_{p}(1)\text{ (}n\rightarrow\infty\text{)}
\]
\[
=\left[MN\left[0,\sigma^{2}_{u}\varrho(0)\int_{\mathbb{R}}K^{2}(u)du\right],\text{ }\varrho(0)\int_{\mathbb{R}}K^{2}(u)du\text{, }\varrho(0)\int_{\mathbb{R}}L(u)du\right]+o_{a.s.}(1)\text{ (}\eta\rightarrow0\text{),}
\]
as required. $\blacksquare$

\bigskip{}

\noindent\textbf{Proof of Theorem \ref{HT-MCLT-LL}}. We begin with
the proof of part (a). Let 
\[
S\sim S_{\alpha_{u}}\left\{ \sigma^{\alpha_{u}}_{u}\int_{\mathbb{R}}\left\vert f(x)\right\vert ^{\alpha_{u}}\varrho(x)dx,\beta^{*}_{u}\frac{\int_{\mathbb{R}}\text{sgn[}f(x)\text{]}\left\vert f(x)\right\vert ^{\alpha_{u}}\varrho(x)dx}{\int_{\mathbb{R}}\left\vert f(x)\right\vert ^{\alpha_{u}}\varrho(x)dx}\right\} .
\]
Noting that $\mathbf{E}\exp\left[i\lambda S\right]\neq0$ ($i=\sqrt{-1}$)
for all $\lambda\in\mathbb{R}$, the result will follow from Jakubowski
(1986), Lemma 1.2 if we show that for each $\lambda$ 
\[
\Psi_{n}(\lambda):=\prod^{n}_{t=1}\left\{ \mathbf{E}\exp\left[i\lambda\frac{1}{n^{1/\alpha_{u}}}f\left(\frac{x_{t-1}}{\gamma_{n}}\right)u_{t}\right]\mid\mathcal{F}_{t-1}\right\} \overset{p}{\rightarrow}\mathbf{E}\exp\left[i\lambda S\right].
\]
Using the same arguments as those used in the proof of Thm 4.1 in
DK2021 we have as $n\rightarrow\infty$
\begin{equation}
\max_{1\leq t\leq n}\left\vert g_{t,n}\right\vert :=\max_{1\leq t\leq n}\left\vert \lambda n^{-1/\alpha_{u}}f\left(\frac{x_{t-1}}{\gamma_{n}}\right)\right\vert =o_{p}(1).\label{proofHT-MCLT-LL-1}
\end{equation}
By \textbf{Assumption HTE}, we have
\[
\ln\mathbf{E(}e^{i\lambda u_{t}}\mid\mathcal{F}_{t-1})=-\sigma^{\alpha_{u}}_{u}\left\vert \lambda\right\vert ^{\alpha_{u}}H_{u}\left(\left\vert \lambda\right\vert ^{-1}\right)\left[1+i\beta^{*}_{u}\text{sgn}(\lambda)\tan\left(\frac{\pi\alpha_{u}}{2}\right)\right]\left[1+R_{t}\left(\lambda\right)\right].
\]
In view of condition (ii) of \textbf{Assumption HTE} and (\ref{proofHT-MCLT-LL-1})
and where the remainder term satisfies
\begin{equation}
\mathrm{sup}_{1\leq t\leq n}R_{t}\left(g_{t,n}\right)=o_{p}(1),\label{proofHT-MCLT-LL-2}
\end{equation}
 as $n\rightarrow\infty$. To see this fix $\varepsilon$,$\delta>0$
and write
\[
P\left(\sup_{1\leq t\leq n}\left|R_{t}\left(g_{t,n}\right)\right|>\varepsilon\right)\leq P\left(\sup_{1\leq t\leq n}\left|g_{t,n}\right|>\delta\right)
\]
\[
+P\left(\sup_{1\leq t\leq n}\left|R_{t}\left(g_{t,n}\right)\right|>\varepsilon,\sup_{1\leq t\leq n}\left|g_{t,n}\right|\leq\delta\right)
\]
\[
\leq P\left(\sup_{1\leq t\leq n}\left|g_{t,n}\right|>\delta\right)+P\left(\sup_{\left|\lambda\right|\leq\delta}\sup_{1\leq t\leq n}\left|R_{t}\left(\lambda\right)\right|>\varepsilon\right).
\]
Then (\ref{proofHT-MCLT-LL-2}) follows from the fact that
\[
\lim_{n\rightarrow\infty}P\left(\sup_{1\leq t\leq n}\left|g_{t,n}\right|>\delta\right)=0
\]
and 
\[
\lim_{\delta\rightarrow0}\limsup_{n\rightarrow\infty}P\left(\sup_{\left|\lambda\right|\leq\delta}\sup_{1\leq t\leq n}\left|R_{t}\left(\lambda\right)\right|>\varepsilon\right)=0.
\]
In view of the above, as $n\rightarrow\infty$
\begin{eqnarray*}
\ln\Psi_{n}(\lambda) & = & \sum^{n}_{t=1}\left\{ -\sigma^{\alpha_{u}}_{u}\left\vert \lambda n^{-1/\alpha_{u}}f\left(\frac{x_{t-1}}{\gamma_{n}}\right)\right\vert ^{\alpha_{u}}H_{u}\left(\left\vert \lambda n^{-1/\alpha_{u}}f\left(\frac{x_{t-1}}{\gamma_{n}}\right)\right\vert ^{-1}\right)\right.\\
 &  & \left.\cdot\left[1+i\beta^{*}_{u}\text{sgn}\left[\lambda n^{-1/\alpha_{u}}f\left(\frac{x_{t-1}}{\gamma_{n}}\right)\right]\tan\left(\frac{\pi\alpha_{u}}{2}\right)\right]\right\} \left[1+R_{t}\left(g_{t,n}\right)\right]
\end{eqnarray*}
\[
=_{(1)}\sum^{n}_{t=1}\left\{ -\sigma^{\alpha_{u}}_{u}\left\vert \lambda n^{-1/\alpha_{u}}f\left(\frac{x_{t-1}}{\gamma_{n}}\right)\right\vert ^{\alpha_{u}}H_{u}\left(\left\vert g_{t,n}\right\vert ^{-1}\right)\right.
\]
\[
\left.\left[1+i\beta^{*}_{u}\text{sgn}\left[\lambda n^{-1/\alpha_{u}}f\left(\frac{x_{t-1}}{\gamma_{n}}\right)\right]\tan\left(\frac{\pi\alpha_{u}}{2}\right)\right]\left[1+o_{p}(1)\right]\right\} 
\]
\[
=_{(2)}\sum^{n}_{t=1}-\sigma^{\alpha_{u}}_{u}\left\vert \lambda n^{-1/\alpha_{u}}f\left(\frac{x_{t-1}}{\gamma_{n}}\right)\right\vert ^{\alpha_{u}}\left[1+i\beta^{*}_{u}\text{sgn}\left[\lambda n^{-1/\alpha_{u}}f\left(\frac{x_{t-1}}{\gamma_{n}}\right)\right]\tan\left(\frac{\pi\alpha_{u}}{2}\right)\right]\left[1+o_{p}(1)\right]
\]
\[
=_{(3)}-\sigma^{\alpha_{u}}_{u}\int_{\mathbb{R}}\left\vert \lambda f\left(x\right)\right\vert ^{\alpha_{u}}\left[1+i\beta^{*}_{u}\text{sgn}\left[\lambda f\left(x\right)\right]\tan\left(\frac{\pi\alpha_{u}}{2}\right)\right]\varrho(x)dx\left[1+o_{p}(1)\right]
\]
\[
=-\sigma^{\alpha_{u}}_{u}\int_{\mathbb{R}}\left\vert f\left(x\right)\right\vert ^{\alpha_{u}}\varrho(x)dx\left\vert \lambda\right\vert ^{\alpha_{u}}\left[1+\frac{i\beta^{*}_{u}\text{sgn}(\lambda)\int_{\mathbb{R}}\left\vert f\left(x\right)\right\vert ^{\alpha_{u}}\text{sgn}\left[f\left(x\right)\right]\varrho(x)dx\tan\left(\frac{\pi\alpha_{u}}{2}\right)}{\int_{\mathbb{R}}\left\vert f\left(x\right)\right\vert ^{\alpha_{u}}\varrho(x)dx}\right]+o_{p}(1),
\]
where $=_{(1)}$ follows from (\ref{proofHT-MCLT-LL-2}), $=_{(2)}$
is due to condition (iii) of \textbf{Assumption HTE} and (\ref{proofHT-MCLT-LL-1}),
and $=_{(3)}$ is a direct consequence of \textbf{Theorem \ref{MainLL}}.

\noindent The proof of part (b) follows from \textbf{Theorem \ref{MainLL_kernel}}
and proceeds by arguments similar to those used above. $\blacksquare$

\section{Proof of Main Lemmas; Lemma \ref{lem:FSalphaM} \& Weak Convergence}

\textbf{Proof of Lemma \ref{WCon} (weak convergence)}. We will only
prove the result for the case $\alpha\neq1$. The proof for $\alpha=1$
uses arguments similar to those used below. Throughout the proof $\text{\ensuremath{\beta^{\#}}:=}\beta^{*}\tan(\alpha\pi/2)$.

\bigskip{}

\textbf{BFR1}. Set 
\[
G^{+}_{nk}:=-\sigma^{\alpha}\left\vert \frac{\lambda}{\gamma_{n}}a_{k}\right\vert ^{\alpha}H\left(\left\vert \frac{\lambda}{\gamma_{n}}a_{k}\right\vert ^{-1}\right)\left[1+i\beta^{\#}\text{sgn}\left(\frac{\lambda}{\gamma_{n}}a_{k}\right)\right],
\]
and
\[
G^{-}_{nk}:=-\sigma^{\alpha}\left\vert \frac{\lambda}{\gamma_{n}}a^{-}_{k,\lfloor nr\rfloor}\right\vert ^{\alpha}H\left(\left\vert \frac{\lambda}{\gamma_{n}}a^{-}_{k,\lfloor nr\rfloor}\right\vert ^{-1}\right)\left[1+i\beta^{\#}\text{sgn}\left(\frac{\lambda}{\gamma_{n}}a^{-}_{k,\lfloor nr\rfloor}\right)\right].
\]
Note that by (\ref{FMtypeI.3}) for $k$ large enough ($k\geq k_{0}$,
say)
\[
\frac{1}{2}\alpha^{\alpha}k^{-1/\alpha}<a_{k}<\frac{3}{2}\alpha^{\alpha}k^{-1/\alpha}.
\]
Hence, for $k$ large and all $n\geq1$
\[
\text{sgn}\left(\frac{\lambda}{\gamma_{n}}a_{k}\right)=\text{sgn}\left(\frac{\lambda}{\gamma_{n}}\right)=\text{sgn}\left(\lambda\right).
\]
In view of the above,
\[
\sum^{\lfloor\left(n-1\right)r\rfloor}_{k=0}G^{+}_{nk}=_{(1)}-\sigma^{\alpha}\left[1+i\beta^{\#}\text{sgn}\left(\lambda\right)\right]\sum^{\lfloor\left(n-1\right)r\rfloor}_{k=k_{0}}\left\vert \frac{\lambda}{\gamma_{n}}a_{k}\right\vert ^{\alpha}H\left(\left\vert \frac{\lambda}{\gamma_{n}}a_{k}\right\vert ^{-1}\right)+o(1),
\]
\begin{equation}
=_{(2)}-\alpha^{\alpha+1}\sigma^{\alpha}\left\vert \lambda\right\vert ^{\alpha}\left[1+i\beta^{\#}\text{sgn}\left(\lambda\right)\right]+o(1),\label{CF1-FR1}
\end{equation}
where the approximation in $=_{(1)}$ follows from the fact that $\mathrm{sup}_{1\leq k\leq k_{0}}\left|\frac{\lambda}{\gamma_{n}}a_{k}\right|\rightarrow0$
as $n\rightarrow\infty$, whilst $=_{(2)}$ follows from arguments
similar to those used for deriving (\ref{FMtypeI.3_0b}). Note that
the limit term above is independent of $r$\ due to the fact $L(x)$
is SV i.e. 
\begin{equation}
\gamma^{\alpha}_{n}\sim L(\gamma_{n}n^{1/\alpha})\sim L(\gamma_{n}\left[\left(n-1\right)r\right]^{1/\alpha}).\label{CF2-FR1}
\end{equation}
Furthermore, identical arguments as those used above give
\begin{equation}
\sum^{\lfloor\left(n-1\right)r\rfloor}_{k=0}\left\vert G^{+}_{nk}\right\vert =O(1).\label{CF3-FR1}
\end{equation}
Next, we show that for each $0<r\leq1$, 
\begin{equation}
\sum^{\infty}_{k=\lfloor\left(n-1\right)r\rfloor+1}G^{-}_{nk}=-\sigma^{\alpha}\alpha^{\alpha+1}\left\vert \lambda\right\vert ^{\alpha}\left[1+i\beta^{\#}\text{sgn}\left(\lambda\right)\right]+o(1).\label{CF4-FR1}
\end{equation}
Without loss of generality set $r=1$. As before the limit does not
depend on $r$ due to eq. (\ref{CF2-FR1}) above. Fix $t_{0}\geq1$.
Then as $n\rightarrow\infty$ we have
\[
\sum^{\infty}_{k=n}G^{-}_{nk}=\sum^{\infty}_{k=n+t_{0}}G^{-}_{nk}+o(1)
\]
where the approximation above follows from the fact that
\[
\mathrm{sup}_{n+t_{0}\leq k\leq n}\left|\frac{\lambda}{\gamma_{n}}a^{-}_{k,\lfloor nr\rfloor}\right|\rightarrow0,
\]
as $n\rightarrow\infty$. Next, using the same decomposition as the
one in eq. (\ref{FMtypeI.3d}) we get for $t_{0}$ large enough and
as $n\rightarrow\infty$
\[
\sum^{\infty}_{k=n+t_{0}}G^{-}_{nk}=\sum^{2\left(n+t_{0}\right)}_{k=n+t_{0}}G^{-}_{nk}+\sum^{\infty}_{k=2\left(n+t_{0}\right)+1}G^{-}_{nk}
\]
\[
=\sum^{n+2t_{0}}_{k=t_{0}}-\sigma^{\alpha}\left\vert \frac{\lambda}{\gamma_{n}}\sum^{k+n}_{j=k}c_{j}\right\vert ^{\alpha}H\left(\left\vert \frac{\lambda}{\gamma_{n}}\sum^{k+n}_{j=k}c_{j}\right\vert ^{-1}\right)\left[1+i\beta^{\#}\text{sgn}\left(\frac{\lambda}{\gamma_{n}}\sum^{k+n}_{j=k}c_{j}\right)\right]
\]
\[
+\sum^{\infty}_{k=2\left(n+t_{0}\right)+1}-\sigma^{\alpha}\left\vert \frac{\lambda}{\gamma_{n}}\sum^{k}_{j=k-n}c_{j}\right\vert ^{\alpha}H\left(\left\vert \frac{\lambda}{\gamma_{n}}\sum^{k}_{j=k-n}c_{j}\right\vert ^{-1}\right)\left[1+i\beta^{\#}\text{sgn}\left(\frac{\lambda}{\gamma_{n}}\sum^{k}_{j=k-n}c_{j}\right)\right]
\]
\[
=_{(3)}-\sigma^{\alpha}\left[1+i\beta^{\#}\text{sgn}\left(\lambda\right)\right]\sum^{n+2t_{0}}_{k=t_{0}}\left\vert \frac{\lambda}{\gamma_{n}}\sum^{k+n}_{j=k}c_{j}\right\vert ^{\alpha}H\left(\left\vert \frac{\lambda}{\gamma_{n}}\sum^{k+n}_{j=k}c_{j}\right\vert ^{-1}\right)
\]
\[
-\sigma^{\alpha}\left[1+i\beta^{\#}\text{sgn}\left(\lambda\right)\right]\sum^{\infty}_{k=2\left(n+t_{0}\right)+1}\left\vert \frac{\lambda}{\gamma_{n}}\sum^{k}_{j=k-n}c_{j}\right\vert ^{\alpha}H\left(\left\vert \frac{\lambda}{\gamma_{n}}\sum^{k}_{j=k-n}c_{j}\right\vert ^{-1}\right)
\]
\[
=:-\sigma^{\alpha}\left[1+i\beta^{\#}\text{sgn}\left(\lambda\right)\right]\left(D^{\prime}_{n}+D^{\prime\prime}_{n}\right)
\]
\begin{equation}
=_{(4)}-\sigma^{\alpha}\alpha^{\alpha+1}\left\vert \lambda\right\vert ^{\alpha}\left[1+i\beta^{\#}\text{sgn}\left(\lambda\right)\right]+o(1),\label{CF5-FR1}
\end{equation}
where $=_{(3)}$ follows from \textbf{Lemma \ref{PropFR1}},\footnote{Let $I_{n}(k)$ and $Q_{n}(k)$ as in \textbf{Lemma \ref{PropFR1}}.
By the aforementioned lemma for $t_{0}$ large enough we get $0<I_{n}(k)/2<\sum^{k}_{j=k-n}c_{j}<3I_{n}(k)/2$
for all $k\geq2\left(n+t_{0}\right)+1$, and $0<Q_{n}(k)/2<\sum^{k+n}_{j=k}c_{j},<3Q_{n}(k)/2$
for all $n+2t_{0}\geq k\geq t_{0}$ and $n\geq t_{0}$. 

{}} $D^{\prime}_{n}$, $D^{\prime\prime}_{n}$ resemble\footnote{Set $R_{\alpha}(x)=\left\vert x\right\vert ^{\alpha}H(\left\vert x\right\vert ^{-1})$.
Then 
\[
D^{\prime}_{n}=\sum^{n+2t_{0}}_{k=t_{0}}R_{\alpha}\left(\frac{\lambda}{\gamma_{n}}\sum^{k+n}_{j=k}c_{j}\right)\text{, whilst }C^{\prime}_{3,n}=\sum^{n+2t_{0}}_{k=t_{0}}R_{\alpha}\left(\frac{1}{\gamma_{n}}\sum^{k+n}_{j=k}c_{j}\right)
\]
Moreover, 
\[
D^{\prime\prime}_{n}=\sum^{\infty}_{k=2\left(n+t_{0}\right)+1}R_{\alpha}\left(\frac{\lambda}{\gamma_{n}}\sum^{k}_{j=k-n}c_{j}\right)\text{, whilst }C^{\prime\prime}_{3,n}=\sum^{\infty}_{k=2\left(n+t_{0}\right)+1}R_{\alpha}\left(\frac{1}{\gamma_{n}}\sum^{k}_{j=k-n}c_{j}\right).
\]
} the terms $C^{\prime}_{3,n},C^{\prime\prime}_{3,n}$ of (\ref{FMtypeI.3d}),
and $=_{(4)}$ follows from an almost identical argument as the one
used for establishing (\ref{FMtypeI.6d2}). Therefore, (\ref{CF4-FR1})
holds true. Moreover, using arguments similar to those used above
we get 
\begin{equation}
\sum^{\infty}_{k=\lfloor\left(n-1\right)r\rfloor+1}\left\vert G^{-}_{nk}\right\vert =O(1).\label{CF6-FR1}
\end{equation}

Next, consider
\[
\mathbf{E}\left(e^{i\lambda x_{\lfloor nr\rfloor}/\gamma_{n}}\right)=\mathbf{E}\left(e^{i\lambda x^{+}_{\lfloor nr\rfloor}/\gamma_{n}}\right)\cdot\mathbf{E}\left(e^{i\lambda x^{-}_{\lfloor nr\rfloor}/\gamma_{n}}\right).
\]
We shall show that for each $\lambda$ as $n\rightarrow\infty$
\[
\mathbf{E}\left(e^{i\lambda x^{+}_{\lfloor nr\rfloor}/\gamma_{n}}\right)\rightarrow\mathbf{E}\left(e^{i\lambda X^{+}}\right),
\]
and
\[
\mathbf{E}\left(e^{i\lambda x^{-}_{\lfloor nr\rfloor}/\gamma_{n}}\right)\rightarrow\mathbf{E}\left(e^{i\lambda X^{+}}\right).
\]
 We start with the former. It can be easily checked that
\[
\sup_{1\leq k\leq n}\left\vert g^{+}_{nk}\right\vert :=\sup_{1\leq k\leq n}\left\vert \frac{\lambda}{\gamma_{n}}a_{k}\right\vert =o(1).
\]
Next, by virtue of (\ref{INN}) there is some $R(\lambda)$ such that
$R(\lambda)\rightarrow0$ as $\lambda\rightarrow0$, and
\[
\mathbf{E}\left(e^{i\lambda x^{+}_{\lfloor nr\rfloor}/\gamma_{n}}\right)=\prod^{\lfloor\left(n-1\right)r\rfloor}_{k=1}\left\vert \psi_{\xi}\left(\frac{\lambda}{\gamma_{n}}a_{k}\right)\right\vert 
\]
\[
=\prod^{\lfloor\left(n-1\right)r\rfloor}_{k=1}\exp\left\{ G^{+}_{nk}\left[1+R(g^{+}_{nk})\right]\right\} =\exp\left\{ \sum^{\lfloor\left(n-1\right)r\rfloor}_{k=1}G^{+}_{nk}\left[1+R(g^{+}_{nk})\right]\right\} .
\]
Set $R^{+}_{n}:=\sup_{1\leq k\leq n}\left\vert R(g^{+}_{nk})\right\vert .$
It follows from the above that $R^{+}_{n}\rightarrow0$, as $n\rightarrow\infty$.
In view of the latter, (\ref{CF1-FR1}) and (\ref{CF3-FR1}) we have
\[
\mathbf{E}\left(e^{i\lambda x^{+}_{\lfloor nr\rfloor}/\gamma_{n}}\right)=\exp\left\{ \sum^{\lfloor\left(n-1\right)r\rfloor}_{k=1}G^{+}_{nk}\left[1+R(g^{+}_{nk})\right]\right\} =\exp\left\{ \sum^{\lfloor\left(n-1\right)r\rfloor}_{k=1}G^{+}_{nk}\right\} +o(1)
\]
\[
=\exp\left\{ -\sigma^{\alpha}\alpha^{\alpha+1}\left[1+i\beta^{\#}\text{sgn}\left(\lambda\right)\right]\right\} +o(1)=\mathbf{E}\left(e^{i\lambda X^{+}}\right)+o(1),
\]
as required. 

Next, it can be easily seen that
\[
\sup_{1\leq k\leq n,0\leq r\leq1}\left\vert g^{-}_{nkr}\right\vert :=\sup_{1\leq k\leq n,0\leq r\leq1}\left\vert \frac{\lambda}{\gamma_{n}}a^{-}_{k,\lfloor nr\rfloor}\right\vert =o(1).
\]
Using identical arguments to those used above together with (\ref{CF5-FR1})
and (\ref{CF6-FR1}) as $n\rightarrow\infty$ we have 
\[
\mathbf{E}\left(e^{i\lambda x^{-}_{\lfloor nr\rfloor}/\gamma_{n}}\right)=\prod^{\infty}_{k=\lfloor\left(n-1\right)r\rfloor+1}\psi_{\xi}\left(\frac{\lambda}{\gamma_{n}}a^{-}_{k,\lfloor nr\rfloor}\right)
\]
\[
=\exp\left\{ \sum^{\infty}_{k=\lfloor\left(n-1\right)r\rfloor+1}G^{-}_{nk}\left[1+R(g^{-}_{nkr})\right]\right\} 
\]
\[
=\exp\left\{ -\sigma^{\alpha}\alpha^{\alpha+1}\left[1+i\beta^{\#}\text{sgn}\left(\lambda\right)\right]\right\} +o(1)=\mathbf{E}\left(e^{i\lambda X^{+}}\right)+o(1),
\]
as required. 

\bigskip{}

\textbf{BFR2}. In view of \textbf{Lemma \ref{LPA}}, it suffices to
show that
\[
\mathbf{E}\left(e^{i\lambda\tilde{x}_{\lfloor nr\rfloor}/\gamma_{n}}\right)=\exp\left\{ -\sigma^{\alpha}\alpha\left\vert \lambda c_{\infty}\right\vert ^{\alpha}\left[1+i\beta^{\#}\text{sgn}\left(\lambda c_{\infty}\right)\right]\right\} +o(1),
\]
for all $0<r\leq1$, where $c_{\infty}:=\sum^{\infty}_{j=0}c_{j}$.
Set 
\[
G_{nk}:=-\sigma^{\alpha}\left\vert \frac{\lambda c_{\infty}}{\gamma_{n}}k^{-\frac{1}{\alpha}}\right\vert ^{\alpha}H\left(\left\vert \frac{\lambda c_{\infty}}{\gamma_{n}}k^{-\frac{1}{\alpha}}\right\vert ^{-1}\right)\left[1+i\beta^{\#}\text{sgn}\left(\frac{\lambda c_{\infty}}{\gamma_{n}}k^{-\frac{1}{\alpha}}\right)\right],
\]
and
\[
g_{nk}:=\frac{\lambda c_{\infty}}{\gamma_{n}}k^{-\frac{1}{\alpha}}.
\]
In view of eq. (\ref{INN}) the characteristic function of $\gamma^{-1}_{n}\tilde{x}_{\lfloor nr\rfloor}$
is 
\[
\mathbf{E}\left(e^{i\lambda\tilde{x}_{\lfloor nr\rfloor}/\gamma_{n}}\right)=\prod^{\lfloor\left(n-1\right)r\rfloor}_{k=1}\left\vert \psi_{\xi}\left(\frac{\lambda c_{\infty}}{\gamma_{n}}k^{-\frac{1}{\alpha}}\right)\right\vert 
\]
\begin{equation}
=\prod^{\lfloor\left(n-1\right)r\rfloor}_{k=1}\exp\left\{ G_{nk}\left[1+R(g_{nk})\right]\right\} =\exp\left\{ \sum^{\lfloor\left(n-1\right)r\rfloor}_{k=1}G_{nk}\left[1+R(g_{nk})\right]\right\} ,\label{CF1-FR2}
\end{equation}
with $R(g_{nk})\rightarrow0$ as $g_{nk}\rightarrow0$. In particular,
we have $R_{n}:=\sup_{1\leq k\leq n}\left\vert R(g_{nk})\right\vert \rightarrow0$
as $n\rightarrow\infty$, due to the fact $\sup_{1\leq k\leq n}\left\vert g_{nk}\right\vert \rightarrow0$
as $n\rightarrow\infty$. Therefore, as $n\rightarrow\infty$ we have
\[
\left\vert \sum^{\lfloor\left(n-1\right)r\rfloor}_{k=1}G_{nk}R(g_{nk})\right\vert \leq R_{n}\sum^{\lfloor\left(n-1\right)r\rfloor}_{k=1}\left\vert G_{nk}\right\vert =o(1),
\]
if we show that $\sum^{n}_{k=1}\left\vert G_{nk}\right\vert $ is
bounded. It follows from the above that as $n\rightarrow\infty$
\begin{equation}
\mathbf{E}\left(e^{i\lambda\tilde{x}_{\lfloor nr\rfloor}/\gamma_{n}}\right)=\exp\left\{ \sum^{n}_{k=1}G_{nk}\right\} +o(1).\label{CF2-FR2}
\end{equation}
Next, we show that 
\begin{equation}
\sum^{\lfloor\left(n-1\right)r\rfloor}_{k=1}G_{nk}=-\sigma^{\alpha}\alpha\left\vert \lambda c_{\infty}\right\vert ^{\alpha}\left[1+i\beta^{\#}\text{sgn}\left(\lambda c_{\infty}\right)\right]+o(1).\label{CF3-FR2}
\end{equation}
 An almost identical argument to the one provided below shows that
$\sum^{n}_{k=1}\left\vert G_{nk}\right\vert $ is bounded. In view
of (\ref{CF1-FR2})-(\ref{CF3-FR2}), we have
\begin{equation}
\mathbf{E}\left(e^{i\lambda\tilde{x}_{\lfloor nr\rfloor}/\gamma_{n}}\right)=\exp\left\{ -\sigma^{\alpha}\alpha\left\vert \lambda c_{\infty}\right\vert ^{\alpha}\left[1+i\beta^{\#}\text{sgn}\left(\lambda c_{\infty}\right)\right]\right\} +o(1),\label{CF4-FR2}
\end{equation}
which is the requisite result. We finish with the proof of (\ref{CF3-FR2}).
Note that for all $n$,
\begin{equation}
\sum^{\lfloor\left(n-1\right)r\rfloor}_{k=1}G_{nk}=-\sigma^{\alpha}\left[1+i\beta^{\#}\text{sgn}\left(\lambda c_{\infty}\right)\right]\sum^{\lfloor\left(n-1\right)r\rfloor}_{k=1}\left\vert \frac{\lambda c_{\infty}}{\gamma_{n}}k^{-\frac{1}{\alpha}}\right\vert ^{\alpha}H\left(\left\vert \frac{\lambda c_{\infty}}{\gamma_{n}}k^{-\frac{1}{\alpha}}\right\vert ^{-1}\right).\label{CF4b-FR2}
\end{equation}
Next, as $n\rightarrow\infty$, the sum 
\[
\sum^{\lfloor\left(n-1\right)r\rfloor}_{k=1}\left\vert \frac{\lambda c_{\infty}}{\gamma_{n}}k^{-\frac{1}{\alpha}}\right\vert ^{\alpha}H\left(\left\vert \frac{\lambda c_{\infty}}{\gamma_{n}}k^{-\frac{1}{\alpha}}\right\vert ^{-1}\right)=\left\vert \frac{\lambda c_{\infty}}{\gamma_{n}}\right\vert ^{\alpha}\sum^{\lfloor\left(n-1\right)r\rfloor}_{k=1}k^{-1}H\left(\left\vert c_{\infty}\lambda\right\vert ^{-1}\gamma_{n}k^{\frac{1}{\alpha}}\right).
\]
\begin{equation}
=_{\left(1\right)}\alpha\left\vert \frac{\lambda c_{\infty}}{\gamma_{n}}\right\vert ^{\alpha}\int^{\left\vert c_{\infty}\lambda\right\vert ^{-1}\gamma_{n}\left(\lfloor\left(n-1\right)r\rfloor\right)^{\frac{1}{\alpha}}}_{u=\left\vert c_{\infty}\lambda\right\vert ^{-1}\gamma_{n}}\frac{H(u)}{u}du\left[1+o(1)\right]=_{\left(2\right)}\alpha\left\vert \lambda c_{\infty}\right\vert ^{\alpha}+o(1),\label{CF5-FR2}
\end{equation}
where $=_{\left(1\right)}$ is due to\textbf{ Lemma \ref{LemmaSV1}(i)},
$=_{\left(2\right)}$ due to the fact that 
\[
L\left(\gamma_{n}n^{\frac{1}{\alpha}}\right)\sim\gamma^{\alpha}_{n},
\]
and
\[
\int^{\left\vert c_{\infty}\lambda\right\vert ^{-1}\gamma_{n}\left(\lfloor\left(n-1\right)r\rfloor\right)^{\frac{1}{\alpha}}}_{u=\left\vert c_{\infty}\lambda\right\vert ^{-1}\gamma_{n}}H(u)u^{-1}du\sim\int^{\gamma_{n}n^{\frac{1}{\alpha}}}_{u=\gamma_{n}}\frac{H(u)}{u}du\sim L\left(\gamma_{n}n^{\frac{1}{\alpha}}\right).
\]
To see validity of the latter statement write
\[
\frac{L\left(\gamma_{n}n^{\frac{1}{\alpha}}\right)-L\left(\gamma_{n}\right)}{L\left(\gamma_{n}n^{\frac{1}{\alpha}}\right)}=1-\frac{L\left(\gamma_{n}\right)}{L\left(\gamma_{n}n^{\frac{1}{\alpha}}\right)}
\]
\[
\sim1-\frac{L\left(\gamma_{n}\right)}{\gamma^{\alpha}_{n}}\rightarrow1.
\]
Note that (\ref{CF3-FR2}) follows from (\ref{CF4b-FR2}) and (\ref{CF5-FR2}).

\bigskip{}

\textbf{MI-SM}. From \textbf{Lemma FS$_{\alpha}$A} we have $\left\vert \gamma^{-1}_{n}x_{\lfloor nr\rfloor}(n)^{-}\right\vert =o_{p}(1)$
for all $0<r\leq1$. Therefore, it suffices to show that 
\[
\mathbf{E}\exp\left[i\lambda\gamma^{-1}_{n}x_{t}(n)^{+}\right]=\exp\left\{ -\sigma^{\alpha}\alpha^{-1}\left\vert \lambda c_{\infty}\right\vert ^{\alpha}\left[1+i\beta^{\#}\text{sgn}\left(\lambda c_{\infty}\right)\right]\right\} +o(1).
\]
Set
\[
G_{nk}:=-\sigma^{\alpha}\left\vert \frac{\lambda}{\gamma_{n}}a_{k}\right\vert ^{\alpha}H\left(\left\vert \frac{\lambda}{\gamma_{n}}a_{k}\right\vert ^{-1}\right)\left[1+i\beta^{\#}\text{sgn}\left(\frac{\lambda}{\gamma_{n}}a_{k}\right)\right],
\]
where $a_{k}(n):=a_{k}=\sum^{k}_{l=0}\rho^{k-l}_{n}c_{l}$ with $0<\rho^{j}_{n}<1$
for all $j\geq1$ and $n$ large enough. 

In view of the above, as $n\rightarrow\infty$, 
\[
\sum^{\lfloor\left(n-1\right)r\rfloor}_{k=1}G_{nk}=_{\left(1\right)}-\sigma^{\alpha}\left\vert \lambda\right\vert ^{\alpha}\kappa^{-1}_{n}\sum^{\lfloor\left(n-1\right)r\rfloor}_{k=1}\left\vert a_{k}\right\vert ^{\alpha}\left[1+i\beta^{\#}\text{sgn}\left(\frac{\lambda}{\gamma_{n}}a_{k}\right)\right]\left[1+o(1)\right]
\]
\[
=_{\left(2\right)}-\sigma^{\alpha}\left\vert \lambda\right\vert ^{\alpha}\kappa^{-1}_{n}\sum^{\lfloor\left(n-1\right)r\rfloor}_{k=1}\left\vert a_{k}\right\vert ^{\alpha}\left[1+i\beta^{\#}\text{sgn}\left(\lambda\rho^{k}_{n}c_{\infty}\right)\right]\left[1+o(1)\right]
\]
\[
=_{\left(3\right)}-\sigma^{\alpha}\left\vert \lambda\right\vert ^{\alpha}\kappa^{-1}_{n}\sum^{\lfloor\left(n-1\right)r\rfloor}_{k=1}\left\vert a_{k}\right\vert ^{\alpha}\left[1+i\beta^{\#}\text{sgn}\left(\lambda c_{\infty}\right)\right]\left[1+o(1)\right]
\]
\[
=_{(4)}-\sigma^{\alpha}\alpha^{-1}\left\vert \lambda c_{\infty}\right\vert ^{\alpha}\left[1+i\beta^{\#}\text{sgn}\left(\lambda c_{\infty}\right)\right]+o(1),
\]
where $=_{\left(1\right)}$ follows from an argument analogous to
the one used for the first approximation in (\ref{MI000a}); the equalities
$=_{\left(2\right)}$ and $=_{\left(4\right)}$ from \textbf{Lemma
\ref{LemmaMI-2}}, and $=_{\left(3\right)}$ from the fact that, for
$n$ sufficiently large, $\text{sgn}\left(\lambda c_{\infty}\right)=\text{sgn}\left(\lambda\rho^{k}_{n}c_{\infty}\right)$.
Furthermore, it is straightforward to verify that $\sum^{n}_{k=1}\left\vert G_{nk}\right\vert =O(1)$
and $\sup_{k\geq0}\left\vert a_{k}\right\vert <\infty$. Hence, in
view of \textbf{Assumption INN}, there is some $R(x)\rightarrow0$
as $x\rightarrow0$ such that
\[
\mathbf{E}\left(e^{i\lambda x^{+}_{\lfloor nr\rfloor}/\gamma_{n}}\right)=\prod^{\lfloor\left(n-1\right)r\rfloor}_{k=1}\psi_{\xi}\left(\frac{\lambda}{\gamma_{n}}a_{k}\right)
\]
\[
=\prod^{\lfloor\left(n-1\right)r\rfloor}_{k=1}\exp\left\{ G_{nk}\left[1+R\left(\frac{\lambda}{\gamma_{n}}a_{k}\right)\right]\right\} =\exp\left\{ \sum^{\lfloor\left(n-1\right)r\rfloor}_{k=1}G_{nk}\left[1+R\left(\frac{\lambda}{\gamma_{n}}a_{k}\right)\right]\right\} ,
\]
\[
=_{\left(4\right)}\exp\left\{ \sum^{\lfloor\left(n-1\right)r\rfloor}_{k=1}G_{nk}+o(1)\right\} =\exp\left\{ -\sigma^{\alpha}\alpha^{-1}\left\vert \lambda c_{\infty}\right\vert ^{\alpha}\left[1+i\beta^{\#}\text{sgn}\left(\lambda c_{\infty}\right)\right]\right\} +o(1),
\]
where $=_{\left(4\right)}$ follows from the fact that
\[
\sum^{n}_{k=1}\left\vert G_{nk}R\left(\frac{\lambda}{\gamma_{n}}a_{k}\right)\right\vert \leq\sup_{k}\left\vert R\left(\frac{\lambda}{\gamma_{n}}a_{k}\right)\right\vert \sum^{n}_{k=1}\left\vert G_{nk}\right\vert =o(1),
\]
and the proof is complete. \bigskip{}

\textbf{MI-LM} ($1/\alpha<m<1$). In view of (\ref{linprocsplit}),
we can write
\[
\gamma^{-1}_{n}x_{\left\lfloor nr\right\rfloor }=\gamma^{-1}_{n}x^{+}_{\left\lfloor nr\right\rfloor }+\gamma^{-1}_{n}x^{-}_{\left\lfloor nr\right\rfloor },
\]
where $\gamma^{-1}_{n}x^{-}_{\left\lfloor nr\right\rfloor }=o_{p}(1)$
due to \textbf{Lemma \ref{lem:FSalphaM}}. Therefore, it suffices
to show that for each $0<r\leq1$, as $n\rightarrow\infty$,
\[
\mathbf{E}\exp\left(\lambda\gamma^{-1}_{n}x^{+}_{\left\lfloor nr\right\rfloor }\right)\rightarrow\exp\left(-\left\vert \lambda\right\vert ^{\alpha}\sigma^{\alpha}\left|m_{0}\right|^{\alpha}\mathbb{V}_{1}\left[1+i\beta^{\#}\text{sgn}\left(\lambda\right)\right]\right).
\]

Write
\[
G_{nk}:=-\sigma^{\alpha}\left\vert \frac{\lambda}{\gamma_{n}}a_{k}\right\vert ^{\alpha}H\left(\left\vert \frac{\lambda}{\gamma_{n}}a_{k}\right\vert ^{-1}\right)\left[1+i\beta^{\#}\text{sgn}\left(\frac{\lambda}{\gamma_{n}}a_{k}\right)\right],
\]
where as in the previous part i.e. $a_{k}(n):=a_{k}=\sum^{k}_{l=0}\rho^{k-l}_{n}c_{l}$.\textbf{
}Set
\[
q(x):=e^{-x}\int^{x}_{0}e^{y}y^{-m}dy.
\]
By Lemma 7.2 of DK2018, as $n\rightarrow\infty$,
\[
\sup_{0\leq k\leq n}\left\vert \kappa^{-\left(1-m\right)}_{n}\sum^{k}_{j=0}\rho^{j}_{n}c_{k-j}-q(k/\kappa_{n})\right\vert =o(1).
\]
In fact, under the additional requirement $c_{j}=m_{0}\cdot j^{-m}[1+O(j^{-\nu})]$,
$0<\nu<1-m$, it is possible to strengthen this result as follows:
\begin{equation}
\sup_{0\leq k\leq n}\left\vert \kappa^{-\left(1-m\right)}_{n}\sum^{k}_{j=0}\rho^{j}_{n}c_{k-j}-q(k/\kappa_{n})\right\vert =O(\kappa^{-\nu}_{n}),\label{FSM-LM00-1}
\end{equation}
as $n\rightarrow\infty$. Let $m_{n}:=\left\lfloor \kappa^{p}_{n}\right\rfloor $,
for some $0<p<1$. Moreover, set $\delta_{n}:=\left\lfloor \kappa_{n}\varsigma_{n}\right\rfloor $,
$\varsigma_{n}=\left[\ln\left(1+\kappa_{n}\right)\right]^{2}$. By
(\ref{FSM-LM00-1}), for $n$ large enough,
\[
\inf_{m_{n}\leq k\leq\delta_{n}}\kappa^{-\left(1-m\right)}_{n}a_{k}(n)=\inf_{m_{n}\leq k\leq n}q\left(\frac{k}{\kappa_{n}}\right)+R_{n}
\]
\[
\geq_{(1)}q\left(\frac{m_{n}}{\kappa_{n}}\right)\wedge q\left(\frac{\delta_{n}}{\kappa_{n}}\right)+R_{n}\geq_{(2)}q\left(\frac{m_{n}}{\kappa_{n}}\right)+R_{n},
\]
where $R_{n}=O(\kappa^{-\nu}_{n})$, and $\geq_{(1)}$ follows from
the fact that $q(x)$ attains its minimum on $[m_{n},\delta_{n}]$
on a boundary point of that interval.\footnote{Note that the second derivative of $q(x)$ is $q^{\prime\prime}(x)=-mx^{-m-1}-q^{\prime}(x)$.
Hence, if $x_{*}$ is a stationary point of $q(x)$ interior in some
interval $[a,b]$, we have $q^{\prime\prime}(x_{*})=-mx^{-m-1}_{*}<0.$
Hence, $q(x_{*})$ corresponds to a maximum. In view of this, $q(x)$
attains its minimum on $[a,b]$ at a boundary point of the interval. } Moreover, $\geq_{(2)}$ holds from the fact that as $n\rightarrow\infty$,
\[
q\left(\frac{m_{n}}{\kappa_{n}}\right)\cdot q\left(\frac{\delta_{n}}{\kappa_{n}}\right)^{-1}\rightarrow0.
\]
To see this by by \textbf{Proposition \ref{SterlingApprox}(a,b)},
as $n\rightarrow\infty$,
\[
\frac{q\left(\frac{m_{n}}{\kappa_{n}}\right)}{q\left(\frac{\delta_{n}}{\kappa_{n}}\right)}\sim\frac{\frac{1}{1-m}\kappa^{-(1-m)(1-p)}_{n}}{\left(\frac{\delta_{n}}{\kappa_{n}}\right)^{-m}}\sim\frac{\frac{1}{1-m}\kappa^{-(1-m)(1-p)}_{n}}{\varsigma^{-m}_{n}}\rightarrow0.
\]

Next, by \textbf{Proposition \ref{SterlingApprox}(b)} we have
\[
q\left(\frac{m_{n}}{\kappa_{n}}\right)=\frac{1}{1-m}\left(\frac{m_{n}}{\kappa_{n}}\right)^{1-m}+R_{1n},\:R_{1n}=O\left(\left(\frac{m_{n}}{\kappa_{n}}\right)^{2-m}\right).
\]
 Choose $\nu$ such that $1-m>\nu>(1-p)(1-m)$. Then, for all $m_{n}\leq k\leq\delta_{n}$,
and $n$ large enough, we get
\[
\kappa^{-\left(1-m\right)}_{n}a_{k}(n)\geq q\left(\frac{m_{n}}{\kappa_{n}}\right)+R_{n}=\frac{1}{1-m}\left(\frac{m_{n}}{\kappa_{n}}\right)^{1-m}+R_{1n}+R_{n}
\]
\[
=\frac{1}{1-m}\kappa^{-(1-p)\left(1-m\right)}_{n}+O\left(\kappa^{-(1-p)\left(2-m\right)}_{n}\right)+O\left(\kappa^{-\nu}_{n}\right)
\]
\[
=\kappa^{-(1-p)\left(1-m\right)}_{n}\left[\frac{1}{1-m}+O\left(\kappa^{-1}_{n}\right)+O\left(\kappa^{(1-p)\left(1-m\right)-\nu}_{n}\right)\right]
\]
\[
\geq\kappa^{-(1-p)\left(1-m\right)}_{n}\left[\frac{1}{1-m}-\frac{1}{2\left(1-m\right)}\right]=\kappa^{-(1-p)\left(1-m\right)}_{n}\frac{1}{2\left(1-m\right)}>0.
\]
Hence, for all $m_{n}\leq k\leq\delta_{n}$, and $n$ large enough,
\begin{equation}
\textrm{sgn}\left(\lambda\right)=\textrm{sgn}\left(\lambda\frac{\kappa^{-(1-p)\left(1-m\right)}_{n}}{2\left(1-m\right)}\right)=\textrm{sgn}\left(\lambda\kappa^{-\left(1-m\right)}_{n}a_{k}(n)\right)=\text{sgn}\left(\frac{\lambda}{\gamma_{n}}a_{k}(n)\right).\label{FSM-LM00-2}
\end{equation}

In view of the above for each $0<r\leq1$, as $n\rightarrow\infty$,
\[
\sum^{\lfloor\left(n-1\right)r\rfloor}_{k=1}G_{nk}=_{\left(1\right)}-\sigma^{\alpha}\left\vert \lambda\right\vert ^{\alpha}H(m_{n})\gamma^{-\alpha}_{n}\sum^{\lfloor\left(n-1\right)r\rfloor}_{k=1}\left\vert a_{k}\right\vert ^{\alpha}\left[1+i\beta^{\#}\text{sgn}\left(\frac{\lambda}{\gamma_{n}}a_{k}\right)\right]+o(1)
\]
\[
=_{\left(2\right)}-\sigma^{\alpha}\left\vert \lambda\right\vert ^{\alpha}\kappa^{-b}_{n}\sum^{\lfloor\left(n-1\right)r\rfloor}_{k=1}\left\vert a_{k}\right\vert ^{\alpha}\left[1+i\beta^{\#}\text{sgn}\left(\frac{\lambda}{\gamma_{n}}a_{k}\right)\right]\left[1+o(1)\right]
\]
\[
=_{\left(3\right)}-\sigma^{\alpha}\left\vert \lambda\right\vert ^{\alpha}\kappa^{-b}_{n}\sum^{\delta_{n}}_{k=m_{n}}\left\vert a_{k}\right\vert ^{\alpha}\left[1+i\beta^{\#}\text{sgn}\left(\frac{\lambda}{\gamma_{n}}a_{k}\right)\right]\left[1+o(1)\right]
\]
\[
=_{\left(4\right)}-\sigma^{\alpha}\left\vert \lambda\right\vert ^{\alpha}\kappa^{-b}_{n}\sum^{\delta_{n}}_{k=m_{n}}\left\vert a_{k}\right\vert ^{\alpha}\left[1+i\beta^{\#}\text{sgn}\left(\lambda\right)\right]\left[1+o(1)\right]
\]
\begin{equation}
=_{(5)}-\left\vert \lambda\right\vert ^{\alpha}\sigma^{\alpha}\left|m_{0}\right|^{\alpha}\mathbb{V}_{1}\left[1+i\beta^{\#}\text{sgn}\left(\lambda\right)\right]+o(1),\label{FSM-LM00-3}
\end{equation}
where $=_{\left(1\right)}$ follows from analogous arguments as those
used for the proof of (\ref{FSM-LM01}), $=_{\left(2\right)}$ follows
from the fact that $H(m_{n})\gamma^{-\alpha}_{n}\sim\kappa^{-b}_{n}$,
$=_{(3)}$ is due to \textbf{Lemma \ref{LemmaMI1}} and \textbf{Lemma
\ref{LemmaMI3}}, $=_{(4)}$ is a direct consequence of (\ref{FSM-LM00-2}),
and $=_{(5)}$ follows from \textbf{Lemma \ref{LemmaMI3}}. Moreover,
arguments analogous to those used for the proof of (\ref{FSM-LM00-3})
show that 
\[
\sum^{n}_{k=1}\left\vert G_{nk}\right\vert =O(1).
\]

Next, note that $\gamma^{-1}_{n}\sup_{1\leq k\leq n}\left|a_{k}\right|\rightarrow0$,
as $n\rightarrow\infty$, cf. (\ref{FSM-LM0}). In view of \textbf{Assumption
INN}, there is some $R(x)\rightarrow0$ as $x\rightarrow0$ such that
the characteristic function of $x^{+}_{\lfloor nr\rfloor}/\gamma_{n}$
\[
\mathbf{E}\left(e^{i\lambda x^{+}_{\lfloor nr\rfloor}/\gamma_{n}}\right)=\prod^{\lfloor\left(n-1\right)r\rfloor}_{k=1}\psi_{\xi}\left(\frac{\lambda}{\gamma_{n}}a_{k}\right)
\]
\[
=\prod^{\lfloor\left(n-1\right)r\rfloor}_{k=1}\exp\left\{ G_{nk}\left[1+R\left(\frac{\lambda}{\gamma_{n}}a_{k}\right)\right]\right\} =\exp\left\{ \sum^{\lfloor\left(n-1\right)r\rfloor}_{k=1}G_{nk}\left[1+R\left(\frac{\lambda}{\gamma_{n}}a_{k}\right)\right]\right\} ,
\]
\[
=_{\left(6\right)}\exp\left\{ \sum^{\lfloor\left(n-1\right)r\rfloor}_{k=1}G_{nk}+o(1)\right\} =_{(7)}\exp\left\{ -\left\vert \lambda\right\vert ^{\alpha}\sigma^{\alpha}\left|m_{0}\right|^{\alpha}\mathbb{V}_{1}(\alpha,m)\left[1+i\beta^{\#}\text{sgn}\left(\lambda\right)\right]\right\} +o(1),
\]
where $=_{\left(6\right)}$ follows from the fact that as $n\rightarrow\infty$,
$\sup_{1\leq k\leq n}\left\vert R\left(\frac{\lambda}{\gamma_{n}}a_{k}\right)\right\vert =o(1)$,
and
\[
\sum^{n}_{k=1}\left\vert G_{nk}R\left(\frac{\lambda}{\gamma_{n}}a_{k}\right)\right\vert \leq\sup_{k}\left\vert R\left(\frac{\lambda}{\gamma_{n}}a_{k}\right)\right\vert \sum^{n}_{k=1}\left\vert G_{nk}\right\vert =o(1).
\]
Finally, $=_{(7)}$ is a direct consequence of (\ref{FSM-LM00-3}).\bigskip{}

\textbf{MI-NM} ($1<m<1+1/\alpha$). In view of (\ref{linprocsplit}),
we can write
\[
\gamma^{-1}_{n}x_{\left\lfloor nr\right\rfloor }=\gamma^{-1}_{n}x^{+}_{\left\lfloor nr\right\rfloor }+\gamma^{-1}_{n}x^{-}_{\left\lfloor nr\right\rfloor }
\]
where $\gamma^{-1}_{n}x^{-}_{\left\lfloor nr\right\rfloor }=0$, by
\textbf{Assumption MI-NM}. It suffices to show that for each $0<r\leq1$,
as $n\rightarrow\infty$,
\[
\mathbf{E}\exp\left(\lambda\gamma^{-1}_{n}x^{+}_{\left\lfloor nr\right\rfloor }\right)\rightarrow\exp\left(-\left\vert \lambda\right\vert ^{\alpha}\sigma^{\alpha}\left|m_{0}\right|^{\alpha}\mathbb{V}_{2}\left[1+i\beta^{\#}\mathbb{V}_{3}\text{sgn}\left(m_{0}\right)\text{sgn}\left(\lambda\right)\right]\right).
\]

Write
\[
G_{nk}:=-\sigma^{\alpha}\left\vert \frac{\lambda}{\gamma_{n}}a_{k}\right\vert ^{\alpha}H\left(\left\vert \frac{\lambda}{\gamma_{n}}a_{k}\right\vert ^{-1}\right)\left[1+i\beta^{\#}\text{sgn}\left(\frac{\lambda}{\gamma_{n}}a_{k}\right)\right],
\]
where recall that $a_{k}(n):=a_{k}=\sum^{k}_{l=0}\rho^{k-l}_{n}c_{l}$.
Moreover,
\[
p(x):=x^{1-m}-e^{-x}\int^{x}_{y=0}e^{y}y^{1-m}dy,
\]
and by (\ref{Ts_rmk-1}), as $n\rightarrow\infty$,
\begin{equation}
\frac{1}{\kappa_{n}}\sum^{\delta_{n}}_{k=1}\left\vert \frac{1}{\kappa^{1-m}_{n}}a_{k}-\frac{m_{0}p(k/\kappa_{n})}{1-m}\right\vert ^{\alpha}=o(1).\label{MN-NM-01}
\end{equation}

Let $\delta_{n}:=\left\lfloor \kappa_{n}\varsigma_{n}\right\rfloor $,
$\varsigma_{n}=\left[\ln\left(1+\kappa_{n}\right)\right]^{2}$. For
each $0<r\leq1$, as $n\rightarrow\infty$, 
\[
\sum^{\lfloor\left(n-1\right)r\rfloor}_{k=1}G_{nk}=_{(1)}-\sigma^{\alpha}\left\vert \lambda\right\vert ^{\alpha}\gamma^{-\alpha}_{n}\sum^{\lfloor\left(n-1\right)r\rfloor}_{k=1}\left\vert a_{k}\right\vert ^{\alpha}\left[1+i\beta^{\#}\text{sgn}\left(\frac{\lambda}{\gamma_{n}}a_{k}\right)\right]
\]
\[
=_{\left(2\right)}-\sigma^{\alpha}\left\vert \lambda\right\vert ^{\alpha}\kappa^{-b}_{n}\sum^{\lfloor\left(n-1\right)r\rfloor}_{k=1}\left\vert a_{k}\right\vert ^{\alpha}\left[1+i\beta^{\#}\text{sgn}\left(\frac{\lambda}{\gamma_{n}}a_{k}\right)\right]\left[1+o(1)\right]
\]
\[
=_{\left(3\right)}-\sigma^{\alpha}\left\vert \lambda\right\vert ^{\alpha}\kappa^{-b}_{n}\sum^{\delta_{n}}_{k=1}\left\vert a_{k}\right\vert ^{\alpha}\left[1+i\beta^{\#}\text{sgn}\left(\frac{\lambda}{\gamma_{n}}a_{k}\right)\right]\left[1+o(1)\right]
\]
\[
=_{\left(4\right)}-\sigma^{\alpha}\left\vert \lambda\right\vert ^{\alpha}\kappa^{-1}_{n}\sum^{\delta_{n}}_{k=1}\left\vert \frac{m_{0}p(k/\kappa_{n})}{1-m}\right\vert ^{\alpha}\left[1+i\beta^{\#}\text{sgn}\left(\frac{\lambda}{\kappa^{1-m}_{n}}a_{k}\right)\right]\left[1+o(1)\right]
\]
\[
=_{\left(5\right)}-\sigma^{\alpha}\left\vert \lambda\right\vert ^{\alpha}\kappa^{-1}_{n}\sum^{\delta_{n}}_{k=1}\left\vert \frac{m_{0}p(k/\kappa_{n})}{1-m}\right\vert ^{\alpha}\left[1+i\beta^{\#}\text{sgn}\left(\lambda\frac{m_{0}p(k/\kappa_{n})}{1-m}\right)\right]\left[1+o(1)\right]
\]
\[
=_{(6)}-\left\vert \lambda\right\vert ^{\alpha}\sigma^{\alpha}\left|m_{0}\right|^{\alpha}\int^{\infty}_{0}\left\vert \frac{p(x)}{1-m}\right\vert ^{\alpha}\left[1+i\beta^{\#}\text{sgn}\left(\lambda\frac{m_{0}p(x)}{1-m}\right)\right]dx+o(1)
\]
\[
=-\left\vert \lambda\right\vert ^{\alpha}\sigma^{\alpha}\left|m_{0}\right|^{\alpha}\left[\int^{\infty}_{0}\left\vert \frac{p(x)}{1-m}\right\vert ^{\alpha}dx+i\beta^{\#}\int^{\infty}_{0}\left\vert \frac{p(x)}{1-m}\right\vert ^{\alpha}\text{sgn}\left(\lambda\frac{m_{0}p(x)}{1-m}\right)dx\right]+o(1)
\]
\[
=-\left\vert \lambda\right\vert ^{\alpha}\sigma^{\alpha}\left|m_{0}\right|^{\alpha}\int^{\infty}_{0}\left\vert \frac{p(x)}{1-m}\right\vert ^{\alpha}dx\left[1+i\beta^{\#}\frac{\int^{\infty}_{0}\left\vert p(x)\right\vert ^{\alpha}\text{sgn}\left(\lambda\frac{p(x)}{1-m}\right)dx}{\int^{\infty}_{0}\left\vert p(x)\right\vert ^{\alpha}dx}\right]+o(1)
\]
\[
=-\left\vert \lambda\right\vert ^{\alpha}\sigma^{\alpha}\left|m_{0}\right|^{\alpha}\int^{\infty}_{0}\left\vert \frac{p(x)}{1-m}\right\vert ^{\alpha}dx\left[1-i\beta^{\#}\text{sgn}\left(m_{0}\right)\frac{\int^{\infty}_{0}\left\vert p(x)\right\vert ^{\alpha}\text{sgn}\left(p(x)\right)dx}{\int^{\infty}_{0}\left\vert p(x)\right\vert ^{\alpha}dx}\text{sgn}\left(\lambda\right)\right]+o(1)
\]
\begin{equation}
=-\left\vert \lambda\right\vert ^{\alpha}\sigma^{\alpha}\left|m_{0}\right|^{\alpha}\mathbb{V}_{2}\left[1-i\beta^{\#}\text{sgn}\left(m_{0}\right)\mathbb{V}_{3}\text{sgn}\left(\lambda\right)\right]+o(1),\label{MN-NM-02}
\end{equation}
where $=_{\left(1\right)}$ follows from (\ref{MI-NM2}), $=_{\left(2\right)}$
the fact that $\gamma^{\alpha}_{n}\sim\kappa^{b}_{n}\ell_{1}(\kappa_{n})\sim\kappa^{b}_{n}$,
$=_{\left(3\right)}$ is a consequence of \textbf{Lemma \ref{LemmaMI1}},
$=_{\left(4\right)}$ is due to \textbf{Corollary \ref{Ts_corol}},
and $=_{\left(6\right)}$ follows from \textbf{Corollary \ref{pi corol}}.
Finally, we show that $=_{\left(5\right)}$ holds. Set 
\[
K_{n,k}:=\left\{ 1\leq k\leq\delta_{n}:\text{sgn}\left(\frac{1}{\kappa^{1-m}_{n}}a_{k}\right)\neq\text{sgn}\left(p(k/\kappa_{n})\right)\right\} .
\]
As $n\rightarrow\infty$, we have
\[
\kappa^{-1}_{n}\left|\sum^{\delta_{n}}_{k=1}\left\vert p(k/\kappa_{n})\right\vert ^{\alpha}\left[\text{sgn}\left(\frac{1}{\kappa^{1-m}_{n}}a_{k}\right)-\text{sgn}\left(p(k/\kappa_{n})\right)\right]\right|
\]
\[
\leq\kappa^{-1}_{n}\sum^{\delta_{n}}_{k=1}\left\vert p(k/\kappa_{n})\right\vert ^{\alpha}\left|\text{sgn}\left(\frac{1}{\kappa^{1-m}_{n}}a_{k}\right)-\text{sgn}\left(p(k/\kappa_{n})\right)\right|
\]
\[
=\kappa^{-1}_{n}\sum_{k\in K_{n}}\left\vert p(k/\kappa_{n})\right\vert ^{\alpha}\left|\text{sgn}\left(\frac{1}{\kappa^{1-m}_{n}}a_{k}\right)-\text{sgn}\left(p(k/\kappa_{n})\right)\right|
\]
\[
\leq2\kappa^{-1}_{n}\sum^{\delta_{n}}_{k=1}\left\vert p(k/\kappa_{n})-\frac{1}{\kappa^{1-m}_{n}}a_{k}\right\vert ^{\alpha}=o(1),
\]
where the last approximation follows from (\ref{MN-NM-01}).

Next, note that an argument analogous to that used for (\ref{MN-NM-02})
gives
\begin{equation}
\sum^{n}_{k=1}\left|G_{nk}\right|=O(1).\label{MN-NM-03}
\end{equation}
Furthermore, 
\[
\gamma^{-1}_{n}\sup_{1\leq k\leq n}\left|a_{k}\right|=\gamma^{-1}_{n}\sup_{1\leq k\leq n}\left|\sum^{k}_{l=0}\rho^{k-l}_{n}c_{l}\right|\leq\gamma^{-1}_{n}\sum^{\infty}_{l=0}\left|c_{l}\right|\rightarrow0,
\]
as $n\rightarrow\infty$. In view of \textbf{Assumption INN}, there
is some $R(x)\rightarrow0$ as $x\rightarrow0$ such that the characteristic
function of $x^{+}_{\lfloor nr\rfloor}/\gamma_{n}$ 
\[
\mathbf{E}\left(e^{i\lambda x^{+}_{\lfloor nr\rfloor}/\gamma_{n}}\right)=\prod^{\lfloor\left(n-1\right)r\rfloor}_{k=1}\psi_{\xi}\left(\frac{\lambda}{\gamma_{n}}a_{k}\right)
\]
\ 
\[
=\prod^{\lfloor\left(n-1\right)r\rfloor}_{k=1}\exp\left\{ G_{nk}\left[1+R\left(\frac{\lambda}{\gamma_{n}}a_{k}\right)\right]\right\} =\exp\left\{ \sum^{\lfloor\left(n-1\right)r\rfloor}_{k=1}G_{nk}\left[1+R\left(\frac{\lambda}{\gamma_{n}}a_{k}\right)\right]\right\} ,
\]
\[
=\exp\left\{ \sum^{\lfloor\left(n-1\right)r\rfloor}_{k=1}G_{nk}+o(1)\right\} =-\left\vert \lambda\right\vert ^{\alpha}\sigma^{\alpha}\left|m_{0}\right|^{\alpha}\mathbb{V}_{2}\left[1-i\beta^{\#}\text{sgn}\left(m_{0}\right)\mathbb{V}_{3}\text{sgn}\left(\lambda\right)\right]+o(1),
\]
where we have used (\ref{MN-NM-02}) and the fact that $\sup_{1\leq k\leq n}\left\vert R\left(\frac{\lambda}{\gamma_{n}}a_{k}\right)\right\vert =o(1)$.
Notice that the latter together with (\ref{MN-NM-03}) implies that
\[
\sum^{n}_{k=1}\left\vert G_{nk}R\left(\frac{\lambda}{\gamma_{n}}a_{k}\right)\right\vert \leq\sup_{k}\left\vert R\left(\frac{\lambda}{\gamma_{n}}a_{k}\right)\right\vert \sum^{n}_{k=1}\left\vert G_{nk}\right\vert =o(1).
\]
This completes the proof. $\blacksquare$\bigskip{}

\noindent\textbf{Proof of} \textbf{Lemma \ref{lem:FSalphaM}}.\textbf{
(a)} We start with part (a). We shall prove the result first under
\textbf{BFR1}, then under \textbf{BFR2} and finally under \textbf{MI}.
Before we present the proofs, we provide a preliminary result that
will be used throughout. Recall that $h(x)$ SV such that (see Remark
\ref{Rmk1})
\[
H(x)=\left\{ \begin{array}{l}
h(x),\text{ }\alpha\in(0,2)\\
-\int^{x}_{0}s^{2}d\left(\frac{h(s)}{s^{2}}\right),\text{ }\alpha=2
\end{array}\right.,
\]
where in the latter case $h(x)=o\left(H(x)\right)$. To see this note
that
\begin{eqnarray*}
\int^{x}_{0}s^{2}d\left(\frac{h(s)}{s^{2}}\right) & = & h(x)-h(0)-\int^{x}_{0}\frac{h(s)}{s^{2}}ds^{2}\\
 & = & h(x)-h(0)-2\int^{x}_{0}\frac{h(s)}{s}ds,
\end{eqnarray*}
with $\int^{x}_{0}\frac{h(s)}{s}ds$ SV (see e.g. BGT 1987; Proposition
1.5.7). Now by BGT 1987, Propopsition 1.5.9a, $\lim_{x\rightarrow\infty}\int^{x}_{0}\frac{h(s)}{s}ds/h(x)=\infty$.
Hence, $h(x)=o\left(H(x)\right)$. Therefore, for $\alpha\in(0,2]$
\begin{equation}
h(x)=O(H(x)),\:\mathrm{as}\:x\rightarrow\infty.\label{a-moments1ab}
\end{equation}

\bigskip{}

\textbf{BFR1}. Let $\varepsilon\in\left(0,\alpha\right)$ and
\[
A_{n,t}:=\mathbf{E}\left\vert \gamma^{-1}_{n}\sum^{t-1}_{k=0}a_{k}\xi_{t-k}\right\vert ^{\alpha-\varepsilon}\text{ and }B_{n,t}:=\mathbf{E}\left\vert \gamma^{-1}_{n}\sum^{\infty}_{k=t}a^{-}_{k,t}\xi_{t-k}\right\vert ^{\alpha-\varepsilon},
\]
Then 
\[
\mathbf{E}\left\vert \gamma^{-1}_{n}x_{t}\right\vert ^{\alpha-\varepsilon}=\mathbf{E}\left\vert \gamma^{-1}_{n}\sum^{\infty}_{k=0}a^{-}_{k,t}\xi_{t-k}\right\vert ^{\alpha-\varepsilon}\leq2\left(A_{n,t}+B_{n,t}\right).
\]
Therefore, it suffices to show that $\sup_{1\leq t\leq n}A_{n,t},\sup_{1\leq t\leq n}B_{n,t}<\infty$.

First, we consider $A_{n,t}$. Let
\[
S_{n,k}:=\gamma_{n}\left\vert a_{k}\right\vert ^{-1}.
\]
Similarly to Astrauskas (1983), set 
\[
\overline{\xi}_{n,k}:=\xi_{k}1\left\{ \left\vert \xi_{k}\right\vert \geq S_{n,k}\right\} \text{ and }\underline{\xi}_{n,k}:=\xi_{k}1\left\{ \left\vert \xi_{k}\right\vert <S_{n,k}\right\} .
\]
 Without loss of generality, for $\alpha\in(1,2]$ we assume that
$\mathbf{E}\overline{\xi}_{n,k}=\mathbf{E}\underline{\xi}_{n,k}=0$.
This convention can be dispensed with if we instead set
\[
\overline{\xi}_{n,k}=\xi_{k}1\left\{ \left\vert \xi_{k}\right\vert \geq S_{n,k}\right\} -\mathbf{E}\left[\xi_{k}1\left\{ \left\vert \xi_{k}\right\vert \geq S_{n,k}\right\} \right],
\]
\begin{equation}
\underline{\xi}_{n,k}=\xi_{k}1\left\{ \left\vert \xi_{k}\right\vert <S_{n,k}\right\} -\mathbf{E}\left[\xi_{k}1\left\{ \left\vert \xi_{k}\right\vert <S_{n,k}\right\} \right].\label{a-moments1abc}
\end{equation}
For more details see footnote \ref{FtNoteAstrauskas }. Fix $k_{0}$$\geq1$.
Note that by \textbf{Assumption INN} we have $\mathbf{E}\left\vert \xi_{k}\right\vert ^{\alpha-\varepsilon}$<$\infty$
-see e.g. Thm 2.6.4 in Ibragimov and Linnik (1971). Hence, as $n\rightarrow\infty$
we get
\[
A_{n,t}=\mathbf{E}\left\vert \gamma^{-1}_{n}\sum^{t-1}_{k=0}a_{k}\xi_{k}\right\vert ^{\alpha-\varepsilon}\leq2\mathbf{E}\left\vert \gamma^{-1}_{n}\sum^{k_{0}-1}_{k=0}a_{k}\xi_{k}\right\vert ^{\alpha-\varepsilon}+2\mathbf{E}\left\vert \gamma^{-1}_{n}\sum^{t-1}_{k=k_{0}}a_{k}\xi_{k}\right\vert ^{\alpha-\varepsilon}
\]
\[
\leq2k_{0}\gamma^{-1}_{n}\sum^{k_{0}-1}_{k=0}\left\vert a_{k}\right\vert ^{\alpha-\varepsilon}\mathbf{E}\left\vert \xi_{k}\right\vert ^{\alpha-\varepsilon}+2\mathbf{E}\left\vert \gamma^{-1}_{n}\sum^{t-1}_{k=k_{0}}a_{k}\xi_{k}\right\vert ^{\alpha-\varepsilon}=2\mathbf{E}\left\vert \gamma^{-1}_{n}\sum^{t-1}_{k=k_{0}}a_{k}\xi_{k}\right\vert ^{\alpha-\varepsilon}+o(1)
\]
\[
\leq4\mathbf{E}\left\vert \gamma^{-1}_{n}\sum^{t-1}_{k=k_{0}}a_{k}\overline{\xi}_{n,k}\right\vert ^{\alpha-\varepsilon}+4\mathbf{E}\left\vert \gamma^{-1}_{n}\sum^{t-1}_{k=k_{0}}a_{k}\underline{\xi}_{n,k}\right\vert ^{\alpha-\varepsilon}+o(1)
\]
\[
\leq4\mathbf{E}\left\vert \gamma^{-1}_{n}\sum^{t-1}_{k=k_{0}}a_{k}\overline{\xi}_{n,k}\right\vert ^{\alpha-\varepsilon}+4\left\{ \mathbf{E}\left\vert \gamma^{-1}_{n}\sum^{t-1}_{k=k_{0}}a_{k}\underline{\xi}_{n,k}\right\vert ^{\alpha+\varepsilon}\right\} ^{\frac{\alpha-\varepsilon}{\alpha+\varepsilon}}+o(1)
\]
\[
\leq_{\left(1\right)}4\left[\gamma^{-\left(\alpha-\varepsilon\right)}_{n}\sum^{n}_{k=k_{0}}\left\vert a_{k}\right\vert ^{\alpha-\varepsilon}\mathbf{E}\left\vert \overline{\xi}_{n,k}\right\vert ^{\alpha-\varepsilon}+\left[\gamma^{-\left(\alpha+\varepsilon\right)}_{n}\sum^{n}_{k=k_{0}}\left\vert a_{k}\right\vert ^{\alpha+\varepsilon}\mathbf{E}\left\vert \underline{\xi}_{n,k}\right\vert ^{\alpha+\varepsilon}\right]^{\frac{\alpha-\varepsilon}{\alpha+\varepsilon}}\right]
\]
\[
=:4\left[A_{1,n}+A^{\frac{\alpha-\varepsilon}{\alpha+\varepsilon}}_{2,n}\right],
\]
where $\leq_{\left(1\right)}$ follows either from Loève's or the
Bahr-Esseen\footnote{\label{FtNoteAstrauskas }We remark that if we set $\overline{\xi}_{n,k}$
and $\underline{\xi}_{n,k}$ as in (\ref{a-moments1abc}), then
\[
\mathbf{E}\left\vert \overline{\xi}_{n,k}\right\vert ^{\alpha-\varepsilon}\leq2\mathbf{E}\left\vert \xi_{k}1\left\{ \left\vert \xi_{k}\right\vert \geq S_{n,k}\right\} \right\vert ^{\alpha-\varepsilon}\:\textrm{and}\:\mathbf{E}\left\vert \underline{\xi}_{n,k}\right\vert ^{\alpha+\varepsilon}\leq4\mathbf{E}\left\vert \xi_{k}1\left\{ \left\vert \xi_{k}\right\vert <S_{n,k}\right\} \right\vert ^{\alpha+\varepsilon}.
\]
In addition, using Lemma 1 in Astrauskas (1983), as $S_{n,k}\rightarrow\infty$,
\begin{equation}
\mathbf{E}\left\vert \xi_{k}1\left\{ \left\vert \xi_{k}\right\vert \geq S_{n,k}\right\} \right\vert ^{\alpha-\varepsilon}\asymp S^{-\varepsilon}_{n,k}h\left(S_{n,k}\right)\text{ and }\mathbf{E}\left\vert \xi_{k}1\left\{ \left\vert \xi_{k}\right\vert <S_{n,k}\right\} \right\vert ^{\alpha+\varepsilon}\asymp S^{\varepsilon}_{n,k}h\left(S_{n,k}\right).\label{FMtypeI.1-1}
\end{equation}
} inequalities depending on the value of the exponents (see p. 133
in Astrauskas, 1983). By Lemma 1 in Astrauskas (1983), for $\alpha\in(0,2]$
as $S_{n,k}\rightarrow\infty$,
\begin{equation}
\mathbf{E}\left\vert \overline{\xi}_{n,k}\right\vert ^{\alpha-\varepsilon}\asymp S^{-\varepsilon}_{n,k}h\left(S_{n,k}\right)\text{ and }\mathbf{E}\left\vert \underline{\xi}_{n,k}\right\vert ^{\alpha+\varepsilon}\asymp S^{\varepsilon}_{n,k}h\left(S_{n,k}\right).\label{FMtypeI.1}
\end{equation}
Note that under \textbf{BFR1},
\[
\left\vert a_{k}\right\vert =\left\vert \sum^{k}_{j=0}c_{k-j}\right\vert \leq\sum^{\infty}_{j=0}\left\vert c_{j}\right\vert =:\epsilon<\infty.
\]
Hence, as $n\rightarrow\infty$ 
\[
\inf_{0\leq k\leq n}S_{n,k}=\gamma_{n}\inf_{0\leq k\leq n}\left\vert a_{k}\right\vert ^{-1}=\gamma_{n}\frac{1}{\sup_{0\leq k\leq n}\left\vert a_{k}\right\vert }\geq\gamma_{n}\epsilon^{-1}\rightarrow\infty.
\]
Therefore, the approximations in (\ref{FMtypeI.1}) hold uniformly
in $0\leq k\leq n$.

Set $R_{-\alpha}\left(x\right):=\left(x\right)^{-\alpha}H\left(x\right)$,
$x>0$. In view of the above and (\ref{a-moments1ab}), as $n\rightarrow\infty$
\[
A_{1,n}=\gamma^{-\left(\alpha-\varepsilon\right)}_{n}\sum^{n}_{k=k_{0}}\left\vert a_{k}\right\vert ^{\alpha-\varepsilon}\mathbf{E}\left\vert \overline{\xi}_{n,k}\right\vert ^{\alpha-\varepsilon}=C_{1}\gamma^{-\left(\alpha-\varepsilon\right)}_{n}\sum^{n}_{k=k_{0}}\left\vert a_{k}\right\vert ^{\alpha-\varepsilon}S^{-\varepsilon}_{n,k}h\left(S_{n,k}\right)\left[1+o(1)\right]
\]
\[
\leq C\gamma^{-\alpha}_{n}\sum^{n}_{k=k_{0}}\left\vert a_{k}\right\vert ^{\alpha}H\left(S_{n,k}\right)\left[1+o(1)\right]=C\gamma^{-\alpha}_{n}\sum^{n}_{k=k_{0}}\left\vert a_{k}\right\vert ^{\alpha}H\left(\gamma_{n}\left\vert a_{k}\right\vert ^{-1}\right)\left[1+o(1)\right]
\]
\[
=C\sum^{n}_{k=k_{0}}R_{-\alpha}\left(\gamma_{n}\left\vert a_{k}\right\vert ^{-1}\right)\left[1+o(1)\right]=:CA_{3,n}\left[1+o(1)\right].
\]
Similarly, 
\[
A_{2,n}=\gamma^{-\left(\alpha+\varepsilon\right)}_{n}\sum^{n}_{k=k_{0}}\left\vert a_{k}\right\vert ^{\alpha+\varepsilon}\mathbf{E}\left\vert \underline{\xi}_{n,k}\right\vert ^{\alpha+\varepsilon}\leq CA_{3,n}\left[1+o(1)\right].
\]
Set $I_{k}:=\int^{\infty}_{k+1}x^{-\left(1+\frac{1}{\alpha}\right)}dx$
$=\alpha\left(k+1\right)^{-1/\alpha}$. Euler summation (see e.g.
Magdalinos, 2012 or \textbf{Lemma \ref{LemmaSV1}} here) gives
\[
\left\vert \sum^{\infty}_{r=k+1}r^{-\left(1+\frac{1}{\alpha}\right)}-\int^{\infty}_{k+1}x^{-\left(1+\frac{1}{\alpha}\right)}dx\right\vert \leq k^{-\left(1+\frac{1}{\alpha}\right)}.
\]
Noting that 
\[
a_{k}=\sum^{k}_{r=0}c_{r}=-\sum^{\infty}_{r=k+1}c_{r},
\]
and $c_{r}\sim r^{-\left(1+\frac{1}{\alpha}\right)}$, for $k$ large
enough and $\eta/2>0$ arbitrarily small we get
\[
\left\vert \sum^{\infty}_{r=k+1}c_{r}-I_{k}\right\vert \leq\sum^{\infty}_{r=k+1}\left\vert c_{r}-r^{-\left(1+\frac{1}{\alpha}\right)}\right\vert +\sum^{\infty}_{r=k+1}\left\vert r^{-\left(1+\frac{1}{\alpha}\right)}-I_{k}\right\vert 
\]
\[
\leq\sum^{\infty}_{r=k+1}\left\vert c_{r}r^{\left(1+\frac{1}{\alpha}\right)}-1\right\vert r^{-\left(1+\frac{1}{\alpha}\right)}+k^{-\left(1+\frac{1}{\alpha}\right)}
\]
\[
\leq\left(\eta/2\right)\sum^{\infty}_{r=k+1}r^{-\left(1+\frac{1}{\alpha}\right)}+k^{-\left(1+\frac{1}{\alpha}\right)}\leq\left(\eta/2\right)\sum^{\infty}_{r=k+1}\left\vert r^{-\left(1+\frac{1}{\alpha}\right)}-I_{k}\right\vert +\left(\eta/2\right)I_{k}+k^{-\left(1+\frac{1}{\alpha}\right)}
\]
\[
\leq\left(\eta/2+1\right)k^{-\left(1+\frac{1}{\alpha}\right)}+\left(\eta/2\right)I_{k}=O\left(k^{-\left(1+\frac{1}{\alpha}\right)}\right)+\left(\eta/2\right)I_{k}.
\]
Hence, as $k\rightarrow\infty$
\[
\left\vert a_{k}/I_{k}-1\right\vert =\left\vert \frac{a_{k}-I_{k}}{I_{k}}\right\vert \leq O(k^{-1})+\eta/2<\eta.
\]
Thus, as $k\rightarrow\infty$ we have
\begin{equation}
\left\vert a_{k}\right\vert \sim\alpha k^{-1/\alpha}\text{ \ and\  }\left\vert a_{k}\right\vert ^{\alpha}\sim\alpha^{\alpha}k^{-1}.\label{FMtypeI.3}
\end{equation}
Let $\eta>0$ arbitrarily small and $k_{0}\geq1$ large enough such
that $\left\vert \frac{\left\vert a_{k}\right\vert ^{\alpha}}{\alpha^{\alpha}k^{-1}}-1\right\vert <\eta$,
for all $k\geq k_{0}$, which possible because of (\ref{FMtypeI.3}).
Hence, as $n\rightarrow\infty$ first and then as $\eta\rightarrow0$
\[
\left\vert A_{4,n}-A_{5,n}\right\vert :=\left\vert \sum^{n}_{k=k_{0}}\gamma^{-\alpha}_{n}\left\vert a_{k}\right\vert ^{\alpha}H\left(\gamma_{n}k^{\frac{1}{\alpha}}\right)-\alpha^{\alpha}\sum^{n}_{k=k_{0}}\gamma^{-\alpha}_{n}k^{-1}H\left(\gamma_{n}k^{\frac{1}{\alpha}}\right)\right\vert 
\]
\[
\leq\alpha^{\alpha}\sum^{n}_{k=k_{0}}\left\vert \frac{\left\vert a_{k}\right\vert ^{\alpha}}{\alpha^{\alpha}k^{-1}}-1\right\vert \gamma^{-\alpha}_{n}k^{-1}H\left(\gamma_{n}k^{\frac{1}{\alpha}}\right)\leq\eta A_{5,n}
\]
\begin{equation}
=_{(1)}\eta\alpha^{\alpha+1}\left(1+o(1)\right)\overset{\eta\rightarrow0}{=}o(1),\label{FMtypeI.3_0a}
\end{equation}
where $=_{(1)}$ follows from (\ref{CF5-FR2}). Next, using (\ref{FMtypeI.3})
again, for $k$ large enough ($k\geq k_{0}$ say) we get $\left\vert k^{-\frac{1}{\alpha}}\left\vert a_{k}\right\vert ^{-1}-\alpha^{-1}\right\vert \leq\delta$,
for all $\delta\in\left(0,\alpha^{-1}\right)$. Therefore, by the
uniform convergence of SV functions (e.g. BGT 1987, Theorem 1.5.2),
for all $\eta>0$ and $k\geq k_{0}$, as $n\rightarrow\infty$
\[
\left\vert \frac{H\left(\gamma_{n}\left\vert a_{k}\right\vert ^{-1}\right)}{H\left(\gamma_{n}k^{\frac{1}{\alpha}}\right)}-1\right\vert \leq\sup_{\left\vert y-\alpha^{-1}\right\vert \leq\delta}\left\vert \frac{H\left(\gamma_{n}k^{\frac{1}{\alpha}}y\right)}{H\left(\gamma_{n}k^{\frac{1}{\alpha}}\right)}-1\right\vert <\eta.
\]
Hence, as $n\rightarrow\infty$ first and then as $\eta\rightarrow0$
\[
\left\vert A_{3,n}-A_{4,n}\right\vert :=\left\vert A_{3,n}-\gamma^{-\alpha}_{n}\sum^{n}_{k=k_{0}}\left\vert a_{k}\right\vert ^{\alpha}H\left(\gamma_{n}k^{\frac{1}{\alpha}}\right)\right\vert 
\]
\[
\leq\sum^{n}_{k=k_{0}}\left\vert \frac{H\left(\gamma_{n}\left\vert a_{k}\right\vert ^{-1}\right)}{H\left(\gamma_{n}k^{\frac{1}{\alpha}}\right)}-1\right\vert \gamma^{-\alpha}_{n}\left\vert a_{k}\right\vert ^{\alpha}H\left(\gamma_{n}k^{\frac{1}{\alpha}}\right)\leq\eta A_{4,n}\overset{\eta\rightarrow0}{\rightarrow}0,
\]
where the last approximation follows from the fact that $A_{4,n}$
is bounded as demonstrated above -cf. (\ref{CF5-FR2}) and (\ref{FMtypeI.3_0a}).
The preceding arguments show that 
\begin{equation}
A_{1,n}=A_{5,n}+o(1)=\alpha^{\alpha+1}+o(1),\label{FMtypeI.3_0b}
\end{equation}
and hence $A_{1,n}$ is bounded. Identical arguments also show that
$A_{2,n}=\alpha^{\alpha+1}+o(1)$. Therefore, we have shown $\sup_{1\leq t\leq n}A_{n,t}<\infty$.

Next, consider $B_{n,t}$. Fix $\eta>0$ and choose $t_{0}$ large
enough such that $\left\vert c_{t}t^{\left(1+\frac{1}{\alpha}\right)}-1\right\vert <\eta$
for $t\geq t_{0}$. Write 
\[
B_{n,t}=\mathbf{E}\left\vert \gamma^{-1}_{n}\sum^{\infty}_{k=t}a^{-}_{k,t}\xi_{k}\right\vert ^{\alpha-\varepsilon}
\]
\[
\leq2\left\{ \mathbf{E}\left\vert \gamma^{-1}_{n}\sum^{t_{0}+t-1}_{k=t}a^{-}_{k,t}\xi_{k}\right\vert ^{\alpha-\varepsilon}+\mathbf{E}\left\vert \gamma^{-1}_{n}\sum^{\infty}_{k=t_{0}+t}a^{-}_{k,t}\xi_{k}\right\vert ^{\alpha-\varepsilon}\right\} :=2\left(B^{\prime}_{n,t}+B^{\prime\prime}_{n,t}\right).
\]
Recall that $a^{-}_{k,t}=\sum^{t-1}_{j=0}c_{k-j}$ (with $k\geq t$),
and $a_{k}=\sum^{k}_{j=0}c_{k-j}$. As $n\rightarrow\infty$, the
first term
\[
B^{\prime}_{n,t}\leq\mathbf{E}\left\vert \gamma^{-1}_{n}\sum^{t_{0}+t-1}_{k=t}a^{-}_{k,t}\xi_{k}\right\vert ^{\alpha-\varepsilon}\leq C\gamma^{-\left(\alpha-\varepsilon\right)}_{n}\sum^{t_{0}+t-1}_{k=t}\left\vert a^{-}_{k,t}\right\vert ^{\alpha-\varepsilon}
\]
\[
\leq\gamma^{-\left(\alpha-\varepsilon\right)}_{n}\sum^{t_{0}+t}_{k=t}\left(\sum^{t}_{j=0}\left\vert c_{k-j}\right\vert \right)^{\alpha-\varepsilon}=\gamma^{-\left(\alpha-\varepsilon\right)}_{n}\sum^{t_{0}+t}_{k=t}\left(\sum^{k}_{j=k-t}\left\vert c_{j}\right\vert \right)^{\alpha-\varepsilon}
\]
\begin{equation}
\leq t_{0}\gamma^{-\left(\alpha-\varepsilon\right)}_{n}\left(\sum^{\infty}_{j=0}\left\vert c_{j}\right\vert \right)^{\alpha-\varepsilon}=o(1),\label{FMtypeI.3_0c}
\end{equation}
uniformly in $1\leq t\leq n$.

Next, we consider $B^{\prime\prime}_{n,t}$. Recall that $t_{0}\in\mathbb{N}$.
First, we shall show that as $n\rightarrow\infty$ 
\begin{equation}
\sup_{1\leq t\leq t_{0}}B^{\prime\prime}_{n,t}=O(1),\label{FMtypeI.3a}
\end{equation}
and then 
\begin{equation}
\sup_{t_{0}\leq t\leq n}B^{\prime\prime}_{n,t}=O(1).\label{FMtypeI.3b}
\end{equation}
We start with the proof of (\ref{FMtypeI.3a}). For $t\leq t_{0}$,
by the Loève or the Bahr-Esseen inequality and Fatou's Lemma
\[
B^{\prime\prime}_{n,t}=\gamma^{-\left(\alpha-\varepsilon\right)}_{n}\mathbf{E}\mathrm{lim}_{M\rightarrow\infty}\left|\sum^{M}_{k=t_{0}+t}\left(\sum^{t-1}_{j=0}c_{k-j}\right)\xi_{k}\right|^{\alpha-\varepsilon}
\]
\[
\leq\gamma^{-\left(\alpha-\varepsilon\right)}_{n}\mathrm{lim\,inf}_{M\rightarrow\infty}\mathbf{E}\left|\sum^{M}_{k=t_{0}+t}\left(\sum^{t-1}_{j=0}c_{k-j}\right)\xi_{k}\right|^{\alpha-\varepsilon}
\]
\[
\leq\gamma^{-\left(\alpha-\varepsilon\right)}_{n}\mathrm{lim\,inf}_{M\rightarrow\infty}\sum^{M}_{k=t_{0}+t}\left|\sum^{t-1}_{j=0}c_{k-j}\right|^{\alpha-\varepsilon}\mathbf{E}\left|\xi_{k}\right|^{\alpha-\varepsilon}
\]
\[
\leq C\gamma^{-\left(\alpha-\varepsilon\right)}_{n}\sum^{\infty}_{k=t_{0}+t}\left|\sum^{t-1}_{j=0}c_{k-j}\right|^{\alpha-\varepsilon}
\]
\[
\leq C\gamma^{-\left(\alpha-\varepsilon\right)}_{n}\sum^{\infty}_{k=0}\sum^{t_{0}}_{j=0}\left\vert c_{k-j}\right\vert ^{\alpha-\varepsilon}\leq Ct_{0}\gamma^{-\left(\alpha-\varepsilon\right)}_{n}\sum^{\infty}_{k=0}\left\vert c_{k}\right\vert ^{\alpha-\varepsilon}=o(1),
\]
which establishes (\ref{FMtypeI.3a}).

Next, we show (\ref{FMtypeI.3b}). Write\footnote{Note that due to \textbf{Assumption INN}, 
\[
\mathbf{E}\left\vert \gamma^{-1}_{t}\sum^{\infty}_{k=t_{0}+t}a^{-}_{k,t}\xi_{t-k}\right\vert ^{\alpha-\varepsilon}=\mathbf{E}\left\vert \gamma^{-1}_{t}\sum^{\infty}_{k=t_{0}+t}a^{-}_{k,t}\xi_{k}\right\vert ^{\alpha-\varepsilon}.
\]
} 

\[
B^{\prime\prime}_{n,t}=\frac{\gamma_{t}}{\gamma_{n}}\mathbf{E}\left\vert \gamma^{-1}_{t}\sum^{\infty}_{k=t_{0}+t}a^{-}_{k,t}\xi_{k}\right\vert ^{\alpha-\varepsilon}=:\frac{\gamma_{t}}{\gamma_{n}}C_{t}.
\]
We shall show $C_{t}$ is bounded. Furthermore, without loss of generality
we assume that $\gamma_{t}$ is increasing. This implies that $\gamma_{t}/\gamma_{n}$
is bounded. Notice that 
\[
\gamma_{t}\sim L(\gamma_{t}t^{1/\alpha})^{1/\alpha},
\]
where $L$ is increasing and $\gamma_{t}t^{1/\alpha}$ is regularly
varying of positive index which has an increasing asymptotic equivalent
(e.g. Theorem 1.5.3 in BGT, 1987). Set $T_{k,t}:=\gamma_{t}\left\vert a^{-}_{k,t}\right\vert ^{-1}$.
Note that as $t\rightarrow\infty$ 
\[
T_{t}:=\inf_{t\leq k}T_{k,t}=\frac{\gamma_{t}}{\sup_{t\leq k}\left\vert a^{-}_{k,t}\right\vert }
\]
\begin{equation}
=\frac{\gamma_{t}}{\sup_{t\leq k}\left\vert \sum^{t-1}_{j=0}c_{k-j}\right\vert }\geq\frac{\gamma_{t}}{\sum^{\infty}_{j=0}\left\vert c_{j}\right\vert }\rightarrow\infty.\label{FM_typeI}
\end{equation}
Hence, as $t\rightarrow\infty$, $T_{k,t}\rightarrow\infty$, uniformly
in $k\geq t$. Furthermore, set 
\[
\overline{\xi}_{k,t}:=\xi_{k}1\left\{ \left\vert \xi_{k}\right\vert \geq T_{k,t}\right\} \text{ and }\underline{\xi}_{k,t}:=\xi_{k}1\left\{ \left\vert \xi_{k}\right\vert <T_{k,t}\right\} .
\]
In view of the above and Lemma 1 in Astrauskas (1983), for $k\geq t$
and $\alpha\in(0,2]$ as $t\rightarrow\infty$ (and hence $T_{k,t}\rightarrow\infty$)
we have
\begin{equation}
\mathbf{E}\left\vert \overline{\xi}_{k,t}\right\vert ^{\alpha-\varepsilon}\asymp T^{-\varepsilon}_{k,t}h\left(T_{k,t}\right)\text{ and }\mathbf{E}\left\vert \underline{\xi}_{k,t}\right\vert ^{\alpha+\varepsilon}\asymp T^{\varepsilon}_{k,t}h\left(T_{k,t}\right).\label{FMtypeI.3da}
\end{equation}
Set $\overline{a}_{k}:=\sum^{k}_{j=0}\left\vert c_{k-j}\right\vert $.
As before, without loss of generality, assume that $\mathbf{E}\overline{\xi}_{k,t}=\mathbf{E}\underline{\xi}_{k,t}=0$.
Consider
\[
2^{-1}C_{t}\leq\mathbf{E}\left\vert \gamma^{-\alpha}_{t}\sum^{\infty}_{k=t+t_{0}}a^{-}_{k,t}\overline{\xi}_{k,t}\right\vert ^{\alpha-\varepsilon}+\left[\mathbf{E}\left\vert \gamma^{-\alpha}_{t}\sum^{\infty}_{k=t+t_{0}}a^{-}_{k,t}\underline{\xi}_{k,t}\right\vert ^{\alpha+\varepsilon}\right]^{\frac{\alpha-\varepsilon}{\alpha+\varepsilon}}
\]
\[
\leq_{\left(2\right)}\gamma^{-\left(\alpha-\varepsilon\right)}_{t}\sum^{\infty}_{k=t+t_{0}}\left\vert a^{-}_{k,t}\right\vert ^{\alpha-\varepsilon}\mathbf{E}\left\vert \overline{\xi}_{k,t}\right\vert ^{\alpha-\varepsilon}+\left[\gamma^{-\left(\alpha+\varepsilon\right)}_{t}\sum^{\infty}_{k=t+t_{0}}\left\vert a^{-}_{k,t}\right\vert ^{\alpha+\varepsilon}\mathbf{E}\left\vert \underline{\xi}_{k,t}\right\vert ^{\alpha+\varepsilon}\right]^{\frac{\alpha-\varepsilon}{\alpha+\varepsilon}}
\]
\[
=:C_{1,t}+C^{\frac{\alpha-\varepsilon}{\alpha+\varepsilon}}_{2,t}.
\]
In view of (\ref{a-moments1ab}) and (\ref{FMtypeI.3da}), as $t\rightarrow\infty$,
the first term
\[
C_{1,t}=\gamma^{-\left(\alpha-\varepsilon\right)}_{t}\sum^{\infty}_{k=t+t_{0}}\left\vert a^{-}_{k,t}\right\vert ^{\alpha-\varepsilon}\mathbf{E}\left\vert \overline{\xi}_{k,t}\right\vert ^{\alpha-\varepsilon}=C\gamma^{-\left(\alpha-\varepsilon\right)}_{t}\sum^{\infty}_{k=t+t_{0}}\left\vert a^{-}_{k,t}\right\vert ^{\alpha-\varepsilon}T^{-\varepsilon}_{k,t}h\left(T_{k,t}\right)\left[1+o(1)\right]
\]
\[
\leq C\sum^{\infty}_{k=t+t_{0}}\left\vert \gamma^{-1}_{t}a^{-}_{k,t}\right\vert ^{\alpha-\varepsilon}\left(\gamma^{-1}_{t}\left\vert a^{-}_{k,t}\right\vert \right)^{\varepsilon}H\left(\gamma_{t}\left\vert a^{-}_{k,t}\right\vert ^{-1}\right)\left[1+o(1)\right]=:C\cdot C_{3,t}\left[1+o(1)\right].
\]
Similarly, as $t\rightarrow\infty$
\[
C_{2,t}=C\cdot C_{3,t}\left[1+o(1)\right].
\]
Write
\[
C_{3,t}=\sum^{\infty}_{k=t+t_{0}}\left\vert \gamma^{-1}_{t}a^{-}_{k,t}\right\vert ^{\alpha}H\left(\gamma_{t}\left\vert a^{-}_{k,t}\right\vert ^{-1}\right)=\sum^{\infty}_{k=t+t_{0}}\left\vert \gamma^{-1}_{t}\sum^{t}_{j=0}c_{k-j}\right\vert ^{\alpha}H\left(\gamma_{t}\left\vert \sum^{t}_{j=0}c_{k-j}\right\vert ^{-1}\right).
\]
Let $0<\delta_{1}\leq1<\delta_{2}\leq2$, and $1\geq\delta:=\delta_{2}-\delta_{1}$.
In addition, $t^{\prime}:=t+t_{0}$. For the purposes of part (b),
we consider the following generalization of $C_{3,t}$ 
\[
C^{\ast}_{3,t}:=\sum^{\infty}_{k=\delta_{1}t+t_{0}}\left\vert \gamma^{-1}_{t}\sum^{\delta_{1}t}_{j=0}c_{k-j}\right\vert ^{\alpha}H\left(\gamma_{t}\left\vert \sum^{\delta_{1}t}_{j=0}c_{k-j}\right\vert ^{-1}\right)
\]
\[
=\left(\sum^{\delta_{2}t^{\prime}}_{k=\delta_{1}t+t_{0}}+\sum^{\infty}_{k=\delta_{2}t^{\prime}+1}\right)\left\vert \gamma^{-1}_{t}\sum^{\delta_{1}t}_{j=0}c_{k-j}\right\vert ^{\alpha}H\left(\gamma_{t}\left\vert \sum^{\delta_{1}t}_{j=0}c_{k-j}\right\vert ^{-1}\right)
\]
\[
=\sum^{\delta t+\delta_{2}t_{0}}_{k=t_{0}}\left\vert \gamma^{-1}_{t}\sum^{k+\delta_{1}t}_{j=k}c_{j}\right\vert ^{\alpha}H\left(\gamma_{t}\left\vert \sum^{k+\delta_{1}t}_{j=k}c_{j}\right\vert ^{-1}\right)
\]
\begin{equation}
+\sum^{\infty}_{k=\delta_{2}t^{\prime}+1}\left\vert \gamma^{-1}_{t}\sum^{k}_{j=k-\delta_{1}t}c_{j}\right\vert ^{\alpha}H\left(\gamma_{t}\left\vert \sum^{k}_{j=k-\delta_{1}t}c_{j}\right\vert ^{-1}\right)=:C^{\prime}_{3,t}+C^{\prime\prime}_{3,t}.\label{FMtypeI.3d}
\end{equation}
Define
\[
Q_{t}(k):=\int^{k+\delta_{1}t}_{k}x^{-\left(1+\frac{1}{\alpha}\right)}dx=\alpha\left[k^{-1/\alpha}-\left(k+\delta_{1}t\right)^{-1/\alpha}\right]
\]
and
\[
I_{t}(k):=\int^{k}_{k-\delta_{1}t}x^{-\left(1+\frac{1}{\alpha}\right)}dx=\alpha\left[\left(k-\delta_{1}t\right)^{-1/\alpha}-k^{-1/\alpha}\right],\text{ }k>\delta_{1}t.
\]
Recall that $R_{-\alpha}(x)=x^{-\alpha}H(x)$, $x>0$ and let \underline{$R$}$_{-\alpha}(x)$
be decreasing such that $R_{-\alpha}(x)\sim$\underline{$R$}$_{-\alpha}(x)$,
see e.g. BGT 1987 Thm1.5.3. Let $\eta\in(0,1)$ and choose $t_{0}$
large enough such that
\[
\left|\frac{Q_{t}(k)}{\sum^{k+\delta_{1}t}_{j=k}c_{j}}-1\right|<\eta,
\]
for all $t_{0}\leq k\leq\delta t+\delta_{2}t_{0}$, which is possible
due to \textbf{Lemma \ref{PropFR1}(aii)}. Using \textbf{Lemma \ref{PropFR1}(b)}
and the uniform convergence of regularly functions (e.g. BGT 1987;
Thm 1.5.2) we have as $t\rightarrow\infty$
\[
\left\vert C^{\prime}_{3,t}-\sum^{\delta t+\delta_{2}t_{0}}_{k=t_{0}}\left|\gamma^{-1}_{t}\sum^{k+\delta_{1}t}_{j=k}c_{j}\right|^{-\alpha}H\left(\gamma_{t}Q_{t}(k)^{-1}\right)\right\vert 
\]
\[
\leq\sum^{\delta t+\delta_{2}t_{0}}_{k=t_{0}}\left\vert \frac{R_{-\alpha}\left(\gamma_{t}\left|\sum^{k+\delta_{1}t}_{j=k}c_{j}\right|^{-1}\right)}{R_{-\alpha}\left(\gamma_{t}Q_{t}(k)^{-1}\right)}-\left[\frac{\gamma_{t}\left|\sum^{k+\delta_{1}t}_{j=k}c_{j}\right|^{-1}}{\gamma_{t}Q_{t}(k)^{-1}}\right]^{-\alpha}\right\vert R_{-\alpha}\left(\gamma_{t}Q_{t}(k)^{-1}\right)
\]
\[
\leq\sup_{1-\eta<y<1+\eta}\left\vert \frac{R_{-\alpha}\left(y\gamma_{t}Q_{t}(k)^{-1}\right)}{R_{-\alpha}\left(\gamma_{t}Q_{t}(k)^{-1}\right)}-y^{-\alpha}\right\vert \sum^{\delta t+\delta_{2}t_{0}}_{k=t_{0}}R_{-\alpha}\left(\gamma_{t}Q_{t}(k)^{-1}\right)
\]
\begin{equation}
=o(1)\cdot\sum^{\delta t+\delta_{2}t_{0}}_{k=t_{0}}R_{-\alpha}\left(\gamma_{t}Q_{t}(k)^{-1}\right).\label{FMtypeI.4a}
\end{equation}
Similarly, using \textbf{Lemma \ref{PropFR1}(aii,b)}, for $t_{0}$
large enough and as $t\rightarrow\infty$ we get
\[
\left\vert \sum^{\delta t+\delta_{2}t_{0}}_{k=t_{0}}\left|\gamma^{-1}_{t}\sum^{k+\delta_{1}t}_{j=k}c_{j}\right|^{-\alpha}H\left(\gamma_{t}Q_{t}(k)^{-1}\right)-\sum^{\delta t+\delta_{2}t_{0}}_{k=t_{0}}R_{-\alpha}\left(\gamma_{t}Q_{t}(k)^{-1}\right)\right\vert 
\]
\[
\leq\sum^{\delta t+\delta_{2}t_{0}}_{k=t_{0}}\left\vert \left(\frac{\gamma_{t}\left|\sum^{k+\delta_{1}t}_{j=k}\right|^{-1}}{\gamma_{t}Q_{t}(k)^{-1}}\right)^{-\alpha}-1\right\vert R_{-\alpha}\left(\gamma_{t}Q_{t}(k)^{-1}\right)
\]
\begin{equation}
\leq o(1)\cdot\sum^{\delta t+\delta_{2}t_{0}}_{k=t_{0}}R_{-\alpha}\left(\gamma_{t}Q_{t}(k)^{-1}\right).\label{FMtypeI.4ab}
\end{equation}
Combining (\ref{FMtypeI.4a})-(\ref{FMtypeI.4ab}),
\[
\left\vert C^{\prime}_{3,t}-\sum^{\delta t+\delta_{2}t_{0}}_{k=t_{0}}R_{-\alpha}\left(\gamma_{t}Q_{t}(k)^{-1}\right)\right\vert \leq o(1)\cdot\sum^{\delta t+\delta_{2}t_{0}}_{k=t_{0}}R_{-\alpha}\left(\gamma_{t}Q_{t}(k)^{-1}\right),
\]
for $t_{0}$ large enough and as $t\rightarrow\infty$. Furthermore,
using the fact that $R_{-\alpha}(x)\sim$\underline{$R$}$_{-\alpha}(x)$
and \textbf{Lemma \ref{PropFR1}(b)}, we get as $t\rightarrow\infty$
\[
\left\vert \sum^{\delta t+\delta_{2}t_{0}}_{k=t_{0}}R_{-\alpha}\left(\gamma_{t}Q_{t}(k)^{-1}\right)-\sum^{\delta t+\delta_{2}t_{0}}_{k=t_{0}}\underline{R}_{-\alpha}\left(\gamma_{t}Q_{t}(k)^{-1}\right)\right\vert \leq o(1)\cdot\sum^{\delta t+\delta_{2}t_{0}}_{k=t_{0}}\underline{R}_{-\alpha}\left(\gamma_{t}Q_{t}(k)^{-1}\right),
\]
which in turn implies that for $t_{0}$ large enough, 
\begin{equation}
\left\vert C^{\prime}_{3,t}-\sum^{\delta t+\delta_{2}t_{0}}_{k=t_{0}}R_{-\alpha}\left(\gamma_{t}Q_{t}(k)^{-1}\right)\right\vert \leq o(1)\cdot\sum^{\delta t+\delta_{2}t_{0}}_{k=t_{0}}\underline{R}_{-\alpha}\left(\gamma_{t}Q_{t}(k)^{-1}\right),\label{FMtypeI.4abc}
\end{equation}
where the $o(1)$ term is negligible as $t\rightarrow\infty$. Set
$t^{\prime}:=t+t_{0}$. Using \textbf{Lemma \ref{PropFR1}(ai,b)}
and a similar argument to the one that yields (\ref{FMtypeI.4abc})
we also get for $t_{0}$ large enough

\begin{equation}
\left\vert C^{\prime\prime}_{3,t}-\sum^{\infty}_{k=\delta_{2}t^{\prime}+1}R_{-\alpha}\left(\gamma_{t}I_{t}(k)^{-1}\right)\right\vert \leq o(1)\cdot\sum^{\infty}_{k=\delta_{2}t^{\prime}+1}\underline{R}_{-\alpha}\left(\gamma_{t}I_{t}(k)^{-1}\right).\label{FMtypeI.4b}
\end{equation}
Next, we show that 
\begin{equation}
\sum^{\delta t+\delta_{2}t_{0}}_{k=t_{0}}\underline{R}_{-\alpha}\left(\gamma_{t}Q_{t}(k)^{-1}\right)=O(1),\label{FMtypeI.5a}
\end{equation}
and 
\begin{equation}
\sum^{\infty}_{k=\delta_{2}t^{\prime}+1}\underline{R}_{-\alpha}\left(\gamma_{t}I_{t}(k)^{-1}\right)=O(1),\label{FMtypeI.5b}
\end{equation}
i.e. the l.h.s. terms of (\ref{FMtypeI.4abc}) and (\ref{FMtypeI.4b})
are $O(1)$. Noting that $\underline{R}_{-\alpha}\left(\gamma_{t}Q_{t}(k)^{-1}\right)$,
$\underline{R}_{-\alpha}\left(\gamma_{t}I_{t}(k)^{-1}\right)$ are
decreasing\footnote{See the last footnote in the proof of \textbf{Lemma \ref{PropFR1}}. 

{}} in $k$, Euler summation gives as $t\rightarrow\infty$
\[
\sum^{\delta t+\delta_{2}t_{0}}_{k=t_{0}}\underline{R}_{-\alpha}\left(\gamma_{t}Q_{t}(k)^{-1}\right)=\int^{\delta t+\delta_{2}t_{0}}_{t_{0}}\underline{R}_{-\alpha}\left(\gamma_{t}Q_{t}(y)^{-1}\right)dy+O\left(\underline{R}_{-\alpha}\left(\gamma_{t}Q_{t}(t_{0})^{-1}\right)\right)
\]
\[
=\int^{\delta t+\delta_{2}t_{0}}_{t_{0}}\underline{R}_{-\alpha}\left(\gamma_{t}Q_{t}(y)^{-1}\right)dy+o\left(1\right)
\]
\begin{equation}
=_{\left(2\right)}\int^{\delta t+\delta_{2}t_{0}}_{t_{0}}R_{-\alpha}\left(\gamma_{t}Q_{t}(y)^{-1}\right)dy\left[1+o(1)\right]+o(1),\label{FMtypeI.6a}
\end{equation}
where $=_{\left(2\right)}$ follows from the fact that that $R_{-\alpha}(x)\sim$\underline{$R$}$_{-\alpha}(x)$.
Similarly, as $t\rightarrow\infty$ 
\[
\sum^{\infty}_{k=\delta_{2}t^{\prime}+1}\underline{R}_{-\alpha}\left(\gamma_{t}I^{-1}_{k,t}\right)=\int^{\infty}_{\delta_{2}t^{\prime}+1}\underline{R}_{-\alpha}\left(\gamma_{t}I_{t}(y)^{-1}\right)dy+O\left(\underline{R}_{-\alpha}\left(\gamma_{t}I_{t}(\delta_{2}t^{\prime}+1)^{-1}\right)\right)
\]
\begin{eqnarray*}
 & = & \int^{\infty}_{\delta_{2}t^{\prime}+1}\underline{R}_{-\alpha}\left(\gamma_{t}I_{t}(y)^{-1}\right)dy+o\left(1\right)\\
 &  & \left[x=y-\delta_{1}t\Leftrightarrow y=x+\delta_{1}t\right]
\end{eqnarray*}
\[
=\int^{\infty}_{\delta t+\delta_{2}t_{0}+1}\underline{R}_{-\alpha}\left(\gamma_{t}Q_{t}(x)^{-1}\right)dx+o\left(1\right)
\]
\begin{equation}
=\int^{\infty}_{\delta t+\delta_{2}t_{0}+1}R_{-\alpha}\left(\gamma_{t}Q_{t}(x)^{-1}\right)dx\left[1+o\left(1\right)\right]+o\left(1\right).\label{FMtypeI.6b}
\end{equation}
In view of (\ref{FMtypeI.6a}) and (\ref{FMtypeI.6b}), (\ref{FMtypeI.5a})
and (\ref{FMtypeI.5b}) will follow if we show that as $t\rightarrow\infty$
\begin{equation}
\int^{\infty}_{t_{0}}R_{-\alpha}\left(\gamma_{t}Q_{t}(x)^{-1}\right)dx=O(1).\label{FMtypeI.6c}
\end{equation}

In fact, the above together with (\ref{FMtypeI.4a})-(\ref{FMtypeI.4b})
and (\ref{FMtypeI.6a})-(\ref{FMtypeI.6b}) imply that as $t\rightarrow\infty$,
\begin{equation}
C^{\prime}_{3,t}=\alpha^{\alpha+1}+o(1)\text{ and }C^{\prime\prime}_{3,t}=o(1).\label{FMtypeI.6d2}
\end{equation}

We proceed with the proof of (\ref{FMtypeI.6c}). We start with the
latter. Consider
\[
\int^{\infty}_{t_{0}}R_{-\alpha}\left(\gamma_{t}Q_{t}(x)^{-1}\right)dx
\]
\[
=\int^{\infty}_{t_{0}}\left\{ \gamma_{t}\alpha^{-1}\left[u^{-1/\alpha}-\left(u+\delta_{1}t\right)^{-1/\alpha}\right]^{-1}\right\} ^{-\alpha}H\left(\gamma_{t}\alpha^{-1}\left[u^{-1/\alpha}-\left(u+\delta_{1}t\right)^{-1/\alpha}\right]^{-1}\right)du
\]
\[
=\int^{\infty}_{t_{0}}\left\{ \gamma_{t}\alpha^{-1}\left(t\delta_{1}\right)^{1/\alpha}\left[\left(\frac{u}{\delta_{1}t}\right)^{-1/\alpha}-\left(\frac{u}{\delta_{1}t}+1\right)^{-1/\alpha}\right]^{-1}\right\} ^{-\alpha}
\]
\[
\times H\left(\gamma_{t}\alpha^{-1}\left(t\delta_{1}\right)^{1/\alpha}\left[\left(\frac{u}{\delta_{1}t}\right)^{-1/\alpha}-\left(\frac{u}{\delta_{1}t}+1\right)^{-1/\alpha}\right]^{-1}\right)du\text{ }\left[z=\frac{u}{\delta_{1}t}\rightarrow du=\delta_{1}tdz\right]
\]
\[
=\delta_{1}t\int^{\infty}_{z=t_{0}/\delta_{1}t}\left\{ \gamma_{t}\alpha^{-1}\left(t\delta_{1}\right)^{1/\alpha}\left[z^{-1/\alpha}-\left(z+1\right)^{-1/\alpha}\right]^{-1}\right\} ^{-\alpha}
\]
\[
\times H\left(\gamma_{t}\alpha^{-1}\left(t\delta_{1}\right)^{1/\alpha}\left[z^{-1/\alpha}-\left(z+1\right)^{-1/\alpha}\right]^{-1}\right)dz
\]
\[
=:\delta_{1}t\int^{\infty}_{z=t_{0}/\delta_{1}t}\left\{ \gamma_{t}\alpha^{-1}\left(t\delta_{1}\right)^{1/\alpha}g\left(z\right)\right\} ^{-\alpha}H\left(\gamma_{t}\alpha^{-1}\left(t\delta_{1}\right)^{1/\alpha}g\left(z\right)\right)dz
\]
\[
=\delta_{1}t\int^{\infty}_{z=t_{0}/\delta_{1}t}R_{-\alpha}\left(\gamma_{t}\alpha^{-1}\left(t\delta_{1}\right)^{1/\alpha}g\left(z\right)\right)dz,
\]
where $g\left(z\right):=\left[z^{-1/\alpha}-\left(z+1\right)^{-1/\alpha}\right]^{-1}$.
Now 
\[
\frac{g\left(z\right)}{z^{1/\alpha}}=z^{-1/\alpha}\left[z^{-1/\alpha}-\left(z+1\right)^{-1/\alpha}\right]^{-1}
\]
\begin{equation}
=\left[1-\left(1+1/z\right)^{-1/\alpha}\right]^{-1}=\frac{1}{1-\frac{1}{\left(1+1/z\right)^{1/\alpha}}}>1,\label{FMtypeI.6d3}
\end{equation}
for all $z>0$. Further, it can be easily seen from the above that
$g\left(z\right)/z^{1/\alpha}\rightarrow1$ as $z\rightarrow0$. Hence,
\begin{equation}
\left|\left\{ g\left(z\right)/z^{1/\alpha}\right\} ^{-\alpha}-1\right|<\varepsilon\text{, for }0<z<\delta_{\varepsilon}.\label{FMtypeI.6d5}
\end{equation}
Without loss of generality choose $\delta_{\varepsilon}<\delta/\delta_{1}$.
Using similar arguments\footnote{Note that 
\begin{eqnarray*}
 &  & t\int^{\delta_{\varepsilon}}_{z=t_{0}/t}R_{-\alpha}\left(\gamma_{t}\alpha^{-1}t^{1/\alpha}z^{1/\alpha}\right)dz\\
 &  & \left[u=\gamma_{t}\alpha^{-1}t^{1/\alpha}z^{1/\alpha}\Rightarrow z=\left(\gamma_{t}\alpha^{-1}t^{1/\alpha}\right)^{-\alpha}u^{\alpha}\Rightarrow dz=\alpha\left(\gamma_{t}\alpha^{-1}t^{1/\alpha}\right)^{-\alpha}u^{\alpha-1}du\right]
\end{eqnarray*}
\[
=\alpha\left(\gamma_{t}\alpha^{-1}t^{1/\alpha}\right)^{-\alpha}t\int^{\gamma_{t}\alpha^{-1}t^{1/\alpha}\delta^{1/\alpha}_{\varepsilon}}_{u=\gamma_{t}\alpha^{-1}t^{1/\alpha}\left(t_{0}/t\right)^{1/\alpha}}R_{-\alpha}\left(u\right)u^{\alpha-1}du
\]
\[
=\alpha\gamma^{-\alpha}_{t}\alpha^{\alpha}\int^{\gamma_{t}\alpha^{-1}t^{1/\alpha}\delta^{1/\alpha}_{\varepsilon}}_{u=\gamma_{t}\alpha^{-1}t^{1/\alpha}\left(t_{0}/t\right)^{1/\alpha}}H\left(u\right)u^{-1}du
\]
\[
=\gamma^{-\alpha}_{t}\alpha^{\alpha+1}\int^{\gamma_{t}\alpha^{-1}t^{1/\alpha}\delta^{1/\alpha}_{\varepsilon}}_{u=\gamma_{t}\alpha^{-1}t^{1/\alpha}_{0}}H\left(u\right)u^{-1}du=\alpha^{\alpha+1}+o(1).
\]
} as those used for (\ref{CF5-FR2}) we have
\begin{equation}
\delta_{1}t\int^{\delta_{\varepsilon}}_{z=t_{0}/\delta_{1}t}R_{-\alpha}\left(\gamma_{t}\alpha^{-1}\left(t\delta_{1}\right)^{1/\alpha}z^{1/\alpha}\right)dz=\alpha^{\alpha+1}+o(1),\label{FMtypeI.6e0}
\end{equation}
as $t\rightarrow\infty$. In view of (\ref{FMtypeI.6e0}), for $\eta>0$
arbitrarily small and $t$ large enough, the difference,
\[
\delta_{1}t\left\vert \int^{\delta_{\varepsilon}}_{z=t_{0}/\delta_{1}t}R_{-\alpha}\left(\gamma_{t}\alpha^{-1}\left(t\delta_{1}\right)^{1/\alpha}g\left(z\right)\right)dz-\left\{ \gamma_{t}\alpha^{-1}\left(t\delta_{1}\right)^{1/\alpha}g\left(z\right)\right\} ^{-\alpha}H\left(\gamma_{t}\alpha^{-1}\left(t\delta_{1}\right)^{1/\alpha}z^{1/\alpha}\right)dz\right\vert 
\]
\[
=\delta_{1}t\left\vert \int^{\delta_{\varepsilon}}_{z=t_{0}/\delta_{1}t}\left\{ \gamma_{t}\alpha^{-1}\left(t\delta_{1}\right)^{1/\alpha}g\left(z\right)\right\} ^{-\alpha}H\left(\gamma_{t}\alpha^{-1}\left(t\delta_{1}\right)^{1/\alpha}g\left(z\right)\right)\right.
\]
\[
\left.-\left\{ \gamma_{t}\alpha^{-1}\left(t\delta_{1}\right)^{1/\alpha}g\left(z\right)\right\} ^{-\alpha}H\left(\gamma_{t}\alpha^{-1}\left(t\delta_{1}\right)^{1/\alpha}z^{1/\alpha}\right)dz\right\vert 
\]
\[
\leq\delta_{1}t\int^{\delta_{\varepsilon}}_{z=t_{0}/\delta_{1}t}\left\vert \left[\frac{\left\{ \gamma_{t}\alpha^{-1}\left(t\delta_{1}\right)^{1/\alpha}g\left(z\right)\right\} ^{-\alpha}H\left(\gamma_{t}\alpha^{-1}\left(t\delta_{1}\right)^{1/\alpha}g\left(z\right)\right)}{\left\{ \gamma_{t}\alpha^{-1}\left(t\delta_{1}\right)^{1/\alpha}z^{1/\alpha}\right\} ^{-\alpha}H\left(\gamma_{t}\alpha^{-1}\left(t\delta_{1}\right)^{1/\alpha}z^{1/\alpha}\right)}-\frac{\left\{ \gamma_{t}\alpha^{-1}\left(t\delta_{1}\right)^{1/\alpha}g\left(z\right)\right\} ^{-\alpha}}{\left\{ \gamma_{t}\alpha^{-1}\left(t\delta_{1}\right)^{1/\alpha}z^{1/\alpha}\right\} ^{-\alpha}}\right]\right.
\]
\[
\times\left.\left\{ \gamma_{t}\alpha^{-1}\left(t\delta_{1}\right)^{1/\alpha}z^{1/\alpha}\right\} ^{-\alpha}H\left(\gamma_{t}\alpha^{-1}\left(t\delta_{1}\right)^{1/\alpha}z^{1/\alpha}\right)\right\vert dz
\]
\[
=\delta_{1}t\int^{\delta_{\varepsilon}}_{z=t_{0}/\delta_{1}t}\left\vert \left[\frac{R_{-\alpha}\left(\gamma_{t}\alpha^{-1}\left(t\delta_{1}\right)^{1/\alpha}z^{1/\alpha}\frac{g\left(z\right)}{z^{1/\alpha}}\right)}{R_{-\alpha}\left(\gamma_{t}\alpha^{-1}\left(t\delta_{1}\right)^{1/\alpha}z^{1/\alpha}\right)}-\frac{\left\{ g\left(z\right)\right\} ^{-\alpha}}{\left\{ z^{1/\alpha}\right\} ^{-\alpha}}\right]R_{-\alpha}\left(\gamma_{t}\alpha^{-1}\left(t\delta_{1}\right)^{1/\alpha}z^{1/\alpha}\right)\right\vert dz
\]
\[
\leq\delta_{1}t\int^{\delta_{\varepsilon}}_{z=t_{0}/\delta_{1}t}\sup_{1\leq y<\infty}\left\vert \frac{R_{-\alpha}\left(\gamma_{t}\alpha^{-1}\left(t\delta_{1}\right)^{1/\alpha}z^{1/\alpha}y\right)}{R_{-\alpha}\left(\gamma_{t}\alpha^{-1}\left(t\delta_{1}\right)^{1/\alpha}z^{1/\alpha}\right)}-y^{-\alpha}\right\vert R_{-\alpha}\left(\gamma_{t}\alpha^{-1}\left(t\delta_{1}\right)^{1/\alpha}z^{1/\alpha}\right)dz.
\]
\begin{equation}
\leq_{\left(3\right)}\eta\delta_{1}t\int^{\delta_{\varepsilon}}_{z=t_{0}/\delta_{1}t}R_{-\alpha}\left(\gamma_{t}\alpha^{-1}\left(t\delta_{1}\right)^{1/\alpha}z^{1/\alpha}\right)dz\leq C\eta,\label{FMtypeI.6e}
\end{equation}
where $\leq_{\left(3\right)}$ follows from (\ref{FMtypeI.6d3}),
the fact that
\[
\gamma_{t}\alpha^{-1}\left(t\delta_{1}\right)^{1/\alpha}\inf_{t_{0}/\delta_{1}t\leq z<\infty}z^{1/\alpha}\rightarrow\infty,
\]
and the uniform convergence of regularly varying functions -see e.g.
BGT 1987, Thm 1.5.2. Next, for $\varepsilon>0$ arbitrarily small
and $t$ large enough, the difference 
\[
\delta_{1}t\left\vert \int^{\delta_{\varepsilon}}_{z=t_{0}/\delta_{1}t}\left[\left\{ \gamma_{t}\alpha^{-1}\left(t\delta_{1}\right)^{1/\alpha}g\left(z\right)\right\} ^{-\alpha}H\left(\gamma_{t}\alpha^{-1}\left(t\delta_{1}\right)^{1/\alpha}z^{1/\alpha}\right)-R_{-\alpha}\left(\gamma_{t}\alpha^{-1}\left(t\delta_{1}\right)^{1/\alpha}z^{1/\alpha}\right)\right]dz\right\vert 
\]
\[
\leq\delta_{1}t\int^{\delta_{\varepsilon}}_{z=t_{0}/\delta_{1}t}\left\vert \frac{\left\{ \gamma_{t}\alpha^{-1}\left(t\delta_{1}\right)^{1/\alpha}g\left(z\right)\right\} ^{-\alpha}}{\left(\gamma_{t}\alpha^{-1}\left(t\delta_{1}\right)^{1/\alpha}z^{1/\alpha}\right)^{-\alpha}}-1\right\vert R_{-\alpha}\left(\gamma_{t}\alpha^{-1}\left(t\delta_{1}\right)^{1/\alpha}z^{1/\alpha}\right)dz
\]
\[
=\delta_{1}t\int^{\delta_{\varepsilon}}_{z=t_{0}/\delta_{1}t}\left\vert \frac{\left\{ g\left(z\right)\right\} ^{-\alpha}}{\left(z^{1/\alpha}\right)^{-\alpha}}-1\right\vert R_{-\alpha}\left(\gamma_{t}\alpha^{-1}\left(t\delta_{1}\right)^{1/\alpha}z^{1/\alpha}\right)dz
\]
\begin{equation}
\leq\varepsilon\delta_{1}t\int^{\delta_{\varepsilon}}_{z=t_{0}/\delta_{1}t}R_{-\alpha}\left(\gamma_{t}\alpha^{-1}\left(t\delta_{1}\right)^{1/\alpha}z^{1/\alpha}\right)dz,\label{FMtypeI.6f}
\end{equation}
where the last inequality follows from (\ref{FMtypeI.6d5}). In view
of (\ref{FMtypeI.6e0})-(\ref{FMtypeI.6f}), as $t\rightarrow\infty$
we get 
\[
\int^{\delta_{\varepsilon}t}_{t_{0}}R_{-\alpha}\left(\gamma_{t}Q_{t}(x)^{-1}\right)dx
\]
\begin{equation}
=\delta_{1}t\int^{\delta_{\varepsilon}}_{z=t_{0}/\delta_{1}t}R_{-\alpha}\left(\gamma_{t}\alpha^{-1}\left(t\delta_{1}\right)^{1/\alpha}g\left(z\right)\right)dx=O(1).\label{FMtypeI.7a}
\end{equation}
Next, we show that
\begin{equation}
\int^{\infty}_{\delta_{\varepsilon}t}R_{-\alpha}\left(\gamma_{t}Q_{t}(x)^{-1}\right)dx=o(1),\label{FMtypeI.7b}
\end{equation}
as $t\rightarrow\infty$ . Note that (\ref{FMtypeI.6c}) follows immediately
from (\ref{FMtypeI.7a}) and (\ref{FMtypeI.7b}).

There is no loss of generality assuming that $R_{-\alpha}(x)$ is
decreasing -see the arguments that lead to (\ref{FMtypeI.6b}) and
Theorem 1.5.3 BGT 1987. Using the mean value theorem 
\[
z^{-1/\alpha}-\left(z+1\right)^{-1/\alpha}=\alpha^{-1}\left(z+\overline{z}\right)^{-1-\frac{1}{\alpha}}=\alpha^{-1}\left(z+\overline{z}\right)^{-\frac{1+\alpha}{\alpha}},
\]
where $0\leq\overline{z}\leq1$. Hence, $g(z)\geq\alpha z^{\frac{1+\alpha}{\alpha}}$.
In view of this and setting $\bar{\delta}:=\alpha\delta^{\frac{1}{\alpha}}_{1}\left(\delta_{\varepsilon}/\delta_{1}\right)^{\frac{1+\alpha}{\alpha}}$
we get as $t\rightarrow\infty$,
\[
\int^{\infty}_{\delta_{\varepsilon}t}R_{-\alpha}\left(\gamma_{t}Q_{t}(x)^{-1}\right)dx=\delta_{1}t\int^{\infty}_{\delta_{\varepsilon}/\delta_{1}}R_{-\alpha}\left(\gamma_{t}\left(t\delta_{1}\right)^{1/\alpha}g(z)\right)dz
\]
\[
\leq\delta_{1}t\int^{\infty}_{\delta_{\varepsilon}/\delta_{1}}R_{-\alpha}\left(\alpha\gamma_{t}\left(t\delta_{1}\right)^{1/\alpha}z^{\frac{1+\alpha}{\alpha}}\right)dz
\]
\[
\left[\begin{gathered}u=\alpha\gamma_{t}\left(t\delta_{1}\right)^{1/\alpha}z^{\frac{1+\alpha}{\alpha}}\rightarrow\left(\alpha^{-1}\gamma^{-1}_{t}\left(t\delta_{1}\right)^{-1/\alpha}u\right)^{\frac{\alpha}{1+\alpha}}=z\\
\rightarrow\frac{\alpha}{1+\alpha}\left(\alpha^{-1}\gamma^{-1}_{t}\left(t\delta_{1}\right)^{-1/\alpha}\right)^{\frac{\alpha}{1+\alpha}}u^{-\frac{1}{\alpha+1}}=dz/du\end{gathered}\right]
\]
\[
=\alpha^{-\frac{\alpha}{1+\alpha}}\frac{\alpha}{1+\alpha}\delta_{1}t\left(\gamma^{-1}_{t}\left(t\delta_{1}\right)^{-1/\alpha}\right)^{\frac{\alpha}{1+\alpha}}\int^{\infty}_{\bar{\delta}\gamma_{t}t^{1/\alpha}}u^{-\frac{1}{\alpha+1}}R_{-\alpha}\left(u\right)du
\]

\[
=\alpha^{-\frac{\alpha}{1+\alpha}}\frac{\alpha}{1+\alpha}\delta_{1}t\gamma^{-\frac{\alpha}{1+\alpha}}_{t}\left(t\delta_{1}\right)^{-\frac{1}{1+\alpha}}\int^{\infty}_{\bar{\delta}\gamma_{t}t^{1/\alpha}}u^{-\left(\frac{\alpha^{2}}{\alpha+1}+1\right)}H\left(u\right)du
\]
\[
\sim_{\left(4\right)}\alpha^{-\frac{\alpha}{1+\alpha}}\frac{\alpha}{1+\alpha}\frac{\alpha+1}{\alpha^{2}}\delta_{1}t\gamma^{-\frac{\alpha}{1+\alpha}}_{t}\left(t\delta_{1}\right)^{-\frac{1}{1+\alpha}}\left(\bar{\delta}\gamma_{t}t^{1/\alpha}\right)^{-\left(\frac{\alpha^{2}}{\alpha+1}\right)}H\left(\bar{\delta}\gamma_{t}t^{1/\alpha}\right)
\]

\[
=Ct\gamma^{-\alpha}_{t}t^{-\frac{1}{1+\alpha}}t^{-\left(\frac{\alpha}{\alpha+1}\right)}H\left(\bar{\delta}\gamma_{t}t^{1/\alpha}\right)=C\gamma^{-\alpha}_{t}H\left(\bar{\delta}\gamma_{t}t^{1/\alpha}\right)=_{\left(5\right)}o(1),
\]
where $\sim_{\left(4\right)}$ and $=_{\left(5\right)}$ follow from
BGT 1987 Proposition 1.5.10 and 1.5.9a respectively.

\bigskip{}

\textbf{BFR2. }Let $\varepsilon\in(0$,$\alpha)$. Consider the term
\[
\mathbf{E}\left\vert \gamma^{-1}_{n}x_{t}\right\vert ^{\alpha-\varepsilon}=\mathbf{E}\left\vert \gamma^{-1}_{n}\sum^{\infty}_{k=0}a_{k,t}\xi_{t-k}\right\vert ^{\alpha-\varepsilon}
\]
\[
\leq2\mathbf{E}\left\vert \gamma^{-1}_{n}\sum^{t-1}_{k=0}a_{k}\xi_{t-k}\right\vert ^{\alpha-\varepsilon}+2\mathbf{E}\left\vert \gamma^{-1}_{n}\sum^{\infty}_{k=t}a_{k,t}\xi_{t-k}\right\vert ^{\alpha-\varepsilon}
\]
\[
=:2\left[A_{n,t}+B_{n,t}\right].
\]
Recall that $a_{k}=\sum^{k}_{j=0}\phi_{j}c_{k-j}$. Set $S_{n,k}:=\gamma_{n}\left\vert a_{k}\right\vert ^{-1}$,
and note that
\[
\left\vert a_{k}\right\vert \leq\sup_{j\geq0}\left\vert \phi_{j}\right\vert \sum^{\infty}_{j=0}\left\vert c_{j}\right\vert =:\epsilon<\infty.
\]
Hence, 
\begin{equation}
\inf_{0\leq k\leq n}S_{n,k}=\gamma_{n}\inf_{0\leq k\leq n}\left\vert a_{k}\right\vert ^{-1}=\gamma_{n}\frac{1}{\sup_{0\leq k\leq n}\left\vert a_{k}\right\vert }\geq\gamma_{n}\epsilon^{-1}\rightarrow\infty.\label{a-moments1}
\end{equation}
Next, set 
\[
\overline{\xi}_{n,k}:=\xi_{k}1\left\{ \left\vert \xi_{k}\right\vert \geq S_{n,k}\right\} \text{ and }\underline{\xi}_{n,k}:=\xi_{k}1\left\{ \left\vert \xi_{k}\right\vert <S_{n,k}\right\} .
\]
As before, and without loss of generality we assume that $\mathbf{E}\overline{\xi}_{n,k}=0\text{ and }\mathbf{E}\underline{\xi}_{n,k}=0$,
when $\alpha\in(1,2]$. Let $k_{0}$ be some fixed integer. As $n\rightarrow\infty$,
we have
\[
A_{n,t}=\mathbf{E}\left\vert \gamma^{-1}_{n}\sum^{t-1}_{k=0}a_{k}\xi_{k}\right\vert ^{\alpha-\varepsilon}\leq2\mathbf{E}\left\vert \gamma^{-1}_{n}\sum^{k_{0}-1}_{k=0}a_{k}\xi_{k}\right\vert ^{\alpha-\varepsilon}+2\mathbf{E}\left\vert \gamma^{-1}_{n}\sum^{t-1}_{k=k_{0}}a_{k}\xi_{k}\right\vert ^{\alpha-\varepsilon}
\]
\[
\leq2\gamma^{-1}_{n}\sum^{k_{0}-1}_{k=0}\left\vert a_{k}\right\vert ^{\alpha-\varepsilon}\mathbf{E}\left\vert \xi_{k}\right\vert ^{\alpha-\varepsilon}+2\mathbf{E}\left\vert \gamma^{-1}_{n}\sum^{t-1}_{k=k_{0}}a_{k}\xi_{k}\right\vert ^{\alpha-\varepsilon}
\]
\[
\leq4\mathbf{E}\left\vert \gamma^{-1}_{n}\sum^{t-1}_{k=k_{0}}a_{k}\overline{\xi}_{n,k}\right\vert ^{\alpha-\varepsilon}+4\left\{ \mathbf{E}\left\vert \gamma^{-1}_{n}\sum^{t-1}_{k=k_{0}}a_{k}\underline{\xi}_{n,k}\right\vert ^{\alpha+\varepsilon}\right\} ^{\frac{\alpha-\varepsilon}{\alpha+\varepsilon}}+o(1)
\]
\[
\leq_{\left(1\right)}C\left\{ \gamma^{-\left(\alpha-\varepsilon\right)}_{n}\sum^{n}_{k=k_{0}}\left\vert a_{k}\right\vert ^{\alpha-\varepsilon}\mathbf{E}\left\vert \overline{\xi}_{n,k}\right\vert ^{\alpha-\varepsilon}+\left[\gamma^{-\left(\alpha+\varepsilon\right)}_{n}\sum^{n}_{k=k_{0}}\left\vert a_{k}\right\vert ^{\alpha+\varepsilon}\mathbf{E}\left\vert \underline{\xi}_{n,k}\right\vert ^{\alpha+\varepsilon}\right]^{\frac{\alpha-\varepsilon}{\alpha+\varepsilon}}\right\} +o(1)
\]
\begin{equation}
=:C\left\{ A_{1,n}+\left[A_{2,n}\right]^{\frac{\alpha-\varepsilon}{\alpha+\varepsilon}}\right\} +o(1).\label{a-moments1abb}
\end{equation}
where $\leq_{\left(1\right)}$ follows either from Loève's or the
Bahr-Esseen inequalities depending on the value of the exponents (see
p. 133 in Astrauskas, 1983). .

By Lemma 1 in Astrauskas (1983), for $\alpha\in(0,2]$ as $S_{n,k}\rightarrow\infty$,
\begin{equation}
\mathbf{E}\left\vert \overline{\xi}_{n,k}\right\vert ^{\alpha-\varepsilon}\asymp S^{-\varepsilon}_{n,k}h\left(S_{n,k}\right)\text{ and }\mathbf{E}\left\vert \underline{\xi}_{n,k}\right\vert ^{\alpha+\varepsilon}\asymp S^{\varepsilon}_{n,k}h\left(S_{n,k}\right).\label{a-moments1a}
\end{equation}
In view of (\ref{a-moments1}) and (\ref{a-moments1a}) as $n\rightarrow\infty$
\[
A_{1,n}=C\gamma^{-\left(\alpha-\varepsilon\right)}_{n}\sum^{n}_{k=k_{0}}\left\vert a_{k}\right\vert ^{\alpha-\varepsilon}S^{-\varepsilon}_{n,k}h\left(S_{n,k}\right)\left[1+o(1)\right]
\]
\[
=C\gamma^{-\alpha}_{n}\sum^{n}_{k=k_{0}}\left\vert a_{k}\right\vert ^{\alpha}h\left(S_{n,k}\right)\left[1+o(1)\right]\leq_{(2)}C\gamma^{-\alpha}_{n}\sum^{n}_{k=k_{0}}\left\vert a_{k}\right\vert ^{\alpha}H\left(S_{n,k}\right)\left[1+o(1)\right]
\]
\begin{equation}
=C\gamma^{-\alpha}_{n}\sum^{n}_{k=k_{0}}\left\vert a_{k}\right\vert ^{\alpha}H\left(\gamma_{n}\left\vert a_{k}\right\vert ^{-1}\right)\left[1+o(1)\right]=:C\cdot A_{3,n}\left[1+o(1)\right],\label{a-moments1abbb}
\end{equation}
where $\leq_{(2)}$ follows from (\ref{a-moments1ab}). Similarly,
as $n\rightarrow\infty$
\begin{equation}
A_{2,n}\leq C\cdot A_{3,n}\left[1+o(1)\right].\label{a-moments1b}
\end{equation}
Without loss of generality set $c_{\infty}=1$. Fix $\varepsilon\in\left(0,1\right)$.
By \textbf{Lemma \ref{cauchy}}, and \textbf{Assumption BFR2} for
$k$ large enough, $k\geq k_{0}$ say, we have 
\begin{equation}
\left\vert c_{\infty}k^{-\frac{1}{\alpha}}\left\vert a_{k}\right\vert ^{-1}-1\right\vert \leq\varepsilon.\label{a-moments1c}
\end{equation}
In view of this, and the uniform convergence of regularly varying
functions (e.g. BGT 1987; Thm 1.5.2) as $n\rightarrow\infty$
\[
\sup_{k\geq k_{0}}\left\vert \frac{\left[\left\vert a_{k}\right\vert ^{-1}k^{-1/\alpha}\right]^{-\alpha}H\left(\gamma_{n}\left\vert a_{k}\right\vert ^{-1}\right)}{H\left(\gamma_{n}k^{\frac{1}{\alpha}}\right)}-\left[\left\vert a_{k}\right\vert ^{-1}k^{-1/\alpha}\right]^{-\alpha}\right\vert 
\]
\begin{equation}
\leq\sup_{k\geq k_{0}}\sup_{y\in[1-\varepsilon,\infty)}\left\vert \frac{\left(\gamma_{n}k^{\frac{1}{\alpha}}y\right)^{-\alpha}H\left(\gamma_{n}k^{\frac{1}{\alpha}}y\right)}{\left(\gamma_{n}k^{\frac{1}{\alpha}}\right)^{-\alpha}H\left(\gamma_{n}k^{\frac{1}{\alpha}}\right)}-y^{-\alpha}\right\vert =o(1).\label{a-moments2}
\end{equation}
Hence,
\[
\left\vert \gamma^{-\alpha}_{n}\sum^{n}_{k=k_{0}}\left\vert a_{k}\right\vert ^{\alpha}H\left(\gamma_{n}\left\vert a_{k}\right\vert ^{-1}\right)-\gamma^{-\alpha}_{n}\sum^{n}_{k=k_{0}}\left\vert a_{k}\right\vert ^{\alpha}H\left(\gamma_{n}k^{\frac{1}{\alpha}}\right)\right\vert 
\]
\[
\leq\sum^{n}_{k=k_{0}}\left\vert \frac{\left\vert a_{k}\right\vert ^{\alpha}H\left(\gamma_{n}\left\vert a_{k}\right\vert ^{-1}\right)}{H\left(\gamma_{n}k^{\frac{1}{\alpha}}\right)}-\left\vert a_{k}\right\vert ^{\alpha}\right\vert \gamma^{-\alpha}_{n}H\left(\gamma_{n}k^{\frac{1}{\alpha}}\right)
\]
\[
=\sum^{n}_{k=k_{0}}\left\vert \frac{\left\vert a_{k}\right\vert ^{\alpha}\left[k^{-1/\alpha}\right]^{-\alpha}H\left(\gamma_{n}\left\vert a_{k}\right\vert ^{-1}\right)}{H\left(\gamma_{n}k^{\frac{1}{\alpha}}\right)}-\left\vert a_{k}\right\vert ^{\alpha}\left[k^{-1/\alpha}\right]^{-\alpha}\right\vert \left[\gamma_{n}k^{1/\alpha}\right]^{-\alpha}H\left(\gamma_{n}k^{\frac{1}{\alpha}}\right)
\]
\begin{equation}
\leq_{(3)}o(1)\cdot\sum^{n}_{k=k_{0}}\gamma^{-\alpha}_{n}k^{-1}H\left(\gamma_{n}k^{\frac{1}{\alpha}}\right)=_{(4)}o(1)\cdot\left(\alpha+o(1)\right)=o(1),\label{a-moments2a}
\end{equation}
where $\leq_{(3)}$ and $=_{(4)}$ follow from (\ref{a-moments2})
and (\ref{CF5-FR2}) respectively. Next, fix $\varepsilon>0$ and
let $k$ large enough, $k\geq k_{0}$ say, such that $\left\vert \left\vert a_{k}\right\vert ^{\alpha}k-1\right\vert <\varepsilon$,
which is possible to due to \textbf{Assumption BFR2}. As $n\rightarrow\infty$
first and then as $\varepsilon\rightarrow0$,
\[
\left\vert \gamma^{-\alpha}_{n}\sum^{n}_{k=k_{0}}\left\vert a_{k}\right\vert ^{\alpha}H\left(\gamma_{n}k^{\frac{1}{\alpha}}\right)-\gamma^{-\alpha}_{n}\sum^{n}_{k=k_{0}}k^{-1}H\left(\gamma_{n}k^{\frac{1}{\alpha}}\right)\right\vert 
\]
\begin{equation}
\leq C\gamma^{-\alpha}_{n}\sum^{n}_{k=k_{0}}\left\vert \frac{\left\vert a_{k}\right\vert ^{\alpha}}{k^{-1}}-1\right\vert k^{-1}H\left(\gamma_{n}k^{\frac{1}{\alpha}}\right)\leq C\varepsilon\gamma^{-\alpha}_{n}\sum^{n}_{k=k_{0}}k^{-1}H\left(\gamma_{n}k^{\frac{1}{\alpha}}\right).\label{a-moments2b}
\end{equation}
Therefore, combining (\ref{a-moments2a}) and (\ref{a-moments2b})
as $n$$\rightarrow\infty$ we get
\begin{equation}
A_{3,n}=\gamma^{-\alpha}_{n}\sum^{n}_{k=k_{0}}\left\vert a_{k}\right\vert ^{\alpha}H\left(\gamma_{n}\left\vert a_{k}\right\vert ^{-1}\right)\leq C\gamma^{-\alpha}_{n}\sum^{n}_{k=k_{0}}k^{-1}H\left(\gamma_{n}k^{1/\alpha}\right)+o(1).\label{a-moments2b2}
\end{equation}

Using (\ref{a-moments1abb}), (\ref{a-moments1abbb}), (\ref{a-moments1b}),
and (\ref{a-moments2b2}) we get as $n\rightarrow\infty$
\[
\sup_{1\leq t\leq n}A_{n,t}\leq o(1)+C\gamma^{-\alpha}_{n}\sum^{n}_{k=k_{0}}k^{-1}H\left(\gamma_{n}k^{\frac{1}{\alpha}}\right)\left[1+o(1)\right]
\]
\[
+C\left\{ \gamma^{-\alpha}_{n}\sum^{n}_{k=k_{0}}k^{-1}H\left(\gamma_{n}k^{\frac{1}{\alpha}}\right)\left[1+o(1)\right]\right\} ^{\frac{\alpha-\varepsilon}{\alpha+\varepsilon}}
\]
\[
=_{(5)}C\left(\alpha+\alpha^{\frac{\alpha-\varepsilon}{\alpha+\varepsilon}}\right)+o(1),
\]
where $=_{(5)}$ follows from (\ref{CF5-FR2}).

Next, we show that $B_{n,t}$ is uniformly bounded for $1\leq t\leq n$.
Write
\[
B_{n,t}=\mathbf{E}\left\vert \gamma^{-\alpha}_{n}\sum^{\infty}_{k=t}a_{k,t}\xi_{t-k}\right\vert ^{\alpha-\varepsilon}=\mathbf{E}\left\vert \gamma^{-\alpha}_{n}\sum^{\infty}_{k=t}a_{k,t}\xi_{k}\right\vert ^{\alpha-\varepsilon}.
\]
Recall that $a_{k,t}=\sum^{\left(k-1\right)\wedge t}_{j=0}\phi_{j}c_{k-j}$
and $a_{k}=\sum^{k}_{j=0}\phi_{j}c_{k-j}$. Furthermore, set $\overline{a}_{k}=\sum^{k}_{j=0}\left\vert \phi_{j}c_{k-j}\right\vert $,
and note that $\left\vert a_{k}\right\vert ,\left\vert a_{k,t}\right\vert \leq\overline{a}_{k}$.
Consider
\[
\mathbf{E}\left\vert \gamma^{-\alpha}_{n}\sum^{\infty}_{k=t}a_{k,t}\xi_{k}\right\vert ^{\alpha-\varepsilon}\leq2\mathbf{E}\left\vert \gamma^{-\alpha}_{n}\sum^{\infty}_{k=t}a_{k,t}\overline{\xi}_{n,k}\right\vert ^{\alpha-\varepsilon}+2\left[\mathbf{E}\left\vert \gamma^{-\alpha}_{n}\sum^{\infty}_{k=t}a_{k,t}\underline{\xi}_{n,k}\right\vert ^{\alpha+\varepsilon}\right]^{\frac{\alpha-\varepsilon}{\alpha+\varepsilon}}
\]
\[
\leq_{\left(6\right)}2\gamma^{-\left(\alpha-\varepsilon\right)}_{n}\sum^{\infty}_{k=t}\left\vert a_{k,t}\right\vert ^{\alpha-\varepsilon}\mathbf{E}\left\vert \overline{\xi}_{n,k}\right\vert ^{\alpha-\varepsilon}+2\left[\gamma^{-\left(\alpha+\varepsilon\right)}_{n}\sum^{\infty}_{k=t}\left\vert a_{k,t}\right\vert ^{\alpha+\varepsilon}\mathbf{E}\left\vert \underline{\xi}_{n,k}\right\vert ^{\alpha+\varepsilon}\right]^{\frac{\alpha-\varepsilon}{\alpha+\varepsilon}},
\]
where $\leq_{\left(6\right)}$ follows from a similar argument to
the one used for $\leq_{\left(1\right)}$. Recall that (Assumption
\textbf{SM}) $\sum^{\infty}_{j=0}\left\vert c_{j}\right\vert ^{\nu}<\infty$
where $\nu=1\wedge\tau$, $\tau\in\left(0,\alpha\right)$. Choose
$\varepsilon>0$ such that $\alpha-\nu-\varepsilon>0$. Once more,
using (\ref{a-moments1}), (\ref{a-moments1a}) and (\ref{a-moments1ab}),
as $n\rightarrow\infty$ the first term above is
\[
\gamma^{-\left(\alpha-\varepsilon\right)}_{n}\sum^{\infty}_{k=t}\left\vert a_{k,t}\right\vert ^{\alpha-\varepsilon}\mathbf{E}\left\vert \overline{\xi}_{n,k}\right\vert ^{\alpha-\varepsilon}=C\gamma^{-\left(\alpha-\varepsilon\right)}_{n}\sum^{\infty}_{k=t}\left\vert a_{k,t}\right\vert ^{\alpha-\varepsilon}S^{-\varepsilon}_{n,k}h\left(S_{n,k}\right)\left[1+o(1)\right]
\]
\[
\leq C\gamma^{-\left(\alpha-\varepsilon\right)}_{n}\sum^{\infty}_{k=t}\left\vert a_{k,t}\right\vert ^{\alpha-\varepsilon}\gamma^{-\varepsilon}_{n}\left\vert a_{k}\right\vert ^{\varepsilon}H\left(\gamma_{n}\left\vert a_{k}\right\vert ^{-1}\right)\left[1+o(1)\right]
\]
\[
=C\gamma^{-\alpha}_{n}\sum^{\infty}_{k=t}\left\vert a_{k,t}\right\vert ^{\alpha-\nu-\varepsilon}\left\vert a_{k,t}\right\vert ^{\nu}\left\vert a_{k}\right\vert ^{\varepsilon}H\left(\gamma_{n}\left\vert a_{k}\right\vert ^{-1}\right)\left[1+o(1)\right]
\]
\[
\leq C\gamma^{-\alpha}_{n}\sum^{\infty}_{k=t}\left(\overline{a}_{k}\right)^{\alpha-\nu-\varepsilon}\left\vert a_{k,t}\right\vert ^{\nu}\left(\overline{a}_{k}\right)^{\varepsilon}H\left(\gamma_{n}\left\vert a_{k}\right\vert ^{-1}\right)\left[1+o(1)\right]
\]
\[
=C\gamma^{-\alpha}_{n}\sum^{\infty}_{k=t}\left(\overline{a}_{k}\right)^{\alpha-\nu}\left\vert a_{k,t}\right\vert ^{\nu}H\left(\gamma_{n}\left\vert a_{k}\right\vert ^{-1}\right)\left[1+o(1)\right]=:C\cdot B_{1,n}(t)\left[1+o(1)\right].
\]
Similarly, 
\[
\gamma^{-\left(\alpha+\varepsilon\right)}_{n}\sum^{\infty}_{k=t}\left\vert a_{k,t}\right\vert ^{\alpha+\varepsilon}\mathbf{E}\left\vert \underline{\xi}_{n,k}\right\vert ^{\alpha+\varepsilon}\leq C\cdot B_{1,n}(t)\left[1+o(1)\right].
\]
Next, note that 
\[
B_{1,n}(t)\leq\gamma^{-\alpha}_{n}\sum^{\infty}_{k=t}\left(\overline{a}_{k}\right)^{\alpha-\nu}H\left(\gamma_{n}\left\vert a_{k}\right\vert ^{-1}\right)\left(\sum^{t-1}_{j=0}\left\vert \phi_{j}c_{k-j}\right\vert \right)^{\nu}
\]
\[
\leq\gamma^{-\alpha}_{n}\sum^{t-1}_{j=0}\left\vert \phi_{j}\right\vert ^{\nu}\sum^{\infty}_{k=t}\left(\overline{a}_{k}\right)^{\alpha-\nu}H\left(\gamma_{n}\left\vert a_{k}\right\vert ^{-1}\right)\left\vert c_{k-j}\right\vert ^{\nu},
\]
where recall that by definition $0<\nu\leq1$.

Next, choose $0<\delta<\alpha-\nu$. Without loss of generality assume
that $\sum^{\infty}_{j=0}\left\vert c_{j}\right\vert =1$ (and recall
that we have already assumed that $\sum^{\infty}_{j=0}c_{j}=1$).
For $0<\eta<1$ arbitrarily small let $t_{0}\in\mathbb{N}$ large
enough such that
\[
\left\vert \frac{a_{k}}{k^{-1/\alpha}}-1\right\vert ,\left\vert \frac{k^{-1/\alpha}}{a_{k}}-1\right\vert ,\left\vert \frac{\overline{a}_{k}}{k^{-1/\alpha}}-1\right\vert <\eta,
\]
for all $k\geq t_{0}$, which is possible due to \textbf{Lemma \ref{cauchy}}.
Then as $n\rightarrow\infty$ we have
\[
\sup_{1\leq t\leq t_{0}}B_{1,n}(t)\leq\gamma^{-\alpha}_{n}\sum^{t_{0}}_{j=0}\left\vert \phi_{j}\right\vert ^{\nu}\sum^{\infty}_{k=1}\left(\overline{a}_{k}\right)^{\alpha-\nu}H\left(\gamma_{n}\left\vert a_{k}\right\vert ^{-1}\right)\left\vert c_{k-j}\right\vert ^{\nu}
\]
\[
\leq_{_{\left(7\right)}}C\gamma^{-\alpha}_{n}H\left(\gamma_{n}\right)\sum^{t_{0}}_{j=0}\sum^{\infty}_{k=1}\left(\overline{a}_{k}\right)^{\alpha-\nu}\left\vert c_{k-j}\right\vert ^{\nu}\max\left\{ \left\vert a_{k}\right\vert ^{\delta},\left\vert a_{k}\right\vert ^{-\delta}\right\} 
\]
\[
=_{\left(8\right)}C\gamma^{-\alpha}_{n}H\left(\gamma_{n}\right)\sum^{t_{0}}_{j=0}\sum^{\infty}_{k=t_{0}+1}\left(\overline{a}_{k}\right)^{\alpha-\nu}\left\vert c_{k-j}\right\vert ^{\nu}\max\left\{ \left\vert a_{k}\right\vert ^{\delta},\left\vert a_{k}\right\vert ^{-\delta}\right\} +o(1)
\]
\[
\leq\left(1+\eta^{\delta}\right)C\gamma^{-\alpha}_{n}H\left(\gamma_{n}\right)\sum^{t_{0}}_{j=0}\sum^{\infty}_{k=t_{0}+1}\left(\overline{a}_{k}\right)^{\alpha-\nu}\left\vert c_{k-j}\right\vert ^{\nu}\max\left\{ k^{-\delta/\alpha},k^{\delta/\alpha}\right\} +o(1)
\]
\[
\leq\left(1+\eta^{\delta}\right)\left(1+\eta^{\alpha-\nu}\right)C\gamma^{-\alpha}_{n}H\left(\gamma_{n}\right)\sum^{t_{0}}_{j=0}\sum^{\infty}_{k=t_{0}+1}\left(k^{-\left(\alpha-\nu\right)/\alpha}\right)\left\vert c_{k-j}\right\vert ^{\nu}k^{\delta/\alpha}+o(1)
\]
\[
=C\gamma^{-\alpha}_{n}H\left(\gamma_{n}\right)\sum^{t_{0}}_{j=0}\sum^{\infty}_{k=t_{0}+1}\left(k^{\frac{-\left(\alpha-\nu\right)+\delta}{\alpha}}\right)\left\vert c_{k-j}\right\vert ^{\nu}+o(1)
\]
\[
=C\gamma^{-\alpha}_{n}H\left(\gamma_{n}\right)\sum^{t_{0}}_{j=0}\sum^{\infty}_{k=t_{0}+1-j}\left(k+j\right)^{\frac{-\left(\alpha-\nu\right)+\delta}{\alpha}}\left\vert c_{k}\right\vert ^{\nu}+o(1)
\]
\[
\leq C\gamma^{-\alpha}_{n}H\left(\gamma_{n}\right)\sum^{\infty}_{k=1}k^{\frac{-\left(\alpha-\nu\right)+\delta}{\alpha}}\left\vert c_{k}\right\vert ^{\nu}+o(1)
\]
\[
\leq C\gamma^{-\alpha}_{n}H\left(\gamma_{n}\right)+o(1)=_{\left(9\right)}o(1),
\]
where $\leq_{_{\left(7\right)}}$ follows from Potter's Theorem (see
e.g. BGT 1987; Thm 1.5.6), and $=_{\left(8\right)}$, $=_{\left(9\right)}$
because $\gamma^{\alpha}_{n}\sim L(n^{1/\alpha}\gamma_{n})$ and $H\left(\gamma_{n}\right)/L(n^{1/\alpha}\gamma_{n})\rightarrow0$.
To see the validity of the last statement recall that $\lim_{x\rightarrow\infty}L(x)/H(x)=\infty$
-cf. BGT 1987; Proposition 1.5.9a. In view of this as $n\rightarrow\infty$
\[
\frac{H\left(\gamma_{n}\right)}{L(n^{1/\alpha}\gamma_{n})}=\frac{L(\gamma_{n})}{L(n^{1/\alpha}\gamma_{n})}\frac{H\left(\gamma_{n}\right)}{L\left(\gamma_{n}\right)}\leq\frac{H\left(\gamma_{n}\right)}{L\left(\gamma_{n}\right)}\rightarrow0,
\]
where we have used the fact that $L(x)$ is increasing. In view of
the above, as $n\rightarrow\infty$ 
\begin{equation}
\sup_{1\leq t\leq t_{0}}B_{1,n}(t)=o(1).\label{a-moments2c}
\end{equation}

Next, we show that 
\begin{equation}
\sup_{t_{0}\leq t\leq n}B_{1,n}(t)=O(1).\label{a-moments2c0}
\end{equation}
 Set $R_{\nu-\alpha}\left(x\right):=x^{\nu-\alpha}H(x)$, $x>0$.
Note that $R_{\nu-\alpha}\left(x\right)$ is regularly varying of
index $\nu-\alpha<0$. Further, without loss of generality we assume
that $R_{\nu-\alpha}\left(x\right)$ is decreasing\footnote{Note that by BGT 1987 (Thm 1.5.3) regularly varying functions have
monotone asympotic equivalents, i.e. there a decreasing $\underline{R}_{\nu-\alpha}\left(x\right)$
such that $R_{\nu-\alpha}\left(x\right)\sim\underline{R}_{\nu-\alpha}\left(x\right),$
as $x\rightarrow\infty$.}. Recall that
\begin{equation}
B_{1,n}(t)\leq\gamma^{-\alpha}_{n}\sum^{t-1}_{j=0}\left\vert \phi_{j}\right\vert ^{\nu}\sum^{\infty}_{k=t}\left(\overline{a}_{k}\right)^{\alpha-\nu}H\left(\gamma_{n}\left\vert a_{k}\right\vert ^{-1}\right)\left\vert c_{k-j}\right\vert ^{\nu}=:B_{2,n}(t).\label{a-moments2ca}
\end{equation}
It is easy to see that by \textbf{Lemma \ref{cauchy}}, for $k$ large
enough $\overline{a}_{k}/\left\vert a_{k}\right\vert $ is bounded.
Hence, for $t\geq t_{0}$ as $n\rightarrow\infty$
\[
\gamma^{-\alpha}_{n}\sum^{t-1}_{j=0}\left\vert \phi_{j}\right\vert ^{\nu}\sum^{\infty}_{k=t}\left(\overline{a}_{k}\right)^{\alpha-\nu}\left[H\left(\gamma_{n}\left\vert a_{k}\right\vert ^{-1}\right)-H\left(\gamma_{n}k^{1/\alpha}\right)\right]\left\vert c_{k-j}\right\vert ^{\nu}
\]
\[
\leq\gamma^{-\alpha}_{n}\sum^{t-1}_{j=0}\left\vert \phi_{j}\right\vert ^{\nu}\sum^{\infty}_{k=t}\left\vert c_{l}\right\vert ^{\nu}\left\vert \frac{\left(\overline{a}^{-1}_{k}\right)^{-\left(\alpha-\nu\right)}H\left(\gamma_{n}\left\vert a_{k}\right\vert ^{-1}\right)}{\left[k^{1/\alpha}\right]^{-\left(\alpha-\nu\right)}H\left(\gamma_{n}k^{1/\alpha}\right)}-\frac{\left(\overline{a}^{-1}_{k}\right)^{-\left(\alpha-\nu\right)}}{\left(k^{1/\alpha}\right)^{-\left(\alpha-\nu\right)}}\right\vert \left[k^{1/\alpha}\right]^{-\left(\alpha-\nu\right)}H\left(\gamma_{n}k^{1/\alpha}\right)
\]
\[
=\gamma^{-\alpha}_{n}\sum^{t-1}_{j=0}\left\vert \phi_{j}\right\vert ^{\nu}\sum^{\infty}_{k=t}\left\vert c_{l}\right\vert ^{\nu}\frac{\left(\overline{a}^{-1}_{k}\right)^{-\left(\alpha-\nu\right)}}{\left(\left\vert a^{-1}_{k}\right\vert \right)^{-\left(\alpha-\nu\right)}}\left\vert \frac{\left(\left\vert a^{-1}_{k}\right\vert \right)^{-\left(\alpha-\nu\right)}h\left(\gamma_{n}\left\vert a_{k}\right\vert ^{-1}\right)}{\left[k^{1/\alpha}\right]^{-\left(\alpha-\nu\right)}h\left(\gamma_{n}k^{1/\alpha}\right)}-\frac{\left(\left\vert a^{-1}_{k}\right\vert \right)^{-\left(\alpha-\nu\right)}}{\left(k^{1/\alpha}\right)^{-\left(\alpha-\nu\right)}}\right\vert 
\]
\[
\cdot\left[k^{1/\alpha}\right]^{-\left(\alpha-\nu\right)}H\left(\gamma_{n}k^{1/\alpha}\right)
\]
\[
\leq C\gamma^{-\alpha}_{n}\sum^{t-1}_{j=0}\left\vert \phi_{j}\right\vert ^{\nu}\sum^{\infty}_{k=t}\left\vert c_{l}\right\vert ^{\nu}\left\vert \frac{\left(\left\vert a^{-1}_{k}\right\vert \right)^{-\left(\alpha-\nu\right)}H\left(\gamma_{n}\left\vert a_{k}\right\vert ^{-1}\right)}{\left[k^{1/\alpha}\right]^{-\left(\alpha-\nu\right)}H\left(\gamma_{n}k^{1/\alpha}\right)}-\frac{\left(\left\vert a^{-1}_{k}\right\vert \right)^{-\left(\alpha-\nu\right)}}{\left(k^{1/\alpha}\right)^{-\left(\alpha-\nu\right)}}\right\vert \left[k^{1/\alpha}\right]^{-\left(\alpha-\nu\right)}H\left(\gamma_{n}k^{1/\alpha}\right)
\]
\[
\leq C\sup_{k\geq t_{0}}\sup_{\left|y-1\right|<\eta}\left\vert \frac{\left(\gamma_{n}k^{1/\alpha}y\right)^{-\left(\alpha-\nu\right)}H\left(\gamma_{n}k^{1/\alpha}y\right)}{\left(\gamma_{n}k^{1/\alpha}\right)^{-\left(\alpha-\nu\right)}H\left(\gamma_{n}k^{1/\alpha}\right)}-y^{-\left(\alpha-\nu\right)}\right\vert \gamma^{-\alpha}_{n}\sum^{t-1}_{j=0}\left\vert \phi_{j}\right\vert ^{\nu}\sum^{\infty}_{k=t}\left\vert c_{l}\right\vert ^{\nu}\left[k^{1/\alpha}\right]^{-\left(\alpha-\nu\right)}H\left(\gamma_{n}k^{1/\alpha}\right)
\]
\begin{equation}
\leq_{\left(10\right)}o(1)\cdot\gamma^{-\alpha}_{n}\sum^{t-1}_{j=0}\left\vert \phi_{j}\right\vert ^{\nu}\sum^{\infty}_{k=t}\left\vert c_{l}\right\vert ^{\nu}\left[k^{1/\alpha}\right]^{-\left(\alpha-\nu\right)}H\left(\gamma_{n}k^{1/\alpha}\right),\label{a-moments3a}
\end{equation}
where $\leq_{\left(10\right)}$ follows from an almost identical argument
as the one used in (\ref{a-moments2}). Next, arguments similar to
those used in (\ref{a-moments2b}) show that for $t\geq t_{0}$,
\[
\gamma^{-\alpha}_{n}\sum^{t-1}_{j=0}\left\vert \phi_{j}\right\vert ^{\nu}\sum^{\infty}_{k=t}\left\vert \left(k^{-1/\alpha}\right)^{\alpha-\nu}-\left(\overline{a}_{k}\right)^{\alpha-\nu}\right\vert H\left(\gamma_{n}k^{1/\alpha}\right)\left\vert c_{k-j}\right\vert ^{\nu}
\]
\begin{equation}
\leq\eta\gamma^{-\alpha}_{n}\sum^{t-1}_{j=0}\left\vert \phi_{j}\right\vert ^{\nu}\sum^{\infty}_{k=t}\left(k^{-1/\alpha}\right)^{\alpha-\nu}H\left(\gamma_{n}k^{1/\alpha}\right)\left\vert c_{k-j}\right\vert ^{\nu}=:\eta B_{3,n}(t).\label{a-moments3b}
\end{equation}
Combining (\ref{a-moments3a})-(\ref{a-moments3b}) we have for $t\geq t_{0}$,
\begin{equation}
\left\vert B_{2,n}(t)-B_{3,n}(t)\right\vert \leq\eta B_{3,n}(t).\label{a-moments3c}
\end{equation}
In view of this and (\ref{a-moments2ca}), (\ref{a-moments2c0}) will
follow if we show that 
\[
\sup_{t_{0}\leq t\leq n}B_{3,n}(t)=O(1),
\]
which what we set out doing next.

In view of (\ref{CF5-FR2}) and the fact that $R_{\nu-\alpha}\left(x\right)$
is decreasing and $\phi_{j}\sim j^{-1/\alpha}$, it follows that,
for all $1\leq t\leq n$ as $n\rightarrow\infty$, 
\[
B_{3,n}(t)=\gamma^{-\alpha}_{n}\sum^{t-1}_{j=0}\left\vert \phi_{j}\right\vert ^{\nu}\sum^{\infty}_{l=t-j}\left(\left(l+j\right)^{-1/\alpha}\right)^{\alpha-\nu}H\left(\gamma^{1/\alpha}_{n}\left(l+j\right)^{1/\alpha}\right)\left\vert c_{l}\right\vert ^{\nu}
\]
\[
=\gamma^{-\nu}_{n}\sum^{t-1}_{j=0}\left\vert \phi_{j}\right\vert ^{\nu}\sum^{\infty}_{l=t-j}R_{\nu-\alpha}\left(\gamma_{n}\left(l+j\right)^{1/\alpha}\right)\left\vert c_{l}\right\vert ^{\nu}
\]
\begin{equation}
\leq C\gamma^{-\nu}_{n}\sum^{t-1}_{j=0}\left\vert \phi_{j}\right\vert ^{\nu}R_{\nu-\alpha}\left(\gamma_{n}j^{1/\alpha}\right)\sum^{\infty}_{l=t-j}\left\vert c_{l}\right\vert ^{\nu}\label{a-moments3d}
\end{equation}
\[
\leq C\gamma^{-\nu}_{n}\sum^{t-1}_{j=0}\left\vert \phi_{j}\right\vert ^{\nu}R_{\nu-\alpha}\left(\gamma_{n}j^{1/\alpha}\right)\sum^{\infty}_{l=0}\left\vert c_{l}\right\vert ^{\nu}
\]
\[
\leq C\gamma^{-\nu}_{n}\sum^{t-1}_{j=1}j^{-\nu/\alpha}R_{\nu-\alpha}\left(\gamma_{n}j^{1/\alpha}\right)\left[1+o(1)\right]
\]
\begin{equation}
=C\gamma^{-\alpha}_{n}\sum^{n}_{j=1}j^{-1}H\left(\gamma_{n}j^{1/\alpha}\right)\left[1+o(1)\right]=C\alpha+o(1),\label{a-moments4}
\end{equation}
as required.

$\bigskip$

\textbf{MI-SM}. Recall that
\[
a_{k,t}=\sum^{(t-1)\wedge k}_{j=0}\phi_{j}c_{k-j}=\begin{cases}
\sum^{k}_{j=0}\phi_{j}c_{k-j}=:a_{k} & \text{if }0\leq k\leq t-1,\\
\sum^{t-1}_{j=0}\phi_{j}c_{k-j}=:a^{-}_{k,t} & \text{if }k\geq t.
\end{cases}
\]
with
\[
\phi_{j}:=\phi_{j}(n)=\rho^{j}_{n}.
\]
Next, we introduce some notation. Let $\epsilon>0$ such that 
\[
\alpha-\epsilon>1,\:\textrm{if \:\ensuremath{\alpha\in}(1,2]},
\]
and
\[
\alpha-\epsilon\geq\nu,\:\textrm{if \:\ensuremath{\alpha\in}(0,1]},
\]
with $\nu$ given in \textbf{Assumption MI-SM}. In addition, set $\alpha^{*}:=\alpha-\epsilon$,
Therefore, by \textbf{Assumption MI-SM}
\[
\sum^{\infty}_{l=0}\left\vert c_{l}\right\vert ^{\alpha^{*}}<\infty.
\]

Choose $n$ large enough such that $0<\rho^{j}_{n}<1$ for all $j\geq1$.
Hence, for $n$ large enough,
\[
\left\vert a_{k}\right\vert \leq\sum^{k}_{j=0}\left\vert \phi_{j}c_{k-j}\right\vert \leq\sum^{\infty}_{j=0}\left\vert c_{j}\right\vert .
\]
which gives
\[
\inf_{k\geq1}\left\vert a^{-1}_{k}\right\vert \geq1/\sum^{\infty}_{j=0}\left\vert c_{j}\right\vert =:\eta>0,
\]
for $n$ large. Set 
\[
A_{1,n}:=\gamma^{-\alpha}_{n}\sum^{n}_{k=1}\left\vert a_{k}\right\vert ^{\alpha}H\left(\gamma_{n}\left\vert a_{k}\right\vert ^{-1}\right).
\]
First, we shall show that
\begin{equation}
\sup_{n\geq1,1\leq t\leq n}A_{n}(t):=\sup_{n\geq1,1\leq t\leq n}\mathbf{E}\left\vert \gamma^{-1}_{n}\sum^{t-1}_{k=0}a_{k}\xi_{t-k}\right\vert ^{\alpha-\varepsilon}<\infty,\label{alpha01}
\end{equation}
for all $\varepsilon>0$ such that $\alpha-\varepsilon>0$. Note that
(\ref{alpha01}) together with the Lyapunov inequality imply that
\[
\sup_{n\geq1,1\leq t\leq n}\mathbf{E}\left\vert \gamma^{-1}_{n}\sum^{t-1}_{k=0}a_{k}\xi_{t-k}\right\vert ^{\alpha^{\prime}}<\infty,
\]
for all $\alpha^{\prime}\leq\alpha-\varepsilon$.

In view of the fact that $\gamma_{n}\inf_{k\geq1}\left\vert a^{-1}_{k}\right\vert \rightarrow\infty$
and, using the same arguments as those used for establishing (\ref{a-moments1b}),
we get 
\begin{equation}
\sup_{1\leq t\leq n}A_{n}(t):=\sup_{1\leq t\leq n}\mathbf{E}\left\vert \gamma^{-1}_{n}\sum^{t-1}_{k=0}a_{k}\xi_{t-k}\right\vert ^{\alpha-\varepsilon}\leq C\left[A_{1,n}+A^{\frac{\alpha-\varepsilon}{\alpha+\varepsilon}}_{1,n}\right]\left[1+o(1)\right],\label{alpha001}
\end{equation}
for $n$ large enough. Next, we show that $A_{1,n}$ is bounded.

Set $\beta_{k}:=a^{-1}_{k}$. As $n\rightarrow\infty$, the difference
\[
\left\vert \sum^{n}_{k=1}\left\vert \gamma_{n}\beta_{k}\right\vert ^{-\alpha}H\left(\left\vert \gamma_{n}\beta_{k}\right\vert \right)-\sum^{n}_{k=1}\left\vert \gamma_{n}\beta_{k}\right\vert ^{-\alpha}H\left(\gamma_{n}\right)\right\vert 
\]
\[
\leq\sum^{n}_{k=1}\left\vert \frac{\left\vert \gamma_{n}\beta_{k}\right\vert ^{-\epsilon}H\left(\left\vert \gamma_{n}\beta_{k}\right\vert \right)}{\gamma^{-\epsilon}_{n}H\left(\gamma_{n}\right)}-\left\vert \beta_{k}\right\vert ^{-\epsilon}\right\vert H\left(\gamma_{n}\right)\gamma^{-\alpha}_{n}\left\vert \beta_{k}\right\vert ^{-\alpha+\epsilon}
\]
\[
\leq\sup_{1\leq k\leq n}\left\vert \frac{\left\vert \gamma_{n}\beta_{k}\right\vert ^{-\epsilon}H\left(\left\vert \gamma_{n}\beta_{k}\right\vert \right)}{\gamma^{-\epsilon}_{n}H\left(\gamma_{n}\right)}-\left\vert \beta_{k}\right\vert ^{-\epsilon}\right\vert H\left(\gamma_{n}\right)\sum^{n}_{k=1}\gamma^{-\alpha}_{n}\left\vert \beta_{k}\right\vert ^{-\alpha+\epsilon}
\]
\[
\leq\sup_{\eta\leq y<\infty}\left\vert \frac{\left(\gamma_{n}y\right)^{-\epsilon}H\left(\left\vert \gamma_{n}y\right\vert \right)}{\gamma^{-\epsilon}_{n}H\left(\gamma_{n}\right)}-y^{-\epsilon}\right\vert H\left(\gamma_{n}\right)\sum^{n}_{k=1}\gamma^{-\alpha}_{n}\left\vert \beta_{k}\right\vert ^{-\alpha+\epsilon}
\]
\[
=_{\left(1\right)}o(1)\cdot H\left(\gamma_{n}\right)\sum^{n}_{k=1}\gamma^{-\alpha}_{n}\left\vert \beta_{k}\right\vert ^{-\alpha+\epsilon}\sim_{\left(2\right)}o(1)\cdot\kappa^{-1}_{n}\sum^{n}_{k=1}\left\vert \beta_{k}\right\vert ^{-\alpha+\epsilon}=_{(3)}o(1),
\]

where $=_{\left(1\right)}$ follows from BGT 1987, Thm 1.5.2, $\sim_{\left(2\right)}$
from \textbf{Definition 1 (MI)} with $b=1$, and $=_{(3)}$ from the
fact that 
\begin{equation}
S_{n}:=\kappa^{-1}_{n}\sum^{n}_{k=1}\left\vert \beta_{k}\right\vert ^{-\alpha+\epsilon}=O(1).\label{MI0000a}
\end{equation}
A proof of (\ref{MI0000a}) is provided below. In view the above and
\textbf{Lemma \ref{LemmaMI-2}}, as $n\rightarrow\infty$ we get
\[
A_{1,n}=H\left(\gamma_{n}\right)\sum^{n}_{k=1}\left\vert \gamma^{-1}_{n}a_{k}\right\vert ^{\alpha}+o(1)
\]
\begin{equation}
=\kappa^{-1}_{n}\sum^{n}_{k=1}\left\vert a_{k}\right\vert ^{\alpha}[1+o(1)]+o(1)=\alpha^{-1}\left\vert \sum^{\infty}_{l=0}c_{l}\right\vert ^{\alpha}+o(1),\label{MI000a}
\end{equation}
where we have used the fact that $H\left(\gamma_{n}\right)/\gamma^{\alpha}_{n}\sim\kappa^{-1}_{n}$.
Note that (\ref{MI000a}) implies that $\sup_{1\leq t\leq n}A_{n}(t)=O(1)$. 

Next, we show that $S_{n}$ is bounded. First, suppose that $\alpha\in(0,1]$.
Recall that in this case $\alpha^{*}=\alpha-\epsilon\geq\nu$, which
implies that $\sum^{\infty}_{l=0}\left\vert c_{l}\right\vert ^{\alpha^{*}}<\infty.$
Using the fact that {[}\textbf{Check again: need to replace }$a_{k}=\sum^{k}_{l=1}\rho^{k-l}_{n}c_{l}$\textbf{,
by }$a_{k}=\sum^{k}_{l=0}\rho^{k-l}_{n}c_{l}$\textbf{{]}} $a_{k}=\sum^{k}_{l=1}\rho^{k-l}_{n}c_{l}$
we get
\[
S_{n}=\kappa^{-1}_{n}\sum^{n}_{k=1}\left\vert a_{k}\right\vert ^{\alpha^{*}}\leq\kappa^{-1}_{n}\sum^{n}_{k=1}\left(\sum^{k}_{l=0}\left\vert \rho^{k-l}_{n}c_{l}\right\vert \right)^{\alpha^{*}}
\]
\[
\leq\kappa^{-1}_{n}\sum^{n}_{k=1}\sum^{k}_{l=0}\left\vert \rho^{k-l}_{n}c_{l}\right\vert ^{\alpha^{*}}=:S^{^{\prime}}_{n}.
\]
Next, 
\[
S^{^{\prime}}_{n}=\kappa^{-1}_{n}\sum^{n}_{l=0}\sum^{n}_{k=l\lor1}\rho^{\alpha^{*}\left(k-l\right)}_{n}\left\vert c_{l}\right\vert ^{\alpha^{*}}\leq\kappa^{-1}_{n}\sum^{n}_{l=0}\rho^{-\alpha^{*}l}_{n}\left\vert c_{l}\right\vert ^{\alpha^{*}}\sum^{n}_{k=l}\rho^{\alpha^{*}k}_{n}
\]
\[
=\kappa^{-1}_{n}\sum^{n}_{l=0}\rho^{-\alpha^{*}l}_{n}\left\vert c_{l}\right\vert ^{\alpha^{*}}\left[\sum^{n}_{k=1}\rho^{\alpha^{*}k}_{n}-\sum^{l-1}_{k=1}\rho^{\alpha^{*}k}_{n}\right]
\]
\[
=\kappa^{-1}_{n}\sum^{n}_{l=0}\rho^{-\alpha^{*}l}_{n}\left\vert c_{l}\right\vert ^{\alpha^{*}}\left[\rho^{\alpha^{*}}_{n}\frac{1-\rho^{\alpha^{*}n}_{n}}{1-\rho^{\alpha^{*}}_{n}}-\rho^{\alpha^{*}}_{n}\frac{1-\rho^{\alpha^{*}\left(l-1\right)}_{n}}{1-\rho^{\alpha^{*}}_{n}}\right]
\]
\[
=\frac{\rho^{\alpha^{*}}_{n}}{1-\rho^{\alpha^{*}}_{n}}\kappa^{-1}_{n}\sum^{n}_{l=0}\rho^{-\alpha^{*}l}_{n}\left\vert c_{l}\right\vert ^{\alpha^{*}}\left[\rho^{\alpha^{*}\left(l-1\right)}_{n}-\rho^{\alpha^{*}n}_{n}\right]
\]
\[
=\frac{1}{1-\rho^{\alpha^{*}}_{n}}\kappa^{-1}_{n}\sum^{n}_{l=0}\left\vert c_{l}\right\vert ^{\alpha^{*}}-\frac{\rho^{\alpha^{*}}_{n}}{1-\rho^{\alpha^{*}}_{n}}\kappa^{-1}_{n}\sum^{n}_{l=0}\rho^{\alpha^{*}\left(n-l\right)}_{n}\left\vert c_{l}\right\vert ^{\alpha^{*}}
\]
\[
=\frac{1}{\alpha^{*}}\left[\sum^{n}_{l=0}\left\vert c_{l}\right\vert ^{\alpha^{*}}-\sum^{n}_{l=0}\rho^{\alpha^{*}\left(n-l\right)}_{n}\left\vert c_{l}\right\vert ^{\alpha^{*}}\right]\left[1+o(1)\right],
\]
\begin{equation}
\leq\frac{1}{\alpha^{*}}\left[\sum^{\infty}_{l=0}\left\vert c_{l}\right\vert ^{\alpha^{*}}+\sum^{n}_{l=0}\rho^{\alpha^{*}\left(n-l\right)}_{n}\left\vert c_{l}\right\vert ^{\alpha^{*}}\right]\left[1+o(1)\right]=O(1),\label{MI00a}
\end{equation}
where we have used the fact that $1/\left(1-\rho^{\alpha^{*}}_{n}\right)\kappa_{n}\rightarrow1/\alpha^{*}$.\footnote{By the mean value theorem, 
\[
f(x)=1+f^{\prime}(\bar{x})x,\text{ }\left\vert \bar{x}\right\vert \leq\left\vert x\right\vert 
\]
with $f(x)=(1-x)^{\alpha^{*}}\Rightarrow f^{\prime}(x)=-\alpha^{*}(1-x)^{\alpha^{*}-1}$.
Hence,
\[
1-f(x)=-f^{\prime}(\bar{x})x=\alpha^{*}(1-\bar{x})^{\alpha^{*}-1}x=\alpha^{*}\frac{1}{(1-\bar{x})^{1-\alpha^{*}}}x
\]
Therefore,
\[
\frac{x}{1-f(x)}=\frac{x(1-\bar{x})^{1-\alpha^{*}}}{\alpha^{*}x}\overset{\kappa^{-1}_{n}=x}{=}\frac{\kappa^{-1}_{n}(1-\bar{x})^{1-\alpha^{*}}}{\alpha^{*}\kappa^{-1}_{n}}\rightarrow\frac{1}{\alpha^{*}},
\]
where we have used the fact that $\bar{x}\leq\kappa^{-1}_{n}\rightarrow0$. 

{}}\textbf{\label{ft11}} It is straightforward to verify that the second
sum in (\ref{MI00a}) is bounded; however, for the subsequent analysis
we prove the stronger claim that it is asymptotically negligible.
Let $m_{n}$ an increasing integer valued sequence such that $m^{-1}_{n}+m_{n}\kappa^{-1}_{n}\rightarrow0$.
Noting that for $n$ large enough $0<\rho_{n}<1$, the term
\[
\sum^{n}_{l=0}\rho^{\alpha^{*}\left(n-l\right)}_{n}\left\vert c_{l}\right\vert ^{\alpha^{*}}=\left(\sum^{m_{n}}_{l=0}+\sum^{n}_{l=m_{n}+1}\right)\rho^{\alpha^{*}\left(n-l\right)}_{n}\left\vert c_{l}\right\vert ^{\alpha^{*}}
\]
\[
\leq\sum^{m_{n}}_{l=0}\rho^{\alpha^{*}\left(n-m_{n}\right)}_{n}\left\vert c_{l}\right\vert ^{\alpha^{*}}+\sum^{n}_{l=m_{n}+1}\rho^{\alpha^{*}\left(n-l\right)}_{n}\left\vert c_{l}\right\vert ^{\alpha^{*}}
\]
\begin{equation}
\leq\rho^{\alpha^{*}\left(n-m_{n}\right)}_{n}\sum^{\infty}_{l=0}\left\vert c_{l}\right\vert ^{\alpha^{*}}+\sum^{n}_{l=m_{n}+1}\left\vert c_{l}\right\vert ^{\alpha^{*}}\rightarrow0,\label{MI00b}
\end{equation}
as required. The above show that (\ref{MI0000a}) holds for $\alpha\in(0,1].$

Next, suppose that $\alpha\in(1,2]$. Recall that $\alpha^{*}=\alpha-\epsilon$.
Hence, in this case $\alpha^{*}>1$. Set $A_{k,n}:=\sum^{k}_{l=1}\left\vert \rho^{k-l}_{n}c_{l}\right\vert $,
$\overline{A}_{k,n}:=A_{k,n}1\left\{ A_{k,n}\geq1\right\} $ and $\underline{A}_{k,n}:=A_{k,n}1\left\{ A_{k,n}<1\right\} $.
Hence,
\[
S_{n}=\kappa^{-1}_{n}\sum^{n}_{k=1}\left\vert a_{k}\right\vert ^{\alpha^{*}}\leq\kappa^{-1}_{n}\sum^{n}_{k=1}\left(\sum^{k}_{l=0}\left\vert \rho^{k-l}_{n}c_{l}\right\vert \right)^{\alpha^{*}}.
\]
\[
=\kappa^{-1}_{n}\sum^{n}_{k=1}A^{\alpha^{*}}_{k,n}\leq2\kappa^{-1}_{n}\sum^{n}_{k=1}\left(\overline{A}^{\alpha^{*}}_{k,n}+\underline{A}^{\alpha^{*}}_{k,n}\right)\leq2\kappa^{-1}_{n}\sum^{n}_{k=1}\left(\overline{A}^{2}_{k,n}+\underline{A}_{k,n}\right)
\]
\begin{equation}
\leq2\kappa^{-1}_{n}\sum^{n}_{k=1}\left(A^{2}_{k,n}+A_{k,n}\right).\label{MI0a}
\end{equation}
Now, the arguments that lead to (\ref{MI00a}), also give $\kappa^{-1}_{n}\sum^{n}_{k=1}A_{k,n}=O(1)$.
Furthermore, $\kappa^{-1}_{n}\sum^{n}_{k=1}A^{2}_{k,n}=O(1)$ by DK2018
eq. (46), (47). Hence, (\ref{MI0000a}) also holds for $\alpha\in(1,2].$ 

The arguments shown above imply that for all $\alpha^{\prime}\in(0,\alpha)$
\[
\sup_{n\geq1,1\leq t\leq n}\mathbf{E}\left\vert \gamma^{-1}_{n}x^{+}_{t}(n)\right\vert ^{\alpha^{\prime}}<\infty.
\]

Next, we show that
\[
\sup_{n\geq1,1\leq t\leq n}\mathbf{E}\left\vert \gamma^{-1}_{n}x^{-}_{t}(n)\right\vert ^{\alpha^{\prime}}<\infty.
\]
Notice that (recall that $a^{-}_{k,t}=\sum^{t-1}_{l=1}\rho^{k-l}_{n}c_{l}$
$=\sum^{t}_{j=0}\rho^{j}_{n}c_{k-j}$, $k\geq t$) 
\[
\inf_{1\leq t\leq n,k\geq t}\left\vert a^{-}_{k,t}\right\vert ^{-1}=\frac{1}{\sup_{1\leq t\leq n,k\geq t}\left\vert a^{-}_{k,t}\right\vert }\geq\frac{1}{\sup_{1\leq t\leq n,k\geq t}\sum^{t}_{j=0}\left\vert \rho^{j}_{n}c_{k-j}\right\vert }
\]
\[
\geq1/\sum^{\infty}_{j=0}\left\vert c_{j}\right\vert =:\eta>0.
\]
Hence, $\gamma_{n}\inf_{1\leq t\leq n,k\geq t}\left\vert a^{-}_{k,t}\right\vert ^{-1}\rightarrow\infty$.
Consequently, by arguments analogous to those used to establish (\ref{a-moments1b}),
we obtain, as $n\rightarrow\infty$,
\begin{equation}
\sup_{1\leq t\leq n}\mathbf{E}\left\vert \gamma^{-1}_{n}x^{-}_{t}(n)\right\vert ^{\alpha-\varepsilon}\leq\left[\sup_{1\leq t\leq n}B_{n}(t)+\left(\sup_{1\leq t\leq n}B_{n}(t)\right)^{\frac{\alpha-\varepsilon}{\alpha+\varepsilon}}\right]\left[1+o(1)\right],\label{MI6}
\end{equation}
where
\[
B_{n}(t):=\sum^{\infty}_{k=t}\left\vert \frac{\lambda}{\gamma_{n}}a^{-}_{k,t}\right\vert ^{\alpha}H\left(\left\vert \gamma^{-1}_{n}a^{-}_{k,t}\right\vert ^{-1}\right).
\]
Therefore, it suffices to show that $\sup_{1\leq t\leq n}B_{n}(t)=O(1)$.
Set $\beta^{-1}_{k,t}:=a^{-}_{k,t}$. As $n\rightarrow\infty$,
\[
\left\vert B_{n}(t)-B_{n}(t)^{\prime}\right\vert 
\]
\[
:=\left\vert B_{n}(t)-\sum^{\infty}_{k=t}\left\vert \gamma^{-1}_{n}a^{-}_{k,t}\right\vert ^{\alpha}H\left(\gamma_{n}\right)\right\vert 
\]
\[
=\left\vert \sum^{\infty}_{k=t}\left\vert \gamma^{-1}_{n}a^{-}_{k,t}\right\vert ^{\alpha}H\left(\left\vert \gamma^{-1}_{n}a^{-}_{k,t}\right\vert ^{-1}\right)-\sum^{\infty}_{k=t}\left\vert \gamma^{-1}_{n}a^{-}_{k,t}\right\vert ^{\alpha}H\left(\gamma_{n}\right)\right\vert 
\]
\[
\leq\sum^{\infty}_{k=t}\left\vert \left\vert \gamma_{n}\beta_{k,t}\right\vert ^{-\alpha}H\left(\left\vert \gamma_{n}\beta_{k,t}\right\vert ^{-1}\right)-\left\vert \gamma_{n}\beta_{k,t}\right\vert ^{-\alpha}H\left(\gamma_{n}\right)\right\vert 
\]
\[
\leq\sum^{\infty}_{k=t}\left\vert \frac{\left\vert \gamma_{n}\beta_{k,t}\right\vert ^{-\epsilon}H\left(\left\vert \gamma_{n}\beta_{k,t}\right\vert \right)}{\left\vert \gamma_{n}\right\vert ^{-\epsilon}H\left(\gamma_{n}\right)}-\left\vert \beta_{k,t}\right\vert ^{-\epsilon}\right\vert \gamma^{-\alpha}_{n}\left\vert \beta_{k,t}\right\vert ^{-\alpha+\epsilon}H\left(\gamma_{n}\right)
\]
\[
\leq\sup_{\eta\leq y<\infty}\left\vert \frac{\left\vert \gamma_{n}y\right\vert ^{-\epsilon}H\left(\left\vert \gamma_{n}y\right\vert \right)}{\left\vert \gamma_{n}\right\vert ^{-\epsilon}H\left(\gamma_{n}\right)}-y^{-\epsilon}\right\vert \gamma^{-\alpha}_{n}\sum^{\infty}_{k=t}\left\vert \beta_{k,t}\right\vert ^{-\alpha+\epsilon}H\left(\gamma_{n}\right)
\]
\begin{equation}
\leq o(1)\cdot\gamma^{-\alpha}_{n}\sum^{\infty}_{k=t}\left\vert \beta_{k,t}\right\vert ^{-\alpha+\epsilon}H\left(\gamma_{n}\right)\sim o(1)\cdot\kappa^{-1}_{n}\sum^{\infty}_{k=t}\left\vert \beta_{k,t}\right\vert ^{-\alpha+\epsilon}=:o(1)\cdot T_{n}(t).\label{MI7}
\end{equation}
In view of (\ref{MI6})-(\ref{MI7}), Next, we show that $\sup_{1\leq t\leq n}T_{n}(t)=O(1)$.
First, suppose that $\alpha\in(0,1]$. We have (recall that $\alpha^{*}=\alpha-\epsilon$)
\[
T_{n}(t)=\kappa^{-1}_{n}\sum^{\infty}_{k=t}\left\vert a^{-}_{k,t}\right\vert ^{\alpha^{*}}\leq\kappa^{-1}_{n}\sum^{\infty}_{k=t}\sum^{t-1}_{j=0}\left\vert \rho^{j}_{n}c_{k-j}\right\vert ^{\alpha^{*}}
\]
\[
\leq\kappa^{-1}_{n}\sum^{\infty}_{k=0}\sum^{k\wedge n}_{j=0}\left\vert \rho^{j}_{n}c_{k-j}\right\vert ^{\alpha^{*}}=\kappa^{-1}_{n}\sum^{n-1}_{k=0}\sum^{k}_{j=0}\left\vert \rho^{j}_{n}c_{k-j}\right\vert ^{\alpha^{*}}+\kappa^{-1}_{n}\sum^{\infty}_{k=n}\sum^{n}_{j=0}\left\vert \rho^{j}_{n}c_{k-j}\right\vert ^{\alpha^{*}}
\]
\[
=S^{^{\prime}}_{n}+\kappa^{-1}_{n}\sum^{\infty}_{k=n}\left(\sum^{[n/2]}_{j=0}+\sum^{n}_{j=[n/2]+1}\right)\left\vert \rho^{j}_{n}c_{k-j}\right\vert ^{\alpha^{*}}
\]
\begin{equation}
\leq O(1)+\sum^{\infty}_{k=[n/2]}\left\vert c_{k}\right\vert ^{\alpha^{*}}\kappa^{-1}_{n}\sum^{[n/2]}_{j=0}\rho^{\alpha^{*}j}_{n}+\kappa^{-1}_{n}\sum^{n}_{j=[n/2]+1}\rho^{\alpha^{*}j}_{n}\sum^{\infty}_{k=0}\left\vert c_{k}\right\vert ^{\alpha^{*}}=O(1),\label{MI8}
\end{equation}
where we have used the facts that $S^{^{\prime}}_{n}=O(1)$, $\kappa^{-1}_{n}\sum^{[n/2]}_{j=0}\rho^{\alpha^{*}j}_{n}=O(1)$,
and $\kappa^{-1}_{n}\sum^{n}_{j=[n/2]+1}\rho^{\alpha^{*}j}_{n}=o(1)$
(for the latter two see Lemma 6.4 in DK2018).

Next, consider the case $\alpha\in(1,2]$. First, note that (\ref{MI8})
holds true if we set $\alpha^{*}=1$. Set $A_{tkn}:=\sum^{t-1}_{j=0}\left\vert \rho^{j}_{n}c_{k-j}\right\vert $
and $\overline{A}_{tkn}:=A_{tkn}1\left\{ A_{tkn}>1\right\} $, \underline{$A$}$_{tkn}=:A_{tkn}1\left\{ A_{tkn}\leq1\right\} $.
The term
\[
T_{n}(t)\leq\kappa^{-1}_{n}\sum^{\infty}_{k=t}\left(A_{tkn}\right)^{\alpha^{*}}=\kappa^{-1}_{n}\sum^{\infty}_{k=t}\left(\overline{A}_{tkn}+\underline{A}_{tkn}\right)^{\alpha^{*}}
\]
\[
\leq2\kappa^{-1}_{n}\sum^{\infty}_{k=t}\left[\left(\overline{A}_{tkn}\right)^{\alpha^{*}}+\left(\underline{A}_{tkn}\right)^{\alpha^{*}}\right]\leq2\kappa^{-1}_{n}\sum^{\infty}_{k=t}\left[\left(\overline{A}_{tkn}\right)^{2}+\underline{A}_{tkn}\right]
\]
\begin{equation}
\leq2\kappa^{-1}_{n}\sum^{\infty}_{k=t}\left[A^{2}_{tkn}+A_{tkn}\right].\label{MI8b}
\end{equation}
It follows from (\ref{MI8}), with $\alpha^{*}=1$, that $\kappa^{-1}_{n}\sup_{1\leq t\leq n}\sum^{\infty}_{k=t}A_{tkn}=O(1)$.
Furthermore, $\kappa^{-1}_{n}\sup_{1\leq t\leq n}\sum^{\infty}_{k=t}A^{2}_{tkn}=O(1)$
by DK2018 (proof of Proposition 2.2).

Combining (\ref{MI7})-(\ref{MI8b}) we have
\begin{equation}
\sup_{1\leq t\leq n}B_{n}(t)=\sup_{1\leq t\leq n}B_{n}(t)^{\prime}+o(1)=O(1),\label{MI8c}
\end{equation}
where the last equality follows from identical arguments as those
used in (\ref{MI8}) and (\ref{MI8b}) with $\alpha^{*}$ replaced
by $\alpha$.

\bigskip{}

\textbf{MI-LM/NM. }Next, we show the finite sub-$\alpha$ moments
property under \textbf{MI-LM} and \textbf{MI-NM}. We start with some
notation. Set
\begin{equation}
m_{n}:=\gamma_{n}\kappa^{-\left(1-m\right)}_{n}\left(\Longleftrightarrow m_{n}/\gamma_{n}=\kappa^{-1+m}_{n}\right)\label{MI8d}
\end{equation}
In addition, recall that (\textbf{Definition 1})
\[
\gamma_{n}=\ell_{1}\left(\kappa_{n}\right)^{1/\alpha}\kappa^{b/\alpha}_{n}\:\textrm{where}\:\ell_{1}(x)\sim H\left(x^{1/\alpha}\ell_{1}(x)^{1/\alpha}\right),
\]
and
\[
a_{k,t}=\sum^{(t-1)\wedge k}_{j=0}\phi_{j}c_{k-j}=\begin{cases}
\sum^{k}_{j=0}\phi_{j}c_{k-j}=:a_{k} & \text{if }0\leq k\leq t-1,\\
\sum^{t-1}_{j=0}\phi_{j}c_{k-j}=:a^{-}_{k,t} & \text{if }k\geq t.
\end{cases}
\]
with
\[
\phi_{j}:=\phi_{j}(n)=\rho^{j}_{n}.
\]

\begin{remark} \textbf{\label{RmkMI-LM/NM(1)}}
\begin{itemize}
\item[(a)]  Note that the above imply that $m_{n}\sim\ell_{1}\left(\kappa_{n}\right)^{1/\alpha}\kappa^{1/\alpha}_{n}$.
To see this write
\end{itemize}
\[
m_{n}=\gamma_{n}\kappa^{-\left(1-m\right)}_{n}\sim\ell_{1}\left(\kappa_{n}\right)^{1/\alpha}\kappa^{b/\alpha}_{n}\kappa^{-\left(1-m\right)}_{n}
\]
\[
=\ell_{1}\left(\kappa_{n}\right)^{1/\alpha}\kappa^{\frac{1}{\alpha}\left(1+\alpha-\alpha m\right)}_{n}\kappa^{-\left(1-m\right)}_{n}=\ell_{1}\left(\kappa_{n}\right)^{1/\alpha}\kappa^{\left(\frac{1}{\alpha}+1-m\right)}_{n}\kappa^{-\left(1-m\right)}_{n}
\]
\[
=\ell_{1}\left(\kappa_{n}\right)^{1/\alpha}\kappa^{1/\alpha}_{n}.
\]

\begin{itemize}
\item[(b)]  In addition, we have the following relationship that will be utilized
later. As $n\rightarrow\infty$,
\[
H\left(m_{n}\right)\sim H\left(\kappa^{1/\alpha}_{n}\ell_{1}\left(\kappa_{n}\right)^{1/\alpha}\right)\sim\ell_{1}\left(\kappa_{n}\right).
\]
\end{itemize}
\end{remark} 

\bigskip{}

\textbf{MI-LM ($m\in(1/\alpha,1)$). }Set
\[
q(x):=e^{-x}\int^{x}_{0}e^{y}y^{-m}dy.
\]
By Lemma 7.2 of DK2018, 
\begin{equation}
\sup_{k_{0}\leq k\leq n}\left\vert \kappa^{-\left(1-m\right)}_{n}\sum^{k}_{j=0}\rho^{j}_{n}c_{k-j}-q(k/\kappa_{n})\right\vert =o(1),\label{FSM-LM00}
\end{equation}
as $n\rightarrow\infty$. Hence,
\begin{equation}
\sup_{0\leq k\leq n}\kappa^{-\left(1-m\right)}_{n}\left\vert \sum^{k}_{j=0}\rho^{j}_{n}c_{k-j}\right\vert =\sup_{0\leq k\leq n}q(k/\kappa_{n})+o(1).\label{FSM-LM000}
\end{equation}
In view of the fact\footnote{Notice that 
\[
0\leq q(x)=\frac{\int^{x}_{0}e^{y}y^{-m}dy}{\int^{x}_{0}e^{y}dy+1}=\frac{\int^{1}_{0}e^{y}y^{-m}dy+\int^{x}_{1}e^{y}y^{-m}dy}{\int^{x}_{0}e^{y}dy+1}\leq\frac{e\int^{1}_{0}y^{-m}dy}{\int^{x}_{0}e^{y}dy+1}+1<\infty.
\]
} that $\sup_{x\geq0}q(x)\leq C_{0}$, for some $C_{0}\in(0,\infty)$,
we get
\[
\sup_{0\leq k\leq n}q(k/\kappa_{n})\leq C_{0}.
\]
Next, note that (see \textbf{Remark \ref{RmkMI-LM/NM(1)}}) $m_{n}=\gamma_{n}\kappa^{-\left(1-m\right)}_{n}\sim\kappa^{1/\alpha}_{n}\ell_{1}\left(\kappa_{n}\right)^{1/\alpha}\rightarrow\infty$.
Set $S_{n,k}:=\gamma_{n}\left\vert a_{k}\right\vert ^{-1}$, and $S^{\ast}_{n,k}:=m^{-1}_{n}S_{n,k}$.
In view of the above, 
\[
\inf_{1\leq k\leq n}S_{n,k}=m_{n}\inf_{1\leq k\leq n}\left(\kappa^{1-m}_{n}\gamma^{-1}_{n}S_{n,k}\right)=m_{n}\inf_{1\leq k\leq n}S^{\ast}_{n,k}
\]
\begin{equation}
=m_{n}\frac{1}{\kappa^{-\left(1-m\right)}_{n}\sup_{1\leq k\leq n}\left\vert a_{k}\right\vert }\geq m_{n}\left(2C_{0}\right)^{-1}\rightarrow\infty,\label{FSM-LM0}
\end{equation}
as $n\rightarrow\infty$.

Set
\begin{equation}
A(n):=\frac{1}{\gamma^{\alpha}_{n}}\sum^{n-1}_{k=1}\left\vert a_{k}\right\vert ^{\alpha}H\left(\left\vert \gamma_{n}a^{-1}_{k}\right\vert \right)=\sum^{n-1}_{k=1}\left\vert S_{n,k}\right\vert ^{-\alpha}H\left(\left\vert S_{n,k}\right\vert \right).\label{FSM-LM001}
\end{equation}

In view of (\ref{FSM-LM0}), and using the same arguments as those
used for establishing (\ref{alpha001}) we get 
\[
\sup_{1\leq t\leq n}A_{n,t}:=\sup_{1\leq t\leq n}\mathbf{E}\left\vert \gamma^{-1}_{n}\sum^{t-1}_{k=0}a_{k}\xi_{t-k}\right\vert ^{\alpha-\varepsilon}\leq C\left[A(n)+A(n)^{\frac{\alpha-\varepsilon}{\alpha+\varepsilon}}\right]\left[1+o(1)\right],
\]
for all $\varepsilon>0$ such that $\alpha-\varepsilon>0$ and $n$
large enough. Hence, it suffices to show that $A(n)=O(1)$, which
is what we set out to do next.

Let $\epsilon>0$ and $b_{\epsilon}:=1+\left(\alpha-\epsilon\right)-\left(\alpha-\epsilon\right)m$
such that (\ref{epsi3}) is satisfied. For $n$ large enough the difference
\[
\left\vert A(n)-\sum^{n-1}_{k=1}\left\vert S_{n,k}\right\vert ^{-\alpha}H\left(m_{n}\right)\right\vert =\left\vert \sum^{n-1}_{k=1}\left\vert S_{n,k}\right\vert ^{-\alpha}H\left(S_{n,k}\right)-\sum^{n-1}_{k=1}\left\vert S_{n,k}\right\vert ^{-\alpha}H\left(m_{n}\right)\right\vert 
\]
\[
\leq\sum^{n-1}_{k=1}\left\{ \left\vert \frac{\left\vert m_{n}\kappa^{1-m}_{n}\gamma^{-1}_{n}S_{n,k}\right\vert ^{-\epsilon}H\left(\left\vert m_{n}\kappa^{1-m}_{n}\gamma^{-1}_{n}S_{n,k}\right\vert \right)}{m^{-\epsilon}_{n}H\left(m_{n}\right)}-\left\vert \kappa^{1-m}_{n}\gamma^{-1}_{n}S_{n,k}\right\vert ^{-\epsilon}\right\vert \right.
\]
\[
\cdot\left.\kappa^{\epsilon\left(1-m\right)}_{n}H\left(m_{n}\right)\gamma^{-\alpha}_{n}\left\vert \gamma^{-1}_{n}S_{n,k}\right\vert ^{-\alpha+\epsilon}\right\} 
\]
\[
\leq\sum^{n-1}_{k=1}\left\vert \frac{\left\vert m_{n}S^{\ast}_{n,k}\right\vert ^{-\epsilon}H\left(\left\vert m_{n}S^{\ast}_{n,k}\right\vert \right)}{m^{-\epsilon}_{n}H\left(m_{n}\right)}-\left\vert S^{\ast}_{n,k}\right\vert ^{-\epsilon}\right\vert \kappa^{\epsilon\left(1-m\right)}_{n}H\left(m_{n}\right)\gamma^{-\alpha}_{n}\left\vert \gamma^{-1}_{n}S_{n,k}\right\vert ^{-\alpha+\epsilon}
\]
\[
\leq\sup_{\left(2C_{0}\right)^{-1}\leq y<\infty}\left\vert \frac{\left(m_{n}y\right)^{-\epsilon}H\left(\left\vert m_{n}y\right\vert \right)}{m^{-\epsilon}_{n}H\left(m_{n}\right)}-y^{-\epsilon}\right\vert \kappa^{\epsilon\left(1-m\right)}_{n}H\left(m_{n}\right)\sum^{n-1}_{k=1}\gamma^{-\alpha}_{n}\left\vert a^{-1}_{k}\right\vert ^{-\alpha+\epsilon}
\]
\[
=_{\left(1\right)}o(1)\cdot H\left(m_{n}\right)\kappa^{\epsilon\left(1-m\right)}_{n}\sum^{n-1}_{k=1}\gamma^{-\alpha}_{n}\left\vert a_{k}\right\vert ^{\alpha-\epsilon}
\]
\begin{equation}
\sim o(1)\cdot\frac{H\left(m_{n}\right)}{H\left(m_{n}\right)}\frac{1}{\kappa^{b_{\epsilon}}_{n}}\sum^{n-1}_{k=1}\left\vert a_{k}\right\vert ^{\alpha-\epsilon}=_{(2)}o(1),\label{FSM-LM01}
\end{equation}
where $=_{\left(1\right)}$ follows from BGT 1987, Thm 1.5.2. In addition,
for $=_{(2)}$ we have used the fact that (cf. \textbf{Remark \ref{RmkMI-LM/NM(1)}})
\[
\kappa^{\epsilon\left(1-m\right)}_{n}\gamma^{-\alpha}_{n}\sim\kappa^{\epsilon\left(1-m\right)}_{n}\kappa^{-b}_{n}H\left(m_{n}\right)^{-1}=\kappa^{\epsilon\left(1-m\right)}_{n}\kappa^{-\left(1+\alpha-\alpha m\right)}_{n}H\left(m_{n}\right)^{-1}
\]
\[
=\kappa^{-\left(1+\left(\alpha-\epsilon\right)-\left(\alpha-\epsilon\right)m\right)}_{n}H\left(m_{n}\right)^{-1}=\kappa^{-b_{\epsilon}}_{n}H\left(m_{n}\right)^{-1},
\]
and
\[
\frac{1}{\kappa^{b_{\epsilon}}_{n}}\sum^{n-1}_{k=1}\left\vert a_{k}\right\vert ^{\alpha-\epsilon}=\left\vert m_{0}\right\vert ^{\alpha-\epsilon}\mathbb{V}_{1}\left(\alpha-\epsilon\right)+o(1),
\]
with $\mathbb{V}_{1}$ given in\textbf{ Lemma \ref{LemmaMI3}} (see
\textbf{Lemma \ref{LemmaMI1} }and\textbf{ Lemma \ref{LemmaMI3}}
for more details). In view the above, as $n\rightarrow\infty$
\[
A(n)=H\left(m_{n}\right)\gamma^{-\alpha}_{n}\sum^{n}_{k=1}\left\vert a_{k}\right\vert ^{\alpha}+o(1)
\]
\begin{equation}
\sim\frac{H\left(m_{n}\right)}{H\left(m_{n}\right)}\frac{1}{\kappa^{b}_{n}}\sum^{n}_{k=1}\left\vert a_{k}\right\vert ^{\alpha}=\left\vert m_{0}\right\vert ^{\alpha}\mathbb{V}_{1}(\alpha)+o(1).\label{FSM-LM02}
\end{equation}

Let
\begin{equation}
B_{n,t}:=\mathbf{E}\left\vert \gamma^{-1}_{n}\sum^{\infty}_{k=t}a^{-}_{k,t}\xi_{t-k}\right\vert ^{\alpha-\varepsilon}=\mathbf{E}\left\vert \gamma^{-1}_{n}\sum^{\infty}_{k=t}a^{-}_{k,t}\xi_{k}\right\vert ^{\alpha-\varepsilon},\label{FSM-LM02b}
\end{equation}
for all $\varepsilon>0$ such that $\alpha-\varepsilon>0$. We shall
show that $\sup_{1\leq t\leq n}B_{n,t}=O(1)$, which will conclude
the proof this Lemma under \textbf{Assumption MI-LM}. Set
\[
S_{n,n\wedge k}:=\left(\sum^{n\wedge k}_{j=0}\left\vert \rho^{j}_{n}c_{k-j}\right\vert \right)^{-1}.
\]
The following calculations are similar to those that lead to (\ref{FSM-LM0}).
First, notice that
\[
\kappa^{m-1}_{n}\sup_{k\geq0}\sum^{n\wedge k}_{j=0}\left\vert \rho^{j}_{n}c_{k-j}\right\vert \leq\kappa^{m-1}_{n}\sup_{0\leq k\leq n}\sum^{k}_{j=0}\left\vert \rho^{j}_{n}c_{k-j}\right\vert +\kappa^{m-1}_{n}\sup_{k>n}\sum^{n}_{j=0}\left\vert \rho^{j}_{n}c_{k-j}\right\vert .
\]
By (\ref{FSM-LM00}), the first term above satisfies (\ref{FSM-LM000}).
Let $\delta_{n}>0$ be an integer valued sequence such that $\delta^{-1}_{n}\kappa_{n}\ln(\kappa^{m}_{n})=o(1)$.
This implies that $\kappa^{m}_{n}\rho^{\delta_{n}}_{n}=o(1)$ (cf.
DK2018, Lemma 6.4). In view if this, the second term

\[
\kappa^{m-1}_{n}\sup_{k>n}\sum^{n}_{j=0}\left\vert \rho^{j}_{n}c_{k-j}\right\vert \leq\kappa^{m-1}_{n}\left[\sup_{k>n}\sum^{\delta_{n}}_{j=0}\left\vert \rho^{j}_{n}c_{k-j}\right\vert +\sup_{k>n}\sum^{n}_{j=\delta_{n}}\left\vert \rho^{j}_{n}c_{k-j}\right\vert \right]
\]
\[
=\kappa^{m-1}_{n}\left[\sup_{k>n}\sum^{k}_{j=k-\delta_{n}}\left\vert \rho^{k-j}_{n}c_{j}\right\vert +\sup_{k\geq0}\left\vert c_{k}\right\vert \sum^{n}_{j=\delta_{n}}\left\vert \rho^{j}_{n}\right\vert \right]
\]
\[
\leq\kappa^{m-1}_{n}\left[\sup_{j>n-\delta_{n}}\left\vert c_{j}\right\vert \sup_{k>n}\sum^{k}_{j=k-\delta_{n}}\left\vert \rho^{k-j}_{n}\right\vert +\sup_{k\geq0}\left\vert c_{k}\right\vert \sum^{n}_{j=\delta_{n}}\left\vert \rho^{j}_{n}\right\vert \right]
\]
\[
=\kappa^{m-1}_{n}\left[\sup_{j>n-\delta_{n}}\left\vert c_{j}\right\vert \sum^{\delta_{n}}_{j=0}\left\vert \rho^{j}_{n}\right\vert +C\sum^{n}_{j=\delta_{n}}\left\vert \rho^{j}_{n}\right\vert \right]=\kappa^{m-1}_{n}\left[\sup_{j>n-\delta_{n}}\left\vert c_{j}\right\vert O\left(\kappa_{n}\right)+O\left(\kappa_{n}\rho^{\delta_{n}}_{n}\right)\right]
\]
\[
=\sup_{j>n-\delta_{n}}\left\vert c_{j}\right\vert O\left(\kappa^{m}_{n}\right)+O\left(\kappa^{m}_{n}\rho^{\delta_{n}}_{n}\right)\leq C\left(n-\delta_{n}\right)^{-m}\kappa^{m}_{n}+o\left(1\right)=o(1).
\]

Hence, similarly to (\ref{FSM-LM0}) as $n\rightarrow\infty$ 
\begin{equation}
\gamma_{n}\inf_{k\geq0}S_{n,n\wedge k}=\frac{\gamma_{n}}{\sup_{k\geq0}\sum^{n\wedge k}_{j=0}\left\vert \rho^{j}_{n}c_{k-j}\right\vert }\geq\left(2C_{0}\right)^{-1}m_{n}\rightarrow\infty,\label{FSM-LM1}
\end{equation}
for some $C_{0}<\infty$. Next, set
\[
\overline{\xi}_{n,k}:=\xi_{k}1\left\{ \left\vert \xi_{k}\right\vert \geq\gamma_{n}S_{n,n\wedge k}\right\} \text{ \ and \ }\underline{\xi}_{n,k}:=\xi_{k}1\left\{ \left\vert \xi_{k}\right\vert <\gamma_{n}S_{n,n\wedge k}\right\} .
\]
As before, without loss of generality we will use the convention that
$\mathbf{E}\overline{\xi}_{n,k}=\mathbf{E}\underline{\xi}_{n,k}=0$
for $1<\alpha\leq2$. It follows from (\ref{FSM-LM1}) and Astrauskas
(1983; Lemma 1) that as $n\rightarrow\infty$
\begin{equation}
\mathbf{E}\left\vert \overline{\xi}_{n,k}\right\vert ^{\alpha-\varepsilon}\asymp S^{-\varepsilon}_{n,k}h(S_{n,n\wedge k})\text{ \ and \ }\mathbf{E}\left\vert \underline{\xi}_{n,k}\right\vert ^{\alpha+\varepsilon}\asymp S^{\varepsilon}_{n,k}h(S_{n,n\wedge k}),\label{FSM-LM1b}
\end{equation}
uniformly in $k\geq0$. Similar arguments as those used before, see
e.g. (\ref{a-moments1b}), give
\[
B_{n,t}\leq C\left[B_{1n,t}+B^{\frac{\alpha-\varepsilon}{\alpha+\varepsilon}}_{2n,t}\right],
\]
with 
\[
B_{1n,t}:=\gamma^{-\left(\alpha-\varepsilon\right)}_{n}\sum^{\infty}_{k=t}\left\vert a^{-}_{k,t}\right\vert ^{\alpha-\varepsilon}\mathbf{E}\left\vert \overline{\xi}_{n,k}\right\vert ^{\alpha-\varepsilon},
\]
and
\[
B_{2n,t}:=\gamma^{-\left(\alpha+\varepsilon\right)}_{n}\sum^{\infty}_{k=t}\left\vert a^{-}_{k,t}\right\vert ^{\alpha+\varepsilon}\mathbf{E}\left\vert \underline{\xi}_{n,k}\right\vert ^{\alpha+\varepsilon}.
\]
Using (\ref{a-moments1ab}), (\ref{FSM-LM1}), and (\ref{FSM-LM1b})
we get
\[
\sup_{1\leq t\leq n}B_{1n,t}\leq C\gamma^{-\left(\alpha-\varepsilon\right)}_{n}\sup_{1\leq t\leq n}\sum^{\infty}_{k=t}\left\vert a^{-}_{k,t}\right\vert ^{\alpha-\varepsilon}S^{-\varepsilon}_{n,n\wedge k}H(S_{n,n\wedge k})\left[1+o(1)\right]
\]
\[
=C\gamma^{-\left(\alpha-\varepsilon\right)}_{n}\sup_{1\leq t\leq n}\sum^{\infty}_{k=t}\left\vert \sum^{t-1}_{j=0}\rho^{j}_{n}c_{k-j}\right\vert ^{\alpha-\varepsilon}S^{-\varepsilon}_{n,n\wedge k}H(S_{n,n\wedge k})\left[1+o(1)\right]
\]
\[
\leq C\gamma^{-\left(\alpha-\varepsilon\right)}_{n}\sup_{1\leq t\leq n}\sum^{\infty}_{k=t}\left(\sum^{t-1}_{j=0}\left\vert \rho^{j}_{n}c_{k-j}\right\vert \right)^{\alpha-\varepsilon}S^{-\varepsilon}_{n,n\wedge k}H(S_{n,n\wedge k})\left[1+o(1)\right]
\]
\[
=C\sup_{1\leq t\leq n}\sum^{\infty}_{k=t}\left(\gamma^{-1}_{n}\sum^{t-1}_{j=0}\left\vert \rho^{j}_{n}c_{k-j}\right\vert \right)^{\alpha-\varepsilon}S^{-\varepsilon}_{n,n\wedge k}H(S_{n,n\wedge k})\left[1+o(1)\right]
\]
\[
\leq C\sup_{1\leq t\leq n}\sum^{\infty}_{k=t}S^{-\alpha+\varepsilon}_{n,n\wedge k}S^{-\varepsilon}_{n,n\wedge k}H(S_{n,n\wedge k})\left[1+o(1)\right]
\]
\begin{equation}
=C\sup_{1\leq t\leq n}\sum^{\infty}_{k=t}S^{-\alpha}_{n,n\wedge k}H(S_{n,n\wedge k})\left[1+o(1)\right]\leq C\sum^{\infty}_{k=0}S^{-\alpha}_{n,n\wedge k}H(S_{n,n\wedge k})\left[1+o(1)\right].\label{FSM-LM2}
\end{equation}
The following arguments are similar to those that lead to (\ref{FSM-LM01})
above. Let $\epsilon$ and $b_{\epsilon}$ as in (\ref{FSM-LM01}).
Then as $n\rightarrow\infty$ the difference 
\[
\left\vert \sum^{\infty}_{k=0}S^{-\alpha}_{n,n\wedge k}H(S_{n,n\wedge k})-\sum^{\infty}_{k=0}S^{-\alpha}_{n,n\wedge k}H\left(m_{n}\right)\right\vert 
\]
\[
\leq\sum^{\infty}_{k=0}\left\{ \left\vert \frac{\left(m_{n}\kappa^{1-m}_{n}\gamma^{-1}_{n}S_{n,n\wedge k}\right)^{-\epsilon}H(m_{n}\kappa^{1-m}_{n}\gamma^{-1}_{n}S_{n,n\wedge k})}{m^{-\epsilon}_{n}H(m_{n})}-\left(\kappa^{1-m}_{n}\gamma^{-1}_{n}S_{n,n\wedge k}\right)^{-\epsilon}\right\vert \right.
\]

\[
\left.\kappa^{\epsilon\left(1-m\right)}_{n}\left(\gamma^{-1}_{n}S_{n,n\wedge k}\right)^{-\alpha+\epsilon}\gamma^{-\alpha}_{n}H\left(m_{n}\right)\right\} 
\]
\[
\leq\sup_{\left(2C_{0}\right)^{-1}\leq y<\infty}\left\vert \frac{\left(m_{n}y\right)^{-\epsilon}H(m_{n}y)}{m^{-\epsilon}_{n}H(m_{n})}-y^{-\epsilon}\right\vert \sum^{\infty}_{k=0}\kappa^{\epsilon\left(1-m\right)}_{n}\left(\gamma^{-1}_{n}S_{n,n\wedge k}\right)^{-\alpha+\epsilon}\gamma^{-\alpha}_{n}H\left(m_{n}\right)
\]

\[
=\sup_{\left(2C_{0}\right)^{-1}\leq y<\infty}\left\vert \frac{\left(m_{n}y\right)^{-\epsilon}H(m_{n}y)}{m^{-\epsilon}_{n}H(m_{n})}-y^{-\epsilon}\right\vert \frac{H\left(m_{n}\right)}{\kappa^{b_{\epsilon}}_{n}H\left(m_{n}\right)}\sum^{\infty}_{k=0}\left(\sum^{n\wedge k}_{j=0}\left\vert \rho^{j}_{n}c_{k-j}\right\vert \right)^{-\alpha+\epsilon}\left[1+o(1)\right]
\]
\begin{equation}
=o(1)\cdot\frac{1}{\kappa^{b_{\epsilon}}_{n}}\sum^{\infty}_{k=0}\left(\sum^{n\wedge k}_{j=0}\left\vert \rho^{j}_{n}c_{k-j}\right\vert \right)^{\alpha-\epsilon}\left[1+o(1)\right]=_{(3)}o(1),\label{FSM-LM3}
\end{equation}
where $=_{(3)}$ is due to the fact that 
\[
\frac{1}{\kappa^{b_{\epsilon}}_{n}}\sum^{\infty}_{k=0}\left(\sum^{n\wedge k}_{j=0}\left\vert \rho^{j}_{n}c_{k-j}\right\vert \right)^{\alpha-\epsilon}\leq\frac{1}{\kappa^{b_{\epsilon}}_{n}}\sum^{n}_{k=0}\left(\sum^{k}_{j=0}\left\vert \rho^{j}_{n}c_{k-j}\right\vert \right)^{\alpha-\epsilon}
\]
\begin{equation}
+\frac{1}{\kappa^{b_{\epsilon}}_{n}}\sum^{\infty}_{k=n}\left(\sum^{n}_{j=0}\left\vert \rho^{j}_{n}c_{k-j}\right\vert \right)^{\alpha-\epsilon}=_{(4)}\mathbb{V}_{1}\left(\alpha-\epsilon\right)+o(1)=O(1).\label{FSM-LM4}
\end{equation}
We note that $=_{(4)}$ is a direct consequence of \textbf{Lemma \ref{LemmaMI1}
}and\textbf{ Lemma \ref{LemmaMI3}}. Combining (\ref{FSM-LM2})-(\ref{FSM-LM4})
we have 
\[
\sup_{1\leq t\leq n}B_{1n,t}\leq C\sum^{\infty}_{k=0}S^{-\alpha}_{n,n\wedge k}H(S_{n,n\wedge k})\left[1+o(1)\right]
\]
 as
\[
=C\sum^{\infty}_{k=0}S^{-\alpha}_{n,n\wedge k}H\left(m_{n}\right)\left[1+o(1)\right]\leq\frac{1}{\kappa^{b}_{n}}\sum^{n}_{k=0}\left(\sum^{k}_{j=0}\left\vert \rho^{j}_{n}c_{k-j}\right\vert \right)^{\alpha}
\]
\begin{equation}
+\frac{1}{\kappa^{b}_{n}}\sum^{\infty}_{k=n}\left(\sum^{n}_{j=0}\left\vert \rho^{j}_{n}c_{k-j}\right\vert \right)^{\alpha}=_{(5)}\left\vert m_{0}\right\vert ^{\alpha}\mathbb{V}_{1}\left(\alpha\right)+o(1)=O(1),\label{FSM-LM5}
\end{equation}
where $=_{(5)}$ is again a consequence of \textbf{Lemma \ref{LemmaMI1}
}and\textbf{ Lemma \ref{LemmaMI3}}. 

Finally, identical arguments to those used above show that 
\[
\sup_{1\leq t\leq n}B_{2n,t}=O(1).
\]
Hence, the proof is complete. $\blacksquare$\bigskip{}

\textbf{MI-NM ($m\in(1,1+1/\alpha)$).} Noting that $0<\rho_{n}<1$
for all $n$ large enough we have as $n\rightarrow\infty$
\begin{equation}
\gamma_{n}\inf_{1\leq k\leq n}\left\vert a_{k}\right\vert ^{-1}=\gamma_{n}\frac{1}{\sup_{1\leq k\leq n}\left\vert a_{k}\right\vert }\geq\frac{\gamma_{n}}{\sum^{\infty}_{j=0}\left\vert c_{j}\right\vert }\rightarrow\infty.\label{MI-NM1}
\end{equation}

As in the previous part, set
\[
A(n):=\frac{1}{\gamma^{\alpha}_{n}}\sum^{n-1}_{k=1}\left\vert a_{k}\right\vert ^{\alpha}H\left(\left\vert \gamma_{n}a^{-1}_{k}\right\vert \right).
\]

Noting that $\gamma_{n}\inf_{1\leq k\leq n}\left\vert a^{-1}_{k}\right\vert \rightarrow\infty$
and, using the same arguments as those used for in the previous parts
(see e.g., (\ref{a-moments1abb}) and (\ref{a-moments1b})), we get
\[
\sup_{1\leq t\leq n}A_{n,t}:=\sup_{1\leq t\leq n}\mathbf{E}\left\vert \gamma^{-1}_{n}\sum^{t-1}_{k=0}a_{k}\xi_{t-k}\right\vert ^{\alpha-\varepsilon}\leq C\left[A(n)+A(n)^{\frac{\alpha-\varepsilon}{\alpha+\varepsilon}}\right]\left[1+o(1)\right],
\]
for all $\varepsilon>0$ such that $\alpha-\varepsilon>0$. Therefore,
it suffices to show that $A(n)=O(1)$, which is what we set out to
do next.

As $n\rightarrow\infty$, the difference 
\[
\left\vert A(n)-\frac{1}{\gamma^{\alpha}_{n}}\sum^{n-1}_{k=1}\left\vert a_{k}\right\vert ^{\alpha}\right\vert 
\]
\[
\leq\sup_{1\leq k\leq n}\left\vert H\left(\left\vert \gamma_{n}a^{-1}_{k}\right\vert \right)-1\right\vert \frac{1}{\gamma^{\alpha}_{n}}\sum^{n-1}_{k=1}\left\vert a_{k}\right\vert ^{\alpha}=_{\left(1\right)}o(1)\cdot\frac{1}{\gamma^{\alpha}_{n}}\sum^{n-1}_{k=1}\left\vert a_{k}\right\vert ^{\alpha}
\]
\begin{equation}
=_{\left(2\right)}o(1)\cdot\left\{ \left\vert m_{0}\right\vert ^{\alpha}\mathbb{V}_{2}(\alpha,m)\left[1+o(1)\right]\right\} =o(1),\label{MI-NM2}
\end{equation}
where $=_{\left(1\right)}$ is due to \textbf{Assumption NM}, and
(\ref{MI-NM1}), whilst $=_{\left(2\right)}$ follows from \textbf{Lemma
\ref{LemmaMI4}}. Hence, 
\[
A(n)=\frac{1}{\gamma^{\alpha}_{n}}\sum^{n-1}_{k=1}\left\vert a_{k}\right\vert ^{\alpha}+o(1)=\left\vert m_{0}\right\vert ^{\alpha}\mathbb{V}_{2}(\alpha,m)+o(1),
\]
and the proof is complete. 

\bigskip{}

\textbf{(b) }Next, we show part (b). Recall that under \textbf{BFR1}
$a^{-}_{k,t}=\sum^{t-1}_{j=0}c_{k-j}$. Write
\[
x^{-}_{t}=\sum^{\infty}_{k=t}a^{-}_{k,t}\xi_{t-k}=\sum^{\infty}_{k=t}\sum^{t-1}_{j=0}c_{k-j}\xi_{t-k}
\]
\[
\left[-i=t-k\Longleftrightarrow k=t+i\right],\:\left[l=t-j\Longleftrightarrow j=t-l\right]
\]
\[
=\sum^{\infty}_{i=0}\sum^{t-1}_{j=0}c_{t+i-j}\xi_{-i}=\sum^{\infty}_{i=0}\sum^{t}_{l=0}c_{i+l}\xi_{-i}.
\]
Hence,

\[
x^{-}_{t}-x^{-}_{s}=\sum^{\infty}_{i=0}\sum^{t}_{l=0}c_{i+l}\xi_{-i}-\sum^{\infty}_{i=0}\sum^{s}_{l=0}c_{i+l}\xi_{-i}
\]
\[
=\sum^{\infty}_{i=0}\sum^{t}_{l=s+1}c_{i+l}\xi_{-i}=\sum^{\infty}_{i=0}\sum^{i+t}_{j=i+s+1}c_{j}\xi_{-i}\overset{d}{=}\sum^{\infty}_{i=0}\sum^{i+t}_{j=i+s+1}c_{j}\xi_{i}
\]
\[
\left[i=k-t\Longleftrightarrow k=i+t\right]
\]
\[
=\sum^{\infty}_{k=t}\sum^{k}_{j=k-(t-s)+1}c_{j}\xi_{k-t}\overset{d}{=}\sum^{\infty}_{k=t}\sum^{k}_{j=k-(t-s)+1}c_{j}\xi_{k}.
\]
Set $\varepsilon>0$ such that $\alpha-\varepsilon>0$. Further, $0<s<r\leq1$.
In view of the above
\[
\mathbf{E}\left\vert \gamma^{-1}_{n}\left(x^{-}_{\lfloor nr\rfloor}-x^{-}_{\lfloor ns\rfloor}\right)\right\vert ^{\alpha-\varepsilon}=\mathbf{E}\left\vert \gamma^{-1}_{n}\sum^{\infty}_{k=\lfloor nr\rfloor}\sum^{k}_{j=k-\left(\lfloor nr\rfloor-\lfloor ns\rfloor\right)+1}c_{j}\xi_{k}\right\vert ^{\alpha-\varepsilon}
\]
\[
=:\mathbf{E}\left\vert \gamma^{-1}_{n}\sum^{\infty}_{k=\lfloor nr\rfloor}b_{k,n}\xi_{k}\right\vert ^{\alpha-\varepsilon}.
\]
Next, set $T_{n,t}:=\gamma_{n}\left\vert b_{k,n}\right\vert ^{-1}$.
Note that as $n\rightarrow\infty$ 
\[
T_{n}:=\inf_{k\geq1}T_{k,n}=\frac{\gamma_{n}}{\sup_{1\leq k}\left\vert b_{k,n}\right\vert }\geq\frac{\gamma_{n}}{\sum^{\infty}_{j=0}\left\vert c_{j}\right\vert }\rightarrow\infty.
\]
Hence, $T_{k,n}\rightarrow\infty$, uniformly in $k\geq1$, as $n\rightarrow\infty$.
Furthermore, set 
\[
\overline{\xi}_{k,n}:=\xi_{k}1\left\{ \left\vert \xi_{k}\right\vert \geq T_{k,n}\right\} \text{ and }\underline{\xi}_{k,n}:=\xi_{k}1\left\{ \left\vert \xi_{k}\right\vert <T_{k,n}\right\} .
\]
By Lemma 1 in Astrauskas (1983), for $k\geq1$ and $\alpha\in(0,2]$
as $n\rightarrow\infty$ (and hence as $T_{k,n}\rightarrow\infty$)
we have
\[
\mathbf{E}\left\vert \overline{\xi}_{k,n}\right\vert ^{\alpha-\varepsilon}\asymp T^{-\varepsilon}_{k,n}h\left(T_{k,n}\right)\text{ and }\mathbf{E}\left\vert \underline{\xi}_{k,n}\right\vert ^{\alpha+\varepsilon}\asymp T^{\varepsilon}_{k,n}h\left(T_{k,n}\right).
\]
In view of the above, (\ref{a-moments1ab}), and using similar arguments
as those used in the proof of part (a) we get as $n\rightarrow\infty$,
\[
\mathbf{E}\left\vert \gamma^{-1}_{n}\sum^{\infty}_{k=\lfloor nr\rfloor}b_{k,n}\xi_{k}\right\vert ^{\alpha-\varepsilon}\leq2\gamma^{-\left(\alpha-\varepsilon\right)}_{n}\sum^{\infty}_{k=\lfloor nr\rfloor}\left\vert b_{k,n}\right\vert ^{\alpha-\varepsilon}\mathbf{E}\left\vert \overline{\xi}_{k,n}\right\vert ^{\alpha-\varepsilon}
\]
\[
+2\left[\gamma^{-\left(\alpha+\varepsilon\right)}_{n}\sum^{\infty}_{k=\lfloor nr\rfloor}\left\vert b_{k,n}\right\vert ^{\alpha+\varepsilon}\mathbf{E}\left\vert \underline{\xi}_{k,n}\right\vert ^{\alpha+\varepsilon}\right]^{\frac{\alpha-\varepsilon}{\alpha+\varepsilon}}
\]
\[
\leq2\gamma^{-\alpha}_{n}\sum^{\infty}_{k=\lfloor nr\rfloor}\left\vert b_{k,n}\right\vert ^{\alpha}H\left(\gamma_{n}\left\vert b_{k,n}\right\vert ^{-1}\right)\left[1+o(1)\right]
\]
\[
+2\left[\gamma^{-\alpha}_{n}\sum^{\infty}_{k=\lfloor nr\rfloor}\left\vert b_{k,n}\right\vert ^{\alpha}H\left(\gamma_{n}\left\vert b_{k,n}\right\vert ^{-1}\right)\left[1+o(1)\right]\right]^{\frac{\alpha-\varepsilon}{\alpha+\varepsilon}},
\]
Next, note that 
\[
\gamma^{-\alpha}_{n}\sum^{\infty}_{k=\lfloor nr\rfloor}\left\vert b_{k,n}\right\vert ^{\alpha}H\left(\gamma_{n}\left\vert b_{k,n}\right\vert ^{-1}\right)
\]
is the same form as the term $C^{\prime\prime}_{3,n}$ of (\ref{FMtypeI.3d})
which is $o(1)$ -see also (\ref{FMtypeI.6d2}). This completes the
proof of part (b).

\bigskip{}

\textbf{(c) }Finally, we show part (c) of this Lemma under \textbf{BFR2}
and \textbf{MI}-\textbf{SM/LM}. The result under \textbf{MI}-\textbf{NM}
holds trivially.

$\bigskip$

\textbf{BFR2}. The result under \textbf{Assumption BFR2} follows directly
from eq. (\ref{LPA1}) cf. proof of \textbf{Lemma \ref{LPA}}.

$\bigskip$

\textbf{MI-SM}. Without loss of generality set $r=1$. It follows
from (\ref{MI6}), (\ref{MI7}) and (\ref{MI8c}) that 
\begin{equation}
B_{n}(n)^{\prime}:=\frac{H(\gamma_{n})}{\gamma^{\alpha}_{n}}\sum^{\infty}_{k=n}\left\vert a^{-}_{k,n}\right\vert ^{\alpha}=o(1),\label{MI9}
\end{equation}
is sufficient for the requisite result -the notation used above is
the same as that in the aforementioned equation.

We first prove (\ref{MI9}) for $\alpha\in(0,1]$. In view of the
fact that, $\gamma^{-\alpha}_{n}H(\gamma_{n})\sim\kappa^{-1}_{n}$
write
\[
B_{n}(n)^{\prime}\leq C\kappa^{-1}_{n}\sum^{\infty}_{k=n}\sum^{n}_{j=0}\left\vert \rho^{j}_{n}c_{k-j}\right\vert ^{\alpha}=C\kappa^{-1}_{n}\sum^{\infty}_{k=n}\left(\sum^{\lfloor n/2\rfloor}_{j=1}+\sum^{n}_{j=\lfloor n/2\rfloor+1}\right)\left\vert \rho^{j}_{n}c_{k-j}\right\vert ^{\alpha}
\]
\begin{equation}
\leq C\kappa^{-1}_{n}\sum^{\lfloor n/2\rfloor}_{j=1}\left\vert \rho^{\alpha j}_{n}\right\vert \sum^{\infty}_{k=\lfloor n/2\rfloor}\left\vert c_{k}\right\vert ^{\alpha}+C\kappa^{-1}_{n}\sum^{n}_{j=\lfloor n/2\rfloor+1}\left\vert \rho^{\alpha j}_{n}\right\vert \sum^{\infty}_{k=0}\left\vert c_{k}\right\vert ^{\alpha}=o(1),\label{MI10}
\end{equation}
where we have used the fact that $\kappa^{-1}_{n}\sum^{\lfloor n/2\rfloor}_{j=1}\left\vert \rho^{\alpha j}_{n}\right\vert =O(1)$
and $\kappa^{-1}_{n}\sum^{n}_{j=\lfloor n/2\rfloor+1}\left\vert \rho^{\alpha j}_{n}\right\vert =o(1)$
cf. DK2018, Lemma 6.4.

Next, suppose that $\alpha\in(1,2]$. It follows from (\ref{MI8b})
that
\[
B_{n}(n)^{\prime}\leq2C\kappa^{-1}_{n}\sum^{\infty}_{k=n}\left[\left(\sum^{n}_{j=0}\left\vert \rho^{j}_{n}c_{k-j}\right\vert \right)^{2}+\sum^{n}_{j=0}\left\vert \rho^{j}_{n}c_{k-j}\right\vert \right]=o(1),
\]
where the negligibility of the first term above follows from Lemma
LVAR in DK2018, whilst the negligibility of the second term follows
from (\ref{MI10}) with $\alpha=1$.

$\bigskip$

\textbf{MI-LM}. Without loss of generality set $r=1$. Let $B_{n,n}$
as in (\ref{FSM-LM02b}). In view of From (\ref{FSM-LM2})-(\ref{FSM-LM5})
we have
\[
B_{n,n}=\frac{1}{\kappa^{b}_{n}}\sum^{\infty}_{k=n}\left(\sum^{n}_{j=0}\left\vert \rho^{j}_{n}c_{k-j}\right\vert \right)^{\alpha}\left[1+o(1)\right]=o(1).
\]
The requisite result follows easily from the above. $\blacksquare$

\bigskip{}
\newpage{}

\section{Technical Lemmas}

\begin{lemma} \label{LemmaSV1}\textit{ The following hold.}

\textbf{\textit{(i)}}\textit{ Let $\varsigma(x)>0$ SV and locally
bounded on $[x_{0}-1,\infty)$ for some $x_{0}-1>0$. In addition,
for some $\alpha>0$
\[
S(x,x_{0}):=\alpha\int^{x}_{x_{0}}\frac{\varsigma(u)}{u}du\rightarrow\infty
\]
as $x\rightarrow\infty$. Moreover, suppose that $\tau\geq0$, }$m_{n}$
\textit{is a positive sequence such that $m^{-1}_{n}+m^{\tau/\alpha}_{n}n^{-1}\rightarrow0$,
and $S(m_{n}n^{1/\alpha},m^{1+\tau/\alpha}_{n}x^{1/\alpha}_{0})\sim S(m_{n}n^{1/\alpha},x_{0})$
as $n\rightarrow\infty$. Then as $n\rightarrow\infty$,
\[
\sum^{n}_{k=x_{0}\left\lfloor m^{\tau}_{n}\right\rfloor }\varsigma(m_{n}k^{1/\alpha})/k\sim S(m_{n}n^{1/\alpha},m^{1+\tau/\alpha}_{n}x^{1/\alpha}_{0})=\alpha\int^{m_{n}n^{1/\alpha}}_{m^{1+\tau/\alpha}_{n}x^{1/\alpha}_{0}}\frac{\varsigma(x)}{x}dx.
\]
}

\textbf{\textit{(ii)}}\textit{ Let $l>-1$, $m(x)$ SV, $q>0$ and
$q_{n}\rightarrow\infty$. As $n\rightarrow\infty$ }
\[
\frac{1}{n^{1+l}m\left(q_{n}n^{q}\right)}\sum^{n}_{j=1}j^{l}m\left(q_{n}j^{q}\right)=\int^{1}_{0}x^{l}dx+o(1).
\]

\end{lemma} \bigskip{}

\noindent\textbf{Proof of Lemma \ref{LemmaSV1}}.

\bigskip{}

\textbf{(i)} By Theorem 1.3.3 and Theorem 1.5.3 of BGT 1987 $\varsigma(x^{1/\alpha})/x$\
can be approximated (in the ``$\sim$''\
sense) by a regularly varying function of index $-1$ that is continuous
on $[x^{\prime}_{0},\infty)$ with $x^{\prime}_{0}>0$ and decreasing.
In view of this, there is no loss of generality if in the subsequent
derivations we assume that $\varsigma(x^{1/\alpha})/x$ is continuous
and decreasing for $x\geq x^{\prime}_{0}$. For $n$ large enough,
such that $m_{n}\left(x_{0}-1\right)^{1/\alpha}\geq x^{\prime}_{0}$,
write
\[
\int^{n}_{x_{0}\left\lfloor m^{\tau}_{n}\right\rfloor }\frac{\varsigma(m_{n}x^{1/\alpha})}{x}dx=\int^{m_{n}n^{1/\alpha}}_{u=m_{n}\left\lfloor m^{\tau/\alpha}_{n}\right\rfloor x^{1/\alpha}_{0}}\frac{\varsigma(u)}{m^{-\alpha}_{n}u^{\alpha}}m^{-\alpha}_{n}\alpha u^{\alpha-1}du
\]
\[
\left[u=m_{n}x^{1/\alpha}\Leftrightarrow x=m^{-\alpha}_{n}u^{\alpha}\Leftrightarrow dx=m^{-\alpha}_{n}\alpha u^{\alpha-1}du\right]
\]
\[
=\alpha\int^{m_{n}n^{1/\alpha}}_{u=m_{n}\left\lfloor m^{\tau/\alpha}_{n}\right\rfloor x^{1/\alpha}_{0}}\frac{\varsigma(u)}{u}du\sim_{(1)}\alpha\int^{m_{n}n^{1/\alpha}}_{u=m^{1+\tau/\alpha}_{n}x^{1/\alpha}_{0}}\frac{\varsigma(u)}{u}du
\]
\[
=S(m_{n}n^{1/\alpha},m_{n}x^{1/\alpha}_{0})=:S_{n},
\]
where note that under our assumptions we have $S_{n}\rightarrow\infty$.
In addition, $\sim_{(1)}$ follows from the fact that $S(x,x_{0})$
is SV (cf. BGT 1987; Proposition 1.5.9a). To see this, first note
that it be easily checked that
\begin{equation}
S\left(m^{1+\tau/\alpha}_{n}x^{1/\alpha}_{0},x_{0}\right)\sim S\left(m_{n}\left\lfloor m^{\tau/\alpha}_{n}\right\rfloor x^{1/\alpha}_{0},x_{0}\right).\label{LemSV1.0}
\end{equation}
 Next, write
\[
S\left(m_{n}n^{1/\alpha},m_{n}\left\lfloor m^{\tau/\alpha}_{n}\right\rfloor x^{1/\alpha}_{0}\right)=S\left(m_{n}n^{1/\alpha},x_{0}\right)-S\left(m_{n}\left\lfloor m^{\tau/\alpha}_{n}\right\rfloor x^{1/\alpha}_{0},x_{0}\right).
\]
Hence,
\[
\frac{S\left(m_{n}n^{1/\alpha},m_{n}\left\lfloor m^{\tau/\alpha}_{n}\right\rfloor x^{1/\alpha}_{0}\right)}{S\left(m_{n}n^{1/\alpha},x_{0}\right)}=1-\frac{S\left(m_{n}\left\lfloor m^{\tau/\alpha}_{n}\right\rfloor x^{1/\alpha}_{0},x_{0}\right)}{S\left(m_{n}n^{1/\alpha},x_{0}\right)}
\]
\[
=1-\frac{S\left(m^{1+\tau/\alpha}_{n}x^{1/\alpha}_{0},x_{0}\right)}{S\left(m_{n}n^{1/\alpha},x_{0}\right)}\frac{S\left(m_{n}\left\lfloor m^{\tau/\alpha}_{n}\right\rfloor x^{1/\alpha}_{0},x_{0}\right)}{S\left(m^{1+\tau/\alpha}_{n}x^{1/\alpha}_{0},x_{0}\right)}
\]
\[
=_{(2)}1-\frac{S\left(m^{1+\tau/\alpha}_{n}x^{1/\alpha}_{0},x_{0}\right)}{S\left(m_{n}n^{1/\alpha},x_{0}\right)}\left[1+o(1)\right]=_{(3)}1+o(1),
\]
where $=_{(2)}$ follows from (\ref{LemSV1.0}). For $=_{(3)}$ recall
that $S(m_{n}n^{1/\alpha},x_{0})\sim S\left(m_{n}n^{1/\alpha},m^{1+\tau/\alpha}_{n}x^{1/\alpha}_{0}\right)$.
Therefore, 
\[
\frac{S\left(m^{1+\tau/\alpha}_{n}x^{1/\alpha}_{0},x_{0}\right)}{S\left(m_{n}n^{1/\alpha},x_{0}\right)}=\frac{S(m_{n}n^{1/\alpha},x_{0})-S\left(m_{n}n^{1/\alpha},m^{1+\tau/\alpha}_{n}x^{1/\alpha}_{0}\right)}{S\left(m_{n}n^{1/\alpha},x_{0}\right)}=o(1).
\]

Next, note that because $\varsigma(x^{1/\alpha})/x$ is continuous
we can write
\[
S^{-1}_{n}\sum^{n}_{k=x_{0}}\varsigma(m_{n}k^{1/\alpha})/k=S^{-1}_{n}\int^{n}_{x_{0}-1}\frac{\varsigma(m_{n}x^{1/\alpha})}{x}d\lfloor x\rfloor
\]
\[
=_{(4)}S^{-1}_{n}\int^{n}_{x_{0}}\frac{\varsigma(m_{n}x^{1/\alpha})}{x}d\lfloor x\rfloor+o(1),
\]
where $=_{(4)}$ follows from the fact that $\varsigma(x)=o\left(S(x,x_{0})\right)$
(cf. BGT 1987; Proposition 1.5.9a) i.e.
\[
S^{-1}_{n}\int^{x_{0}}_{x_{0}-1}\frac{\varsigma(m_{n}x^{1/\alpha})}{x}d\lfloor x\rfloor=S^{-1}_{n}\frac{\varsigma(m_{n}x^{1/\alpha}_{0})}{x_{0}}
\]
\[
\sim S(m_{n}n^{1/\alpha},x_{0})^{-1}\frac{\varsigma(m_{n}x^{1/\alpha}_{0})}{x_{0}}
\]
\begin{equation}
\leq S(m_{n}x^{1/\alpha}_{0},x_{0})^{-1}\frac{\varsigma(m_{n}x^{1/\alpha}_{0})}{x_{0}}=o(1).\label{LemSV1.1}
\end{equation}
 In addition, using the fact that $\varsigma(x^{1/\alpha})/x$ is
decreasing the total variation of $\varsigma(m_{n}x^{1/\alpha})x^{-1}$
on $[x_{0},n]$ is 
\[
\mathbf{V}\left(\frac{\varsigma(m_{n}x^{1/\alpha})}{x}\right)^{n}_{x_{0}}=m^{\alpha}_{n}\mathbf{V}\left(\frac{\varsigma(m_{n}x^{1/\alpha})}{m^{\alpha}_{n}x}\right)^{n}_{x_{0}}
\]
\[
\leq m^{\alpha}_{n}\frac{\varsigma(m_{n}x^{1/\alpha}_{0})}{m^{\alpha}_{n}x_{0}}=\frac{\varsigma(m_{n}x^{1/\alpha}_{0})}{x_{0}}=o\left(S_{n}\right),
\]
where the last equality above follows from (\ref{LemSV1.1}). In view
of the above,
\[
\left\vert S^{-1}_{n}\sum^{n}_{k=x_{0}}\varsigma(m_{n}k^{1/\alpha})/k-1\right\vert =S^{-1}_{n}\left\vert \sum^{n}_{k=x_{0}}\varsigma(m_{n}k^{1/\alpha})/k-\int^{n}_{x_{0}}\frac{\varsigma(m_{n}x^{1/\alpha})}{x}dx\right\vert 
\]
\[
=S^{-1}_{n}\left\vert \int^{n}_{x_{0}}\frac{\varsigma(m_{n}x^{1/\alpha})}{x}d\lfloor x\rfloor-\int^{n}_{x_{0}}\frac{\varsigma(m_{n}x^{1/\alpha})}{x}dx\right\vert +o(1)
\]
\[
=S^{-1}_{n}\left\vert \int^{n}_{x_{0}}\frac{\varsigma(m_{n}x^{1/\alpha})}{x}d\left(x-\lfloor x\rfloor\right)\right\vert +o(1)
\]
\[
=S^{-1}_{n}\left|\int^{n}_{x_{0}}\left(x-\lfloor x\rfloor\right)d\left(\frac{\varsigma(m_{n}x^{1/\alpha})}{x}\right)\right|+o(1)
\]
\[
\leq S^{-1}_{n}\int^{n}_{x_{0}}\left|x-\lfloor x\rfloor\right|d\mathbf{V}\left(\frac{\varsigma(m_{n}x^{1/\alpha})}{x}\right)+o(1)
\]
\[
\leq S^{-1}_{n}\int^{n}_{x_{0}}d\mathbf{V}\left(\frac{\varsigma(m_{n}x^{1/\alpha})}{x}\right)+o(1)
\]
\[
\leq S^{-1}_{n}\mathbf{V}\left(\frac{\varsigma(m_{n}x^{1/\alpha})}{x}\right)^{n}_{x_{0}}+o(1)\leq S^{-1}_{n}\frac{\varsigma(m_{n}x^{1/\alpha}_{0})}{x_{0}}=o(1),
\]
as required.\bigskip{}

\textbf{(ii)} Let $l_{2}>-1$ with $l_{1}:=l-l_{2}>0$, and set $R_{r}(x):=x^{r}m\left(x\right)$,
$r\in\mathbb{R}$. Then
\[
\left\vert \frac{1}{n^{1+l}m\left(q_{n}n^{q}\right)}\sum^{n}_{j=1}j^{l}m\left(q_{n}j^{q}\right)-\frac{1}{n^{1+l}}\sum^{n}_{j=1}j^{l}\right\vert \leq\frac{1}{n}\sum^{n}_{j=1}\left\vert \frac{\left(q_{n}j^{q}\right)^{l/q}m\left(q_{n}j^{q}\right)}{\left(q_{n}n^{q}\right)^{l/q}m\left(q_{n}n^{q}\right)}-\left(\frac{j^{q}}{n^{q}}\right)^{l/q}\right\vert 
\]
\[
=\frac{1}{n}\sum^{n}_{j=1}\left\vert \frac{\left(q_{n}j^{q}\right)^{l_{1}/q}m\left(q_{n}j^{q}\right)}{\left(q_{n}n^{q}\right)^{l_{1}/q}m\left(q_{n}n^{q}\right)}-\left(\frac{j^{q}}{n^{q}}\right)^{l_{1}/q}\right\vert \left(\frac{j^{q}}{n^{q}}\right)^{l_{2}/q}=\frac{1}{n}\sum^{n}_{j=1}\left\vert \frac{R_{l_{1}/q}\left(q_{n}j^{q}\right)}{R_{l_{1}/q}\left(q_{n}n^{q}\right)}-\left(\frac{j^{q}}{n^{q}}\right)^{l_{1}/q}\right\vert \left(\frac{j}{n}\right)^{l_{2}}
\]
\[
\leq\sup_{0<y\leq1}\left\vert \frac{R_{l_{1}/q}\left(q_{n}n^{q}y\right)}{R_{l_{1}/q}\left(q_{n}n^{q}\right)}-y^{l_{1}/q}\right\vert \frac{1}{n}\sum^{n}_{j=1}\left(\frac{j}{n}\right)^{l_{2}}
\]
\[
=\sup_{0<y\leq1}\left\vert \frac{R_{l_{1}/q}\left(q_{n}n^{q}y\right)}{R_{l_{1}/q}\left(q_{n}n^{q}\right)}-y^{l_{1}/q}\right\vert \left[\int^{1}_{0}x^{l_{2}}dx+o(1)\right]=o(1),
\]
where the last equality above is a consequence of the uniform convergence
of regularly varying functions -see e.g. Theorem 1.5.2. in BGT 1987.
In view of the above, the result follows from the fact that $n^{-\left(1+l\right)}\sum^{n}_{j=1}j^{l}$
$=\int^{1}_{0}x^{l}dx+o(1)$. $\blacksquare$

\bigskip{}

\begin{lemma} \label{LemmaMI} Let $0<m_{n}$ an integer valued sequence
such that $m^{-1}_{n}+m_{n}\kappa^{-1}_{n}\rightarrow0$. Under \textbf{Assumption
MI-SM}, for all $r>0$,
\[
\kappa^{-1}_{n}\sum^{\lfloor\left(n-1\right)r\rfloor}_{k=1}\left\vert a_{k}\right\vert ^{\alpha}=\kappa^{-1}_{n}\sum^{\lfloor\left(n-1\right)r\rfloor}_{k=m_{n}}\left\vert a_{k}\right\vert ^{\alpha}+o(1)
\]
\[
=\alpha^{-1}\left\vert \sum^{\infty}_{l=0}c_{l}\right\vert ^{\alpha}+o(1).
\]
\end{lemma}

\bigskip{}

\begin{lemma} \label{LemmaMI-2} Under \textbf{Assumption MI-SM},
the following hold:

\textbf{(a)} As $n\rightarrow\infty$,
\[
\kappa^{-1}_{n}\sum^{n}_{k=1}\left\vert a_{k}-\rho^{k}_{n}c_{\infty}\right\vert ^{\alpha}=o(1).
\]

\textbf{(b)} For each $0<r\leq1$, as $n\rightarrow\infty$,
\[
\kappa^{-1}_{n}\sum^{\lfloor nr\rfloor}_{k=1}\left\vert a_{k}\right\vert ^{\alpha}=\kappa^{-1}_{n}\sum^{\lfloor nr\rfloor}_{k=1}\left\vert \rho^{k}_{n}c_{\infty}\right\vert ^{\alpha}+o(1)=\alpha^{-1}\left\vert \sum^{\infty}_{l=0}c_{l}\right\vert ^{\alpha}+o(1).
\]

\textbf{(c)} For each $0<r\leq1$, as $n\rightarrow\infty$,
\[
\kappa^{-1}_{n}\sum^{\lfloor nr\rfloor}_{k=1}\left\vert a_{k}\right\vert ^{\alpha}\text{sgn}\left(\frac{1}{\gamma_{n}}a_{k}\right)=\kappa^{-1}_{n}\sum^{\lfloor nr\rfloor}_{k=1}\left\vert a_{k}\right\vert ^{\alpha}\text{sgn}\left(\rho^{k}_{n}c_{\infty}\right)+o(1).
\]

\end{lemma}

\bigskip{}

\noindent\textbf{Proof of Lemma \ref{LemmaMI}}. Note that for $n$
large enough, $0<\rho^{j}_{n}\leq1$, for all $j\geq0$. Hence,
\[
\kappa^{-1}_{n}\sum^{m_{n}}_{k=1}\left\vert \sum^{k}_{l=0}\rho^{k-l}_{n}c_{l}\right\vert ^{\alpha}\leq\kappa^{-1}_{n}\sum^{m_{n}}_{k=1}\left(\sum^{k}_{l=0}\rho^{k-l}_{n}\left\vert c_{l}\right\vert \right)^{\alpha}
\]
\begin{equation}
\leq\kappa^{-1}_{n}\sum^{m_{n}}_{k=1}\left(\sum^{\infty}_{l=0}\left\vert c_{l}\right\vert \right)^{\alpha}=O\left(\frac{m_{n}}{\kappa_{n}}\right)\rightarrow0.\label{MI1b}
\end{equation}
Next, consider
\[
\kappa^{-1}_{n}\sum^{\lfloor\left(n-1\right)r\rfloor}_{k=m_{n}}\left\vert a_{k}\right\vert ^{\alpha}=\kappa^{-1}_{n}\sum^{\lfloor\left(n-1\right)r\rfloor}_{k=m_{n}}\left\vert \sum^{k}_{l=0}\rho^{k-l}_{n}c_{l}\right\vert ^{\alpha}
\]

\begin{equation}
=\kappa^{-1}_{n}\sum^{\lfloor\left(n-1\right)r\rfloor}_{k=m_{n}}\left\vert \left(\sum^{m_{n}}_{l=0}+\sum^{k}_{l=m_{n}+1}\right)\rho^{k-l}_{n}c_{l}\right\vert ^{\alpha}.\label{MI2}
\end{equation}
Using the fact $\lim_{n\rightarrow\infty}\rho^{-l}_{n}=1$ for each
$l$ and $\sum^{\infty}_{l=0}\left\vert c_{l}\right\vert <\infty$
dominated convergence gives $\sum^{m_{n}}_{l=0}\left\vert \rho^{-l}_{n}c_{l}-c_{l}\right\vert \rightarrow0$.
It follows from DK2018, Lemma 6.4, that
\[
\kappa^{-1}_{n}\sum^{\lfloor\left(n-1\right)r\rfloor}_{k=m_{n}}\left\vert \sum^{m_{n}}_{l=0}\rho^{k-l}_{n}c_{l}\right\vert ^{\alpha}=\kappa^{-1}_{n}\sum^{\lfloor\left(n-1\right)r\rfloor}_{k=m_{n}}\left\vert \rho^{\alpha k}_{n}\right\vert \left\vert \sum^{m_{n}}_{l=0}\rho^{-l}_{n}c_{l}\right\vert ^{\alpha}
\]
\[
=\kappa^{-1}_{n}\sum^{\lfloor\left(n-1\right)r\rfloor}_{k=m_{n}}\left\vert \rho^{\alpha k}_{n}\right\vert \left\vert \sum^{m_{n}}_{l=0}c_{l}\right\vert ^{\alpha}\left[1+o(1)\right]
\]
\[
=\int^{\lfloor\left(n-1\right)r\rfloor/\kappa_{n}}_{m_{n}/\kappa_{n}}e^{-\alpha x}dx\left\vert \sum^{\infty}_{l=0}c_{l}\right\vert ^{\alpha}+o(1)
\]

\begin{equation}
=\int^{\infty}_{0}e^{-\alpha x}dx\left\vert \sum^{\infty}_{l=0}c_{l}\right\vert ^{\alpha}+o(1)=\alpha^{-1}\left\vert \sum^{\infty}_{l=0}c_{l}\right\vert ^{\alpha}.\label{MI3}
\end{equation}
\bigskip{}
Next, set
\[
R_{n}(\alpha):=\kappa^{-1}_{n}\sum^{\lfloor\left(n-1\right)r\rfloor}_{k=m_{n}}\left\vert \sum^{k}_{l=m_{n}+1}\rho^{k-l}_{n}c_{l}\right\vert ^{\alpha}.
\]
We shall show that
\begin{equation}
R_{n}(\alpha)=o(1).\label{MI4}
\end{equation}
First, we consider the $\alpha\in(0,1]$. As $n\rightarrow\infty$
\[
R_{n}(\alpha)\leq\kappa^{-1}_{n}\sum^{n}_{k=m_{n}}\sum^{k}_{l=m_{n}+1}\rho^{\alpha\left(k-l\right)}_{n}\left\vert c_{l}\right\vert ^{\alpha}
\]
\[
=\kappa^{-1}_{n}\sum^{n}_{l=m_{n}+1}\sum^{n}_{k=m_{n}\vee l}\rho^{\alpha\left(k-l\right)}_{n}\left\vert c_{l}\right\vert ^{\alpha}\leq\kappa^{-1}_{n}\sum^{n}_{l=m_{n}}\rho^{-\alpha l}_{n}\left\vert c_{l}\right\vert ^{\alpha}\sum^{n}_{k=l}\rho^{\alpha k}_{n}
\]
\[
=\kappa^{-1}_{n}\sum^{n}_{l=m_{n}}\rho^{-\alpha l}_{n}\left\vert c_{l}\right\vert ^{\alpha}\left[\frac{1-\rho^{\alpha\left(n+1\right)}_{n}}{1-\rho^{\alpha}_{n}}-\frac{1-\rho^{\alpha l}_{n}}{1-\rho^{\alpha}_{n}}\right]
\]
\[
=\frac{1}{1-\rho^{\alpha}_{n}}\kappa^{-1}_{n}\sum^{n}_{l=m_{n}}\rho^{-\alpha l}_{n}\left\vert c_{l}\right\vert ^{\alpha}\left[\rho^{\alpha l}_{n}-\rho^{\alpha\left(n+1\right)}_{n}\right]
\]
\[
=\frac{1}{1-\rho^{\alpha}_{n}}\kappa^{-1}_{n}\sum^{n}_{l=m_{n}}\left\vert c_{l}\right\vert ^{\alpha}-\frac{\rho^{\alpha\left(n+1\right)}_{n}}{1-\rho^{\alpha}_{n}}\kappa^{-1}_{n}\sum^{n}_{l=m_{n}}\left\vert c_{l}\right\vert ^{\alpha}
\]
\[
=_{(1)}\left[1+o(1)\right]C\sum^{\infty}_{l=m_{n}}\left\vert c_{l}\right\vert ^{\alpha}=o(1),
\]
where $=_{(1)}$ follows from the fact that $\rho^{\alpha\left(n+1\right)}_{n}=o(1)$
and the fact that $\left(1-\rho^{\alpha}_{n}\right)^{-1}\kappa^{-1}_{n}\rightarrow1/\alpha$.
For the latter, see Footnote 11. Moreover, an argument similar to
the one shows that (\ref{MI0a}), and $R_{n}(\alpha)=o(1)$ for $\alpha\in(1,2]$.
Next, we show that 
\[
R^{\prime}_{n}:=\left\vert \kappa^{-1}_{n}\sum^{\lfloor\left(n-1\right)r\rfloor}_{k=m_{n}}\left\vert a_{k}\right\vert ^{\alpha}-\kappa^{-1}_{n}\sum^{\lfloor\left(n-1\right)r\rfloor}_{k=m_{n}}\left\vert \sum^{m_{n}}_{l=0}\rho^{k-l}_{n}c_{l}\right\vert ^{\alpha}\right\vert =o(1)
\]
First, suppose that $0<\alpha\leq1$. In view of (\ref{MI2}), for
\[
R^{\prime}_{n}=\kappa^{-1}_{n}\sum^{\lfloor\left(n-1\right)r\rfloor}_{k=m_{n}}\left\vert \left\vert \left(\sum^{m_{n}}_{l=0}+\sum^{k}_{l=m_{n}+1}\right)\rho^{k-l}_{n}c_{l}\right\vert ^{\alpha}-\left\vert \sum^{m_{n}}_{l=0}\rho^{k-l}_{n}c_{l}\right\vert ^{\alpha}\right\vert 
\]
\begin{equation}
\leq_{\left(2\right)}R_{n}(\alpha)=_{(3)}o(1),\label{MI5}
\end{equation}
where $\leq_{\left(2\right)}$ follows from Loève's inequality and
$=_{(3)}$ from the fact that $R_{n}(\alpha)=o(1)$. Next, suppose
that $1<\alpha\leq2$. Using the fact that, see (\ref{MVT}),
\[
\left\vert \left|y+x\right|^{\alpha}-\left|y\right|^{\alpha}\right\vert \leq C\left(\left\vert y\right\vert ^{\alpha-1}\left\vert x\right\vert +\left\vert x\right\vert ^{\alpha}\right),
\]
for $n$ large enough,
\[
R^{\prime}_{n}\leq\kappa^{-1}_{n}\sum^{n-1}_{k=m_{n}}\left\vert \sum^{m_{n}}_{l=0}\rho^{k-l}_{n}c_{l}\right\vert ^{\alpha-1}\left\vert \sum^{k}_{l=m_{n}+1}\rho^{k-l}_{n}c_{l}\right\vert +\kappa^{-1}_{n}\sum^{n-1}_{k=m_{n}}\left\vert \sum^{k}_{l=m_{n}+1}\rho^{k-l}_{n}c_{l}\right\vert ^{\alpha}
\]
\[
\leq_{\left(4\right)}R_{n}(1)\sum^{m_{n}}_{l=0}\left\vert c_{l}\right\vert +R_{n}(\alpha)\leq R_{n}(1)\sum^{\infty}_{l=0}\left\vert c_{l}\right\vert +R_{n}(\alpha)=o(1),
\]
where $\leq_{\left(4\right)}$ follows from the fact that $\left\vert \rho^{k-l}_{n}\right\vert \leq1$
for $n$ large and all $k\geq l$, the last approximation follows
from the fact that $R_{n}(\alpha)=o(1)$, $\alpha\in(0,2]$. Therefore,
for all $\alpha\in(0,2]$
\begin{equation}
R^{\prime}_{n}=o(1).\label{MI5a}
\end{equation}
Using (\ref{MI1b})-(\ref{MI4}) and (\ref{MI5a}) for all $\alpha\in(0,2]$
we get
\[
\kappa^{-1}_{n}\sum^{\lfloor\left(n-1\right)r\rfloor}_{k=1}\left\vert a_{k}\right\vert ^{\alpha}=\kappa^{-1}_{n}\sum^{\lfloor\left(n-1\right)r\rfloor}_{k=m_{n}}\left\vert a_{k}\right\vert ^{\alpha}+o(1)
\]
\[
=\kappa^{-1}_{n}\sum^{\lfloor\left(n-1\right)r\rfloor}_{k=m_{n}}\left\vert \sum^{m_{n}}_{l=0}\rho^{k-l}_{n}c_{l}\right\vert ^{\alpha}+o(1)=\alpha^{-1}\left\vert \sum^{\infty}_{l=0}c_{l}\right\vert ^{\alpha}+o(1),
\]
as required. $\blacksquare$\bigskip{}

\noindent\textbf{Proof of Lemma \ref{LemmaMI-2}}. \textbf{(a)} Write
\[
\kappa^{-1}_{n}\sum^{n}_{k=1}\left|a_{k}-\rho^{k}_{n}c_{\infty}\right|^{\alpha}
\]
\[
\leq C\kappa^{-1}_{n}\sum^{n}_{k=1}\left\vert \sum^{k}_{l=0}\rho^{k-l}_{n}c_{l}-\sum^{k}_{l=0}\rho^{k}_{n}c_{l}\right\vert ^{\alpha}+C\kappa^{-1}_{n}\sum^{n}_{k=1}\left\vert \sum^{\infty}_{l=k+1}\rho^{k}_{n}c_{l}\right\vert ^{\alpha}
\]
\[
=:C\left[R_{1n}(\alpha)+R_{2n}(\alpha)\right].
\]
We shall show that as $n\rightarrow\infty$
\begin{equation}
R_{1n}(\alpha)+R_{2n}(\alpha)=o(1).\label{LemmaMI-2.1}
\end{equation}

\noindent First, suppose that $\alpha\in(0,1]$. Note that for $n$
large enough, $0<\rho^{j}_{n}\leq1$, for all $j\geq0$. In addition,
for each $j\geq0$ fixed $\rho^{j}_{n}\rightarrow1$, as $n\rightarrow\infty$.
In view of this, as $n\rightarrow\infty$, the first term above is
\[
R_{1n}(\alpha)=\kappa^{-1}_{n}\sum^{n}_{k=1}\left\vert \sum^{k}_{l=0}\rho^{k-l}_{n}\left(1-\rho^{l}_{n}\right)c_{l}\right\vert ^{\alpha}\leq\kappa^{-1}_{n}\sum^{n}_{k=1}\sum^{k}_{l=0}\rho^{\alpha\left(k-l\right)}_{n}\left(1-\rho^{l}_{n}\right)^{\alpha}\left|c_{l}\right|^{\alpha}
\]
\[
=\kappa^{-1}_{n}\sum^{n}_{l=0}\left(1-\rho^{l}_{n}\right)^{\alpha}\left|c_{l}\right|^{\alpha}\sum^{n}_{k=l}\rho^{\alpha\left(k-l\right)}_{n}=\kappa^{-1}_{n}\sum^{n}_{l=0}\left(1-\rho^{l}_{n}\right)^{\alpha}\left|c_{l}\right|^{\alpha}\sum^{n}_{k=l}\rho^{\alpha\left(k-l\right)}_{n}
\]
\[
=\kappa^{-1}_{n}\sum^{n}_{l=0}\left(1-\rho^{l}_{n}\right)^{\alpha}\left|c_{l}\right|^{\alpha}\sum^{n-l}_{k=0}\rho^{\alpha k}_{n}\leq\kappa^{-1}_{n}\sum^{n}_{l=0}\left(1-\rho^{l}_{n}\right)^{\alpha}\left|c_{l}\right|^{\alpha}\sum^{n}_{k=0}\rho^{\alpha k}_{n}
\]
\[
\leq_{(1)}C\sum^{\infty}_{l=0}\left(1-\rho^{l}_{n}\right)^{\alpha}\left|c_{l}\right|^{\alpha}=_{(2)}o(1),
\]
where $\leq_{(1)}$ follows from the fact that $\kappa^{-1}_{n}\sum^{n}_{k=0}\rho^{\alpha k}_{n}=O(1)$,
and $=_{(2)}$ from dominated convergence. Next, let $m_{n}>0$ an
integer valued sequence such that $m^{-1}_{n}+m_{n}\kappa^{-1}_{n}\rightarrow0$.
As $n\rightarrow\infty$, the second term
\[
R_{2n}(\alpha)\leq\kappa^{-1}_{n}\sum^{n}_{k=1}\sum^{\infty}_{l=k+1}\rho^{\alpha k}_{n}\left|c_{l}\right|^{\alpha}
\]
\[
=\kappa^{-1}_{n}\left(\sum^{n}_{k=m_{n}+1}+\sum^{m_{n}}_{k=1}\right)\rho^{\alpha k}_{n}\sum^{\infty}_{l=k+1}\left|c_{l}\right|^{\alpha}
\]
\[
\leq\kappa^{-1}_{n}\sum^{n}_{k=m_{n}+1}\rho^{\alpha k}_{n}\sum^{\infty}_{l=m_{n}+1}\left|c_{l}\right|^{\alpha}+\kappa^{-1}_{n}\sum^{m_{n}}_{k=1}\rho^{\alpha k}_{n}\sum^{\infty}_{l=0}\left|c_{l}\right|^{\alpha}
\]
\[
=C\sum^{\infty}_{l=m_{n}+1}\left|c_{l}\right|^{\alpha}+O(m_{n}/\kappa_{n})=o(1).
\]
Hence, (\ref{LemmaMI-2.1}) holds.

\noindent Next, suppose that $1<\alpha\leq2$. First, consider $R_{1n}$($\alpha$).
As $n\rightarrow\infty$
\[
\kappa^{-1}_{n}\sum^{n}_{k=1}\left\vert \sum^{k}_{l=0}\rho^{k-l}_{n}c_{l}-\sum^{k}_{l=0}\rho^{k}_{n}c_{l}\right\vert ^{\alpha}=\kappa^{-1}_{n}\sum^{n}_{k=1}\left\vert \sum^{k}_{l=0}\rho^{k-l}_{n}\left(1-\rho^{l}_{n}\right)c_{l}\right\vert ^{\alpha}
\]
\[
=\kappa^{-1}_{n}\sum^{n}_{k=1}\left\vert \sum^{k}_{l=0}\rho^{k-l}_{n}\left(1-\rho^{l}_{n}\right)c_{l}\right\vert ^{\alpha-1}\left\vert \sum^{k}_{l=0}\rho^{k-l}_{n}\left(1-\rho^{l}_{n}\right)c_{l}\right\vert 
\]
\[
\leq\kappa^{-1}_{n}\sum^{n}_{k=1}\left\vert \sum^{\infty}_{l=0}2\left|c_{l}\right|\right\vert ^{\alpha-1}\left\vert \sum^{k}_{l=0}\rho^{k-l}_{n}\left(1-\rho^{l}_{n}\right)c_{l}\right\vert 
\]
\[
=C\kappa^{-1}_{n}\sum^{n}_{k=1}\left\vert \sum^{k}_{l=0}\rho^{k-l}_{n}\left(1-\rho^{l}_{n}\right)c_{l}\right\vert =CR_{1n}(1)=o(1).
\]
Finally, as $n\rightarrow\infty$
\[
R_{2n}(\alpha)=\kappa^{-1}_{n}\sum^{n}_{k=1}\left\vert \sum^{\infty}_{l=k+1}\rho^{k}_{n}c_{l}\right\vert ^{\alpha-1}\left\vert \sum^{\infty}_{l=k+1}\rho^{k}_{n}c_{l}\right\vert 
\]
\[
\leq\left(\sum^{\infty}_{l=k+1}\left|c_{l}\right|\right)^{\alpha-1}\kappa^{-1}_{n}\sum^{n}_{k=1}\left\vert \sum^{\infty}_{l=k+1}\rho^{k}_{n}c_{l}\right\vert 
\]
\[
=CR_{2n}(1)=o(1),
\]
as required. The proof for $0<\alpha\leq1$ is analogous for the one
provided above and therefore omitted. 

\noindent\textbf{(b)} We will prove the result for the case where
$1<\alpha\leq2$ only. Recall that by (\ref{MVT}), 
\[
\kappa^{-1}_{n}\sum^{n}_{k=1}\left\vert \left\vert a_{k}\right\vert ^{\alpha}-\left\vert \rho^{k}_{n}c_{\infty}\right\vert ^{\alpha}\right\vert \leq C\cdot\left[\kappa^{-1}_{n}\sum^{n}_{k=1}\left\vert \rho^{k}_{n}c_{\infty}\right\vert ^{\alpha-1}\left\vert a_{k}-\rho^{k}_{n}c_{\infty}\right\vert +\kappa^{-1}_{n}\sum^{n}_{k=1}\left\vert a_{k}-\rho^{k}_{n}c_{\infty}\right\vert ^{\alpha}\right]\text{.}
\]
\[
\leq C\cdot\left[\kappa^{-1}_{n}\sum^{n}_{k=1}\left\vert a_{k}-\rho^{k}_{n}c_{\infty}\right\vert +\kappa^{-1}_{n}\sum^{n}_{k=1}\left\vert a_{k}-\rho^{k}_{n}c_{\infty}\right\vert ^{\alpha}\right]=o(1),
\]
where the last approximation follows from part (a) above. This shows
the first approximation of part (b).

Finally, we show the second approximation. By Lemma 6.4 in DK2018
we get
\[
\kappa^{-1}_{n}\sum^{\lfloor\left(n-1\right)r\rfloor}_{k=1}\left\vert \rho^{k}_{n}c_{\infty}\right\vert ^{\alpha}=\left\vert c_{\infty}\right\vert ^{\alpha}\kappa^{-1}_{n}\sum^{\lfloor\left(n-1\right)r\rfloor}_{k=1}\rho^{\alpha k}_{n}
\]
\[
=\left\vert c_{\infty}\right\vert ^{\alpha}\kappa^{-1}_{n}\left[\sum^{\kappa_{n}\ln(\kappa_{n})}_{k=1}\rho^{\alpha k}_{n}+\sum^{\lfloor\left(n-1\right)r\rfloor}_{k=\kappa_{n}\ln(\kappa_{n})}\rho^{\alpha k}_{n}\right]=\left\vert c_{\infty}\right\vert ^{\alpha}\kappa^{-1}_{n}\sum^{\kappa_{n}\ln(\kappa_{n})}_{k=1}\rho^{\alpha k}_{n}+o(1)
\]
\[
=\left\vert c_{\infty}\right\vert ^{\alpha}\kappa^{-1}_{n}\sum^{\kappa_{n}\ln(\kappa_{n})}_{k=1}e^{-\alpha\frac{k}{\kappa_{n}}}+o(1)=\left\vert c_{\infty}\right\vert ^{\alpha}\int^{\kappa_{n}\ln(\kappa_{n})/\kappa_{n}}_{1/\kappa_{n}}e^{-\alpha x}dx+o(1)
\]
\[
=\left\vert c_{\infty}\right\vert ^{\alpha}\int^{\infty}_{0}e^{-\alpha x}dx+o(1)=\alpha^{-1}\left\vert c_{\infty}\right\vert ^{\alpha}+o(1),
\]
as required.

\textbf{(c)} Set 
\[
K_{n,k}(r)=\left\{ 1\leq k\leq\lfloor nr\rfloor:\text{sgn}\left(\frac{1}{\gamma_{n}}a_{k}\right)\neq\text{sgn}\left(\rho^{k}_{n}c_{\infty}\right)\right\} .
\]
Next, note that as $n\rightarrow\infty$,
\[
\kappa^{-1}_{n}\sum^{\lfloor nr\rfloor}_{k=1}\left\vert a_{k}\right\vert ^{\alpha}\left|\text{sgn}\left(\frac{1}{\gamma_{n}}a_{k}\right)-\text{sgn}\left(\rho^{k}_{n}c_{\infty}\right)\right|
\]
\[
=\kappa^{-1}_{n}\sum_{k\in K_{n,k}(r)}\left\vert a_{k}\right\vert ^{\alpha}\left|\text{sgn}\left(\frac{1}{\gamma_{n}}a_{k}\right)-\text{sgn}\left(\rho^{k}_{n}c_{\infty}\right)\right|
\]
\[
\leq_{(2)}2\kappa^{-1}_{n}\sum^{n}_{k=1}\left|a_{k}-\rho^{k}_{n}c_{\infty}\right|^{\alpha}=_{(3)}o(1),
\]
where $\leq_{(2)}$ follows from the fact that
\[
\left|X\right|^{\alpha}\left|\text{sgn}\left(X\right)-\text{sgn}\left(Y\right)\right|\leq2\left|X-Y\right|^{\alpha},
\]
while $=_{(3)}$ follows from part (a). $\blacksquare$

\bigskip{}

\begin{lemma}\textbf{\label{PropFR1}}\textit{Let $\delta_{1}>0$,
$\delta_{2}>1$, and $\delta:=\delta_{2}-\delta_{1}>0$. Suppose that
}\textbf{\textit{Assumption BFR1}}\textit{ holds. Define 
\[
I_{t}(k):=\int^{k}_{k-\delta_{1}t}x^{-\left(1+\frac{1}{\alpha}\right)}dx=\alpha\left[\left(k-\delta_{1}t\right)^{-1/\alpha}-k^{-1/\alpha}\right],
\]
and
\[
Q_{t}(k):=\int^{k+\delta_{1}t}_{k}x^{-\left(1+\frac{1}{\alpha}\right)}dx=\alpha\left[k^{-1/\alpha}-\left(k+\delta_{1}t\right)^{-1/\alpha}\right].
\]
(a) Fix $\eta\in\left(0,1\right)$. Then for some $t_{0}\in\mathbb{N}$
large enough, depending on $\eta$, we have}

\bigskip{}

(i) $\left\vert \sum^{k}_{j=k-\delta_{1}t}c_{j}/I_{t}(k)-1\right\vert <\eta$\textit{,
for }$k\geq\delta_{2}\left(t+t_{0}\right)+1$;\bigskip{}

(ii) $\left\vert \sum^{k+\delta_{1}t}_{j=k}c_{j}/Q_{t}(k)-1\right\vert <\eta$\textit{,
for }$\delta t+\delta_{2}t_{0}\geq k\geq t_{0}$\textit{ and }$t\geq t_{0}$;

\bigskip{}

\noindent\textit{(b) Let }$t_{0}\in\mathbb{N}$. \textit{As} $t\rightarrow\infty$,
\[
\inf_{k\geq\delta_{2}\left(t+t_{0}\right)+1}\gamma_{t}I_{t}(k)^{-1},\inf_{\delta t+\delta_{2}t_{0}\geq k\geq t_{0}}\gamma_{t}Q_{t}(k)^{-1}\rightarrow\infty\text{.}
\]

\end{lemma}

\bigskip{}

\textbf{Proof of Lemma \ref{PropFR1}. }Let $0<\eta<1$ and choose
$t_{0}$ large enough such that 
\begin{equation}
\left\vert c_{j}j^{\left(1+\frac{1}{\alpha}\right)}-1\right\vert ,\frac{3}{2}\left\{ \alpha j\left[1-\left(\frac{2\delta_{2}-\delta_{1}}{2\delta_{2}}\right)^{1/\alpha}\right]\right\} ^{-1},\frac{3}{2\delta_{1}}\left(\frac{\delta_{2}}{\delta_{2}-\delta_{1}}\right)^{^{\left(1+\frac{1}{\alpha}\right)}}j^{-1}<\eta/2,\label{FMtypeI.3c}
\end{equation}
\ for $j\geq t_{0}$. Further, set $t^{\prime}:=t+t_{0}$. Euler
summation yields 
\begin{equation}
\left\vert \sum^{k}_{j=k-\delta_{1}t}j^{-\left(1+\frac{1}{\alpha}\right)}-I_{t}(k)\right\vert \leq\left(k-\delta_{1}t\right)^{-\left(1+\frac{1}{\alpha}\right)}\text{ and }\left\vert \sum^{k+\delta_{1}t}_{j=k}j^{-\left(1+\frac{1}{\alpha}\right)}-Q_{t}(k)\right\vert \leq k^{-\left(1+\frac{1}{\alpha}\right)}.\label{PropFR1.1}
\end{equation}
In view of (\ref{FMtypeI.3c}), and (\ref{PropFR1.1}) for $k\geq\delta_{2}t^{\prime}+1$
we have
\[
\left\vert \sum^{k}_{j=k-\delta_{1}t}c_{j}-I_{t}(k)\right\vert \leq\sum^{k}_{j=k-\delta_{1}t}\left\vert c_{j}-j^{-\left(1+\frac{1}{\alpha}\right)}\right\vert +\left\vert \sum^{k}_{j=k-\delta_{1}t}j^{-\left(1+\frac{1}{\alpha}\right)}-I_{t}(k)\right\vert 
\]
\[
\leq\sum^{k}_{j=k-\delta_{1}t}\left\vert c_{j}j^{\left(1+\frac{1}{\alpha}\right)}-1\right\vert j^{-\left(1+\frac{1}{\alpha}\right)}+\left(k-\delta_{1}t\right)^{-\left(1+\frac{1}{\alpha}\right)}\leq\left(\eta/2\right)\sum^{k}_{j=k-\delta_{1}t}j^{-\left(1+\frac{1}{\alpha}\right)}+\left(k-\delta_{1}t\right)^{-\left(1+\frac{1}{\alpha}\right)}
\]
\[
=\left(\eta/2\right)\left(\sum^{k}_{j=k-\delta_{1}t}j^{-\left(1+\frac{1}{\alpha}\right)}-I_{t}(k)\right)+\left(\eta/2\right)I_{t}(k)+\left(k-\delta_{1}t\right)^{-\left(1+\frac{1}{\alpha}\right)}
\]
\[
\leq\left(1+\eta/2\right)\left(k-\delta_{1}t\right)^{-\left(1+\frac{1}{\alpha}\right)}+\left(\eta/2\right)I_{t}(k)
\]
Therefore, for $k\geq\delta_{2}t^{\prime}+1$,
\[
\left\vert \sum^{k}_{j=k-\delta_{1}t}c_{j}/I_{t}(k)-1\right\vert \leq\left(1+\eta/2\right)\left(k-\delta_{1}t\right)^{-\left(1+\frac{1}{\alpha}\right)}/I_{t}(k)+\eta/2
\]
\begin{equation}
=\left(1+\eta/2\right)\frac{\left(k-\delta_{1}t\right)^{-\left(1+\frac{1}{\alpha}\right)}}{\alpha\left[\left(k-\delta_{1}t\right)^{-1/\alpha}-k^{-1/\alpha}\right]}+\eta/2.\label{PropFR1.2}
\end{equation}
Consider the term
\[
\frac{\left(k-\delta_{1}t\right)^{-\left(1+\frac{1}{\alpha}\right)}}{\alpha\left[\left(k-\delta_{1}t\right)^{-1/\alpha}-k^{-1/\alpha}\right]}=\frac{k^{-\left(1+\frac{1}{\alpha}\right)}\left(1-\delta_{1}t/k\right)^{-\left(1+\frac{1}{\alpha}\right)}}{\alpha k^{-1/\alpha}\left[\left(1-\delta_{1}t/k\right)^{-1/\alpha}-1\right]}
\]
\[
=\frac{k^{-1}\left(1-\delta_{1}t/k\right)^{-\left(1+\frac{1}{\alpha}\right)}}{\alpha\left[\left(1-\delta_{1}t/k\right)^{-1/\alpha}-1\right]}\leq_{\left(1\right)}k^{-1}\frac{\alpha\left(1-\delta_{1}t/k\right)^{-\left(1+\frac{1}{\alpha}\right)}}{\alpha\delta_{1}t/k}
\]
\[
=\frac{1}{\left(1-\delta_{1}t/k\right)^{\left(1+\frac{1}{\alpha}\right)}\delta_{1}t}\leq\frac{1}{\left[1-\delta_{1}t/\left(\delta_{2}t^{\prime}+1\right)\right]^{\left(1+\frac{1}{\alpha}\right)}\delta_{1}t}
\]
\begin{equation}
\leq\frac{1}{\left[1-\delta_{1}t/\delta_{2}t\right]^{\left(1+\frac{1}{\alpha}\right)}\delta_{1}t}\leq\left(\frac{\delta_{2}}{\delta_{2}-\delta_{1}}\right)^{^{\left(1+\frac{1}{\alpha}\right)}}\frac{1}{\delta_{1}t},\label{PropFR1.3}
\end{equation}
where $\leq_{\left(1\right)}$ follows from the mean value theorem\footnote{Let $f(y)=(1-y)^{-1/\alpha}$. Then for $0\leq\bar{y}\leq y<1$, 
\[
f(y)-f(0)=\frac{1}{\alpha}\frac{y}{\left(1-\bar{y}\right)^{1+\frac{1}{\alpha}}}\geq\frac{1}{\alpha}y.
\]
Setting $y=t/k$, the above yields
\[
(1-t/k)^{-1/\alpha}-1\geq\frac{1}{\alpha}t/k.
\]
}. In view of (\ref{PropFR1.2})-(\ref{PropFR1.3}), for $k\geq\delta_{2}t^{\prime}+1$
we get
\[
\left\vert \sum^{k}_{j=k-\delta_{1}t}c_{j}/I_{t}(k)-1\right\vert \leq\frac{1}{\delta_{1}t}\left(\frac{\delta_{2}}{\delta_{2}-\delta_{1}}\right)^{^{\left(1+\frac{1}{\alpha}\right)}}\left(1+\eta/2\right)+\eta/2
\]
\[
\leq\frac{1}{\delta_{1}}\left(\frac{\delta_{2}}{\delta_{2}-\delta_{1}}\right)^{^{\left(1+\frac{1}{\alpha}\right)}}\left(1+\eta/2\right)/t_{0}+\eta/2\leq\frac{\eta}{3}\left(1+\eta/2\right)+\eta/2\leq\eta,
\]
and this shows part (a.i).

Next, we show part (ii). In view of (\ref{FMtypeI.3c}) and (\ref{PropFR1.1})
for $\delta t+\delta_{2}t_{0}\geq k\geq t_{0}$ we have
\[
\left\vert \sum^{k+\delta_{2}t}_{j=k}c_{j}-Q_{t}(k)\right\vert \leq\sum^{k+\delta_{2}t}_{j=k}\left\vert c_{j}-j^{-\left(1+\frac{1}{\alpha}\right)}\right\vert +\left\vert \sum^{k+\delta_{2}t}_{j=k}j^{-\left(1+\frac{1}{\alpha}\right)}-Q_{t}(k)\right\vert 
\]
\[
\leq\sum^{k+\delta_{2}t}_{j=k}\left\vert c_{j}j^{\left(1+\frac{1}{\alpha}\right)}-1\right\vert j^{-\left(1+\frac{1}{\alpha}\right)}+k^{-\left(1+\frac{1}{\alpha}\right)}\leq\left(\eta/2\right)\sum^{k+\delta_{2}t}_{j=k}j^{-\left(1+\frac{1}{\alpha}\right)}+k^{-\left(1+\frac{1}{\alpha}\right)}
\]
\[
=\left(\eta/2\right)\left(\sum^{k+\delta_{2}t}_{j=k}j^{-\left(1+\frac{1}{\alpha}\right)}-Q_{t}(k)\right)+\left(\eta/2\right)Q_{t}(k)+k^{-\left(1+\frac{1}{\alpha}\right)}
\]
\[
\leq\left(1+\eta/2\right)k^{-\left(1+\frac{1}{\alpha}\right)}+\left(\eta/2\right)Q_{t}(k).
\]
Hence,
\[
\left\vert \frac{\sum^{k+\delta_{2}t}_{j=k}c_{j}}{Q_{t}(k)}-1\right\vert =\left\vert \frac{\sum^{k+\delta_{2}t}_{j=k}c_{j}-Q_{t}(k)}{Q_{t}(k)}\right\vert 
\]
\begin{equation}
\leq\frac{\left(1+\eta/2\right)k^{-\left(1+\frac{1}{\alpha}\right)}+\left(\eta/2\right)Q_{t}(k)}{Q_{t}(k)}\leq\frac{3k^{-\left(1+\frac{1}{\alpha}\right)}}{2Q_{t}(k)}+\eta/2.\label{PropFR1.4}
\end{equation}
Next, for $\delta t+\delta_{2}t_{0}\geq k\geq t_{0}$ and $t\geq t_{0}$
\[
\frac{k^{-\left(1+\frac{1}{\alpha}\right)}}{Q_{t}(k)}=\frac{k^{-\left(1+\frac{1}{\alpha}\right)}}{\alpha\left[k^{-1/\alpha}-\left(k+\delta_{1}t\right)^{-1/\alpha}\right]}=\frac{k^{-\left(1+\frac{1}{\alpha}\right)}}{\alpha k^{-1/\alpha}\left[1-\left(1+\delta_{1}t/k\right)^{-1/\alpha}\right]}
\]
\[
=\frac{1}{\alpha k\left[1-\left\{ 1+\delta_{1}t/k\right\} ^{-1/\alpha}\right]}
\]
\[
\overset{k=\delta t+\delta_{2}t_{0}}{\leq}\frac{1}{\alpha k\left[1-\left\{ 1+\delta_{1}t/\left(\delta t+\delta_{2}t_{0}\right)\right\} ^{-1/\alpha}\right]}\overset{k,t=t_{0}}{\leq}\frac{1}{\alpha t_{0}\left[1-\left\{ 1+\delta_{1}t_{0}/\left(\delta t_{0}+\delta_{2}t_{0}\right)\right\} ^{-1/\alpha}\right]}
\]
\[
=\frac{1}{\alpha t_{0}\left[1-\left\{ 1+\delta_{1}/\left(\delta+\delta_{2}\right)\right\} ^{-1/\alpha}\right]}=\frac{1}{\alpha t_{0}\left[1-\left\{ \frac{\delta+\delta_{2}}{\delta+\delta_{2}}+\frac{\delta_{1}}{\delta+\delta_{2}}\right\} ^{-1/\alpha}\right]}
\]
\begin{equation}
=\frac{1}{\alpha t_{0}\left[1-\left(\frac{\delta+\delta_{2}+\delta_{1}}{\delta+\delta_{2}}\right)^{-1/\alpha}\right]}=\frac{1}{\alpha t_{0}\left[1-\left(\frac{2\delta_{2}-\delta_{1}}{2\delta_{2}}\right)^{1/\alpha}\right]}<\eta/3,\label{PropFR1.5}
\end{equation}
where the last inequality follows from (\ref{FMtypeI.3c}).

Using (\ref{PropFR1.4})-(\ref{PropFR1.5}), for $\delta t+\delta_{2}t_{0}\geq k\geq t_{0}$
and $t\geq t_{0}$

\[
\left\vert \frac{\sum^{k+\delta_{1}t}_{j=k}c_{j}}{Q_{t}(k)}-1\right\vert <\eta,
\]
as required.

Finally, we show part (b). First, consider $I_{t}(k)^{-1}$, $k\geq\delta_{2}t^{\prime}+1$.
We have
\[
I_{t}(k)=\alpha\left[\left(k-\delta_{1}t\right)^{-1/\alpha}-k^{-1/\alpha}\right]\leq\alpha k^{-1/\alpha}\left[\left(1-t/k\right)^{-1/\alpha}-1\right]
\]
\[
\leq_{\left(2\right)}k^{-1/\alpha}(1-t/k)^{-(1+1/\alpha)}\delta_{1}t/k\leq(1-t/k)^{-(1+1/\alpha)}t/k^{1+1/\alpha}
\]
\[
\leq\left(1-t/\left(\delta_{2}t\right)\right){}^{-(1+1/\alpha)}t/\left(\delta_{2}t\right)^{1+1/\alpha}
\]
\[
=(1-1/\delta_{2})^{-(1+1/\alpha)}(1/\delta_{2})^{(1+1/\alpha)}/t^{1/\alpha}=\left[\delta_{2}(1-1/\delta_{2})\right]^{-(1+1/\alpha)}/t^{1/\alpha},
\]
where $\leq_{\left(2\right)}$ follows from the mean value theorem\footnote{Set $f(u)=(1-u)^{-1/\alpha}$. Then for $0<\bar{u}\leq u<1$, 
\[
0<(1-u)^{-1/\alpha}-1=f(u)-f(0)=f^{\prime}(\bar{u})u=\frac{1}{\alpha}(1-\bar{u})^{-(1+1/\alpha)}u
\]
\[
\leq\frac{1}{\alpha}(1-u)^{-(1+1/\alpha)}u.
\]

{}}. In view of this, as $t\rightarrow\infty$
\[
\inf_{k\geq\delta_{2}t^{\prime}+1}\gamma_{t}I_{t}(k)^{-1}\geq\left[\delta_{2}(1-1/\delta_{2})\right]^{(1+1/\alpha)}\gamma_{t}t^{1/\alpha}\rightarrow\infty\text{.}
\]

Next, consider $Q_{t}(k)^{-1}$, $\delta t+\delta_{2}t_{0}\geq k\geq t_{0}$.
Notice that $Q_{t}(k)^{-1}$ is increasing in $k$.\footnote{Note that
\[
\frac{\partial Q_{t}(k)}{\partial k}=\alpha\frac{\partial\left(k^{-1/\alpha}-\left(k+\delta_{1}t\right)^{-1/\alpha}\right)}{\partial k}=\left(k+\delta_{1}t\right)^{-\frac{1}{\alpha}\left(\alpha+1\right)}-k^{-\frac{1}{\alpha}\left(\alpha+1\right)}<0.
\]
} Hence, as $t\rightarrow\infty$ 
\[
Q_{t}(k)^{-1}\geq Q_{t}(t_{0})^{-1}=\left[\frac{\alpha}{t^{-1/\alpha}_{0}-\left(t_{0}+t\right)^{-1/\alpha}}\right]^{-1}
\]
\[
=\frac{t^{-1/\alpha}_{0}-\left(t_{0}+t\right)^{-1/\alpha}}{\alpha}\rightarrow t^{-1/\alpha}_{0}/\alpha.
\]
Therefore, as $t\rightarrow\infty$,
\[
\inf_{\delta t+\delta_{2}t_{0}\geq k\geq t_{0}}\gamma_{t}Q_{t}(k)^{-1}\geq\gamma_{t}Q_{t}(t_{0})^{-1}\rightarrow\infty,
\]
as required. $\blacksquare$

\bigskip{}
\newpage{}

\begin{lemma} \label{bracketing}\textit{Suppose that }\textbf{\textit{Assumption
BFR2}}\textit{ holds and $a_{k}=\sum^{k}_{j=0}\phi_{j}c_{k-j}$. There
is some $t_{0}>0$ and constants $0<\underline{D}\leq\overline{D}<\infty$\ such
that}
\[
\underline{D}\leq\inf_{t_{0}\leq k}k^{1/\alpha}\left\vert a_{k}\right\vert \leq\sup_{0\leq k}k^{1/\alpha}\left\vert a_{k}\right\vert \leq\overline{D}.
\]
\end{lemma}

\bigskip{}

\textbf{Proof of Lemma \ref{bracketing}.} It follows directly from
\textbf{Lemma \ref{cauchy}} and the fact that $\sum^{\infty}_{j=0}c_{j}\neq0$.
$\blacksquare$

\bigskip{}

\begin{lemma} \label{cauchy}\textit{The following hold.}

\textit{Suppose that }$\phi_{j}\sim cons.\cdot j^{-1/\alpha}$\textit{,
$\alpha\in(0,2]$,} $\sum^{\infty}_{j=0}\left\vert c_{j}\right\vert <\infty$,
$\sum^{\infty}_{j=0}c_{j}\neq0$, \textit{and one of the following
holds:}
\begin{itemize}
\item[(a)] $\lim\sup_{j\rightarrow\infty}\bar{c}_{j}\left\vert c_{j}\right\vert <\infty$,
\textit{where $\bar{c}_{j}:=\varsigma(j)j^{\frac{1}{\alpha\wedge1}}$
with $\varsigma(x)$ SV such that}
\end{itemize}
$\begin{array}{l}
\lim_{x\rightarrow\infty}\varsigma(x)^{-1}=0,\:\mathit{if}\:\alpha<1;\\
\lim_{x\rightarrow\infty}\varsigma(x)^{-1}\ln(x)=0,\:\mathit{if}\:\alpha=1;\\
\varsigma(x)=1,\:\mathit{if}\:\alpha\in(1,2].
\end{array}$ 
\begin{itemize}
\item[(b)] $\sum^{\infty}_{j=1}j^{1/\alpha}\left\vert c_{j}\right\vert <\infty$.
\end{itemize}
\noindent\textit{Then as} $s\rightarrow\infty$

\noindent
\[
\sum^{s}_{j=0}\phi_{j}c_{s-j}=cons.\cdot s^{-1/\alpha}\sum^{\infty}_{j=0}c_{j}+o\left(s^{-1/\alpha}\right).
\]
\end{lemma}

\noindent\textbf{Remark}. \textbf{Lemma \ref{cauchy}} implies that
as $s\rightarrow\infty$
\[
\sum^{s}_{j=0}\phi_{j}c_{s-j}\sim cons.\cdot s^{-1/\alpha}\sum^{\infty}_{j=0}c_{j}\text{.}
\]

\bigskip{}

\textbf{Proof of Lemma \ref{cauchy}}. Without loss of generality,
suppose that $cons.=1$. Set $\lambda\left(s\right):=s^{-1/\alpha}$.
Note that for all $\varepsilon>0$, there some $N_{\varepsilon}$
such that $\left\vert \phi_{s-l}/\lambda\left(s-l\right)-1\right\vert <\varepsilon$
when $s-l>N_{\varepsilon}$ . We have

\begin{eqnarray}
\sum^{s}_{j=0}\phi_{j}c_{s-j} & = & \sum^{s}_{l=0}\phi_{s-l}c_{l}=\sum^{s-1-N_{\varepsilon}}_{l=0}\lambda\left(s-l\right)c_{l}+S_{\varepsilon}(s)\label{Cauchy1}\\
 &  & \left[l=s-j\right]\nonumber 
\end{eqnarray}
where 
\[
S_{\varepsilon}(s):=\sum^{s-1-N_{\varepsilon}}_{l=0}\left[\phi_{s-l}/\lambda\left(s-l\right)-1\right]\lambda\left(s-l\right)c_{l}+\sum^{s}_{l=s-N_{\varepsilon}}\phi_{s-l}c_{l}
\]
We shall show that as $s\rightarrow\infty$
\begin{equation}
\sum^{s-1-N_{\varepsilon}}_{l=0}\lambda\left(s-l\right)c_{l}=\lambda\left(s\right)\sum^{\infty}_{j=0}c_{j}+o\left(\lambda\left(s\right)\right),\label{Cauchy2}
\end{equation}
and as $s\rightarrow\infty$ first and then $\varepsilon\rightarrow0$
\begin{equation}
\lambda^{-1}\left(s\right)S_{\varepsilon}(s)=o\left(1\right).\label{Cauchy3}
\end{equation}
The requisite result follows from (\ref{Cauchy1})-(\ref{Cauchy3}).

We start with the proof of (\ref{Cauchy3}). As $s\rightarrow\infty$
first and then $\varepsilon\rightarrow0$, the first component of
$S_{\varepsilon}(s)$ is
\[
\sum^{s-1-N_{\varepsilon}}_{l=0}\left\vert \phi_{s-l}/\lambda\left(s-l\right)-1\right\vert \left\vert \lambda\left(s-l\right)c_{l}\right\vert \leq\varepsilon\sum^{s-1-N_{\varepsilon}}_{l=0}\left\vert \lambda\left(s-l\right)c_{l}\right\vert 
\]
\[
=_{(1)}\varepsilon\left[\left\vert \lambda\left(s\right)\right\vert \sum^{\infty}_{l=0}\left\vert c_{l}\right\vert +o\left(\lambda\left(s\right)\right)\right]=o\left(\lambda\left(s\right)\right),
\]
where $=_{(1)}$follows from (\ref{Cauchy2}) applied $\left|c_{j}\right|$
in place of $c_{j}$.

Next, we consider the second component of $S_{\varepsilon}(s)$. First,
suppose that condition (a) holds. Then
\[
\left\vert \sum^{s}_{l=s-N_{\varepsilon}}\phi_{s-l}c_{l}\right\vert \leq\sup_{0\leq j\leq N_{\varepsilon}}\left\vert \phi_{j}\right\vert \sum^{s}_{l=s-N_{\varepsilon}}\left\vert c_{l}\right\vert \leq C\sup_{0\leq j\leq N_{\varepsilon}}\left\vert \phi_{j}\right\vert \sum^{s}_{l=s-N_{\varepsilon}}\left\vert \bar{c}_{l}\right\vert ^{-1}
\]
\[
=C\sum^{s}_{l=s-N_{\varepsilon}}\left\vert \varsigma(l)l^{\frac{1}{\alpha\wedge1}}\right\vert ^{-1}.
\]
Set $\bar{R}\left(l\right)=\left[\varsigma(l)l^{\frac{1}{\alpha\wedge1}}\right]^{-1}$.
For any $\eta>0$ and $s$ large enough, the sum term above is 
\[
\sum^{s}_{l=s-N_{\varepsilon}}\left\vert \varsigma(l)l^{^{\frac{1}{\alpha\wedge1}}}\right\vert ^{-1}=\bar{R}\left(s\right)\sum^{s}_{l=s-N_{\varepsilon}}\frac{\bar{R}\left(l\right)}{\bar{R}\left(s\right)}
\]
\[
=\bar{R}\left(s\right)\sum^{s}_{l=s-N_{\varepsilon}}\left[\frac{\bar{R}\left(\frac{l}{s}s\right)}{\bar{R}\left(s\right)}-\left(\frac{l}{s}\right)^{-\frac{1}{\alpha\wedge1}}\right]+\bar{R}\left(s\right)\sum^{s}_{l=s-N_{\varepsilon}}\left(\frac{l}{s}\right)^{-\frac{1}{\alpha\wedge1}}
\]
\[
\leq\bar{R}\left(s\right)\sum^{s}_{l=s-N_{\varepsilon}}\sup_{1-\frac{N_{\varepsilon}}{s}\leq y\leq1}\left[\frac{\bar{R}\left(ys\right)}{\bar{R}\left(s\right)}-y^{-\frac{1}{\alpha\wedge1}}\right]+N_{\varepsilon}\left(\frac{s}{s-N_{\varepsilon}}\right)^{1/\alpha}\bar{R}\left(s\right)
\]
\[
\leq_{\left(2\right)}\bar{R}\left(s\right)N_{\varepsilon}\left[\eta+\left(\frac{s}{s-N_{\varepsilon}}\right)^{\frac{1}{\alpha\wedge1}}\right]=C\frac{\bar{R}\left(s\right)}{\lambda(s)}\lambda(s)
\]
\[
=\frac{C}{\varsigma(s)}\lambda(s)I\left\{ \alpha\in(0,1]\right\} +\frac{Cs^{1/\alpha}}{\varsigma(s)s}\lambda(s)I\left\{ \alpha\in(1,2]\right\} =o(\lambda(s)),
\]
as required. We note that $\leq_{\left(2\right)}$ follows from the
uniform convergence of regularly varying functions. Next, suppose
that condition (b) holds. In this case, for $s$ large enough, the
second component of $S_{\varepsilon}(s)$ is 
\[
\left\vert \sum^{s}_{l=s-N_{\varepsilon}}\phi_{s-l}c_{l}\right\vert \leq C\sum^{s}_{l=s-N_{\varepsilon}}\left\vert c_{l}\right\vert \leq C\left(s-N_{\varepsilon}\right)^{-1/\alpha}\sum^{\infty}_{l=s-N_{\varepsilon}}l^{1/\alpha}\left\vert c_{l}\right\vert 
\]
\[
=C\left(\frac{s}{s-N_{\varepsilon}}\right)^{1/\alpha}\lambda(s)\sum^{\infty}_{l=s-N_{\varepsilon}}l^{1/\alpha}\left\vert c_{l}\right\vert =o\left(\lambda(s)\right).
\]
Therefore, (\ref{Cauchy3}) holds either under condition (a) or (b).

Finally, we show (\ref{Cauchy2}). Set $m_{s}:=\lfloor s\delta\rfloor$
where $1/2\leq\delta<1$. Consider

\[
\sum^{s-1-N_{\varepsilon}}_{l=0}\lambda\left(s-l\right)c_{l}=\left(\sum^{m_{s}}_{l=0}+\sum^{s-1-N_{\varepsilon}}_{l=m_{s}+1}\right)\lambda\left(s-l\right)c_{l}=:T_{1}(s)+T_{2}(s)\text{.}
\]
As $s\rightarrow\infty$, the first term 
\[
\left\vert T_{1}(s)-\lambda(s)\sum^{m_{s}}_{l=0}c_{l}\right\vert \leq\sum^{m_{s}}_{l=0}\left\vert \lambda\left(s\right)-\lambda\left(s-l\right)\right\vert \left\vert c_{l}\right\vert 
\]
\[
\leq\sum^{m_{s}}_{l=0}\left\vert \frac{\lambda\left(s-l\right)}{\lambda\left(s\right)}-1\right\vert \lambda\left(s\right)\left\vert c_{l}\right\vert =\lambda\left(s\right)\sum^{m_{s}}_{l=0}q_{l,s}\left\vert c_{l}\right\vert 
\]
where $q_{l,s}:=\left\vert \left(\frac{s-l}{s}\right)^{-1/\alpha}-1\right\vert $.
Next, we show that $\sum^{m_{s}}_{l=0}q_{l,s}\left\vert c_{l}\right\vert =o(1)$
as $s\rightarrow\infty$. Note that for $\alpha\in(0,1)$ the mean
value theorem yields
\[
q_{l,s}=\frac{s^{1/\alpha}-\left(s-l\right)^{1/\alpha}}{\left(s-l\right)^{1/\alpha}}\leq\frac{1}{\alpha}\frac{s^{\frac{1}{\alpha}-1}l}{\left(s-l\right)^{1/\alpha}},
\]
whilst for $\alpha\in[1,2]$ by Loève's inequality
\[
q_{l,s}=\frac{s^{1/\alpha}-\left(s-l\right)^{1/\alpha}}{\left(s-l\right)^{1/\alpha}}\leq\frac{l^{1/\alpha}}{\left(s-l\right)^{1/\alpha}}.
\]
Therefore, as $s\rightarrow\infty$ 
\[
\sum^{m_{s}}_{l=0}q_{l,s}\left\vert c_{l}\right\vert \leq C\sum^{m_{s}}_{l=0}\left\vert c_{l}\right\vert \left[\frac{l^{1/\alpha}}{\left(s-l\right)^{1/\alpha}}1\left\{ \alpha\in[1,2]\right\} +\frac{s^{\frac{1}{\alpha}-1}l}{\left(s-l\right)^{1/\alpha}}1\left\{ \alpha\in(0,1)\right\} \right]
\]
\[
\leq\frac{1}{\left(s-m_{s}\right)^{1/\alpha}}\sum^{m_{s}}_{l=0}l^{1/\alpha}\left\vert c_{l}\right\vert +\frac{Cs^{\frac{1}{\alpha}-1}}{\left(s-m_{s}\right)^{1/\alpha}}\sum^{m_{s}}_{l=0}l\left\vert c_{l}\right\vert =_{(4)}o(1),
\]
where $=_{(4)}$ follows from Kronecker's lemma (e.g. p. 31 in Hall
and Heyde, 1980).

Next, we consider $T_{2}\left(s\right)$. First, assume that condition
(a) holds. Given that $\bar{c}_{l}$ is regularly varying of a positive
index, there is no loss of generality if we assume it is increasing
-see e.g. BGT 1987, Thm 1.5.3. Hence, as $s\rightarrow\infty$, we
have
\[
\left\vert T_{2}\left(s\right)\right\vert \leq\sum^{s-1}_{l=m_{s}}\lambda\left(s-l\right)\bar{c}_{l}{}^{-1}\bar{c}_{l}\left\vert c_{l}\right\vert 
\]
\[
\leq C\varsigma(m_{s})^{-1}m^{-\frac{1}{\alpha\wedge1}}_{s}\sum^{s-1}_{l=m_{s}}\lambda\left(s-l\right)=C\varsigma(m_{s})^{-1}m^{-\frac{1}{\alpha\wedge1}}_{s}\sum^{s-m_{s}+1}_{l=1}\lambda\left(l\right)
\]
\[
\leq C\varsigma(s)^{-1}s^{-\frac{1}{\alpha\wedge1}}\sum^{s-m_{s}+1}_{l=1}\lambda\left(l\right)=:CQ(s).
\]
The result follows if we show that $Q(s)=o(\lambda(s))$. Note that
for $\alpha\in(0,1)$, $\varsigma(s)^{-1}s^{-\frac{1}{\alpha\wedge1}}=o(\lambda(s))$
and $\sum^{\infty}_{l=1}\lambda\left(l\right)<\infty$. Hence, $Q(s)=o(\lambda(s))$.
For $\alpha=1$, $\sum^{s-m_{s}+1}_{l=1}\lambda\left(l\right)\sim\ln(s)$.
Hence,
\[
Q(s)\sim\frac{\lambda(s)}{\varsigma(s)}\ln(s)=_{(5)}o(\lambda(s)),
\]
where $=_{(5)}$ follows due to the requirement of condition (a) for
$\alpha=1$. For $\alpha\in(1,2]$, using Euler summation (see e.g.
\textbf{Lemma \ref{LemmaSV1}}) we get as $s\rightarrow\infty$ first,
and then as $\delta\rightarrow1$
\[
Q(s)=s^{-1}\sum^{s-m_{s}+1}_{l=1}l^{-1/\alpha}=\lambda(s)s^{-1}\sum^{s-m_{s}+1}_{l=1}\left(\frac{l}{s}\right)^{-1/\alpha}
\]
\[
=\lambda(s)\int^{1-\delta}_{0}x^{-1/\alpha}dx\left(1+o(1)\right)=o\left(\lambda(s)\right).
\]

Finally, we consider $T_{2}(s)$ under condition (b). As $s\rightarrow\infty$,

\[
\left\vert T_{2}(s)\right\vert \leq\sup_{l>N_{\varepsilon}}\lambda\left(l\right)\sum^{s-1-N_{\varepsilon}}_{l=m_{s}}\left\vert c_{l}\right\vert \leq C\sum^{\infty}_{l=m_{s}}l^{1/\alpha}
\]
\[
\leq C\lambda(m_{s})\sum^{\infty}_{l=m_{s}}l^{1/\alpha}\left\vert c_{l}\right\vert =o(\lambda(s)),
\]
as required. $\blacksquare$

\bigskip{}

\newpage{}

\noindent\begin{lemma}\label{LPA}\textbf{ }(linear process approximation).
\textit{Suppose that }\textbf{\textit{Assumptions BFR2}}\textit{ and
}\textbf{\textit{INN}}\textit{ hold. Then for all }$\alpha^{\prime}\in(0,\alpha)$\textit{,
as }$n\rightarrow\infty$ 
\[
\sup_{0\leq r\leq1}\left\Vert \gamma^{-1}_{n}x_{\lfloor nr\rfloor}-\gamma^{-1}_{n}\tilde{x}_{\lfloor nr\rfloor}\right\Vert _{L_{\alpha^{\prime}}}=o(1),
\]
\textit{where} 
\[
\tilde{x}_{n}:=\left(\sum^{\infty}_{j=0}c_{j}\right)\sum^{n-1}_{k=1}k^{-\frac{1}{\alpha}}\xi_{n-k}.
\]

\end{lemma}

\bigskip{}

\noindent\textbf{Proof of Lemma \ref{LPA}}. Let $\varepsilon>0$
such that $\alpha-\varepsilon>0$. We shall show that for each $0<r\leq1$
\begin{equation}
\mathbf{E}\left\vert \gamma^{-1}_{n}x^{-}_{\lfloor nr\rfloor}\right\vert ^{\alpha-\varepsilon}\rightarrow0\label{LPA1}
\end{equation}
and
\begin{equation}
\mathbf{E}\left\vert \gamma^{-1}_{n}x^{+}_{\lfloor nr\rfloor}-\gamma^{-1}_{n}\tilde{x}_{\lfloor nr\rfloor}\right\vert ^{\alpha-\varepsilon}\rightarrow0.\label{LPA2}
\end{equation}
We start with the proof of (\ref{LPA1}). Without loss of generality
set $r=1$, and assume that $R_{v-\alpha}(x)=x^{\nu-\alpha}H(x)$
is decreasing, where $\nu=1$ for $\alpha\in(1,2)$ and $\nu\in(0,\alpha)$
for and $\alpha\in(0,1]$. Moreover, recall that under our assumptions
$\sum^{\infty}_{l=0}\left\vert c_{l}\right\vert ^{\nu}<\infty$. In
addition, to simplify the notation assume that $\sum^{\infty}_{l=0}c_{l}=1$.
Using arguments similar to those used for establishing (\ref{a-moments2ca})
and (\ref{a-moments3d}), we have
\[
\mathbf{E}\left\vert \gamma^{-1}_{n}x^{-}_{n}\right\vert ^{\alpha-\varepsilon}\leq C\gamma^{-\nu}_{n}\sum^{n}_{j=0}\left\vert \phi_{j}\right\vert ^{\nu}R_{\nu-\alpha}\left(\gamma_{n}j^{1/\alpha}\right)\sum^{\infty}_{l=n-j}\left\vert c_{l}\right\vert ^{\nu}
\]
\[
\leq C\gamma^{-\alpha}_{n}\sum^{n}_{j=1}j^{-1}H\left(\gamma_{n}j^{1/\alpha}\right)\sum^{\infty}_{l=n-j}\left\vert c_{l}\right\vert ^{\nu}
\]
\[
=C\gamma^{-\alpha}_{n}\left(\sum^{n}_{j=\lfloor n/2\rfloor}+\sum^{\lfloor n/2\rfloor-1}_{j=1}\right)j^{-1}H\left(\gamma_{n}j^{1/\alpha}\right)\sum^{\infty}_{l=n-j}\left\vert c_{l}\right\vert ^{\nu}
\]
\[
\leq C\gamma^{-\alpha}_{n}\sum^{n}_{j=\lfloor n/2\rfloor}j^{-1}H\left(\gamma_{n}j^{1/\alpha}\right)\sum^{\infty}_{l=0}\left\vert c_{l}\right\vert ^{\nu}
\]
\[
+C\gamma^{-\alpha}_{n}\sum^{n}_{j=1}j^{-1}H\left(\gamma_{n}j^{1/\alpha}\right)\sum^{\infty}_{l=\lfloor n/2\rfloor}\left\vert c_{l}\right\vert ^{\nu}
\]
\[
=_{(1)}C\gamma^{-\alpha}_{n}\left[\int^{\gamma_{n}n^{\frac{1}{\alpha}}}_{\gamma_{n}\lfloor n/2\rfloor^{\frac{1}{\alpha}}}\frac{H\left(x\right)}{x}dx+\int^{\gamma_{n}n^{\frac{1}{\alpha}}}_{1}\frac{H\left(x\right)}{x}dx\sum^{\infty}_{l=\lfloor n/2\rfloor}\left\vert c_{l}\right\vert ^{\nu}\right]\left[1+o(1)\right]=_{(2)}o(1),
\]
where $=_{(1)}$ follows from \textbf{Lemma \ref{LemmaSV1}(i)}, and
$=_{(2)}$ from the fact that $\gamma^{-\alpha}_{n}\int^{\gamma_{n}n^{\frac{1}{\alpha}}}_{\gamma_{n}\lfloor n/2\rfloor^{\frac{1}{\alpha}}}\frac{H\left(x\right)}{x}dx=o(1)$.

Next, we show (\ref{LPA2}). Again, without loss of generality set
$r=1$. It can be easily checked that as $n\rightarrow\infty$
\[
\inf_{1\leq k\leq n}S_{n,k}:=\inf_{1\leq k\leq n}\gamma_{n}\left(k^{-1/\alpha}\right)^{-1}\rightarrow\infty.
\]
Set
\[
\overline{\xi}_{n,k}:=\xi_{k}1\left\{ \left\vert \xi_{k}\right\vert \geq S_{n,k}\right\} \text{ \ and \ }\underline{\xi}_{n,k}:=\xi_{k}1\left\{ \left\vert \xi_{k}\right\vert <S_{n,k}\right\} .
\]
Hence, by Lemma 1 in Astrauskas (1983) 
\[
\mathbf{E}\left\vert \overline{\xi}_{n,k}\right\vert ^{\alpha-\varepsilon}\asymp S^{-\varepsilon}_{n,k}h\left(S_{n,k}\right)\text{ \ and \ }\mathbf{E}\left\vert \underline{\xi}_{n,k}\right\vert ^{\alpha+\varepsilon}\asymp S^{\varepsilon}_{n,k}h\left(S_{n,k}\right),
\]
as $n\rightarrow\infty$. Let $r_{k}:=\left(\frac{a_{k}}{k^{-1/\alpha}}-1\right)$.
Due to \textbf{Lemma \ref{cauchy}}, for $\epsilon>0$ arbitrarily
small we have $\left\vert r_{k}\right\vert <\epsilon$ for all $k$
large enough, $k\geq N_{\epsilon}$ say. Consider
\[
\mathbf{E}\left\vert \gamma^{-1}_{n}x^{+}_{n}-\gamma^{-1}_{n}\tilde{x}_{n}\right\vert ^{\alpha-\varepsilon}=\mathbf{E}\left\vert \gamma^{-1}_{n}\sum^{n-1}_{k=0}a_{k}\xi_{n-k}-\gamma^{-1}_{n}\sum^{n-1}_{k=1}k^{-\frac{1}{\alpha}}\xi_{n-k}\right\vert ^{\alpha-\varepsilon}
\]
\begin{equation}
\leq\left\vert \gamma^{-1}_{n}\sum^{n-1}_{k=N_{\epsilon}}\left(\frac{a_{k}}{k^{-\frac{1}{\alpha}}}-1\right)k^{-\frac{1}{\alpha}}\xi_{n-k}\right\vert ^{\alpha-\varepsilon}=\mathbf{E}\left\vert \gamma^{-1}_{n}\sum^{n-1}_{k=N_{\epsilon}}r_{k}k^{-\frac{1}{\alpha}}\xi_{n-k}\right\vert ^{\alpha-\varepsilon}+o(1).\label{LPA3}
\end{equation}
Next, arguments similar to those used for (\ref{a-moments1abb}) give
\[
\mathbf{E}\left\vert \gamma^{-1}_{n}\sum^{n-1}_{k=N_{\epsilon}}r_{k}k^{-1/\alpha}\xi_{n-k}\right\vert ^{\alpha-\varepsilon}\leq2\gamma^{-\left(\alpha-\varepsilon\right)}_{n}\sum^{n-1}_{k=N_{\epsilon}}\left(\left\vert r_{k}\right\vert k^{-1/\alpha}\right)^{\alpha-\varepsilon}\mathbf{E}\left\vert \overline{\xi}_{n,k}\right\vert ^{\alpha-\varepsilon}
\]
\[
+2\left[\gamma^{-\left(\alpha+\varepsilon\right)}_{n}\sum^{n-1}_{k=N_{\epsilon}}\left(\left\vert r_{k}\right\vert k^{-1/\alpha}\right)^{\alpha+\varepsilon}\mathbf{E}\left\vert \underline{\xi}_{n,k}\right\vert ^{\alpha+\varepsilon}\right]^{\frac{\alpha-\varepsilon}{\alpha+\varepsilon}}
\]

\[
\leq2\left\{ \epsilon^{\alpha-\varepsilon}\sum^{n-1}_{k=N_{\epsilon}}S^{-\alpha+\varepsilon}_{n,k}S^{-\varepsilon}_{n,k}h\left(S_{n,k}\right)+\left[\epsilon^{\alpha+\varepsilon}\sum^{n-1}_{k=N_{\epsilon}}S^{-\alpha-\varepsilon}_{n,k}S^{\varepsilon}_{n,k}h\left(S_{n,k}\right)\right]^{\frac{\alpha-\varepsilon}{\alpha+\varepsilon}}\right\} \left[1+o(1)\right]
\]
\[
\leq2\epsilon^{\alpha-\varepsilon}\left\{ \sum^{n-1}_{k=N_{\epsilon}}S^{-\alpha}_{n,k}h\left(S_{n,k}\right)+\left[\sum^{n-1}_{k=N_{\epsilon}}S^{-\alpha}_{n,k}h\left(S_{n,k}\right)\right]^{\frac{\alpha-\varepsilon}{\alpha+\varepsilon}}\right\} \left[1+o(1)\right]
\]
\begin{equation}
\leq2\epsilon^{\alpha-\varepsilon}\left\{ \sum^{n}_{k=1}S^{-\alpha}_{n,k}H\left(S_{n,k}\right)+\left[\sum^{n}_{k=1}S^{-\alpha}_{n,k}H\left(S_{n,k}\right)\right]^{\frac{\alpha-\varepsilon}{\alpha+\varepsilon}}\right\} \left[1+o(1)\right]\leq_{(3)}C\epsilon^{\alpha-\varepsilon}\label{LPA4}
\end{equation}
where $\leq_{(3)}$ follows from the fact that (cf. \textbf{Lemma
\ref{LemmaSV1}})
\[
\sum^{n}_{k=1}S^{-\alpha}_{n,k}H\left(S_{n,k}\right)=\gamma^{-\alpha}_{n}\sum^{n}_{k=1}k^{-1}H\left(\gamma_{n}k^{1/\alpha}\right)=1+o(1).
\]
The result follows from (\ref{LPA3}), (\ref{LPA4}) and the fact
that $\epsilon$ can be made arbitrarily small. $\blacksquare$

\bigskip{}

\section{Verification of HL4}

\begin{lemma} \label{MartingApprox}S\textit{uppose the }\textbf{\textit{Assumption
INN}}\textit{ and that $x_{t}$ satisfiews either }\textbf{\textit{Assumption
BFR}}\textit{ or }\textbf{\textit{MI}}\textit{ holds. Then, for every
bounded and Lipschitz function $g:\mathbb{R}\rightarrow\mathbb{R}$,
as} $n\rightarrow\infty$
\[
\frac{1}{n}\sum^{n}_{t=1}g\left(\gamma^{-1}_{n}x_{t}\right)=\frac{1}{n}\sum^{n}_{t=1}\mathbf{E}\left[g\left(\gamma^{-1}_{n}x_{t}\right)\mid\mathcal{F}_{0}\right]+o_{p}(1).
\]
\end{lemma}

\bigskip{}

\textbf{Proof of Lemma \ref{MartingApprox}}. Let $\mathbf{E}_{t}\left[.\right]:=\mathbf{E}\left[.\mid\mathcal{F}_{t}\right]$
and as in DK2021, write 
\[
\sum^{n}_{t=1}\left[g\left(\gamma^{-1}_{n}x_{t}\right)-\mathbf{E}_{0}g\left(\gamma^{-1}_{n}x_{t}\right)\right]=\sum^{n-1}_{s=0}M_{n,s},
\]
where $M_{n,s}=\sum^{n}_{t=s+1}\left[\mathbf{E}_{t-s}g\left(\gamma^{-1}_{n}x_{t}\right)-\mathbf{E}_{\left(t-s\right)-1}g\left(\gamma^{-1}_{n}x_{t}\right)\right]$
-see DK2021, Proof of Proposition B.1(iii). Note that the summand
in $M_{n,s}$ forms a martingale difference w.r.t. $\mathcal{F}_{t-s}$. 

Let
\begin{equation}
\alpha_{\varepsilon}:=\alpha+\varepsilon\:\textrm{with }\left\{ \begin{array}{l}
\varepsilon=0,\text{ if }\alpha=2\\
\varepsilon>0\text{, if }\alpha\in(0,2)\text{ such that }\alpha_{\varepsilon}\in(1,2]
\end{array}\right..\label{MartinApprox0}
\end{equation}
Note that in the following calculations all relevant moments exist
because $g$ is bounded. Therefore, in view of this and the fact that
$\alpha_{\varepsilon}>1$ we have
\[
I_{n}:=\frac{1}{n}\mathbf{E}\left\vert \sum^{n}_{t=1}\left[g\left(\gamma^{-1}_{n}x_{t}\right)-\mathbf{E}_{0}g\left(\gamma^{-1}_{n}x_{t}\right)\right]\right\vert =\frac{1}{n}\mathbf{E}\left\vert \sum^{n-1}_{s=0}M_{n,s}\right\vert 
\]
\begin{equation}
\leq\frac{1}{n}\sum^{n-1}_{s=0}\mathbf{E}\left\vert M_{n,s}\right\vert \leq\frac{1}{n}\sum^{n-1}_{s=0}\left[\mathbf{E}\left\vert M_{n,s}\right\vert ^{\alpha_{\varepsilon}}\right]^{1/\alpha_{\varepsilon}}.\label{MartinApprox1}
\end{equation}

Next, for $0\leq s\leq n-1$ set 
\[
x^{\ast}_{t,s}:=x^{\ast}_{t}:=\sum^{\infty}_{k=0,k\neq s}a_{k,t}\xi_{t-k}+a_{s,t}\xi^{\ast}
\]
with $\xi^{\ast}\overset{d}{=}\xi_{1}$\ and $\xi^{\ast}\perp\xi_{1}$.
Then $x^{\ast}_{t}\overset{d}{=}x_{t}$. In fact, $\left\{ \mathbf{E}_{t-s}g\left(\gamma^{-1}_{n}x_{t}\right),\mathbf{E}_{\left(t-s\right)-1}g\left(\gamma^{-1}_{n}x_{t}\right)\right\} $
$\overset{d}{=}\left\{ \mathbf{E}_{t-s}g\left(\gamma^{-1}_{n}x_{t}\right),\mathbf{E}_{t-s}g\left(\gamma^{-1}_{n}x^{\ast}_{t}\right)\right\} $.
Set $\Delta_{s,n}:=\gamma^{-1}_{n}\left\vert a_{s}\right\vert \left\vert \xi_{1}-\xi^{\ast}\right\vert $.
Recall that the coefficients $a_{k,t}$, $a_{k}$ are defined in (\ref{lincoefs}).
Due to the fact that $g$ is bounded and Lipschitz, we have $\left\vert g(x)-g(z)\right\vert \leq C_{0}\left(1\wedge\left\vert x-z\right\vert \right)$
for some $C_{0}\in\left(0,\infty\right)$. In view of the above,
\[
\mathbf{E}\left\vert M_{n,s}\right\vert ^{\alpha_{\varepsilon}}=\mathbf{E}\left\vert \sum^{n}_{t=s+1}\mathbf{E}_{t-s}g\left(\gamma^{-1}_{n}x_{t}\right)-\mathbf{E}_{\left(t-s\right)-1}g\left(\gamma^{-1}_{n}x_{t}\right)\right\vert ^{\alpha_{\varepsilon}}
\]
\[
\leq2\mathbf{E}\sum^{n}_{t=s+1}\left\vert \mathbf{E}_{t-s}g\left(\gamma^{-1}_{n}x_{t}\right)-\mathbf{E}_{\left(t-s\right)-1}g\left(\gamma^{-1}_{n}x_{t}\right)\right\vert ^{\alpha_{\varepsilon}}=2\mathbf{E}\sum^{n}_{t=s+1}\left\vert \mathbf{E}_{t-s}g\left(\gamma^{-1}_{n}x_{t}\right)-\mathbf{E}_{t-s}g\left(\gamma^{-1}_{n}x^{\ast}_{t}\right)\right\vert ^{\alpha_{\varepsilon}}
\]
\[
\leq2\sum^{n}_{t=s+1}\mathbf{E}\left[\mathbf{E}_{t-s}\left\vert g\left(\gamma^{-1}_{n}x_{t}\right)-g\left(\gamma^{-1}_{n}x^{\ast}_{t}\right)\right\vert \right]^{\alpha_{\varepsilon}}\leq2\sum^{n}_{t=s+1}\mathbf{E}\left\{ \left[\mathbf{E}_{t-s}\left\vert g\left(\gamma^{-1}_{n}x_{t}\right)-g\left(\gamma^{-1}_{n}x^{\ast}_{t}\right)\right\vert ^{\alpha_{\varepsilon}}\right]^{1/\alpha_{\varepsilon}}\right\} ^{\alpha_{\varepsilon}}
\]
\[
=2\sum^{n}_{t=s+1}\mathbf{EE}_{t-s}\left\vert g\left(\gamma^{-1}_{n}x_{t}\right)-g\left(\gamma^{-1}_{n}x^{\ast}_{t}\right)\right\vert ^{\alpha_{\varepsilon}}=2\sum^{n}_{t=s+1}\mathbf{E}\left\vert g\left(\gamma^{-1}_{n}x_{t}\right)-g\left(\gamma^{-1}_{n}x^{\ast}_{t}\right)\right\vert ^{\alpha_{\varepsilon}}
\]
\[
\leq2C^{\alpha_{\varepsilon}}_{0}\sum^{n}_{t=s+1}\mathbf{E}\left\{ \left[\gamma^{-1}_{n}\left\vert a_{s,t}\right\vert \left\vert \xi_{t-s}-\xi^{\ast}\right\vert \wedge1\right]\right\} ^{\alpha_{\varepsilon}}=_{(1)}C\sum^{n}_{t=s+1}\mathbf{E}\left\{ \left[\gamma^{-1}_{n}\left\vert a_{s}\right\vert \left\vert \xi_{1}-\xi^{\ast}\right\vert \wedge1\right]\right\} ^{\alpha_{\varepsilon}}
\]
\[
=C\sum^{n}_{t=s+1}\mathbf{E}\left[\left(\Delta_{s,n}\wedge1\right)^{\alpha_{\varepsilon}}\right]=C\sum^{n}_{t=s+1}\mathbf{E}\left\{ \left[\left(\Delta_{s,n}I\left\{ \Delta_{s,n}\leq1\right\} +I\left\{ \Delta_{s,n}>1\right\} \right)\right]^{\alpha_{\varepsilon}}\right\} 
\]
\begin{equation}
\leq C\sum^{n}_{t=s+1}\mathbf{E}\left[\Delta^{\alpha_{\varepsilon}}_{s,n}I\left\{ \Delta_{s,n}\leq1\right\} +I\left\{ \Delta_{s,n}>1\right\} \right],\label{MartinApprox2}
\end{equation}
where $=_{(1)}$ holds for $0\leq s\leq t-1$. Recall that (see the
proof of \textbf{Lemma \ref{lem:FSalphaM}})
\begin{equation}
\lim_{n\rightarrow\infty}\gamma_{n}\inf_{0\leq s\leq n}\left\vert a_{s}(n)\right\vert ^{-1}=\infty,\label{MartinApprox2b}
\end{equation}
where $a_{s}:=a_{s}(n)$.

In view of \textbf{Assumption INN}, (\ref{MartinApprox2b}), and the
fact that $\xi_{1}-\xi^{\ast}$\
is a (symmetric) random variable in the domain of attraction of an
$\alpha$-stable law (see Remark 1), as $n\rightarrow\infty$
\[
\mathbf{E}\left(I\left\{ \Delta_{s,n}>1\right\} \right)=P\left(\gamma^{-1}_{n}\left\vert a_{s}\right\vert \left\vert \xi_{1}-\xi^{\ast}\right\vert >1\right)
\]
\begin{equation}
=C\left(\frac{\left\vert a_{s}\right\vert }{\gamma_{n}}\right)^{\alpha}h\left(\frac{\gamma_{n}}{\left\vert a_{s}\right\vert }\right)\left[1+o(1)\right],\label{MartinApprox3}
\end{equation}
uniformly in $0\leq s\leq n$. Next, using (\ref{MartinApprox2b})
once more and Lemma 1 of Astrauskas (1983) as $n\rightarrow\infty$
we have
\[
\mathbf{E}\left(\Delta^{\alpha_{\varepsilon}}_{s,n}I\left\{ \Delta_{s,n}\leq1\right\} \right)=\gamma^{-\alpha_{\varepsilon}}_{n}\left\vert a_{s}\right\vert ^{\alpha_{\varepsilon}}\mathbf{E}\left(\left\vert \xi_{1}-\xi^{\ast}\right\vert ^{\alpha_{\varepsilon}}I\left\{ \left\vert \xi_{1}-\xi^{\ast}\right\vert \leq\gamma_{n}\left\vert a_{s}\right\vert ^{-1}\right\} \right)
\]
\begin{equation}
=C\left(\frac{\left\vert a_{s}\right\vert }{\gamma_{n}}\right)^{\alpha}\left\{ H\left(\frac{\gamma_{n}}{\left\vert a_{s}\right\vert }\right)1\left\{ \alpha=2\right\} +h\left(\frac{\gamma_{n}}{\left\vert a_{s}\right\vert }\right)1\left\{ \alpha\in(0,2)\right\} \right\} \left[1+o(1)\right],\label{MartinApprox4}
\end{equation}
uniformly in $0\leq s\leq n$. Since $h(x)=O\left(H(x)\right)$ as
$x\rightarrow\infty$, the contributions involving $h$ can be bounded
by the corresponding expressions involving $H$, after an adjustment
of the multiplicative constant. We shall use these bounds throughout
the remainder of the proof. In view of this, and (\ref{MartinApprox1})-(\ref{MartinApprox4})
we get
\[
I_{n}\leq\frac{C}{n}\sum^{n-1}_{s=0}\left(\sum^{n}_{t=s+1}\left(\frac{\left\vert a_{s}\right\vert }{\gamma_{n}}\right)^{\alpha}H\left(\frac{\gamma_{n}}{\left\vert a_{s}\right\vert }\right)\right)^{1/\alpha_{\varepsilon}}\left[1+o(1)\right]
\]
\begin{equation}
\leq\frac{Cn^{1/\alpha_{\varepsilon}}}{n}\sum^{n}_{s=0}\left(\frac{\left\vert a_{s}\right\vert }{\gamma_{n}}\right)^{\alpha/\alpha_{\varepsilon}}H\left(\frac{\gamma_{n}}{\left\vert a_{s}\right\vert }\right)^{1/\alpha_{\varepsilon}}\left[1+o(1)\right]=:CI^{\prime}_{n}(\varepsilon,\alpha)\left[1+o(1)\right].\label{MartinApprox4b}
\end{equation}
as $n\rightarrow\infty$. Next, we shall show that, under each assumption
(i.e. under \textbf{BFR} and \textbf{MI}), there is some $\varepsilon>0$
such that 
\begin{equation}
I^{\prime}_{n}:=I^{\prime}_{n}(\varepsilon,\alpha)=o(1),\:\textrm{for}\:\alpha\in(0,2)\:\textrm{and}\:\alpha_{\varepsilon}>1,\label{MartinApprox4c-1}
\end{equation}
and
\begin{equation}
I^{\prime}_{n}(0,2)=o(1),\label{MartinApprox4c-2}
\end{equation}
 as $n\rightarrow\infty$.

\bigskip{}

\textbf{BFR (}$\alpha\in(0,2]$\textbf{)}. Using the approximation
$a_{s}\asymp s^{-1/\alpha}$ (see e.g. (\ref{FMtypeI.3})), and similar
arguments as those used for the proof of \textbf{Lemma \ref{lem:FSalphaM}},
(see e.g. (\ref{FMtypeI.3_0a})) together with \textbf{Lemma \ref{LemmaSV1}(ii)}
we have as $n\rightarrow\infty$
\[
I^{\prime}_{n}=\frac{n^{1/\alpha_{\varepsilon}}}{\gamma^{\alpha/\alpha_{\varepsilon}}_{n}n}\sum^{n}_{s=0}\left\vert a_{s}\right\vert ^{\alpha/\alpha_{\varepsilon}}H\left(\frac{\gamma_{n}}{\left\vert a_{s}\right\vert }\right)^{1/\alpha_{\varepsilon}}=\frac{n^{1/\alpha_{\varepsilon}}}{\gamma^{\alpha/\alpha_{\varepsilon}}_{n}n}\sum^{n}_{s=1}s^{-1/\alpha_{\varepsilon}}H\left(s^{1/\alpha}\gamma_{n}\right)^{1/\alpha_{\varepsilon}}\left[1+o(1)\right]
\]
\[
=\left[\frac{H\left(n^{1/\alpha}\gamma_{n}\right)}{\gamma^{\alpha}_{n}}\right]^{1/\alpha_{\varepsilon}}\int^{1}_{0}x^{-1/\alpha_{\varepsilon}}dx\left[1+o(1)\right]=_{(1)}o(1),
\]
where for $=_{(1)}$ we have used the fact that $\gamma^{\alpha}_{n}\sim L\left(n^{1/\alpha}\gamma_{n}\right)$
with $H\left(x\right)=o\left(L(x)\right)$, as $x\rightarrow\infty$.
Hence, (\ref{MartinApprox4c-1}) and (\ref{MartinApprox4c-2}) hold.\bigskip{}

\textbf{MI-SM (}$\alpha\in(0,2]$\textbf{)}. We start with a proof
under \textbf{Assumption SM(ii)}. In addition to $(\ref{MartinApprox0})$,
let
\begin{equation}
0<\delta<1-\nu\text{ and }0\leq\varepsilon<\alpha\left(\frac{1}{\delta+\nu}-1\right),\label{MartinApprox0b}
\end{equation}
where recall that $\nu\in(0,\alpha\wedge1)$ -cf. \textbf{Assumption
SM}. The former restriction in (\ref{MartinApprox0b}) ensures that
$\left[0,\alpha\left(\frac{1}{\delta+\nu}-1\right)\right)$ is not
empty, while the latter restriction ensures that $\sum^{\infty}_{s=0}\left\vert c_{s}\right\vert ^{\frac{\alpha}{\alpha_{\varepsilon}}-\delta}<\infty$.\footnote{By  (\ref{MartinApprox0b}) to $\frac{\alpha}{\alpha_{\varepsilon}}-\delta>\nu$.
Therefore, $\sum^{\infty}_{s=0}\left\vert c_{s}\right\vert ^{\nu}<\infty$
(cf. \textbf{Assumption SM}) implies $\sum^{\infty}_{s=0}\left\vert c_{s}\right\vert ^{\frac{\alpha}{\alpha_{\varepsilon}}-\delta}<\infty$. } Next, by (\ref{MartinApprox2b}), for all $\delta>0$ as $n\rightarrow\infty$
\[
H\left(\frac{\gamma_{n}}{\left\vert a_{s}\right\vert }\right)^{1/\alpha_{\varepsilon}}\left(\frac{\gamma_{n}}{\left\vert a_{s}\right\vert }\right)^{-\delta}\rightarrow0,
\]
uniformly in $s\geq1$. In view of this and (\ref{MartinApprox2b})
we get as $n\rightarrow\infty$
\[
I^{\prime}_{n}=\frac{n^{1/\alpha_{\varepsilon}}}{\gamma^{\alpha/\alpha_{\varepsilon}}_{n}n}\sum^{n}_{s=0}\left\vert a_{s}\right\vert ^{\alpha/\alpha_{\varepsilon}}H\left(\frac{\gamma_{n}}{\left\vert a_{s}\right\vert }\right)^{1/\alpha_{\varepsilon}}
\]
\[
\leq\frac{n^{1/\alpha_{\varepsilon}}}{\gamma^{\frac{\alpha}{\alpha_{\varepsilon}}-\delta}_{n}n}\sum^{n}_{s=0}\left\vert a_{s}\right\vert ^{\frac{\alpha}{\alpha_{\varepsilon}}-\delta}=\frac{n^{1/\alpha_{\varepsilon}}}{\gamma^{\frac{\alpha}{\alpha_{\varepsilon}}-\delta}_{n}n}\sum^{n}_{s=0}\left\vert \sum^{s}_{l=0}\rho^{l}_{n}c_{s-l}\right\vert ^{\frac{\alpha}{\alpha_{\varepsilon}}-\delta}
\]
\[
\leq\frac{n^{1/\alpha_{\varepsilon}}}{\gamma^{\frac{\alpha}{\alpha_{\varepsilon}}-\delta}_{n}n}\sum^{n}_{s=0}\sum^{s}_{l=0}\left\vert \rho^{l}_{n}c_{s-l}\right\vert ^{\frac{\alpha}{\alpha_{\varepsilon}}-\delta}=\frac{n^{1/\alpha_{\varepsilon}}}{\gamma^{\frac{\alpha}{\alpha_{\varepsilon}}-\delta}_{n}n}\sum^{n}_{l=0}\sum^{n}_{s=l}\left\vert \rho^{l}_{n}c_{s-l}\right\vert ^{\frac{\alpha}{\alpha_{\varepsilon}}-\delta}
\]
\begin{equation}
\leq\frac{n^{1/\alpha_{\varepsilon}}}{\gamma^{\frac{\alpha}{\alpha_{\varepsilon}}-\delta}_{n}n}\sum^{n}_{l=0}\rho^{\left(\frac{\alpha}{\alpha_{\varepsilon}}-\delta\right)l}_{n}\sum^{\infty}_{s=0}\left\vert c_{s}\right\vert ^{\frac{\alpha}{\alpha_{\varepsilon}}-\delta}\leq_{(1)}C\frac{n^{1/\alpha_{\varepsilon}}\kappa_{n}}{\gamma^{\frac{\alpha}{\alpha_{\varepsilon}}-\delta}_{n}n},\label{MartinApprox4c}
\end{equation}
where $\leq_{(1)}$ follows from Lemma 6.4 of DK2018, and the fact
that $\frac{\alpha}{\alpha_{\varepsilon}}-\delta>\nu>0$. Recall that
under the current assumptions we have $\kappa_{n}=n^{\alpha_{\kappa}}\ell_{\kappa}(n)\rightarrow\infty$
with $\alpha_{\kappa}\in[0,1)$, $\ell_{\kappa}$ SV, and $\gamma_{n}\sim\kappa^{1/\alpha}_{n}\ell_{1}(\kappa_{n})^{1/\alpha}$.
Hence, as $n\rightarrow\infty$ the last term in (\ref{MartinApprox4c})
is 
\[
\frac{n^{1/\alpha_{\varepsilon}}n^{\alpha_{\kappa}}\ell_{\kappa}(n)}{\left[\kappa^{\frac{1}{\alpha}}_{n}\ell_{1}(\kappa_{n})^{\frac{1}{\alpha}}\right]^{\frac{\alpha}{\alpha_{\varepsilon}}-\delta}n}=\frac{n^{1/\alpha_{\varepsilon}}n^{\alpha_{\kappa}}\ell_{\kappa}(n)}{\left[n^{\frac{\alpha_{\kappa}}{\alpha}}\ell_{\kappa}(n)^{\frac{1}{\alpha}}\ell_{1}(\kappa_{n})^{\frac{1}{\alpha}}\right]^{\frac{\alpha}{\alpha_{\varepsilon}}-\delta}n}
\]
\[
=\frac{\ell_{\kappa}(n)n^{1/\alpha_{\varepsilon}}n^{\frac{\alpha_{\kappa}}{\alpha}\left(\alpha-\frac{\alpha}{\alpha_{\varepsilon}}+\delta\right)}}{\left[\ell_{\kappa}(n)^{\frac{1}{\alpha}}\ell_{1}(\kappa_{n})^{\frac{1}{\alpha}}\right]^{\frac{\alpha}{\alpha_{\varepsilon}}-\delta}n}=:\lambda_{n}n^{\frac{1}{\alpha_{\varepsilon}}+\frac{\alpha_{\kappa}}{\alpha}\left(\alpha-\frac{\alpha}{\alpha_{\varepsilon}}+\delta\right)-1}=:\lambda_{n}n^{\eta}.
\]
Notice that the term 
\[
\lambda_{n}=\frac{\ell_{\kappa}(n)}{\left[\ell_{\kappa}(n)^{\frac{1}{\alpha}}\ell_{1}(\kappa_{n})^{\frac{1}{\alpha}}\right]^{\frac{\alpha}{\alpha_{\varepsilon}}-\delta}}=\frac{\ell_{\kappa}(n)}{\left[\ell_{\kappa}(n)^{\frac{1}{\alpha}}\ell_{1}(n^{\alpha_{\kappa}}\ell_{\kappa}(n))^{\frac{1}{\alpha}}\right]^{\frac{\alpha}{\alpha_{\varepsilon}}-\delta}}
\]
is SV in $n$ (see e.g. BGT 1987, Proposition 1.5.7). In view of this,
the result will follow if we show that $\eta<0$. Indeed, for $\alpha_{\kappa}=0$
it can be easily checked that $\eta<0$. For $\alpha_{\kappa}\in(0,1)$,
the following additional restrictions on $\delta$ and $\varepsilon$
ensure that $\eta$ is strictly negative.

\bigskip{}

\begin{itemize}
\item $1<\alpha\leq2$: In addition to the restrictions of (\ref{MartinApprox0})
and (\ref{MartinApprox0b}), $\delta$ satisfies 
\[
\delta<\left(\alpha-1\right)\left(\frac{1}{\alpha_{\kappa}}-1\right)\wedge\left(1-\nu\right).
\]
\item $0<\alpha\leq1$: In addition to the restrictions of (\ref{MartinApprox0})
and (\ref{MartinApprox0b}), $\delta$ and $\varepsilon$ satisfy
\[
\delta<\left(\alpha-\nu\right)\left(1-\alpha_{\kappa}\right)\text{\ and \ }\frac{\left(1-\alpha\right)\left(1-\alpha_{\kappa}\right)+\alpha_{\kappa}\delta}{1-\alpha_{\kappa}\left(1+\frac{\delta}{\alpha}\right)}<\varepsilon.
\]
\end{itemize}
Finally, we show the result under \textbf{Assumption SM(i)}. We only
need to consider the case $\alpha\in(1,2]$. Note that \textbf{Assumption
SM(i)} and \textbf{Assumption SM(ii)} are identical for $\alpha\in(0,1]$.
Let $\Delta_{s,n}$ as above. Using a similar argument to the one
of Astrauskas (1983, Lemma) together with (\ref{MartinApprox2b})
we get
\[
\mathbf{E}\left\vert \xi_{1}-\xi^{\ast}\right\vert ^{\alpha}I\left\{ \left\vert \xi_{1}-\xi^{\ast}\right\vert \leq\gamma_{n}\left\vert a_{s}\right\vert ^{-1}\right\} \asymp S\left(\frac{\gamma_{n}}{\left\vert a_{s}\right\vert }\right),
\]
where $S$ a bounded SV function such that
\[
S(x)=\int^{x}_{1}\frac{h(u)}{u}du.
\]
Hence, as $n\rightarrow\infty$
\[
\mathbf{E}\Delta^{\alpha}_{s,n}I\left\{ \Delta_{s,n}\leq1\right\} =C\left(\frac{\left\vert a_{s}\right\vert }{\gamma_{n}}\right)^{\alpha}S\left(\frac{\gamma_{n}}{\left\vert a_{s}\right\vert }\right)\left[1+o(1)\right].
\]
In view of this and using similar arguments as those used in (\ref{MartinApprox4b})
as $n\rightarrow\infty$ we get
\[
I_{n}\leq\frac{1}{n}\sum^{n-1}_{s=0}\left[\mathbf{E}\left\vert M_{n,s}\right\vert ^{\alpha}\right]^{1/\alpha}\leq C\frac{n^{1/\alpha}}{\gamma_{n}n}\sum^{n}_{s=0}\left\vert a_{s}\right\vert S\left(\frac{\gamma_{n}}{\left\vert a_{s}\right\vert }\right)^{1/\alpha}
\]
\[
\leq C\frac{n^{1/\alpha}}{\gamma_{n}n}\sum^{n}_{s=0}\left\vert a_{s}\right\vert \leq\frac{n^{1/\alpha}}{\gamma_{n}n}\sum^{n}_{l=0}\rho^{l}_{n}\sum^{\infty}_{s=0}\left\vert c_{s}\right\vert 
\]
\[
\leq C\frac{n^{1/\alpha}\kappa_{n}}{\gamma_{n}n}=C\frac{n^{1/\alpha}\kappa_{n}}{\kappa^{1/\alpha}_{n}\ell_{1}(\kappa_{n})^{\frac{1}{\alpha}}n}
\]
\[
=\frac{C}{\ell_{1}(\kappa_{n})^{\frac{1}{\alpha}}}\left(\frac{n}{\kappa_{n}}\right)^{\frac{1}{\alpha}-1}=\frac{C}{\ell_{1}(n^{\alpha_{\kappa}}\ell_{\kappa}(n))^{\frac{1}{\alpha}}}\left(\frac{n}{n^{\alpha_{\kappa}}\ell_{\kappa}(n)}\right)^{\frac{1}{\alpha}-1}=o(1),
\]
as required. 

\bigskip{}

\textbf{MI-LM (}$\alpha\in(1,2]$\textbf{, $m\in(1/\alpha,1)$)}.
Similarly to the previous part, let
\begin{equation}
0<\delta<1\text{ and }0\leq\varepsilon<\alpha\left(\frac{1}{\delta}-1\right).\label{MartinApprox0c}
\end{equation}
 Then as $n\rightarrow\infty$
\[
I^{\prime}_{n}=\frac{n^{1/\alpha_{\varepsilon}}}{\gamma^{\alpha/\alpha_{\varepsilon}}_{n}n}\sum^{n}_{s=0}\left\vert a_{s}\right\vert ^{\alpha/\alpha_{\varepsilon}}H\left(\frac{\gamma_{n}}{\left\vert a_{s}\right\vert }\right)^{1/\alpha_{\varepsilon}}
\]
\[
\leq\frac{n^{1/\alpha_{\varepsilon}}}{\gamma^{\frac{\alpha}{\alpha_{\varepsilon}}-\delta}_{n}n}\sum^{n}_{s=0}\left\vert a_{s}\right\vert ^{\frac{\alpha}{\alpha_{\varepsilon}}-\delta}=\frac{n^{1/\alpha_{\varepsilon}}}{\gamma^{\frac{\alpha}{\alpha_{\varepsilon}}-\delta}_{n}n}\sum^{n}_{s=0}\left\vert \sum^{s}_{l=0}\rho^{l}_{n}c_{s-l}\right\vert ^{\frac{\alpha}{\alpha_{\varepsilon}}-\delta}
\]
\[
\leq\frac{n^{1/\alpha_{\varepsilon}}}{\gamma^{\frac{\alpha}{\alpha_{\varepsilon}}-\delta}_{n}n}\sum^{n}_{s=0}\sum^{s}_{l=0}\left\vert \rho^{l}_{n}c_{s-l}\right\vert ^{\frac{\alpha}{\alpha_{\varepsilon}}-\delta}=\frac{n^{1/\alpha_{\varepsilon}}}{\gamma^{\frac{\alpha}{\alpha_{\varepsilon}}-\delta}_{n}n}\sum^{n}_{l=0}\sum^{n}_{s=l}\left\vert \rho^{l}_{n}c_{s-l}\right\vert ^{\frac{\alpha}{\alpha_{\varepsilon}}-\delta}
\]
\[
\leq\frac{n^{1/\alpha_{\varepsilon}}}{\gamma^{\frac{\alpha}{\alpha_{\varepsilon}}-\delta}_{n}n}\sum^{n}_{l=0}\left\vert \rho^{l}_{n}\right\vert ^{\frac{\alpha}{\alpha_{\varepsilon}}-\delta}\sum^{n}_{s=0}\left\vert c_{s}\right\vert ^{\frac{\alpha}{\alpha_{\varepsilon}}-\delta}\leq_{(2)}C\frac{n^{1/\alpha_{\varepsilon}}\kappa_{n}}{\gamma^{\frac{\alpha}{\alpha_{\varepsilon}}-\delta}_{n}n}\sum^{n}_{s=1}s^{-m\left(\frac{\alpha}{\alpha_{\varepsilon}}-\delta\right)}\left[1+o(1)\right]
\]
\[
=C\frac{n^{\frac{1}{\alpha_{\varepsilon}}-m\left(\frac{\alpha}{\alpha_{\varepsilon}}-\delta\right)}\kappa_{n}}{\gamma^{\frac{\alpha}{\alpha_{\varepsilon}}-\delta}_{n}}\frac{1}{n}\sum^{n}_{s=1}\left(\frac{s}{n}\right)^{-m\left(\frac{\alpha}{\alpha_{\varepsilon}}-\delta\right)}\left[1+o(1)\right]
\]
\[
=C\frac{n^{\frac{1}{\alpha_{\varepsilon}}-m\left(\frac{\alpha}{\alpha_{\varepsilon}}-\delta\right)}\kappa_{n}}{\gamma^{\frac{\alpha}{\alpha_{\varepsilon}}-\delta}_{n}}\int^{1}_{0}x^{-m\left(\frac{\alpha}{\alpha_{\varepsilon}}-\delta\right)}dx\left[1+o(1)\right],
\]
where $\leq_{(2)}$ follows from DK2018 (Lemma 6.4), and the fact
that $\frac{\alpha}{\alpha_{\varepsilon}}-\delta>0$, under (\ref{MartinApprox0c}).
Moreover, notice that the integral above is finite because $m\left(\frac{\alpha}{\alpha_{\varepsilon}}-\delta\right)<1$.
Recall that $\gamma_{n}\sim\kappa^{b/\alpha}_{n}\ell_{1}\left(\kappa_{n}\right)^{1/\alpha}$,
$b=1+\alpha-\alpha m$ with $m\in\left(1/\alpha,1\right)$, and $\kappa_{n}=n^{\alpha_{\kappa}}\ell_{\kappa}(n)$.
Write 
\[
\frac{n^{\frac{1}{\alpha_{\varepsilon}}-m\left(\frac{\alpha}{\alpha_{\varepsilon}}-\delta\right)}\kappa_{n}}{\gamma^{\frac{\alpha}{\alpha_{\varepsilon}}-\delta}_{n}}=\frac{n^{\frac{1}{\alpha_{\varepsilon}}-m\left(\frac{\alpha}{\alpha_{\varepsilon}}-\delta\right)}n^{\alpha_{\kappa}}\ell_{\kappa}(n)}{\left[\kappa^{b/\alpha}_{n}\ell_{1}\left(\kappa_{n}\right)^{1/\alpha}\right]^{\left(\frac{\alpha}{\alpha_{\varepsilon}}-\delta\right)}}=\frac{n^{\frac{1}{\alpha_{\varepsilon}}-m\left(\frac{\alpha}{\alpha_{\varepsilon}}-\delta\right)}n^{\alpha_{\kappa}}\ell_{\kappa}(n)}{\left[\left(n^{\alpha_{\kappa}}\ell_{\kappa}(n)\right)^{b/\alpha}\ell_{1}\left(\kappa_{n}\right)^{1/\alpha}\right]^{\left(\frac{\alpha}{\alpha_{\varepsilon}}-\delta\right)}}
\]
\[
=\frac{\ell_{\kappa}(n)}{\left[\ell_{\kappa}(n)^{b/\alpha}\ell_{1}\left(\kappa_{n}\right)^{1/\alpha}\right]^{\left(\frac{\alpha}{\alpha_{\varepsilon}}-\delta\right)}}\frac{n^{\frac{1}{\alpha_{\varepsilon}}+\alpha_{\kappa}-m\left(\frac{\alpha}{\alpha_{\varepsilon}}-\delta\right)}}{n^{\frac{\alpha_{\kappa}}{\alpha}b\left(\frac{\alpha}{\alpha_{\varepsilon}}-\delta\right)}}=:\lambda_{n}n^{\eta},
\]
where as in the previous part $\lambda_{n}$ is SV in $n$. Therefore,
is suffices to show that $\eta<0$. Lengthy but straightforward calculations
show that the latter is true if $\delta$ and $\varepsilon$ satisfy
the additional requirements
\[
\delta<\frac{\left(\alpha m-1\right)\left(1-\alpha_{\kappa}\right)}{\alpha_{\kappa}b+\alpha m}\text{\ and\ }\varepsilon<\frac{\left(\alpha m-1\right)\left(1-\alpha_{\kappa}\right)-\left[\alpha_{\kappa}b+\alpha m\right]\delta}{m\delta+\alpha_{\kappa}+\alpha_{\kappa}\frac{b}{\alpha}\delta}.
\]

\bigskip{}

\textbf{MI-NM (}$\alpha\in(1,2]$\textbf{, $m\in(1,1+1/\alpha)$)}.
Similarly to the previous part, assume that (\ref{MartinApprox0})
holds. In addition, suppose that 
\begin{equation}
\delta<\frac{\left(\alpha-1\right)\left[1-\alpha_{\kappa}\right]}{\alpha_{\kappa}}1\left\{ \alpha_{\kappa}\neq0\right\} +1\left\{ \alpha_{\kappa}=0\right\} \wedge\frac{m-1}{m},\label{MartinApprox5a}
\end{equation}
and

\begin{equation}
\varepsilon<\frac{\alpha m\left(1-\delta\right)-\alpha}{1+m\delta}\wedge\alpha\left(\delta^{-1}-1\right).\label{MartinApprox5b}
\end{equation}
We remark that the requirement $\delta<\frac{m-1}{m}$ ensures that
the upper bound $\frac{\alpha m\left(1-\delta\right)-\alpha}{1+m\delta}$
in (\ref{MartinApprox5b}) is strictly positive. Then as $n\rightarrow\infty$
\[
I^{\prime}_{n}=\frac{n^{1/\alpha_{\varepsilon}}}{\gamma^{\alpha/\alpha_{\varepsilon}}_{n}n}\sum^{n}_{s=0}\left\vert a_{s}\right\vert ^{\alpha/\alpha_{\varepsilon}}H\left(\frac{\gamma_{n}}{\left\vert a_{s}\right\vert }\right)^{1/\alpha_{\varepsilon}}
\]
\[
\leq_{\left(1\right)}\frac{n^{1/\alpha_{\varepsilon}}}{\gamma^{\frac{\alpha}{\alpha_{\varepsilon}}-\delta}_{n}n}\sum^{n}_{s=0}\left\vert a_{s}\right\vert ^{\frac{\alpha}{\alpha_{\varepsilon}}-\delta}=\frac{n^{1/\alpha_{\varepsilon}}}{\gamma^{\frac{\alpha}{\alpha_{\varepsilon}}-\delta}_{n}n}\sum^{n}_{s=0}\left\vert \sum^{s}_{l=0}\rho^{l}_{n}c_{s-l}\right\vert ^{\frac{\alpha}{\alpha_{\varepsilon}}-\delta}
\]
\[
=_{(2)}\frac{n^{1/\alpha_{\varepsilon}}}{\gamma^{\frac{\alpha}{\alpha_{\varepsilon}}-\delta}_{n}n}\sum^{\delta_{n}}_{s=0}\left\vert \sum^{s}_{l=0}\rho^{l}_{n}c_{s-l}\right\vert ^{\frac{\alpha}{\alpha_{\varepsilon}}-\delta}\left[1+o(1)\right]
\]
\[
=\frac{n^{1/\alpha_{\varepsilon}}\kappa^{1+\frac{\alpha}{\alpha_{\varepsilon}}-\delta-m\left(\frac{\alpha}{\alpha_{\varepsilon}}-\delta\right)}_{n}}{\gamma^{\frac{\alpha}{\alpha_{\varepsilon}}-\delta}_{n}n}\frac{1}{\kappa_{n}}\sum^{\delta_{n}}_{s=0}\left\vert \frac{1}{\kappa_{n}}\sum^{s}_{l=1}\rho^{s-l}_{n}\left(\frac{l}{^{\kappa_{n}}}\right)^{-m}\right\vert ^{\frac{\alpha}{\alpha_{\varepsilon}}-\delta}\left[1+o(1)\right]
\]
\[
=_{\left(3\right)}C\int^{\infty}_{x=0}\left\vert p(x)\right\vert ^{\frac{\alpha}{\alpha_{\varepsilon}}-\delta}dx\cdot\frac{n^{1/\alpha_{\varepsilon}}\kappa^{1+\frac{\alpha}{\alpha_{\varepsilon}}-\delta-m\left(\frac{\alpha}{\alpha_{\varepsilon}}-\delta\right)}_{n}}{\gamma^{\frac{\alpha}{\alpha_{\varepsilon}}-\delta}_{n}n}\left[1+o(1)\right]
\]
\[
=C\frac{n^{1/\alpha_{\varepsilon}}\kappa^{1+\frac{\alpha}{\alpha_{\varepsilon}}-\delta-m\left(\frac{\alpha}{\alpha_{\varepsilon}}-\delta\right)}_{n}}{\gamma^{\frac{\alpha}{\alpha_{\varepsilon}}-\delta}_{n}n}=C\frac{n^{1/\alpha_{\varepsilon}}\kappa^{1+\frac{\alpha}{\alpha_{\varepsilon}}-\delta-m\left(\frac{\alpha}{\alpha_{\varepsilon}}-\delta\right)}_{n}}{\left[\kappa^{b/\alpha}_{n}\ell_{1}\left(\kappa_{n}\right)^{1/\alpha}\right]^{\left(\frac{\alpha}{\alpha_{\varepsilon}}-\delta\right)}n}
\]
\[
=C\frac{n^{1/\alpha_{\varepsilon}}\kappa^{1+\frac{\alpha}{\alpha_{\varepsilon}}-\delta-m\left(\frac{\alpha}{\alpha_{\varepsilon}}-\delta\right)}_{n}}{\ell_{1}\left(\kappa_{n}\right)^{\frac{1}{\alpha}\left(\frac{\alpha}{\alpha_{\varepsilon}}-\delta\right)}\kappa^{\left(\frac{1}{\alpha}+1-m\right)\left(\frac{\alpha}{\alpha_{\varepsilon}}-\delta\right)}_{n}n}=C\frac{n^{1/\alpha_{\varepsilon}}\kappa^{1+\frac{\alpha}{\alpha_{\varepsilon}}-\delta}_{n}}{\ell_{1}\left(\kappa_{n}\right)^{\frac{1}{\alpha}\left(\frac{\alpha}{\alpha_{\varepsilon}}-\delta\right)}\kappa^{\left(\frac{1}{\alpha}+1\right)\left(\frac{\alpha}{\alpha_{\varepsilon}}-\delta\right)}_{n}n}
\]
\[
=C\frac{n^{\frac{1}{\alpha_{\varepsilon}}-1}\kappa^{1-\frac{1}{\alpha_{\varepsilon}}+\frac{\delta}{\alpha}}_{n}}{\ell_{1}\left(\kappa_{n}\right)^{\frac{1}{\alpha}\left(\frac{\alpha}{\alpha_{\varepsilon}}-\delta\right)}}=C=:\lambda_{n}\cdot n^{\frac{1}{\alpha_{\varepsilon}}-1}n^{\alpha_{\kappa}\left(1-\frac{1}{\alpha_{\varepsilon}}+\frac{\delta}{\alpha}\right)}
\]
\[
=:C\lambda_{n}n^{\eta},
\]
where $\leq_{\left(1\right)}$ follows from (\ref{MartinApprox2b}).
Moreover, $=_{\left(2\right)}$ follows from \textbf{Lemma \ref{LemmaMI0}}
and the fact that $m\left(\frac{\alpha}{\alpha_{\varepsilon}}-\delta\right)>1$.
Moreover, $=_{(3)}$ is due to \textbf{Lemma \ref{MI-NM Scale}}.\footnote{Note that (\ref{MartinApprox5b}) implies that $m\left(\frac{\alpha}{\alpha_{\varepsilon}}-\delta\right)>1$.}
Lengthy but straightforward calculations show that under the current
assumptions, i.e. (\ref{MartinApprox0}), (\ref{MartinApprox5a})
and (\ref{MartinApprox5b}), $\eta<0$. In view of this, the requisite
result follows from the fact that $\lambda_{n}$ is SV in $n$. $\blacksquare$

\section{Additional Workings for the MI-LM and MI-NM Cases}

We start with some notation. For some $m,\theta>0$ set 
\[
b_{\theta}:=1+\theta-\theta m
\]

\bigskip{}
We have the following results.\bigskip{}

\begin{lemma} \label{Ts} Let $m\in\left(1,2\right)$ and $\theta>0$
such that $1/m<\theta<1/(m-1)$. Moreover, suppose that $\rho_{n}$
satisfies \textbf{Assumption MI}. Set 
\begin{equation}
T_{1nk}:=-\kappa^{-1+m}_{n}\sum^{\infty}_{j=k+1}c_{j}\text{\ and\ }T_{2nk}:=\frac{1}{\kappa^{2-m}_{n}}\sum^{k}_{l=1}\rho^{k-l}_{n}\sum^{\infty}_{j=l}c_{j}.\label{Ts0}
\end{equation}
Finally, let $\varsigma(.)$ SV such that $\lim_{x\rightarrow\infty}\varsigma(x)=\infty$,
and set $\delta_{n}:=\kappa_{n}\cdot\varsigma(\kappa_{n})$, $\mu_{n}:=\mu\left(\kappa^{b_{\theta}/2}_{n}\right)$
with $\mu\left(j\right):=\left|\frac{c_{j}}{j^{-m}}-1\right|=O(j^{-\nu})$
for some $\nu>0$.

Then as $n\rightarrow\infty$ 
\begin{equation}
\frac{1}{\kappa_{n}}\sum^{\delta_{n}}_{k=1}\left\vert \Delta_{1nk}\right\vert ^{\theta}:=\frac{1}{\kappa_{n}}\sum^{\delta_{n}}_{k=1}\left\vert T_{1nk}+\int^{\infty}_{\frac{k}{\kappa_{n}}}x^{-m}dx\right\vert ^{\theta}=o(1),\label{Ts1}
\end{equation}
and 
\begin{equation}
\frac{1}{\kappa_{n}}\sum^{\delta_{n}}_{k=1}\left\vert \Delta_{2nk}\right\vert ^{\theta}:=\frac{1}{\kappa_{n}}\sum^{\delta_{n}}_{k=1}\left\vert T_{2nk}-\frac{1}{1-m}e^{-\frac{k}{\kappa_{n}}}\int^{k/\kappa_{n}}_{0}e^{y}y^{1-m}dy\right\vert ^{\theta}=o(1).\label{Ts2}
\end{equation}
 \end{lemma}

\bigskip{}

\begin{corollary} \label{Ts_corol} Suppose that \textbf{Assumption
MI-NM} holds. Let $T_{1nk}$ and $T_{2nk}$ as in (\ref{Ts0}). Then
as $n\rightarrow\infty$ 
\[
\frac{1}{\kappa_{n}}\sum^{\delta_{n}}_{k=1}\left\vert \Delta_{1nk}\right\vert ^{\alpha}:=\frac{1}{\kappa_{n}}\sum^{\delta_{n}}_{k=1}\left\vert T_{1nk}+\int^{\infty}_{\frac{k}{\kappa_{n}}}x^{-m}dx\right\vert ^{\alpha}=o(1),
\]
and 
\[
\frac{1}{\kappa_{n}}\sum^{\delta_{n}}_{k=1}\left\vert \Delta_{2nk}\right\vert ^{\alpha}:=\frac{1}{\kappa_{n}}\sum^{\delta_{n}}_{k=1}\left\vert T_{2nk}+\frac{1}{1-m}e^{-\frac{k}{\kappa_{n}}}\int^{k/\kappa_{n}}_{0}e^{y}y^{1-m}dy\right\vert ^{\alpha}=o(1),
\]
where $\delta_{n}$ is given in\textbf{ Lemma \ref{Ts}}. \end{corollary}\bigskip{}

\begin{remark} Note that under the assumptions of \textbf{Lemma \ref{Ts}},
for all $a,b>0$ 
\begin{equation}
\mu^{a}_{n}\varsigma(\kappa_{n})^{b}\rightarrow0,\label{Ts_rmk}
\end{equation}
as $n\rightarrow\infty$.

\end{remark}

For the subsequent results recall that 
\[
p(x):=x^{1-m}-e^{-x}\int^{x}_{0}e^{y}y^{1-m}dy.
\]
\bigskip{}

\begin{remark} \textbf{(a)} It can be shown that -see (\ref{LemmaMI4.3})
below- that under \textbf{MI-NM},
\[
\kappa^{1-m}_{n}a_{k}(n)=\kappa^{1-m}_{n}\sum^{k}_{l=0}\rho^{k-l}_{n}c_{l}=T_{1nk}+T_{2nk}.
\]
Moreover,
\[
-\int^{\infty}_{\frac{k}{\kappa_{n}}}x^{-m}dx-\frac{1}{1-m}e^{-\frac{k}{\kappa_{n}}}\int^{k/\kappa_{n}}_{0}e^{y}y^{1-m}dy
\]
\[
=\frac{\left(k/\kappa_{n}\right)^{1-m}}{1-m}-\frac{1}{1-m}e^{-\frac{k}{\kappa_{n}}}\int^{k/\kappa_{n}}_{0}e^{y}y^{1-m}dy=\frac{p(k/\kappa_{n})}{1-m}.
\]

\textbf{(b)} The following fact is utilized in the proof of\textbf{
Lemma \ref{WCon}} and is a direct consequence of \textbf{Lemma \ref{Ts_corol}}
\begin{equation}
\frac{1}{\kappa_{n}}\sum^{\delta_{n}}_{k=1}\left\vert \kappa^{m-1}_{n}a_{k}(n)-\frac{p(k/\kappa_{n})}{1-m}\right\vert ^{\alpha}=o(1),\label{Ts_rmk-1}
\end{equation}
as $n\rightarrow\infty$.

\end{remark}\bigskip{}

\begin{lemma} \label{pi} Let $m\in\left(1,2\right)$ and $\theta>0$
such that $1/m<\theta<1/(m-1)$. Then as $n\rightarrow\infty$
\[
\frac{1}{\kappa_{n}}\sum^{\delta_{n}}_{k=1}\left\vert p\left(k/\kappa_{n}\right)\right\vert ^{\theta}=\int^{\infty}_{0}\left\vert p\left(x\right)\right\vert ^{\theta}dx+o(1),
\]
with $\delta_{n}$ given in\textbf{ Lemma \ref{Ts}}. $\allowbreak$
\end{lemma}

\bigskip{}

\begin{corollary} \label{pi corol} Suppose that \textbf{Assumption
MI-NM} holds. Then as $n\rightarrow\infty$
\[
\frac{1}{\kappa_{n}}\sum^{\delta_{n}}_{k=1}\left\vert p\left(k/\kappa_{n}\right)\right\vert ^{\alpha}=\int^{\infty}_{0}\left\vert p\left(x\right)\right\vert ^{\alpha}dx+o(1),
\]
and
\[
\frac{1}{\kappa_{n}}\sum^{\delta_{n}}_{k=1}p\left(k/\kappa_{n}\right)^{\left\langle \alpha\right\rangle }=\int^{\infty}_{0}p\left(x\right)^{\left\langle \alpha\right\rangle }dx+o(1),
\]
where $\delta_{n}$ is given in\textbf{ Lemma \ref{Ts}}.$\allowbreak$
\end{corollary}

\bigskip{}

\begin{lemma} \label{MI-NM Scale} Let $m\in\left(1,2\right)$ and
$\theta>0$ such that $1/m<\theta<1/(m-1)$. Then as $n\rightarrow\infty$
\[
\frac{1}{\kappa_{n}}\sum^{\delta_{n}}_{k=1}\left\vert T_{1nk}+T_{2nk}\right\vert ^{\theta}=\frac{1}{\left(m-1\right)^{\theta}}\int^{\infty}_{0}\left\vert p\left(x\right)\right\vert ^{\theta}dx+o(1),
\]
with $\delta_{n}$ given in\textbf{ Lemma \ref{Ts}}.\end{lemma}

\bigskip{}

\begin{corollary} \label{MI-NM Scale corol} Suppose that \textbf{Assumption
MI-NM} holds. Then as $n\rightarrow\infty$
\[
\frac{1}{\kappa_{n}}\sum^{\delta_{n}}_{k=1}\left\vert T_{1nk}+T_{2nk}\right\vert ^{\alpha}=\frac{1}{\left(m-1\right)^{\alpha}}\int^{\infty}_{0}\left\vert p\left(x\right)\right\vert ^{\alpha}dx+o(1),
\]
with $\delta_{n}$ given in\textbf{ Lemma \ref{Ts}}. \end{corollary}

\bigskip{}

\textbf{Proof of Lemma \ref{Ts}}. Notice that under the assumptions
of \textbf{Lemma \ref{Ts}} we have
\begin{equation}
\theta m>1\text{ \ and \ }b_{\theta}:=1+\theta-\theta m>0.\label{T1.00}
\end{equation}

\begin{itemize}
\item $T_{1nk}$: We start with the proof of (\ref{Ts1}). Set
\[
T_{1nk}:=-\kappa^{m-1}_{n}\sum^{\infty}_{j=k+1}c_{j},\text{ \ }T^{\left(1\right)}_{1nk}:=-\kappa^{-1+m}_{n}\sum^{\infty}_{j=k+1}j^{-m}
\]
\begin{equation}
T^{\left(2\right)}_{1nk}:=-\int^{\infty}_{\frac{k+1}{\kappa_{n}}}x^{-m}dx,\text{ \ }T^{\left(3\right)}_{1nk}:=-\int^{\infty}_{\frac{k}{\kappa_{n}}}x^{-m}dx.\label{T1.001}
\end{equation}
\end{itemize}
In addition, let
\[
\Delta^{\left(1\right)}_{1nk}:=T_{1nk}-T^{\left(1\right)}_{1nk},\text{\ }\Delta^{\left(2\right)}_{1nk}:=T^{\left(1\right)}_{1nk}-T^{\left(2\right)}_{1nk},
\]
\[
\Delta^{\left(3\right)}_{1nk}:=T^{\left(2\right)}_{1nk}-T^{\left(3\right)}_{1nk},\text{\ and \ }\Delta_{1nk}:=T_{1nk}-T^{\left(3\right)}_{1nk}.
\]
First, note that 
\begin{equation}
\frac{1}{\kappa_{n}}\sum^{\kappa^{b_{\theta}/2}_{n}}_{k=1}\left\vert T_{1nk}\right\vert ^{\theta}\leq\frac{\kappa^{b_{\theta}/2}_{n}}{\kappa^{b_{\theta}}_{n}}\left(\sum^{\infty}_{k=1}\left\vert c_{k}\right\vert \right)^{\theta}=O\left(\kappa^{-b_{\theta}/2}_{n}\right)=o\left(1\right).\label{T1.0}
\end{equation}
Write 
\begin{equation}
\frac{1}{\kappa_{n}}\sum^{\delta_{n}}_{k=\kappa^{b_{\theta}/2}_{n}}\left\vert \Delta_{1nk}\right\vert ^{\theta}\leq\frac{C}{\kappa_{n}}\sum^{\delta_{n}}_{k=\kappa^{b_{\theta}/2}_{n}}\sum^{3}_{i=1}\left\vert \Delta^{\left(i\right)}_{1nk}\right\vert ^{\theta}.\label{T1.0b}
\end{equation}
We shall show that the terms on the r.h.s. above are $o(1)$, as $n\rightarrow\infty$.

As $n\rightarrow\infty$ the first term
\[
\frac{1}{\kappa_{n}}\sum^{\delta_{n}}_{k=\kappa^{b_{\theta}/2}_{n}}\left\vert \Delta^{\left(1\right)}_{1nk}\right\vert ^{\theta}=\frac{1}{\kappa_{n}}\sum^{\delta_{n}}_{k=\kappa^{b_{\theta}/2}_{n}}\left\vert \kappa^{-1+m}_{n}\sum^{\infty}_{j=k+1}c_{j}-\kappa^{-1+m}_{n}\sum^{\infty}_{j=k+1}j^{-m}\right\vert ^{\theta}
\]
\[
\leq\frac{1}{\kappa_{n}}\sum^{\delta_{n}}_{k=\kappa^{b_{\theta}/2}_{n}}\left\vert \kappa^{-1+m}_{n}\sum^{\infty}_{j=k+1}\left\vert \frac{c_{j}}{j^{-m}}-1\right\vert j^{-m}\right\vert ^{\theta}\leq\mu^{\theta}_{n}\frac{C}{\kappa_{n}}\sum^{\delta_{n}}_{k=1}\left\vert \kappa^{-1+m}_{n}\sum^{\infty}_{j=k+1}j^{-m}\right\vert ^{\theta}
\]
\begin{equation}
=\mu^{\theta}_{n}\frac{C}{\kappa_{n}}\sum^{\delta_{n}}_{k=1}\left\vert T^{\left(1\right)}_{1nk}\right\vert ^{\theta}=_{\left(1\right)}o(1),\label{T1.1}
\end{equation}
where $=_{\left(1\right)}$ will be proved later.

Next, consider
\[
\frac{1}{\kappa_{n}}\sum^{\delta_{n}}_{k=1}\left\vert \Delta^{\left(2\right)}_{1nk}\right\vert ^{\theta}=\frac{1}{\kappa_{n}}\sum^{\delta_{n}}_{k=1}\left\vert \kappa^{-1+m}_{n}\sum^{\infty}_{j=k+1}j^{-m}-\int^{\infty}_{\frac{k+1}{\kappa_{n}}}x^{-m}dx\right\vert ^{\theta}
\]
\[
\leq\frac{1}{\kappa_{n}}\sum^{\delta_{n}}_{k=1}\left[\frac{1}{\kappa_{n}}\mathbf{V}\left(x^{-m}\right)^{\infty}_{\frac{k+1}{\kappa_{n}}}\right]^{\theta}\leq\frac{1}{\kappa_{n}}\sum^{\delta_{n}}_{k=1}\left[\frac{1}{\kappa_{n}}\left(\frac{k+1}{\kappa_{n}}\right)^{-m}\right]^{\theta}
\]

\begin{equation}
=\frac{1}{\kappa^{b_{\theta}}_{n}}\sum^{\delta_{n}}_{k=1}\left(k+1\right)^{-\theta m}\leq\frac{1}{\kappa^{b_{\theta}}_{n}}\sum^{\infty}_{k=1}\left(k+1\right)^{-\theta m}=O\left(\kappa^{-b_{\theta}}_{n}\right)=o(1),\label{T1.1b}
\end{equation}
where recall that $b_{\theta}>0$.

The third term 
\begin{equation}
\frac{1}{\kappa_{n}}\sum^{\delta_{n}}_{k=1}\left\vert \Delta^{\left(3\right)}_{1nk}\right\vert ^{\theta}=\frac{C}{\kappa_{n}}\sum^{\delta_{n}}_{k=1}\left\vert \left(\frac{k+1}{\kappa_{n}}\right)^{1-m}-\left(\frac{k}{\kappa_{n}}\right)^{1-m}\right\vert ^{\theta}.\label{T1.2}
\end{equation}
The summand in (\ref{T1.2}) is

\[
\left(\frac{\kappa_{n}}{k}\right)^{m-1}-\left(\frac{\kappa_{n}}{k+1}\right)^{m-1}=\left(\frac{\kappa_{n}}{k}\right)^{m-1}-\left(\frac{\kappa_{n}}{k}+\frac{\kappa_{n}}{k+1}-\frac{\kappa_{n}}{k}\right)^{m-1}
\]
\[
=\left(\frac{\kappa_{n}}{k}\right)^{m-1}-\left(\frac{\kappa_{n}}{k}+\frac{k\kappa_{n}-\left(k+1\right)\kappa_{n}}{\left(k+1\right)k}\right)^{m-1}
\]
\[
=\left(\frac{\kappa_{n}}{k}\right)^{m-1}-\left(\frac{\kappa_{n}}{k}-\frac{\kappa_{n}}{\left(k+1\right)k}\right)^{m-1}
\]
\[
=\left(\frac{\kappa_{n}}{\left(k+1\right)k}+\frac{\kappa_{n}}{k}-\frac{\kappa_{n}}{\left(k+1\right)k}\right)^{m-1}-\left(\frac{\kappa_{n}}{k}-\frac{\kappa_{n}}{\left(k+1\right)k}\right)^{m-1}
\]
\[
\leq\left(\frac{\kappa_{n}}{\left(k+1\right)k}\right)^{m-1}+\left(\frac{\kappa_{n}}{k}-\frac{\kappa_{n}}{\left(k+1\right)k}\right)^{m-1}-\left(\frac{\kappa_{n}}{k}-\frac{\kappa_{n}}{\left(k+1\right)k}\right)^{m-1}
\]
\[
=\left[\frac{\kappa_{n}}{\left(k+1\right)k}\right]^{m-1}.
\]
Hence,
\begin{equation}
\frac{1}{\kappa_{n}}\sum^{\delta_{n}}_{k=1}\left\vert \Delta^{\left(3\right)}_{1nk}\right\vert ^{\theta}\leq\frac{1}{\kappa_{n}}\sum^{\delta_{n}}_{k=1}\left[\frac{\kappa_{n}}{\left(k+1\right)k}\right]^{\theta\left(m-1\right)}=_{\left(2\right)}o(1).\label{T1.2b}
\end{equation}
To see why $=_{\left(2\right)}$ holds, first assume that $2\theta\left(m-1\right)>1$.\footnote{By definition, $2\theta\left(m-1\right)>0$.}
In this case
\[
\frac{1}{\kappa_{n}}\sum^{\delta_{n}}_{k=1}\left[\frac{\kappa_{n}}{\left(k+1\right)k}\right]^{\theta\left(m-1\right)}=\frac{1}{\kappa^{b_{\theta}}_{n}}\sum^{\delta_{n}}_{k=1}\left[\frac{1}{\left(k+1\right)k}\right]^{\theta\left(m-1\right)}\leq\frac{1}{\kappa^{b_{\theta}}_{n}}\sum^{\infty}_{k=1}k^{-2\theta\left(m-1\right)}=o(1).
\]
Next, suppose that $0<2\theta\left(m-1\right)\leq1$. Then
\[
\frac{1}{\kappa_{n}}\sum^{\delta_{n}}_{k=1}\left[\frac{\kappa_{n}}{\left(k+1\right)k}\right]^{\left(m-1\right)\theta}\leq\frac{1}{\kappa_{n}}\sum^{\kappa_{n}}_{k=1}\left(\frac{\kappa_{n}}{k^{2}}\right)^{\left(m-1\right)\theta}+\frac{1}{\kappa_{n}}\sum^{\delta_{n}}_{k=\kappa_{n}}\left(\frac{\kappa_{n}}{k^{2}}\right)^{\left(m-1\right)\theta}
\]
\[
\leq\frac{1}{\kappa^{1+\left(m-1\right)\theta}_{n}}\sum^{\kappa_{n}}_{k=1}\left(\frac{\kappa^{2}_{n}}{k^{2}}\right)^{\left(m-1\right)\theta}+\frac{1}{\kappa_{n}}\sum^{\delta_{n}}_{k=\kappa_{n}}\left(\frac{\kappa_{n}}{\kappa^{2}_{n}}\right)^{\left(m-1\right)\theta}
\]
\[
\leq\frac{1}{\kappa^{1+\left(m-1\right)\theta}_{n}}\sum^{\kappa_{n}}_{k=1}\left(\frac{\kappa^{2}_{n}}{k^{2}}\right)^{\left(m-1\right)\theta}+\frac{\delta_{n}}{\kappa^{1+\left(m-1\right)\theta}_{n}}
\]
\[
=\frac{1}{\kappa^{\left(m-1\right)\theta}_{n}}\left[\int^{1}_{1/\kappa_{n}}x^{-2\left(m-1\right)\theta}dx+o(1)\right]+\frac{\delta_{n}}{\kappa^{1+\left(m-1\right)\theta}_{n}}=o(1),
\]
where we have used the fact that $\delta_{n}=\kappa_{n}\varsigma(\kappa_{n})$
with $\varsigma$ slowly varying, and
\[
\int^{1}_{1/\kappa_{n}}x^{-2\left(m-1\right)\theta}dx=\left\{ \begin{array}{l}
\text{bounded, for }0<2\theta\left(m-1\right)<1\\
\ln(\kappa_{n})\text{, for }2\theta\left(m-1\right)=1
\end{array}\right..
\]

Finally, we show that $=_{\left(1\right)}$ in (\ref{T1.1}) holds
true. Write
\[
\frac{\mu^{\theta}_{n}}{\kappa_{n}}\sum^{\delta_{n}}_{k=1}\left\vert T^{\left(1\right)}_{1nk}\right\vert ^{\theta}\leq C\left(\frac{\mu^{\theta}_{n}}{\kappa_{n}}\sum^{\delta_{n}}_{k=1}\left\vert T^{\left(1\right)}_{1nk}-T^{\left(3\right)}_{1nk}\right\vert ^{\theta}+\frac{\mu^{\theta}_{n}}{\kappa_{n}}\sum^{\delta_{n}}_{k=1}\left\vert T^{\left(3\right)}_{1nk}\right\vert ^{\theta}\right)
\]
\[
\leq C\left[\frac{\mu^{\theta}_{n}}{\kappa_{n}}\sum^{\delta_{n}}_{k=1}\left(\left\vert \Delta^{\left(2\right)}_{1nk}\right\vert ^{\theta}+\left\vert \Delta^{\left(3\right)}_{1nk}\right\vert ^{\theta}\right)+\frac{\mu^{\theta}_{n}}{\kappa_{n}}\sum^{\delta_{n}}_{k=1}\left\vert T^{\left(3\right)}_{1nk}\right\vert ^{\theta}\right]
\]
\[
=_{\left(3\right)}o\left(\mu^{\theta}_{n}\right)+C\frac{\mu^{\theta}_{n}}{\kappa_{n}}\sum^{\delta_{n}}_{k=1}\left\vert T^{\left(3\right)}_{1nk}\right\vert ^{\theta}
\]
\[
=_{\left(4\right)}o\left(1\right)+C\mu^{\theta}_{n}\left[\frac{1}{\left\vert 1-m\right\vert ^{\theta}}\int^{\delta_{n}/\kappa_{n}}_{1/\kappa_{n}}x^{-\theta\left(m-1\right)}dx+o(1)\right]
\]
\[
=o\left(1\right)+C\mu^{\theta}_{n}\left[x^{-\theta m+\theta+1}\right]^{\delta_{n}/\kappa_{n}}_{1/\kappa_{n}}=o\left(1\right)+C\mu^{\theta}_{n}\left[\left(\delta_{n}/\kappa_{n}\right)^{b_{\theta}}-\left(1/\kappa_{n}\right)^{b_{\theta}}\right]=o(1),
\]
\[
\leq o\left(1\right)+C\mu^{\theta}_{n}\left[\left(\kappa_{n}\varsigma(\kappa_{n})/\kappa_{n}\right)^{b_{\theta}}\right]=O(\mu^{\theta}_{n}\varsigma(\kappa_{n})^{b_{\theta}})=_{(5)}o(1),
\]
where $=_{\left(3\right)}$ follows from (\ref{T1.2}) and (\ref{T1.2b}),
whilst $=_{\left(4\right)}$ follows from Euler summation, see, e.g.
\textbf{Lemma \ref{pi}}. Finally, $=_{(5)}$ is due to (\ref{Ts_rmk}).
Hence, (\ref{Ts1}) is a direct consequence of (\ref{T1.0b}), (\ref{T1.1}),
(\ref{T1.1b}), and (\ref{T1.2b}).

\bigskip{}

\begin{itemize}
\item $T_{2nk}$: Next, we prove (\ref{Ts2}). Let
\[
p>\frac{b_{\theta}}{2-m}\lor1,
\]
where note that $0<2-m<1.$ Hence, 
\begin{equation}
0<b_{\theta}/p<2-m<1.\label{T2.00}
\end{equation}
Set
\[
T_{2nk}:=T^{(0)}_{2nk}:=\frac{1}{\kappa^{2-m}_{n}}\sum^{k}_{l=1}\rho^{k-l}_{n}\sum^{\infty}_{j=l}c_{j},\text{ \ }T^{\left(1\right)}_{2nk}:=\frac{1}{\kappa^{2-m}_{n}}\sum^{k}_{l=\kappa^{b_{\theta}/p}_{n}}\rho^{k-l}_{n}\sum^{\infty}_{j=l}c_{j},
\]
\end{itemize}
\[
T^{\left(2\right)}_{2nk}:=\frac{1}{\kappa^{2-m}_{n}}\sum^{k}_{l=\kappa^{b_{\theta}/p}_{n}}\rho^{k-l}_{n}\sum^{\infty}_{j=l}j^{-m},\text{ \ }T^{\left(3\right)}_{2nk}:=\frac{1}{\kappa_{n}}\sum^{k}_{l=\kappa^{b_{\theta}/p}_{n}}\rho^{k-l}_{n}\int^{\infty}_{l/\kappa_{n}}y^{-m}dy,
\]
\[
T^{\left(4\right)}_{2nk}:=\frac{1}{\kappa_{n}}\sum^{k}_{l=\kappa^{b_{\theta}/p}_{n}}e^{-\frac{k-l}{\kappa_{n}}}\int^{\infty}_{l/\kappa_{n}}y^{-m}dy=-\frac{1}{1-m}\kappa^{-1}_{n}e^{-\frac{k}{\kappa_{n}}}\sum^{k}_{l=\kappa^{b_{\theta}/p}_{n}}e^{\frac{l}{\kappa_{n}}}\left(\frac{l}{\kappa_{n}}\right)^{1-m},
\]
\[
T^{\left(5\right)}_{2nk}:=-\frac{1}{1-m}e^{-\frac{k}{\kappa_{n}}}\int^{k/\kappa_{n}}_{\kappa^{b_{\theta}/p}_{n}/\kappa_{n}}e^{y}y^{1-m}dy,\text{ and}
\]
\begin{equation}
T^{\left(6\right)}_{2nk}:=-\frac{1}{1-m}e^{-\frac{k}{\kappa_{n}}}\int^{k/\kappa_{n}}_{0}e^{y}y^{1-m}dy.\label{T2.001}
\end{equation}
Moreover, for $i=1,...,6$ set
\[
\Delta^{\left(i\right)}_{2nk}:=T^{\left(i\right)}_{2nk}-T^{(i-1)}_{2nk},\text{\ and \ }\Delta_{2nk}:=T_{2nk}-T^{\left(6\right)}_{2nk}.
\]

Consider

\[
\frac{1}{\kappa_{n}}\sum^{\kappa^{b_{\theta}/p}_{n}}_{k=1}\left\vert T_{2nk}\right\vert ^{\theta}=\frac{1}{\kappa_{n}}\sum^{\kappa^{b_{\theta}/p}_{n}}_{k=1}\left\vert \frac{1}{\kappa^{2-m}_{n}}\sum^{k}_{l=1}\rho^{k-l}_{n}\sum^{\infty}_{j=l}c_{j}\right\vert ^{\theta}
\]
\begin{equation}
\leq\frac{1}{\kappa_{n}}\sum^{\kappa^{b_{\theta}/p}_{n}}_{k=1}\left(\frac{1}{\kappa^{2-m}_{n}}\sum^{\kappa^{b_{\theta}/p}_{n}}_{l=1}\sum^{\infty}_{j=0}\left\vert c_{j}\right\vert \right)^{\theta}\leq C\frac{\kappa^{b_{\theta}/p}_{n}}{\kappa_{n}}\left(\frac{\kappa^{b_{\theta}/p}_{n}}{\kappa^{2-m}_{n}}\right)^{\theta}=o(1),\label{T2.0}
\end{equation}
where the last equality follows from (\ref{T2.00}) above. Next, we
shall show that
\begin{equation}
\frac{1}{\kappa_{n}}\sum^{\delta_{n}}_{k=\kappa^{b_{\theta}/p}_{n}}\left\vert \Delta_{2nk}\right\vert ^{\theta}\leq C\frac{1}{\kappa_{n}}\sum^{\delta_{n}}_{k=\kappa^{b_{\theta}/p}_{n}}\sum^{6}_{i=1}\left\vert \Delta^{\left(i\right)}_{2nk}\right\vert ^{\theta}=o(1).\label{T2.0b}
\end{equation}

\begin{itemize}
\item[(i)] $\Delta^{(1)}_{2nk}$: Consider the first summand on the r.h.s. of
(\ref{T2.0b}) 
\[
\frac{1}{\kappa_{n}}\sum^{\delta_{n}}_{k=\kappa^{b_{\theta}/p}_{n}}\left\vert \Delta^{(1)}_{2nk}\right\vert ^{\theta}=\frac{1}{\kappa_{n}}\sum^{\delta_{n}}_{k=\kappa^{b_{\theta}/p}_{n}}\left\vert T_{2nk}-T^{\left(1\right)}_{2nk}\right\vert ^{\theta}
\]
\[
=\frac{1}{\kappa_{n}}\sum^{\delta_{n}}_{k=\kappa^{b_{\theta}/p}_{n}}\left\vert \frac{1}{\kappa^{2-m}_{n}}\sum^{\kappa^{b_{\theta}/p}_{n}-1}_{l=1}\rho^{k-l}_{n}\sum^{\infty}_{j=l}c_{j}\right\vert ^{\theta}
\]
\[
\leq\frac{C}{\kappa_{n}}\sum^{\delta_{n}}_{k=\kappa^{b_{\theta}/p}_{n}}\left(\frac{1}{\kappa^{2-m}_{n}}\sum^{\kappa^{b_{\theta}/p}_{n}-1}_{l=1}\sum^{\infty}_{j=l}\left\vert c_{j}\right\vert \right)^{\theta}\leq\left(\frac{\kappa^{b_{\theta}/p}_{n}}{\kappa^{2-m}_{n}}\right)^{\theta}\frac{C\delta_{n}}{\kappa_{n}}\left(\sum^{\infty}_{j=0}\left\vert c_{j}\right\vert \right)^{\theta}
\]
\[
=C\left(\frac{\kappa^{b_{\theta}/p}_{n}}{\kappa^{2-m}_{n}}\right)^{\theta}\frac{\delta_{n}}{\kappa_{n}}=C\left(\frac{\kappa^{b_{\theta}/p}_{n}}{\kappa^{2-m}_{n}}\right)^{\theta}\varsigma\left(\kappa_{n}\right)=o(1),
\]
where we have used the fact that $b_{\theta}/p<2-m$ and that $\varsigma\left(.\right)$
is slowly varying.
\item[(ii)] $\Delta^{(2)}_{2nk}$:
\[
\frac{1}{\kappa_{n}}\sum^{\delta_{n}}_{k=\kappa^{b_{\theta}/p}_{n}}\left\vert \Delta^{\left(2\right)}_{2nk}\right\vert ^{\theta}=\frac{1}{\kappa_{n}}\sum^{\delta_{n}}_{k=\kappa^{b_{\theta}/p}_{n}}\left\vert \frac{1}{\kappa^{2-m}_{n}}\sum^{k}_{l=\kappa^{b_{\theta}/p}_{n}}\rho^{k-l}_{n}\sum^{\infty}_{j=l}j^{-m}-\frac{1}{\kappa^{2-m}_{n}}\sum^{k}_{l=\kappa^{b_{\theta}/p}_{n}}\rho^{k-l}_{n}\sum^{\infty}_{j=l}c_{j}\right\vert ^{\theta}
\]
\[
\leq\frac{1}{\kappa_{n}}\sum^{\delta_{n}}_{k=\kappa^{b_{\theta}/p}_{n}}\left\vert \frac{1}{\kappa^{2-m}_{n}}\sum^{k}_{l=\kappa^{b_{\theta}/p}_{n}}\rho^{k-l}_{n}\sum^{\infty}_{j=l}j^{-m}\left\vert \frac{c_{j}}{j^{-m}}-1\right\vert \right\vert ^{\theta}
\]
\[
\leq\frac{\mu^{\theta}_{n}}{\kappa_{n}}\sum^{\delta_{n}}_{k=\kappa^{b_{\theta}/p}_{n}}\left\vert \frac{1}{\kappa^{2-m}_{n}}\sum^{k}_{l=\kappa^{b_{\theta}/p}_{n}}\rho^{k-l}_{n}\sum^{\infty}_{j=l}j^{-m}\right\vert ^{\theta}
\]
\begin{equation}
=\frac{\mu^{\theta}_{n}}{\kappa_{n}}\sum^{\delta_{n}}_{k=\kappa^{b_{\theta}/p}_{n}}\left\vert T^{\left(2\right)}_{2nk}\right\vert ^{\theta}=_{\left(5\right)}o(1),\label{Delta2.1}
\end{equation}
where $=_{\left(5\right)}$ will be established in detail later.
\item[(iii)] $\Delta^{\left(3\right)}_{2nk}$: Next, noting that
\[
\left\vert \frac{1}{\kappa_{n}}\sum^{\infty}_{j=l}\left(\frac{j}{\kappa_{n}}\right)^{-m}-\int^{\infty}_{l/\kappa_{n}}x^{-m}dx\right\vert 
\]
\begin{equation}
\leq\kappa^{-1}_{n}\left[\mathbf{V}\left(x^{-m}\right)\right]^{\infty}_{l/\kappa_{n}}\leq\kappa^{-1}_{n}\left(l/\kappa_{n}\right)^{-m},\label{Delta2.1b}
\end{equation}
\end{itemize}
we get
\[
\frac{1}{\kappa_{n}}\sum^{\delta_{n}}_{k=\kappa^{b_{\theta}/p}_{n}}\left\vert \Delta^{\left(3\right)}_{2nk}\right\vert ^{\theta}=\frac{1}{\kappa_{n}}\sum^{\delta_{n}}_{k=\kappa^{b_{\theta}/p}_{n}}\left\vert T^{\left(3\right)}_{2nk}-T^{\left(2\right)}_{2nk}\right\vert ^{\theta}
\]
\[
\leq\frac{1}{\kappa_{n}}\sum^{\delta_{n}}_{k=\kappa^{b_{\theta}/p}_{n}}\left\vert \frac{1}{\kappa_{n}}\sum^{k}_{l=\kappa^{b_{\theta}/p}_{n}}\rho^{k-l}_{n}\int^{\infty}_{l/\kappa_{n}}y^{-m}dy-\frac{1}{\kappa_{n}}\sum^{k}_{l=\kappa^{b_{\theta}/p}_{n}}\rho^{k-l}_{n}\frac{1}{\kappa_{n}}\sum^{\infty}_{j=l}\left(\frac{j}{\kappa_{n}}\right)^{-m}\right\vert ^{\theta}
\]
\begin{equation}
\leq_{_{\left(6\right)}}\frac{C}{\kappa_{n}}\sum^{\delta_{n}}_{k=1}\left\vert \frac{1}{\kappa^{2}_{n}}\sum^{k}_{l=1}\rho^{k-l}_{n}\left(l/\kappa_{n}\right)^{-m}\right\vert ^{\theta}=_{\left(7\right)}o(1),\label{Delta2.2}
\end{equation}
where the $\leq_{\left(6\right)}$ follows from (\ref{Delta2.1b})
whilst $=_{\left(7\right)}$ follows from the fact that 
\[
\frac{C}{\kappa_{n}}\sum^{\delta_{n}}_{k=1}\left\vert \frac{1}{\kappa^{2}_{n}}\sum^{k}_{l=1}\rho^{k-l}_{n}\left(l/\kappa_{n}\right)^{-m}\right\vert ^{\theta}
\]
\[
\leq\frac{C}{\kappa_{n}\kappa^{\theta\left(2-m\right)}_{n}}\sum^{\delta_{n}}_{k=1}\left(\sum^{\infty}_{l=1}l^{-m}\right)^{\theta}=\frac{C\delta_{n}}{\kappa_{n}\kappa^{\theta\left(2-m\right)}_{n}}=o(1),
\]
where the last equality is due to the fact that $2-m>0$ and $\delta_{n}=\kappa_{n}\varsigma(\kappa_{n})$,
with $\varsigma(.)$ SV.
\begin{itemize}
\item[(iv)] $\Delta^{(4)}_{2nk}$:
\[
\frac{1}{\kappa_{n}}\sum^{\delta_{n}}_{k=\kappa^{b_{\theta}/p}_{n}}\left\vert \Delta^{\left(4\right)}_{2nk}\right\vert ^{\theta}
\]
\[
=\frac{1}{\kappa_{n}}\sum^{\delta_{n}}_{k=\kappa^{b_{\theta}/p}_{n}}\left\vert \frac{1}{\kappa_{n}}\sum^{k}_{l=\kappa^{b_{\theta}/p}_{n}}e^{-\frac{k-l}{\kappa_{n}}}\int^{\infty}_{l/\kappa_{n}}y^{-m}dy-\frac{1}{\kappa_{n}}\sum^{k}_{l=\kappa^{b_{\theta}/p}_{n}}\rho^{k-l}_{n}\int^{\infty}_{l/\kappa_{n}}y^{-m}dy\right\vert ^{\theta}
\]
\[
\leq\frac{1}{\kappa_{n}}\sum^{\delta_{n}}_{k=1}\left[\frac{1}{\kappa_{n}}\sum^{k}_{l=\kappa^{b_{\theta}/p}_{n}}\left\vert e^{-\frac{k-l}{\kappa_{n}}}-\rho^{k-l}_{n}\right\vert \int^{\infty}_{l/\kappa_{n}}y^{-m}dy\right]^{\theta}
\]
\[
=\frac{C}{\kappa_{n}}\sum^{\delta_{n}}_{k=1}\left[\frac{1}{\kappa_{n}}\sum^{k}_{l=\kappa^{b_{\theta}/p}_{n}}\left\vert e^{-\frac{k-l}{\kappa_{n}}}-\rho^{k-l}_{n}\right\vert \left(l/\kappa_{n}\right)^{1-m}\right]^{\theta}
\]
\[
\leq_{\left(8\right)}\frac{C}{\kappa_{n}}\sum^{\delta_{n}}_{k=1}\left\vert \kappa^{-1}_{n}\sum^{k}_{l=\kappa^{b_{\theta}/p}_{n}}e^{-\frac{k-l}{\kappa_{n}}}\frac{k-l}{\kappa^{2}_{n}}\left(\frac{l}{\kappa_{n}}\right)^{1-m}\right\vert ^{\theta}
\]
\[
\leq C\left(\frac{\delta_{n}}{\kappa^{2}_{n}}\right)^{\theta}\frac{1}{\kappa_{n}}\sum^{\delta_{n}}_{k=1}\left\vert \kappa^{-1}_{n}\sum^{k}_{l=\kappa^{b_{\theta}/p}_{n}}e^{-\frac{k-l}{\kappa_{n}}}\left(\frac{l}{\kappa_{n}}\right)^{1-m}\right\vert ^{\theta}
\]
\begin{equation}
=C\left(\frac{\delta_{n}}{\kappa^{2}_{n}}\right)^{\theta}\frac{1}{\kappa_{n}}\sum^{\delta_{n}}_{k=1}\left\vert T^{\left(4\right)}_{2nk}\right\vert ^{\theta}=_{\left(9\right)}o(1),\label{Delta2.3}
\end{equation}
where $\leq_{\left(8\right)}$ follows from DK2018, Lemma 6.4 whilst
$=_{(9)}$ will be proved in detail later.
\item[(v)] $\Delta^{(5)}_{2nk}$: Set $\omega\left(y\right):=e^{y}y^{1-m}$.
We have
\[
\frac{1}{\kappa_{n}}\sum^{\delta_{n}}_{k=\kappa^{b_{\theta}/p}_{n}}\left\vert \Delta^{\left(5\right)}_{2nk}\right\vert ^{\theta}=\frac{1}{\kappa_{n}}\sum^{\delta_{n}}_{k=\kappa^{b_{\theta}/p}_{n}}\left\vert T^{\left(5\right)}_{2nk}-T^{\left(4\right)}_{2nk}\right\vert ^{\theta}
\]
\[
=\left\vert \frac{1}{1-m}\right\vert ^{\theta}\frac{1}{\kappa_{n}}\sum^{\delta_{n}}_{k=\kappa^{b_{\theta}/p}_{n}}\left\vert e^{-\frac{k}{\kappa_{n}}}\int^{k/\kappa_{n}}_{\kappa^{b_{\theta}/p}_{n}/\kappa_{n}}e^{y}y^{1-m}dy-\frac{e^{-\frac{k}{\kappa_{n}}}}{\kappa_{n}}\sum^{k}_{l=\kappa^{b_{\theta}/p}_{n}}e^{\frac{l}{\kappa_{n}}}\left(\frac{l}{\kappa_{n}}\right)^{1-m}\right\vert ^{\theta}
\]
\end{itemize}
\begin{equation}
\leq\frac{C}{\kappa_{n}}\sum^{\delta_{n}}_{k=1}\left[\frac{e^{-\frac{k}{\kappa_{n}}}}{\kappa_{n}}\mathbf{V}\left[\omega\left(y\right)\right]^{k/\kappa_{n}}_{\kappa^{b_{\theta}/p}_{n}/\kappa_{n}}\right]^{\theta}\leq\frac{C}{\kappa_{n}}\sum^{\delta_{n}}_{k=1}\left[\frac{e^{-\frac{k}{\kappa_{n}}}}{\kappa_{n}}\mathbf{V}\left[\omega\left(y\right)\right]^{k/\kappa_{n}}_{1/\kappa_{n}}\right]^{\theta},\label{Delta3.1}
\end{equation}
where recall that $\mathbf{V}\left[\omega\left(y\right)\right]^{k/\kappa_{n}}_{1/\kappa_{n}}$
is the total variation of $\omega\left(y\right)$ on $[1/\kappa_{n},k/\kappa_{n}]$.
It follows from the proof of \textbf{Proposition \ref{TotalVar}}
that 
\[
e^{-\frac{k}{\kappa_{n}}}\mathbf{V}\left[\omega\left(y\right)\right]^{k/\kappa_{n}}_{1/\kappa_{n}}\leq2e^{-\frac{k}{\kappa_{n}}}\left[\omega(k/\kappa_{n})+\omega(m-1)+\omega(1/\kappa_{n})\right].
\]
Hence,
\[
\frac{1}{\kappa_{n}}\sum^{\delta_{n}}_{k=1}\left[\frac{e^{-\frac{k}{\kappa_{n}}}}{\kappa_{n}}\mathbf{V}\left[\omega\left(y\right)\right]^{k/\kappa_{n}}_{1/\kappa_{n}}\right]^{\theta}\leq\frac{C}{\kappa_{n}}\sum^{\delta_{n}}_{k=1}\left[\frac{e^{-\frac{k}{\kappa_{n}}}}{\kappa_{n}}\left[\omega(k/\kappa_{n})+\omega(m-1)+\omega(1/\kappa_{n})\right]\right]^{\theta}
\]
\[
=\frac{C}{\kappa_{n}}\sum^{\delta_{n}}_{k=1}\left[\frac{1}{\kappa_{n}}\left[(k/\kappa_{n})^{1-m}+e^{-\frac{k}{\kappa_{n}}}\omega(m-1)+e^{-\frac{k}{\kappa_{n}}}\omega(1/\kappa_{n})^{1-m}\right]\right]^{\theta}
\]
\[
\leq\frac{C}{\kappa_{n}}\sum^{\delta_{n}}_{k=1}\left[\frac{1}{\kappa_{n}}(k/\kappa_{n})^{1-m}\right]^{\theta}+\frac{C}{\kappa_{n}}\sum^{\delta_{n}}_{k=1}\left[\frac{1}{\kappa_{n}}e^{-\frac{k}{\kappa_{n}}}\omega(m-1)\right]^{\theta}
\]
\[
+\frac{C}{\kappa_{n}}\sum^{\delta_{n}}_{k=1}\left[\frac{1}{\kappa_{n}}e^{-\frac{k-1}{\kappa_{n}}}(1/\kappa_{n})^{1-m}\right]^{\theta}=:R_{1n}+R_{2n}+R_{3n}.
\]
We shall show that
\[
R_{in}=o(1),\,i=1,2,3,
\]
which in view of (\ref{Delta3.1}) implies that 
\begin{equation}
\frac{1}{\kappa_{n}}\sum^{\delta_{n}}_{k=\kappa^{b_{\theta}/p}_{n}}\left\vert \Delta^{\left(5\right)}_{2nk}\right\vert ^{\theta}=o(1).\label{Delta3.2}
\end{equation}
The first term
\[
R_{1n}=\frac{C}{\kappa^{\theta+b_{\theta}}_{n}}\sum^{\delta_{n}}_{k=1}k^{\theta\left(1-m\right)}=\frac{C}{\kappa^{\theta+b_{\theta}}_{n}}\sum^{\delta_{n}}_{k=1}k^{-\theta\left(m-1\right)}
\]
\[
=\frac{C}{\kappa^{\theta+b_{\theta}}_{n}}\delta^{1+\theta\left(1-m\right)}_{n}\frac{1}{\delta_{n}}\sum^{\delta_{n}}_{k=1}\left(\frac{k}{\delta_{n}}\right)^{\theta\left(1-m\right)}=\frac{C}{\kappa^{\theta+b_{\theta}}_{n}}\delta^{b_{\theta}}_{n}\left[\int^{1}_{0}x^{-\theta\left(m-1\right)}dx+o(1)\right]
\]
\[
\leq\frac{C\varsigma\left(\kappa_{n}\right)^{b_{\theta}}}{\kappa^{\theta}_{n}}\left[\int^{1}_{0}x^{-\theta\left(m-1\right)}dx+o(1)\right]=o(1),
\]
where we have used the fact that $\theta(m-1)<1$. The second term
\[
R_{2n}=\frac{C}{\kappa_{n}}\sum^{\delta_{n}}_{k=1}\left[\frac{1}{\kappa_{n}}e^{-\frac{k}{\kappa_{n}}}\right]^{\theta}
\]
\[
=\frac{C}{\kappa^{\theta}_{n}}\left[\int^{\infty}_{0}e^{-\theta x}dx+o(1)\right]=O\left(\kappa^{-\theta}_{n}\right)=o(1).
\]
Finally, 
\[
R_{3n}=\frac{C}{\kappa_{n}}\sum^{\delta_{n}}_{k=1}\left[\frac{1}{\kappa_{n}}e^{-\frac{k-1}{\kappa_{n}}}(1/\kappa_{n})^{1-m}\right]^{\theta}
\]
\[
=\frac{C}{\kappa^{\theta\left(2-m\right)}_{n}}\left[\int^{\infty}_{0}e^{-\theta x}dx+o(1)\right]=O\left(\kappa^{\theta\left(2-m\right)}_{n}\right)=o(1).
\]
Hence, (\ref{Delta3.2}) holds.
\begin{itemize}
\item[(vi)]  Consider
\[
\frac{1}{\kappa_{n}}\sum^{\delta_{n}}_{k=\kappa^{b_{\theta}/p}_{n}}\left\vert \Delta^{\left(6\right)}_{2nk}\right\vert ^{\theta}=\frac{1}{\kappa_{n}}\sum^{\delta_{n}}_{k=\kappa^{b_{\theta}/p}_{n}}\left\vert T^{\left(6\right)}_{2nk}-T^{\left(5\right)}_{2nk}\right\vert ^{\theta}
\]
\[
=C\frac{1}{\kappa_{n}}\sum^{\delta_{n}}_{k=\kappa^{b_{\theta}/p}_{n}}\left\vert e^{-\frac{k}{\kappa_{n}}}\int^{\kappa^{b_{\theta}/p}_{n}/\kappa_{n}}_{0}e^{y}y^{1-m}dy\right\vert ^{\theta}\leq Ce^{\theta\frac{\kappa^{b_{\theta}/p}_{n}}{\kappa_{n}}}\left(\int^{\kappa^{b_{\theta}/p}_{n}/\kappa_{n}}_{0}y^{1-m}dy\right)\cdot\frac{1}{\kappa_{n}}\sum^{\delta_{n}}_{k=\kappa^{b_{\theta}/p}_{n}}e^{-\theta\frac{k}{\kappa_{n}}}
\]
\[
=C\left(\frac{\kappa^{b_{\theta}/p}_{n}}{\kappa_{n}}\right)^{2-m}\left[\int^{\delta_{n}/\kappa_{n}}_{\kappa^{b_{\theta}/p}_{n}/\kappa_{n}}e^{-\theta x}dx+o(1)\right]=o(1)\cdot\left[\int^{\infty}_{0}e^{-\theta x}dx+o(1)\right]=o(1),
\]
\end{itemize}
\noindent where we have used the facts that $\beta_{\theta}/p\in(0,1)$
and $2-m>0$.

Finally, we establish (\ref{Delta2.1}) and (\ref{Delta2.3}). We
start with the latter. Write 
\[
\left(\frac{\delta_{n}}{\kappa^{2}_{n}}\right)^{\theta}\frac{1}{\kappa_{n}}\sum^{\delta_{n}}_{k=\kappa^{b_{\theta}/p}_{n}}\left\vert T^{\left(4\right)}_{2nk}\right\vert ^{\theta}\leq\left(\frac{\delta_{n}}{\kappa^{2}_{n}}\right)^{\theta}\frac{C}{\kappa_{n}}\sum^{\delta_{n}}_{k=\kappa^{b_{\theta}/p}_{n}}\left[\left\vert T^{\left(5\right)}_{2nk}-T^{\left(4\right)}_{2nk}\right\vert ^{\theta}+\left\vert T^{\left(6\right)}_{2nk}-T^{\left(5\right)}_{2nk}\right\vert ^{\theta}\right]
\]
\[
+\left(\frac{\delta_{n}}{\kappa^{2}_{n}}\right)^{\theta}\frac{C}{\kappa_{n}}\sum^{\delta_{n}}_{k=\kappa^{b_{\theta}/p}_{n}}\left\vert T^{\left(6\right)}_{2nk}\right\vert ^{\theta}=_{\left(10\right)}o\left(\left(\frac{\delta_{n}}{\kappa^{2}_{n}}\right)^{\theta}\right)+C\left(\frac{\delta_{n}}{\kappa^{2}_{n}}\right)^{\theta}\left[\int^{\delta_{n}/\kappa_{n}}_{\kappa^{b_{\theta}/p}_{n}/\kappa_{n}}\left[e^{-x}\int^{x}_{0}y^{1-m}dy\right]^{\theta}dx+o(1)\right]
\]
\[
\leq C\left(\frac{\delta_{n}}{\kappa^{2}_{n}}\right)^{\theta}\int^{\delta_{n}/\kappa_{n}}_{1/\kappa_{n}}\left[e^{-x}\int^{x}_{0}e^{y}y^{1-m}dy\right]^{\theta}dx+o(1)
\]
\[
\leq_{\left(11\right)}C\left(\frac{\delta_{n}}{\kappa^{2}_{n}}\right)^{\theta}\int^{\delta_{n}/\kappa_{n}}_{1}x^{-\theta\left(m-1\right)}dx+o\left(1\right)
\]
\[
=C\left(\frac{\delta_{n}}{\kappa^{2}_{n}}\right)^{\theta}\left(\delta_{n}/\kappa_{n}\right)^{-\theta\left(m-1\right)+1}+o\left(1\right)=C\left(\frac{\delta_{n}}{\kappa^{2}_{n}}\right)^{\theta}\left(\delta_{n}/\kappa_{n}\right)^{b_{\theta}}+o\left(1\right)
\]
\[
=C\frac{\varsigma\left(\kappa_{n}\right)^{\theta+b_{\theta}}}{\kappa^{\theta}_{n}}+o\left(1\right)=o\left(1\right),
\]
where $=_{\left(10\right)}$ follows from parts (v)-(vi) above and
an Euler summation argument (see e.g. \textbf{Lemma \ref{pi}}), whilst
$\leq_{\left(11\right)}$ from the fact that $e^{-x}\int^{x}_{0}e^{y}y^{1-m}dy\sim x^{1-m}$
(cf. \textbf{Proposition \ref{SterlingApprox}}).
\begin{description}
\item [{Next,}] we conclude with the proof of $=_{\left(5\right)}$ in
(\ref{Delta2.1}). Using similar arguments as those used above for
the proof of (\ref{Delta2.3}) we get
\[
\frac{\mu^{\theta}_{n}}{\kappa_{n}}\sum^{\delta_{n}}_{k=\kappa^{b_{\theta}/p}_{n}}\left\vert T^{\left(2\right)}_{2nk}\right\vert ^{\theta}\leq C\left\{ \frac{\mu^{\theta}_{n}}{\kappa_{n}}\sum^{\delta_{n}}_{k=\kappa^{b_{\theta}/p}_{n}}\sum^{6}_{i=3}\left\vert T^{\left(i\right)}_{2nk}-T^{\left(i-1\right)}_{2nk}\right\vert ^{\theta}+\frac{\mu^{\theta}_{n}}{\kappa_{n}}\sum^{\delta_{n}}_{k=\kappa^{b_{\theta}/p}_{n}}\left\vert T^{\left(6\right)}_{2nk}\right\vert \right\} 
\]
\[
=C\left\{ \frac{\mu^{\theta}_{n}}{\kappa_{n}}\sum^{\delta_{n}}_{k=\kappa^{b_{\theta}/p}_{n}}\sum^{6}_{i=3}\left\vert \Delta^{\left(i\right)}_{2nk}\right\vert ^{\theta}+\frac{\mu^{\theta}_{n}}{\kappa_{n}}\sum^{\delta_{n}}_{k=\kappa^{b_{\theta}/p}_{n}}\left\vert T^{\left(6\right)}_{2nk}\right\vert ^{\theta}\right\} 
\]
\[
=o\left(\mu^{\theta}_{n}\right)+\frac{\mu^{\theta}_{n}}{\kappa_{n}}\sum^{\delta_{n}}_{k=\kappa^{b_{\theta}/p}_{n}}\left\vert T^{\left(6\right)}_{2nk}\right\vert ^{\theta}
\]
\[
=o\left(1\right)+\mu^{\theta}_{n}\left[\int^{\delta_{n}/\kappa_{n}}_{\kappa^{b_{\theta}/p}_{n}/\kappa_{n}}\left[e^{-x}\int^{x}_{0}y^{1-m}dy\right]^{\theta}dx+o(1)\right]
\]
\[
\leq_{(11)}C\mu^{\theta}_{n}\int^{\delta_{n}/\kappa_{n}}_{1}x^{-\theta\left(m-1\right)}dx+o\left(1\right)=C\mu^{\theta}_{n}\left(\delta_{n}/\kappa_{n}\right)^{b_{\theta}}+o\left(1\right)
\]
\[
\leq C\mu^{\theta}_{n}\varsigma(\kappa_{n})^{b_{\theta}}+o\left(1\right)=o\left(1\right),
\]
where $\leq_{(11)}$ follows from \textbf{Proposition \ref{SterlingApprox}}.
This completes the proof. $\blacksquare$
\end{description}
\bigskip{}

\textbf{Proof of Corollary \ref{Ts_corol}}. The result follows directly
from \textbf{Lemma \ref{Ts}} if we set $\theta=\alpha$ with $\alpha$
as in \textbf{Assumption MI-NM}. Notice in this case the conditions
$m\in\left(1,2\right)$ and $1/m<\theta<1/(m-1)$ of \textbf{Lemma
\ref{Ts}} are satisfied. Recall that under \textbf{Assumption MI-NM}
$\alpha\in(1,2]$ and $m\in(1,1+1/\alpha)$. Hence, we do have $1/m<\alpha<1/(m-1)$
and $m\in\left(1,2\right)$. $\blacksquare$

\bigskip{}

\textbf{Proof of Lemma \ref{pi}}. Set
\[
\mathcal{P}(x)=\left\vert p\left(x\right)\right\vert ^{\theta}.
\]
We shall show that the total variation of $\mathcal{P}(x)$ on $[1/\kappa_{n},\delta_{n}/\kappa_{n}]$
satisfies
\begin{equation}
\frac{1}{\kappa_{n}}\mathbf{V}\left[\mathcal{P}(x)\right]^{\delta_{n}/\kappa_{n}}_{1/\kappa_{n}}\rightarrow0,\label{pi1}
\end{equation}
as $n\rightarrow\infty$. Noting that 
\begin{equation}
\left\vert \frac{1}{\kappa_{n}}\sum^{\delta_{n}}_{k=1}\left\vert p\left(k/\kappa_{n}\right)\right\vert ^{\theta}-\int^{\delta_{n}/\kappa_{n}}_{1/\kappa_{n}}\left\vert p\left(x\right)\right\vert ^{\theta}dx\right\vert \leq\frac{1}{\kappa_{n}}\mathbf{V}\left[\mathcal{P}(x)\right]^{\delta_{n}/\kappa_{n}}_{1/\kappa_{n}},\label{pi2}
\end{equation}
the requisite result follows from (\ref{pi1}), (\ref{pi2}) and dominated
convergence.

We set out to prove (\ref{pi1}). By \textbf{Proposition \ref{pi-roots}},
$p\left(x\right)$ has a single root, $x_{\ast}$, such that $x_{\ast}>2-m$,
and
\[
\left|p(x)\right|=\left\{ \begin{array}{l}
p\left(x\right)\text{, }0<x\leq x_{\ast}\\
-p\left(x\right)\text{, }x>x_{\ast}
\end{array}\right..
\]
Hence,
\[
\mathcal{P}(x)=\left\{ \begin{array}{l}
p\left(x\right)^{\theta}\text{, }0<x\leq x_{\ast}\\
\left(-p\left(x\right)\right)^{\theta}\text{, }x>x_{\ast}
\end{array}\right..
\]
Therefore, for $n$ large enough (such that $1/\kappa_{n}<x_{\ast}<\delta_{n}/\kappa_{n}$)
we can write
\[
\mathbf{V}\left[\mathcal{P}(x)\right]^{\delta_{n}/\kappa_{n}}_{1/\kappa_{n}}=\mathbf{V}\left[\mathcal{P}(x)\right]^{x_{\ast}}_{1/\kappa_{n}}+\mathbf{V}\left[\mathcal{P}(x)\right]^{\delta_{n}/\kappa_{n}}_{x_{\ast}}
\]
\begin{equation}
=\int^{x_{\ast}}_{1/\kappa_{n}}\left\vert \frac{\mathcal{\partial P}(x)}{\partial x}\right\vert dx+\int^{\delta_{n}/\kappa_{n}}_{x_{\ast}}\left\vert \frac{\mathcal{\partial P}(x)}{\partial x}\right\vert dx.\label{pi3}
\end{equation}
Denote the derivative of $p(x)$ as $p(x)^{\prime}$. In view of the
above,
\[
\mathbf{V}\left[\mathcal{P}(x)\right]^{\delta_{n}/\kappa_{n}}_{1/\kappa_{n}}=\int^{x_{\ast}}_{1/\kappa_{n}}\left\vert \theta p\left(x\right)^{\theta-1}p(x)^{\prime}\right\vert dx+\int^{\delta_{n}/\kappa_{n}}_{x_{\ast}}\left\vert \left[-\theta p\left(x\right)\right]^{\theta-1}p(x)^{\prime}\right\vert dx
\]
\[
=\int^{x_{\ast}}_{1/\kappa_{n}}\left|p\left(x\right)\right|^{\theta-1}\left\vert \theta p(x)^{\prime}\right\vert dx+\int^{\delta_{n}/\kappa_{n}}_{x_{\ast}}\left|p\left(x\right)\right|^{\theta-1}\left\vert \theta p(x)^{\prime}\right\vert dx.
\]
Next, consider
\[
p(x)^{\prime}=(1-m)x^{-m}+e^{-x}\int^{x}_{0}e^{y}y^{1-m}dy-x^{1-m}.
\]
First, suppose that $x\geq x_{\ast}$. As $x\rightarrow\infty$ we
have
\[
p(x)^{\prime}=(1-m)x^{-m}+e^{-x}\int^{x_{\ast}}_{0}e^{y}y^{1-m}dy+e^{-x}\int^{x}_{x_{\ast}}e^{y}y^{1-m}dy-x^{1-m}
\]
\[
=(1-m)x^{-m}+O\left(e^{-x}\right)+e^{-x}\left[\left[e^{y}y^{1-m}\right]^{x}_{x_{\ast}}-\left(1-m\right)\int^{x}_{x_{\ast}}e^{y}y^{-m}dy\right]-x^{1-m}
\]
\[
=(1-m)x^{-m}+O\left(e^{-x}\right)+e^{-x}\left[e^{x}x^{1-m}-e^{x_{\ast}}x^{1-m}_{\ast}-\left(1-m\right)\int^{x}_{x_{\ast}}e^{y}y^{-m}dy\right]-x^{1-m}
\]
\[
=(1-m)x^{-m}+O\left(e^{-x}\right)-\left(1-m\right)e^{-x}\int^{x}_{x_{\ast}}e^{y}y^{-m}dy=o(x^{-m}),
\]
where we have used the fact that $e^{-x}\int^{x}_{x_{\ast}}e^{y}y^{-m}dy\sim x^{-m}$
(cf. \textbf{Proposition \ref{SterlingApprox}}). Furthermore, it
follows from the proof of \textbf{Lemma \ref{LemmaMI2}} that 
\[
p(x)\asymp x^{-m}.
\]
Hence, as $n\rightarrow\infty$
\[
\frac{1}{\kappa_{n}}\int^{\delta_{n}/\kappa_{n}}_{x_{\ast}}\left\vert \theta p\left(x\right)^{\theta-1}p(x)^{\prime}\right\vert dx\leq\frac{C}{\kappa_{n}}\int^{\delta_{n}/\kappa_{n}}_{x_{\ast}}\left(x^{-m}\right)^{\theta-1}x^{-m}dx
\]
\begin{equation}
=\frac{C}{\kappa_{n}}\int^{\delta_{n}/\kappa_{n}}_{x_{\ast}}x^{-\theta m}dx\leq\frac{C}{\kappa_{n}}\int^{\infty}_{x_{\ast}}x^{-\theta m}dx=o(1).\label{pi4}
\end{equation}
Next, note that for $x>0$
\[
\left\vert p\left(x\right)\right\vert ^{\theta-1}\leq C\left(x^{1-m}+x^{2-m}\right)^{\theta-1}=C\left[x^{1-m}\left(1+x\right)\right]^{\theta-1}
\]
\[
\leq_{\left(1\right)}C\cdot\left\{ \begin{array}{c}
x^{\left(1-m\right)\left(\theta-1\right)}+x^{\left(2-m\right)\left(\theta-1\right)},\text{ }\theta\geq1\\
x^{\left(m-2\right)\left(1-\theta\right)},\text{ }0<\theta<1
\end{array}\right.\leq C\cdot\left[x^{\left(1-m\right)\left(\theta-1\right)}+x^{\left(2-m\right)\left(\theta-1\right)}\right].
\]
where $\leq_{\left(1\right)}$ follows from standard arguments\footnote{In particular, for $\theta<1$ note that
\[
\left(x^{1-m}+x^{2-m}\right)^{\theta-1}=\left[x^{1-m}\left(1+x\right)\right]^{\theta-1}=\left[x^{-\left(m-1\right)}\left(1+x\right)\right]^{\theta-1}
\]
\[
=\left[\frac{1+x}{x^{m-1}}\right]^{\theta-1}=\left[\frac{x^{m-1}}{1+x}\right]^{1-\theta}\leq\left[\frac{x^{m-1}}{x}\right]^{1-\theta}=\left[x^{m-2}\right]^{1-\theta}.
\]
}. In addition, for $x>0$
\[
\left\vert p(x)^{\prime}\right\vert \leq C\left(x^{-m}+x^{1-m}+x^{2-m}\right).
\]
Hence,
\begin{equation}
\frac{1}{\kappa_{n}}\int^{x_{\ast}}_{1/\kappa_{n}}\left|p\left(x\right)\right|^{\theta-1}\left\vert \theta p(x)^{\prime}\right\vert dx\label{pi5}
\end{equation}
\[
\leq\frac{C}{\kappa_{n}}\int^{x_{\ast}}_{1/\kappa_{n}}\left[x^{\left(\theta-1\right)\left(1-m\right)}+x^{\left(\theta-1\right)\left(2-m\right)}\right]\left(x^{-m}+x^{1-m}+x^{2-m}\right)dx.
\]
The leading term above is 
\[
\frac{C}{\kappa_{n}}\int^{x_{\ast}}_{1/\kappa_{n}}\left[x^{\left(\theta-1\right)\left(1-m\right)-m}+x^{\left(\theta-1\right)\left(2-m\right)-m}\right]dx
\]
\[
=\frac{C}{\kappa_{n}}\int^{x_{\ast}}_{1/\kappa_{n}}\left[x^{\theta-\theta m-1+m-m}+x^{2\theta-\theta m-2+m-m}\right]dx
\]
\[
=\frac{C}{\kappa_{n}}\int^{x_{\ast}}_{1/\kappa_{n}}\left[x^{\theta-\theta m-1}+x^{2\theta-\theta m-2}\right]dx=\frac{C}{\kappa_{n}}\left[\left(\frac{1}{\kappa_{n}}\right)^{\theta-\theta m}+\left(\frac{1}{\kappa_{n}}\right)^{2\theta-\theta m-1}\right]+O(\kappa^{-1}_{n})
\]
\[
=C\left[\left(\frac{1}{\kappa_{n}}\right)^{1+\theta-\theta m}+\left(\frac{1}{\kappa_{n}}\right)^{\theta\left(2-m\right)}\right]+o(1)
\]
\begin{equation}
=C\left[\left(\frac{1}{\kappa_{n}}\right)^{b_{\theta}}+\left(\frac{1}{\kappa_{n}}\right)^{\theta\left(2-m\right)}\right]+o(1)=o(1),\label{pi6}
\end{equation}
where we have used the fact that $b_{\theta}:=1+\theta-\theta m>0$,
cf. (\ref{T1.00}), and $\theta\left(2-m\right)>0$. Equation (\ref{pi1})
follows directly from (\ref{pi3})-(\ref{pi6}). $\blacksquare$

\bigskip{}

\textbf{Proof of Lemma \ref{MI-NM Scale}}. We start with some preliminary
facts. Recall that (cf. (\ref{MVT})) for all $\theta>0$ and some
$C_{\theta}\in(0,\infty)$
\[
\left\vert \left\vert Y\right\vert ^{\theta}-\left\vert X\right\vert ^{\theta}\right\vert \leq C_{\theta}\cdot\left\{ \begin{array}{l}
\left(\left\vert X\right\vert ^{\theta-1}\left\vert Y-X\right\vert +\left\vert Y-X\right\vert ^{\theta}\right),\text{ }\theta>1\\
\left\vert Y-X\right\vert ^{\theta},\text{ }0<\theta\leq1
\end{array}\right.\text{.}
\]
In addition, for any two sequences $X_{k}$ and $Y_{k}$, such that
$\frac{1}{\kappa_{n}}\sum^{\delta_{n}}_{k=1}\left\vert Y_{k}-X_{k}\right\vert ^{\theta}=o(1)$
and $\frac{1}{\kappa_{n}}\sum^{\delta_{n}}_{k=1}\left\vert X_{k}\right\vert ^{\theta}=O(1)$,
Hölder's inequality gives
\[
\frac{1}{\kappa_{n}}\sum^{\delta_{n}}_{k=1}\left\vert X_{k}\right\vert ^{\theta-1}\left\vert Y_{k}-X_{k}\right\vert \leq\left[\frac{1}{\kappa_{n}}\sum^{\delta_{n}}_{k=1}\left(\left\vert X_{k}\right\vert ^{\theta-1}\right)^{\frac{\theta}{\theta-1}}\right]^{\frac{\theta-1}{\theta}}\left[\frac{1}{\kappa_{n}}\sum^{\delta_{n}}_{k=1}\left\vert Y_{k}-X_{k}\right\vert ^{\theta}\right]^{1/\theta}
\]
\begin{equation}
=\left[\frac{1}{\kappa_{n}}\sum^{\delta_{n}}_{k=1}\left\vert X_{k}\right\vert ^{\theta}\right]^{\frac{\theta-1}{\theta}}\left[\frac{1}{\kappa_{n}}\sum^{\delta_{n}}_{k=1}\left\vert Y_{k}-X_{k}\right\vert ^{\theta}\right]^{1/\theta}=o(1).\label{MVT2}
\end{equation}

In view of the above, set \textbf{{[}check the sign{]}}
\[
Y_{k}=T_{1nk}+T_{2nk}\text{\ and \ }X_{k}=T^{\left(3\right)}_{1nk}+T^{(6)}_{2nk}=\frac{1}{1-m}p(k/\kappa_{n}),
\]
where $T_{1nk}$, $T_{2nk}$ are given \textbf{Lemma \ref{Ts}}, $T^{\left(3\right)}_{1nk}$
and $T^{(6)}_{2nk}$ are defined in (\ref{T1.001}) and (\ref{T2.001}),
respectively, and $p(x)$\
is defined in \textbf{Lemma \ref{pi}}. Hence, by \textbf{Lemma \ref{pi}},
\begin{equation}
\frac{1}{\kappa_{n}}\sum^{\delta_{n}}_{k=1}\left\vert X_{k}\right\vert ^{\theta}=\frac{1}{\left(m-1\right)^{\theta}}\int^{\infty}_{0}\left\vert p\left(x\right)\right\vert ^{\theta}dx+o(1)=O(1).\label{MVT3}
\end{equation}
Moreover, by \textbf{Lemma \ref{Ts}} as $n\rightarrow\infty$ we
have
\begin{equation}
\frac{1}{\kappa_{n}}\sum^{\delta_{n}}_{k=1}\left\vert Y_{k}-X_{k}\right\vert ^{\theta}\leq\frac{2}{\kappa_{n}}\sum^{\delta_{n}}_{k=1}\left[\left\vert \Delta_{1nk}\right\vert ^{\theta}+\left\vert \Delta_{2nk}\right\vert ^{\theta}\right]=o(1),\label{MVT4}
\end{equation}
Hence,
\[
\frac{1}{\kappa_{n}}\sum^{\delta_{n}}_{k=1}\left\vert T_{1nk}-T_{2nk}\right\vert ^{\theta}=\frac{1}{\kappa_{n}}\sum^{\delta_{n}}_{k=1}\left\vert T^{(3)}_{1nk}-T^{(6)}_{2nk}\right\vert ^{\theta}
\]
\[
+\frac{1}{\kappa_{n}}\sum^{\delta_{n}}_{k=1}\left\vert T_{1nk}-T_{2nk}\right\vert ^{\theta}-\frac{1}{\kappa_{n}}\sum^{\delta_{n}}_{k=1}\left\vert T^{(3)}_{1nk}-T^{(6)}_{2nk}\right\vert ^{\theta}
\]
\[
=_{(1)}\frac{1}{\kappa_{n}}\sum^{\delta_{n}}_{k=1}\left\vert T^{(3)}_{1nk}-T^{(6)}_{2nk}\right\vert ^{\theta}+o(1)=_{(2)}\frac{1}{\left(m-1\right)^{\theta}}\int^{\infty}_{0}\left\vert p\left(x\right)\right\vert ^{\theta}dx+o(1),
\]
where $=_{(1)}$ follows from (\ref{MVT2})-(\ref{MVT4}) and $=_{(2)}$
from \textbf{Lemma \ref{pi}}. $\blacksquare$

\bigskip{}

\begin{proposition} \label{TotalVar} Suppose that $m\in(1,1+1/\alpha)$,
$\alpha\in(1,2]$ (i.e. \textbf{NM} holds). Set $\omega(y):=e^{y}y^{1-m}$.
The total variation of $\omega(y)$ on $[1/\kappa_{n},k/\kappa_{n}]$,
$k\in\mathbb{N}$ satisfies
\[
\mathbf{V}\left[\omega(x)\right]^{k/\kappa_{n}}_{1/\kappa_{n}}\leq2\left[\omega(k/\kappa_{n})+\omega(m-1)+\omega(1/\kappa_{n})\right].
\]
\end{proposition}

\bigskip{}

\textbf{Proof of Proposition \ref{TotalVar}}. First, note that 
\[
\left\vert \omega(y)^{\prime}\right\vert =\left\vert e^{y}y^{1-m}+\left(1-m\right)e^{y}y^{-m}\right\vert 
\]
\[
=e^{y}y^{-m}\left\vert y+1-m\right\vert =\left\{ \begin{array}{c}
e^{y}y^{-m}\left(y+1-m\right),\text{ }y\geq m-1\\
-e^{y}y^{-m}\left(y+1-m\right),\text{ }y<m-1
\end{array}\right.=\left\{ \begin{array}{c}
\omega(y)^{\prime},\text{ }y\geq m-1\\
-\omega(y)^{\prime},\text{ }y<m-1
\end{array}\right..
\]
In view of the above,
\[
\mathbf{V}\left[\omega(x)\right]^{k/\kappa_{n}}_{1/\kappa_{n}}=\int^{k/\kappa_{n}}_{1/\kappa_{n}}\left\vert \omega(y)^{\prime}\right\vert dy
\]
\[
=\left[\int^{k/\kappa_{n}}_{m-1}\omega(y)^{\prime}dy-\int^{m-1}_{1/\kappa_{n}}\omega(y)^{\prime}dy\right]I_{\left\{ k/\kappa_{n}\geq m-1\right\} }-\int^{k/\kappa_{n}}_{1/\kappa_{n}}\omega(y)^{\prime}dyI_{\left\{ k/\kappa_{n}<m-1\right\} }
\]
\[
=\left[\left[\omega(y)\right]^{k/\kappa_{n}}_{m-1}-\left[\omega(y)\right]^{m-1}_{1/\kappa_{n}}\right]I_{\left\{ k/\kappa_{n}\geq m-1\right\} }-\left[\omega(y)\right]^{k/\kappa_{n}}_{1/\kappa_{n}}I_{\left\{ k/\kappa_{n}<m-1\right\} }
\]
\[
=\left\{ \left[\omega(k/\kappa_{n})-\omega(m-1)\right]-\left[\omega(m-1)-\omega(1/\kappa_{n})\right]\right\} I_{\left\{ k/\kappa_{n}\geq m-1\right\} }-\left[\omega(k/\kappa_{n})-\omega(1/\kappa_{n})\right]I_{\left\{ k/\kappa_{n}<m-1\right\} }
\]
\[
=\left\{ \omega(k/\kappa_{n})-2\omega(m-1)+\omega(1/\kappa_{n})\right\} I_{\left\{ k/\kappa_{n}\geq m-1\right\} }-\left[\omega(k/\kappa_{n})-\omega(1/\kappa_{n})\right]I_{\left\{ k/\kappa_{n}<m-1\right\} }
\]
\[
=\left\{ \omega(k/\kappa_{n})-2\omega(m-1)\right\} I_{\left\{ k/\kappa_{n}\geq m-1\right\} }-\omega(k/\kappa_{n})I_{\left\{ k/\kappa_{n}<m-1\right\} }+\omega(1/\kappa_{n})
\]
\[
\leq\left\{ \omega(k/\kappa_{n})+2\omega(m-1)\right\} I_{\left\{ k/\kappa_{n}\geq m-1\right\} }+\omega(k/\kappa_{n})I_{\left\{ k/\kappa_{n}<m-1\right\} }+\omega(1/\kappa_{n})
\]
\[
\leq2\left[\omega(k/\kappa_{n})+\omega(m-1)+\omega(1/\kappa_{n})\right],
\]
as required. $\blacksquare$\bigskip{}

\subsection{Technical Lemmas}

Set $b:=1+\alpha-\alpha m$.

\bigskip{}

\textbf{LM} ($m\in(1/\alpha,1)$). For a given $\alpha\in(1,2]$ set
\[
\eta_{1}:=\frac{1}{2}(\alpha-1)\text{ and }\eta_{2}:=\frac{1}{2}\frac{\alpha m-1}{m}
\]
It can be easily checked that $\eta_{1},\eta_{2}>0$. Moreover, let
\begin{equation}
\epsilon_{-}\in\lbrack0,\eta_{1}\wedge\eta_{2}].\label{epsi}
\end{equation}
Define $\alpha_{\epsilon_{-}}:=\alpha-\epsilon_{-}$ and $b_{\epsilon_{-}}:=1+\alpha_{\epsilon_{-}}-\alpha_{\epsilon_{-}}m$.
Notice that by definition
\begin{equation}
b_{\epsilon_{-}}>1\text{, \ }\alpha_{\epsilon_{-}}>1\text{ \ and\ \ }\alpha_{\epsilon_{-}}m>1.\label{epsi.a}
\end{equation}
To see this first note that
\[
b_{\epsilon_{-}}>1\Leftrightarrow1+\left(\alpha-\epsilon_{-}\right)m>1\Leftrightarrow\left(\alpha-\epsilon_{-}\right)m>0
\]
\[
\Leftrightarrow\alpha>\epsilon_{-},
\]
which is true because $\alpha>\eta_{1}\geq\epsilon_{-}$. For the
latter two results in (\ref{epsi.a}) note that 
\[
\alpha_{\epsilon_{-}}\geq\alpha-\eta_{1}=\alpha-\frac{1}{2}(\alpha-1)=\frac{\alpha}{2}+\frac{1}{2}>1.
\]
In addition, 
\[
\alpha_{\epsilon_{-}}m\geq(\alpha-\eta_{2})m
\]
and
\[
(\alpha-\eta_{2})m>1\Leftrightarrow\alpha m-\eta_{2}m>1
\]
\[
\Leftrightarrow\alpha m-1>\eta_{2}m\Leftrightarrow\frac{\alpha m-1}{m}>\eta_{2},
\]
which is true.

\bigskip{}

\textbf{NM} ($m\in(1,1+1/\alpha)$). Set 
\[
\eta_{3}:=b/2,
\]
and
\[
\eta_{4}:=\frac{1}{2}(2-\alpha),
\]
and note that $\eta_{3}>0$, $\eta_{4}\geq0$. Further, let
\begin{equation}
\epsilon_{+}\in\left\{ \begin{array}{l}
\left[0,\eta_{3}\wedge\eta_{4}\right]\text{, }\alpha\in(1,2)\\{}
[0,\eta_{3}]\text{, }\alpha=2
\end{array}\right.,\label{epsi2}
\end{equation}
and $\alpha_{\epsilon_{+}}=\alpha+\epsilon_{+}$, $b_{\epsilon_{+}}:=1+\alpha_{\epsilon_{+}}-\alpha_{\epsilon_{+}}m$.
Notice that
\begin{equation}
\alpha_{\epsilon_{+}}m>1\text{, \ }b_{\epsilon_{+}}>0\text{ \ and \ }\alpha_{\epsilon_{+}}m-\epsilon_{+}>1.\label{epsi2.a}
\end{equation}
The first fact is obvious. For the second, note that 
\[
b_{\epsilon_{+}}=1+\alpha+\epsilon_{+}-\alpha m-\epsilon_{+}m>0
\]
\[
\Leftrightarrow1+\alpha+\epsilon_{+}\left(1-m\right)>\alpha m\Leftrightarrow\epsilon_{+}\left(1-m\right)>\alpha m-\alpha-1
\]
\[
\Leftrightarrow\epsilon_{+}<\frac{\alpha m-\alpha-1}{1-m}=\frac{b}{m-1},
\]
which conforms with (\ref{epsi2}), because 
\[
\epsilon_{+}\leq\eta_{3}=\frac{b}{2}=\frac{m-1}{2}\frac{b}{m-1}<\frac{1}{2\alpha}\frac{b}{m-1}<\frac{b}{m-1}.
\]

In addition,
\[
\alpha_{\epsilon_{+}}m-\epsilon_{+}>1\Leftrightarrow(\alpha+\epsilon_{+})m-\epsilon_{+}>1
\]
\[
\Leftrightarrow\alpha m+\left(m-1\right)\epsilon_{+}>1,
\]
which is true.

Notice that the intervals of (\ref{epsi}) and (\ref{epsi2}) have
strictly positive lengths. Next, we introduce some further notation.
Set 
\begin{equation}
\left[\epsilon,\alpha_{\epsilon},b_{\epsilon}\right]:=\left\{ \begin{array}{c}
\left[\epsilon_{-},\alpha_{\epsilon_{-}},b_{\epsilon_{-}}\right]\text{, under \textbf{LM}}\\
\left[\epsilon_{+},\alpha_{\epsilon_{+}},b_{\epsilon_{+}}\right]\text{, under \textbf{NM}}
\end{array}\right.\label{epsi3}
\end{equation}

\bigskip{}

For any integer valued sequences $a_{n}\leq b_{n}$, set 
\[
V_{\epsilon}\left(a_{n},b_{n}\right):=\kappa^{-b_{\epsilon}}_{n}\sum^{b_{n}}_{k=a_{n}}\left\vert \sum^{k}_{l=0}\rho^{l}_{n}c_{k-l}\right\vert ^{\alpha_{\epsilon}},
\]
and
\[
Q_{\epsilon}\left(a_{n},b_{n}\right):=\kappa^{-b_{\epsilon}}_{n}\sum^{\infty}_{k=b_{n}}\left\vert \sum^{a_{n}-1}_{l=0}\rho^{l}_{n}c_{k-l}\right\vert ^{\alpha_{\epsilon}}.
\]

Moreover, for some $\theta>0$ define
\[
b_{\theta}:=1+\theta-\theta m,
\]
\[
\tilde{V}_{\theta}\left(a_{n},b_{n}\right):=\kappa^{-b_{\theta}}_{n}\sum^{b_{n}}_{k=a_{n}}\left\vert \sum^{k}_{l=0}\rho^{l}_{n}c_{k-l}\right\vert ^{\theta},
\]
and
\[
\tilde{Q}_{\theta}\left(a_{n},b_{n}\right):=\kappa^{-b_{\theta}}_{n}\sum^{\infty}_{k=b_{n}}\left\vert \sum^{a_{n}-1}_{l=0}\rho^{l}_{n}c_{k-l}\right\vert ^{\theta}.
\]

\begin{lemma} \textbf{\label{LemmaMI0}} \textit{Set }$\delta_{n}:=\lfloor\varsigma_{n}\kappa_{n}\rfloor$,
\textit{with }$\varsigma_{n}:=\left[\ln(1+\kappa_{n})\right]^{2}$\textit{,
and }$t_{n}=n\vee\delta_{n}$\textit{. Suppose that }$c_{k}\sim m_{0}\cdot k^{-m}$,
$m_{0}\neq0$ such that
\begin{itemize}
\item[\textbf{(i)}] $\theta\in(0,3]$;
\item[\textbf{(ii)}] $\theta m>1.$ 
\end{itemize}
\noindent\textit{Then as} $n\rightarrow\infty$
\begin{itemize}
\item[\textbf{(a)}] $\tilde{V}_{\theta}\left(\delta_{n},t_{n}\right)=o(1)$\textit{;}
\item[\textbf{(b)}] $\tilde{Q}_{\theta}\left(t_{n},t_{n}\right)=\tilde{Q}_{\theta}\left(t_{n},\delta_{n}\right)+o(1)$\textit{;}
\item[\textbf{(c)}] $\tilde{Q}_{\theta}\left(t_{n},\delta_{n}\right)=o(1).$ 
\end{itemize}
\end{lemma}

\bigskip{}
\textbf{Remark}. The proof provided for \textbf{Lemma \ref{LemmaMI0}}
in fact shows that
\[
\kappa^{-b_{\theta}}_{n}\sum^{\infty}_{k=t_{n}}\left(\sum^{t_{n}-1}_{l=0}\left|\rho^{l}_{n}c_{k-l}\right|\right)^{\theta}=\kappa^{-b_{\theta}}_{n}\sum^{\infty}_{k=t_{n}}\left(\sum^{\delta_{n}-1}_{l=0}\left|\rho^{l}_{n}c_{k-l}\right|\right)^{\theta}=o(1),
\]
which implies that $\tilde{Q}_{\theta}\left(t_{n},t_{n}\right)=o(1).$

\bigskip{}

The following result is a direct consequence of \textbf{Lemma \ref{LemmaMI0}}.\ 

\bigskip{}

\begin{lemma} \textbf{\label{LemmaMI1}} \textit{Set }$\delta_{n}:=\lfloor\varsigma_{n}\kappa_{n}\rfloor$,
\textit{with }$\varsigma_{n}:=\left[\ln(1+\kappa_{n})\right]^{2}$\textit{,
and }$t_{n}=n\vee\delta_{n}$\textit{. Suppose that} 
\[
c_{k}\sim m_{0}\cdot k^{-m},\text{ }m_{0}\neq0,
\]
\textit{and either condition }\textbf{\textit{MI-NM}}\textit{ or }\textbf{\textit{MI-LM}}\textit{
holds.}

\noindent\textit{Then as} $n\rightarrow\infty$
\begin{itemize}
\item[\textbf{(a)}] $V_{\epsilon}\left(\delta_{n},t_{n}\right)=o(1)$\textit{;}
\item[\textbf{(b)}] $Q_{\epsilon}\left(t_{n},t_{n}\right)=Q_{\epsilon}\left(t_{n},\delta_{n}\right)+o(1)$\textit{;}
\item[\textbf{(c)}] $Q_{\epsilon}\left(t_{n},\delta_{n}\right)=o(1).$ 
\end{itemize}
\end{lemma}

\bigskip{}

Consider the ``\textit{reversed''} Gamma function 
\[
\gamma\left(m;a,x\right):=\int^{x}_{y=a}e^{y}y^{1-m}dy\text{, }m\in\mathbb{R}\text{, }0\leq a<x.
\]
Moreover, set
\[
\gamma\left(m;x\right):=\gamma\left(m;0,x\right).
\]

\begin{remark} Notice that for $m>2$, we need $a>0$ to ensure that
$\gamma\left(m;a,x\right)<\infty$. \end{remark}

\bigskip{}

\begin{lemma} \textbf{\label{LemmaMI2}} The following hold.
\begin{itemize}
\item[\textbf{(i)}] \textit{Let} $\theta_{1}>0$\textit{, and }$\theta_{2}\in(0,1)$\textit{
such that }$\theta_{1}\theta_{2}>1$\textit{. Then}
\begin{equation}
\int^{\infty}_{x=0}e^{-\theta_{1}x}\left(\int^{x}_{y=0}e^{y}y^{-\theta_{2}}dy\right)^{\theta_{1}}dx<\infty.\label{MI2a}
\end{equation}
\item[\textbf{(ii)}] \textit{Let }$\theta_{1}>0$\textit{ and }$\theta_{2}\in(1,2)$\textit{
such that }$1/\theta_{2}<\theta_{1}<1/\left(\theta_{2}-1\right)$\textit{.
Then}
\begin{equation}
\int^{\infty}_{x=0}e^{-\theta_{1}x}\left(\int^{x}_{y=1}e^{y}y^{-\theta_{2}}dy\right)^{\theta_{1}}dx<\infty\text{,}\label{MI2b}
\end{equation}
\textit{and}
\begin{equation}
\int^{\infty}_{x=0}\left\vert x^{1-\theta_{2}}-e^{-x}\gamma\left(\theta_{2};x\right)\right\vert ^{\theta_{1}}dx<\infty.\label{MI2c}
\end{equation}
\end{itemize}
\end{lemma}

\bigskip{}

\begin{corollary} \textbf{\label{CorollaryMI2} }Let\textit{\ }$\alpha_{\epsilon}$\textit{\
as in (\ref{epsi3}). The following hold.}
\begin{itemize}
\item[\textbf{(i)}] Under condition \textbf{MI-LM}
\[
\mathbb{V}_{1}(\alpha_{\epsilon},m):=\int^{\infty}_{x=0}e^{-\alpha_{\epsilon}x}\left(\int^{x}_{y=0}e^{y}y^{-m}dy\right)^{\alpha_{\epsilon}}dx<\infty.
\]
\item[\textbf{(ii)}] Under condition \textbf{MI-NM}
\[
\mathbb{V}^{\prime}_{1}(\alpha_{\epsilon},m):=\int^{\infty}_{x=0}e^{-\alpha_{\epsilon}x}\left(\int^{x}_{y=1}e^{y}y^{-m}dy\right)^{\alpha_{\epsilon}}dx<\infty\text{,}
\]
and
\[
\mathbb{V}_{2}(\alpha_{\epsilon},m):=\left\vert 1-m\right\vert ^{-\alpha_{\epsilon}}\int^{\infty}_{x=0}\left\vert x^{1-m}-e^{-x}\gamma\left(m;x\right)\right\vert ^{\alpha_{\epsilon}}dx<\infty.
\]
\end{itemize}
\end{corollary}

\bigskip{}

\noindent\textbf{Proof of Corollary \ref{CorollaryMI2}.} It follows
directly from \textbf{Lemma \ref{LemmaMI2}} with $\alpha_{\epsilon}=\theta_{1}$,
and $m=\theta_{2}$. $\blacksquare$\bigskip{}

\begin{remark} \label{RemarkMI1}

\textbf{(a)} Notice that under condition \textbf{MI-LM} $\sum^{\infty}_{l=0}\left\vert c_{l}\right\vert ^{\alpha_{\epsilon}}<\infty$,
whilst under condition \textbf{MI-NM} $\sum^{\infty}_{l=0}\left(\left\vert c_{l}\right\vert +\left\vert c_{l}\right\vert ^{\alpha_{\epsilon}}\right)<\infty$.

\textbf{(b) }Due to the approximation $c_{j}\sim m_{0}\cdot j^{-m}$,
there is no loss of generality if we assume that $\left\vert c_{j}\right\vert $
is decreasing.

\textbf{(c)} Recall that $b:=1+\alpha-\alpha m$. By definition, $b>0$.
In particular, 
\[
b\in\left\{ \begin{array}{l}
(0,1)\text{, under condition MI-NM}\\
(1,\alpha)\text{, under condition MI-LM}
\end{array}\right.
\]

\textbf{(d)} Under condition \textbf{MI-NM} or \textbf{MI-LM} we have
$\alpha m>1$.

\textbf{(e)} For $\alpha=2$, it can be shown (cf. DK2018) that 
\[
\mathbb{V}_{1}(\alpha)=B\left(1-m,2m-1\right)\Gamma(2-2m),
\]
with $B$ and $\Gamma$\ being beta and gamma functions respectively.
\end{remark}

\bigskip{}

\begin{lemma} \textbf{\label{LemmaMI3} }Suppose that the assumptions
of \textbf{Lemma \ref{LemmaMI1}} hold with $m\in\left(1/\alpha,1\right)$,
$\alpha>1$. 

\textbf{(a)} As $n\rightarrow\infty$,

\[
V_{\epsilon}\left(0,\delta_{n}\right)=\left\vert m_{0}\right\vert ^{\alpha_{\epsilon}}\mathbb{V}_{1}(\alpha_{\epsilon},m)+o(1).
\]

\textbf{(b)} Let $(\delta^{\prime}_{n})^{-1}+\delta^{\prime}_{n}\kappa^{-1}_{n}\rightarrow0$.
As $n\rightarrow\infty$,
\[
V_{\epsilon}\left(0,\delta^{\prime}_{n}\right)=o(1).
\]
\end{lemma}

\bigskip{}

\begin{lemma} \textbf{\label{LemmaMI4} }Suppose that the assumptions
of \textbf{Lemma \ref{LemmaMI1}} hold with $m\in\left(1,1+1/\alpha\right)$,
$\alpha>1$. Then as $n\rightarrow\infty$
\[
V_{\epsilon}\left(0,\delta_{n}\right)=\left\vert m_{0}\right\vert ^{\alpha_{\epsilon}}\mathbb{V}_{2}(\alpha_{\epsilon},m)+o(1).
\]
\end{lemma}

\bigskip{}

\begin{proposition} \label{SterlingApprox} The following hold

\textbf{(a)} Let
\[
q(x):=e^{-x}\int^{x}_{1}e^{y}y^{-q_{0}}dy,\:q_{0}>0.
\]
For all integer $N\geq1$, as $x\rightarrow\infty$
\[
q(x)=x^{-q_{0}}\left[\sum^{N-1}_{j=0}\frac{\left(q_{0}\right)_{j}}{x^{j}}+O\left(x^{-N}\right)\right],
\]
where $\left(.\right)_{j}$ Pochhammer\textquoteright s symbol for
the forward factorial function, i.e. $\left(q_{0}\right)_{j}=q_{0}\left(q_{0}+1\right)...\left(q_{0}+j-1\right).$ 

\textbf{(b)} Let
\[
q(x):=e^{-x}\int^{x}_{0}e^{y}y^{-q_{0}}dy,\:0<q_{0}<1.
\]
For all integer $N\geq1$, as $x\rightarrow0$
\[
q(x)=\frac{x^{1-q_{0}}}{1-q_{0}}\left[\sum^{N-1}_{j=0}\frac{\left(-x\right)^{j}}{\left(2-q_{0}\right)_{j}}+O\left(x^{-N}\right)\right],
\]

\bigskip{}

\textbf{Remark}. \textbf{(a)} It follows from \textbf{Proposition
\ref{SterlingApprox}(a)} that as $x\rightarrow\infty$

\[
q(x)\sim x^{-q_{0}}.
\]
In addition,
\[
\frac{q(x)}{x^{-q_{0}}}-1\sim\frac{q_{0}}{x}.
\]

\textbf{(b)} It follows from \textbf{Proposition \ref{SterlingApprox}(b)}
that as that as $x\rightarrow0$

\[
q(x)\sim\frac{x^{1-q_{0}}}{1-q_{0}}.
\]
In addition,
\[
\frac{q(x)}{x^{1-q_{0}}/(1-q_{0})}-1\sim-\frac{x^{2-q_{0}}}{2-q_{0}}.
\]
\end{proposition}

\bigskip{}

For the following proof recall that

\[
V_{\theta}\left(\delta_{n},t_{n}\right)=\kappa^{-b_{\theta}}_{n}\sum^{t_{n}}_{k=\delta_{n}}\left\vert \sum^{k}_{l=0}\rho^{l}_{n}c_{k-l}\right\vert ^{\theta}
\]
with $\theta\in(0,3]$, $b_{\theta}:=1+\theta-\theta m$ and $c_{k}\sim m_{0}k^{-m}$,
$m_{0}\neq0$ such that $\theta m>1$.\bigskip{}

\textbf{Proof of Lemma \ref{LemmaMI0}}. Without loss of generality,
we assume that $m_{0}=1$.

\textbf{(a)} We provide a proof of the case $\theta\in(2,3]$. The
proof for $\theta\in(0,2]$ is similar but simpler. Set $\bar{\delta}_{n}=\lfloor\delta_{n}/2\rfloor$
and recall that $\delta_{n}=\kappa_{n}\left[\ln(1+\kappa_{n})\right]^{2}$.
Hence, in view of this and the fact that $\kappa_{n}\rightarrow\infty$,
\begin{equation}
\lim_{n\rightarrow\infty}\ln\left(\kappa^{q}_{n}\right)\kappa_{n}\bar{\delta}^{-1}_{n}=0,\text{ for all }q\in\mathbb{R}.\label{VarMI-0}
\end{equation}
Note that (\ref{VarMI-0}) together with Lemma 6.4 in DK2018 imply
that
\begin{equation}
\lim_{n\rightarrow\infty}\kappa^{q}_{n}\rho^{\vartheta\bar{\delta}_{n}}_{n}=0,\text{ for all }\vartheta>0.\label{VarMI-01}
\end{equation}
In addition, under our assumptions $\sum^{\infty}_{k=0}\left\vert c_{k}\right\vert ^{\theta}<\infty$.
In view of these facts write
\[
V_{\theta}\left(\delta_{n},t_{n}\right)=\kappa^{-b_{\theta}}_{n}\sum^{t_{n}}_{k=\delta_{n}}\left\vert \sum^{k}_{l=0}\rho^{l}_{n}c_{k-l}\right\vert ^{\theta}\leq\kappa^{-b_{\theta}}_{n}\sum^{t_{n}}_{k=\delta_{n}}\left(\sum^{k}_{l=0}\left\vert \rho^{l}_{n}c_{k-l}\right\vert \right)^{\theta}
\]
\[
=\kappa^{-b_{\theta}}_{n}\sum^{t_{n}}_{k=\delta_{n}}\left(\sum^{k}_{l_{1}=0}\left\vert \rho^{l_{1}}_{n}c_{k-l_{1}}\right\vert \right)^{\theta-2}\left(\sum^{k}_{l_{2}=0}\left\vert \rho^{l_{2}}_{n}c_{k-l_{2}}\right\vert \right)\left(\sum^{k}_{l_{3}=0}\left\vert \rho^{l_{3}}_{n}c_{k-l_{3}}\right\vert \right)
\]
\[
\leq\kappa^{-b_{\theta}}_{n}\sum^{t_{n}}_{k=\delta_{n}}\left[\left(\sum^{\bar{\delta}_{n}}_{l_{1}=0}\left\vert \rho^{l_{1}}_{n}c_{k-l_{1}}\right\vert \right)^{\theta-2}+\left(\sum^{k}_{l_{1}=\bar{\delta}_{n}}\left\vert \rho^{l_{1}}_{n}c_{k-l_{1}}\right\vert \right)^{\theta-2}\right]\left(\sum^{k}_{l_{2}=0}\left\vert \rho^{l_{2}}_{n}c_{k-l_{2}}\right\vert \right)\left(\sum^{k}_{l_{3}=0}\left\vert \rho^{l_{3}}_{n}c_{k-l_{3}}\right\vert \right)
\]
\begin{equation}
=:A_{n}+B_{n}.\label{VarMI-1}
\end{equation}
Noting that $c_{k}\sim m_{0}k^{-m}$ where the exponent $-m$ is strictly
negative, there is no loss of generality if we assume that $c_{k}$
is decreasing. In view of this convention, the first term is
\[
A_{n}\leq\kappa^{-b_{\theta}}_{n}\sum^{t_{n}}_{k=\delta_{n}}\left(\sum^{\bar{\delta}_{n}}_{l_{1}=0}\left\vert \rho^{l_{1}}_{n}c_{k-\bar{\delta}_{n}}\right\vert \right)^{\theta-2}\left(\sum^{k}_{l_{2}=0}\left\vert \rho^{l_{2}}_{n}c_{k-l_{2}}\right\vert \right)\left(\sum^{k}_{l_{3}=0}\left\vert \rho^{l_{3}}_{n}c_{k-l_{3}}\right\vert \right)
\]
\[
=\kappa^{-b_{\theta}}_{n}\left(\sum^{\bar{\delta}_{n}}_{l_{1}=0}\rho^{l_{1}}_{n}\right)^{\theta-2}\sum^{t_{n}}_{k=\delta_{n}}\left\vert c_{k-\bar{\delta}_{n}}\right\vert ^{\theta-2}\left(\sum^{k}_{l_{2}=0}\left\vert \rho^{l_{2}}_{n}c_{k-l_{2}}\right\vert \right)\left(\sum^{k}_{l_{3}=0}\left\vert \rho^{l_{3}}_{n}c_{k-l_{3}}\right\vert \right)
\]
\[
=\kappa^{-b_{\theta}}_{n}\left(\sum^{\bar{\delta}_{n}}_{l_{1}=0}\rho^{l_{1}}_{n}\right)^{\theta-2}\sum^{t_{n}}_{l_{2}=0}\left\vert \rho^{l_{2}}_{n}\right\vert \sum^{t_{n}}_{l_{3}=0}\left\vert \rho^{l_{3}}_{n}\right\vert \sum^{t_{n}}_{k=\delta_{n}\vee l_{2}\vee l_{3}}\left\vert c_{k-\bar{\delta}_{n}}\right\vert ^{\theta-2}\left\vert c_{k-l_{2}}\right\vert \left\vert c_{k-l_{3}}\right\vert 
\]
\[
=\kappa^{-b_{\theta}}_{n}\left(\sum^{\bar{\delta}_{n}}_{l_{1}=0}\rho^{l_{1}}_{n}\right)^{\theta-2}\prod^{3}_{i=2}\left(\sum^{\bar{\delta}_{n}}_{l_{i}=0}\left\vert \rho^{l_{i}}_{n}\right\vert +\sum^{t_{n}}_{l_{i}=\bar{\delta}_{n}+1}\left\vert \rho^{l_{i}}_{n}\right\vert \right)\sum^{t_{n}}_{k=\delta_{n}\vee l_{2}\vee l_{3}}\left\vert c_{k-\bar{\delta}_{n}}\right\vert ^{\theta-2}\left\vert c_{k-l_{2}}\right\vert \left\vert c_{k-l_{3}}\right\vert 
\]
\[
\leq_{(1)}C\kappa^{\theta-2-b_{\theta}}_{n}\cdot\prod^{3}_{i=2}\left(\sum^{\bar{\delta}_{n}}_{l_{i}=0}\left\vert \rho^{l_{i}}_{n}\right\vert +\sum^{t_{n}}_{l_{i}=\bar{\delta}_{n}+1}\left\vert \rho^{l_{i}}_{n}\right\vert \right)\sum^{t_{n}}_{k=\delta_{n}\vee l_{2}\vee l_{3}}\left\vert c_{k-\bar{\delta}_{n}}\right\vert ^{\theta-2}\left\vert c_{k-l_{2}}\right\vert \left\vert c_{k-l_{3}}\right\vert 
\]
\begin{equation}
=:C\kappa^{\theta-2-b_{\theta}}_{n}\cdot A^{\prime}_{n},\label{VarMI-2}
\end{equation}
where $\leq_{(1)}$ follows from Lemma 6.4 of DK2018. Consider the
term
\[
A^{\prime}_{n}=\prod^{3}_{i=2}\left(\sum^{\bar{\delta}_{n}}_{l_{i}=0}\left\vert \rho^{l_{i}}_{n}\right\vert +\sum^{t_{n}}_{l_{i}=\bar{\delta}_{n}+1}\left\vert \rho^{l_{i}}_{n}\right\vert \right)\sum^{t_{n}}_{k=\delta_{n}\vee l_{2}\vee l_{3}}\left\vert c_{k-\bar{\delta}_{n}}\right\vert ^{\theta-2}\left\vert c_{k-l_{2}}\right\vert \left\vert c_{k-l_{3}}\right\vert 
\]
\[
=\sum^{\bar{\delta}_{n}}_{l_{2}=0}\left\vert \rho^{l_{2}}_{n}\right\vert \sum^{\bar{\delta}_{n}}_{l_{3}=0}\left\vert \rho^{l_{3}}_{n}\right\vert \sum^{t_{n}}_{k=\delta_{n}\vee l_{2}\vee l_{3}}\left\vert c_{k-\bar{\delta}_{n}}\right\vert ^{\theta-2}\left\vert c_{k-l_{2}}\right\vert \left\vert c_{k-l_{3}}\right\vert 
\]
\[
+\sum^{t_{n}}_{l_{2}=\bar{\delta}_{n}+1}\left\vert \rho^{l_{2}}_{n}\right\vert \sum^{t_{n}}_{l_{3}=\bar{\delta}_{n}+1}\left\vert \rho^{l_{3}}_{n}\right\vert \sum^{t_{n}}_{k=\delta_{n}\vee l_{2}\vee l_{3}}\left\vert c_{k-\bar{\delta}_{n}}\right\vert ^{\theta-2}\left\vert c_{k-l_{2}}\right\vert \left\vert c_{k-l_{3}}\right\vert 
\]
\[
+2\sum^{\bar{\delta}_{n}}_{l_{2}=0}\left\vert \rho^{l_{2}}_{n}\right\vert \sum^{t_{n}}_{l_{3}=\bar{\delta}_{n}+1}\left\vert \rho^{l_{3}}_{n}\right\vert \sum^{t_{n}}_{k=\delta_{n}\vee l_{2}\vee l_{3}}\left\vert c_{k-\bar{\delta}_{n}}\right\vert ^{\theta-2}\left\vert c_{k-l_{2}}\right\vert \left\vert c_{k-l_{3}}\right\vert 
\]
\[
\leq\sum^{\bar{\delta}_{n}}_{l_{2}=0}\left\vert \rho^{l_{2}}_{n}\right\vert \sum^{\bar{\delta}_{n}}_{l_{3}=0}\left\vert \rho^{l_{3}}_{n}\right\vert \sum^{t_{n}}_{k=\delta_{n}\vee l_{2}\vee l_{3}}\left\vert c_{k-\bar{\delta}_{n}}\right\vert ^{\theta-2}\left\vert c_{k-\bar{\delta}_{n}}\right\vert \left\vert c_{k-\bar{\delta}_{n}}\right\vert 
\]
\[
+\sum^{t_{n}}_{l_{2}=\bar{\delta}_{n}+1}\left\vert \rho^{l_{2}}_{n}\right\vert \sum^{t_{n}}_{l_{3}=\bar{\delta}_{n}+1}\left\vert \rho^{l_{3}}_{n}\right\vert \sum^{t_{n}}_{k=\delta_{n}\vee l_{2}\vee l_{3}}\left\vert c_{k-\left(\delta_{n}\vee l_{2}\vee l_{3}\right)}\right\vert ^{\theta-2}\left\vert c_{k-\left(\delta_{n}\vee l_{2}\vee l_{3}\right)}\right\vert \left\vert c_{k-\left(\delta_{n}\vee l_{2}\vee l_{3}\right)}\right\vert 
\]
\[
+2\sum^{\bar{\delta}_{n}}_{l_{2}=0}\left\vert \rho^{l_{2}}_{n}\right\vert \sum^{t_{n}}_{l_{3}=\bar{\delta}_{n}+1}\left\vert \rho^{l_{3}}_{n}\right\vert \sum^{t_{n}}_{k=\delta_{n}\vee l_{2}\vee l_{3}}\left\vert c_{k-\left(\delta_{n}\vee l_{2}\vee l_{3}\right)}\right\vert ^{\theta-2}\left\vert c_{k-\left(\delta_{n}\vee l_{2}\vee l_{3}\right)}\right\vert \left\vert c_{k-\left(\delta_{n}\vee l_{2}\vee l_{3}\right)}\right\vert 
\]
\[
\leq\sum^{\bar{\delta}_{n}}_{l_{2}=0}\left\vert \rho^{l_{2}}_{n}\right\vert \sum^{\bar{\delta}_{n}}_{l_{3}=0}\left\vert \rho^{l_{3}}_{n}\right\vert \sum^{\infty}_{k=\bar{\delta}_{n}}\left\vert c_{k}\right\vert ^{\theta}+\sum^{t_{n}}_{l_{2}=\bar{\delta}_{n}+1}\left\vert \rho^{l_{2}}_{n}\right\vert \sum^{t_{n}}_{l_{3}=\bar{\delta}_{n}+1}\left\vert \rho^{l_{3}}_{n}\right\vert \sum^{\infty}_{k=0}\left\vert c_{k}\right\vert ^{\theta}
\]
\[
+2\sum^{\bar{\delta}_{n}}_{l_{2}=0}\left\vert \rho^{l_{2}}_{n}\right\vert \sum^{t_{n}}_{l_{3}=\bar{\delta}_{n}+1}\left\vert \rho^{l_{3}}_{n}\right\vert \sum^{\infty}_{k=0}\left\vert c_{k}\right\vert ^{\theta}
\]
\[
\leq_{(2)}C\left[\kappa^{2}_{n}\bar{\delta}^{1-m\theta}_{n}\frac{1}{\bar{\delta}_{n}}\sum^{\infty}_{k=\bar{\delta}_{n}}\left(\frac{k}{\bar{\delta}_{n}}\right)^{-m\theta}+\kappa^{2}_{n}\rho^{2\bar{\delta}_{n}}_{n}\sum^{\infty}_{k=0}\left\vert c_{k}\right\vert ^{\theta}+\kappa^{2}_{n}\rho^{\bar{\delta}_{n}}_{n}\sum^{\infty}_{k=0}\left\vert c_{k}\right\vert ^{\theta}\right]
\]
\[
=_{(3)}C\left[\kappa^{2}_{n}\bar{\delta}^{1-m\theta}_{n}\int^{\infty}_{1}x^{-m\theta}dx+\kappa^{2}_{n}\rho^{2\bar{\delta}_{n}}_{n}+\kappa^{2}_{n}\rho^{\bar{\delta}_{n}}_{n}\right],
\]
where $\leq_{(2)}$ follows from DK2018 (Lemma 6.4), and $=_{(3)}$
from Euler summation (see e.g. \ref{LemmaSV1}). Combining (\ref{VarMI-2})-(\ref{VarMI-3})
we get 
\[
A_{n}\leq C\kappa^{\theta-2-b_{\theta}}_{n}\left[\kappa^{2}_{n}\bar{\delta}^{1-m\theta}_{n}+\kappa^{2}_{n}\rho^{2\bar{\delta}_{n}}_{n}+\kappa^{2}_{n}\rho^{\bar{\delta}_{n}}_{n}\right]
\]
\begin{equation}
=C\left[\kappa^{\theta-b_{\theta}}_{n}\bar{\delta}^{1-m\theta}_{n}+\kappa^{\theta-b_{\theta}}_{n}\rho^{2\bar{\delta}_{n}}_{n}+\kappa^{\theta-b_{\theta}}_{n}\rho^{\bar{\delta}_{n}}_{n}\right]=_{(4)}C\frac{\kappa^{\theta}_{n}}{\kappa^{\theta}_{n}}\frac{\bar{\delta}^{1-m\theta}_{n}}{\kappa^{1-m\theta}_{n}}+o(1)=_{(5)}o(1),\label{VarMI-3}
\end{equation}
where $=_{(4)}$ is due to (\ref{VarMI-01}), and $=_{(5)}$ due to
the fact that $\kappa_{n}=o\left(\bar{\delta}_{n}\right)$.

Next, using (\ref{VarMI-01}) once more and DK2018 (Lemma 6.4) the
term $B_{n}$ 
\[
B_{n}\leq\kappa^{-b_{\theta}}_{n}\sum^{t_{n}}_{k=\delta_{n}}\sum^{k}_{l_{1}=\bar{\delta}_{n}}\left\vert \rho^{l_{1}}_{n}c_{k-l_{1}}\right\vert ^{\theta-2}\sum^{k}_{l_{2}=0}\left\vert \rho^{l_{2}}_{n}c_{k-l_{2}}\right\vert \sum^{k}_{l_{3}=0}\left\vert \rho^{l_{3}}_{n}c_{k-l_{3}}\right\vert 
\]
\[
=\kappa^{-b_{\theta}}_{n}\sum^{t_{n}}_{l_{1}=\bar{\delta}_{n}}\rho^{\left(\theta-2\right)l_{1}}_{n}\sum^{t_{n}}_{l_{2}=0}\rho^{l_{2}}_{n}\sum^{t_{n}}_{l_{3}=0}\rho^{l_{3}}_{n}\sum^{t_{n}}_{k=\delta_{n}\vee l_{1}\vee l_{2}\vee l_{3}}\left\vert c_{k-l_{1}}\right\vert ^{\theta-2}\left\vert c_{k-l_{2}}\right\vert \left\vert c_{k-l_{3}}\right\vert 
\]
\[
\leq\kappa^{-b_{\theta}}_{n}\sum^{t_{n}}_{l_{1}=\bar{\delta}_{n}}\rho^{\left(\theta-2\right)l_{1}}_{n}\sum^{t_{n}}_{l_{2}=0}\rho^{l_{2}}_{n}\sum^{t_{n}}_{l_{3}=0}\rho^{l_{3}}_{n}\sum^{t_{n}}_{k=\delta_{n}\vee l_{1}\vee l_{2}\vee l_{3}}\left\vert c_{k-\left(\delta_{n}\vee l_{1}\vee l_{2}\vee l_{3}\right)}\right\vert ^{\theta}
\]
\begin{equation}
\leq\kappa^{-b_{\theta}}_{n}\sum^{t_{n}}_{l_{1}=\bar{\delta}_{n}}\rho^{\left(\theta-2\right)l_{1}}_{n}\sum^{t_{n}}_{l_{2}=0}\rho^{l_{2}}_{n}\sum^{t_{n}}_{l_{3}=0}\rho^{l_{3}}_{n}\sum^{\infty}_{k=0}\left\vert c_{k}\right\vert ^{\theta}\leq C\kappa^{3-b_{\theta}}_{n}\rho^{\left(\theta-2\right)\bar{\delta}_{n}}_{n}\sum^{\infty}_{k=0}\left\vert c_{k}\right\vert ^{\theta}=o(1).\label{VarMI-4}
\end{equation}
The requisite result follows from (\ref{VarMI-1}), (\ref{VarMI-3})
and (\ref{VarMI-4}).

\textbf{(b)} Write 
\[
\left\vert Q_{\theta}\left(t_{n},t_{n}\right)-Q_{\theta}\left(t_{n},\delta_{n}\right)\right\vert 
\]
\[
\leq\kappa^{-b_{\theta}}_{n}\sum^{\infty}_{k=t_{n}}\left\vert \left\vert \sum^{t_{n}-1}_{l=0}\rho^{l}_{n}c_{k-l}\right\vert ^{\theta}-\left\vert \sum^{\delta_{n}-1}_{l=0}\rho^{l}_{n}c_{k-l}\right\vert ^{\theta}\right\vert 
\]
\[
\leq C\kappa^{-b_{\theta}}_{n}\sum^{\infty}_{k=t_{n}}\left[\left\vert \sum^{t_{n}-1}_{l=\delta_{n}}\rho^{l}_{n}c_{k-l}\right\vert ^{\theta}+\left\vert \sum^{\delta_{n}-1}_{l=0}\rho^{l}_{n}c_{k-l}\right\vert ^{\theta-1}\left\vert \sum^{t_{n}-1}_{l=\delta_{n}}\rho^{l}_{n}c_{k-l}\right\vert \right]=:C\left(\Gamma_{1n}+\Gamma_{2n}\right),
\]
where the last inequality follows from (\ref{MVT}). As before, without
loss of generality we assume that $\left\vert c_{k}\right\vert $
is decreasing. In view of (\ref{VarMI-01}), consider 
\[
\Gamma_{1n}=\kappa^{-b_{\theta}}_{n}\sum^{\infty}_{k=t_{n}}\left\vert \sum^{t_{n}-1}_{l=\delta_{n}}\rho^{l}_{n}c_{k-l}\right\vert ^{\theta}\leq\kappa^{-b_{\theta}}_{n}\sum^{\infty}_{k=t_{n}}\left(\sum^{t_{n}-1}_{l=\delta_{n}}\left\vert \rho^{l}_{n}c_{k-t_{n}}\right\vert \right)^{\theta}
\]
\[
=\kappa^{-b_{\theta}}_{n}\sum^{\infty}_{k=t_{n}}\left\vert c_{k-t_{n}}\right\vert ^{\theta}\left(\sum^{t_{n}-1}_{l=\delta_{n}}\left\vert \rho^{l}_{n}\right\vert \right)^{\theta}\leq C\kappa^{\theta}_{n}\rho^{\theta\delta_{n}}_{n}\kappa^{-b_{\theta}}_{n}\sum^{\infty}_{k=0}\left\vert c_{k}\right\vert ^{\theta}
\]
\[
=C\rho^{\theta\delta_{n}}_{n}\kappa^{-\theta m+1}_{n}=o(\kappa^{-\theta m+1}_{n})=o(1).
\]
Next,
\[
\Gamma_{2n}\leq\kappa^{-b_{\theta}}_{n}\sum^{\infty}_{k=t_{n}}\left(\sum^{\delta_{n}-1}_{l=0}\left\vert \rho^{l}_{n}c_{k-l}\right\vert \right)^{\theta-1}\sum^{t_{n}-1}_{l=\delta_{n}}\left\vert \rho^{l}_{n}c_{k-l}\right\vert 
\]
\[
\leq\kappa^{-b_{\theta}}_{n}\sum^{\infty}_{k=t_{n}}\left(\sum^{\delta_{n}-1}_{l=0}\left\vert \rho^{l}_{n}c_{k-t_{n}}\right\vert \right)^{\theta-1}\sum^{t_{n}-1}_{l=\delta_{n}}\left\vert \rho^{l}_{n}c_{k-t_{n}}\right\vert \leq\kappa^{-b_{\theta}}_{n}\left(\sum^{\delta_{n}-1}_{l=0}\rho^{l}_{n}\right)^{\theta-1}\sum^{t_{n}-1}_{l=\delta_{n}}\rho^{l}_{n}\sum^{\infty}_{k=0}\left\vert c_{k}\right\vert ^{\theta}
\]
\[
\leq C\kappa^{-b_{\theta}}_{n}\left(\sum^{\delta_{n}-1}_{l=0}\rho^{l}_{n}\right)^{\theta-1}\sum^{t_{n}-1}_{l=\delta_{n}}\rho^{l}_{n}=C\kappa^{\theta-b_{\theta}}_{n}\rho^{\delta_{n}}_{n}=o\left(1\right).
\]
as required.

\textbf{(c)} As $n\rightarrow\infty$
\[
Q_{\theta}(t_{n},\delta_{n})=\kappa^{-b_{\theta}}_{n}\sum^{\infty}_{k=t_{n}}\left\vert \sum^{\delta_{n}-1}_{l=0}\rho^{l}_{n}c_{k-l}\right\vert ^{\theta}=\kappa^{-b_{\theta}}_{n}\sum^{\infty}_{k=t_{n}}\left\vert \sum^{k}_{l=k-\left(\delta_{n}-1\right)}\rho^{k-l}_{n}c_{l}\right\vert ^{\theta}
\]
\[
=_{(1)}\frac{1}{\kappa_{n}}\sum^{\infty}_{k=t_{n}}\left\vert \frac{1}{\kappa_{n}}\sum^{k}_{l=k-\left(\delta_{n}-1\right)}e^{-\frac{k-l}{\kappa_{n}}}\left(\frac{l}{\kappa_{n}}\right)^{-m}\right\vert ^{\theta}\left[1+o(1)\right]
\]
\[
=_{(2)}\int^{\infty}_{x=t_{n}/\kappa_{n}}\left\vert \int^{x}_{y=x-\left(\delta_{n}-1\right)/\kappa_{n}}e^{-x+y}y^{-m}dy\right\vert ^{\theta}dx\left[1+o(1)\right]
\]
\[
\leq\int^{\infty}_{x=t_{n}/\kappa_{n}}\left(\int^{x}_{y=(t_{n}-\delta_{n})/\kappa_{n}}e^{-x+y}y^{-m}dy\right)^{\theta}dx\left[1+o(1)\right]
\]
\[
\leq_{(3)}\int^{\infty}_{x=t_{n}/\kappa_{n}}\left(\int^{x}_{y=1}e^{-x+y}y^{-m}dy\right)^{\theta}dx\left[1+o(1)\right]=_{(4)}o(1),
\]
where $=_{(1)}$ follows from Lemma 6.4 in DK2018 and the fact that
$c_{l}\sim l^{-m}$, $=_{(2)}$ from Euler summation, and $\leq_{(3)}$
from the fact that $(t_{n}-\delta_{n})/\kappa_{n}\rightarrow\infty$.
Finally, $=_{(4)}$ follows from Lemma \textbf{\ref{LemmaMI2}. }$\blacksquare$

\bigskip{}

\textbf{Proof of Lemma \ref{LemmaMI1}}. \textbf{Lemma \ref{LemmaMI1}}
follows directly from \textbf{Lemma \ref{LemmaMI0}} if we set $\theta=\alpha_{\epsilon}$,
$b_{\theta}=b_{\epsilon}$ with $\epsilon,\alpha_{\epsilon}$, and
$b_{\epsilon}$ defined in (\ref{epsi3}) -see also (\ref{epsi.a})
and (\ref{epsi2.a}). $\blacksquare$

\bigskip{}

\textbf{Proof of Lemma \ref{LemmaMI2}}. \textbf{(i)} Recall that
$\theta_{1}>0$ and $\theta_{2}\in(1,2)$ such that $\theta_{1}\theta_{2}>1$.
Let $N>0$ and write
\[
\int^{\infty}_{x=0}e^{-\theta_{1}x}\left(\int^{x}_{y=0}e^{y}y^{-\theta_{2}}dy\right)^{\theta_{1}}dx\leq\int^{N}_{x=0}e^{-\theta_{1}x}\left(\int^{x}_{y=0}e^{y}y^{-\theta_{2}}dy\right)^{\theta_{1}}dx
\]
\[
+C\int^{\infty}_{x=N}e^{-\theta_{1}x}\left(\int^{N}_{y=0}e^{y}y^{-\theta_{2}}dy\right)^{\theta_{1}}dx+C\int^{\infty}_{x=N}e^{-\theta_{1}x}\left(\int^{x}_{y=N}e^{y}y^{-\theta_{2}}dy\right)^{\theta_{1}}dx
\]
\begin{equation}
=:I_{1}+CI_{2}+CI_{3}.\label{ProofLemmaMI2.2}
\end{equation}
 The first term
\[
I_{1}\leq\int^{N}_{x=0}e^{-\theta_{1}x}\left(\int^{N}_{y=0}e^{y}y^{-\theta_{2}}dy\right)^{\theta_{1}}dx
\]
\[
\leq C\left(\int^{N}_{y=0}y^{-\theta_{2}}dy\right)^{\theta_{1}}<\infty.
\]
Similarly, the second term
\[
I_{2}\leq\int^{\infty}_{x=N}e^{-\theta_{1}x}dx\left(\int^{N}_{y=0}e^{y}y^{-\theta_{2}}dy\right)^{\theta_{1}}<\infty.
\]
Finally, let $N$ large enough such that 
\[
\left|\frac{e^{-x}\int^{x}_{y=N}e^{y}y^{-\theta_{2}}dy}{x^{-\theta_{2}}}-1\right|<1,
\]
which is possible due to \textbf{Proposition \ref{SterlingApprox}}.
The third term 
\[
I_{3}=\int^{\infty}_{x=N}\left(e^{-x}\int^{x}_{y=N}e^{y}y^{-\theta_{2}}dy\right)^{\theta_{1}}dx
\]
\[
=\int^{\infty}_{x=N}\left(e^{-x}\int^{x}_{y=N}e^{y}y^{-\theta_{2}}dy-x^{-\theta_{2}}+x^{-\theta_{2}}\right)^{\theta_{1}}dx
\]
\begin{equation}
\leq\int^{\infty}_{x=N}\left(2x^{-\theta_{2}}\right)^{\theta_{1}}dx=C\int^{\infty}_{x=N}x^{-\theta_{1}\theta_{2}}dx<\infty,\label{ProofLemmaMI2.1}
\end{equation}
as required. 

\textbf{(ii)} It can be easily checked that (\ref{MI2b}) follows
from identical arguments to those used for (\ref{ProofLemmaMI2.1}).
We conclude the proof by establishing the validity of (\ref{MI2c}).
Recall that
\[
\gamma(\theta_{2};1,x):=\int^{x}_{y=1}e^{y}y^{1-\theta_{2}}dy.
\]
By \textbf{Proposition \ref{SterlingApprox}} as $x\rightarrow\infty$
\[
e^{-x}\gamma(\theta_{2};1,x)-x^{1-\theta_{2}}=O\left(x^{-\theta_{2}}\right).
\]
Therefore, for $x$ large enough, $x\geq x_{0}$ say, 
\[
\int^{\infty}_{x=x_{0}\lor1}\left|x^{1-\theta_{2}}-e^{-x}\gamma(\theta_{2};1,x)\right|^{\theta_{1}}dx<C\int^{\infty}_{x=x_{0}\lor1}x^{-\theta_{1}\theta_{2}}dx<\infty.
\]
Finally, noting that under our assumptions $\left(2-\theta_{2}\right)\theta_{1}>0$
and $\left(1-\theta_{2}\right)\theta_{1}+1>0$
\[
\int^{x_{0}\lor1}_{x=0}\left|x^{1-\theta_{2}}-e^{-x}\gamma(\theta_{2};1,x)\right|^{\theta_{1}}dx
\]
\[
\leq C\left[\int^{x_{0}\lor1}_{x=0}x^{\left(1-\theta_{2}\right)\theta_{1}}dx+\int^{x_{0}\lor1}_{x=0}\left(\left[y^{2-\theta_{2}}\right]^{x}_{0}\right)^{\theta_{1}}dx\right]
\]
\[
=C\left[\int^{x_{0}\lor1}_{x=0}x^{\left(1-\theta_{2}\right)\theta_{1}}dx+\int^{x_{0}\lor1}_{x=0}x^{\left(2-\theta_{2}\right)\theta_{1}}dx\right]<\infty,
\]
as required. $\blacksquare$

\bigskip{}

\textbf{Proof of Lemma} \textbf{\ref{LemmaMI3}}. Without loss of
generality assume that $m_{0}=1$. In view of the fact that $\alpha_{\epsilon}>1$
and (\ref{MVT}), the difference between $V_{\epsilon}\left(0,\delta_{n}\right)$
and $\kappa^{-b_{\epsilon}}_{n}\sum^{\delta_{n}}_{k=0}\left(\sum^{k}_{l=0}e^{-\frac{l}{\kappa_{n}}}c_{k-l}\right)^{\alpha_{\epsilon}}$
is
\[
\left\vert \kappa^{-b_{\epsilon}}_{n}\sum^{\delta_{n}}_{k=0}\left\vert \sum^{k}_{l=0}\rho^{l}_{n}c_{k-l}\right\vert ^{\alpha_{\epsilon}}-\kappa^{-b_{\epsilon}}_{n}\sum^{\delta_{n}}_{k=0}\left\vert \sum^{k}_{l=0}e^{-\frac{l}{\kappa_{n}}}c_{k-l}\right\vert ^{\alpha_{\epsilon}}\right\vert 
\]
\[
\leq C\kappa^{-b_{\epsilon}}_{n}\sum^{\delta_{n}}_{k=0}\left\vert \sum^{k}_{l=0}\left(\rho^{l}_{n}-e^{-\frac{l}{\kappa_{n}}}\right)c_{k-l}\right\vert ^{\alpha_{\epsilon}}+C\kappa^{-b_{\epsilon}}_{n}\sum^{\delta_{n}}_{k=0}\left\vert \sum^{k}_{l=0}e^{-\frac{l}{\kappa_{n}}}c_{k-l}\right\vert ^{\alpha_{\epsilon}-1}\left\vert \sum^{k}_{l=0}\left(\rho^{l}_{n}-e^{-\frac{l}{\kappa_{n}}}\right)c_{k-l}\right\vert 
\]
\[
\leq C\kappa^{-b_{\epsilon}}_{n}\sum^{\delta_{n}}_{k=0}\left[\sum^{k}_{l=0}\left\vert \left(\rho^{l}_{n}-e^{-\frac{l}{\kappa_{n}}}\right)c_{k-l}\right\vert \right]^{\alpha_{\epsilon}}
\]
\begin{equation}
+C\kappa^{-b_{\epsilon}}_{n}\sum^{\delta_{n}}_{k=0}\left[\sum^{k}_{l=0}\left\vert e^{-\frac{l}{\kappa_{n}}}c_{k-l}\right\vert \right]^{\alpha_{\epsilon}-1}\sum^{k}_{l=0}\left\vert \left(\rho^{l}_{n}-e^{-\frac{l}{\kappa_{n}}}\right)c_{k-l}\right\vert :=C\left(B_{1n}+B_{2n}\right).\label{VarMI-LM6}
\end{equation}
Using the approximation $c_{l}\sim l^{-m}$, the first term
\[
B_{1n}=\kappa^{-b_{\epsilon}}_{n}\sum^{\delta_{n}}_{k=0}\left[\sum^{k}_{l=0}\left\vert \left(\rho^{l}_{n}-e^{-\frac{l}{\kappa_{n}}}\right)c_{k-l}\right\vert \right]^{\alpha_{\epsilon}}\leq_{\left(1\right)}C\left(\frac{\delta_{n}}{\kappa^{2}_{n}}\right)^{\alpha_{\epsilon}}\kappa^{-b_{\epsilon}}_{n}\sum^{\delta_{n}}_{k=0}\left[\sum^{k}_{l=0}\left\vert c_{k-l}\right\vert \right]^{\alpha_{\epsilon}}
\]
\[
=C\left(\frac{\delta_{n}}{\kappa^{2}_{n}}\right)^{\alpha_{\epsilon}}\kappa^{-b_{\epsilon}}_{n}\sum^{\delta_{n}}_{k=0}\left[\sum^{k}_{l=0}\left\vert c_{l}\right\vert \right]^{\alpha_{\epsilon}}=\alpha_{\epsilon}\left(\frac{\delta_{n}}{\kappa^{2}_{n}}\right)^{\alpha_{\epsilon}}\kappa^{-b_{\epsilon}}_{n}\sum^{\delta_{n}}_{k=0}\left[\sum^{k}_{l=1}l^{-m}\right]^{\alpha_{\epsilon}}\left[1+o(1)\right]
\]
\[
=C\left(\frac{\delta_{n}}{\kappa^{2}_{n}}\right)^{\alpha_{\epsilon}}\kappa^{-b_{\epsilon}}_{n}\delta^{b_{\epsilon}}_{n}\frac{1}{\delta_{n}}\sum^{\delta_{n}}_{k=0}\left[\frac{1}{\delta_{n}}\sum^{k}_{l=0}\left(\frac{l}{\delta_{n}}\right)^{-m}\right]^{\alpha_{\epsilon}}\left[1+o(1)\right]
\]
\[
=C\left(\frac{\delta_{n}}{\kappa^{2}_{n}}\right)^{\alpha_{\epsilon}}\frac{\delta^{1+\alpha_{\epsilon}-\alpha_{\epsilon}m}_{n}}{\kappa^{1+\alpha_{\epsilon}-\alpha_{\epsilon}m}_{n}}\int^{1}_{0}\left(\int^{x}_{0}y^{-m}dy\right)^{\alpha_{\epsilon}}dx\left[1+o(1)\right]
\]
\begin{equation}
=C\left(\frac{\kappa_{n}\varsigma_{n}}{\kappa^{2}_{n}}\right)^{\alpha_{\epsilon}}\varsigma^{b_{\epsilon}}_{n}\int^{1}_{0}\left(\int^{x}_{0}y^{-m}dy\right)^{\alpha_{\epsilon}}dx\left[1+o(1)\right]=_{(2)}o(1),\label{VarMI-LM7}
\end{equation}
where $\leq_{\left(1\right)}$ follows from DK2018 (Lemma 6.4), and
$=_{(2)}$ from the fact that $\varsigma_{n}$ is SV in $\kappa_{n}$.
Similarly, the second term
\[
B_{2n}\leq C\left(\frac{\delta_{n}}{\kappa^{2}_{n}}\right)\kappa^{-b_{\epsilon}}_{n}\sum^{\delta_{n}}_{k=0}\left[\sum^{k}_{l=0}\left\vert e^{-\frac{l}{\kappa_{n}}}c_{k-l}\right\vert \right]^{\alpha_{\epsilon}-1}\sum^{k}_{l=0}\left\vert c_{k-l}\right\vert 
\]
\[
\leq C\left(\frac{\delta_{n}}{\kappa^{2}_{n}}\right)\kappa^{-b_{\epsilon}}_{n}\left[\sum^{\delta_{n}}_{k=0}\left\{ \left[\sum^{k}_{l=0}\left\vert e^{-\frac{l}{\kappa_{n}}}c_{k-l}\right\vert \right]^{\alpha_{\epsilon}-1}\right\} ^{\frac{\alpha_{\epsilon}}{\alpha_{\epsilon}-1}}\right]^{\frac{\alpha_{\epsilon}-1}{\alpha_{\epsilon}}}\left[\sum^{\delta_{n}}_{k=0}\left\{ \sum^{k}_{l=0}\left\vert c_{k-l}\right\vert \right\} ^{\alpha_{\epsilon}}\right]^{\frac{1}{\alpha_{\epsilon}}}
\]
\[
=C\left(\frac{\delta_{n}}{\kappa^{2}_{n}}\right)\delta^{b_{\epsilon}/\alpha_{\epsilon}}_{n}\kappa^{-b_{\epsilon}/\alpha_{\epsilon}}_{n}
\]
\begin{equation}
\cdot\left[\kappa^{-b_{\epsilon}}_{n}\sum^{\delta_{n}}_{k=0}\left\{ \sum^{k}_{l=0}\left\vert e^{-\frac{l}{\kappa_{n}}}c_{k-l}\right\vert \right\} ^{\alpha_{\epsilon}}\right]^{\frac{\alpha_{\epsilon}-1}{\alpha_{\epsilon}}}\left[\delta^{-b_{\epsilon}}_{n}\sum^{\delta_{n}}_{k=0}\left\{ \sum^{k}_{l=0}\left\vert c_{l}\right\vert \right\} ^{\alpha_{\epsilon}}\right]^{\frac{1}{\alpha_{\epsilon}}}.\label{VarMI-LM8}
\end{equation}
Euler summation together with the approximation $c_{k}\sim k^{-m}$
gives 
\begin{equation}
\kappa^{-b_{\epsilon}}_{n}\sum^{\delta_{n}}_{k=0}\left\{ \sum^{k}_{l=0}\left\vert e^{-\frac{l}{\kappa_{n}}}c_{k-l}\right\vert \right\} ^{\alpha_{\epsilon}}=\int^{\infty}_{0}e^{-\alpha_{\epsilon}x}\left(\int^{x}_{0}e^{y}y^{-m}dy\right)^{\alpha_{\epsilon}}dx+o(1)\label{VarMI-LM9}
\end{equation}
and
\begin{equation}
\delta^{-b_{\epsilon}}_{n}\sum^{\delta_{n}}_{k=0}\left\{ \sum^{k}_{l=0}\left\vert c_{l}\right\vert \right\} ^{\alpha_{\epsilon}}=\int^{1}_{0}\left(\int^{x}_{0}y^{-m}dy\right)^{\alpha_{\epsilon}}dx+o(1),\label{VarMI-LM10}
\end{equation}
where the integral term of (\ref{VarMI-LM9}) is finite by \textbf{Corollary
\ref{CorollaryMI2}}. Hence, in view of (\ref{VarMI-LM8})-(\ref{VarMI-LM10})
we have
\[
B_{2n}\leq C\left(\frac{\delta_{n}}{\kappa^{2}_{n}}\right)\delta^{b_{\epsilon}/\alpha_{\epsilon}}_{n}\kappa^{-\beta_{\epsilon}/\alpha_{\epsilon}}_{n}=\frac{\varsigma^{1+b_{\epsilon}/\alpha_{\epsilon}}_{n}}{\kappa_{n}}=o(1),
\]
where we have used the fact that $\varsigma_{n}$ is SV in $\kappa_{n}$.
Once more, using Euler summation together with (\ref{VarMI-LM7})
and (\ref{VarMI-LM10}) we have 
\[
V_{\epsilon}(0,\delta_{n})=\kappa^{-b_{\epsilon}}_{n}\sum^{\delta_{n}}_{k=0}\left\vert \sum^{k}_{l=0}e^{-\frac{l}{\kappa_{n}}}c_{k-l}\right\vert ^{\alpha_{\epsilon}}+o(1)
\]
\[
=\int^{\infty}_{0}e^{-\alpha_{\epsilon}x}\left(\int^{x}_{0}e^{y}y^{-m}dy\right)^{\alpha_{\epsilon}}dx+o(1),
\]
as required.

(b) Analogous arguments to those used for the proof of part (a) together
with Euler summation give
\[
V_{\epsilon}\left(0,\delta^{\prime}_{n}\right)=\kappa^{-b_{\epsilon}}_{n}\sum^{\delta^{\prime}_{n}}_{k=0}\left\vert \sum^{k}_{l=0}e^{-\frac{l}{\kappa_{n}}}c_{k-l}\right\vert ^{\alpha_{\epsilon}}+o(1)
\]
\[
=\int^{\delta^{\prime}_{n}/\kappa_{n}}_{0}e^{-\alpha_{\epsilon}x}\left(\int^{x}_{0}e^{y}y^{-m}dy\right)^{\alpha_{\epsilon}}dx+o(1)=o(1),
\]
as required. $\blacksquare$

\bigskip{}

\textbf{Proof of Lemma \ref{LemmaMI4}}. Without loss of generality
we assume that $m_{0}=1$. In addition, note that for $n$ large enough,
$0<\rho^{j}_{n}\leq1$, for all $j\geq0$. Set

\[
\tilde{c}_{s}:=\sum^{\infty}_{j=s}c_{j}.
\]
Write
\begin{equation}
a_{k}:=a_{k}(n)=\sum^{k}_{l=0}\rho^{k-l}_{n}c_{l}=\sum^{k}_{l=0}\rho^{k-l}_{n}\left(\tilde{c}_{l}-\tilde{c}_{l+1}\right).\label{LemmaMI4.1}
\end{equation}
Summation by parts gives
\[
\sum^{k}_{l=0}\rho^{-l}_{n}\left(\tilde{c}_{l}-\tilde{c}_{l+1}\right)
\]
\[
=\rho^{-0}_{n}\left(\tilde{c}_{0}-\tilde{c}_{0+1}\right)+\rho^{-1}_{n}\left(\tilde{c}_{1}-\tilde{c}_{2}\right)+\rho^{-2}_{n}\left(\tilde{c}_{2}-\tilde{c}_{3}\right)+...
\]
\[
...+\rho^{-k-1}_{n}\left(\tilde{c}_{k-1}-\tilde{c}_{k}\right)+\rho^{-k}_{n}\left(\tilde{c}_{k}-\tilde{c}_{k+1}\right)
\]
\[
=\tilde{c}_{0}-\rho^{-k}_{n}\tilde{c}_{k+1}-\tilde{c}_{1}+\rho^{-1}_{n}\tilde{c}_{1}-\rho^{-1}_{n}\tilde{c}_{2}+\rho^{-2}_{n}\tilde{c}_{2}-\rho^{-2}_{n}\tilde{c}_{3}+...
\]
\[
...+\rho^{-k-1}_{n}\tilde{c}_{k-1}-\rho^{-k-1}_{n}\tilde{c}_{k}+\rho^{-k}_{n}\tilde{c}_{k}
\]
\[
=\tilde{c}_{0}-\rho^{-k}_{n}\tilde{c}_{k+1}-\sum^{k}_{l=1}\tilde{c}_{l}\left(\rho^{-l+1}_{n}-\rho^{-l}_{n}\right)
\]
\[
=\tilde{c}_{0}-\rho^{-k}_{n}\tilde{c}_{k+1}-\sum^{k}_{l=1}\rho^{-l}_{n}\tilde{c}_{l}\left(\rho_{n}-1\right)
\]
\begin{equation}
=\tilde{c}_{0}-\rho^{-k}_{n}\tilde{c}_{k+1}+\kappa^{-1}_{n}\sum^{k}_{l=1}\tilde{c}_{l}\rho^{-l}_{n}.\label{LemmaMI4.2}
\end{equation}
\bigskip{}

In view of (\ref{LemmaMI4.1}) and (\ref{LemmaMI4.2}),
\[
\kappa^{-\left(1-m\right)}_{n}a_{k}=\kappa^{-\left(1-m\right)}_{n}\sum^{k}_{l=0}\rho^{k-l}_{n}c_{l}
\]
\begin{equation}
=T_{1nk}+T_{2nk}.\label{LemmaMI4.3}
\end{equation}
Finally, by (\ref{LemmaMI4.3})
\[
V_{\epsilon}\left(0,\delta_{n}\right)=\kappa^{-b_{\epsilon}}_{n}\sum^{\delta_{n}}_{k=1}\left\vert a_{k}\right\vert ^{\alpha_{\epsilon}}
\]
\[
=\kappa^{-\left(1+\alpha_{\epsilon}-\alpha_{\epsilon}m\right)}_{n}\sum^{\delta_{n}}_{k=1}\left\vert \sum^{k}_{l=0}\rho^{k-l}_{n}c_{l}\right\vert ^{\alpha_{\epsilon}}=\kappa^{-1}_{n}\sum^{\delta_{n}}_{k=1}\left\vert \kappa^{-\left(1-m\right)}_{n}\sum^{k}_{l=0}\rho^{k-l}_{n}c_{l}\right\vert ^{\alpha_{\epsilon}}
\]
\[
=\kappa^{-1}_{n}\sum^{\delta_{n}}_{k=1}\left\vert T_{1nk}+T_{2nk}\right\vert ^{\alpha_{\epsilon}}=\left\vert 1-m\right\vert ^{-\alpha_{\epsilon}}\int^{\infty}_{x=0}\left\vert p(x)\right\vert ^{\alpha_{\epsilon}}dx+o(1),
\]
where the last approximation above is a direct consequence of \textbf{Corollary
\ref{Ts_corol}}. $\blacksquare$

\bigskip{}

\textbf{Proof of Proposition \ref{SterlingApprox}}. \textbf{(a)}
Integration by parts gives
\[
q(x)=e^{-x}\int^{x}_{1}e^{y}y^{-q_{0}}dy
\]
\[
=e^{-x}\left\{ \left[e^{y}y^{-q_{0}}\right]^{x}_{1}+q_{0}\int^{x}_{1}e^{y}y^{-q_{0}-1}dy\right\} 
\]
\begin{equation}
=x^{-q_{0}}+e^{-x}+q_{0}e^{-x}\int^{x}_{1}e^{y}y^{-q_{0}-1}dy.\label{Sterling1}
\end{equation}
Consider the term
\[
e^{-x}\int^{x}_{1}e^{y}y^{-q_{0}-1}dy=e^{-x}\int^{x/2}_{1}e^{y}y^{-q_{0}-1}dy+e^{-x}\int^{x}_{x/2}e^{y}y^{-q_{0}-1}dy
\]
\[
\leq e^{-x/2}\int^{x/2}_{1}y^{-q_{0}-1}dy+\left(\frac{x}{2}\right)^{-q_{0}-1}e^{-x}\int^{x}_{x/2}e^{y}dy
\]
\begin{equation}
\leq C\left[e^{-x/2}+\left(\frac{x}{2}\right)^{-q_{0}-1}\right]=O\left(x^{-q_{0}-1}\right).\label{Sterling2}
\end{equation}
The result for $N=1$ follows directly from (\ref{Sterling1}) and
(\ref{Sterling2}). The general case follows from arguments similar
to those used above and induction.

\textbf{(b)} Recall that $0<q_{0}<1$. Integration by parts gives
\[
q(x)=e^{-x}\int^{x}_{0}e^{y}y^{-q_{0}}dy=e^{-x}\left[e^{y}\int y^{-q_{0}}dy\right]^{x}_{0}-e^{-x}\int^{x}_{0}\int y^{-q_{0}}dye^{y}dy
\]
\[
=e^{-x}\left[e^{y}\int y^{-q_{0}}dy\right]^{x}_{0}-e^{-x}\int^{x}_{0}y^{1-q_{0}}e^{y}dy=e^{-x}\left[e^{x}\frac{x^{1-q_{0}}}{1-q_{0}}\right]-e^{-x}\int^{x}_{0}y^{1-q_{0}}e^{y}dy
\]
\[
=\frac{x^{1-q_{0}}}{1-q_{0}}-e^{-x}\int^{x}_{0}y^{1-q_{0}}e^{y}dy.
\]
Next, note that as $x\rightarrow0$
\[
\frac{e^{-x}\int^{x}_{0}y^{1-q_{0}}e^{y}dy}{x^{1-q_{0}}}\leq\frac{\int^{x}_{0}y^{1-q_{0}}dy}{x^{1-q_{0}}}=\frac{1}{1-q_{0}}\frac{x^{2-q_{0}}}{x^{1-q_{0}}}\rightarrow0.
\]
In view of the above
\[
q(x)\sim\frac{x^{1-q_{0}}}{1-q_{0}}.
\]
Therefore, the requisite result holds for $N=1$. The general case
follows from arguments analogous to those used above and induction.
$\blacksquare$

\bigskip{}

\[
\text{\textbf{Reversed Gamma Function}}
\]
For $a<x$ set.
\[
\gamma(m;a,x):=\int^{x}_{y=a}e^{y}y^{1-m}dy.
\]
Then integration by parts gives
\[
\gamma(m;a,x)=\left[e^{x}x^{1-m}-a^{1-m}e^{a}\right]-\left(1-m\right)\int^{x}_{y=a}e^{y}y^{-m}dy
\]
\[
=\left[e^{x}x^{1-m}-a^{1-m}e^{a}\right]+\left(m-1\right)\gamma(m+1;a,x).
\]
Hence,
\[
x^{1-m}-e^{-x}\gamma(m,2-m,x)=x^{1-m}-e^{-x}\int^{x}_{y=2-m}e^{y}y^{1-m}dy
\]
\[
=x^{1-m}-e^{-x}\left\{ \left[e^{y}y^{1-m}\right]^{x}_{2-m}+\left(m-1\right)\int^{x}_{y=2-m}e^{y}y^{-m}dy\right\} 
\]
\[
=x^{1-m}-e^{-x}\left\{ \left[e^{x}x^{1-m}-\left(2-m\right)^{1-m}e^{2-m}\right]+\left(m-1\right)\int^{x}_{y=2-m}e^{y}y^{-m}dy\right\} 
\]
\[
=e^{-x}\left[\left(2-m\right)^{1-m}e^{2-m}-\left(m-1\right)\int^{x}_{y=2-m}e^{y}y^{-m}dy\right].
\]
For $x>2-m$ the above imply
\[
p(x)=x^{1-m}-e^{-x}\gamma(m,0,x)=x^{1-m}-e^{-x}\gamma(m,2-m,x)-e^{-x}\gamma(m,0,2-m)
\]
\[
=e^{-x}\left[\left(2-m\right)^{1-m}e^{2-m}-\left(m-1\right)\int^{x}_{y=2-m}e^{y}y^{-m}dy\right]-e^{-x}\gamma(m,0,2-m)
\]
\[
=e^{-x}\left[\left(2-m\right)^{1-m}e^{2-m}-\left(m-1\right)\int^{x}_{y=2-m}e^{y}y^{-m}dy-\int^{2-m}_{y=0}e^{y}y^{1-m}dy\right].
\]

\bigskip{}
\[
\text{\textbf{Roots of }}p(x)
\]
Recall that
\[
p(x)=x^{1-m}-e^{-x}\int^{x}_{0}e^{y}y^{1-m}dy.
\]

\begin{proposition} \label{pi-roots} \emph{$p(x)$ has a single
root $x_{*}$ on $(0,\infty)$ such that $x_{*}>2-m$.}

\end{proposition} 

\textbf{Case }$x\in(0,2-m)$:
\[
p(x)=x^{1-m}-e^{-x}\int^{x}_{0}e^{y}y^{1-m}dy\geq x^{1-m}-\int^{x}_{0}y^{1-m}dy
\]
\[
\geq x^{1-m}-\frac{x^{2-m}}{2-m}:=p_{1}(x)
\]
Now 
\[
p_{1}(x)>0\Longleftrightarrow x^{1-m}\geq\frac{x^{2-m}}{2-m}\Longleftrightarrow\left(1-m\right)\ln x>(2-m)\ln(x)-\ln(2-m)
\]
\[
\Longleftrightarrow\ln(2-m)>(2-m)\ln(x)-\left(1-m\right)\ln\left(x\right)
\]
\[
\Longleftrightarrow\ln(2-m)>\ln\left(x\right)\Longleftrightarrow(2-m)>x.
\]
Hence,
\[
p(x)>0\text{ for }0<x<2-m
\]

\bigskip{}
\textbf{Case} $x\geq2-m$:
\[
p(x)=e^{-x}\left[\left[2-m\right]^{1-m}e^{2-m}-\left(m-1\right)\int^{x}_{y=2-m}e^{y}y^{-m}dy-\int^{2-m}_{0}e^{y}y^{1-m}dy\right].
\]
The last term above
\[
\int^{2-m}_{0}e^{y}y^{1-m}dy<e^{2-m}\int^{2-m}_{0}y^{1-m}dy=e^{2-m}\frac{\left(2-m\right)^{2-m}}{2-m}=e^{2-m}\left(2-m\right)^{1-m}.
\]
Hence,
\[
p(x)=e^{-x}\left[\underset{=:p_{2}}{\underbrace{\left[2-m\right]^{1-m}e^{2-m}-\int^{2-m}_{0}e^{y}y^{1-m}dy}}-\left(m-1\right)\int^{x}_{y=2-m}e^{y}y^{-m}dy\right],
\]
 which implies that 
\[
p(2-m)=e^{-\left(2-m\right)}p_{2}>0.
\]
Moreover, in view of this we can write 
\[
p(x)=e^{-x}\left[p_{2}-\left(m-1\right)\int^{x}_{y=2-m}e^{y}y^{-m}dy\right],
\]
where $p_{2}=\left[2-m\right]^{1-m}e^{2-m}-\int^{2-m}_{0}e^{y}y^{1-m}dy>0$.
Noting that $\int^{x}_{y=2-m}e^{y}y^{-m}dy$ is continuous in $x$,
strictly increasing and diverging, by the intermediate value Theorem
there is a some $x^{\ast}>2-m$ such that $p(x^{\ast})=0$. Further,
$x^{\ast}$ is a unique root due to the monotonicity of the aforementioned
term.

\bigskip{}

\section{Uniformly Bounded Density Results}

We start with some notation and preliminary results. First, we will
expand the definition of the sequence $m_{n}$ introduced in (\ref{MI8d}).
Recall that $m\in(1/\alpha,1)$, $m\in(1,1+1/\alpha)$ under \textbf{LM}
and \textbf{NM} respectively. For \textbf{SM}, we shall use the convention
$m=1$. Define
\begin{equation}
m_{n}:=\gamma_{n}\kappa^{-\left(1-m\right)}_{n}\text{, }\text{ with \ensuremath{m}\ensuremath{\in\left\{  \begin{array}{l}
(1/\alpha,1)\text{, under \textbf{LM}, }\alpha\in(1,2];\\
\left\{ 1\right\} \text{, under \textbf{SM}, }\alpha\in(0,2];\\
(1,1+1/\alpha)\text{, under \textbf{NM}, }\alpha\in(1,2].
\end{array}\right.}}\label{m_n(1)}
\end{equation}
In addition, recall that
\[
\gamma_{n}\sim\left\{ \begin{array}{l}
L\left(\gamma_{n}n^{1/\alpha}\right)^{1/\alpha}\text{, under \textbf{BFR}}\\
\ell_{1}\left(\kappa_{n}\right)^{1/\alpha}\kappa^{b/\alpha}_{n}\text{, }\text{ under \textbf{MI}}
\end{array}\right.,
\]
where $b=1+\alpha-\alpha m$ and $\ell_{1}(x)\sim H\left(x^{1/\alpha}\ell_{1}(x)^{1/\alpha}\right)$.
By construction we have
\[
m_{n}\rightarrow\infty.
\]
In particular, (see Remark \ref{RmkMI-LM/NM(1)})
\begin{equation}
m_{n}\sim\ell_{1}\left(\kappa_{n}\right)^{1/\alpha}\kappa^{1/\alpha}_{n}\rightarrow\infty.\label{m_n(1b)}
\end{equation}

\bigskip{}

Finally, set 
\begin{equation}
\left[\eta_{n},\overline{q}_{n},\underline{q}_{n}\right]:=\left\{ \begin{array}{l}
\left[n^{1/\alpha},n-1,\left\lfloor \gamma_{n}\right\rfloor \right]\text{, under \textbf{BFR}}\\
\left[1,\left\lfloor \delta\kappa_{n}\right\rfloor ,\left\lfloor \delta\kappa_{n}/2\right\rfloor \right]\text{, }\delta\in(0,1]\text{ under \textbf{MI}}
\end{array}\right.,\label{m_n(2)}
\end{equation}
and 
\[
\delta_{n}:=\left\{ \begin{array}{l}
1\text{, under \textbf{BFR}}\\
\kappa^{1-m}_{n},\text{ under \textbf{MI}}
\end{array}\right.
\]

Since the modulus of a characteristic function does not exceed one,
if Assumption INN(iii) holds for some $\theta\geq1$, it also holds
for every $\theta^{\prime}\geq\theta$. Hence, without loss of generality,
throughout this section we take $\theta\in\mathbb{N}$\textcolor{red}{.}

\bigskip{}
We have the following technical result that will be utilized in this
section.\bigskip{}

\begin{lemma} \label{BoundCoefGeneral} Suppose that \textbf{Assumption
INN}, and \textbf{Assumption BFR} or \textbf{MI} hold. 

(i) There are $n_{0},t_{0}\geq1$, sequences $\left[\eta_{n},\overline{q}_{n},\underline{q}_{n}\right]$
that conform with (\ref{m_n(2)}), and constants $0<\mathcal{D}_{0}\leq\mathcal{D}_{1}<\infty$
such that
\begin{equation}
\mathcal{D}_{0}<\inf_{\underline{q}_{t}\leq k\leq\overline{q}_{t}}\delta^{-1}_{t}\left\vert \eta_{k}a_{k}(n)\right\vert \leq\sup_{\underline{q}_{t}\leq k\leq\overline{q}_{t}}\delta^{-1}_{t}\left\vert \eta_{k}a_{k}(n)\right\vert \leq\mathcal{D}_{1},\label{bracket1}
\end{equation}
for all $n\geq n_{0}$ and $t_{0}\leq t\leq n$.

(ii) There is an $n_{0}\geq1$ and some $\zeta\in(0,1)$ such that
as $t\rightarrow\infty$
\begin{equation}
\sup_{0\leq k\leq\overline{q}_{t},n\geq n_{0}}\left\vert \frac{a_{k}(n)}{\gamma^{\zeta}_{t}}\right\vert =o(1).\label{bracket2-1}
\end{equation}

\end{lemma}\bigskip{}

\begin{remark} \label{RmkUBD3} Notice that as $n\rightarrow\infty$,
$n\geq\overline{q}_{n}>\underline{q}_{n}>0$, and $q_{n}:=\overline{q}_{n}-\underline{q}_{n}\rightarrow\infty$. 

\end{remark} 

\bigskip{}

For all $\lambda\in\mathbb{R}$ define
\begin{equation}
G^{\ast}(x):=x^{-\alpha}H(x),\quad x>0,\qquad G(\lambda):=\left\{ \begin{array}{l}
\left\vert \lambda\right\vert ^{\alpha}H(\left\vert \lambda\right\vert ^{-1}),\text{ }\lambda\neq0\\
0,\text{ }\lambda=0
\end{array}\right..\label{Gee}
\end{equation}

\bigskip{}

\begin{lemma} \label{BoundCoefGeneral-1} 

(i) Suppose that \textbf{Assumption INN}, and \textbf{Assumption BFR}
or \textbf{MI} hold, and $C_{0}\in(0,\infty)$. There are $n_{0},t_{0}\geq1$
\[
\inf_{n\geq n_{0},t_{0}\leq t\leq n}\sum^{\overline{q}_{t}}_{k=\underline{q}_{t}+\theta}G\left(C_{0}\frac{a_{k}(n)}{\gamma_{t}}\right)>0.
\]

(ii) Let $r_{1},r_{2}>0$. As $n\rightarrow\infty$,
\[
\exp\left(-r_{1}\gamma^{r_{2}}_{n}\right)\gamma^{2}_{n}\delta^{-2}_{n}\eta^{2}_{\underline{q}_{n}+\theta}\rightarrow0.
\]

\end{lemma}

\bigskip{}

Under \textbf{Assumption INN}, for every $\eta\in(0,1)$ there is
an $N^{\prime}_{\eta}>0$ such that 
\begin{equation}
\left\vert \psi_{\xi}(\lambda)\right\vert \leq e^{-\eta\left\vert \lambda\right\vert ^{\alpha}H(\left\vert \lambda\right\vert ^{-1})}\:\:\text{for all \:}\:0<\left\vert \lambda\right\vert \leq N^{\prime}_{\eta}.\label{UBD1}
\end{equation}
In addition, there is some $N>0$ such that $\left\vert \lambda\right\vert ^{\alpha}H(\left\vert \lambda\right\vert ^{-1})<1/2$
for all $0<\left\vert \lambda\right\vert \leq N^{\prime}$ (see e.g.
BGT 1987; Proposition 1.3.6(v)). Set $N_{\eta}:=N^{\prime}_{\eta}\wedge N^{\prime}$.
Moreover, because $\xi$ has an absolutely continuous distribution,
for each $\delta>0$ there some $\rho_{\delta}$, $0<\rho_{\delta}<1$
such that (see e.g. Feller, 1971; p. 516)
\begin{equation}
\left\vert \psi_{\xi}(\lambda)\right\vert \leq\rho_{\delta}\text{,}\:\left\vert \lambda\right\vert >\delta.\label{UBD1b}
\end{equation}
Combining (\ref{UBD1})-(\ref{UBD1b}) we have 
\[
\left\vert \psi_{\xi}(\lambda)\right\vert \leq e^{-\eta G(\lambda)}1_{\left\{ \left\vert \lambda\right\vert \leq N_{\eta}\right\} }+e^{\ln\rho_{N_{\eta}}}1_{\left\{ \left\vert \lambda\right\vert >N_{\eta}\right\} }
\]
\[
=:e^{-\eta G(\lambda)}1_{\left\{ \left\vert \lambda\right\vert \leq N_{\eta}\right\} }+e^{-\vartheta_{\eta}}1_{\left\{ \left\vert \lambda\right\vert >N_{\eta}\right\} }
\]
\[
\leq e^{-\left(\eta\wedge\vartheta_{\eta}\right)G(\lambda)}1_{\left\{ \left\vert \lambda\right\vert \leq N_{\eta}\right\} }+e^{-\left(\eta\wedge\vartheta_{\eta}\right)}1_{\left\{ \left\vert \lambda\right\vert >N_{\eta}\right\} }
\]
\begin{equation}
=:e^{-\eta_{\ast}G(\lambda)}1_{\left\{ \left\vert \lambda\right\vert \leq N_{\eta}\right\} }+e^{-\eta_{\ast}}1_{\left\{ \left\vert \lambda\right\vert >N_{\eta}\right\} }\label{UBD1c}
\end{equation}
\begin{equation}
\leq e^{-\eta_{\ast}\left[G(\lambda)\wedge1\right]}1_{\left\{ \left\vert \lambda\right\vert \leq N_{\eta}\right\} }+e^{-\eta_{\ast}\left[G(\lambda)\wedge1\right]}1_{\left\{ \left\vert \lambda\right\vert >N_{\eta}\right\} }=e^{-\eta_{\ast}\left[G(\lambda)\wedge1\right]}.\label{UBD1d}
\end{equation}

\bigskip{}

Noting that $G^{\ast}(.)$ (cf. eq. (\ref{Gee})) is regularly varying
of index $-\alpha$, there no loss of generality (see e.g. Theorem
1.5.3 in BGT 1987) in the subsequent derivations, if we assume that
$G^{\ast}(.)$ is decreasing ($G(.)$ increasing). To see this note
that for large $\left\vert \lambda\right\vert ^{-1}$, there is some
decreasing regularly varying function of index $-\alpha$, $\underline{G}^{\ast}(.)$
say, such that $G^{\ast}(\mu)\sim\underline{G}^{\ast}(\mu)$, $\mu\rightarrow\infty$.
Hence, for $\mu$ large enough, $\mu\geq1/N_{\eta}$ say 
\[
\left|G^{\ast}(\mu)/\underline{G}^{\ast}(\mu)-1\right|\leq1/2.
\]
Hence,
\[
G^{\ast}(\mu)=\left[G^{\ast}(\mu)/\underline{G}^{\ast}(\mu)-1\right]\underline{G}^{\ast}(\mu)+\underline{G}^{\ast}(\mu)\geq-\frac{1}{2}\underline{G}^{\ast}(\mu)+\underline{G}^{\ast}(\mu)=\frac{1}{2}\underline{G}^{\ast}(\mu).
\]
Setting $\mu=\left\vert \lambda\right\vert ^{-1}$, $\left\vert \lambda\right\vert \leq N_{\eta}$
we get 
\[
\left\vert \psi_{\xi}(\lambda)\right\vert \leq e^{-\eta G(\lambda)}=e^{-\eta G^{*}(|\lambda|^{-1})}\leq e^{-\eta\frac{1}{2}\underline{G}^{\ast}(|\lambda|^{-1})}.
\]
Therefore, there is no loss of generality if we use the convention
$G(\lambda)=\frac{1}{2}\underline{G}^{\ast}(|\lambda|^{-1})$.

\bigskip{}

\begin{remark} \textbf{\label{RmkUBD1}}In view of the the monotonicity
convention for $G(.)$, we have
\[
\sup_{\left\vert \lambda\right\vert \geq N_{\eta}}e^{-\eta_{\ast}\left[G(\lambda)\wedge1\right]}\leq e^{-\eta_{\ast}\left[G(N_{\eta})\wedge1\right]}=e^{-\eta_{\ast}G(N_{\eta})}.
\]

Moreover, by construction, $G(\lambda)<1/2$ for all $0<\left\vert \lambda\right\vert \leq N_{\eta}$.
Hence,
\[
\sup_{0<\left\vert \lambda\right\vert \leq N_{\eta}}e^{-\eta_{\ast}\left[G(\lambda)\wedge1\right]}\leq\sup_{0<\left\vert \lambda\right\vert \leq N_{\eta}}e^{-\eta_{\ast}G(\lambda)}.
\]

\end{remark} 

\bigskip{}

In addition, we shall use the following facts:

\begin{remark} \textbf{\label{RmkUBD4}}
\begin{itemize}
\item[\textbf{(i)}] Fix $\epsilon>0$ and note that $\left\vert \lambda\right\vert ^{-\epsilon}H(\left\vert \lambda\right\vert ^{-1})\rightarrow\infty$
as $\lambda\rightarrow0$. Hence, for all $\left|\lambda\right|$
($\lambda\neq0$) small enough
\[
\left\vert \lambda\right\vert ^{\alpha}H(\left\vert \lambda\right\vert ^{-1})=\left\vert \lambda\right\vert ^{\alpha+\epsilon}\left\vert \lambda\right\vert ^{-\epsilon}H(\left\vert \lambda\right\vert ^{-1})\geq\left\vert \lambda\right\vert ^{\alpha+\epsilon}.
\]
\item[\textbf{(ii)}]  We remark that (\ref{bracket2-1}) implies that there is $0<\mathcal{D}_{2}<\infty$
such that 
\begin{equation}
\sup_{t\geq1}\sup_{0\leq k\leq\overline{q}_{t},n\geq n_{0}}\left\vert \frac{a_{k}(n)}{\gamma^{\zeta}_{t}}\right\vert \leq\mathcal{D}_{2}.\label{bracket2}
\end{equation}
In addition, as $t\rightarrow\infty$
\begin{equation}
\inf_{0\leq k\leq\overline{q}_{t},n\geq n_{0}}\left\vert \frac{a_{k}(n)}{\gamma^{\zeta}_{t}}\right\vert ^{-1}\rightarrow\infty.\label{bracket3}
\end{equation}
\item[\textbf{(iii)}]  Fix $A\in\left(1,\infty\right)$ and $\epsilon\in(0,\alpha)$. In
view of (\ref{bracket3}) and Potter's Theorem for all $\lambda\in[1,\mathcal{D}^{-1}_{2}\gamma^{1-\zeta}_{t}]$,
$\underline{q}_{t}\leq k\leq\overline{q}_{t}$, and $t$ large enough
we have
\[
\frac{H\left(\left\vert \frac{\gamma_{t}}{a_{k}(n)}\right\vert \right)}{H\left(\left\vert \frac{\gamma_{t}}{a_{k}(n)\lambda}\right\vert \right)}\leq A\left[\left(\frac{\left\vert \frac{\gamma_{t}}{a_{k}(n)}\right\vert }{\left\vert \frac{\gamma_{t}}{a_{k}(n)\lambda}\right\vert }\right)^{\epsilon}\vee\left(\frac{\left\vert \frac{\gamma_{t}}{a_{k}(n)}\right\vert }{\left\vert \frac{\gamma_{t}}{a_{k}(n)\lambda}\right\vert }\right)^{-\epsilon}\right]
\]
\[
=A\left[\left\vert \lambda\right\vert ^{\epsilon}\vee\left\vert \lambda\right\vert ^{-\epsilon}\right]=A\left\vert \lambda\right\vert ^{\epsilon}.
\]
Hence, for all $\lambda\in[1,\mathcal{D}^{-1}_{2}\gamma^{1-\zeta}_{t}]$,
$\underline{q}_{t}\leq k\leq\overline{q}_{t}$, and $t$ large enough
\[
G\left(\frac{a_{k}(n)}{\gamma_{t}}\lambda\right)=\left\vert \frac{a_{k}(n)}{\gamma_{t}}\lambda\right\vert ^{\alpha}H\left(\left\vert \frac{\gamma_{t}}{a_{k}(n)\lambda}\right\vert \right)
\]
\begin{equation}
\geq A^{-1}\left\vert \lambda\right\vert ^{-\epsilon}\left\vert \frac{a_{k}(n)}{\gamma_{t}}\lambda\right\vert ^{\alpha}H\left(\left\vert \frac{\gamma_{t}}{a_{k}(n)}\right\vert \right)=A^{-1}\left\vert \lambda\right\vert ^{\alpha-\epsilon}G\left(\frac{a_{k}(n)}{\gamma_{t}}\right).\label{UBDeq2b}
\end{equation}
\end{itemize}
 Let $\zeta\in(0,1)$. Fix $\epsilon\in\left(0,\alpha\right)$ and
$A>1$. Then there is some $t_{1}>0$ such that for all $t\geq t_{1}$
the following hold
\[
\frac{H\left(\left\vert \frac{a_{k}(n)}{\mathcal{D}_{2}\gamma_{t}}\right\vert ^{-1}\right)}{H\left(\left\vert \frac{a_{k}(n)}{\mathcal{D}_{2}\gamma^{\zeta}_{t}}\right\vert ^{-1}\right)}\leq A\left[\left(\frac{\left\vert \frac{a_{k}(n)}{\mathcal{D}_{2}\gamma_{t}}\right\vert ^{-1}}{\left\vert \frac{a_{k}(n)}{\mathcal{D}_{2}\gamma^{\zeta}_{t}}\right\vert ^{-1}}\right)^{\epsilon}\vee\left(\frac{\left\vert \frac{a_{k}(n)}{\mathcal{D}_{2}\gamma_{t}}\right\vert ^{-1}}{\left\vert \frac{a_{k}(n)}{\mathcal{D}_{2}\gamma^{\zeta}_{t}}\right\vert ^{-1}}\right)^{-\epsilon}\right]
\]
\[
=C\left[\left(\frac{\gamma_{t}}{\gamma^{\zeta}_{t}}\right)^{\epsilon}\vee\left(\frac{\gamma_{t}}{\gamma^{\zeta}_{t}}\right)^{-\epsilon}\right]=A\left[\gamma^{\epsilon\left(1-\zeta\right)}_{t}\vee\gamma^{-\epsilon\left(1-\zeta\right)}_{t}\right]=A\gamma^{\epsilon\left(1-\zeta\right)}_{t},
\]
where the inequality above is due to Potter's Theorem. A simple rearrangement
of the terms above gives
\[
G\left(\frac{a_{k}(n)}{\mathcal{D}_{2}\gamma^{\zeta}_{t}}\right)=\left\vert \frac{a_{k}(n)}{\mathcal{D}_{2}\gamma^{\zeta}_{t}}\right\vert ^{\alpha}H\left(\left\vert \frac{a_{k}(n)}{\mathcal{D}_{2}\gamma^{\zeta}_{t}}\right\vert ^{-1}\right)
\]
\[
\geq\left\vert \frac{a_{k}(n)}{\mathcal{D}_{2}\gamma^{\zeta}_{t}}\right\vert ^{\alpha}H\left(\left\vert \frac{a_{k}(n)}{\mathcal{D}_{2}\gamma_{t}}\right\vert ^{-1}\right)\cdot\frac{1}{A\gamma^{\epsilon\left(1-\zeta\right)}_{t}}
\]
\[
=\frac{\gamma^{\alpha}_{t}}{A\gamma^{\epsilon\left(1-\zeta\right)}_{t}\gamma^{\alpha\zeta}_{t}}\left\vert \frac{a_{k}(n)}{\mathcal{D}_{2}\gamma_{t}}\right\vert ^{\alpha}H\left(\left\vert \frac{a_{k}(n)}{\mathcal{D}_{2}\gamma_{t}}\right\vert ^{-1}\right)
\]
\begin{equation}
=A^{-1}\gamma^{\left(\alpha-\epsilon\right)\left(1-\zeta\right)}_{t}G\left(\frac{a_{k}(n)}{\mathcal{D}_{2}\gamma_{t}}\right).\label{UBDeq2}
\end{equation}

\end{remark}

\bigskip{}

\begin{lemma} \label{UBD} Suppose that the conditions of \textbf{Lemma
\ref{BoundCoefGeneral}} hold. There are some $n_{0}\geq t_{0}\geq1$
such that
\[
\sup_{n\geq n_{0},t_{0}\leq t\leq n}\int_{\mathbb{R}}|\lambda|\left|\mathbf{E}\left(e^{i\lambda x^{+}_{t}(n)/\gamma_{t}}\right)\right|d\lambda<\infty.
\]
and
\[
\sup_{n\geq n_{0},t_{0}\leq t\leq n}\int_{\mathbb{R}}\left\vert \lambda\right\vert \left\vert \mathbf{E}\left(e^{i\lambda\frac{x_{t}(n)}{\gamma_{t}}}\right)\right\vert d\lambda<\infty.
\]
\end{lemma}

\bigskip{}

\textbf{Proof of Lemma \ref{UBD}}. Fix $\eta\in(0,1)$ and let $N_{\eta}$,
$\eta^{*}$ as in (\ref{UBD1c}). Set $\mathcal{N}=N_{\eta}\mathcal{D}^{-1}_{2}$,
with $\mathcal{D}_{2}$ given in (\ref{bracket2}). For some $\zeta\in(0,1)$,\textbf{
}since $x^{+}_{t}(n)$ and $x^{-}_{t}(n)$ are independent and the
modulus of the characteristic function of $x^{-}_{t}(n)$ is at most
one, it suffices to bound the characteristic function of $x^{+}_{t}(n)$.
Write
\[
\int_{\mathbb{R}}\left\vert \lambda\right\vert \left\vert \mathbf{E}\left(e^{i\lambda\frac{x_{t}(n)}{\gamma_{t}}}\right)\right\vert d\lambda\leq\int_{\mathbb{R}}\left\vert \lambda\right\vert \left\vert \mathbf{E}\left(e^{i\lambda\gamma^{-1}_{t}\sum^{t-1}_{k=0}a_{k}(n)\xi_{t-k}}\right)\right\vert d\lambda
\]
\[
=\int_{\mathbb{R}}\left\vert \lambda\right\vert \prod^{t-1}_{k=0}\left\vert \psi\left(\frac{\lambda}{\gamma_{t}}a_{k}(n)\right)\right\vert d\lambda
\]
\[
=\left(\int_{\left\vert \lambda\right\vert <1}+\int_{1\leq\left\vert \lambda\right\vert \leq\gamma^{1-\zeta}_{t}\mathcal{N}}+\int_{\left\vert \lambda\right\vert >\gamma^{1-\zeta}_{t}\mathcal{N}}\right)\left\vert \lambda\right\vert \prod^{t-1}_{k=0}\left\vert \psi\left(\frac{\lambda}{\gamma_{t}}a_{k}(n)\right)\right\vert d\lambda=:I_{0n}(t)+I_{1n}(t)+I_{2n}(t).
\]
Clearly, $I_{0n}(t)\leq2$.

Next, let $\epsilon,A$ as in (\ref{UBDeq2b}). For $t$ large enough
the second term
\[
I_{1n}(t)\leq\int_{1\leq\left\vert \lambda\right\vert \leq\gamma^{1-\zeta}_{t}\mathcal{N}}\left\vert \lambda\right\vert \prod^{t-1}_{k=0}\left\vert \psi\left(\frac{\lambda}{\gamma_{t}}a_{k}(n)\right)\right\vert d\lambda
\]
\[
\leq\int_{1\leq\left\vert \lambda\right\vert \leq\gamma^{1-\zeta}_{t}\mathcal{N}}\left\vert \lambda\right\vert \prod^{\overline{q}_{t}}_{k=\underline{q}_{t}}\left\vert \psi\left(\frac{\lambda}{\gamma_{t}}a_{k}(n)\right)\right\vert d\lambda\leq_{(1)}\int_{\left\vert \lambda\right\vert \leq\gamma^{1-\zeta}_{t}\mathcal{N}}\left\vert \lambda\right\vert \exp\left[-\eta_{\ast}\sum^{\overline{q}_{t}}_{k=\underline{q}_{t}}G\left(\left\vert \frac{\lambda}{\gamma_{t}}a_{k}(n)\right\vert \right)\right]d\lambda
\]
\[
\leq\int_{\left\vert \lambda\right\vert \leq\gamma^{1-\zeta}_{t}\mathcal{N}}\left\vert \lambda\right\vert \exp\left[-A^{-1}\eta_{\ast}\left\vert \lambda\right\vert ^{\alpha-\epsilon}\sum^{\overline{q}_{t}}_{k=\underline{q}_{t}}G\left(\frac{a_{k}(n)}{\gamma_{t}}\right)\right]d\lambda
\]
\[
\left[\mu=\left[\sum^{\overline{q}_{t}}_{k=\underline{q}_{t}}G\left(\frac{a_{k}(n)}{\gamma_{t}}\right)\right]^{\frac{1}{\alpha-\epsilon}}\lambda\right]
\]
\[
\leq\left[\inf_{n\geq n_{0},t_{0}\leq t\leq n}\sum^{\overline{q}_{t}}_{k=\underline{q}_{t}}G\left(\frac{a_{k}(n)}{\gamma_{t}}\right)\right]^{-\frac{2}{\alpha-\epsilon}}\int_{\mathbb{R}}\left\vert \mu\right\vert \exp\left[-A^{-1}\eta_{\ast}\left\vert \mu\right\vert ^{\alpha-\epsilon}\right]d\mu<_{(3)}\infty,
\]
where $\leq_{(1)}$ and $\leq_{(2)}$ follows from (\ref{UBD1c})
and (\ref{UBDeq2b}) respectively. Moreover, $<_{(3)}$ is a direct
consequence of \textbf{Lemma \ref{BoundCoefGeneral-1}}. 

Next, in view of \textbf{Remark \ref{RmkUBD3}}, 
\[
I_{2n}(t)\leq\int_{\left\vert \lambda\right\vert >\gamma^{1-\zeta}_{t}\mathcal{N}}\left\vert \lambda\right\vert \prod^{\overline{q}_{t}}_{k=\underline{q}_{t}}\left\vert \psi\left(\frac{\lambda}{\gamma_{t}}a_{k}(n)\right)\right\vert d\lambda
\]
\[
=\int_{\left\vert \lambda\right\vert >\gamma^{1-\zeta}_{t}\mathcal{N}}\left\vert \lambda\right\vert \prod^{\underline{q}_{t}+\theta-1}_{k=\underline{q}_{t}}\left\vert \psi\left(\frac{\lambda}{\gamma_{t}}a_{k}(n)\right)\right\vert \prod^{\overline{q}_{t}}_{k=\underline{q}_{t}+\theta}\left\vert \psi\left(\frac{\lambda}{\gamma_{t}}a_{k}(n)\right)\right\vert d\lambda.
\]
For $\left\vert \lambda\right\vert >\gamma^{1-\zeta}_{t}\mathcal{N}$
and $\underline{q}_{t}+\theta\leq k\leq\overline{q}_{t}$ we have
\[
\left\vert \psi\left(\frac{\lambda}{\gamma_{t}}a_{k}(n)\right)\right\vert \leq_{(4)}e^{-\eta_{\ast}\left[G\left(\frac{\lambda}{\gamma_{t}}a_{k}(n)\right)\wedge1\right]}
\]
\[
\leq_{(5)}\exp\left\{ -\eta_{\ast}\left[G\left(\frac{\gamma^{1-\zeta}_{t}\mathcal{N}}{\gamma_{t}}a_{k}(n)\right)\wedge1\right]\right\} \leq_{(6)}\exp\left\{ -\eta_{\ast}G\left(\gamma^{-\zeta}_{t}\mathcal{N}a_{k}(n)\right)\right\} ,
\]
where $\leq_{(4)}$ follows from (\ref{UBD1d}), $\leq_{(5)}$ from
the monotonicity of $G$, and $\leq_{(6)}$ from the fact that $G\left(\left|\lambda\right|\right)\leq1/2$,
for $\left|\lambda\right|\leq N_{\eta}$. In view of the above, (\ref{UBDeq2}),
and \textbf{Lemma \ref{BoundCoefGeneral-1}} we have
\[
\sup_{\left\vert \lambda\right\vert >\gamma^{1-\zeta}_{t}\mathcal{N}}\prod^{\overline{q}_{t}}_{k=\underline{q}_{t}+\theta}\left\vert \psi\left(\frac{\lambda}{\gamma_{t}}a_{k}(n)\right)\right\vert \leq\prod^{\overline{q}_{t}}_{k=\underline{q}_{t}+\theta}\exp\left\{ -\eta_{\ast}G\left(\gamma^{-\zeta}_{t}\mathcal{N}a_{k}(n)\right)\right\} 
\]
\[
\leq\exp\left\{ -\eta_{\ast}A^{-1}\gamma^{\left(\alpha-\epsilon\right)\left(1-\zeta\right)}_{t}\sum^{\overline{q}_{t}}_{k=\underline{q}_{t}+\theta}G\left(\mathcal{N}\frac{a_{k}(n)}{\gamma_{t}}\right)\right\} 
\]
\[
\leq\exp\left\{ -\eta_{\ast}C_{0}\gamma^{\left(\alpha-\epsilon\right)\left(1-\zeta\right)}_{t}\right\} =:Q_{t},
\]
for some $C_{0}\in(0,\infty)$. Hence,
\[
I_{2n}(t)\leq Q_{t}\int_{\left\vert \lambda\right\vert >\gamma^{1-\zeta}_{t}\mathcal{N}}\left\vert \lambda\right\vert \prod^{\underline{q}_{t}+\theta-1}_{k=\underline{q}_{t}}\left\vert \psi\left(\frac{\lambda}{\gamma_{t}}a_{k}(n)\right)\right\vert d\lambda
\]
\[
\leq_{(7)}Q_{t}\prod^{\underline{q}_{t}+\theta-1}_{k=\underline{q}_{t}}\left(\int_{\left\vert \lambda\right\vert >\gamma^{1-\zeta}_{t}\mathcal{N}}\left\vert \lambda\right\vert \left\vert \psi\left(\frac{\lambda}{\gamma_{t}}a_{k}(n)\right)\right\vert ^{\theta}d\lambda\right)^{1/\theta}
\]
\[
\left[\mu=\frac{\lambda}{\gamma_{t}}a_{k}(n)\right]
\]
\[
\leq Q_{t}\gamma^{2}_{t}\prod^{\underline{q}_{t}+\theta-1}_{k=\underline{q}_{t}}\left(\left\vert a_{k}(n)\right\vert ^{-2}\int_{\mathbb{R}}\left\vert \mu\right\vert \left\vert \psi\left(\mu\right)\right\vert ^{\theta}d\mu\right)^{1/\theta}
\]
\[
\leq Q_{t}\gamma^{2}_{t}\sup_{\underline{q}_{t}\leq k\leq\underline{q}_{t}+\theta-1}\left\vert a_{k}(n)\right\vert ^{-2}\int_{\mathbb{R}}\left\vert \mu\right\vert \left\vert \psi\left(\mu\right)\right\vert ^{\theta}d\mu
\]
\[
=Q_{t}\left(\inf_{\underline{q}_{t}\leq k\leq\underline{q}_{t}+\theta-1}\frac{1}{\gamma_{t}}\left\vert a_{k}(n)\right\vert \right)^{-2}\int_{\mathbb{R}}\left\vert \mu\right\vert \left\vert \psi\left(\mu\right)\right\vert ^{\theta}d\mu
\]
\[
=Q_{t}\gamma^{2}_{t}\delta^{-2}_{t}\left(\delta^{-1}_{t}\inf_{\underline{q}_{t}\leq k\leq\underline{q}_{t}+\theta-1}\eta^{-1}_{k}\eta_{k}\left\vert a_{k}(n)\right\vert \right)^{-2}\int_{\mathbb{R}}\left\vert \mu\right\vert \left\vert \psi\left(\mu\right)\right\vert ^{\theta}d\mu
\]
\[
\leq Q_{t}\gamma^{2}_{t}\delta^{-2}_{t}\left(\left[\inf_{\underline{q}_{t}\leq k\leq\underline{q}_{t}+\theta-1}\eta^{-1}_{k}\delta^{-1}_{t}\inf_{\underline{q}_{t}\leq k\leq\underline{q}_{t}+\theta-1}\eta_{k}\left\vert a_{k}(n)\right\vert \right]\right)^{-2}\int_{\mathbb{R}}\left\vert \mu\right\vert \left\vert \psi\left(\mu\right)\right\vert ^{\theta}d\mu
\]
\[
\leq_{(5)}\mathcal{D}^{-2}_{0}Q_{t}\gamma^{2}_{t}\delta^{-2}_{t}\left(\left[\inf_{\underline{q}_{t}\leq k\leq\underline{q}_{t}+\theta-1}\eta^{-1}_{k}\right]\right)^{-2}\int_{\mathbb{R}}\left\vert \mu\right\vert \left\vert \psi\left(\mu\right)\right\vert ^{\theta}d\mu
\]
\[
\leq\mathcal{D}^{-2}_{0}Q_{t}\gamma^{2}_{t}\delta^{-2}_{t}\eta^{2}_{\underline{q}_{t}+\theta}\int_{\mathbb{R}}\left\vert \mu\right\vert \left\vert \psi\left(\mu\right)\right\vert ^{\theta}d\mu<C<\infty,
\]
where $\leq_{(7)}$ is a consequence of H$\ddot{\textrm{o}}$lder's
inequality (see e.g. Jeganathan, 2008; eq. (35)), $\leq_{(5)}$ follows
from \textbf{Lemma \ref{BoundCoefGeneral}}. Finally, the last inequality
follows from the fact that $Q_{t}\gamma^{2}_{t}\delta^{-2}_{t}\eta^{2}_{\underline{q}_{t}+\theta}$
is real valued sequence such that $Q_{t}\gamma^{2}_{t}\delta^{-2}_{t}\eta^{2}_{\underline{q}_{t}+\theta}\rightarrow0$
as $t\rightarrow\infty$ -cf. \textbf{Lemma \ref{BoundCoefGeneral-1}}.
$\blacksquare$

\bigskip{}

\textbf{Proof of Lemma \ref{BoundCoefGeneral}}. 

(i) We start with the proof of (\ref{bracket1}). We shall show the
result under \textbf{MI-NM} (i.e. $1<m<1+1/\alpha$). Without loss
of generality set $m_{0}=1$. For the other cases see DK2018 (Lemma
6.2). Recall that for $0\leq k\leq t-1$, we have $a_{k,t}=a_{k}$.
Set $0\leq k\leq t-1$. For some $\delta>0$, that will be specified
in detail later, set $\underline{q}_{t}=\lfloor\delta\kappa_{t}/2\rfloor$
and $\overline{q}_{t}=\lfloor\delta\kappa_{t}\rfloor$. Using the
fact that $m_{t}/\gamma_{t}\sim\kappa^{m-1}_{t}$ -cf. (\ref{m_n(1)})-
and the mean value theorem (see e.g. DK2018; Lemma 6.4(i)), as $t(\leq n)\rightarrow\infty$
we have uniformly in $\underline{q}_{t}\leq k\leq\overline{q}_{t}$,
\[
\frac{m_{t}}{\gamma_{t}}\left\vert a_{k}-\sum^{k}_{j=0}e^{-\frac{k-j}{\kappa_{n}}}c_{j}\right\vert =\frac{m_{t}}{\gamma_{t}}\left\vert \sum^{k}_{j=0}\rho^{k-j}_{n}c_{j}-\sum^{k}_{j=0}e^{-\frac{k-j}{\kappa_{n}}}c_{j}\right\vert 
\]
\[
\leq C\kappa^{m-1}_{t}\sum^{k}_{j=0}\left\vert \rho^{k-j}_{n}-e^{-\frac{k-j}{\kappa_{n}}}\right\vert \left\vert c_{j}\right\vert \leq\frac{2C}{^{\kappa^{2}_{n}}}\kappa^{m-1}_{t}\sum^{k}_{j=0}e^{-\frac{k-j}{\kappa_{n}}}\left(k-j\right)\left\vert c_{j}\right\vert 
\]
\[
\leq\frac{2C}{^{\kappa^{2}_{n}}}\kappa^{m-1}_{t}\sup_{1\leq l\leq n}\left(e^{-\frac{l}{\kappa_{n}}}l\right)\sum^{\infty}_{j=0}\left\vert c_{j}\right\vert \leq\frac{2Ce^{-1}}{^{\kappa_{n}}}\kappa^{m-1}_{t}\sum^{\infty}_{j=0}\left\vert c_{j}\right\vert 
\]
\begin{equation}
\leq\frac{2C}{^{\kappa_{n}}}\kappa^{m-1}_{n}\sum^{\infty}_{j=0}\left\vert c_{j}\right\vert =o(1).\label{BoundCoeffMI_NM1}
\end{equation}
Next, in view of (\ref{BoundCoeffMI_NM1}) as $t\rightarrow\infty$

\[
\frac{m_{t}a_{k}(n)}{\gamma_{t}}=\frac{m_{t}\sum^{k}_{j=0}\rho^{k-j}_{n}c_{j}}{\gamma_{t}}
\]
\begin{equation}
\sim\kappa^{m-1}_{t}\sum^{k}_{j=0}\rho^{k-j}_{n}c_{j}=\kappa^{m-1}_{t}\sum^{k}_{j=0}e^{-\frac{k-j}{\kappa_{n}}}c_{j}+o(1),\label{BoundCoeffMI_NM1b}
\end{equation}
uniformly in $\underline{q}_{t}\leq k\leq\overline{q}_{t}$. Fix $\epsilon>0$
arbitrarily small and let $\left\vert c_{j}/j^{-m}-1\right\vert <\epsilon$
for all $j\geq N_{\epsilon}$. Note also that $e^{x}-1\leq xe^{x}$.
Therefore, for $\underline{q}_{t}\leq k\leq\overline{q}_{t}$, as
$t\left(\leq n\right)\rightarrow\infty$, the difference
\[
\left\vert R_{nk}\right\vert :=\left\vert \sum^{k}_{j=0}e^{-\frac{k-j}{\kappa_{n}}}c_{j}-\sum^{k}_{j=0}e^{-\frac{k}{\kappa_{n}}}c_{j}\right\vert \leq e^{-\frac{k}{\kappa_{n}}}\left\vert \sum^{k}_{j=0}e^{\frac{j}{\kappa_{n}}}c_{j}-\sum^{k}_{j=0}c_{j}\right\vert 
\]
\[
\leq Ce^{-\frac{k}{\kappa_{n}}}\sum^{k}_{j=0}\frac{j}{\kappa_{n}}\left\vert c_{j}\right\vert \leq Ce^{-\frac{k}{\kappa_{n}}}\sum^{k}_{j=1}\frac{j}{\kappa_{n}}j^{-m}+Ce^{-\frac{k}{\kappa_{n}}}\sum^{k}_{j=1}\left\vert \frac{\left\vert c_{j}\right\vert }{j^{-m}}-1\right\vert \frac{j}{\kappa_{n}}j^{-m}
\]
\[
\leq C\kappa^{-1}_{n}e^{-\frac{k}{\kappa_{n}}}\left[\sum^{k}_{j=1}j^{1-m}+\epsilon\sum^{k}_{j=N_{\epsilon}}j^{1-m}+N_{\epsilon}\sup_{0\leq j\leq N_{\epsilon}}\left\vert jc_{j}\right\vert \right]
\]
\[
\leq C_{0}\kappa^{-1}_{n}e^{-\frac{k}{\kappa_{n}}}\sum^{k}_{j=1}j^{1-m}
\]
\[
=C_{0}\kappa^{-1}_{n}e^{-\frac{k}{\kappa_{n}}}k^{2-m}\int^{1}_{0}x^{1-m}dx\left[1+o(1)\right]+o(1)
\]
\begin{equation}
\leq C_{1}\kappa^{-1}_{n}e^{-\frac{k}{\kappa_{n}}}k^{2-m}\left[1+o(1)\right],\label{BoundCoeffMI_NM2}
\end{equation}
where $C_{0},C_{1}\in[1,\infty)$. Set $C_{\infty}:=\sup_{n\geq1,1\leq t\leq n}\frac{\kappa_{t}}{\kappa_{n}}\geq1$
and note that by assumption $C_{\infty}<\infty$. Hence, $1\leq C_{\infty}<\infty$.
Choose $\delta:=1\wedge\frac{1}{8\left(m-1\right)C_{1}C_{\infty}}$.
In view of (\ref{BoundCoeffMI_NM1})-(\ref{BoundCoeffMI_NM2}), for\,\underline{$%
q$}$_{t}\leq k\leq\overline{q}_{t}$ and $t\left(\leq n\right)\rightarrow\infty$
\[
\left\vert \frac{m_{t}a_{k}(n)}{\gamma_{t}}\right\vert \geq\kappa^{m-1}_{t}e^{-\frac{k}{\kappa_{n}}}\left\vert \sum^{k}_{j=0}c_{j}\right\vert -\kappa^{m-1}_{t}\left\vert R_{nk}\right\vert +o(1)
\]
\[
=\kappa^{m-1}_{t}e^{-\frac{k}{\kappa_{n}}}k^{1-m}\int^{\infty}_{1}x^{-m}dx\left[1+o(1)\right]-\kappa^{m-1}_{t}\left\vert R_{nk}\right\vert +o(1)
\]
\[
\geq\kappa^{m-1}_{t}e^{-\frac{k}{\kappa_{n}}}k^{1-m}\int^{\infty}_{1}x^{-m}dx\left[1+o(1)\right]-C_{1}\kappa^{m-1}_{t}\kappa^{-1}_{n}e^{-\frac{k}{\kappa_{n}}}k^{2-m}\left[1+o(1)\right]
\]
\[
=\kappa^{m-1}_{t}e^{-\frac{k}{\kappa_{n}}}k^{1-m}\left[\frac{1}{m-1}\left[1+o(1)\right]-C_{1}k\kappa^{-1}_{n}\left[1+o(1)\right]\right]+o(1)
\]
\[
\geq\kappa^{m-1}_{t}e^{-\frac{k}{\kappa_{n}}}k^{1-m}\left[\frac{1}{2\left(m-1\right)}-2C_{1}\lfloor\delta\kappa_{t}\rfloor\kappa^{-1}_{n}\right]+o(1)
\]
\[
\geq\kappa^{m-1}_{t}e^{-\frac{k}{\kappa_{n}}}k^{1-m}\left[\frac{1}{2\left(m-1\right)}-2C_{1}C_{\infty}\delta\right]+o(1)
\]
\[
\geq\left(\frac{\kappa_{t}}{k}\right)^{m-1}e^{-\frac{k}{\kappa_{n}}}\frac{1}{4\left(m-1\right)}+o(1)
\]
\[
\geq\left(\frac{\kappa_{t}}{\delta\kappa_{t}}\right)^{m-1}e^{-C_{\infty}}\frac{1}{4\left(m-1\right)}+o(1)\geq\delta^{1-m}e^{-C_{\infty}}\frac{1}{8\left(m-1\right)}>0.
\]

Similarly, using (\ref{BoundCoeffMI_NM1})-(\ref{BoundCoeffMI_NM2})
for $\underline{q}_{t}\leq k\leq\overline{q}_{t}$ as $t\rightarrow\infty$
\[
\left\vert \frac{m_{t}a_{k}(n)}{\gamma_{t}}\right\vert =\left\vert \kappa^{m-1}_{t}\sum^{k}_{j=0}e^{-\frac{k-j}{\kappa_{n}}}c_{j}\right\vert +o(1)
\]
\[
\leq C\left[\kappa^{m-1}_{t}e^{-\frac{k}{\kappa_{n}}}\left\vert \sum^{k}_{j=0}c_{j}\right\vert +\kappa^{m-1}_{t}\left\vert R_{nk}\right\vert \right]+o(1)
\]
\[
=C\left[\kappa^{m-1}_{t}e^{-\frac{k}{\kappa_{n}}}k^{1-m}\int^{\infty}_{1}x^{-m}dx\left(1+o(1)\right)+\kappa^{m-1}_{t}\left\vert R_{nk}\right\vert \right]+o(1)
\]
\[
\leq C\kappa^{m-1}_{t}e^{-\frac{k}{\kappa_{n}}}k^{1-m}\left[\frac{2}{m-1}+C_{1}k\kappa^{-1}_{n}\right]+o(1)
\]
\[
\leq C\left(\frac{\kappa_{t}}{\delta\kappa_{t}/2}\right)^{m-1}\left[\frac{2}{m-1}+C_{1}C_{\infty}\delta\right]+o(1)
\]
\[
\leq2C\left(\frac{2}{\delta}\right)^{m-1}\left[\frac{2}{m-1}+C_{1}C_{\infty}\delta\right]<\infty,
\]
as required.\bigskip{}

(ii) Next, we show that (\ref{bracket2-1}) holds under \textbf{BFR1},\textbf{
BFR2} and \textbf{MI}. We start with the\textbf{ BFR1} case.
\begin{itemize}
\item[\textbf{BFR1}] As $t\rightarrow\infty$ 
\[
\gamma^{-\zeta}_{t}\sup_{0\leq k\leq\overline{q}_{t}}\left\vert a_{k}(n)\right\vert =\gamma^{-\zeta}_{t}\sup_{0\leq k\leq\overline{q}_{t}}\left\vert a_{k}\right\vert 
\]
\[
=\gamma^{-\zeta}_{t}\sup_{0\leq k\leq\overline{q}_{t}}\left\vert \sum^{k}_{j=0}\phi_{j}c_{k-j}\right\vert \leq\gamma^{-\zeta}_{t}\sum^{\infty}_{j=0}\left\vert c_{j}\right\vert \rightarrow0,
\]
for all $\zeta\in(0,1)$. Note that we have used the fact that $\phi_{j}=1$. 
\item[\textbf{BFR2}]  By \textbf{Lemma \ref{cauchy},} as $k\rightarrow\infty$ we have
$\sum^{k}_{j=0}\phi_{j}c_{k-j}\sim k^{-1/\alpha}\sum^{\infty}_{j=0}c_{j}$.
Hence,
\[
\gamma^{-\zeta}_{t}\sup_{0\leq k\leq\overline{q}_{t}}\left\vert a_{k}(n)\right\vert =\gamma^{-\zeta}_{t}\sup_{0\leq k\leq\overline{q}_{t}}\left\vert a_{k}\right\vert 
\]
\[
\leq\gamma^{-\zeta}_{t}\sup_{k\geq0}\left\vert \sum^{k}_{j=0}\phi_{j}c_{k-j}\right\vert \leq C\gamma^{-\zeta}_{t}\rightarrow0,
\]
\end{itemize}
\bigskip{}

\begin{itemize}
\item[\textbf{MI}] Let $n$ large enough, $n\geq n_{0}$ say, such that $\left\vert \rho_{n}\right\vert <1$.
First, suppose that either \textbf{Assumptions MI-SM }or\textbf{ MI-NM}
hold. Noting that $\phi_{j}=\rho^{j}_{n}$, as $t\rightarrow\infty$
we have 
\[
\gamma^{-\zeta}_{t}\sup_{0\leq k\leq\overline{q}_{t},n\geq n_{0}}\left\vert a_{k}(n)\right\vert =\gamma^{-\zeta}_{t}\sup_{0\leq k\leq\overline{q}_{t},n\geq n_{0}}\left\vert \sum^{k}_{j=0}\phi_{j}c_{k-j}\right\vert \leq\gamma^{-\zeta}_{t}\sup_{j\geq0,n\geq n_{0}}\left\vert \rho^{j}_{n}\right\vert \sum^{\infty}_{j=0}\left\vert c_{j}\right\vert \leq C\gamma^{-\zeta}_{t}\rightarrow0,
\]
for all $\zeta\in(0,1)$.

Finally, suppose that \textbf{Assumption MI-LM} holds. Recall that
in this case $|c_{j}|\asymp j^{-m}$ with $m\in(1/\alpha,1)$. Moreover,
note that as $k\rightarrow\infty$, $\frac{1}{k}\sum^{k}_{j=1}\left(\frac{j}{k}\right)^{-m}\rightarrow\frac{1}{1-m}$.\footnote{Notice that 
\[
\frac{1}{k}\sum^{k}_{j=1}\left(\frac{j}{k}\right)^{-m}=\int^{1}_{1/k}x^{-m}dx+o(1)=\frac{1}{1-m}+o(1).
\]
} Therefore, the aforementioned term is bounded. Write
\[
\gamma^{-\zeta}_{t}\sup_{0\leq k\leq\overline{q}_{t},n\geq n_{0}}\left\vert a_{k}(n)\right\vert =\gamma^{-\zeta}_{t}\sup_{0\leq k\leq\overline{q}_{t},n\geq n_{0}}\left\vert \sum^{k}_{j=0}\rho^{k-j}_{n}c_{j}\right\vert \leq C\gamma^{-\zeta}_{t}\left(\left\vert c_{0}\right\vert +\sup_{1\leq k\leq\overline{q}_{t}}\sum^{k}_{j=1}\left\vert c_{j}\right\vert \right)
\]
\[
\leq C\gamma^{-\zeta}_{t}\left(1+\sup_{1\leq k\leq\overline{q}_{t}}\sum^{k}_{j=1}j^{-m}\right)\leq C\gamma^{-\zeta}_{t}\sup_{1\leq k\leq\overline{q}_{t}}k^{1-m}\frac{1}{k}\sum^{k}_{j=1}\left(\frac{j}{k}\right)^{-m}+o(1)
\]
\[
\leq C\gamma^{-\zeta}_{t}\overline{q}^{1-m}_{t}+o(1)=\frac{C\kappa^{1-m}_{t}}{\ell_{1}\left(\kappa_{t}\right)^{\zeta/\alpha}\kappa^{\zeta b/\alpha}_{t}}+o(1)
\]
\[
=\frac{C\kappa^{1-m}_{t}}{\ell_{1}\left(\kappa_{t}\right)^{\zeta/\alpha}\kappa^{\zeta\frac{\left(1+\alpha-\alpha m\right)}{\alpha}}_{t}}+o(1)=\frac{C\kappa^{1-m}_{t}}{\ell_{1}\left(\kappa_{t}\right)^{\zeta/\alpha}\kappa^{\zeta\left(\frac{1}{\alpha}+1-m\right)}_{t}}\rightarrow0,
\]
for all $\zeta\in\left(\frac{1-m}{\frac{1}{\alpha}+1-m},1\right)$.
Notice that for such \textbf{$\zeta$} the exponent \textbf{$1-m-\zeta\left(\frac{1}{\alpha}+1-m\right)$
}is negative; since \textbf{$\ell_{1}(\kappa_{t})^{-\zeta/\alpha}$
}is slowly varying, the negative power of $\kappa_{t}$ dominates
and the preceding term converges to zero. $\blacksquare$
\end{itemize}
\bigskip{}

\textbf{Proof of Lemma \ref{BoundCoefGeneral-1}}. 

\textbf{(i)} Without loss of generality assume that $G(.)$ is increasing.
In addition, set $C_{0}=1$. Let $n_{0}$ and $t_{0}$ as in \textbf{Lemma
\ref{BoundCoefGeneral}}. Using the aforementioned lemma, for all
$n\geq n_{0}$ and $t_{0}\leq t\leq n$
\[
\sum^{\overline{q}_{t}}_{k=\underline{q}_{t}+\theta}G\left(\frac{a_{k}(n)}{\gamma_{t}}\right)=\sum^{\overline{q}_{t}}_{k=\underline{q}_{t}+\theta}G\left(\frac{\delta^{-1}_{t}\eta_{k}a_{k}(n)}{\delta^{-1}_{t}\eta_{k}\gamma_{t}}\right)\geq\sum^{\overline{q}_{t}}_{k=\underline{q}_{t}+\theta}G\left(\frac{\mathcal{D}_{0}}{\delta^{-1}_{t}\eta_{k}\gamma_{t}}\right).
\]

Next, we show that the last term above is bounded away from zero.
\bigskip{}

\textbf{BFR}. Suppose that \textbf{Assumption BFR} holds (i.e. $\delta_{t}=1$,
$\eta_{t}=t^{1/\alpha}$, $\overline{q}_{t}=t-1$, $\underline{q}_{t}=\left\lfloor \gamma_{t}\right\rfloor $).
As $t\rightarrow\infty$
\[
\sum^{\overline{q}_{t}}_{k=\underline{q}_{t}+\theta}G\left(\frac{1}{\delta^{-1}_{t}\eta_{k}\gamma_{t}}\mathcal{D}_{0}\right)=\sum^{t-1}_{k=\gamma_{t}+\theta}G\left(\mathcal{D}_{0}\frac{k^{-1/\alpha}}{\gamma_{t}}\right)
\]
\[
=\sum^{t-1}_{k=\gamma_{t}+\theta}\left\vert \mathcal{D}_{0}\frac{k^{-1/\alpha}}{\gamma_{t}}\right\vert ^{\alpha}H\left(\left\vert \mathcal{D}_{0}\frac{k^{-1/\alpha}}{\gamma_{t}}\right\vert ^{-1}\right)\geq\mathcal{D}^{\alpha}_{0}\gamma^{-\alpha}_{t}\sum^{t-1}_{k=2\gamma_{t}}k^{-1}H\left(\left\vert \mathcal{D}_{0}\frac{k^{-1/\alpha}}{\gamma_{t}}\right\vert ^{-1}\right)
\]
\[
=_{(1)}\mathcal{D}^{\alpha}_{0}\gamma^{-\alpha}_{t}\alpha\int^{\mathcal{D}^{-1}_{0}\gamma_{t}(t-1)^{1/\alpha}}_{2^{1/\alpha}\mathcal{D}^{-1}_{0}\gamma^{1+1/\alpha}_{t}}\frac{H\left(x\right)}{x}dx\left[1+o(1)\right]
\]
\[
=\alpha\mathcal{D}^{\alpha}_{0}\gamma^{-\alpha}_{t}\left[L\left(\mathcal{D}^{-1}_{0}\gamma_{t}(t-1)^{1/\alpha}\right)-L\left(2^{1/\alpha}\mathcal{D}^{-1}_{0}\gamma^{1+1/\alpha}_{t}\right)\right]\left[1+o(1)\right]=\alpha\mathcal{D}^{\alpha}_{0}+o(1),
\]
where $=_{(1)}$ follows from \textbf{Lemma \ref{LemmaSV1}}, and
the last approximation is due to the facts that $L(.)$ is SV, $L\left(\gamma^{1+1/\alpha}_{t}\right)$
is SV in $\gamma_{t}$ (see BGT 1987; Proposition 1.3.6, Proposition
1.5.7), and $\gamma^{\alpha}_{t}\sim L\left(\gamma_{t}t^{1/\alpha}\right)$.\textbf{
}Since $(t-1)^{1/\alpha}\sim t^{1/\alpha}$ and $L$ is slowly varying,
the fixed multiplicative constants and the replacement of $t^{1/\alpha}$
by $(t-1)^{1/\alpha}$ do not affect the preceding asymptotic equivalence.
The result follows easily from the observation that $\alpha\mathcal{D}^{\alpha}_{0}$
is bounded away from zero. 

\bigskip{}

\textbf{MI}. Suppose that \textbf{Assumption MI} holds (i.e. $\delta_{t}=\kappa^{1-m}_{t}$,
$\eta_{t}=1$, $\overline{q}_{t}=\left\lfloor \delta\kappa_{t}\right\rfloor $,
$\underline{q}_{t}=\left\lfloor \delta\kappa_{t}/2\right\rfloor $
for some $0<\delta\leq1$). Recall that $q_{t}=\overline{q}_{t}-\underline{q}_{t}$.
In view of (\ref{m_n(1b)}), as $t\rightarrow\infty$
\[
\sum^{\overline{q}_{t}}_{k=\underline{q}_{t}+\theta}G\left(\frac{1}{\delta^{-1}_{t}\eta_{k}\gamma_{t}}\mathcal{D}_{0}\right)=\sum^{\overline{q}_{t}}_{k=\underline{q}_{t}+\theta}G\left(\mathcal{D}_{0}\frac{\kappa^{1-m}_{t}}{\gamma_{t}}\right)
\]
\[
=\left[\overline{q}_{t}-\left(\underline{q}_{t}+\theta\right)+1\right]G\left(\frac{\mathcal{D}_{0}}{m_{t}}\right)\sim q_{t}G\left(\frac{\mathcal{D}_{0}}{m_{t}}\right)
\]
\[
=q_{t}\left(\frac{m_{t}}{\mathcal{D}_{0}}\right)^{-\alpha}H\left(\frac{m_{t}}{\mathcal{D}_{0}}\right)\sim\mathcal{D}^{\alpha}_{0}\frac{q_{t}}{\ell_{1}\left(\kappa_{t}\right)\kappa_{t}}H\left(\ell_{1}\left(\kappa_{t}\right)^{1/\alpha}\kappa^{1/\alpha}_{t}\right)
\]
\[
\sim\mathcal{D}^{\alpha}_{0}\frac{q_{t}}{\kappa_{t}}=\delta\mathcal{D}^{\alpha}_{0}/2+o(1),
\]
where the limit term is bounded away from zero.\bigskip{}

\textbf{(ii)} We start with the \textbf{BFR} case. As $n\rightarrow\infty$
\[
\exp\left(-r_{1}\gamma^{r_{2}}_{n}\right)\gamma^{2}_{n}\delta^{-2}_{n}\eta^{2}_{\underline{q}_{n}+\theta}=\exp\left(-r_{1}\gamma^{r_{2}}_{n}\right)\gamma^{2}_{n}\left(\gamma_{n}+\theta\right)^{2/\alpha}\rightarrow0.
\]

Next, suppose that \textbf{Assumption MI} holds. Once more, using
(\ref{m_n(1b)}) as $n\rightarrow\infty$
\[
\exp\left(-r_{1}\gamma^{r_{2}}_{n}\right)\gamma^{2}_{n}\delta^{-2}_{n}\eta^{2}_{\underline{q}_{n}+\theta}=\exp\left(-r_{1}\gamma^{r_{2}}_{n}\right)\gamma^{2}_{n}\left(\kappa^{1-m}_{n}\right)^{-2}
\]
\[
=\exp\left(-r_{1}\gamma^{r_{2}}_{n}\right)m^{2}_{n}\leq C\exp\left(-r_{1}\gamma^{r_{2}}_{n}\right)\left[\ell_{1}\left(\kappa_{n}\right)^{1/\alpha}\kappa^{1/\alpha}_{n}\right]^{2}
\]
\[
=C\exp\left(-r_{1}\left(\ell_{1}\left(\kappa_{n}\right)^{1/\alpha}\kappa^{b/\alpha}_{n}\right)^{r_{2}}+2\ln\left[\ell_{1}\left(\kappa_{n}\right)^{1/\alpha}\kappa^{1/\alpha}_{n}\right]\right)\rightarrow0,
\]
as required. $\blacksquare$\bigskip{}

\section{Proofs for Section \ref{sec:Limit-Theory-SNPs}}

\bigskip{}

\subsection{Supporting Results for Section \ref{sec:Limit-Theory-SNPs}}

First, we provide some preliminary technical results. Their proofs
appear at the end of this Section.\bigskip{}

\begin{lemma}\label{SNP} \textit{Suppose that the conditions of
}\textbf{\textit{Lemma \ref{lem:SSalphaM}}}\textit{ hold. Then, as
$n\to\infty$,} 
\[
\frac{1}{n}\sum^{n}_{t=1}f\left(\frac{x_{t}}{\gamma^{*}_{n}}\right)\xrightarrow{d}\int^{1}_{0}f\bigl(\Lambda_{\alpha,d}(r)\bigr)dr.
\]
\end{lemma}

\bigskip{}

\begin{proposition}\label{BoundSmoothAppox} \emph{Let $\tilde{X}_{t,n}(n):=\tilde{X}_{t,n}:=\lambda^{-1}_{n}x_{t}:=\lambda^{-1}_{n}x_{t}(n),$
where $\lambda_{n}>0$ is some deterministic sequence, possibly $\lambda_{n}=1$
for all $n\geq1$. Suppose that the following hold:}
\begin{itemize}
\item[\foreignlanguage{american}{\textbf{(i)}}]  \emph{$f:\mathbb{R}\rightarrow\mathbb{R}$ is locally integrable
such that $\left|f(x)\right|=O$$\left(\left|x\right|^{\lambda^{\prime}}\right),$
as $\left|x\right|\rightarrow\infty$, for some $0<\lambda^{\prime}<\alpha$; }
\item[\foreignlanguage{american}{\textbf{(ii)}}] \emph{ $\sup_{n\geq1,1\leq t\leq n}\mathbf{E}\left|\tilde{X}_{t,n}\right|^{\alpha^{\prime}}<\infty$,
for all $0<\alpha^{\prime}<\alpha$;}
\item[\textbf{(iii)}]  \emph{either }
\begin{itemize}
\item[\textbf{(a)}] \emph{ $\tilde{X}_{t,t}$ has density $\tilde{D}_{n,t}(x)$ such
that $\sup_{n\geq n_{0},1\leq t\leq n}\sup_{x}\tilde{D}_{n,t}(x)<\infty,$
for some $n_{0}\geq1$; or}
\item[\textbf{(b)}] \emph{ $\tilde{X}_{t,t}$ has density $\tilde{D}_{n,t}(x)$ such
that $\sup_{n\geq n_{0},t_{0}\leq t\leq n}\sup_{x}\tilde{D}_{n,t}(x)<\infty,$
and $\frac{1}{n}\sum^{t^{\prime}}_{t=1}f(\tilde{X}_{t})=o_{p}(1)$,
for all fixed $t^{\prime}\geq1$. }
\end{itemize}
\item[\textbf{(iv)}]  $\sup_{n\geq1}\frac{1}{n}\sum^{n}_{t=1}\frac{\lambda_{n}}{\lambda_{t}}<\infty.$
\end{itemize}
\noindent\emph{Then for each $\varepsilon>0$, there is a there is
a bounded Lipschitz continuous $f_{\varepsilon}$ such that as $n\rightarrow\infty$
first and then $\varepsilon\rightarrow0$
\begin{equation}
\frac{1}{n}\sum^{n}_{t=1}f(\tilde{X}_{t,n})=\frac{1}{n}\sum^{n}_{t=1}f_{\varepsilon}(\tilde{X}_{t,n})+o_{p}(1).\label{BoundSmoothAppox1}
\end{equation}
}

\noindent\emph{Suppose that in addition to (i)-(iv), the following
holds:}
\begin{itemize}
\item[\textbf{(v)}]  \emph{$\tilde{X}_{t,n}=X_{t,n}+R_{t,n}$, where $R_{t,n}$ such
that for each $\epsilon>0$, $\lim_{n\rightarrow\infty}\sup_{1\leq t\leq n}P\left(\left|R_{t,n}\right|>\epsilon\right)=0.$}
\end{itemize}
\noindent\emph{Then for each $\varepsilon>0$, there is a there is
a bounded Lipschitz continuous $f_{\varepsilon}$ such that as $n\rightarrow\infty$
first and then $\varepsilon\rightarrow0$
\begin{equation}
\frac{1}{n}\sum^{n}_{t=1}f(\tilde{X}_{t,n})=\frac{1}{n}\sum^{n}_{t=1}f_{\varepsilon}(X_{t,n})+o_{p}(1).\label{BoundSmoothAppox2}
\end{equation}
}

\end{proposition}

\bigskip{}

Let $\left(x\right)_{+}=xI_{\left\{ x>0\right\} }$. Moreover, for
$0\leq r\leq1$ and $0<\mathcal{H<}1$ let
\[
\varLambda(r):=\varLambda_{\alpha,d}(r)=\left\{ \begin{array}{l}
\sigma_{0}\int_{\mathbb{R}}\left[\left(r-s\right)^{\mathcal{H}-\frac{1}{\alpha}}_{+}-(-s)^{\mathcal{H}-\frac{1}{\alpha}}_{+}\right]dZ_{\alpha}(s),\;\mathcal{H}\neq\frac{1}{\alpha}\\
\sigma_{0}Z_{\alpha}(r),\;\mathcal{H}=\frac{1}{\alpha}
\end{array}\right.,
\]
where $Z_{\alpha}(r)$ is an $\alpha$-stable Lévy motion with scale
parameter $1$ and skewness parameter $\beta_{*}=\beta I_{\left\{ \alpha\neq1\right\} }$,
$\left|\beta\right|\leq1$. Recall the connection between the the
Hurst index $\mathcal{H}$ and the memory parameter is $\mathcal{H}=d-1+1/\alpha$.

\bigskip{}

\begin{proposition}\label{IntegralExist}\emph{ }\textbf{\emph{(a)}}\emph{
For $0<r\leq1$, $\varLambda(r)$ has an $L_{1}$ integrable characteristic
function and possesses a density $p_{\varLambda(r)}$ such that}
\[
\int^{1}_{0}\sup_{x}\left|p_{\varLambda(r)}(x)\right|dr<\infty.
\]

\textbf{\emph{(b) }}\emph{Suppose that $f:\mathbb{R}\rightarrow\overline{\mathbb{R}}$
is locally integrable and satisfies the growth condition $\left|f(x)\right|=O\left(\left|x\right|^{\alpha'}\right)$,
$0<\alpha'<\alpha$, $\left|x\right|\rightarrow\infty$. Then 
\[
\int^{1}_{0}\left|f(\varLambda(r))\right|dr<\infty\;a.s.
\]
}

\end{proposition}

\subsection{Proof of Lemma \ref{lem:SSalphaM} and Theorem \ref{SNP-1}}

\bigskip{}
\textbf{Proof of Lemma \ref{lem:SSalphaM}. (a)} First, we present
some preliminary results and introduce some notation. Recall that
for $t=1,...,n$ we can write
\[
x_{t}=\sum^{\infty}_{k=0}a_{k,t}\xi_{t-k},
\]
where
\[
a_{k,t}:=\sum^{\left(t-1\right)\wedge k}_{j=0}\varphi_{j}c_{k-j}=\left\{ \begin{array}{c}
\sum^{k}_{j=0}\varphi_{j}c_{k-j}=:a_{k},\:k\leq t-1\\
\sum^{t-1}_{j=0}\varphi_{j}c_{k-j}=:a^{-}_{k,t},\:k>t-1
\end{array}\right.,
\]
where under \textbf{SNP} and \textbf{SNP{*}} $\varphi_{j}=1$ for
all $j$. Moreover, decompose
\[
x_{t}=x^{+}_{t}+x^{-}_{t}:=\sum^{t-1}_{k=0}a_{k}\xi_{t-k}+\sum^{\infty}_{k=t}a^{-}_{k,t}\xi_{t-k}.
\]
We will use the convention $x^{+}_{0},x^{-}_{0}=0$.

Let $m^{\ast}_{n}:=\gamma^{\ast}_{n}n^{-\left(1-m\right)}$ and note
that
\begin{equation}
m^{\ast}_{n}\sim n^{1/\alpha}\ell_{1}\left(n\right)^{1/\alpha},\text{ }H\left(m^{\ast}_{n}\right)\sim\ell_{1}\left(n\right)\text{.}\label{SNP-NM -0.1}
\end{equation}
Let
\[
\overline{\xi}_{n,k}:=\xi_{k}1\left\{ \left\vert \xi_{k}\right\vert \geq m^{\ast}_{n}\right\} -\mathbf{E}\left[\xi_{k}1\left\{ \left\vert \xi_{k}\right\vert \geq m^{\ast}_{n}\right\} \right],
\]
\[
\underline{\xi}_{n,k}:=\xi_{k}1\left\{ \left\vert \xi_{k}\right\vert <m^{\ast}_{n}\right\} -\mathbf{E}\left[\xi_{k}1\left\{ \left\vert \xi_{k}\right\vert <m^{\ast}_{n}\right\} \right].
\]
Set
\[
0<\varepsilon<\left(\frac{1}{m-1}-\alpha\right)\wedge\left(\alpha-\frac{1}{m}\right).
\]
For $\alpha<2$, we will impose the additional restriction, $\varepsilon<2-\alpha$,
and if $\alpha=2$ we will use the convention $\alpha+\varepsilon=2$.
These restrictions on ensure that for all $\delta>0$ the following
integrals are finite:
\[
\int^{1}_{0}\left(\int^{\infty}_{x}y^{-m}dy\right)^{\alpha-\varepsilon}dx,\int^{\delta}_{0}\left(\int^{\infty}_{x}y^{-m}dy\right)^{\alpha\pm\varepsilon}dx,
\]
\[
\int^{1}_{0}\left(\int^{\delta+x}_{x}y^{-m}dy\right)^{\alpha\pm\varepsilon}dx,\int^{\infty}_{0}\left[\int^{1+x}_{x}\left[y^{-m}-\left(\delta+y\right)^{-m}\right]dy\right]^{\alpha\pm\varepsilon}dx.
\]
Finally, recall that 
\[
b:=1+\alpha\left(1-m\right)\:\textrm{and\:\ensuremath{b_{\varepsilon^{\pm}}}=1+\ensuremath{\left(\alpha\pm\varepsilon\right)\left(1-m\right)}}.
\]

We provide a proof only under\textbf{ Assumption SNP{*}-NM} $(m\in(1,1+1/\alpha);\ \varphi_{j}=1)$.
\bigskip{}

\textbf{(1.a)} $x^{+}_{t}$: Consider
\[
A_{n,t}:=\mathbf{E}\left\vert \gamma^{*-1}_{n}x^{+}_{t}\right\vert ^{\alpha-\varepsilon}=\mathbf{E}\left\vert \gamma^{*-1}_{n}\sum^{t-1}_{k=0}a_{k}\xi_{t-k}\right\vert ^{\alpha-\varepsilon}
\]
\[
\leq_{\left(1\right)}C\gamma^{-\left(\alpha-\varepsilon\right)}_{n}\sum^{n}_{k=0}\left\vert a_{k}\right\vert ^{\alpha-\varepsilon}\mathbf{E}\left\vert \overline{\xi}_{n,k}\right\vert ^{\alpha-\varepsilon}+C\left[\gamma^{-\left(\alpha+\varepsilon\right)}_{n}\sum^{n}_{k=0}\left\vert a_{k}\right\vert ^{\alpha+\varepsilon}\mathbf{E}\left\vert \underline{\xi}_{n,k}\right\vert ^{\alpha+\varepsilon}\right]^{\frac{\alpha-\varepsilon}{\alpha+\varepsilon}}
\]
\begin{equation}
=:C\left(A_{1n}+A^{\frac{\alpha-\varepsilon}{\alpha+\varepsilon}}_{2n}\right),\label{SNP-NM -1.0}
\end{equation}
where $\leq_{\left(1\right)}$ follows either from the Loève or the
Bahr-Esseen inequalities. Moreover, by Astrauskas (1983; Lemma 1)
\begin{equation}
\mathbf{E}\left\vert \overline{\xi}_{n,k}\right\vert ^{\alpha-\varepsilon}\leq4\mathbf{E}\left[\left\vert \xi_{k}\right\vert ^{\alpha-\varepsilon}1\left\{ \left\vert \xi_{k}\right\vert \geq m^{\ast}_{n}\right\} \right]\asymp\left(m^{\ast}_{n}\right)^{-\varepsilon}h(m^{\ast}_{n}),\label{SNP-NM -0.01}
\end{equation}
and 
\begin{equation}
\mathbf{E}\left\vert \underline{\xi}_{n,k}\right\vert ^{\alpha+\varepsilon}\leq8\mathbf{E}\left[\left\vert \xi_{k}\right\vert ^{\alpha+\varepsilon}1\left\{ \left\vert \xi_{k}\right\vert <m^{\ast}_{n}\right\} \right]\asymp\left(m^{\ast}_{n}\right)^{\varepsilon}h(m^{\ast}_{n}).\label{SNP-NM -0.02}
\end{equation}

In view of the above, as $n\rightarrow\infty$,
\[
A_{1n}\leq C\gamma^{-\left(\alpha-\varepsilon\right)}_{n}\sum^{n}_{k=0}\left\vert a_{k}\right\vert ^{\alpha-\varepsilon}\left(m^{\ast}_{n}\right)^{-\varepsilon}h(m^{\ast}_{n})
\]
\[
\leq C\gamma^{-\left(\alpha-\varepsilon\right)}_{n}\sum^{n}_{k=0}\left\vert a_{k}\right\vert ^{\alpha-\varepsilon}\left(m^{\ast}_{n}\right)^{-\varepsilon}H(m^{\ast}_{n})\sim_{(1)}C\frac{\left[n\ell_{1}(n)\right]^{-\frac{\varepsilon}{\alpha}}\ell_{1}\left(n\right)}{\left[n^{b}\ell_{1}(n)\right]^{\frac{\alpha-\varepsilon}{\alpha}}}\sum^{n}_{k=0}\left\vert a_{k}\right\vert ^{\alpha-\varepsilon}
\]
\[
=C\frac{n^{-\frac{\varepsilon}{\alpha}}\ell_{1}\left(n\right)^{1-\frac{\varepsilon}{\alpha}}}{n^{b\frac{\alpha-\varepsilon}{\alpha}}\ell_{1}(n)^{\frac{\alpha-\varepsilon}{\alpha}}}\sum^{n}_{k=0}\left\vert a_{k}\right\vert ^{\alpha-\varepsilon}=C\frac{n^{-\frac{\varepsilon}{\alpha}}}{n^{b\frac{\alpha-\varepsilon}{\alpha}}}\sum^{n}_{k=0}\left\vert a_{k}\right\vert ^{\alpha-\varepsilon}
\]
\[
=_{\left(2\right)}C\frac{1}{n^{b_{\varepsilon^{-}}}}\sum^{n}_{k=0}\left\vert a_{k}\right\vert ^{\alpha-\varepsilon}=C\frac{1}{n^{b_{\varepsilon^{-}}}}\sum^{n}_{k=0}\left\vert \sum^{k}_{j=1}c_{j}\right\vert ^{\alpha-\varepsilon}
\]
\[
=C\frac{1}{n^{b_{\varepsilon^{-}}}}\sum^{n}_{k=0}\left\vert \sum^{\infty}_{j=k+1}c_{j}\right\vert ^{\alpha-\varepsilon}=_{(3)}C\frac{1}{n^{b_{\varepsilon^{-}}}}\sum^{n}_{k=0}\left\vert \sum^{\infty}_{j=k+1}j^{-m}\right\vert ^{\alpha-\varepsilon}\left[1+o(1)\right]
\]
\[
=C\frac{1}{n}\sum^{n}_{k=0}\left(\frac{1}{n}\sum^{\infty}_{j=k+1}\left(\frac{j}{n}\right)^{-m}\right)^{\alpha-\varepsilon}\left[1+o(1)\right]
\]
\begin{equation}
=_{(4)}C\int^{1}_{0}\left(\int^{\infty}_{x}y^{-m}dy\right)^{\alpha-\varepsilon}dx\left[1+o(1)\right]<C<\infty,\label{SNP-NM -0.2}
\end{equation}
where $\sim_{(1)}$ follows from (\ref{SNP-NM -0.1}), and $=_{\left(2\right)}$
follows from the fact that
\[
b\frac{\alpha-\varepsilon}{\alpha}=\left(1+\alpha-\alpha m\right)\left(1-\frac{\varepsilon}{\alpha}\right)=1+\alpha-\alpha m-\left(1+\alpha-\alpha m\right)\frac{\varepsilon}{\alpha}
\]
\[
=1+\alpha-\alpha m-\left(\frac{\varepsilon}{\alpha}+\varepsilon-\varepsilon m\right)=1+\alpha-\alpha m-\frac{\varepsilon}{\alpha}-\varepsilon+\varepsilon m
\]
\[
=1+\left(\alpha-\varepsilon\right)-\left(\alpha-\varepsilon\right)m-\frac{\varepsilon}{\alpha}=b_{\varepsilon^{-}}-\frac{\varepsilon}{\alpha}.
\]
Moreover, $=_{(3)}$ and $=_{(4)}$ follow, respectively, from the
fact that $c_{j}\asymp j^{-m}$ and Euler summation. Similarly, we
get
\begin{equation}
A_{2n}\leq C\frac{1}{n^{b_{\varepsilon^{+}}}}\sum^{n}_{k=0}\left\vert a_{k}\right\vert ^{\alpha+\varepsilon}=C\int^{1}_{0}\left(\int^{\infty}_{x}y^{-m}dy\right)^{\alpha+\varepsilon}dx\left[1+o(1)\right]<C<\infty.\label{SNP-NM -0.3}
\end{equation}
Therefore, using (\ref{SNP-NM -1.0})-(\ref{SNP-NM -0.3}) we get
\[
\sup_{n\geq1,1\leq t\leq n}\mathbf{E}\left\vert \gamma^{*-1}_{n}x^{+}_{t}\right\vert ^{\alpha-\varepsilon}=\sup_{n\geq1,1\leq t\leq n}A_{n,t}<\infty,
\]
as required.

\bigskip{}

\textbf{(1.b)} $x^{+}_{t}-x^{+}_{s}$: Let $0<\delta\leq1$, $t,s=1,...,n$,
$0\leq t-s\leq n\delta$ and consider 
\[
x^{+}_{t}-x^{+}_{s}=\sum^{t}_{l=s+1}a_{t-l}\xi_{l}-\sum^{s}_{l=1}\left(a_{t-l}-a_{s-l}\right)\xi_{l},
\]
where we are using the convention $\sum^{b}_{j=a}(.)=0$ for $a>b$.

Using analogous arguments to those used for (\ref{SNP-NM -0.2}) and
(\ref{SNP-NM -0.3}) as $n\rightarrow\infty$ first and then $\delta\rightarrow0$
we get
\[
A_{n,t,s}:=\mathbf{E}\left\vert \gamma^{*-1}_{n}\sum^{t}_{l=s+1}a_{t-l}\xi_{l}\right\vert ^{\alpha-\varepsilon}=\mathbf{E}\left\vert \gamma^{*-1}_{n}\sum^{t-\left(s+1\right)}_{k=0}a_{k}\xi_{t-k}\right\vert ^{\alpha-\varepsilon}
\]
\[
\leq C\gamma^{-\left(\alpha-\varepsilon\right)}_{n}\sum^{t-\left(s+1\right)}_{k=0}\left\vert a_{k}\right\vert ^{\alpha-\varepsilon}\mathbf{E}\left\vert \overline{\xi}_{n,k}\right\vert ^{\alpha-\varepsilon}+C\left[\gamma^{-\left(\alpha+\varepsilon\right)}_{n}\sum^{t-\left(s+1\right)}_{k=0}\left\vert a_{k}\right\vert ^{\alpha+\varepsilon}\mathbf{E}\left\vert \underline{\xi}_{n,k}\right\vert ^{\alpha+\varepsilon}\right]^{\frac{\alpha-\varepsilon}{\alpha+\varepsilon}}
\]
\[
\leq C\left\{ \int^{\delta}_{0}\left(\int^{\infty}_{x}y^{-m}dy\right)^{\alpha-\varepsilon}dx+\left[\int^{\delta}_{0}\left(\int^{\infty}_{x}y^{-m}dy\right)^{\alpha+\varepsilon}dx\right]^{\frac{\alpha-\varepsilon}{\alpha+\varepsilon}}\right\} \left[1+o(1)\right]=o(1).
\]

Next, consider
\[
B_{n,t,s}:=\mathbf{E}\left\vert \gamma^{*-1}_{n}\sum^{s}_{l=1}\left(a_{t-l}-a_{s-l}\right)\xi_{l}\right\vert ^{\alpha-\varepsilon}=\mathbf{E}\left\vert \gamma^{*-1}_{n}\sum^{s-1}_{k=0}\left(a_{t-s+k}-a_{k}\right)\xi_{s-k}\right\vert ^{\alpha-\varepsilon}
\]
\[
=\mathbf{E}\left\vert \gamma^{*-1}_{n}\sum^{s-1}_{k=0}\left(a_{t-s+k}-a_{k}\right)\xi_{k}\right\vert ^{\alpha-\varepsilon}
\]
\[
\leq C\gamma^{-\left(\alpha-\varepsilon\right)}_{n}\sum^{n}_{k=0}\left\vert \left(a_{t-s+k}-a_{k}\right)\right\vert ^{\alpha-\varepsilon}\mathbf{E}\left\vert \overline{\xi}_{n,k}\right\vert ^{\alpha-\varepsilon}
\]
\[
+C\left[\gamma^{-\left(\alpha+\varepsilon\right)}_{n}\sum^{n}_{k=0}\left\vert \left(a_{t-s+k}-a_{k}\right)\right\vert ^{\alpha+\varepsilon}\mathbf{E}\left\vert \underline{\xi}_{n,k}\right\vert ^{\alpha+\varepsilon}\right]^{\frac{\alpha-\varepsilon}{\alpha+\varepsilon}}
\]
\[
=:C\left(B_{1n,t,s}+B^{\frac{\alpha-\varepsilon}{\alpha+\varepsilon}}_{2n,t,s}\right),
\]

Once more, using (\ref{SNP-NM -0.01}) and (\ref{SNP-NM -0.02}) as
$n\rightarrow\infty$ first, and then $\delta\rightarrow0$, we have
\[
B_{1n,t,s}\leq C\gamma^{-\left(\alpha-\varepsilon\right)}_{n}\sum^{n}_{k=0}\left\vert a_{t-s+k}-a_{k}\right\vert ^{\alpha-\varepsilon}\left(m^{\ast}_{n}\right)^{-\varepsilon}h(m^{\ast}_{n})
\]
\[
\leq C\frac{1}{n^{b_{\varepsilon^{-}}}}\sum^{n}_{k=0}\left\vert a_{t-s+k}-a_{k}\right\vert ^{\alpha-\varepsilon}=C\frac{1}{n^{b_{\varepsilon^{-}}}}\sum^{n}_{k=0}\left\vert \sum^{t-s+k}_{j=0}c_{j}-\sum^{k}_{j=0}c_{j}\right\vert ^{\alpha-\varepsilon}
\]
\[
=C\frac{1}{n^{b_{\varepsilon^{-}}}}\sum^{n}_{k=0}\left\vert \sum^{t-s+k}_{j=k+1}c_{j}\right\vert ^{\alpha-\varepsilon}\leq C\frac{1}{n^{b_{\varepsilon^{-}}}}\sum^{n}_{k=0}\left(\sum^{t-s+k}_{j=k+1}j^{-m}\right)^{\alpha-\varepsilon}\left[1+o(1)\right]
\]
\[
\leq C\frac{1}{n^{b_{\varepsilon^{-}}}}\sum^{n}_{k=0}\left(\sum^{\lfloor n\delta\rfloor+k}_{j=k+1}j^{-m}\right)^{\alpha-\varepsilon}\left[1+o(1)\right]=C\frac{1}{n}\sum^{n}_{k=0}\left[\frac{1}{n}\sum^{\lfloor n\delta\rfloor+k}_{j=k+1}\left(\frac{j}{n}\right)^{-m}\right]^{\alpha-\varepsilon}\left[1+o(1)\right]
\]
\[
=C\int^{1}_{0}\left[\int^{\delta+x}_{x}y^{-m}dy\right]^{\alpha-\varepsilon}dx\left[1+o(1)\right]=o(1),
\]
where the last equality follows from dominated convergence. Similarly,
as $n\rightarrow\infty$ first and then $\delta\rightarrow0$
\[
B_{2n,t,s}\leq C\int^{1}_{0}\left[\int^{\delta+x}_{x}y^{-m}dy\right]^{\alpha+\varepsilon}dx=o(1).
\]
In view of the above, 
\[
\lim_{\delta\rightarrow0}\limsup_{n\rightarrow\infty}\sup_{\begin{array}{c}
0\leq r,s\leq1\\
0\leq r-s\leq\delta
\end{array}}\mathbf{E}\left\vert \gamma^{*-1}_{n}x^{+}_{\left\lfloor nr\right\rfloor }-\gamma^{*-1}_{n}x^{+}_{\left\lfloor ns\right\rfloor }\right\vert ^{\alpha-\varepsilon}
\]
\[
=\lim_{\delta\rightarrow0}\limsup_{n\rightarrow\infty}\sup_{\begin{array}{c}
1\leq t,s\leq n\\
0\leq t-s\leq\delta n
\end{array}}\mathbf{E}\left\vert \gamma^{*-1}_{n}x^{+}_{t}-\gamma^{*-1}_{n}x^{+}_{s}\right\vert ^{\alpha-\varepsilon}=0.
\]

\bigskip{}

\textbf{(2a)} $x^{-}_{t}$: Write
\begin{eqnarray*}
x^{-}_{t} & = & \sum^{\infty}_{k=t}a^{-}_{k,t}\xi_{t-k}=\sum^{0}_{l=-\infty}a^{-}_{t-l,t}\xi_{l}.\\
 &  & \left[l=t-k\Leftrightarrow k=t-l\right]
\end{eqnarray*}
Using arguments analogous to those used below for the increment process,
it can be shown that 
\[
\sup_{n\geq1,1\leq t\leq n}\mathbf{E}\left\vert \gamma^{*-1}_{n}x^{-}_{t}\right\vert ^{\alpha-\varepsilon}<\infty.
\]

\textbf{(2b)} $x^{-}_{t}-x^{-}_{s}$: Let $0<\delta\leq1$, $t,s=1,...,n$,
$0\leq t-s\leq n\delta$. The increment 
\begin{eqnarray*}
x^{-}_{t}-x^{-}_{s} & = & \sum^{0}_{l=-\infty}a^{-}_{t-l,t}\xi_{l}-\sum^{0}_{l=-\infty}a^{-}_{s-l,s}\xi_{l}=\sum^{0}_{l=-\infty}\left(a^{-}_{t-l,t}-a^{-}_{s-l,s}\right)\xi_{l}\\
 &  & \left[k=s-l-j\Rightarrow j=s-l-k\right]
\end{eqnarray*}
\begin{eqnarray*}
 & = & \sum^{0}_{l=-\infty}\left[\sum^{t-1}_{j=0}c_{t-l-j}-\sum^{s-1}_{j=0}c_{s-l-j}\right]\xi_{l}=\sum^{0}_{l=-\infty}\left[\sum^{s-l}_{k=s-l-\left(t-1\right)}c_{t-s+k}-\sum^{s-l}_{k=s-l-\left(s-1\right)}c_{k}\right]\xi_{l}\\{}
[r & = & -l]
\end{eqnarray*}
\[
=\sum^{0}_{l=-\infty}\left[\sum^{s-l}_{k=s-t-l+1}c_{t-s+k}-\sum^{s-l}_{k=-l+1}c_{k}\right]\xi_{l}=\sum^{\infty}_{r=0}\left[\sum^{s+r}_{k=s-t+r+1}c_{t-s+k}-\sum^{s+r}_{k=r+1}c_{k}\right]\xi_{-r}
\]
\[
=\sum^{\infty}_{r=0}\left[\sum^{s+r}_{k=-\left(t-s\right)+r+1}c_{t-s+k}-\sum^{s+r}_{k=r+1}c_{k}\right]\xi_{-r}=\sum^{\infty}_{r=0}\left[\sum^{r}_{k=-\left(t-s\right)+r+1}c_{t-s+k}+\sum^{s+r}_{k=r+1}c_{t-s+k}-\sum^{s+r}_{k=r+1}c_{k}\right]\xi_{-r}
\]
\[
=\sum^{\infty}_{r=0}\left[\sum^{r}_{k=-\left(t-s\right)+r+1}c_{t-s+k}+\sum^{s+r}_{k=r+1}c_{t-s+k}-\sum^{s+r}_{k=r+1}c_{k}\right]\xi_{-r}
\]
\[
=\sum^{\infty}_{k=0}\left[\sum^{k}_{j=-\left(t-s\right)+k+1}c_{t-s+j}\right]\xi_{-k}+\sum^{\infty}_{k=0}\left[\sum^{s+k}_{j=k+1}c_{t-s+j}-\sum^{s+k}_{j=k+1}c_{j}\right]\xi_{-k}
\]

First, consider
\[
C_{n,t,s}:=\mathbf{E}\left\vert \gamma^{*-1}_{n}\sum^{\infty}_{k=0}\left(\sum^{k}_{j=-\left(t-s\right)+k+1}c_{t-s+j}\right)\xi_{-k}\right\vert ^{\alpha-\varepsilon}
\]
\[
\leq C\gamma^{-\left(\alpha-\varepsilon\right)}_{n}\sum^{\infty}_{k=0}\left\vert \sum^{k}_{j=-\left(t-s\right)+k+1}c_{t-s+j}\right\vert ^{\alpha-\varepsilon}\mathbf{E}\left\vert \overline{\xi}_{n,k}\right\vert ^{\alpha-\varepsilon}
\]
\[
+C\left[\gamma^{-\left(\alpha+\varepsilon\right)}_{n}\sum^{\infty}_{k=0}\left\vert \sum^{k}_{j=-\left(t-s\right)+k+1}c_{t-s+j}\right\vert ^{\alpha+\varepsilon}\mathbf{E}\left\vert \underline{\xi}_{n,k}\right\vert ^{\alpha+\varepsilon}\right]^{\frac{\alpha-\varepsilon}{\alpha+\varepsilon}}
\]
\[
=:C\left(C_{1n,t,s}+C^{\frac{\alpha-\varepsilon}{\alpha+\varepsilon}}_{2n,t,s}\right).
\]
The first term above, as $n\rightarrow\infty$ first and then $\delta\rightarrow0$
\[
C_{1n,t,s}:=C\frac{1}{n^{b_{\varepsilon^{-}}}}\sum^{\infty}_{k=0}\left\vert \sum^{k}_{j=-\left(t-s\right)+k+1}c_{t-s+j}\right\vert ^{\alpha-\varepsilon}\leq C\frac{1}{n^{b_{\varepsilon^{-}}}}\sum^{\infty}_{k=0}\left\vert \sum^{k}_{j=-\left(t-s\right)+k+1}c_{t-s+j}\right\vert ^{\alpha-\varepsilon}
\]
\[
=C\frac{1}{n^{b_{\varepsilon^{-}}}}\sum^{\infty}_{k=0}\left\vert \sum^{t-s+k}_{j=t-s+\left[-\left(t-s\right)+k+1\right]}c_{j}\right\vert ^{\alpha-\varepsilon}=C\frac{1}{n^{b_{\varepsilon^{-}}}}\sum^{\infty}_{k=0}\left\vert \sum^{t-s+k}_{j=k+1}c_{j}\right\vert ^{\alpha-\varepsilon}
\]
\[
\leq C\frac{1}{n^{b_{\varepsilon^{-}}}}\sum^{\infty}_{k=0}\left(\sum^{\lfloor n\delta\rfloor+k}_{j=k+1}j^{-m}\right)^{\alpha-\varepsilon}\left[1+o(1)\right]=C\frac{1}{n}\sum^{\infty}_{k=0}\left(\frac{1}{n}\sum^{\lfloor n\delta\rfloor+k}_{j=k+1}\left(\frac{j}{n}\right)^{-m}\right)^{\alpha-\varepsilon}\left[1+o(1)\right]
\]
\[
=C\int^{\infty}_{0}\left(\int^{\delta+x}_{x}y^{-m}dy\right)^{\alpha-\varepsilon}dx=o(1).
\]
Similarly,
\[
C_{2n,t,s}\leq C\int^{\infty}_{0}\left(\int^{\delta+x}_{x}y^{-m}dy\right)^{\alpha+\varepsilon}dx=o(1).
\]
Finally, consider 
\[
D_{n,t,s}:=\mathbf{E}\left\vert \gamma^{*-1}_{n}\sum^{\infty}_{k=0}\left(\sum^{s+k}_{j=k+1}c_{t-s+j}-\sum^{s+k}_{j=k+1}c_{j}\right)\xi_{-k}\right\vert ^{\alpha-\varepsilon}
\]
\[
\leq C\gamma^{-\left(\alpha-\varepsilon\right)}_{n}\sum^{\infty}_{k=0}\left\vert \sum^{s+k}_{j=k+1}c_{t-s+j}-\sum^{s+k}_{j=k+1}c_{j}\right\vert ^{\alpha-\varepsilon}\mathbf{E}\left\vert \overline{\xi}_{n,k}\right\vert ^{\alpha-\varepsilon}
\]
\[
+C\left[\gamma^{-\left(\alpha+\varepsilon\right)}_{n}\sum^{n}_{k=0}\left\vert \sum^{s+k}_{j=k+1}c_{t-s+j}-\sum^{s+k}_{j=k+1}c_{j}\right\vert ^{\alpha+\varepsilon}\mathbf{E}\left\vert \underline{\xi}_{n,k}\right\vert ^{\alpha+\varepsilon}\right]^{\frac{\alpha-\varepsilon}{\alpha+\varepsilon}}
\]
\[
\leq C\left(D_{1n,t,s}+D^{\frac{\alpha-\varepsilon}{\alpha+\varepsilon}}_{2n,t,s}\right).
\]
As $n\rightarrow\infty$ first and then $\delta\rightarrow0$, the
first term
\[
D_{1n,t,s}\leq C\frac{1}{n^{b_{\varepsilon^{-}}}}\sum^{\infty}_{k=0}\left\vert \sum^{s+k}_{j=k+1}c_{t-s+j}-\sum^{s+k}_{j=k+1}c_{j}\right\vert ^{\alpha-\varepsilon}
\]
\[
=C\frac{1}{n^{b_{\varepsilon^{-}}}}\sum^{\infty}_{k=0}\left\vert \sum^{s+k}_{j=k+1}\left(t-s+j\right)^{-m}-\sum^{s+k}_{j=k+1}j^{-m}\right\vert ^{\alpha-\varepsilon}\left[1+o(1)\right]
\]
\[
=C\frac{1}{n^{b_{\varepsilon^{-}}}}\sum^{\infty}_{k=0}\left\vert \sum^{s+k}_{j=k+1}j^{-m}-\sum^{s+k}_{j=k+1}\left(t-s+j\right)^{-m}\right\vert ^{\alpha-\varepsilon}\left[1+o(1)\right]
\]
\[
\leq C\frac{1}{n^{b_{\varepsilon^{-}}}}\sum^{\infty}_{k=0}\left[\sum^{s+k}_{j=k+1}j^{-m}-\sum^{s+k}_{j=k+1}\left(\lfloor n\delta\rfloor+j\right)^{-m}\right]^{\alpha-\varepsilon}\left[1+o(1)\right]
\]
\[
\leq C\frac{1}{n}\sum^{\infty}_{k=0}\left[\frac{1}{n}\sum^{n+k}_{j=k+1}\left\{ \left(\frac{j}{n}\right)^{-m}-\left(\frac{\lfloor n\delta\rfloor+j}{n}\right)^{-m}\right\} \right]^{\alpha-\varepsilon}\left[1+o(1)\right]
\]
\[
\leq C\int^{\infty}_{0}\left[\int^{1+x}_{x}\left[y^{-m}-\left(\delta+y\right)^{-m}\right]dy\right]^{\alpha-\varepsilon}dx\left[1+o(1)\right]=o(1).
\]
Finally, using identical arguments as those used above we get
\[
D_{2n,t,s}\leq C\int^{\infty}_{0}\left[\int^{1+x}_{x}\left[y^{-m}-\left(\delta+y\right)^{-m}\right]dy\right]^{\alpha+\varepsilon}dx=o(1).
\]
Hence, 
\[
\lim_{\delta\rightarrow0}\limsup_{n\rightarrow\infty}\sup_{\begin{array}{c}
0\leq r,s\leq1\\
0\leq r-s\leq\delta
\end{array}}\mathbf{E}\left\vert \gamma^{*-1}_{n}x^{-}_{\left\lfloor nr\right\rfloor }-\gamma^{*-1}_{n}x^{-}_{\left\lfloor ns\right\rfloor }\right\vert ^{\alpha-\varepsilon}
\]
\[
=\lim_{\delta\rightarrow0}\limsup_{n\rightarrow\infty}\sup_{\begin{array}{c}
1\leq t,s\leq n\\
0\leq t-s\leq\delta n
\end{array}}\mathbf{E}\left\vert \gamma^{*-1}_{n}x^{-}_{t}-\gamma^{*-1}_{n}x^{-}_{s}\right\vert ^{\alpha-\varepsilon}=0.
\]
The proof of part (a) is complete.

\textbf{(b)} From part (a) we have
\begin{equation}
\lim_{\delta\rightarrow0}\limsup_{n\rightarrow\infty}\sup_{\begin{array}{c}
0\leq r,s\leq1\\
0\leq r-s\leq\delta
\end{array}}\mathbf{E}\left\vert \gamma^{*-1}_{n}x_{\left\lfloor nr\right\rfloor }-\gamma^{*-1}_{n}x_{\left\lfloor ns\right\rfloor }\right\vert ^{\alpha-\varepsilon}=0\label{SNP-NM -0.4}
\end{equation}
Set $p=\alpha-\varepsilon$, $\Delta:=\gamma^{*-1}_{n}x_{\left\lfloor nr\right\rfloor }-\gamma^{*-1}_{n}x_{\left\lfloor ns\right\rfloor }$
and $X:=\gamma^{*-1}_{n}x_{\left\lfloor ns\right\rfloor }$.

Recall that $f$ is Lipschitz and bounded i.e. $\left|f(X)\right|\leq M<\infty$.
First, assume that $p\leq1$. Then we have 
\[
\left|f(X+\Delta)-f(X)\right|\leq C_{Lip}\left|\Delta\right|\land2M\leq C^{p}_{Lip}\left(2M\right)^{1-p}\left|\Delta\right|^{p}
\]
To see this, first suppose that $C_{Lip}\left|\Delta\right|\leq2M$.
Then
\[
C_{Lip}\left|\Delta\right|\land2M=C_{Lip}\left|\Delta\right|\leq\left(C_{Lip}\left|\Delta\right|\right)^{p}\left(2M\right)^{1-p}=C^{p}_{Lip}\left(2M\right)^{1-p}\left|\Delta\right|^{p}.
\]
Next, suppose that $C_{Lip}\left|\Delta\right|>2M$. Then
\[
C_{Lip}\left|\Delta\right|\land2M=2M=\left(2M\right)^{p}\left(2M\right)^{1-p}\leq\left(C_{Lip}\left|\Delta\right|\right)^{p}\left(2M\right)^{1-p}=C^{p}_{Lip}\left(2M\right)^{1-p}\left|\Delta\right|^{p}.
\]
Hence, 
\begin{equation}
\mathbf{E}\left|f(X+\Delta)-f(X)\right|\leq C^{p}_{Lip}\left(2M\right)^{1-p}\mathbf{E}\left|\Delta\right|^{p}\label{SNP-NM -0.5}
\end{equation}
Next, suppose the $p>1$. The by Lyapunov's inequality
\begin{equation}
\mathbf{E}\left|f(X+\Delta)-f(X)\right|\leq C_{Lip}\mathbf{E}\left|\Delta\right|\leq C_{Lip}\left\{ \mathbf{E}\left|\Delta\right|^{p}\right\} ^{1/p}.\label{SNP-NM -0.6}
\end{equation}
The result follows from (\ref{SNP-NM -0.4})-(\ref{SNP-NM -0.6}).
$\blacksquare$

\bigskip{}
\textbf{Proof of Theorem \ref{SNP-1}}. First, using analogous arguments
to those used for \textbf{Lemma \ref{UBD}} -see also Jeganathan (2008;
Section 3) it can be shown that there is some $t_{0}\geq1$ such that
$\left(\gamma^{*}_{t}\right)^{-1}x_{t}$ has density $\tilde{D}_{t}(x)$
such that 
\begin{equation}
\sup_{x\geq t_{0}}\sup_{x}\tilde{D}_{t}(x)<\infty.\label{proofSNP1}
\end{equation}
If in addition \textbf{Assumption INN} holds with $\theta=1$, then
we can set $t_{0}=1$. By \textbf{Lemma \ref{lem:SSalphaM}}, 
\begin{equation}
\sup_{n\geq1}\sup_{1\leq t\leq n}\mathbf{E}\left|\frac{x_{t}}{\gamma^{*}_{n}}\right|^{\alpha^{'}}<\infty,\label{proofSNP2}
\end{equation}
for all $0<\alpha^{'}<\alpha$. Moreover, $\gamma^{*}_{n}$ is regularly
varying of index $0<\mathcal{H}<1$. Hence, standard arguments e.g.
DK2021 (proof of Theorem 3.2), show that 
\begin{equation}
\frac{1}{n}\sum^{n}_{t=1}\frac{\gamma^{*}_{n}}{\gamma^{*}_{t}}<C<\infty\label{proofSNP3}
\end{equation}
 We will only prove the result under condition (iv.a). Fix $\varepsilon>0$.
In view (\ref{proofSNP1})-(\ref{proofSNP3}) and \textbf{Proposition
\ref{BoundSmoothAppox}}, there some bounded Lipschitz continuous
$f_{\varepsilon}$ such that 
\[
\int_{\left|x\right|\leq\varepsilon^{-1}}\left|f_{\varepsilon}\left(x\right)-f\left(x\right)\right|dx<\varepsilon,
\]
 $f_{\varepsilon}(x)=0$ for $\left|x\right|>\varepsilon^{-1}$, and
as $n\rightarrow\infty$ first and then $\varepsilon\rightarrow0$\emph{
\begin{equation}
\frac{1}{n}\sum^{n}_{t=1}f\left(\frac{x_{t}}{\gamma^{*}_{n}}\right)=_{(1)}\frac{1}{n}\sum^{n}_{t=t_{0}}f\left(\frac{x_{t}}{\gamma^{*}_{n}}\right)+o_{p}(1)=\frac{1}{n}\sum^{n}_{t=1}f_{\varepsilon}\left(\frac{x_{t}}{\gamma^{*}_{n}}\right)+o_{p}(1),\label{BoundSmoothAppox2-1}
\end{equation}
}where $=_{(1)}$ follows from condition (iv.a). Next, by \textbf{Lemma
\ref{SNP}} we get 
\begin{equation}
\frac{1}{n}\sum^{n}_{t=1}f_{\varepsilon}\left(\frac{x_{t}}{\gamma^{*}_{n}}\right)\stackrel{d}{\rightarrow}\int^{1}_{0}f_{\varepsilon}\left(\varLambda(r)\right)dr.\label{BoundSmoothAppox2-1b}
\end{equation}
Next, as $\varepsilon\rightarrow0$ we get
\[
\mathbf{E}\left|\int^{1}_{0}f_{\varepsilon}\left(\varLambda(r)\right)dr-\int^{1}_{0}f\left(\varLambda(r)\right)dr\right|\leq\int^{1}_{0}\mathbf{E}\left|f_{\varepsilon}\left(\varLambda(r)\right)-f\left(\varLambda(r)\right)\right|dr
\]
\[
=\int^{1}_{0}\int_{\mathbb{R}}\left|f_{\varepsilon}\left(x\right)-f\left(x\right)\right|p_{\Lambda(r)}(x)dxdr
\]
\[
=\int^{1}_{0}\int_{\left|x\right|>\varepsilon^{-1}}\left|f_{\varepsilon}\left(x\right)-f\left(x\right)\right|p_{\Lambda(r)}(x)dxdr+\int^{1}_{0}\int_{\left|x\right|\leq\varepsilon^{-1}}\left|f_{\varepsilon}\left(x\right)-f\left(x\right)\right|p_{\Lambda(r)}(x)dxdr
\]
\[
=\int^{1}_{0}\int_{\left|x\right|>\varepsilon^{-1}}\left|f\left(x\right)\right|p_{\Lambda(r)}(x)dxdr+\int^{1}_{0}\sup_{x}p_{\Lambda(r)}(x)dr\int_{\left|x\right|\leq\varepsilon^{-1}}\left|f_{\varepsilon}\left(x\right)-f\left(x\right)\right|dx
\]
\[
\leq C\int^{1}_{0}\int_{\left|x\right|>\varepsilon^{-1}}\left|x\right|^{\alpha'}p_{\Lambda(r)}(x)dxdr+\varepsilon\int^{1}_{0}\sup_{x}p_{\Lambda(r)}(x)dr
\]
\begin{equation}
\leq_{(1)}C\int^{1}_{0}r^{\alpha\mathcal{H}}dr\int_{\left|x\right|>\varepsilon^{-1}}\left|x\right|^{\alpha'}\left|x\right|^{-1-\alpha}dx+\varepsilon\int^{1}_{0}\sup_{x}p_{\Lambda(r)}(x)dr\rightarrow0.\label{BoundSmoothAppox2-2}
\end{equation}
Combining (\ref{BoundSmoothAppox2-1}) - (\ref{BoundSmoothAppox2-2})
we get
\[
\frac{1}{n}\sum^{n}_{t=1}f\left(\frac{x_{t}}{\gamma^{*}_{n}}\right)\stackrel{d}{\rightarrow}\int^{1}_{0}f\left(\varLambda(r)\right)dr,
\]
as required.

Finally, suppose that there is $\tilde{f}:\mathbb{R}\rightarrow\mathbb{\bar{R}}$
such that $\tilde{f}=f$ a.e. and 
\[
\frac{1}{n}\sum^{t_{0}-1}_{t=1}\tilde{f}\left(\frac{x_{t}}{\gamma^{*}_{n}}\right)=o_{p}(1).
\]
Using the limit result established in the first part of this proof,
as $n\rightarrow\infty$, we get
\[
\frac{1}{n}\sum^{n}_{t=1}\tilde{f}\left(\frac{x_{t}}{\gamma^{*}_{n}}\right)=\frac{1}{n}\sum^{n}_{t=t_{0}}\tilde{f}\left(\frac{x_{t}}{\gamma^{*}_{n}}\right)+o_{p}(1)\stackrel{a.s.}{=}_{(2)}\frac{1}{n}\sum^{n}_{t=t_{0}}f\left(\frac{x_{t}}{\gamma^{*}_{n}}\right)\stackrel{d}{\rightarrow}\int^{1}_{0}f\left(\varLambda(r)\right)dr,
\]
where $\stackrel{a.s.}{=}_{(2)}$ follows from (\ref{proofSNP1}).
$\blacksquare$

\subsection{Proofs of Supporting Results for Section \ref{sec:Limit-Theory-SNPs}}

\bigskip{}
\textbf{Proof of} \textbf{Lemma \ref{SNP}.} By Astrauskas (1983)
and Kasahara and Maejima (1988), we have 
\[
\frac{1}{\gamma^{*}_{n}}x_{\left\lfloor nr\right\rfloor }\stackrel{f.d.d.}{\rightarrow}\varLambda_{\alpha,d}(r).
\]
 In view of the above and \textbf{Lemma \ref{lem:SSalphaM}}, the
result follows directly from Theorem 1 of Gikhman and Skorokhod (1969,
p. 485). $\blacksquare$

\bigskip{}
\foreignlanguage{american}{\textbf{P}}\textbf{roof of Proposition}
\textbf{\ref{BoundSmoothAppox}}. Suppose that $f:\mathbb{R}\rightarrow\mathbb{R}$
is locally integrable. It follows from Folland (1999) that for each
$\varepsilon>0$, there is a there is a Lipschitz continuous $f_{\varepsilon}$
such that
\[
\int_{\left|x\right|\leq\varepsilon^{-1}}\left|f(x)-f_{\varepsilon}(x)\right|dx<\varepsilon,\;\textrm{and}\;f_{\varepsilon}(x)=0,\;\left|x\right|>\varepsilon^{-1}.
\]
Write
\[
f(x)=f(x)1_{\left\{ \left|x\right|\leq\varepsilon^{-1}\right\} }+f(x)1_{\left\{ \left|x\right|>\varepsilon^{-1}\right\} }=:f_{1,\varepsilon}(x)+f_{2,\varepsilon}(x).
\]
Choose $\varepsilon$ small enough such that $\left|f_{2,\varepsilon}(x)\right|$$\leq C\left|x\right|^{\lambda^{\prime}}$.
Choose $\alpha^{\prime}$ such that $\lambda^{\prime}<\alpha^{\prime}<\alpha$.
By condition (iii), $\left|\tilde{X}_{t}(n)\right|^{\lambda^{\prime}}$
is uniformly integrable. Therefore, as $\varepsilon\rightarrow0$
\begin{equation}
\frac{1}{n}\mathbf{E}\sum^{n}_{t=1}\left|f_{2,\varepsilon}(\tilde{X}_{t,n}(n))\right|\leq C\sup_{n\geq1,1\leq t\leq n}\mathbf{E}\left|\tilde{X}_{t,n}(n)\right|^{\lambda^{\prime}}1_{\left\{ \left|\tilde{X}_{t,n}(n)\right|^{\lambda^{\prime}}>\varepsilon^{-\lambda^{\prime}}\right\} }=o(1).\label{BoundSmoothAppox3}
\end{equation}

Next, suppose that condition (iii.b) holds. An analogous derivation
holds under (iii.a). As $n\rightarrow\infty$ 
\begin{equation}
\frac{1}{n}\sum^{n}_{t=1}f(\tilde{X}_{t,n}(n))=\frac{1}{n}\sum^{n}_{t=t_{0}}f(\tilde{X}_{t,n}(n))+o_{p}(1).\label{BoundSmoothAppox4}
\end{equation}
In view of (\ref{BoundSmoothAppox3}) and (\ref{BoundSmoothAppox4}),
as $n\rightarrow\infty$ first and then $\epsilon\rightarrow0$, 
\[
\frac{1}{n}\sum^{n}_{t=1}f(\tilde{X}_{t,n}(n))=\frac{1}{n}\sum^{n}_{t=t_{0}}f_{1,\varepsilon}(\tilde{X}_{t,n}(n))+\frac{1}{n}\sum^{n}_{t=t_{0}}f_{2,\varepsilon}(\tilde{X}_{t,n}(n))+o_{p}(1)
\]
\[
=\frac{1}{n}\sum^{n}_{t=t_{0}}f_{1,\varepsilon}(\tilde{X}_{t,n}(n))+o_{p}(1).
\]
Next, note that $\tilde{D}_{n,t}(x)$ is the density of $\tilde{X}_{t,t}(n)=\lambda^{-1}_{t}x_{t}(n)$.
Hence,
\[
\frac{1}{n}\sum^{n}_{t=t_{0}}\mathbf{E}\left|f_{1,\varepsilon}(\tilde{X}_{t,n}(n))-f_{\varepsilon}(\tilde{X}_{t,n}(n))\right|=\frac{1}{n}\sum^{n}_{t=t_{0}}\int_{\mathbb{R}}\left|f_{1,\varepsilon}\left(\frac{\lambda_{t}}{\lambda_{n}}x\right)-f_{\varepsilon}\left(\frac{\lambda_{t}}{\lambda_{n}}x\right)\right|\tilde{D}_{n,t}(x)dx
\]
\[
\leq\sup_{n\geq n_{0},t_{0}\leq t\leq n}\sup_{x}\tilde{D}_{n,t}(x)\int_{\mathbb{R}}\left|f_{1,\varepsilon}\left(x\right)-f_{\varepsilon}\left(x\right)\right|dx\frac{1}{n}\sum^{n}_{t=1}\frac{\lambda_{n}}{\lambda_{t}}\leq C\varepsilon\rightarrow0,
\]
as required for (\ref{BoundSmoothAppox1}). 

Next, we show (\ref{BoundSmoothAppox2}). Using the fact that $f_{\varepsilon}$
is Lipschitz and bounded, as $n\rightarrow\infty$, 
\[
\frac{1}{n}\sum^{n}_{t=1}\mathbf{E}\left|f_{\varepsilon}(\tilde{X}_{t,n})-f_{\varepsilon}(X_{t,n})\right|\leq C\sup_{1\leq t\leq n}\mathbf{E}\left(1\land\left|R_{t,n}\right|\right)\rightarrow0.
\]
In view of the above and $(\ref{BoundSmoothAppox1})$, we get $(\ref{BoundSmoothAppox2})$. 

Finally, the result for $\tilde{f}:\mathbb{R}\rightarrow\overline{\mathbb{R}}$
follows easily from the fact that $\tilde{f}\left(\tilde{X}_{t,n}\right)=f\left(\tilde{X}_{t,n}\right)$
a.s. for all $1\leq t\leq n$. $\blacksquare$\bigskip{}

\noindent\textbf{Proof of Proposition \ref{IntegralExist}}. \textbf{(a)}
We only consider the case $\mathcal{H}\neq\frac{1}{\alpha}$. Moreover,
without loss of generality set $\sigma_{0}=1$. Let
\[
\mathcal{K}_{r}(s):=\left(r-s\right)^{\mathcal{H}-\frac{1}{\alpha}}_{+}-(-s)^{\mathcal{H}-\frac{1}{\alpha}}_{+}.
\]
Under our Assumptions, it follows from Proposition 3.4.1 in Samorodnitsky
and Taqqu (1994), that the characteristic function of $\varLambda(r)$
is
\[
\mathbf{E}e^{i\lambda\varLambda(r)}=\exp\left\{ -\int_{\mathbb{R}}\left|\lambda\mathcal{K}_{r}(s)\right|^{\alpha}\left\{ 1-i\tilde{\beta}\mathrm{sgn}\left[\lambda\mathcal{K}_{r}(s)\right]\right\} ds\right\} ,
\]
where 
\[
\tilde{\beta}=\beta_{*}\tan\left(\frac{\pi\alpha}{2}\right),
\]
$\beta_{*}:=\beta I_{\left\{ \alpha\neq1\right\} }$ and $\left|\beta\right|\leq1$
the skewness of the Lévy motion $Z_{\alpha}(r)$ , see e.g., Example
3.3 in in Samorodnitsky and Taqqu (1994). We will prove the result
for $\alpha\neq1$. The case $\alpha=1$ follows analogously, since
the modulus of the characteristic function has the same form required
below. 

We start with the case $\mathcal{H}\neq\frac{1}{\alpha}$. Consider
\[
\int_{\mathbb{R}}\left|\mathcal{K}_{r}(s)\right|^{\alpha}ds=\int_{\mathbb{R}}\left|\left(r-s\right)^{\mathcal{H}-\frac{1}{\alpha}}_{+}-(-s)^{\mathcal{H}-\frac{1}{\alpha}}_{+}\right|^{\alpha}ds
\]
\[
=r^{\mathcal{\alpha H}}\int_{\mathbb{R}}\left[\left(1-u\right)^{\mathcal{H}-\frac{1}{\alpha}}_{+}-(-u)^{\mathcal{H}-\frac{1}{\alpha}}_{+}\right]^{\alpha}du=C_{\mathcal{H},\alpha}r^{\mathcal{\alpha H}}
\]

where $0<C_{\mathcal{H},\alpha}<\infty$. In view of the above for
$r>0$
\[
\int_{\mathbb{R}}\left|\mathbf{E}e^{i\lambda\varLambda(r)}\right|d\lambda=\int_{\mathbb{R}}\left|\exp\left\{ -\int_{\mathbb{R}}\left|\lambda\mathcal{K}_{r}(s)\right|^{\alpha}\left\{ 1-i\tilde{\beta}\mathrm{sgn}\left[\lambda\mathcal{K}_{r}(s)\right]\right\} ds\right\} \right|d\lambda
\]
\[
\leq\int_{\mathbb{R}}\exp\left\{ -\int_{\mathbb{R}}\left|\lambda\mathcal{K}_{r}(s)\right|^{\alpha}ds\right\} d\lambda=\int_{\mathbb{R}}\exp\left\{ -C_{\mathcal{H},\alpha}r^{\mathcal{\alpha H}}\left|\lambda\right|^{\alpha}\right\} d\lambda
\]
\[
=2\int^{\infty}_{0}\exp\left\{ -C_{\mathcal{H},\alpha}r^{\mathcal{\alpha H}}\lambda^{\alpha}\right\} d\lambda
\]

\[
\left[\mu=C_{\mathcal{H},\alpha}r^{\mathcal{\alpha H}}\lambda^{\alpha}\Longleftrightarrow\lambda=\left(C_{\mathcal{H},\alpha}r^{\mathcal{\alpha H}}\right)^{-1/\alpha}\mu^{1/\alpha}\Longleftrightarrow\frac{d\lambda}{d\mu}=\frac{1}{\alpha}\left(C_{\mathcal{H},\alpha}r^{\mathcal{\alpha H}}\right)^{-1/\alpha}\mu^{1/\alpha-1}\right]
\]

\[
=\frac{2}{\alpha}\left(C_{\mathcal{H},\alpha}r^{\mathcal{\alpha H}}\right)^{-1/\alpha}\int^{\infty}_{0}\mu^{1/\alpha-1}e^{-\mu}d\mu=\frac{2}{\alpha}C_{\mathcal{H},\alpha}r^{-\mathcal{H}}\Gamma(\alpha^{-1})<\infty,
\]
 where $\Gamma$ denotes the gamma function.

Next, in view of the fact that the Hurst index is $0<\mathcal{H}<1$
we have
\[
\int^{1}_{0}r^{-\mathcal{H}}dr=\frac{1}{1-\mathcal{H}}<\infty.
\]
Hence, 
\[
\int^{1}_{0}\int_{\mathbb{R}}\left|\mathbf{E}e^{i\lambda\varLambda(r)}\right|d\lambda dr<\infty.
\]
 In view of the integrability of $\mathbf{E}e^{i\lambda\varLambda(r)}$,
it follows from Fourier inversion that $\varLambda(r)$ possesses
a density, $p_{\varLambda(r)}(x)$ say, such that
\[
p_{\varLambda(r)}(x)=\frac{1}{2\pi}\int_{\mathbb{R}}e^{-i\lambda x}\mathbf{E}e^{i\lambda\varLambda(r)}d\lambda.
\]
Hence,
\[
\int^{1}_{0}\sup_{x}\left|p_{\varLambda(r)}(x)\right|dr\leq\frac{1}{2\pi}\int^{1}_{0}\int_{\mathbb{R}}\left|\mathbf{E}e^{i\lambda\varLambda(r)}\right|d\lambda dr<\infty,
\]
as required.

Finally, for $\mathcal{H}=\frac{1}{\alpha}$ note that
\[
\left|\mathbf{E}e^{i\lambda\varLambda(r)}\right|=\exp\left\{ -\int^{r}_{0}\left|\lambda\right|^{\alpha}ds\right\} 
\]
\[
=e^{-r\left|\lambda\right|^{\alpha}}=e^{-r^{\alpha\mathcal{H}}\left|\lambda\right|^{\alpha}}.
\]
Hence, the result follows from analogous arguments as those used above. 

\bigskip{}

\textbf{(b)} We will show the result for $\mathcal{H}\neq\frac{1}{\alpha}$.
The proof for the case $\mathcal{H}=\frac{1}{\alpha}$ is analogous
to the one provided below and therefore omitted. Recall that
\[
\mathbf{E}e^{i\lambda\varLambda(r)}=\exp\left\{ -\sigma^{\alpha}_{0}\int_{\mathbb{R}}\left|\lambda\mathcal{K}_{r}(s)\right|^{\alpha}\left\{ 1-i\tilde{\beta}\mathrm{sgn}\left[\lambda\mathcal{K}_{r}(s)\right]\right\} ds\right\} .
\]
A change of variables gives 
\[
\int_{\mathbb{R}}\left|\mathcal{K}_{r}(s)\right|^{\alpha}ds=\int_{\mathbb{R}}\left|\left(r-s\right)^{\mathcal{H}-\frac{1}{\alpha}}_{+}-(-s)^{\mathcal{H}-\frac{1}{\alpha}}_{+}\right|^{\alpha}ds
\]
\[
=r^{\mathcal{\alpha H}}\int_{\mathbb{R}}\left[\left(1-u\right)^{\mathcal{H}-\frac{1}{\alpha}}_{+}-(-u)^{\mathcal{H}-\frac{1}{\alpha}}_{+}\right]^{\alpha}du
\]
\[
=r^{\mathcal{\alpha H}}\int_{\mathbb{R}}\left|\mathcal{K}_{1}(s)\right|^{\alpha}ds=:C_{\mathcal{H},\alpha}r^{\mathcal{\alpha H}}.
\]
Similarly, 
\[
\int_{\mathbb{R}}\left|\mathcal{K}_{r}(s)\right|^{\alpha}\mathrm{sgn}\left[\mathcal{K}_{r}(s)\right]ds=r^{\mathcal{\alpha H}}\int_{\mathbb{R}}\left|\mathcal{K}_{1}(s)\right|^{\alpha}\mathrm{sgn}\left[\mathcal{K}_{1}(s)\right]ds=:D_{\mathcal{H},\alpha}r^{\mathcal{\alpha H}}.
\]
In view of the above,
\[
\mathbf{E}e^{i\lambda\varLambda(r)}=\exp\left\{ -\sigma^{\alpha}_{0}\int_{\mathbb{R}}\left|\lambda\mathcal{K}_{r}(s)\right|^{\alpha}ds-i\tilde{\beta}\int_{\mathbb{R}}\left|\lambda\mathcal{K}_{r}(s)\right|^{\alpha}\mathrm{sgn}\left[\lambda\mathcal{K}_{r}(s)\right]ds\right\} 
\]
\[
=\exp\left\{ -\sigma^{\alpha}_{0}\int_{\mathbb{R}}\left|\lambda\mathcal{K}_{r}(s)\right|^{\alpha}ds\left[1-i\tilde{\beta}\frac{\int_{\mathbb{R}}\left|\lambda\mathcal{K}_{r}(s)\right|^{\alpha}\mathrm{sgn}\left[\lambda\mathcal{K}_{r}(s)\right]ds}{\int_{\mathbb{R}}\left|\lambda\mathcal{K}_{r}(s)\right|^{\alpha}ds}\right]\right\} 
\]
\[
=\exp\left\{ -\sigma^{\alpha}_{0}\int_{\mathbb{R}}\left|\mathcal{K}_{r}(s)\right|^{\alpha}ds\left|\lambda\right|^{\alpha}\left[1-i\tilde{\beta}\mathrm{sgn}\left(\lambda\right)\frac{\int_{\mathbb{R}}\left|\mathcal{K}_{r}(s)\right|^{\alpha}\mathrm{sgn}\left[\mathcal{K}_{r}(s)\right]ds}{\int_{\mathbb{R}}\left|\mathcal{K}_{r}(s)\right|^{\alpha}ds}\right]\right\} 
\]
\begin{equation}
=\exp\left\{ -\sigma^{\alpha}_{0}C_{\mathcal{H},\alpha}r^{\mathcal{\alpha H}}\left|\lambda\right|^{\alpha}\left[1-i\beta_{*}\frac{D_{\mathcal{H},\alpha}}{C_{\mathcal{H},\alpha}}\tan\left(\frac{\pi\alpha}{2}\right)\mathrm{sgn}\left(\lambda\right)\right]\right\} .\label{integralexist1}
\end{equation}
Using Ibragimov and Linnik (1971) Theorems 2.4.1-2.4.3 together with
(\ref{integralexist1}), for $0<\alpha<2$ we get 
\[
p_{\varLambda(r)}(x)\sim K_{\alpha}\sigma^{\alpha}_{0}\left(C_{\mathcal{H},\alpha}\pm\beta_{*}D_{\mathcal{H},\alpha}\right)r^{\mathcal{\alpha H}}\left|x\right|^{-1-\alpha},\quad\pm x\rightarrow\infty.
\]
Moreover, for $\alpha=2$, the density $p_{\varLambda_{\alpha}(r)}(x)$
is Gaussian. Hence, there is some $\delta_{1}>0$ and $0<C_{1}<\infty$
such that
\[
p_{\varLambda(r)}(x)\leq C_{1}\left|x\right|^{-1-\alpha}
\]
for all $\left|x\right|>\delta_{1}$. In addition, due to the assumption
$\left|f(x)\right|=O\left(\left|x\right|^{\alpha'}\right)$, $0<\alpha'<\alpha$,
there is some $\delta_{2}>0$ and $0<C_{2}<\infty$ such that
\[
\left|f(x)\right|\leq C_{2}\left|x\right|^{\alpha'},
\]
for all $\left|x\right|>\delta_{2}$. Set $C_{*}=C_{1}\lor C_{2}$
and $\delta_{*}=\delta_{1}\lor\delta_{2}$. In view of the above,
\[
\mathbf{E}\int^{1}_{0}\left|f\left(\varLambda(r)\right)\right|dr=\int_{\mathbb{R}}\int^{1}_{0}\left|f\left(x\right)\right|p_{\varLambda(r)}(x)drdx
\]
\[
=\int_{\left|x\right|\leq\delta_{*}}\int^{1}_{0}\left|f\left(x\right)\right|p_{\varLambda(r)}(x)drdx+\int_{\left|x\right|>\delta_{*}}\int^{1}_{0}\left|f\left(x\right)\right|p_{\varLambda(r)}(x)drdx
\]
\[
\leq\int_{\left|x\right|\leq\delta_{*}}\left|f\left(x\right)\right|dx\int^{1}_{0}\sup_{x}p_{\varLambda(r)}(x)dr+C_{*}\int_{\left|x\right|>\delta_{*}}\left|x\right|^{\alpha'}\left|x\right|^{-1-\alpha}dx\int^{1}_{0}dr<\infty.
\]
Hence, 
\[
\int^{1}_{0}\left|f\left(\varLambda(r)\right)\right|dr<\infty\;a.s.,
\]
as required. $\blacksquare$ 

\bigskip{}

\section{Proofs for Sections \ref{sec:Kernel-Regression}-\ref{sec:IV-Methods}}

\noindent\textbf{Proof of Theorem \ref{MainNonParam}}. In view of
\textbf{Theorem \ref{MCLT-LLker}}, the result follows by arguments
analogous to those used in the proofs of Thm 3.1 in Wang and Phillips
(2009), and Thm. 4.2 in DK2021. $\blacksquare$

\bigskip{}

\noindent\textbf{Proof of Theorem \ref{NonParam_t-test}}. In view
of \textbf{Theorem \ref{MCLT-LLker}}, the result follows by arguments
analogous to those used in the proof of Thm. 5.1 in DK2021. $\blacksquare$\bigskip{}

\noindent\textbf{Proof of Theorem \ref{MainNonParam-HTE}}. We start
with some preliminary results. Write
\[
\hat{m}(x)-m(x)=\frac{\sum^{n}_{t=1}K_{n,t-1}\left(x\right)u_{t}}{\sum^{n}_{t=1}K_{n,t-1}\left(x\right)}+\frac{\sum^{n}_{t=1}K_{n,t-1}\left(x\right)\left[m(x_{t-1})-m(x)\right]}{\sum^{n}_{t=1}K_{n,t-1}\left(x\right)}
\]
\begin{equation}
=:S_{n}+B_{n}.\label{MainNonParam-HTE-0}
\end{equation}

\noindent By \textbf{Theorem \ref{MainLL_kernel} }and \textbf{Theorem
\ref{HT-MCLT-LL}} it follows that the rescaled score term 
\[
\left(\frac{\gamma_{n}}{h_{n}n}\right)^{\frac{1}{\alpha_{u}}-1}S_{n}=\frac{\left(\frac{\gamma_{n}}{h_{n}n}\right)^{1/\alpha_{u}}\sum^{n}_{t=1}K_{n,t-1}\left(x\right)u_{t}}{\left(\frac{\gamma_{n}}{h_{n}n}\right)\sum^{n}_{t=1}K_{n,t-1}\left(x\right)}
\]
\[
\stackrel{d}{\rightarrow}\frac{1}{\varrho(0)\int_{\mathbb{R}}K(x)dx}\cdot S_{\alpha_{u}}\left\{ \sigma^{\alpha_{u}}_{u}\varrho(0)\int_{\mathbb{R}}\left\vert K(x)\right\vert ^{\alpha_{u}}dx,\beta^{\ast}_{u}\frac{\int_{\mathbb{R}}K(x)^{\left\langle \alpha_{u}\right\rangle }dx}{\int_{\mathbb{R}}\left\vert K(x)\right\vert ^{\alpha_{u}}dx}\right\} 
\]
\begin{equation}
=_{d}S_{\alpha_{u}}\left\{ \frac{\sigma^{\alpha_{u}}_{u}\varrho(0)^{1-\alpha_{u}}\int_{\mathbb{R}}\left\vert K(x)\right\vert ^{\alpha_{u}}dx}{\left(\int_{\mathbb{R}}K(x)dx\right)^{\alpha_{u}}},\beta^{\ast}_{u}\frac{\int_{\mathbb{R}}K(x)^{\left\langle \alpha_{u}\right\rangle }dx}{\int_{\mathbb{R}}\left\vert K(x)\right\vert ^{\alpha_{u}}dx}\right\} =S_{\alpha_{u}}\left\{ \tau_{1},\beta^{\ast}_{u}\tau_{2}\right\} .\label{MainNonParam-HTE-1}
\end{equation}
Next, recall that under our assumptions, $\gamma^{-1}_{t}x_{t}$ as
density function $\mathcal{D}_{n,t}(u)$ (cf. \textbf{Lemma \ref{UBD}})
such that 
\[
\sup_{n\geq1,1\leq t\leq n}\sup_{u}\mathcal{D}_{n,t}(u)<\infty
\]
and $\frac{1}{n}\sum^{n}_{t=1}\frac{\gamma_{n}}{\gamma_{t}}<C<\infty$.
For the latter see the proof of \textbf{Theorem \ref{HT-MCLT-LL}}.
In view of these facts, consider
\[
\mathbf{E}\left|\text{\ensuremath{\sum^{n}_{t=1}K_{n,t-1}\left(x\right)\left[m(x_{t-1})-m(x)\right]}}\right|\leq\sum^{n}_{t=1}\mathbf{E}\left|K\left(\frac{x_{t}-x}{h_{n}}\right)\left[m(x_{t-1})-m(x)\right]\right|
\]
\[
=\sum^{n}_{t=1}\int_{\mathbb{R}}\left|K\left(\frac{\gamma_{t}u-x}{h_{n}}\right)\left[m(\gamma_{t}u)-m(x)\right]\right|\mathcal{D}_{n,t}(u)du
\]

\noindent
\[
\left[z=\frac{\gamma_{t}u-x}{h_{n}}\Longleftrightarrow u=\gamma^{-1}_{t}\left(h_{n}z+x\right)\right]
\]
\[
=\sum^{n}_{t=1}\frac{h_{n}}{\gamma_{t}}\int_{\mathbb{R}}\left|K\left(z\right)\left[m(h_{n}z+x)-m(x)\right]\right|\mathcal{D}_{t}(\gamma^{-1}_{t}\left(h_{n}z+x\right))dz
\]
\begin{equation}
\leq Ch^{1+\nu}_{n}\int_{\mathbb{R}}\left|K\left(z\right)\right|\bar{m}(x,z)dz\sum^{n}_{t=1}\gamma^{-1}_{t}\leq Ch^{1+\nu}_{n}\sum^{n}_{t=1}\gamma^{-1}_{t}.\label{MainNonParam-HTE-2}
\end{equation}
Using \textbf{Theorem \ref{MainLL_kernel}},
\[
\left|B_{n}\right|=\left|\frac{\frac{\gamma_{n}}{nh_{n}}\sum^{n}_{t=1}K_{n,t-1}\left(x\right)\left[m(x_{t-1})-m(x)\right]}{\frac{\gamma_{n}}{nh_{n}}\sum^{n}_{t=1}K_{n,t-1}\left(x\right)}\right|
\]
\[
=\left|\frac{\frac{\gamma_{n}}{nh_{n}}\sum^{n}_{t=1}K_{n,t-1}\left(x\right)\left[m(x_{t-1})-m(x)\right]}{\varrho(0)\int_{\mathbb{R}}K(x)dx}\right|+o_{p}(1)\leq C\frac{h^{1+\nu}_{n}\gamma_{n}}{nh_{n}}\sum^{n}_{t=1}\gamma^{-1}_{t}+o_{p}(1)
\]
\begin{equation}
=Ch^{\nu}_{n}\frac{1}{n}\sum^{n}_{t=1}\frac{\gamma_{n}}{\gamma_{t}}+o_{p}(1)=O_{p}\left(h^{\nu}_{n}\right)=o_{p}(1).\label{MainNonParam-HTE-2b}
\end{equation}
\textbf{(a)} Next, we show part (a). Suppose that $\alpha_{u}>1$.
Consider the numerator of the 'bias' term $B_{n}$Combining (\ref{MainNonParam-HTE-0})-(\ref{MainNonParam-HTE-2}),
\[
\hat{m}(x)-m(x)=O_{p}\left(h^{\nu}_{n}\right)+O_{p}\left(\left(\frac{\gamma_{n}}{h_{n}n}\right)^{1-1/\alpha_{u}}\right)=o_{p}(1),
\]
as required. 

\noindent Next, we show the weak limit result of part (a). Using (\ref{MainNonParam-HTE-2}),
\[
\left(\frac{\gamma_{n}}{nh_{n}}\right)^{\frac{1}{\alpha_{u}}}\sum^{n}_{t=1}\mathbf{E}\left|K\left(\frac{x_{t}-x}{h_{n}}\right)\left[m(x_{t-1})-m(x)\right]\right|\leq C\left(\frac{\gamma_{n}}{nh_{n}}\right)^{\frac{1}{\alpha_{u}}}\sum^{n}_{t=1}\frac{h^{1+\nu}_{n}}{\gamma_{t}}
\]
\[
=Ch^{1+\nu}_{n}\left(\frac{\gamma_{n}}{nh_{n}}\right)^{\frac{1}{\alpha_{u}}}\gamma^{-1}_{n}\sum^{n}_{t=1}\frac{\gamma_{n}}{\gamma_{t}}=Ch^{1+\nu}_{n}\left(\frac{\gamma_{n}}{nh_{n}}\right)^{\frac{1}{\alpha_{u}}}\frac{n}{\gamma_{n}}\frac{1}{n}\sum^{n}_{t=1}\frac{\gamma_{n}}{\gamma_{t}}
\]
\[
=Ch^{\nu}_{n}\left(\frac{\gamma_{n}}{nh_{n}}\right)^{\frac{1}{\alpha_{u}}}\frac{nh_{n}}{\gamma_{n}}\frac{1}{n}\sum^{n}_{t=1}\frac{\gamma_{n}}{\gamma_{t}}=Ch^{\nu}_{n}\left(\frac{\gamma_{n}}{nh_{n}}\right)^{\frac{1}{\alpha_{u}}}\frac{nh_{n}}{\gamma_{n}}\frac{1}{n}\sum^{n}_{t=1}\frac{\gamma_{n}}{\gamma_{t}}\leq Ch^{\nu}_{n}\left(\frac{\gamma_{n}}{nh_{n}}\right)^{\frac{1}{\alpha_{u}}-1}.
\]
Therefore, using the upper bound shown above together with \textbf{Theorem
\ref{MCLT-LLker}} and (\ref{MainNonParam-HTE-1}) we get
\[
\left(\frac{nh_{n}}{\gamma_{n}}\right)^{1-\frac{1}{\alpha_{u}}}\left[\hat{m}(x)-m(x)\right]=\frac{\left(\frac{\gamma_{n}}{nh_{n}}\right)^{\frac{1}{\alpha_{u}}}\sum^{n}_{t=1}K_{n,t-1}\left(x\right)\left[m(x_{t-1})-m(x)\right]}{\frac{\gamma_{n}}{nh_{n}}\sum^{n}_{t=1}K_{n,t-1}\left(x\right)}
\]
\[
+\frac{\left(\frac{\gamma_{n}}{nh_{n}}\right)^{\frac{1}{\alpha_{u}}}\sum^{n}_{t=1}K_{n,t-1}\left(x\right)u_{t}}{\frac{\gamma_{n}}{nh_{n}}\sum^{n}_{t=1}K_{n,t-1}\left(x\right)}=O_{p}\left(h^{\nu}_{n}\left(\frac{\gamma_{n}}{nh_{n}}\right)^{\frac{1}{\alpha_{u}}-1}\right)+\frac{\left(\frac{\gamma_{n}}{nh_{n}}\right)^{\frac{1}{\alpha_{u}}}\sum^{n}_{t=1}K_{n,t-1}\left(x\right)u_{t}}{\frac{\gamma_{n}}{nh_{n}}\sum^{n}_{t=1}K_{n,t-1}\left(x\right)}
\]
\[
=\frac{\left(\frac{\gamma_{n}}{nh_{n}}\right)^{\frac{1}{\alpha_{u}}}\sum^{n}_{t=1}K_{n,t-1}\left(x\right)u_{t}}{\frac{\gamma_{n}}{nh_{n}}\sum^{n}_{t=1}K_{n,t-1}\left(x\right)}+o_{p}(1)\stackrel{d}{\rightarrow}S_{\alpha_{u}}\left\{ \tau_{1},\beta^{\ast}_{u}\tau_{2}\right\} ,
\]
as required.

\noindent\textbf{(b)} Suppose that $\alpha_{u}=1$. Using (\ref{MainNonParam-HTE-1})
together with (\ref{MainNonParam-HTE-2b})
\[
\hat{m}(x)-m(x)=O_{p}\left(h^{\nu}_{n}\right)+\frac{\left(\frac{\gamma_{n}}{h_{n}n}\right)\sum^{n}_{t=1}K_{n,t-1}\left(x\right)u_{t}}{\left(\frac{\gamma_{n}}{h_{n}n}\right)\sum^{n}_{t=1}K_{n,t-1}\left(x\right)}\stackrel{d}{\rightarrow}S_{\alpha_{u}}\left\{ \tau_{1},\beta^{\ast}_{u}\tau_{2}\right\} ,
\]
as required. 

\noindent\textbf{(c)} Suppose that $\alpha_{u}<1$. Using (\ref{MainNonParam-HTE-1})
together with (\ref{MainNonParam-HTE-2b})
\[
\left(\frac{nh_{n}}{\gamma_{n}}\right)^{1-\frac{1}{\alpha_{u}}}\left[\hat{m}(x)-m(x)\right]=\frac{\left(\frac{\gamma_{n}}{nh_{n}}\right)^{\frac{1}{\alpha_{u}}}\sum^{n}_{t=1}K_{n,t-1}\left(x\right)\left[m(x_{t-1})-m(x)\right]}{\frac{\gamma_{n}}{nh_{n}}\sum^{n}_{t=1}K_{n,t-1}\left(x\right)}
\]
\[
+\frac{\left(\frac{\gamma_{n}}{nh_{n}}\right)^{\frac{1}{\alpha_{u}}}\sum^{n}_{t=1}K_{n,t-1}\left(x\right)u_{t}}{\frac{\gamma_{n}}{nh_{n}}\sum^{n}_{t=1}K_{n,t-1}\left(x\right)}=O_{p}\left(h^{\nu}_{n}\left(\frac{\gamma_{n}}{nh_{n}}\right)^{\frac{1}{\alpha_{u}}-1}\right)+\frac{\left(\frac{\gamma_{n}}{nh_{n}}\right)^{\frac{1}{\alpha_{u}}}\sum^{n}_{t=1}K_{n,t-1}\left(x\right)u_{t}}{\frac{\gamma_{n}}{nh_{n}}\sum^{n}_{t=1}K_{n,t-1}\left(x\right)}
\]
\[
\stackrel{d}{\rightarrow}S_{\alpha_{u}}\left\{ \tau_{1},\beta^{\ast}_{u}\tau_{2}\right\} ,
\]
as required. $\blacksquare$\bigskip{}

\noindent\textbf{Proof of Theorem \ref{WIV_SP}}. \textbf{(a) }Recall
 that\textbf{ $x_{t}$} which is linear process of the form 
\[
x_{t}=\sum^{\infty}_{s=0}c_{s}\xi_{t-s},
\]
which exists a.s. In particular, $x_{t}$ is a causal measurable tranformation
of the i.i.d. sequence $\left\{ \xi_{j}\right\} ^{t}_{-\infty}$ .
Therefore, $x_{t}$ is strictly stationary and ergodic; see e.g. Straumann
(2005, Corollary 2.1.3). Let $g$ be a $\mathcal{B}(\mathbb{\overline{R}})$-measurable
function. Then, $\left\{ g(x_{t})\right\} ^{n}_{t=1}$ is also strictly
stationary and ergodic; see, e.g., Straumann (2005, Corollary 2.1.1).
Moreover, if $\mathbf{E}\left|g(x_{t})\right|<\infty$, then by ergodic
theorem (see, e.g., Shiryaev, 1996, p. 409-410)
\begin{equation}
\frac{1}{n}\sum^{n}_{t=1}g(x_{t})\stackrel{a.s.}{\rightarrow}\mathbf{E}g(x_{t}),\label{Proof_WIV_SP_1}
\end{equation}
as $n\rightarrow\infty$. Next, consider the WIV estimator
\begin{equation}
\sqrt{n}\left(\mathbf{\tilde{\mathbf{b}}}_{WIV}-\mathbf{b}\right)=\left[\frac{1}{n}\sum^{n}_{t=1}\mathbf{w}_{t-1}\mathbf{x}^{\prime}_{t-1}\right]^{-1}\frac{1}{\sqrt{n}}\sum^{n}_{t=1}\mathbf{w}_{t-1}u_{t}.\label{Proof_WIV_SP2}
\end{equation}
Note that $f,w_{0},w_{1},w_{0}f,w_{1}f,w^{2}_{0},w^{2}_{1}$ all satisfy
the conditions imposed on $g$; that is, they are $\mathcal{B}(\mathbb{\overline{R}})$-measurable
and have finite absolute moments. Let $\mathbf{j}\in\mathbb{R}^{2}$,
$\mathbf{j}\neq\mathbf{0}$. Consequently, 
\begin{equation}
\left[\frac{1}{n}\sum^{n}_{t=1}\mathbf{A}(x_{t}),\;\frac{1}{n}\sum^{n}_{t=1}\left(\mathbf{j}^{\prime}\mathbf{w}_{t-1}\right)^{2}\right]\stackrel{p}{\rightarrow}\left[\mathbf{E}\left(\mathbf{A}(x_{t})\right),\;\mathbf{E}\left(\mathbf{j}^{\prime}\mathbf{w}_{t-1}\right)^{2}\right],\label{Proof_WIV_SP3}
\end{equation}
and
\begin{equation}
\mathbf{E}_{t-1}\left(\frac{1}{\sqrt{n}}\sum^{n}_{t=1}\mathbf{j}^{\prime}\mathbf{w}_{t-1}u_{t}\right)^{2}=\frac{\sigma^{2}_{u}}{n}\sum^{n}_{t=1}\left(\mathbf{j}^{\prime}\mathbf{w}_{t-1}\right)^{2}\stackrel{p}{\rightarrow}\sigma^{2}_{u}\mathbf{E}\left(\mathbf{j}^{\prime}\mathbf{w}_{t-1}\right)^{2}.\label{Proof_WIV_SP4}
\end{equation}
To verify the corresponding Lindeberg condition, fix $\epsilon>0$.
Moreover, set $J_{t}:=\mathbf{j}^{\prime}\mathbf{w}_{t-1}$. Since
$\mathbf{E}\left\Vert \mathbf{w}_{t-1}\right\Vert ^{2+\delta}<\infty$,
for some $\delta>0$. we have $\mathbf{E}J^{2+\delta}_{t}<\infty$.
Therefore, dominated convergence yields 
\[
\mathbf{E}\mathbf{E}_{t-1}\left(\frac{1}{\sqrt{n}}\sum^{n}_{t=1}J_{t}u_{t}I_{\left\{ \left|J_{t}\right|>\sqrt{n}\epsilon\right\} }\right)^{2}=\frac{\sigma^{2}_{u}}{n}\mathbf{E}\sum^{n}_{t=1}J^{2}_{t}I_{\left\{ \left|J_{t}\right|>\sqrt{n}\epsilon\right\} }
\]
\begin{equation}
=\mathbf{E}J^{2}_{t}I_{\left\{ \left|J_{t}\right|>\sqrt{n}\epsilon\right\} }=o(1),\label{Proof_WIV_SP5}
\end{equation}
as $n\rightarrow\infty$. Using (\ref{Proof_WIV_SP2})-(\ref{Proof_WIV_SP5})
together with the Cramér--Wold device and Hall and Heyde (1980, Corollary
3.1) we get
\begin{equation}
\sqrt{n}\left(\mathbf{\tilde{\mathbf{b}}}_{WIV}-\mathbf{b}\right)\stackrel{d}{\rightarrow}\left[\mathbf{E}\left(\mathbf{A}(x_{t})\right)\right]^{-1}N\left[\mathbf{0},\sigma^{2}_{u}\mathbf{E}\left(\mathbf{W}(x_{t})\right)\right],\label{Proof_WIV_SP5a}
\end{equation}
as required. 

\noindent\textbf{(b)} Next, we show the consistency of 
\[
\tilde{\sigma}^{2}_{w}=\frac{\sum^{n}_{t=1}\tilde{u}^{2}_{t}w_{u}(x_{t-1})}{\sum^{n}_{t=1}w_{u}(x_{t-1})}.
\]
In view of the fact that $w_{u}(x_{t-1})$ is strictly stationary
and ergodic with a finite mean we obtain
\begin{equation}
\frac{1}{n}\sum^{n}_{t=1}w_{u}(x_{t-1})\stackrel{p}{\rightarrow}\mathbf{E}\left[w_{u}(x_{t-1})\right].\label{Proof_WIV_SP5b}
\end{equation}
Next, consider
\[
\frac{1}{n}\sum^{n}_{t=1}\tilde{u}^{2}_{t}w_{u}(x_{t-1})=\frac{1}{n}\sum^{n}_{t=1}u^{2}_{t}w_{u}(x_{t-1})
\]
\begin{equation}
+\frac{1}{n}\sum^{n}_{t=1}\left[\left(\mathbf{\tilde{\mathbf{b}}}_{WIV}-\mathbf{b}\right)^{\prime}\mathbf{x}_{t-1}\right]^{2}w_{u}(x_{t-1})-\frac{2}{n}\sum^{n}_{t=1}\left(\mathbf{\tilde{\mathbf{b}}}_{WIV}-\mathbf{b}\right)^{\prime}\mathbf{x}_{t-1}u_{t}w_{u}(x_{t-1})\label{Proof_WIV_SP6}
\end{equation}
By the Bahr-Essen inequality for martingales, see e.g. Pinelis (2015,
Thm 1.1), for some $1<p<2$, the
\[
\mathbf{E}\left|\frac{1}{n}\sum^{n}_{t=1}\left(u^{2}_{t}-\sigma^{2}_{u}\right)w_{u}(x_{t-1})\right|^{p}\leq\frac{1}{n^{p}}\sum^{n}_{t=1}\mathbf{E}\left[\left|u^{2}_{t}-\sigma^{2}_{u}\right|^{p}w_{u}(x_{t-1})^{p}\right]
\]
\begin{equation}
\leq\frac{C}{n^{p}}\sum^{n}_{t=1}\mathbf{E}w_{u}(x_{t-1})^{p}\leq\frac{C}{n^{p-1}}=o(1),\label{Proof_WIV_SP7}
\end{equation}
The second term
\[
\frac{1}{n}\sum^{n}_{t=1}\left[\left(\mathbf{\tilde{\mathbf{b}}}_{WIV}-\mathbf{b}\right)^{\prime}\mathbf{x}_{t-1}\right]^{2}w_{u}(x_{t-1})\leq\left\Vert \mathbf{\tilde{\mathbf{b}}}_{WIV}-\mathbf{b}\right\Vert ^{2}\frac{1}{n}\sum^{n}_{t=1}\left\Vert \mathbf{x}_{t-1}\right\Vert ^{2}w_{u}(x_{t-1})
\]
\begin{equation}
\leq C\left\Vert \mathbf{\tilde{\mathbf{b}}}_{WIV}-\mathbf{b}\right\Vert ^{2}\frac{1}{n}\sum^{n}_{t=1}\left(1+\left|f(x_{t-1})\right|\right)^{2}w_{u}(x_{t-1})=O_{p}(n^{-1}),\label{Proof_WIV_SP8}
\end{equation}
where the last equality above follows from the $\sqrt{n}$-consistency
of the WIV estimator and the fact that $\mathbf{E}f(x_{t-1})^{2}w_{u}(x_{t-1})<\infty$.
Similarly, the third term
\[
\frac{1}{n}\sum^{n}_{t=1}\left|\left(\mathbf{\tilde{\mathbf{b}}}_{WIV}-\mathbf{b}\right)^{\prime}\mathbf{x}_{t-1}u_{t}w_{u}(x_{t-1})\right|
\]
\[
\leq\frac{1}{n}\sum^{n}_{t=1}\left|\left(\mathbf{\tilde{\mathbf{b}}}_{WIV}-\mathbf{b}\right)^{\prime}\mathbf{x}_{t-1}u_{t}w_{u}(x_{t-1})\right|
\]
\begin{equation}
\leq\left\Vert \mathbf{\tilde{\mathbf{b}}}_{WIV}-\mathbf{b}\right\Vert \frac{1}{n}\sum^{n}_{t=1}\left|f(x_{t-1})w_{u}(x_{t-1})u_{t}\right|=O_{p}(n^{-1/2}),\label{Proof_WIV_SP9}
\end{equation}
where we have used the fact that
\[
\frac{1}{n}\sum^{n}_{t=1}\mathbf{E}\left|f(x_{t-1})w_{u}(x_{t-1})u_{t}\right|\leq\frac{1}{n}\sum^{n}_{t=1}\mathbf{E}\left|f(x_{t-1})w_{u}(x_{t-1})\right|\mathbf{E}_{t-1}\left|u_{t}\right|
\]
\[
\leq\frac{1}{n}\sum^{n}_{t=1}\mathbf{E}\left|f(x_{t-1})w_{u}(x_{t-1})\right|\left[\mathbf{E}_{t-1}u^{2}_{t}\right]^{1/2}\leq\frac{\left|\sigma_{u}\right|}{n}\sum^{n}_{t=1}\mathbf{E}\left|f(x_{t-1})w_{u}(x_{t-1})\right|=O(1).
\]
The result follows from (\ref{Proof_WIV_SP5b})-(\ref{Proof_WIV_SP9}).

\noindent Finally, the limit distribution of the t-statistic follows
easily form part (a) and the consistency of $\tilde{\sigma}^{2}_{w}$.
Indeed, set $\mathbf{e}^{\prime}_{2}:=\left[0,1\right]$ and 
\[
Q_{n}:=n^{-1}\mathbf{e}^{\prime}_{2}\left[\sum^{n}_{t=1}\mathbf{A}(x_{t-1})\right]^{-1}\sum^{n}_{t=1}\mathbf{W}(x_{t-1})\left[\sum^{n}_{t=1}\mathbf{A^{\prime}}(x_{t-1})\right]^{-1}\mathbf{e}_{2}.
\]
 Then by the ergodic theorem
\begin{equation}
Q_{n}\stackrel{p}{\rightarrow}\mathbf{e}^{\prime}_{2}\mathbf{E}\left[\mathbf{A}(x_{t-1})\right]\mathbf{E}\left[\mathbf{W}(x_{t-1})\right]\left\{ \mathbf{E}\left[\mathbf{A^{\prime}}(x_{t-1})\right]\right\} ^{-1}\mathbf{e}_{2}.\label{Proof_WIV_SP10}
\end{equation}
Combining (\ref{Proof_WIV_SP5a}), (\ref{Proof_WIV_SP10}) and the
consistency of $\tilde{\sigma}^{2}_{u}$ we get
\[
\tilde{t}_{WIV}=\frac{\sqrt{n}\mathbf{e}^{\prime}_{2}\left(\mathbf{\tilde{\mathbf{b}}}_{WIV}-\mathbf{b}\right)}{\sqrt{\tilde{\sigma}^{2}_{u}Q_{n}}}\stackrel{d}{\rightarrow}N[0,1],
\]
as required. $\blacksquare$

\bigskip{}

\noindent\textbf{Proof of Theorem \ref{MainParametric}}. \textbf{(a)}
Set $C_{n}=\sqrt{n}diag\{1,k_{f}(\gamma_{n})\}$ and $D_{n}=\sqrt{n}diag\{k_{w_{0}}(\gamma_{n}),k_{w_{1}}(\gamma_{n})\}$.
By \textbf{Theorem \ref{MCLT-LL}},
\[
C_{n}\left(\mathbf{\tilde{\mathbf{b}}}_{WIV}-\mathbf{b}\right)=C_{n}\left[\sum^{n}_{t=1}\mathbf{w}_{t-1}\mathbf{x}^{\prime}_{t-1}\right]^{-1}D_{n}D^{-1}_{n}\sum^{n}_{t=1}\mathbf{w}_{t-1}u_{t},
\]
\begin{equation}
=\left\{ D^{-1}_{n}\left[\sum^{n}_{t=1}\mathbf{w}_{t-1}\mathbf{x}^{\prime}_{t-1}\right]C^{-1}_{n}\right\} ^{-1}D_{n}\sum^{n}_{t=1}\mathbf{w}_{t-1}u_{t}=:Q^{-1}_{n}\cdot D_{n}\sum^{n}_{t=1}\mathbf{w}_{t-1}u_{t}\label{Proof_WIV_MainParam_WNP}
\end{equation}
\[
\overset{d}{\rightarrow}\left[\int_{\mathbb{R}}\mathbf{A}_{H}(x)\varrho(x)dx\right]^{-1}\cdot MN\left[\mathbf{0},\sigma^{2}_{u}\int_{\mathbb{R}}\mathbf{W}_{H}(x)\varrho(x)dx\right],
\]
as required. 

\noindent\textbf{(b)} Set $\mathbf{e}^{\prime}_{2}:=\left[0,1\right]$.
Using the limit distribution result of part (a) above, the numerator
of $\tilde{t}_{WLS}$ after scaling by $\sqrt{n}k_{f}(\gamma_{n})$
is
\[
\sqrt{n}k_{f}(\gamma_{n})\left(\mathbf{\tilde{b}}_{WIV}-\mathbf{b}\right)=\mathbf{e}^{\prime}_{2}C_{n}\left(\mathbf{\tilde{b}}_{WIV}-\mathbf{b}\right)
\]
\[
\stackrel{d}{\rightarrow}MN\left[0,\sigma^{2}_{u}\mathbf{e}^{\prime}_{2}\left[\int_{\mathbb{R}}\mathbf{A}_{H}(x)\varrho(x)dx\right]^{-1}\int_{\mathbb{R}}\mathbf{W}_{H}(x)\varrho(x)dx\left[\int_{\mathbb{R}}\mathbf{A}^{\prime}_{H}(x)\varrho(x)dx\right]^{-1}\mathbf{e}_{2}\right].
\]
In view of \textbf{Theorem \ref{MCLT-LL}}, the normalizing term
\[
nk_{f}(\gamma_{n})^{2}\mathbf{e}^{\prime}_{2}\left(\sum^{n}_{t=1}\mathbf{w}_{t-1}\mathbf{x}^{\prime}_{t-1}\right)^{-1}\left(\sum^{n}_{t=1}\mathbf{w}_{t-1}\mathbf{w}^{\prime}_{t-1}\right)\left(\sum^{n}_{t=1}\mathbf{x}_{t-1}\mathbf{w}^{\prime}_{t-1}\right)^{-1}\mathbf{e}_{2}
\]
\[
=nk_{f}(\gamma_{n})^{2}\mathbf{e}^{\prime}_{2}C^{-1}_{n}C_{n}\left(\sum^{n}_{t=1}\mathbf{w}_{t-1}\mathbf{x}^{\prime}_{t-1}\right)^{-1}D_{n}D^{-1}_{n}\left(\sum^{n}_{t=1}\mathbf{w}_{t-1}\mathbf{w}^{\prime}_{t-1}\right)D^{-1}_{n}D_{n}\left(\sum^{n}_{t=1}\mathbf{x}_{t-1}\mathbf{w}^{\prime}_{t-1}\right)^{-1}C_{n}C^{-1}_{n}\mathbf{e}_{2}
\]
\[
=\mathbf{e}^{\prime}_{2}C_{n}\left(\sum^{n}_{t=1}\mathbf{w}_{t-1}\mathbf{x}^{\prime}_{t-1}\right)^{-1}D_{n}D^{-1}_{n}\left(\sum^{n}_{t=1}\mathbf{w}_{t-1}\mathbf{w}^{\prime}_{t-1}\right)D^{-1}_{n}D_{n}\left(\sum^{n}_{t=1}\mathbf{x}_{t-1}\mathbf{w}^{\prime}_{t-1}\right)^{-1}C_{n}\mathbf{e}_{2}
\]
\[
\stackrel{d}{\rightarrow}\mathbf{e}^{\prime}_{2}\left[\int_{\mathbb{R}}\mathbf{A}_{H}(x)\varrho(x)dx\right]^{-1}\int_{\mathbb{R}}\mathbf{W}_{H}(x)\varrho(x)dx\left[\int_{\mathbb{R}}\mathbf{A}^{\prime}_{H}(x)\varrho(x)dx\right]^{-1}\mathbf{e}_{2}.
\]
To complete the proof, it suffices to show that $\tilde{\sigma}^{2}_{u}=\sigma^{2}_{u}+o_{p}(1)$.
This can be established using arguments analogous to those used for
the proof of \textbf{Theorem \ref{MainParametricWLS}(b)}. $\blacksquare$

\bigskip{}
\textbf{Proof of Theorem} \textbf{\ref{MainParametric-HTE}}. Let
$C_{n}$, $D_{n}$ $Q_{n}$ be as in (\ref{Proof_WIV_MainParam_WNP}).
Then by \textbf{Theorem \ref{MCLT-LL}}, as $n\rightarrow\infty$,
\begin{equation}
Q^{-1}_{n}\stackrel{p}{\rightarrow}Q^{-1}:=\left[\int_{\mathbb{R}}\mathbf{A}_{H}(x)\varrho(x)dx\right]^{-1}.\label{ProofMainParametric-HTE-1}
\end{equation}
Next, fix $\mathbf{v}=(v_{0},v_{1})^{\prime}\in\mathbb{R}^{2}$. Since
$w_{0},w_{1}$ are \textbf{HF}, the score term is

\[
\mathbf{v}^{\prime}D^{-1}_{n}\sum^{n}_{t=1}\mathbf{w}_{t-1}u_{t}=n^{-1/\alpha_{u}}\sum^{n}_{t=1}\left[v_{0}\frac{w_{0}(x_{t-1})}{k_{w_{0}}(\gamma_{n})}+v_{1}\frac{w_{1}(x_{t-1})}{k_{w_{1}}(\gamma_{n})}\right]u_{t}
\]
\[
=n^{-1/\alpha_{u}}\sum^{n}_{t=1}\mathbf{v}^{\prime}\mathbf{w}(\gamma^{-1}_{n}x_{t-1})u_{t}.
\]
We verify the conditions on $\mathbf{v}^{\prime}\mathbf{w}(.)$ required
by\textbf{ Theorem \ref{HT-MCLT-LL}}. For all $\alpha_{u}>0$, there
exists $C_{\mathbf{v},\alpha_{u}}<\infty$ such that

\[
|\mathbf{v}^{\prime}\mathbf{w}(x)|^{\alpha_{u}}\le C_{\mathbf{v},\alpha_{u}}\big(|w_{0}(x)|^{\alpha_{u}}+|w_{1}(x)|^{\alpha_{u}}\big)
\]
By condition (iv) and \textbf{Definition 4,} each $|w_{i}|^{\alpha_{u}}$
is locally integrable. In view of the above, $|\mathbf{v}^{\prime}\mathbf{w}(x)|^{\alpha_{u}}$
is locally integrable. Moreover, setting $\lambda=\max\{\lambda_{0},\lambda_{1}\}<\alpha$,
condition (iv) gives

\[
|\mathbf{v}^{\prime}\mathbf{w}(x)|^{\alpha_{u}}=O(|x|^{\lambda}),
\]
as $|x|\rightarrow\infty$. Next, we show condition (iii) of \textbf{Theorem
\ref{HT-MCLT-LL}} i.e. the positivity requirement $\int_{\mathbb{R}}|\mathbf{v}^{\prime}\mathbf{w}(x)|^{\alpha_{u}}\varrho(x)dx\neq0$.
Indeed, if $\int_{\mathbb{R}}|\mathbf{v}^{\prime}\mathbf{w}(x)|^{\alpha_{u}}\varrho(x)dx=0$
for some $\mathbf{v}\neq0$, then $\mathbf{v}^{\prime}\mathbf{w}(x)=0$
for $\varrho(x)dx$-almost every $x$. In view of the fact $\mathbf{A}_{H}(x)=\mathbf{w}(x)(1,f_{H}(x))$,
this implies $\mathbf{v}^{\prime}Q=0$, contradicting condition (v)
of \textbf{Theorem} \textbf{\ref{MainParametric-HTE}.}

Thus, by\textbf{ Theorem \ref{HT-MCLT-LL}} as $n\rightarrow\infty$,

\begin{equation}
\mathbf{v}^{\prime}D^{-1}_{n}\sum^{n}_{t=1}\mathbf{w}_{t-1}u_{t}\stackrel{d}{\rightarrow}S_{\alpha_{u}}\left\{ \sigma^{\alpha_{u}}_{u}\int_{\mathbb{R}}|\mathbf{v}^{\prime}\mathbf{w}(x)|^{\alpha_{u}}\varrho(x)dx,\ \beta^{*}_{u}\frac{\int_{\mathbb{R}}\mathbf{v}^{\prime}\mathbf{w}(x)^{\langle\alpha_{u}\rangle}\varrho(x)dx}{\int_{\mathbb{R}}|\mathbf{v}^{\prime}\mathbf{w}(x)|^{\alpha_{u}}\varrho(x)dx}\right\} .\label{ProofMainParametric-HTE-2}
\end{equation}
The requisite result follows from ($\ref{ProofMainParametric-HTE-1}$)
and $(\ref{ProofMainParametric-HTE-2})$. $\blacksquare$

\bigskip{}

\noindent\textbf{Proof of Theorem \ref{MainParametricWLS}}. \textbf{(a)}
Set $C^{*}_{n}=\sqrt{n}diag\{1,k_{f}(\gamma^{*}_{n})\}$ and $D^{*}_{n}=\sqrt{n}diag\{k_{w_{0}}(\gamma^{*}_{n}),k_{w_{1}}(\gamma^{*}_{n})\}$.
Write
\[
C^{*}_{n}\left(\mathbf{\tilde{b}}_{WIV}-\mathbf{b}\right)=C^{*}_{n}\left[\sum^{n}_{t=1}\mathbf{w}_{t-1}\mathbf{x}^{\prime}_{t-1}\right]^{-1}D^{*}_{n}D^{*-1}_{n}\sum^{n}_{t=1}\mathbf{w}_{t-1}u_{t}
\]
\[
=\left\{ D^{*-1}_{n}\left[\sum^{n}_{t=1}\mathbf{w}_{t-1}\mathbf{x}^{\prime}_{t-1}\right]C^{*-1}_{n}\right\} ^{-1}D^{*-1}_{n}\sum^{n}_{t=1}\mathbf{w}_{t-1}u_{t}.
\]
Furthermore,
\begin{equation}
D^{*-1}_{n}\left[\sum^{n}_{t=1}\mathbf{w}_{t-1}\mathbf{x}^{\prime}_{t-1}\right]C^{*-1}_{n}\overset{d}{\rightarrow}\int^{1}_{0}\mathbf{A}_{H}(\varLambda_{\alpha,d}(r))dr.\label{MainParamWLS0.1}
\end{equation}
Next, for $i=\left\{ 0,1\right\} $, 
\[
\mathbf{E}\left|\frac{1}{\sqrt{n}k_{w_{i}}(\gamma^{*}_{n})}\sum^{n}_{t=1}w_{i}(x_{t-1})u_{t}\right|^{2}\leq2\mathbf{E}\sum^{n}_{t=1}\left|\frac{1}{\sqrt{n}k_{w_{i}}(\gamma^{*}_{n})}w_{i}(x_{t-1})u_{t}\right|^{2}
\]
\[
=\frac{2}{n}\mathbf{E}\sum^{n}_{t=1}\left|\frac{1}{k_{w_{i}}(\gamma^{*}_{n})}w_{i}(x_{t-1})u_{t}\right|^{2}=\frac{2}{n}\mathbf{E}\sum^{n}_{t=1}\left|\frac{1}{k_{w_{i}}(\gamma^{*}_{n})}w_{i}(x_{t-1})\right|^{2}\mathbf{E}_{t-1}u^{2}_{t}
\]
\[
=\frac{2}{n}\mathbf{E}\sum^{n}_{t=1}\left|\frac{1}{k_{w_{i}}(\gamma^{*}_{n})}w_{i}(x_{t-1})\right|^{2}\mathbf{E}_{t-1}u^{2}_{t}\leq\frac{\sigma^{2}}{n}\mathbf{E}\sum^{n}_{t=1}\left|\frac{1}{k_{w_{i}}(\gamma^{*}_{n})}w_{i}(x_{t-1})\right|^{2}
\]
\begin{equation}
\leq_{(1)}\frac{\sigma^{2}}{n}\mathbf{E}\sum^{n}_{t=1}\left|\frac{1}{k_{w_{i}}(\gamma^{*}_{n})}w_{i}(x_{t-1})\right|^{2}=O\left(1\right),\label{MainParamWLS0.2}
\end{equation}
where $\leq_{(1)}$ follows from \textbf{Lemma \ref{lem:SSalphaM}}.
Hence,
\[
C^{*}_{n}\left(\mathbf{\tilde{b}}_{WIV}-\mathbf{b}\right)=O_{p}(1),
\]
as required.

\noindent\textbf{(b)} Next, we show the consistency of 
\[
\tilde{\sigma}^{2}_{w}=\frac{\sum^{n}_{t=1}\tilde{u}^{2}_{t}w_{u}(x_{t-1})}{\sum^{n}_{t=1}w_{u}(x_{t-1})}.
\]
Recall that $\mathbf{x}_{t-1}:=\left[1,f(x_{t-1})\right]^{\prime}$.
Then
\[
\tilde{u}_{t}=y_{t}-\mathbf{x}^{\prime}_{t-1}\mathbf{\tilde{b}}_{WIV}=\mathbf{x}^{\prime}_{t-1}\left(\mathbf{\tilde{b}}_{WIV}-\mathbf{b}\right)+u_{t}
\]
 For a matrix $\mathbf{A}=[a_{ij}]$, set $\left\Vert \mathbf{A}\right\Vert =\sum_{ij}\left|a_{ij}\right|$.
By \textbf{Theorem \ref{MCLT-LL}} as, $n\rightarrow\infty$, 
\begin{equation}
\frac{1}{nk_{w_{u}}(\gamma^{*}_{n})}\sum^{n}_{t=1}\left\Vert \mathbf{x}^{\prime}_{t-1}C^{*-1}_{n}\right\Vert w_{u}(x_{t-1})\stackrel{d}{\rightarrow}\int^{1}_{0}\left[\left(1+\left|f_{H}\right|\right)w_{u}\right]\left(\varLambda_{\alpha,d}(r)\right)dr,\label{MainParamWLS0.3}
\end{equation}
\begin{equation}
\frac{1}{nk_{w_{u}}(\gamma^{*}_{n})}\sum^{n}_{t=1}\left\Vert \mathbf{x}^{\prime}_{t-1}C^{*-1}_{n}\right\Vert ^{2}w_{u}(x_{t-1})\stackrel{d}{\rightarrow}\int^{1}_{0}\left[\left(1+\left|f_{H}\right|\right)^{2}w_{u}\right]\left(\varLambda_{\alpha,d}(r)\right)dr,\label{MainParamWLS0.4}
\end{equation}
and
\begin{equation}
\frac{1}{nk_{w_{u}}(\gamma^{*}_{n})}\sum^{n}_{t=1}w_{u}(x_{t-1})\stackrel{d}{\rightarrow}\int^{1}_{0}w_{u}\left(\varLambda_{\alpha,d}(r)\right)dr.\label{MainParamWLS0.5}
\end{equation}
Moreover, by the Bahr-Essen inequality for martingales, see e.g. Pinelis
(2015, Thm 1.1), for some $1<p<2$,
\[
\mathbf{E}\left|\frac{1}{nk_{w_{u}}(\gamma^{*}_{n})}\sum^{n}_{t=1}\left(u^{2}_{t}-\sigma^{2}_{u}\right)w_{u}(x_{t-1})\right|^{p}\leq\frac{1}{n^{p}k_{w_{u}}(\gamma^{*}_{n})^{p}}\sum^{n}_{t=1}\mathbf{E}\left[\left|u^{2}_{t}-\sigma^{2}_{u}\right|^{p}w_{u}(x_{t-1})^{p}\right]
\]
\begin{equation}
\leq\frac{C}{n^{p}k_{w_{u}}(\gamma^{*}_{n})^{p}}\sum^{n}_{t=1}\mathbf{E}w_{u}(x_{t-1})^{p}\leq_{(2)}\frac{C}{n^{p-1}}=o(1),\label{MainParamWLS0.6}
\end{equation}
where $\leq_{(2)}$ follows from \textbf{Lemma \ref{lem:SSalphaM}}.

\noindent Next, consider the numerator of $\tilde{\sigma}^{2}_{u}$
\[
\frac{1}{nk_{w_{u}}(\gamma^{*}_{n})}\sum^{n}_{t=1}\tilde{u}^{2}_{t}w_{u}(x_{t-1})=\frac{1}{nk_{w_{u}}(\gamma^{*}_{n})}\sum^{n}_{t=1}u^{2}_{t}w_{u}(x_{t-1})+2\frac{1}{nk_{w_{u}}(\gamma^{*}_{n})}\sum^{n}_{t=1}\mathbf{x}^{\prime}_{t}C^{*-1}_{n}C^{*}_{n}\left(\mathbf{\tilde{b}}_{WIV}-\mathbf{b}\right)w_{u}(x_{t-1})
\]
\[
+\frac{1}{nk_{w_{u}}(\gamma^{*}_{n})}\sum^{n}_{t=1}\left[\mathbf{x}^{\prime}_{t}C^{*-1}_{n}C^{*}_{n}\left(\mathbf{\tilde{b}}_{WIV}-\mathbf{b}\right)\right]^{2}w_{u}(x_{t-1})
\]
\[
=:\frac{1}{nk_{w_{u}}(\gamma^{*}_{n})}\sum^{n}_{t=1}u^{2}_{t}w_{u}(x_{t-1})+2R_{1n}+R_{2n}.
\]

\noindent In view of (\ref{MainParamWLS0.3})-(\ref{MainParamWLS0.6}),
we get
\[
R_{1n}\leq\left\Vert \sqrt{n}C^{*}_{n}\left(\mathbf{\tilde{b}}_{WIV}-\mathbf{b}\right)\right\Vert \frac{1}{n^{3/2}k_{w_{u}}(\gamma^{*}_{n})}\sum^{n}_{t=1}\left\Vert \mathbf{x}^{\prime}_{t-1}C^{*-1}_{n}\right\Vert w_{u}(x_{t-1})=O_{p}(n^{-1/2}),
\]
\[
R_{2n}\leq\left\Vert \sqrt{n}C^{*}_{n}\left(\mathbf{\tilde{b}}_{WIV}-\mathbf{b}\right)\right\Vert ^{2}\frac{1}{n^{2}k_{w_{u}}(\gamma^{*}_{n})}\sum^{n}_{t=1}\left\Vert \mathbf{x}^{\prime}_{t-1}C^{*-1}_{n}\right\Vert ^{2}w_{u}(x_{t-1})=O_{p}(n^{-1}),
\]
\[
\frac{1}{nk_{w_{u}}(\gamma^{*}_{n})}\sum^{n}_{t=1}u^{2}_{t}w_{u}(x_{t-1})
\]
\[
=\frac{\sigma^{2}_{u}}{nk_{w_{u}}(\gamma^{*}_{n})}\sum^{n}_{t=1}w_{u}(x_{t-1})+o_{p}(1)\stackrel{d}{\rightarrow}\sigma^{2}_{u}\int^{1}_{0}w_{u}\left(\varLambda_{\alpha,d}(r)\right)dr.
\]
Combining the above we have
\[
\tilde{\sigma}^{2}_{w}\stackrel{p}{\rightarrow}\sigma^{2}_{u},
\]
as required. $\blacksquare$

\section{Proofs for Section \ref{sec:IVX-Inference}}

\subsection{Supporting Results for Section \ref{sec:IVX-Inference}}

\begin{lemma} \emph{\label{Lemma_IVX_support}Suppose that the following
hold:}
\begin{itemize}
\item[\textbf{(i)}]  \textbf{\emph{Assumption INN}}\emph{ holds with $\theta=1$;}
\item[\textbf{(ii)}]  $x_{t}$ \emph{satisfies }\textbf{\emph{Assumption SP}}\emph{;}
\item[\textbf{(iii)}]  \textit{$z_{t}:=z_{t}(n)$ is defined by (\ref{IVXz0}) and $\rho_{z,n}$
by (\ref{IVXz}).}
\end{itemize}
\emph{Then
\[
\sup_{n\geq1,1\leq t\leq n}\int_{\mathbb{R}}\left|\lambda\right|\left|\mathbf{E}\left(e^{i\lambda z_{t}}\right)\right|d\lambda<\infty.
\]
}

\end{lemma}

\subsection{Proofs of Main Results in Section \ref{sec:IVX-Inference}}

\textbf{Proof of Theorem \ref{IVX_SP}}. \textbf{(a)} Let $g$ locally
integrable and suppose that for some $0<\lambda<\alpha$, $g(x)=O\left(\left|x\right|^{\lambda}\right)$,
as $\left|x\right|\rightarrow\infty$. We will show that 
\begin{equation}
\frac{1}{n}\sum^{n}_{t=1}g(z_{t})\stackrel{p}{\rightarrow}\mathbf{E}g(x_{t}),\label{Proof_IVX_SP_0}
\end{equation}
as $n\rightarrow\infty$. Indeed, since the elements of $\mathbf{A}(x)$
and $\mathbf{W}(x)$ satisfy the same local integrability and growth
conditions as $g(x)$, the result follows by arguments analogous to
those used in the proof of \textbf{Theorem \ref{WIV_SP}}. 

First, by \textbf{Lemma \ref{Lemma_IVX_support}} under \textbf{Assumption
SP}, $z_{t}:=z_{t}(n)$ has a density $D_{z_{t}(n)}(x)$ that is uniformly
bounded i.e. density \emph{
\begin{equation}
\sup_{n\geq1,1\leq t\leq n}\sup_{x}D_{z_{t}(n)}(x)<\infty.\label{Proof_IVX_SP_1}
\end{equation}
}Furthermore, the arguments used for the proof of \textbf{Lemma \ref{Lemma_IVX_support}}
show that $x_{t}$ has also a density $D(x)$ say, such that
\begin{equation}
\sup_{x}D(x)<\infty.\label{Proof_IVX_SP_4}
\end{equation}
Next, is it easy is show that for all $0<\alpha^{\prime}<\alpha$,
\begin{equation}
\sup_{t\geq1}\mathbf{E}\left|x_{t}\right|^{\alpha^{\prime}}<\infty.\label{Proof_IVX_SP_2}
\end{equation}
In addition, it follows from \textbf{Lemma \ref{Lemma_kappa_z}} that
as $n\rightarrow\infty$,
\begin{equation}
\sup_{1\leq t\leq n}\mathbf{E}\left|z_{t}(n)-x_{t}\right|^{\alpha^{\prime}}\rightarrow0.\label{Proof_IVX_SP_3}
\end{equation}

Next, let $g$ locally integrable such that for some $0<\lambda<\alpha$,
$g(x)=O\left(\left|x\right|^{\lambda}\right)$, as $\left|x\right|\rightarrow\infty$.
It follows from (\ref{Proof_IVX_SP_1})-(\ref{Proof_IVX_SP_3}) and
\textbf{Proposition \ref{BoundSmoothAppox}} that for each $\epsilon>0$
there is a bounded and Lipschwitz $g_{\epsilon}$ such that
\begin{equation}
g_{\epsilon}(x)=0,\;\left|x\right|>\epsilon^{-1},\label{Proof_IVX_SP_5}
\end{equation}
\begin{equation}
\int_{\left|x\right|\leq\epsilon^{-1}}\left|g_{\epsilon}(x)-g(x)\right|dx<\epsilon,\label{Proof_IVX_SP_6}
\end{equation}
and as $n\rightarrow\infty$ and the $\epsilon\rightarrow0$
\begin{equation}
\frac{1}{n}\sum^{n}_{t=1}g(z_{t})=\frac{1}{n}\sum^{n}_{t=1}g_{\epsilon}(x_{t})+o_{p}(1).\label{Proof_IVX_SP_7}
\end{equation}
Since, $\left\{ g_{\epsilon}(x_{t})\right\} ^{n}_{t=1}$ is strictly
stationary and ergodic, see, e.g. (\ref{Proof_WIV_SP_1}), by the
ergodic theorem we get 
\begin{equation}
\frac{1}{n}\sum^{n}_{t=1}g_{\epsilon}(x_{t})\stackrel{a.s.}{\rightarrow}\mathbf{E}g_{\epsilon}(x_{t}),\label{Proof_IVX_SP_8}
\end{equation}
as $n\rightarrow\infty$. 

Next, in view of (\ref{Proof_IVX_SP_4}) as $\epsilon\rightarrow0$,
\[
\left|\mathbf{E}g_{\epsilon}(x_{t})-\mathbf{E}g(x_{t})\right|\leq\mathbf{E}\left|g_{\epsilon}(x_{t})-g(x_{t})\right|=\int_{\mathbb{R}}\left|g_{\epsilon}(x)-g(x)\right|D(x)dx
\]
\[
=_{(1)}\int_{\left|x\right|\leq\epsilon^{-1}}\left|g_{\epsilon}(x)-g(x)\right|D(x)dx+\int_{\left|x\right|>\epsilon^{-1}}\left|g(x)\right|D(x)dx
\]
\[
\leq\sup_{x}D(x)\int_{\left|x\right|\leq\epsilon^{-1}}\left|g_{\epsilon}(x)-g(x)\right|dx+\mathbf{E}g(x_{t})I_{\left\{ \left|x_{t}\right|>\epsilon^{-1}\right\} }
\]
\[
\leq\sup_{x}D(x)\int_{\left|x\right|\leq\epsilon^{-1}}\left|g_{\epsilon}(x)-g(x)\right|dx+\mathbf{E}g(x_{t})I_{\left\{ \left|x_{t}\right|>\epsilon^{-1}\right\} }
\]
\begin{equation}
\leq_{(2)}\sup_{x}D(x)\cdot\epsilon+C\mathbf{E}\left|x_{t}\right|^{\lambda}I_{\left\{ \left|x_{t}\right|>\epsilon^{-1}\right\} }=_{(3)}0,\label{Proof_IVX_SP_9}
\end{equation}
where $=_{(1)}$, $\leq_{(2)}$ follow from (\ref{Proof_IVX_SP_5})
and (\ref{Proof_IVX_SP_6}) respectively. Finally, $=_{(3)}$ follows
from dominated convergence i.e. it is possible to find $\lambda<\alpha^{\prime}<\alpha$
such that $\mathbf{E}\left|x_{t}\right|^{\alpha^{\prime}}<C<\infty$.
Combining (\ref{Proof_IVX_SP_7})-(\ref{Proof_IVX_SP_9}) we get (\ref{Proof_IVX_SP_0}).

\textbf{(b)} Finally, for the IVX t-statistic, as $n\rightarrow\infty$,
\[
\tilde{t}_{IVX}=\tilde{t}_{WLS}+o_{p}(1)\stackrel{d}{\rightarrow}N\left(0,1\right),
\]
where the first relation follows from part (a), while the weak convergence
is a direct consequence of Theorem (\ref{WIV_SP}). $\blacksquare$

\bigskip{}

\noindent\textbf{Proof of Theorem \ref{IVX_SNP}}. Under our Assumptions,
$z_{t}$ is of the form 
\[
z_{t}=\sum^{t}_{j=1}\rho^{t-j}_{z,n}v_{j},
\]
and therefore exactly a MI process with localizing coefficients $\rho_{z,n}$.
Moreover, setting $\mu_{n}=\gamma_{z,n}$ it follows from \textbf{Lemma
\ref{WCon}}, 
\[
\frac{1}{\mu_{n}}z_{n}\stackrel{d}{\rightarrow}X^{+},
\]
where $X^{+}$ is some $\alpha$-stable variate. Set $D_{n}=\sqrt{n}diag\{k_{w_{0}}(\mu_{n}),k_{w_{1}}(\mu_{n})\}$.
In view of the above, \textbf{Theorem \ref{MainLL}}, we obtain, for
as $n\rightarrow\infty$, 
\[
D^{-1}_{n}\sum^{n}_{t=1}\mathbf{W}(z_{t-1})D^{-1}_{n}\stackrel{p}{\rightarrow}\int_{\mathbb{R}}\mathbf{W}_{H}(x)\varrho(x)dx
\]

\noindent
\begin{equation}
=\int_{\mathbb{R}}\left[\begin{array}{cc}
w^{2}_{0,H} & w_{0,H}w_{1,H}\\
w_{0,H}w_{1,H} & w^{2}_{1,H}
\end{array}\right](x)\varrho(x)dx,\label{IVX_SNP1}
\end{equation}
where $\varrho(x)$ is the density function of $X^{+}$. Moreover,

Under $H_{0}$,
\[
\hat{t}_{IVX}=\frac{1}{\sqrt{\hat{\sigma}^{2}\sum^{n}_{t=1}\overline{w}_{1}(z_{t-1})^{2}}}\sum^{n}_{t=1}w_{1}(z_{t-1})\overline{u}_{t}
\]
\begin{equation}
=\frac{1}{\sqrt{\hat{\sigma}^{2}\sum^{n}_{t=1}\overline{w}_{1}(z_{t-1})^{2}}}\sum^{n}_{t=1}w_{1}(z_{t-1})\left(u_{t}-\frac{\sum^{n}_{j=1}u_{j}w_{0}(z_{j-1})}{\sum^{n}_{j=1}w_{0}(z_{j-1})}\right).\label{IVX_SNP2}
\end{equation}
Set $e^{\prime}_{2}:=\left[0,1\right]$. In view of (\ref{IVX_SNP1}),
\textbf{Theorem \ref{MCLT-LL} }and the Cramér--Wold device, we obtain
\[
e^{\prime}_{2}D^{-1}_{n}\sum^{n}_{t=1}\mathbf{W}(z_{t-1})\bar{u}_{t}
\]
\[
=\frac{1}{\sqrt{n}k_{w_{1}}(\mu_{n})}\sum^{n}_{t=1}w_{1}(z_{t-1})\left(u_{t}-\frac{\sum^{n}_{j=1}u_{j}w_{0}(z_{j-1})}{\sum^{n}_{j=1}w_{0}(z_{j-1})}\right)
\]
\[
=\left[\begin{array}{cc}
-\frac{\frac{1}{nk_{w_{1}}(\mu_{n})}\sum^{n}_{t=1}w_{1}(z_{t-1})}{\frac{1}{nk_{w_{0}}(\mu_{n})}\sum^{n}_{t=1}w_{0}(z_{t-1})}, & 1\end{array}\right]\left[\begin{array}{c}
\frac{1}{\sqrt{n}k_{w_{0}}(\mu_{n})}\sum^{n}_{t=1}w_{0}(z_{t-1})u_{t}\\
\frac{1}{\sqrt{n}k_{w_{1}}(\mu_{n})}\sum^{n}_{t=1}w_{1}(z_{t-1})u_{t}
\end{array}\right]
\]
\[
=:\left[\begin{array}{cc}
-a_{n}, & 1\end{array}\right]\left[\begin{array}{c}
\frac{1}{\sqrt{n}k_{w_{0}}(\mu_{n})}\sum^{n}_{t=1}w_{0}(z_{t-1})u_{t}\\
\frac{1}{\sqrt{n}k_{w_{1}}(\mu_{n})}\sum^{n}_{t=1}w_{1}(z_{t-1})u_{t}
\end{array}\right]
\]
\[
\stackrel{d}{\rightarrow}\left[-\frac{\int_{\mathbb{R}}w_{1,H}(x)\varrho(x)dx}{\int_{\mathbb{R}}w_{0,H}(x)\varrho(x)dx},1\right]\cdot N\left[\mathbf{0},\sigma^{2}_{u}\int_{\mathbb{R}}\mathbf{W}_{H}(x)\varrho(x)dx\right]
\]

\[
=:\left[-a_{H},1\right]\cdot N\left[\mathbf{0},\sigma^{2}_{u}\int_{\mathbb{R}}\mathbf{W}_{H}(x)\varrho(x)dx\right]
\]
\begin{equation}
=N\left[0,\sigma^{2}_{u}\int_{\mathbb{R}}\left[w_{1,H}(x)-a_{H}w_{0,H}\right]^{2}\varrho(x)dx\right].\label{IVX_SNP3}
\end{equation}

Next, note that 
\[
\overline{w}_{1}(z_{t-1})=w_{1}(z_{t-1})-\frac{\sum^{n}_{j=1}w_{1}(z_{j-1})}{\sum^{n}_{j=1}w_{0}(z_{j-1})}w_{1}(z_{t-1}).
\]
Consequently,
\[
\frac{1}{k_{w_{1}}(\mu_{n})}\overline{w}_{1}(z_{t-1})=\frac{1}{k_{w_{1}}(\mu_{n})}w_{1}(z_{t-1})-a_{n}\frac{1}{k_{w_{0}}(\mu_{n})}w_{1}(z_{t-1}).
\]
Therefore,
\[
\frac{1}{nk^{2}_{w_{1}}(\mu_{n})}\sum^{n}_{t=1}\overline{w}_{1}(z_{t-1})^{2}=\left[\begin{array}{cc}
-a_{n}, & 1\end{array}\right]D^{-1}_{n}\sum^{n}_{t=1}\mathbf{W}(z_{t-1})D^{-1}_{n}\left[\begin{array}{cc}
-a_{n}, & 1\end{array}\right]^{\prime}
\]
\begin{equation}
\stackrel{p}{\rightarrow}\int_{\mathbb{R}}\left[w_{1,H}(x)-a_{H}w_{0,H}\right]^{2}\varrho(x)dx>0.\label{IVX_SNP4}
\end{equation}
Finally, in view of (\ref{IVX_SNP2})-(\ref{IVX_SNP4}) and the fact
that $\hat{\sigma}^{2}=\sigma^{2}_{u}+o_{p}(1)$, we obtain 
\[
\hat{t}_{IVX}=\frac{\frac{1}{\sqrt{n}k_{w_{1}}(\mu_{n})}\sum^{n}_{t=1}w_{1}(z_{t-1})\bar{u}_{t}}{\left[\frac{\hat{\sigma}^{2}}{nk^{2}_{w_{1}}(\mu_{n})}\sum^{n}_{t=1}\overline{w}_{1}(z_{t-1})^{2}\right]^{1/2}}\stackrel{d}{\rightarrow}N(0,1),
\]
as required. $\blacksquare$

\bigskip{}

Before we prove \textbf{Theorem \ref{IVX_SNP-2}}, we introduce some
notation.\bigskip{}

\noindent\textbf{Definition 6}. \emph{Suppose that $x_{t}$ satisfies
either }\textbf{\emph{Assumption SNP }}\emph{or }\textbf{\emph{Assumption
MI}}\emph{. Define the random array $\tilde{z}_{t}(n)$ by}
\[
\tilde{z}_{t}:=\tilde{z}_{t}(n):=\left\{ \begin{array}{l}
\sum^{t}_{j=1}\rho^{t-j}_{z,n}v_{j},\:\textrm{\emph{under \textbf{Assumption SNP}}}\\
\textrm{\emph{or \textbf{Assumption MI} with} \ensuremath{\kappa_{n}}/\ensuremath{\kappa_{z,n}\rightarrow\infty}};\\
\\x_{t},\:\textrm{\emph{under \textbf{Assumption MI} with} \ensuremath{\kappa_{z,n}}/\ensuremath{\kappa_{n}\rightarrow\infty}}.
\end{array}\right.
\]

\bigskip{}

\noindent\textbf{Proof of Theorem \ref{IVX_SNP-2}}. Set 
\[
R_{t}(n):=z_{t}-\tilde{z}_{t}.
\]
By \textbf{Lemma \ref{Lemma_kappa_z}}, for all $1<\alpha^{\prime}<\alpha$
\begin{equation}
\lim_{n\rightarrow\infty}\sup_{1\leq t\leq n}\mathbf{E}\left|\frac{1}{\mu_{n}}R_{t}(n)\right|^{\alpha^{\prime}}=0,\label{IVX_SNP1-02a}
\end{equation}
where $\mu_{n}$ is given in \textbf{Definition 5}. Moreover, note
that $\tilde{z}_{t}$ is a MI process. Therefore, by \textbf{Lemma
\ref{WCon}}, 
\[
\frac{1}{\mu_{n}}\tilde{z}_{n}\stackrel{d}{\rightarrow}X^{+},
\]
where $X^{+}$ is some $\alpha$-stable variate. Set $D_{n}=\sqrt{n}diag\{k_{w_{0}}(\mu_{n}),k_{w_{1}}(\mu_{n})\}$.
In view (\ref{IVX_SNP1-02a}), \textbf{Theorem \ref{MainLL}} and
the fact that $\mathbf{W}(x)$ is \textbf{\emph{AHF}}$\left\{ z_{t},\mu_{n}\right\} $
and \textbf{LH}, as $n\rightarrow\infty$ we obtain 
\[
D^{-1}_{n}\sum^{n}_{t=1}\mathbf{W}(z_{t-1})D^{-1}_{n}=\frac{1}{n}\sum^{n}_{t=1}\mathbf{W}_{H}\left(\mu^{-1}_{n}\tilde{z}_{t-1}\right)+o_{p}(1)
\]
\begin{equation}
\stackrel{p}{\rightarrow}\int_{\mathbb{R}}\mathbf{W}_{H}(x)\varrho(x)dx,\label{IVX_SNP1-02}
\end{equation}
where $\varrho(x)$ is the density of $X^{+}$. In view of (\ref{IVX_SNP1-02}),
the remainder of the proof follows by arguments analogous to those
used in the proof of \textbf{Theorem \ref{IVX_SNP}}. $\blacksquare$

\bigskip{}

\subsection{Proofs of Technical Results for Section \ref{sec:IVX-Inference}}

\textbf{Proof of Lemma \ref{Lemma_kappa_z}}. \textbf{Proof of Lemma
\ref{Lemma_kappa_z}.} We start with two decompositions of $z_{t}$.
\begin{itemize}
\item Recall that under assumptions $x_{0}=0$. Summation by parts gives
\[
z_{t}=x_{t}-\rho^{t-1}_{z,n}x_{0}-(1-\rho_{z,n})\sum^{t-1}_{j=1}\rho^{t-j-1}_{z,n}x_{j}
\]
\[
=x_{t}-\rho^{t-1}_{z,n}x_{0}-(1-\left\{ 1-\kappa^{-1}_{z,n}\right\} )\sum^{t-1}_{j=1}\rho^{t-j-1}_{z,n}x_{j}
\]
\begin{equation}
=x_{t}-\frac{1}{\kappa_{z,n}}\sum^{t-1}_{j=1}\rho^{t-j-1}_{z,n}x_{j}.\label{NI-MIcov0.0}
\end{equation}
Hence,
\begin{equation}
\frac{1}{\gamma_{n}}z_{t}=\frac{1}{\gamma_{n}}x_{t}-\frac{1}{\gamma_{n}\kappa_{z,n}}\sum^{t-1}_{j=1}\rho^{t-j-1}_{z,n}x_{j}.\label{NI-MIcov0.1}
\end{equation}
\item Alternatively, 
\[
z_{t}=\sum^{t}_{j=1}\rho^{t-j}_{z,n}\Delta x_{j}=\sum^{t}_{j=1}\rho^{t-j}_{z,n}v_{j}-\frac{1}{k_{n}}\sum^{t}_{j=1}\rho^{t-j}_{z,n}x_{j-1}
\]
\[
=\sum^{t}_{j=1}\rho^{t-j}_{z,n}v_{j}-\frac{1}{k_{n}}\sum^{t-1}_{j=1}\rho^{t-j-1}_{z,n}x_{j}.
\]
Hence,
\end{itemize}
\[
\frac{1}{\gamma_{z,n}}z_{t}=\frac{1}{\gamma_{z,n}}\sum^{t}_{j=1}\rho^{t-j}_{z,n}\Delta x_{j}
\]

\[
=\frac{1}{\gamma_{z,n}}\sum^{t}_{j=1}\rho^{t-j}_{z,n}v_{j}-\frac{1}{\gamma_{z,n}k_{n}}\sum^{t-1}_{j=1}\rho^{t-j-1}_{z,n}x_{j}.
\]
Facts:
\begin{itemize}
\item $\gamma_{n}=R_{\mathcal{H}}(\kappa_{n}):=\kappa^{\mathcal{H}}_{n}\varsigma\left(\kappa_{n}\right)$,
where $\varsigma\left(.\right)$ is slowly varying -i.e. $\gamma_{n}$
is regularly varying \ in $\kappa_{n}$ of (Hurst) index $0<\mathcal{H<}1$.
\item $\gamma_{z,n}=R_{\mathcal{H}}(\kappa_{z,n}):=\kappa^{\mathcal{H}}_{z,n}\varsigma\left(\kappa_{z,n}\right)$
-i.e. $\gamma_{n}$ is regularly varying in $\kappa_{z,n}$ of (Hurst)
index $0<\mathcal{H<}1$. 
\end{itemize}
\bigskip{}

\begin{itemize}
\item[\textbf{(a)}]  Case $\frac{\kappa_{n}}{\kappa_{z,n}}\rightarrow\infty:$ The result
follows analogous arguments to those used in part (b) below.\bigskip{}
\item[\textbf{(b)}]  Case $\frac{\kappa_{z,n}}{\kappa_{n}}\rightarrow\infty:$ Recall
that we can decompose $x_{t}$ as
\[
x_{t}=x^{+}_{t}+x^{-}_{t},
\]
where 
\[
x^{+}_{t}=\sum^{t-1}_{k=0}a_{k}\xi_{t-k},\text{ \ }x^{-}_{t}=\sum^{\infty}_{k=t}a^{-}_{k,t}\xi_{t-k}.
\]
Hence, we can write
\[
\sum^{t-1}_{j=1}\rho^{t-j-1}_{z,n}x_{j}=\sum^{t-1}_{j=1}\rho^{t-j-1}_{z,n}x^{+}_{j}+\sum^{t-1}_{j=1}\rho^{t-j-1}_{z,n}x^{-}_{j}.
\]
\end{itemize}
\bigskip{}

For $1<\alpha\leq2$, let 
\[
0<\epsilon<\left(\alpha-1\right)\land\left[\left(2-\alpha\right)I_{\{\alpha<2\}}+I_{\{\alpha=2\}}\right]
\]
 and set
\[
\alpha_{\epsilon^{-}}:=\alpha-\epsilon\;\textrm{and}
\]
\[
\alpha_{\epsilon^{+}}:=\left\{ \begin{array}{l}
\alpha+\epsilon,\;\alpha<2\\
\alpha,\;\text{\ensuremath{\alpha}=2}
\end{array}\right..
\]

\begin{itemize}
\item $x^{+}_{t}$: Write 
\end{itemize}
\[
\sum^{t-1}_{j=1}\rho^{t-j-1}_{z,n}x^{+}_{j}=\sum^{t-1}_{j=1}\rho^{t-j-1}_{z,n}\sum^{j-1}_{k=0}a_{k}\xi_{j-k}
\]

\[
\left[l=j-k\rightarrow k=j-l\right]
\]
\[
=\sum^{t-1}_{j=1}\rho^{t-j-1}_{z,n}\sum^{j}_{l=1}a_{j-l}\xi_{l}=\sum^{t-1}_{l=1}\left(\sum^{t-1}_{j=l}\rho^{t-j-1}_{z,n}\sum^{j-l}_{i=0}\rho^{j-l-i}_{n}c_{i}\right)\xi_{l}
\]
\[
\left.\begin{array}{c}
1\leq j\leq t-1\\
1\leq l\leq j
\end{array}\right\} \rightarrow\left.\begin{array}{c}
1\leq l\leq t-1\\
l\leq j\leq t-1
\end{array}\right\} 
\]
\[
=\sum^{t-1}_{l=1}\left(\sum^{t-1}_{j=l}\rho^{t-j-1}_{z,n}\sum^{j-l}_{i=0}\rho^{j-l-i}_{n}c_{i}\right)\xi_{l}
\]
\[
=:\sum^{t-1}_{l=1}A_{nl}\xi_{l}=:\mathcal{A}_{n,t}.
\]
Recall that that under \textbf{Assumption MI-SM},
\begin{equation}
\gamma_{n}\sim\kappa^{1/\alpha}_{n}\ell_{1}\left(\kappa_{n}\right)^{1/\alpha},\text{ }H\left(\gamma_{n}\right)\sim\ell_{1}\left(\kappa_{n}\right)\text{.}\label{SNP-NM -0.1-1}
\end{equation}
Let
\[
\overline{\xi}_{n,l}:=\xi_{l}1\left\{ \left\vert \xi_{l}\right\vert \geq\gamma_{n}\right\} -\mathbf{E}\left[\xi_{l}1\left\{ \left\vert \xi_{l}\right\vert \geq\gamma_{n}\right\} \right],
\]
\[
\underline{\xi}_{n,l}:=\xi_{l}1\left\{ \left\vert \xi_{l}\right\vert <\gamma_{n}\right\} -\mathbf{E}\left[\xi_{l}1\left\{ \left\vert \xi_{l}\right\vert <\gamma_{n}\right\} \right].
\]
Note that as $n\rightarrow\infty$, cf. Astrauskas (1983, Lemma 1),
\[
\mathbf{E}\left|\overline{\xi}_{n,l}\right|^{\alpha_{\epsilon^{-}}}\leq Ch(\gamma_{n})\gamma^{-\epsilon}_{n}\left[1+o(1)\right],
\]
and 
\[
\mathbf{E}\left|\underline{\xi}_{n,l}\right|^{\alpha_{\epsilon^{+}}}\leq C\left[h(\gamma_{n})\gamma^{\epsilon}_{n}1\{\alpha\neq2\}+H(\gamma_{n})1\{\alpha=2\}\right]\left[1+o(1)\right],
\]
where recall that $h(\gamma_{n})/H(\gamma_{n})<C$. We explicitly
consider only the case $\alpha\neq2$. The proof for the case $\alpha=2$
follows from analogous arguments to those used below. Write
\[
\mathbf{E}\left|\frac{1}{\gamma_{n}\kappa_{z,n}}\mathcal{A}_{n,t}\right|^{\alpha_{\epsilon^{-}}}\leq\mathbf{E}\left|\frac{1}{\gamma_{n}\kappa_{z,n}}\sum^{t-1}_{l=1}A_{nl}\xi_{l}\right|^{\alpha_{\epsilon^{-}}}
\]
\[
\leq2\mathbf{E}\left|\frac{1}{\gamma_{n}\kappa_{z,n}}\sum^{t-1}_{l=1}A_{nl}\overline{\xi}_{n,l}\right|^{\alpha_{\epsilon^{-}}}+2\mathbf{E}\left|\frac{1}{\gamma_{n}\kappa_{z,n}}\sum^{t-1}_{l=1}A_{nl}\underline{\xi}_{n,l}\right|^{\alpha_{\epsilon^{-}}}
\]
\[
\leq2\mathbf{E}\left|\frac{1}{\gamma_{n}\kappa_{z,n}}\sum^{t-1}_{l=1}A_{nl}\overline{\xi}_{n,l}\right|^{\alpha_{\epsilon^{-}}}+2\left[\mathbf{E}\left|\frac{1}{\gamma_{n}\kappa_{z,n}}\sum^{t-1}_{l=1}A_{nl}\underline{\xi}_{n,l}\right|^{\alpha_{\epsilon^{+}}}\right]^{\frac{\alpha_{\epsilon^{-}}}{\alpha_{\epsilon^{+}}}}
\]
\[
\leq_{(1)}2\left(\frac{1}{\gamma_{n}\kappa_{z,n}}\right)^{\alpha_{\epsilon^{-}}}\sum^{t-1}_{l=1}\left|A_{nl}\right|^{\alpha_{\epsilon^{-}}}\mathbf{E}\left|\overline{\xi}_{n,l}\right|^{\alpha_{\epsilon^{-}}}+2\left[\left(\frac{1}{\gamma_{n}\kappa_{z,n}}\right)^{\alpha_{\epsilon^{+}}}\sum^{t-1}_{l=1}\left|A_{nl}\right|^{\alpha_{\epsilon^{+}}}\mathbf{E}\left|\underline{\xi}_{n,l}\right|^{\alpha_{\epsilon^{+}}}\right]^{\frac{\alpha_{\epsilon^{-}}}{\alpha_{\epsilon^{+}}}}
\]
\[
\leq C\left(\frac{1}{\gamma_{n}\kappa_{z,n}}\right)^{\alpha_{\epsilon^{-}}}\sum^{t-1}_{l=1}\left|A_{nl}\right|^{\alpha_{\epsilon^{-}}}h(\gamma_{n})\gamma^{-\epsilon}_{n}\left[1+o(1)\right]
\]
\[
+C\left[\left(\frac{1}{\gamma_{n}\kappa_{z,n}}\right)^{\alpha_{\epsilon^{+}}}\sum^{t-1}_{l=1}\left|A_{nl}\right|^{\alpha_{\epsilon^{+}}}h(\gamma_{n})\gamma^{\epsilon}_{n}\left[1+o(1)\right]\right]^{\frac{\alpha_{\epsilon^{-}}}{\alpha_{\epsilon^{+}}}}
\]
\[
\leq C\left(\frac{1}{\gamma^{\alpha}_{n}}\right)\sum^{t-1}_{l=1}\frac{1}{\kappa^{\alpha_{\epsilon^{-}}}_{z,n}}\left|A_{nl}\right|^{\alpha_{\epsilon^{-}}}h(\gamma_{n})\left[1+o(1)\right]+C\left[\frac{1}{\gamma^{\alpha}_{n}}\sum^{t-1}_{l=1}\frac{1}{\kappa^{\alpha_{\epsilon^{+}}}_{z,n}}\left|A_{nl}\right|^{\alpha_{\epsilon^{+}}}h(\gamma_{n})\left[1+o(1)\right]\right]^{\frac{\alpha_{\epsilon^{-}}}{\alpha_{\epsilon^{+}}}}
\]
\[
\leq C\frac{1}{\gamma^{\alpha}_{n}}\sum^{t-1}_{l=1}\frac{1}{\kappa^{\alpha_{\epsilon^{-}}}_{z,n}}\left|A_{nl}\right|^{\alpha_{\epsilon^{-}}}H(\gamma_{n})\left[1+o(1)\right]+C\left[\frac{1}{\gamma^{\alpha}_{n}}\sum^{t-1}_{l=1}\frac{1}{\kappa^{\alpha_{\epsilon^{+}}}_{z,n}}\left|A_{nl}\right|^{\alpha_{\epsilon^{+}}}H(\gamma_{n})\left[1+o(1)\right]\right]^{\frac{\alpha_{\epsilon^{-}}}{\alpha_{\epsilon^{+}}}}
\]
\[
\leq C\frac{1}{\kappa_{n}}\sum^{t-1}_{l=1}\frac{1}{\kappa^{\alpha_{\epsilon^{-}}}_{z,n}}\left|A_{nl}\right|^{\alpha_{\epsilon^{-}}}\left[1+o(1)\right]+C\left[\frac{1}{\kappa_{n}}\sum^{t-1}_{l=1}\frac{1}{\kappa^{\alpha_{\epsilon^{+}}}_{z,n}}\left|A_{nl}\right|^{\alpha_{\epsilon^{+}}}\left[1+o(1)\right]\right]^{\frac{\alpha_{\epsilon^{-}}}{\alpha_{\epsilon^{+}}}},
\]
where $\leq_{(1)}$ follows from Bahr-Esseen inequality if the exponents
are greater than one and the Loève inequality when there less or equal
to unity.

Next, set 
\[
\alpha_{\epsilon}:=\{\alpha_{\epsilon^{+}},\alpha_{\epsilon^{-}}\}.
\]
Moreover,
\[
R_{1nt}:=\frac{1}{\kappa_{n}}\sum^{t-1}_{l=1}\frac{1}{\kappa^{\alpha_{\epsilon}}_{z,n}}\left|A_{nl}\right|^{\alpha_{\epsilon}}=\frac{1}{\kappa_{n}}\frac{1}{\kappa^{\alpha_{\epsilon}}_{z,n}}\sum^{t-1}_{l=1}\left\vert \sum^{t-1}_{j=l}\rho^{t-j-1}_{z,n}\sum^{j-l}_{i=0}\rho^{j-l-i}_{n}c_{i}\right\vert ^{\alpha_{\epsilon}}.
\]
The summand
\[
\sum^{t-1}_{j=l}\rho^{t-j-1}_{z,n}\sum^{j-l}_{i=0}\rho^{j-l-i}_{n}c_{i}
\]
\[
\left[k=j-l\rightarrow j=k+l\right]
\]
\begin{equation}
=\sum^{t-1-l}_{k=0}\rho^{t-\left(k+l\right)-1}_{z,n}\sum^{k+l-l}_{i=0}\rho^{k+l-l-i}_{n}c_{i}=\sum^{t-1-l}_{k=0}\rho^{t-1-l-k}_{z,n}\sum^{k}_{i=0}\rho^{k-i}_{n}c_{i}.\label{NI-MIcov0.2}
\end{equation}
Set
\[
B_{n}\left(h\right):=\left\{ \begin{array}{l}
\sum^{h}_{l=0}\rho^{h-l}_{z,n}\rho^{l}_{n},\;h\geq0\\
0,\;h<0
\end{array}\right..
\]
Note that for $p\geq1$, and $n$ large enough such that $0<\rho_{z,n},\rho_{n}<1$
and  $q\in\mathbb{Z}$, 
\[
\left[\sum^{\infty}_{h=0}B_{n}\left(q+h\right)^{p}\right]^{1/p}=\left[\sum^{\infty}_{h=q}B_{n}\left(h\right)^{p}I_{\left\{ h\geq0\right\} }\right]^{1/p}\leq\left[\sum^{\infty}_{h=0}B_{n}\left(h\right)^{p}\right]^{1/p}
\]
\[
=\left[\sum^{\infty}_{h=0}\left(\sum^{h}_{l=0}\rho^{h-l}_{z,n}\rho^{l}_{n}\right)^{p}\right]^{1/p}=\left[\sum^{\infty}_{h=0}\left[\sum^{\infty}_{l=0}\rho^{h-l}_{z,n}\rho^{l}_{n}I_{\left\{ h\geq l\right\} }\right]^{p}\right]^{1/p}
\]
\[
\leq_{(1)}\sum^{\infty}_{l=0}\left[\sum^{\infty}_{h=0}\left[\rho^{h-l}_{z,n}\rho^{l}_{n}I_{\left\{ h\geq l\right\} }\right]^{p}\right]^{1/p}=\sum^{\infty}_{l=0}\rho^{l}_{n}\left[\sum^{\infty}_{h=0}\left[\rho^{p\left(h-l\right)}_{z,n}I_{\left\{ h\geq l\right\} }\right]\right]^{1/p}
\]
\[
\left[j=h-l\right]
\]
\begin{equation}
=\sum^{\infty}_{l=0}\rho^{l}_{n}\left[\sum^{\infty}_{j=-l}\rho^{pj}_{z,n}I_{\left\{ j\geq0\right\} }\right]^{1/p}=\sum^{\infty}_{l=0}\rho^{l}_{n}\left[\sum^{\infty}_{j=0}\rho^{pj}_{z,n}\right]^{1/p}\leq C\kappa_{n}\left(\kappa_{z,n}\right)^{1/p},\label{NI-MIcov0.21}
\end{equation}
where $\leq_{(1)}$ follows from Minkowski's inequality. Moreover,
we have used the fact that 
\[
\sum^{\infty}_{l=0}\rho^{l}_{n}=\frac{1}{1-\rho_{n}}=\kappa_{n},
\]
and that as $n\rightarrow\infty$,
\[
\sum^{\infty}_{j=0}\rho^{pj}_{z,n}=\frac{1}{1-\left(1-\kappa^{-1}_{z,n}\right)^{p}}\sim\frac{\kappa_{z,n}}{p}.
\]
First, assume that $1<\alpha\leq2$. In this case $\alpha_{\epsilon}>1$.
In view of (\ref{NI-MIcov0.2}),
\[
\kappa_{n}\kappa^{\alpha_{\epsilon}}_{z,n}R_{1nt}=\sum^{t-1}_{l=1}\left\vert \sum^{t-1-l}_{k=0}\rho^{t-1-l-k}_{z,n}\sum^{k}_{i=0}\rho^{k-i}_{n}c_{i}\right\vert ^{\alpha_{\epsilon}}
\]
\[
\left[j=t-1-l\right]
\]
\[
=\sum^{t-2}_{j=0}\left\vert \sum^{j}_{k=0}\rho^{j-k}_{z,n}\sum^{k}_{i=0}\rho^{k-i}_{n}c_{i}\right\vert ^{\alpha_{\epsilon}}=\sum^{t-2}_{j=0}\left\vert \sum^{j}_{k=0}\rho^{j-k}_{z,n}\sum^{k}_{i=0}\rho^{k-i}_{n}c_{i}\right\vert ^{\alpha_{\epsilon}}
\]
\[
\left.\begin{array}{c}
0\leq k\leq j\\
0\leq i\leq k
\end{array}\right\} \rightarrow\left.\begin{array}{c}
0\leq i\leq j\\
i\leq k\leq j
\end{array}\right\} 
\]
\[
=\sum^{t-2}_{j=0}\left\vert \sum^{j}_{i=0}c_{i}\sum^{j}_{k=i}\rho^{j-k}_{z,n}\rho^{k-i}_{n}\right\vert ^{\alpha_{\epsilon}}
\]
\[
l=k-i\rightarrow k=l+i
\]
\[
=\sum^{t-2}_{j=0}\left\vert \sum^{j}_{i=0}c_{i}\sum^{j-i}_{k=0}\rho^{j-\left(l+i\right)}_{z,n}\rho^{l}_{n}\right\vert ^{\alpha_{\epsilon}}=\sum^{t-2}_{j=0}\left\vert \sum^{j}_{i=0}c_{i}\sum^{j-i}_{l=0}\rho^{j-i-l}_{z,n}\rho^{l}_{n}\right\vert ^{\alpha_{\epsilon}}
\]
\begin{equation}
=\sum^{t-2}_{j=0}\left\vert \sum^{j}_{i=0}c_{i}B_{n}\left(j-i\right)\right\vert ^{\alpha_{\epsilon}}\leq\sum^{t-2}_{j=0}\left\vert \sum^{j}_{i=0}\left\vert c_{i}\right\vert B_{n}\left(j-i\right)\right\vert ^{\alpha_{\epsilon}}\label{SNP-NM -0.1.1}
\end{equation}
\[
=\left\{ \left[\sum^{t-2}_{j=0}\left\vert \sum^{\infty}_{i=0}\left\vert c_{i}\right\vert B_{n}\left(j-i\right)\right\vert ^{\alpha_{\epsilon}}\right]^{1/\alpha_{\epsilon}}\right\} ^{\alpha_{\epsilon}}\leq_{(2)}\left\{ \sum^{\infty}_{i=0}\left[\sum^{t-2}_{j=0}\left[\left\vert c_{i}\right\vert B_{n}\left(j-i\right)\right]^{\alpha_{\epsilon}}\right]^{1/\alpha_{\epsilon}}\right\} ^{\alpha_{\epsilon}}
\]
\[
\leq\left\{ \sum^{\infty}_{i=0}\left\vert c_{i}\right\vert \left[\sum^{t}_{j=0}B_{n}\left(j-i\right)^{\alpha_{\epsilon}}\right]^{1/\alpha_{\epsilon}}\right\} ^{\alpha_{\epsilon}}\leq\left\{ \sum^{\infty}_{i=0}\left\vert c_{i}\right\vert \left[\sum^{\infty}_{j=0}B_{n}\left(j-i\right)^{\alpha_{\epsilon}}\right]^{1/\alpha_{\epsilon}}\right\} ^{\alpha_{\epsilon}}
\]
\begin{equation}
\leq_{(3)}C\left\{ \sum^{\infty}_{i=0}\left\vert c_{i}\right\vert \kappa_{n}\kappa^{1/\alpha_{\epsilon}}_{z,n}\right\} ^{\alpha_{\epsilon}},\label{NI-MIcov0.22}
\end{equation}
where $\leq_{(2)}$ follows from Minkowski's inequality and (\ref{NI-MIcov0.21}).
In view of (\ref{NI-MIcov0.22}), as $n\rightarrow\infty$, we get
\[
R_{1nt}\leq\frac{1}{\kappa_{n}\kappa^{\alpha_{\epsilon}}_{z,n}}\left\{ \sum^{\infty}_{i=0}\left\vert c_{i}\right\vert \left[\sum^{t}_{j=0}B_{n}\left(j-i\right)^{\alpha_{\epsilon}}\right]^{1/\alpha_{\epsilon}}\right\} ^{\text{\ensuremath{\alpha_{\epsilon}}}}\leq\frac{C}{\kappa_{n}\kappa^{\alpha_{\epsilon}}_{z,n}}\left\{ \sum^{\infty}_{i=0}\left\vert c_{i}\right\vert \kappa_{n}\kappa^{1/\alpha_{\epsilon}}_{z,n}\right\} ^{\alpha_{\epsilon}}
\]
\begin{equation}
\leq\frac{C}{\kappa_{n}\kappa^{\alpha_{\epsilon}}_{z,n}}\kappa^{\alpha_{\epsilon}}_{n}\kappa_{z,n}=C\frac{\kappa^{\alpha_{\epsilon}-1}_{n}}{\kappa^{\alpha_{\epsilon}-1}_{z,n}}\rightarrow0.\label{NI-MIcov0.23}
\end{equation}

Next suppose that $0<\alpha\leq1$. This case $0<\alpha_{\epsilon^{-}}<1$,
while $\alpha_{\epsilon^{+}}$ can be either greater or smaller than
unity. Therefore, we either have $\text{\ensuremath{\alpha_{\epsilon}}}>1$
or $\text{\ensuremath{\alpha_{\epsilon}}}\leq1$. For the former case,
the result follows using the same arguments as those used for establishing
(\ref{NI-MIcov0.23}). Therefore, lets assume that $0<\text{\ensuremath{\alpha_{\epsilon}}}\leq1$.
In view of (\ref{SNP-NM -0.1.1}),
\[
\kappa_{n}\kappa^{\alpha_{\epsilon}}_{z,n}R_{1nt}\leq\sum^{t-2}_{j=0}\left\vert \sum^{j}_{i=0}\left\vert c_{i}\right\vert B_{n}\left(j-i\right)\right\vert ^{\alpha_{\epsilon}}
\]
\[
\leq\sum^{t-2}_{j=0}\sum^{j}_{i=0}\left\vert c_{i}\right\vert ^{\alpha_{\epsilon}}B_{n}\left(j-i\right)^{\alpha_{\epsilon}}\leq\sum^{\infty}_{i=0}\left\vert c_{i}\right\vert ^{\alpha_{\epsilon}}\sum^{\infty}_{h=0}B_{n}\left(h\right)^{\alpha_{\epsilon}}.
\]
Next, note that

\[
\sum^{\infty}_{h=0}B^{\alpha_{\epsilon}}_{n}\left(h\right)=\sum^{\infty}_{h=0}\left[\sum^{h}_{l=0}\rho^{h-l}_{z,n}\rho^{l}_{n}\right]^{\alpha_{\epsilon}}=\sum^{\infty}_{h=0}\left[\rho^{h}_{z,n}\sum^{h}_{l=0}\left(\frac{\rho_{n}}{\rho_{z,n}}\right)^{l}\right]^{\alpha_{\epsilon}}
\]
\[
\leq\sum^{\infty}_{h=0}\left[\rho^{h}_{z,n}\sum^{\infty}_{l=0}\left(\frac{\rho_{n}}{\rho_{z,n}}\right)^{l}\right]^{\alpha_{\epsilon}}\leq\sum^{\infty}_{h=0}\rho^{\alpha_{\epsilon}h}_{z,n}\left[\sum^{\infty}_{l=0}\left(\frac{\rho_{n}}{\rho_{z,n}}\right)^{l}\right]^{\alpha_{\epsilon}}
\]
\[
\leq\sum^{\infty}_{h=0}\rho^{\alpha_{\epsilon}h}_{z,n}\left[\sum^{\infty}_{l=0}\left(\frac{\rho_{n}}{\rho_{z,n}}\right)^{l}\right]^{\alpha_{\epsilon}}=\frac{1}{1-\rho^{\alpha_{\epsilon}}_{z,n}}\left[\frac{1}{1-\frac{\rho_{n}}{\rho_{z,n}}}\right]^{\alpha_{\epsilon}}
\]
\[
=\frac{1}{1-\rho^{\alpha_{\epsilon}}_{z,n}}\left[\frac{\rho_{z,n}}{\rho_{z,n}-\rho_{n}}\right]^{\alpha_{\epsilon}}=\frac{1}{1-\rho^{\alpha_{\epsilon}}_{z,n}}\left[\frac{\rho_{z,n}}{\rho_{z,n}-\rho_{n}}\right]^{\alpha_{\epsilon}}=\frac{1}{1-\rho^{\alpha_{\epsilon}}_{z,n}}\left[\frac{\rho_{z,n}}{\kappa^{-1}_{n}-\kappa^{-1}_{z,n}}\right]^{\alpha_{\epsilon}}
\]
\[
=\frac{1}{1-\rho^{\alpha_{\epsilon}}_{z,n}}\left[\frac{\kappa_{n}\rho_{z,n}}{1-\kappa_{n}\kappa^{-1}_{z,n}}\right]^{\alpha_{\epsilon}}\sim\frac{\kappa_{z,n}}{\alpha_{\epsilon}}\kappa^{\alpha_{\epsilon}}_{n}.
\]
 Hence, 
\[
R_{1nt}\leq C\frac{1}{\kappa_{n}\kappa^{\alpha_{\epsilon}}_{z,n}}\frac{\kappa_{z,n}}{\alpha_{\epsilon}}\kappa^{\alpha_{\epsilon}}_{n}\rightarrow\infty
\]

\begin{itemize}
\item $x^{-}_{t}$: Recall that 
\[
x^{-}_{t}=\sum^{\infty}_{k=t}a^{-}_{k,t}(n)\xi_{t-k}=\sum^{\infty}_{k=t}\left(\sum^{t-1}_{i=0}\rho^{i}_{n}c_{k-i}\right)\xi_{t-k}
\]
Hence, 
\[
\sum^{t-1}_{j=1}\rho^{t-j-1}_{z,n}x^{-}_{j}=\sum^{t-1}_{j=1}\rho^{t-j-1}_{z,n}\sum^{\infty}_{k=j}\left(\sum^{j-1}_{i=0}\rho^{i}_{n}c_{k-i}\right)\xi_{j-k}
\]
\[
\left[-l=j-k\rightarrow l=k-j\rightarrow k=l+j\right]
\]
\[
=\sum^{t-1}_{j=1}\rho^{t-j-1}_{z,n}\sum^{\infty}_{l=0}\left(\sum^{j-1}_{i=0}\rho^{i}_{n}c_{l+j-i}\right)\xi_{-l}
\]
\end{itemize}
\[
=\sum^{\infty}_{l=0}\sum^{t-1}_{j=1}\rho^{t-j-1}_{z,n}\left(\sum^{j-1}_{i=0}\rho^{i}_{n}c_{l+j-i}\right)\xi_{-l}.
\]
Set
\[
R_{2nt}:=\frac{1}{\kappa_{n}\kappa^{\alpha_{\epsilon}}_{z,n}}\sum^{\infty}_{l=0}\left\vert \sum^{t-1}_{j=1}\rho^{t-j-1}_{z,n}\left(\sum^{j-1}_{i=0}\rho^{i}_{n}c_{l+j-i}\right)\right\vert ^{\alpha_{\epsilon}}.
\]
Using arguments analogous to those employed in the proof for $x^{+}_{t}$
it is easily seen that it suffices to establish 
\[
\sup_{1\leq t\leq n}R_{2nt}=o(1),
\]
to obtain the requisite result. Set $I_{l}(r):=I_{\left\{ l+1\leq r\leq l+t-1\right\} }$
and consider
\[
\kappa_{n}\kappa^{\alpha_{\epsilon}}_{z,n}R_{2nt}:=\sum^{\infty}_{l=0}\left\vert \sum^{t-1}_{j=1}\rho^{t-j-1}_{z,n}\left(\sum^{j-1}_{i=0}\rho^{i}_{n}c_{l+j-i}\right)\right\vert ^{\alpha_{\epsilon}}
\]
\[
\left[q=j-i\Longleftrightarrow i=j-q\right]
\]

\[
=\sum^{\infty}_{l=0}\left\vert \sum^{t-1}_{j=1}\rho^{t-j-1}_{z,n}\sum^{j}_{q=1}\rho^{j-q}_{n}c_{l+q}\right\vert ^{\alpha_{\epsilon}}
\]

\[
\left.\begin{array}{c}
1\leq j\leq t-1\\
1\leq q\leq j
\end{array}\right\} \rightarrow\left.\begin{array}{c}
1\leq q\leq t-1\\
q\leq j\leq t-1
\end{array}\right\} 
\]
\[
=\sum^{\infty}_{l=0}\left\vert \sum^{t-1}_{q=1}c_{l+q}\sum^{t-1}_{j=q}\rho^{t-j-1}_{z,n}\rho^{j-q}_{n}\right\vert ^{\alpha_{\epsilon}}
\]
\[
\left[r=j-q\Longleftrightarrow j=r+q\right]
\]

\[
=\sum^{\infty}_{l=0}\left\vert \sum^{t-1}_{q=1}c_{l+q}\sum^{t-1-q}_{r=0}\rho^{t-1-q-r}_{z,n}\rho^{r}_{n}\right\vert ^{\alpha_{\epsilon}}=:\sum^{\infty}_{l=0}\left\vert \sum^{t-1}_{q=1}c_{l+q}B_{n}\left(t-1-q\right)\right\vert ^{\alpha_{\epsilon}}
\]
\[
\left[r=l+q\Longleftrightarrow q=r-l\right]
\]
\[
=\sum^{\infty}_{l=0}\left\vert \sum^{l+t-1}_{r=l+1}c_{r}B_{n}\left(l+t-1-r\right)\right\vert ^{\alpha_{\epsilon}}=\sum^{\infty}_{l=0}\left\vert \sum^{\infty}_{r=0}c_{r}B_{n}\left(l+t-1-r\right)I_{l}(r)\right\vert ^{\alpha_{\epsilon}}
\]
\[
=:\left\{ \left[\sum^{\infty}_{l=0}\left\vert \sum^{\infty}_{r=0}c_{r}B_{n}\left(l+t-1-r\right)I_{l}(r)\right\vert ^{\alpha_{\epsilon}}\right]^{\frac{1}{\alpha_{\epsilon}}}\right\} ^{\alpha_{\epsilon}}\leq_{(4)}\left\{ \sum^{\infty}_{r=0}\left[\sum^{\infty}_{l=0}I_{l}(r)\left\vert c_{r}B_{n}\left(l+t-1-r\right)\right\vert ^{\alpha_{\epsilon}}\right]^{\frac{1}{\alpha_{\epsilon}}}\right\} ^{\alpha_{\epsilon}}
\]
\[
=\left\{ \sum^{\infty}_{r=0}\left|c_{r}\right|\left[\sum^{\infty}_{l=0}I_{l}(r)\left\vert B_{n}\left(l+t-1-r\right)\right\vert ^{\alpha_{\epsilon}}\right]^{\frac{1}{\alpha_{\epsilon}}}\right\} ^{\alpha_{\epsilon}}\leq\left\{ \sum^{\infty}_{r=0}\left|c_{r}\right|\left[\sum^{\infty}_{l=0}\left\vert B_{n}\left(l+t-1-r\right)\right\vert ^{\alpha_{\epsilon}}\right]^{\frac{1}{\alpha_{\epsilon}}}\right\} ^{\alpha_{\epsilon}}
\]

\begin{equation}
\leq_{(5)}C\left\{ \sum^{\infty}_{r=0}\left|c_{r}\right|\kappa_{n}\left(\kappa_{z,n}\right)^{\frac{1}{\alpha_{\epsilon}}}\right\} ^{\alpha_{\epsilon}},\label{NI-MIcov0.3}
\end{equation}
where $\leq_{(4)}$ follows from Mikowski's inequality and $\leq_{(5)}$
from (\ref{NI-MIcov0.21}). In view of (\ref{NI-MIcov0.3}), as $n\rightarrow\infty$,
we get
\[
R_{2nt}\leq C\kappa^{-1}_{n}\kappa^{-\alpha_{\epsilon}}_{z,n}\left\{ \sum^{\infty}_{r=0}\left|c_{r}\right|\kappa_{n}\left(\kappa_{z,n}\right)^{\frac{1}{\alpha_{\epsilon}}}\right\} ^{\alpha_{\epsilon}}\leq C\kappa^{-1}_{n}\kappa^{-\alpha_{\epsilon}}_{z,n}\left\{ \kappa_{n}\left(\kappa_{z,n}\right)^{\frac{1}{\alpha_{\epsilon}}}\right\} ^{\alpha_{\epsilon}}
\]
\[
=C\kappa^{-1}_{n}\kappa^{-\alpha_{\epsilon}}_{z,n}\kappa^{\alpha_{\epsilon}}_{n}\kappa_{z,n}=C\kappa^{1-\alpha_{\epsilon}}_{z,n}\kappa^{\alpha_{\epsilon}-1}_{n}=o(1).
\]

\bigskip{}

\textbf{(c)} First, recall that $\gamma_{z,n}=R_{\mathcal{H}}(\kappa_{z,n}):=\kappa^{\mathcal{H}}_{z,n}\varsigma\left(\kappa_{z,n}\right)$
-i.e. $\gamma_{n}$ is regularly varying in $\kappa_{z,n}$ of (Hurst)
index $0<\mathcal{H<}1$. Notice that by \textbf{Lemma FS}$_{\alpha}$\textbf{A},
for all $0<\alpha^{\prime}<\alpha$,
\[
\sup_{n\geq1,1\leq t\leq n}\mathbf{E}\left\vert \frac{1}{\gamma_{z,n}}\sum^{t}_{j=1}\rho^{t-j}_{z,n}x_{j}\right\vert ^{\alpha^{\prime}}=\sup_{n\geq1,1\leq t\leq n}\mathbf{E}\left\vert \frac{1}{\gamma_{z,n}}\sum^{t}_{j=1}\rho^{t-j}_{z,n}v_{j}\right\vert ^{\alpha^{\prime}}<\infty.
\]
In view of this, as $n\rightarrow\infty$,
\[
\sup_{1\leq t\leq n}\mathbf{E}\left\vert \frac{1}{\kappa_{z,n}}\sum^{t}_{j=1}\rho^{t-j}_{z,n}x_{j}\right\vert ^{\alpha^{\prime}}=\left(\frac{\gamma_{z,n}}{\kappa_{z,n}}\right)^{\alpha^{\prime}}\sup_{1\leq t\leq n}\mathbf{E}\left\vert \frac{1}{\gamma_{z,n}}\sum^{t}_{j=1}\rho^{t-j}_{z,n}x_{j}\right\vert ^{\alpha^{\prime}}
\]
\[
\leq C\left(\frac{\gamma_{z,n}}{\kappa_{z,n}}\right)^{\alpha^{\prime}}=C\left(\frac{R_{\mathcal{H}}\left(\kappa_{z,n}\right)}{\kappa_{z,n}}\right)^{\alpha^{\prime}}=o(1).
\]
In view of the above and (\ref{NI-MIcov0.0}), as $n\rightarrow\infty$,
we get
\[
\sup_{1\leq t\leq n}\mathbf{E}\left|z_{t}-x_{t}\right|^{\alpha^{\prime}}=\sup_{1\leq t\leq n}\mathbf{E}\left|\frac{1}{\kappa_{z,n}}\sum^{t-1}_{j=1}\rho^{t-j-1}_{z,n}x_{j}\right|^{\alpha^{\prime}}=o(1).
\]
The proof is complete. $\blacksquare$

\subsection{Proofs of Supporting Results for Section \ref{sec:IVX-Inference}}

\textbf{Proof of Lemma \ref{Lemma_IVX_support}}. Recall that $\{\xi_{t}\}$
is i.i.d. with characteristic function
\[
\psi_{\xi}(\lambda):=\mathbf{E}\left(e^{i\lambda\xi_{t}}\right).
\]
Moreover, we can write
\begin{equation}
z_{t}=x_{t}-\frac{1}{\kappa_{z,n}}\sum^{t-1}_{j=1}\rho^{t-j-1}_{z,n}x_{j}=:x_{t}-z^{-}_{t}.\label{Lemma_IVX_support-0}
\end{equation}
Since
\[
x_{t}=\sum^{\infty}_{k=0}c_{k}\xi_{t-k},
\]
the characteristic function of the first term is 
\begin{equation}
\phi_{x_{t}}(\lambda):=\mathbf{E}\left(e^{i\lambda x_{t}}\right)=\prod^{\infty}_{k=0}\psi_{\xi}(\lambda c_{k}).\label{Lemma_IVX_support-1}
\end{equation}

Next, consider $z^{-}_{t}$. Write
\[
z^{-}_{t}:=\frac{1}{\kappa_{z,n}}\sum^{t-1}_{j=1}\rho^{t-j-1}_{z,n}x_{j}=\frac{1}{\kappa_{z,n}}\sum^{t-1}_{j=1}\rho^{t-j-1}_{z}\sum^{\infty}_{i=0}c_{i}\xi_{j-i}
\]
\[
=\frac{1}{\kappa_{z,n}}\sum^{t-2}_{h=0}\rho^{h}_{z,n}\sum^{\infty}_{i=0}c_{i}\xi_{t-h-1-i}=\frac{1}{\kappa_{z,n}}\sum^{\infty}_{i=0}\sum^{t-2}_{h=0}\rho^{h}_{z,n}c_{i}\xi_{t-h-1-i}
\]
\[
\left[k=h+1+i\Longleftrightarrow h=k-i-1\right]
\]
\[
=\frac{1}{\kappa_{z,n}}\sum^{\infty}_{i=0}\sum^{t+i-1}_{k=i+1}\rho^{k-i-1}_{z,n}c_{i}\xi_{t-k}=\frac{1}{\kappa_{z,n}}\sum^{\infty}_{i=0}\sum^{t+i-1}_{k=i+1}\rho^{k-i-1}_{z,n}c_{i}\xi_{t-k}
\]
\[
\left[\begin{array}{c}
0\leq i\leq\infty\\
i+1\leq k\leq t-2+i+1
\end{array}\right]\Longleftrightarrow\left[\begin{array}{c}
1\leq k\leq\infty\\
0\vee\left(k+1-t\right)\leq i\leq k-1
\end{array}\right]
\]
\begin{equation}
=\frac{1}{\kappa_{z,n}}\sum^{\infty}_{k=1}\sum^{k-1}_{i=0\vee\left(k+1-t\right)}\rho^{k-i-1}_{z,n}c_{i}\xi_{t-k}=:\frac{1}{\kappa_{z,n}}\sum^{\infty}_{k=1}b_{t,k}\xi_{t-k}.\label{Lemma_IVX_support-1b}
\end{equation}
Hence, the characteristic function of $z^{-}_{t}$ is 
\[
\phi_{z^{-}_{t}}(\lambda)=\prod^{\infty}_{k=1}\psi_{\xi}\left(\lambda\frac{b_{t,k}}{\kappa_{z,n}}\right).
\]
Specifically, after a change of variables we can write the coefficients
$b_{t,k}$ as
\[
b_{t,k}=\sum^{k-1}_{i=0\vee\left(k+1-t\right)}\rho^{k-i-1}_{z,n}c_{i}
\]
\[
\left[l=k-i-1\Longleftrightarrow i=k-l-1\right]
\]
\[
=\sum^{k-1-\left(0\vee\left(k+1-t\right)\right)}_{l=0}\rho^{l}_{z,n}c_{k-1-l}=\sum^{\left(k-1-0\right)\land\left(k-1-\left(k+1-t\right)\right)}_{l=0}\rho^{l}_{z,n}c_{k-1-l}
\]
\begin{equation}
=\sum^{\left(k-1\right)\land\left(t-2\right)}_{l=0}\rho^{l}_{z,n}c_{k-1-l}.\label{Lemma_IVX_support-2}
\end{equation}

Next, we find an expression for the characteristic function of\textbf{
}$z_{t}$. Using (\ref{Lemma_IVX_support-0}) and (\ref{Lemma_IVX_support-1b})
we obtain
\[
z_{t}=c_{0}\xi_{t}+\sum^{\infty}_{k=1}\left(c_{k}-\frac{1}{\kappa_{z,n}}b_{t,k}\right)\xi_{t-k}.
\]
Consequently,
\begin{equation}
\phi_{z_{t}}(\lambda)=\psi_{\xi}(\lambda c_{0})\prod^{\infty}_{k=1}\psi_{\xi}\left(\lambda\left[c_{k}-\frac{1}{\kappa_{z,n}}b_{t,k}\right]\right).\label{Lemma_IVX_support-3}
\end{equation}
For $k\geq0$, write
\[
\tilde{c}_{k,t}:=\left\{ \begin{array}{l}
c_{0},\;k=0,\\
c_{k}-\frac{1}{\kappa_{z,n}}b_{t,k},\;k\geq1
\end{array}\right..
\]
In particular, using (\ref{Lemma_IVX_support-2}) note that for $k\geq1$,
\[
\tilde{c}_{k,t}=c_{k}-\frac{1}{\kappa_{z,n}}b_{t,k}=c_{k}-\frac{1}{\kappa_{z,n}}\sum^{\left(k-1\right)\land\left(t-2\right)}_{l=0}\rho^{l}_{z,n}c_{k-1-l}.
\]
Under \textbf{Assumption SP}, there is at least one coefficient in
the sequence $\left\{ c_{k}\right\} _{k\geq0}$ that is nonzero. Let
$k_{*}\geq0$ the minimal index for which $c_{k_{*}}\neq0$. In view
of the above we have,
\[
\tilde{c}_{k_{*},t}=c_{k_{*}}\neq0.
\]
By virtue of (\ref{Lemma_IVX_support-3}), we get
\[
\int_{\mathbb{R}}\left|\lambda\right|\left|\phi_{z_{t}}(\lambda)\right|d\lambda=\int_{\mathbb{R}}\left|\lambda\right|\left|\prod^{\infty}_{k=0}\psi_{\xi}\left(\lambda\tilde{c}_{k,t}\right)\right|d\lambda
\]
\[
\leq\int_{\mathbb{R}}\left|\lambda\right|\left|\psi_{\xi}\left(\lambda\tilde{c}_{k_{*},t}\right)\right|d\lambda=\left|\tilde{c}_{k_{*},t}\right|^{-2}\int_{\mathbb{R}}\left|\lambda\right|\left|\psi_{\xi}\left(\lambda\right)\right|d\lambda<\infty,
\]
uniformly in $1\leq t\leq n$, $n\geq1$, as required. $\blacksquare$

\clearpage{}

\section{Simulation Results}\label{sec:Simulation-Results}

This section collects the simulation evidence discussed in Section
\ref{sec:Simulations}. The four sets of size tables use the eight
stability parameters in the order 
\[
\alpha=2,1.99,1.8,1.6,1.4,1.2,1,0.8.
\]

Tables \ref{tab:sim-dgp1-linear-w05-alpha-2p00}, \ref{tab:sim-dgp1-linear-w05-alpha-1p99},
\ref{tab:sim-dgp1-linear-w05-alpha-1p80}, \ref{tab:sim-dgp1-linear-w05-alpha-1p60},
\ref{tab:sim-dgp1-linear-w05-alpha-1p40}, \ref{tab:sim-dgp1-linear-w05-alpha-1p20},
\ref{tab:sim-dgp1-linear-w05-alpha-1p00}, and \ref{tab:sim-dgp1-linear-w05-alpha-0p80}
report empirical size under DGP1 for a linear regression function
with the fixed weight $w(x)=x^{\langle0.5\rangle}$.

Tables \ref{tab:sim-dgp1-linear-weights-alpha-2p00}, \ref{tab:sim-dgp1-linear-weights-alpha-1p99},
\ref{tab:sim-dgp1-linear-weights-alpha-1p80}, \ref{tab:sim-dgp1-linear-weights-alpha-1p60},
\ref{tab:sim-dgp1-linear-weights-alpha-1p40}, \ref{tab:sim-dgp1-linear-weights-alpha-1p20},
\ref{tab:sim-dgp1-linear-weights-alpha-1p00}, and \ref{tab:sim-dgp1-linear-weights-alpha-0p80}
report empirical size under DGP1 for a linear regression function
with alternative weight functions.

Tables \ref{tab:sim-dgp1-nonlinear-alpha-2p00}, \ref{tab:sim-dgp1-nonlinear-alpha-1p99},
\ref{tab:sim-dgp1-nonlinear-alpha-1p80}, \ref{tab:sim-dgp1-nonlinear-alpha-1p60},
\ref{tab:sim-dgp1-nonlinear-alpha-1p40}, \ref{tab:sim-dgp1-nonlinear-alpha-1p20},
\ref{tab:sim-dgp1-nonlinear-alpha-1p00}, and \ref{tab:sim-dgp1-nonlinear-alpha-0p80}
report empirical size under DGP1 for the nonlinear regression specification.

Tables \ref{tab:sim-dgp2-linear-weights-alpha-2p00}, \ref{tab:sim-dgp2-linear-weights-alpha-1p99},
\ref{tab:sim-dgp2-linear-weights-alpha-1p80}, \ref{tab:sim-dgp2-linear-weights-alpha-1p60},
\ref{tab:sim-dgp2-linear-weights-alpha-1p40}, \ref{tab:sim-dgp2-linear-weights-alpha-1p20},
\ref{tab:sim-dgp2-linear-weights-alpha-1p00}, and \ref{tab:sim-dgp2-linear-weights-alpha-0p80}
report empirical size under DGP2 for a linear regression function
with alternative weight functions.

\begin{landscape} 

\begin{figure}[htbp]
\caption{Asymptotic Power (nominal size 5\%): $[f(x),w(x)]=[x,x^{\langle0.5\rangle}]$,
$\delta=-0.95$, $n=1000$}
\label{fig:sim-power-linear-w05}\centering \captionsetup[subfigure]{position=top}

\begin{subfigure}[t]{0.48\linewidth} \centering \caption{$\alpha=2$, $d_{*}=0.5$}
\label{fig:sim-power-linear-w05-alpha-2p00} \includegraphics[width=0.49\linewidth]{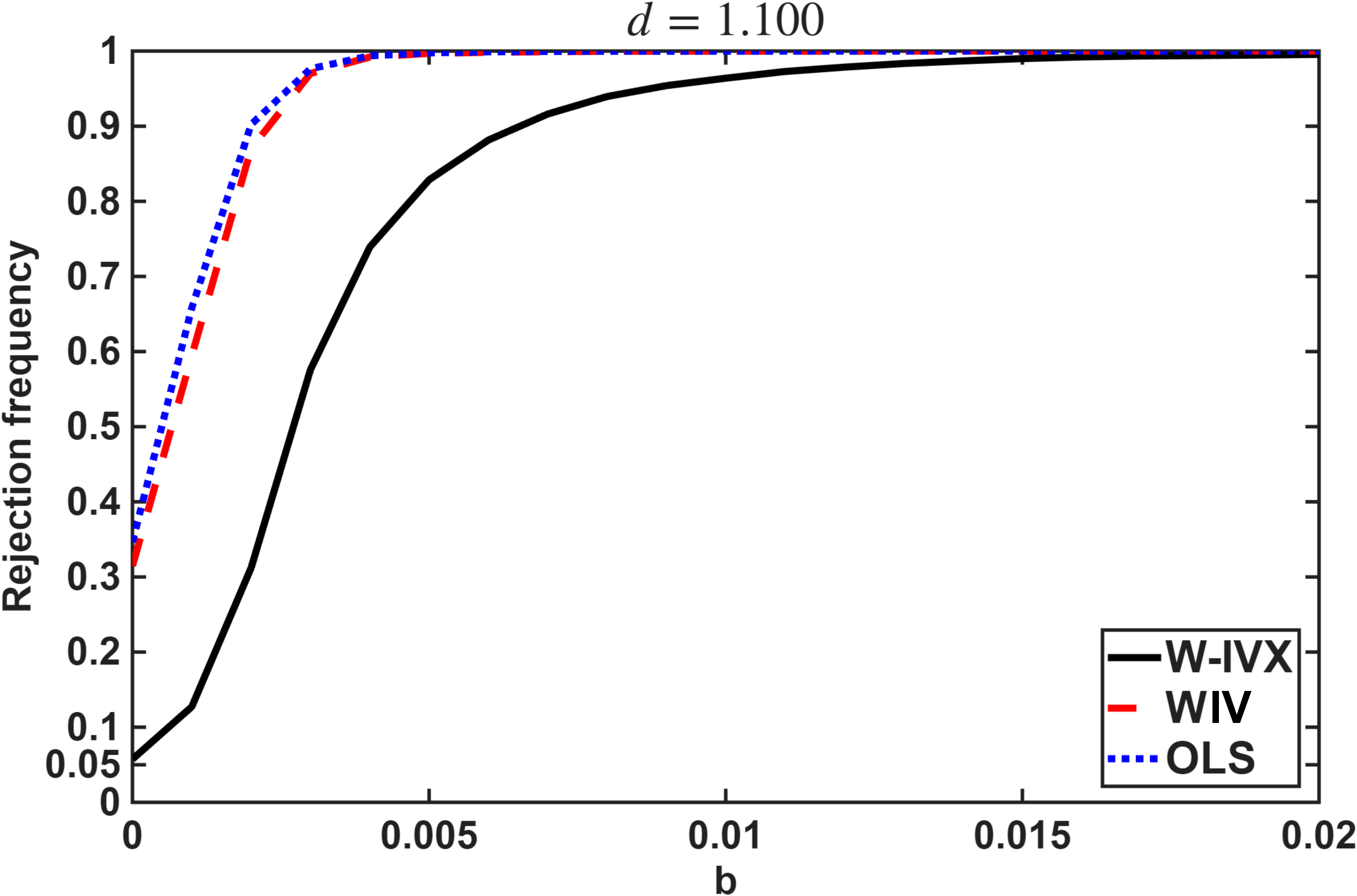}
\includegraphics[width=0.49\linewidth]{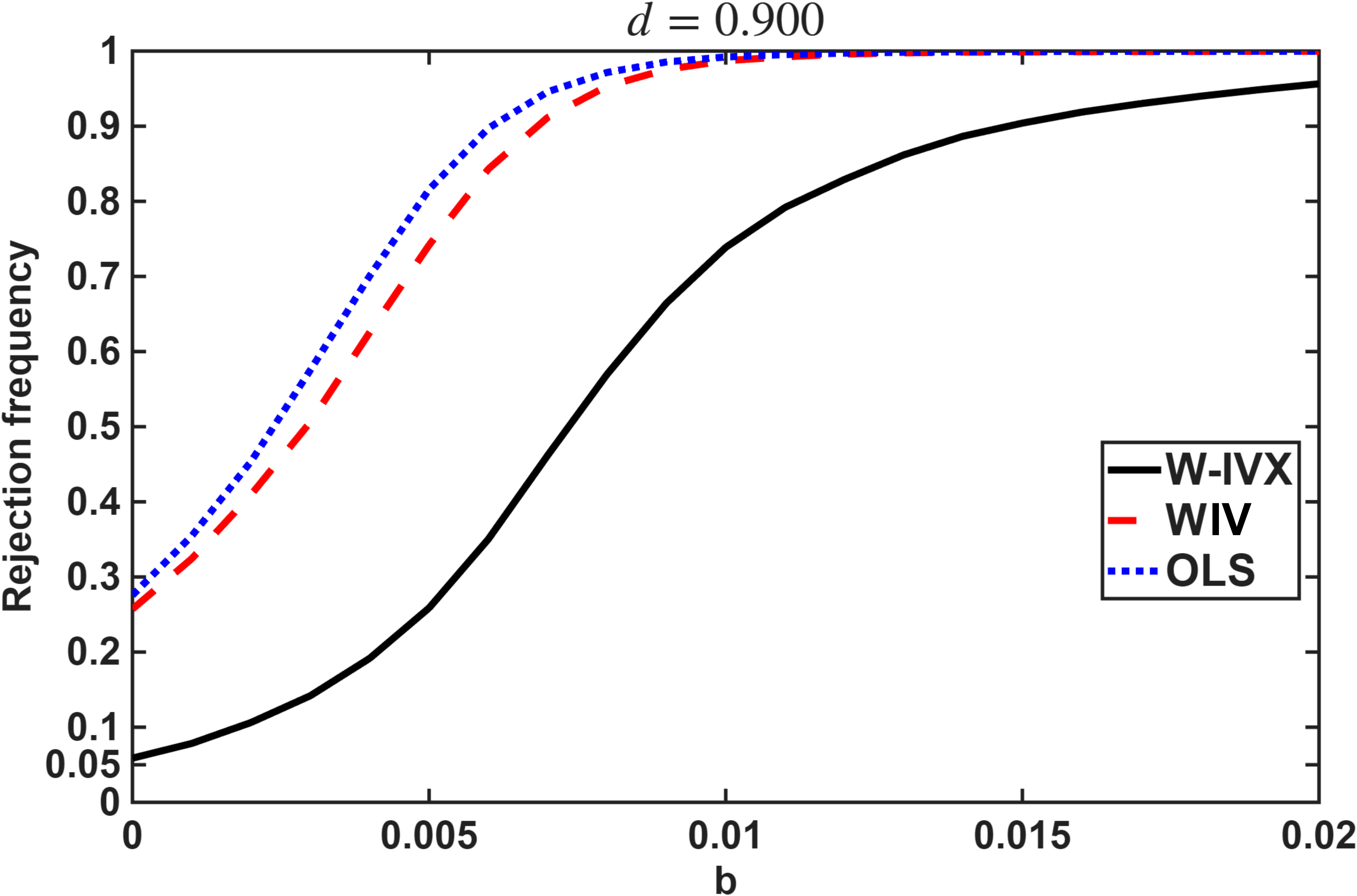}\\[0.6ex]
\includegraphics[width=0.49\linewidth]{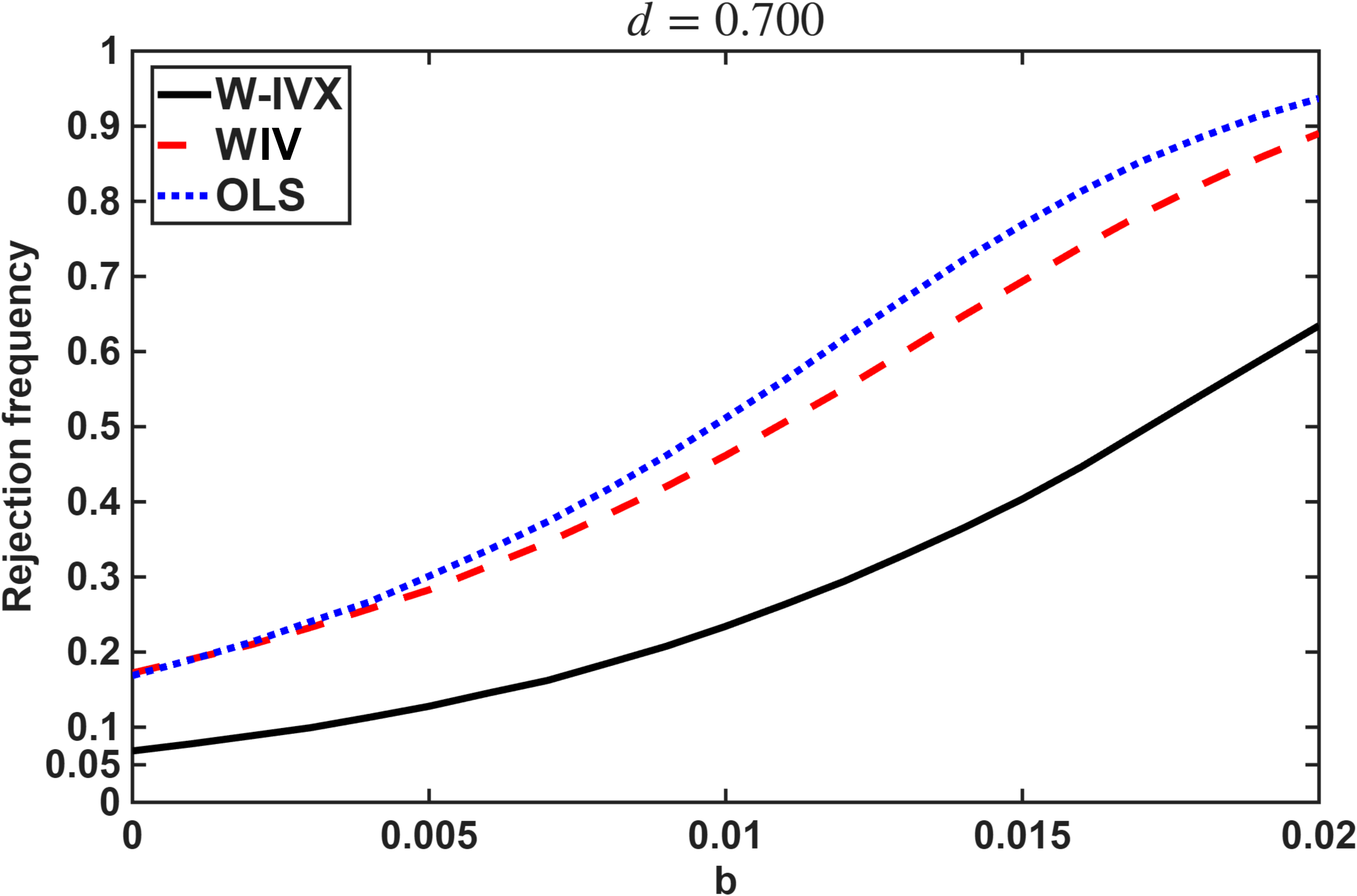}
\includegraphics[width=0.49\linewidth]{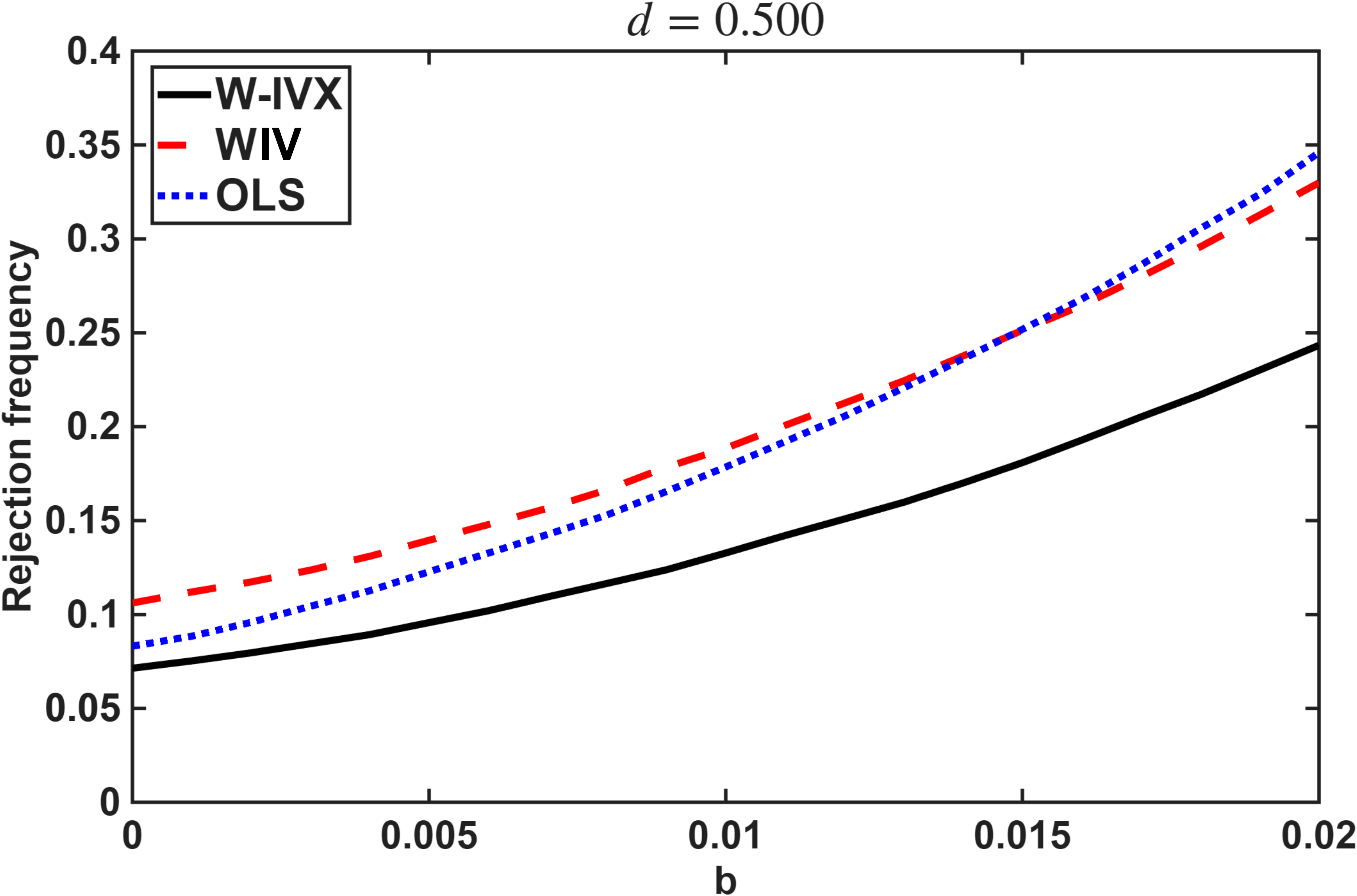}\\[0.6ex]
\includegraphics[width=0.49\linewidth]{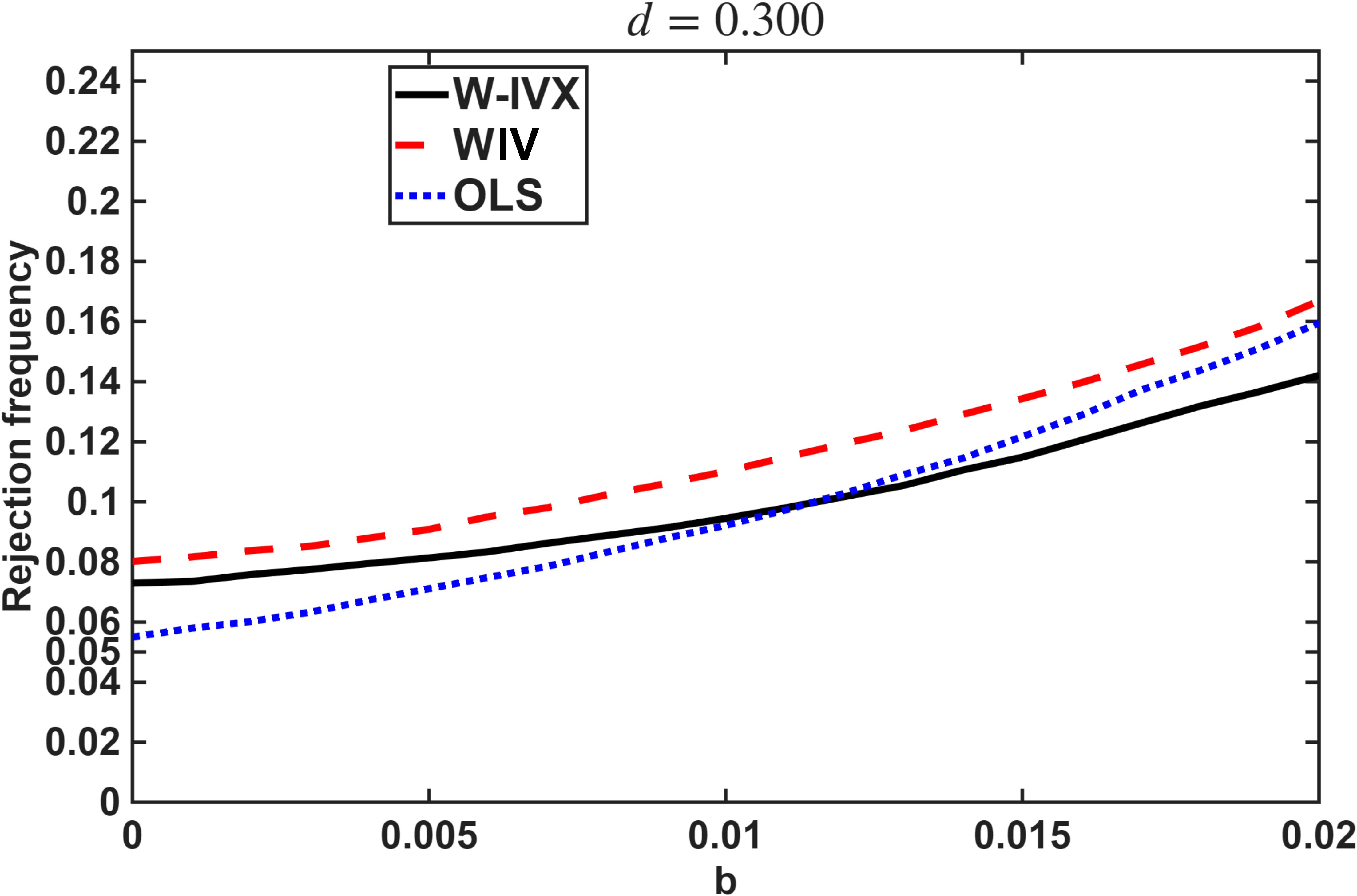}
\includegraphics[width=0.49\linewidth]{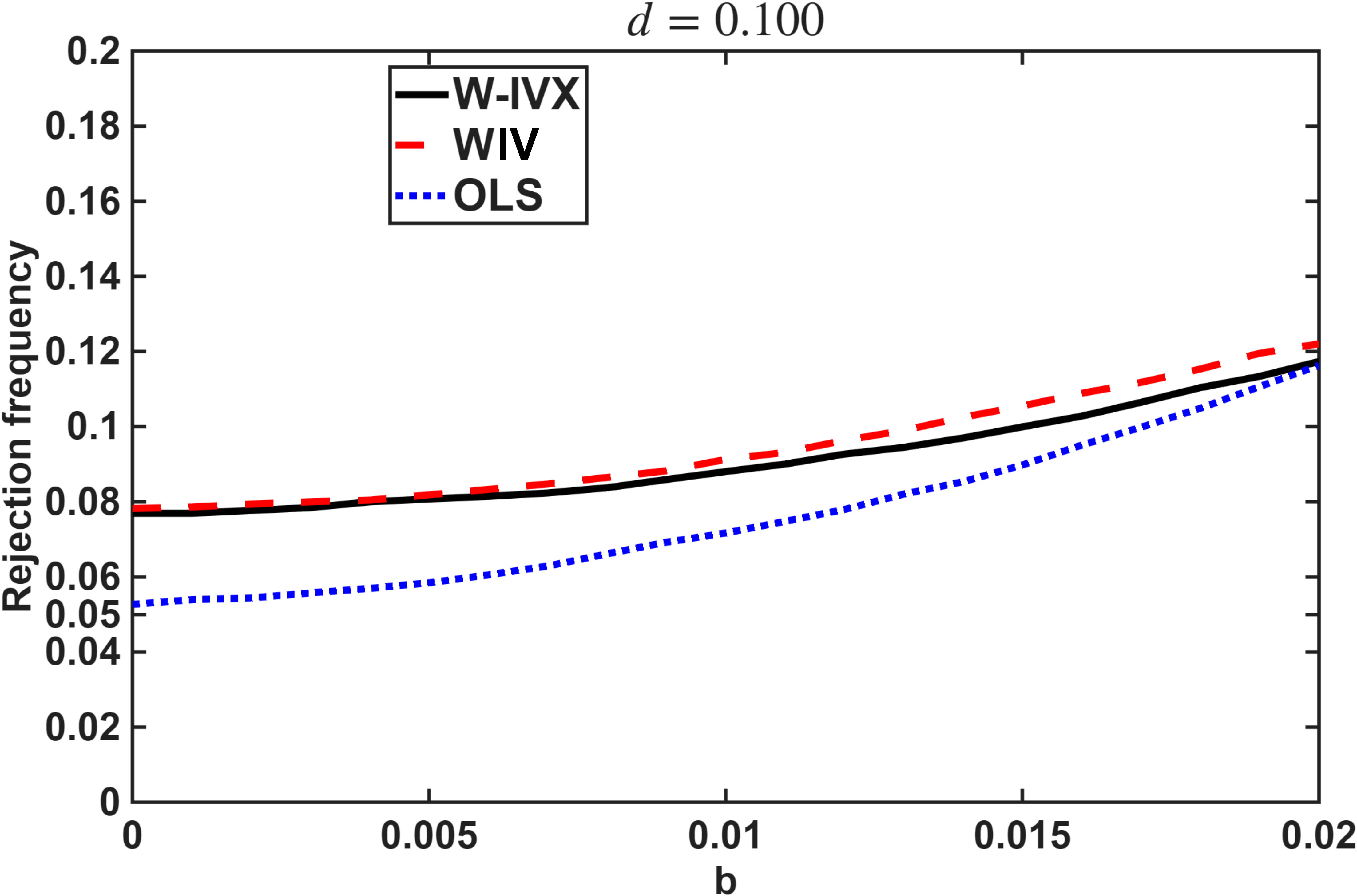}
\end{subfigure} \hfill{}
\begin{subfigure}[t]{0.48\linewidth} \centering \caption{$\alpha=1.8$, $d_{*}=0.444$.}
\label{fig:sim-power-linear-w05-alpha-1p80} \includegraphics[width=0.49\linewidth]{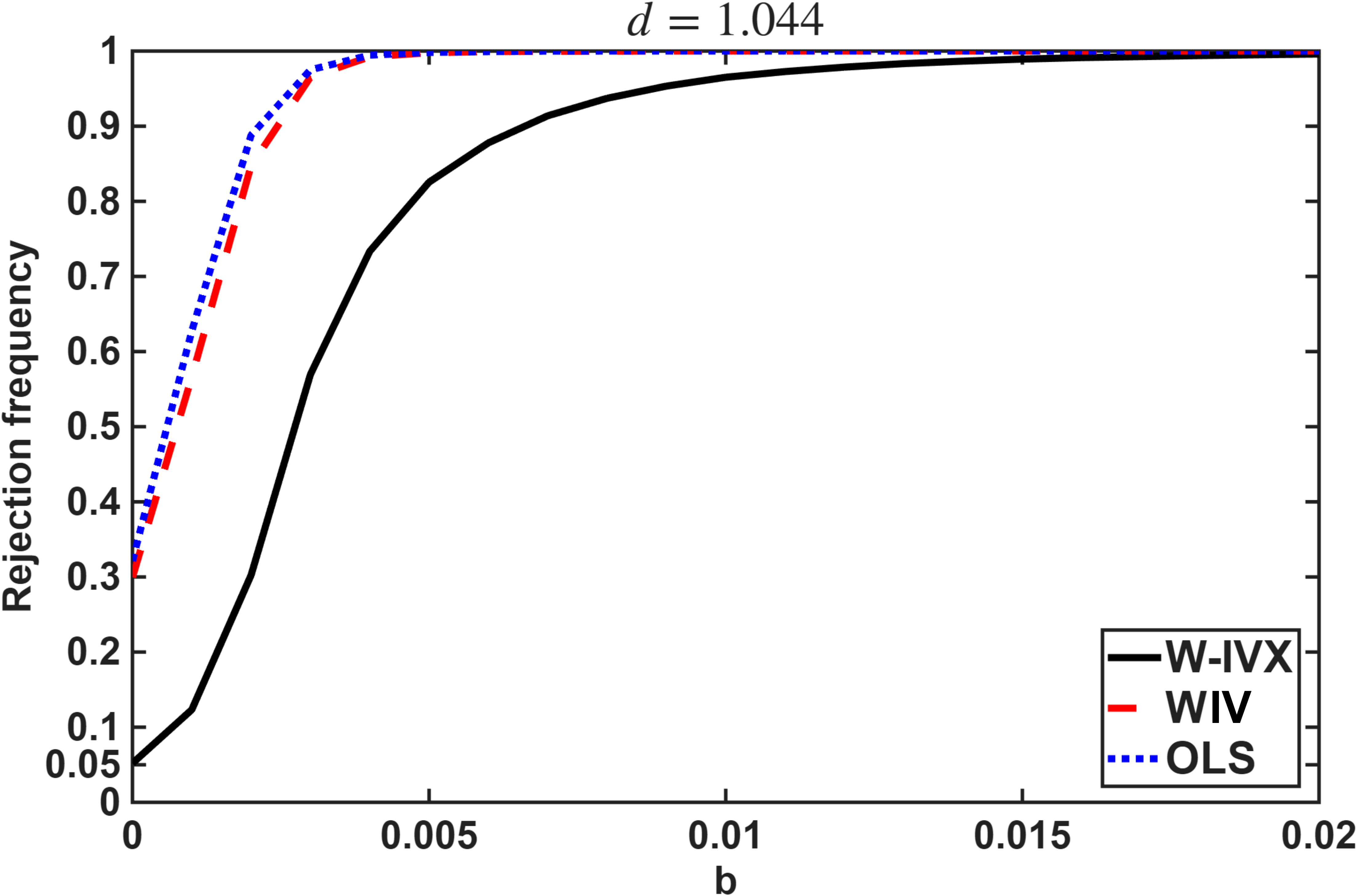}
\includegraphics[width=0.49\linewidth]{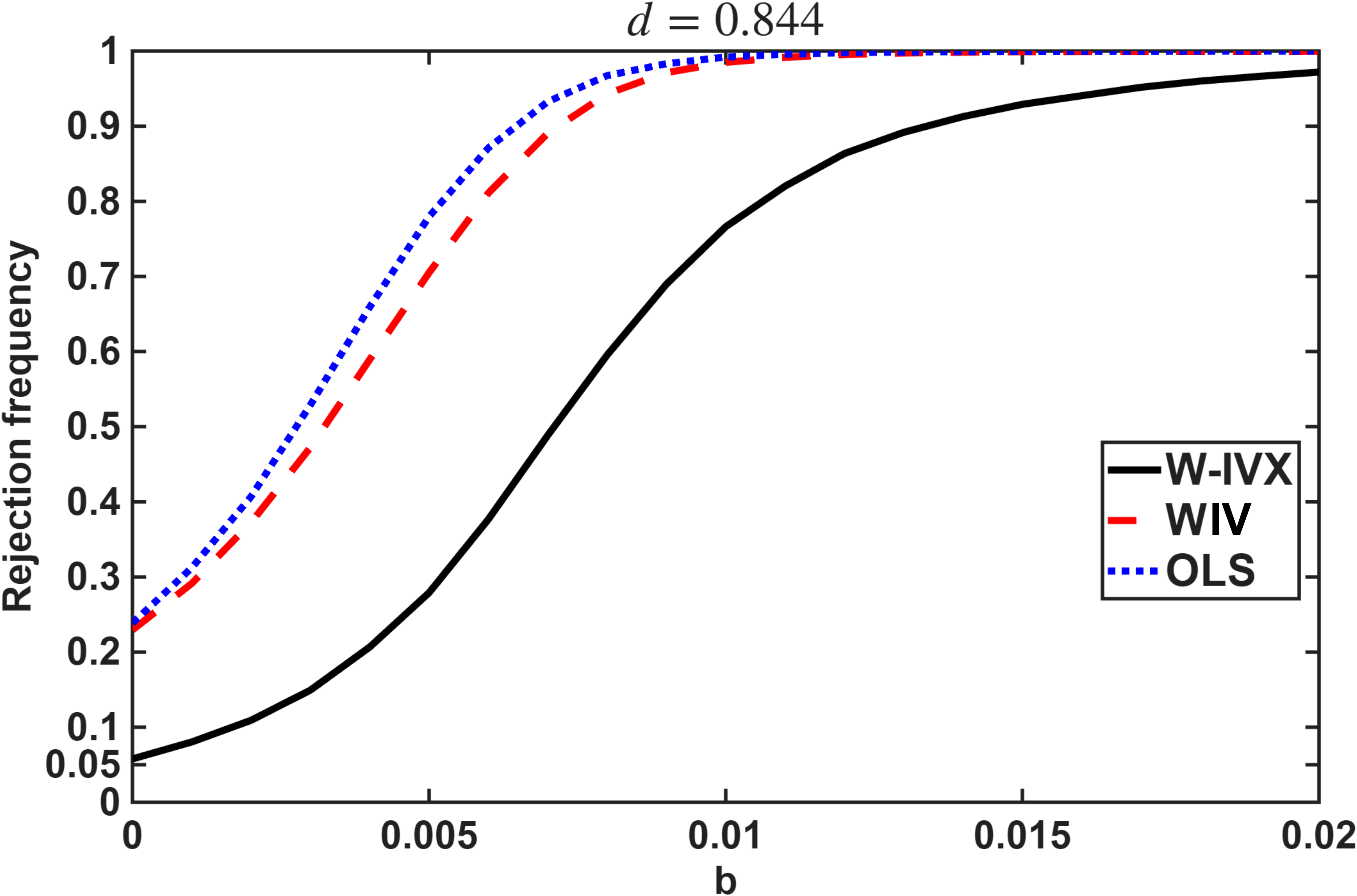}\\[0.6ex]
\includegraphics[width=0.49\linewidth]{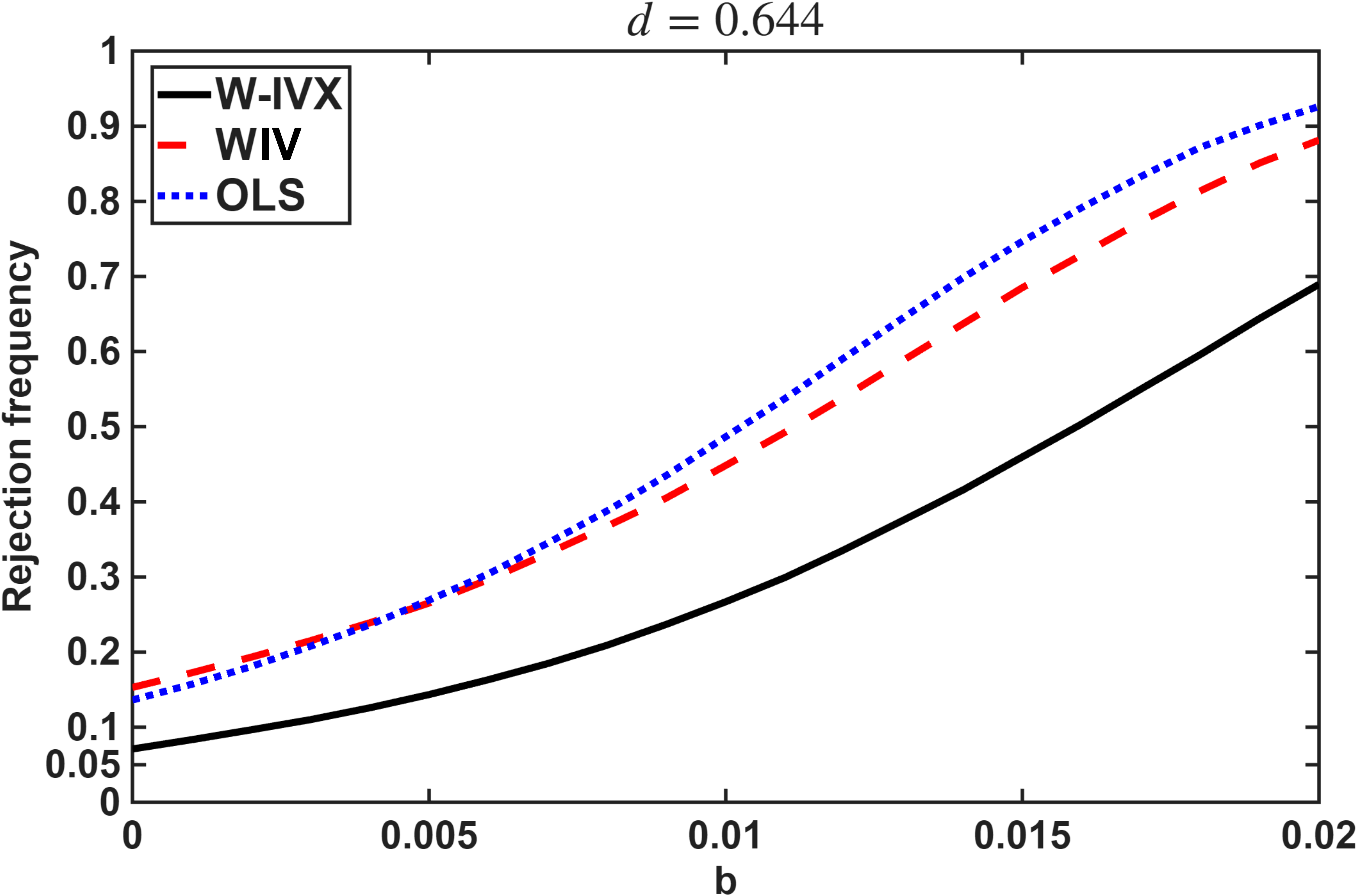}
\includegraphics[width=0.49\linewidth]{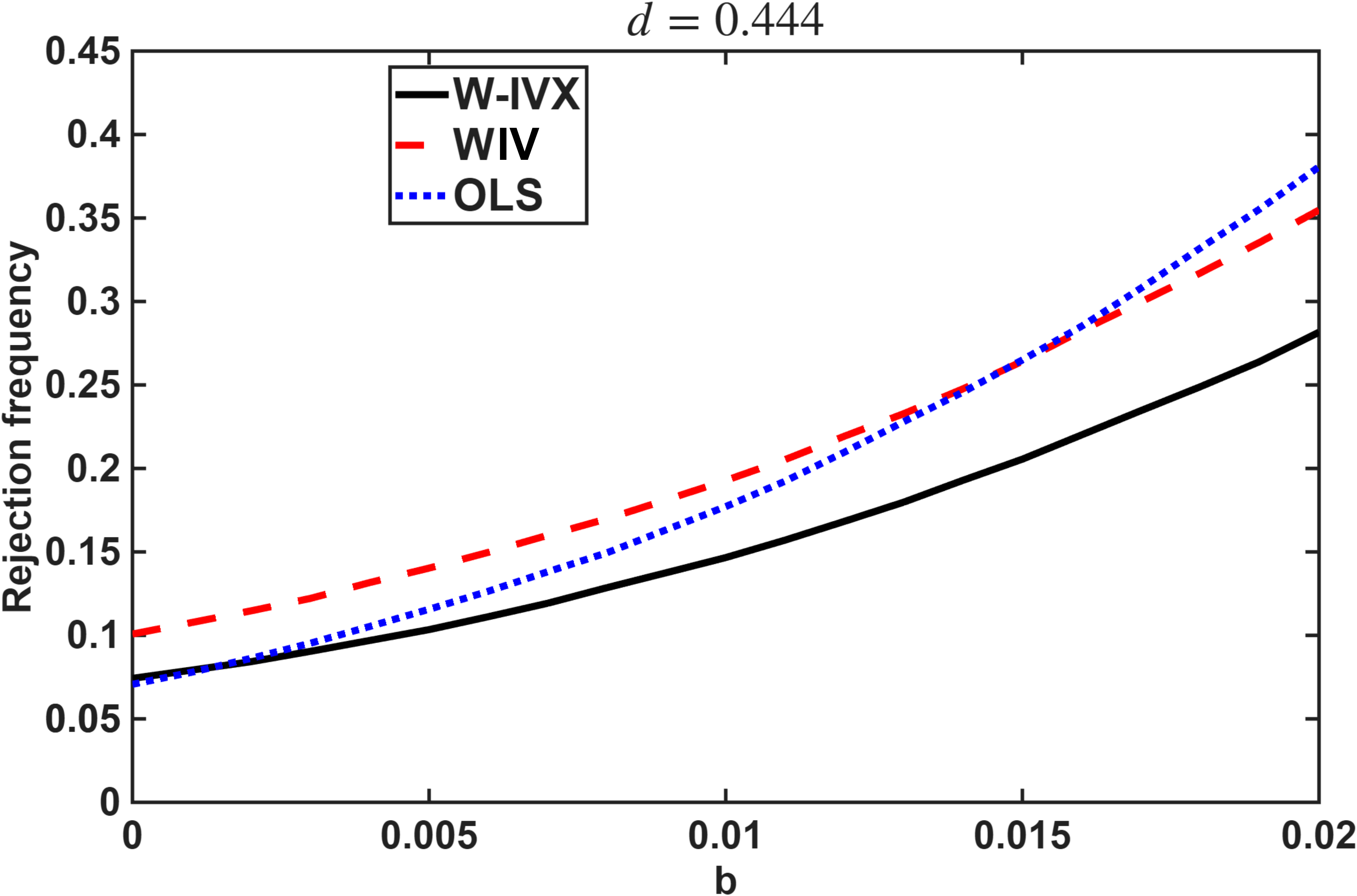}\\[0.6ex]
\includegraphics[width=0.49\linewidth]{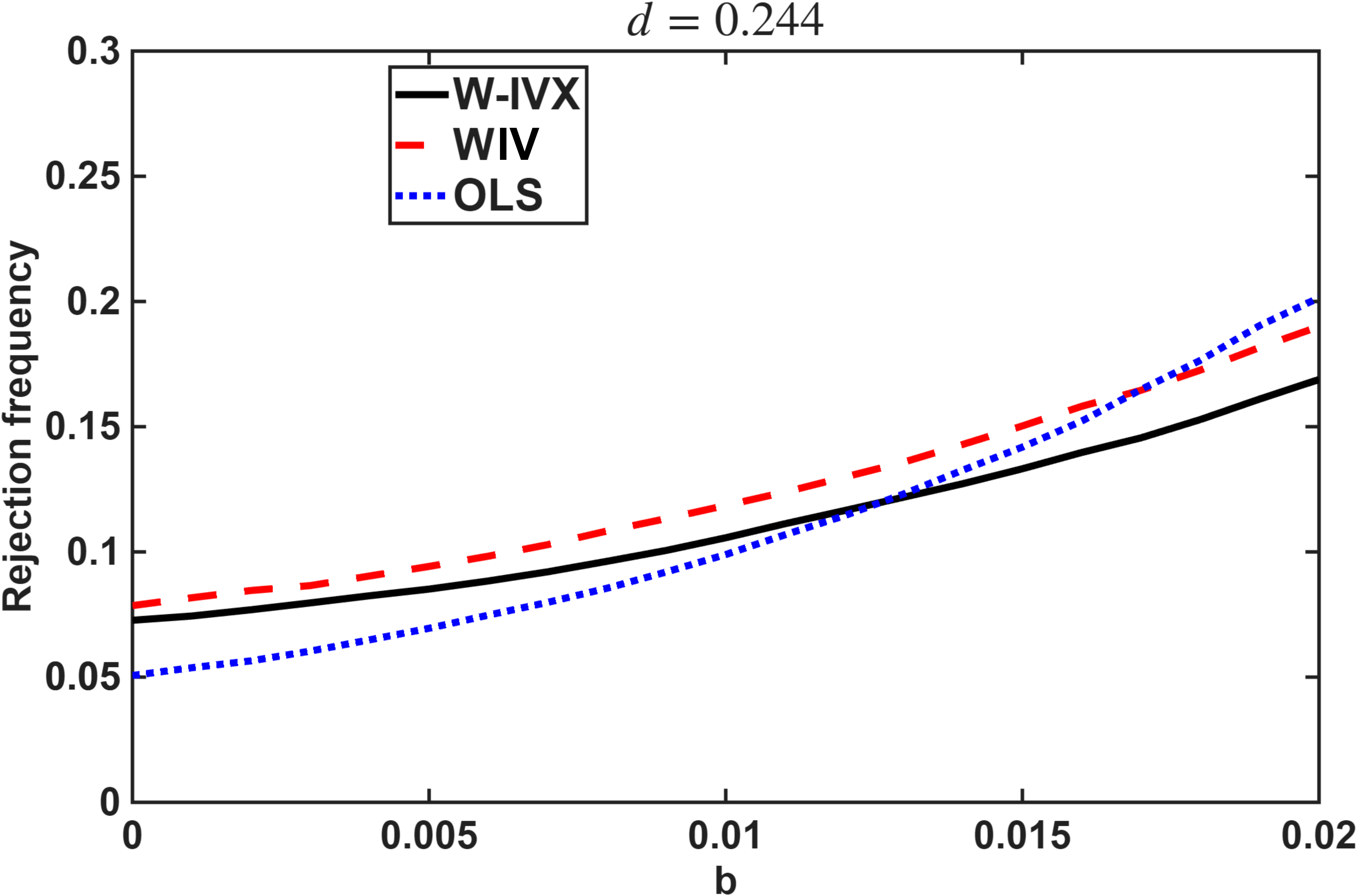}
\includegraphics[width=0.49\linewidth]{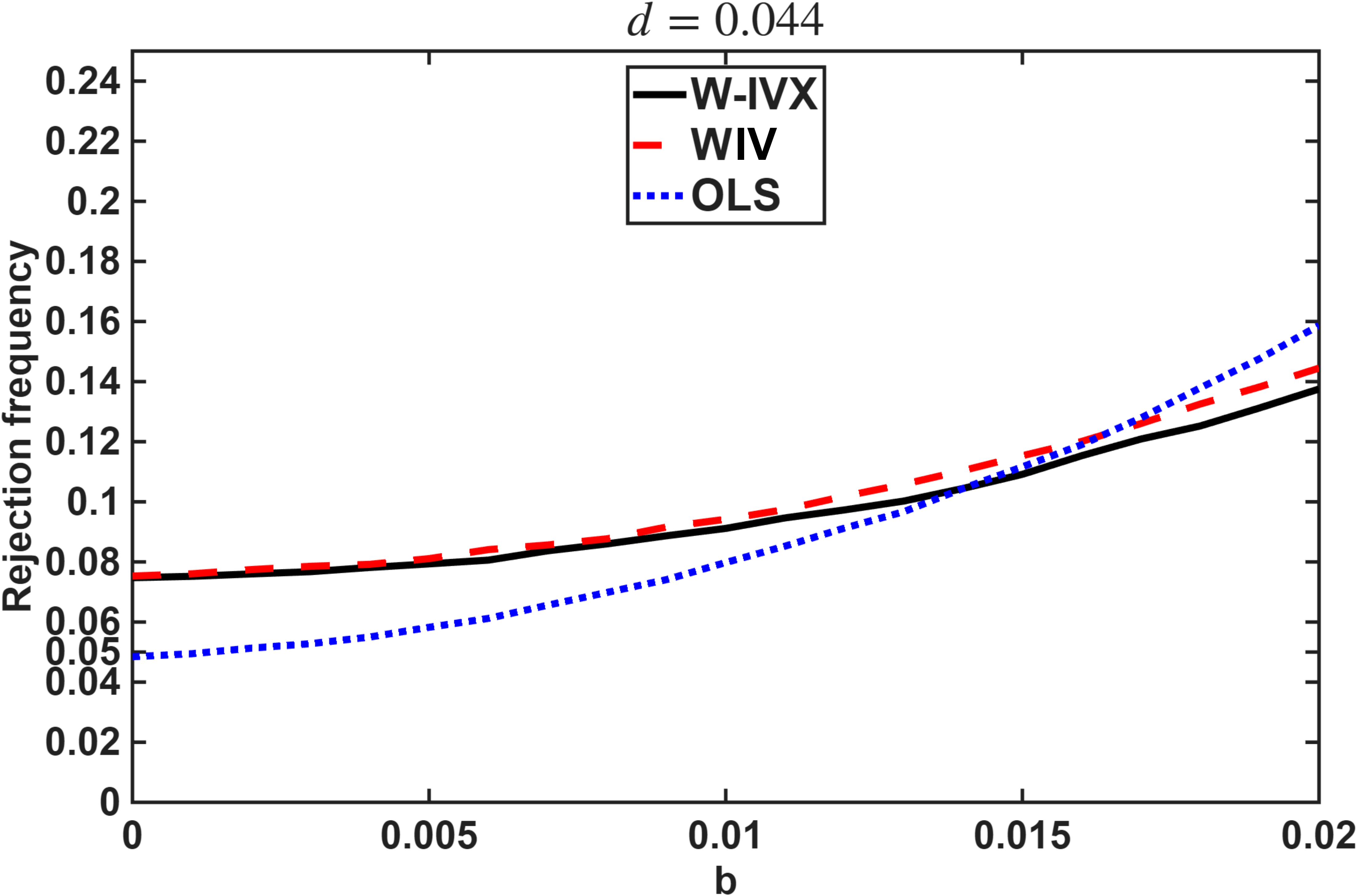}
\end{subfigure}

\end{figure}

\end{landscape}

\begin{landscape} 

\begin{figure}[htbp]
\caption{Asymptotic Power (nominal size 5\%): $[f(x),w(x)]=[x,x^{\langle0.75\rangle}]$,
$\delta=-0.95$, $n=1000$}
\label{fig:sim-power-linear-w075}\centering \captionsetup[subfigure]{position=top}

\begin{subfigure}[t]{0.48\linewidth} \centering \caption{$\alpha=2$, $d_{*}=0.5$}
\label{fig:sim-power-linear-w075-alpha-2p00} \includegraphics[width=0.49\linewidth]{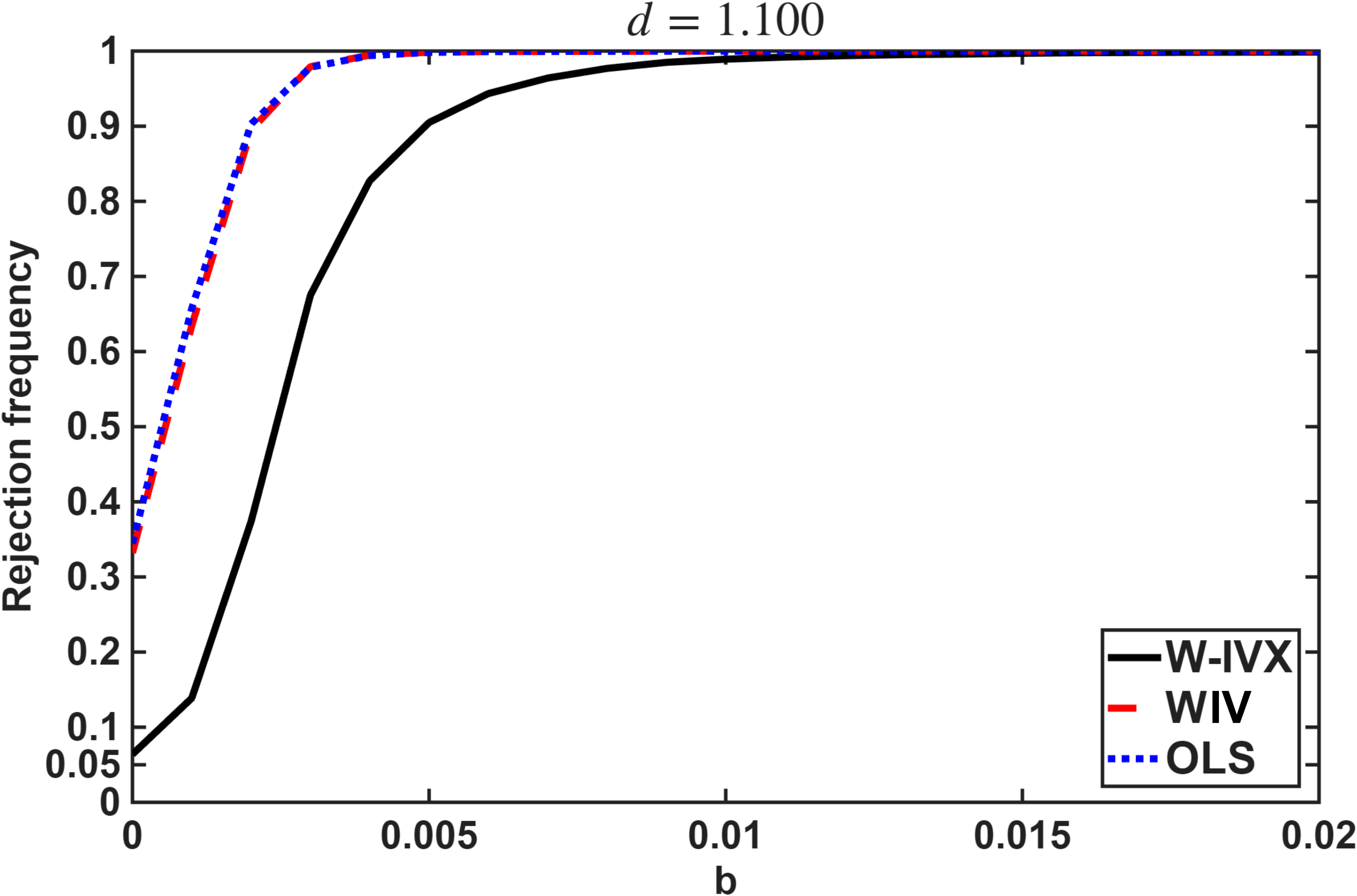}
\includegraphics[width=0.49\linewidth]{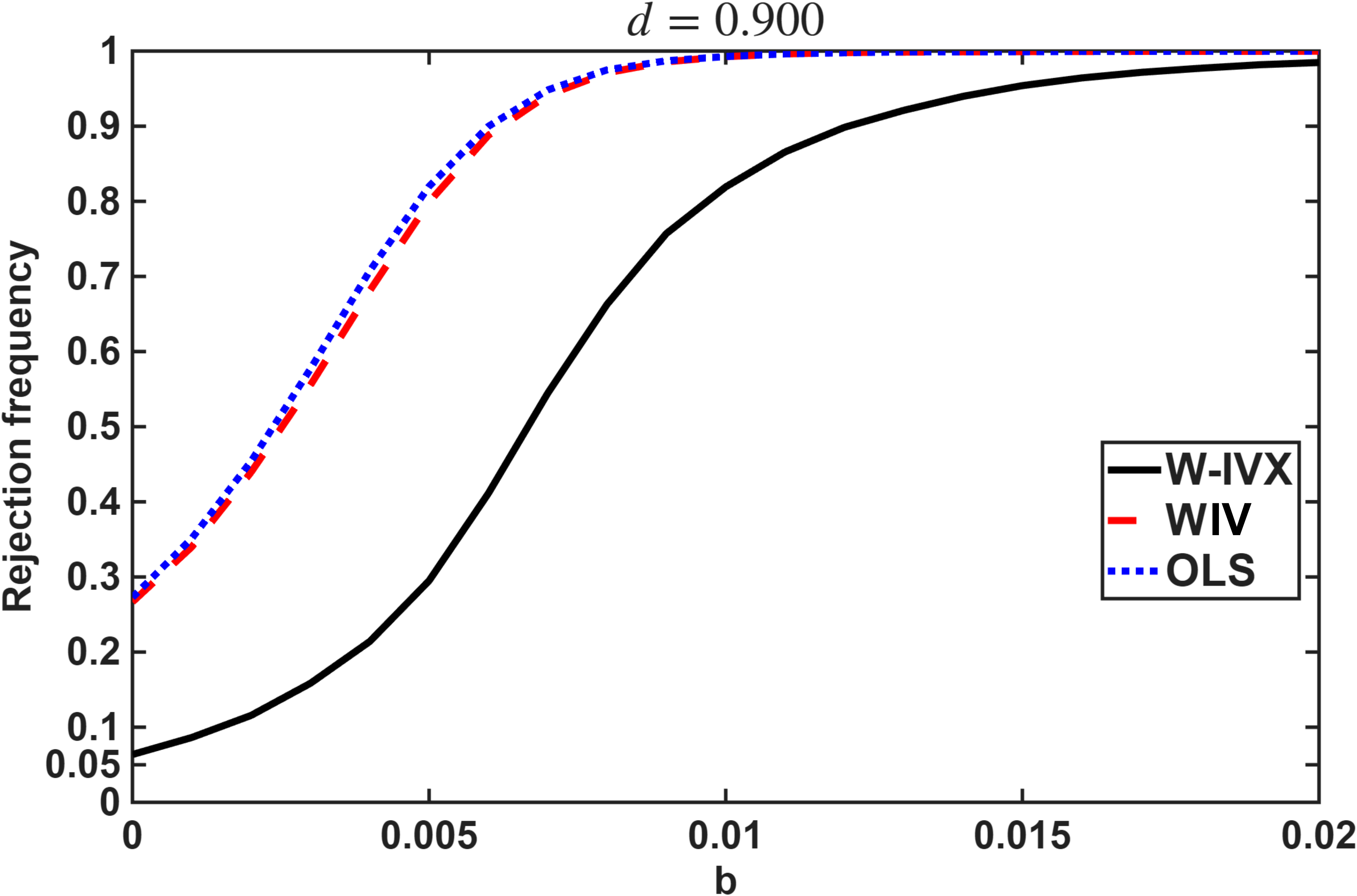}\\[0.6ex]
\includegraphics[width=0.49\linewidth]{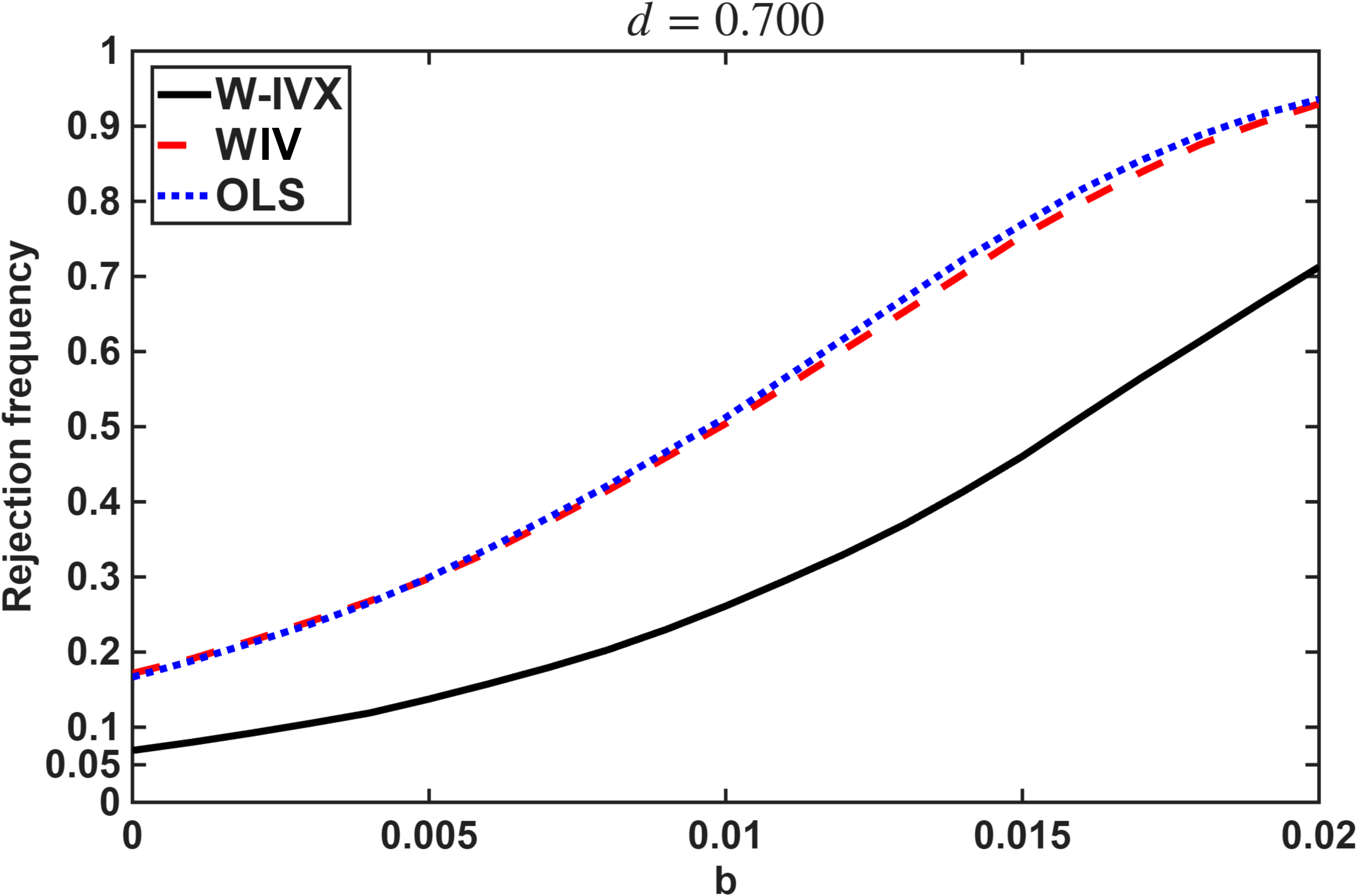}
\includegraphics[width=0.49\linewidth]{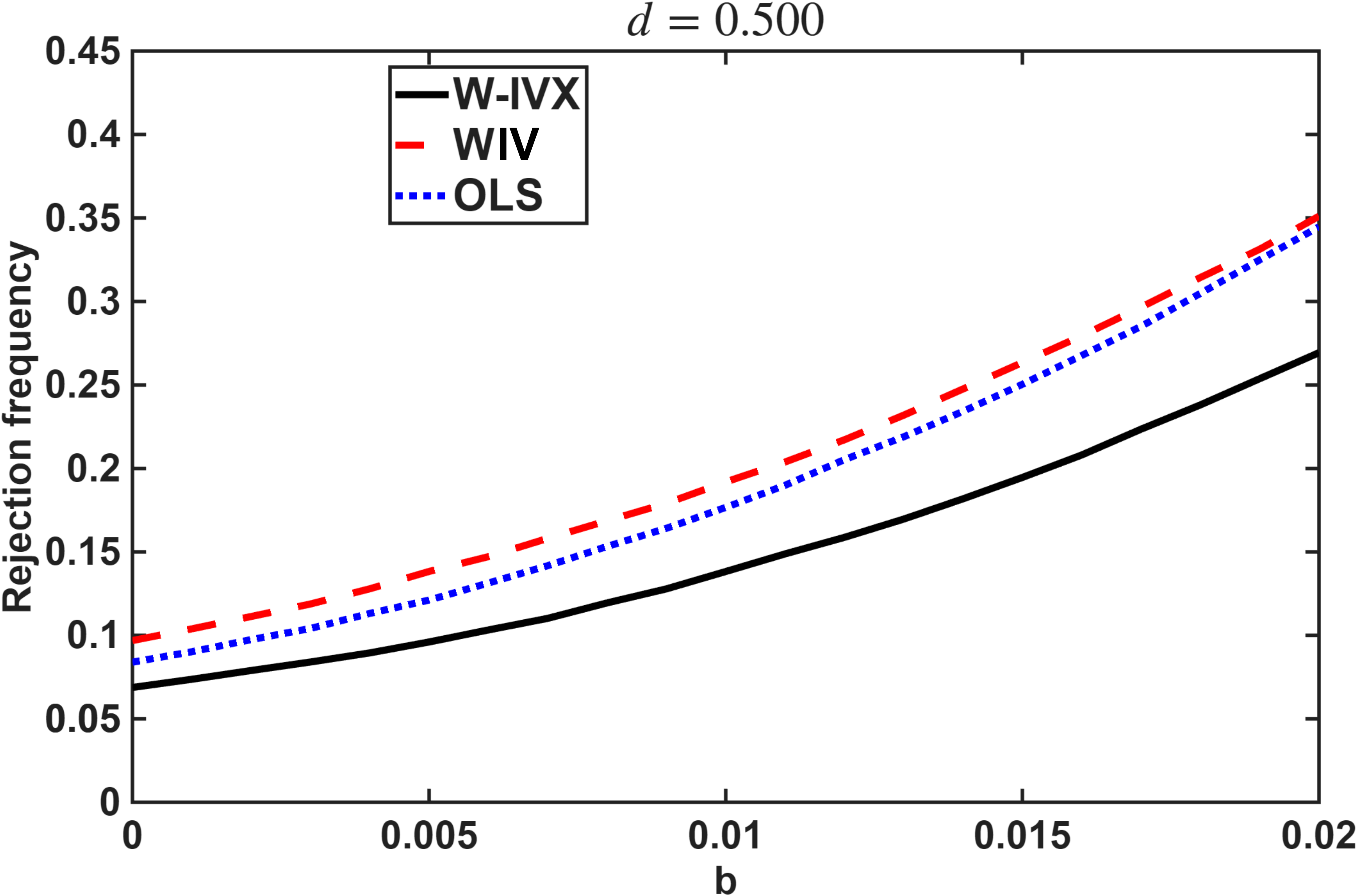}\\[0.6ex]
\includegraphics[width=0.49\linewidth]{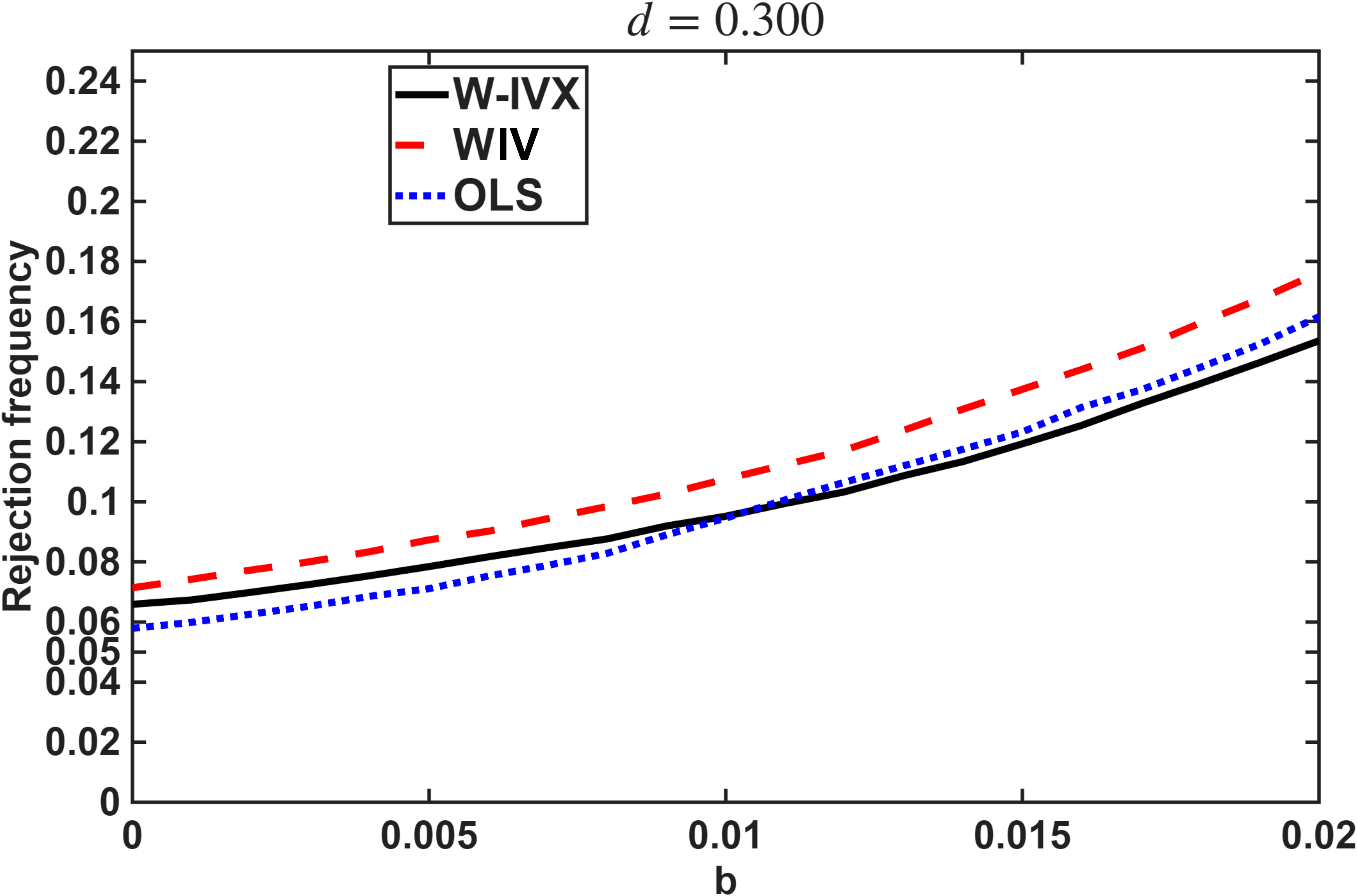}
\includegraphics[width=0.49\linewidth]{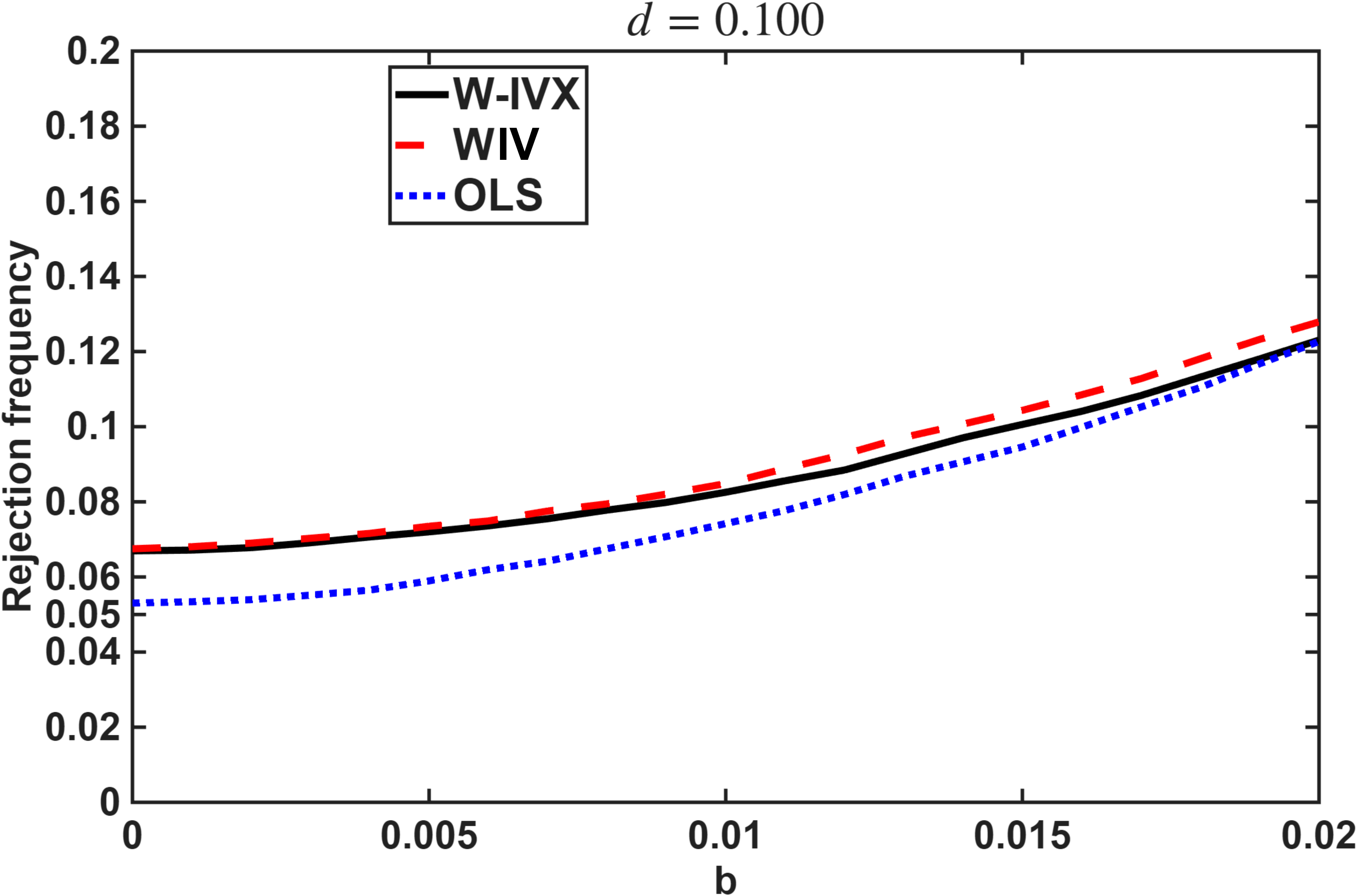}
\end{subfigure} \hfill{}
\begin{subfigure}[t]{0.48\linewidth} \centering \caption{$\alpha=1.8$, $d_{*}=0.444$.}
\label{fig:sim-power-linear-w075-alpha-1p80} \includegraphics[width=0.49\linewidth]{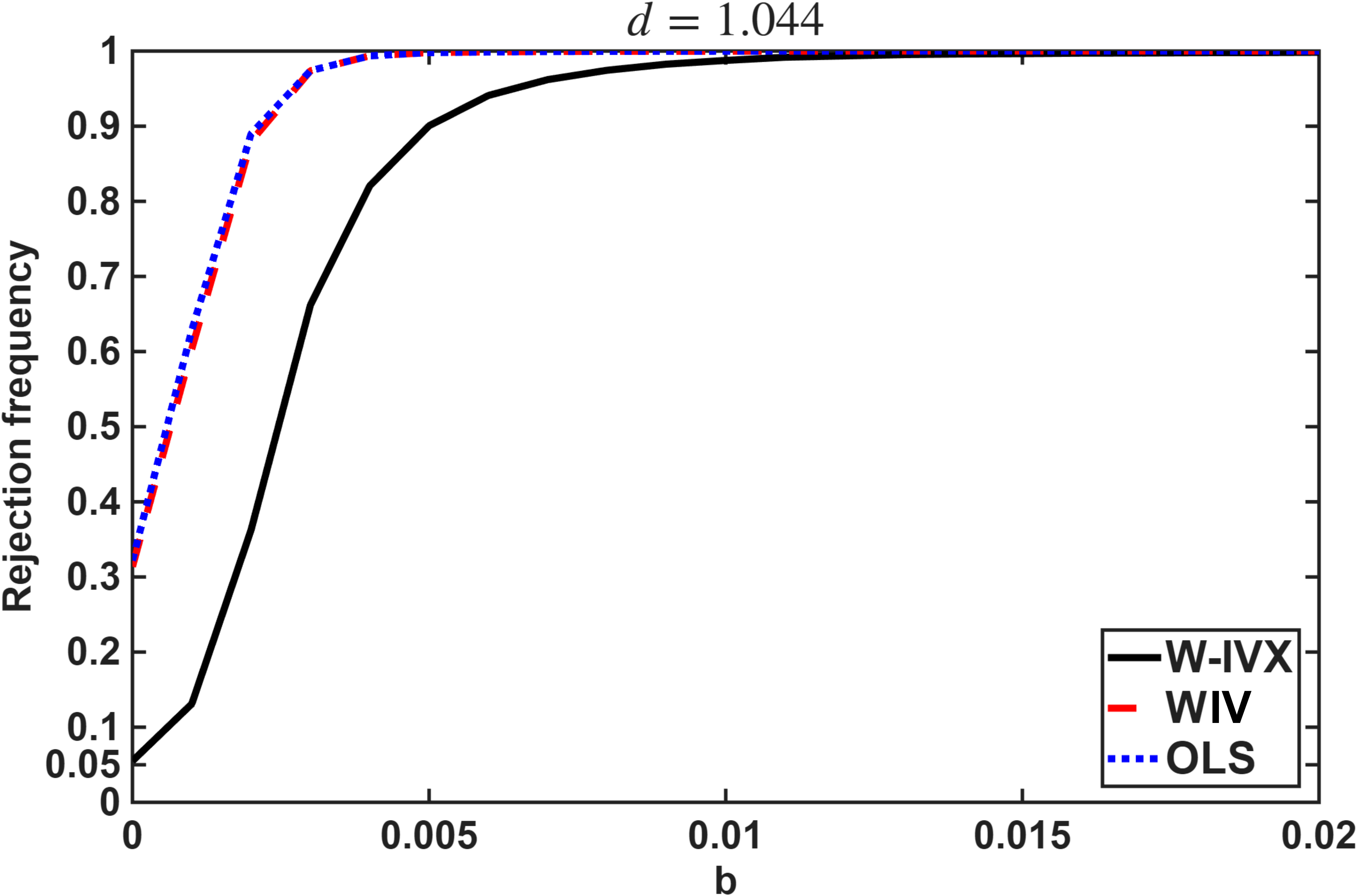}
\includegraphics[width=0.49\linewidth]{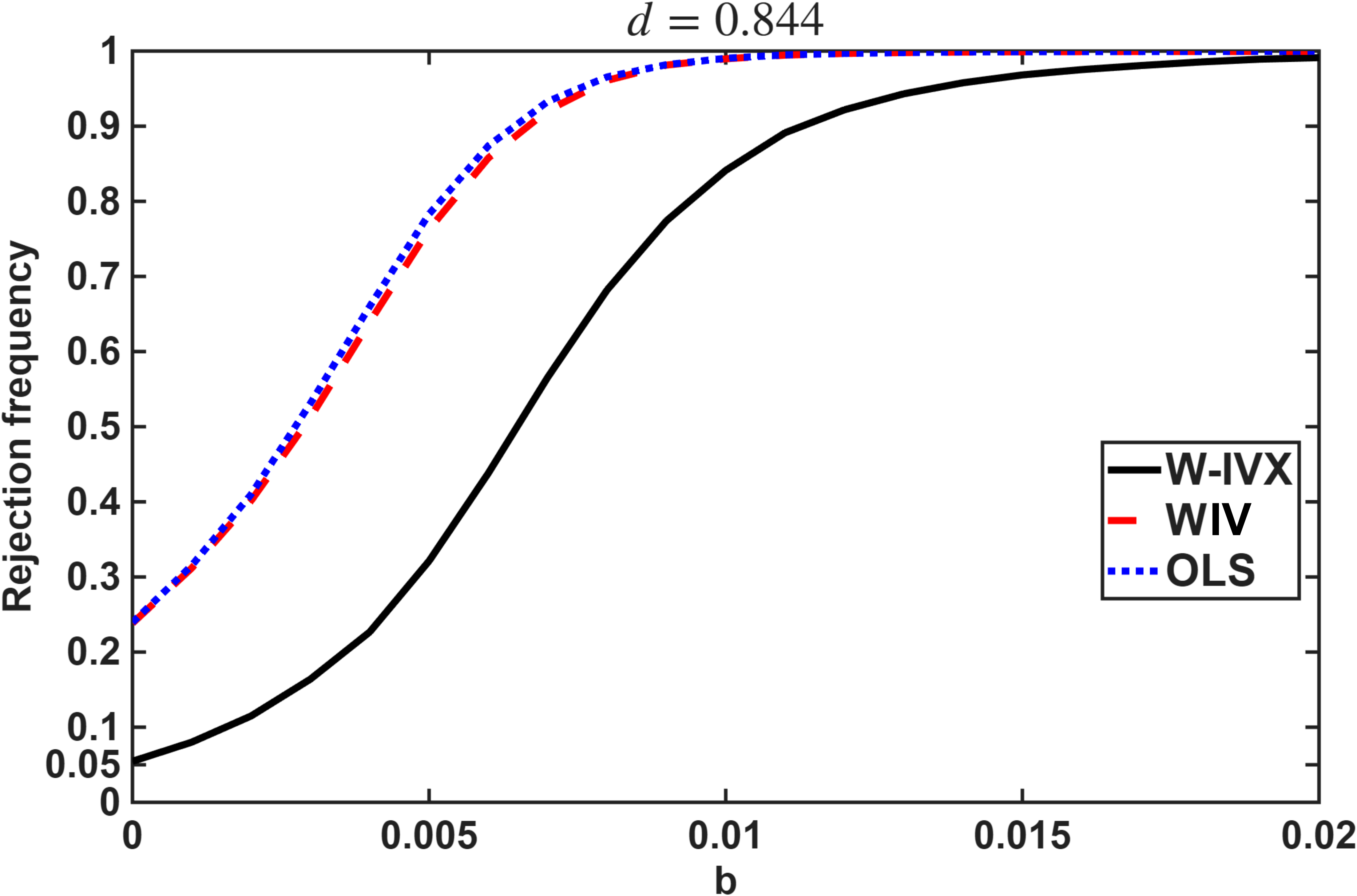}\\[0.6ex]
\includegraphics[width=0.49\linewidth]{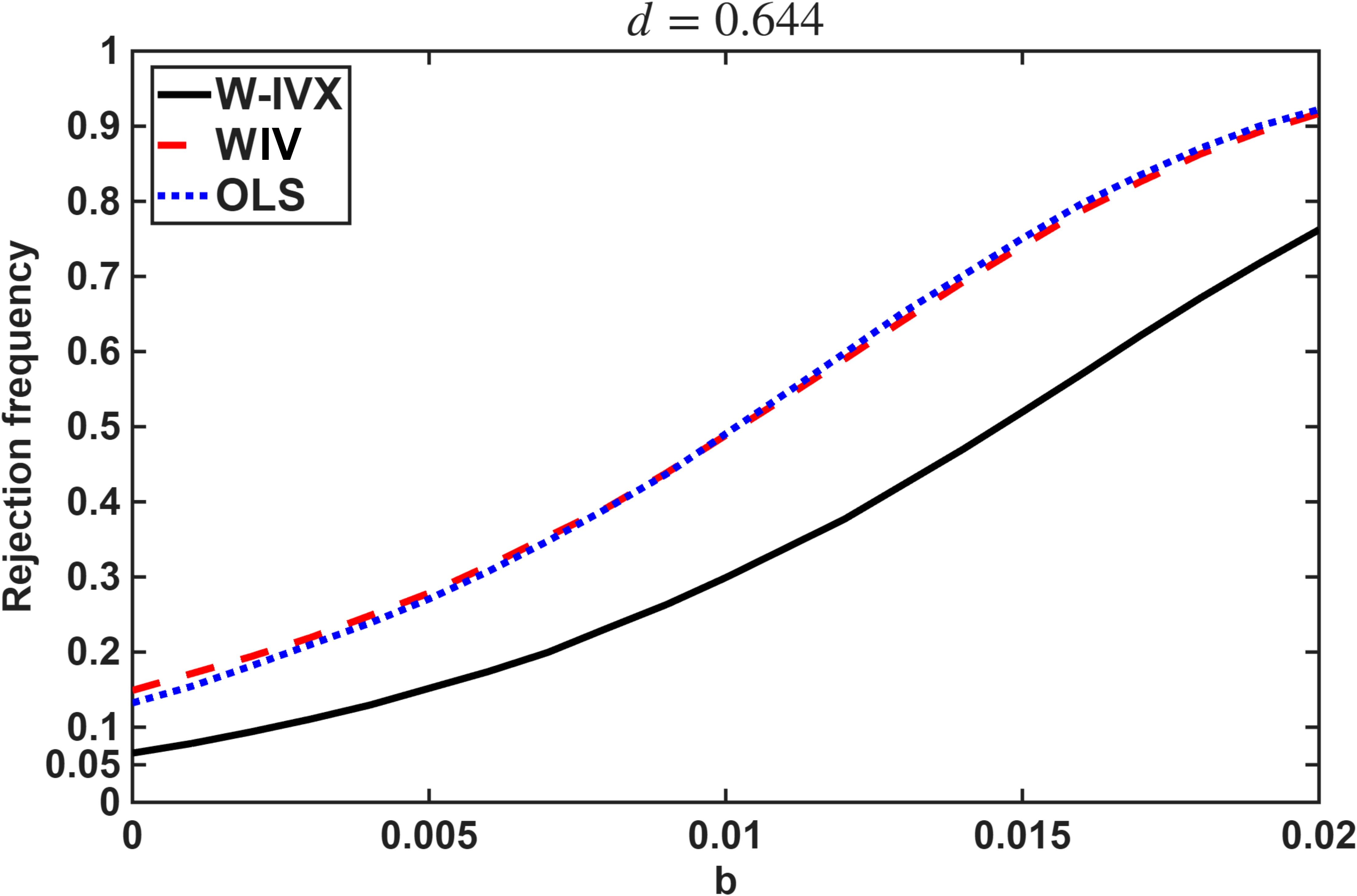}
\includegraphics[width=0.49\linewidth]{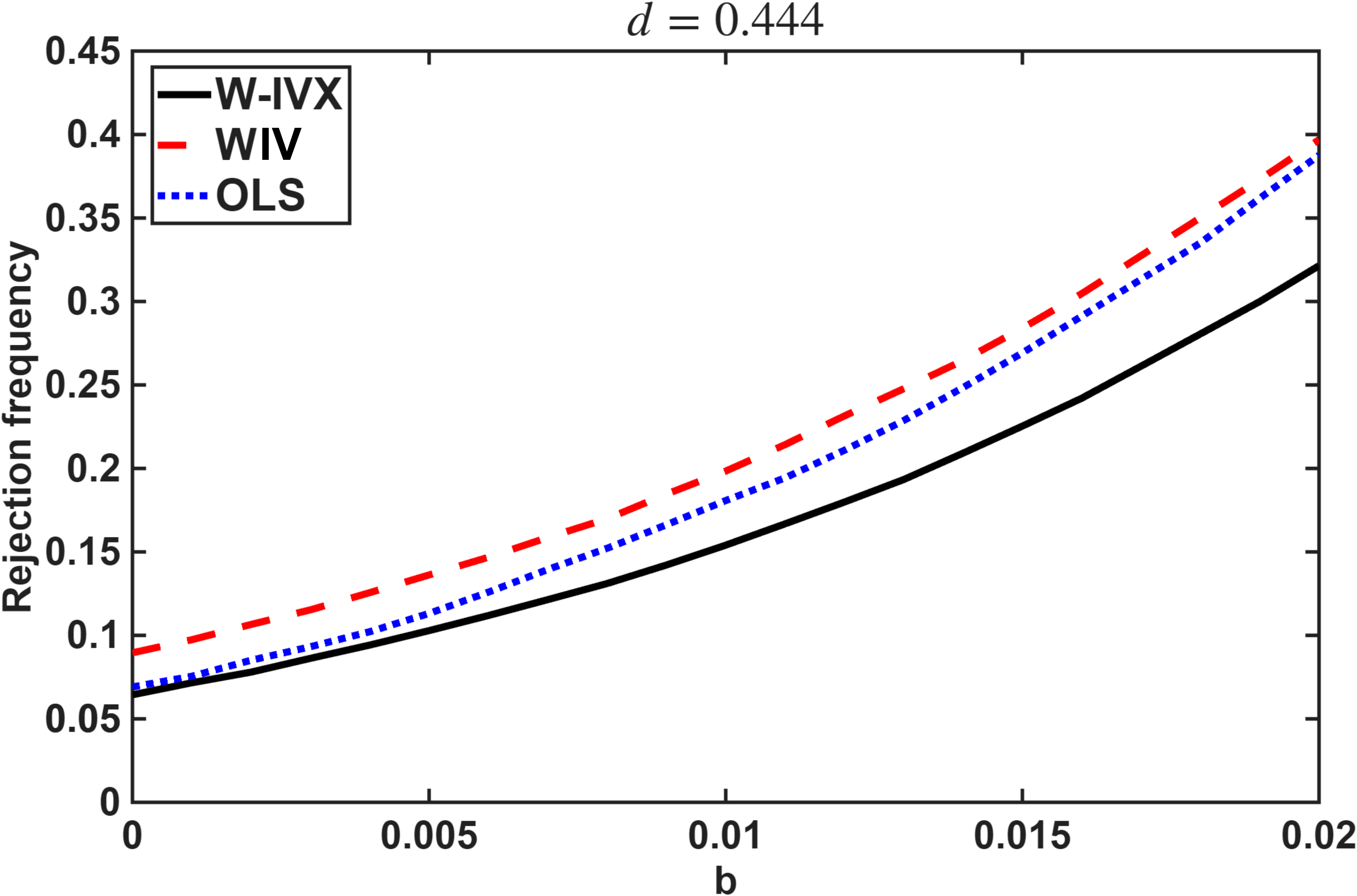}\\[0.6ex]
\includegraphics[width=0.49\linewidth]{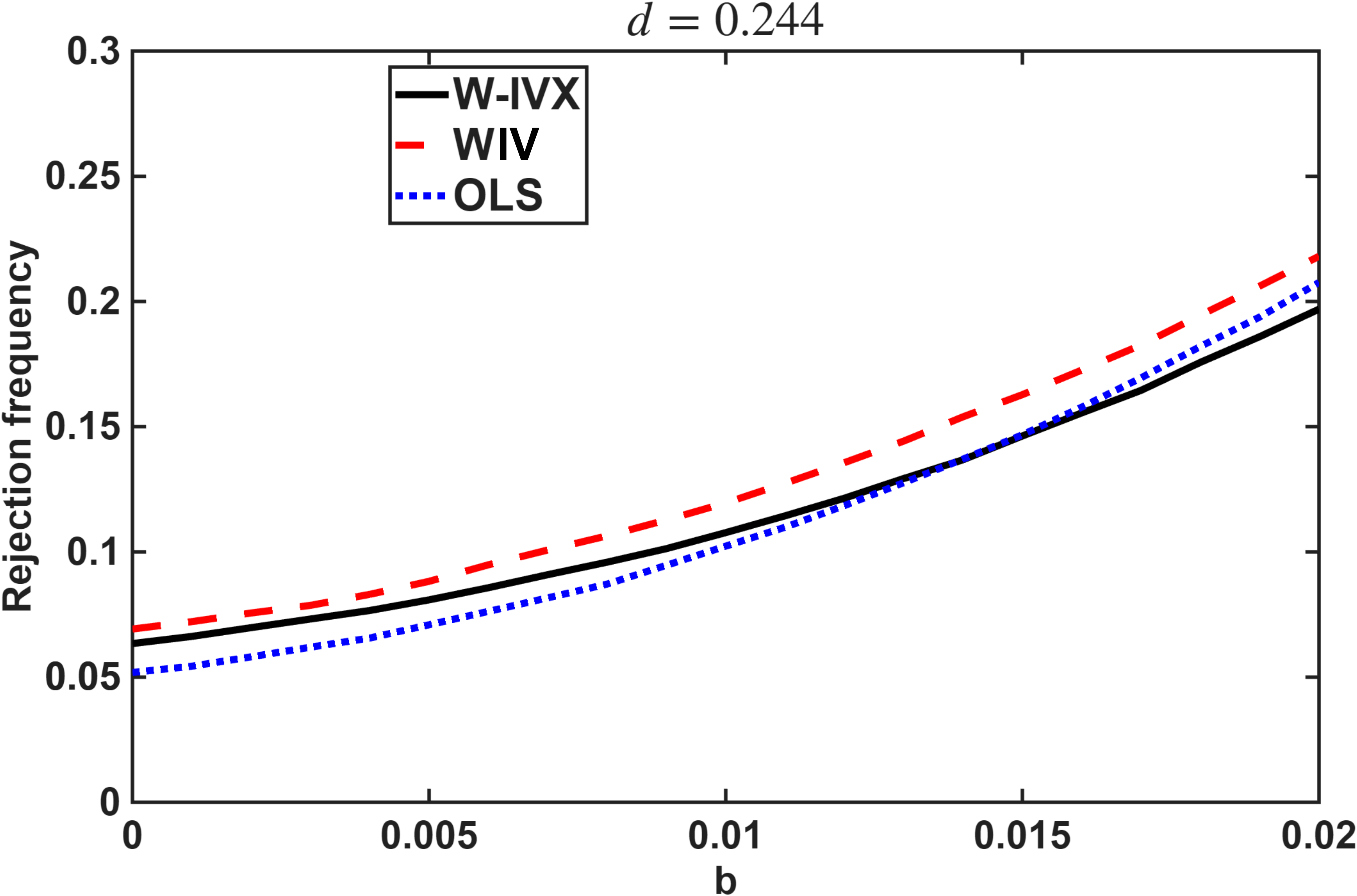}
\includegraphics[width=0.49\linewidth]{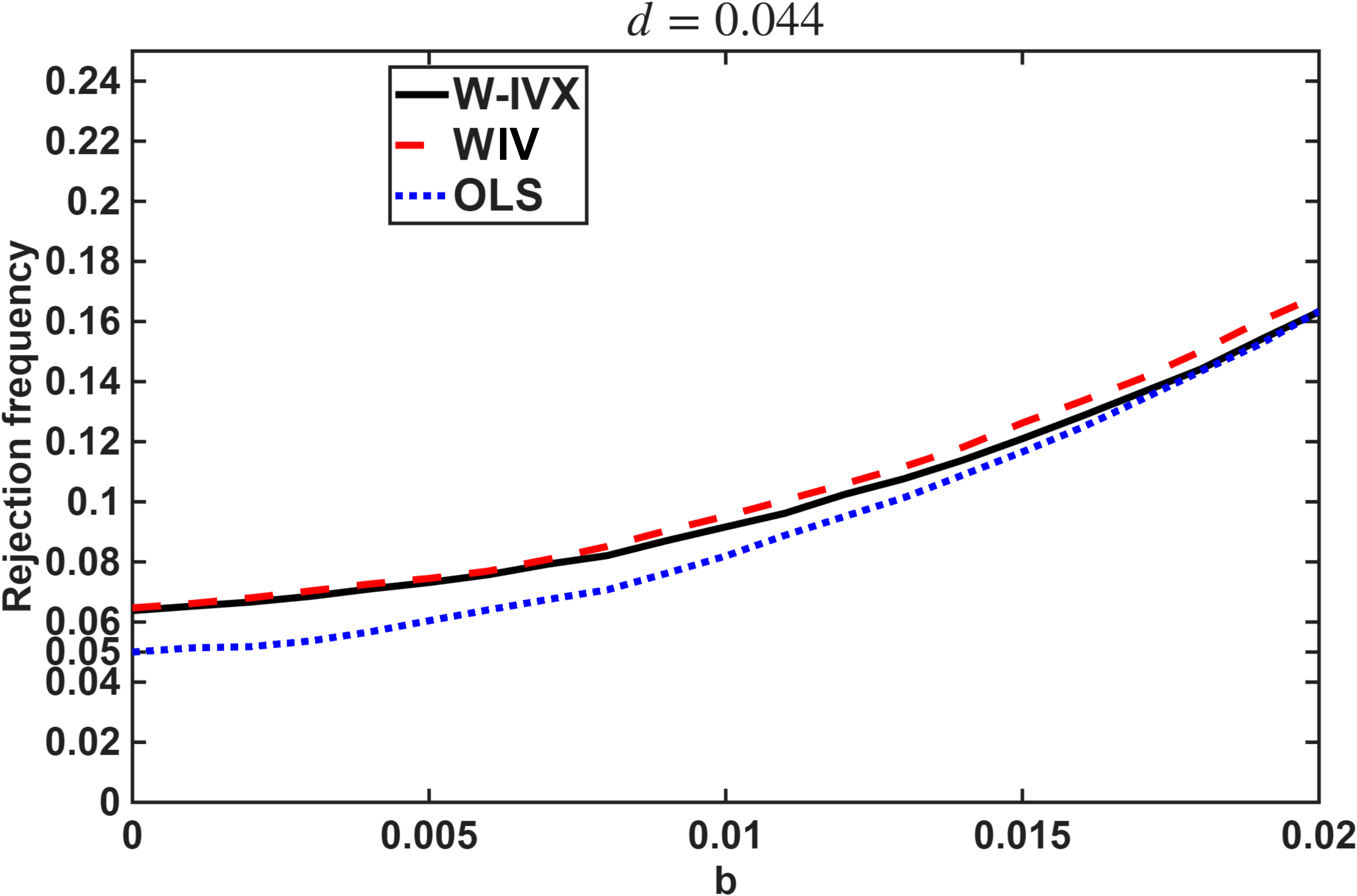}
\end{subfigure}

\end{figure}

\end{landscape}

\begin{landscape} 

\begin{figure}[htbp]
\caption{Asymptotic Power (nominal size 5\%): $[f(x),w(x)]=[x,x]$, $\delta=-0.95$,
$n=1000$}
\label{fig:sim-power-linear-w10}\centering \captionsetup[subfigure]{position=top}

\begin{subfigure}[t]{0.48\linewidth} \centering \caption{$\alpha=2$, $d_{*}=0.5$}
\label{fig:sim-power-linear-w10-alpha-2p00} \includegraphics[width=0.49\linewidth]{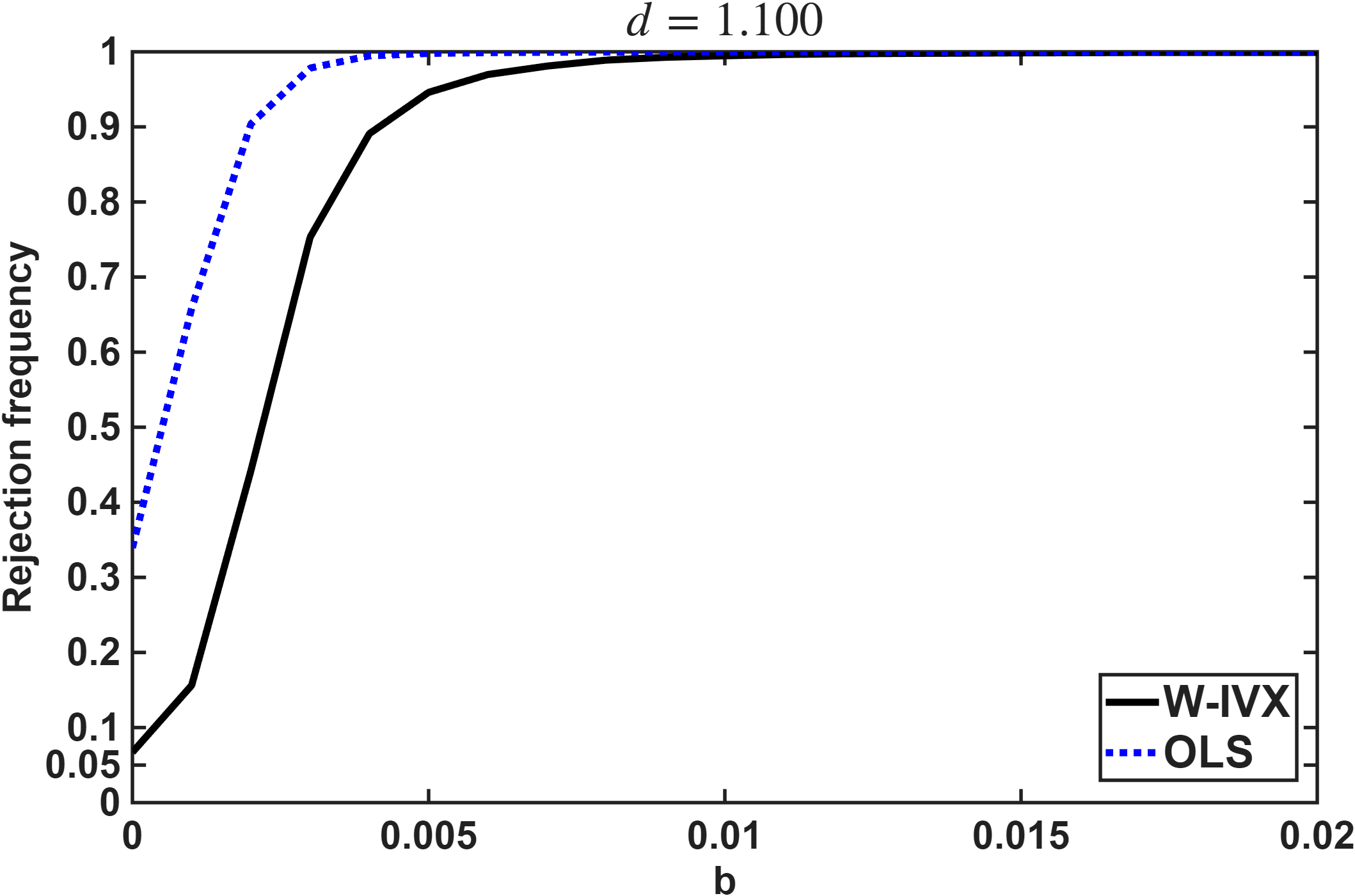}
\includegraphics[width=0.49\linewidth]{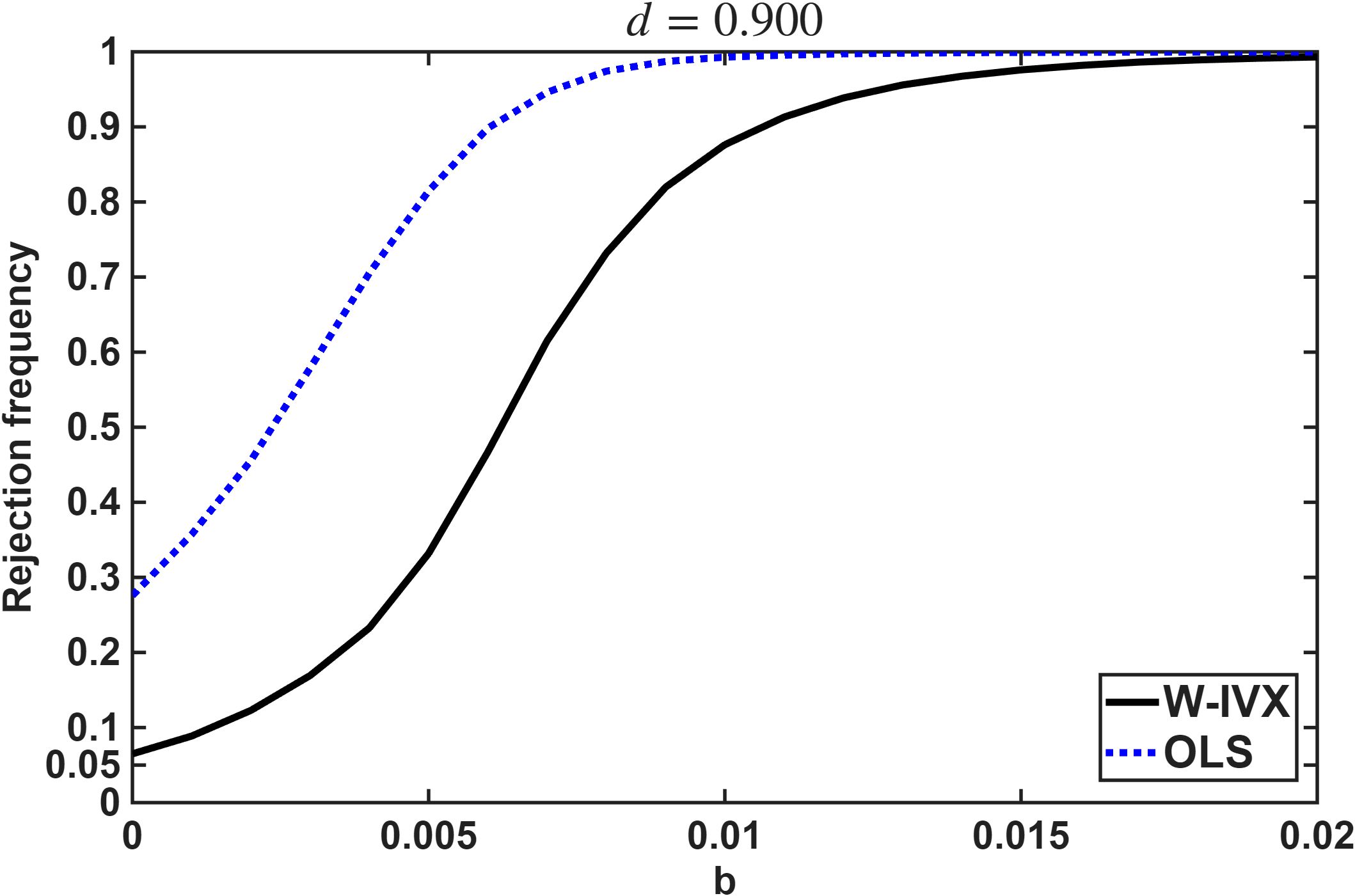}\\[0.6ex]
\includegraphics[width=0.49\linewidth]{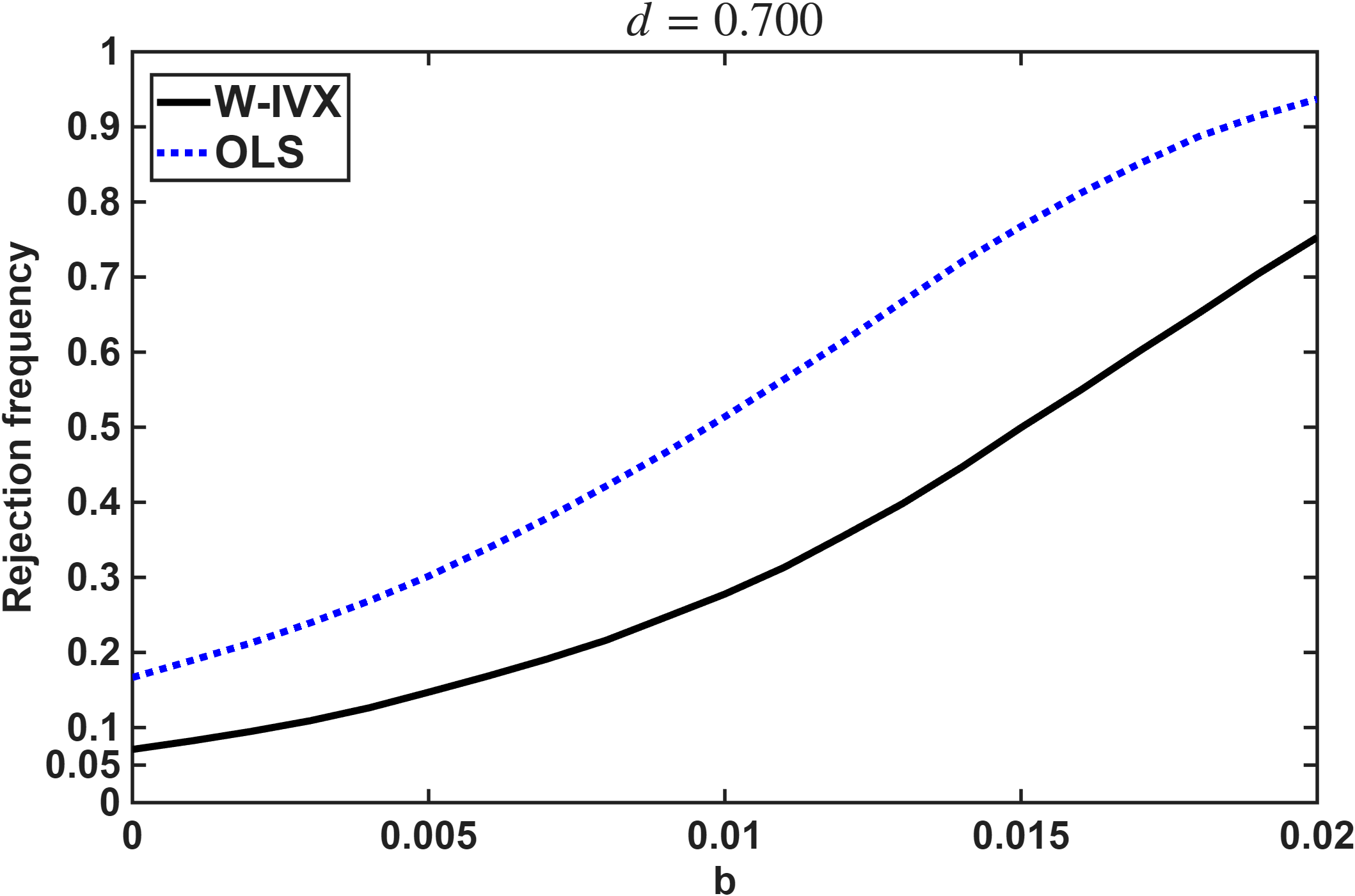}
\includegraphics[width=0.49\linewidth]{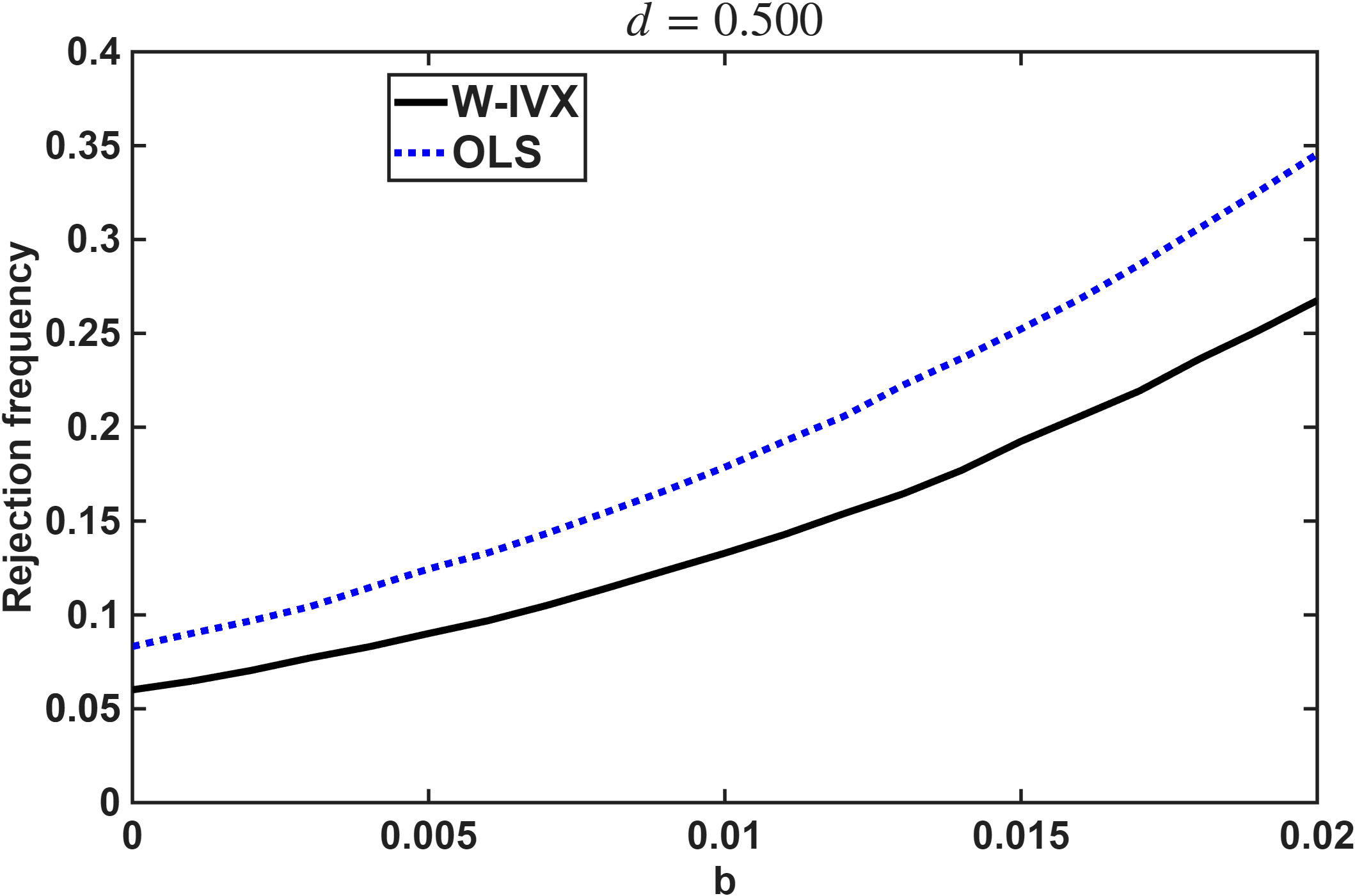}\\[0.6ex]
\includegraphics[width=0.49\linewidth]{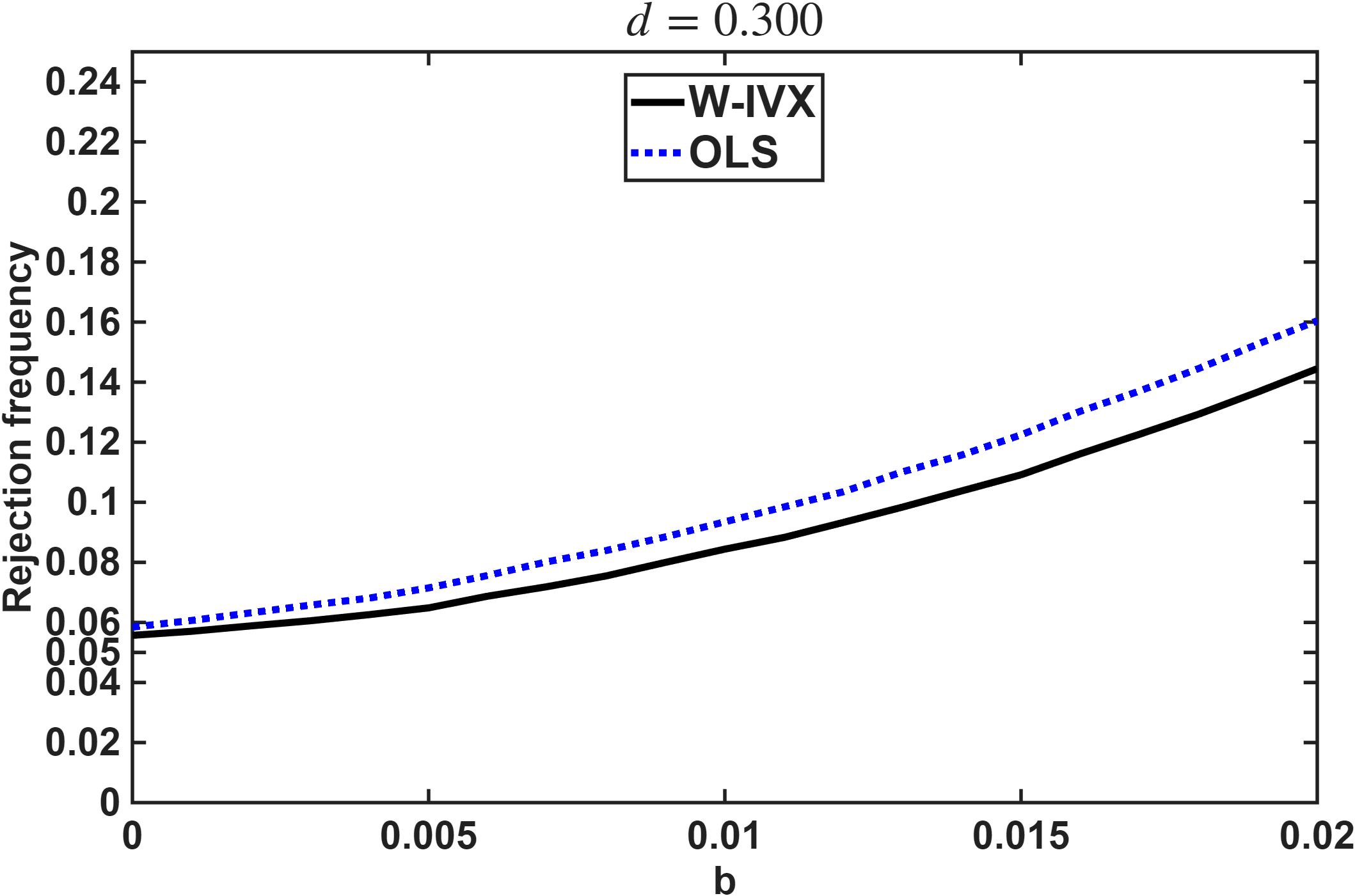}
\includegraphics[width=0.49\linewidth]{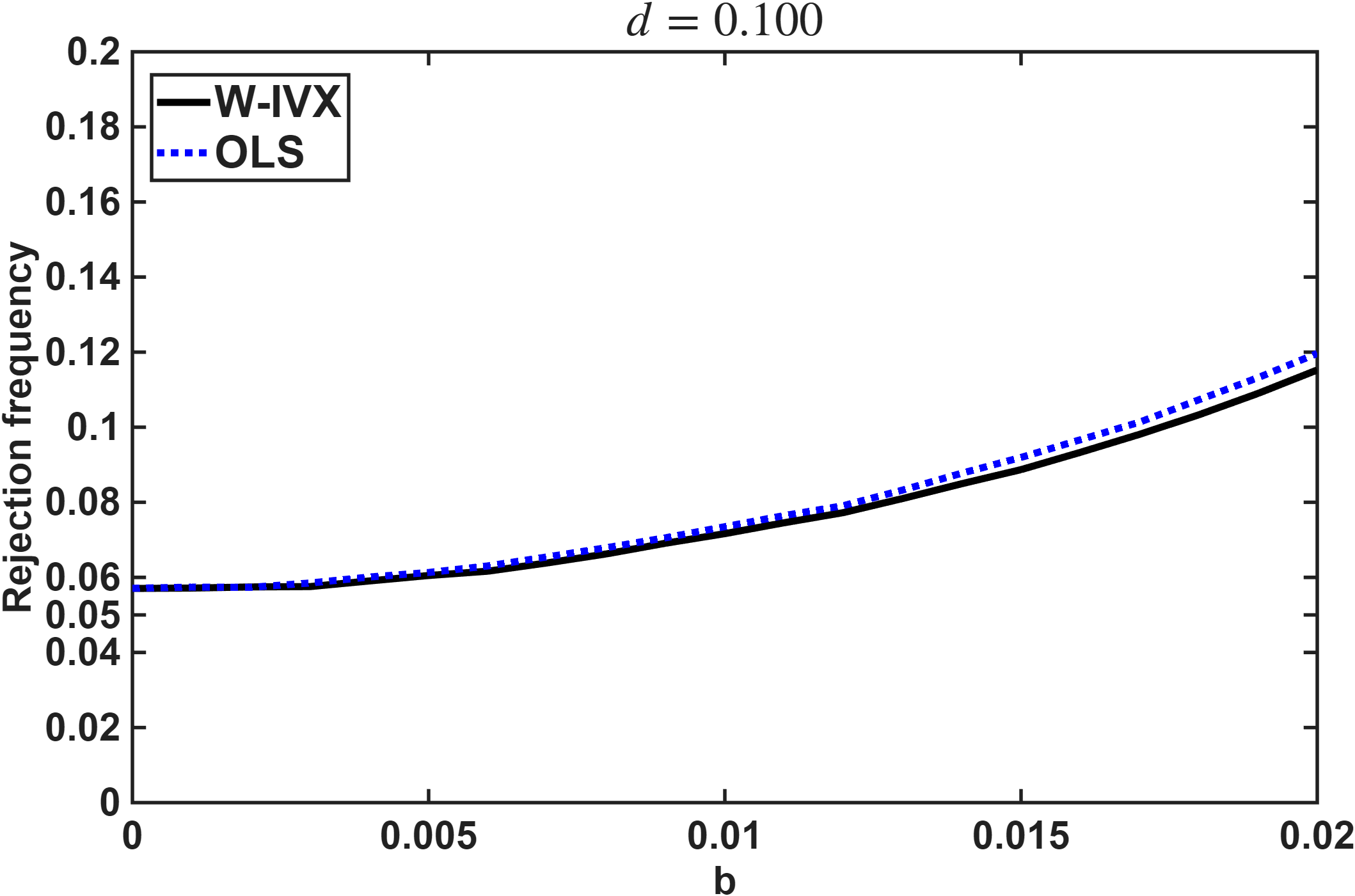}
\end{subfigure} \hfill{}
\begin{subfigure}[t]{0.48\linewidth} \centering \caption{$\alpha=1.8$, $d_{*}=0.444$.}
\label{fig:sim-power-linear-w10-alpha-1p80} \includegraphics[width=0.49\linewidth]{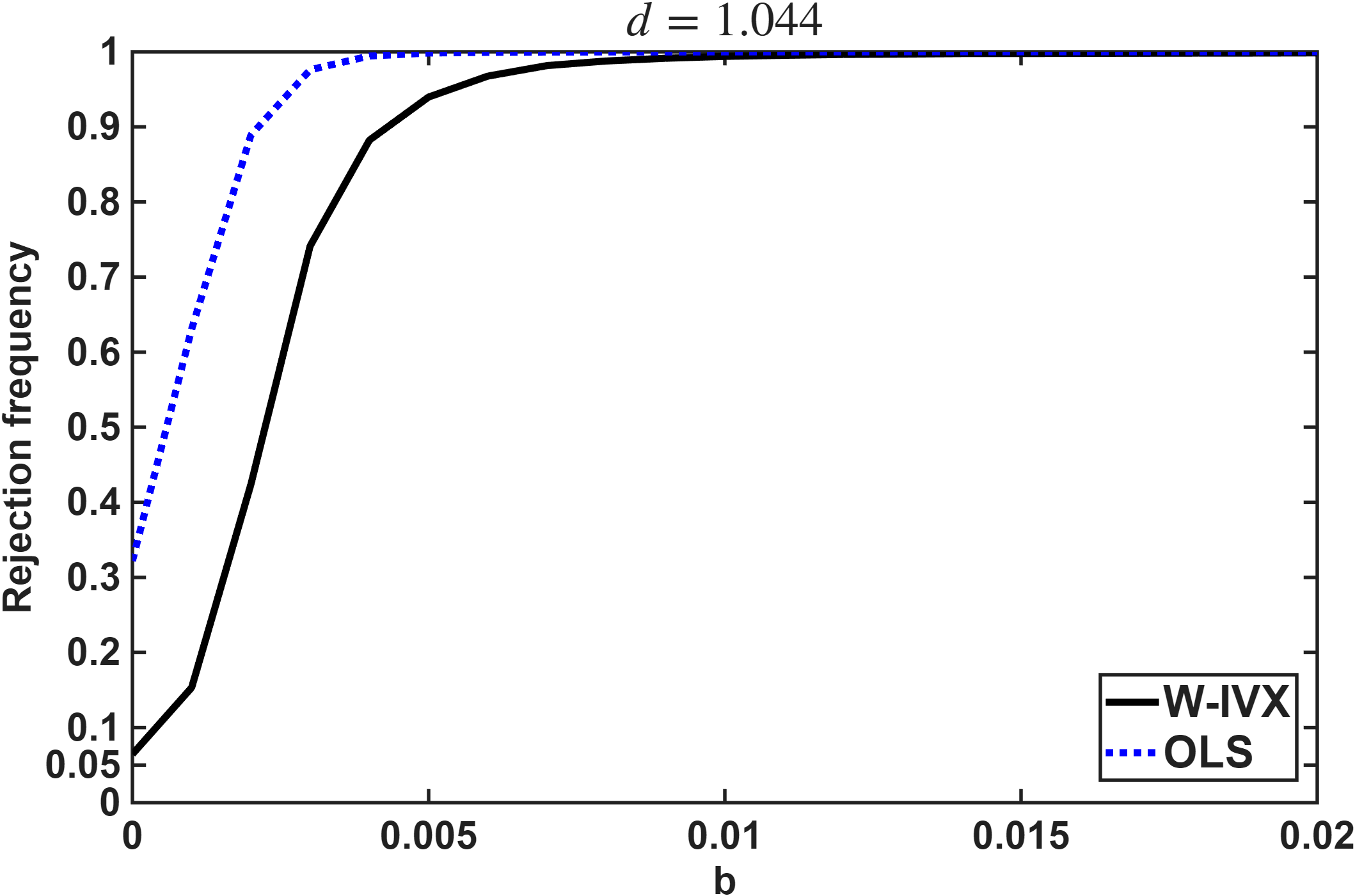}
\includegraphics[width=0.49\linewidth]{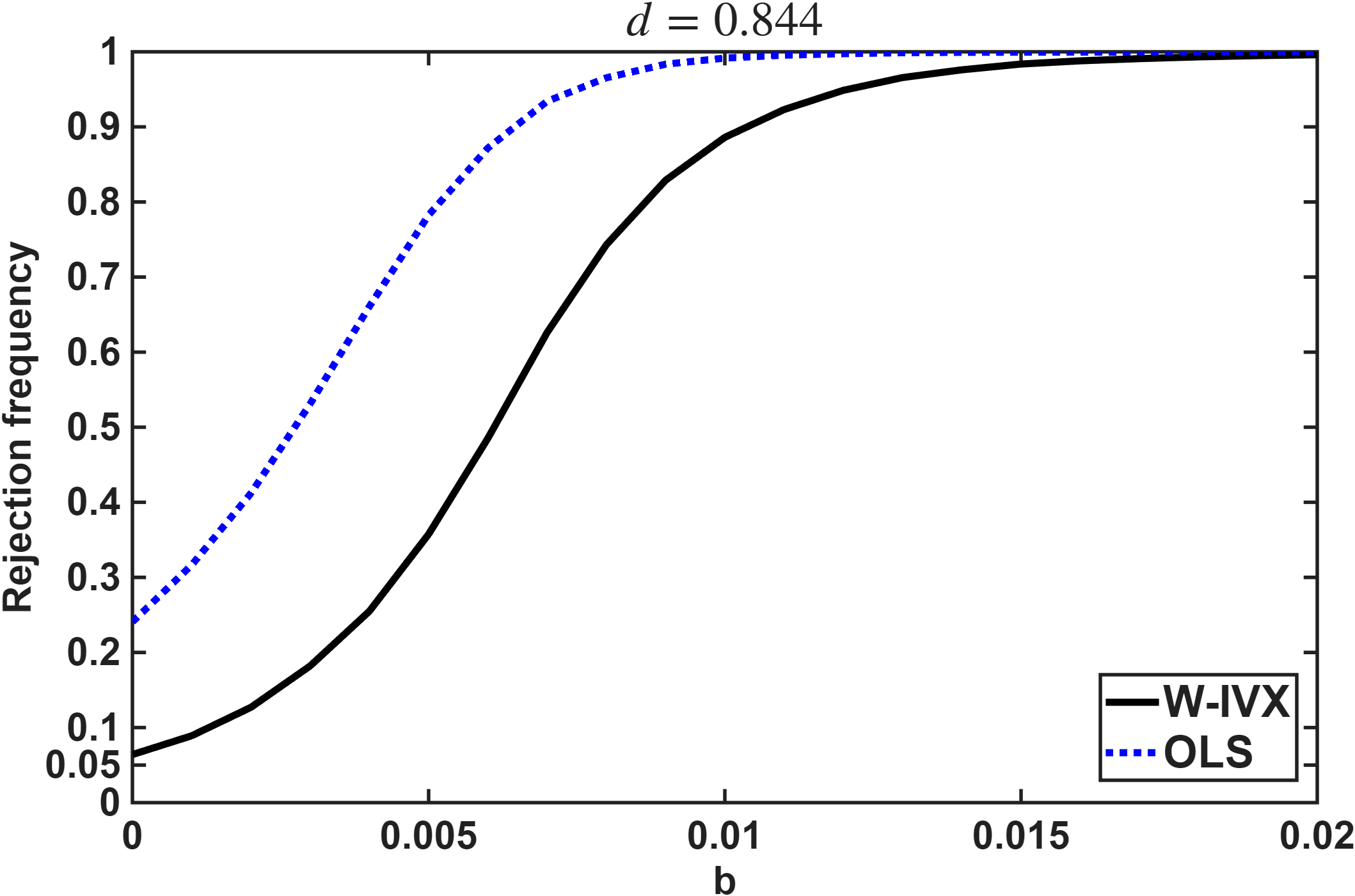}\\[0.6ex]
\includegraphics[width=0.49\linewidth]{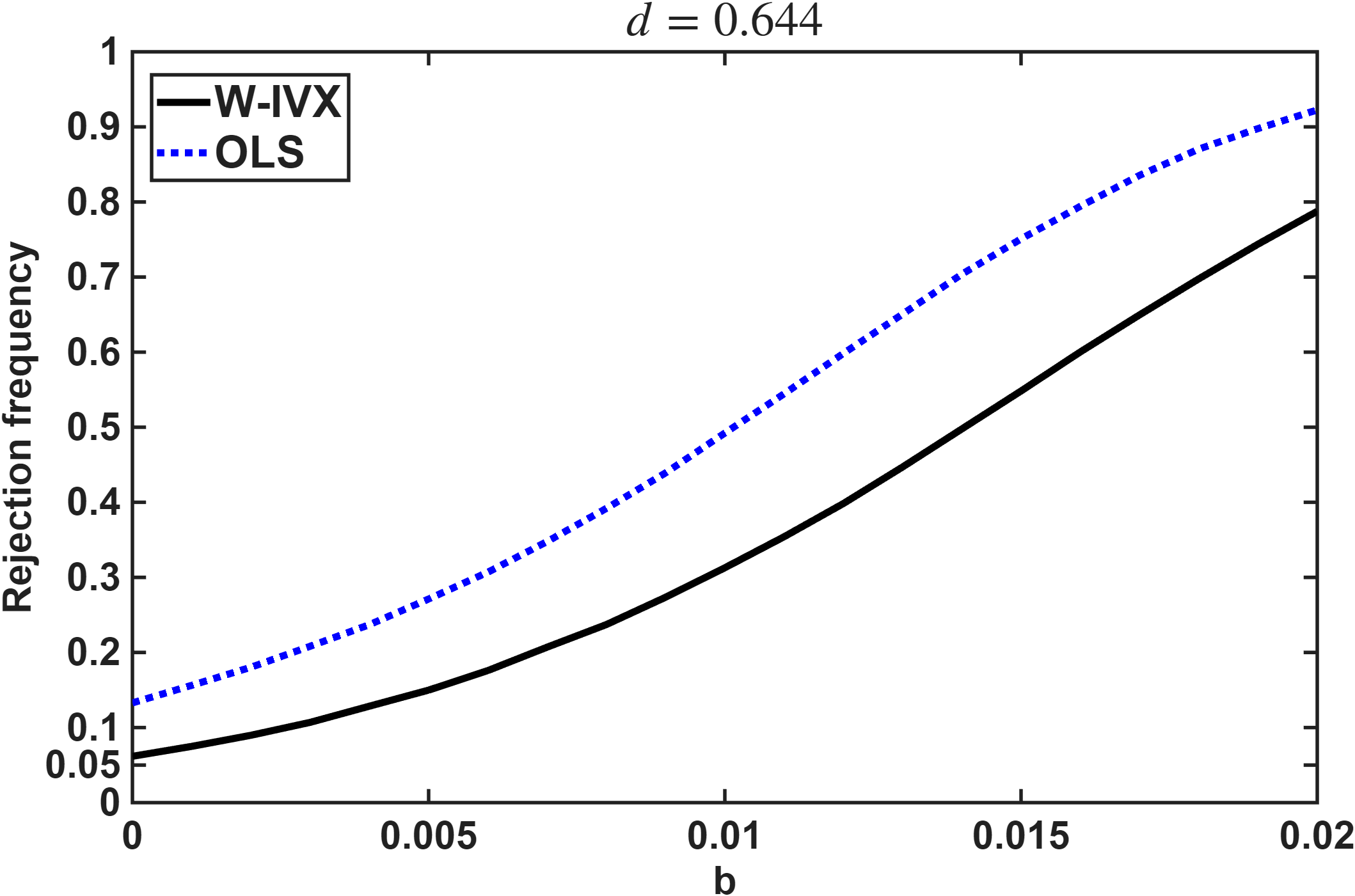}
\includegraphics[width=0.49\linewidth]{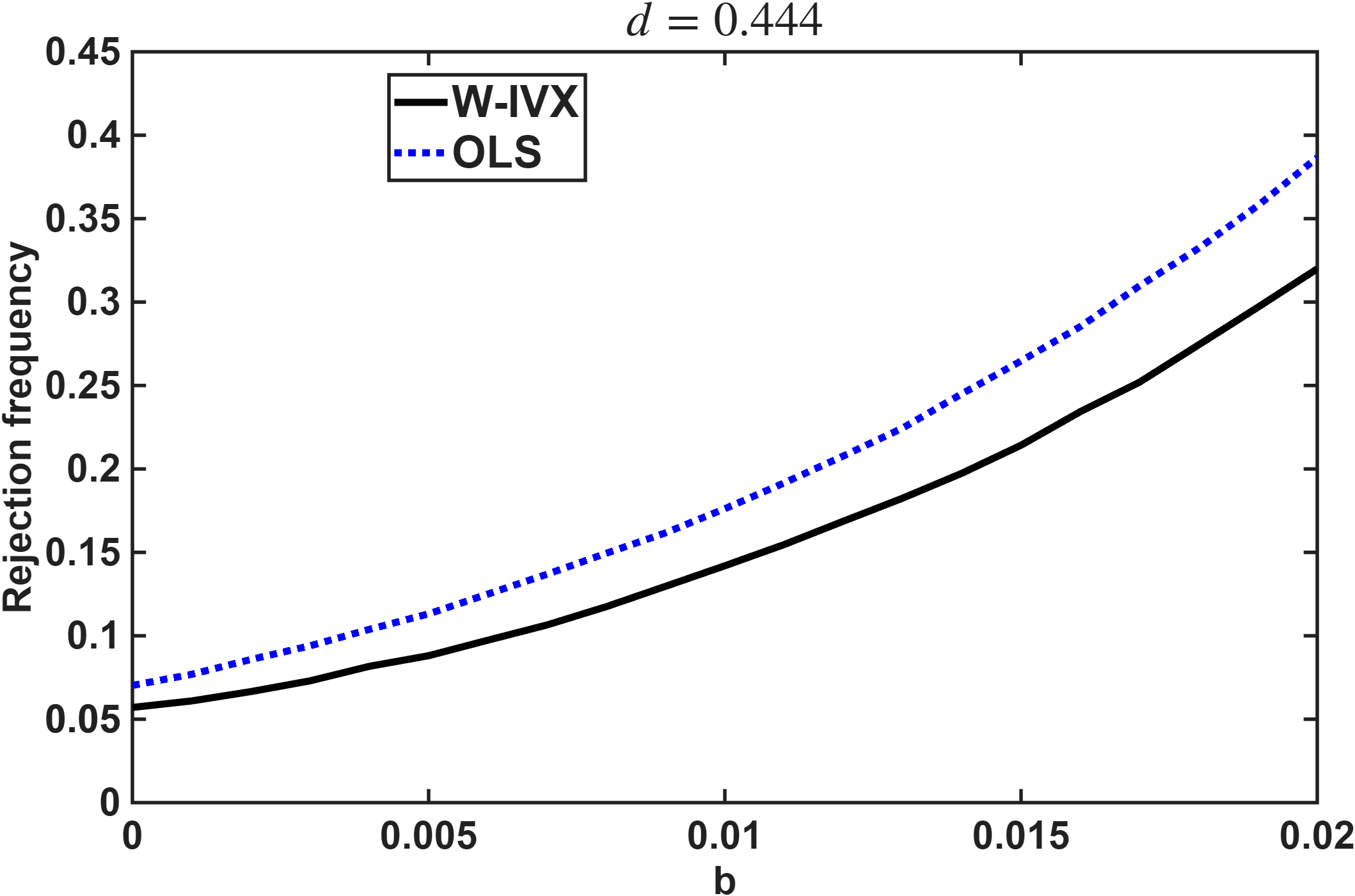}\\[0.6ex]
\includegraphics[width=0.49\linewidth]{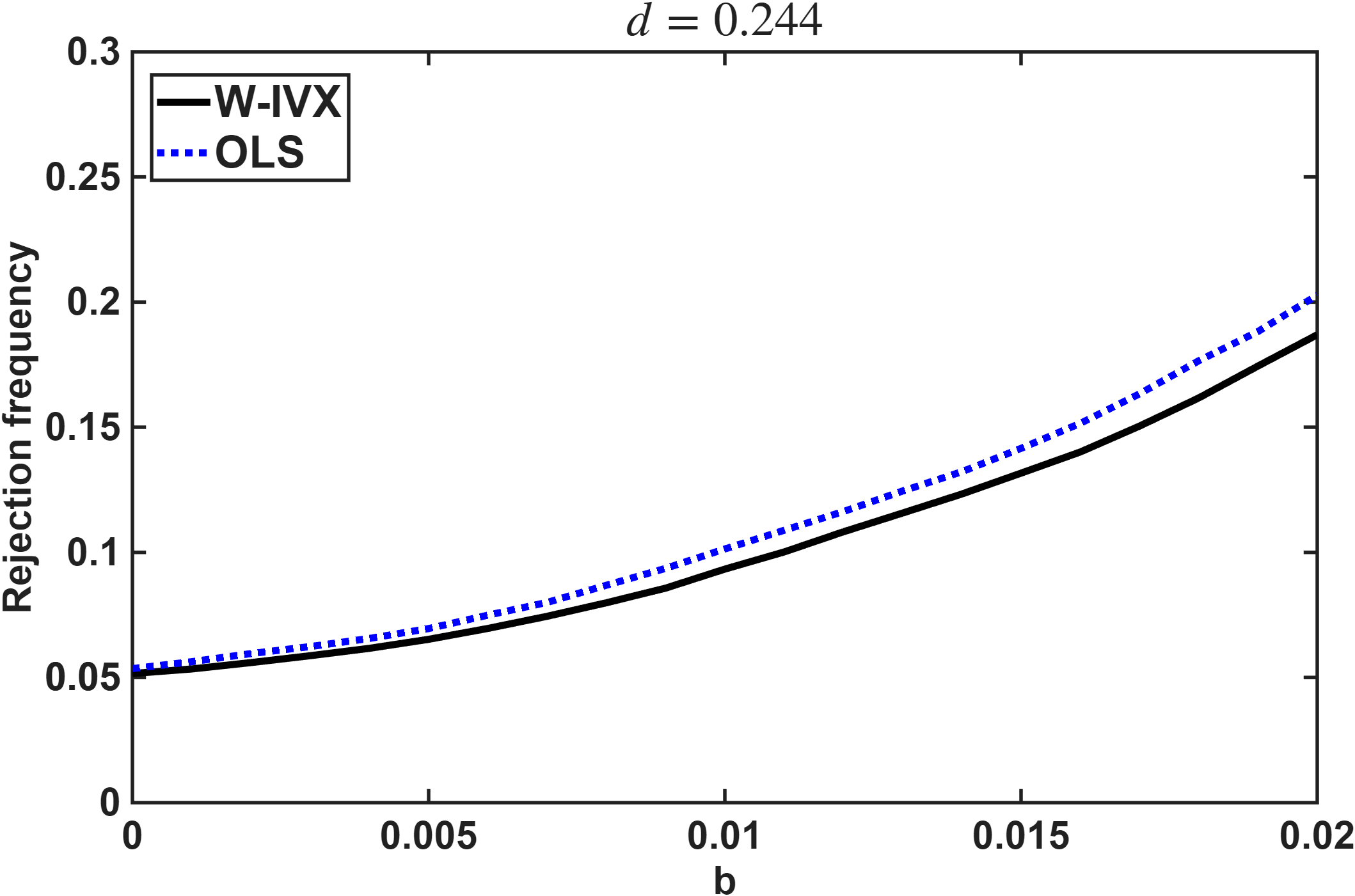}
\includegraphics[width=0.49\linewidth]{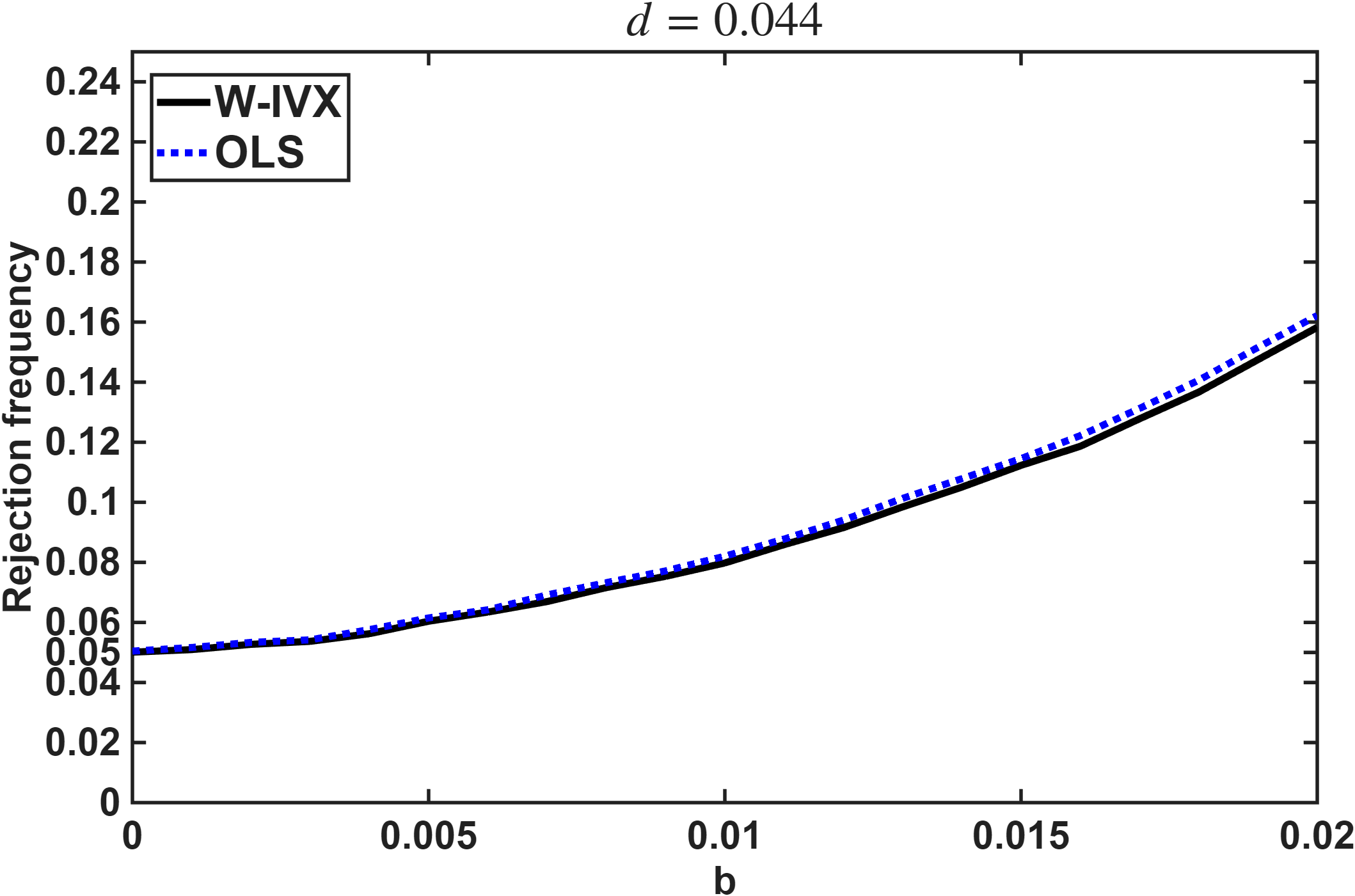}
\end{subfigure}

\end{figure}

\end{landscape}

Figures \ref{fig:sim-power-linear-w00}, \ref{fig:sim-power-linear-w05},
\ref{fig:sim-power-linear-w075}, and \ref{fig:sim-power-linear-w10}
report empirical power for the linear regression specification with
weight exponents $0,0.5,0.75,1$, respectively. Figures \ref{fig:sim-sap-linear-w00},
\ref{fig:sim-sap-linear-w05}, \ref{fig:sim-sap-linear-w075}, and
\ref{fig:sim-sap-linear-w10} report the corresponding size-adjusted
power. Figure \ref{fig:sim-sap-nonlinear-f05-w05} reports size-adjusted
power for the nonlinear specification. Each size-adjusted-power figure
continues on a second page.

\begin{landscape} 

\subsection*{Empirical Size}

\begin{table}[H]
\centering \caption{Empirical size (nominal size 5\%): $[\alpha,\nu]=[2,2.1]$; \textbf{DGP1};
$[f(x),\ w(x)]=[x,\ x^{\left\langle 0.5\right\rangle }]$; $d_{*}=0.5$ }
\begin{adjustbox}{max totalsize={\linewidth}{0.65\textwidth},center}
\begin{tabular}{rrrrrrrrrrrrrrrrr}
 &  & \multicolumn{5}{c}{W-IVX} & \multicolumn{5}{c}{WLS} & \multicolumn{5}{c}{OLS}\tabularnewline
\midrule 
$n$  & $d$  & $\delta$=-0.99  & -0.75  & 0.00  & 0.75  & 0.99  & $\delta$=-0.99  & -0.75  & 0.00  & 0.75  & 0.99  & $\delta$=-0.99  & -0.75  & 0.00  & 0.75  & 0.99\tabularnewline
\midrule
\midrule 
250  & -0.10  & 0.084  & 0.082  & 0.080  & 0.082  & 0.082  & 0.083  & 0.083  & 0.081  & 0.081  & 0.083  & 0.062  & 0.063  & 0.065  & 0.064  & 0.061\tabularnewline
500  & -0.10  & 0.080  & 0.078  & 0.083  & 0.078  & 0.081  & 0.080  & 0.077  & 0.084  & 0.078  & 0.082  & 0.055  & 0.060  & 0.062  & 0.060  & 0.058\tabularnewline
1000  & -0.10  & 0.077  & 0.075  & 0.077  & 0.073  & 0.076  & 0.077  & 0.076  & 0.077  & 0.073  & 0.077  & 0.054  & 0.059  & 0.058  & 0.052  & 0.055\tabularnewline
\midrule 
250  & 0.10  & 0.081  & 0.082  & 0.082  & 0.084  & 0.083  & 0.082  & 0.084  & 0.082  & 0.085  & 0.084  & 0.062  & 0.065  & 0.065  & 0.066  & 0.062\tabularnewline
500  & 0.10  & 0.077  & 0.077  & 0.081  & 0.079  & 0.082  & 0.078  & 0.078  & 0.081  & 0.080  & 0.083  & 0.056  & 0.059  & 0.063  & 0.058  & 0.058\tabularnewline
1000  & 0.10  & 0.079  & 0.078  & 0.076  & 0.076  & 0.077  & 0.079  & 0.078  & 0.076  & 0.077  & 0.078  & 0.054  & 0.057  & 0.058  & 0.055  & 0.054\tabularnewline
\bottomrule
\end{tabular}\end{adjustbox} \label{tab:sim-dgp1-linear-w05-alpha-2p00} 
\end{table}

\end{landscape} \clearpage{}

\begin{landscape}

\begin{table}[htbp]
\centering \caption{Empirical size (nominal size 5\%): $[\alpha,\nu]=[1.99,2.1]$; \textbf{DGP1};
$[f(x),\ w(x)]=[x,\ x^{\left\langle 0.5\right\rangle }]$; $d_{*}=0.497$ }
\begin{adjustbox}{max totalsize={\linewidth}{0.77\textwidth},center}
\begin{tabular}{rrrrrrrrrrrrrrrrr}
 &  & \multicolumn{5}{c}{W-IVX} & \multicolumn{5}{c}{WLS} & \multicolumn{5}{c}{OLS}\tabularnewline
\midrule 
$n$  & $d$  & $\delta$=-0.99  & -0.75  & 0.00  & 0.75  & 0.99  & $\delta$=-0.99  & -0.75  & 0.00  & 0.75  & 0.99  & $\delta$=-0.99  & -0.75  & 0.00  & 0.75  & 0.99\tabularnewline
\midrule
\midrule 
250  & -0.103  & 0.080  & 0.081  & 0.083  & 0.081  & 0.082  & 0.081  & 0.082  & 0.083  & 0.082  & 0.083  & 0.060  & 0.062  & 0.067  & 0.066  & 0.060\tabularnewline
500  & -0.103  & 0.081  & 0.079  & 0.080  & 0.080  & 0.076  & 0.081  & 0.079  & 0.080  & 0.079  & 0.078  & 0.056  & 0.057  & 0.063  & 0.060  & 0.054\tabularnewline
1000  & -0.103  & 0.076  & 0.075  & 0.074  & 0.075  & 0.079  & 0.076  & 0.075  & 0.074  & 0.075  & 0.078  & 0.053  & 0.055  & 0.056  & 0.057  & 0.053\tabularnewline
\midrule 
250  & 0.097  & 0.079  & 0.080  & 0.084  & 0.082  & 0.082  & 0.081  & 0.082  & 0.085  & 0.084  & 0.084  & 0.060  & 0.063  & 0.068  & 0.068  & 0.059\tabularnewline
500  & 0.097  & 0.078  & 0.078  & 0.078  & 0.077  & 0.076  & 0.079  & 0.078  & 0.078  & 0.077  & 0.077  & 0.057  & 0.057  & 0.062  & 0.060  & 0.054\tabularnewline
1000  & 0.097  & 0.075  & 0.076  & 0.074  & 0.077  & 0.078  & 0.076  & 0.077  & 0.075  & 0.078  & 0.079  & 0.052  & 0.054  & 0.056  & 0.057  & 0.052\tabularnewline
\bottomrule
\end{tabular}\end{adjustbox} \label{tab:sim-dgp1-linear-w05-alpha-1p99} 
\end{table}

\end{landscape} \clearpage{}

\begin{landscape}

\begin{table}[htbp]
\centering \caption{Empirical size (nominal size 5\%): $[\alpha,\nu]=[1.8,2.1]$; \textbf{DGP1};
$[f(x),\ w(x)]=[x,\ x^{\left\langle 0.5\right\rangle }]$; $d_{*}=0.444$ }
\begin{adjustbox}{max totalsize={\linewidth}{0.77\textwidth},center}
\begin{tabular}{rrrrrrrrrrrrrrrrr}
 &  & \multicolumn{5}{c}{W-IVX} & \multicolumn{5}{c}{WLS} & \multicolumn{5}{c}{OLS}\tabularnewline
\midrule 
$n$  & $d$  & $\delta$=-0.99  & -0.75  & 0.00  & 0.75  & 0.99  & $\delta$=-0.99  & -0.75  & 0.00  & 0.75  & 0.99  & $\delta$=-0.99  & -0.75  & 0.00  & 0.75  & 0.99\tabularnewline
\midrule
\midrule 
250  & -0.156  & 0.074  & 0.077  & 0.081  & 0.077  & 0.081  & 0.075  & 0.078  & 0.081  & 0.079  & 0.082  & 0.055  & 0.060  & 0.064  & 0.059  & 0.055\tabularnewline
500  & -0.156  & 0.078  & 0.073  & 0.078  & 0.077  & 0.076  & 0.078  & 0.075  & 0.079  & 0.078  & 0.077  & 0.051  & 0.053  & 0.057  & 0.058  & 0.052\tabularnewline
1000  & -0.156  & 0.073  & 0.071  & 0.077  & 0.074  & 0.072  & 0.073  & 0.071  & 0.078  & 0.074  & 0.073  & 0.050  & 0.050  & 0.057  & 0.050  & 0.050\tabularnewline
\midrule 
250  & 0.044  & 0.078  & 0.078  & 0.079  & 0.079  & 0.079  & 0.080  & 0.080  & 0.080  & 0.080  & 0.082  & 0.056  & 0.060  & 0.061  & 0.059  & 0.056\tabularnewline
500  & 0.044  & 0.077  & 0.074  & 0.079  & 0.079  & 0.075  & 0.079  & 0.076  & 0.080  & 0.079  & 0.077  & 0.051  & 0.054  & 0.057  & 0.057  & 0.053\tabularnewline
1000  & 0.044  & 0.074  & 0.073  & 0.078  & 0.072  & 0.073  & 0.076  & 0.073  & 0.079  & 0.072  & 0.074  & 0.051  & 0.050  & 0.058  & 0.048  & 0.049\tabularnewline
\bottomrule
\end{tabular}\end{adjustbox} \label{tab:sim-dgp1-linear-w05-alpha-1p80} 
\end{table}

\end{landscape} \clearpage{}

\begin{landscape}

\begin{table}[htbp]
\centering \caption{Empirical size (nominal size 5\%): $[\alpha,\nu]=[1.6,2.1]$; \textbf{DGP1};
$[f(x),\ w(x)]=[x,\ x^{\left\langle 0.5\right\rangle }]$; $d_{*}=0.375$ }
\begin{adjustbox}{max totalsize={\linewidth}{0.77\textwidth},center}
\begin{tabular}{rrrrrrrrrrrrrrrrr}
 &  & \multicolumn{5}{c}{W-IVX} & \multicolumn{5}{c}{WLS} & \multicolumn{5}{c}{OLS}\tabularnewline
\midrule 
$n$  & $d$  & $\delta$=-0.99  & -0.75  & 0.00  & 0.75  & 0.99  & $\delta$=-0.99  & -0.75  & 0.00  & 0.75  & 0.99  & $\delta$=-0.99  & -0.75  & 0.00  & 0.75  & 0.99\tabularnewline
\midrule
\midrule 
250  & -0.225  & 0.077  & 0.073  & 0.074  & 0.075  & 0.076  & 0.078  & 0.075  & 0.075  & 0.077  & 0.078  & 0.052  & 0.053  & 0.056  & 0.054  & 0.051\tabularnewline
500  & -0.225  & 0.071  & 0.071  & 0.075  & 0.072  & 0.071  & 0.073  & 0.072  & 0.076  & 0.072  & 0.073  & 0.051  & 0.051  & 0.054  & 0.053  & 0.048\tabularnewline
1000  & -0.225  & 0.067  & 0.068  & 0.071  & 0.074  & 0.070  & 0.068  & 0.068  & 0.071  & 0.074  & 0.071  & 0.047  & 0.045  & 0.048  & 0.050  & 0.046\tabularnewline
\midrule 
250  & -0.025  & 0.078  & 0.076  & 0.077  & 0.076  & 0.076  & 0.081  & 0.078  & 0.079  & 0.078  & 0.078  & 0.053  & 0.053  & 0.057  & 0.056  & 0.052\tabularnewline
500  & -0.025  & 0.071  & 0.073  & 0.076  & 0.072  & 0.072  & 0.073  & 0.074  & 0.076  & 0.075  & 0.074  & 0.050  & 0.050  & 0.056  & 0.052  & 0.050\tabularnewline
1000  & -0.025  & 0.070  & 0.068  & 0.074  & 0.072  & 0.071  & 0.072  & 0.069  & 0.075  & 0.074  & 0.073  & 0.047  & 0.046  & 0.048  & 0.049  & 0.048\tabularnewline
\bottomrule
\end{tabular}\end{adjustbox} \label{tab:sim-dgp1-linear-w05-alpha-1p60} 
\end{table}

\end{landscape} \clearpage{}

\begin{landscape}

\begin{table}[htbp]
\centering \caption{Empirical size (nominal size 5\%): $[\alpha,\nu]=[1.4,2.1]$; \textbf{DGP1};
$[f(x),\ w(x)]=[x,\ x^{\left\langle 0.5\right\rangle }]$; $d_{*}=0.286$ }
\begin{adjustbox}{max totalsize={\linewidth}{0.77\textwidth},center}
\begin{tabular}{rrrrrrrrrrrrrrrrr}
 &  & \multicolumn{5}{c}{W-IVX} & \multicolumn{5}{c}{WLS} & \multicolumn{5}{c}{OLS}\tabularnewline
\midrule 
$n$  & $d$  & $\delta$=-0.99  & -0.75  & 0.00  & 0.75  & 0.99  & $\delta$=-0.99  & -0.75  & 0.00  & 0.75  & 0.99  & $\delta$=-0.99  & -0.75  & 0.00  & 0.75  & 0.99\tabularnewline
\midrule
\midrule 
250  & -0.314  & 0.065  & 0.069  & 0.075  & 0.071  & 0.068  & 0.067  & 0.070  & 0.076  & 0.073  & 0.070  & 0.045  & 0.050  & 0.051  & 0.049  & 0.043\tabularnewline
500  & -0.314  & 0.066  & 0.069  & 0.072  & 0.069  & 0.063  & 0.067  & 0.070  & 0.073  & 0.071  & 0.065  & 0.048  & 0.050  & 0.051  & 0.049  & 0.048\tabularnewline
1000  & -0.314  & 0.061  & 0.064  & 0.068  & 0.068  & 0.065  & 0.062  & 0.065  & 0.069  & 0.069  & 0.066  & 0.043  & 0.042  & 0.044  & 0.045  & 0.044\tabularnewline
\midrule 
250  & -0.114  & 0.065  & 0.069  & 0.074  & 0.071  & 0.065  & 0.069  & 0.072  & 0.076  & 0.073  & 0.069  & 0.046  & 0.050  & 0.050  & 0.050  & 0.044\tabularnewline
500  & -0.114  & 0.065  & 0.068  & 0.071  & 0.069  & 0.065  & 0.067  & 0.071  & 0.073  & 0.072  & 0.066  & 0.047  & 0.049  & 0.050  & 0.049  & 0.047\tabularnewline
1000  & -0.114  & 0.062  & 0.063  & 0.067  & 0.068  & 0.065  & 0.064  & 0.065  & 0.069  & 0.069  & 0.067  & 0.042  & 0.042  & 0.043  & 0.044  & 0.044\tabularnewline
\bottomrule
\end{tabular}\end{adjustbox} \label{tab:sim-dgp1-linear-w05-alpha-1p40} 
\end{table}

\end{landscape} \clearpage{}

\begin{landscape}

\begin{table}[htbp]
\centering \caption{Empirical size (nominal size 5\%): $[\alpha,\nu]=[1.2,2.1]$; \textbf{DGP1};
$[f(x),\ w(x)]=[x,\ x^{\left\langle 0.5\right\rangle }]$; $d_{*}=0.167$ }
\begin{adjustbox}{max totalsize={\linewidth}{0.77\textwidth},center}
\begin{tabular}{rrrrrrrrrrrrrrrrr}
 &  & \multicolumn{5}{c}{W-IVX} & \multicolumn{5}{c}{WLS} & \multicolumn{5}{c}{OLS}\tabularnewline
\midrule 
$n$  & $d$  & $\delta$=-0.99  & -0.75  & 0.00  & 0.75  & 0.99  & $\delta$=-0.99  & -0.75  & 0.00  & 0.75  & 0.99  & $\delta$=-0.99  & -0.75  & 0.00  & 0.75  & 0.99\tabularnewline
\midrule
\midrule 
250  & -0.433  & 0.066  & 0.065  & 0.067  & 0.067  & 0.061  & 0.067  & 0.066  & 0.069  & 0.069  & 0.063  & 0.046  & 0.047  & 0.050  & 0.049  & 0.046\tabularnewline
500  & -0.433  & 0.062  & 0.064  & 0.068  & 0.060  & 0.058  & 0.062  & 0.065  & 0.068  & 0.061  & 0.060  & 0.042  & 0.046  & 0.046  & 0.044  & 0.044\tabularnewline
1000  & -0.433  & 0.060  & 0.061  & 0.062  & 0.060  & 0.058  & 0.061  & 0.062  & 0.062  & 0.060  & 0.058  & 0.042  & 0.043  & 0.041  & 0.041  & 0.040\tabularnewline
\midrule 
250  & -0.233  & 0.063  & 0.063  & 0.067  & 0.067  & 0.061  & 0.065  & 0.065  & 0.069  & 0.069  & 0.064  & 0.045  & 0.045  & 0.050  & 0.048  & 0.044\tabularnewline
500  & -0.233  & 0.061  & 0.062  & 0.067  & 0.059  & 0.057  & 0.063  & 0.066  & 0.069  & 0.061  & 0.060  & 0.041  & 0.045  & 0.045  & 0.042  & 0.042\tabularnewline
1000  & -0.233  & 0.060  & 0.059  & 0.060  & 0.060  & 0.058  & 0.062  & 0.061  & 0.060  & 0.061  & 0.060  & 0.041  & 0.042  & 0.040  & 0.040  & 0.039\tabularnewline
\bottomrule
\end{tabular}\end{adjustbox} \label{tab:sim-dgp1-linear-w05-alpha-1p20} 
\end{table}

\end{landscape} \clearpage{}

\begin{landscape}

\begin{table}[htbp]
\centering \caption{Empirical size (nominal size 5\%): $[\alpha,\nu]=[1.0,2.1]$; \textbf{DGP1};
$[f(x),\ w(x)]=[x,\ x^{\left\langle 0.5\right\rangle }]$; $d_{*}=0$ }
\begin{adjustbox}{max totalsize={\linewidth}{0.77\textwidth},center}
\begin{tabular}{rrrrrrrrrrrrrrrrr}
 &  & \multicolumn{5}{c}{W-IVX} & \multicolumn{5}{c}{WLS} & \multicolumn{5}{c}{OLS}\tabularnewline
\midrule 
$n$  & $d$  & $\delta$=-0.99  & -0.75  & 0.00  & 0.75  & 0.99  & $\delta$=-0.99  & -0.75  & 0.00  & 0.75  & 0.99  & $\delta$=-0.99  & -0.75  & 0.00  & 0.75  & 0.99\tabularnewline
\midrule
\midrule 
250  & -0.60  & 0.054  & 0.061  & 0.063  & 0.057  & 0.054  & 0.055  & 0.061  & 0.064  & 0.058  & 0.055  & 0.043  & 0.044  & 0.050  & 0.045  & 0.042\tabularnewline
500  & -0.60  & 0.051  & 0.056  & 0.057  & 0.056  & 0.051  & 0.052  & 0.057  & 0.058  & 0.057  & 0.051  & 0.040  & 0.041  & 0.044  & 0.043  & 0.039\tabularnewline
1000  & -0.60  & 0.055  & 0.051  & 0.056  & 0.056  & 0.049  & 0.055  & 0.051  & 0.056  & 0.056  & 0.049  & 0.043  & 0.037  & 0.043  & 0.042  & 0.040\tabularnewline
\midrule 
250  & -0.40  & 0.052  & 0.061  & 0.063  & 0.058  & 0.053  & 0.055  & 0.062  & 0.064  & 0.060  & 0.056  & 0.042  & 0.043  & 0.049  & 0.045  & 0.041\tabularnewline
500  & -0.40  & 0.049  & 0.056  & 0.056  & 0.058  & 0.050  & 0.052  & 0.058  & 0.058  & 0.059  & 0.051  & 0.039  & 0.040  & 0.044  & 0.043  & 0.038\tabularnewline
1000  & -0.40  & 0.055  & 0.051  & 0.055  & 0.057  & 0.049  & 0.057  & 0.052  & 0.056  & 0.058  & 0.050  & 0.041  & 0.037  & 0.042  & 0.041  & 0.039\tabularnewline
\bottomrule
\end{tabular}\end{adjustbox} \label{tab:sim-dgp1-linear-w05-alpha-1p00} 
\end{table}

\end{landscape} \clearpage{}

\begin{landscape}

\begin{table}[htbp]
\centering \caption{Empirical size (nominal size 5\%): $[\alpha,\nu]=[0.80,2.1]$; \textbf{DGP1};
$[f(x),\ w(x)]=[x,\ x^{\left\langle 0.5\right\rangle }]$; $d_{*}=-0.250$ }
\begin{adjustbox}{max totalsize={\linewidth}{0.77\textwidth},center}
\begin{tabular}{rrrrrrrrrrrrrrrrr}
 &  & \multicolumn{5}{c}{W-IVX} & \multicolumn{5}{c}{WLS} & \multicolumn{5}{c}{OLS}\tabularnewline
\midrule 
$n$  & $d$  & $\delta$=-0.99  & -0.75  & 0.00  & 0.75  & 0.99  & $\delta$=-0.99  & -0.75  & 0.00  & 0.75  & 0.99  & $\delta$=-0.99  & -0.75  & 0.00  & 0.75  & 0.99\tabularnewline
\midrule
\midrule 
250  & -0.85  & 0.051  & 0.055  & 0.056  & 0.054  & 0.049  & 0.051  & 0.055  & 0.056  & 0.054  & 0.049  & 0.044  & 0.047  & 0.044  & 0.045  & 0.043\tabularnewline
500  & -0.85  & 0.048  & 0.053  & 0.053  & 0.047  & 0.048  & 0.049  & 0.053  & 0.053  & 0.047  & 0.048  & 0.039  & 0.043  & 0.045  & 0.041  & 0.042\tabularnewline
1000  & -0.85  & 0.046  & 0.049  & 0.047  & 0.048  & 0.047  & 0.045  & 0.049  & 0.047  & 0.048  & 0.047  & 0.037  & 0.041  & 0.040  & 0.041  & 0.037\tabularnewline
\midrule 
250  & -0.65  & 0.051  & 0.054  & 0.052  & 0.055  & 0.050  & 0.052  & 0.055  & 0.053  & 0.055  & 0.051  & 0.043  & 0.046  & 0.043  & 0.044  & 0.042\tabularnewline
500  & -0.65  & 0.048  & 0.051  & 0.052  & 0.049  & 0.048  & 0.049  & 0.052  & 0.053  & 0.049  & 0.049  & 0.038  & 0.042  & 0.043  & 0.040  & 0.041\tabularnewline
1000  & -0.65  & 0.046  & 0.048  & 0.048  & 0.048  & 0.046  & 0.046  & 0.048  & 0.048  & 0.048  & 0.047  & 0.037  & 0.039  & 0.039  & 0.040  & 0.037\tabularnewline
\bottomrule
\end{tabular}\end{adjustbox} \label{tab:sim-dgp1-linear-w05-alpha-0p80} 
\end{table}

\end{landscape} \clearpage{}

\begin{landscape}

\begin{table}[htbp]
\centering \caption{Empirical size (nominal size 5\%): W-IVX; $[\alpha,\nu]=[2,2.1]$;
$f(x)=x$; $d_{*}=0.5$; \textbf{DGP1}}
\begin{adjustbox}{max totalsize={\linewidth}{0.77\textwidth},center}
\begin{tabular}{rrrrrrrrrrrrrr}
$w(x)$  &  & \multicolumn{3}{c}{$x^{\langle0\rangle}$} & \multicolumn{3}{c}{$x^{\langle0.5\rangle}$} & \multicolumn{3}{c}{$x^{\langle0.75\rangle}$} & \multicolumn{3}{c}{$x^{\langle1\rangle}$}\tabularnewline
\midrule 
$n$  & $d$  & $\delta$=-0.99  & 0.00  & 0.99  & $\delta$=-0.99  & 0.00  & 0.99  & $\delta$=-0.99  & 0.00  & 0.99  & $\delta$=-0.99  & 0.00  & 0.99\tabularnewline
\midrule
\midrule 
250  & -0.10  & 0.074  & 0.075  & 0.072  & 0.084  & 0.080  & 0.082  & 0.062  & 0.067  & 0.059  & 0.058  & 0.070  & 0.061\tabularnewline
500  & -0.10  & 0.070  & 0.069  & 0.067  & 0.080  & 0.083  & 0.081  & 0.053  & 0.064  & 0.058  & 0.057  & 0.064  & 0.056\tabularnewline
1000  & -0.10  & 0.066  & 0.067  & 0.068  & 0.077  & 0.077  & 0.076  & 0.053  & 0.059  & 0.051  & 0.054  & 0.058  & 0.054\tabularnewline
\midrule 
250  & 0.10  & 0.073  & 0.074  & 0.073  & 0.081  & 0.082  & 0.083  & 0.060  & 0.067  & 0.058  & 0.059  & 0.070  & 0.061\tabularnewline
500  & 0.10  & 0.071  & 0.071  & 0.070  & 0.077  & 0.081  & 0.082  & 0.054  & 0.064  & 0.056  & 0.056  & 0.063  & 0.055\tabularnewline
1000  & 0.10  & 0.062  & 0.067  & 0.068  & 0.079  & 0.076  & 0.077  & 0.055  & 0.058  & 0.051  & 0.055  & 0.058  & 0.054\tabularnewline
\bottomrule
\end{tabular}\end{adjustbox} \label{tab:sim-dgp1-linear-weights-alpha-2p00} 
\end{table}

\end{landscape} \clearpage{}

\begin{landscape}

\begin{table}[htbp]
\centering \caption{Empirical size (nominal size 5\%): W-IVX; $[\alpha,\nu]=[1.99,2.1]$;
$f(x)=x$; $d_{*}=0.497$; \textbf{DGP1}}
\begin{adjustbox}{max totalsize={\linewidth}{0.77\textwidth},center}
\begin{tabular}{rrrrrrrrrrrrrr}
$w(x)$  &  & \multicolumn{3}{c}{$x^{\langle0\rangle}$} & \multicolumn{3}{c}{$x^{\langle0.5\rangle}$} & \multicolumn{3}{c}{$x^{\langle0.75\rangle}$} & \multicolumn{3}{c}{$x^{\langle1\rangle}$}\tabularnewline
\midrule 
$n$  & $d$  & $\delta$=-0.99  & 0.00  & 0.99  & $\delta$=-0.99  & 0.00  & 0.99  & $\delta$=-0.99  & 0.00  & 0.99  & $\delta$=-0.99  & 0.00  & 0.99\tabularnewline
\midrule
\midrule 
250  & -0.103  & 0.069  & 0.071  & 0.071  & 0.080  & 0.083  & 0.082  & 0.057  & 0.067  & 0.057  & 0.057  & 0.065  & 0.058\tabularnewline
500  & -0.103  & 0.070  & 0.071  & 0.069  & 0.081  & 0.080  & 0.076  & 0.054  & 0.063  & 0.056  & 0.058  & 0.061  & 0.053\tabularnewline
1000  & -0.103  & 0.066  & 0.062  & 0.065  & 0.076  & 0.074  & 0.079  & 0.056  & 0.057  & 0.051  & 0.053  & 0.056  & 0.052\tabularnewline
\midrule 
250  & 0.097  & 0.071  & 0.071  & 0.072  & 0.079  & 0.084  & 0.082  & 0.057  & 0.069  & 0.055  & 0.056  & 0.064  & 0.058\tabularnewline
500  & 0.097  & 0.069  & 0.068  & 0.069  & 0.078  & 0.078  & 0.076  & 0.054  & 0.062  & 0.055  & 0.057  & 0.060  & 0.053\tabularnewline
1000  & 0.097  & 0.066  & 0.062  & 0.062  & 0.075  & 0.074  & 0.078  & 0.057  & 0.056  & 0.053  & 0.053  & 0.056  & 0.052\tabularnewline
\bottomrule
\end{tabular}\end{adjustbox} \label{tab:sim-dgp1-linear-weights-alpha-1p99} 
\end{table}

\end{landscape} \clearpage{}

\begin{landscape}

\begin{table}[htbp]
\centering \caption{Empirical size (nominal size 5\%): W-IVX; $[\alpha,\nu]=[1.8,2.1]$;
$f(x)=x$; $d_{*}=0.444$; \textbf{DGP1}}
\begin{adjustbox}{max totalsize={\linewidth}{0.77\textwidth},center}
\begin{tabular}{rrrrrrrrrrrrrr}
$w(x)$  &  & \multicolumn{3}{c}{$x^{\langle0\rangle}$} & \multicolumn{3}{c}{$x^{\langle0.5\rangle}$} & \multicolumn{3}{c}{$x^{\langle0.75\rangle}$} & \multicolumn{3}{c}{$x^{\langle1\rangle}$}\tabularnewline
\midrule 
$n$  & $d$  & $\delta$=-0.99  & 0.00  & 0.99  & $\delta$=-0.99  & 0.00  & 0.99  & $\delta$=-0.99  & 0.00  & 0.99  & $\delta$=-0.99  & 0.00  & 0.99\tabularnewline
\midrule
\midrule 
250  & -0.156  & 0.072  & 0.072  & 0.071  & 0.074  & 0.081  & 0.081  & 0.055  & 0.063  & 0.055  & 0.054  & 0.065  & 0.054\tabularnewline
500  & -0.156  & 0.067  & 0.066  & 0.069  & 0.078  & 0.078  & 0.076  & 0.051  & 0.058  & 0.049  & 0.052  & 0.057  & 0.054\tabularnewline
1000  & -0.156  & 0.065  & 0.066  & 0.065  & 0.073  & 0.077  & 0.072  & 0.050  & 0.057  & 0.048  & 0.049  & 0.055  & 0.047\tabularnewline
\midrule 
250  & 0.044  & 0.073  & 0.068  & 0.072  & 0.078  & 0.079  & 0.079  & 0.053  & 0.063  & 0.055  & 0.055  & 0.064  & 0.054\tabularnewline
500  & 0.044  & 0.067  & 0.069  & 0.070  & 0.077  & 0.079  & 0.075  & 0.053  & 0.057  & 0.047  & 0.052  & 0.058  & 0.053\tabularnewline
1000  & 0.044  & 0.066  & 0.066  & 0.064  & 0.074  & 0.078  & 0.073  & 0.050  & 0.057  & 0.048  & 0.049  & 0.053  & 0.048\tabularnewline
\bottomrule
\end{tabular}\end{adjustbox} \label{tab:sim-dgp1-linear-weights-alpha-1p80} 
\end{table}

\end{landscape} \clearpage{}

\begin{landscape}

\begin{table}[htbp]
\centering \caption{Empirical size (nominal size 5\%): W-IVX; $[\alpha,\nu]=[1.6,2.1]$;
$f(x)=x$; $d_{*}=0.375$; \textbf{DGP1}}
\begin{adjustbox}{max totalsize={\linewidth}{0.77\textwidth},center}
\begin{tabular}{rrrrrrrrrrrrrr}
$w(x)$  &  & \multicolumn{3}{c}{$x^{\langle0\rangle}$} & \multicolumn{3}{c}{$x^{\langle0.5\rangle}$} & \multicolumn{3}{c}{$x^{\langle0.75\rangle}$} & \multicolumn{3}{c}{$x^{\langle1\rangle}$}\tabularnewline
\midrule 
$n$  & $d$  & $\delta$=-0.99  & 0.00  & 0.99  & $\delta$=-0.99  & 0.00  & 0.99  & $\delta$=-0.99  & 0.00  & 0.99  & $\delta$=-0.99  & 0.00  & 0.99\tabularnewline
\midrule
\midrule 
250  & -0.225  & 0.071  & 0.069  & 0.072  & 0.077  & 0.074  & 0.076  & 0.048  & 0.060  & 0.050  & 0.049  & 0.057  & 0.052\tabularnewline
500  & -0.225  & 0.067  & 0.069  & 0.065  & 0.071  & 0.075  & 0.071  & 0.049  & 0.055  & 0.049  & 0.047  & 0.052  & 0.046\tabularnewline
1000  & -0.225  & 0.068  & 0.062  & 0.067  & 0.067  & 0.071  & 0.070  & 0.046  & 0.048  & 0.047  & 0.046  & 0.052  & 0.044\tabularnewline
\midrule 
250  & -0.025  & 0.067  & 0.070  & 0.070  & 0.078  & 0.077  & 0.076  & 0.050  & 0.060  & 0.050  & 0.048  & 0.056  & 0.052\tabularnewline
500  & -0.025  & 0.063  & 0.066  & 0.065  & 0.071  & 0.076  & 0.072  & 0.049  & 0.056  & 0.049  & 0.046  & 0.052  & 0.047\tabularnewline
1000  & -0.025  & 0.067  & 0.066  & 0.064  & 0.070  & 0.074  & 0.071  & 0.048  & 0.047  & 0.047  & 0.046  & 0.052  & 0.045\tabularnewline
\bottomrule
\end{tabular}\end{adjustbox} \label{tab:sim-dgp1-linear-weights-alpha-1p60} 
\end{table}

\end{landscape} \clearpage{}

\begin{landscape}

\begin{table}[htbp]
\centering \caption{Empirical size (nominal size 5\%): W-IVX; $[\alpha,\nu]=[1.4,2.1]$;
$f(x)=x$; $d_{*}=0.286$; \textbf{DGP1}}
\begin{adjustbox}{max totalsize={\linewidth}{0.77\textwidth},center}
\begin{tabular}{rrrrrrrrrrrrrr}
$w(x)$  &  & \multicolumn{3}{c}{$x^{\langle0\rangle}$} & \multicolumn{3}{c}{$x^{\langle0.5\rangle}$} & \multicolumn{3}{c}{$x^{\langle0.75\rangle}$} & \multicolumn{3}{c}{$x^{\langle1\rangle}$}\tabularnewline
\midrule 
$n$  & $d$  & $\delta$=-0.99  & 0.00  & 0.99  & $\delta$=-0.99  & 0.00  & 0.99  & $\delta$=-0.99  & 0.00  & 0.99  & $\delta$=-0.99  & 0.00  & 0.99\tabularnewline
\midrule
\midrule 
250  & -0.314  & 0.068  & 0.067  & 0.069  & 0.065  & 0.075  & 0.068  & 0.050  & 0.056  & 0.046  & 0.050  & 0.054  & 0.049\tabularnewline
500  & -0.314  & 0.066  & 0.065  & 0.066  & 0.066  & 0.072  & 0.063  & 0.044  & 0.051  & 0.046  & 0.048  & 0.049  & 0.044\tabularnewline
1000  & -0.314  & 0.067  & 0.067  & 0.064  & 0.061  & 0.068  & 0.065  & 0.045  & 0.045  & 0.044  & 0.039  & 0.048  & 0.042\tabularnewline
\midrule 
250  & -0.114  & 0.067  & 0.069  & 0.068  & 0.065  & 0.074  & 0.065  & 0.048  & 0.055  & 0.046  & 0.051  & 0.052  & 0.049\tabularnewline
500  & -0.114  & 0.063  & 0.065  & 0.068  & 0.065  & 0.071  & 0.065  & 0.044  & 0.049  & 0.047  & 0.048  & 0.048  & 0.044\tabularnewline
1000  & -0.114  & 0.065  & 0.065  & 0.064  & 0.062  & 0.067  & 0.065  & 0.044  & 0.044  & 0.044  & 0.040  & 0.046  & 0.042\tabularnewline
\bottomrule
\end{tabular}\end{adjustbox} \label{tab:sim-dgp1-linear-weights-alpha-1p40} 
\end{table}

\end{landscape} \clearpage{}

\begin{landscape}

\begin{table}[htbp]
\centering \caption{Empirical size (nominal size 5\%): W-IVX; $[\alpha,\nu]=[1.2,2.1]$;
$f(x)=x$; $d_{*}=0.167$; \textbf{DGP1}}
\begin{adjustbox}{max totalsize={\linewidth}{0.77\textwidth},center}
\begin{tabular}{rrrrrrrrrrrrrr}
$w(x)$  &  & \multicolumn{3}{c}{$x^{\langle0\rangle}$} & \multicolumn{3}{c}{$x^{\langle0.5\rangle}$} & \multicolumn{3}{c}{$x^{\langle0.75\rangle}$} & \multicolumn{3}{c}{$x^{\langle1\rangle}$}\tabularnewline
\midrule 
$n$  & $d$  & $\delta$=-0.99  & 0.00  & 0.99  & $\delta$=-0.99  & 0.00  & 0.99  & $\delta$=-0.99  & 0.00  & 0.99  & $\delta$=-0.99  & 0.00  & 0.99\tabularnewline
\midrule
\midrule 
250  & -0.433  & 0.067  & 0.067  & 0.066  & 0.066  & 0.067  & 0.061  & 0.045  & 0.052  & 0.044  & 0.046  & 0.051  & 0.044\tabularnewline
500  & -0.433  & 0.062  & 0.068  & 0.065  & 0.062  & 0.068  & 0.058  & 0.040  & 0.047  & 0.046  & 0.044  & 0.046  & 0.046\tabularnewline
1000  & -0.433  & 0.064  & 0.064  & 0.064  & 0.060  & 0.062  & 0.058  & 0.040  & 0.043  & 0.040  & 0.039  & 0.041  & 0.045\tabularnewline
\midrule 
250  & -0.233  & 0.068  & 0.070  & 0.066  & 0.063  & 0.067  & 0.061  & 0.045  & 0.051  & 0.042  & 0.046  & 0.049  & 0.044\tabularnewline
500  & -0.233  & 0.063  & 0.069  & 0.064  & 0.061  & 0.067  & 0.057  & 0.041  & 0.047  & 0.044  & 0.045  & 0.044  & 0.044\tabularnewline
1000  & -0.233  & 0.062  & 0.065  & 0.063  & 0.060  & 0.060  & 0.058  & 0.039  & 0.042  & 0.040  & 0.037  & 0.040  & 0.043\tabularnewline
\bottomrule
\end{tabular}\end{adjustbox} \label{tab:sim-dgp1-linear-weights-alpha-1p20} 
\end{table}

\end{landscape} \clearpage{}

\begin{landscape}

\begin{table}[htbp]
\centering \caption{Empirical size (nominal size 5\%): W-IVX; $[\alpha,\nu]=[1,2.1]$;
$f(x)=x$; $d_{*}=0$; \textbf{DGP1}}
\begin{adjustbox}{max totalsize={\linewidth}{0.77\textwidth},center}
\begin{tabular}{rrrrrrrrrrrrrr}
$w(x)$  &  & \multicolumn{3}{c}{$x^{\langle0\rangle}$} & \multicolumn{3}{c}{$x^{\langle0.5\rangle}$} & \multicolumn{3}{c}{$x^{\langle0.75\rangle}$} & \multicolumn{3}{c}{$x^{\langle1\rangle}$}\tabularnewline
\midrule 
$n$  & $d$  & $\delta$=-0.99  & 0.00  & 0.99  & $\delta$=-0.99  & 0.00  & 0.99  & $\delta$=-0.99  & 0.00  & 0.99  & $\delta$=-0.99  & 0.00  & 0.99\tabularnewline
\midrule
\midrule 
250  & -0.60  & 0.063  & 0.066  & 0.066  & 0.054  & 0.063  & 0.054  & 0.045  & 0.048  & 0.037  & 0.043  & 0.045  & 0.042\tabularnewline
500  & -0.60  & 0.063  & 0.064  & 0.066  & 0.051  & 0.057  & 0.051  & 0.040  & 0.043  & 0.042  & 0.043  & 0.044  & 0.044\tabularnewline
1000  & -0.60  & 0.063  & 0.062  & 0.061  & 0.055  & 0.056  & 0.049  & 0.038  & 0.044  & 0.038  & 0.042  & 0.042  & 0.039\tabularnewline
\midrule 
250  & -0.40  & 0.063  & 0.066  & 0.064  & 0.052  & 0.063  & 0.053  & 0.043  & 0.046  & 0.038  & 0.043  & 0.043  & 0.041\tabularnewline
500  & -0.40  & 0.062  & 0.065  & 0.065  & 0.049  & 0.056  & 0.050  & 0.039  & 0.042  & 0.042  & 0.042  & 0.043  & 0.042\tabularnewline
1000  & -0.40  & 0.062  & 0.061  & 0.059  & 0.055  & 0.055  & 0.049  & 0.038  & 0.042  & 0.038  & 0.041  & 0.041  & 0.039\tabularnewline
\bottomrule
\end{tabular}\end{adjustbox} \label{tab:sim-dgp1-linear-weights-alpha-1p00} 
\end{table}

\end{landscape} \clearpage{}

\begin{landscape}

\begin{table}[htbp]
\centering \caption{Empirical size (nominal size 5\%): W-IVX; $[\alpha,\nu]=[0.8,2.1]$;
$f(x)=x$; $d_{*}=-0.25$; \textbf{DGP1}}
\begin{adjustbox}{max totalsize={\linewidth}{0.77\textwidth},center}
\begin{tabular}{rrrrrrrrrrrrrr}
$w(x)$  &  & \multicolumn{3}{c}{$x^{\langle0\rangle}$} & \multicolumn{3}{c}{$x^{\langle0.5\rangle}$} & \multicolumn{3}{c}{$x^{\langle0.75\rangle}$} & \multicolumn{3}{c}{$x^{\langle1\rangle}$}\tabularnewline
\midrule 
$n$  & $d$  & $\delta$=-0.99  & 0.00  & 0.99  & $\delta$=-0.99  & 0.00  & 0.99  & $\delta$=-0.99  & 0.00  & 0.99  & $\delta$=-0.99  & 0.00  & 0.99\tabularnewline
\midrule
\midrule 
250  & -0.85  & 0.065  & 0.066  & 0.064  & 0.051  & 0.056  & 0.049  & 0.042  & 0.046  & 0.043  & 0.042  & 0.044  & 0.040\tabularnewline
500  & -0.85  & 0.062  & 0.065  & 0.064  & 0.048  & 0.053  & 0.048  & 0.041  & 0.040  & 0.043  & 0.041  & 0.041  & 0.041\tabularnewline
1000  & -0.85  & 0.062  & 0.062  & 0.064  & 0.046  & 0.047  & 0.047  & 0.039  & 0.040  & 0.038  & 0.040  & 0.039  & 0.039\tabularnewline
\midrule 
250  & -0.65  & 0.059  & 0.061  & 0.062  & 0.051  & 0.052  & 0.050  & 0.042  & 0.045  & 0.043  & 0.041  & 0.043  & 0.040\tabularnewline
500  & -0.65  & 0.060  & 0.061  & 0.059  & 0.048  & 0.052  & 0.048  & 0.039  & 0.040  & 0.043  & 0.041  & 0.040  & 0.041\tabularnewline
1000  & -0.65  & 0.057  & 0.055  & 0.058  & 0.046  & 0.048  & 0.046  & 0.038  & 0.040  & 0.038  & 0.039  & 0.039  & 0.039\tabularnewline
\bottomrule
\end{tabular}\end{adjustbox} \label{tab:sim-dgp1-linear-weights-alpha-0p80} 
\end{table}

\end{landscape} \clearpage{}

\begin{landscape}

\begin{table}[htbp]
\centering \caption{Empirical size (nominal size 5\%): \textbf{DGP1}; $\alpha=2.0$; $[f(x),\ w(x)]=[x^{\left\langle 0.5\right\rangle },\ x^{\left\langle 0.5\right\rangle }]$;
$d_{*}=0.5$ }
\begin{adjustbox}{max totalsize={\linewidth}{0.77\textwidth},center}
\begin{tabular}{rrrrrrrrrrrr}
 &  & \multicolumn{5}{c}{W-IVX} & \multicolumn{5}{c}{OLS}\tabularnewline
\midrule 
$n$  & $d$  & $\delta$=-0.99  & -0.75  & 0.00  & 0.75  & 0.99  & $\delta$=-0.99  & -0.75  & 0.00  & 0.75  & 0.99\tabularnewline
\midrule
\midrule 
250  & -0.30  & 0.081  & 0.079  & 0.082  & 0.082  & 0.082  & 0.081  & 0.080  & 0.083  & 0.083  & 0.083\tabularnewline
500  & -0.30  & 0.078  & 0.077  & 0.078  & 0.076  & 0.075  & 0.079  & 0.078  & 0.078  & 0.077  & 0.075\tabularnewline
1000  & -0.30  & 0.075  & 0.075  & 0.076  & 0.076  & 0.076  & 0.075  & 0.075  & 0.076  & 0.076  & 0.076\tabularnewline
\midrule 
250  & -0.10  & 0.080  & 0.080  & 0.084  & 0.085  & 0.084  & 0.082  & 0.081  & 0.084  & 0.087  & 0.085\tabularnewline
500  & -0.10  & 0.081  & 0.077  & 0.080  & 0.078  & 0.075  & 0.081  & 0.077  & 0.079  & 0.077  & 0.075\tabularnewline
1000  & -0.10  & 0.076  & 0.074  & 0.075  & 0.078  & 0.077  & 0.076  & 0.074  & 0.075  & 0.079  & 0.078\tabularnewline
\bottomrule
\end{tabular}\end{adjustbox} \label{tab:sim-dgp1-nonlinear-alpha-2p00} 
\end{table}

\end{landscape} \clearpage{}

\begin{landscape}

\begin{table}[htbp]
\centering \caption{Empirical size (nominal size 5\%): \textbf{DGP1}; $\alpha=1.99$;
$[f(x),\ w(x)]=[x^{\left\langle 0.5\right\rangle },\ x^{\left\langle 0.5\right\rangle }]$;
$d_{*}=0.497$ }
\begin{adjustbox}{max totalsize={\linewidth}{0.77\textwidth},center}
\begin{tabular}{rrrrrrrrrrrr}
 &  & \multicolumn{5}{c}{W-IVX} & \multicolumn{5}{c}{OLS}\tabularnewline
\midrule 
$n$  & $d$  & $\delta$=-0.99  & -0.75  & 0.00  & 0.75  & 0.99  & $\delta$=-0.99  & -0.75  & 0.00  & 0.75  & 0.99\tabularnewline
\midrule
\midrule 
250  & -0.303  & 0.081  & 0.083  & 0.084  & 0.081  & 0.079  & 0.082  & 0.084  & 0.084  & 0.081  & 0.079\tabularnewline
500  & -0.303  & 0.080  & 0.078  & 0.080  & 0.077  & 0.076  & 0.080  & 0.079  & 0.080  & 0.078  & 0.076\tabularnewline
1000  & -0.303  & 0.075  & 0.075  & 0.076  & 0.074  & 0.072  & 0.075  & 0.075  & 0.076  & 0.075  & 0.072\tabularnewline
\midrule 
250  & -0.103  & 0.082  & 0.083  & 0.084  & 0.082  & 0.080  & 0.082  & 0.083  & 0.084  & 0.082  & 0.081\tabularnewline
500  & -0.103  & 0.079  & 0.081  & 0.080  & 0.079  & 0.078  & 0.080  & 0.080  & 0.080  & 0.078  & 0.079\tabularnewline
1000  & -0.103  & 0.075  & 0.074  & 0.076  & 0.074  & 0.073  & 0.076  & 0.074  & 0.076  & 0.075  & 0.074\tabularnewline
\bottomrule
\end{tabular}\end{adjustbox} \label{tab:sim-dgp1-nonlinear-alpha-1p99} 
\end{table}

\end{landscape} \clearpage{}

\begin{landscape}

\begin{table}[htbp]
\centering \caption{Empirical size (nominal size 5\%): \textbf{DGP1}; $\alpha=1.80$;
$[f(x),\ w(x)]=[x^{\left\langle 0.5\right\rangle },\ x^{\left\langle 0.5\right\rangle }]$;
$d_{*}=0.444$ }
\begin{adjustbox}{max totalsize={\linewidth}{0.77\textwidth},center}
\begin{tabular}{rrrrrrrrrrrr}
 &  & \multicolumn{5}{c}{W-IVX} & \multicolumn{5}{c}{OLS}\tabularnewline
\midrule 
$n$  & $d$  & $\delta$=-0.99  & -0.75  & 0.00  & 0.75  & 0.99  & $\delta$=-0.99  & -0.75  & 0.00  & 0.75  & 0.99\tabularnewline
\midrule
\midrule 
250  & -0.356  & 0.076  & 0.078  & 0.083  & 0.074  & 0.074  & 0.076  & 0.079  & 0.083  & 0.075  & 0.076\tabularnewline
500  & -0.356  & 0.075  & 0.078  & 0.074  & 0.078  & 0.074  & 0.075  & 0.079  & 0.074  & 0.078  & 0.074\tabularnewline
1000  & -0.356  & 0.077  & 0.074  & 0.076  & 0.073  & 0.072  & 0.077  & 0.073  & 0.076  & 0.073  & 0.072\tabularnewline
\midrule 
250  & -0.156  & 0.077  & 0.078  & 0.084  & 0.076  & 0.076  & 0.078  & 0.079  & 0.084  & 0.077  & 0.078\tabularnewline
500  & -0.156  & 0.075  & 0.077  & 0.075  & 0.078  & 0.075  & 0.077  & 0.078  & 0.075  & 0.078  & 0.076\tabularnewline
1000  & -0.156  & 0.079  & 0.073  & 0.076  & 0.073  & 0.072  & 0.080  & 0.073  & 0.076  & 0.073  & 0.072\tabularnewline
\bottomrule
\end{tabular}\end{adjustbox} \label{tab:sim-dgp1-nonlinear-alpha-1p80} 
\end{table}

\end{landscape} \clearpage{}

\begin{landscape}

\begin{table}[htbp]
\centering \caption{Empirical size (nominal size 5\%): \textbf{DGP1}; $\alpha=1.60$;
$[f(x),\ w(x)]=[x^{\left\langle 0.5\right\rangle },\ x^{\left\langle 0.5\right\rangle }]$;
$d_{*}=0.375$ }
\begin{adjustbox}{max totalsize={\linewidth}{0.77\textwidth},center}
\begin{tabular}{rrrrrrrrrrrr}
 &  & \multicolumn{5}{c}{W-IVX} & \multicolumn{5}{c}{OLS}\tabularnewline
\midrule 
$n$  & $d$  & $\delta$=-0.99  & -0.75  & 0.00  & 0.75  & 0.99  & $\delta$=-0.99  & -0.75  & 0.00  & 0.75  & 0.99\tabularnewline
\midrule
\midrule 
250  & -0.425  & 0.070  & 0.075  & 0.080  & 0.070  & 0.074  & 0.071  & 0.075  & 0.080  & 0.071  & 0.075\tabularnewline
500  & -0.425  & 0.072  & 0.070  & 0.074  & 0.073  & 0.071  & 0.072  & 0.071  & 0.075  & 0.073  & 0.071\tabularnewline
1000  & -0.425  & 0.067  & 0.073  & 0.071  & 0.069  & 0.073  & 0.067  & 0.073  & 0.071  & 0.068  & 0.073\tabularnewline
\midrule 
250  & -0.225  & 0.070  & 0.074  & 0.078  & 0.072  & 0.073  & 0.071  & 0.076  & 0.078  & 0.074  & 0.076\tabularnewline
500  & -0.225  & 0.071  & 0.071  & 0.073  & 0.073  & 0.070  & 0.072  & 0.072  & 0.073  & 0.073  & 0.071\tabularnewline
1000  & -0.225  & 0.070  & 0.072  & 0.073  & 0.068  & 0.072  & 0.070  & 0.073  & 0.074  & 0.068  & 0.073\tabularnewline
\bottomrule
\end{tabular}\end{adjustbox} \label{tab:sim-dgp1-nonlinear-alpha-1p60} 
\end{table}

\end{landscape} \clearpage{}

\begin{landscape}

\begin{table}[htbp]
\centering \caption{Empirical size (nominal size 5\%): \textbf{DGP1}; $\alpha=1.40$;
$[f(x),\ w(x)]=[x^{\left\langle 0.5\right\rangle },\ x^{\left\langle 0.5\right\rangle }]$;
$d_{*}=0.286$ }
\begin{adjustbox}{max totalsize={\linewidth}{0.77\textwidth},center}
\begin{tabular}{rrrrrrrrrrrr}
 &  & \multicolumn{5}{c}{W-IVX} & \multicolumn{5}{c}{OLS}\tabularnewline
\midrule 
$n$  & $d$  & $\delta$=-0.99  & -0.75  & 0.00  & 0.75  & 0.99  & $\delta$=-0.99  & -0.75  & 0.00  & 0.75  & 0.99\tabularnewline
\midrule
\midrule 
250  & -0.514  & 0.067  & 0.070  & 0.072  & 0.069  & 0.069  & 0.068  & 0.071  & 0.073  & 0.070  & 0.069\tabularnewline
500  & -0.514  & 0.067  & 0.066  & 0.072  & 0.067  & 0.065  & 0.067  & 0.066  & 0.072  & 0.067  & 0.066\tabularnewline
1000  & -0.514  & 0.068  & 0.067  & 0.068  & 0.070  & 0.064  & 0.068  & 0.067  & 0.068  & 0.070  & 0.064\tabularnewline
\midrule 
250  & -0.314  & 0.067  & 0.070  & 0.072  & 0.067  & 0.067  & 0.069  & 0.071  & 0.072  & 0.069  & 0.069\tabularnewline
500  & -0.314  & 0.066  & 0.066  & 0.070  & 0.066  & 0.065  & 0.068  & 0.067  & 0.072  & 0.067  & 0.066\tabularnewline
1000  & -0.314  & 0.066  & 0.066  & 0.067  & 0.069  & 0.063  & 0.067  & 0.067  & 0.067  & 0.070  & 0.064\tabularnewline
\bottomrule
\end{tabular}\end{adjustbox} \label{tab:sim-dgp1-nonlinear-alpha-1p40} 
\end{table}

\end{landscape} \clearpage{}

\begin{landscape}

\begin{table}[htbp]
\centering \caption{Empirical size (nominal size 5\%): \textbf{DGP1}; $\alpha=1.20$;
$[f(x),\ w(x)]=[x^{\left\langle 0.5\right\rangle },\ x^{\left\langle 0.5\right\rangle }]$;
$d_{*}=0.167$ }
\begin{adjustbox}{max totalsize={\linewidth}{0.77\textwidth},center}
\begin{tabular}{rrrrrrrrrrrr}
 &  & \multicolumn{5}{c}{W-IVX} & \multicolumn{5}{c}{OLS}\tabularnewline
\midrule 
$n$  & $d$  & $\delta$=-0.99  & -0.75  & 0.00  & 0.75  & 0.99  & $\delta$=-0.99  & -0.75  & 0.00  & 0.75  & 0.99\tabularnewline
\midrule
\midrule 
250  & -0.633  & 0.066  & 0.067  & 0.069  & 0.064  & 0.065  & 0.066  & 0.067  & 0.068  & 0.064  & 0.066\tabularnewline
500  & -0.633  & 0.060  & 0.065  & 0.067  & 0.063  & 0.058  & 0.060  & 0.065  & 0.067  & 0.063  & 0.058\tabularnewline
1000  & -0.633  & 0.059  & 0.061  & 0.065  & 0.065  & 0.062  & 0.059  & 0.061  & 0.065  & 0.065  & 0.062\tabularnewline
\midrule 
250  & -0.433  & 0.064  & 0.067  & 0.069  & 0.064  & 0.064  & 0.066  & 0.068  & 0.070  & 0.065  & 0.066\tabularnewline
500  & -0.433  & 0.061  & 0.064  & 0.067  & 0.063  & 0.057  & 0.062  & 0.065  & 0.068  & 0.064  & 0.059\tabularnewline
1000  & -0.433  & 0.059  & 0.061  & 0.064  & 0.064  & 0.060  & 0.060  & 0.062  & 0.064  & 0.064  & 0.061\tabularnewline
\bottomrule
\end{tabular}\end{adjustbox} \label{tab:sim-dgp1-nonlinear-alpha-1p20} 
\end{table}

\end{landscape} \clearpage{}

\begin{landscape}

\begin{table}[htbp]
\centering \caption{Empirical size (nominal size 5\%): \textbf{DGP1}; $\alpha=1.0$; $[f(x),\ w(x)]=[x^{\left\langle 0.5\right\rangle },\ x^{\left\langle 0.5\right\rangle }]$;
$d_{*}=0$ }
\begin{adjustbox}{max totalsize={\linewidth}{0.77\textwidth},center}
\begin{tabular}{rrrrrrrrrrrr}
 &  & \multicolumn{5}{c}{W-IVX} & \multicolumn{5}{c}{OLS}\tabularnewline
\midrule 
$n$  & $d$  & $\delta$=-0.99  & -0.75  & 0.00  & 0.75  & 0.99  & $\delta$=-0.99  & -0.75  & 0.00  & 0.75  & 0.99\tabularnewline
\midrule
\midrule 
250  & -0.80  & 0.056  & 0.058  & 0.059  & 0.058  & 0.052  & 0.056  & 0.058  & 0.059  & 0.058  & 0.053\tabularnewline
500  & -0.80  & 0.053  & 0.057  & 0.059  & 0.055  & 0.053  & 0.053  & 0.058  & 0.059  & 0.055  & 0.053\tabularnewline
1000  & -0.80  & 0.051  & 0.054  & 0.055  & 0.054  & 0.052  & 0.051  & 0.054  & 0.055  & 0.054  & 0.052\tabularnewline
\midrule 
250  & -0.60  & 0.055  & 0.057  & 0.059  & 0.058  & 0.053  & 0.056  & 0.058  & 0.059  & 0.059  & 0.054\tabularnewline
500  & -0.60  & 0.054  & 0.057  & 0.060  & 0.054  & 0.054  & 0.055  & 0.058  & 0.060  & 0.054  & 0.055\tabularnewline
1000  & -0.60  & 0.050  & 0.054  & 0.055  & 0.055  & 0.053  & 0.050  & 0.054  & 0.055  & 0.055  & 0.053\tabularnewline
\bottomrule
\end{tabular}\end{adjustbox} \label{tab:sim-dgp1-nonlinear-alpha-1p00} 
\end{table}

\end{landscape} \clearpage{}

\begin{landscape}

\begin{table}[htbp]
\centering \caption{Empirical size (nominal size 5\%): \textbf{DGP1}; $\alpha=0.80$;
$[f(x),\ w(x)]=[x^{\left\langle 0.5\right\rangle },\ x^{\left\langle 0.5\right\rangle }]$;
$d_{*}=-0.25$ }
\begin{adjustbox}{max totalsize={\linewidth}{0.77\textwidth},center}
\begin{tabular}{rrrrrrrrrrrr}
 &  & \multicolumn{5}{c}{W-IVX} & \multicolumn{5}{c}{OLS}\tabularnewline
\midrule 
$n$  & $d$  & $\delta$=-0.99  & -0.75  & 0.00  & 0.75  & 0.99  & $\delta$=-0.99  & -0.75  & 0.00  & 0.75  & 0.99\tabularnewline
\midrule
\midrule 
250  & -1.05  & 0.051  & 0.050  & 0.054  & 0.052  & 0.049  & 0.051  & 0.051  & 0.054  & 0.052  & 0.049\tabularnewline
500  & -1.05  & 0.050  & 0.049  & 0.049  & 0.052  & 0.051  & 0.050  & 0.049  & 0.049  & 0.052  & 0.051\tabularnewline
1000  & -1.05  & 0.045  & 0.048  & 0.049  & 0.049  & 0.044  & 0.044  & 0.048  & 0.049  & 0.049  & 0.044\tabularnewline
\midrule 
250  & -0.85  & 0.050  & 0.051  & 0.052  & 0.051  & 0.048  & 0.051  & 0.051  & 0.052  & 0.051  & 0.048\tabularnewline
500  & -0.85  & 0.050  & 0.048  & 0.050  & 0.052  & 0.051  & 0.050  & 0.048  & 0.050  & 0.052  & 0.051\tabularnewline
1000  & -0.85  & 0.045  & 0.048  & 0.049  & 0.048  & 0.046  & 0.045  & 0.048  & 0.050  & 0.048  & 0.045\tabularnewline
\bottomrule
\end{tabular}\end{adjustbox} \label{tab:sim-dgp1-nonlinear-alpha-0p80} 
\end{table}

\end{landscape} \clearpage{}

\begin{landscape}

\begin{table}[htbp]
\centering \caption{Empirical size (nominal size 5\%): W-IVX; $[\alpha,\nu]=[2,2.1]$;
$f(x)=x$; $d_{*}=0.5$; \textbf{DGP2}}
\begin{adjustbox}{max totalsize={\linewidth}{0.77\textwidth},center}
\begin{tabular}{rrrrrrrrrrrrrr}
$w(x)$  &  & \multicolumn{3}{c}{$x^{\langle0\rangle}$} & \multicolumn{3}{c}{$x^{\langle0.5\rangle}$} & \multicolumn{3}{c}{$x^{\langle0.75\rangle}$} & \multicolumn{3}{c}{$x^{\langle1\rangle}$}\tabularnewline
\midrule 
$n$  & $d$  & $\delta$=-0.99  & 0.00  & 0.99  & $\delta$=-0.99  & 0.00  & 0.99  & $\delta$=-0.99  & 0.00  & 0.99  & $\delta$=-0.99  & 0.00  & 0.99\tabularnewline
\midrule
\midrule 
250  & -0.10  & 0.052  & 0.049  & 0.049  & 0.052  & 0.052  & 0.052  & 0.050  & 0.053  & 0.050  & 0.052  & 0.051  & 0.049\tabularnewline
500  & -0.10  & 0.048  & 0.052  & 0.054  & 0.050  & 0.050  & 0.052  & 0.052  & 0.050  & 0.052  & 0.052  & 0.050  & 0.051\tabularnewline
1000  & -0.10  & 0.050  & 0.053  & 0.051  & 0.047  & 0.050  & 0.051  & 0.054  & 0.050  & 0.050  & 0.051  & 0.050  & 0.051\tabularnewline
\midrule 
250  & 0.10  & 0.052  & 0.052  & 0.050  & 0.051  & 0.052  & 0.052  & 0.050  & 0.052  & 0.050  & 0.050  & 0.051  & 0.049\tabularnewline
500  & 0.10  & 0.052  & 0.053  & 0.051  & 0.050  & 0.049  & 0.050  & 0.051  & 0.050  & 0.052  & 0.051  & 0.051  & 0.051\tabularnewline
1000  & 0.10  & 0.050  & 0.053  & 0.049  & 0.049  & 0.051  & 0.052  & 0.054  & 0.049  & 0.048  & 0.051  & 0.051  & 0.052\tabularnewline
\bottomrule
\end{tabular}\end{adjustbox} \label{tab:sim-dgp2-linear-weights-alpha-2p00} 
\end{table}

\end{landscape} \clearpage{}

\begin{landscape}

\begin{table}[htbp]
\centering \caption{Empirical size (nominal size 5\%): W-IVX; $[\alpha,\nu]=[1.99,2.1]$;
$f(x)=x$; $d_{*}=0.497$; \textbf{DGP2}}
\begin{adjustbox}{max totalsize={\linewidth}{0.77\textwidth},center}
\begin{tabular}{rrrrrrrrrrrrrr}
$w(x)$  &  & \multicolumn{3}{c}{$x^{\langle0\rangle}$} & \multicolumn{3}{c}{$x^{\langle0.5\rangle}$} & \multicolumn{3}{c}{$x^{\langle0.75\rangle}$} & \multicolumn{3}{c}{$x^{\langle1\rangle}$}\tabularnewline
\midrule 
$n$  & $d$  & $\delta$=-0.99  & 0.00  & 0.99  & $\delta$=-0.99  & 0.00  & 0.99  & $\delta$=-0.99  & 0.00  & 0.99  & $\delta$=-0.99  & 0.00  & 0.99\tabularnewline
\midrule
\midrule 
250  & -0.103  & 0.050  & 0.054  & 0.050  & 0.051  & 0.051  & 0.051  & 0.051  & 0.051  & 0.052  & 0.051  & 0.055  & 0.055\tabularnewline
500  & -0.103  & 0.054  & 0.054  & 0.050  & 0.051  & 0.053  & 0.050  & 0.049  & 0.056  & 0.051  & 0.054  & 0.053  & 0.049\tabularnewline
1000  & -0.103  & 0.050  & 0.051  & 0.051  & 0.049  & 0.051  & 0.051  & 0.050  & 0.051  & 0.051  & 0.048  & 0.049  & 0.052\tabularnewline
\midrule 
250  & 0.097  & 0.050  & 0.054  & 0.049  & 0.052  & 0.052  & 0.050  & 0.051  & 0.051  & 0.052  & 0.051  & 0.054  & 0.054\tabularnewline
500  & 0.097  & 0.051  & 0.052  & 0.053  & 0.051  & 0.053  & 0.048  & 0.050  & 0.056  & 0.050  & 0.051  & 0.052  & 0.048\tabularnewline
1000  & 0.097  & 0.051  & 0.052  & 0.049  & 0.050  & 0.051  & 0.049  & 0.048  & 0.050  & 0.051  & 0.049  & 0.048  & 0.050\tabularnewline
\bottomrule
\end{tabular}\end{adjustbox} \label{tab:sim-dgp2-linear-weights-alpha-1p99} 
\end{table}

\end{landscape} \clearpage{}

\begin{landscape}

\begin{table}[htbp]
\centering \caption{Empirical size (nominal size 5\%): W-IVX; $[\alpha,\nu]=[1.8,2.1]$;
$f(x)=x$; $d_{*}=0.444$; \textbf{DGP2}}
\begin{adjustbox}{max totalsize={\linewidth}{0.77\textwidth},center}
\begin{tabular}{rrrrrrrrrrrrrr}
$w(x)$  &  & \multicolumn{3}{c}{$x^{\langle0\rangle}$} & \multicolumn{3}{c}{$x^{\langle0.5\rangle}$} & \multicolumn{3}{c}{$x^{\langle0.75\rangle}$} & \multicolumn{3}{c}{$x^{\langle1\rangle}$}\tabularnewline
\midrule 
$n$  & $d$  & $\delta$=-0.99  & 0.00  & 0.99  & $\delta$=-0.99  & 0.00  & 0.99  & $\delta$=-0.99  & 0.00  & 0.99  & $\delta$=-0.99  & 0.00  & 0.99\tabularnewline
\midrule
\midrule 
250  & -0.156  & 0.051  & 0.050  & 0.050  & 0.050  & 0.052  & 0.051  & 0.051  & 0.050  & 0.051  & 0.055  & 0.053  & 0.049\tabularnewline
500  & -0.156  & 0.050  & 0.052  & 0.052  & 0.050  & 0.051  & 0.053  & 0.047  & 0.051  & 0.051  & 0.053  & 0.047  & 0.053\tabularnewline
1000  & -0.156  & 0.048  & 0.049  & 0.049  & 0.053  & 0.050  & 0.050  & 0.051  & 0.051  & 0.053  & 0.050  & 0.049  & 0.050\tabularnewline
\midrule 
250  & 0.044  & 0.053  & 0.052  & 0.048  & 0.051  & 0.053  & 0.050  & 0.051  & 0.050  & 0.049  & 0.054  & 0.054  & 0.050\tabularnewline
500  & 0.044  & 0.053  & 0.050  & 0.053  & 0.051  & 0.049  & 0.052  & 0.046  & 0.051  & 0.049  & 0.050  & 0.049  & 0.052\tabularnewline
1000  & 0.044  & 0.051  & 0.049  & 0.051  & 0.052  & 0.051  & 0.049  & 0.051  & 0.050  & 0.052  & 0.052  & 0.050  & 0.050\tabularnewline
\bottomrule
\end{tabular}\end{adjustbox} \label{tab:sim-dgp2-linear-weights-alpha-1p80} 
\end{table}

\end{landscape} \clearpage{}

\begin{landscape}

\begin{table}[htbp]
\centering \caption{Empirical size (nominal size 5\%): W-IVX; $[\alpha,\nu]=[1.6,2.1]$;
$f(x)=x$; $d_{*}=0.375$; \textbf{DGP2}}
\begin{adjustbox}{max totalsize={\linewidth}{0.77\textwidth},center}
\begin{tabular}{rrrrrrrrrrrrrr}
$w(x)$  &  & \multicolumn{3}{c}{$x^{\langle0\rangle}$} & \multicolumn{3}{c}{$x^{\langle0.5\rangle}$} & \multicolumn{3}{c}{$x^{\langle0.75\rangle}$} & \multicolumn{3}{c}{$x^{\langle1\rangle}$}\tabularnewline
\midrule 
$n$  & $d$  & $\delta$=-0.99  & 0.00  & 0.99  & $\delta$=-0.99  & 0.00  & 0.99  & $\delta$=-0.99  & 0.00  & 0.99  & $\delta$=-0.99  & 0.00  & 0.99\tabularnewline
\midrule
\midrule 
250  & -0.225  & 0.051  & 0.050  & 0.050  & 0.052  & 0.051  & 0.051  & 0.050  & 0.050  & 0.049  & 0.051  & 0.053  & 0.051\tabularnewline
500  & -0.225  & 0.048  & 0.052  & 0.052  & 0.048  & 0.050  & 0.048  & 0.051  & 0.052  & 0.050  & 0.050  & 0.052  & 0.049\tabularnewline
1000  & -0.225  & 0.049  & 0.051  & 0.049  & 0.050  & 0.050  & 0.051  & 0.051  & 0.052  & 0.051  & 0.049  & 0.050  & 0.051\tabularnewline
\midrule 
250  & -0.025  & 0.053  & 0.053  & 0.053  & 0.051  & 0.051  & 0.049  & 0.051  & 0.051  & 0.050  & 0.049  & 0.052  & 0.052\tabularnewline
500  & -0.025  & 0.047  & 0.052  & 0.051  & 0.050  & 0.052  & 0.050  & 0.048  & 0.052  & 0.051  & 0.049  & 0.052  & 0.050\tabularnewline
1000  & -0.025  & 0.048  & 0.049  & 0.048  & 0.051  & 0.051  & 0.052  & 0.051  & 0.052  & 0.050  & 0.048  & 0.049  & 0.049\tabularnewline
\bottomrule
\end{tabular}\end{adjustbox} \label{tab:sim-dgp2-linear-weights-alpha-1p60} 
\end{table}

\end{landscape} \clearpage{}

\begin{landscape}

\begin{table}[htbp]
\centering \caption{Empirical size (nominal size 5\%): W-IVX; $[\alpha,\nu]=[1.4,2.1]$;
$f(x)=x$; $d_{*}=0.286$; \textbf{DGP2}}
\begin{adjustbox}{max totalsize={\linewidth}{0.77\textwidth},center}
\begin{tabular}{rrrrrrrrrrrrrr}
$w(x)$  &  & \multicolumn{3}{c}{$x^{\langle0\rangle}$} & \multicolumn{3}{c}{$x^{\langle0.5\rangle}$} & \multicolumn{3}{c}{$x^{\langle0.75\rangle}$} & \multicolumn{3}{c}{$x^{\langle1\rangle}$}\tabularnewline
\midrule 
$n$  & $d$  & $\delta$=-0.99  & 0.00  & 0.99  & $\delta$=-0.99  & 0.00  & 0.99  & $\delta$=-0.99  & 0.00  & 0.99  & $\delta$=-0.99  & 0.00  & 0.99\tabularnewline
\midrule
\midrule 
250  & -0.314  & 0.049  & 0.050  & 0.052  & 0.052  & 0.056  & 0.051  & 0.052  & 0.050  & 0.049  & 0.048  & 0.050  & 0.052\tabularnewline
500  & -0.314  & 0.050  & 0.052  & 0.050  & 0.051  & 0.052  & 0.048  & 0.046  & 0.048  & 0.050  & 0.050  & 0.050  & 0.051\tabularnewline
1000  & -0.314  & 0.048  & 0.050  & 0.050  & 0.052  & 0.049  & 0.050  & 0.049  & 0.049  & 0.049  & 0.053  & 0.051  & 0.051\tabularnewline
\midrule 
250  & -0.114  & 0.049  & 0.050  & 0.051  & 0.051  & 0.055  & 0.052  & 0.052  & 0.049  & 0.049  & 0.048  & 0.049  & 0.051\tabularnewline
500  & -0.114  & 0.050  & 0.052  & 0.051  & 0.051  & 0.053  & 0.048  & 0.048  & 0.048  & 0.051  & 0.050  & 0.051  & 0.050\tabularnewline
1000  & -0.114  & 0.049  & 0.051  & 0.050  & 0.052  & 0.050  & 0.050  & 0.050  & 0.049  & 0.050  & 0.052  & 0.049  & 0.052\tabularnewline
\bottomrule
\end{tabular}\end{adjustbox} \label{tab:sim-dgp2-linear-weights-alpha-1p40} 
\end{table}

\end{landscape} \clearpage{}

\begin{landscape}

\begin{table}[htbp]
\centering \caption{Empirical size (nominal size 5\%): W-IVX; $[\alpha,\nu]=[1.2,2.1]$;
$f(x)=x$; $d_{*}=0.167$; \textbf{DGP2}}
\begin{adjustbox}{max totalsize={\linewidth}{0.77\textwidth},center}
\begin{tabular}{rrrrrrrrrrrrrr}
$w(x)$  &  & \multicolumn{3}{c}{$x^{\langle0\rangle}$} & \multicolumn{3}{c}{$x^{\langle0.5\rangle}$} & \multicolumn{3}{c}{$x^{\langle0.75\rangle}$} & \multicolumn{3}{c}{$x^{\langle1\rangle}$}\tabularnewline
\midrule 
$n$  & $d$  & $\delta$=-0.99  & 0.00  & 0.99  & $\delta$=-0.99  & 0.00  & 0.99  & $\delta$=-0.99  & 0.00  & 0.99  & $\delta$=-0.99  & 0.00  & 0.99\tabularnewline
\midrule
\midrule 
250  & -0.433  & 0.049  & 0.052  & 0.051  & 0.049  & 0.050  & 0.052  & 0.049  & 0.055  & 0.051  & 0.047  & 0.051  & 0.051\tabularnewline
500  & -0.433  & 0.052  & 0.054  & 0.051  & 0.054  & 0.050  & 0.049  & 0.050  & 0.048  & 0.048  & 0.052  & 0.049  & 0.049\tabularnewline
1000  & -0.433  & 0.051  & 0.051  & 0.050  & 0.051  & 0.049  & 0.053  & 0.048  & 0.052  & 0.050  & 0.051  & 0.050  & 0.050\tabularnewline
\midrule 
250  & -0.233  & 0.051  & 0.051  & 0.052  & 0.050  & 0.050  & 0.051  & 0.049  & 0.053  & 0.051  & 0.046  & 0.052  & 0.050\tabularnewline
500  & -0.233  & 0.049  & 0.052  & 0.049  & 0.052  & 0.051  & 0.048  & 0.050  & 0.048  & 0.047  & 0.051  & 0.050  & 0.049\tabularnewline
1000  & -0.233  & 0.049  & 0.050  & 0.051  & 0.052  & 0.049  & 0.052  & 0.049  & 0.052  & 0.050  & 0.051  & 0.050  & 0.049\tabularnewline
\bottomrule
\end{tabular}\end{adjustbox} \label{tab:sim-dgp2-linear-weights-alpha-1p20} 
\end{table}

\end{landscape} \clearpage{}

\begin{landscape}

\begin{table}[htbp]
\centering \caption{Empirical size (nominal size 5\%): W-IVX; $[\alpha,\nu]=[1,2.1]$;
$f(x)=x$; $d_{*}=0$; \textbf{DGP2}}
\begin{adjustbox}{max totalsize={\linewidth}{0.77\textwidth},center}
\begin{tabular}{rrrrrrrrrrrrrr}
$w(x)$  &  & \multicolumn{3}{c}{$x^{\langle0\rangle}$} & \multicolumn{3}{c}{$x^{\langle0.5\rangle}$} & \multicolumn{3}{c}{$x^{\langle0.75\rangle}$} & \multicolumn{3}{c}{$x^{\langle1\rangle}$}\tabularnewline
\midrule 
$n$  & $d$  & $\delta$=-0.99  & 0.00  & 0.99  & $\delta$=-0.99  & 0.00  & 0.99  & $\delta$=-0.99  & 0.00  & 0.99  & $\delta$=-0.99  & 0.00  & 0.99\tabularnewline
\midrule
\midrule 
250  & -0.60  & 0.051  & 0.054  & 0.050  & 0.051  & 0.051  & 0.051  & 0.046  & 0.053  & 0.048  & 0.046  & 0.052  & 0.048\tabularnewline
500  & -0.60  & 0.052  & 0.053  & 0.050  & 0.051  & 0.050  & 0.050  & 0.051  & 0.051  & 0.051  & 0.049  & 0.048  & 0.046\tabularnewline
1000  & -0.60  & 0.049  & 0.049  & 0.051  & 0.050  & 0.053  & 0.049  & 0.049  & 0.051  & 0.051  & 0.047  & 0.049  & 0.045\tabularnewline
\midrule 
250  & -0.40  & 0.049  & 0.053  & 0.049  & 0.052  & 0.052  & 0.052  & 0.047  & 0.052  & 0.047  & 0.046  & 0.049  & 0.047\tabularnewline
500  & -0.40  & 0.050  & 0.051  & 0.049  & 0.050  & 0.051  & 0.050  & 0.051  & 0.051  & 0.050  & 0.050  & 0.047  & 0.046\tabularnewline
1000  & -0.40  & 0.050  & 0.051  & 0.050  & 0.050  & 0.053  & 0.049  & 0.050  & 0.053  & 0.050  & 0.048  & 0.047  & 0.045\tabularnewline
\bottomrule
\end{tabular}\end{adjustbox} \label{tab:sim-dgp2-linear-weights-alpha-1p00} 
\end{table}

\end{landscape} \clearpage{}

\begin{landscape}

\begin{table}[htbp]
\centering \caption{Empirical size (nominal size 5\%): W-IVX; $[\alpha,\nu]=[0.8,2.1]$;
$f(x)=x$; $d_{*}=-0.25$; \textbf{DGP2}}
\begin{adjustbox}{max totalsize={\linewidth}{0.77\textwidth},center}
\begin{tabular}{rrrrrrrrrrrrrr}
$w(x)$  &  & \multicolumn{3}{c}{$x^{\langle0\rangle}$} & \multicolumn{3}{c}{$x^{\langle0.5\rangle}$} & \multicolumn{3}{c}{$x^{\langle0.75\rangle}$} & \multicolumn{3}{c}{$x^{\langle1\rangle}$}\tabularnewline
\midrule 
$n$  & $d$  & $\delta$=-0.99  & 0.00  & 0.99  & $\delta$=-0.99  & 0.00  & 0.99  & $\delta$=-0.99  & 0.00  & 0.99  & $\delta$=-0.99  & 0.00  & 0.99\tabularnewline
\midrule
\midrule 
250  & -0.85  & 0.053  & 0.053  & 0.054  & 0.052  & 0.053  & 0.050  & 0.047  & 0.050  & 0.043  & 0.045  & 0.045  & 0.043\tabularnewline
500  & -0.85  & 0.053  & 0.051  & 0.051  & 0.050  & 0.051  & 0.050  & 0.047  & 0.047  & 0.047  & 0.045  & 0.045  & 0.044\tabularnewline
1000  & -0.85  & 0.050  & 0.051  & 0.052  & 0.049  & 0.051  & 0.051  & 0.047  & 0.048  & 0.047  & 0.044  & 0.047  & 0.044\tabularnewline
\midrule 
250  & -0.65  & 0.052  & 0.052  & 0.052  & 0.051  & 0.052  & 0.050  & 0.047  & 0.050  & 0.043  & 0.046  & 0.045  & 0.043\tabularnewline
500  & -0.65  & 0.052  & 0.052  & 0.051  & 0.050  & 0.052  & 0.049  & 0.047  & 0.048  & 0.048  & 0.045  & 0.044  & 0.043\tabularnewline
1000  & -0.65  & 0.049  & 0.053  & 0.050  & 0.048  & 0.051  & 0.051  & 0.047  & 0.049  & 0.048  & 0.045  & 0.047  & 0.044\tabularnewline
\bottomrule
\end{tabular}\end{adjustbox} \label{tab:sim-dgp2-linear-weights-alpha-0p80} 
\end{table}

\end{landscape} \clearpage{}

\begin{landscape} 

\subsection*{Empirical Power}


\begin{figure}[H]
\caption{Asymptotic Power (nominal size 5\%): $[f(x),w(x)]=[x,x^{\langle0\rangle}]$,
$\delta=-0.95$, $n=1000$}
\label{fig:sim-power-linear-w00}\centering \captionsetup[subfigure]{position=top}

\begin{subfigure}[t]{0.45\linewidth} \centering \caption{$\alpha=2$, $d_{*}=0.5$}
\label{fig:sim-power-linear-w00-alpha-2p00} \includegraphics[width=0.49\linewidth]{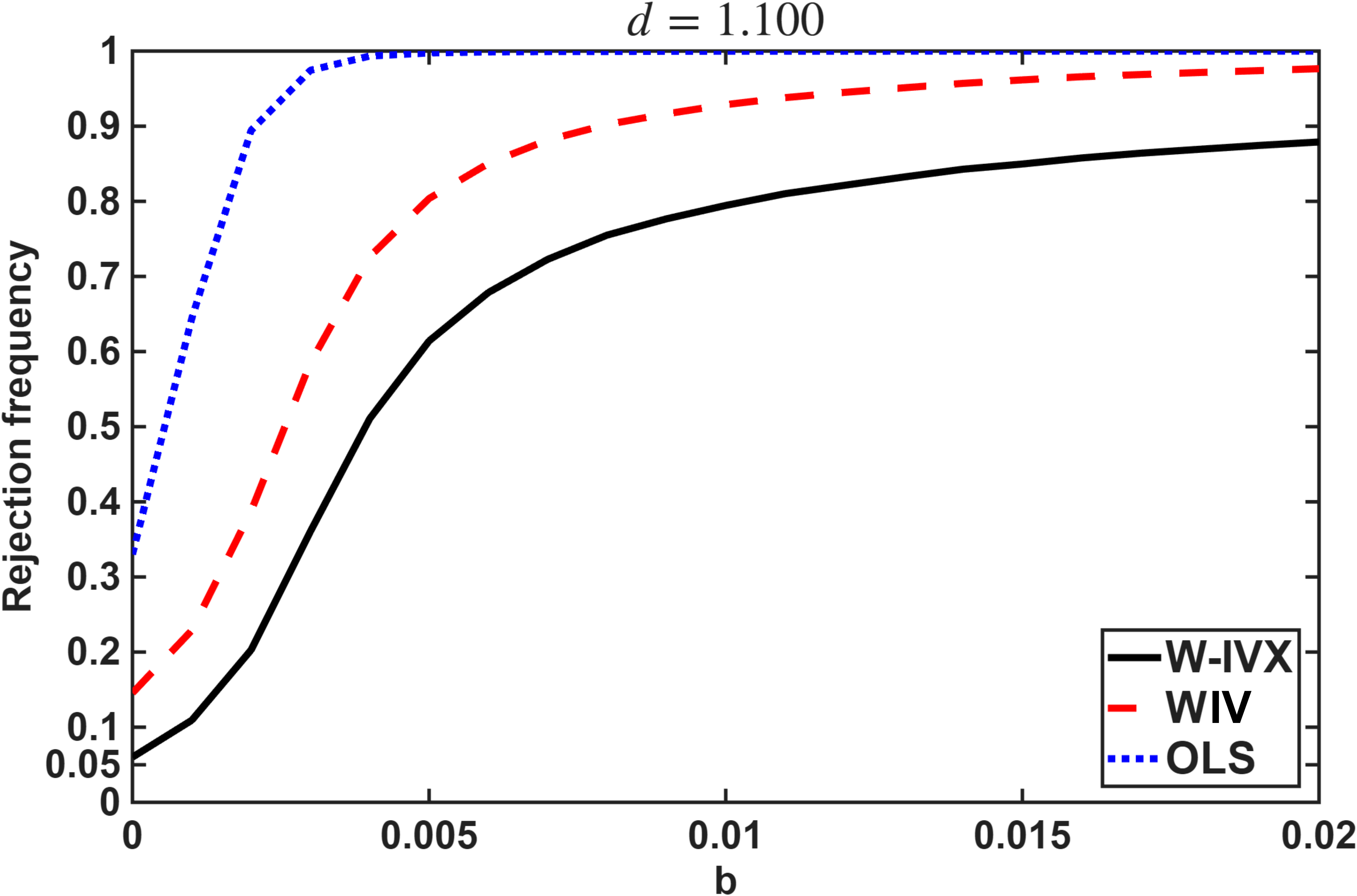}
\includegraphics[width=0.49\linewidth]{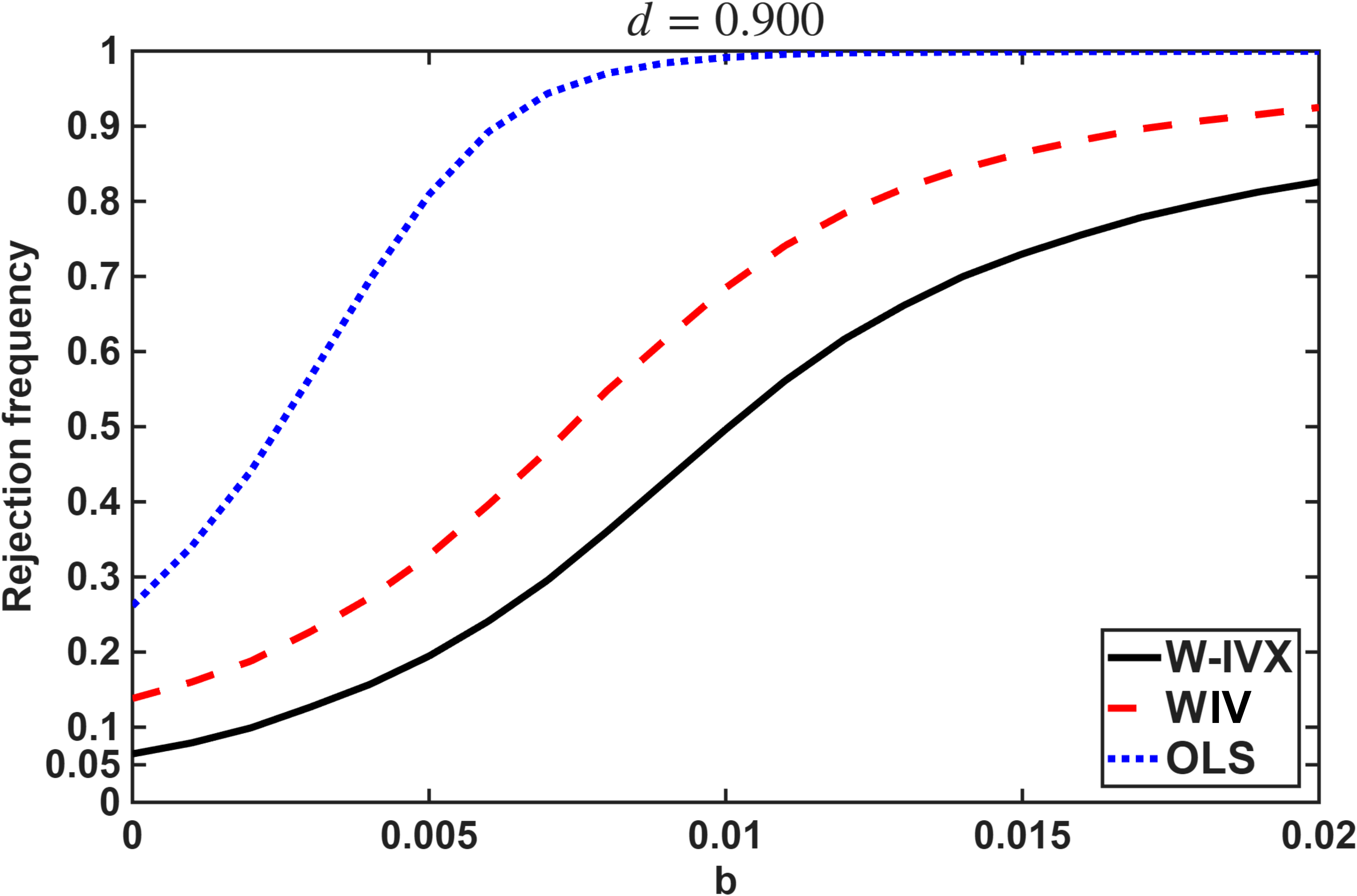}\\[0.6ex]
\includegraphics[width=0.49\linewidth]{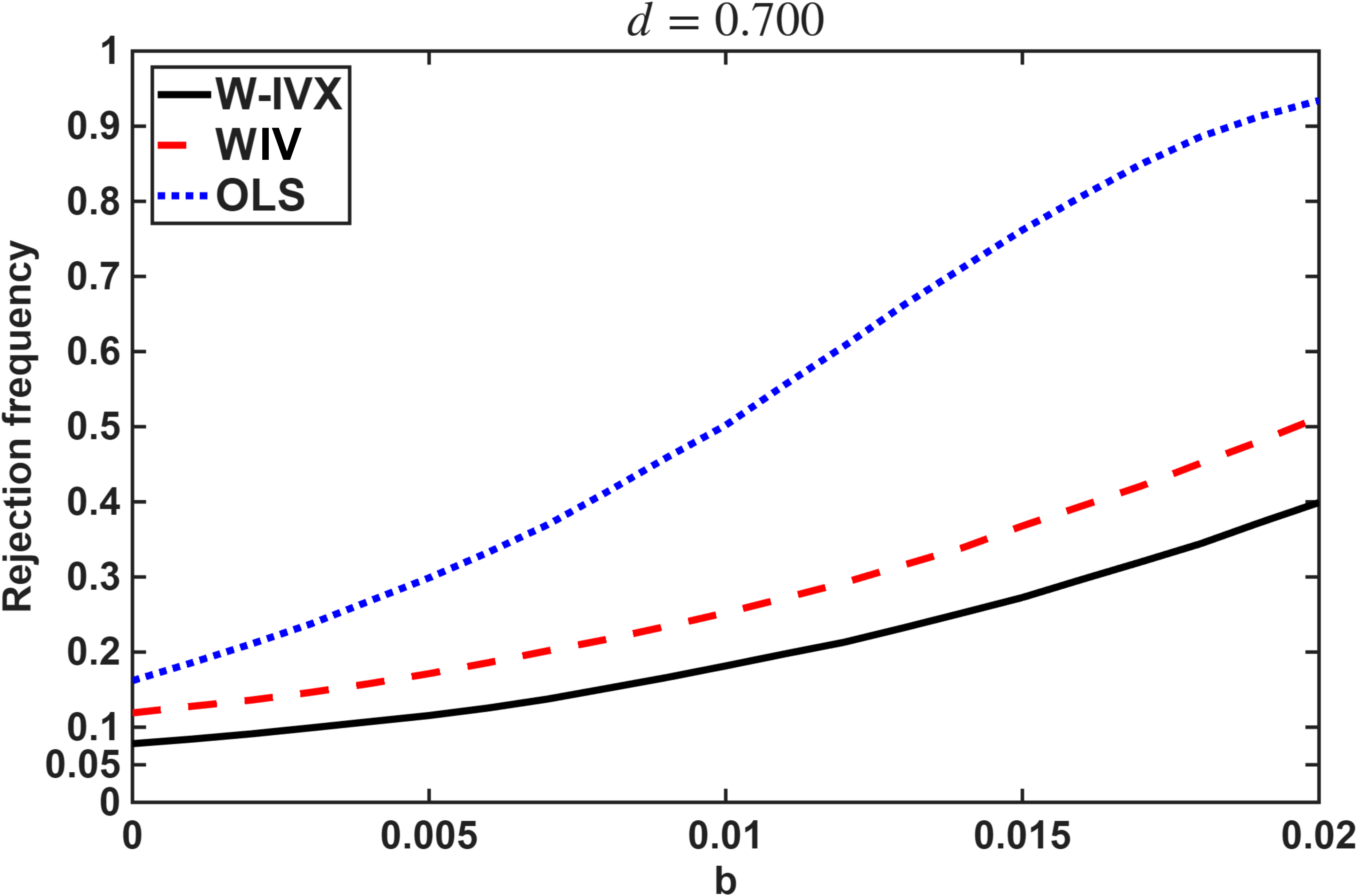}
\includegraphics[width=0.49\linewidth]{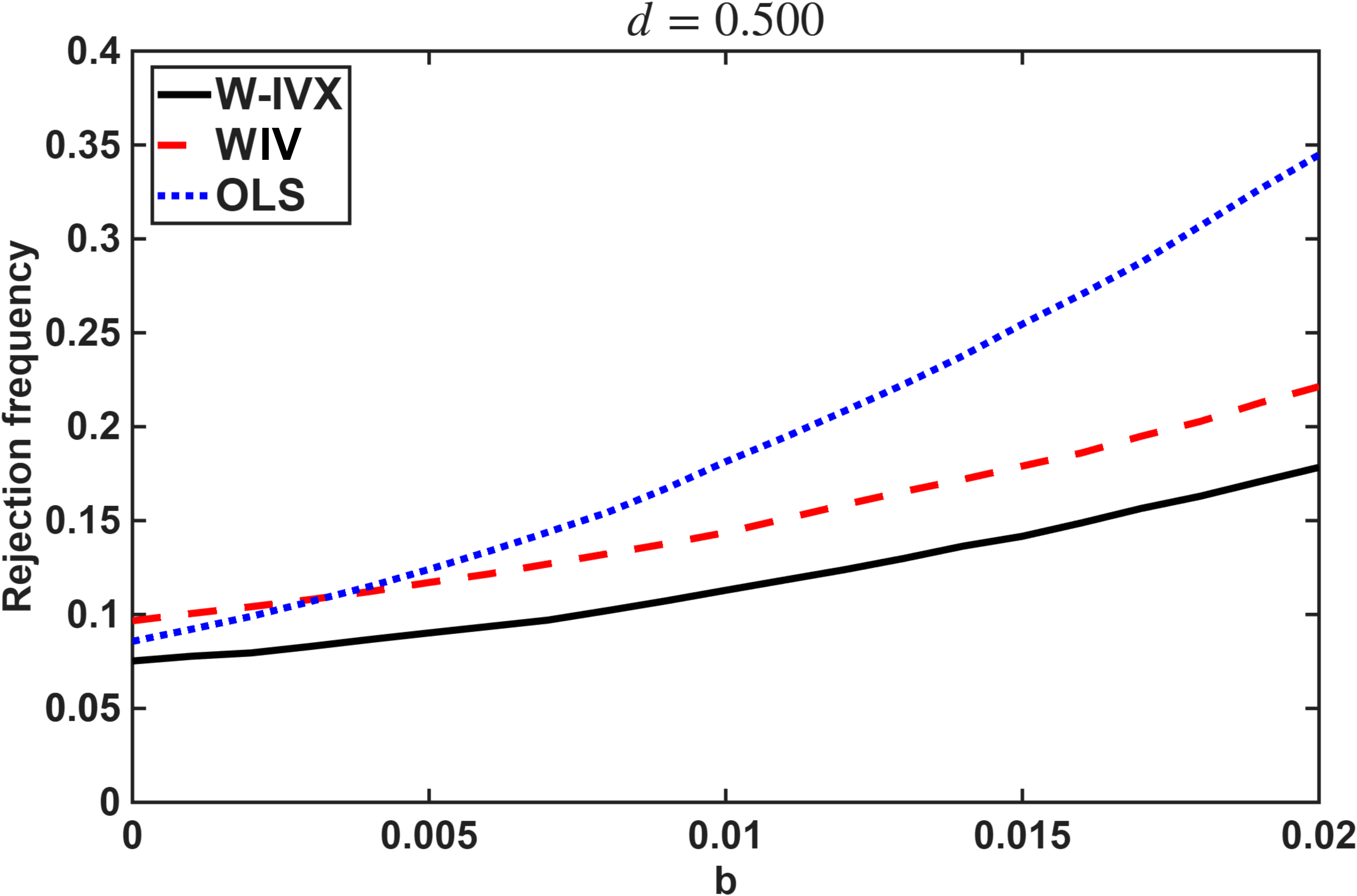}\\[0.6ex]
\includegraphics[width=0.49\linewidth]{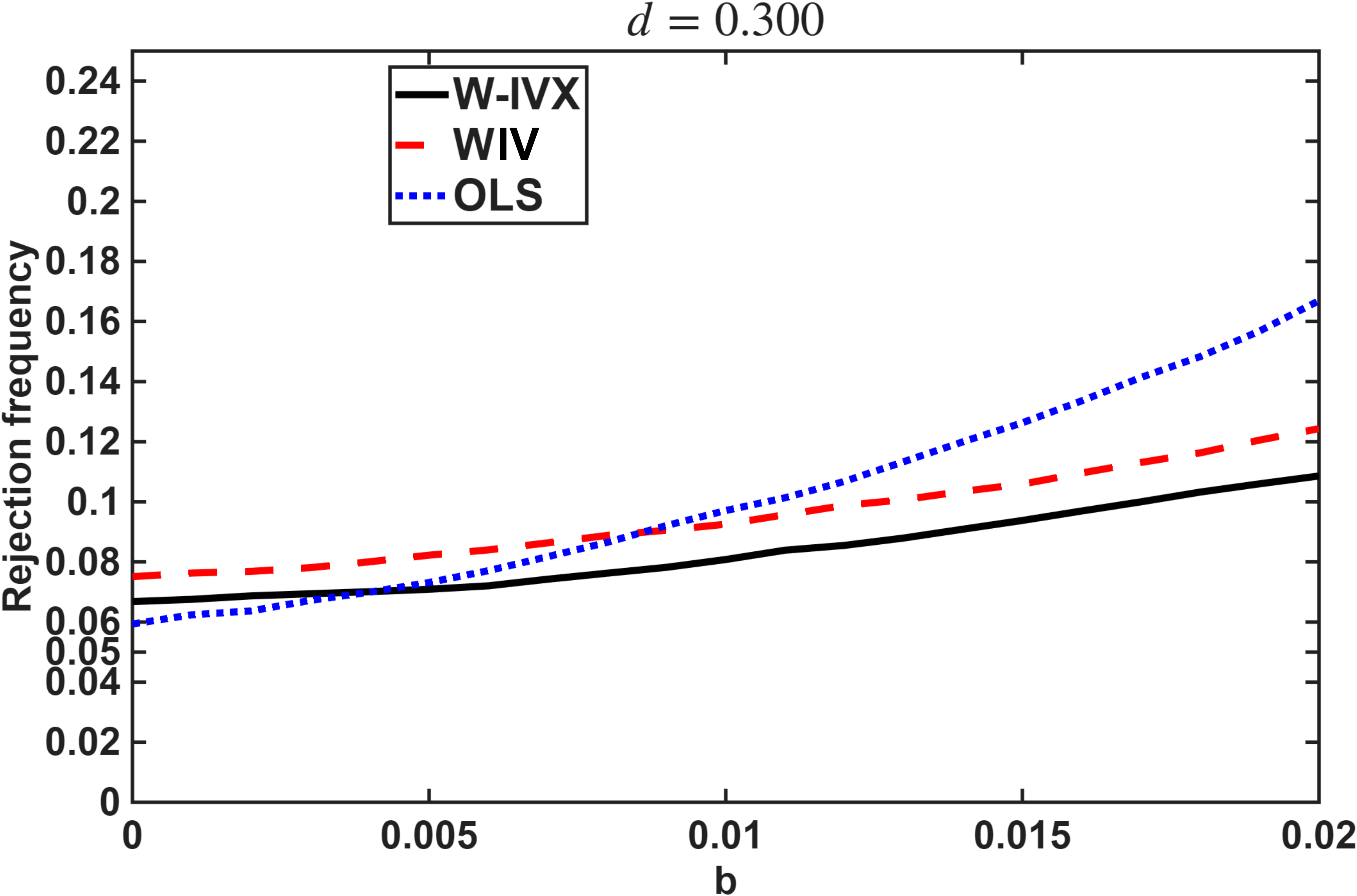}
\includegraphics[width=0.49\linewidth]{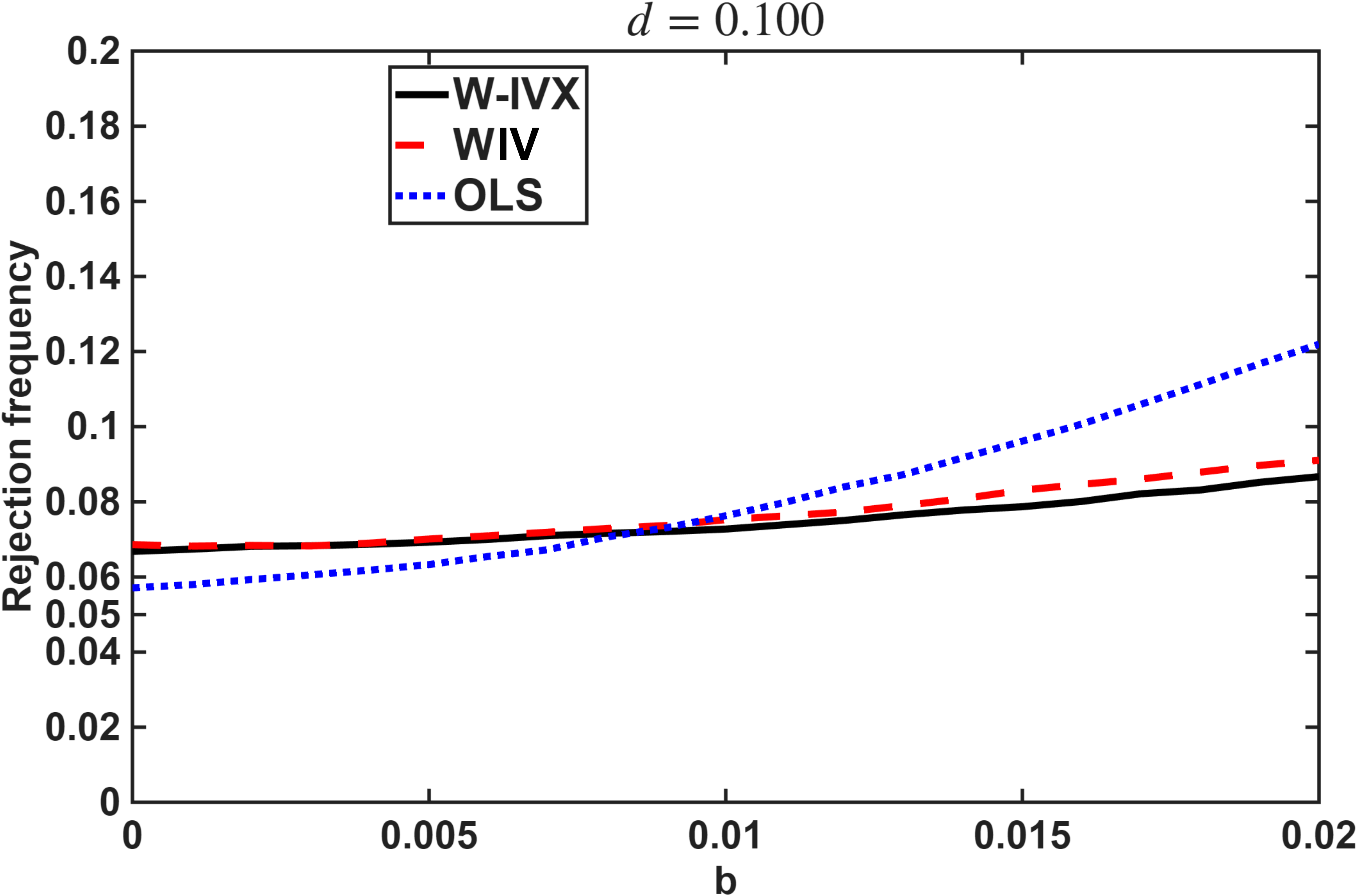}
\end{subfigure} \hfill{}
\begin{subfigure}[t]{0.45\linewidth} \centering \caption{$\alpha=1.8$, $d_{*}=0.444$.}
\label{fig:sim-power-linear-w00-alpha-1p80} \includegraphics[width=0.49\linewidth]{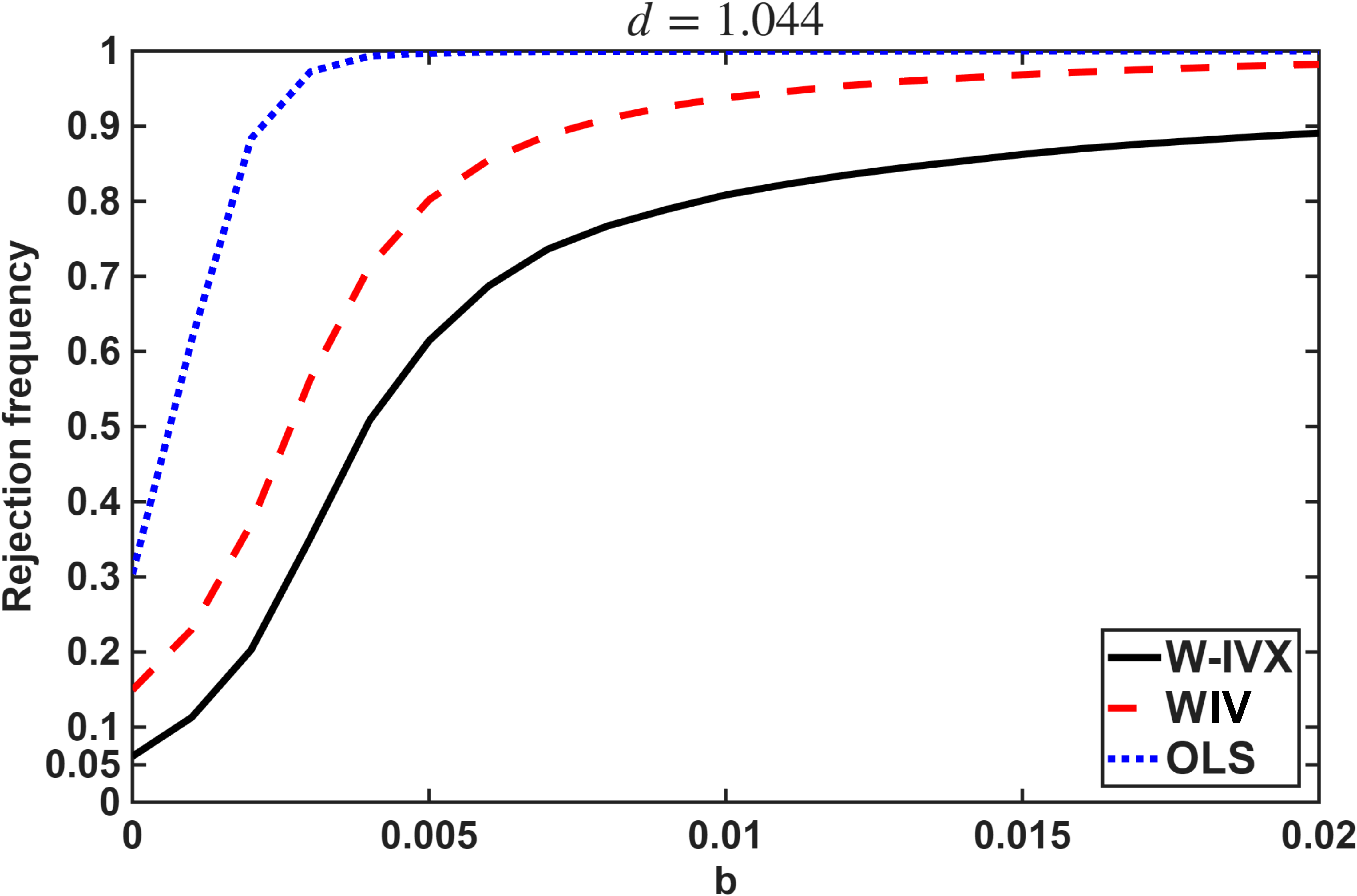}
\includegraphics[width=0.49\linewidth]{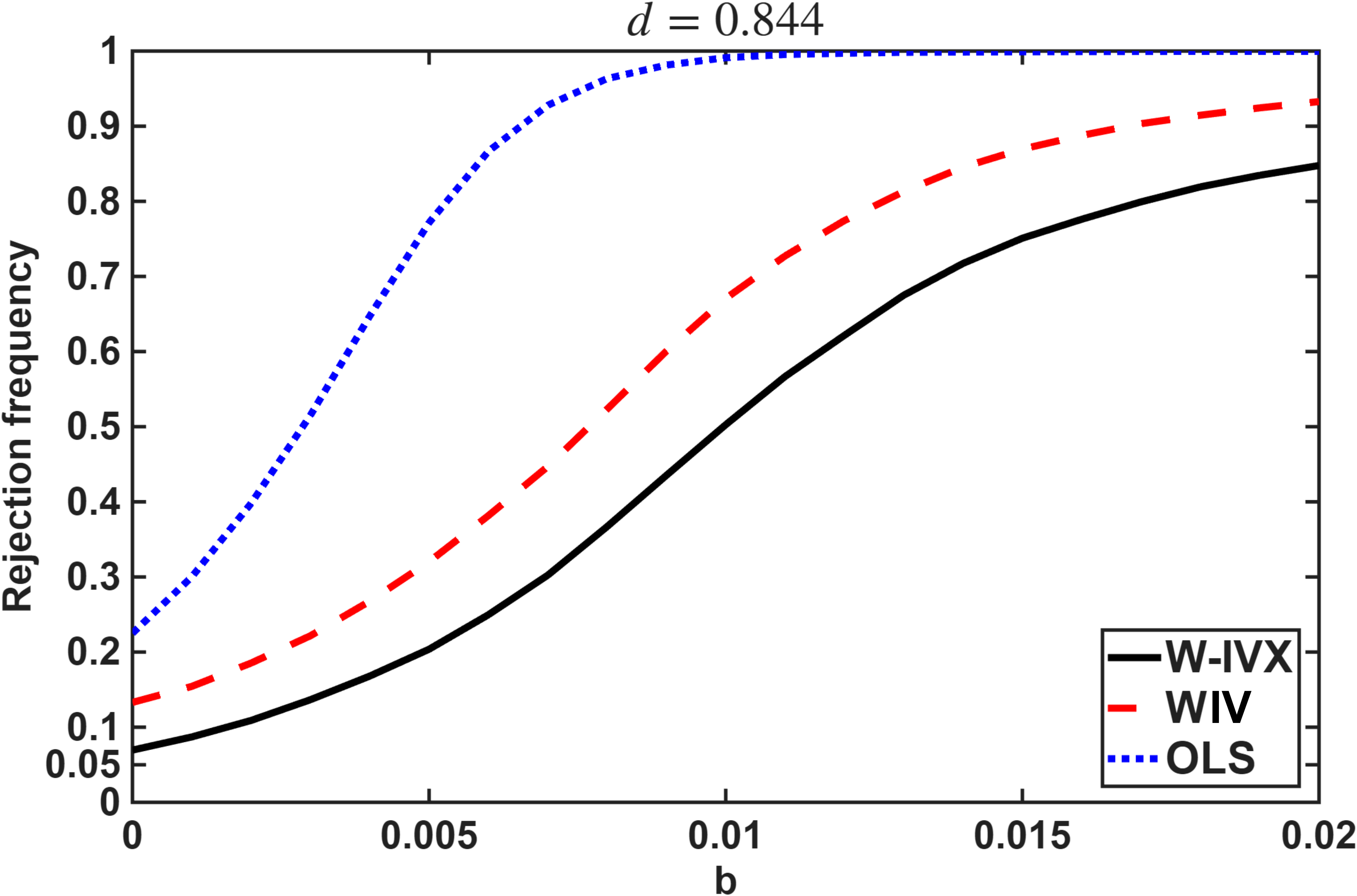}\\[0.6ex]
\includegraphics[width=0.49\linewidth]{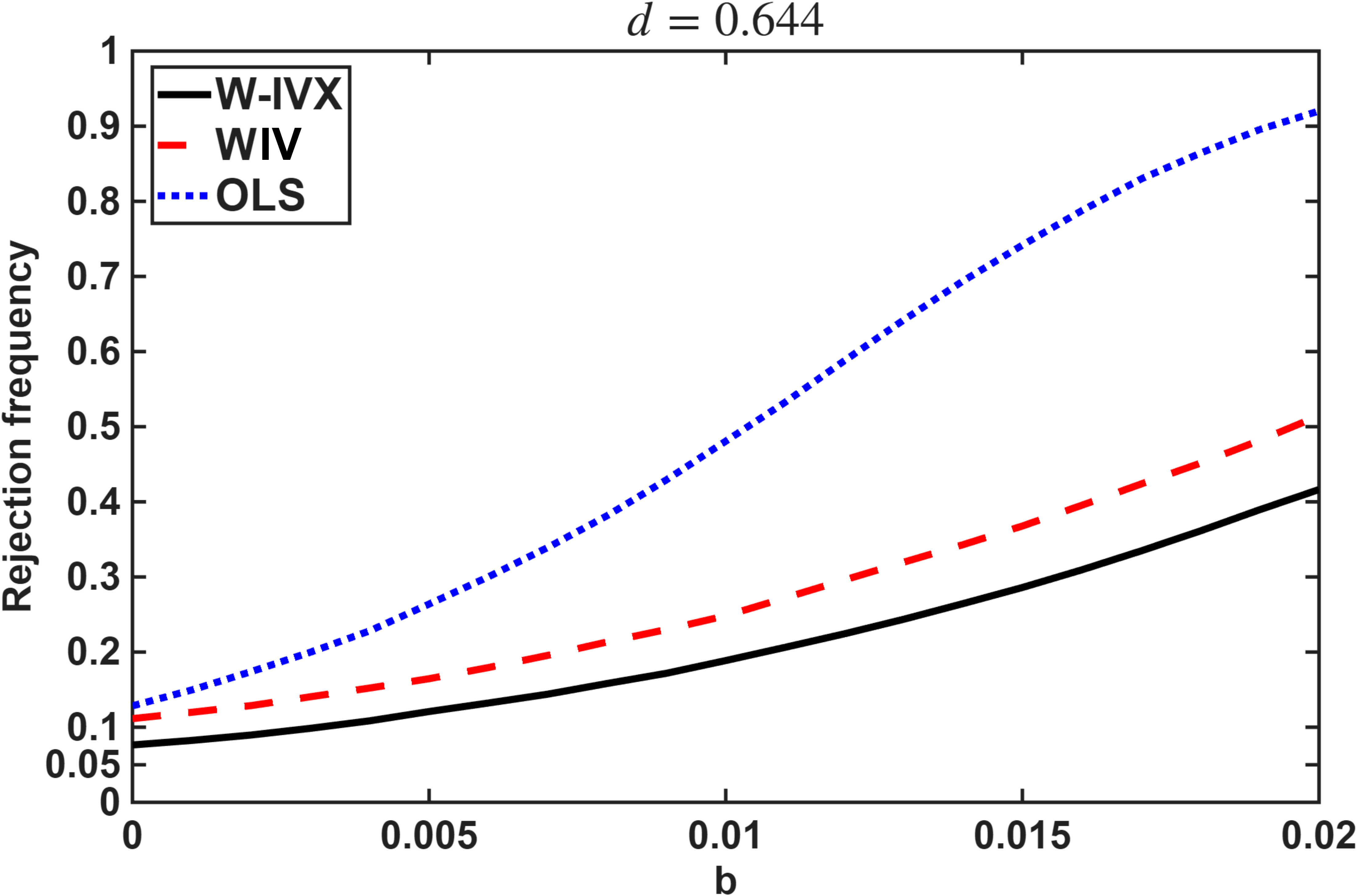}
\includegraphics[width=0.49\linewidth]{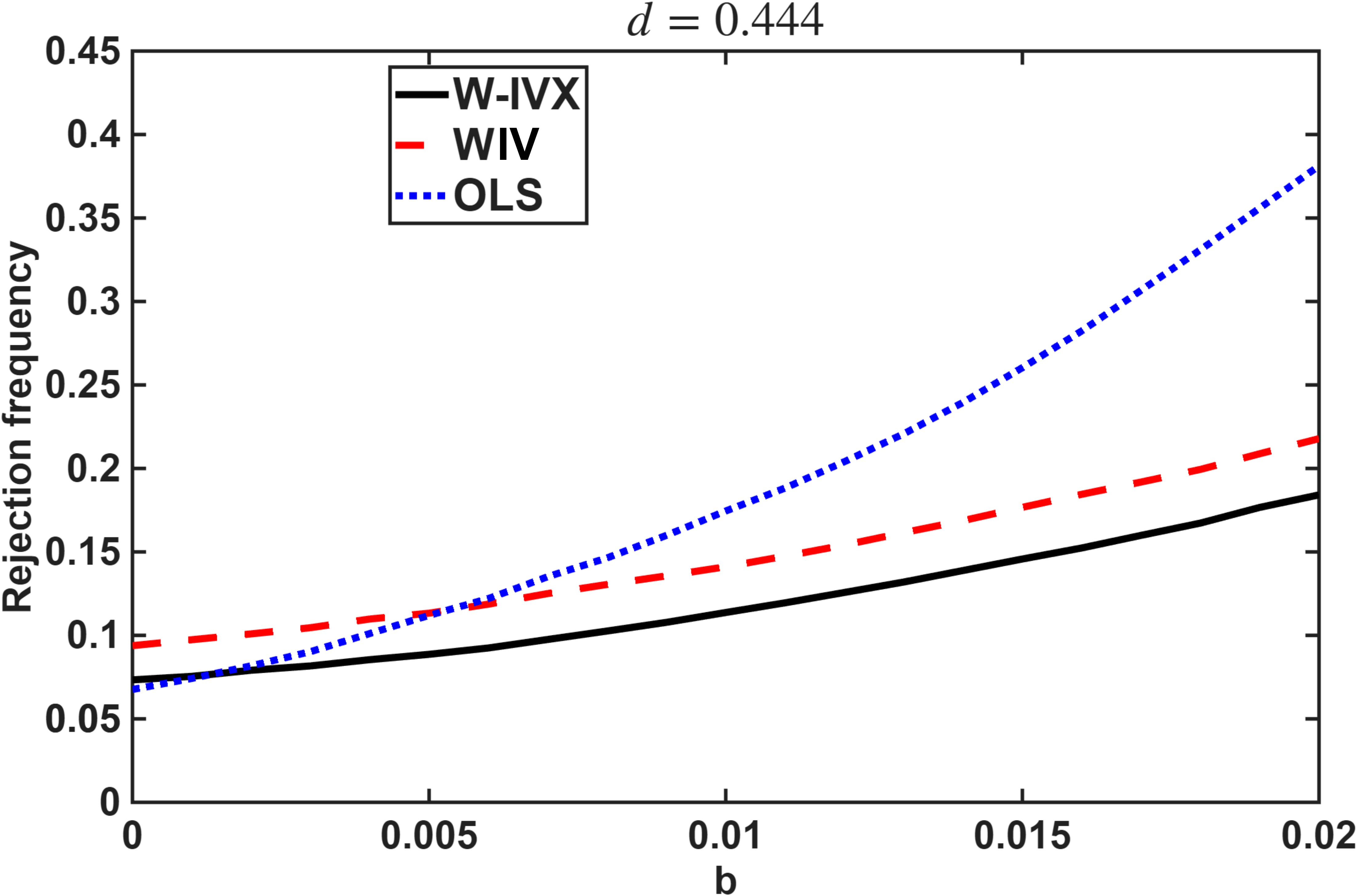}\\[0.6ex]
\includegraphics[width=0.49\linewidth]{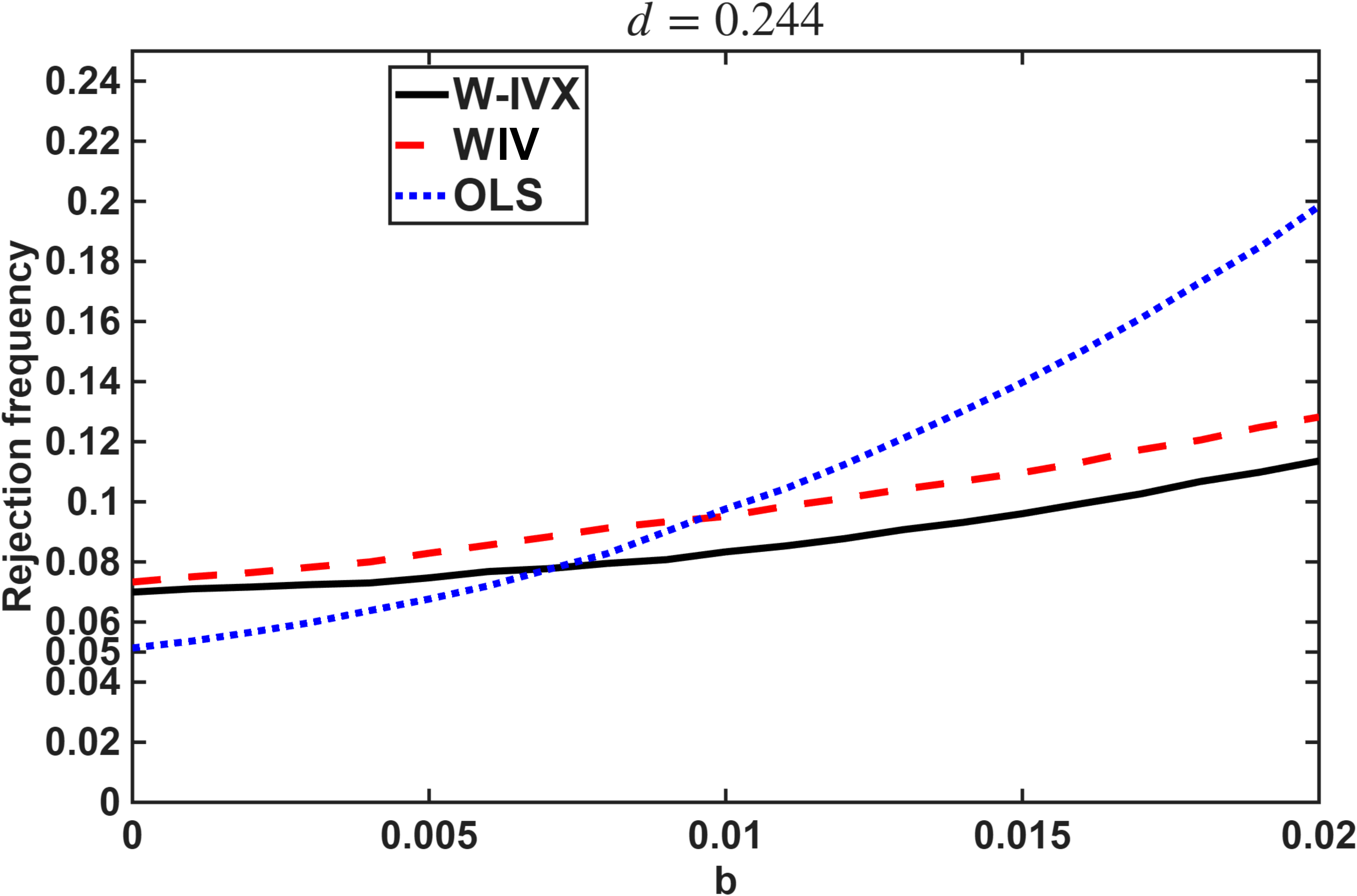}
\includegraphics[width=0.49\linewidth]{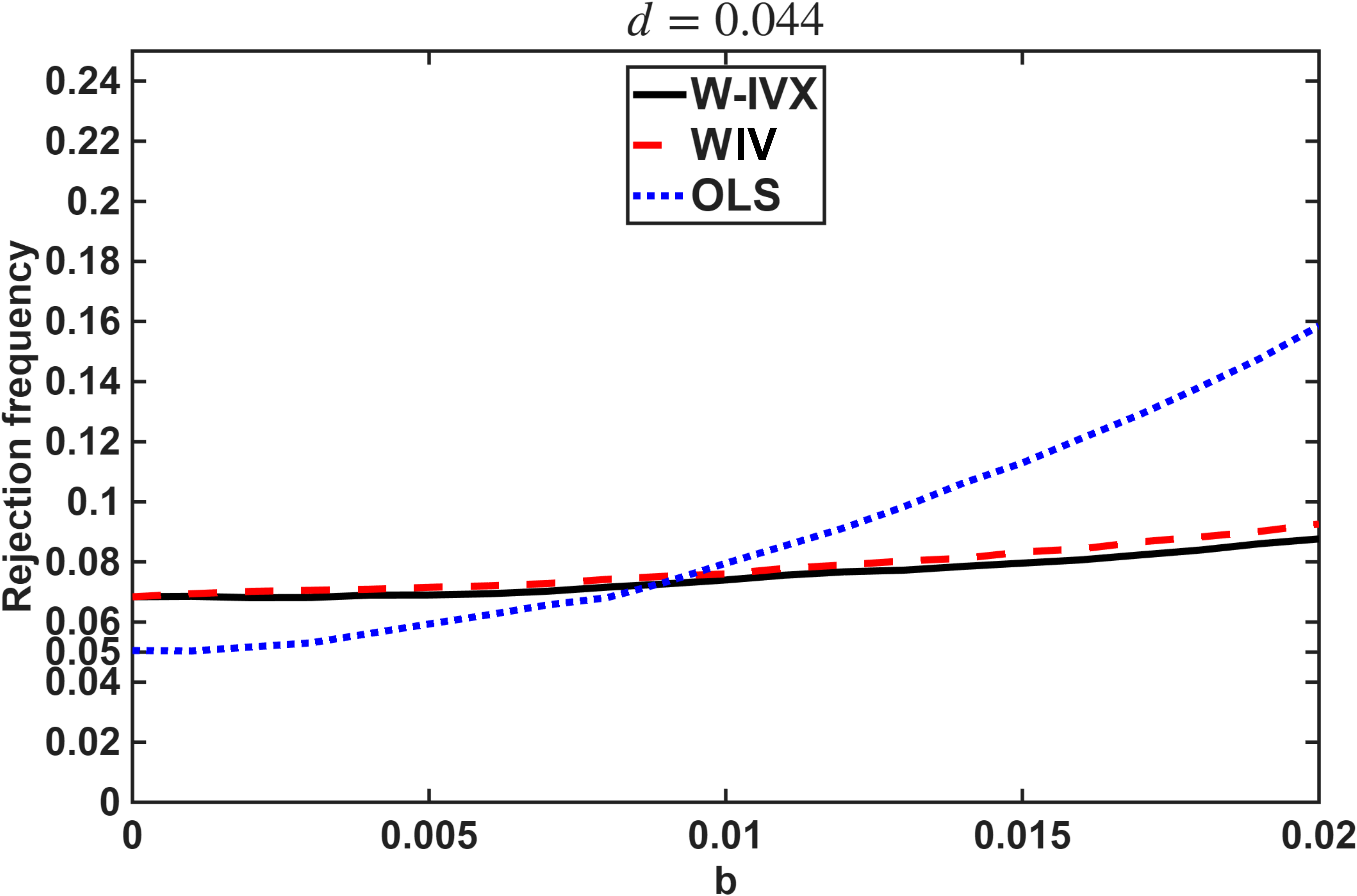}
\end{subfigure}

\end{figure}

\end{landscape} \clearpage{}

 \clearpage{}

 \clearpage{}

 \clearpage{}

\begin{landscape}

\begin{figure}[p]
\caption{Size Adjusted Power (nominal size 5\%): \textbf{DGP1}, $[f(x),w(x)]=[x,x^{\langle0\rangle}]$,
$\delta=-0.95$, $n=1000$}
\label{fig:sim-sap-linear-w00}\centering \captionsetup[subfigure]{position=top}

\begin{subfigure}[t]{0.48\linewidth} \centering \caption{$\alpha=2$, $d_{*}=0.5$}
\label{fig:sim-sap-linear-w00-alpha-2p00} \includegraphics[width=0.49\linewidth]{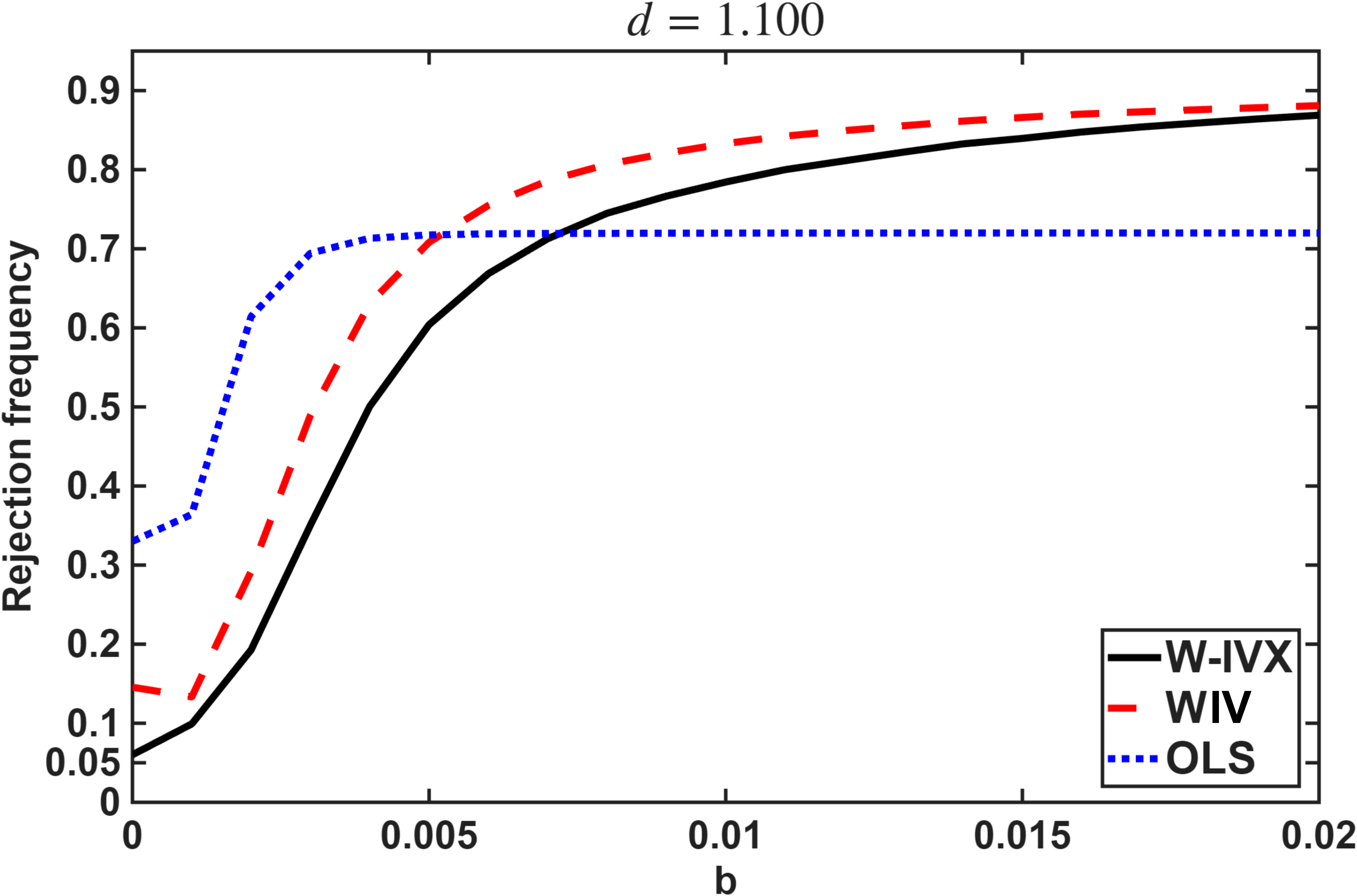}
\includegraphics[width=0.49\linewidth]{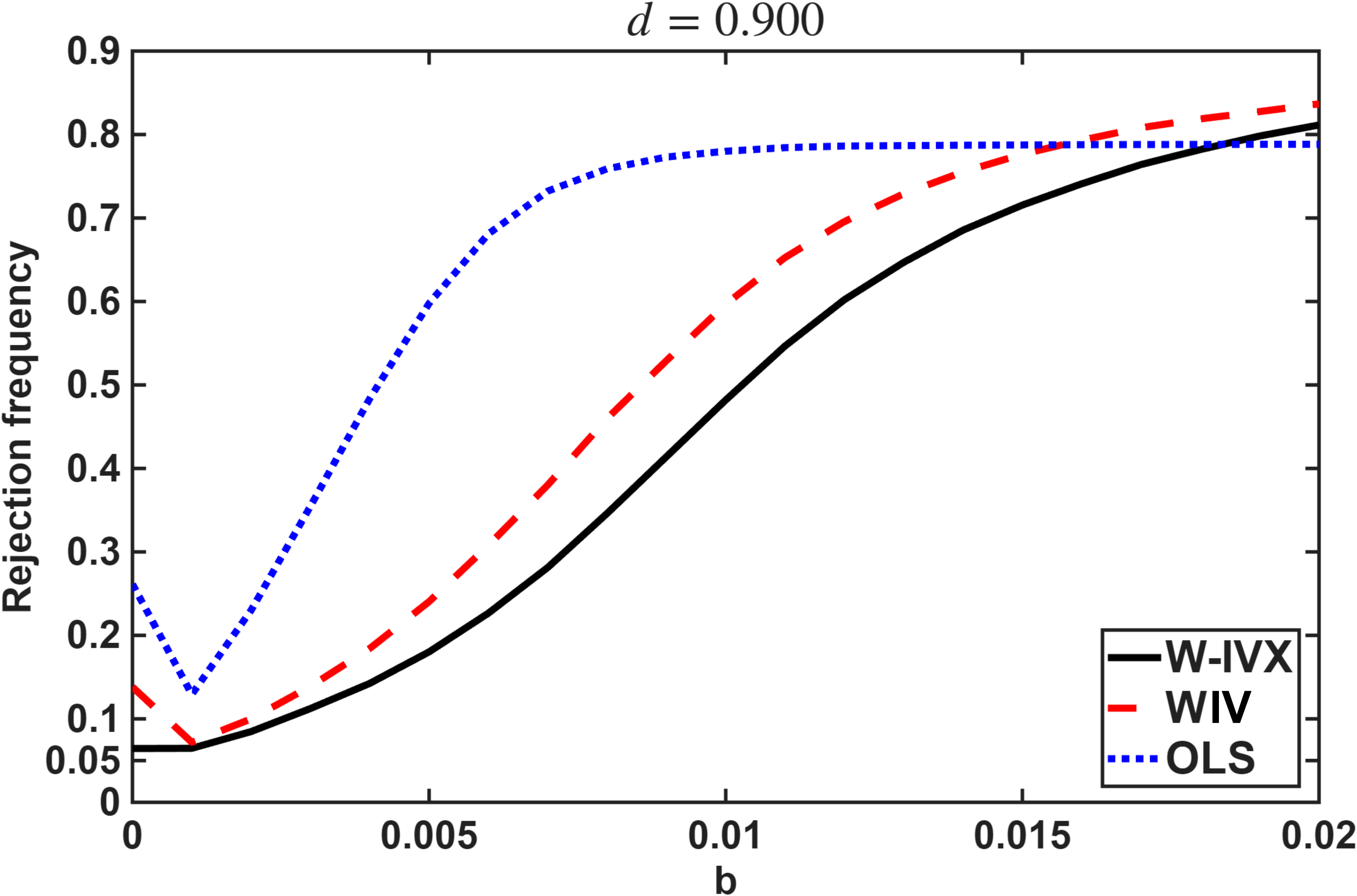}\\[0.6ex]
\includegraphics[width=0.49\linewidth]{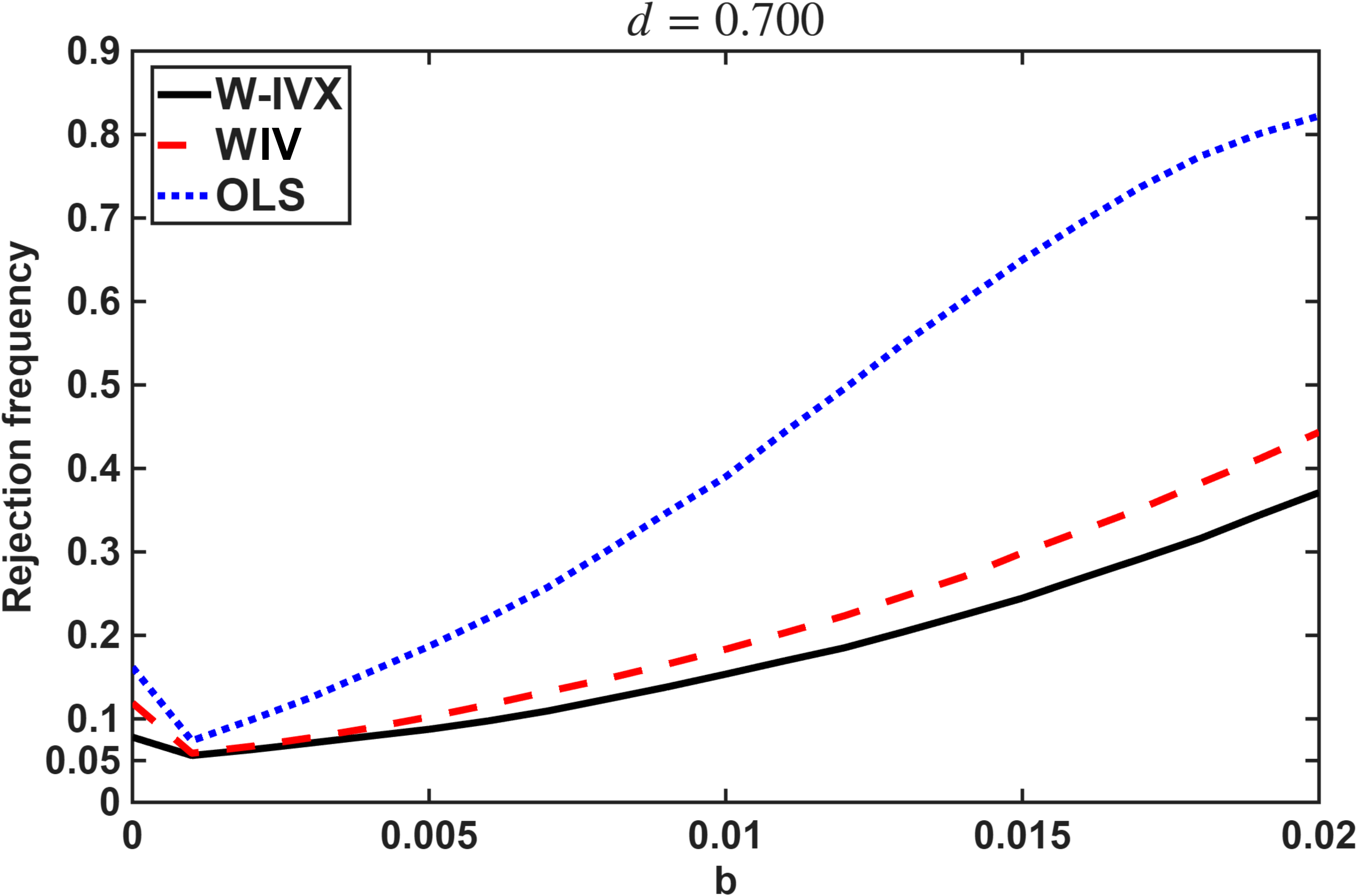}
\includegraphics[width=0.49\linewidth]{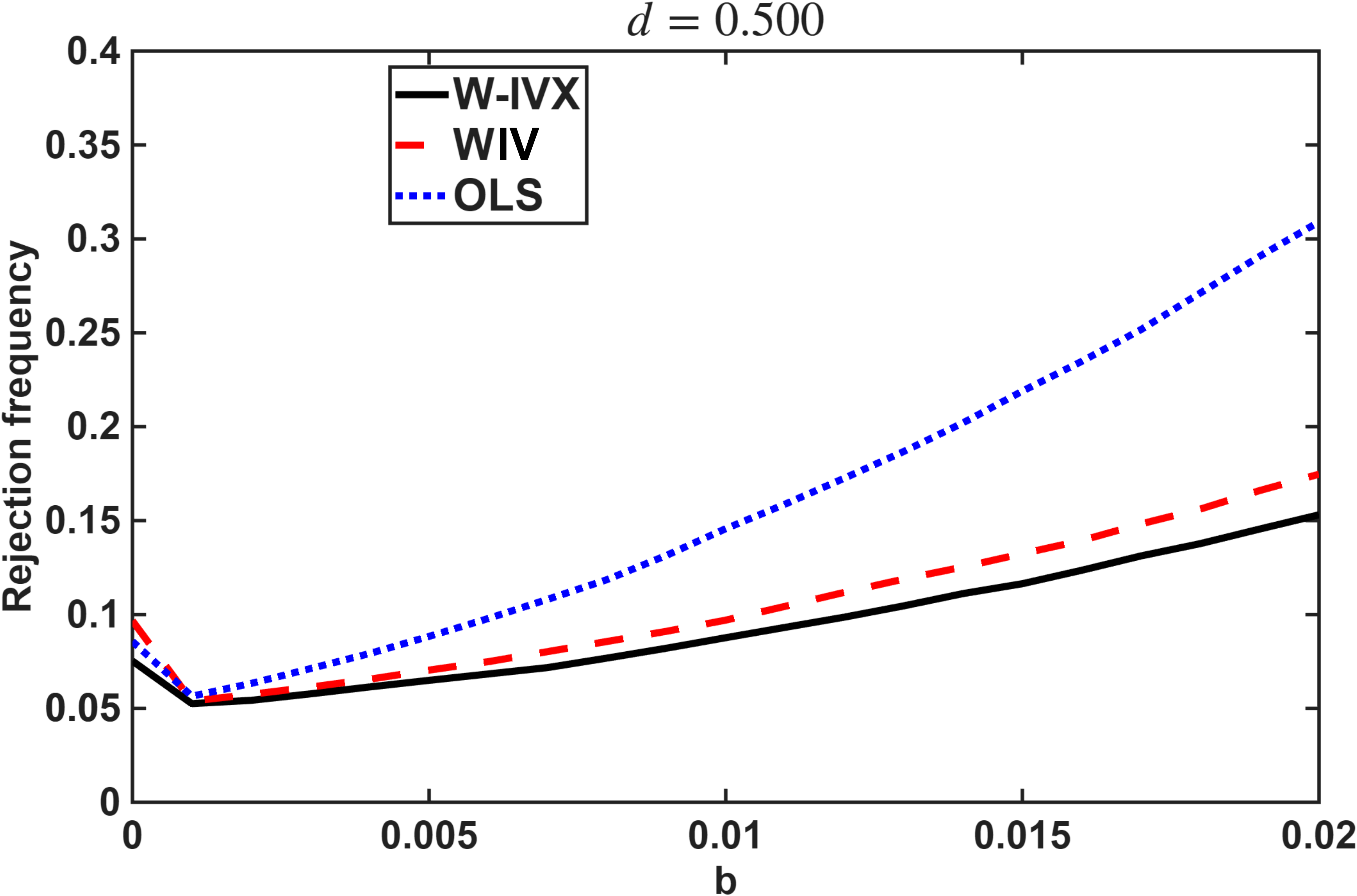}\\[0.6ex]
\includegraphics[width=0.49\linewidth]{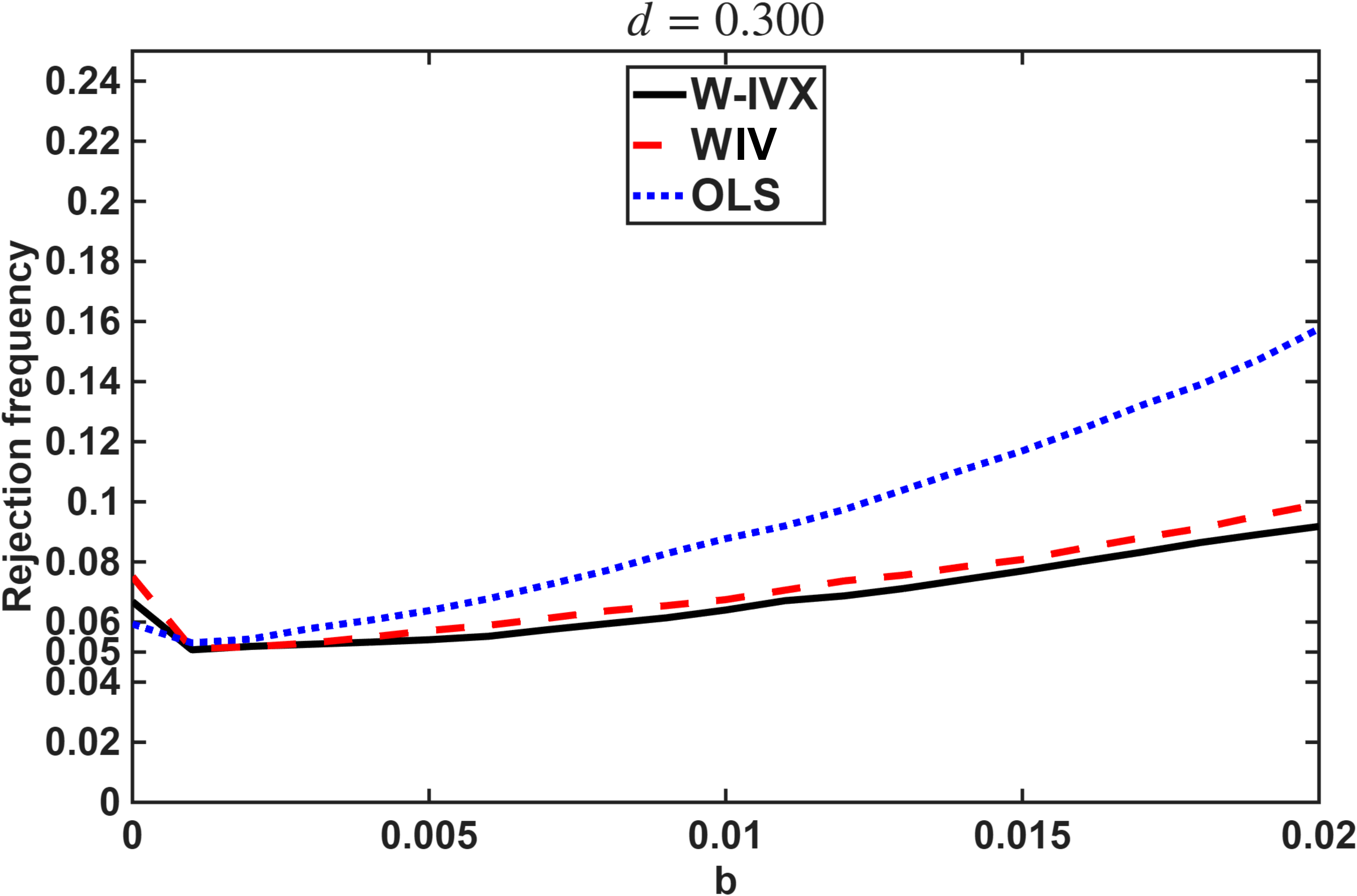}
\includegraphics[width=0.49\linewidth]{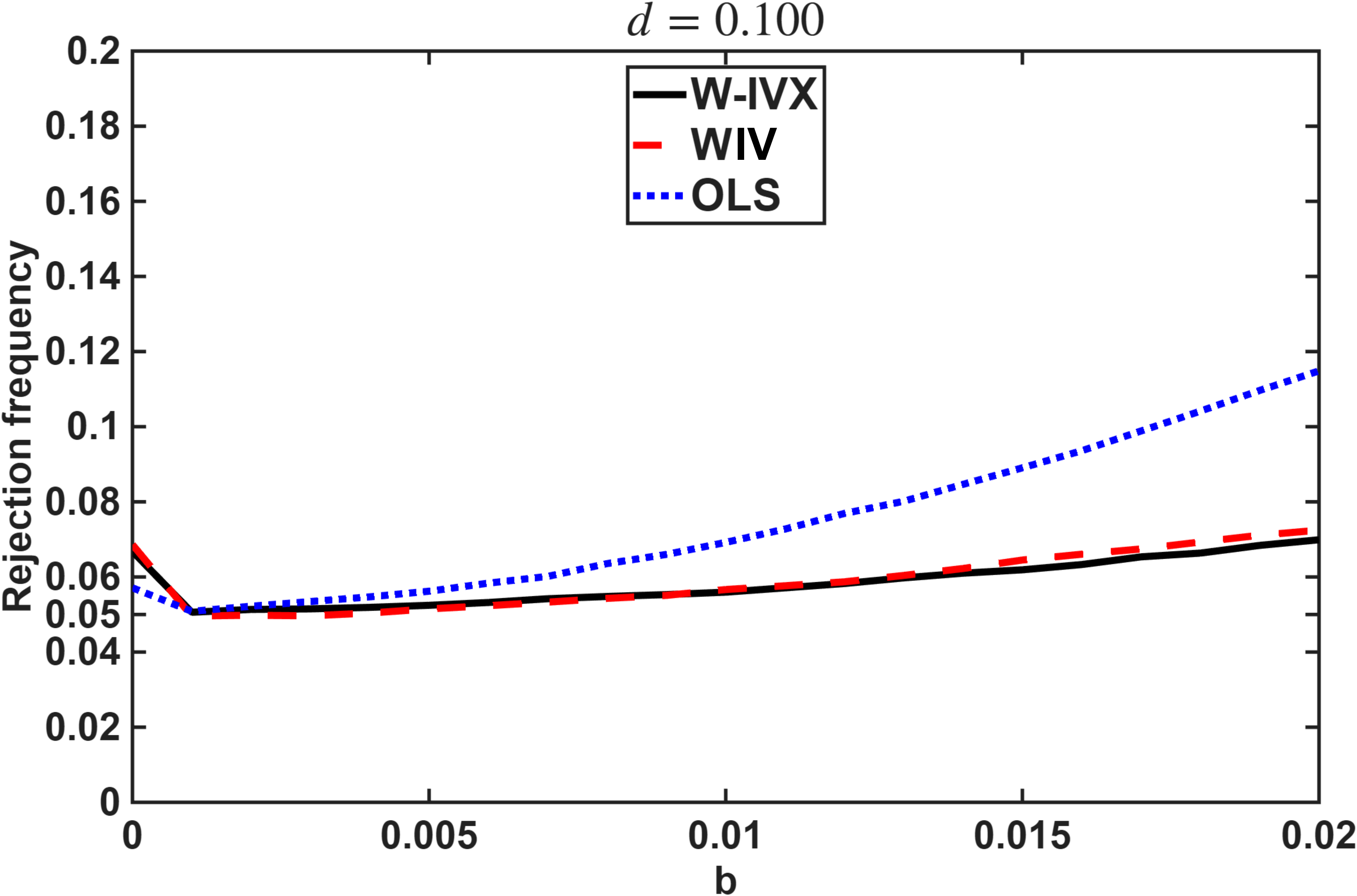}
\end{subfigure} \hfill{}
\begin{subfigure}[t]{0.48\linewidth} \centering \caption{$\alpha=1.8$, $d_{*}=0.444$}
\label{fig:sim-sap-linear-w00-alpha-1p80} \includegraphics[width=0.49\linewidth]{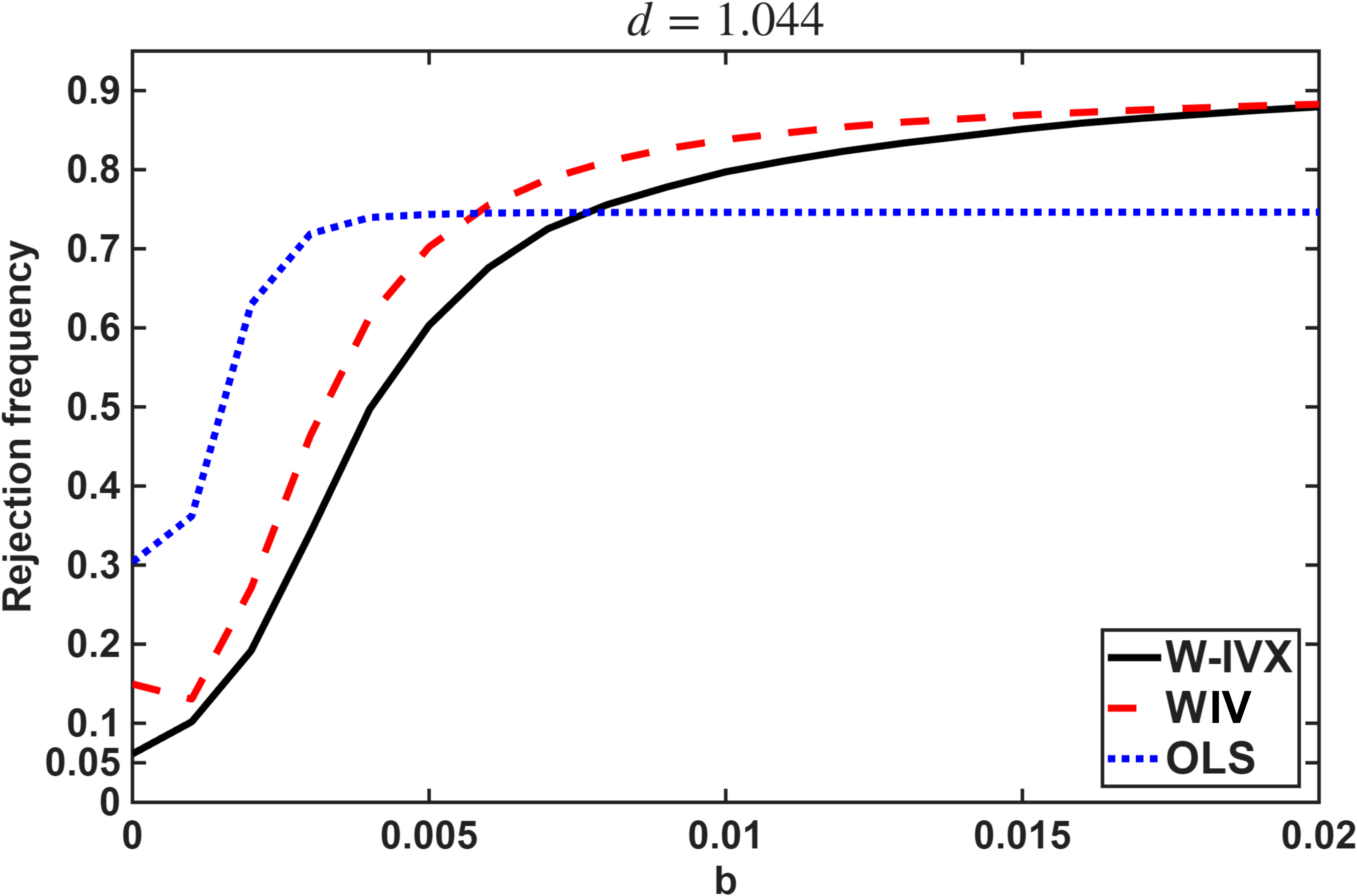}
\includegraphics[width=0.49\linewidth]{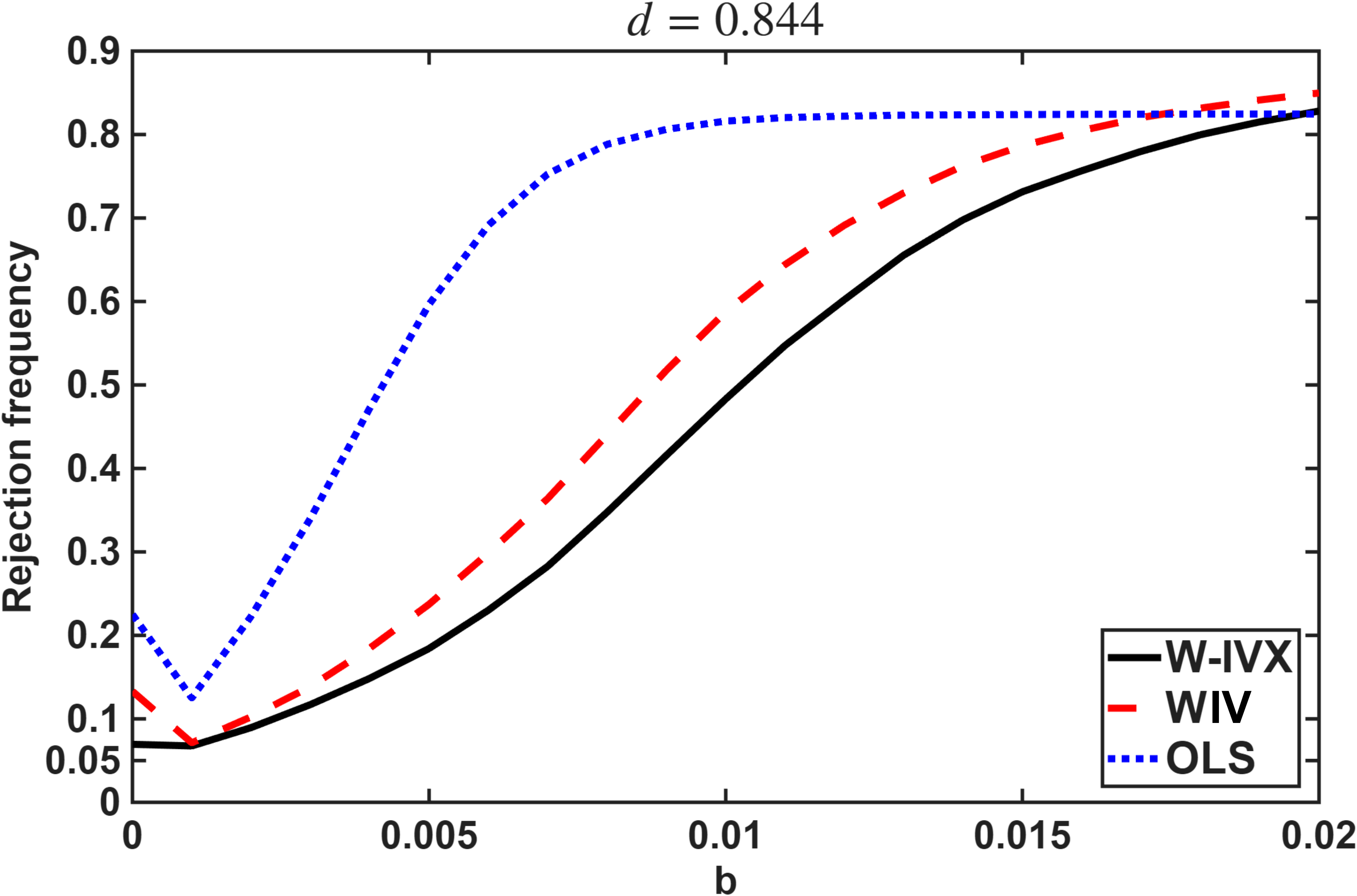}\\[0.6ex]
\includegraphics[width=0.49\linewidth]{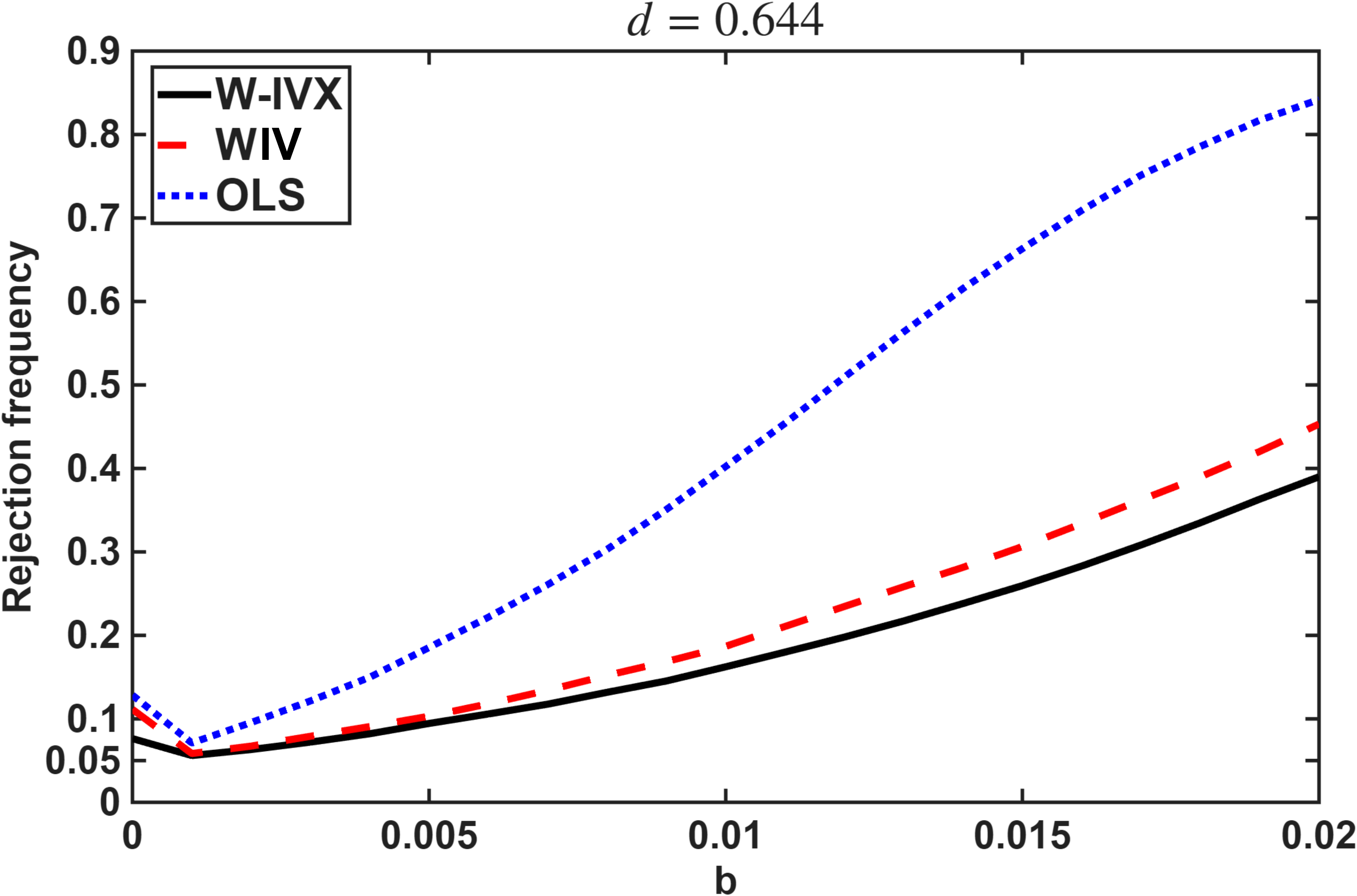}
\includegraphics[width=0.49\linewidth]{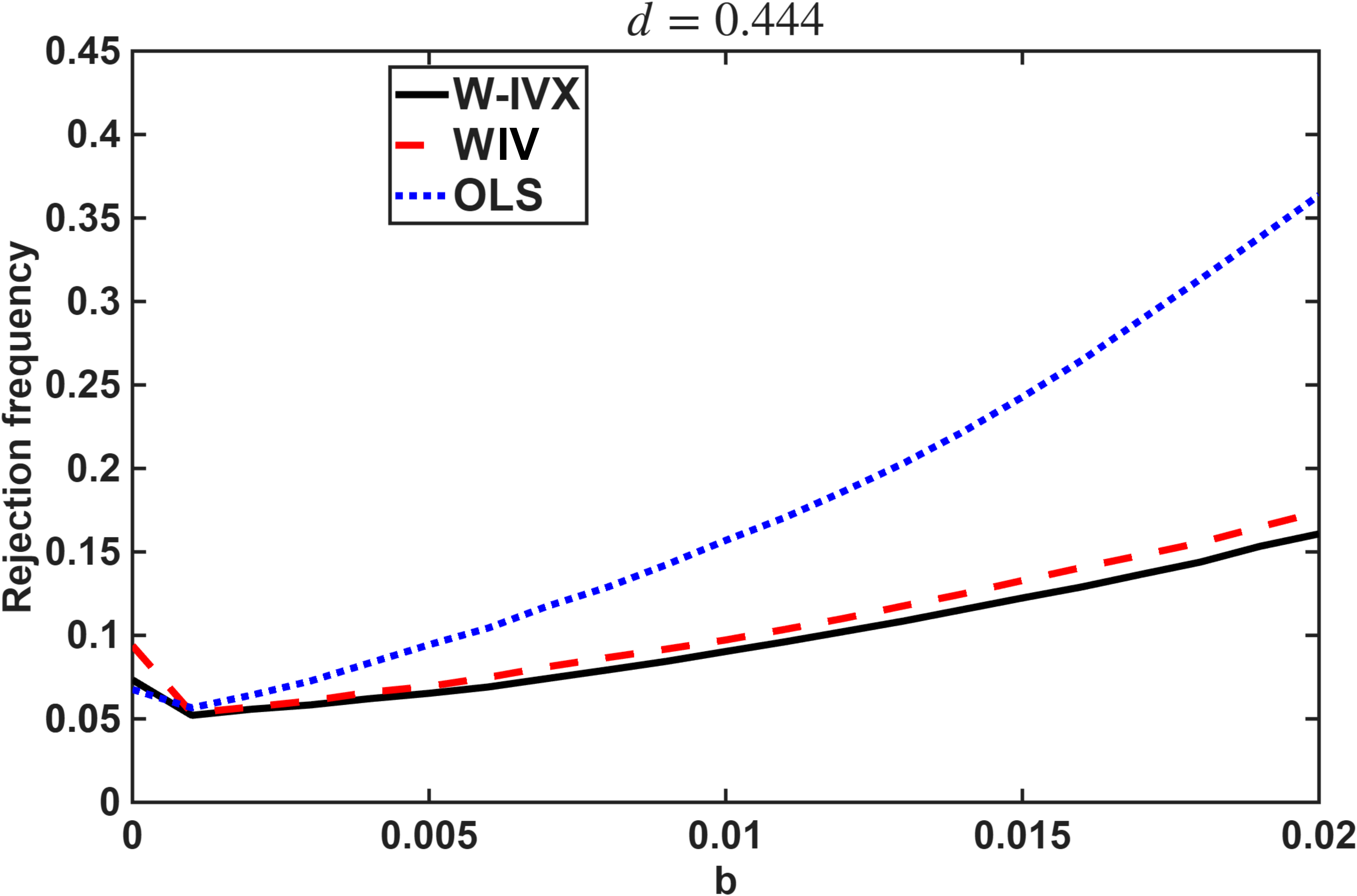}\\[0.6ex]
\includegraphics[width=0.49\linewidth]{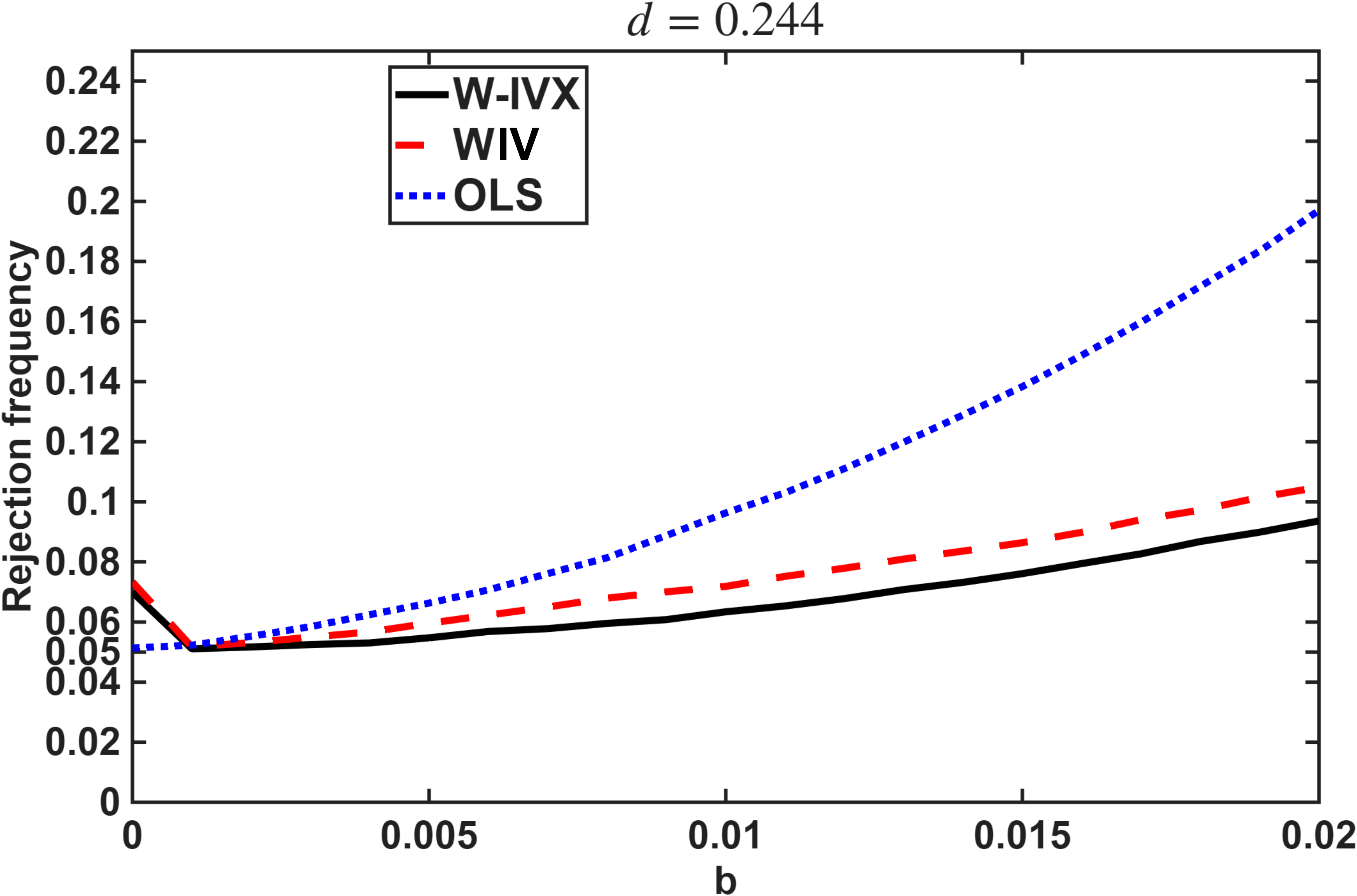}
\includegraphics[width=0.49\linewidth]{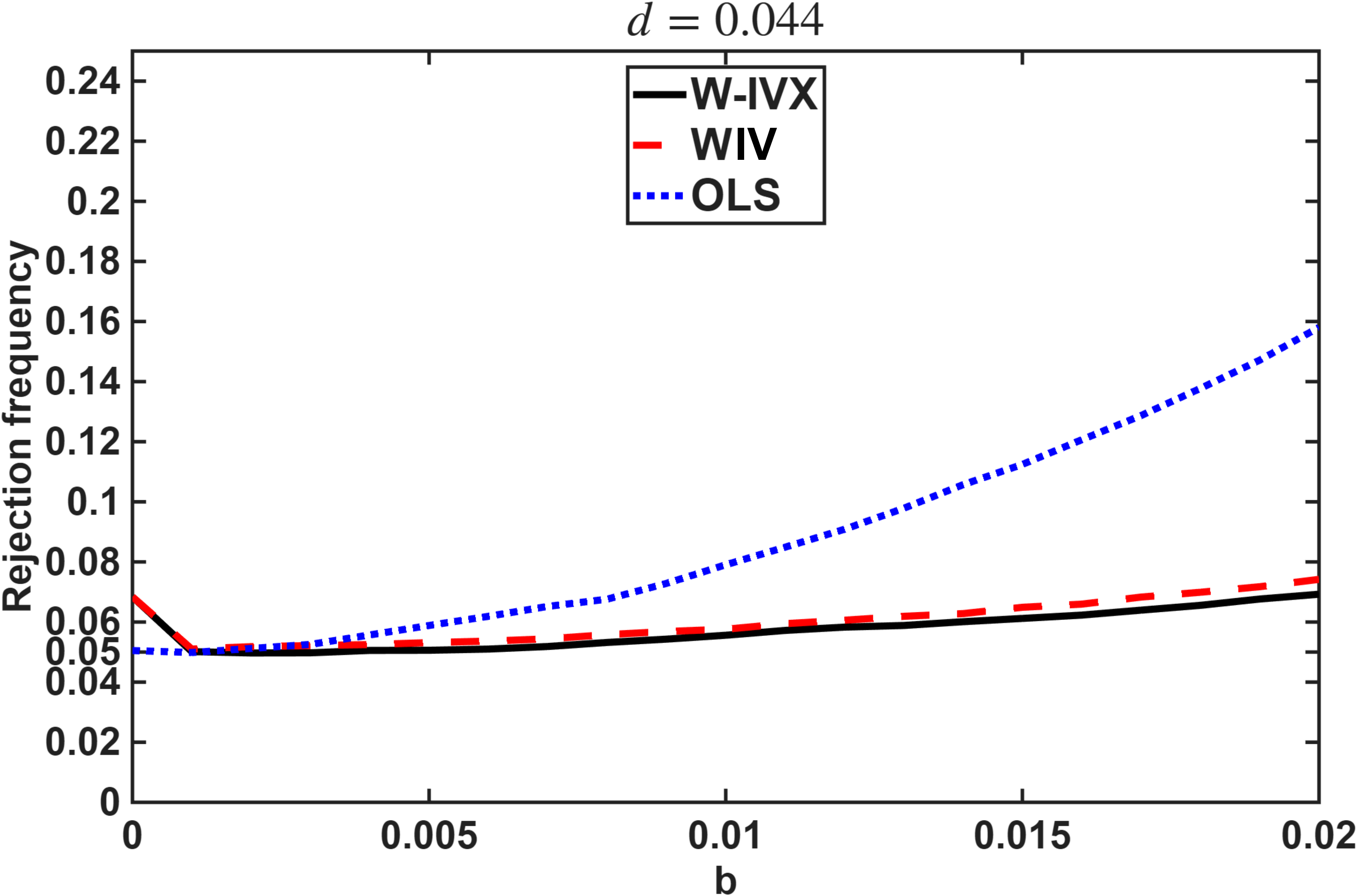}

\end{subfigure}

\end{figure}

\begin{figure}[p]
\ContinuedFloat \centering \captionsetup[subfigure]{position=top}

\setcounter{subfigure}{2}

\caption[]{Size Adjusted Power (nominal size 5\%): $[f(x),w(x)]=[x,x^{\langle0\rangle}]$,
$\delta=-0.95$, $n=1000$ (continued)}
\label{fig:sim-sap-linear-w00-continued}

\begin{subfigure}[t]{0.48\linewidth} \centering \caption{$\alpha=1.6$, $d_{*}=0.375$}
\label{fig:sim-sap-linear-w00-alpha-1p60}\includegraphics[width=0.49\linewidth]{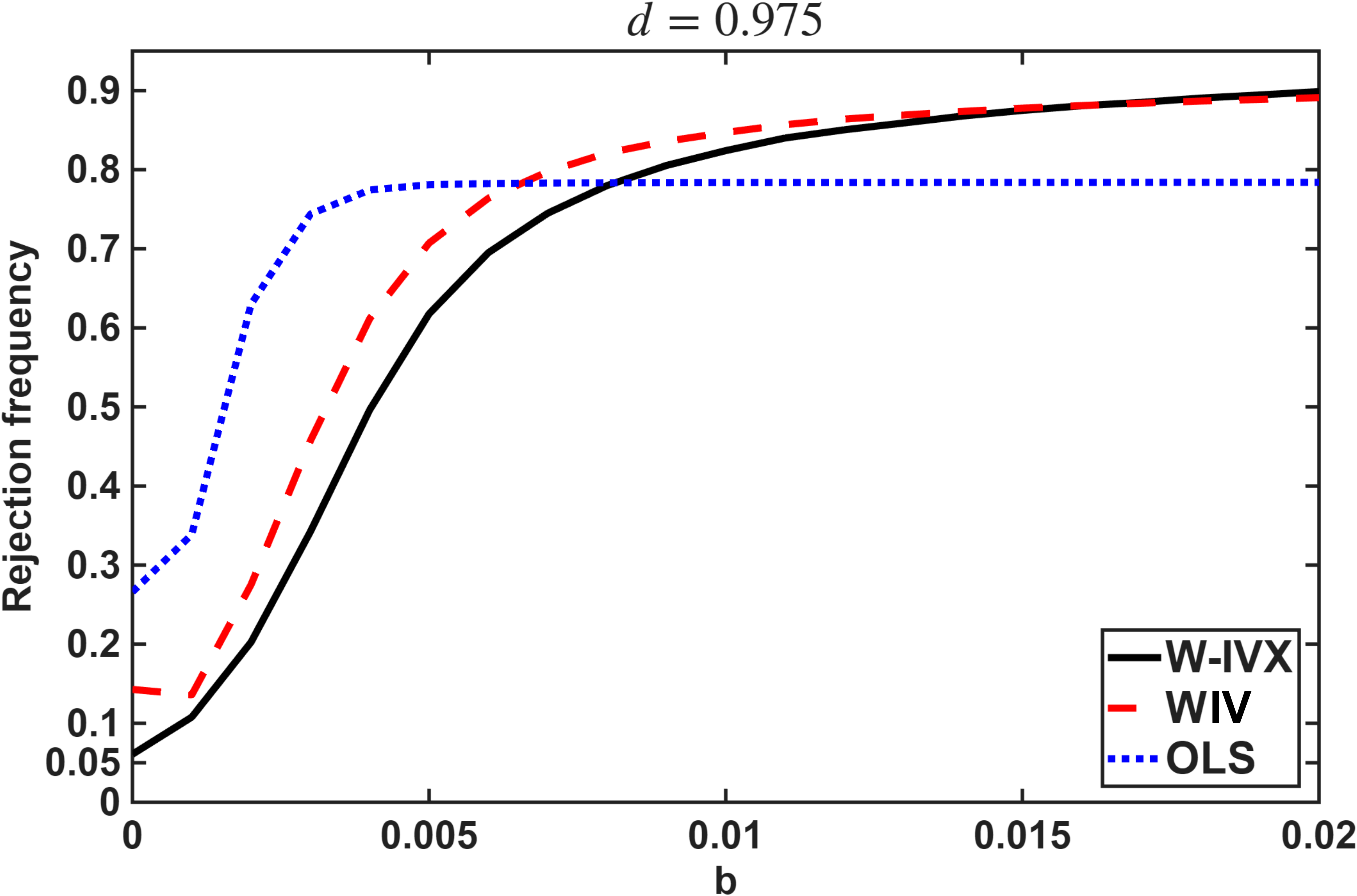}
\includegraphics[width=0.49\linewidth]{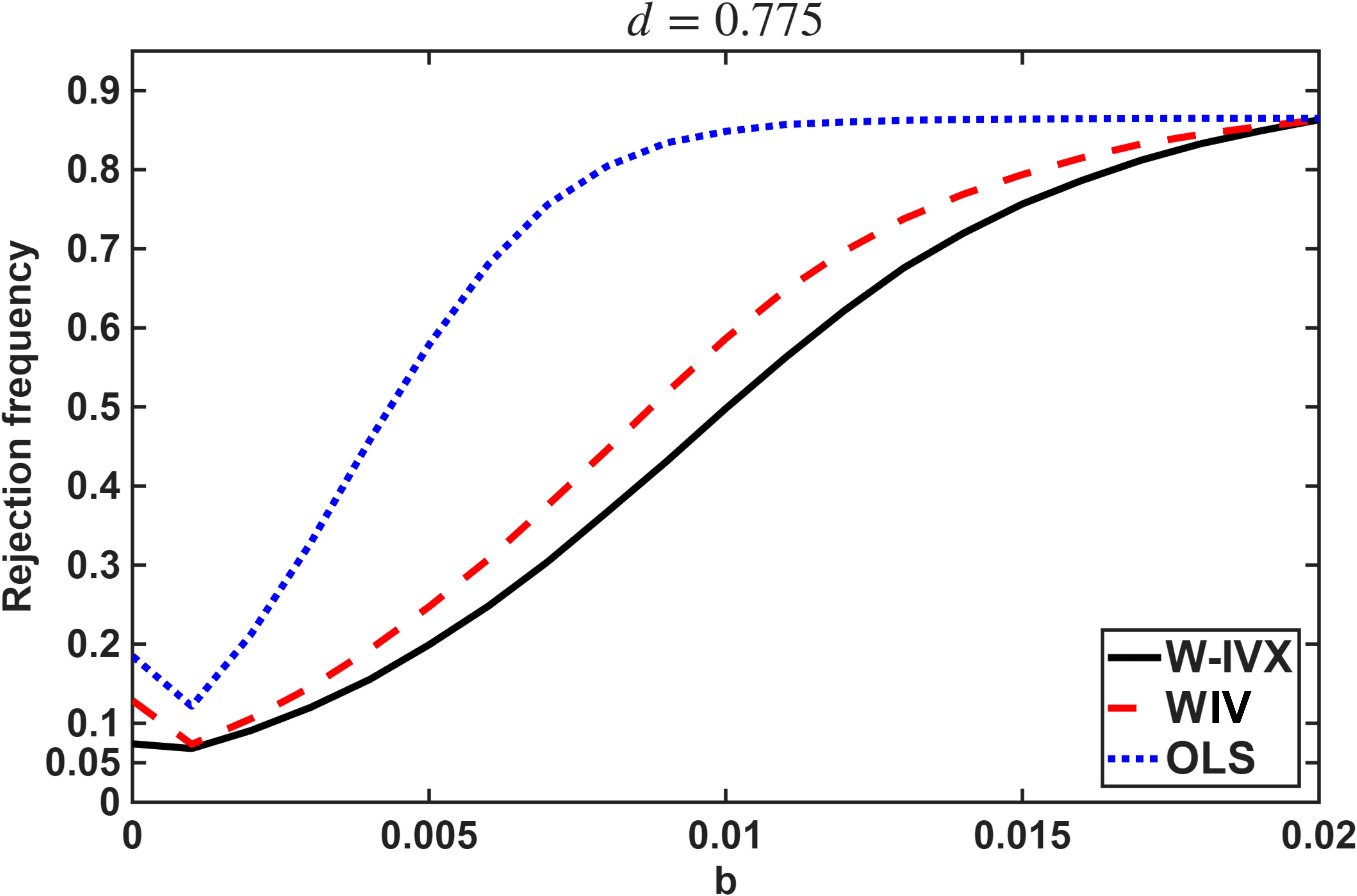}\\[0.6ex]
\includegraphics[width=0.49\linewidth]{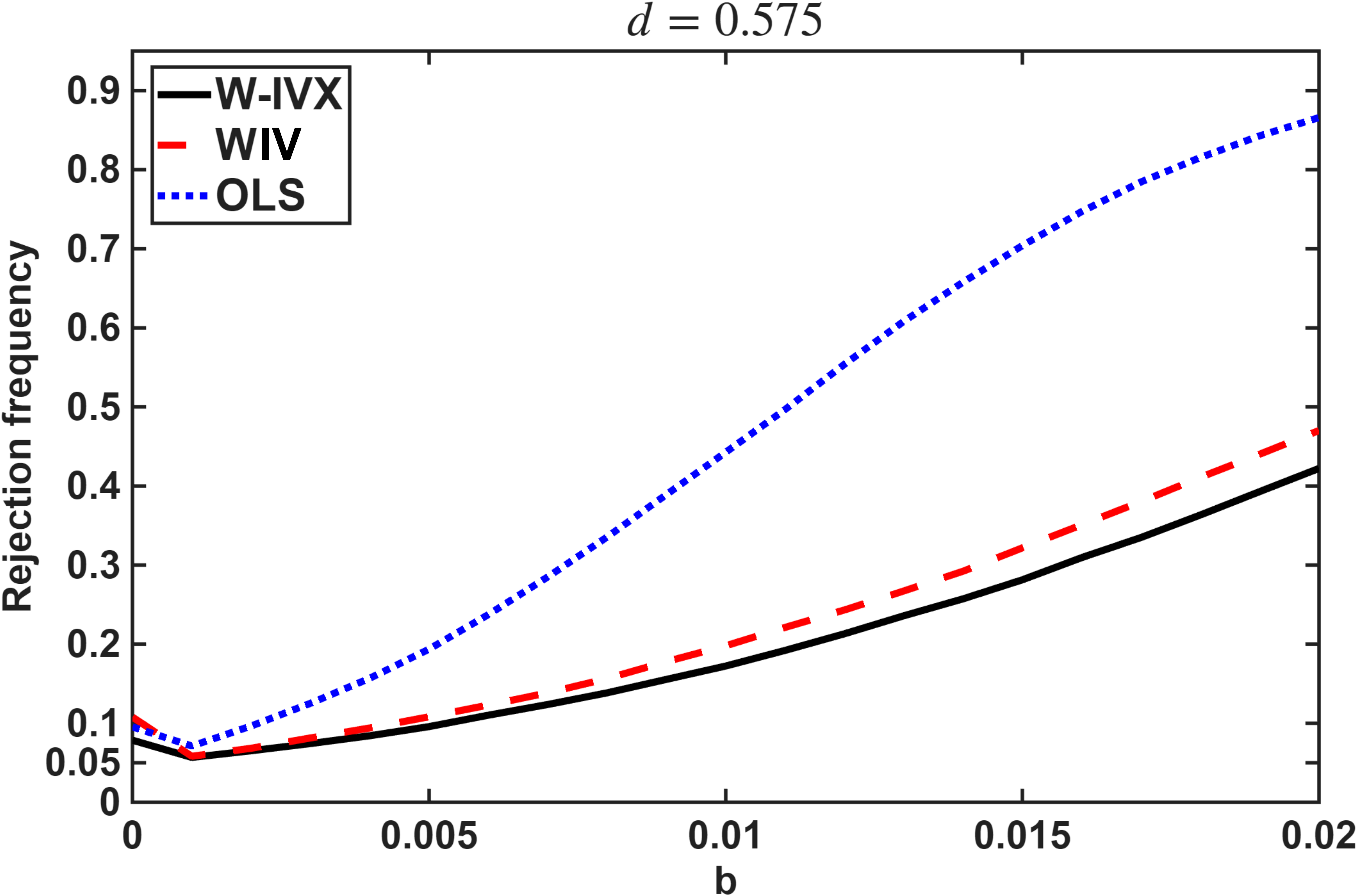}
\includegraphics[width=0.49\linewidth]{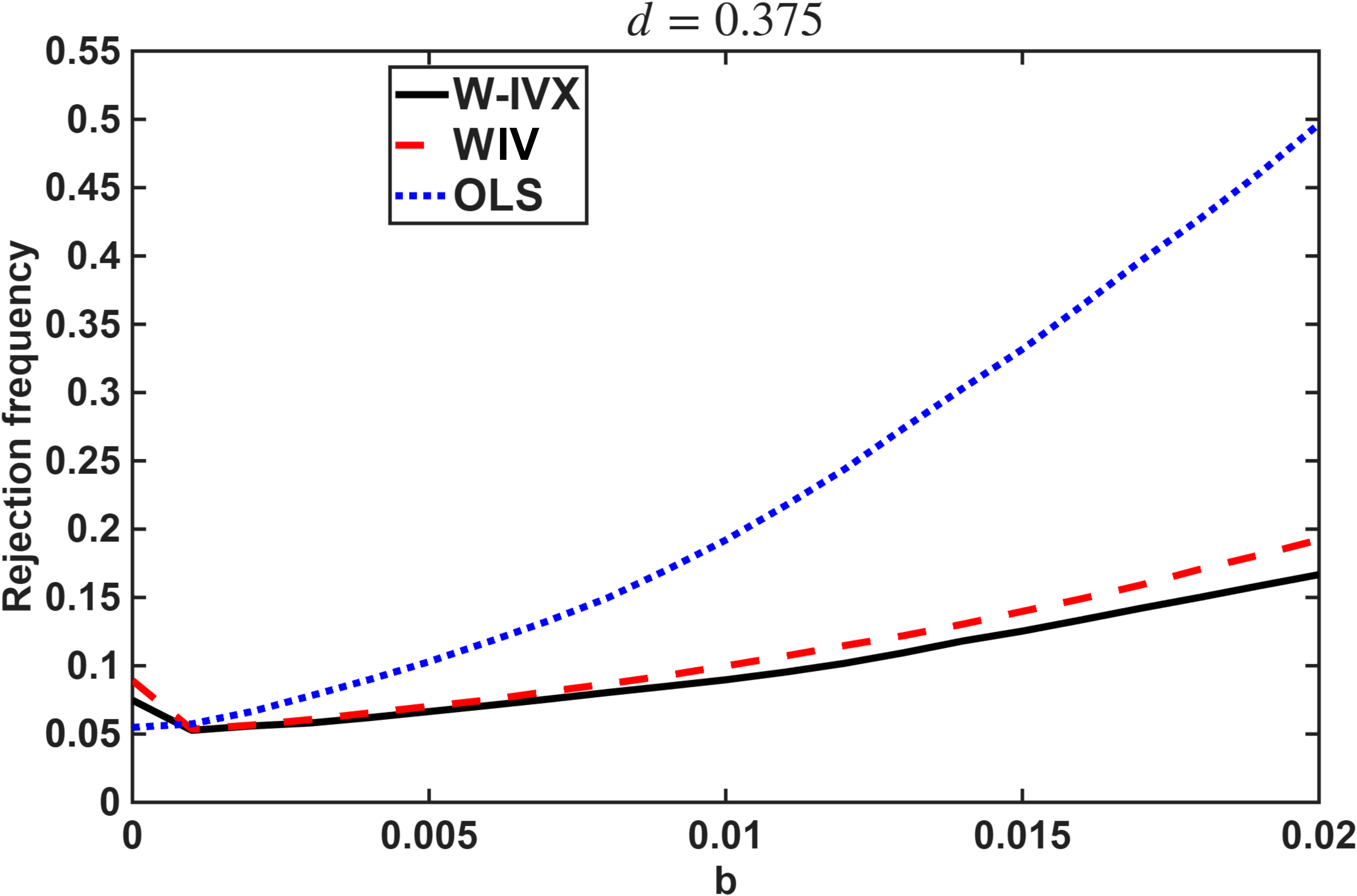}\\[0.6ex]
\includegraphics[width=0.49\linewidth]{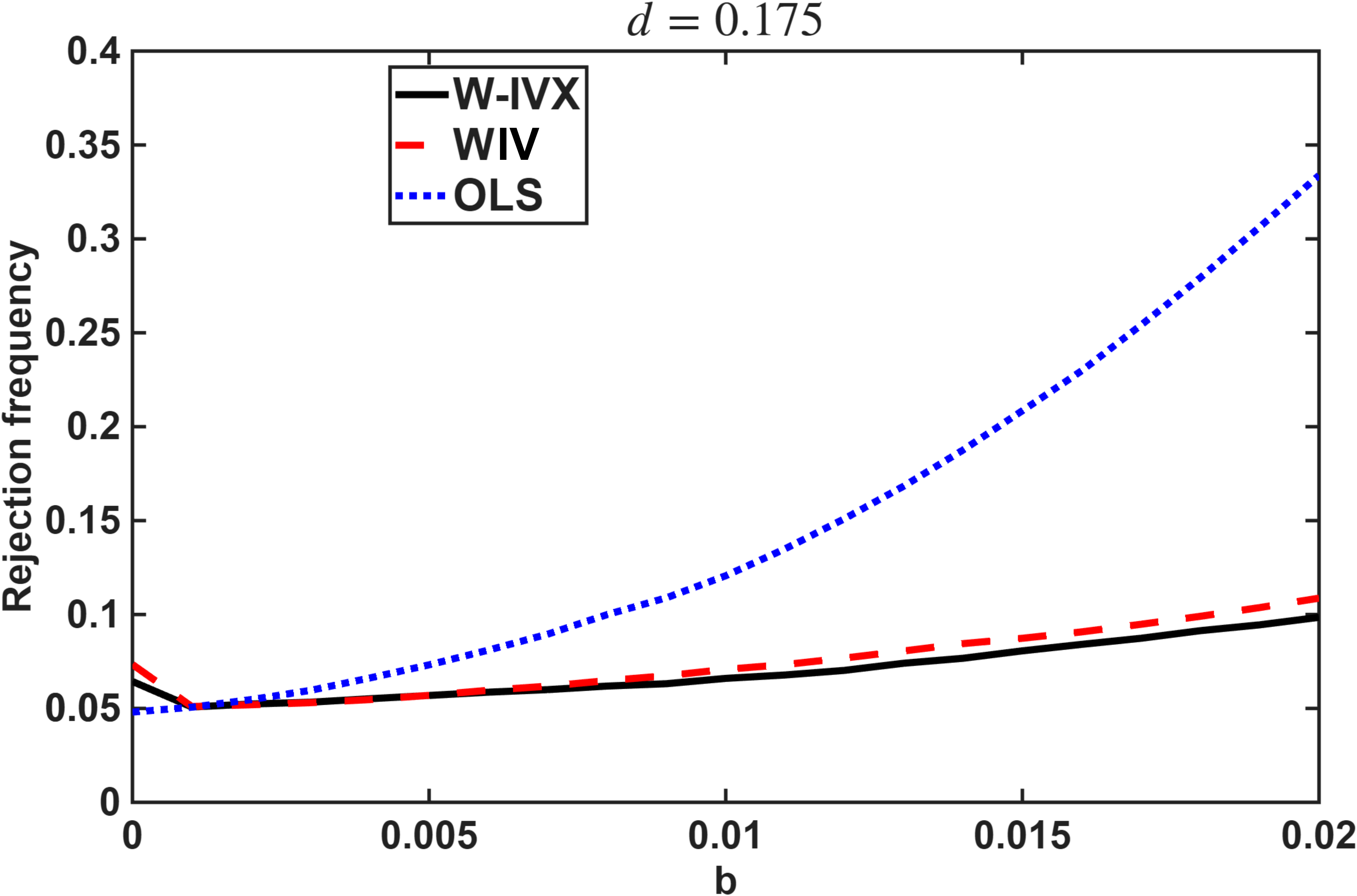}
\includegraphics[width=0.49\linewidth]{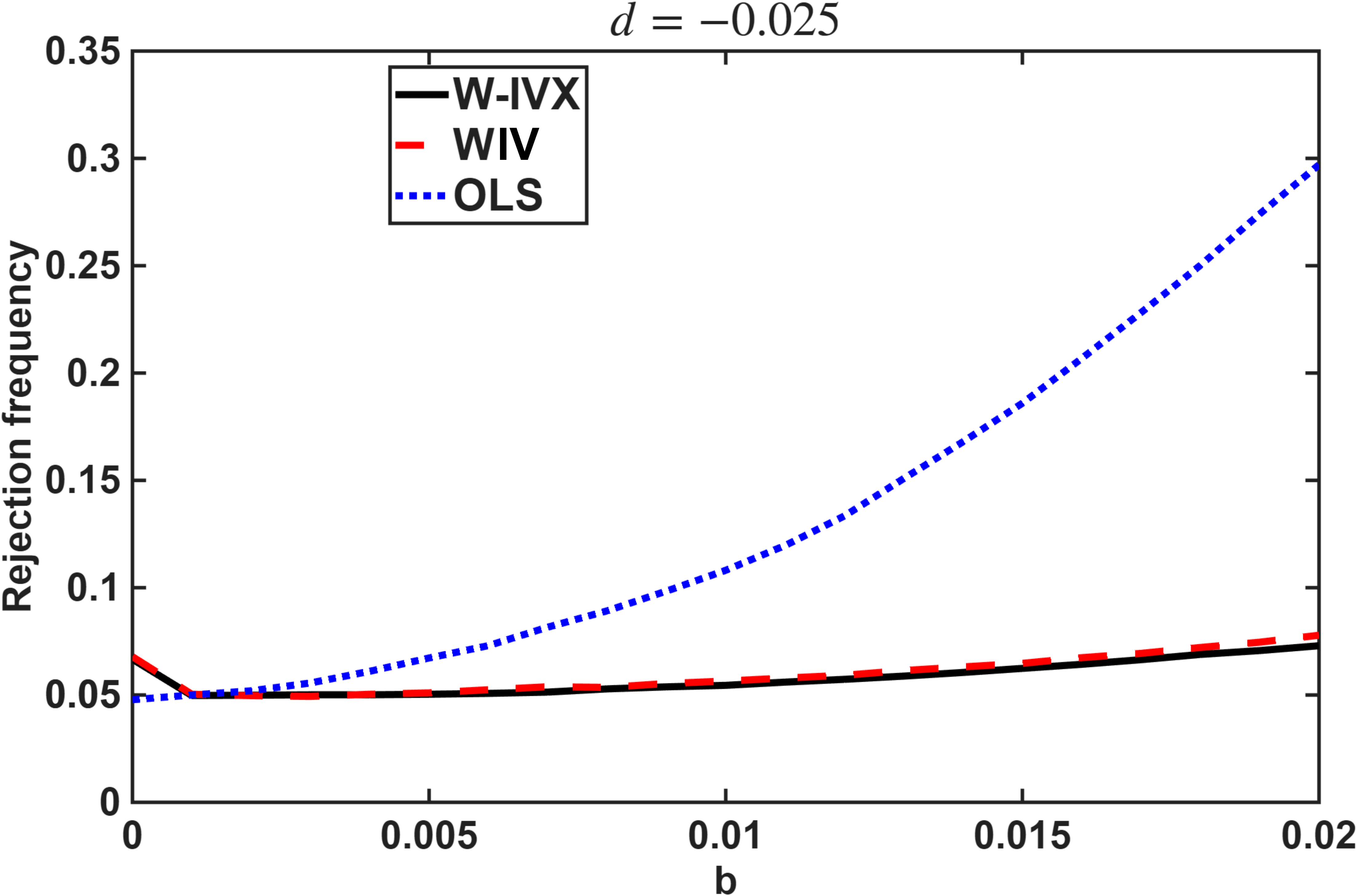}
\end{subfigure} \hfill{}
\begin{subfigure}[t]{0.48\linewidth} \centering \caption{$\alpha=1.4$, $d_{*}=0.286$}
\label{fig:sim-sap-linear-w00-alpha-1p40}\includegraphics[width=0.49\linewidth]{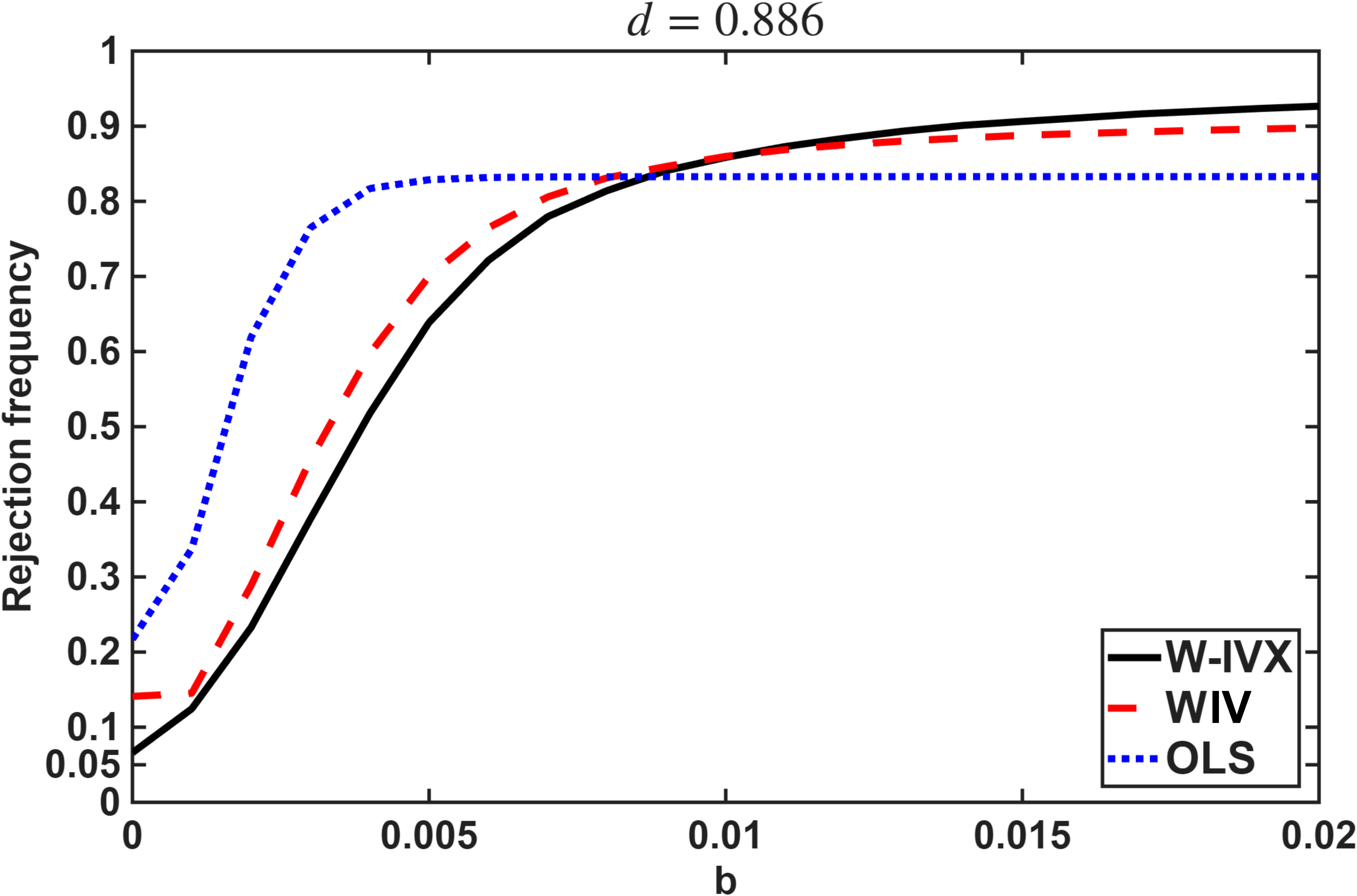}
\includegraphics[width=0.49\linewidth]{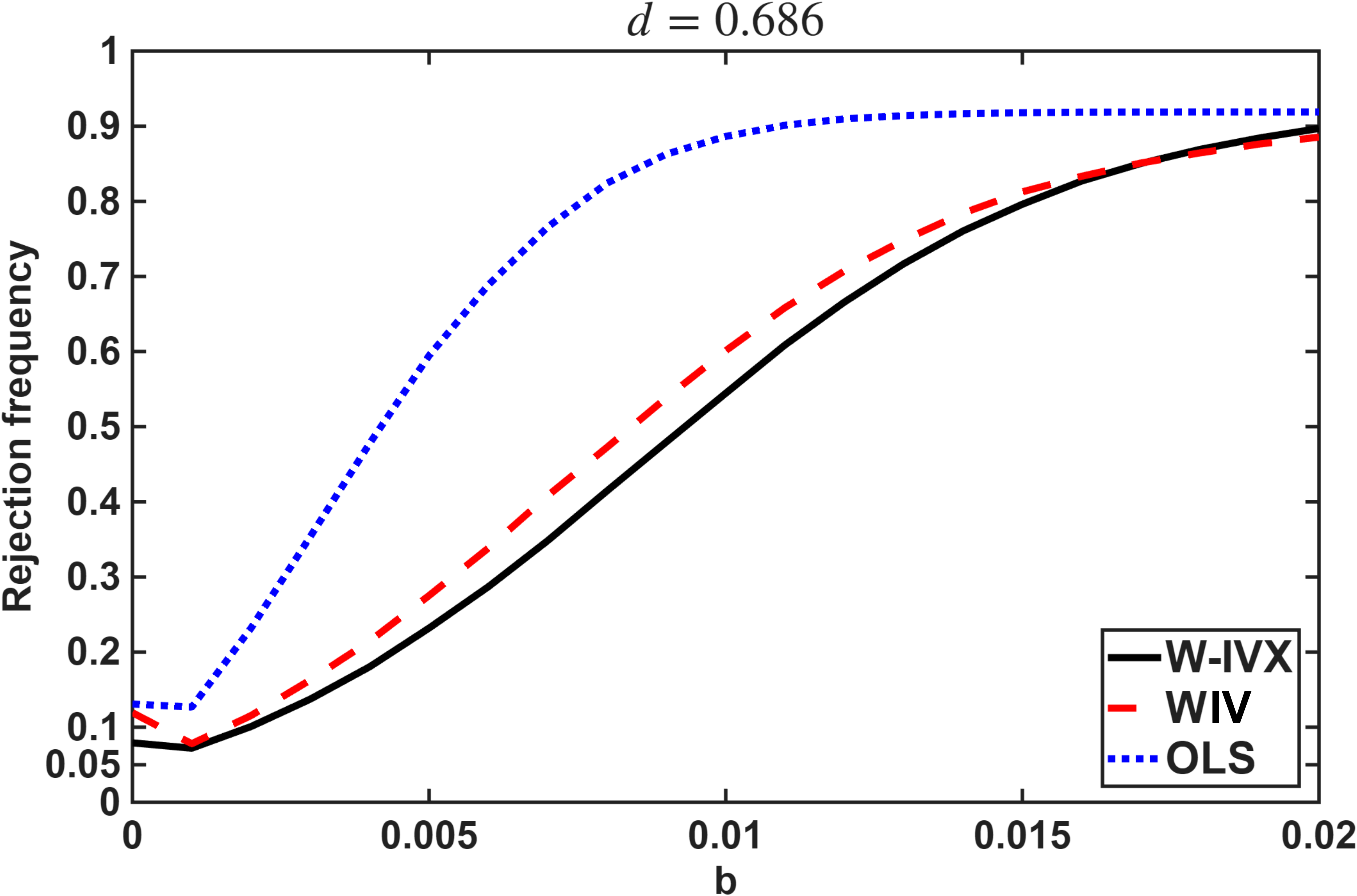}\\[0.6ex]
\includegraphics[width=0.49\linewidth]{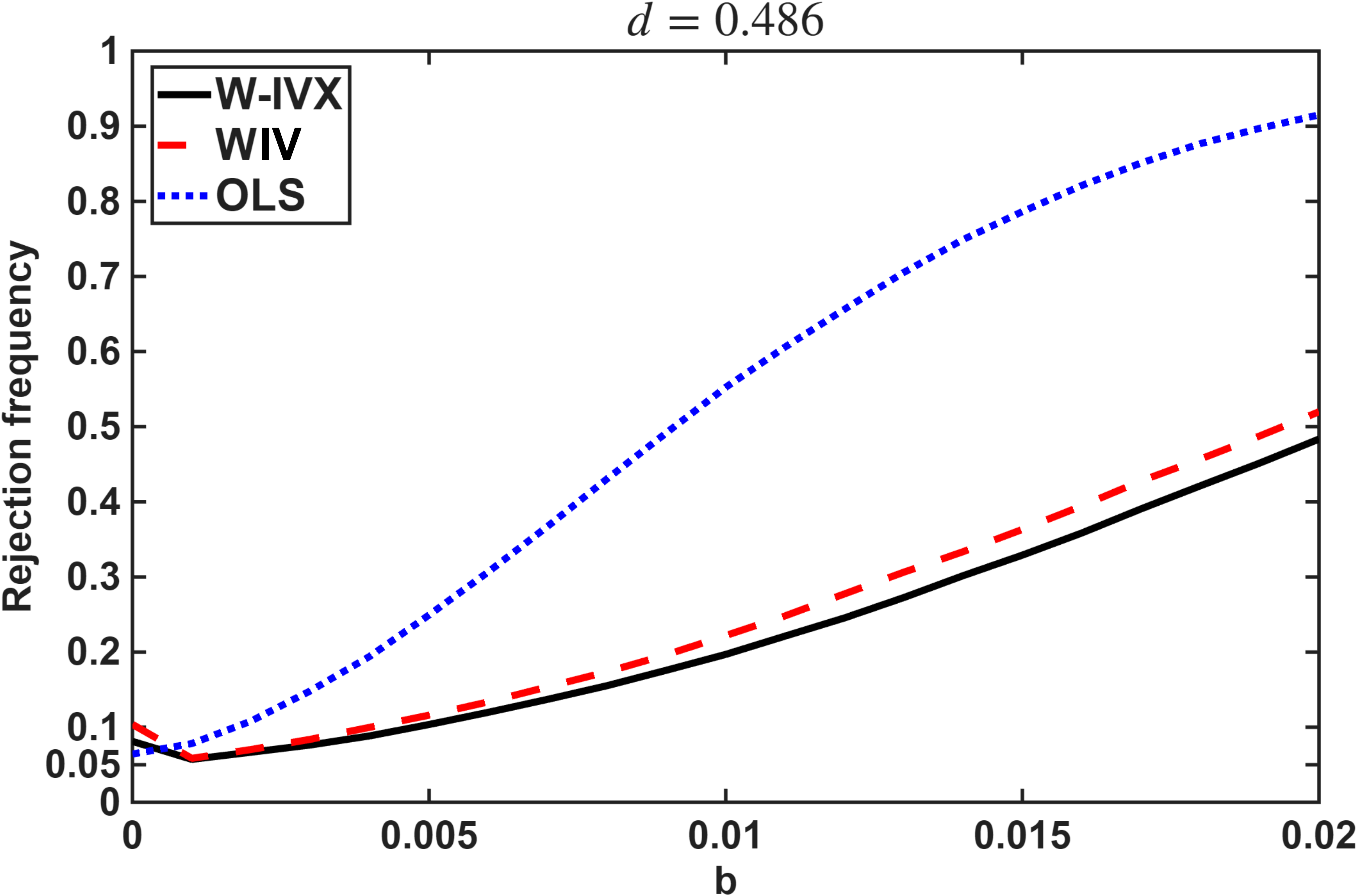}
\includegraphics[width=0.49\linewidth]{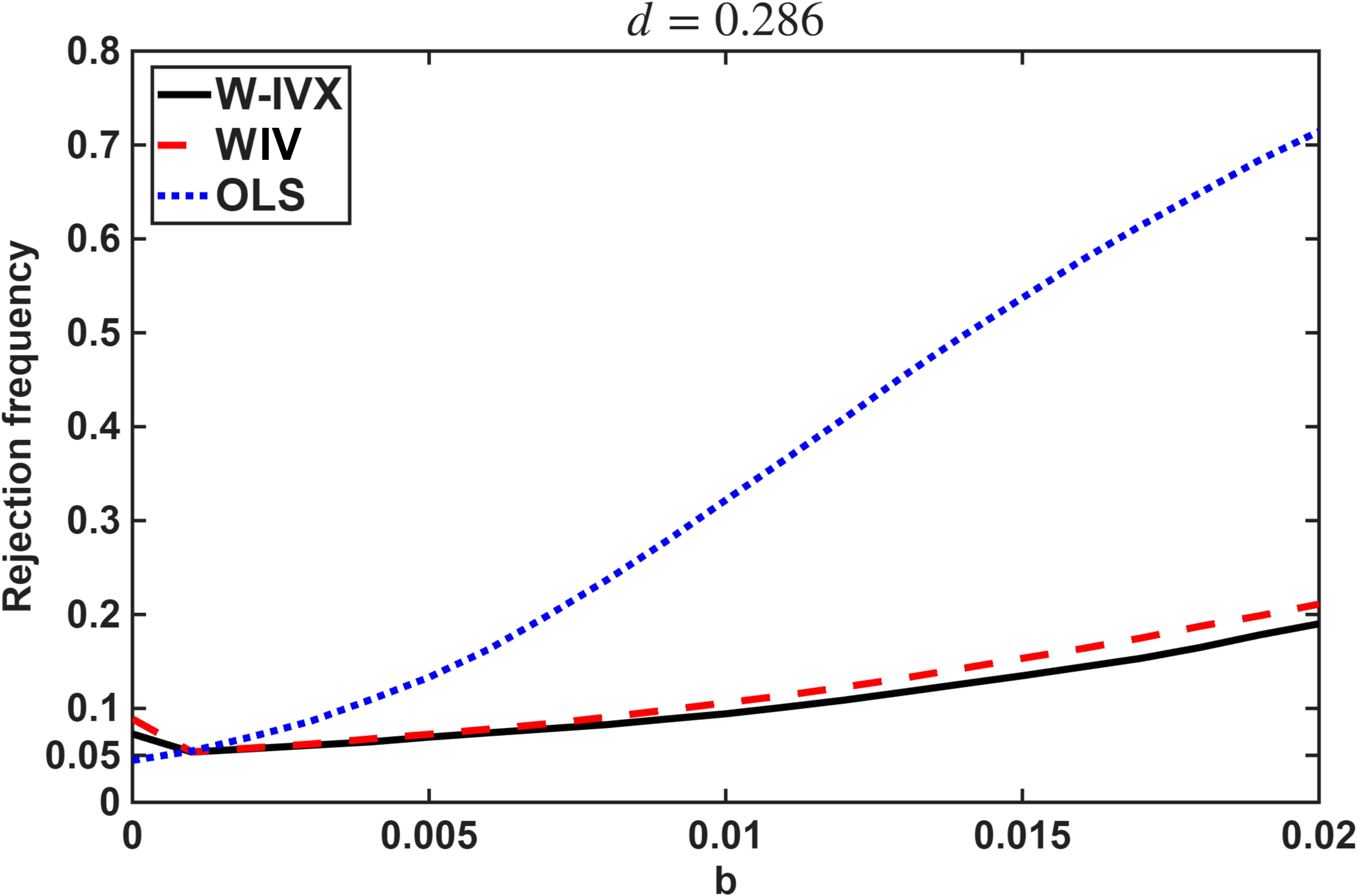}\\[0.6ex]
\includegraphics[width=0.49\linewidth]{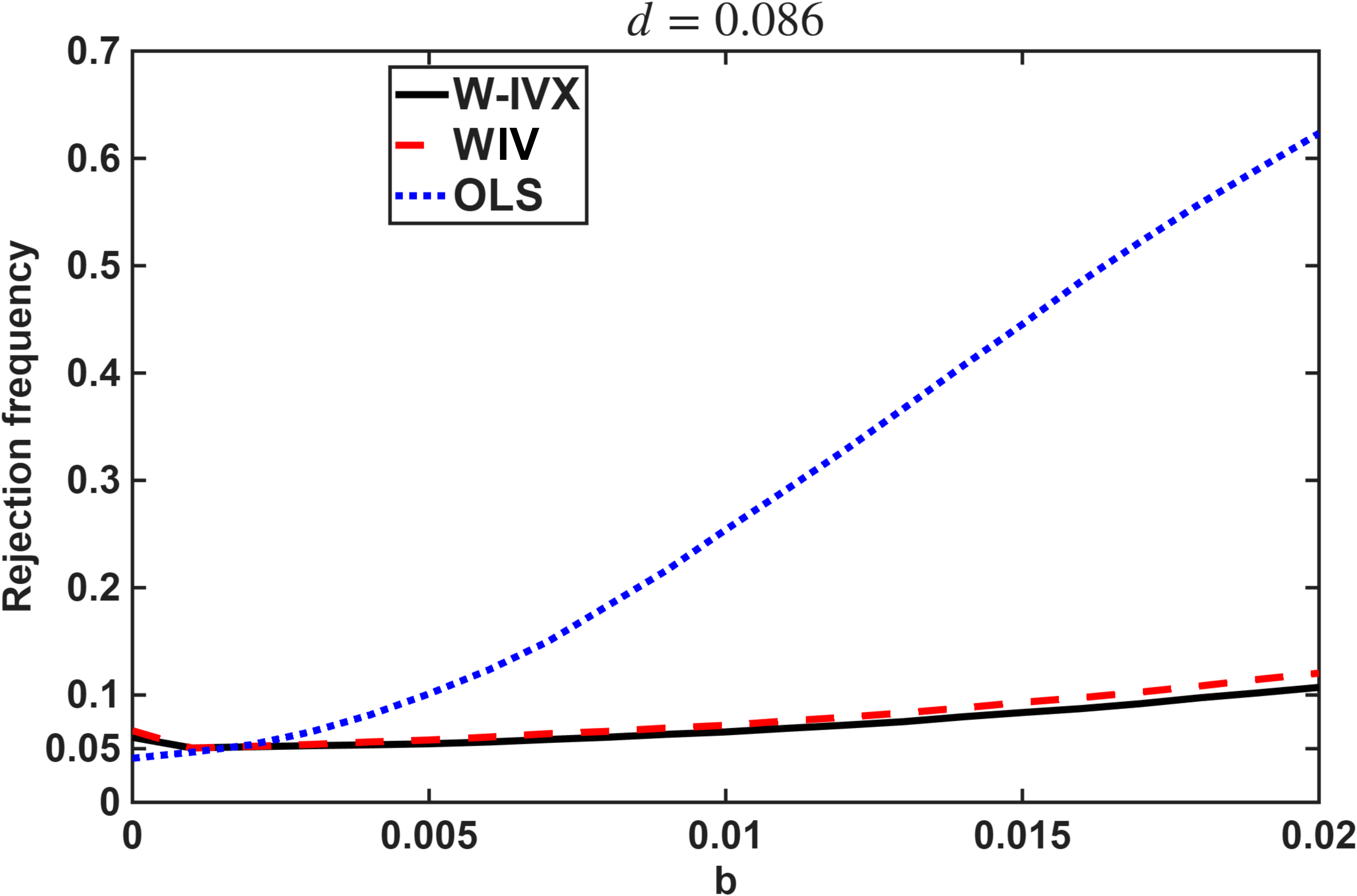}
\includegraphics[width=0.49\linewidth]{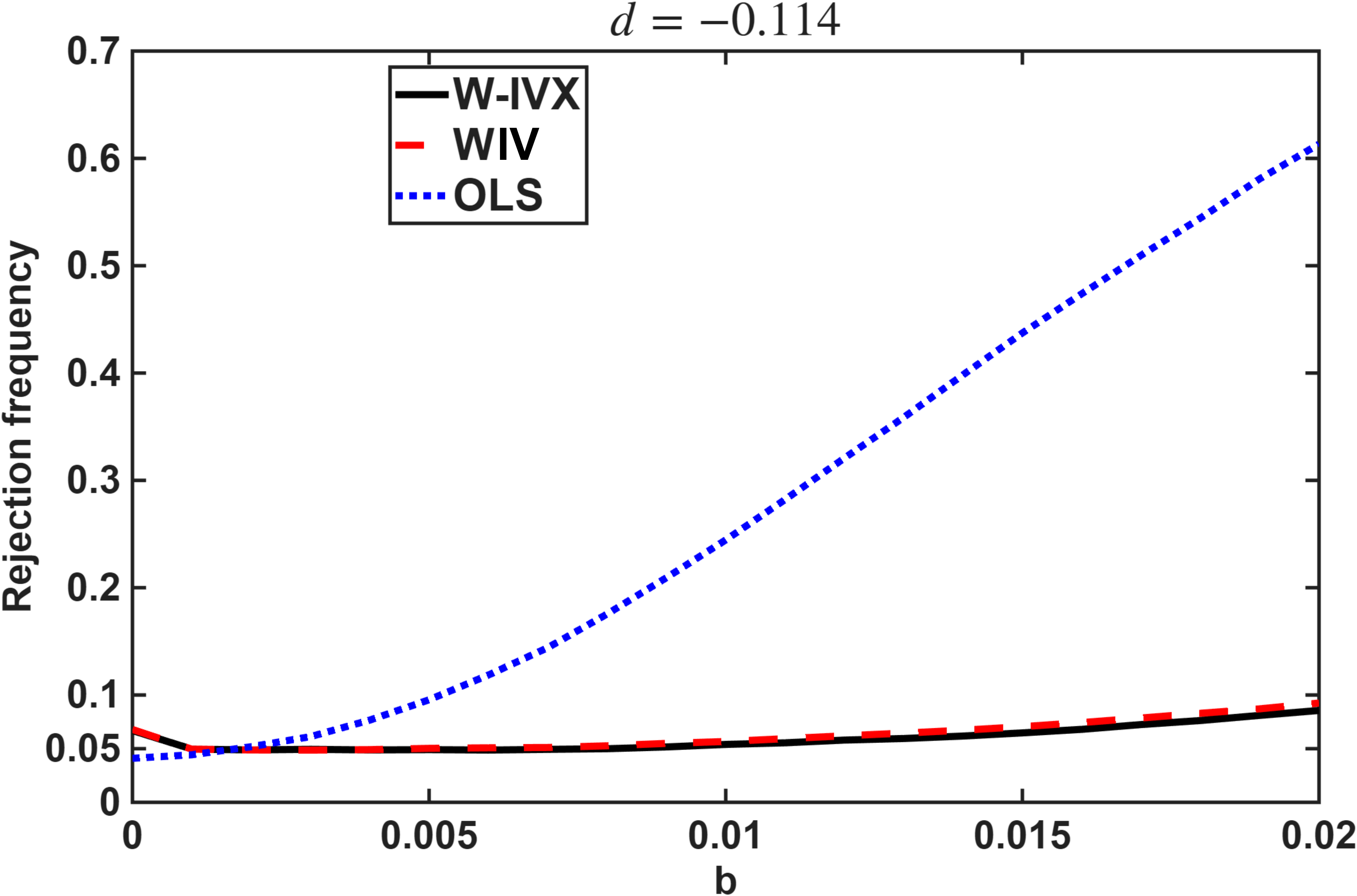}
\end{subfigure}
\end{figure}

\end{landscape} \clearpage{}

\begin{landscape}

\begin{figure}[p]
\caption{Size Adjusted Power (nominal size 5\%): \textbf{DGP1}, $[f(x),w(x)]=[x,x^{\langle0.5\rangle}]$,
$\delta=-0.95$, $n=1000$}
\label{fig:sim-sap-linear-w05}\centering \captionsetup[subfigure]{position=top}

\begin{subfigure}[t]{0.48\linewidth} \centering \caption{$\alpha=2$, $d_{*}=0.5$}
\label{fig:sim-sap-linear-w05-alpha-2p00} \includegraphics[width=0.49\linewidth]{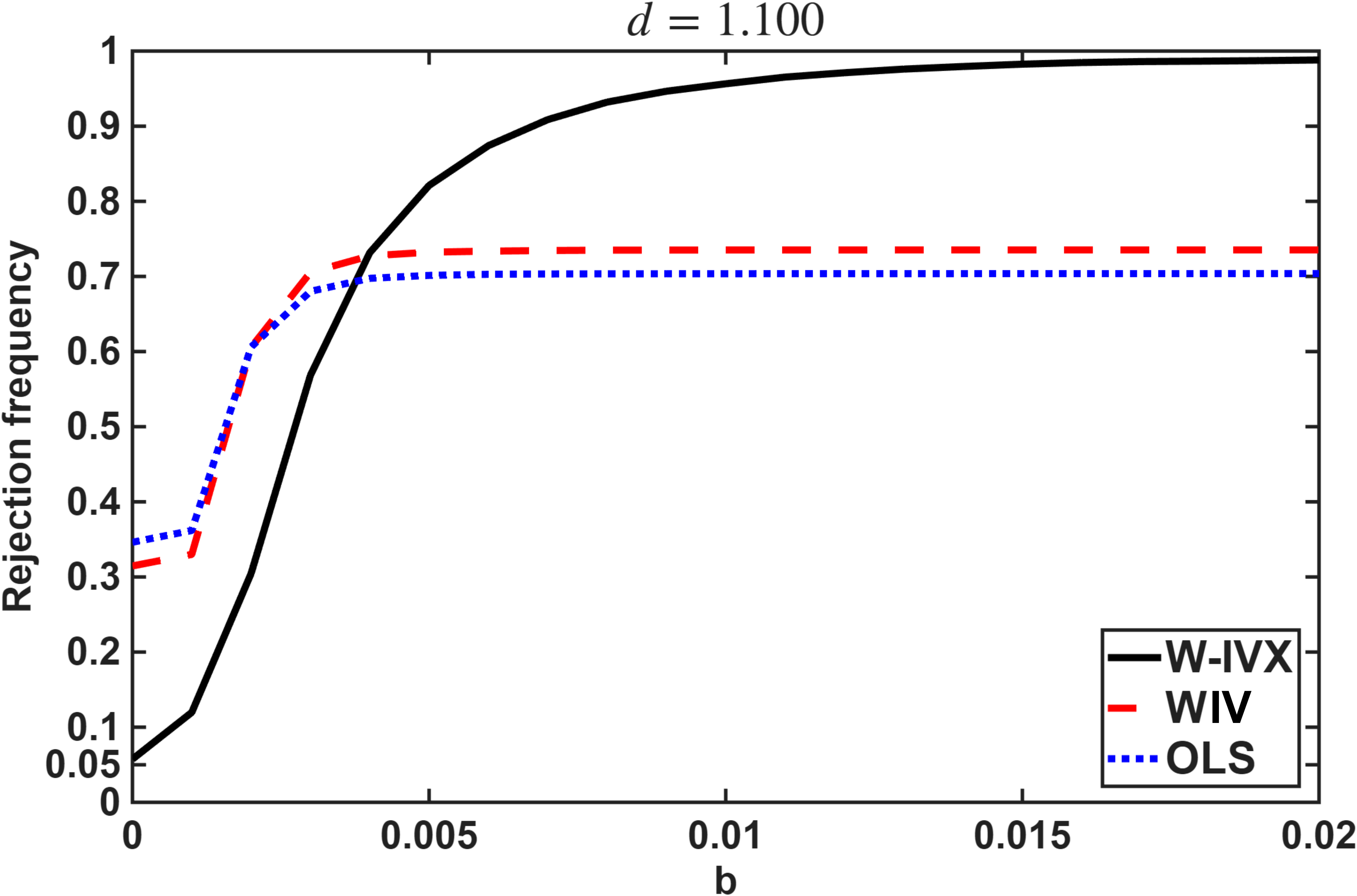}
\includegraphics[width=0.49\linewidth]{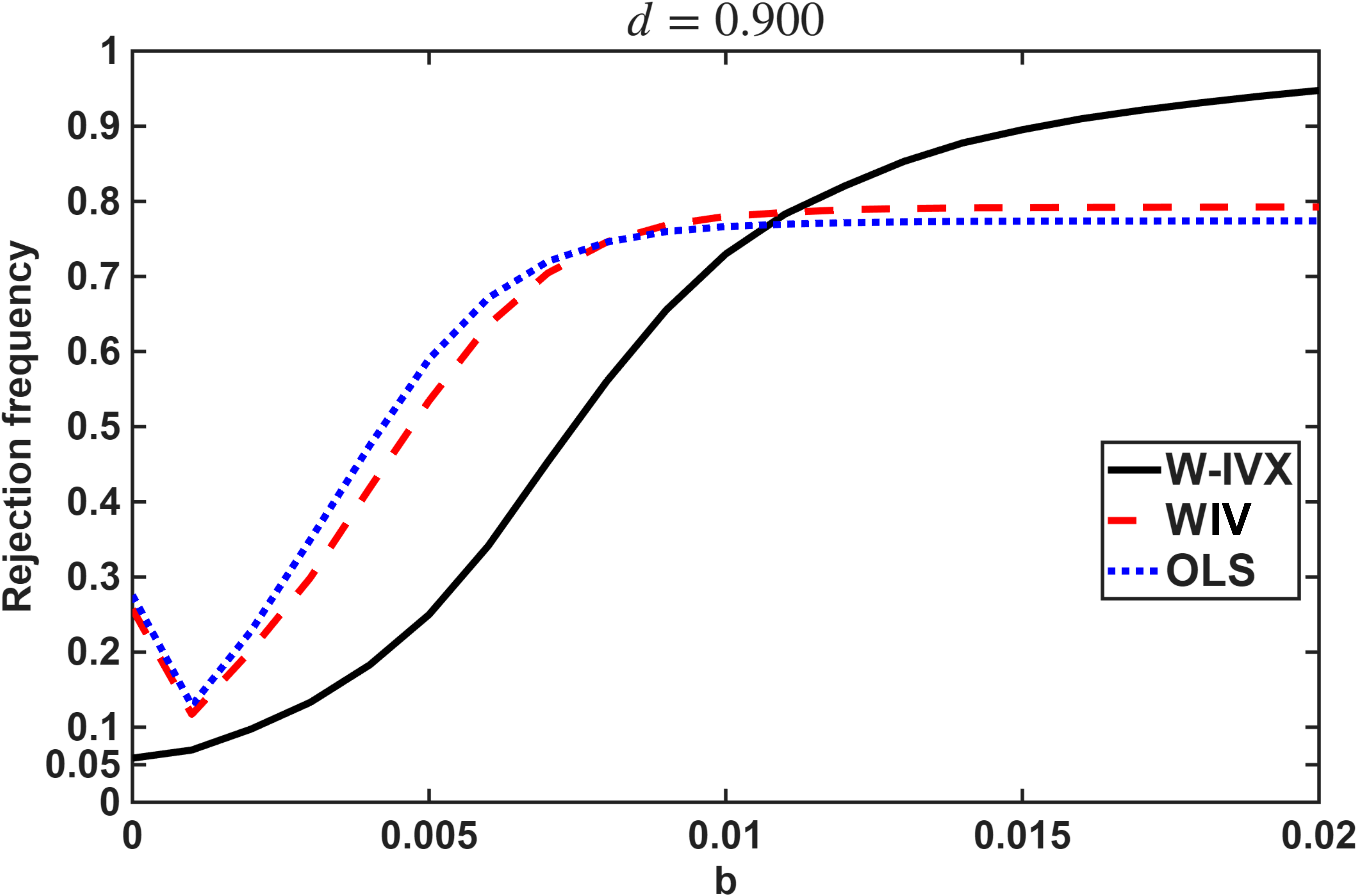}\\[0.6ex]
\includegraphics[width=0.49\linewidth]{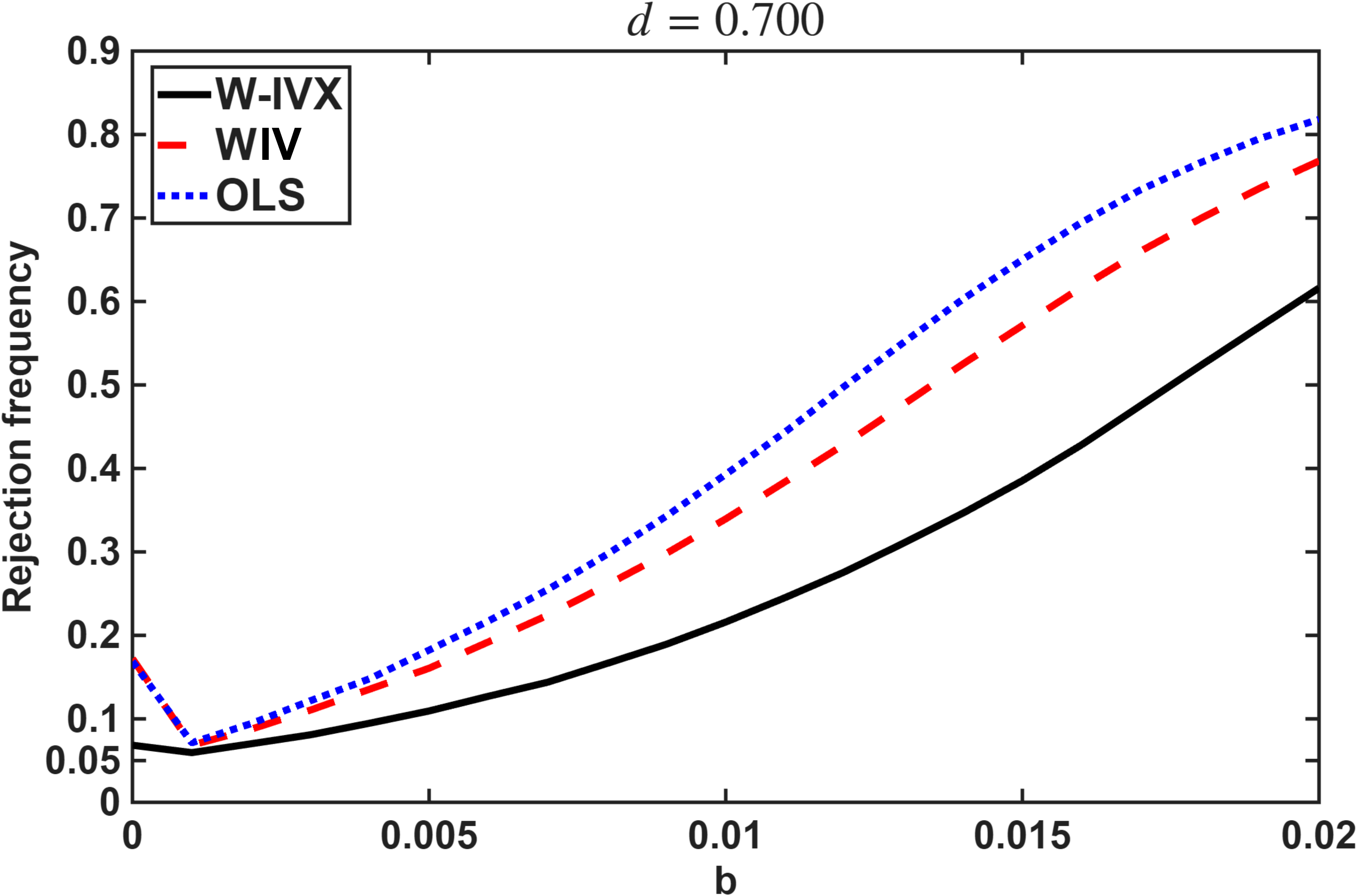}
\includegraphics[width=0.49\linewidth]{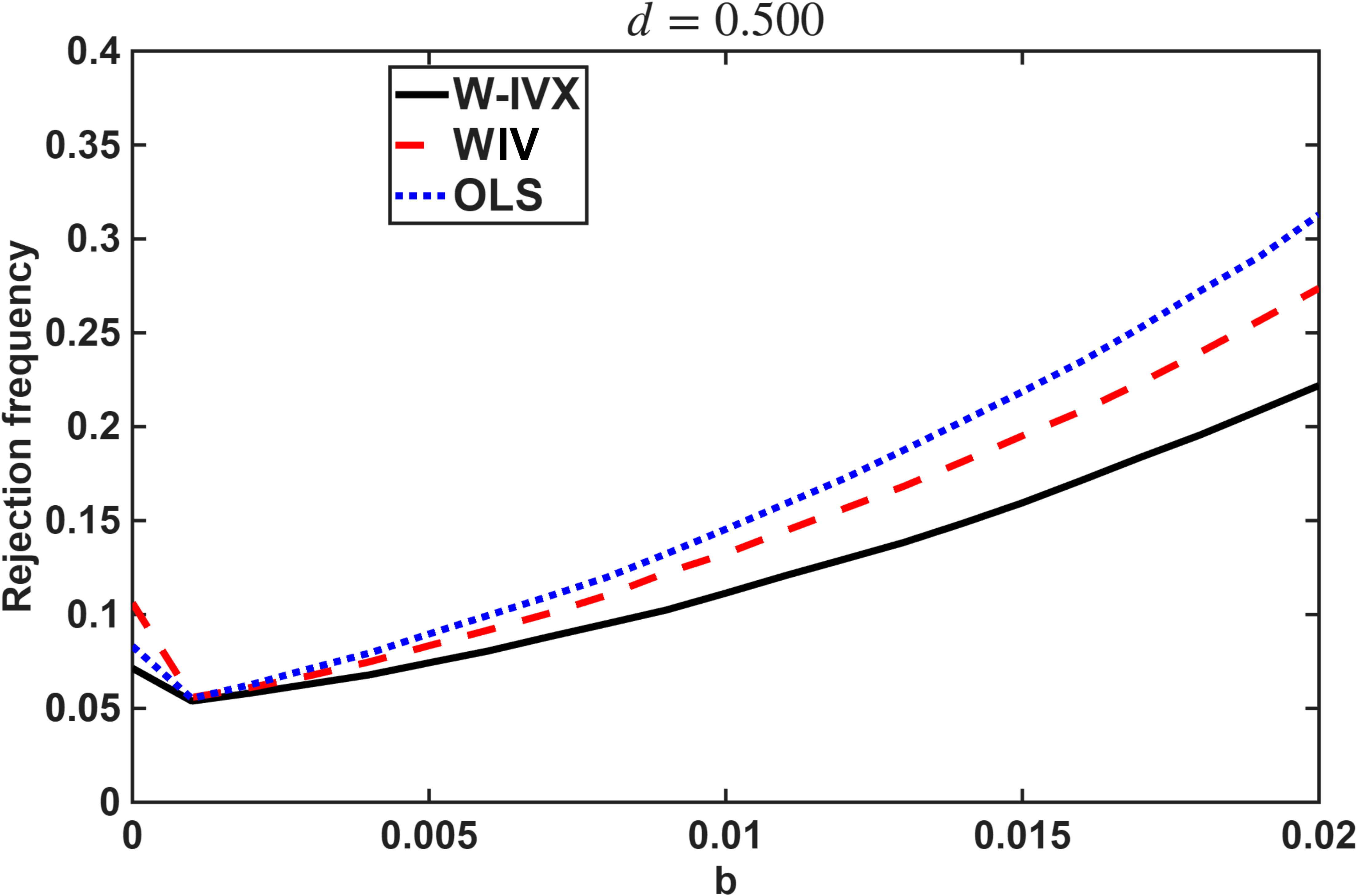}\\[0.6ex]
\includegraphics[width=0.49\linewidth]{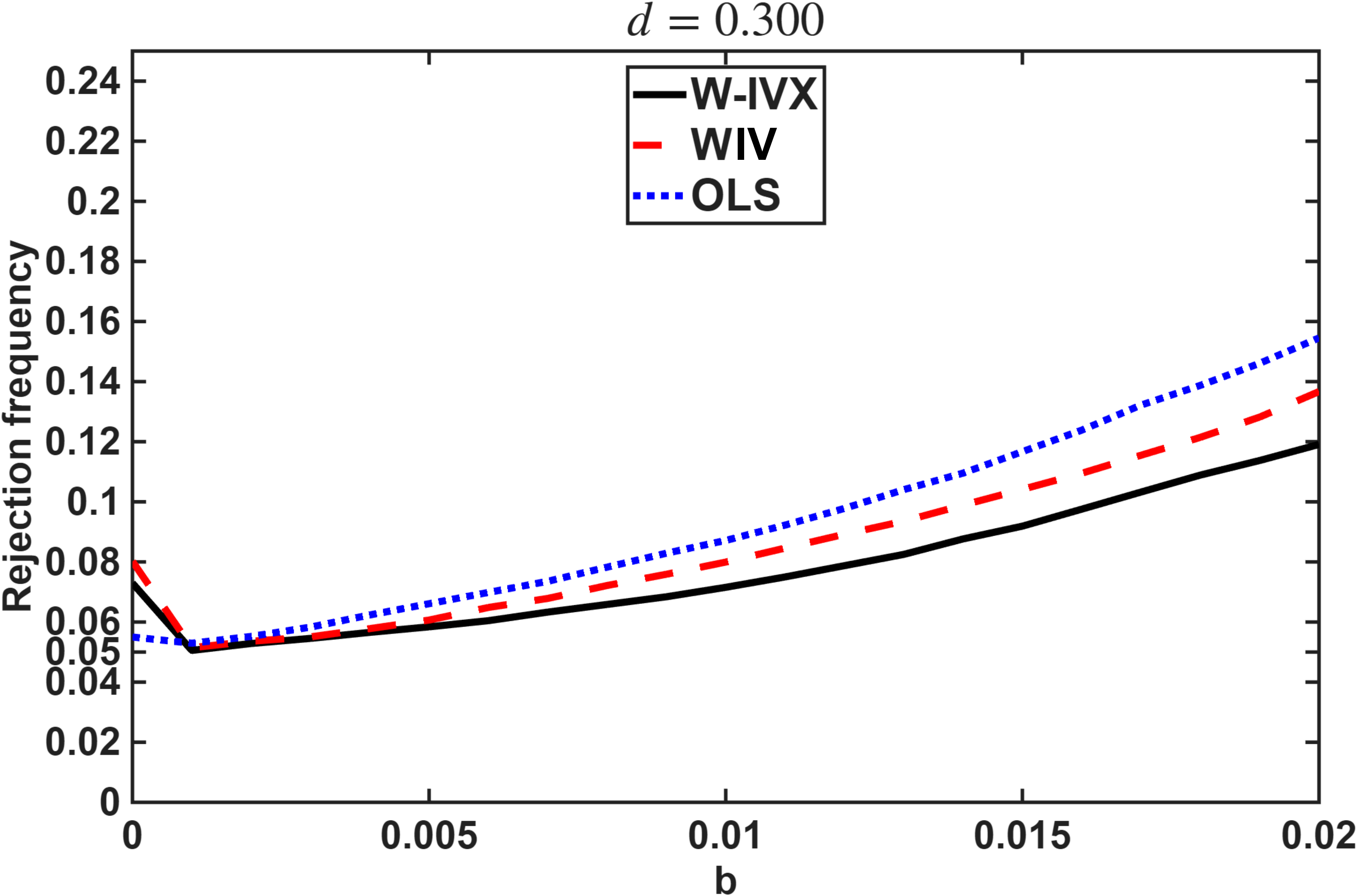}
\includegraphics[width=0.49\linewidth]{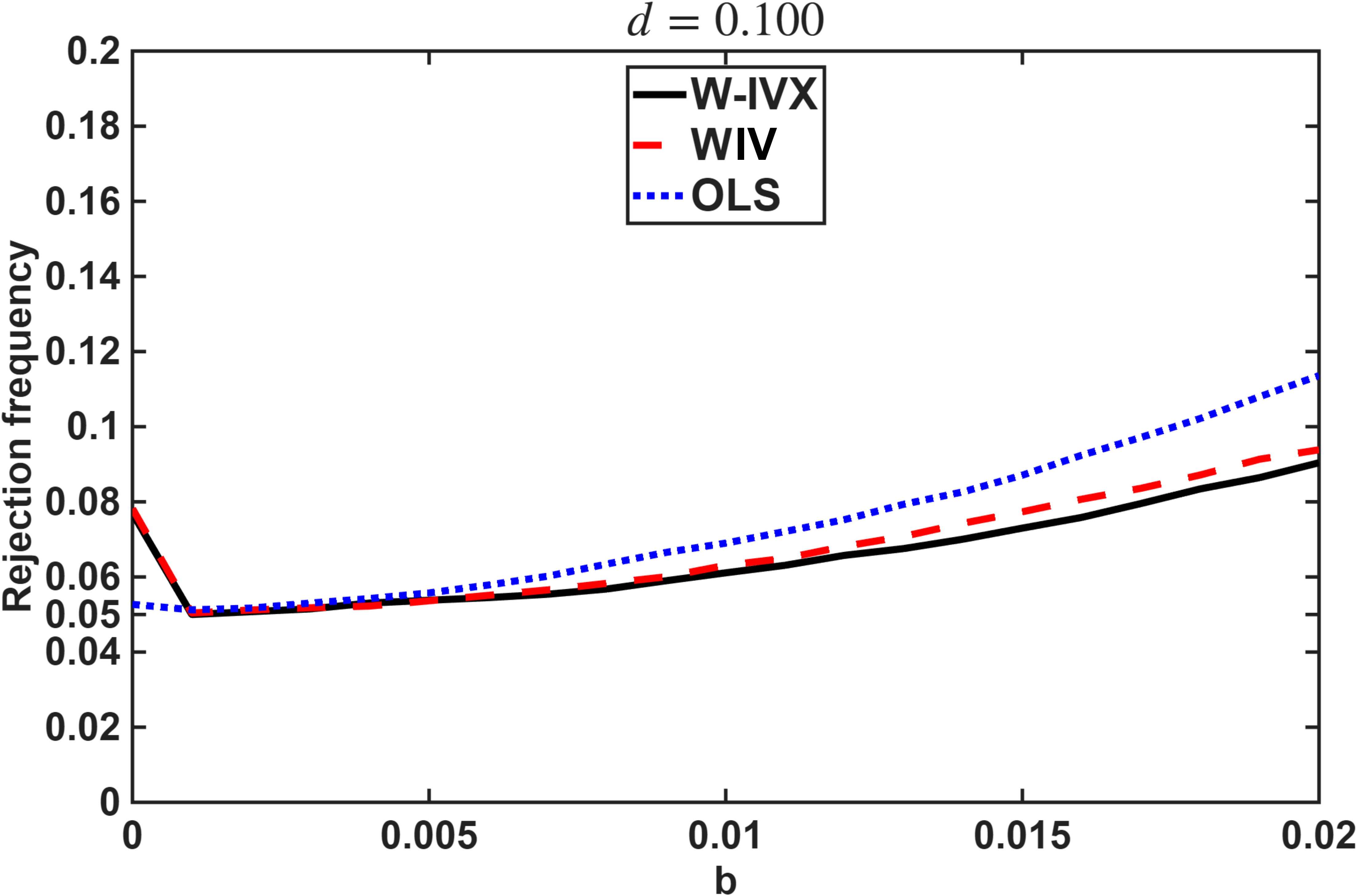}
\end{subfigure} \hfill{}
\begin{subfigure}[t]{0.48\linewidth} \centering \caption{$\alpha=1.8$, $d_{*}=0.444$}
\label{fig:sim-sap-linear-w05-alpha-1p80} \includegraphics[width=0.49\linewidth]{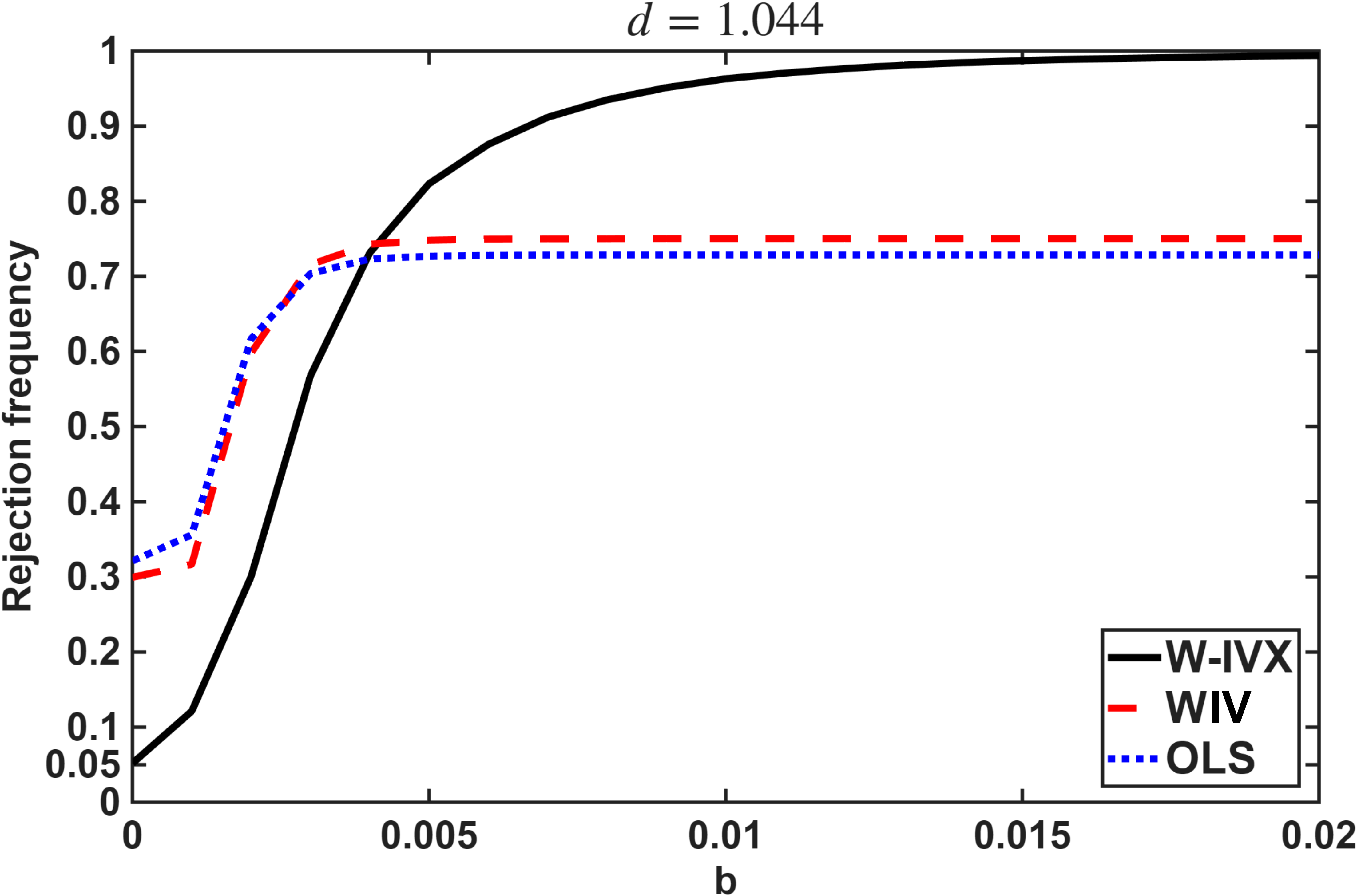}
\includegraphics[width=0.49\linewidth]{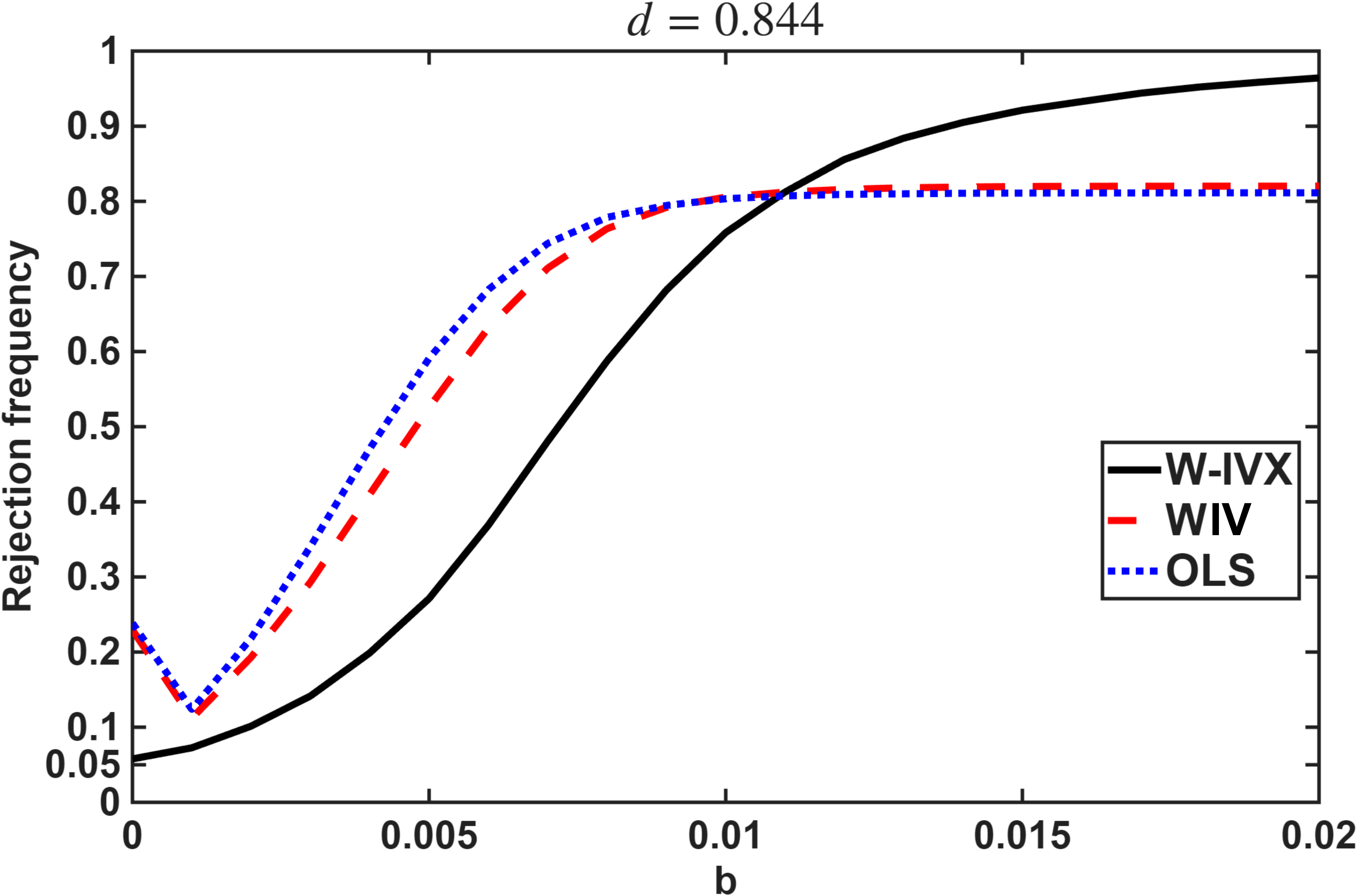}\\[0.6ex]
\includegraphics[width=0.49\linewidth]{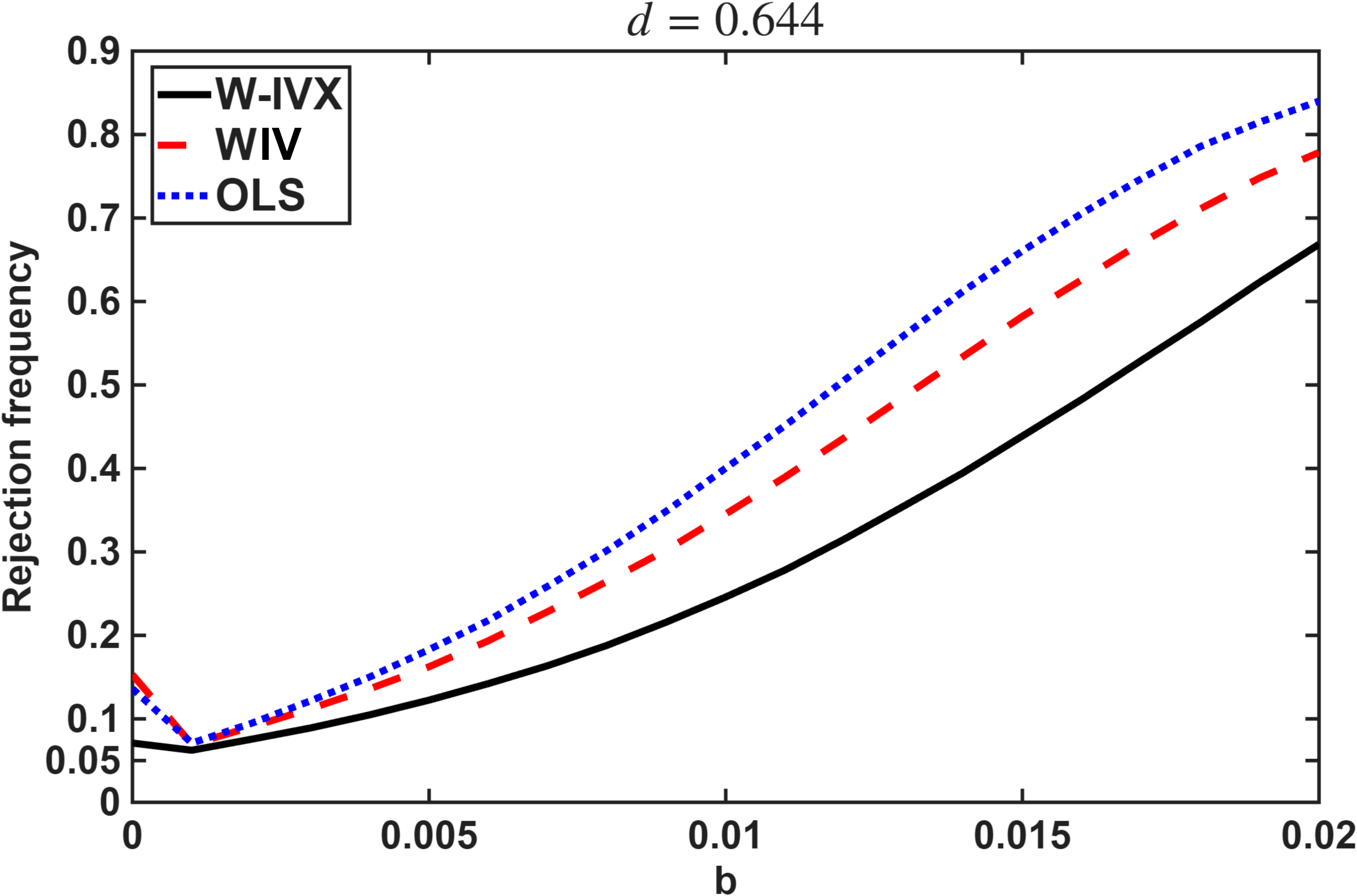}
\includegraphics[width=0.49\linewidth]{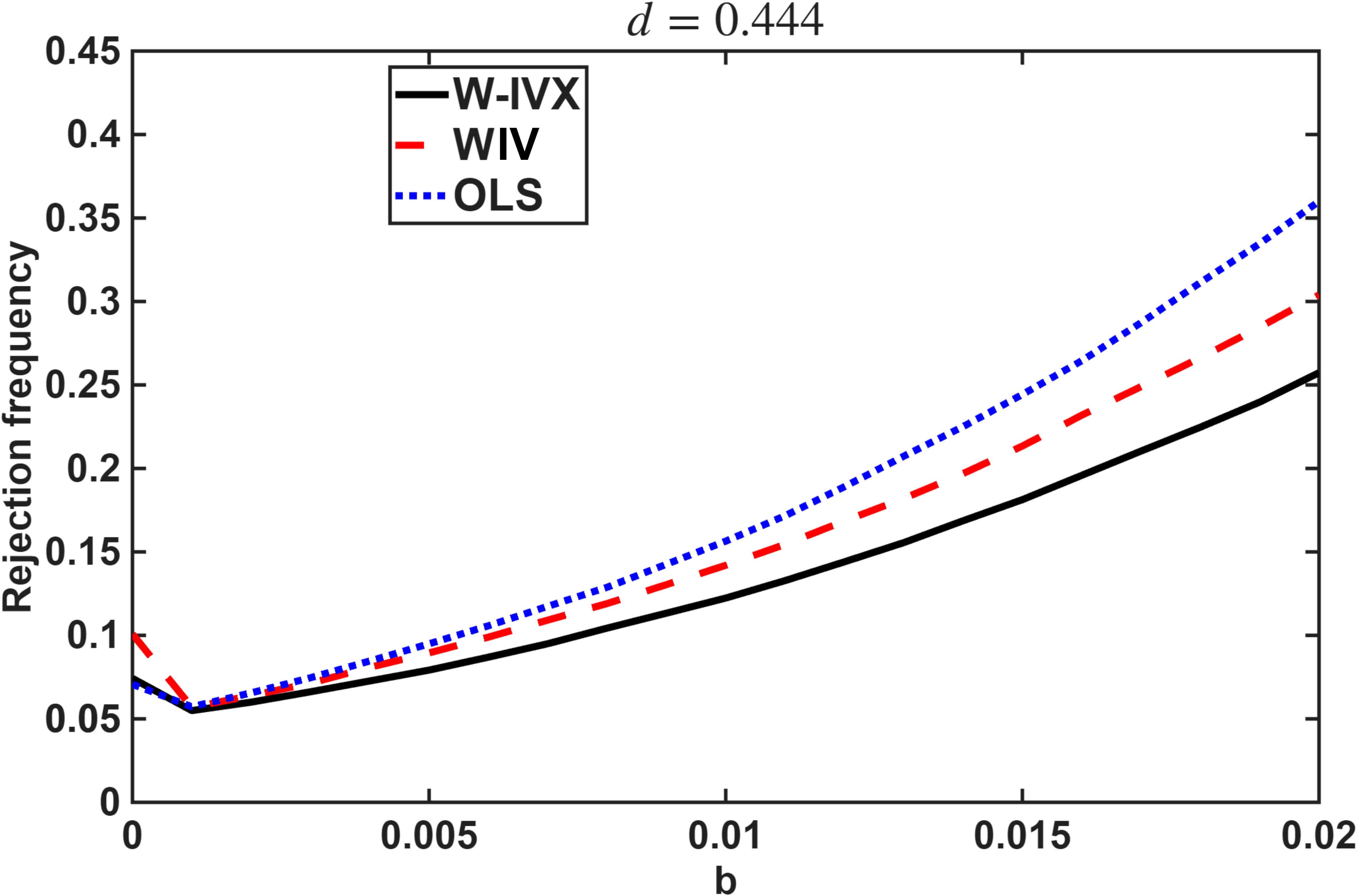}\\[0.6ex]
\includegraphics[width=0.49\linewidth]{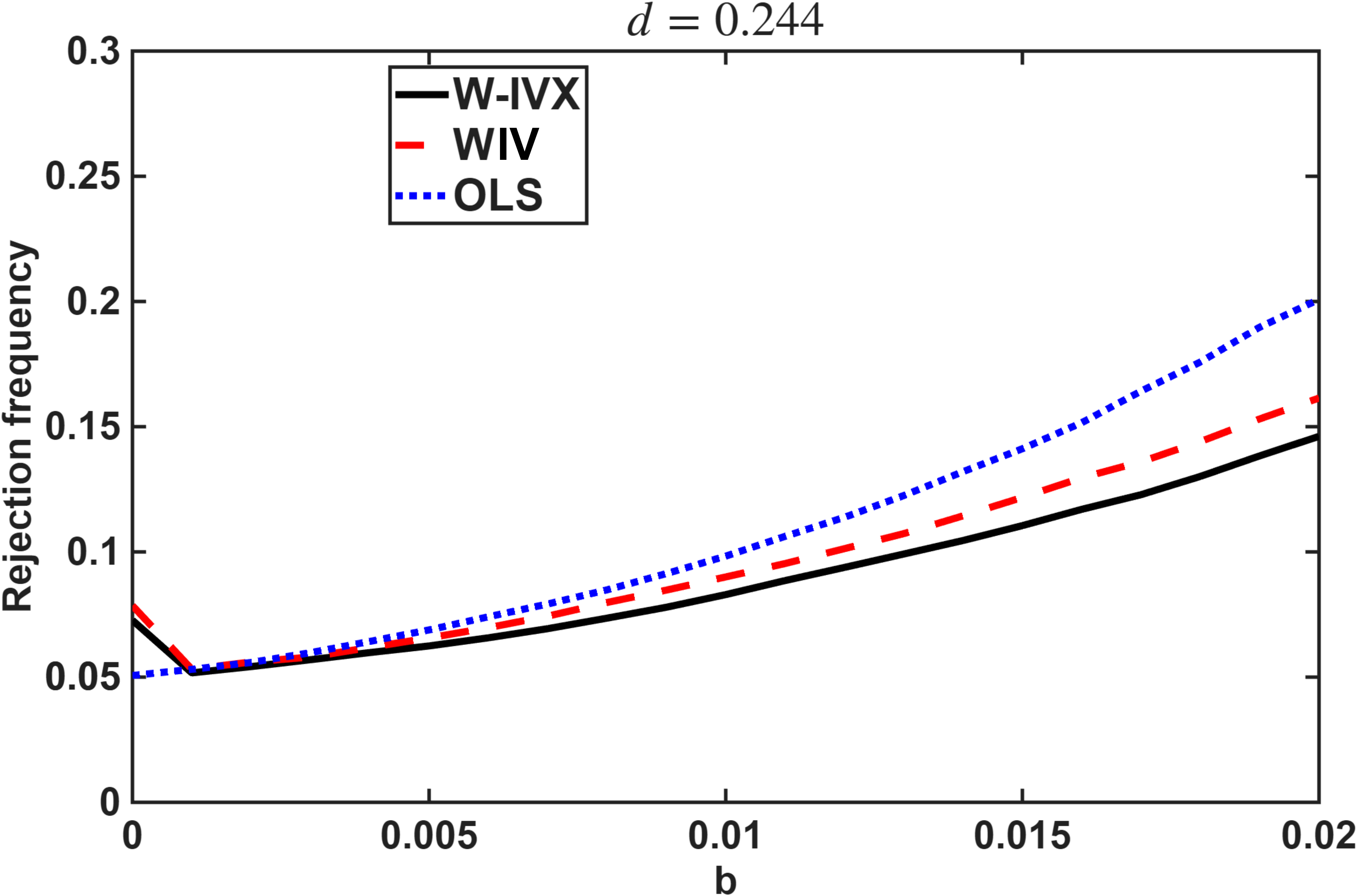}
\includegraphics[width=0.49\linewidth]{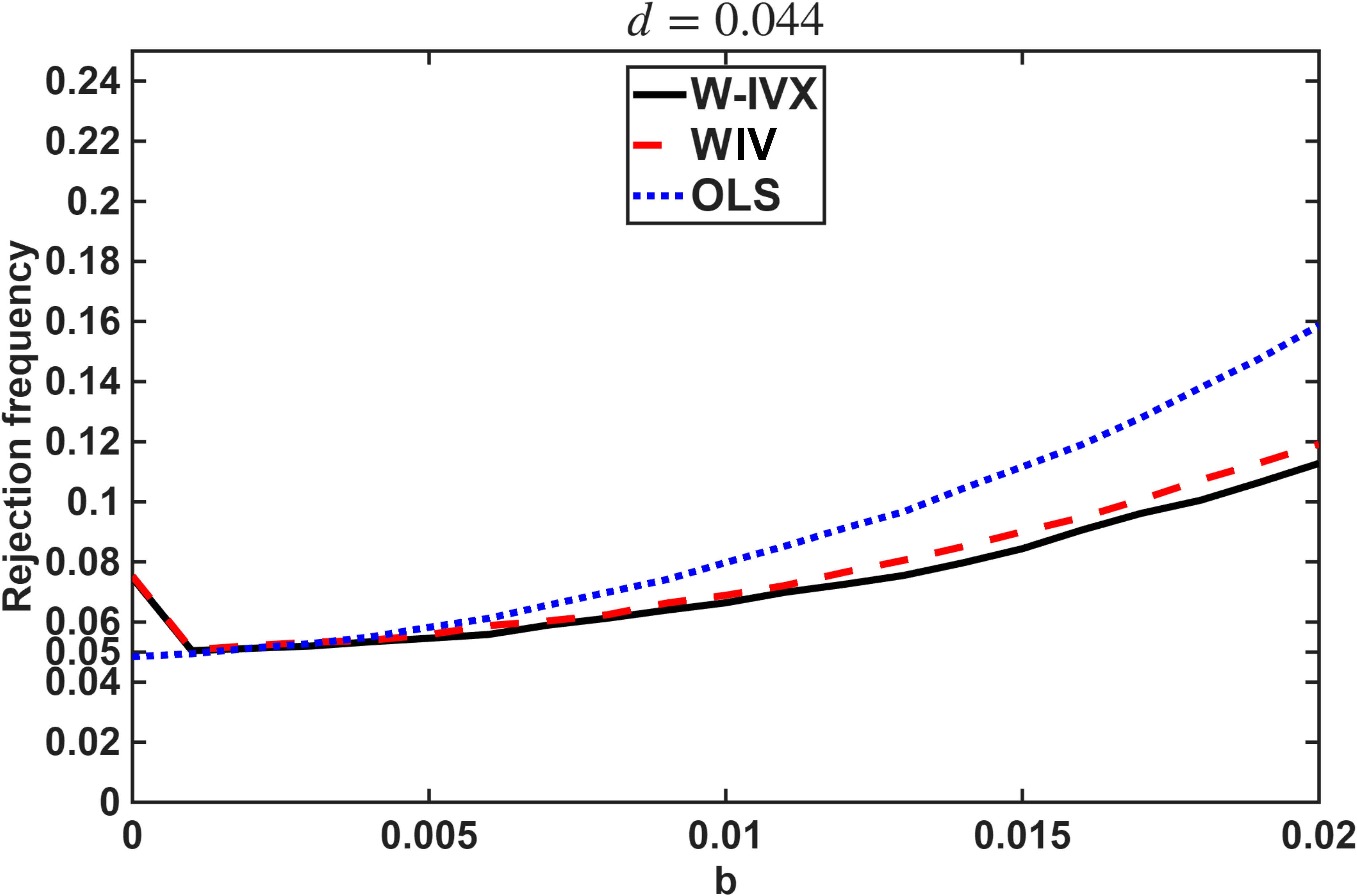}
\end{subfigure}

\end{figure}

\begin{figure}[p]
\ContinuedFloat \centering \captionsetup[subfigure]{position=top}

\setcounter{subfigure}{2}

\caption[]{Size Adjusted Power (nominal size 5\%): $[f(x),w(x)]=[x,x^{\langle0.5\rangle}]$,
$\delta=-0.95$, $n=1000$ (continued)}
\label{fig:sim-sap-linear-w05-continued}

\begin{subfigure}[t]{0.48\linewidth} \centering \caption{$\alpha=1.6$, $d_{*}=0.5$}
\label{fig:sim-sap-linear-w05-alpha-1p60}\includegraphics[width=0.49\linewidth]{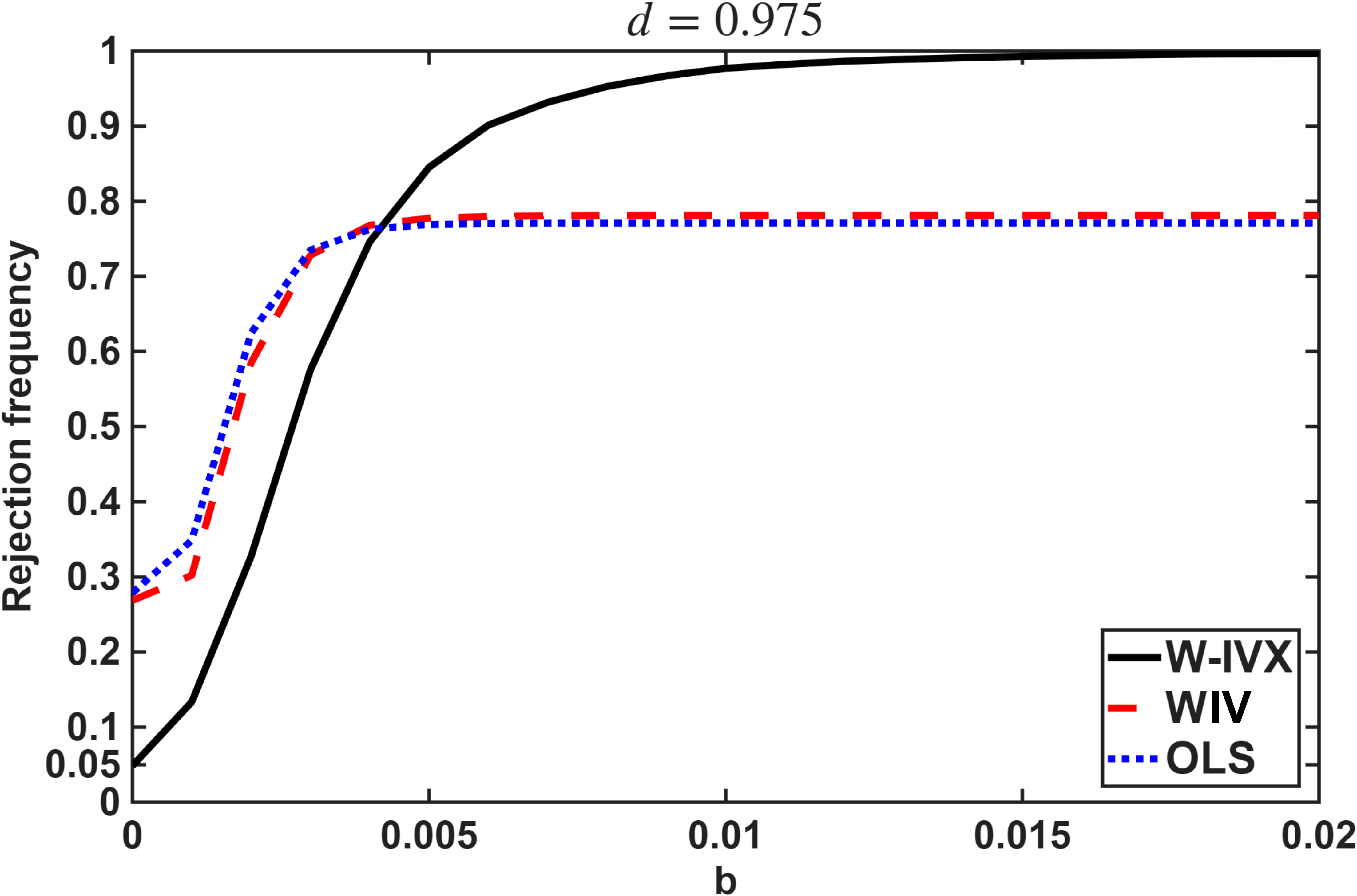}
\includegraphics[width=0.49\linewidth]{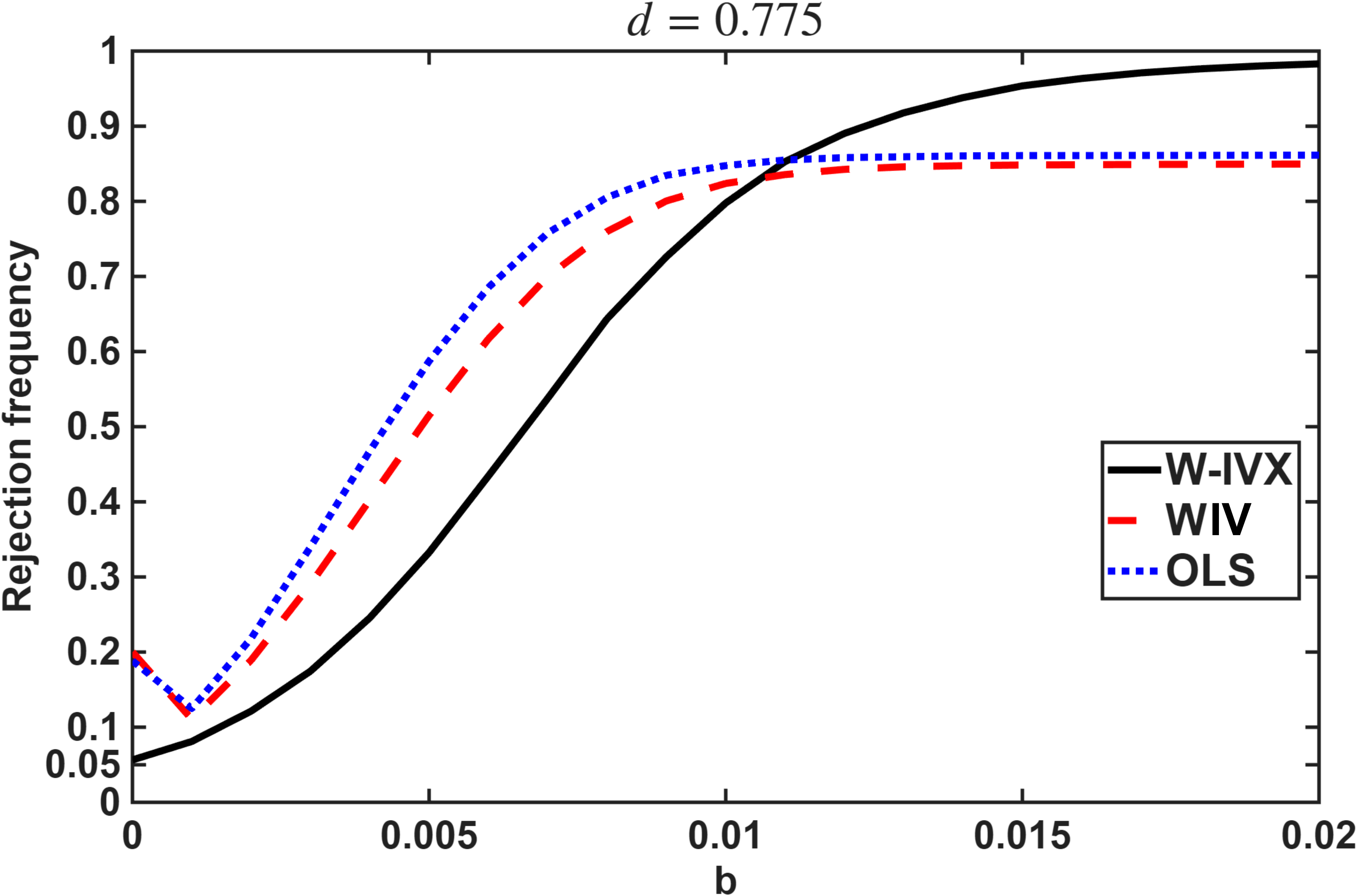}\\[0.6ex]
\includegraphics[width=0.49\linewidth]{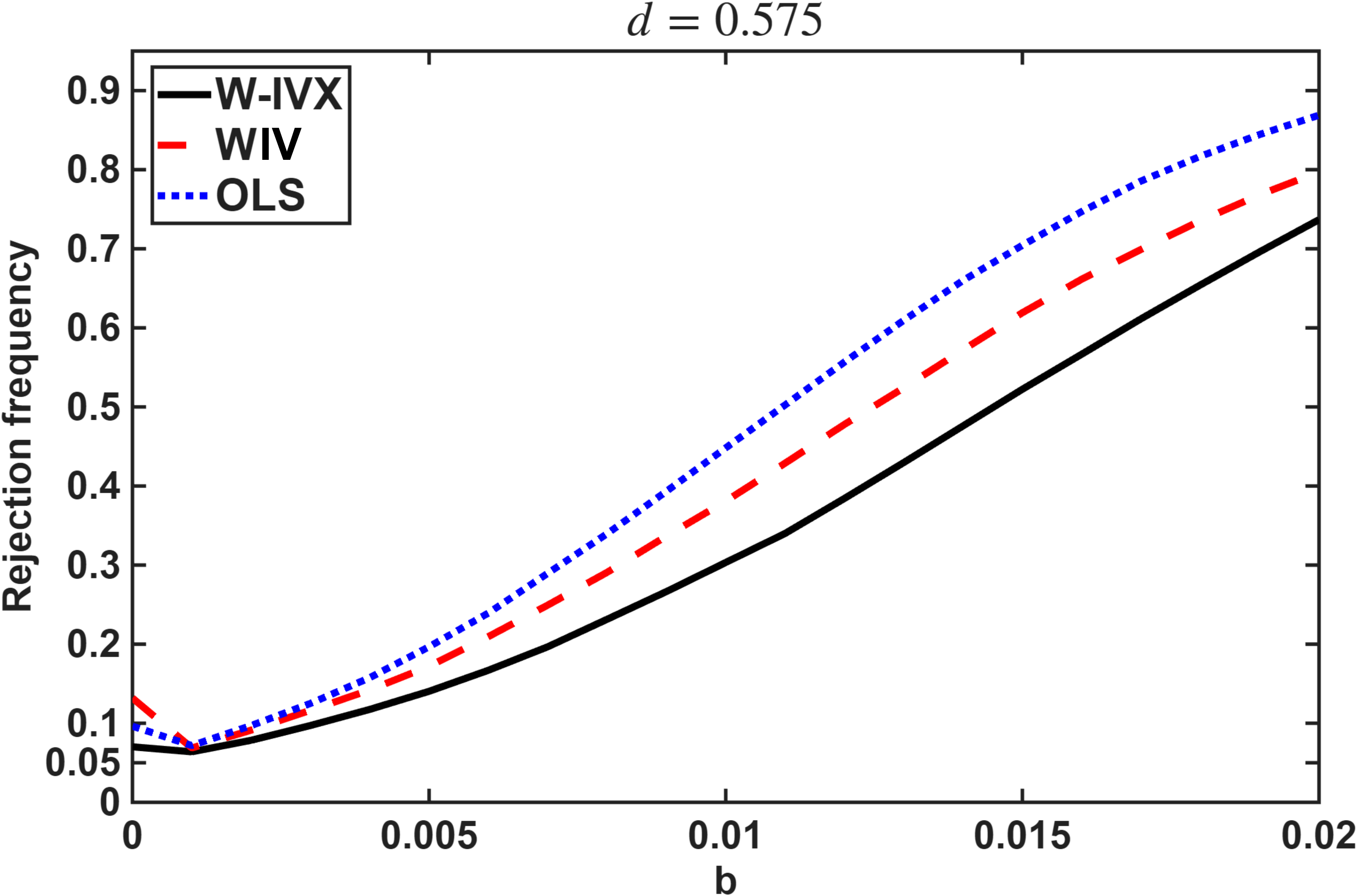}
\includegraphics[width=0.49\linewidth]{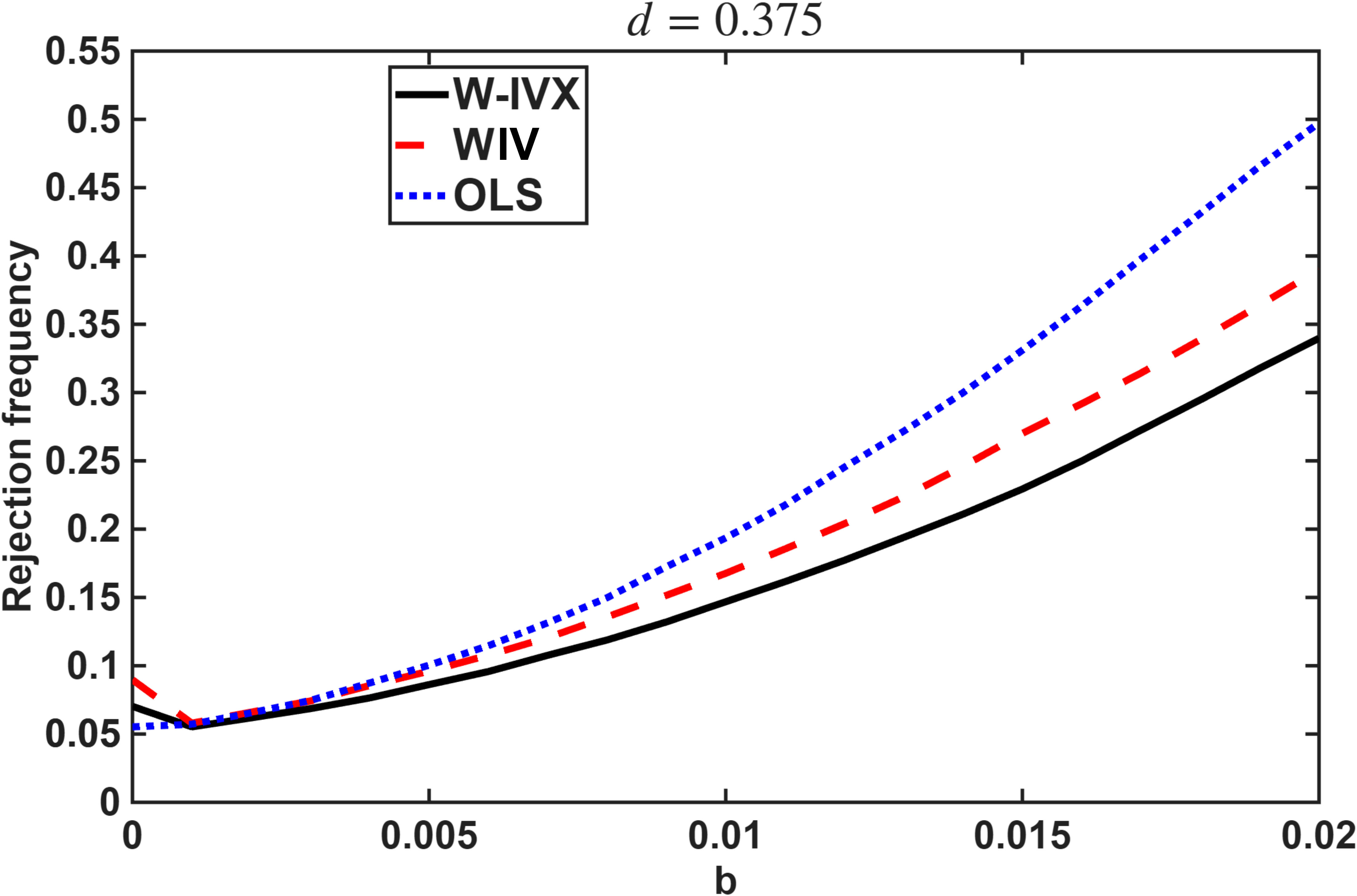}\\[0.6ex]
\includegraphics[width=0.49\linewidth]{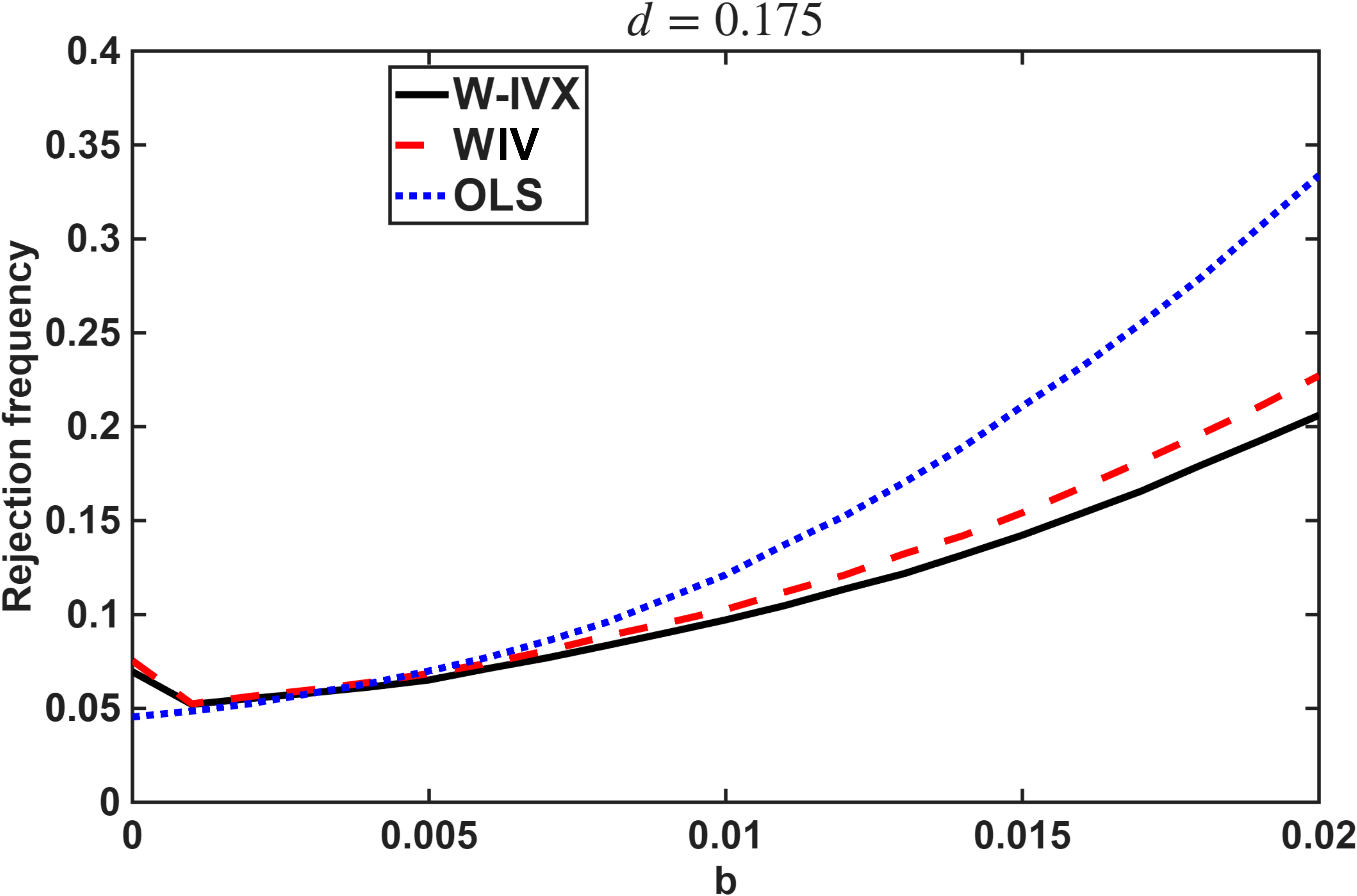}
\includegraphics[width=0.49\linewidth]{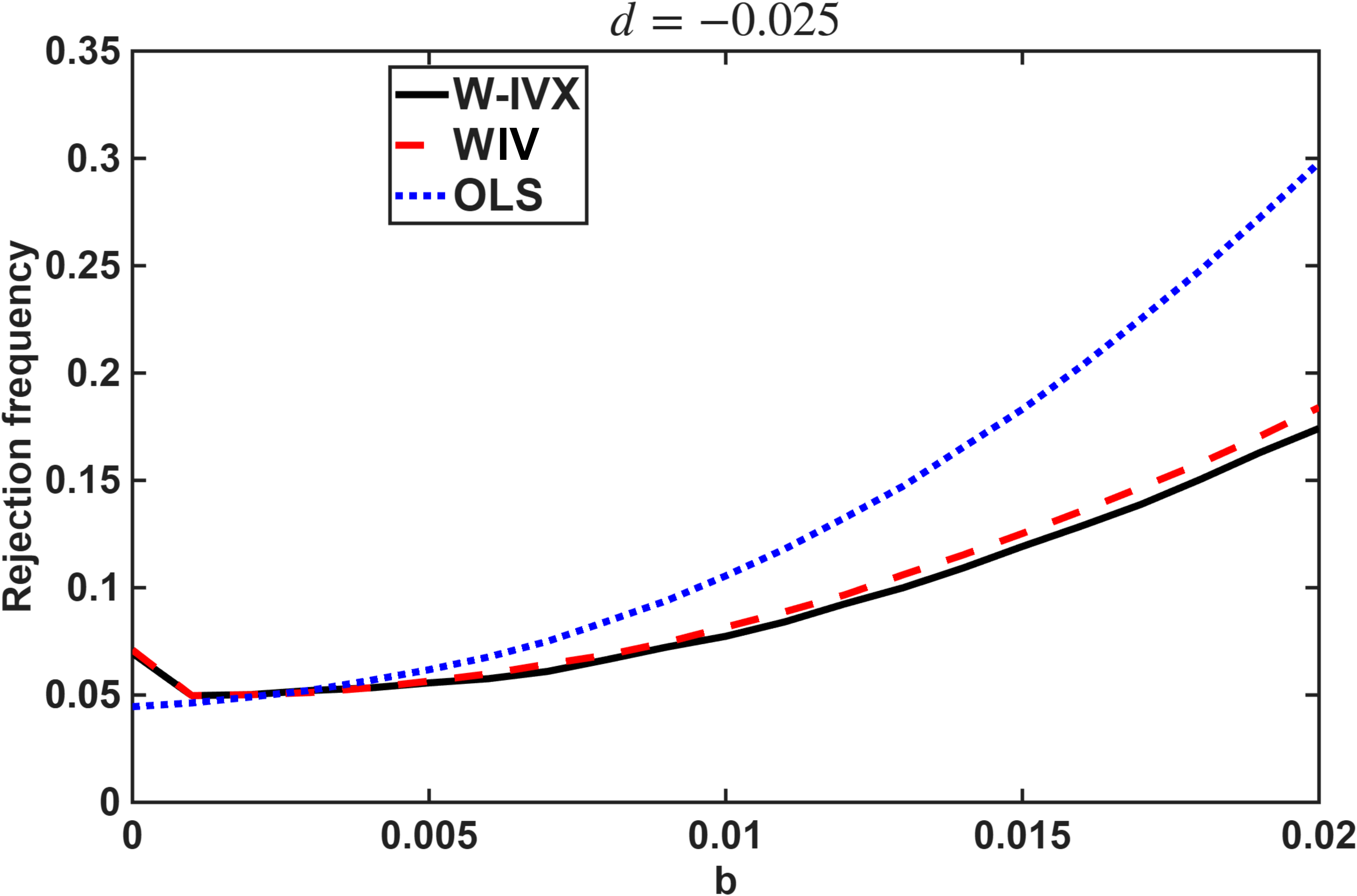}
\end{subfigure} \hfill{}
\begin{subfigure}[t]{0.48\linewidth} \centering \caption{$\alpha=1.4$, $d_{*}=0.444$}
\label{fig:sim-sap-linear-w05-alpha-1p40}\includegraphics[width=0.49\linewidth]{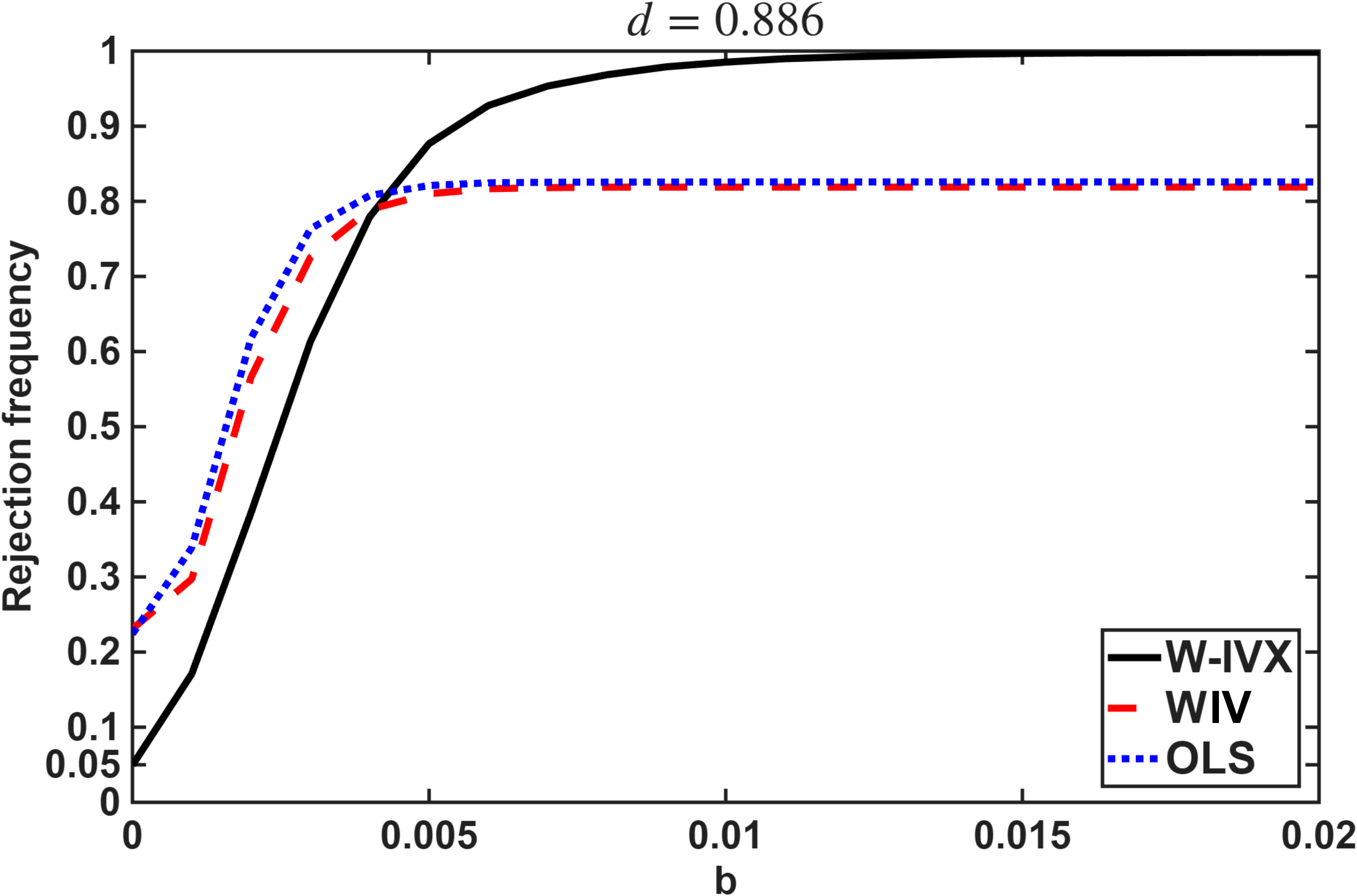}
\includegraphics[width=0.49\linewidth]{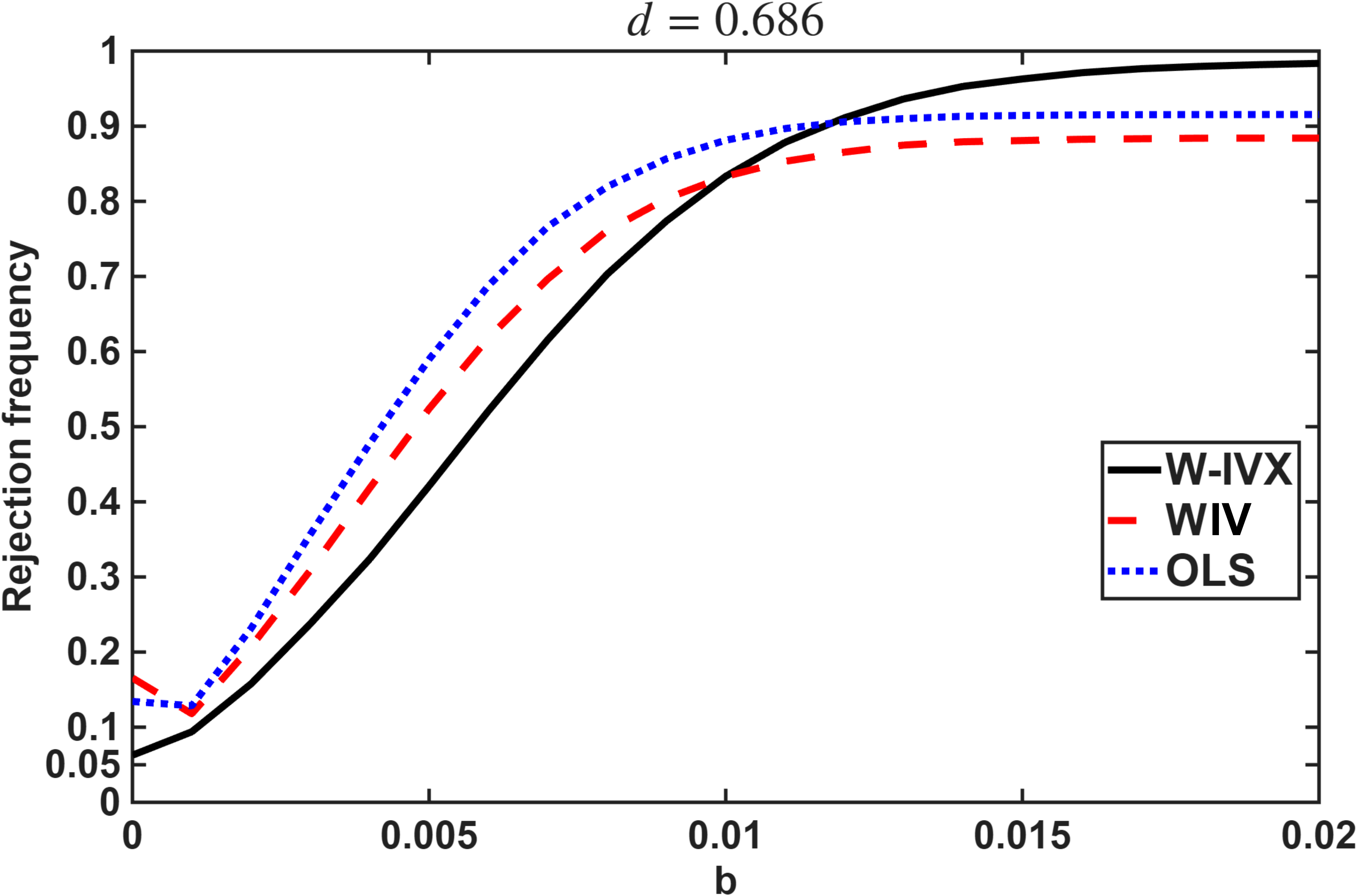}\\[0.6ex]
\includegraphics[width=0.49\linewidth]{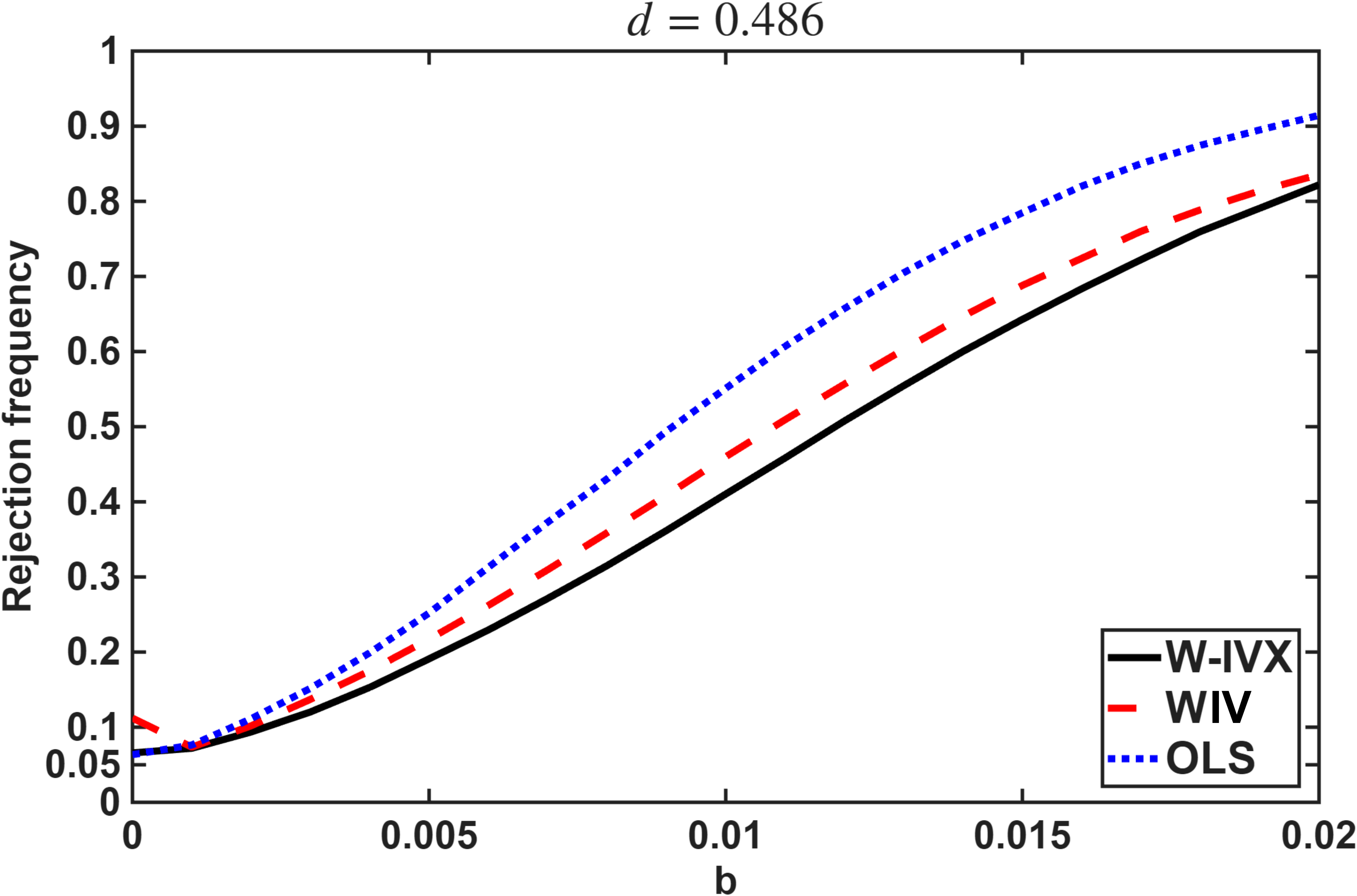}
\includegraphics[width=0.49\linewidth]{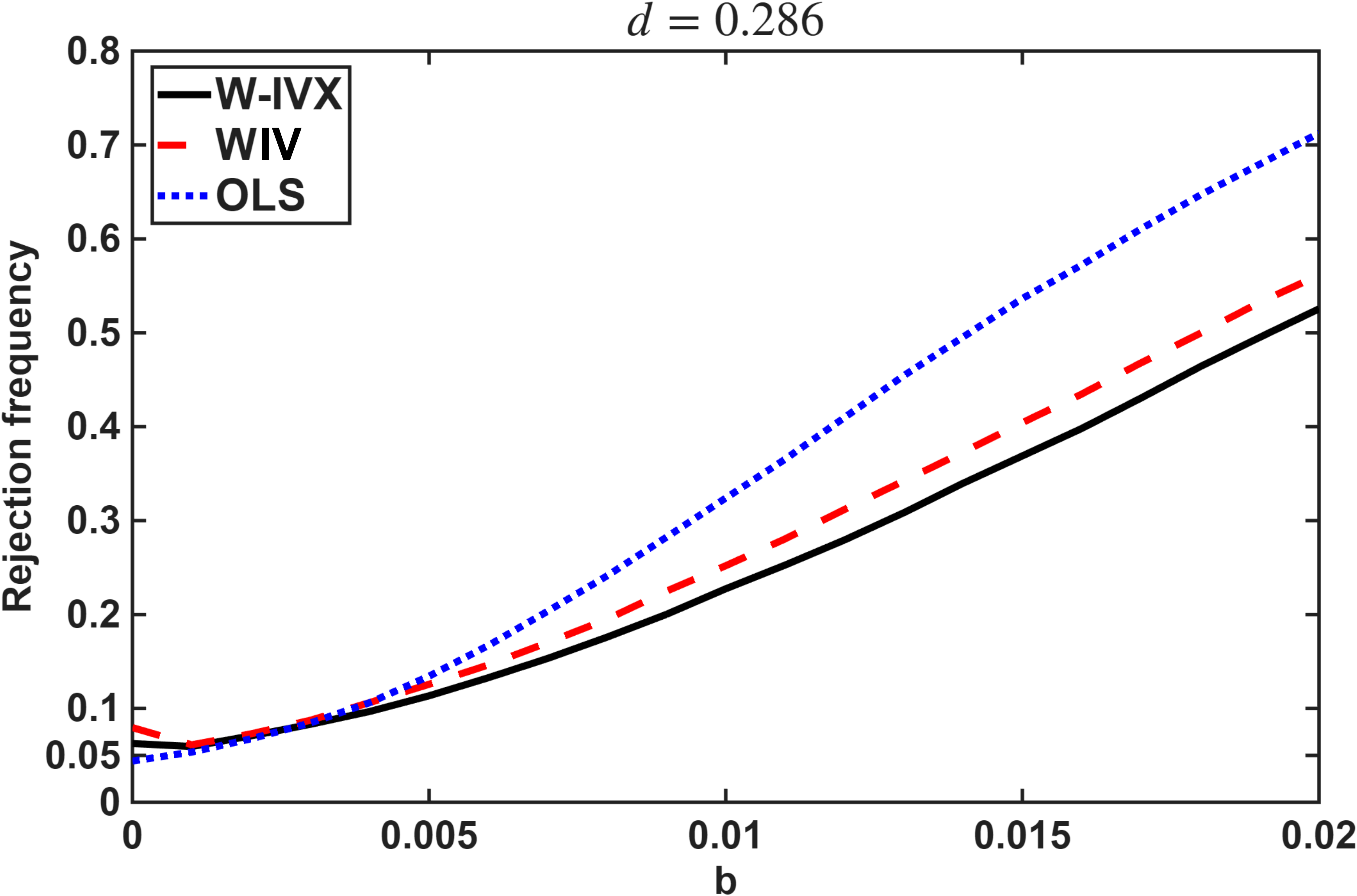}\\[0.6ex]
\includegraphics[width=0.49\linewidth]{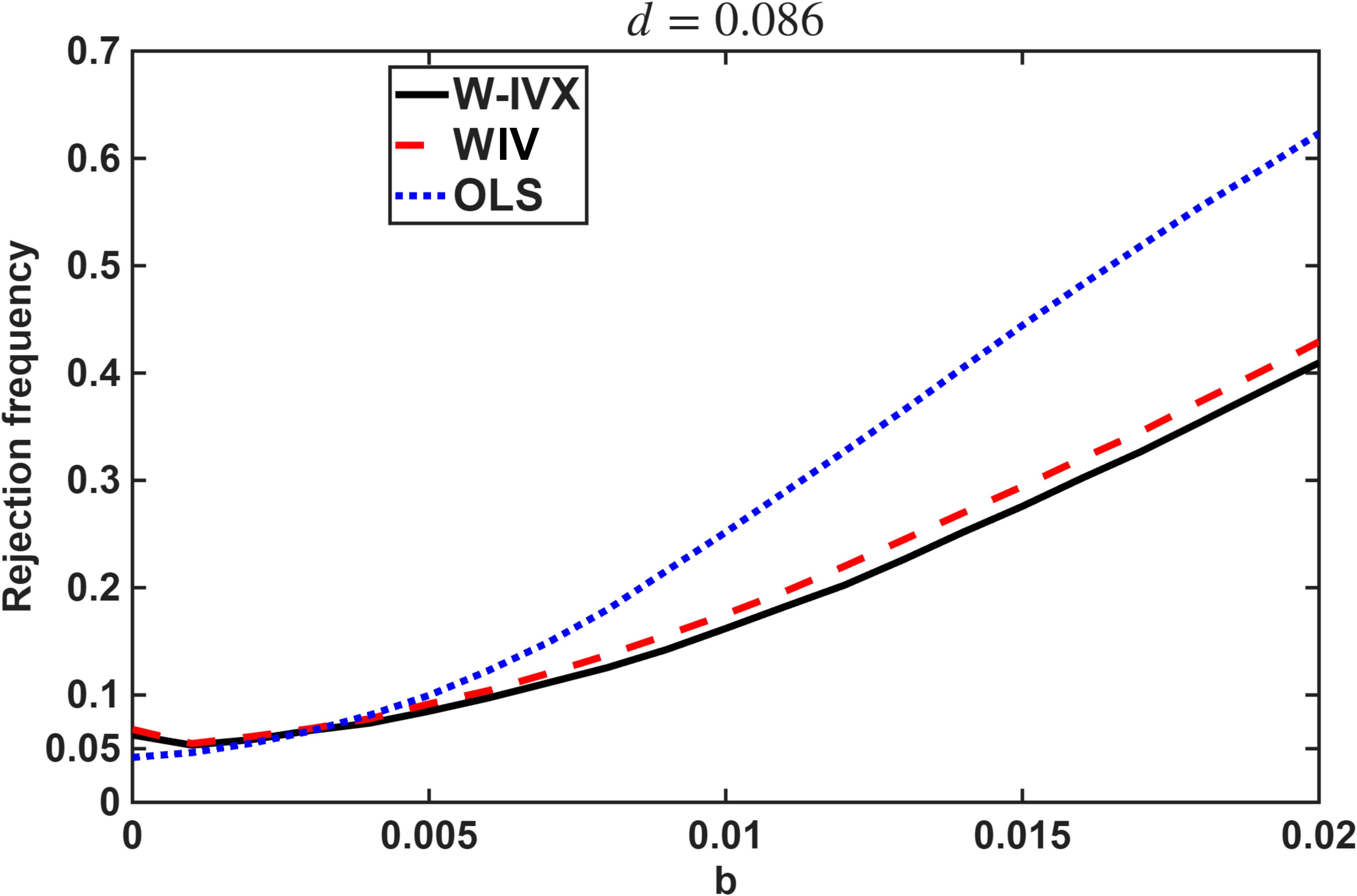}
\includegraphics[width=0.49\linewidth]{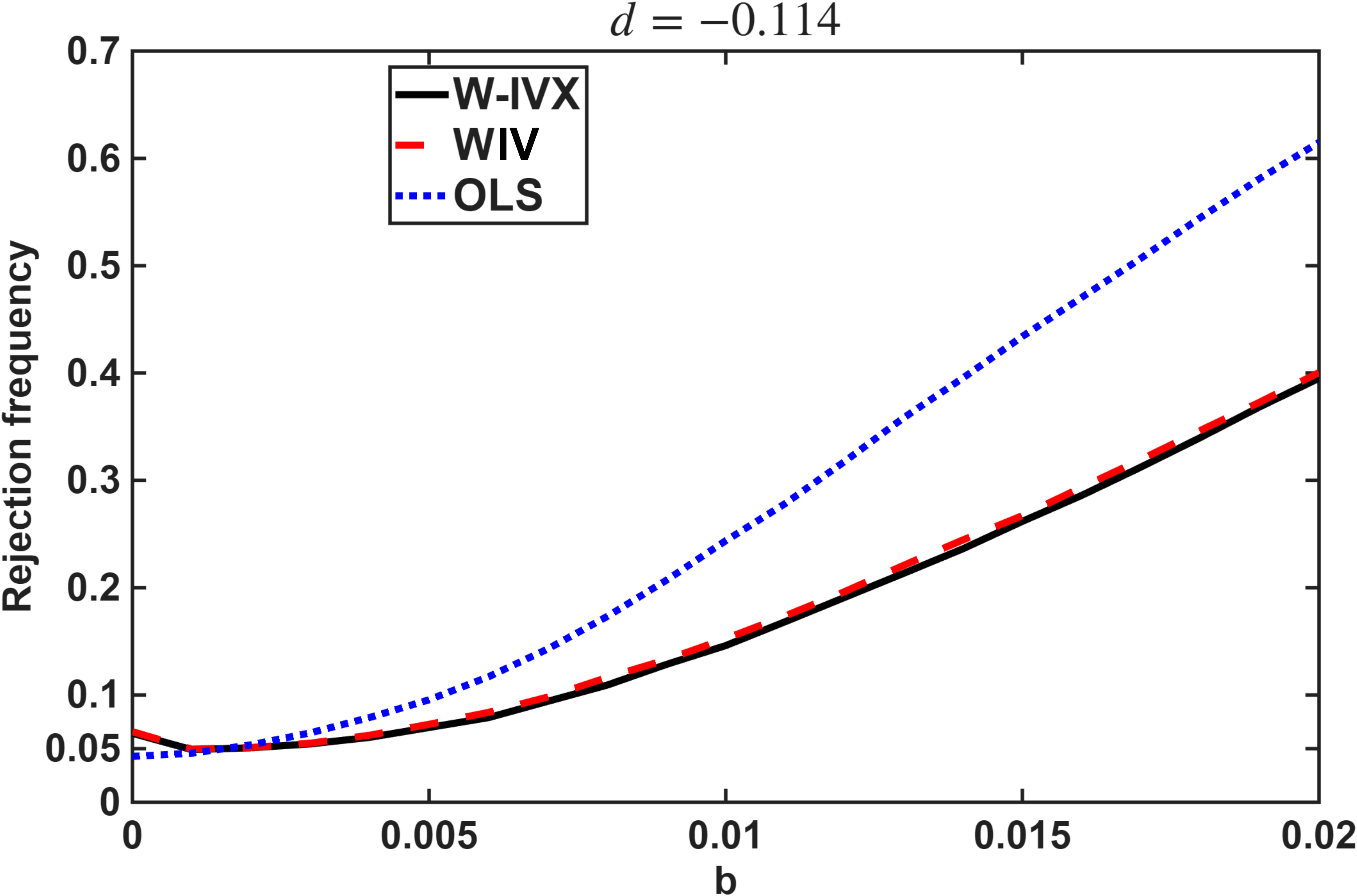}
\end{subfigure}
\end{figure}

\end{landscape} \clearpage{}

\begin{landscape}

\begin{figure}[p]
\caption{Size Adjusted Power (nominal size 5\%): \textbf{DGP1}, $[f(x),w(x)]=[x,x^{\langle0.75\rangle}]$,
$\delta=-0.95$, $n=1000$}
\label{fig:sim-sap-linear-w075}\centering \captionsetup[subfigure]{position=top}

\begin{subfigure}[t]{0.48\linewidth} \centering \caption{$\alpha=2$, $d_{*}=0.5$}
\label{fig:sim-sap-linear-w075-alpha-2p00} \includegraphics[width=0.49\linewidth]{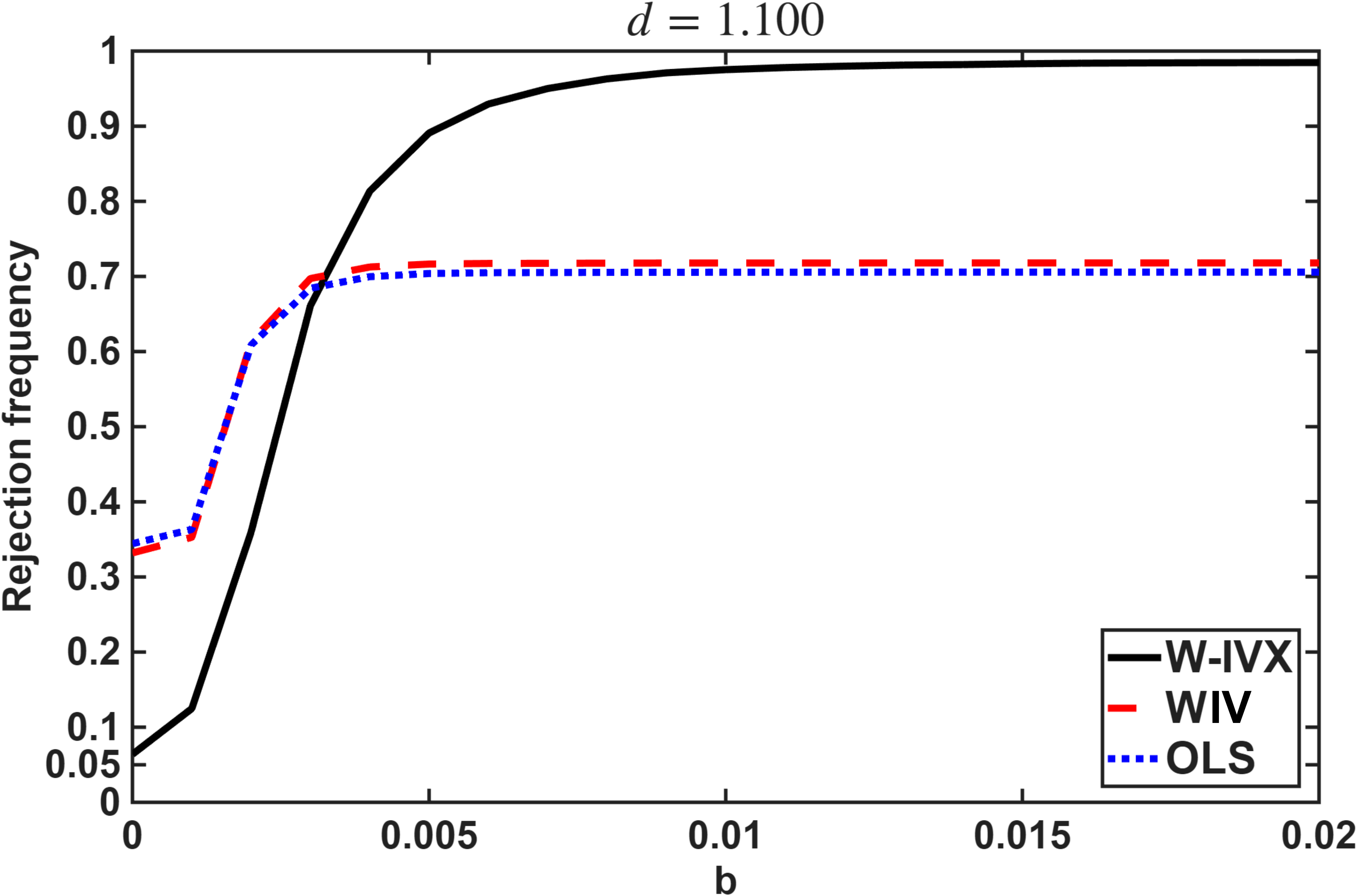}
\includegraphics[width=0.49\linewidth]{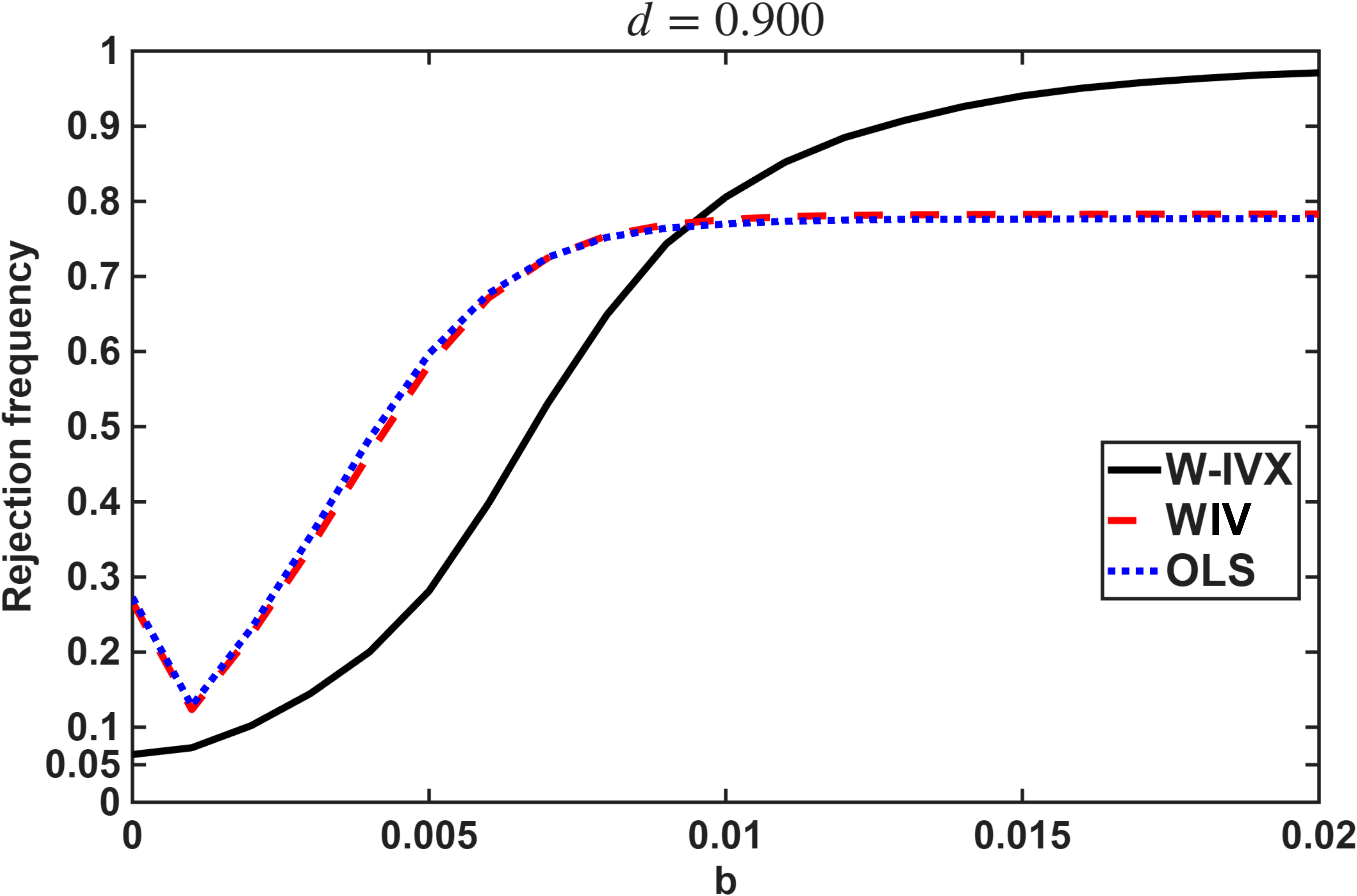}\\[0.6ex]
\includegraphics[width=0.49\linewidth]{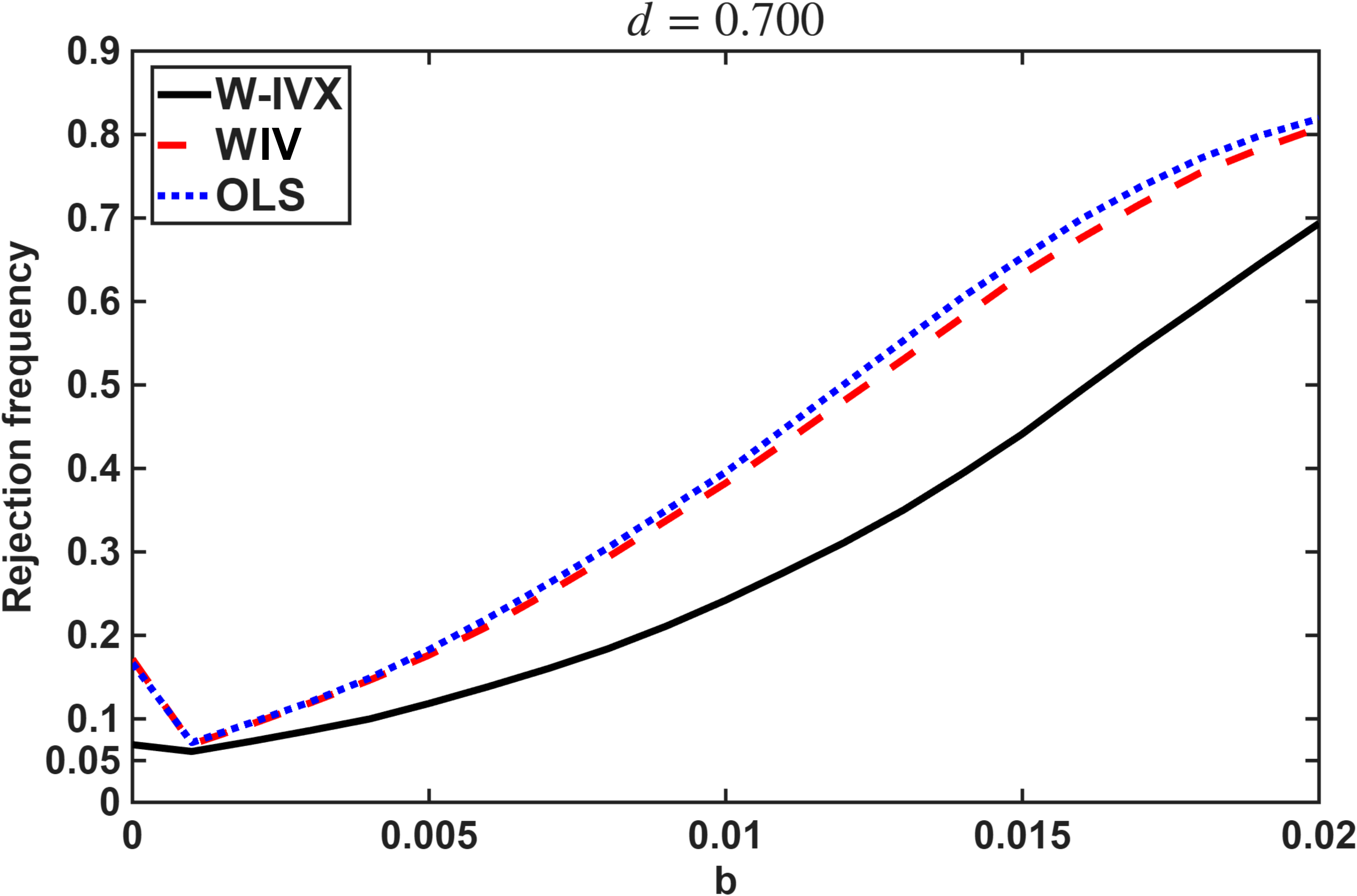}
\includegraphics[width=0.49\linewidth]{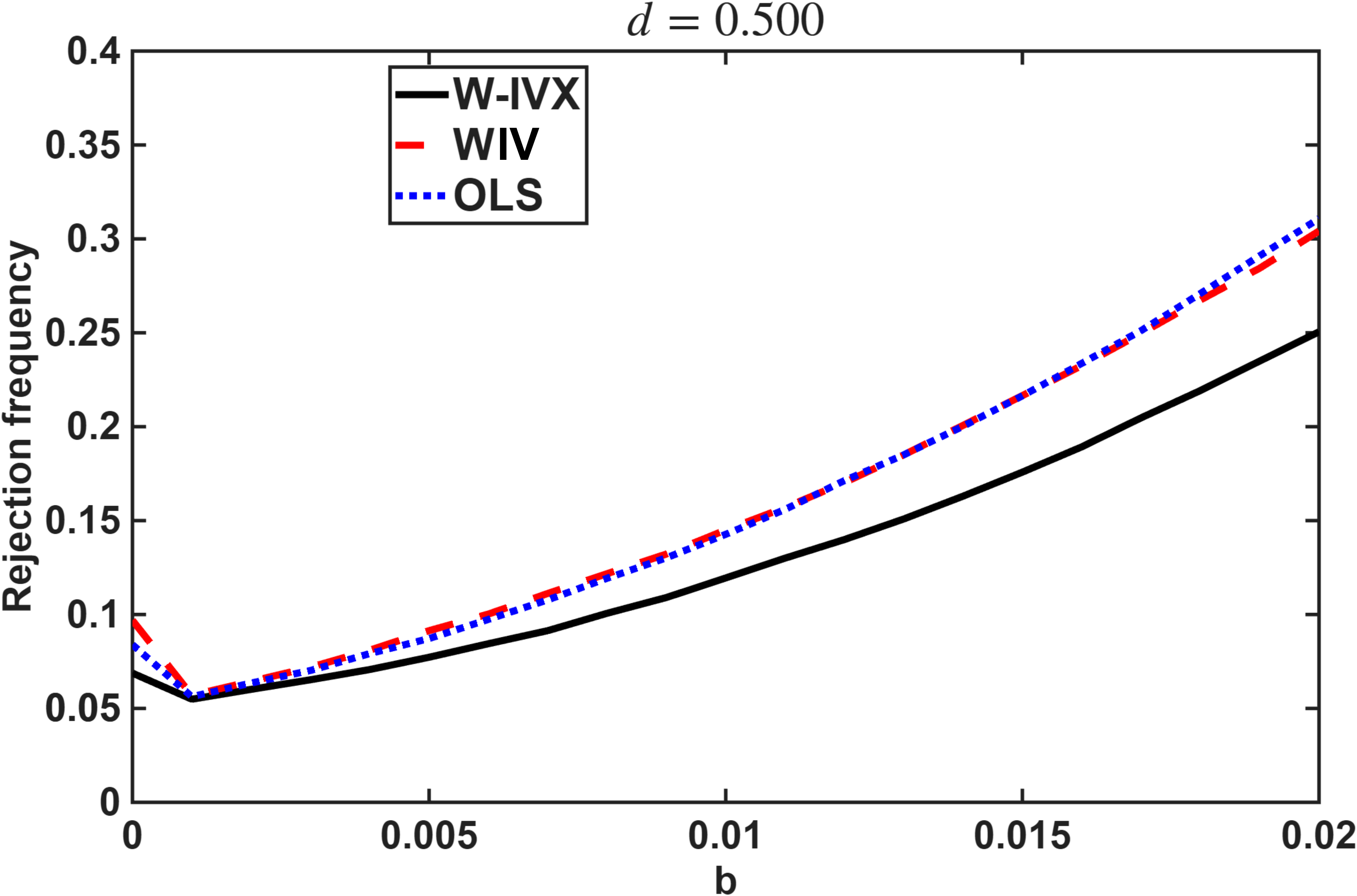}\\[0.6ex]
\includegraphics[width=0.49\linewidth]{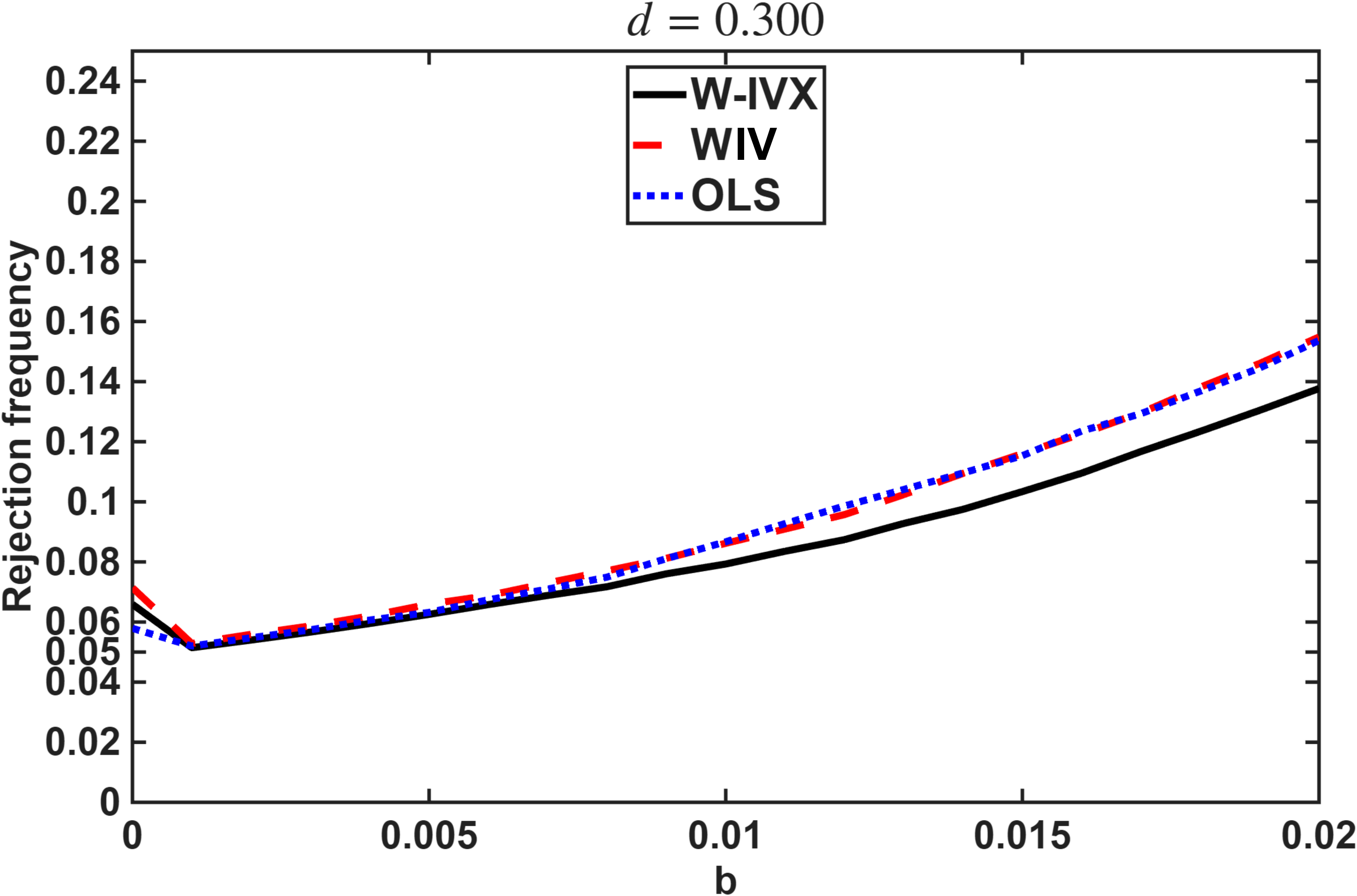}
\includegraphics[width=0.49\linewidth]{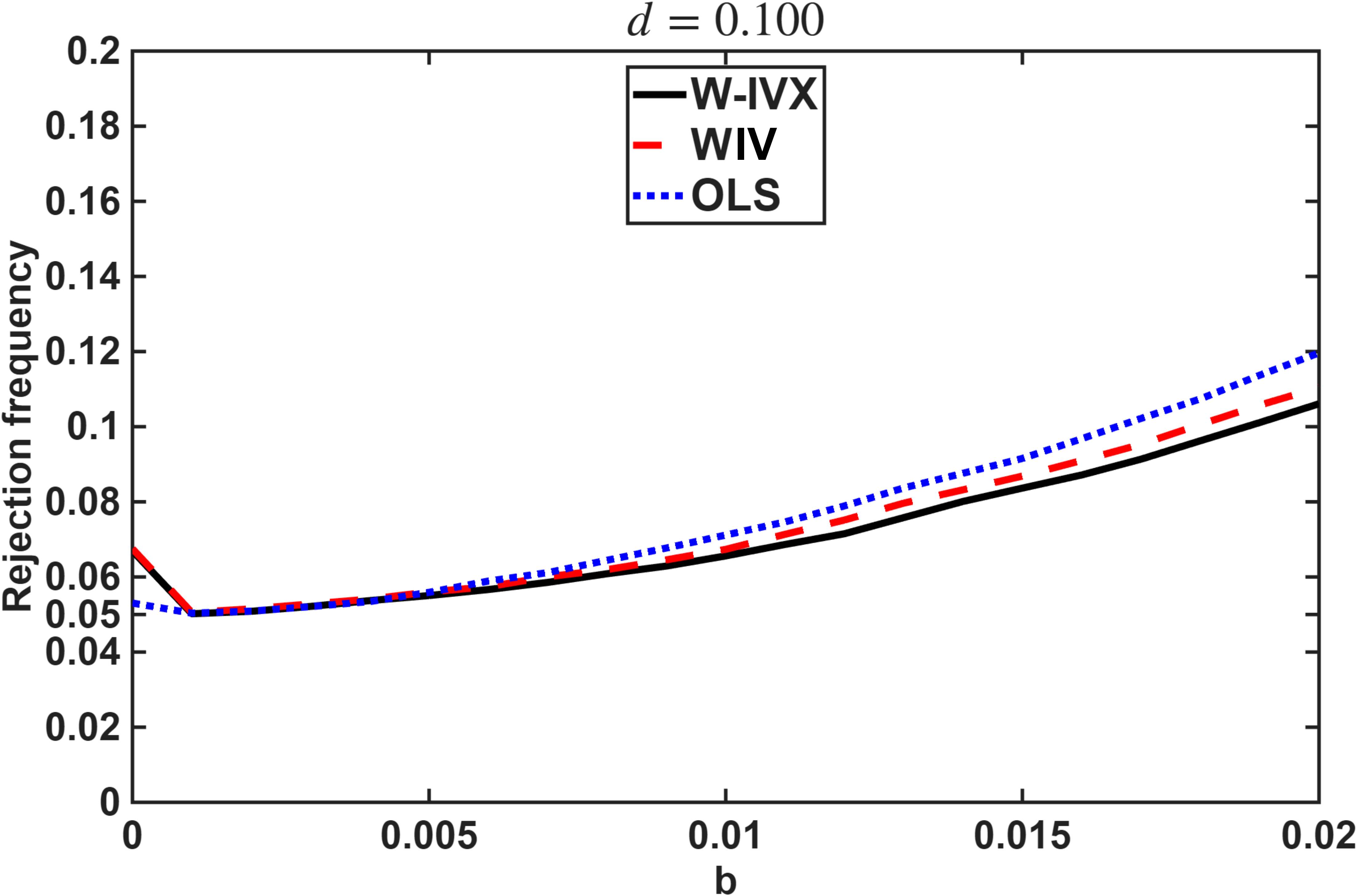}
\end{subfigure} \hfill{}
\begin{subfigure}[t]{0.48\linewidth} \centering \caption{$\alpha=1.8$, $d_{*}=0.444$}
\label{fig:sim-sap-linear-w075-alpha-1p80} \includegraphics[width=0.49\linewidth]{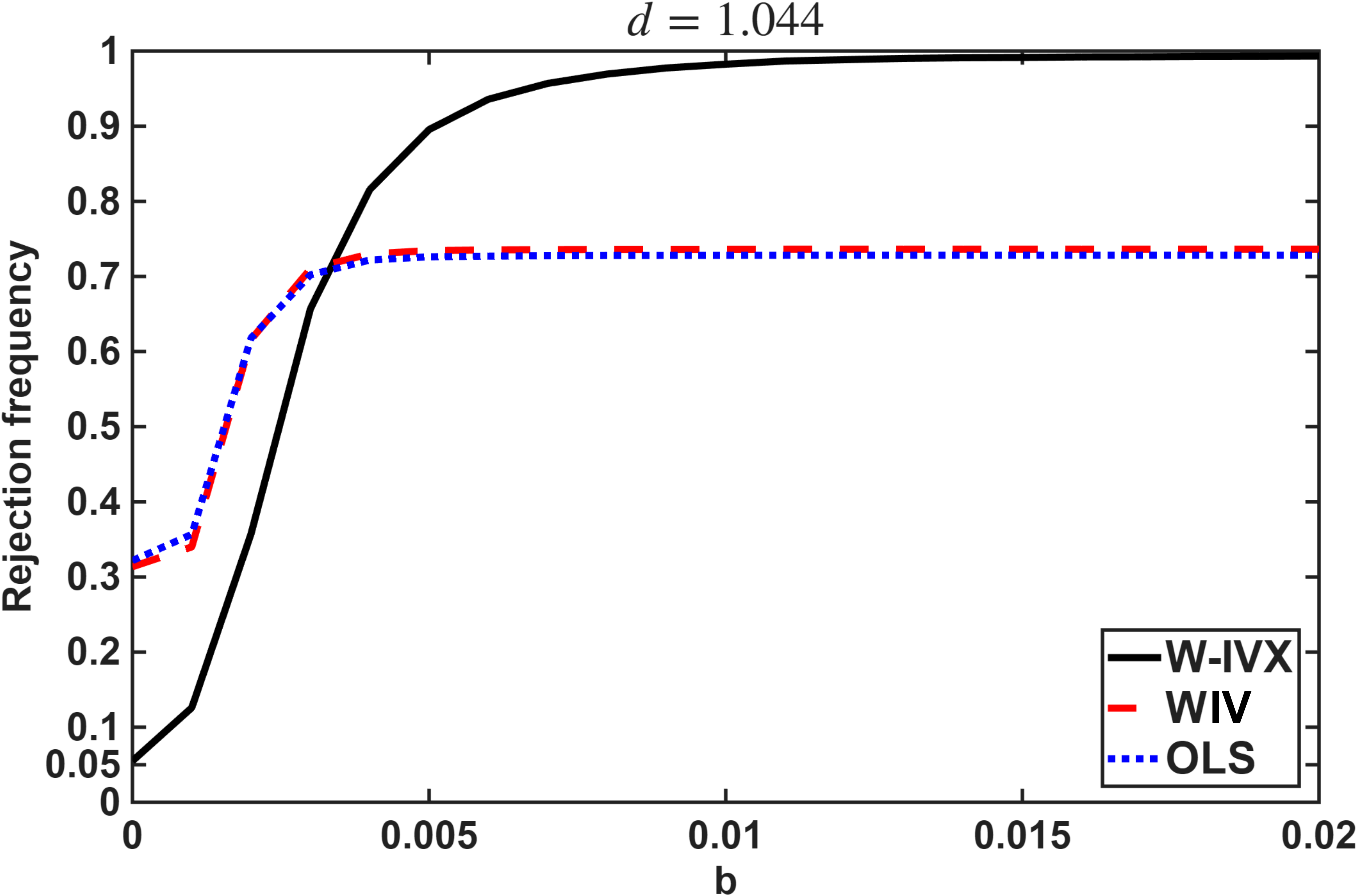}
\includegraphics[width=0.49\linewidth]{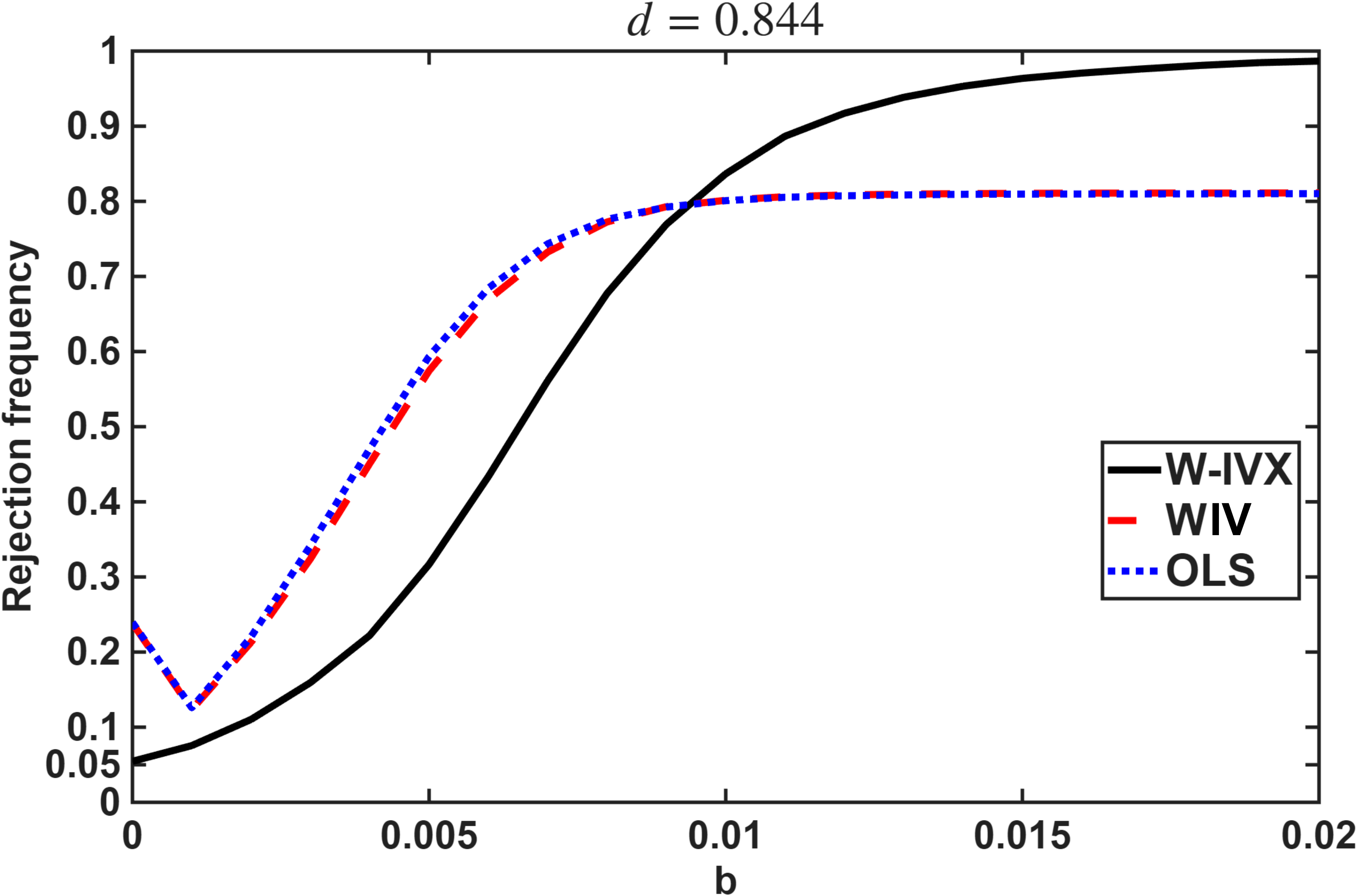}\\[0.6ex]
\includegraphics[width=0.49\linewidth]{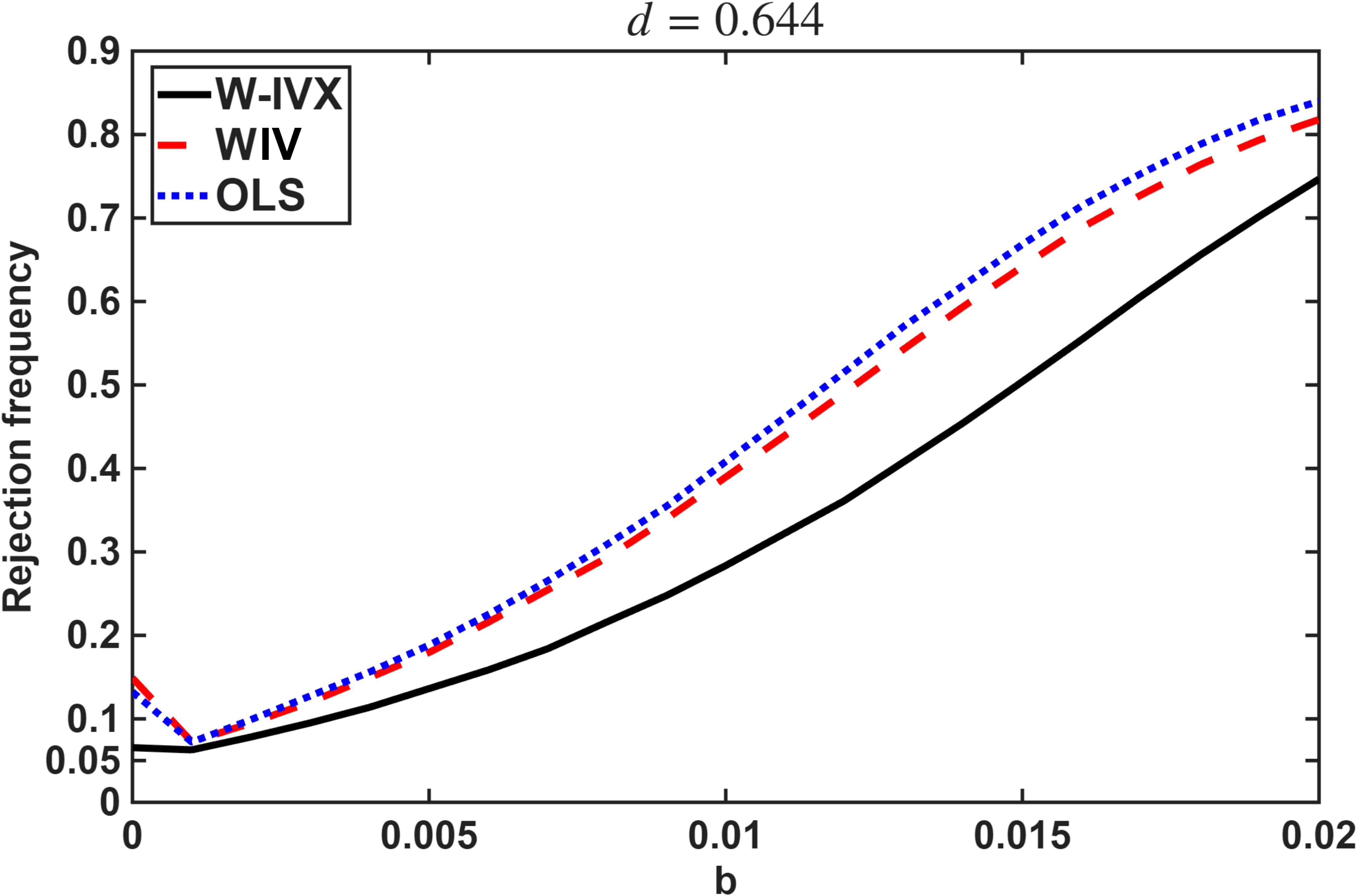}
\includegraphics[width=0.49\linewidth]{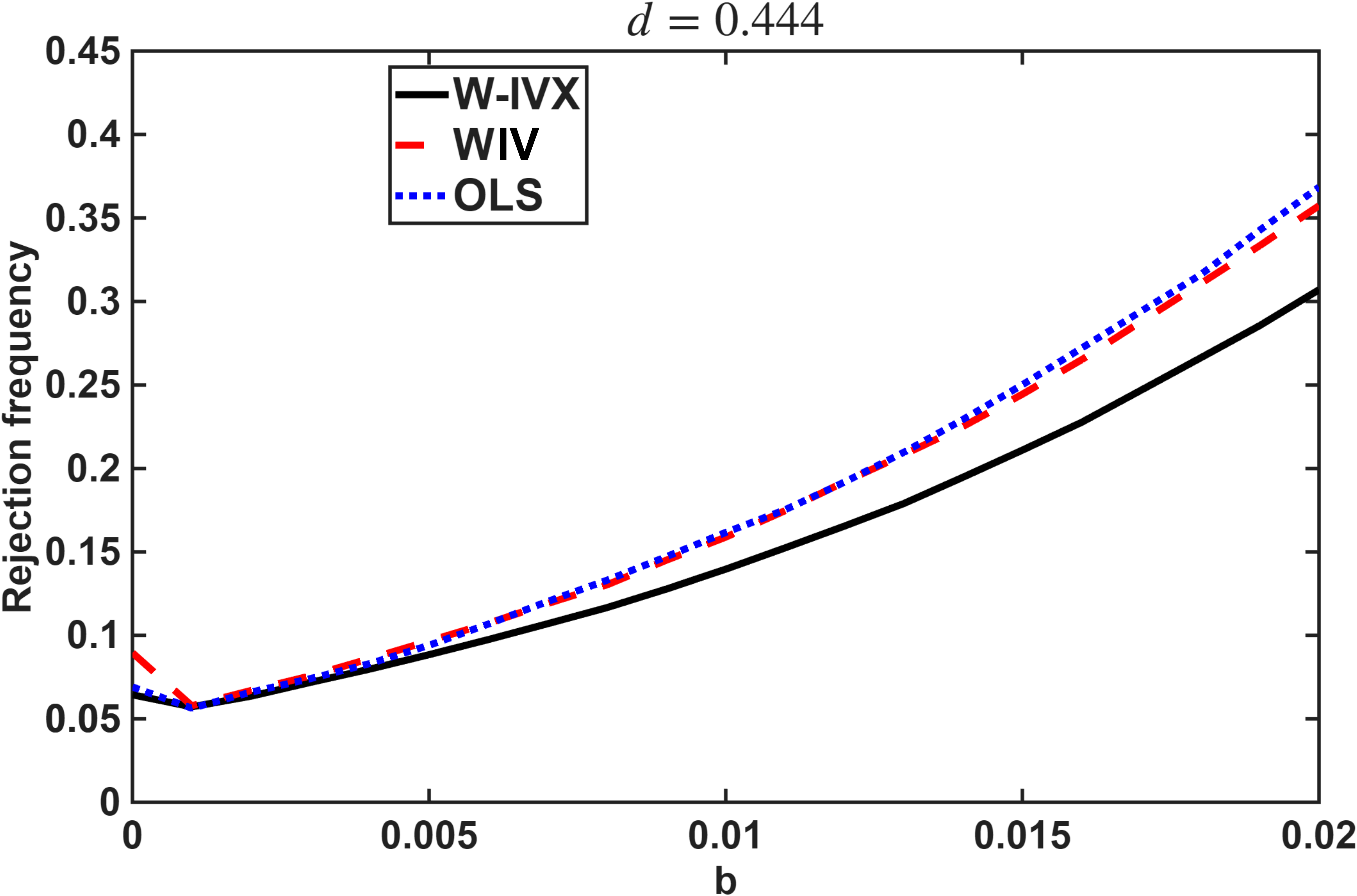}\\[0.6ex]
\includegraphics[width=0.49\linewidth]{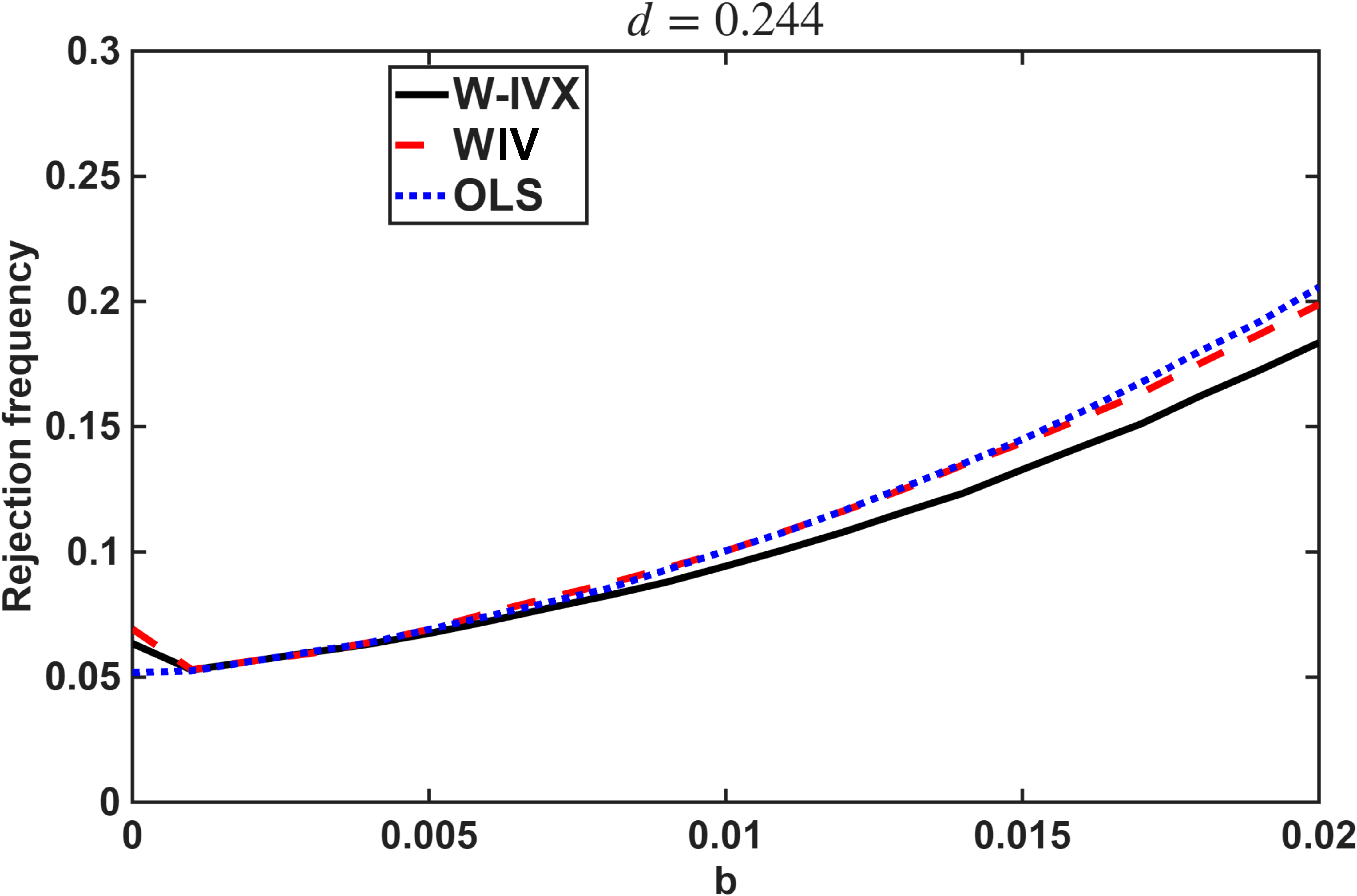}
\includegraphics[width=0.49\linewidth]{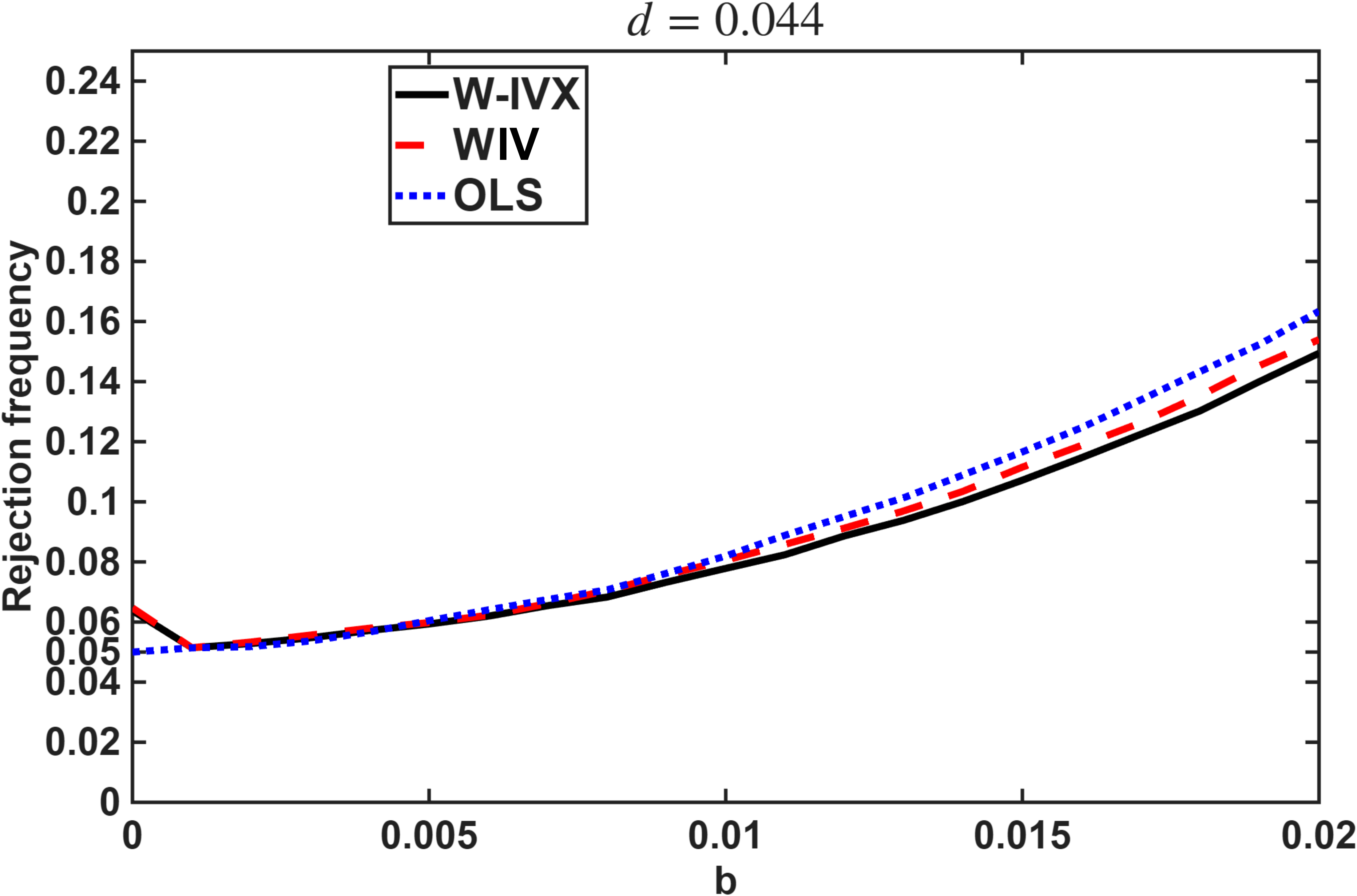}
\end{subfigure}

\end{figure}

\begin{figure}[p]
\ContinuedFloat \centering \captionsetup[subfigure]{position=top}

\setcounter{subfigure}{2}

\caption[]{Size Adjusted Power (nominal size 5\%): $[f(x),w(x)]=[x,x^{\langle0.75\rangle}]$,
$\delta=-0.95$, $n=1000$ (continued)}
\label{fig:sim-sap-linear-w075-continued}

\begin{subfigure}[t]{0.48\linewidth} \centering \caption{$\alpha=1.6$, $d_{*}=0.375$}
\label{fig:sim-sap-linear-w075-alpha-1p60}\includegraphics[width=0.49\linewidth]{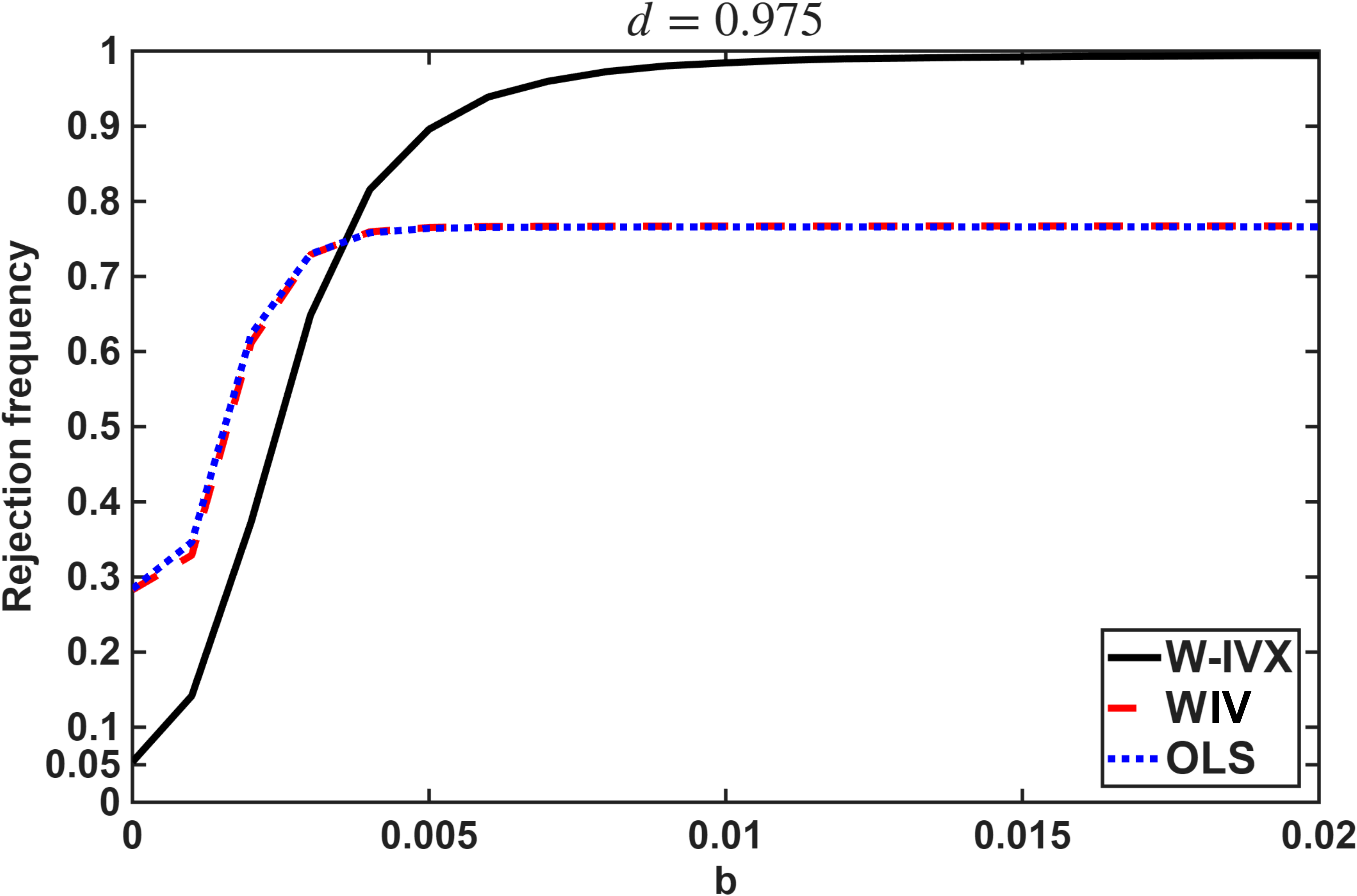}
\includegraphics[width=0.49\linewidth]{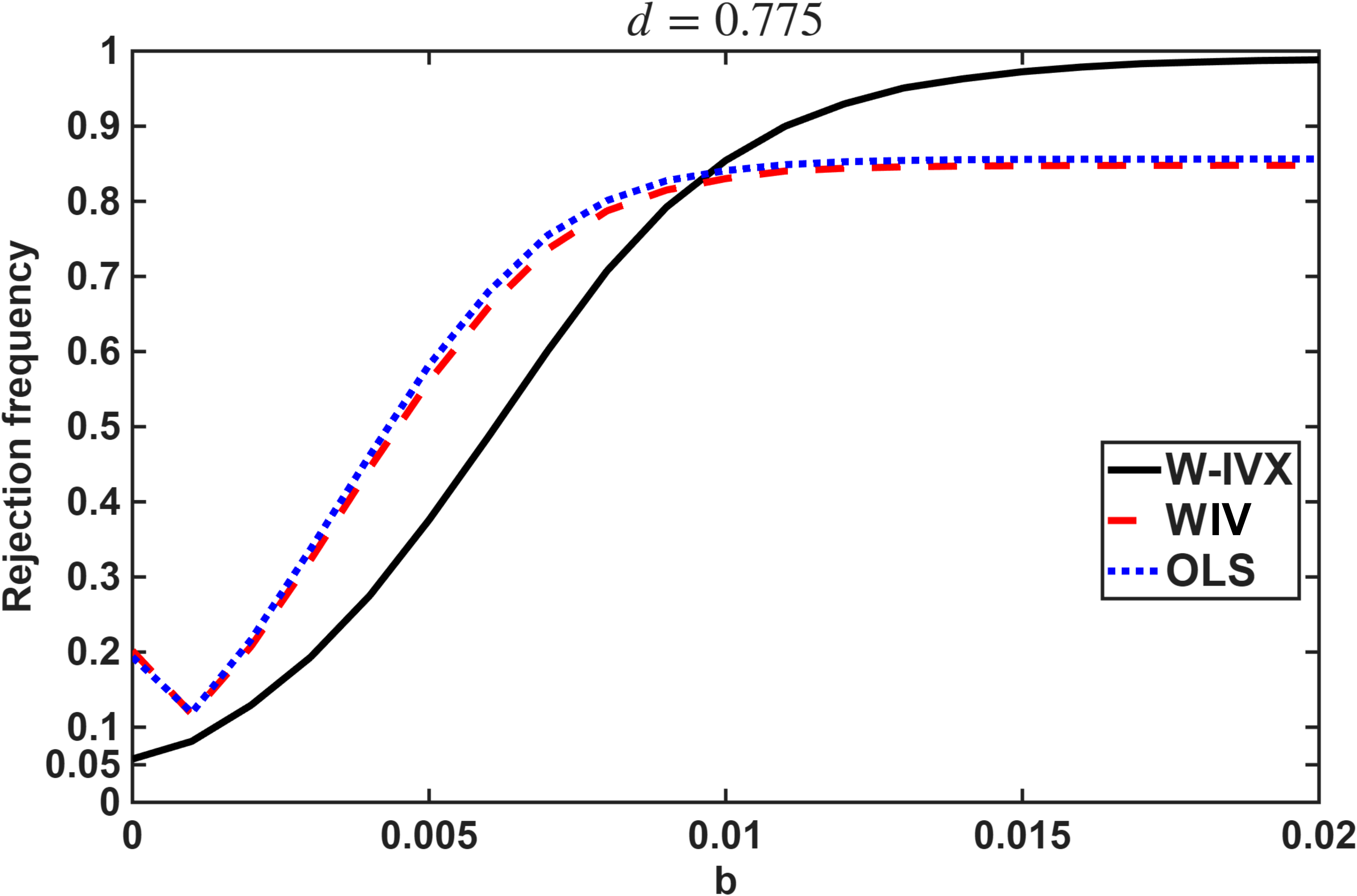}\\[0.6ex]
\includegraphics[width=0.49\linewidth]{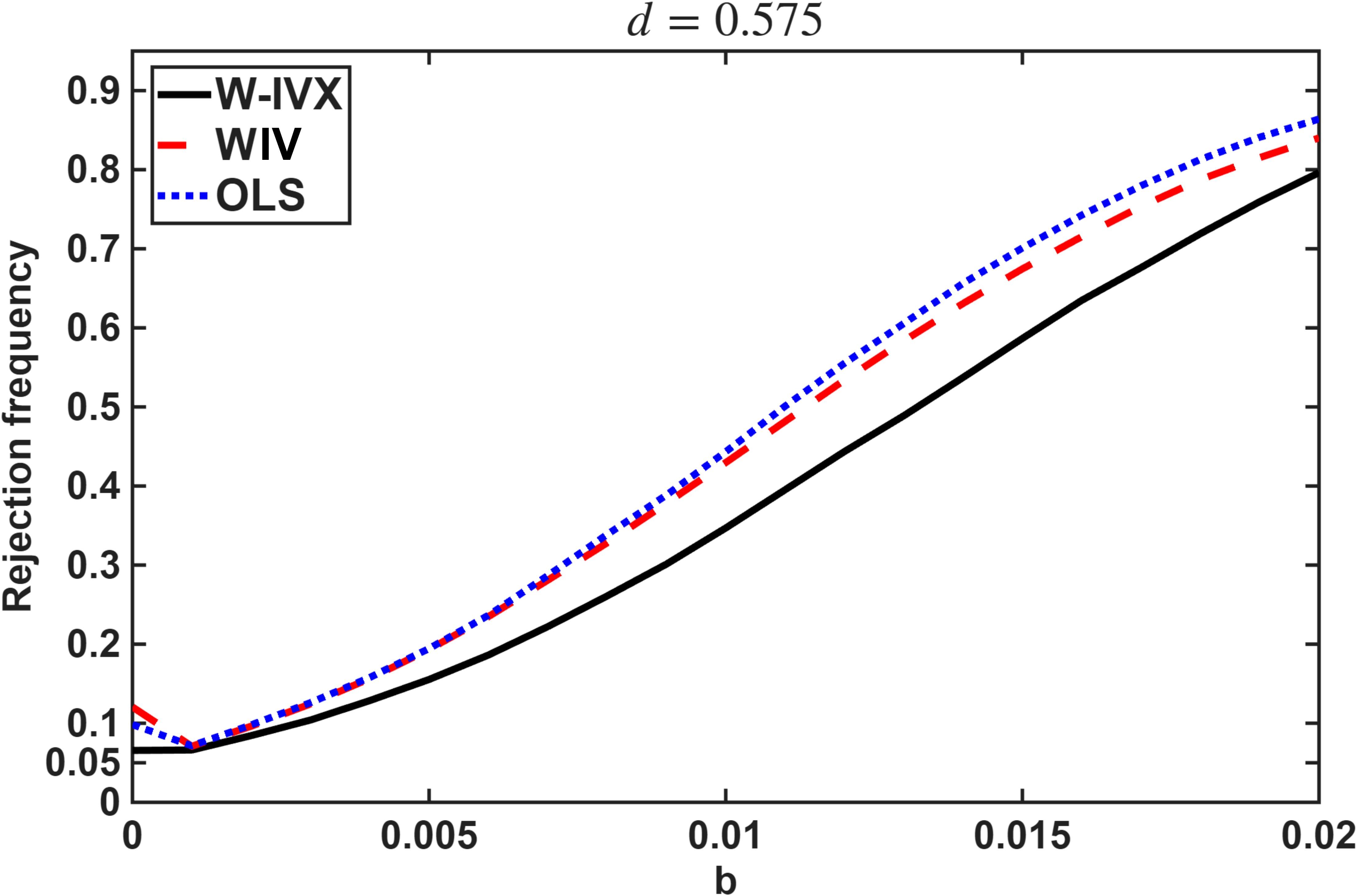}
\includegraphics[width=0.49\linewidth]{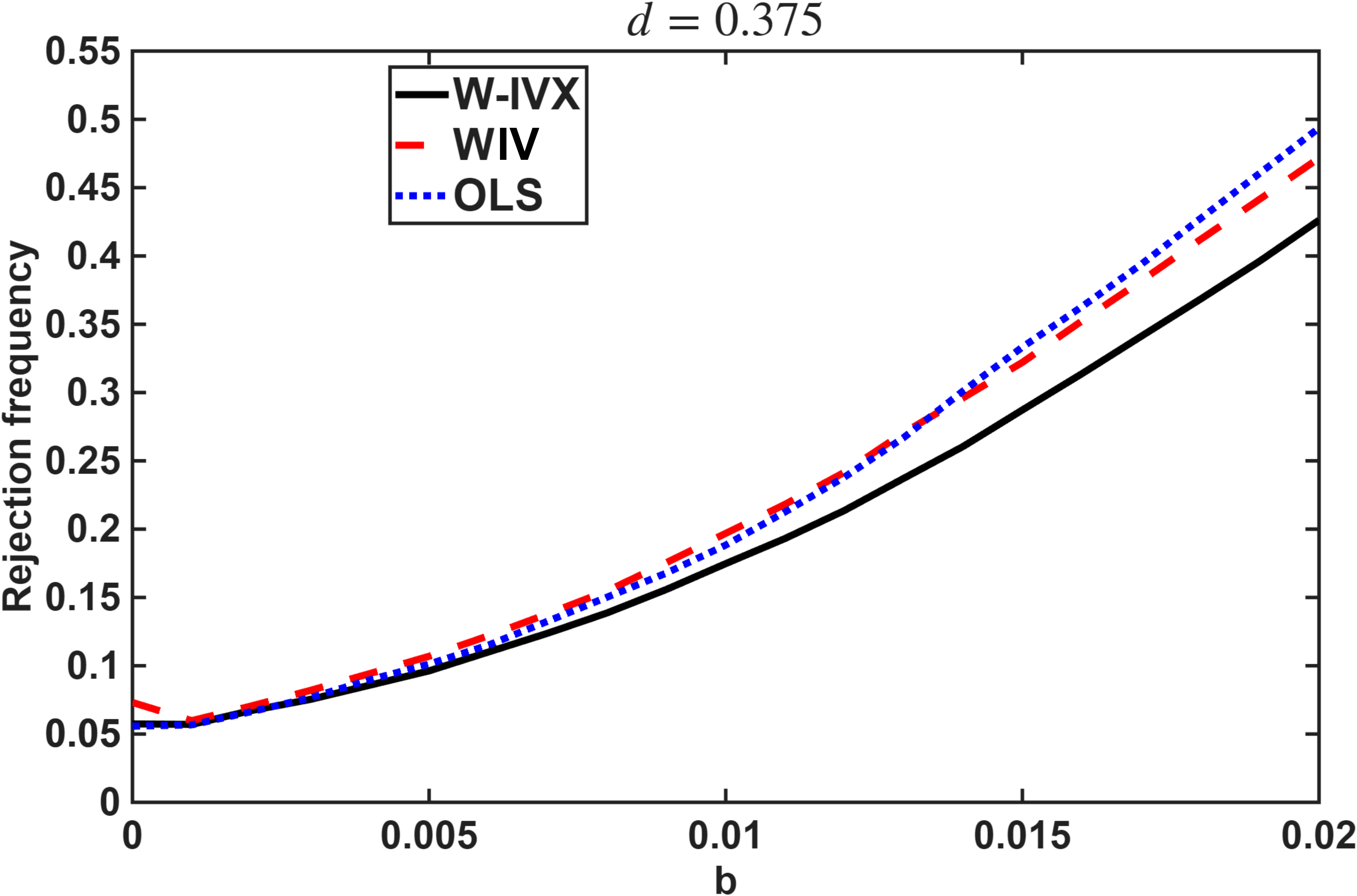}\\[0.6ex]
\includegraphics[width=0.49\linewidth]{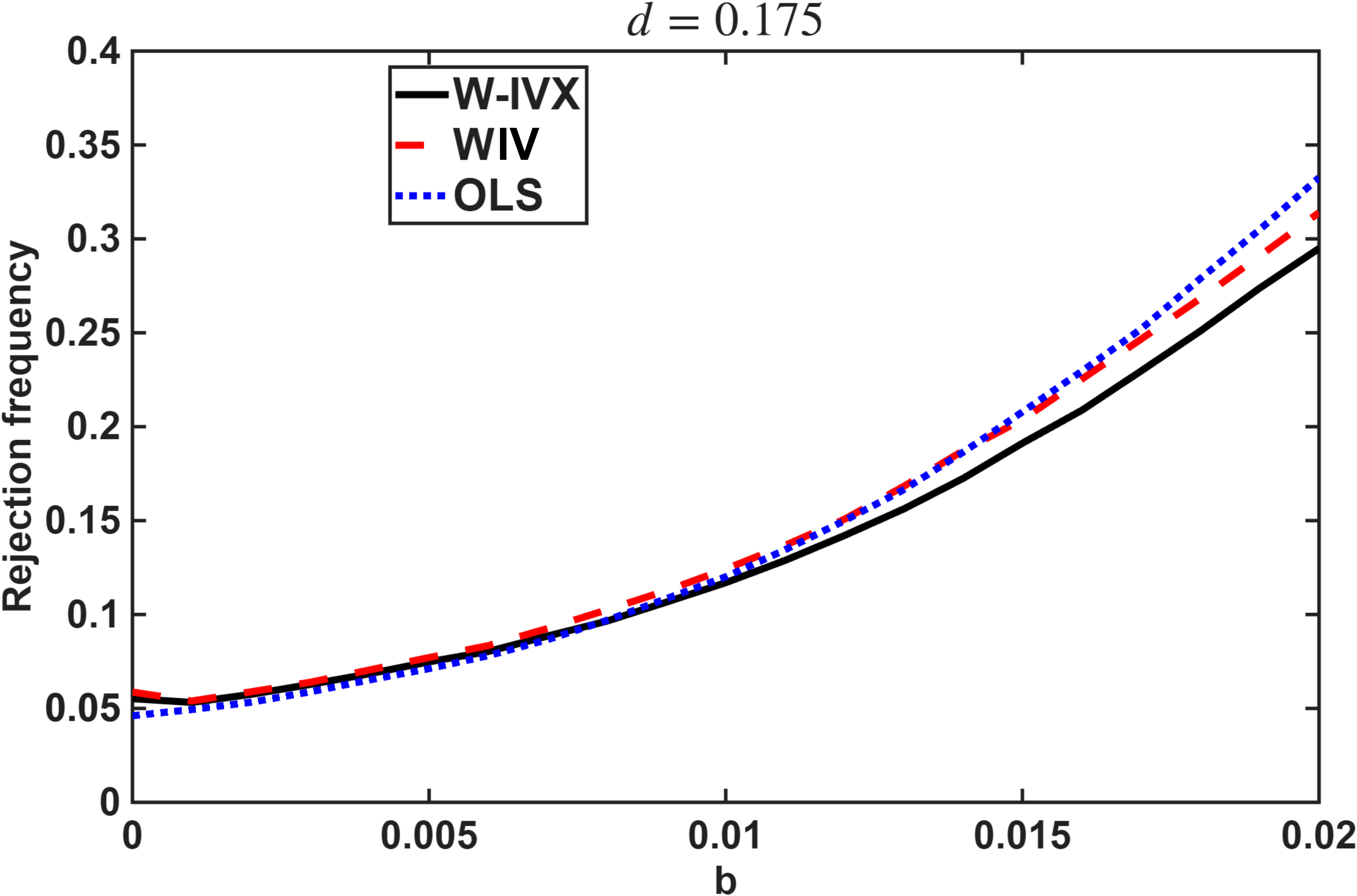}
\includegraphics[width=0.49\linewidth]{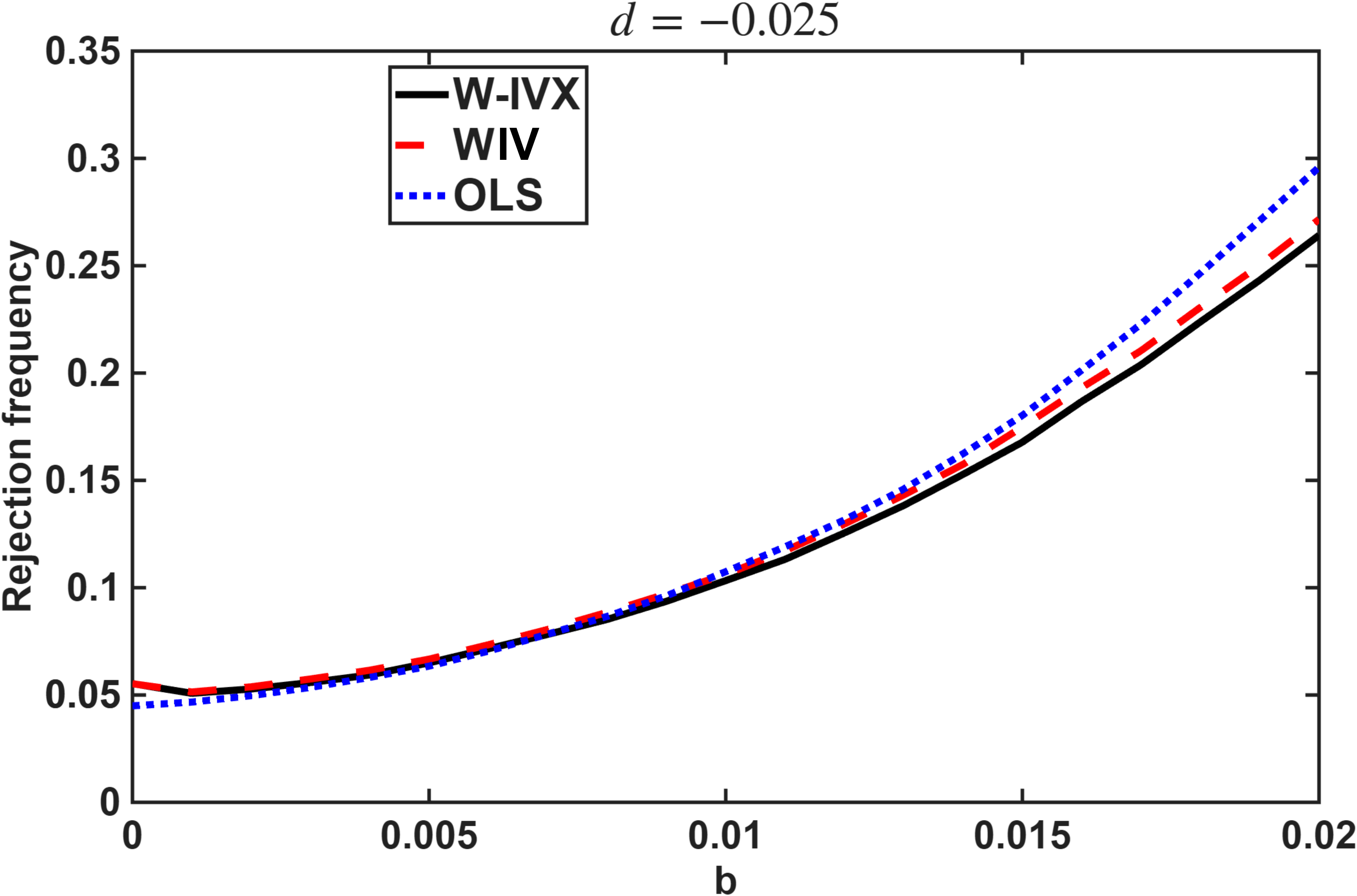}
\end{subfigure} \hfill{}
\begin{subfigure}[t]{0.48\linewidth} \centering \caption{$\alpha=1.4$, $d_{*}=0.286$}
\label{fig:sim-sap-linear-w075-alpha-1p40}\includegraphics[width=0.49\linewidth]{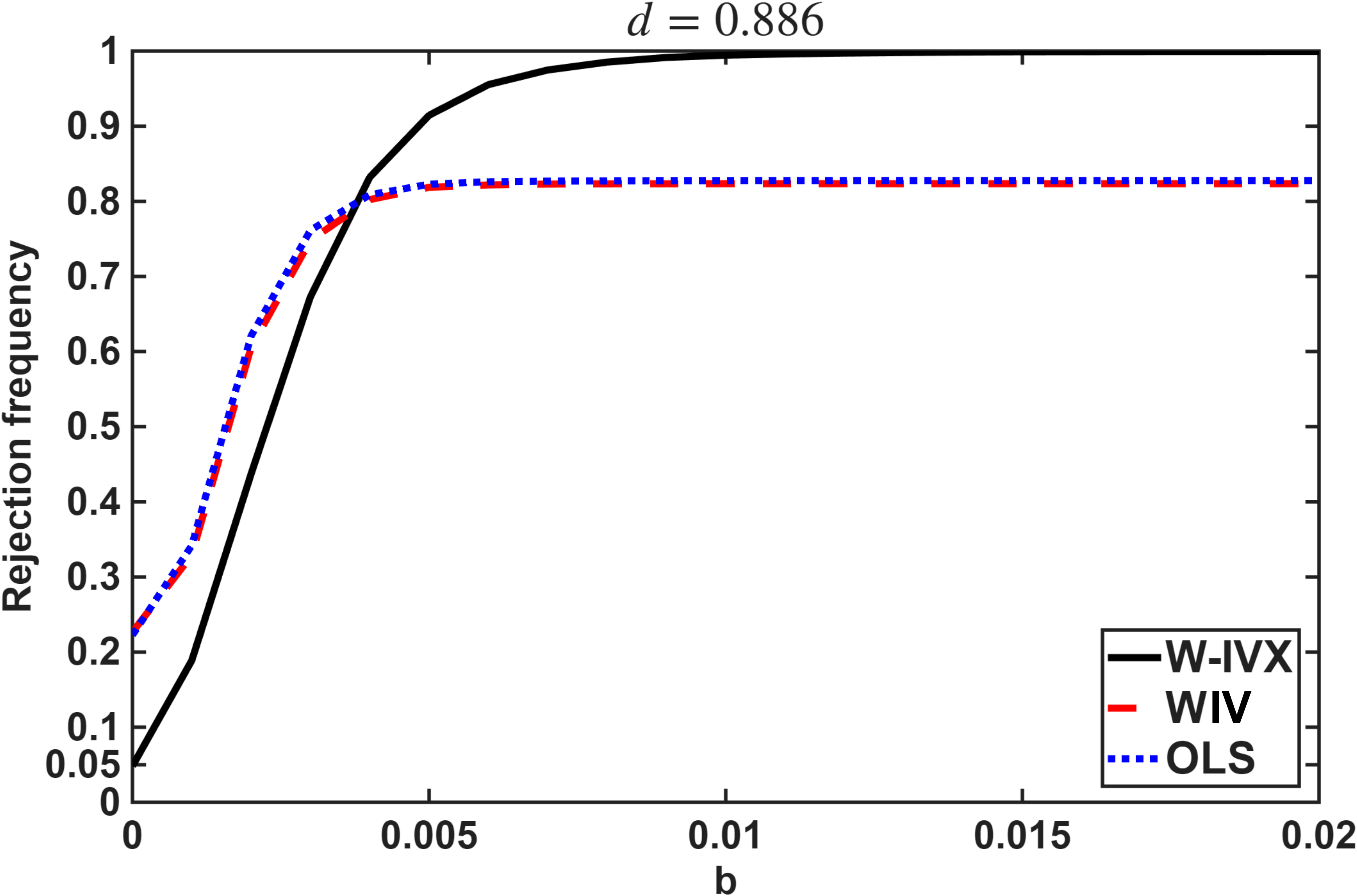}
\includegraphics[width=0.49\linewidth]{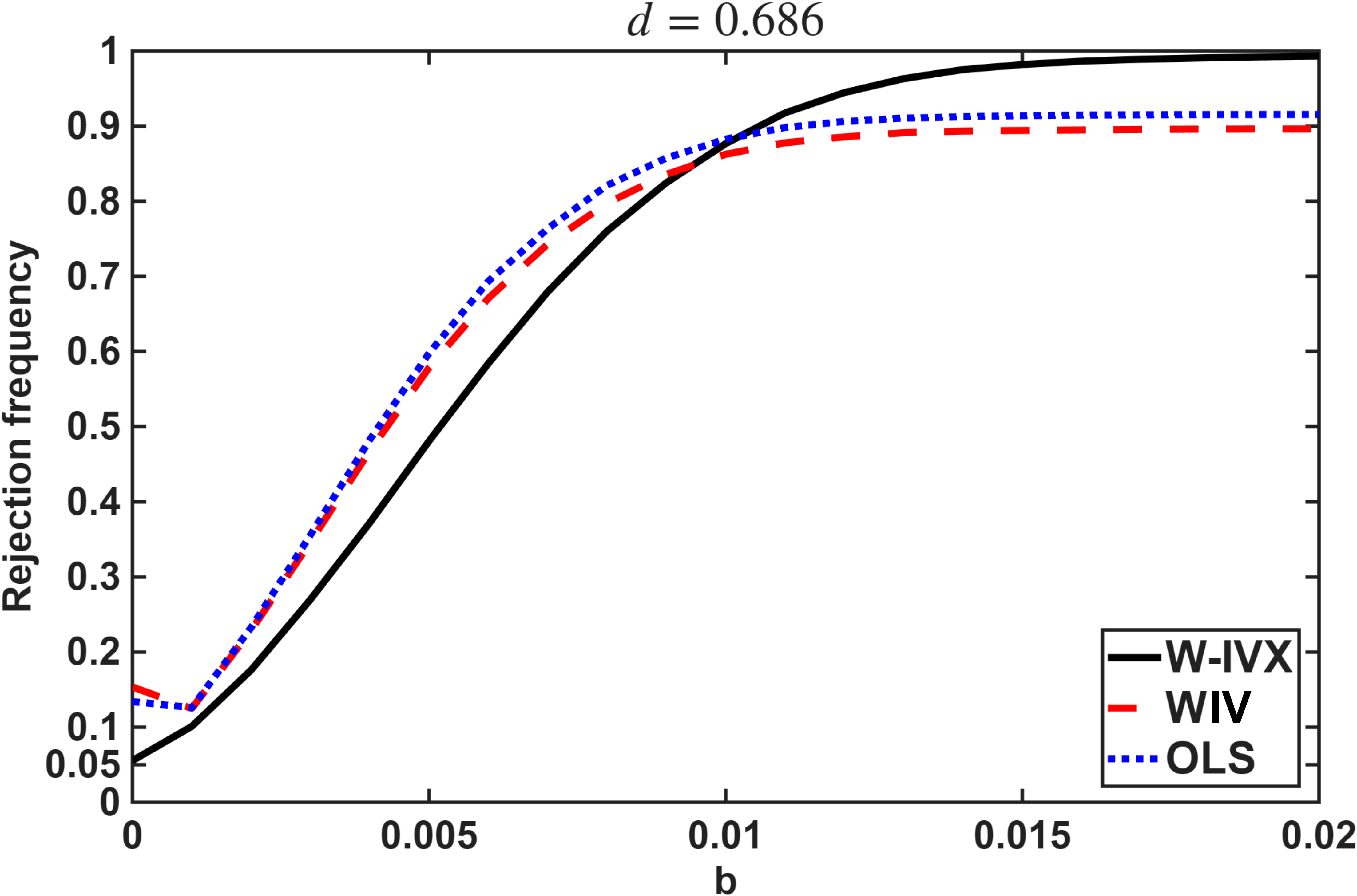}\\[0.6ex]
\includegraphics[width=0.49\linewidth]{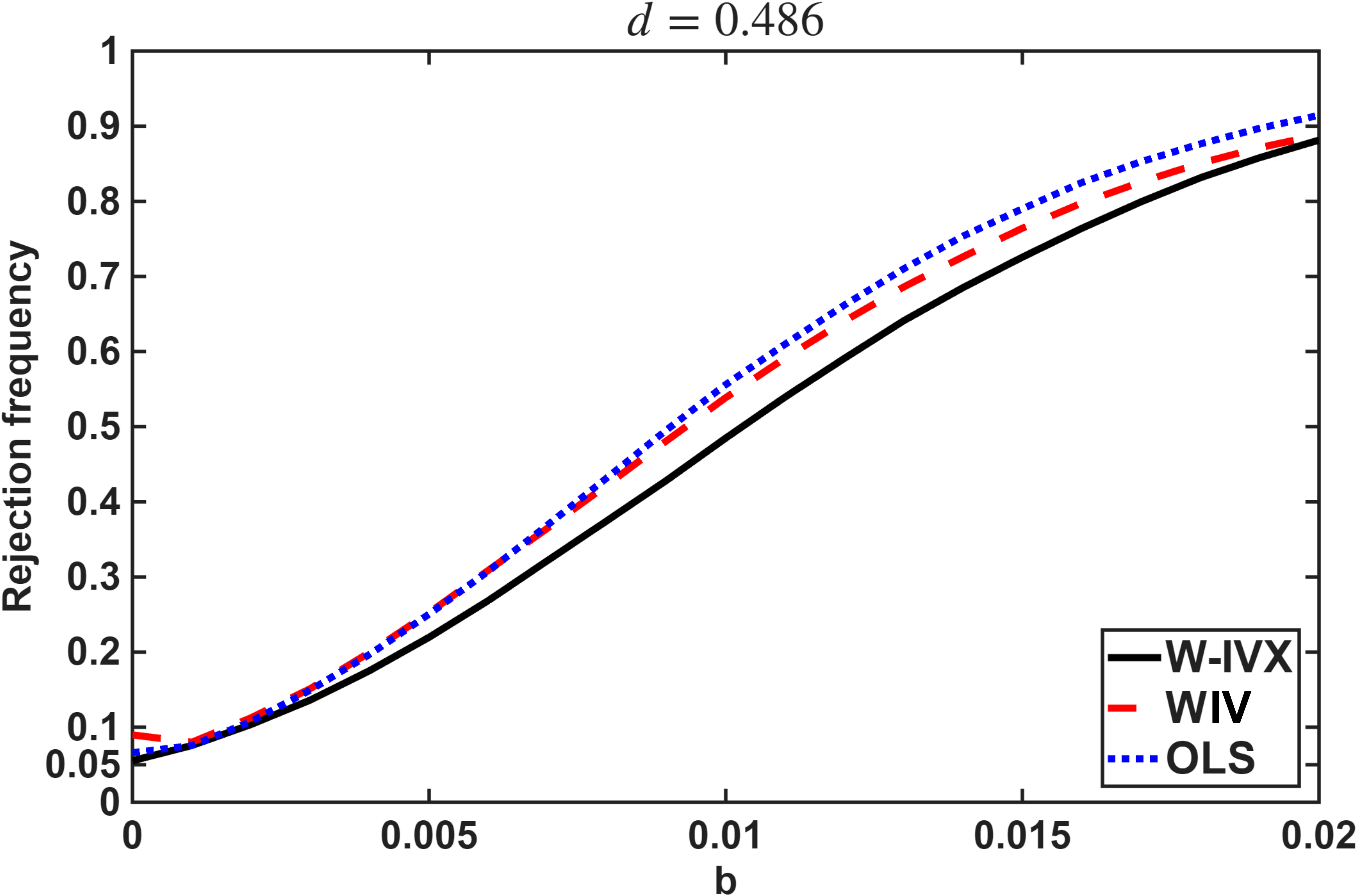}
\includegraphics[width=0.49\linewidth]{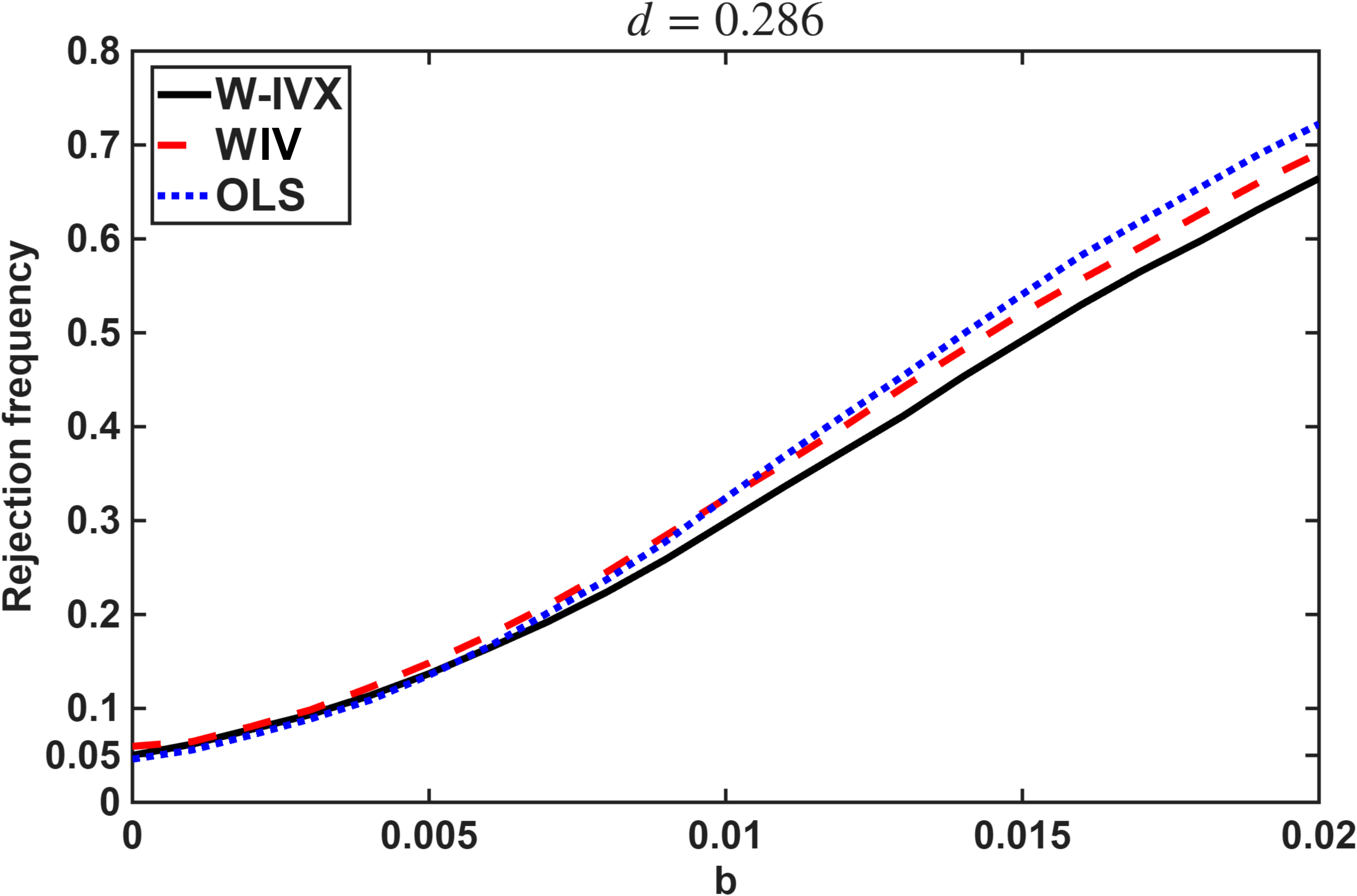}\\[0.6ex]
\includegraphics[width=0.49\linewidth]{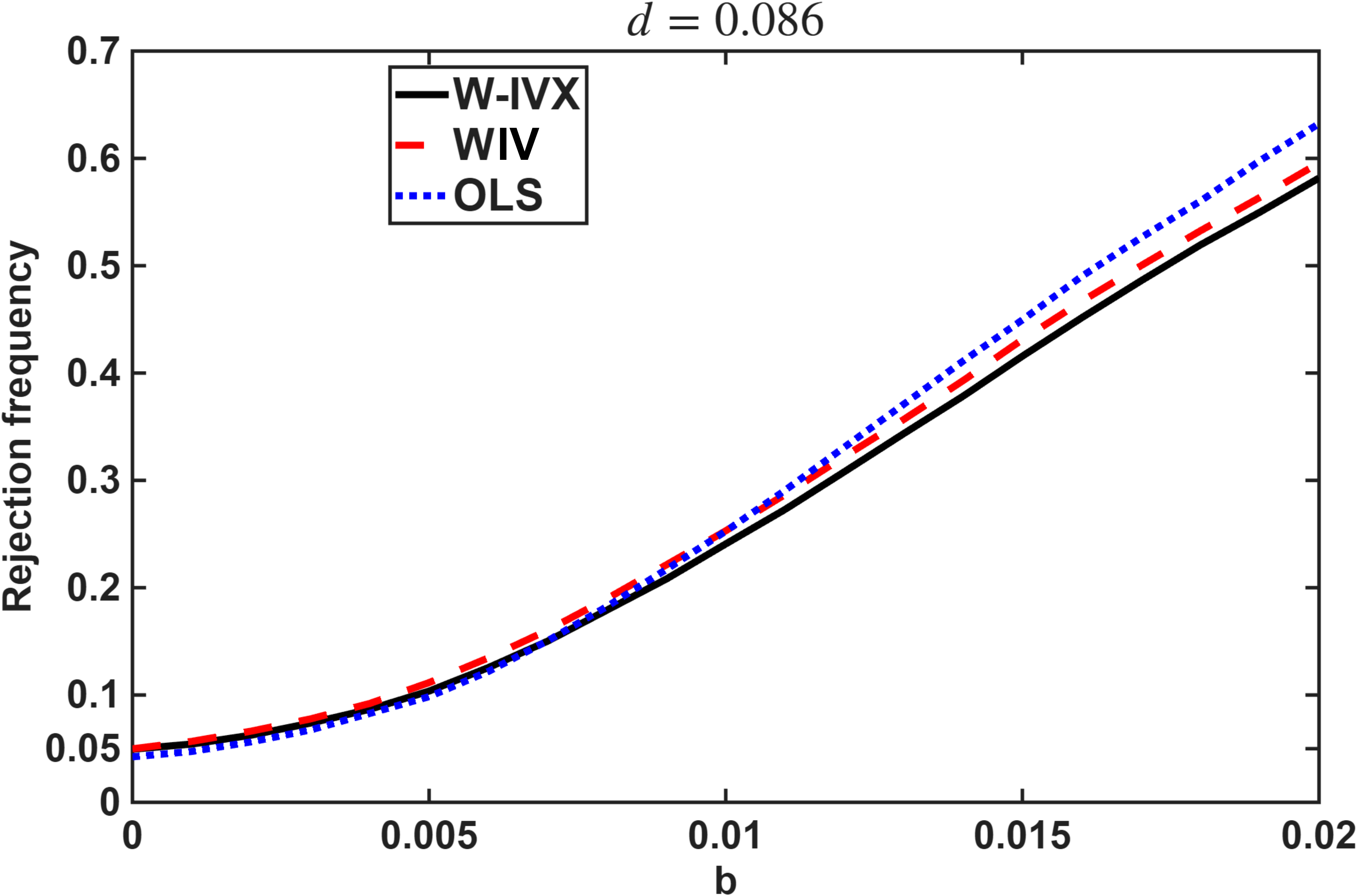}
\includegraphics[width=0.49\linewidth]{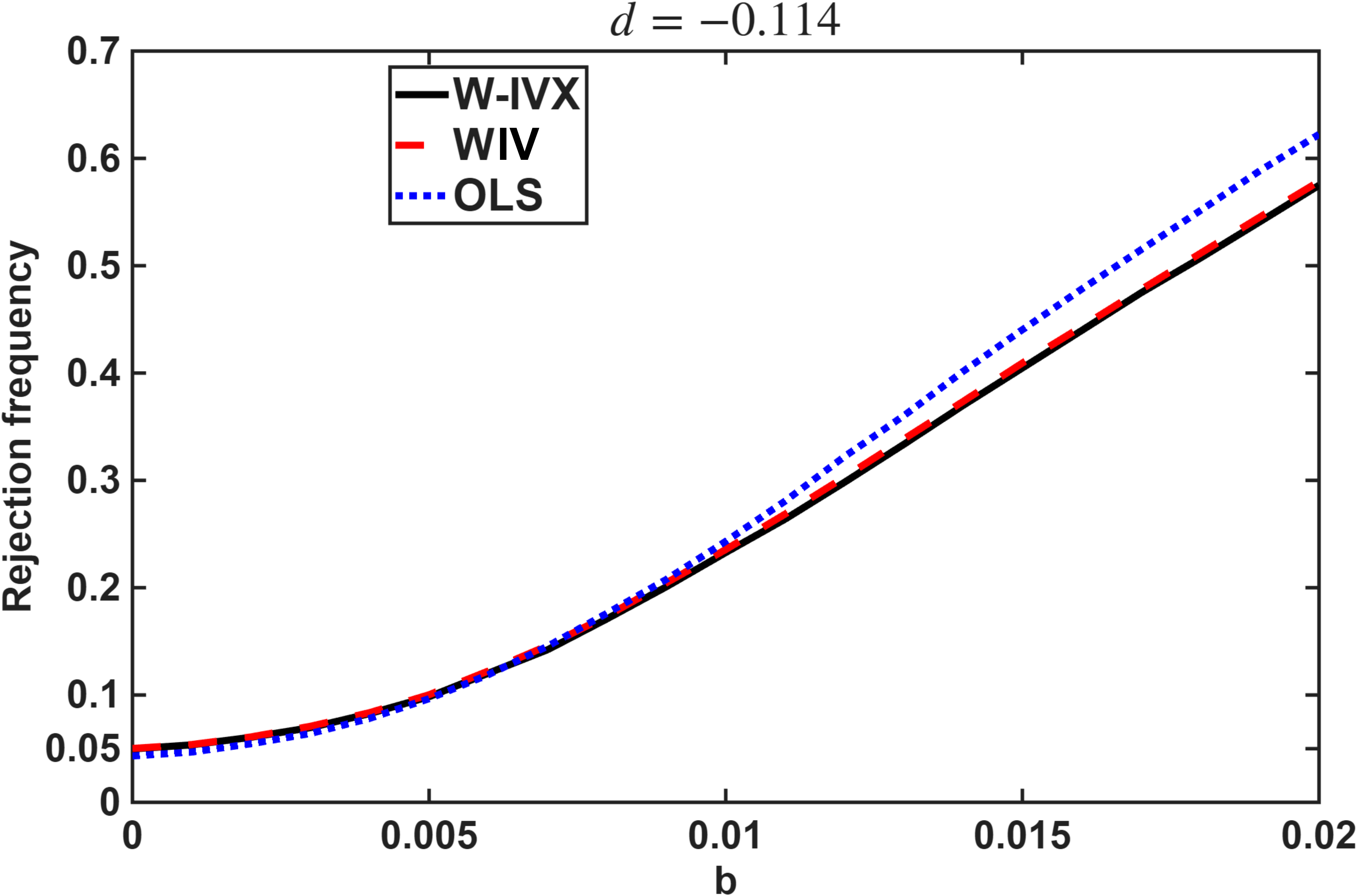}
\end{subfigure}
\end{figure}

\end{landscape} \clearpage{}

\begin{landscape}

\begin{figure}[p]
\caption{Size Adjusted Power (nominal size 5\%): \textbf{DGP1}, $[f(x),w(x)]=[x,x]$,
$\delta=-0.95$, $n=1000$}
\label{fig:sim-sap-linear-w10}\centering \captionsetup[subfigure]{position=top}

\begin{subfigure}[t]{0.48\linewidth} \centering \caption{$\alpha=2$, $d_{*}=0.5$}
\label{fig:sim-sap-linear-w10-alpha-2p00} \includegraphics[width=0.49\linewidth]{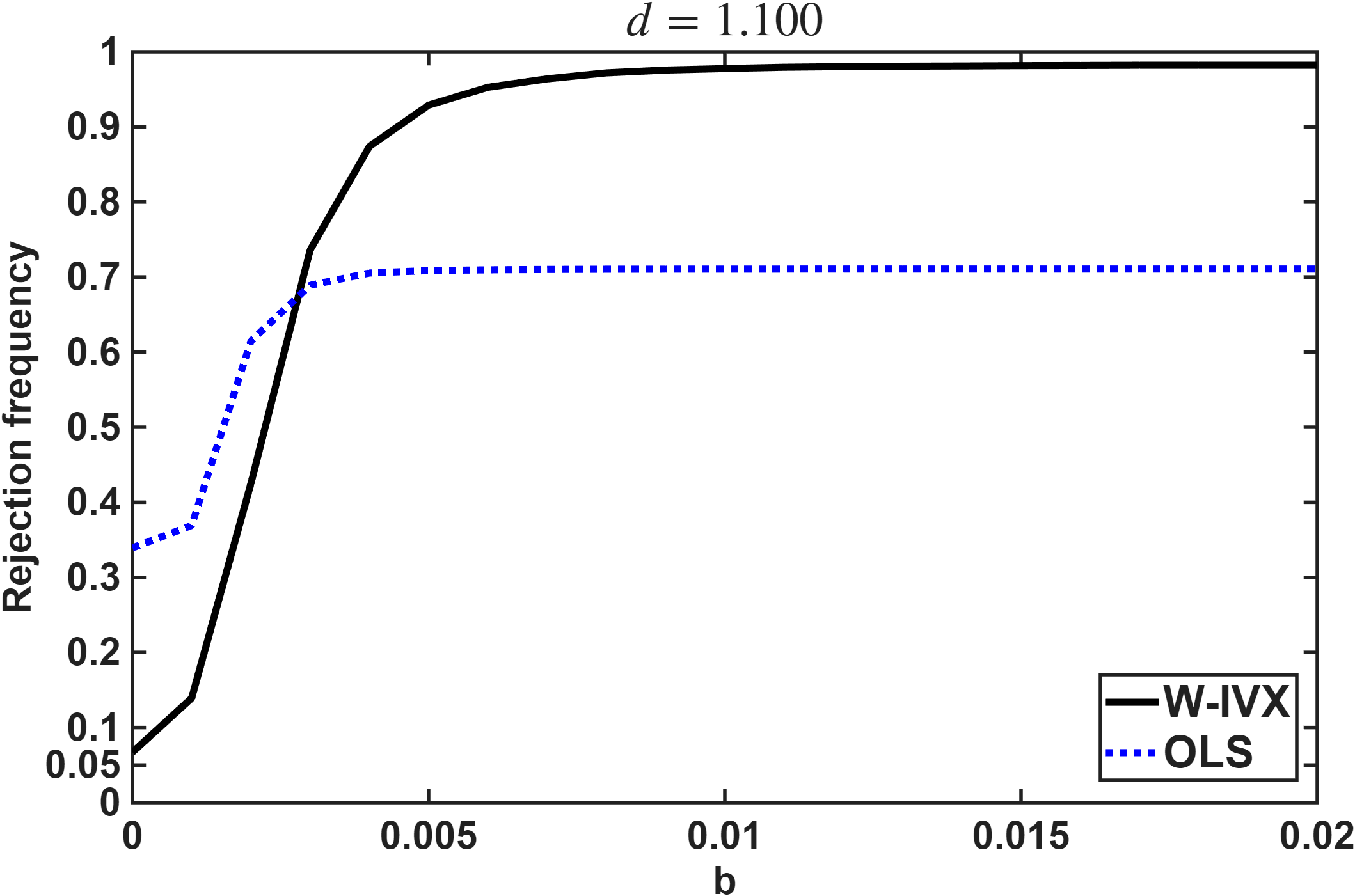}
\includegraphics[width=0.49\linewidth]{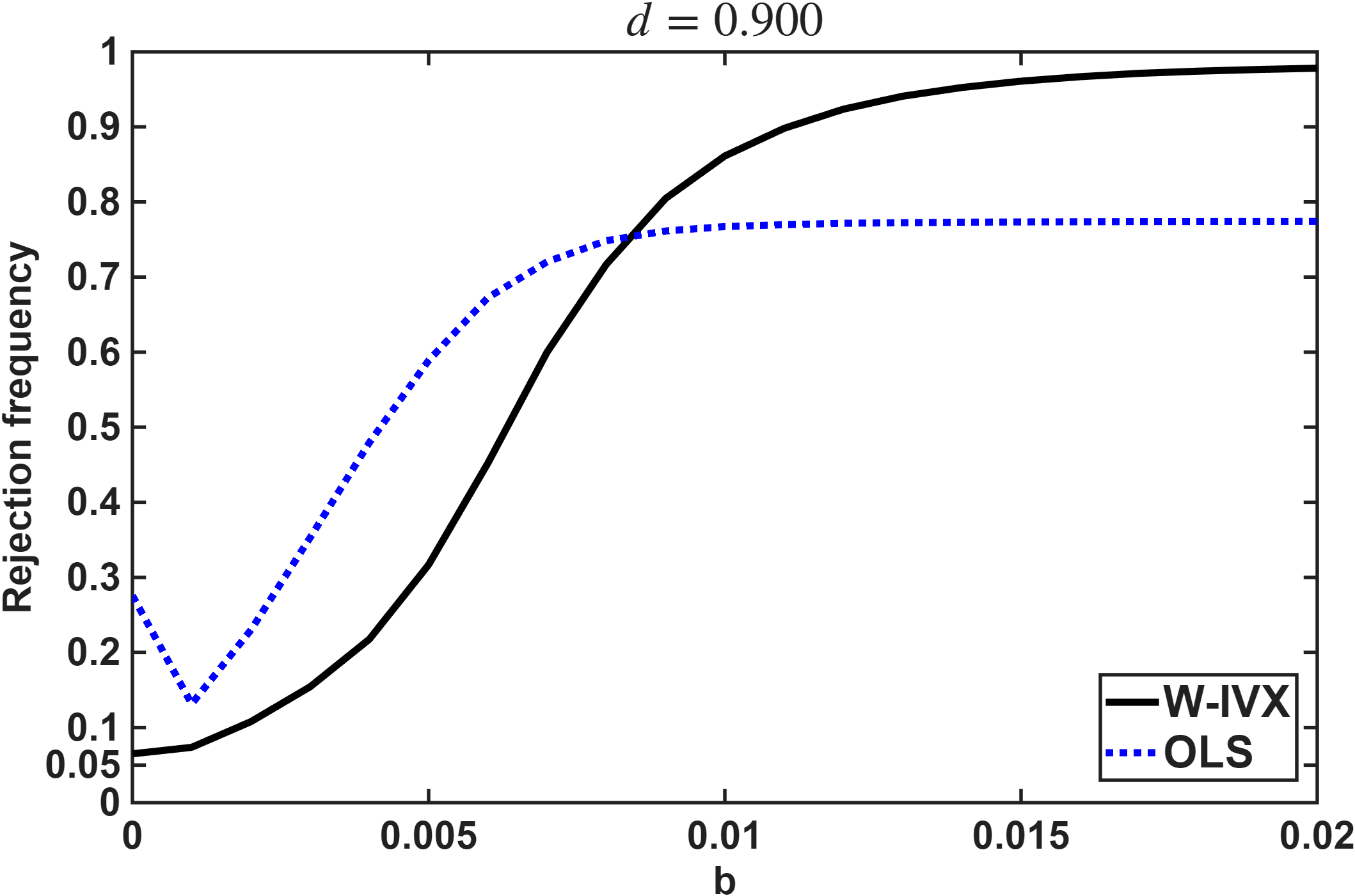}\\[0.6ex]
\includegraphics[width=0.49\linewidth]{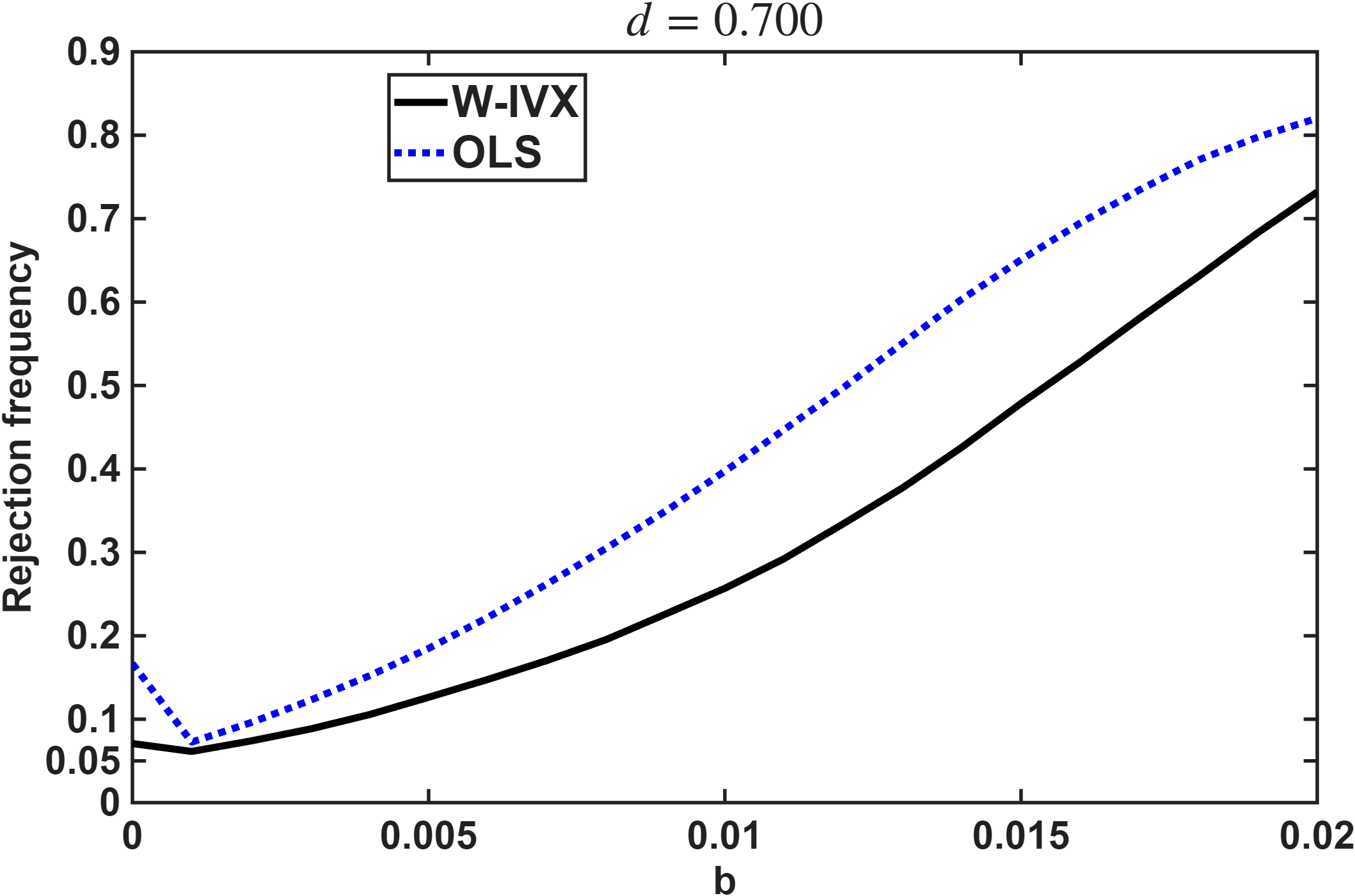}
\includegraphics[width=0.49\linewidth]{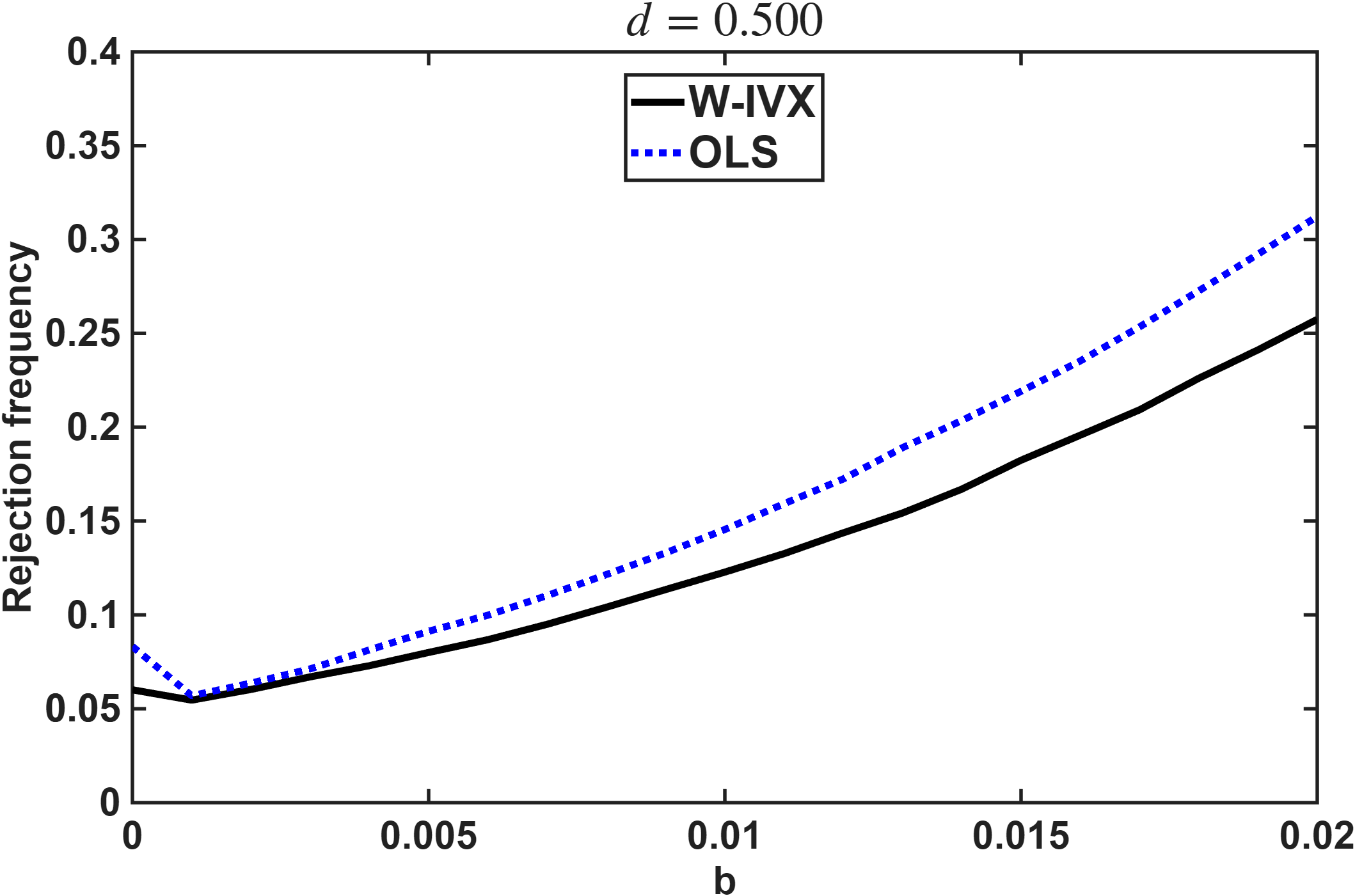}\\[0.6ex]
\includegraphics[width=0.49\linewidth]{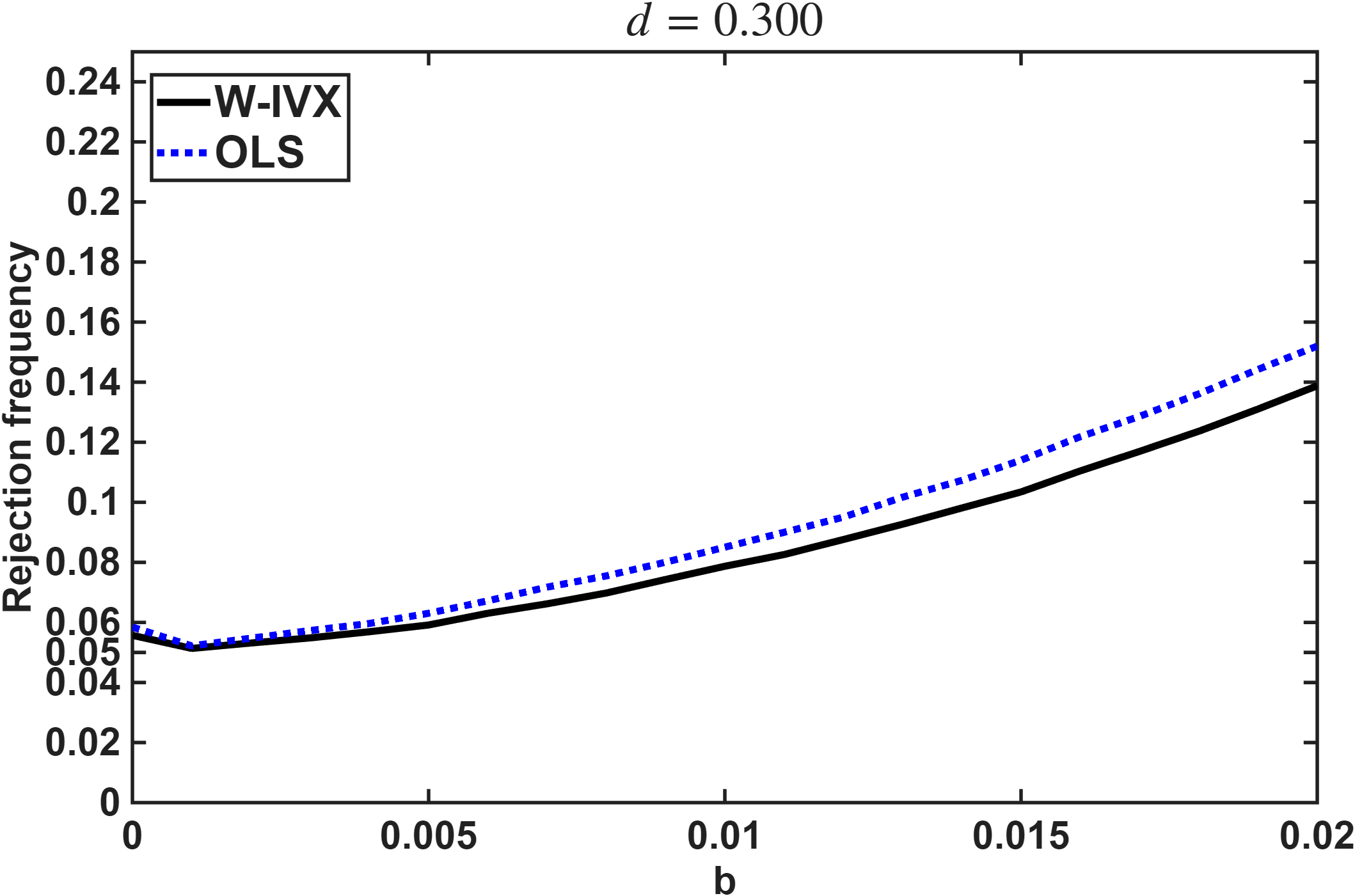}
\includegraphics[width=0.49\linewidth]{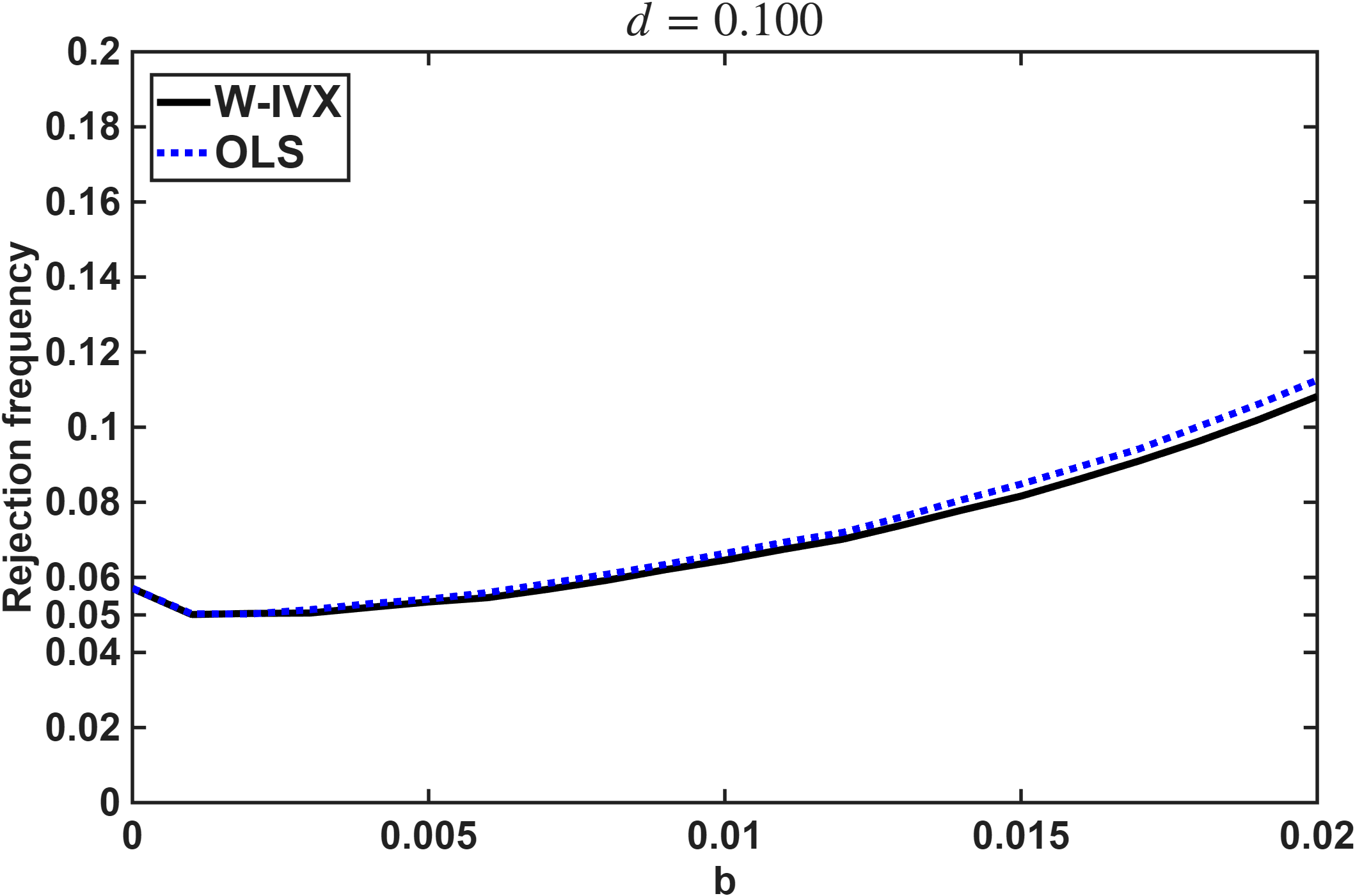}
\end{subfigure} \hfill{}
\begin{subfigure}[t]{0.48\linewidth} \centering \caption{$\alpha=1.8$, $d_{*}=0.444$}
\label{fig:sim-sap-linear-w10-alpha-1p80} \includegraphics[width=0.49\linewidth]{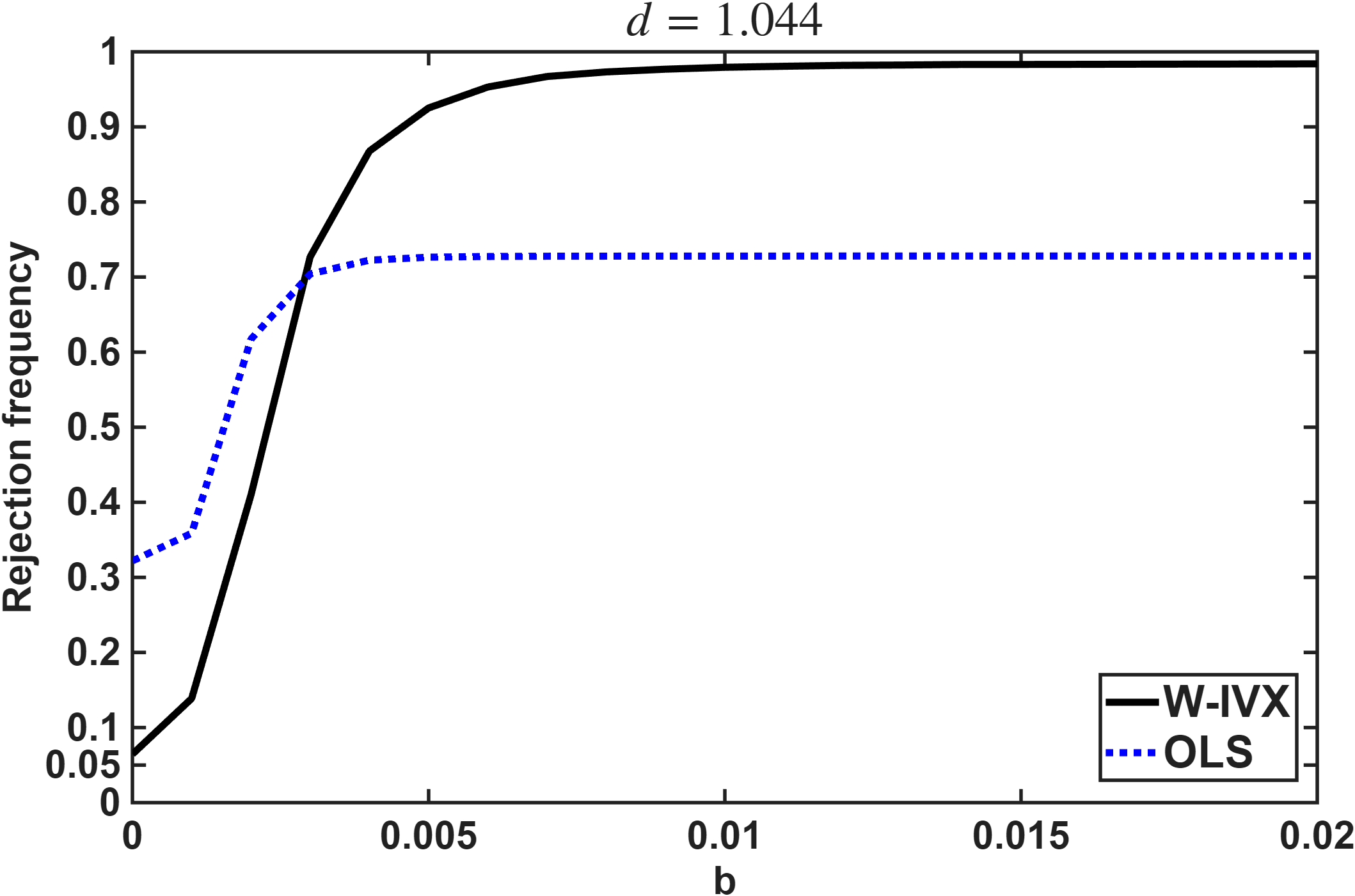}
\includegraphics[width=0.49\linewidth]{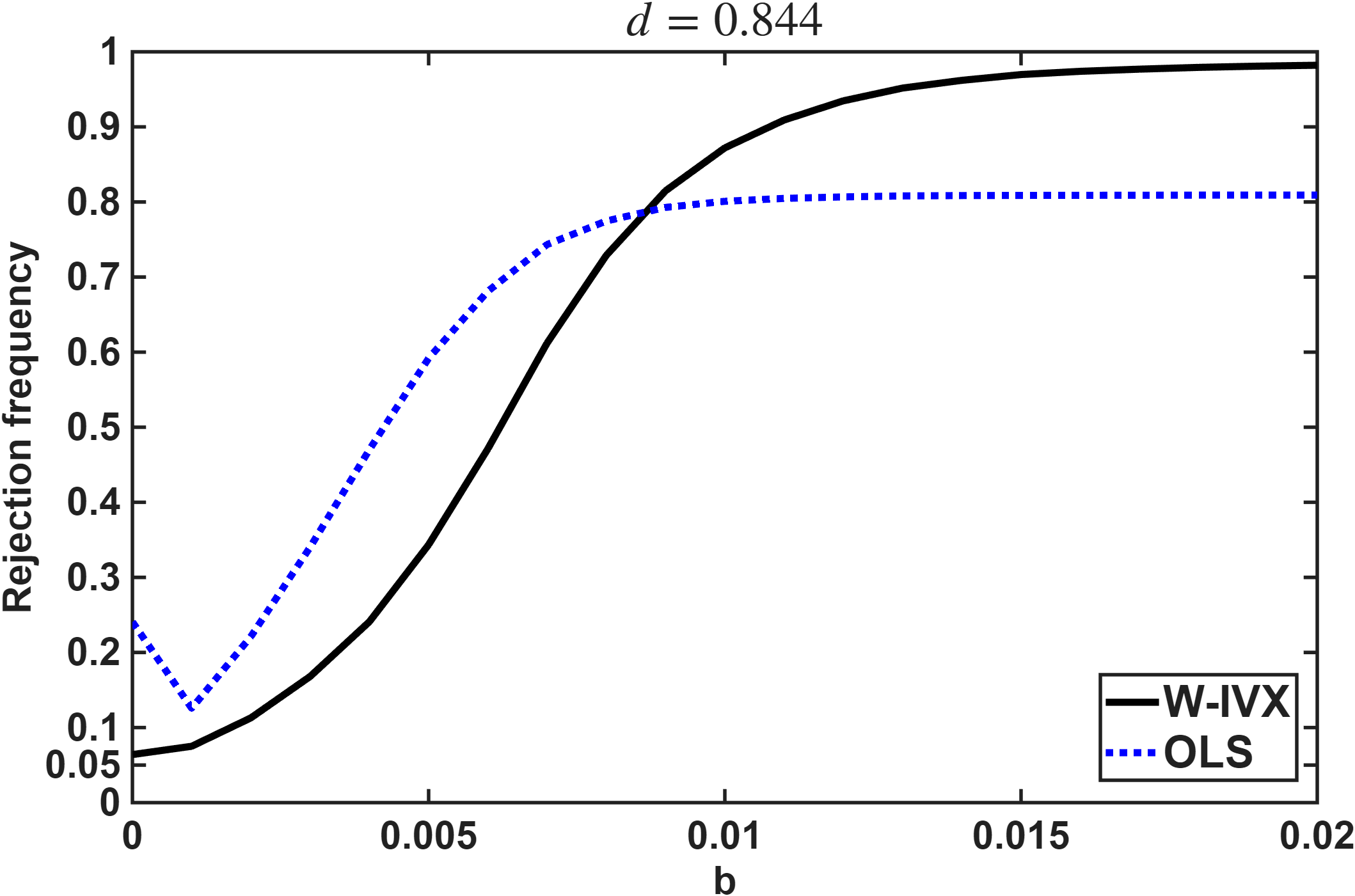}\\[0.6ex]
\includegraphics[width=0.49\linewidth]{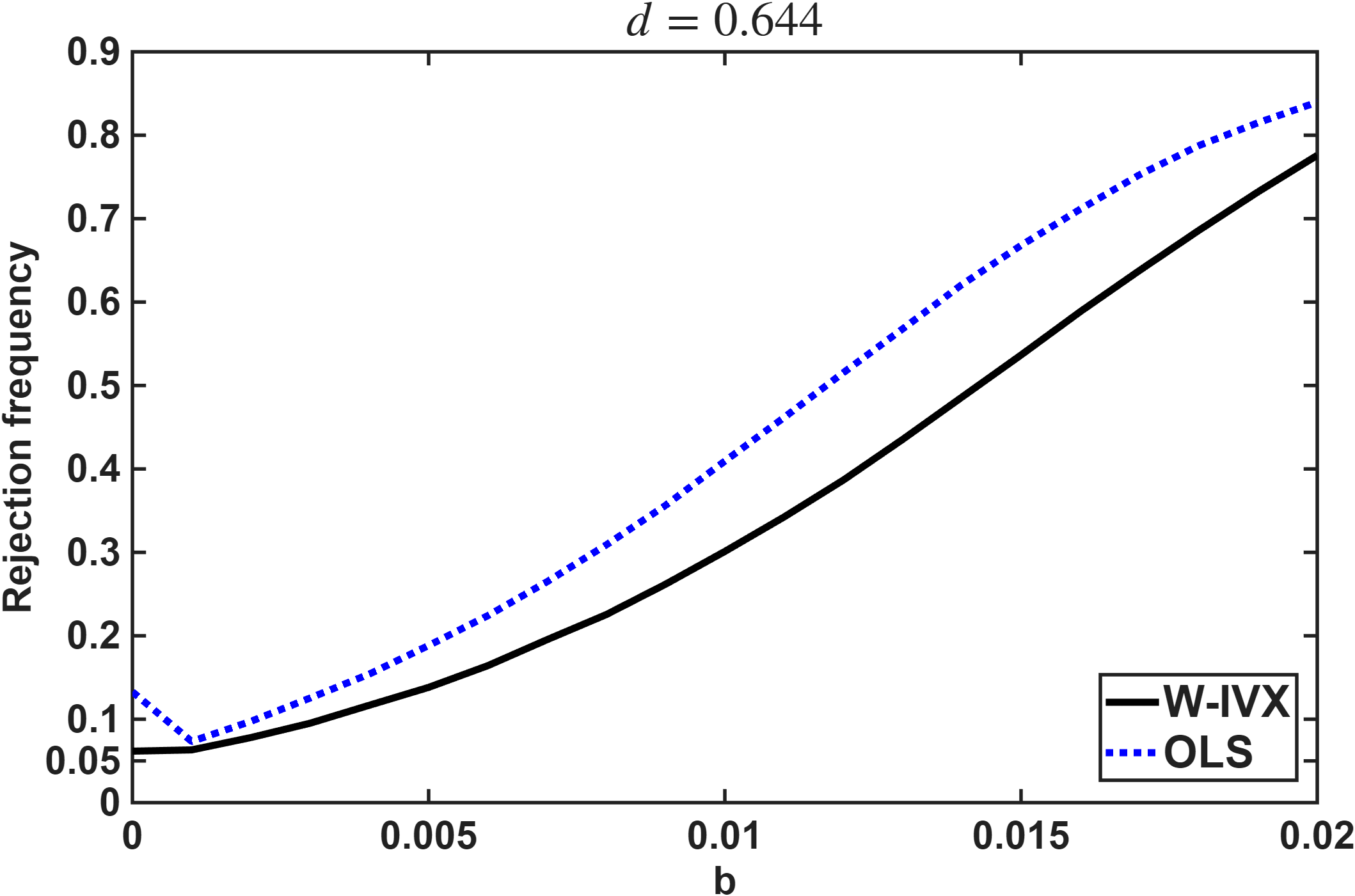}
\includegraphics[width=0.49\linewidth]{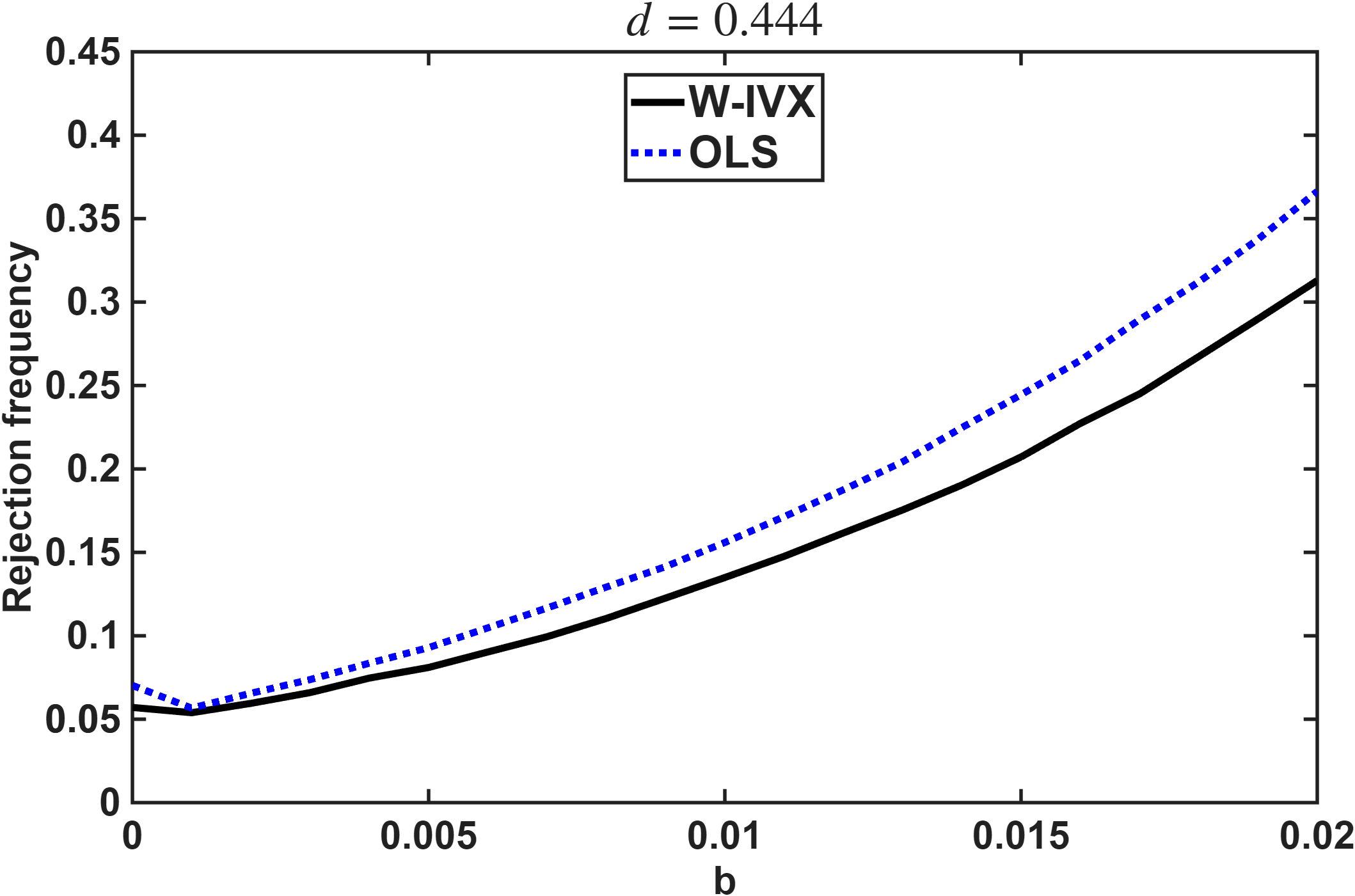}\\[0.6ex]
\includegraphics[width=0.49\linewidth]{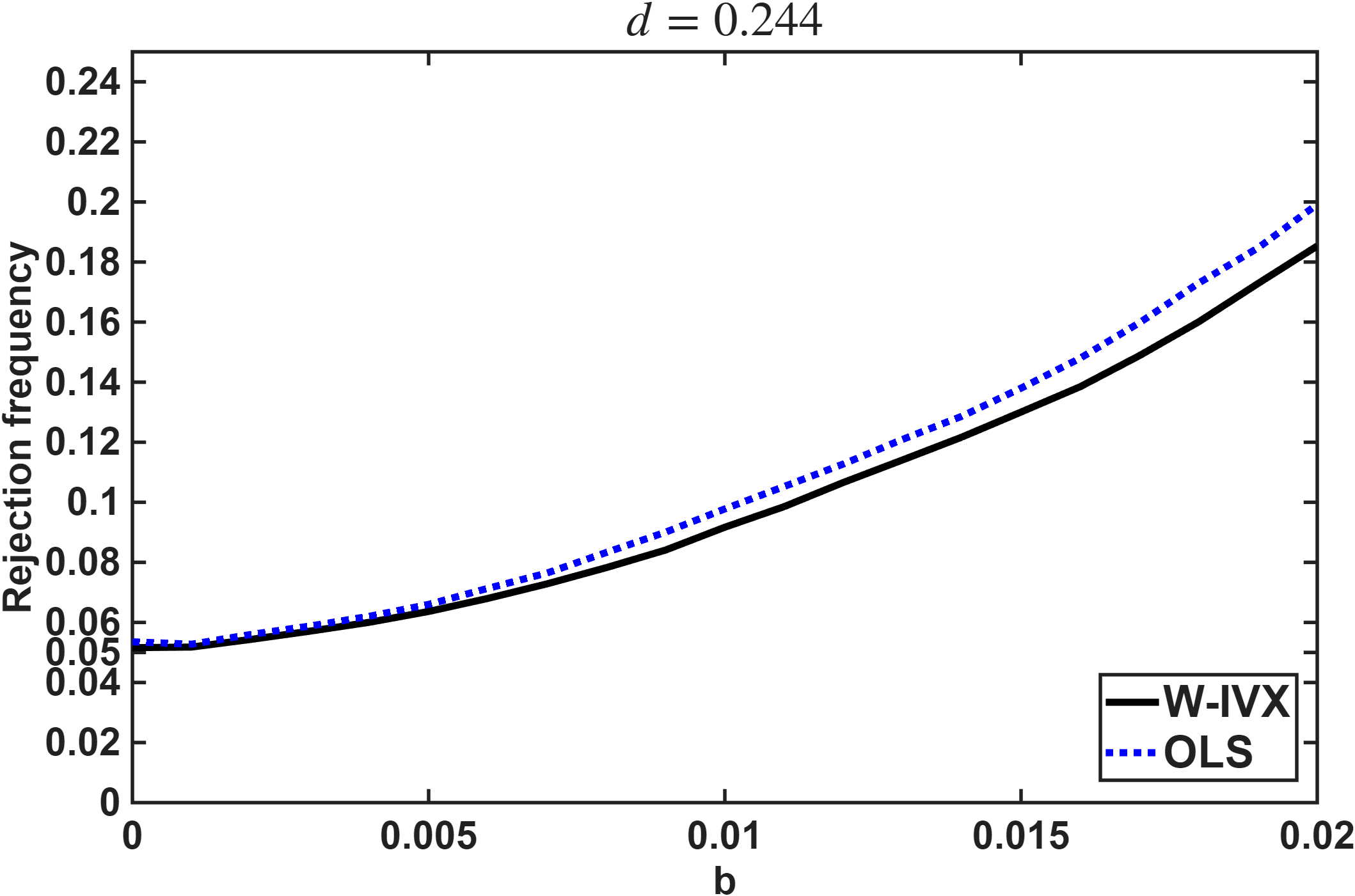}
\includegraphics[width=0.49\linewidth]{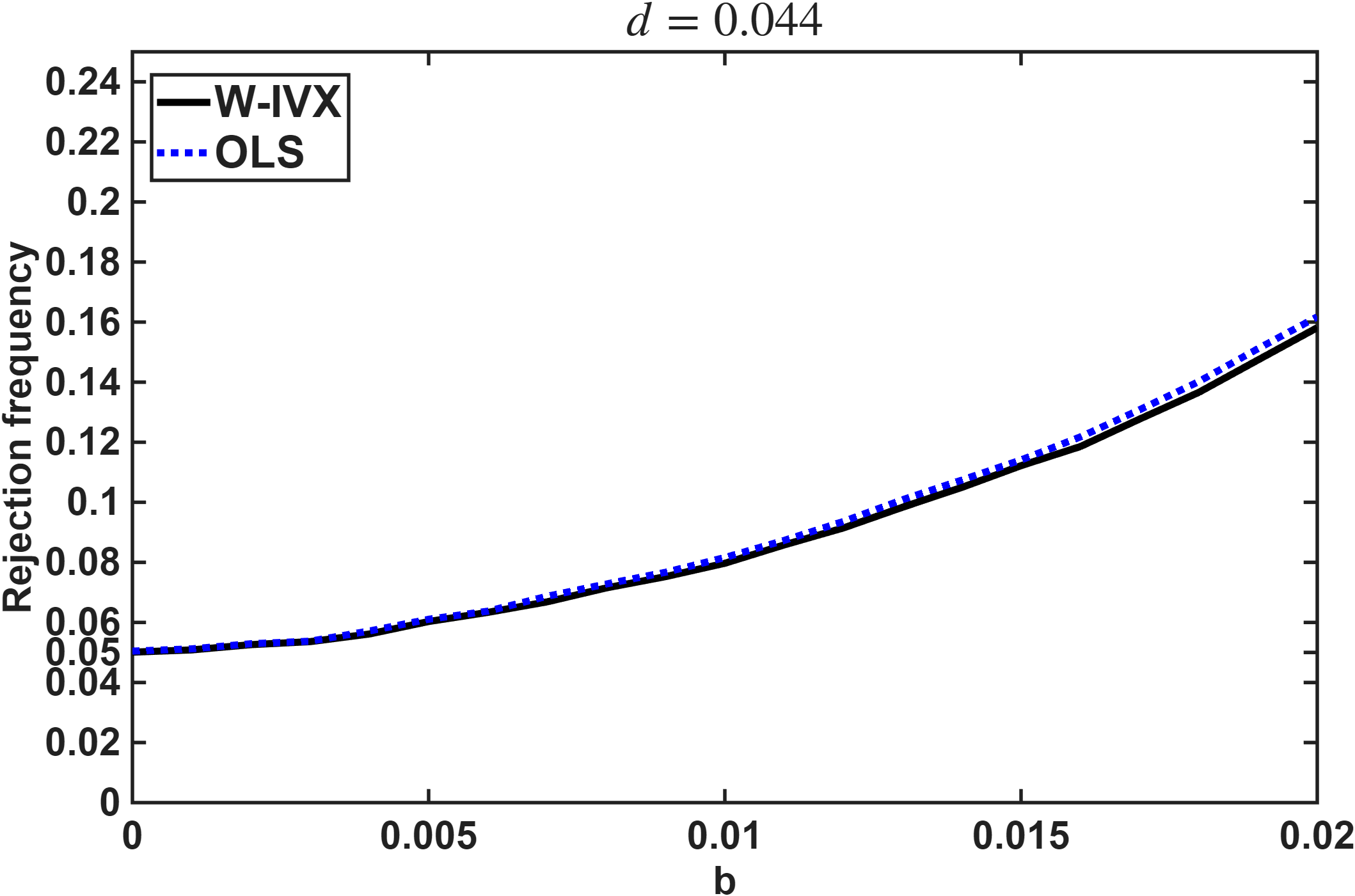}
\end{subfigure}

\end{figure}

\begin{figure}[p]
\ContinuedFloat \centering \captionsetup[subfigure]{position=top}

\setcounter{subfigure}{2}

\caption[]{Size Adjusted Power (nominal size 5\%): $[f(x),w(x)]=[x,x]$, $\delta=-0.95$,
$n=1000$ (continued)}
\label{fig:sim-sap-linear-w10-continued}

\begin{subfigure}[t]{0.48\linewidth} \centering \caption{$\alpha=1.6$, $d_{*}=0.375$}
\label{fig:sim-sap-linear-w10-alpha-1p60}\includegraphics[width=0.49\linewidth]{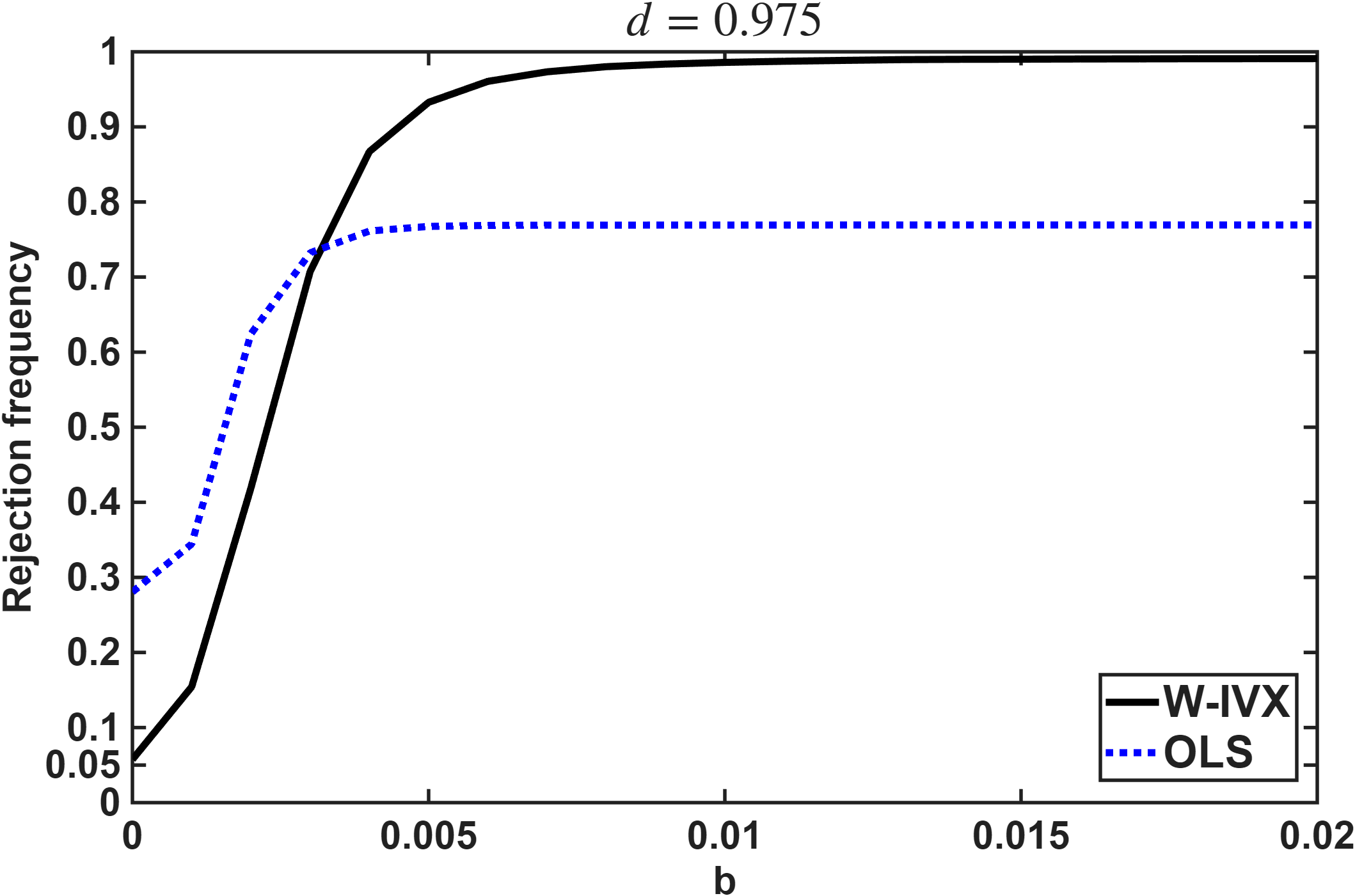}
\includegraphics[width=0.49\linewidth]{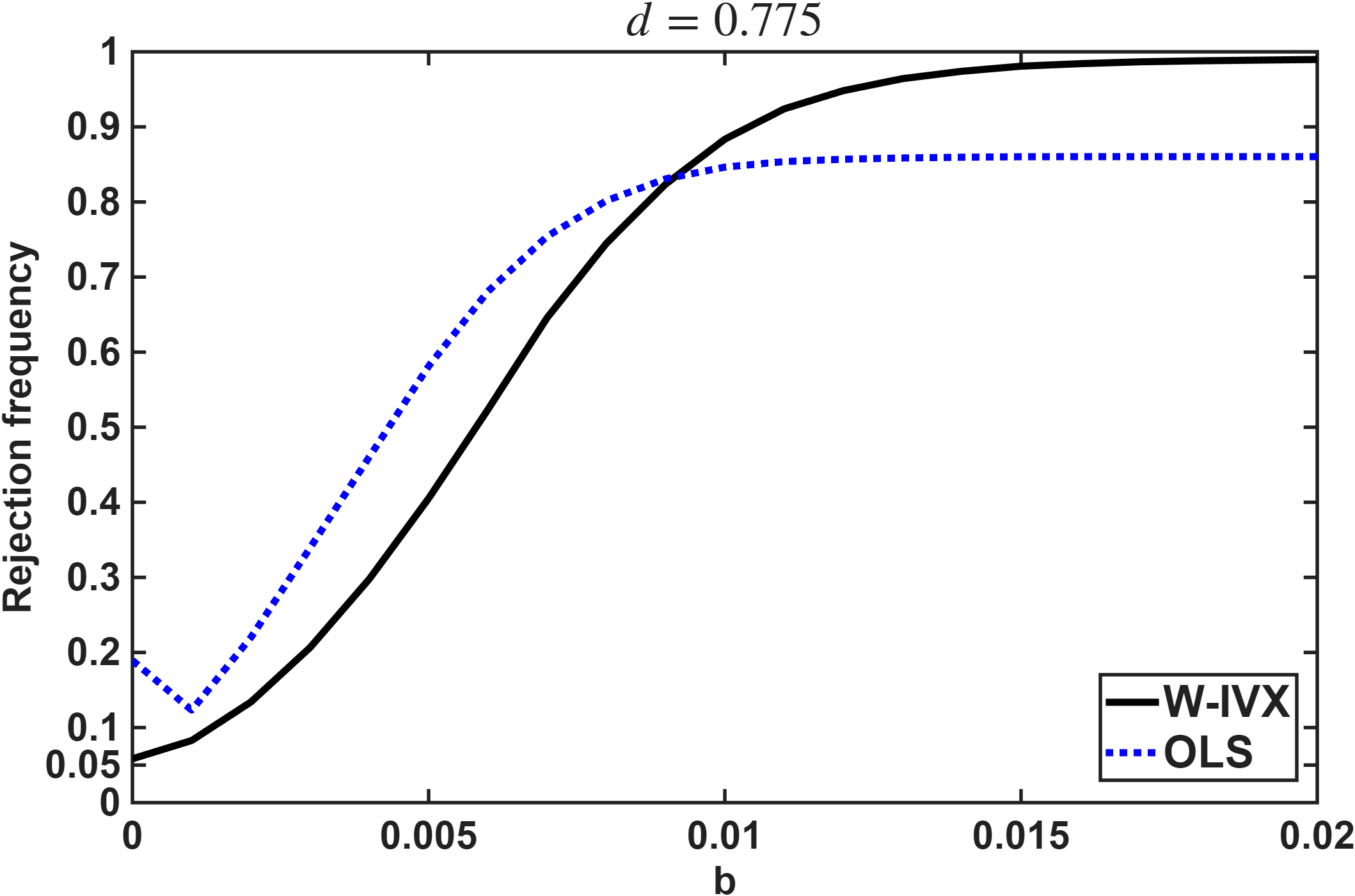}\\[0.6ex]
\includegraphics[width=0.49\linewidth]{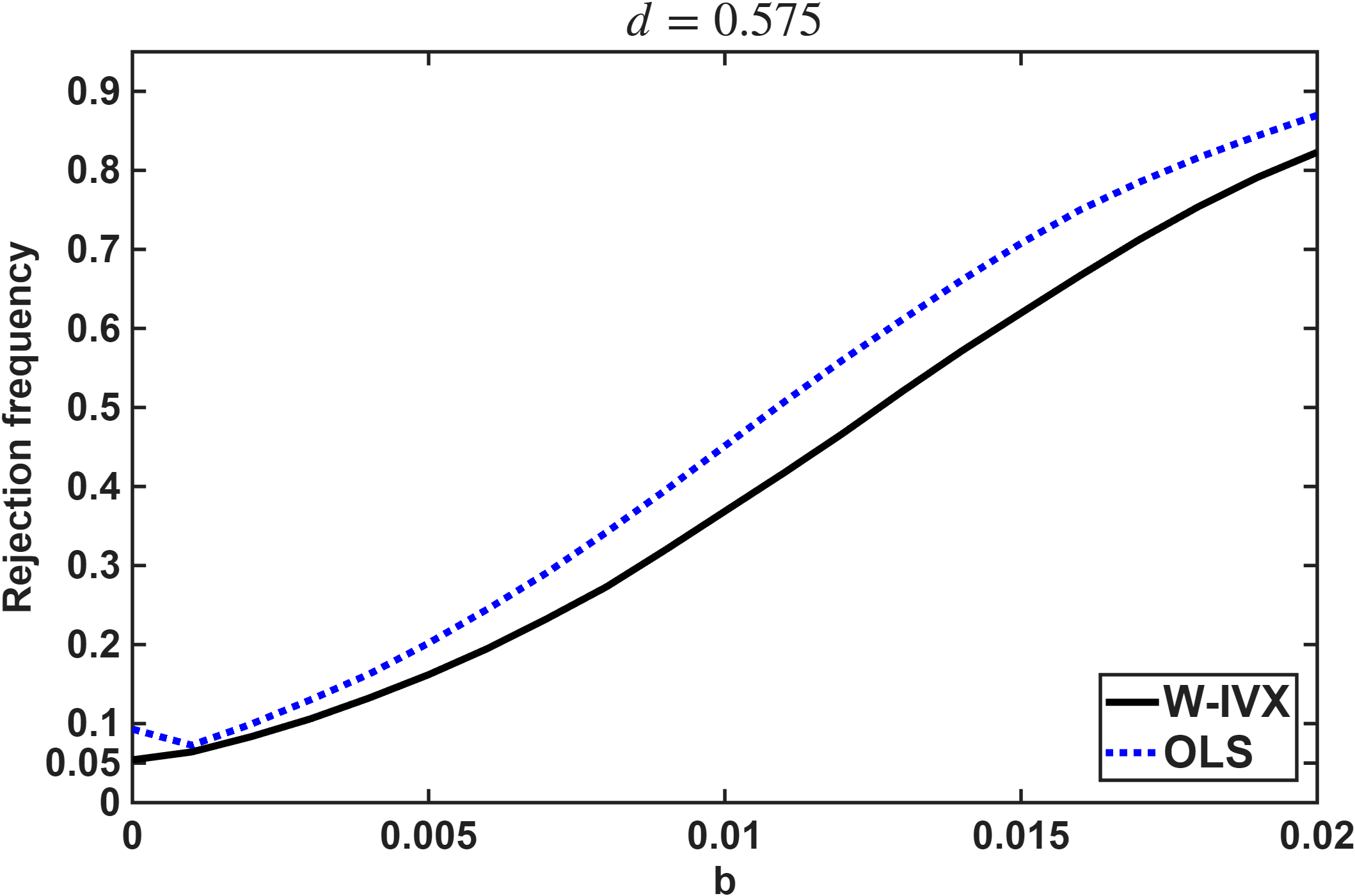}
\includegraphics[width=0.49\linewidth]{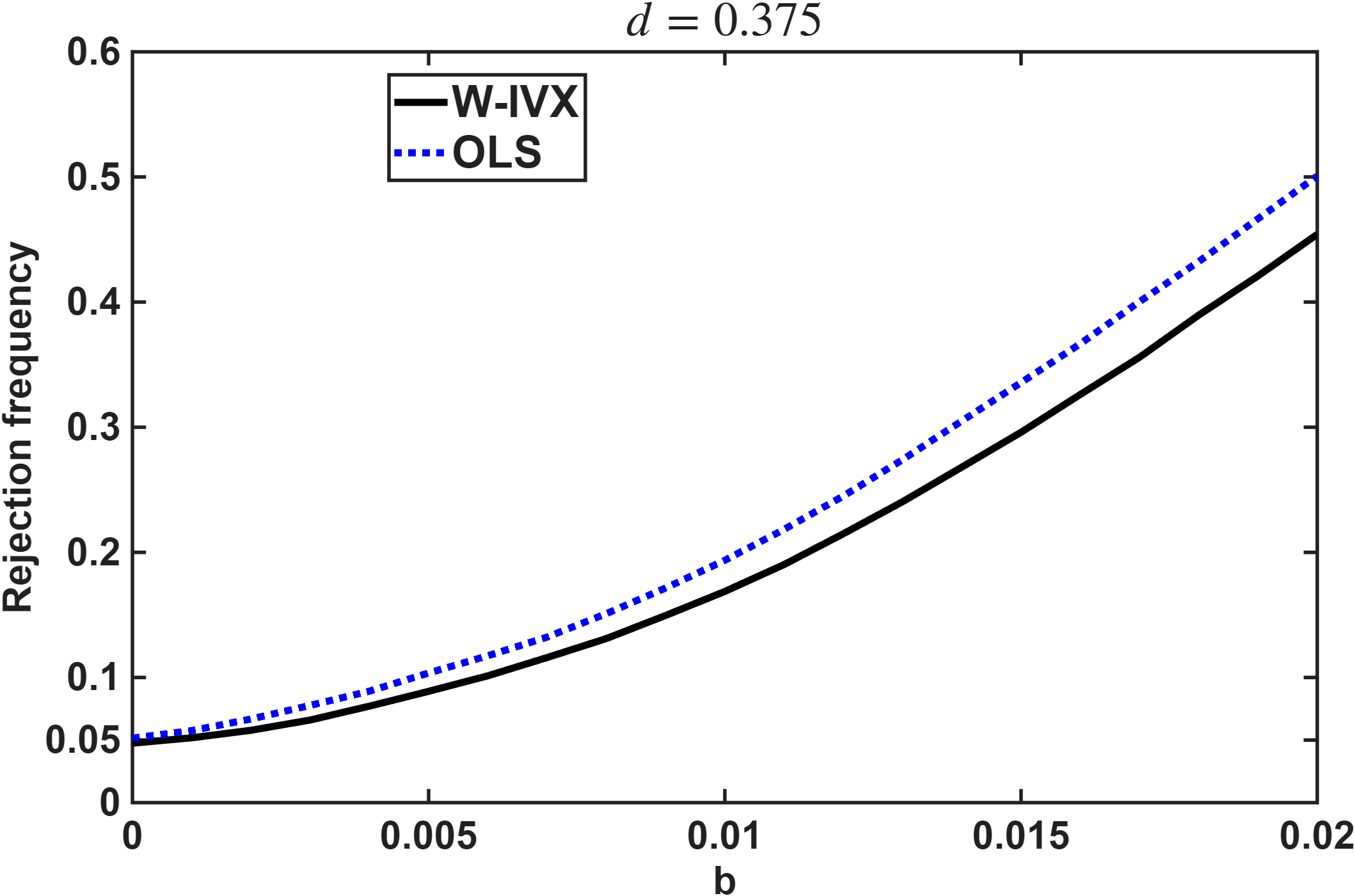}\\[0.6ex]
\includegraphics[width=0.49\linewidth]{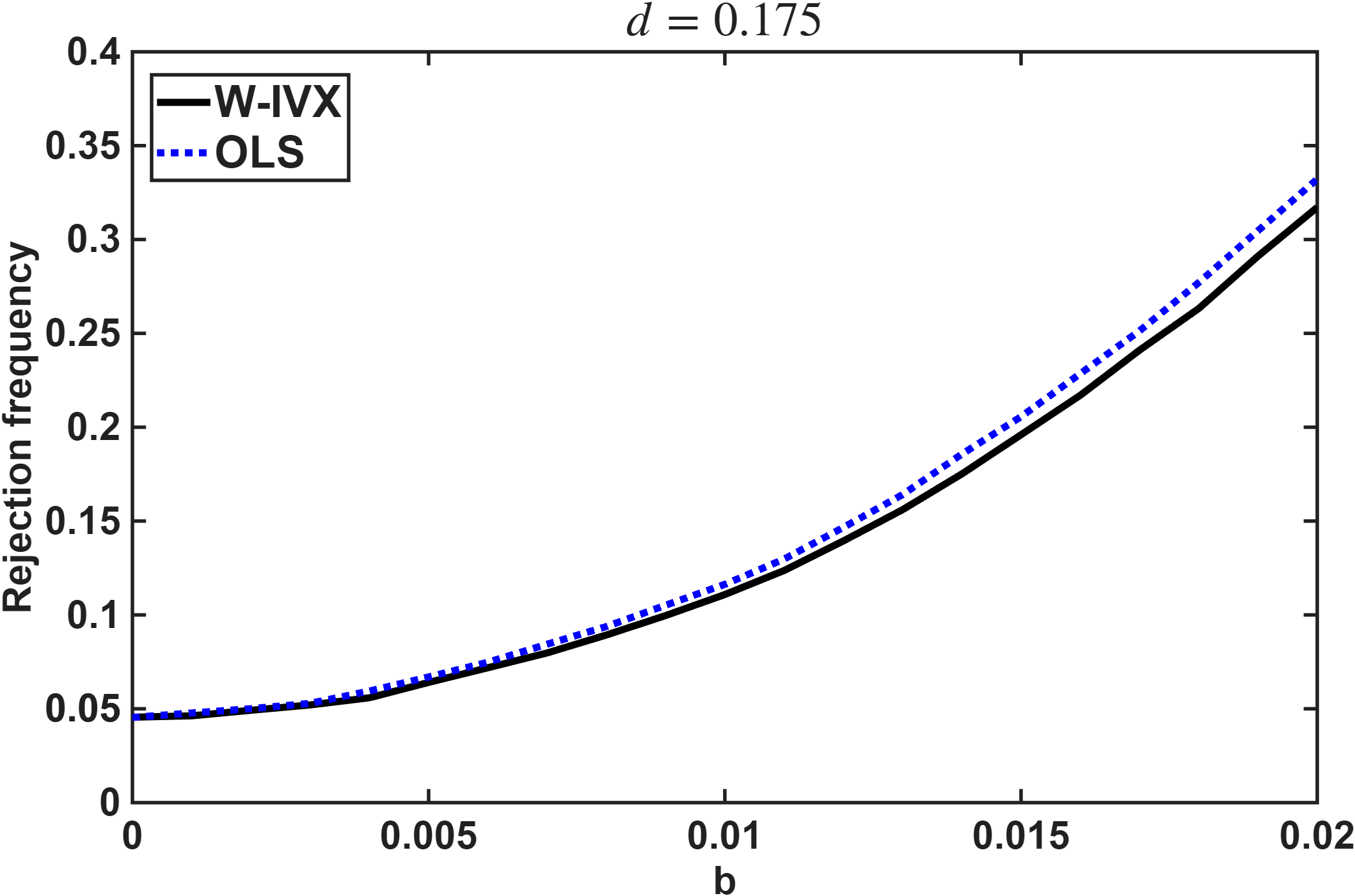}
\includegraphics[width=0.49\linewidth]{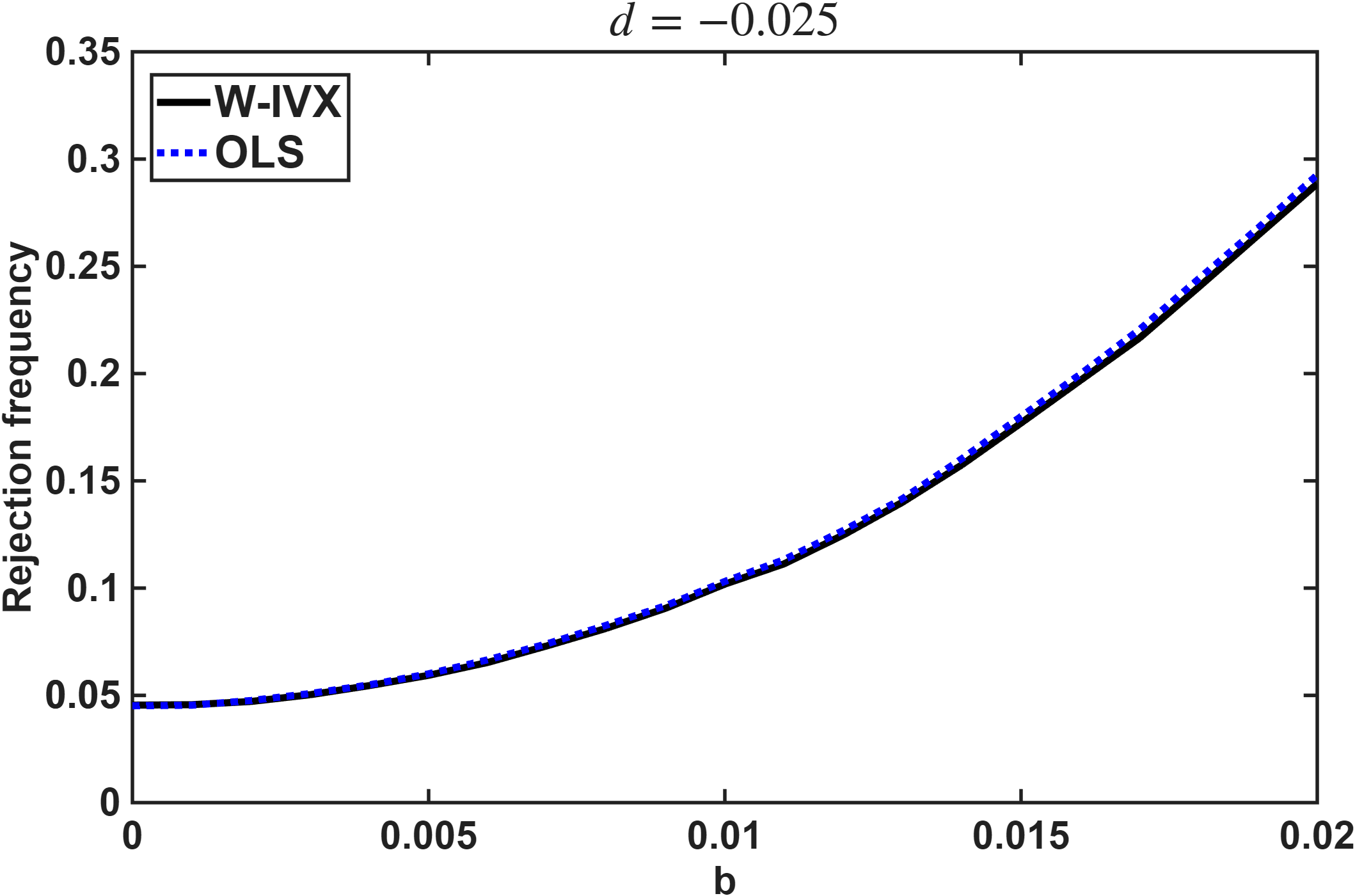}
\end{subfigure} \hfill{}
\begin{subfigure}[t]{0.48\linewidth} \centering \caption{$\alpha=1.4$, $d_{*}=0.286$}
\label{fig:sim-sap-linear-w10-alpha-1p40}\includegraphics[width=0.49\linewidth]{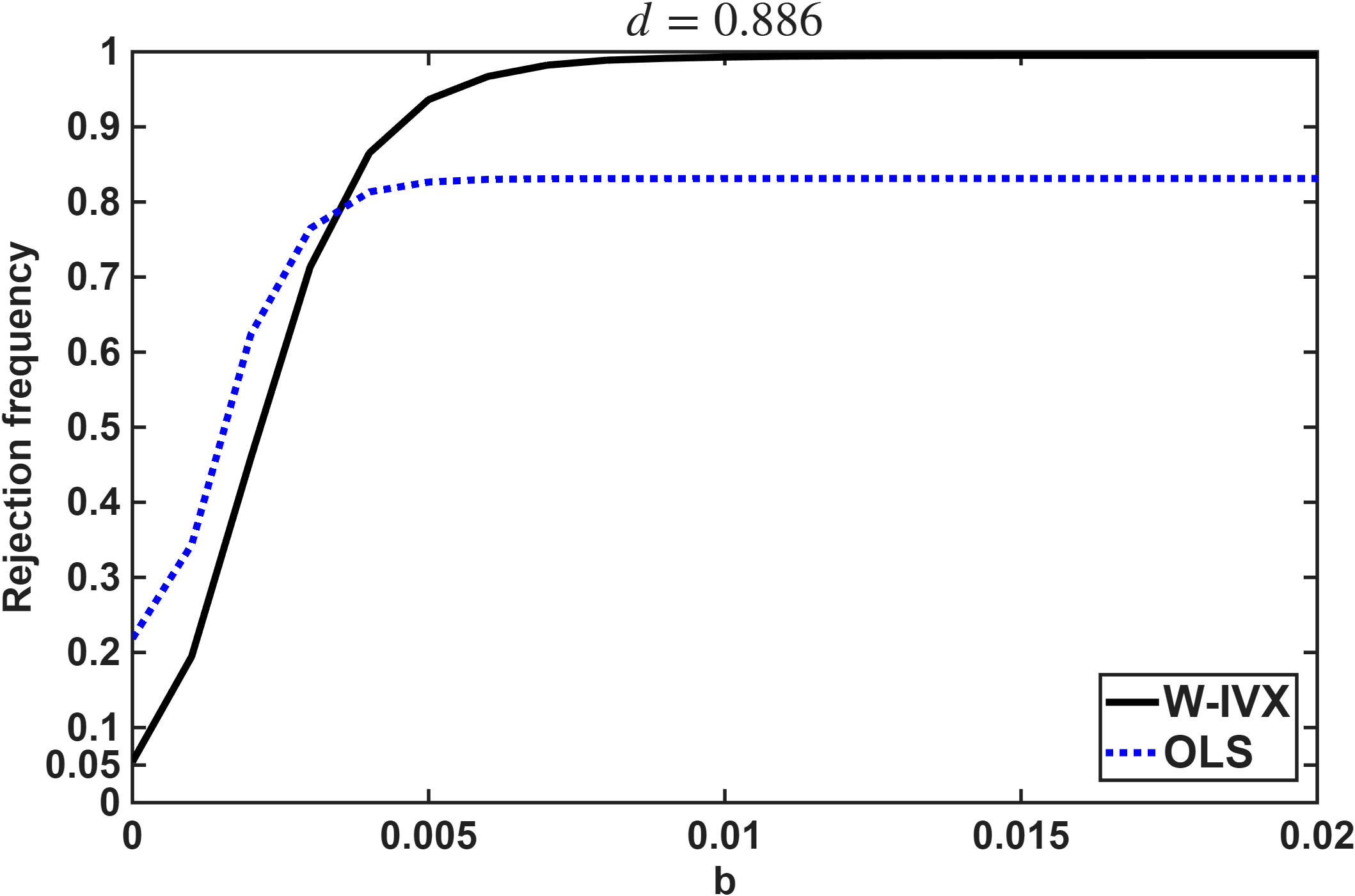}
\includegraphics[width=0.49\linewidth]{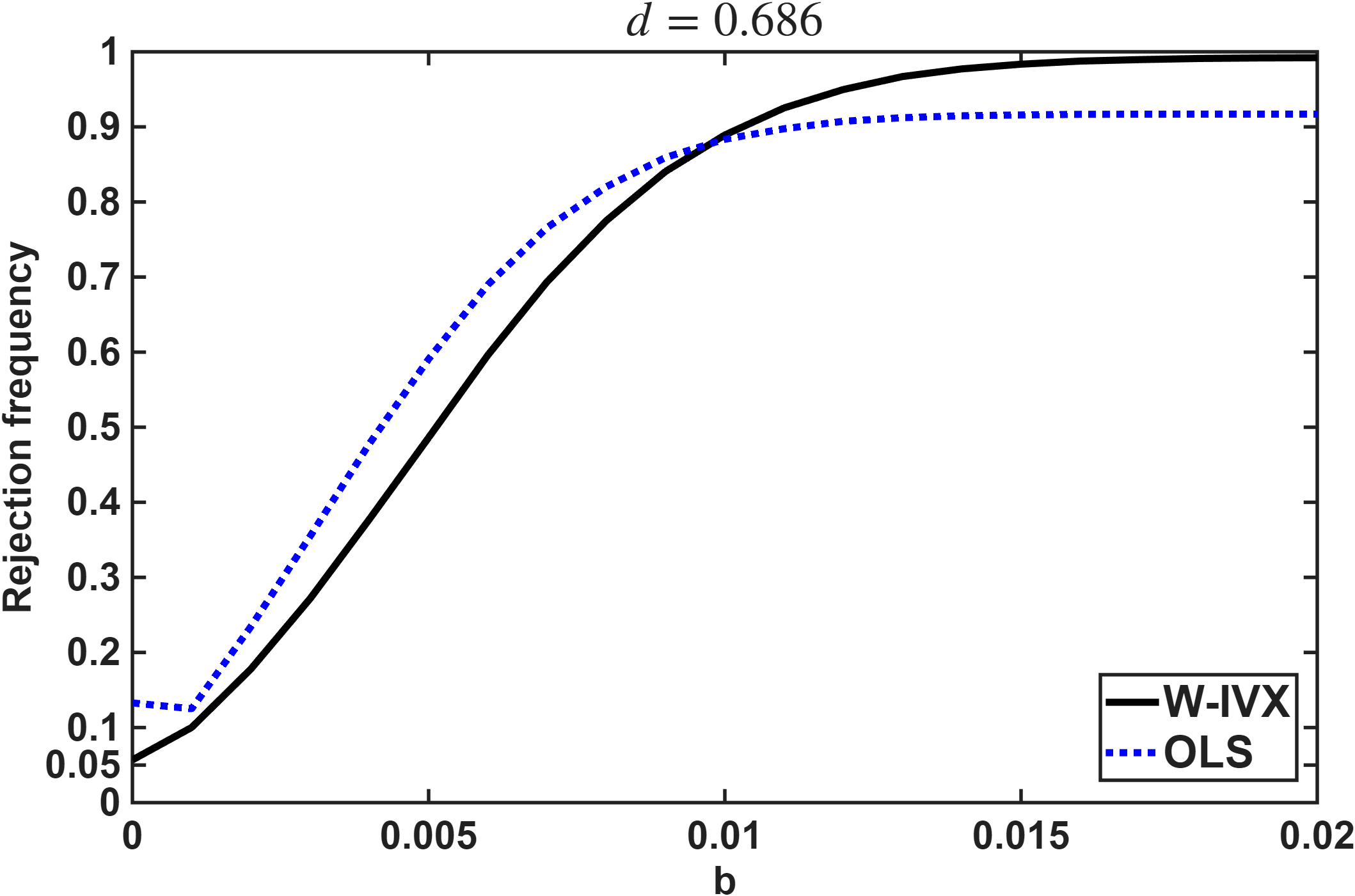}\\[0.6ex]
\includegraphics[width=0.49\linewidth]{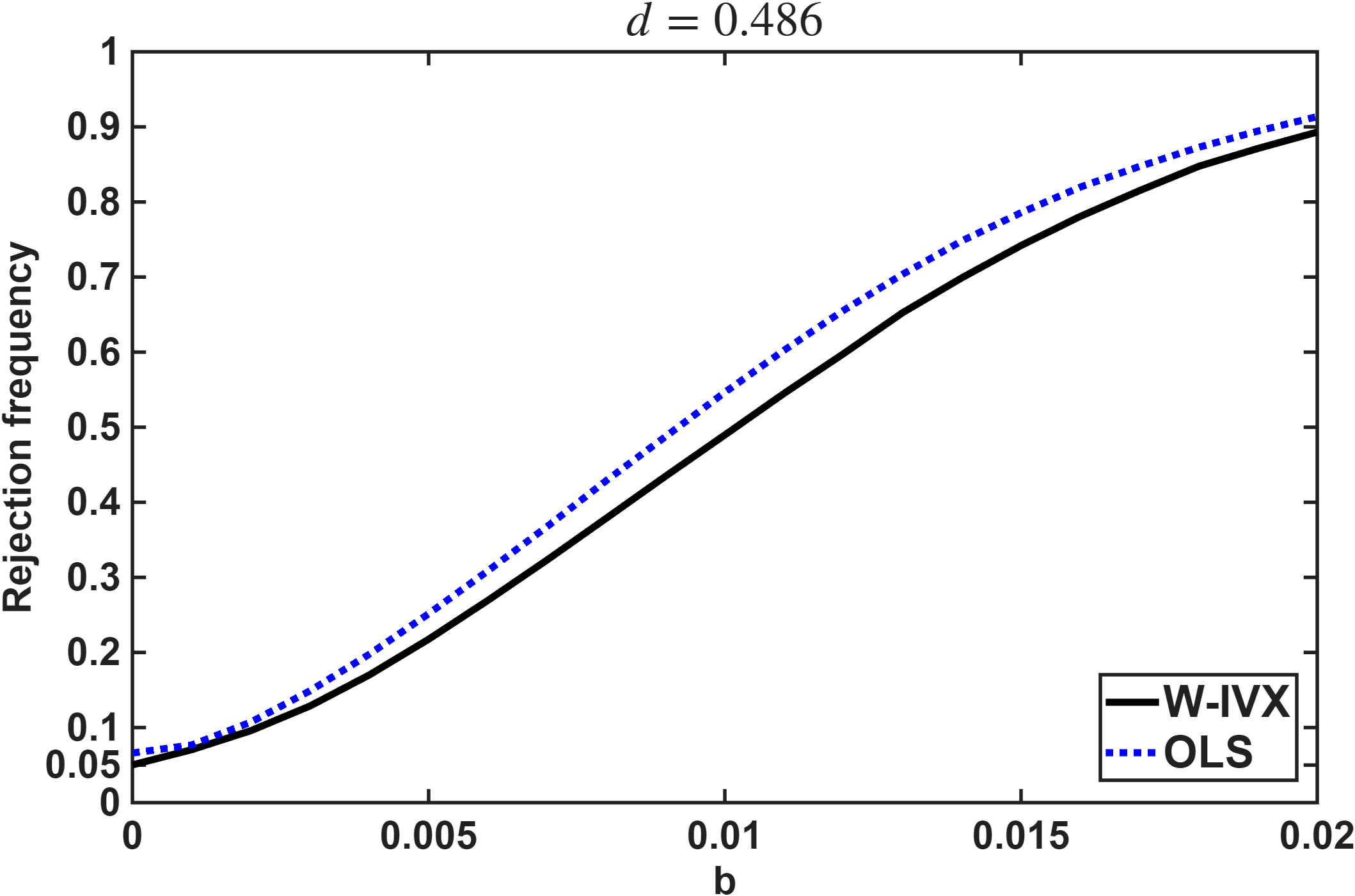}
\includegraphics[width=0.49\linewidth]{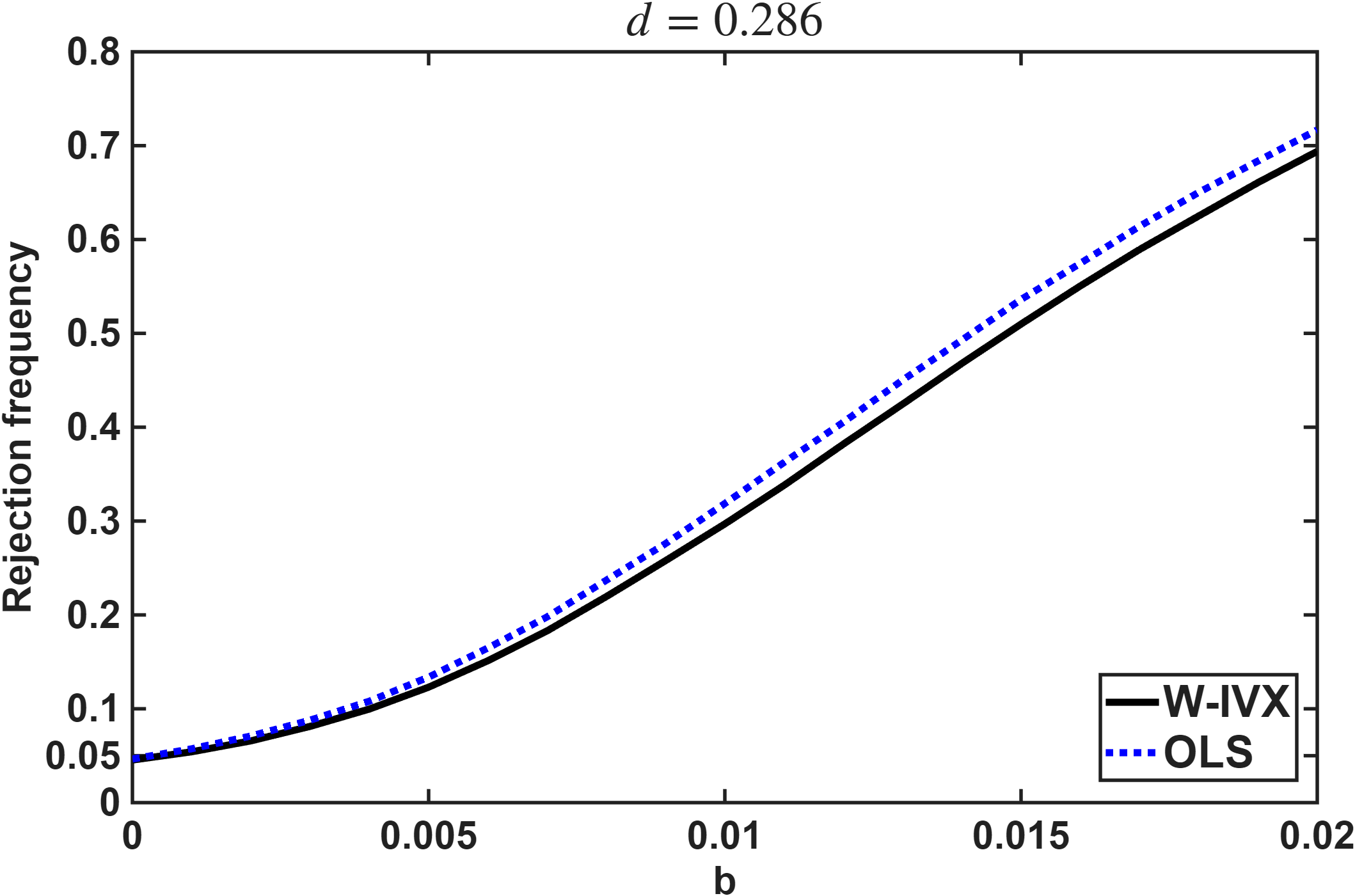}\\[0.6ex]
\includegraphics[width=0.49\linewidth]{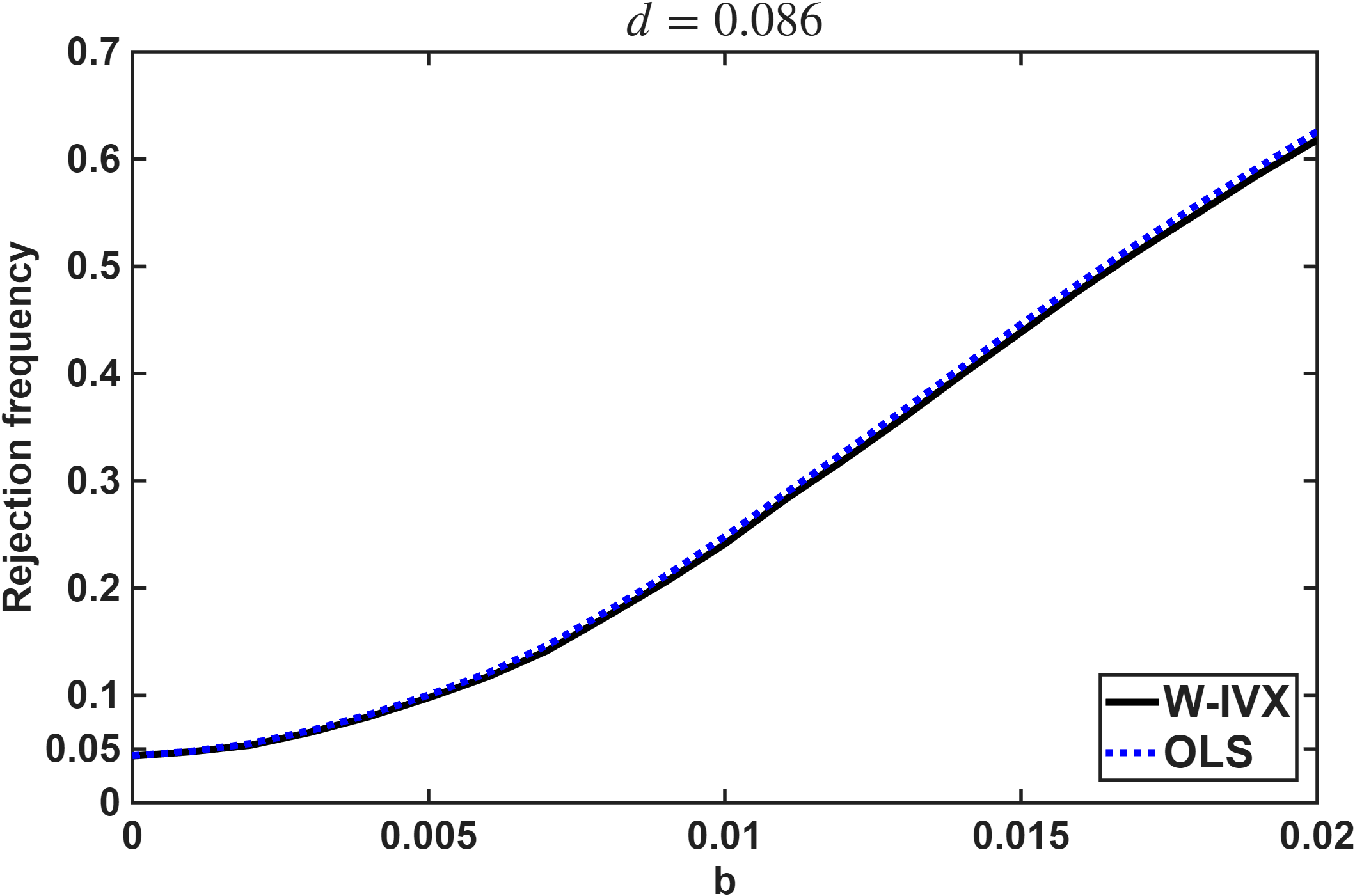}
\includegraphics[width=0.49\linewidth]{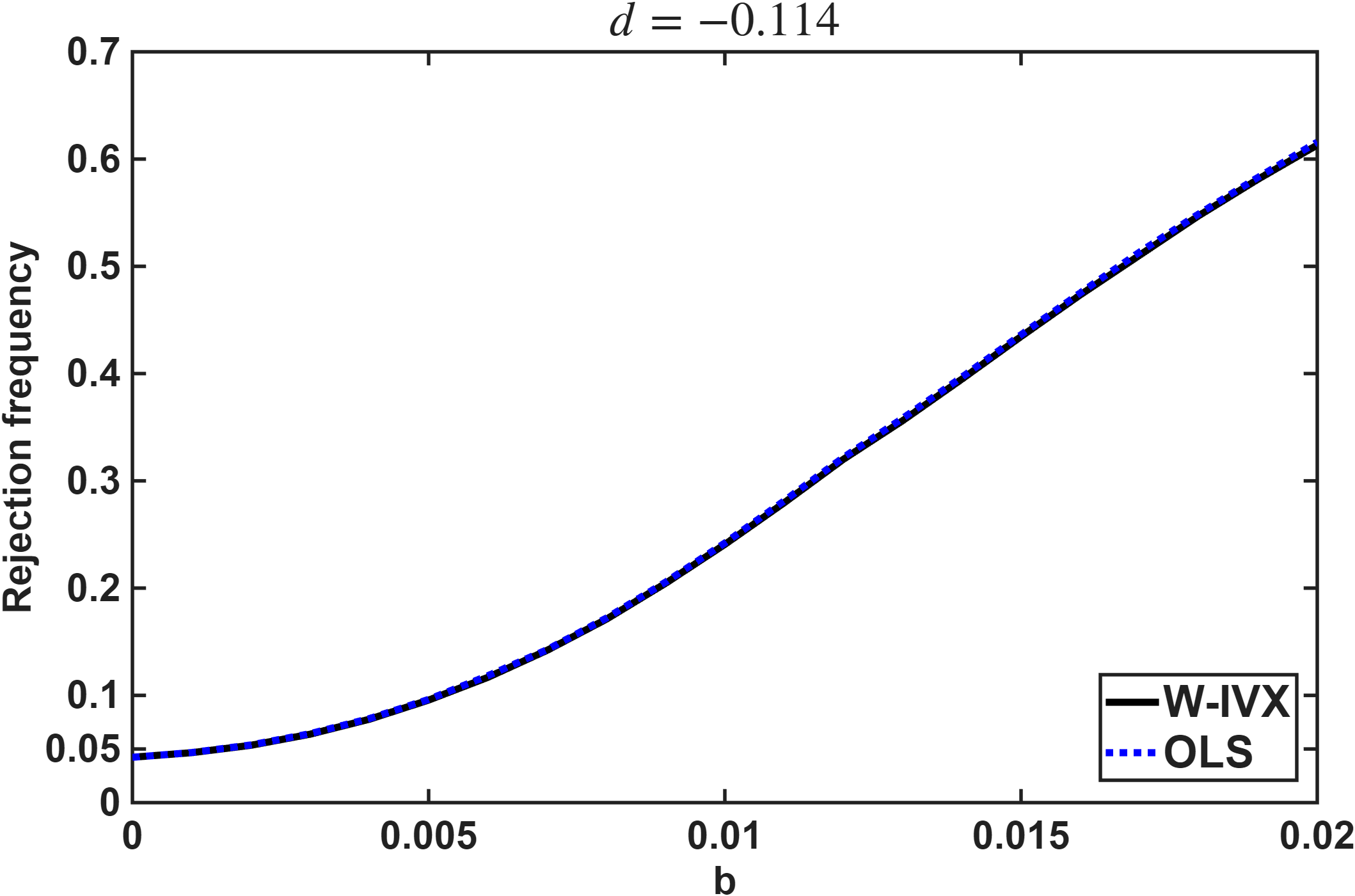}
\end{subfigure}
\end{figure}

\end{landscape} \clearpage{}

\begin{landscape}

\begin{figure}[p]
\caption{Size Adjusted Power (nominal size 5\%): \textbf{DGP1}, $[f(x),w(x)]=[x^{\langle0.5\rangle},x^{\langle0.5\rangle}]$,
$\delta=-0.95$, $n=1000$}
\label{fig:sim-sap-nonlinear-f05-w05}\centering \captionsetup[subfigure]{position=top}

\begin{subfigure}[t]{0.48\linewidth} \centering \caption{$\alpha=2$, $d_{*}=0.5$}
\label{fig:sim-sap-nonlinear-f05-w05-alpha-2p00} \includegraphics[width=0.49\linewidth]{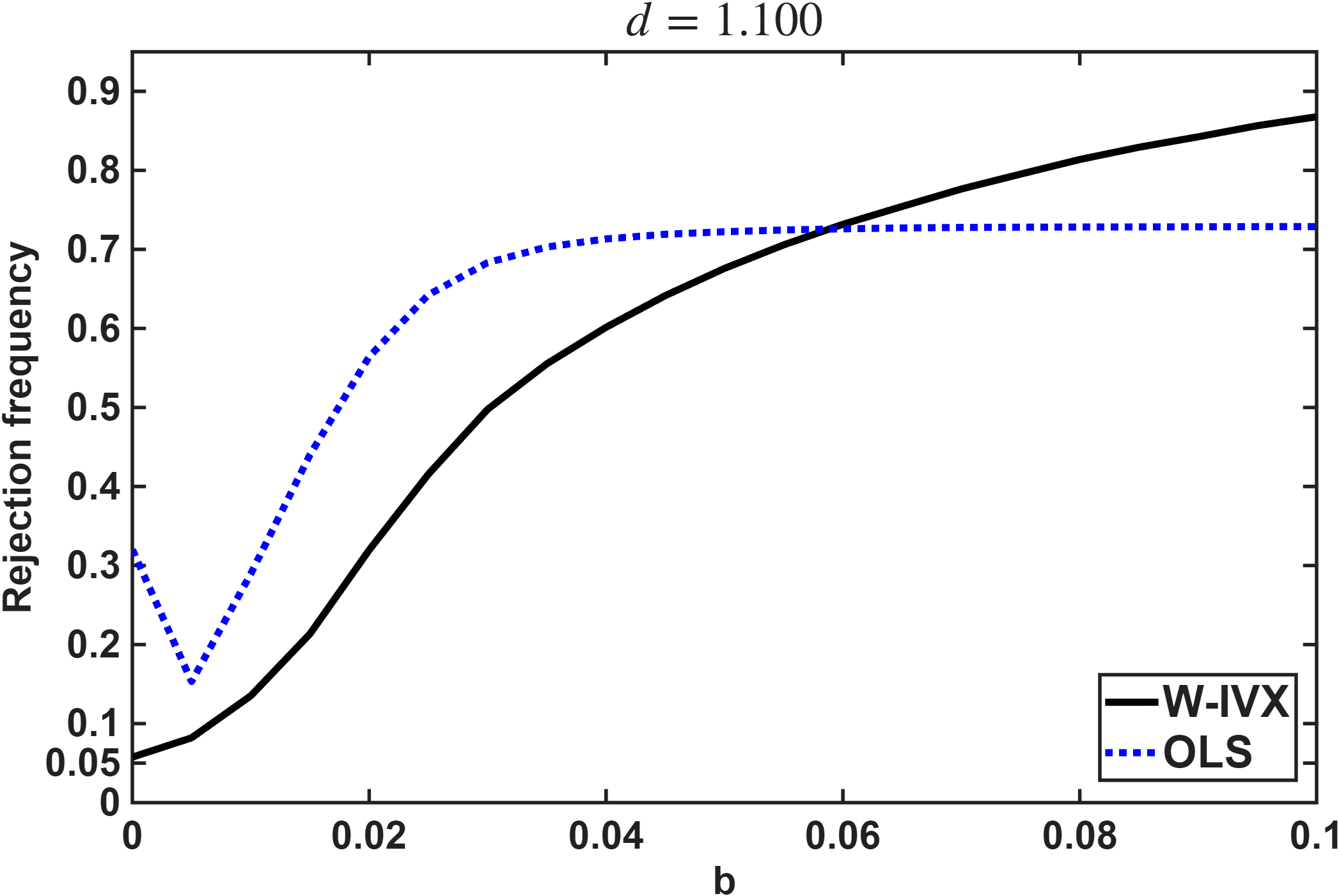}
\includegraphics[width=0.49\linewidth]{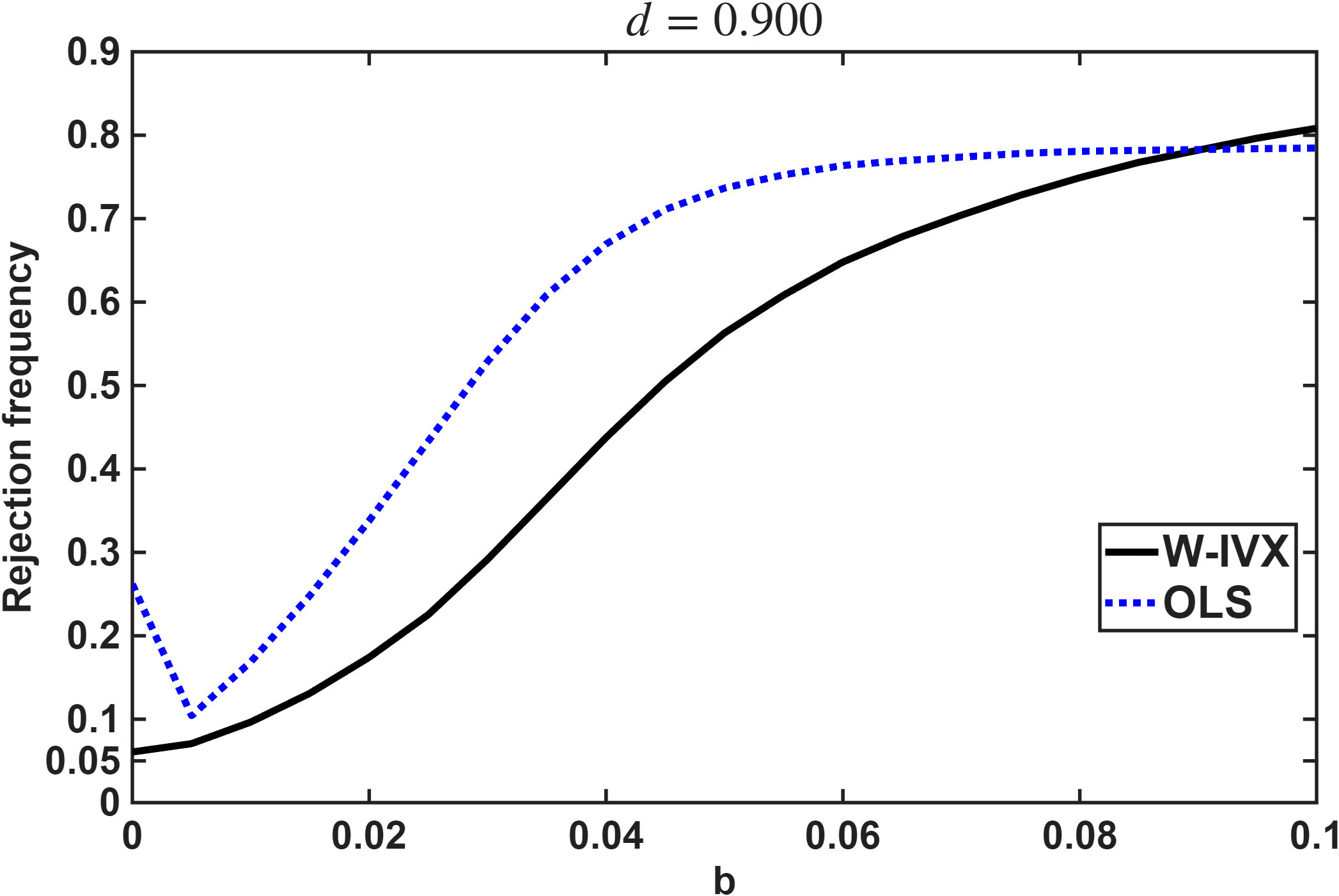}\\[0.6ex]
\includegraphics[width=0.49\linewidth]{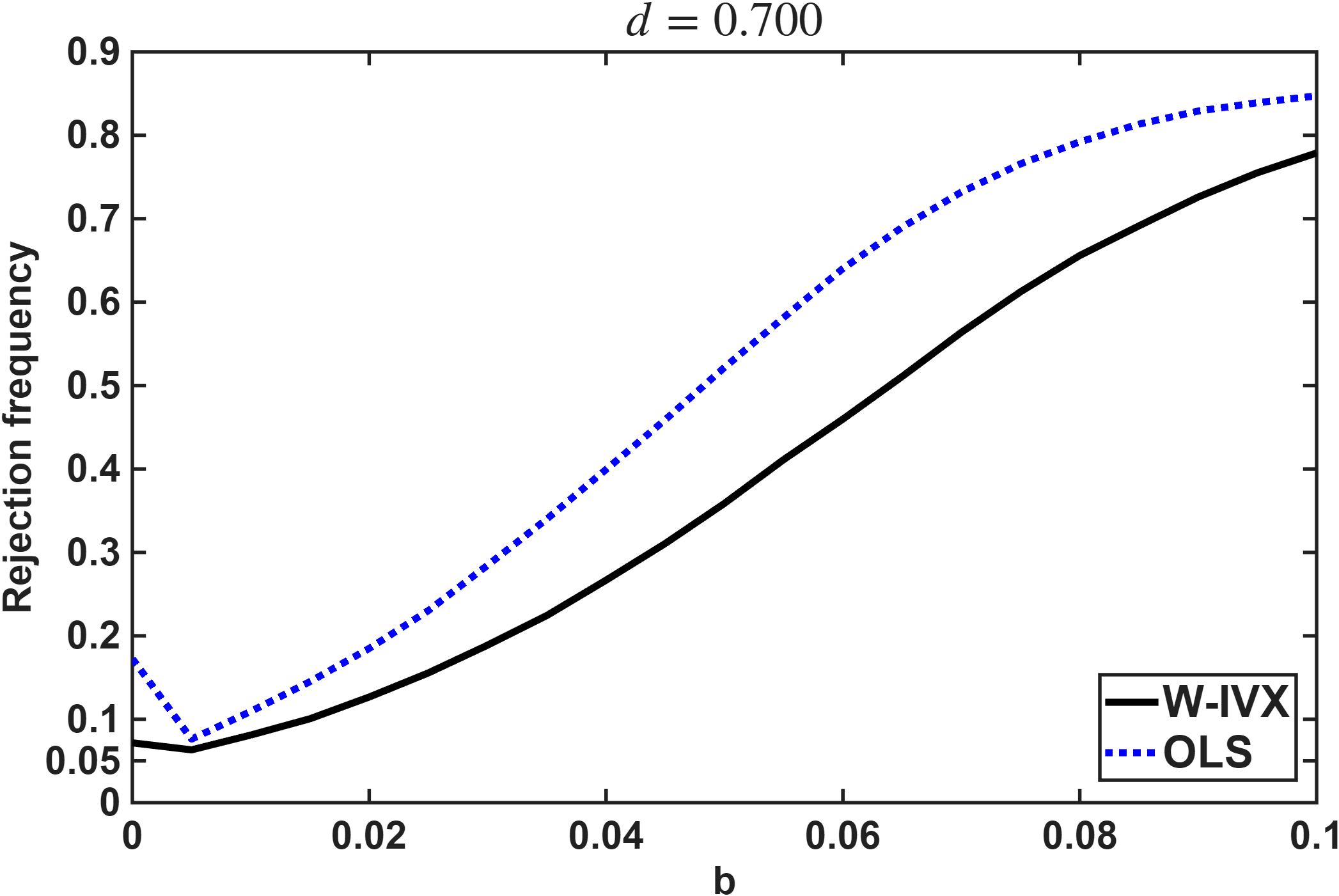}
\includegraphics[width=0.49\linewidth]{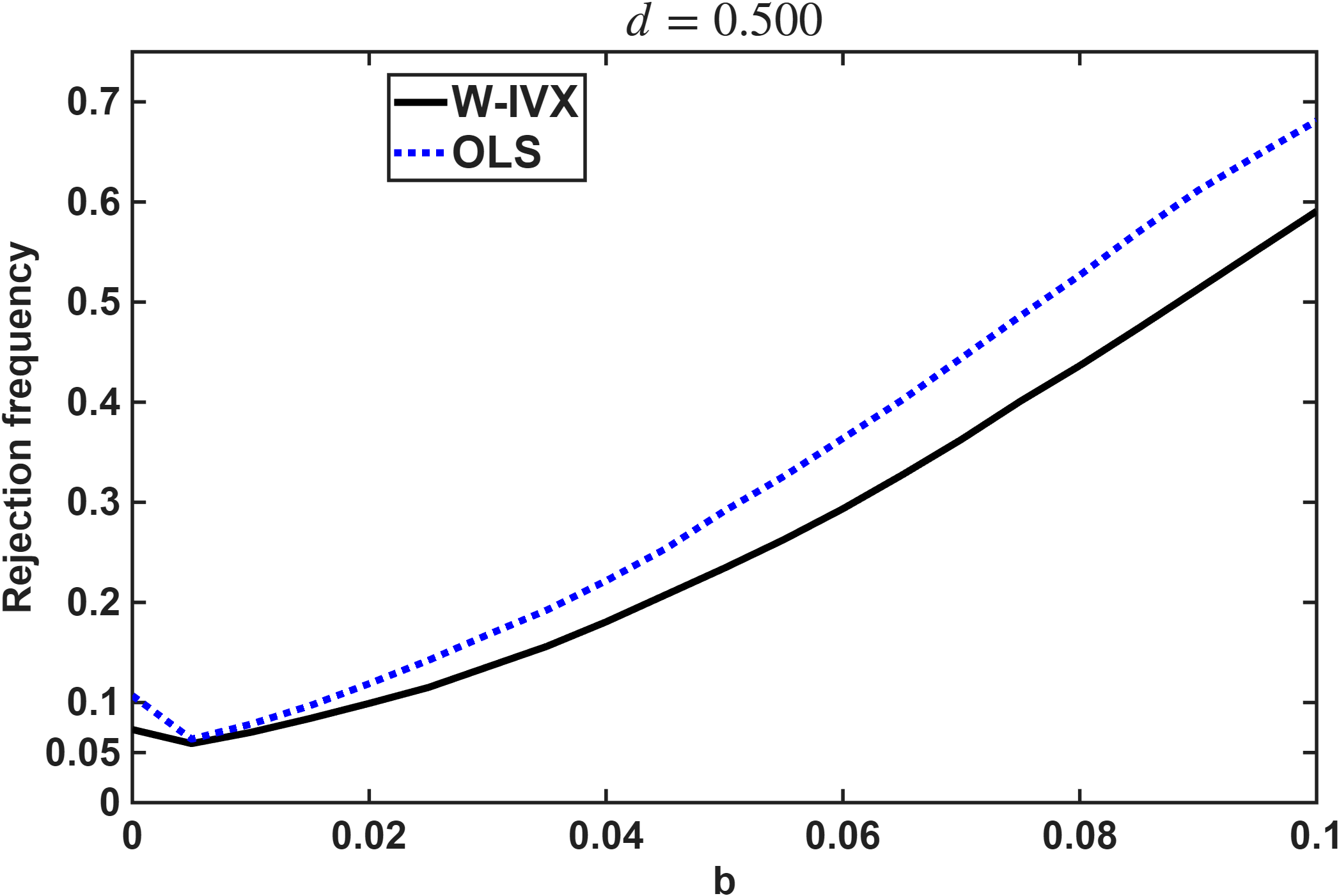}\\[0.6ex]
\includegraphics[width=0.49\linewidth]{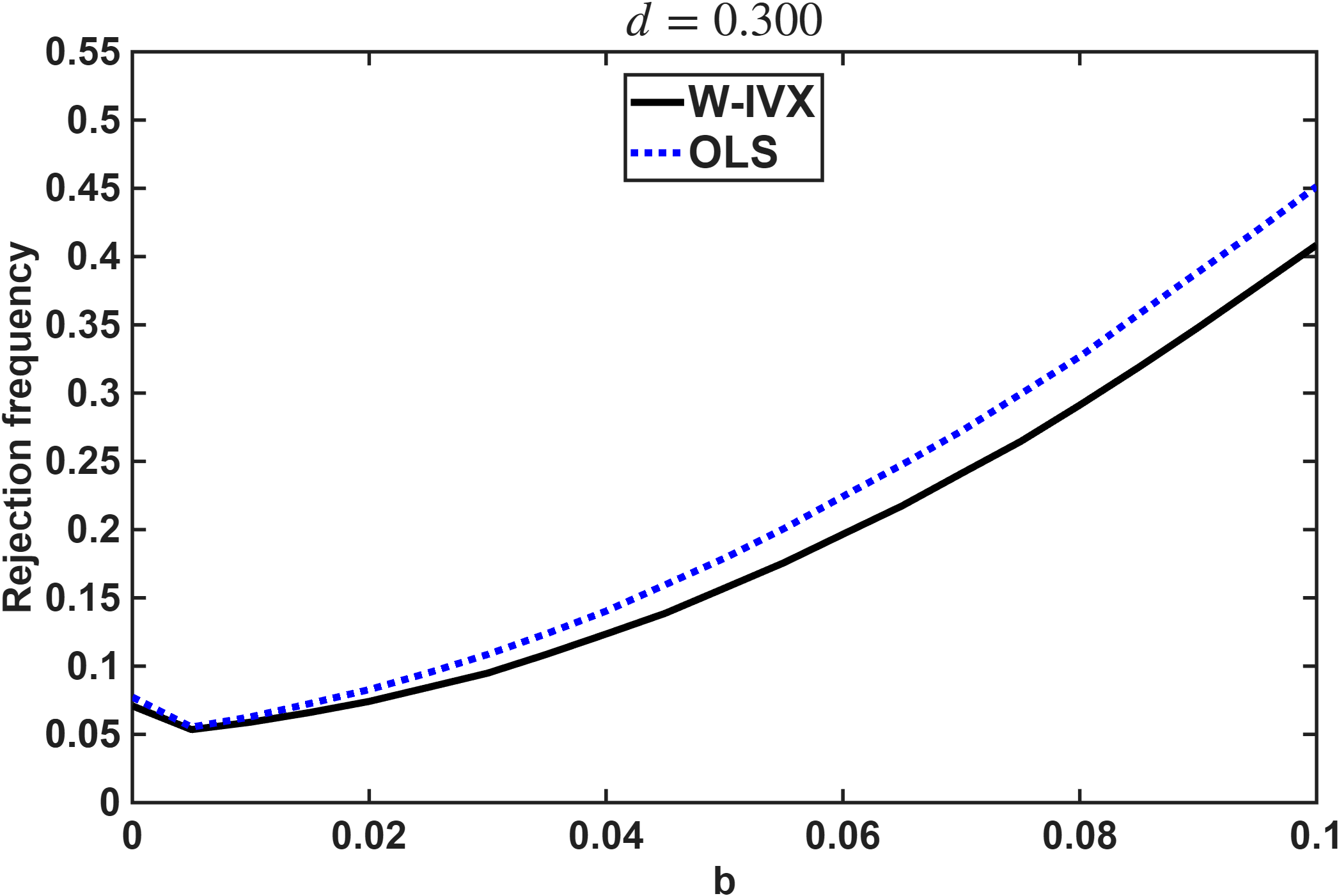}
\includegraphics[width=0.49\linewidth]{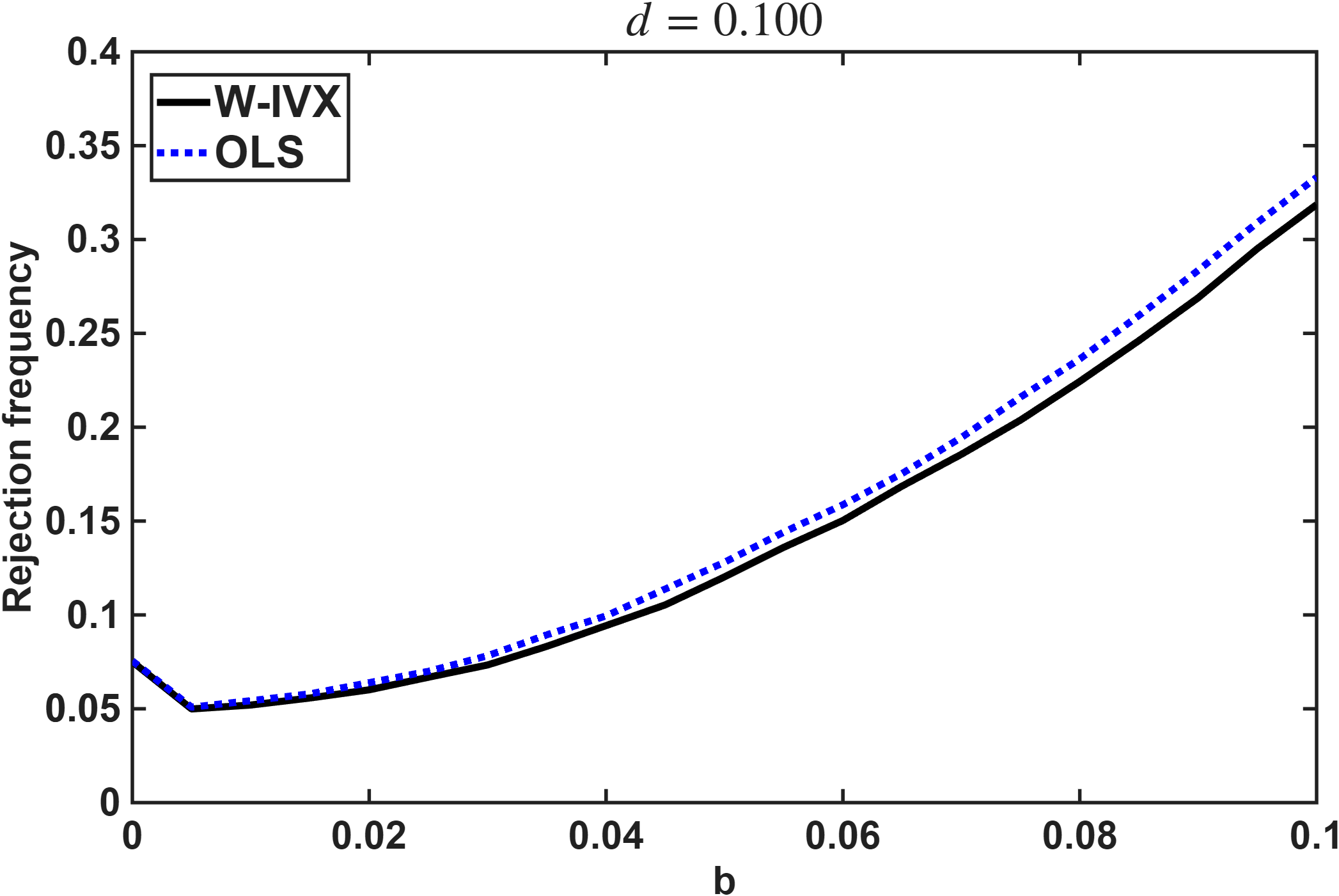}
\end{subfigure} \hfill{}
\begin{subfigure}[t]{0.48\linewidth} \centering \caption{$\alpha=1.8$, $d_{*}=0.444$}
\label{fig:sim-sap-nonlinear-f05-w05-alpha-1p80} \includegraphics[width=0.49\linewidth]{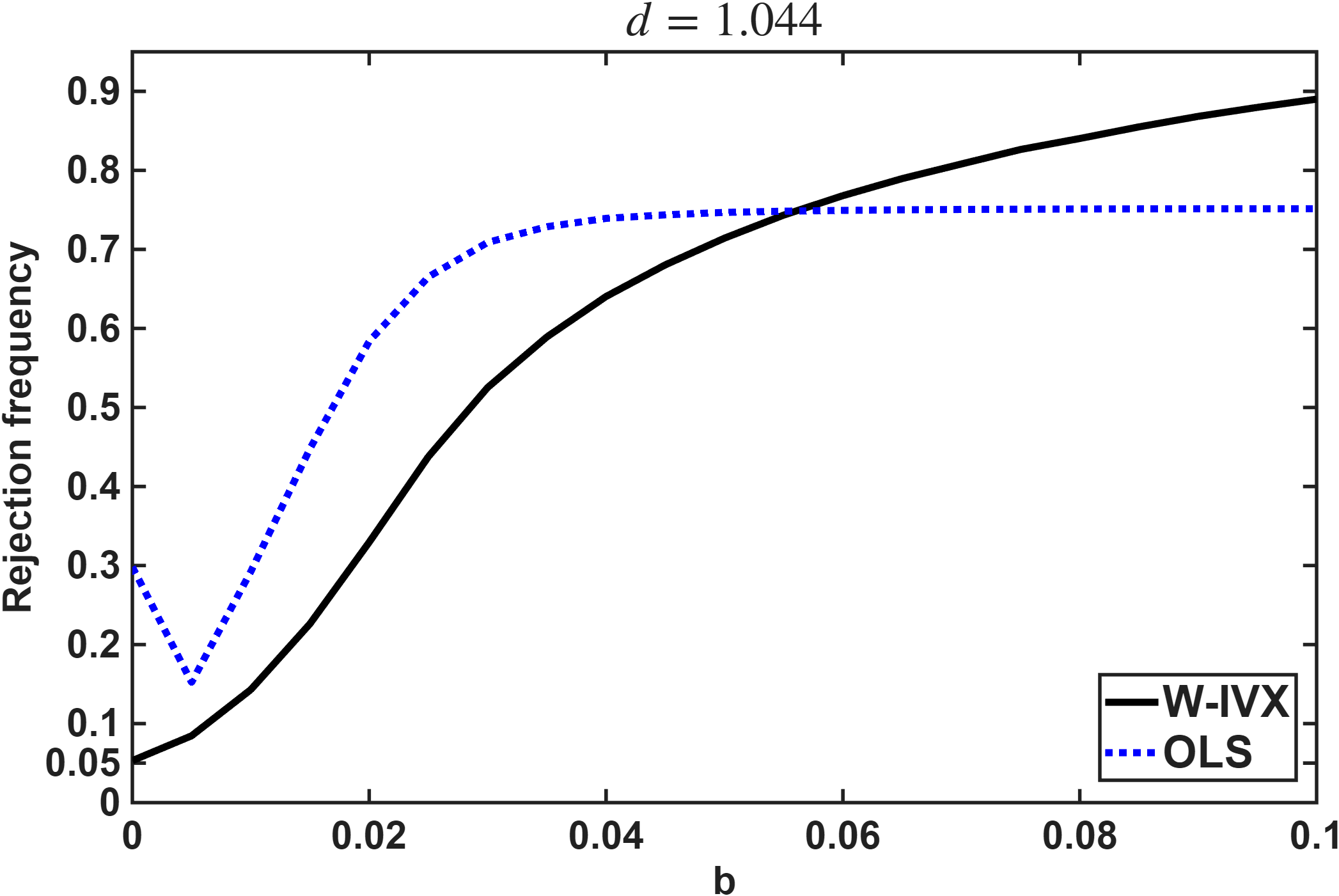}
\includegraphics[width=0.49\linewidth]{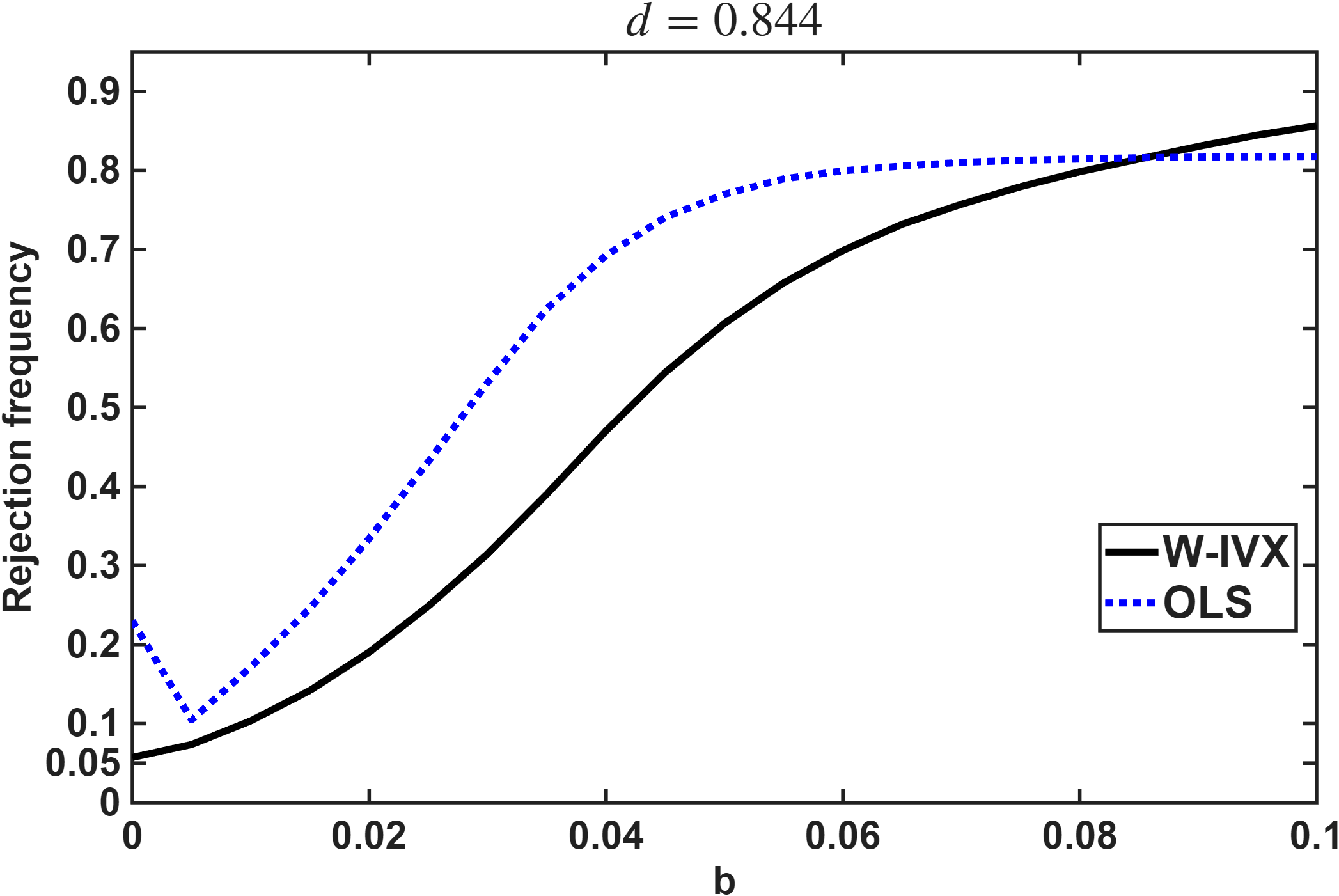}\\[0.6ex]
\includegraphics[width=0.49\linewidth]{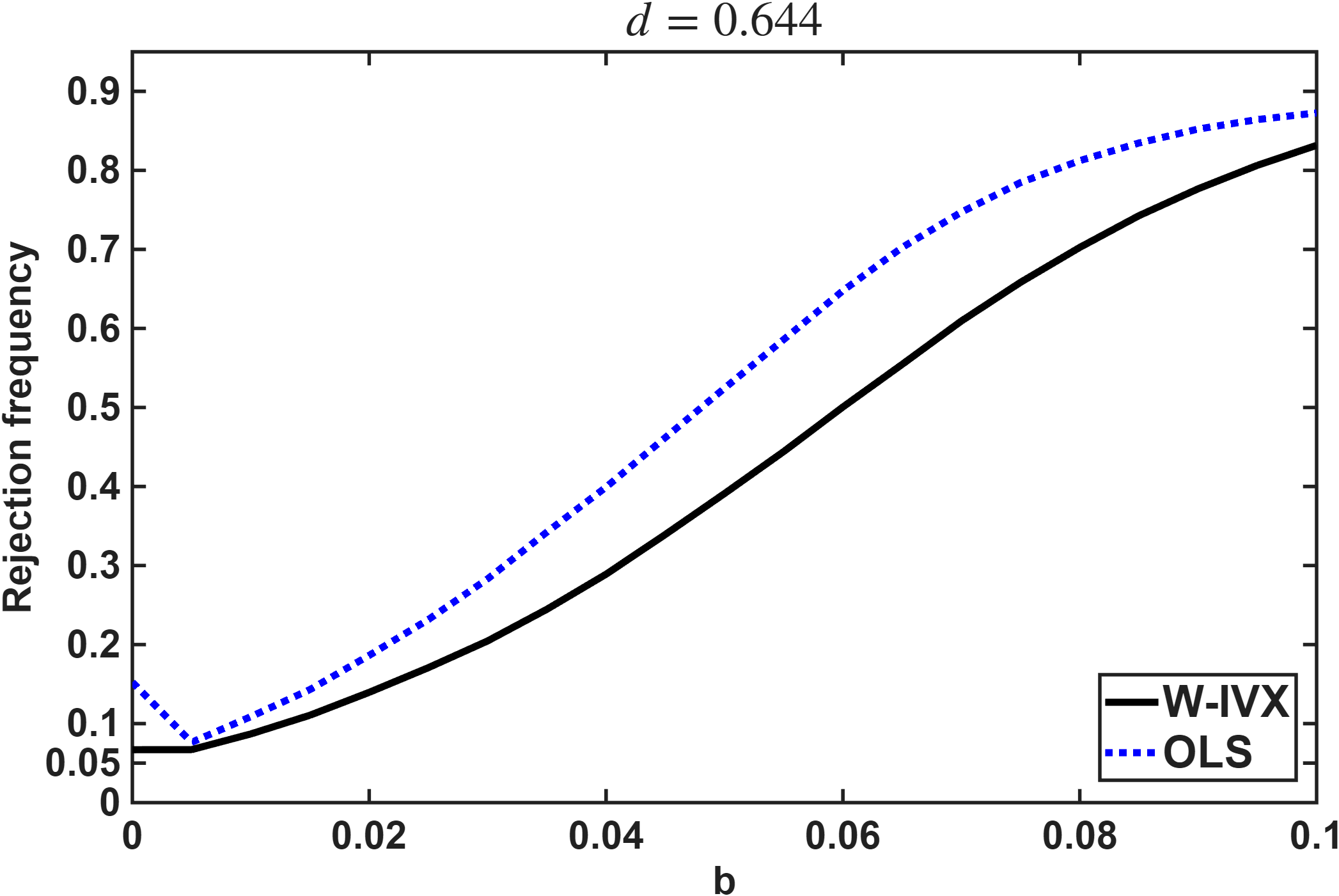}
\includegraphics[width=0.49\linewidth]{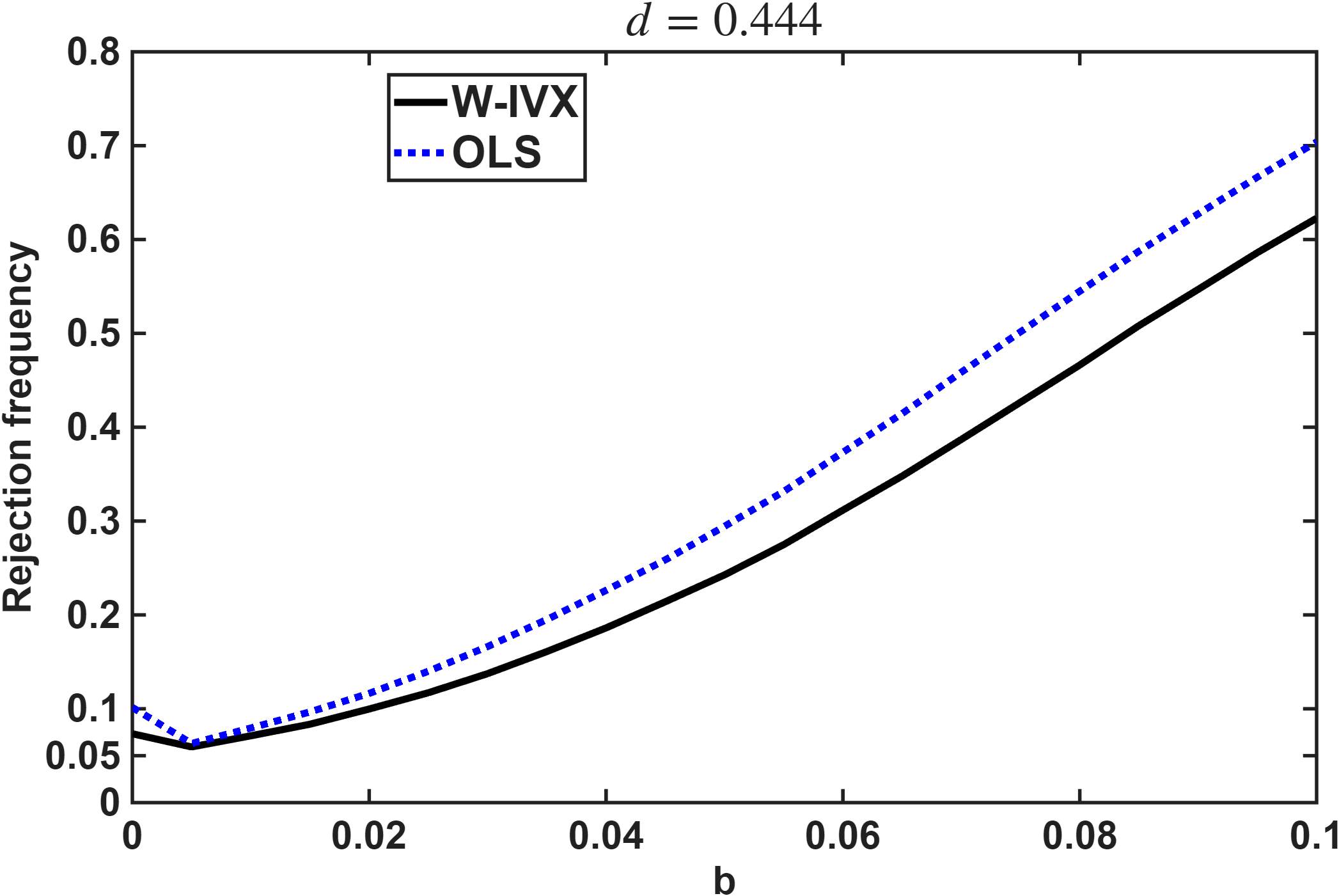}\\[0.6ex]
\includegraphics[width=0.49\linewidth]{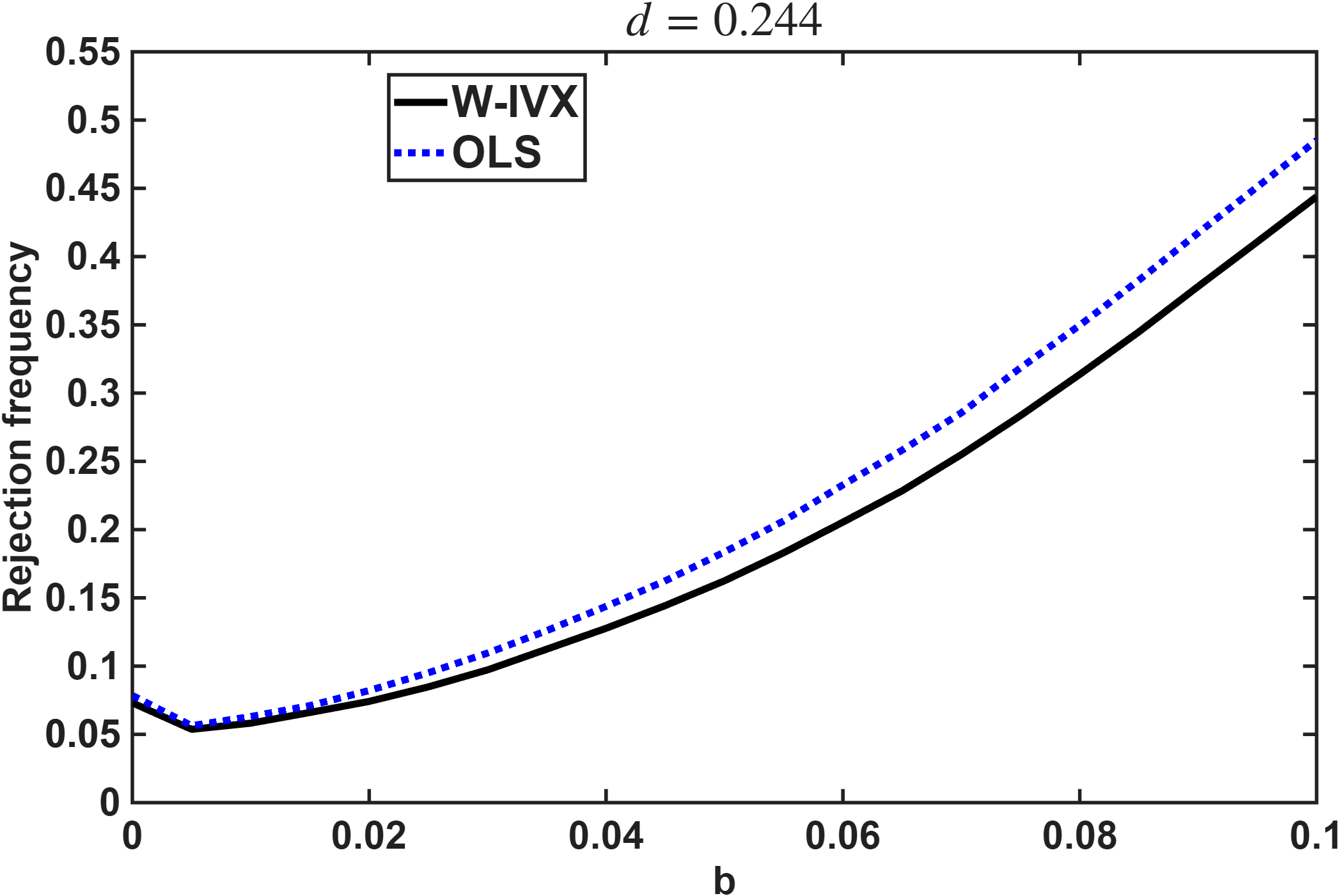}
\includegraphics[width=0.49\linewidth]{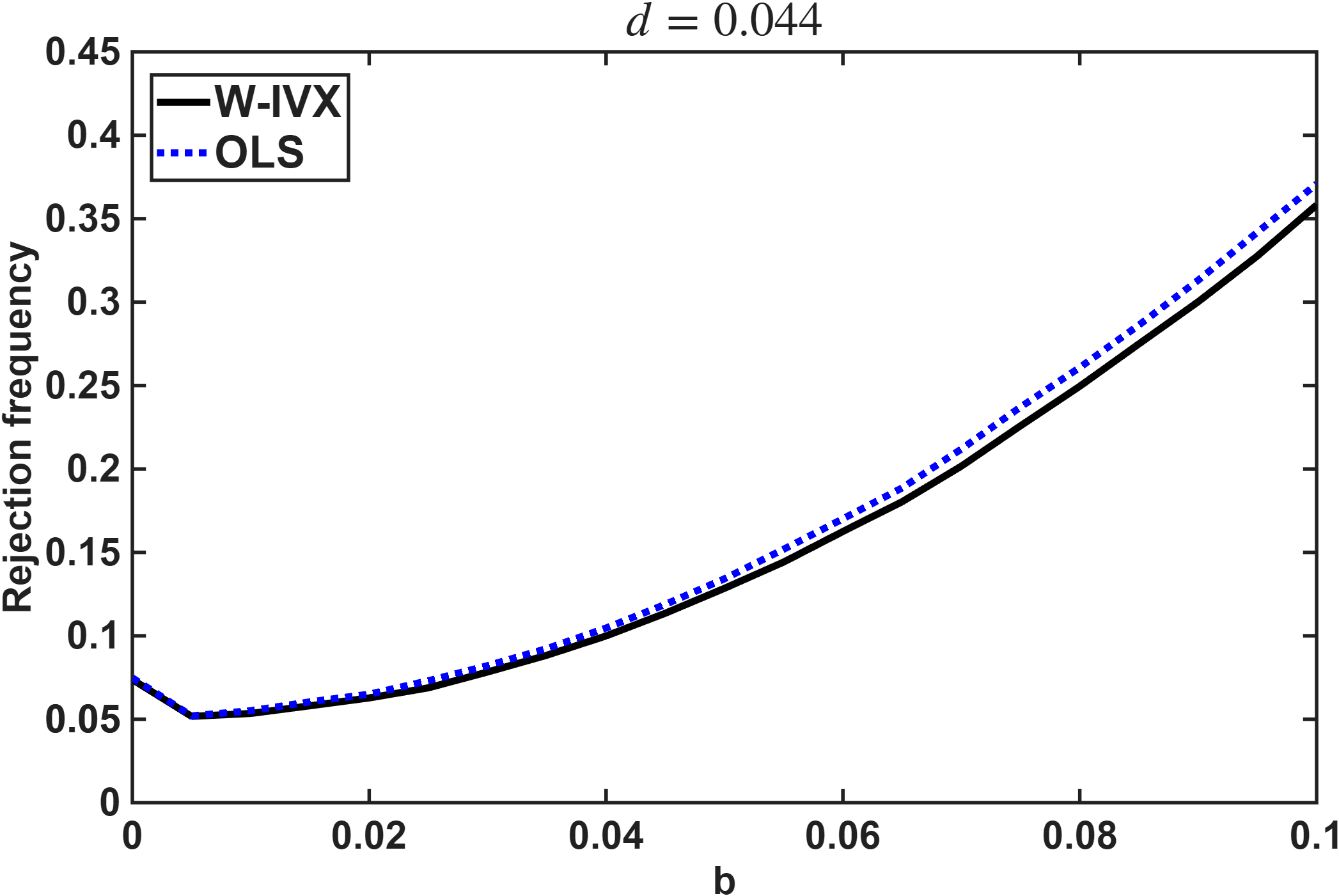}
\end{subfigure}

\end{figure}

\begin{figure}[p]
\ContinuedFloat \centering \captionsetup[subfigure]{position=top}

\setcounter{subfigure}{2}

\caption[]{Size Adjusted Power (nominal size 5\%): $[f(x),w(x)]=[x^{\langle0.5\rangle},x^{\langle0.5\rangle}]$,
$\delta=-0.95$, $n=1000$ (continued)}
\label{fig:sim-sap-nonlinear-f05-w05-continued}

\begin{subfigure}[t]{0.48\linewidth} \centering \caption{$\alpha=1.6$, $d_{*}=0.375$}
\label{fig:sim-sap-nonlinear-f05-w05-alpha-1p60}\includegraphics[width=0.49\linewidth]{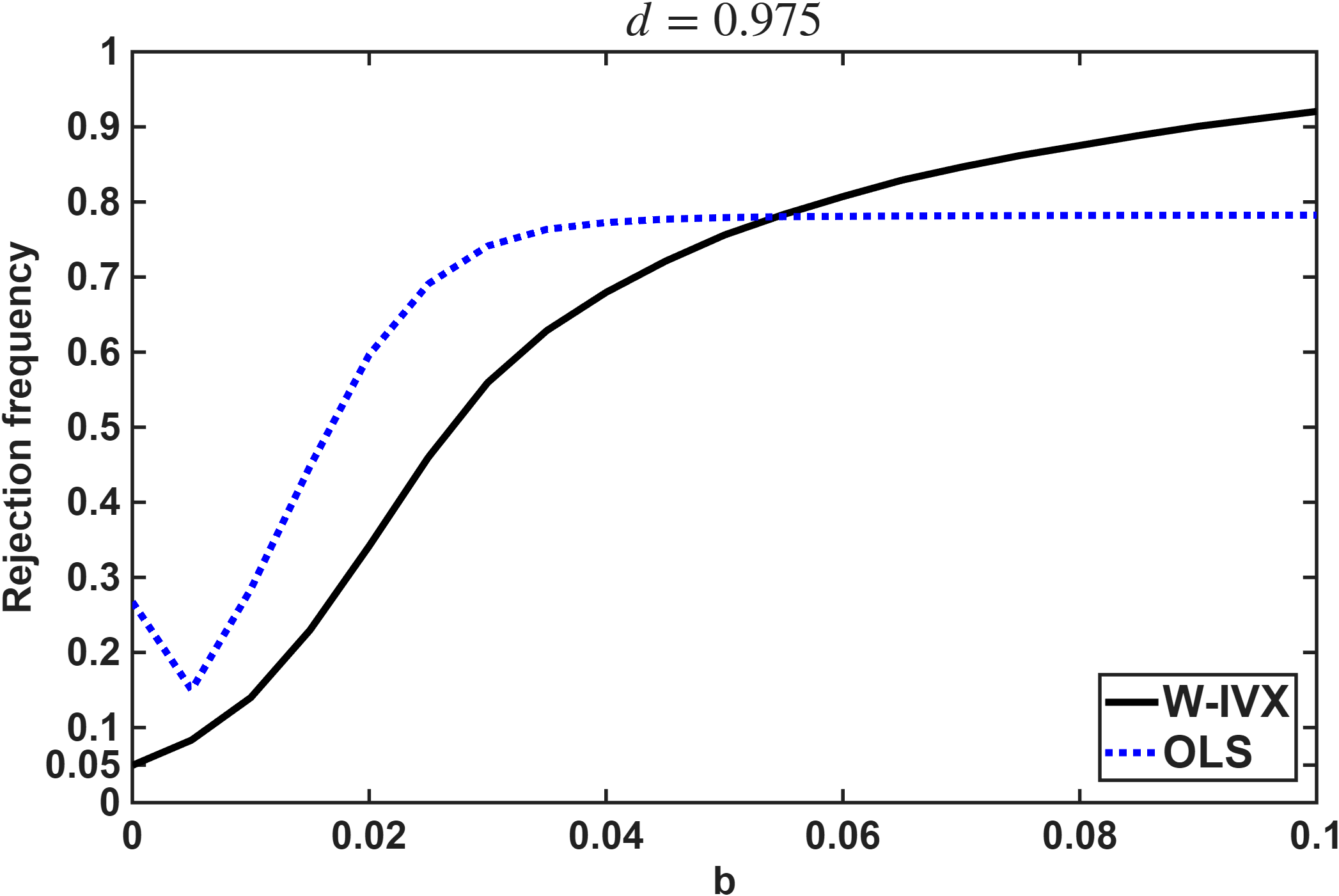}
\includegraphics[width=0.49\linewidth]{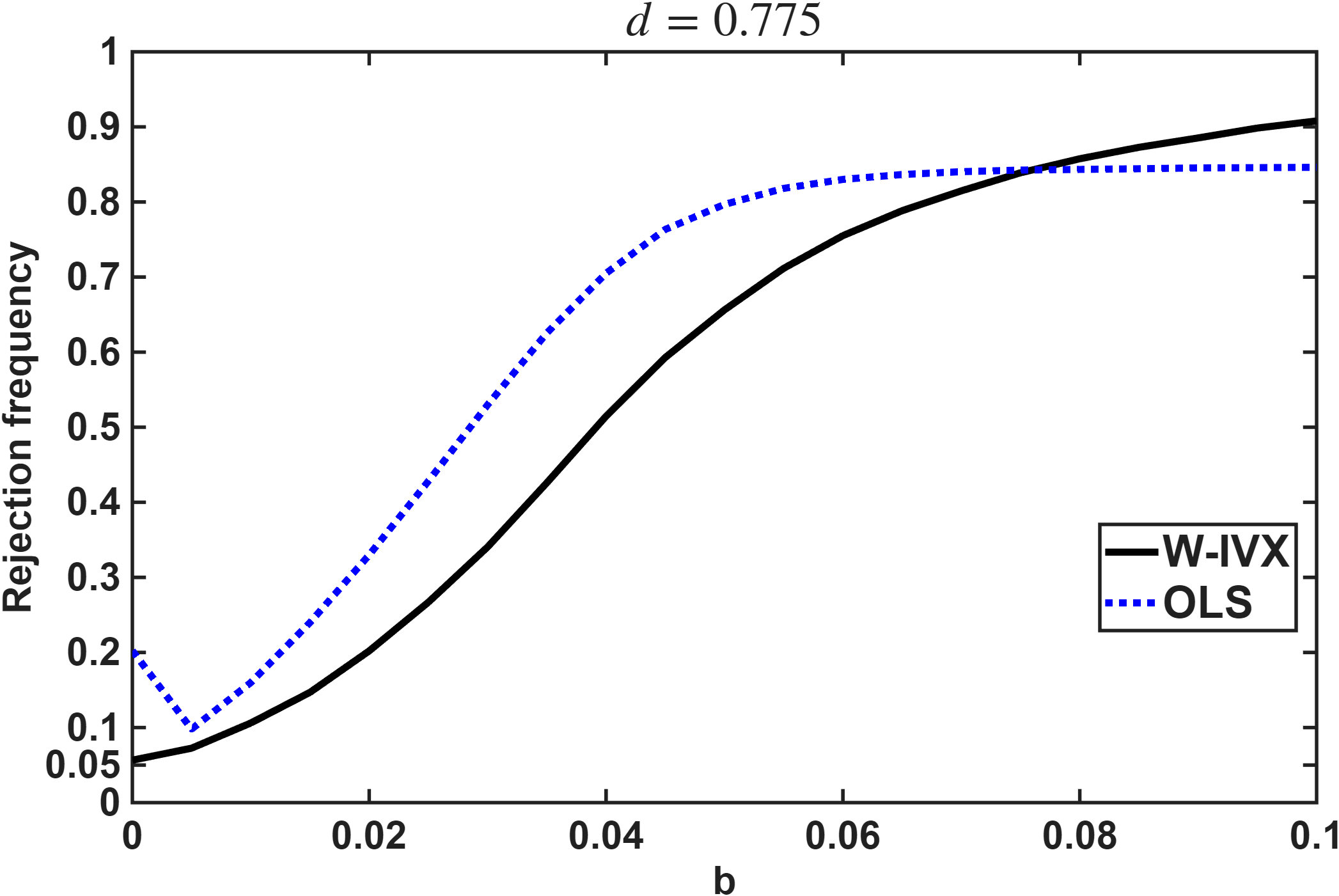}\\[0.6ex]
\includegraphics[width=0.49\linewidth]{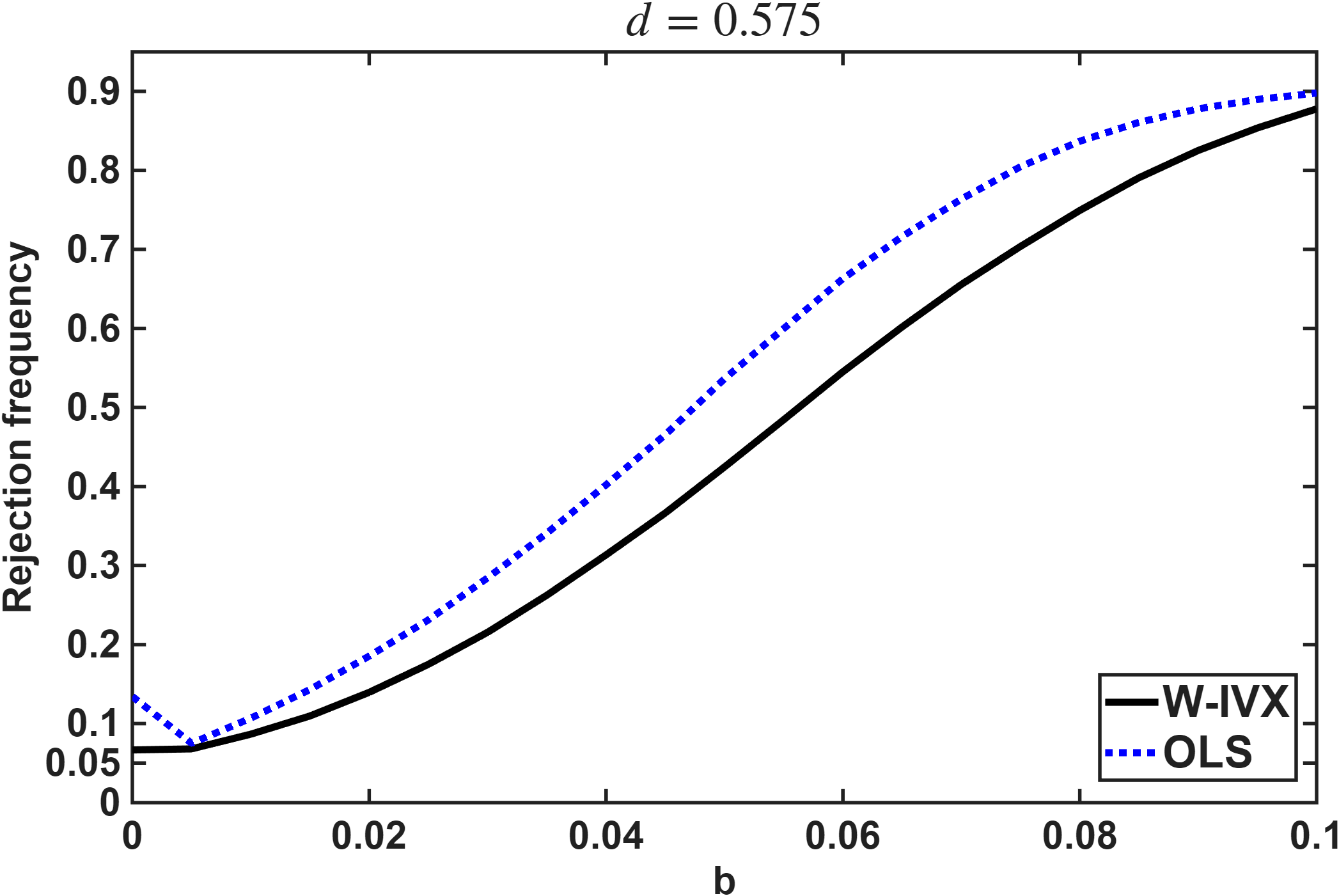}
\includegraphics[width=0.49\linewidth]{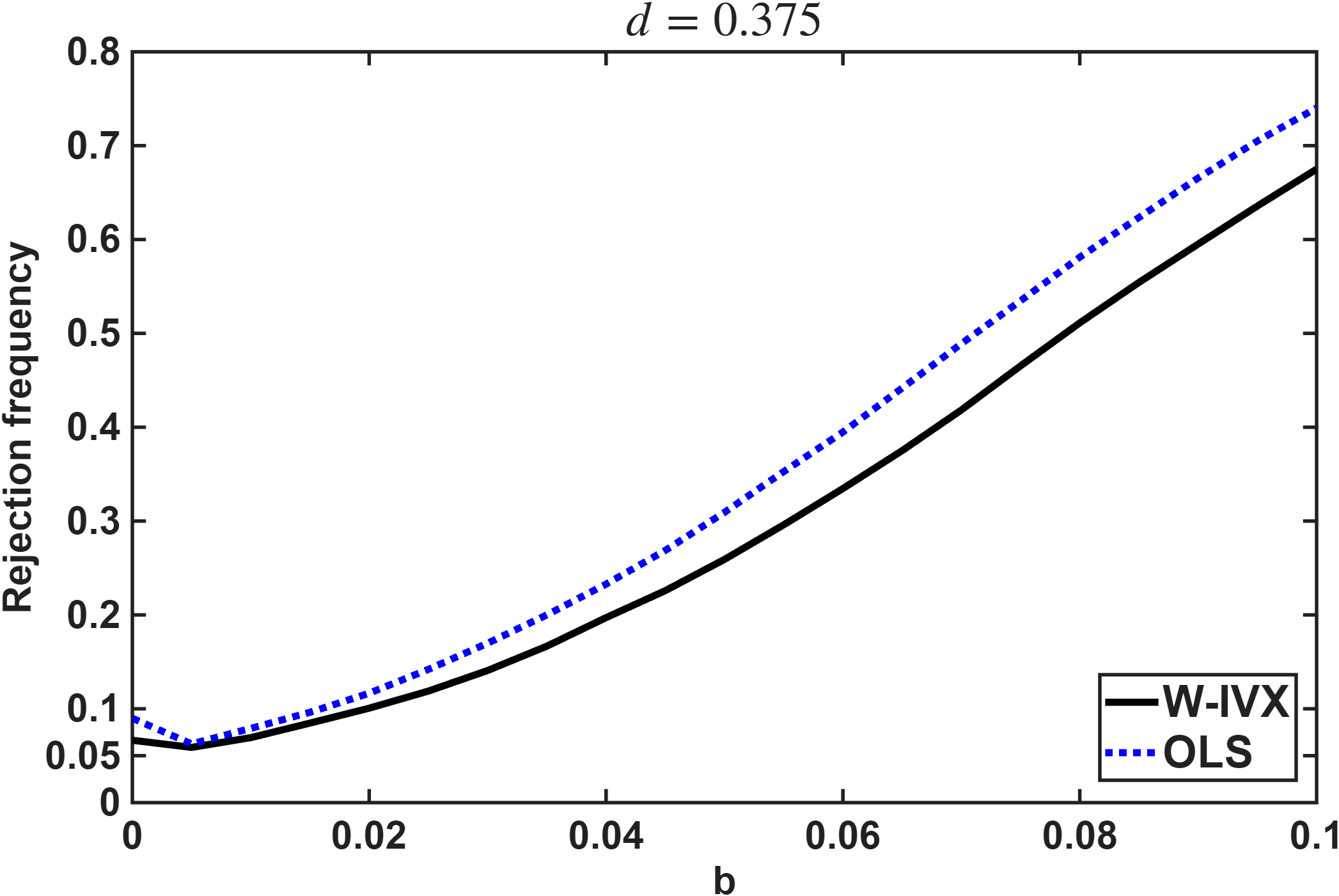}\\[0.6ex]
\includegraphics[width=0.49\linewidth]{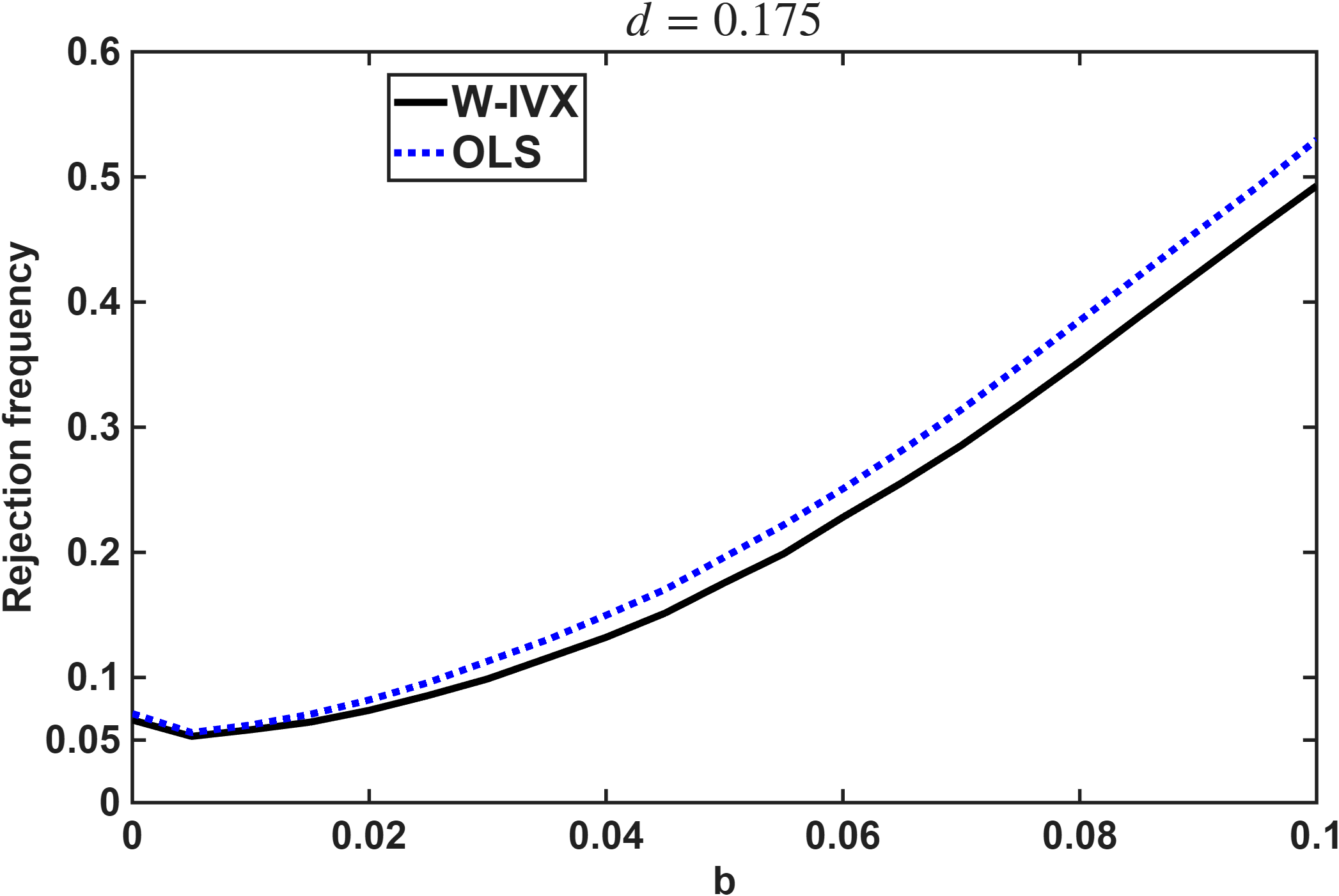}
\includegraphics[width=0.49\linewidth]{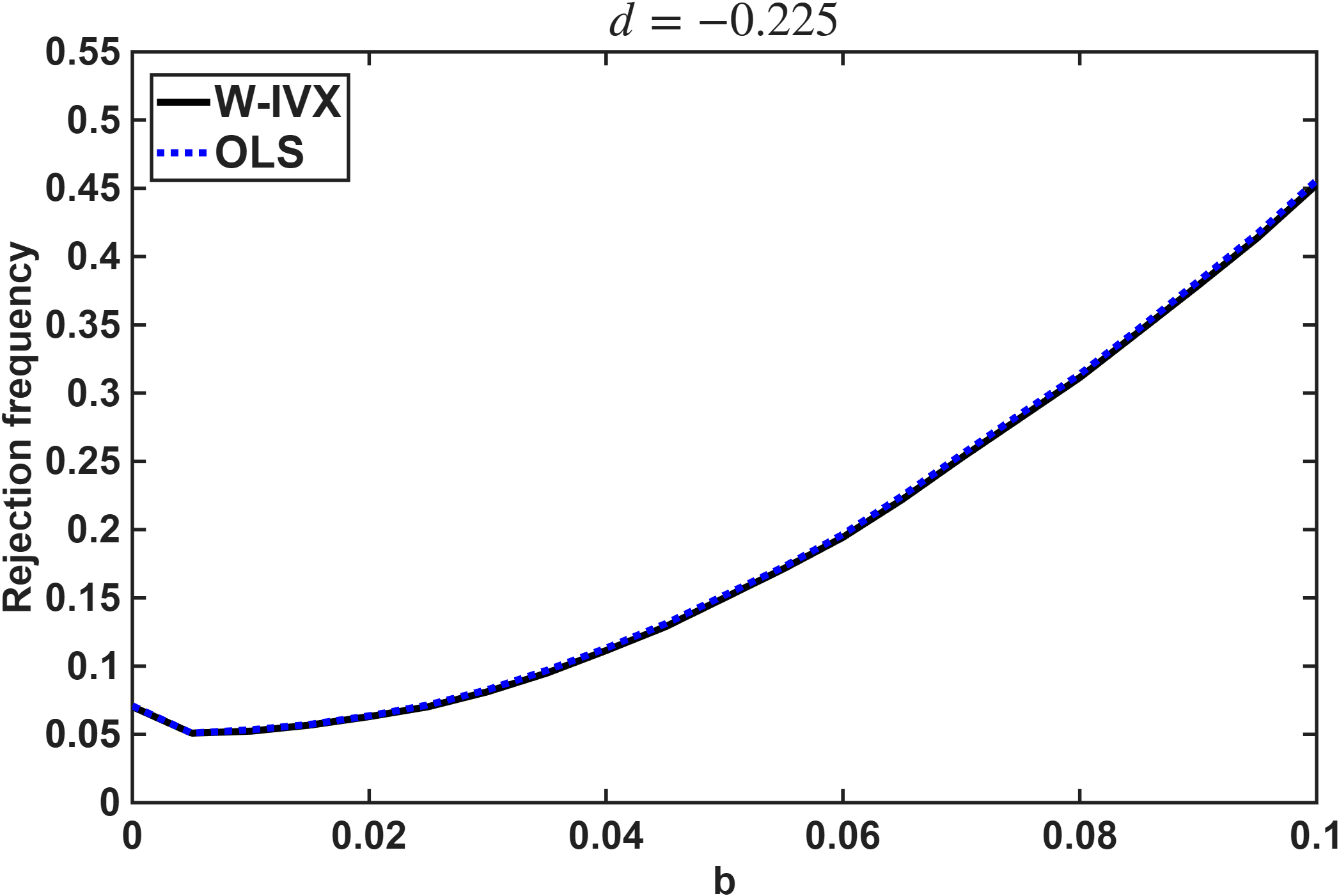}
\end{subfigure} \hfill{}
\begin{subfigure}[t]{0.48\linewidth} \centering \caption{$\alpha=1.4$, $d_{*}=0.286$}
\label{fig:sim-sap-nonlinear-f05-w05-alpha-1p40}\includegraphics[width=0.49\linewidth]{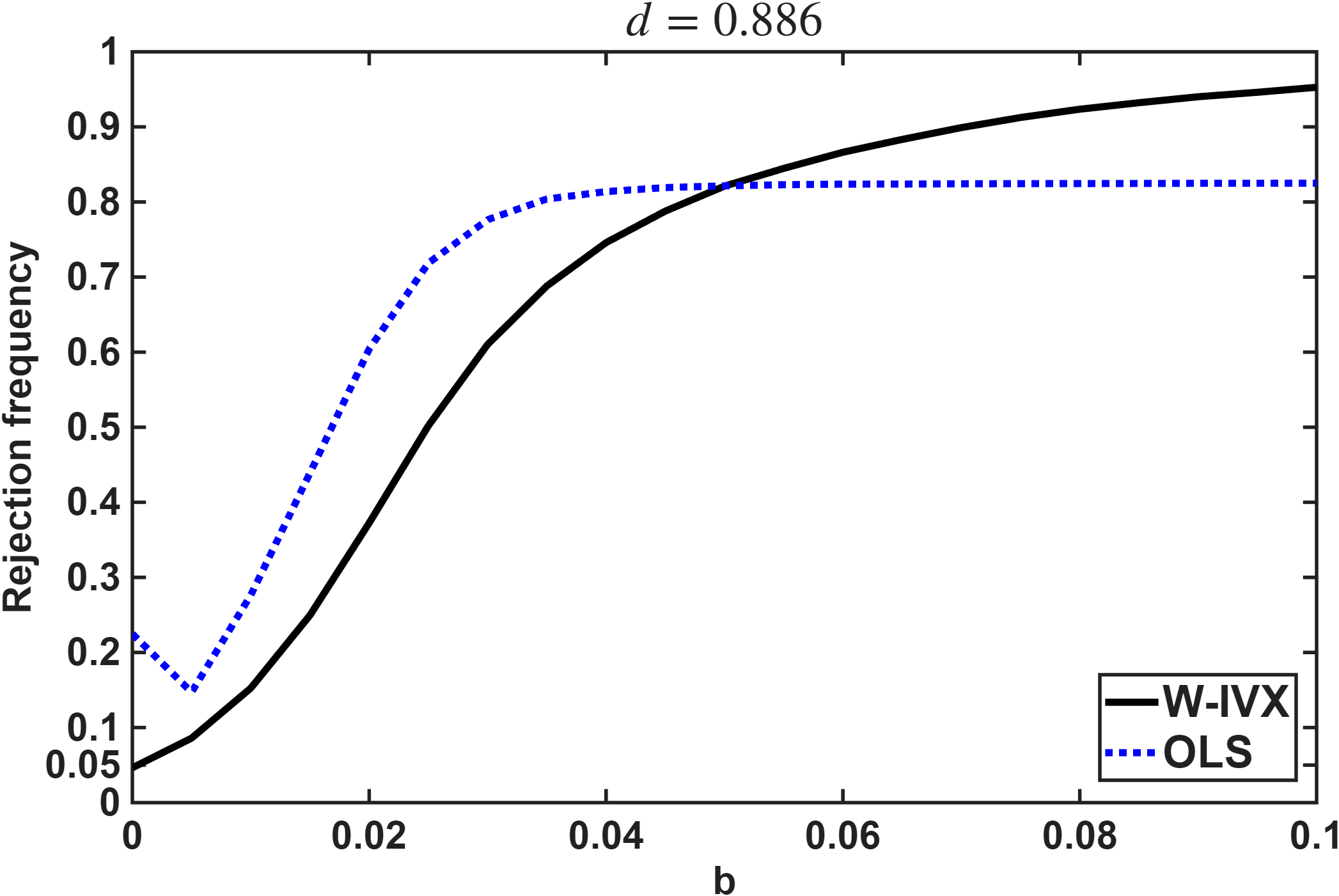}
\includegraphics[width=0.49\linewidth]{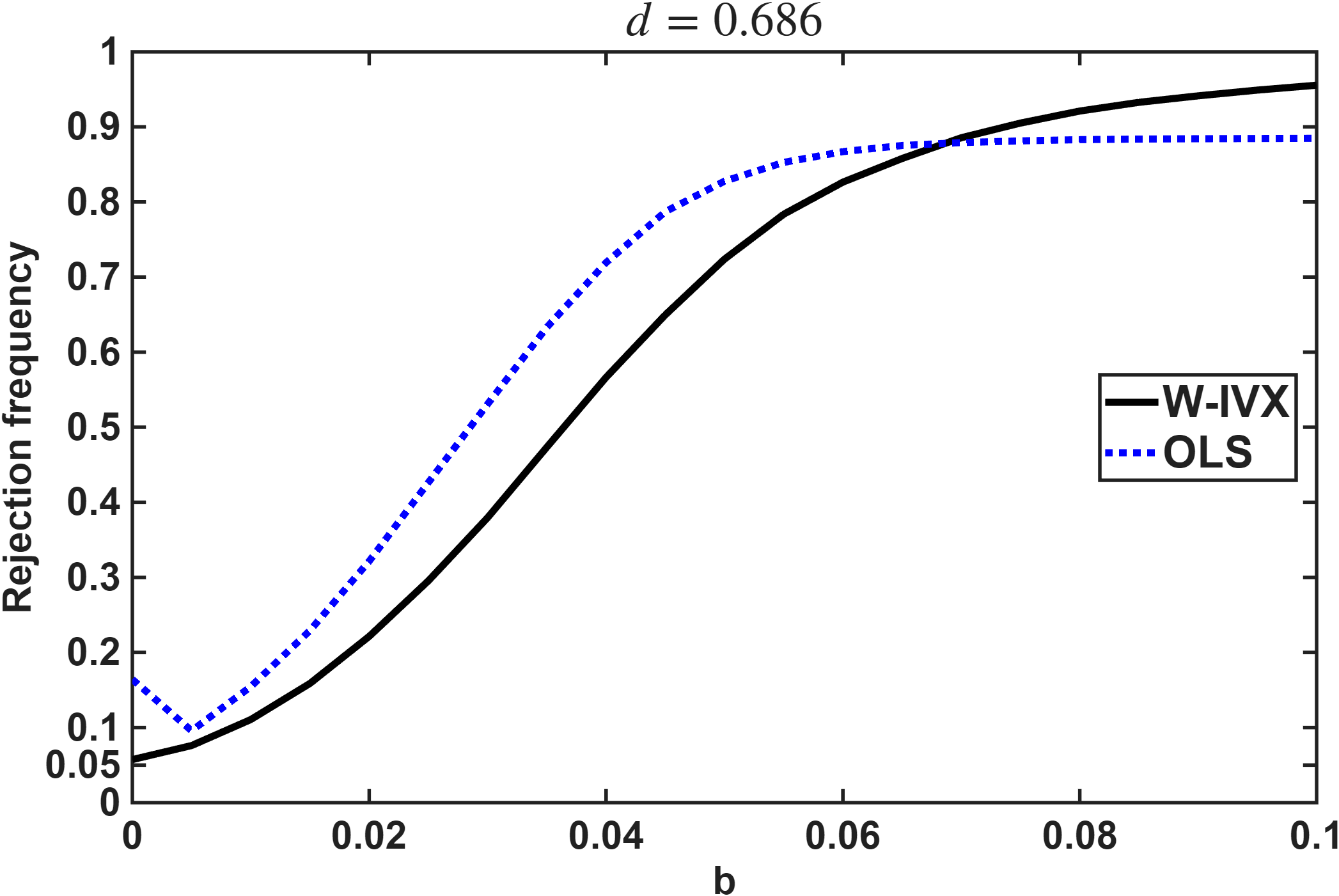}\\[0.6ex]
\includegraphics[width=0.49\linewidth]{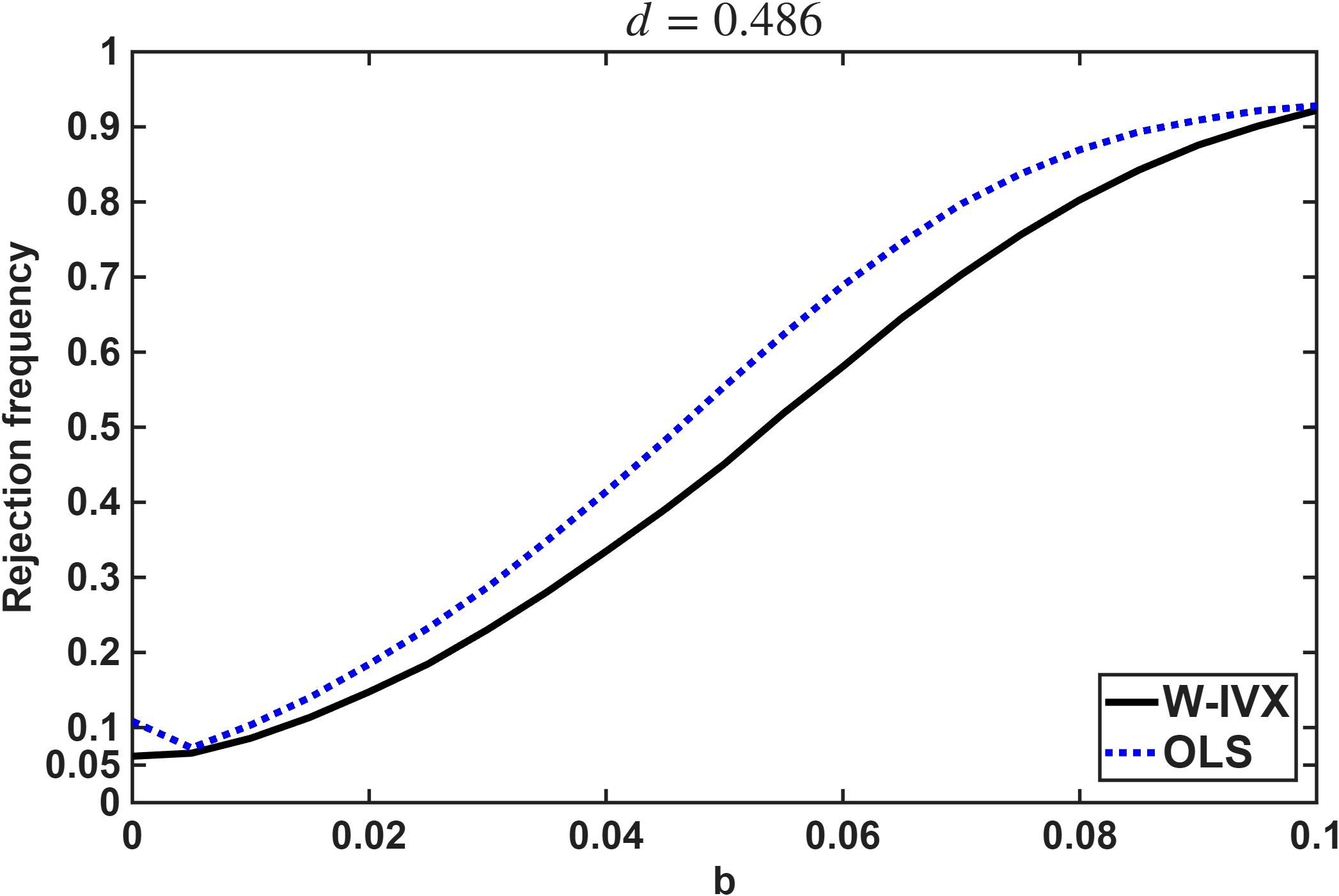}
\includegraphics[width=0.49\linewidth]{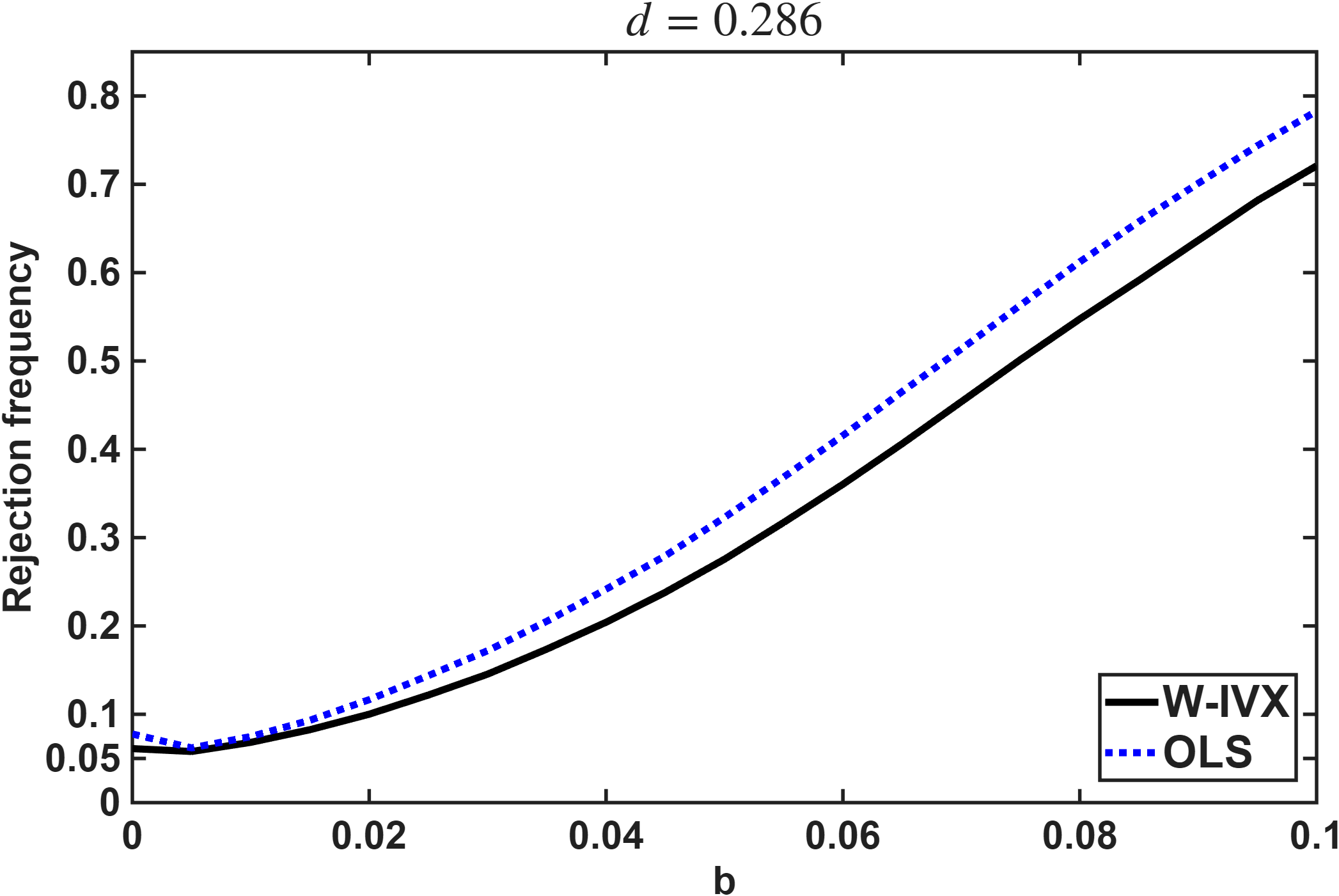}\\[0.6ex]
\includegraphics[width=0.49\linewidth]{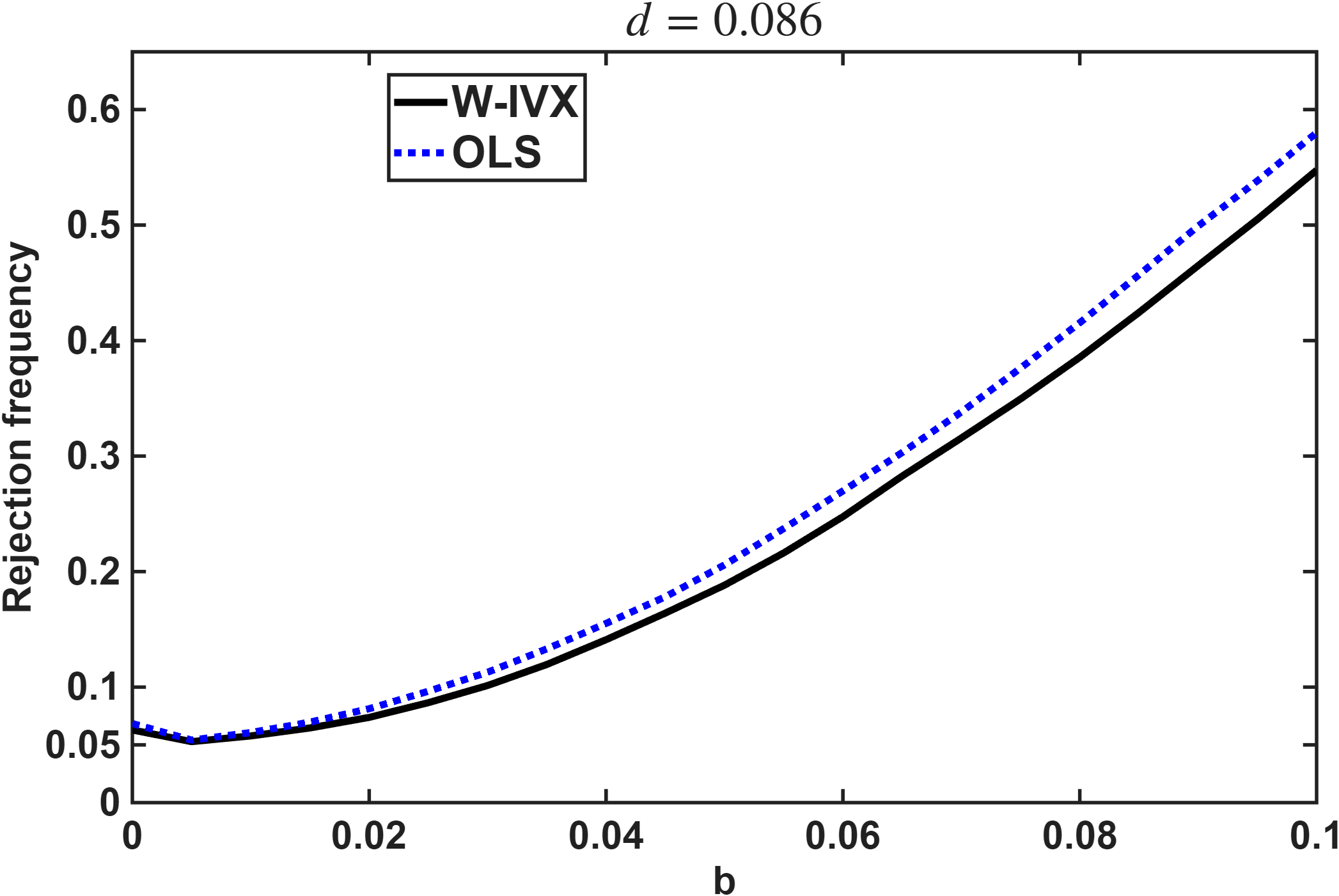}
\includegraphics[width=0.49\linewidth]{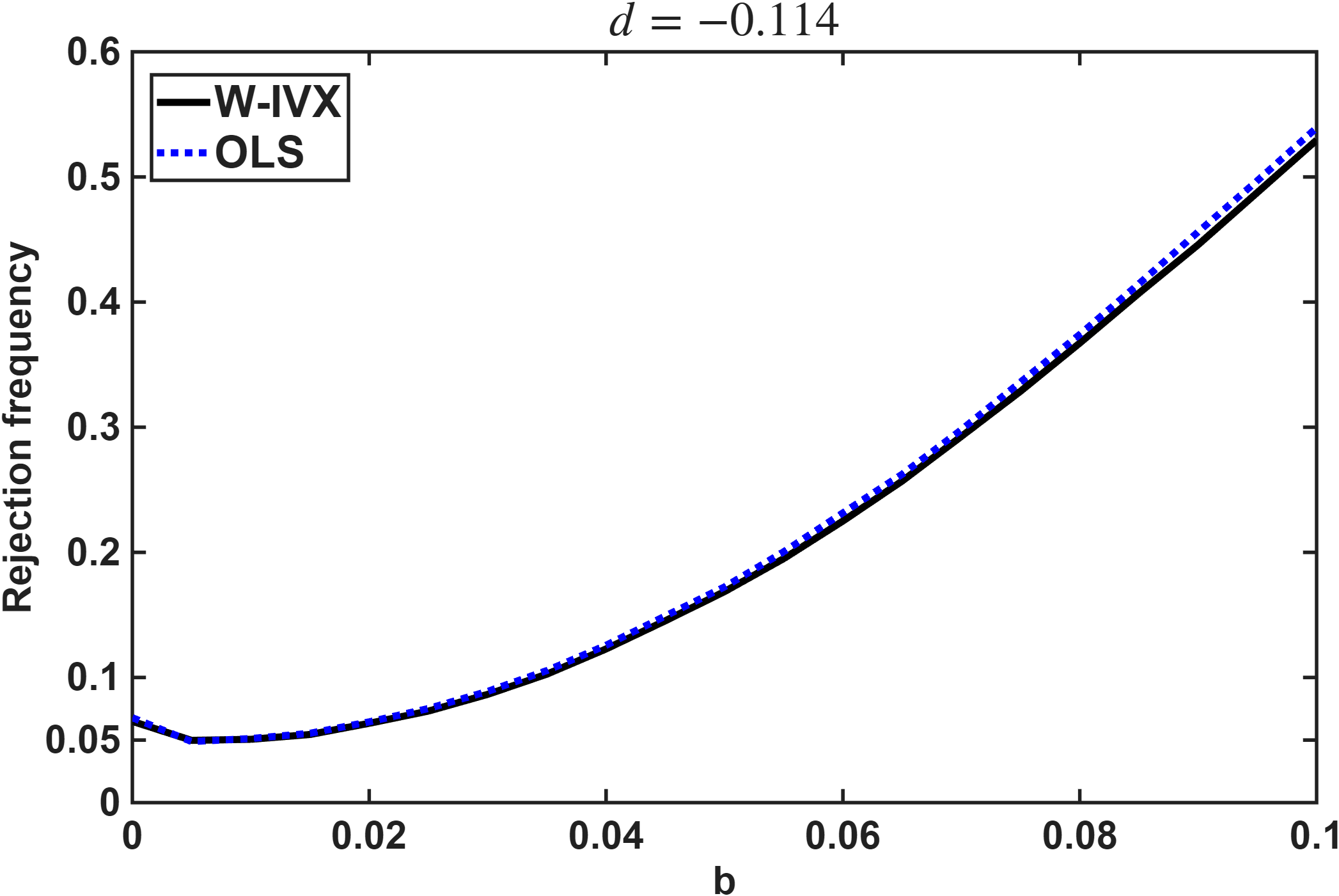}
\end{subfigure}
\end{figure}

\end{landscape} \clearpage{}

\end{document}